\documentclass[12pt,a4paper,leqno]{report}

\usepackage{amsmath,amscd,amsfonts,amsthm,amssymb,dsfont,bbm,manfnt}
\usepackage{hyperref,url,color,appendix,eucal}
\usepackage[usenames,dvipsnames]{xcolor}
\usepackage[nottoc,numbib]{tocbibind}
\usepackage{multicol}
\usepackage[OT2,T1]{fontenc}
\usepackage[utf8]{inputenc}
\usepackage[francais]{babel}
\usepackage[all]{xy}
\usepackage{lscape,slashbox,rotating} 
\usepackage{graphicx}

\usepackage{tocloft}

\def\got{\mathfrak}

\newenvironment{pf}
{\medskip\noindent {\it D\'emonstration.  }}
{\hfill\nobreak $\Box$ \par\bigbreak}

\newcommand{\isomo}{\overset{\sim}{\rightarrow}}
\newcommand{\isomol}{\overset{\sim}{\leftarrow}}
\newcommand{\lisomo}{\overset{\sim}{\longrightarrow}}
\newcommand{\Hom}{\mathrm{Hom}}
\newcommand{\pn}{\par \noindent}

\newcommand{\GL}{\mathrm{GL}}
\newcommand{\SL}{\mathrm{SL}}
\newcommand{\ps}{\par \smallskip}
  
\newcommand{\N}{\mathbb{N}}
\newcommand{\Z}{\mathbb{Z}}

\newcommand{\Q}{\mathbb{Q}}

\newcommand{\R}{\mathbb{R}}    
\newcommand{\C}{\mathbb{C}}

\newcommand{\AAA}{\mathbb{A}}
\newcommand{\F}{{\mathbb{F}}}

\newcommand{\PGL}{\mathrm{PGL}}   

\newcommand{\DET}{\mathrm{d\acute{e}t}}
\newcommand{\cbil}{{\rm B}}

\newtheorem{thm}[subsection]{Th\'eor\`eme}
\newtheorem*{thm*}{Th\'eor\`eme}
\newtheorem{lemme}[subsection]{Lemme}
\newtheorem{lemmev}[subsection]{Lemme{\color{green}${}^\star$}}
\newtheorem{critere}[subsection]{Crit\`ere}
\newtheorem*{lemme*}{Lemme}
\newtheorem{remark}[subsection]{Remarque}
\newtheorem{remarque}[subsection]{Remarque}
\newtheorem{cor}[subsection]{Corollaire} 

\newtheorem{corv}[subsection]{Corollaire{\color{green}${}^\star$}}
\newtheorem{prop}[subsection]{Proposition}
\newtheorem{propv}[subsection]{Proposition{\color{green}${}^\star$}}
\newtheorem{example}[subsection]{Exemple}

\newtheorem{conj}[subsection]{Conjecture}

\newtheorem{definition}[subsection]{D\'efinition}
\newtheorem{defprop}[subsection]{Proposition-D\'efinition}

\newtheorem{thmv}[subsection]{Th\'eor\`eme{\color{green}${}^\star$}}

\newtheorem{scholie}[subsection]{Scholie}
\newtheorem{defpropv}[subsection]{Proposition-D\'efinition{\color{green}${}^\star$}}
\newtheorem{thmintro}{Th\'eor\`eme}

\usepackage{titlesec}

\titleformat{\chapter}
{\bf \large}
{\Roman{chapter}.}
{1 em}
{}

\titleformat{\section}
{\bf \normalsize}
{\arabic{section}.}
{1 em}
{}

\titleformat{\subsection}
{}
{\bf \arabic{section}.\arabic{subsection}.}
{0.5 em}
{\normalsize}

\titleformat{\subsubsection}[runin]
{\small \bf}
{}
{}
{}

\makeatletter

\@addtoreset{equation}{section}
\makeatother

\hypersetup{colorlinks=true,linkcolor=black,citecolor=black,urlcolor=Periwinkle,pdfborderstyle={/S/U/W 1}}

\begin{document}

\title{Formes automorphes et voisins de Kneser des r\'eseaux de Niemeier}
\bigskip
\bigskip

\author{Ga\"etan Chenevier\thanks{Centre de Math\'ematiques Laurent
Schwartz, \'Ecole Polytechnique, 91128 Palaiseau Cedex, France (avant 10/2014). Universit\'e Paris-Sud, D\'epartement de math\'ematiques, B\^atiment 425, 91405 Orsay Cedex (\`a partir de 10/2014). Pendant la r\'edaction de ce m\'emoire, G.  Chenevier a \'et\'e financ\'e par le
C.N.R.S. , et il a re\c{c}u le soutien des projets ANR-10-BLAN 0114 et ANR-14-CE25.} \& Jean
Lannes\thanks{Universit\'e Paris 7, UFR de Math\'ematiques, b\^atiment
Sophie Germain, avenue de France, 75013 Paris, France.}}


\pagenumbering{roman}

\maketitle

{\setlength{\baselineskip}{1.2\baselineskip}
\tableofcontents\par}

\parindent=0.5cm
\renewcommand{\titlepage}{\relax}
\abstract{{}\ps Ce m\'emoire porte sur les r\'eseaux unimodulaires pairs de
dimension $\leq 24$ et sur les formes automorphes pour les groupes
orthogonaux, symplectiques ou lin\'eaires, auxquelles ces r\'eseaux sont reli\'es.  Notre
fil conducteur est la question, d'apparence \'el\'ementaire, de d\'eterminer
le nombre des $p$-voisinages, au sens de M.  Kneser, entre deux classes
d'isom\'etrie de tels r\'eseaux.  En utilisant notamment les travaux
r\'ecents d'Arthur \cite{arthur} sur les formes automorphes des groupes
classiques, nous d\'emontrons une formule explicite pour ces nombres, dans
laquelle interviennent certaines formes modulaires de Siegel de genre $1$ et
$2$.  \ps

Ces formules offrent un point de vue int\'eressant sur des probl\`emes
classiques de ce sujet.  Elles permettent par exemple de comprendre et
d'ordonner un ensemble de constructions du r\'eseau de Leech et de
d\'eterminer le graphe des $p$-voisins de Kneser des r\'eseaux de Niemeier
pour tout nombre premier $p$ (le cas $p=2$ \'etait connu de Borcherds
\cite{borcherdsthese}).  Elles ont d'autres applications que nous
d\'evelopperons dans le texte, comme la d\'emonstration d'une conjecture de
Nebe-Venkov~\cite{nebevenkov} portant sur les combinaisons lin\'eaires de
s\'eries th\^eta de genre sup\'erieur des r\'eseaux de Niemeier, la
d\'etermination de valeurs propres pour les op\'erateurs de Hecke de
certaines formes de genre $2$, ou encore la d\'emonstration d'une conjecture
de Harder~\cite{harder}.  \ps

Par un jeu de traductions assez classiques, le probl\`eme du comptage des
$p$-voisinages entre classes d'isom\'etrie de r\'eseaux de Niemeier se
ram\`ene \`a celui de d\'eterminer les repr\'esentations automorphes $\pi$
d'une certaine forme enti\`ere du groupe orthogonal euclidien ${\rm
O}_{24}(\R)$ qui sont soumises aux deux conditions suivantes : la composante
archim\'edienne $\pi_\infty$ est triviale et pour tout premier $p$ la
composante $\pi_p$ est ``non ramifi\'ee''.  Les travaux de Niemeier~\cite{Ni} et
Nebe-Venkov montrent qu'il y a $24$ telles repr\'esentations.  Les
r\'esultats sus-cit\'es d'Arthur sugg\`erent des angles d'attaque pour cette
question.  Une motivation importante de ce travail \'etait de comprendre ce
qu'ils disent dans ce contexte, et plus g\'en\'eralement sur les formes
automorphes partout non ramifi\'ees des groupes classiques.  Nous en donnons
dans ce m\'emoire plusieurs autres applications, par exemple \`a l'\'etude
des formes modulaires de Siegel de poids $\leq 12$ pour ${\rm Sp}_{2g}(\Z)$.  \ps}

\pagenumbering{arabic}

\chapter{Introduction} \label{chapintro}

\section{R\'eseaux unimodulaires pairs} Fixons un entier $n\geq 1$ et consid\'erons l'espace euclidien 
$\R^n$ muni de son produit scalaire standard $(x_i)\cdot (y_i)=\sum_i x_i y_i$. Un r\'eseau unimodulaire pair de rang $n$ est un r\'eseau $L \subset \R^n$ de covolume $1$ tel que $x \cdot x$ est pair pour tout $x$ dans $L$. L'ensemble $\mathcal{L}_n$ de ces r\'eseaux est muni d'une action du groupe orthogonal euclidien ${\rm O}(\R^n)$, et nous notons $${\rm X}_n \,\,:= \,\,{\rm O}(\R^n) \backslash \mathcal{L}_n$$
l'ensemble des classes d'isom\'etrie de r\'eseaux unimodulaires pairs de rang $n$. \`A chaque $L$ dans $\mathcal{L}_n$, on 
peut associer une forme quadratique $${\rm q}_L : L \rightarrow \Z, \, \, \, x  \mapsto \frac{x \cdot x}{2},$$
dont la forme bilin\'eaire sous-jacente $x\cdot y$ est de d\'eterminant $1$. L'application $L \mapsto {\rm q}_L$ induit alors une bijection entre ${\rm X}_n$ et l'ensemble des classes d'isomorphisme de formes quadratiques de rang $n$, d\'efinies positives, et de d\'eterminant $1$,
sur l'anneau des entiers $\Z$. \ps \ps

Il est bien connu que l'ensemble ${\rm X}_n$ est fini, et qu'il est non vide si, et seulement si, on a $n \equiv 0 \bmod 8$. Un exemple standard d'\'el\'ement de $\mathcal{L}_n$ est $${\rm E}_n \,\,:=\,\, {\rm D}_n \,+ \,\Z e$$
o\`u ${\rm D}_n = \{ (x_i) \in \Z^n, \sum_i x_i \equiv 0 \bmod 2\}$, $e  = \frac{1}{2} (1,1,\dots,1)$, et $n \equiv 0 \bmod 8$. Expliquons cette notation. \`A tout \'el\'ement $L$ de $\mathcal{L}_n$ est associ\'e un syst\`eme de racines (de type {\rm ADE}) $${\rm R}(L)\,\,:=\,\,\{ x \in L, x \cdot x =2\},$$ de rang $\leq n$. Le syst\`eme de racines ${\rm R}({\rm E}_8)$ est alors de type ${\bf E}_8$ et engendre $\Z$-lin\'eairement ${\rm E}_8$. Pour $n>8$, ${\rm R}({\rm E}_n)$ est de type ${\bf D}_n$ et engendre ${\rm D}_n$. \ps \ps

L'ensemble ${\rm X}_n$ n'a \'et\'e d\'etermin\'e qu'en dimension $n \leq 24$. Mordell et Witt ont respectivement d\'emontr\'e 
$${\rm X}_8=\{ {\rm E}_8 \}  \, \, \, {\rm et}\, \, \, {\rm X}_{\rm 16}=\{ {\rm E}_{16}, {\rm E}_8 \oplus {\rm E}_8\}.$$
Les deux r\'eseaux ${\rm E}_{16}$ et ${\rm E}_8 \oplus {\rm E}_8$ joueront un grand r\^ole dans ce m\'emoire. Ils sont \`a la fois faciles et difficiles \`a distinguer : leurs syst\`emes de racines sont diff\'erents, mais ils repr\'esentent chaque entier exactement le m\^eme nombre de fois. Cette derni\`ere propri\'et\'e, fameuse, conduit par exemple aux tores isospectraux de Milnor. \ps \ps

Les \'el\'ements de ${\rm X}_{24}$ ont \'et\'e classifi\'es par Niemeier \cite{Ni}, qui a notamment d\'emontr\'e $|{\rm X}_{24}|=24$. Avant d'en dire plus sur ces derniers, mentionnons que pour $n \geq 32$ la formule de masse Minkowski-Siegel-Smith montre que la taille de ${\rm X}_n$ explose. Par exemple, on a $|{\rm X}_{32}|> 8 \cdot 10^6$ \cite{serre} ; il y a m\^eme plus d'un milliard d'\'el\'ements dans ${\rm X}_{32}$ d'apr\`es King \cite{king}. \ps \ps

Le plus fameux des r\'eseaux de Niemeier (on appelle ainsi les \'el\'ements de $\mathcal{L}_{24}$) est le r\'eseau de Leech. \`A isom\'etrie pr\`es, c'est le seul \'el\'ement $L$ de ${\mathcal{L}}_{24}$ tel que ${\rm R}(L) = \emptyset$ (Conway). Un fait remarquable est que si $L$ est un r\'eseau de Niemeier non isomorphe au r\'eseau de Leech, alors ${\rm R}(L)$ est de rang $24$ et tous ses constituants irr\'eductibles ont m\^eme nombre de Coxeter ; une d\'emonstration simple en a \'et\'e donn\'ee par Venkov~\cite{venkov}.  Le miracle est alors que $L \mapsto {\rm R}(L)$ induit une bijection entre ${\rm X}_{24}- \{{\rm Leech}\}$ et l'ensemble des classes d'isomorphisme de syst\`emes de racines $R$ de rang $24$ dont les constituants irr\'eductibles sont de type ADE et tous de m\^eme nombre de Coxeter ${\rm h}(R)$. La d\'emonstration est une v\'erification fastidieuse, au cas par cas.\ps\ps

\begin{table}[htp]

\renewcommand{\arraystretch}{1.5}

\begin{center}
{\scriptsize
\begin{tabular}{|c|c|c|c|c|c|c|c|c|}

\hline ${\rm R}$ &  ${\bf D}_{24}$ &  ${\bf D}_{16}{\bf E}_{8}$ &  $3 {\bf E}_8$ &  ${\bf A}_{24}$ &  $2 {\bf D}_{12}$ &  ${\bf A}_{17}\,{\bf E}_7$ &  ${\bf D}_{10}\, 2 {\bf E}_7$  &  ${\bf A}_{15}\, {\bf D}_9$ \\

\hline ${\rm h(R)}$ &  $46$ &  $30$ &  $30$ &  $25$ &  $22$ &  $18$ &  $18$ &  $16$ \\

\hline  ${\rm R}$ & $3 {\bf D}_8$ &  $2 {\bf A}_{12}$ &  ${\bf A}_{11}\,{\bf D}_7\,{\bf E}_6$ &  $4 {\bf E}_6$  &  $2 {\bf A}_9 {\bf D}_{6}$ &  $4 {\bf D}_{6}$ &  $3 {\bf A}_8$ &  $2 {\bf A}_{7}\,2 {\bf D}_5$ \\

\hline ${\rm h(R)}$ &  $14$ &  $13$ &  $12$ &  $12$  &  $10$ &  $10$ &  $9$ &  $8$  \\

\hline ${\rm R}$ &  $4 {\bf A}_{6}$ &  $4 {\bf A}_{5} \,{\bf D}_4$ &  $6 {\bf D}_4$ &  $6 {\bf A}_{4}$ &  $8 {\bf A}_3$ &  $12 {\bf A}_2$ &  $24 {\bf A}_1$  &    \\

\hline ${\rm h(R)}$&   $7$ &  $6$ &  $6$ &  $5$ &  $4$ &  $3$ &  $2$ &  \\

\hline
\end{tabular}

\caption{{\small Les $23$ syst\`emes de racines \'equi-Coxeter de type ADE et de rang $24$}}
\label{listeniemeier}
}
\end{center}
\end{table}

Les r\'esultats mentionn\'es dans ce paragraphe sont expos\'es dans le chapitre II. Il s'agit pour l'essentiel d'un chapitre de rappels. Nous d\'eveloppons d'abord des pr\'erequis d'alg\`ebre bilin\'eaire et quadratique. Ils sont n\'ecessaires \`a  la compr\'ehension des constructions de formes quadratiques auxquelles nous avons fait allusion ci-dessus, et d'autres dont nous aurons besoin au cours du m\'emoire. En particulier, nous rappelons la th\'eorie de Venkov et r\'eexpliquons la construction de certains des r\'eseaux de Niemeier. Nous en profitons \'egalement pour faire quelques rappels sur les (sch\'emas en) groupes classiques sur $\Z$, qui seront utilis\'es dans la suite du m\'emoire. L'appendice ${\rm B}$ contient notamment une variante des r\'esultats du chapitre II : nous y \'etudions les r\'eseaux pairs de $\R^n$ de d\'eterminant $2$, et plus g\'en\'eralement la th\'eorie des formes quadratiques sur $\Z$ (de dimension impaire) correspondante.

\section{Voisins \`a la Kneser} \label{secknesserintro}

Soit $p$ un nombre premier. La notion de $p$-voisins a \'et\'e introduite par M. Kneser, elle peut \^etre vue comme un proc\'ed\'e de construction de r\'eseaux \`a partir de r\'eseaux connus. Dans le chapitre III, nous \'etudions plusieurs variations de cette notion et nous en donnons de nombreux exemples. \ps \ps

Suivant Kneser, deux r\'eseaux $L, M$ dans $\mathcal{L}_n$ sont dits {\it $p$-voisins} si $L\cap M$ est d'indice $p$ dans $L$ (et donc dans $M$). Il est ais\'e de fabriquer tous les $p$-voisins d'un r\'eseau $L$ appartenant \`a $\mathcal{L}_n$ donn\'e. En effet, \`a chaque droite isotrope $\ell$ de $L \otimes \F_p$, disons engendr\'ee par un \'el\'ement $x$ de $L$ v\'erifiant ${\rm q}_L(x) \equiv 0 \bmod p^2$, on peut associer le r\'eseau unimodulaire pair
$${\rm vois}_p(L ; \ell)\,\,{:=}\,\, H + \Z\,\, \frac{x}{p},$$
avec $H= \{ y \in L, x \cdot y \equiv 0 \bmod p\}$ (ce r\'eseau ne d\'epend pas du choix de $x$). L'application $\ell \mapsto {\rm vois}_p(L ; \ell)$ induit une bijection entre ${\rm C}_L(\F_p)$ et l'ensemble des $p$-voisins de $L$, o\`u ${\rm C}_L$ d\'esigne la quadrique projective (et lisse) sur $\Z$ d\'efinie par ${\rm q}_L=0$. Cette quadrique \'etant hyperbolique sur $\F_p$ pour tout premier $p$, le nombre de $p$-voisins de $L$ est exactement $|{\rm C}_L(\F_p)| = 1+p+p^2+\cdots +
p^{n-2} + p^{\frac{n}{2}-1}$, une quantit\'e que l'on note ${\rm c}_n(p)$. \ps \ps

Par exemple, consid\'erons l'\'el\'ement $\rho=(0,1,2,\dots,23)$ de ${\rm E}_{24}$. Il engendre une droite isotrope dans ${\rm E}_{24} \otimes \mathbb{F}_{47}$ car $\sum_{i=1}^{23} i^2 \equiv 0 \bmod 47$. Il n'est pas tr\`es dif\-ficile de v\'erifier que ${\rm vois}_{47}( {\rm E}_{24} ; \rho)$ ne poss\`ede pas de racine, d'o\`u une isom\'etrie 
$${\rm vois}_{47}( {\rm E}_{24} ; \rho) \simeq {\rm Leech}.$$
Cette construction particuli\`erement simple du r\'eseau de Leech est attribu\'ee \`a Thompson dans \cite{conwaysloane}, et nous y reviendrons plus loin. Elle illustre le slogan affirmant que de nombreuses constructions de r\'eseaux sont des cas particuliers de constructions de voisins. \ps \ps

Revenant au cadre g\'en\'eral, on dispose pour tout $L$ dans $\mathcal{L}_n$ d'une partition de la quadrique ${\rm C}_L(\F_p)$ donn\'ee par la classe d'isom\'etrie du $p$-voisin associ\'e. Un des objets de ce m\'emoire est d'\'etudier cette partition en dimension $n \leq 24$. Par exemple, peut-on d\'eterminer le nombre ${\rm N}_p(L,M)$ de $p$-voisins de $L$ isom\'etriques \`a un $M \in \mathcal{L}_n$ donn\'e ? Le premier cas int\'eres\-sant est bien entendu celui de la dimension $n=16$. Il sera commode, pour formuler le r\'esultat, d'introduire l'op\'erateur lin\'eaire ${\rm T}_p : \Z[{\rm X}_n] \rightarrow \Z[{\rm X}_n]$ d\'efini par ${\rm T}_p \,[L]= \sum [M]$, la somme portant sur les $p$-voisins $M$ de $L$.\ps \ps

\begin{thmintro}\label{thmintro16} Dans la base \,\,$\mathrm{E}_8\oplus\mathrm{E}_8,\, \,  \mathrm{E}_{16}$\, \, la matrice de $\mathrm{T}_p$ est 
$$\mathrm{c}_{16}(p) \begin{bmatrix}1 & 0 \\ 0 &
1\end{bmatrix} + (1+p+p^2+p^3)\hspace{2pt}\frac{1+p^{11}-\tau(p)}{691}
\begin{bmatrix} -405 & 286 \\ 405 & -286\end{bmatrix},$$
o\`u $\tau$ est la fonction de Ramanujan, d\'efinie par $q\underset{m\geq 1}{\prod}(1-q^m)^{24}=\sum_{n\geq 1} \tau(n) q^n$.
\end{thmintro}
 Par exemple, on a la formule ${\rm N}_p(\mathrm{E}_8\oplus\mathrm{E}_8,\,
\mathrm{E}_{16})=\frac{405}{691}(1+p^{11}-\tau(p))\frac{p^4-1}{p-1}$.  Ce
th\'eor\`eme \'etait sans doute connu des sp\'ecialistes, mais nous n'en
avons pas trouv\'e de trace sous cette forme dans la litt\'erature.  Nous en
donnerons plusieurs d\'emonstrations, sur lesquelles nous reviendrons plus
loin.  D'apr\`es la th\'eorie des formes modulaires pour ${\rm SL}_2(\Z)$ et
des s\'eries th\^eta, la pr\'esence de $\tau(n)$ dans les questions de
r\'eseaux para\^it de prime abord assez classique.  Par exemple, si l'on pose
$r_L(n)=|\{ x \in L, x \cdot x = 2n\}|$, on montre facilement $r_{\rm
Leech}(p)= \frac{65520}{691} (1+p^{11}-\tau(p))$ pour tout premier $p$, une
formule d'apparence similaire \`a celle du th\'eor\`eme.  N\'eanmoins,
l'apparition de $\tau(p) \frac{p^4-1}{p-1}$ dans le probl\`eme de comptage
de $p$-voisins \'evoqu\'e plus haut est sensiblement plus subtile ; elle
s'av\`erera \'equivalente \`a un cas non trivial de fonctorialit\'e
d'Arthur-Langlands.\footnote{La comparaison du th\'eor\`eme \ref{thmintro16} et de la formule pour $r_{\rm
Leech}(p)$ donn\'ee ci-dessus conduit d'ailleurs \`a la relation ``purement quadratique'' suivante, que nous ne savons
pas d\'emontrer directement : 
\,\,\,${\rm N}_p({\rm E}_8\oplus {\rm E}_8,{\rm E}_{16}) = \frac{9}{1456}
\cdot r_{\rm Leech}(p) \cdot \frac{p^4-1}{p-1}$.} \ps \ps

Notre th\'eor\`eme principal est un \'enonc\'e du m\^eme type que le
th\'eor\`eme \ref{thmintro16} mais concernant les r\'eseaux de Niemeier.  Il
est possible d'en donner un \'enonc\'e dans le style du th\'eor\`eme
pr\'ec\'edent, qui est toutefois tr\`es indigeste.  Chose remarquable, il
fait intervenir les formes modulaires paraboliques pour ${\rm SL}_2(\Z)$ de
poids $k\leq 22$ ainsi que $4$ formes modulaires de Siegel pour ${\rm
Sp}_4(\Z)$ \`a valeurs vectorielles.  Les formules explicites contiennent de
si gros coefficients rationnels qu'il est tout \`a fait \'etonnant en les
inspectant que ${\rm N}_p(L,M)$ soit un nombre entier !  Nous \'enoncerons
une version plus conceptuelle (et \'equivalente) de notre r\'esultat au \S
\ref{fautintro} (Th\'eor\`eme \ref{thmintro24}).  Expliquons en revanche d\`es \`a pr\'esent
certaines cons\'equences concernant les r\'eseaux de Niemeier d'une analyse
de nos formules.\ps \ps

Consid\'erons le graphe ${\rm K}_n(p)$ d'ensemble de sommets ${\rm X}_n$, o\`u les classes de deux r\'eseaux non isomorphes $L$ et $M$ sont reli\'ees par une ar\^ete si, et seulement  si, on a ${\rm N}_p(L,M) \neq 0$. Kneser a d\'emontr\'e que ${\rm K}_n(p)$ est connexe, quels que soient $n$ et $p$, comme cons\'equence de son c\'el\`ebre th\'eor\`eme d'approximation forte. Ce joli \'enonc\'e montre que l'on peut th\'eoriquement reconstruire ${\rm X}_n$ \`a partir du r\'eseau ${\rm E}_n$ et d'un nombre premier. Ce fut d'ailleurs utile \`a Niemeier pour comprendre ${\rm X}_{24}$ par $2$-voisinages. \ps \ps

Le graphe ${\rm K}_{16}(p)$ est le graphe connexe \`a deux sommets, d'apr\`es Kneser. C'est bien s\^ur compatible avec l'estimation $|\tau(p)|< 2 p^{\frac{11}{2}}$ (Deligne-Ramanujan)  et la formule pour ${\rm N}_p({\rm E}_8 \oplus {\rm E}_8, {\rm E}_{16})$ donn\'ee par le th\'eor\`eme \ref{thmintro16}. Par contre, le graphe ${\rm K}_{24}(2)$, d\'etermin\'e par
Borcherds~\cite{conwaysloane}, n'est pas du tout trivial. Il est de diam\`etre $5$ et cette page Wikipedia \url{http://en.wikipedia.org/wiki/Niemeier_lattice} en donne une jolie repr\'esentation, aussi due \`a Borcherds. Nos r\'esultats permettent par exemple de d\'eterminer ${\rm K}_{24}(p)$ pour tout nombre premier $p$ (\S X.2).\footnote{Une liste de ces graphes est disponible \`a l'adresse  \url{http://gaetan.chenevier.perso.math.cnrs.fr/niemeier/niemeier.html}.}

\begin{thmintro}\label{thm24pop} \begin{itemize} \item[(i)] Soit $L$ un r\'eseau de Niemeier avec racines. Alors $L$ est un $p$-voisin du r\'eseau de Leech si, et seulement si, $p \geq {\rm h}({\rm R}(L))$. \ps \ps
\item[(ii)] Le graphe ${\rm K}_{24}(p)$ est complet si, et seulement si, $p \geq 47$. \ps \ps
\end{itemize}
\end{thmintro}

Faisons quelques remarques sur cet \'enonc\'e. Le point (i) concerne les cons\-tructions du r\'eseau de Leech comme $p$-voisin d'un r\'eseau de Niemeier avec racines. Par exemple, on observe sur le graphe de Borcherds ${\rm K}_{24}(2)$ que ${\rm Leech}$ est \`a distance $5$ de ${\rm E}_{24}$ et qu'il n'est  reli\'e qu'au r\'eseau de syst\`eme de racines $24 {\bf A}_1$, le r\'eseau de Niemeier avec racines le plus d\'elicat \`a construire (il n\'ecessite le code de Golay, \S II.3). Cette derni\`ere propri\'et\'e est en fait assez simple \`a comprendre : si le r\'eseau de Leech est un $2$-voisin du r\'eseau de Niemeier $L$ (avec racines), alors un sous-groupe d'indice $2$ de $L$ ne  contient  pas de racine. En particulier, ${\rm R}(L)$ a la propri\'et\'e que la somme de deux racines n'est jamais 
une racine, de sorte que ses constituants irr\'eductibles sont de rang $1$, {\it i.e.} ${\rm R}(L)=24 {\bf A}_1$. Parmi ceux de la table \ref{listeniemeier}, ce syst\`eme de racines est aussi le seul dont le nombre de Coxeter est $2$, en accord avec le (i). 

\ps \ps La partie la plus \'el\'ementaire du th\'eor\`eme \ref{thm24pop},
d\'emontr\'ee au \S III.4 et g\'en\'eralisant l'observation ci-dessus,
consiste \`a v\'erifier que l'on a $p \geq {\rm h}({\rm R}(L))$ si ${\rm
Leech}$ est $p$-voisin de $L$.  C'est un analogue formel d'un r\'esultat de
Kostant~\cite{kostantsl2} affirmant que l'ordre minimal d'un \'el\'ement
r\'egulier d'ordre fini dans un groupe de Lie compact connexe adjoint
co\"incide avec le nombre de Coxeter de son syst\`eme de racines.  En
revanche, la d\'emonstration des autres assertions n\'ecessite l'utilisation
du th\'eor\`eme \ref{thmintro24} ainsi qu'un ensemble d'in\'egalit\'es de
type Ramanujan, elle ne sera compl\'et\'ee qu'au chapitre X 
(\S \ref{parfinniemeier} et \S \ref{calcfinvpgenre2}).  \ps \ps

\begin{figure} 
\includegraphics{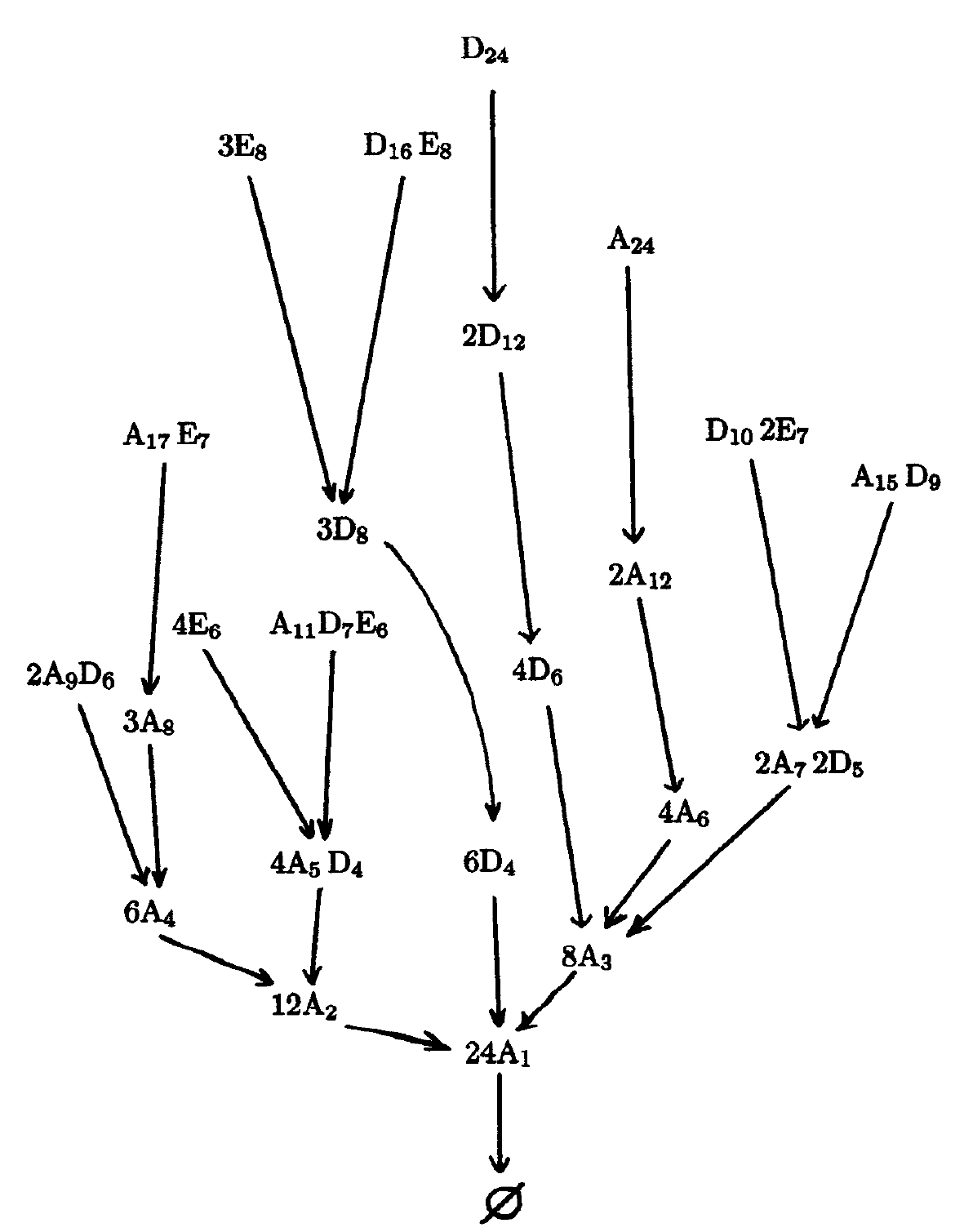} 
\caption{Syst\`eme de racines
du $2$-voisin associ\'e \`a un vecteur de Weyl, d’un r\'eseau unimodulaire
pair de dimension $24$, en fonction du syst\`eme de racines de ce dernier.}
\label{2voisrho}
\end{figure}

Dans le chapitre III, nous \'etudions \'egalement les cas limites de l'assertion (i) du th\'eor\`eme ci-dessus (de mani\`ere directe, {\it i.e.} sans faire usage du th\'eor\`eme \ref{thmintro24}). Nous proc\'edons pour cela \`a une analyse d\'etaill\'ee des \'el\'ements $c$ de ${\rm C}_{\rm L}(\F_p)$ v\'erifiant ${\rm vois}_p(L; c) \simeq {\rm Leech}$, $L$ \'etant un r\'eseau de Niemeier de syst\`eme de racines $R={\rm R}(L)$ non vide. Il est n\'ecessaire \`a la pertinence des \'enonc\'es d'\'etudier plus g\'en\'eralement les $d$-voisins de $L$, o\`u $d \geq 1$ est un entier non n\'ecessairement premier (\S III.1). Nous d\'emontrons que si $\rho$ est un vecteur de Weyl de $R$, et si $h={\rm h}(R)$, on a des isom\'etries (Th\'eor\`eme III.4.2.10)
\begin{equation} \label{holyconstvois} {\rm vois}_h(L ; \rho) \simeq {\rm vois}_{h+1}(L; \rho) \simeq {\rm Leech}. \end{equation}
Cela a un sens car on a $\rho \in L$ (Borcherds) et ${\rm q}_L(\rho)=h(h+1)$ (Venkov). Cet \'enonc\'e contient par exemple l'observation de Thompson susmentionn\'ee. En fait, ces $23$ (ou $46$) constructions du r\'eseau de Leech ne sont autres que les fameuses {\it holy constructions} de Conway et Sloane \cite{conwaysloane23leech}. Nous donnons cependant une nouvelle d\'emonstration des isom\'etries \eqref{holyconstvois} par des m\'ethodes de voisins, et montrons les identit\'es 
\begin{equation} \label{nbholyconstvois} {\rm N}_h(L,{\rm Leech}) = \frac{|W|}{\varphi(h) g}\, \, \, \, {\rm et}\, \, \, \, {\rm N}_{h+1}(L, {\rm Leech}) = \frac{|W|}{h+1},\end{equation}
o\`u $W$ d\'esigne le groupe de Weyl de $R$, et $g^2$ son indice de connexion au sens de Bourbaki. Nous terminons par une \'etude de ${\rm vois}_2(L; \rho)$ inspir\'ee de r\'esultats de Borcherds (Figure \ref{2voisrho}). \ps \ps

\section{S\'eries th\^eta et formes modulaires de Siegel} \label{introthetasiegel}
	
	Revenons \`a la question de d\'eterminer l'op\'erateur ${\rm T}_p$
sur $\Z[{\rm X}_n]$. Quelques observations simples s'imposent : les ${\rm T}_p$ commutent et sont auto-adjoints sur
$\R[{\rm X}_n]$ pour un produit scalaire ad\'equat \cite{nebevenkov} (\S III.2). Il s'agit donc d'en d\'eterminer
une base de vecteurs propres communs ainsi que les collections de valeurs propres des ${\rm
T}_p$ associ\'ees.  La seule droite propre \'evidente, engendr\'ee par
$\sum_{L \in {\rm X}_n} \frac{[L]}{|{\rm O}(L)|}$, est de valeur propre ${\rm
c}_n(p)$ sous ${\rm T}_p$.  \ps \ps

	Nous sommes en fait en pr\'esence d'un probl\`eme d\'eguis\'e de
th\'eorie spectrale des formes automorphes.  En effet, si $G={\rm O}_n$
d\'esigne le (sch\'ema en) groupe orthogonal sur $\Z$ d\'efini par la forme
quadratique ${\rm q}_{{\rm E}_n}$, et $\AAA$ l'anneau des ad\`eles de $\Q$,
des arguments de th\'eorie du genre conduisent \`a un isomorphisme de
$G(\R)$-ensembles $\mathcal{L}_n \simeq G(\Q) \backslash
G(\AAA)/G(\widehat{\Z})$\, (\S II.2, \S \ref{ensclass}). Il s'ensuit que le dual de $\R[{\rm X}_n]$
s'identifie \`a l'espace des fonctions \`a valeurs r\'eelles
sur $G(\Q) \backslash G(\AAA)$ qui sont invariantes \`a droite par $G(\R) \times
G(\widehat{\Z})$. Dans cette description, l'op\'erateur ${}^t\, {\rm T}_p$ appara\^it comme un
\'el\'ement particulier de l'anneau ${\rm H}(G)$ des op\'erateurs de Hecke
de $G$.  \ps \ps
	
	Ces observations tout \`a fait classiques sont rappel\'ees dans le
chapitre \ref{chap4} du m\'emoire.  Bien que nous nous int\'eressons principalement
aux formes automorphes du $\Z$-groupe ${\rm O}_n$, nos \'enonc\'es et d\'emonstrations en feront
intervenir diverses d\'eclinaisons (formes automorphes pour ${\rm SO}_n$, ${\rm PGO}_n$ et ${\rm
PGSO}_n$), ainsi que des formes modulaires pour ${\rm SL}_2(\Z)$, des formes
modulaires de Siegel \`a valeurs vectorielles pour ${\rm Sp}_{2g}(\Z)$, et
m\^eme, par le biais des r\'esultats d'Arthur, des formes automorphes pour
${\rm PGL}_n$. Il sera donc n\'ecessaire d'adopter d\`es le d\'epart un
point de vue suffisamment g\'en\'eral englobant
tous ces objets (\S \ref{paysageaut}). Le lecteur trouvera aux \S \ref{ensclass} et \S \ref{corrhecke} un expos\'e
\'el\'ementaire sur les op\'erateurs de Hecke.  L'accent est mis
sur les exemples fournis par les groupes classiques et leurs variantes
(Hecke, Satake, Shimura), ils donnent plus de largeur de vue sur les
$p$-voisins et leurs g\'en\'eralisations.  Les \S \ref{cascompact} et \S \ref{formesiegel} sont
consacr\'es \`a des rappels sur les formes automorphes pour ${\rm O}_n$ et
les formes modulaires
de Siegel. Soulignons que la r\'edaction de ce chapitre s'adresse au
non-sp\'ecialiste, et pr\'etend \`a peu d'originalit\'e.  \ps \ps 

Une approche pour \'etudier le ${\rm H}({\rm O}_n)$-module $\Z[{\rm X}_n]$
consiste \`a examiner les s\'eries th\^eta de Siegel de tout genre $g \geq
1$ des \'el\'ements de $\mathcal{L}_n$.  Pour tout $n \equiv 0 \bmod 8$ et
$g \geq 1$, elles fournissent une application lin\'eaire $$\vartheta_g :
\C[{\rm X}_n] \rightarrow {\rm M}_{\frac{n}{2}}({\mathrm{Sp}}_{2g}(\Z)),\,\, [L]
	\mapsto \vartheta_g(L),$$ o\`u
$\mathrm{M}_k({\mathrm{Sp}}_{2g}(\Z))$ d\'esigne l'espace des formes
modulaires de Siegel de poids $k \in \Z$ pour le groupe ${\rm Sp}_{2g}(\Z)$
(\S \ref{thetasiegel}).  La pertinence de cette application pour notre probl\`eme vient
des relations de commutation d'Eichler g\'en\'eralis\'ees ; elles montrent
que $\vartheta_g$ entrelace chaque \'el\'ement de ${\rm H}({\rm O}_n)$ avec
un certain \'el\'ement ``explicite'' de ${\rm H}({\rm Sp}_{2g})$ (Eichler,
Freitag, Andrianov, \S \ref{thetasiegel}).  L'application $\vartheta_g$ est trivialement injective pour $g \geq n$.
La question de d\'eterminer la structure du 
${\rm H}({\rm Sp}_{2g})$-module ${\rm M}_{k}({\mathrm{Sp}}_{2g}(\Z))$ est en revanche notoirement difficile, et ce d'autant plus que $g$
est grand. N\'eanmoins, nous d\'evelopperons au chapitre IX une strat\'egie permettant de r\'esoudre des nouveaux cas de cette question, qui repose entre autres sur les r\'esultats r\'ecents d'Arthur \cite{arthur}.   \ps\ps

L'application $\vartheta_g$ a fait l'objet de nombreux travaux. Son noyau, qui d\'ecro\^it quand $g$ grandit, d\'ecrit les relations lin\'eaires entre les s\'eries th\^eta de genre $g$ des \'el\'ements de $\mathcal{L}_n$, et la d\'etermination de son image est un exemple du c\'el\`ebre {\it probl\`eme de la base} d'Eichler. Pr\'ecis\'ement, $\vartheta_g$ induit une application injective 
\begin{equation} \label{pbbaseeichler} {\rm Ker} \, \vartheta_{g-1} / {\rm Ker}\, \vartheta_g \longrightarrow {\rm S}_{\frac{n}{2}}({\rm Sp}_{2g}(\Z)) \end{equation}
o\`u ${\rm S}_k({\rm Sp}_{2g}(\Z)) \subset {\rm M}_k({\rm Sp}_{2g}(\Z))$ d\'esigne le sous-espace des formes paraboliques (voir le \S \ref{thetasiegel} pour la convention sur $\vartheta_0$), et il s'agit de savoir si elle est surjective. Un r\'esultat important de B\"ocherer \cite{bochererbowdoin} donne une condition n\'ecessaire et suffisante pour qu'une forme propre pour ${\rm H}({\rm Sp}_{2g})$ soit dans son image en termes de l'annulation en un certain entier d'une fonction ${\rm L}$ associ\'ee (\S \ref{retourchap4}). \ps \ps
\bigskip

\noindent
{\sc Cas de la dimension $n=16$}.\ps \ps \medskip

Le cas de la dimension $16$ est l'objet d'une histoire fameuse, rappell\'ee au \S \ref{casdim16}. En effet, un r\'esultat classique de Witt et Igusa assure que l'on a
\begin{equation} \label{introidwitt} \vartheta_g(\mathrm{E}_8 \oplus \mathrm{E}_8)=
\vartheta_g(\mathrm{E}_{16}) \, \, \, \, {\rm
si} \, \, \, g\leq 3.\end{equation}
Ces identit\'es remarquables affirment que $\mathrm{E}_8 \oplus
\mathrm{E}_8$ et ${\rm E}_{16}$ repr\'esentent toutes les formes quadratiques enti\`eres positives de rang
$\leq 3$ exactement le m\^eme nombre de fois. L'annulation bien connue ${\rm S}_8({\rm SL}_2(\Z))=0$ traite le cas $g=1$ (et conduit aux tores isospectraux de Milnor d\'ej\`a mentionn\'es). En particulier, ``le'' vecteur propre non trivial de $\Z[{\rm X}_{\rm 16}]$ est $[{\rm E}_{16}]-[{\rm E}_8 \oplus {\rm E}_8]$. La difficult\'e en genres $2$ et $3$ est que l'annulation de ${\rm S}_8({\rm Sp}_{2g}(\Z))$, bien que toujours vraie, est plus d\'elicate \`a d\'emontrer. Dans l'appendice {\rm A}, nous exposons une autre d\'emonstration ing\'enieuse des identit\'es \eqref{introidwitt} due \`a Kneser, qui ne repose pas sur de telles annulations. \ps \ps

La forme $J=\vartheta_4(\mathrm{E}_8 \oplus \mathrm{E}_8)-\vartheta_4(\mathrm{E}_{16})$, qui n'est autre que la fameuse {\it forme de Schottky}, est en revanche non nulle. D'apr\`es Poor et Yuen \cite{py}, on sait m\^eme qu'elle engendre ${\rm S}_8({\rm Sp}_8(\Z))$. Le th\'eor\`eme \ref{thmintro16} r\'esulte alors de la r\'esolution par Ikeda \cite{ikeda1} de la  conjecture de Duke-Imamo\u{g}lu \cite{bkuss}. En effet, appliqu\'ee \`a la forme modulaire $\Delta$ dans ${\rm S}_{12}({\rm SL}_2(\Z))$, cette derni\`ere assure l'existence d'une forme modulaire de Siegel non nulle dans ${\rm S}_8({\rm Sp}_8(\Z))$ qui est propre pour les op\'erateurs de Hecke dans ${\rm H}({\rm Sp}_8)$, de valeurs propres explicitement d\'etermin\'ees par les $\tau(p)$. La d\'emonstration d'Ikeda est difficile ; une des contributions de ce m\'emoire est d'avoir donn\'e une autre d\'emonstration de son th\'eor\`eme, tr\`es diff\'erente, dans le cas particulier ci-dessus. \ps \ps

Le r\'esultat important est le suivant. Pour toute fonction $f : \mathcal{L}_n \rightarrow \C$, on d\'efinit ${\rm T}_p(f) : \mathcal{L}_n \rightarrow \C$ en posant pour tout $L \in \mathcal{L}_n$, ${\rm T}_p(f)(L) = \sum_M f(M)$, la somme portant sur les $p$-voisins de $L$. Si $1 \leq g \leq n/2$, on note ${\rm H}_{d,g}(\R^n)$ l'espace des polyn\^omes $(\R^n)^g \rightarrow \C$ qui sont harmoniques pour le laplacien euclidien de $(\R^n)^g$, et qui satisfont $P \circ \gamma = (\DET \gamma)^d \,P$ pour tout $\gamma \in \GL_g(\C)$ {\rm (\S \ref{i4dpartrialite})}. C'est une repr\'esentation lin\'eaire de ${\rm O}(\R^n)$ de mani\`ere naturelle.

\begin{thmintro}\label{variationtrialite} Soient $q+\sum_{n \geq 2} a_n q^n$ une forme modulaire de poids $k$ pour ${\rm SL}_2(\Z)$ propre pour les op\'erateurs de Hecke, et $d=\frac{k}{2}-2$. Il existe une fonction $f : \mathcal{L}_8 \rightarrow \C$ telle que :\ps \ps \begin{itemize}
\item[(i)] pour tout nombre premier $p$, on a \,\,${\rm T}_p (f) = p^{-d}\, \frac{p^4-1}{p-1}\, a_p \, f$, \ps \ps
\item[(ii)] $f$ engendre sous ${\rm O}(\R^8)$ une repr\'esentation isomorphe \`a ${\rm H}_{d,4}(\R^8)$.
\end{itemize}
\end{thmintro}

Le \S \ref{i4dpartrialite} est principalement consacr\'e \`a la d\'emonstration d'un cas particulier de ce th\'eor\`eme dans le cas $k=12$, qui conduit \`a une d\'emonstration compl\`ete et relativement \'el\'ementaire du th\'eor\`eme \ref{thmintro16}. Le cas g\'en\'eral ci-dessus sera d\'emontr\'e et pr\'ecis\'e au \S \ref{retourchap4}. \ps \ps

Donnons une id\'ee de la d\'emonstration. Nous commen\c{c}ons par r\'ealiser la forme modulaire de d\'epart comme s\'erie th\^eta $\sum_{x \in {\rm E}_8} P(x)\,  q^{\frac{x\cdot x}{2}}$ o\`u $P : \R^8 \rightarrow \C$ est un polyn\^ome harmonique bien choisi. Dans le cas de $\Delta$, ``l'invariant'' harmonique de degr\'e $8$ de ${\rm W}({\bf E}_8)$ convient, et en g\'en\'eral on invoque un r\'esultat de Waldspurger \cite{waldspurger79}. Cette construction d\'efinit un sous-espace des fonctions $\mathcal{L}_8 \rightarrow \C$ qui sont propres pour les op\'erateurs de Hecke dans ${\rm H}({\rm O}_8)$, de valeurs propres reli\'ees aux $a_p$ par les relations d'Eichler, et qui engendrent une repr\'esentation isomorphe \`a ${\rm H}_{8,1}(\R^8)$ sous l'action de ${\rm O}(\R^8)$. L'id\'ee principale consiste \`a leur appliquer \`a la source un automorphisme d'ordre $3$ de $\mathcal{L}_8$ issu de la trialit\'e. Cet automorphisme est construit \`a partir d'une structure d'octonions de Coxeter sur ${\rm E}_8$ et de l'identit\'e $\mathcal{L}_8 \simeq G(\Q) \backslash G(\AAA) / G(\widehat{\Z})$ avec $G={\rm PGSO}_8$. Les fonctions obtenues conviennent : nous renvoyons au \S \ref{i4dpartrialite} pour les d\'emonstrations. \ps \ps

La condition (ii) de l'\'enonc\'e assure que la fonction $f$ donne naissance \`a une s\'erie th\^eta de Siegel de genre $4$ (\`a coefficients ``pluriharmoniques''). Quand elle est non nulle, cette s\'erie th\^eta est un substitut du rel\`evement d'Ikeda de genre $4$ de la forme modulaire de d\'epart. On v\'erifie cette non-annulation quand $k=12$ ; on d\'eduit ais\'ement le th\'eor\`eme \ref{thmintro16}. \ps \ps

Mentionnons enfin que nous d\'emontrerons plus tard (Th\'eor\`eme IX.\ref{thmsiegelsmallweight}) l'annulation ${\rm S}_8({\rm Sp}_{2g}(\Z))=0$ pour tout $g \neq 4$. Quand $g=5, 6$ elle avait d\'ej\`a \'et\'e obtenue par Poor et Yuen \cite{py2}. Par cons\'equent, l'application $\vartheta_g : \C[{\rm X}_{16}] \rightarrow {\rm M}_8({\rm Sp}_{2g}(\Z))$ est surjective pour tout genre $g\geq 1$.

\ps \medskip \noindent
{\sc Cas de la dimension $24$}
\ps  \medskip

Ce cas a fait l'objet de travaux remarquables d'Erokhin \cite{erokhin}, Bor\-cherds-Freitag-Weissauer \cite{bfw} et Nebe-Venkov \cite{nebevenkov} (\S \ref{thetaniemeier}). Erokhin a montr\'e ${\rm Ker} \, \vartheta_{12}=0$, et les trois auteurs de \cite{bfw} ont prouv\'e que ${\rm Ker} \, \, \vartheta_{11}$ est de dimension $1$. Une \'etude fine de toute la filtration de $\Z[{\rm X}_{24}]$ par les ${\rm Ker}\, \vartheta_g$ a \'et\'e entreprise par Nebe et Venkov {\it loc. cit.} Leur point de d\'epart est une explicitation de l'op\'erateur ${\rm T}_2$ sur $\Z[{\rm X}_{24}]$ qu'ils d\'eduisent de r\'esultats de Borcherds (\S III.3.3); cet op\'erateur \'etant \`a valeurs propres distinctes et enti\`eres, ils obtiennent une base de $\Q[{\rm X}_{24}]$ constitu\'ee de vecteurs propres des ${\rm T}_p$ dans $\Q[{\rm X}_{24}]$. Ils proposent \'egalement une conjecture pour la dimension de l'image de l'application \eqref{pbbaseeichler} pour tout entier $1 \leq g \leq 12$, dont ils d\'emontrent des cas particuliers. Nous \'etablissons leur conjecture et montrons que le probl\`eme de la base d'Eichler admet une r\'eponse positive en dimension $n=24$, pour tout genre $1 \leq g \leq 23$ (Th\'eor\`eme IX.\ref{mainthmsiegel12} et Corollaire IX.\ref{dimm12sp2g}) : \ps \bigskip

\begin{thmintro}\label{pbeichler24} L'application $\vartheta_g : \C[{\rm X}_{24}] \rightarrow {\rm M}_{12}({\rm Sp}_{2g}(\Z))$ est surjective, et induit un isomorphisme \,$ {\rm Ker} \, \vartheta_{g-1} / {\rm Ker} \, \vartheta_g \isomo {\rm S}_{12}({\rm Sp}_{2g}(\Z))$, pour tout entier $g \leq 23$. La dimension de ${\rm S}_{12}({\rm Sp}_{2g}(\Z))$ pour $g \neq 24$ est donn\'ee par la table :
\begin{table}[htp]
\renewcommand{\arraystretch}{1.5}
\begin{center}
\begin{tabular}{|c||c|c|c|c|c|c|c|c|c|c|c|c|c|}
\hline $g$ & $1$ & $2$ & $3$ & $4$ & $5$ & $6$ & $7$ & $8$ & $9$ & $10$ & $11$ & $12$ & {\footnotesize $>12$}\cr
\hline $\dim\, {\rm S}_{12}({\rm Sp}_{2g}(\Z))$ & $1$ & $1$ & $1$ & $2$ & $2$ & 
$3$ & $3$ & $4$ & $2$ & $2$ & $1$ & $1$ & $0$ \cr
\hline
\end{tabular}
\end{center}
\end{table}
\end{thmintro}

Nous donnerons une id\'ee de la d\'emonstration au \S \ref{introdem}, le point le plus d\'elicat \'etant la premi\`ere assertion. Il r\'esulte du th\'eor\`eme une description compl\`ete de la filtration des ${\rm Ker} \, \vartheta_g$ sur $\Z[{\rm X}_{24}]$.  Mentionnons que le probl\`eme de la base d'Eichler admet une r\'eponse n\'egative en dimension $n=32$ et genre $g=14$, comme nous l'observons au corollaire VII.\ref{cex32}.
	
\section{Formes automorphes des groupes classiques}\label{fautintro}

	Les formes modulaires de Siegel, ainsi que les espaces de formes
automorphes pour ${\rm O}_n$, peuvent \^etre \'etudi\'es \`a l'aide des
travaux r\'ecents d'Arthur \cite{arthur}.  La formulation m\^eme des
\'enonc\'es n\'ecessite cependant des rappels pr\'eliminaires sur le point
de vue de Langlands sur la th\'eorie des formes automorphes
\cite{langlandspb} \cite{borelcorvallis}, que nous avons regroup\'es dans le
chapitre \ref{chap5}.  \ps \ps
	
	Soit $G$ un sch\'ema en groupes semi-simple\footnote{La discussion
qui suit ne s'y applique pas verbatim \`a certains $\Z$-groupes naturels
ici, comme ${\rm O}_n$ ou ${\rm PGO}_n$, qui sont non connexes.  Nous
pr\'ecisons dans le texte les modifications \`a apporter afin de les
englober, mais nous ignorerons ce d\'etail dans cette introduction.} sur
$\Z$.  On note $\Pi_{\rm disc}(G)$ l'ensemble des sous-repr\'esentations
topologiquement irr\'eductibles de l'espace des fonctions de carr\'e
int\'egrable sur $G(\Q)\backslash G(\AAA)/G(\widehat{\Z})$ pour les actions
naturelles de $G(\R)$ et de l'anneau ${\rm H}(G)$ des op\'erateurs de Hecke
de $G$ (\S \ref{paysageaut}).  L'isomorphisme de Satake associe \`a tout $\pi \in
\Pi_{\rm disc}(G)$ et tout premier $p$ une classe de conjugaison semi-simple
${\rm c}_p(\pi)$ dans $\widehat{G}(\C)$, o\`u $\widehat{G}$ d\'esigne le
groupe semi-simple complexe dual de Langlands de $G_\C$ (\S \ref{rapredchap6},
\S \ref{parametrisationsatake}).  Ce point de vue \'eclairant sur les
valeurs propres des op\'erateurs de Hecke, d\^u \`a Langlands, est
explicit\'e au \S VI.\ref{exemplehecke} dans le cas des groupes classiques
et des op\'erateurs d'int\^eret pour ce m\'emoire, en suivant l'expos\'e de
Gross \cite{grossatake}.  De m\^eme, nous rappelons comment l'isomorphisme
d'Harish-Chandra associe au caract\`ere infinit\'esimal de la composante
archim\'edienne $\pi_\infty$ de $\pi$ une classe de conjugaison semi-simple
${\rm c}_\infty(\pi)$ dans l'alg\`ebre de Lie de $\widehat{G}$ (\S \ref{secisohc}).  \ps \ps
	
	Une conjecture centrale et structurante, initialement due \`a Langlands et
pr\'ecis\'ee par Arthur \cite{arthurunipotent}, est que ces collections de
classes de conjugaison s'expriment toutes en terme des donn\'ees similaires
relatives aux $\Pi_{\rm disc}(\PGL_m)$ pour $m \geq 1$.  Concr\`etement,
fixons $\pi \in \Pi_{\rm disc}(G)$ et $r : \widehat{G}(\C) \rightarrow {\rm
SL}_n(\C)$ une repr\'esentation alg\'ebrique, de sorte que l'on dispose
suivant Langlands d'un produit eul\'erien ${\rm L}(s,\pi,r)=\prod_p \DET(1-
p^{-s} \,r({\rm c}_p(\pi)))^{-1}$, qui est absolument convergent pour tout nombre complexe $s$ de partie r\'eelle assez
grande.  Quand $G$ est le $\Z$-groupe ${\rm PGL}_m$, et que $r$ est la
repr\'esentation tautologique de $\widehat{G}={\rm SL}_m$, on pose simplement
${\rm L}(s,\pi)={\rm L}(s,\pi,r)$.  La conjecture d'Arthur-Langlands pour le couple
$(\pi,r)$ pr\'edit l'existence d'un entier $k\geq 1$, et pour $i=1,\dots,k$
d'une repr\'esentation\!\!\footnote{On note traditionnellement $\Pi_{\rm
cusp}(G) \subset \Pi_{\rm disc}(G)$ le sous-ensemble des repr\'esentations
intervenant dans le sous-espace des formes automorphes paraboliques \cite{GGPS} (\S
\ref{paysageaut}).} $\pi_i\in \Pi_{\rm cusp}(\PGL_{n_i})$ et d'un entier
$d_i\geq 1$, tels que l'on ait l'\'egalit\'e (voir le \S VI.\ref{parconjarthlan} pour une autre formulation)
\begin{equation}\label{idalintro} {\rm L}(s,\pi,r) = \prod_{i=1}^k
\prod_{j=0}^{d_i-1} {\rm L}(s+j-\frac{d_i-1}{2},\pi_i).\end{equation} De
mani\`ere l\'eg\`erement abusive, la collection des classes de conjugaison
$r({\rm c}_v(\pi))$ sera appel\'ee {\it param\`etre de Langlands} du couple
$(\pi,r)$ ; elle sera not\'ee $\psi(\pi,r)$.  Lorsque l'\'egalit\'e \eqref{idalintro} est v\'erifi\'ee, elle sera
not\'ee symboliquement\footnote{Au sens strict, notre notation inclus
l'identit\'e naturelle correspondante sur les caract\`eres infinit\'esimaux
(\S VI.\ref{parconjarthlan}). De plus, le facteur $\pi_i[d_i]$ sera simplement not\'e $[d_i]$ (resp. $\pi_i$) lorsque $n_i=1$ (resp. $d_i=1$). Ces conventions sont utilis\'ees dans la table \ref{table24}.} $$\psi(\pi,r) = \oplus_{i=1}^k \pi_i[d_i].$$
Si $G$ est un {\it $\Z$-groupe classique} (\S VI.4.7, \S VIII.1), alors
$\widehat{G}$ est un groupe classique complexe et vient avec une
repr\'esentation tautologique appel\'ee {\it repr\'esentation standard},
not\'ee ${\rm St}$.  Un r\'esultat important d\'emontr\'e par Arthur
\cite{arthur} est que la conjecture ci-dessus est vraie si l'on suppose que
$G$ est soit ${\rm Sp}_{2g}$, soit un $\Z$-groupe sp\'ecial orthogonal
d\'eploy\'e, et si $r={\rm St}$.  \ps \ps

Au chapitre \ref{chap7}, nous illustrons ces th\'eories en donnant de
nombreux exemples de cas particuliers de la conjecture d'Arthur-Langlands,
concernant les formes automorphes pour ${\rm SO}_n$ ou les formes de Siegel
pour ${\rm Sp}_{2g}(\Z)$.  Ils ne reposent pas sur les travaux d'Arthur,
mais sur des constructions plus classiques de s\'eries th\^eta.  Nous
rappelons le point de vue de Rallis sur les relations d'Eichler (\S
\ref{retoureichler}) ainsi que des r\'esultats importants de
B\"ocherer et d'Ikeda.  Nous d\'emontrons le th\'eor\`eme \ref{variationtrialite} et
donnons d'autres applications de la trialit\'e \`a la construction de
certains \'el\'ements de $\Pi_{\rm disc}({\rm SO}_8)$ (\S
\ref{retourchap4}).  Un ingr\'edient des d\'emonstrations est un l\'eger
raffinement de l'identit\'e de Rallis \`a la paire $({\rm PGO}_n,{\rm
PGSp}_{2g})$ (\S VII.\ref{rallisspin}).  Au final, notre analyse recouvre
suffisamment de constructions pour permettre par exemple de d\'eterminer
$\psi(\pi,{\rm St})$ pour $13$ des $16$ ``premi\`eres'' repr\'esentations $\pi
\in \Pi_{\rm disc}({\rm SO}_8)$ (\S \ref{tableexempleso8}).  \ps \ps

Nous pouvons maintenant \'enoncer l'analogue pour ${\rm X}_{24}$ du
th\'eor\`eme \ref{thmintro16}. Nous renvoyons au \S \ref{repgal24} pour une formulation de
ce th\'eor\`eme en terme de repr\'esentations de ${\rm
Gal}(\overline{\Q}/\Q)$, dans l'esprit de l'annonce \cite{chlannes}. 

\begin{thmintro}\label{thmintro24} Les
param\`etres $\psi(\pi,{\rm St})$ des $24$ repr\'esentations $\pi \in
\Pi_{\rm disc}({\rm O}_{24})$ engendr\'ees par les fonctions ${\rm X}_{24}
\rightarrow \C$ propres pour ${\rm H}({\rm O}_{24})$ sont :
\begin{table}[htp] \hspace{-1cm} \renewcommand{\arraystretch}{1.5} {\footnotesize
\begin{tabular}{cc}

{{}} $[23]\oplus [1]$ & ${\rm Sym}^2 \Delta_{11} \oplus \Delta_{17}[4] \oplus \Delta_{11}[2] \oplus [9]$ \cr

{{}} ${\rm Sym}^2 \Delta_{11} \oplus [21]$ & ${\rm Sym}^2 \Delta_{11} \oplus \Delta_{15}[6] \oplus [9]$ \cr

{{}} $\Delta_{21}[2] \oplus [1] \oplus [19]$ & $\Delta_{15}[8] \oplus [1] \oplus [7]$ \cr

{{}} ${\rm Sym}^2 \Delta_{11} \oplus \Delta_{19}[2] \oplus [17]$ & $\Delta_{21}[2] \oplus \Delta_{17}[2] \oplus \Delta_{11}[4] \oplus [1] \oplus [7]$ \cr

{{}} $\Delta_{21}[2] \oplus \Delta_{17}[2] \oplus [1] \oplus [15]$ & $\Delta_{19}[4] \oplus \Delta_{11}[4] \oplus [1] \oplus [7]$\cr

{{}} $\Delta_{19}[4] \oplus [1] \oplus [15]$ & $\Delta_{21,9}[2] \oplus \Delta_{15}[4] \oplus [1] \oplus [7]$\cr

{{}} ${\rm Sym}^2 \Delta_{11} \oplus \Delta_{19}[2] \oplus \Delta_{15}[2]\oplus [13]$ &${\rm Sym}^2 \Delta_{11} \oplus \Delta_{19}[2] \oplus \Delta_{11}[6] \oplus [5]$\cr

{{}} ${\rm Sym}^2 \Delta_{11} \oplus \Delta_{17}[4] \oplus [13]$ & ${\rm Sym}^2 \Delta_{11} \oplus \Delta_{19,7}[2] \oplus \Delta_{15}[2] \oplus \Delta_{11}[2] \oplus [5]$\cr

{{}} $\Delta_{17}[6] \oplus [1] \oplus [11]$ & $\Delta_{21}[2] \oplus \Delta_{11}[8] \oplus [1] \oplus [3]$\cr

{{}} $\Delta_{21}[2] \oplus \Delta_{15}[4] \oplus [1] \oplus[11]$ & $\Delta_{21,5}[2] \oplus \Delta_{17}[2] \oplus \Delta_{11}[4] \oplus [1] \oplus [3]$\cr

{{}} $\Delta_{21,13}[2] \oplus \Delta_{17}[2] \oplus [1] \oplus [11]$ & ${\rm Sym}^2 \Delta_{11} \oplus \Delta_{11}[10] \oplus [1]$\cr

{{}} ${\rm Sym}^2 \Delta_{11} \oplus \Delta_{19}[2] \oplus \Delta_{15}[2] \oplus \Delta_{11}[2] \oplus [9]$ & $\Delta_{11}[12]$ \cr
{{}}
\end{tabular} \ps
}
\caption{{\small Les param\`etres standards des $\pi \in \Pi_{\rm disc}({\rm O}_{24})$ telles que $\pi_\infty=1$. }}
\label{table24}
\end{table}
\end{thmintro}

Pr\'ecisons les notations intervenant dans l'\'enonc\'e ci-dessus. La repr\'esentation $\Delta_w$, pour $w \in \{11,15,17,19,21\}$, d\'esigne l'\'el\'ement de
$\Pi_{\rm cusp}({\rm PGL}_2)$ engendr\'e par la droite ${\rm S}_{w+1}({\rm
SL}_2(\Z))$ des formes modulaires paraboliques de poids $w+1$ pour ${\rm SL}_2(\Z)$. La repr\'esentation ${\rm Sym}^2 \Delta_w$, caract\'eris\'ee par
l'identit\'e ${\rm c}_v({\rm Sym}^2 \Delta_w)\,=\,{\rm Sym}^2 \, {\rm c}_v(\Delta_w)$ pour toute place $v$ de $\Q$, est le
carr\'e sym\'etrique de Gelbart-Jacquet de $\Delta_w$ \cite{GJ}. C'est un \'el\'ement de $\Pi_{\rm cusp}(\PGL_3)$. Enfin, les $4$ repr\'esentations
$\Delta_{w,v}$ sont des \'el\'ements de $\Pi_{\rm cusp}(\PGL_4)$ d\'efinis et \'etudi\'es en d\'etail au \S \ref{prelimsp4}. Entre autres, le caract\`ere infinit\'esimal ${\rm c}_\infty(\Delta_{w,v})$, qui est par d\'efinition la classe de conjugaison d'un \'el\'ement semi-simple dans ${\rm M}_4(\C)$, a pour valeurs propres $\{\pm \frac{w}{2}, \pm \frac{v}{2}\}$, ce qui explique notre notation. De m\^eme, ${\rm c}_\infty(\Delta_w)$ a pour valeurs propres $\pm \frac{w}{2}$. \ps\ps

Soulignons que dans un travail remarquable \cite{ikeda2}, Ikeda \'etait parvenu \`a d\'eterminer
$20$ des $24$ param\`etres ci-dessus, \`a savoir ceux ne contenant pas les repr\'esentations $\Delta_{w,v}$. \ps\ps

 \'Etant donn\'e l'importance du r\^ole jou\'e par les $\Delta_{w,v}$ dans ce m\'emoire, il s'impose d'indiquer bri\`evement leur origine. Soit $(j,k)$ l'un des $4$ couples $(6,8), (4,10), (8,8)$ ou $(12,6)$. Une formule de dimension due \`a R. Tsushima \cite{tsushima} montre que l'espace des formes modulaires de Siegel pour ${\rm Sp}_4(\Z)$ qui sont paraboliques et \`a coefficients dans ${\rm Sym}^j \otimes \DET^k$ est de dimension $1$. Nous en donnons un g\'en\'erateur concret par une construction de s\'erie th\^eta \`a coefficients ``pluriharmoniques'' bas\'ee sur le r\'eseau ${\rm E}_8$. Si $\pi_{j,k}$ d\'esigne l'\'el\'ement de $\Pi_{\rm cusp}({\rm PGSp}_4)$ engendr\'e par cette forme propre, on a alors la relation $\psi(\pi_{j,k},{\rm St})=\Delta_{w,v}$ avec $(w,v)=(2j+k-3,j+1)$. Remarquons que ${\rm PGSp}_4$ est isomorphe au $\Z$-groupe classique d\'eploy\'e ${\rm SO}_{3,2}$, de groupe dual le groupe complexe ${\rm Sp}_4$, de sorte que la th\'eorie d'Arthur s'applique \`a $(\pi_{j,k},{\rm St})$. \ps \ps
 
Le th\'eor\`eme \ref{thmintro24} sera d\'emontr\'e au \S X.\ref{demothmE}, par une m\'ethode que nous d\'ecrirons au \S \ref{introdem}.
N\'eanmoins, nous en donnerons d'abord deux autres d\'emonstrations {\it conditionnelles} aux \S X.\ref{preuve1thm24} et \S X.\ref{preuve2thm24}. Ces d\'emonstrations, obtenues en appliquant la {\it formule de multiplicit\'e} d'Arthur \cite{arthur}, seront \`a terme les plus naturelles, mais elles d\'ependent actuellement de certains \'enonc\'es pr\'ecisant ceux d'Arthur, qui bien qu'attendus ne sont pas encore disponibles. \ps \ps

Au chapitre \ref{classarthur}, nous revenons donc sur les r\'esultats g\'en\'eraux
d'Arthur \cite{arthur}, que nous sp\'ecifions au cas des groupes classiques $G$ sur $\Z$ et de
leurs formes automorphes ``non ramifi\'ees partout'' ; une telle analyse a
d\'ej\`a \'et\'e men\'ee \`a bien dans \cite[\S 3]{chrenard2}, que nous
reprenons et compl\'etons. L'essentiel du chapitre est consacr\'e \`a
l'explicitation de la fameuse formule de multiplicit\'e sus-mentionn\'ee. Il
s'agit d'une condition n\'ecessaire et suffisante pour qu'une collection
donn\'ee de $(\pi_i,d_i)$ ``provienne'' d'un $\pi \in \Pi_{\rm disc}(G)$, au sens o\`u $\psi(\pi,{\rm
St})=\oplus_{i=1}^k \pi_i[d_i]$, tout en ayant
une composante $\pi_\infty$ prescrite (\S \ref{paramarthur}).  Nous nous limitons au cas
o\`u $\pi_\infty$ est une s\'erie discr\`ete de $G(\R)$ et expliquons de
mani\`ere concr\`ete la param\'etrisation de ces derni\`eres par Shelstad (\S \ref{paramdisc}).
C'est cette param\'etrisation qui intervient dans la formule d'Arthur.  La version de cette formule que nous donnons n'est \`a pr\'esent d\'emontr\'ee que si $G$ est d\'eploy\'e sur $\Z$ et si
tous les entiers $d_i$ sont \'egaux \`a $1$.  Nous la discutons cependant en
g\'en\'eral, en pr\'ecisant les conjectures desquelles divers cas
particuliers d\'ependent (\S VIII.\ref{parconjajenonce}), car elle \'eclaire grandement les constructions particuli\`eres \'etudi\'ees tout au long de ce m\'emoire.  Dans le cas
$G={\rm SO}_n$ ou des formes de Siegel pour ${\rm Sp}_{2g}(\Z)$, des
formulaires concrets sont d\'egag\'es au \S \ref{formulairesAMF}.  Nous
v\'erifions par exemple qu'ils sont compatibles aux r\'esultats du chapitre
\ref{chap7} et aux r\'esultats de B\"ocherer sur l'image de l'application
\eqref{pbbaseeichler} (\S \ref{comptheta}). Comme promis, nous montrons enfin au \S \ref{so24etnv} que ces formules conduisent \`a une d\'emonstration conditionnelle simple, bien qu'assez miraculeuse, du th\'eor\`eme \ref{thmintro24}.  \ps \ps

\section{Repr\'esentations automorphes alg\'ebriques de petit poids}\label{introclass}	

Soit $m\geq 1$ un entier. Nous disons qu'une repr\'esentation $\pi \in \Pi_{\rm
cusp}(\PGL_m)$ est {\it alg\'ebrique} si les valeurs propres de ${\rm
c}_\infty(\pi)$ sont dans $\frac{1}{2}\Z$ et de diff\'erences deux \`a
deux dans $\Z$ (\S VIII.\ref{paralgreg}). Le double de la plus grande valeur propre de ${\rm c}_\infty(\pi)$, not\'e ${\rm w}(\pi)$, est alors appell\'e {\it poids motivique} de
$\pi$ ; c'est un entier $\geq 0$. Ces repr\'esentations alg\'ebriques ont un int\'er\^et propre, car ce sont pr\'ecis\'ement celles qui sont reli\'ees aux repr\'esentations $\ell$-adiques ``g\'eom\'etriques'' du groupe de Galois absolu de $\Q$, par le yoga de Fontaine-Mazur et Langlands (\S VIII.\ref{repgal}). Leur int\'er\^et dans les questions qui nous pr\'eoccupent ici r\'eside plut\^ot dans l'observation suivante. \ps\ps

Soient $G$ un $\Z$-groupe semi-simple, $\pi \in \Pi_{\rm disc}(G)$ telle que $\pi_\infty$ a m\^eme caract\`ere infinit\'esimal qu'une repr\'esentation alg\'ebrique $V$ de dimension finie de $G(\C)$, et $r : \widehat{G} \rightarrow {\rm SL}_n(\C)$ une repr\'esentation alg\'ebrique. Supposons $\psi(\pi,r)=\oplus_{i=1}^k \pi_i[d_i]$ suivant Arthur et Langlands. Les repr\'esentations $\pi_i$ sont alors alg\'ebriques, de poids motivique born\'e en fonction des plus hauts poids de $V$ et $r$ (\S \ref{repselfdual}). Par exemple, si $G={\rm Sp}_{2g}$ 
et si $\pi \in \Pi_{\rm cusp}({\rm Sp}_{2g})$ est engendr\'ee par une forme modulaire de Siegel propre de poids $k$ pour
${\rm Sp}_{2g}(\Z)$ (avec disons $k>g$, mais cette condition peut \^etre affaiblie), on peut \'ecrire $\psi(\pi,{\rm St}) = \oplus_{i=1}^k \pi_i[d_i]$ d'apr\`es Arthur, et les $\pi_i$ sont alg\'ebriques de poids motivique $\leq 2k-2$. Un ingr\'edient important de nos d\'emonstrations est alors l'\'enonc\'e de classification suivant, d'int\'er\^et ind\'ependant, d\'emontr\'e au \S \ref{mrw}.\ps\ps

\begin{thmintro}\label{introw22} Soient $n\geq 1$ et $\pi \in \Pi_{\rm cusp}(\PGL_n)$
alg\'ebrique de poids motivique $\leq 22$. Alors $\pi$ appartient \`a la liste des $11$ repr\'esentations suivantes :
$$1, \hspace{4pt}\Delta_{11}, \hspace{4pt}\Delta_{15},\hspace{4pt} \Delta_{17},\hspace{4pt}\Delta_{19},\hspace{4pt}\Delta_{19,7},\hspace{4pt}\Delta_{21},\hspace{4pt}\Delta_{21,5},\hspace{4pt}\Delta_{21,9},\hspace{7pt}\Delta_{21,13},\hspace{7pt}{\rm Sym}^2 \Delta_{11}.$$
\end{thmintro}
 \ps\ps
En poids motivique $<11$, ce th\'eor\`eme affirme que l'on a $n=1$ et que $\pi$ est la repr\'esentation triviale, un r\'esultat \'etait d\'ej\`a connu de Mestre et Serre (au langage pr\`es, voir \cite[\S III, Remarque 1]{mestre}). Dans ce cas tr\`es particulier, il fournit entre autres une explication ``automorphe'' du th\'eor\`eme classique de Minkowski affirmant que tout corps de nombres distinct de $\Q$ est ramifi\'e en au moins un nombre premier (cas ${\rm w}(\pi)=0$), ou encore de la conjecture de Shaffarevic, d\'emontr\'ee ind\'ependamment par Abrashkin et Fontaine,  selon laquelle il n'y a pas de vari\'et\'e ab\'elienne sur $\Z$ (cas ${\rm w}(\pi)=1$). \`A notre connaissance, l'\'enonc\'e du th\'eor\`eme \ref{introw22} est d\'ej\`a nouveau dans le cas particulier ${\rm w}(\pi)=11$. Soulignons que l'on n'y fait pas d'hypoth\`ese sur l'entier $n$, et qu'il entra\^ine $n \leq 4$. \ps \ps

Notre d\'emonstration de ce th\'eor\`eme, dans le prolongement des travaux de Stark, Odlyzko et Serre sur les
minorations des discriminants de corps de nombres, repose sur un analogue dans le cadre des fonctions ${\rm L}$ automorphes des {\it formules explicites} de Riemann et Weil en th\'eorie des nombres premiers. Cet analogue a \'et\'e d\'evelopp\'e par Mestre \cite{mestre}, puis appliqu\'e par Fermigier aux fonctions
${\rm L}(s,\pi)$ pour $\pi \in \Pi_{\rm cusp}(\PGL_n)$ \cite{fermigier}. Nous l'appliquons plus g\'en\'eralement \`a la fonction ${\rm L}$ ``de Rankin-Selberg'' d'une paire arbitraire $\{\pi, \pi'\}$ de repr\'esentations automorphes cuspidales alg\'ebriques de $\PGL_{n}$ et $\PGL_{n'}$ (Jacquet, Piatetski-Shapiro, Shalika). \ps\ps

Dans le cas particulier o\`u $\pi'$ est la duale de $\pi$, cette m\'ethode avait d\'ej\`a port\'e ses fruits dans un travail de Miller \cite{miller} ; cependant notre \'etude comprend quelques nouveaut\'es qui m\'eritent d'\^etre signal\'ees. Tout d'abord, nous avons d\'ecouvert que certaines formes bilin\'eaires sym\'etriques, qui sont d\'efinies sur l'anneau de Grothendieck ${\rm K}_\infty$ du groupe de Weil de $\R$ et \`a valeurs r\'eelles, et qui interviennent de mani\`ere naturelle dans l'\'ecriture des formules explicites, sont d\'efinies positives sur des sous-groupes assez gros de ${\rm K}_\infty$. C'est cet \'enonc\'e qui est responsable de la finitude de la liste apparaissant dans le th\'eor\`eme \ref{introw22}.  De plus, nous \'etablissons des crit\`eres simples, portant par exemple uniquement sur $\pi_\infty$ et $\pi'_\infty$, interdisant l'existence simultan\'ee de $\pi$ et de $\pi'$. Nous renvoyons au \S \ref{mrw} pour des \'enonc\'es pr\'ecis. \ps \ps

\section{D\'emonstration des th\'eor\`emes \ref{pbeichler24} et \ref{thmintro24}}\label{introdem}

	Donnons les grandes lignes de la d\'emonstration du th\'eor\`eme \ref{thmintro24} (\S IX.\ref{demothmE}). Soit $\pi \in \Pi_{\rm disc}({\rm O}_{24})$ telle que $\pi_\infty$ est la repr\'esentation triviale. Les r\'esultats d'Erokhin et Borcherds-Freitag-Weissauer rappel\'es au \S \ref{introthetasiegel} montrent que $\pi$ admet un ``$\vartheta$-correspondant'' dans $\Pi_{\rm cusp}({\rm Sp}_{2g})$, engendr\'e par une forme de Siegel de poids $12$ et de genre $g  \leq 11$ (\S \ref{retoureichler}), \`a l'exception pr\`es d'un $\pi$ d\'ej\`a d\'etermin\'e par Ikeda \cite{ikeda1} et v\'erifiant $\psi(\pi,{\rm St})=\Delta_{11}[12]$. Le th\'eor\`eme d'Arthur appliqu\'e \`a ce $\vartheta$-correspondant, et le point de vue de Rallis sur les relations d'Eichler, entra\^inent que le couple $(\pi,{\rm St})$ v\'erifie la conjecture d'Arthur-Langlands. Une \'etude combinatoire simple reposant uniquement sur le th\'eor\`eme \ref{introw22} montre qu'il y a au plus $24$ possibilit\'es pour $\psi(\pi,{\rm St})$ (celles de la table \ref{table24}). Mais il y a \'egalement au moins $24$ possibilit\'es pour $\psi(\pi,{\rm St})$, car d'apr\`es Nebe et Venkov l'op\'erateur ${\rm T}_2$ admet des valeurs propres distinctes sur $\C[{\rm X}_{24}]$, ce qui termine la d\'emonstration. \ps\ps

	Cette m\'ethode permet d'\'etudier plus g\'en\'eralement les repr\'esentations dans $\Pi_{\rm cusp}({\rm Sp}_{2g})$ engendr\'ees par une forme modulaire de Siegel de poids $k\leq 12$ pour le groupe ${\rm Sp}_{2g}(\Z)$ (\S \ref{fsiegelpoids12}). Le th\'eor\`eme \ref{pbeichler24} est le fruit de cette \'etude dans le cas particulier $k=12$. On trouve en tout $23$ formes modulaires de Siegel pour ${\rm Sp}_{2g}(\Z)$
qui sont propres pour ${\rm H}({\rm Sp}_{2g})$, de poids $12$, et de genre $g \leq
23$, dont nous d\'ecrivons au passage les param\`etres standards  (Table \ref{tablek=12}). Dans le cas des formes de poids $k \leq 11$, nous d\'emontrons le th\'eor\`eme suivant, g\'en\'eralisant des r\'esultats de \cite{dukei1} et \cite{py2} (Th\'eor\`eme IX.\ref{thmsiegelsmallweight}). 
	
\begin{thmintro}\label{introk11} Soient $g\geq 1$ et $k \in \Z$. \ps \ps
\begin{itemize} \item[(i)] Si $k \leq 10$ alors ${\rm S}_k({\rm Sp}_{2g}(\Z))=0$, \`a moins que $(k,g)$ ne soit parmi $$(8,4), \hspace{.6 cm}(10,2), \hspace{.6 cm}(10,4), \hspace{.6 cm}(10,6), \hspace{.6 cm}(10,8),$$ auquel cas $\dim {\rm S}_k({\rm Sp}_{2g}(\Z)) =1$. Les param\`etres standards des $5$ \'el\'ements de $\Pi_{\rm disc}({\rm Sp}_{2g})$ engendr\'es par ces espaces sont respectivement 
{\small $$\Delta_{11}[4]\oplus[1],\hspace{.3 cm}\Delta_{17}[2]\oplus [1],\hspace{.3 cm}\Delta_{15}[4]\oplus [1],\hspace{.3 cm}\Delta_{17}[2] \oplus \Delta_{11}[4] \oplus [1]\hspace{.3 cm} {\rm et}\hspace{.3 cm} \Delta_{11}[8] \oplus [1].$$ }
\item[(ii)] Si $k=11$ alors ${\rm S}_k({\rm Sp}_{2g}(\Z))=0$, sauf peut-\^etre si $g=6$. \ps\ps
\end{itemize}
\end{thmintro}
\ps\ps
	Indiquons quelques difficult\'es dans les d\'emonstrations des th\'eor\`emes  \ref{pbeichler24} et \ref{introk11} qui n'apparaissent pas dans celle du th\'eor\`eme \ref{thmintro24}. Soit $\pi$ un \'el\'ement de $\Pi_{\rm cusp}({\rm Sp}_{2g})$ engendr\'e par une forme modulaire de Siegel de poids $k \leq 12$ et de genre $g < 12+k$. Le th\'eor\`eme \ref{introw22} permet de montrer qu'il n'y a qu'une liste finie explicite de possibilit\'es pour $\psi(\pi,{\rm St})$. En revanche, contrairement \`a la situation du th\'eor\`eme \ref{thmintro24}, certains \'el\'ements de cette liste ne devraient pas intervenir, comme le montre une inspection de la formule de multiplicit\'e. Nous contournons l'utilisation de cette formule par un recours \`a des r\'esultats de B\"ocherer \cite{bochererbowdoin} et d'Ikeda \cite{ikeda1,ikeda3}, ainsi qu'\`a un ensemble de constructions de s\'eries th\^eta. L'annulation de ${\rm S}_{11}({\rm Sp}_{12}(\Z))$ est attendue, mais nous ne pouvons la d\'emontrer inconditionellement.  Les cas $g \geq k$ sont plus d\'elicats (nous ne savons m\^eme pas comment expliciter la formule de multiplicit\'e d'Arthur dans ce cas). Nous les excluons de mani\`ere {\it ad hoc} en utilisant les travaux de S. Mizumoto \cite{mizumoto} sur les p\^oles de la fonction ${\rm L}(s,\pi,{\rm St})$  (\S \ref{complboc}). \ps\ps

\section{Quelques applications}\label{introappl}

	D'apr\`es le th\'eor\`eme \ref{thmintro24}, la question initiale de
d\'eterminer les ${\rm N}_p(L,M)$ pour $L,M$ dans ${\rm X}_{24}$ et $p$ premier
devient \'equivalente \`a celle de d\'eterminer les valeurs propres des op\'erateurs de
Hecke agissant sur les $4$ formes modulaires de Siegel \`a valeurs vectorielles pour ${\rm Sp}_4(\Z)$ mentionn\'ees au \S \ref{fautintro}.  Nous exposons au \S
\ref{calcfinvpgenre2} une m\'ethode que nous avons d\'ecouverte pour y
parvenir, utilisant l'analyse des $p$-voisins du r\'eseau de Leech
effectu\'ee au \S III.4.  \ps\ps

Soit $(j,k)$ l'un des $4$ couples consid\'er\'es au \S \ref{fautintro}, c'est-\`a-dire $(6,8)$, $(4,10)$, $(8,8)$ ou $(12,6)$. Notons $(w,v)$ le couple $(2j+k-3,j+1)$ correspondant. Si $q$ est un entier de la forme $p^k$ avec $p$ premier et $k$ entier $\geq 1$, nous posons $$\tau_{j,k}(q)\, \, =\, \, q^{\frac{w}{2}} \, {\rm trace} \,\,{\rm c}_p(\Delta_{w,v})^k,$$ c'est un nombre complexe qui est en fait dans $\Z$. \ps \ps

\begin{thmintro}\label{introtjk} Soit $(j,k)$ l'un des $4$ couples $(6,8), (4,10), (8,8)$ ou $(12,6)$. Les entiers $\tau_{j,k}(p)$ avec $p$ premier $\leq 113$, et $\tau_{j,k}(p^2)$ avec $p$ premier $\leq 29$, sont respectivement donn\'es par les tables \ref{taujk} et \ref{taujkp2}.
\end{thmintro} 
\ps\ps

Ces r\'esultats confirment et \'etendent des calculs
ant\'erieurs de Faber et van der Geer \cite{fabervandergeerIetII}  \cite[\S
25]{vandergeer}, par des
m\'ethodes tr\`es diff\'erentes. Au final, nous disposons donc de la valeur exacte de
${\rm N}_p(L,M)$ pour tout $L,M$ dans ${\rm X}_{24}$ et tout nombre premier $p
\leq 113$.   \ps \ps
	
	Le th\'eor\`eme~\ref{introw22} montre que la question du calcul des $\tau_{j,k}(q)$ est peut-\^etre moins futile qu'elle ne le para\^it. En effet, compte tenu des conjectures de Langlands, ce th\'eor\`eme sugg\`ere une classification parall\`ele, encore \`a d\'emontrer du c\^ot\'e $\ell$-adique, des  motifs purs et effectifs sur $\Q$, ayant bonne r\'eduction partout, et de poids motivique $\leq 22$. Par exemple, il impose une contrainte conjecturale remarquable sur la fonction z\^eta de Hasse-Weil du champs de Deligne-Mumford $\overline{\mathcal{M}_{{g,n}}}$ classifiant les courbes stables de genre $g$ munies de $n$ points marqu\'es, avec $g\geq 2$, $n\geq 0$ et $3g-3+n \leq 22$ : elle devrait s'exprimer uniquement en termes des coefficients des formes modulaires paraboliques normalis\'ees de poids $\leq 22$ pour ${\rm SL}_{2}(\Z)$ et des $\tau_{j,k}(q)$. Cela confirme certains r\'esultats (resp. conjectures) de Bergstr\"om, Faber et van der Geer \cite{fabervandergeerIetII,fabermgn,BFVdG} lorsque $g=2$ (resp. $g=3$).  \ps\ps

	Nous donnons au \S \ref{parfinharder} des applications du
th\'eor\`eme \ref{thmintro24} aux con\-gruences satisfaites par les entiers $\tau_{j,k}(p)$ avec $p$ premier.  Ces congruences sont obtenues par une \'etude des vecteurs propres de ${\rm T}_2$ dans la base naturelle de $\Z[{\rm X}_{24}]$ et par des arguments de repr\'esentations galoisiennes. Entre autres, nous d\'emontrons la congruence conjectur\'ee par Harder dans \cite{harder}.\ps \ps

\begin{thmintro} \label{introconjhar}{\rm (Conjecture de Harder)} Pour tout premier $p$ on a la congruence 
$$\tau_{4,10}(p)  \equiv \tau_{22}(p) + p^{13} + p^{8} \bmod 41,$$ o\`u $\tau_{22}(p)$ d\'esigne le $p$-i\`eme coefficient de la forme modulaire parabolique normalis\'ee de poids $22$ pour le groupe ${\rm SL}_2(\Z)$.
\end{thmintro}\ps\ps

	Revenons enfin sur la d\'emonstration du th\'eor\`eme \ref{thmintro24} esquiss\'ee au \S \ref{introdem}. Elle utilise l'\'egalit\'e $|{\rm X}_{24}|=24$, cons\'equence de la classification de Niemeier. Cependant, nous expliquons au \S \ref{nou2425} comment la combinaison des id\'ees ci-dessus et des formules de multiplicit\'e d'Arthur (incluant les conjectures \ref{conjaj} et \ref{conjaj2} formul\'ees au \S VIII) permet en retour de se passer de cette \'egalit\'e, et m\^eme de la red\'emontrer ``sans aucun calcul de r\'eseaux euclidiens''. Mieux, nous retrouvons non seulement qu'il existe exactement $24$ r\'eseaux de Niemeier \`a isom\'etries pr\`es, mais \'egalement qu'un seul d'entre eux ne poss\`ede pas d'isom\'etrie de d\'eterminant $-1$. \ps\ps

Est-il raisonnable d'esp\'erer pouvoir estimer finement le cardinal de ${\rm X}_{32}$ par une telle m\'ethode  ? La question reste enti\`ere, mais l'exemple de la dimension $24$ montre que cette approche, ch\`ere au premier auteur, est moins farfelue qu'elle n'en a l'air. Un ingr\'edient n\'ecessaire \`a ce projet est la connaissance des repr\'esentations alg\'ebriques (disons ``autoduales, r\'eguli\`eres'') dans $\Pi_{\rm cusp}(\PGL_n)$ dont le poids motivique est $\leq 30$ : des progr\`es dans cette direction ont \'et\'e r\'ecemment obtenus dans \cite{chrenard2} et \cite{taibisiegel}.\par

\bigskip{
\begin{center}
*\par
*\hspace*{3ex}*
\end{center}
}

Pour conclure cette introduction, disons un mot sur l'utilisation
dans ce m\'emoire des r\'esultats r\'ecents d'Arthur.  Ces derniers
reposent sur un ensemble impressionnant de travaux difficiles,
certaines d\'emonstrations \'etant encore toutes fra\^iches
(voir \cite{arthur}, \cite{moewal2, stablewal} et la discussion du \S \ref{enonceparamst}).  C'est pourquoi il
nous a sembl\'e utile d'indiquer, par une \'etoile$^{\color{green}
\ast}$, dans le corps du m\'emoire, les \'enonc\'es d\'ependant des
r\'esultats du livre \cite{arthur}.  Dans cette introduction, cela
concerne les d\'emonstrations des th\'eor\`emes \ref{thm24pop},
\ref{pbeichler24}, \ref{thmintro24}, \ref{introw22}\footnote{En ce qui concerne le th\'eor\`eme \ref{introw22}, nous en d\'emontrons en fait une variante quasiment aussi forte sans aucun recours \`a la th\'eorie d'Arthur : voir le th\'eor\`eme IX.\ref{classpoids22}.}, \ref{introk11}, \ref{introtjk}  et \ref{introconjhar}.  En revanche,
pr\'ecisons que contrairement aux annonces que nous avions faites
initialement dans \cite{chlannes} et \cite{chhdr}, nos d\'emonstrations ne reposent
plus sur les r\'esultats annonc\'es par Arthur {\it loc.  cit.} dans
son dernier chapitre sur les formes int\'erieures, ni sur les
propri\'et\'es conjecturales des paquets d'Arthur du type de ceux
\'etudi\'es par Adams et Johnson.\ps \ps

\newpage

\begin{center} {\it Remerciements} \end{center}\ps\ps\ps

\noindent Les auteurs tiennent \`a remercier le Centre de Math\'ematiques de l'\'Ecole
Polytechnique, dont les structures ont permis la collaboration improbable
entre un arithm\'eticien et un topologue amateurs de formes quadratiques.
Ils remercient \'egalement chaleureusement Colette Moeglin, David Renard et
Olivier Ta\"ibi pour des discussions utiles \`a ce travail, ainsi que Richard Borcherds et Thomas M\'egarban\'e pour leurs remarques. Enfin, les
auteurs remercient respectivement Valeria et Martine pour leur soutien tout au long de ce
projet qui a parfois paru sans fin. \ps \ps

\parindent=0cm


\chapter{Algèbre bilinéaire et quadratique}

\section{Généralités sur les formes bilinéaires et quadratiques}

\vspace{0,2cm}
Soit $A$ un anneau commutatif unitaire.

\bigskip
Nous appellerons {\em $\mathrm{b}$-module} sur
$A$, un $A$-module projectif de type fini $L$ muni
d'une forme bilinéaire symétrique non dégénérée, c'est-à-dire telle que\linebreak l'homomorphisme induit $L\to\mathrm{Hom}_{A}(L,A)$ est un isomorphisme~; quand $A$ sera un corps nous remplacerons évidemment
``$\mathrm{b}$-module'' par ``{\em $\mathrm{b}$-espace vectoriel}\hspace{2pt}''. Le plus souvent la forme bilinéaire
symétrique $L\times L\to A$ sera notée
$(x,y)\mapsto x.y$\hspace{3pt}. Dans ce mémoire les anneaux $A$ auxquels nous aurons affaire seront des anneaux principaux ou des corps si bien que nous aurions pu remplacer dans la définition ci-dessus ``$A$-module projectif de type fini'' par ``$A$-module libre de dimension finie''.

\medskip
Soit $S$ une matrice symétrique de taille $(n,n)$ à coefficients dans $A$, la notation $\langle S\rangle$ désignera le $A$-module $A^{n}$ muni de la forme bilinéaire dont la matrice dans la base canonique est $S$~; il est clair que $\langle S\rangle$ est un $\mathrm{b}$-module si et seulement si l'on a $\mathop{\mathrm{d\acute{e}t}}S\in A^{\times}$. Si $S$ est une matrice diagonale de coefficients diagonaux $a_{1},a_{2},\ldots,a_{n}$ alors $\langle S\rangle$ sera plutôt notée $\langle a_{1},a_{2},\ldots,a_{n}\rangle$.

\bigskip
Nous supposons maintenant que $A$ est un anneau de Dedekind~; on note $K$ son corps des fractions. (En fait, étant donné les anneaux que nous avons en tête, à savoir $\mathbb{Z}$ et $\mathbb{Z}_{p}$, nous pourrions remplacer ``anneau de Dedekind'' par ``anneau principal''.)

\bigskip
Soit $V$ un $K$-espace vectoriel de dimension finie. On appelle {\em réseau} de $V$ (par rapport à $A$) un sous-$A$-module $L$ de $V$, qui engendre $V$ comme $K$-espace vectoriel et qui est de type fini sur $A$~; un tel $L$ est un $A$-module projectif dont le rang est égal à $\dim_{K}V$.

\bigskip
Soit $V$ un $\mathrm{b}$-espace vectoriel sur
$K$. Soit $L$ un réseau de $V$. Le sous-$A$-module
de $V$ formé des éléments $y$ tels
que $x.y$ appartient à $A$ pour tout $x$ dans
$L$ est un réseau de $V$ que l'on appelle le
{\em dual} de $L$ et que l'on note
$L^{\sharp}$.

\bigskip
Nous dirons que le réseau $L$ est
{\em entier} si $x.y$ appartient à
$A$ pour tous $x$ et $y$ dans $L$, en d'autres
termes si l'on a $L\subset L^{\sharp}$.

\medskip
Soit $L$ un réseau entier de $V$.
Considérons le quotient $L^{\sharp}/L$~:
\begin{itemize}
\item [--] Le $A$-module $L^{\sharp}/L$ est de
torsion et de type fini.
\item [--] La forme bilinéaire symétrique
dont $V$ est muni induit sur $L^{\sharp}/L$ une
``forme'' bilinéaire symétrique
non dégénérée à valeurs dans
$K/A$. Ici ``non dégénérée''
signifie que l'homomorphisme, induit par la forme
en question, $L^{\sharp}/L\to
\mathrm{Hom}_{A}(L^{\sharp}/L,K/A)$ est un
isomorphisme.
\end{itemize}

\medskip
Nous appelons ce type d'objet un {\em $\mathrm{e}$-module} sur $A$. Un $\mathrm{e}$-module sur $A$ est donc un $A$-module de torsion et de type fini muni d'une forme bilinéaire symétrique non dégénérée à valeurs dans $K/A$ ({\em forme
d'enlacement}). Nous appelons {\em résidu} de $L$ le $\mathrm{e}$-module $L^{\sharp}/L$. La considération du quotient $L^{\sharp}/L$ est bien sûr fort ancienne, mais notre terminologie n'est pas classique, par exemple $L^{\sharp}/L$ est appelé ``dual quotient group'' ou ``glue group'' dans \cite{conwaysloane} et ``conoyau de $L$'' dans \cite{BLLV}~; nous le noterons souvent
$\mathrm{r\acute{e}s}\hspace{1pt}L$. 

\medskip
Dans la mesure du possible nous étendrons aux
$\mathrm{e}$-modules les notations et la
terminologie utilisées pour les modules munis
d'une forme bilinéaire symétrique à valeurs dans $A$. Voici quelques exemples. Soit $C$ un $\mathrm{e}$-module :

\smallskip
-- la forme bilinéaire symétrique $C\times C\to K/A$ sera le plus souvent notée $(x,y)\mapsto x.y$~;

\smallskip
-- nous dirons qu'un sous-module $I$ de $C$ est  {\em isotrope} si l'on a $x.y=0$ pour tous éléments $x$ et $y$ de $I$, autrement dit si l'on a $I\subset I^{\perp}$, $I^{\perp}$ désignant l'orthogonal de $I$~; 

\smallskip
-- nous dirons qu'un sous-module $I$ de $C$ est un {\em lagrangien} si l'on a $I=I^{\perp}$.

\bigskip
Nous appellerons {\em $\widetilde{\mathrm{b}}$-module} sur $A$ (qui est toujours supposé de Dedekind) un $A$-module projectif de type fini $L$ muni d'une forme bilinéaire symétrique telle que l'homomorphisme induit $L\to\mathrm{Hom}_{A}(L,A)$ est injectif (ou, ce qui revient au même, telle que la forme bilinéaire induite sur $K\otimes_{A}L$ est non dégénérée). Un réseau entier d'un $\mathrm{b}$-espace vectoriel sur $K$ est le type même d'un tel
objet. Réciproquement tout $\widetilde{\mathrm{b}}$-module $L$ sur $A$ peut être vu comme un réseau
entier du $\mathrm{b}$-espace vectoriel
$K\otimes_{A}L$. Un $\widetilde{\mathrm{b}}$-module
$L$ possède donc un résidu $\mathrm{r\acute{e}s}\hspace{1pt}L$ qui est un $\mathrm{e}$-module~; le $A$-module sous-jacent à $\mathrm{r\acute{e}s}\hspace{1pt}L$ s'identifie au conoyau de l'injection
$L\hookrightarrow\mathrm{Hom}_{A}(L,A)$.

\bigskip
En remplaçant, dans les définitions précédentes, formes bilinéaires symé\-triques par formes bilinéaire alternées, on obtient, \textit{mutatis mutandis}, les définitions respectives de {\em $\mathrm{a}$-module}, {\em $\mathrm{ae}$-module}, {\em $\widetilde{\mathrm{a}}$-module} et {\em résidu} d'un $\widetilde{\mathrm{a}}$-module. Par exemple, un $\mathrm{ae}$-module sur $A$ est un $A$-module de torsion et de type fini $C$ muni d'une forme bilinéaire alternée $C\times C\to K/A$ telle que l'homomorphisme induit $C\to\mathrm{Hom}_{A}(C,K/A)$ est un isomorphisme. On observera qu'un $\mathrm{a}$-module sur $A$ est toujours de rang pair et qu'il en est de même pour un $\widetilde{\mathrm{a}}$-module si $A$ est de Dedekind.

\bigskip
Rappelons qu'une application $f:M\to N$ entre deux $A$-modules est dite {\em quadratique} si elle vérifie les deux propriétés suivantes~:

\smallskip
-- on a $f(ax)=a^{2}f(x)$ pour tout $a$ dans $A$ et tout $x$ dans $M$~;

\smallskip
-- l'application $M\times M\to N,(x,y)\mapsto f(x+y)-f(x)-f(y)$ est bilinéaire.

\smallskip
Dans le cas $N=A$, on dit que $f$ est une forme quadratique sur $M$.

\bigskip
On obtient, \textit{mutatis mutandis}, les définitions respectives de {\em $\mathrm{q}$-module},\linebreak {\em $\mathrm{qe}$-module} et {\em $\widetilde{\mathrm{q}}$-module}, en remplaçant dans les définitions précédentes formes bilinéaires symétriques par formes quadratiques. Par exemple un $\mathrm{q}$-module sur $A$ est un $A$-module projectif de type fini $L$ muni d'une forme quadratique $\mathrm{q}:L\to A$ telle que la forme bilinéaire symétrique (dite {\em associée})
$$
L\times L\to A
\hspace{24pt},\hspace{24pt}
(x,y)
\hspace{4pt}\mapsto\hspace{4pt}
\mathrm{q}(x+y)-\mathrm{q}(x)-\mathrm{q}(y)
$$
est non dégénérée. Un $\mathrm{qe}$-module sur $A$ est un $A$-module de torsion et de type fini muni d'une
forme quadratique non dégénérée à valeurs dans $K/A$ ({\em forme quadratique d'enlacement})~; le {\em résidu} d'un $\widetilde{\mathrm{q}}$-module est maintenant un\linebreak
$\mathrm{qe}$-module. Un sous-module $I$ d'un $\mathrm{qe}$-module est {\em isotrope} si l'on a $\mathrm{q}(I)=0$ (condition qui implique $I\subset I^{\perp}$)~; c'est un {\em lagrangien} si l'on a $\mathrm{q}(I)=0$ en plus de la condition $I=I^{\perp}$. Un réseau $L$ d'un $\mathrm{q}$-espace vectoriel est {\em entier} si $\mathrm{q}(x)$ appartient à $A$ pour tout $x$ dans $L$ (condition qui implique que $x.y$ appartient à $A$ pour tous $x$ et $y$ dans $L$).

\bigskip
Si $2$ n'est pas diviseur de zéro alors un $\mathrm{q}$-module sur $A$ n'est rien d'autre qu'un $\mathrm{b}$-module {\em pair} c'est-à-dire un $\mathrm{b}$-module $L$ tel que $x.x$ est divisible par $2$ pour tout $x$ dans $L$~: dans ce cas la forme quadratique est déterminée par l'égalité $x.x=2\mathrm{q}(x)$. Si $2$ est inversible dans $A$ alors les notions de $\mathrm{q}$-module et $\mathrm{b}$-module coïncident~; on se souviendra cependant de ce que formes quadratiques et formes bilinéaires symétriques associées sont reliées par l'égalité précédente.

\bigskip
La proposition suivante est évidente\ldots ce qui ne l'empêche pas d'être très utile.

\bigskip
\textbf{Proposition 1.1.} {\em Soient $V$ un $\mathrm{b}$-espace vectoriel (resp. $\mathrm{a}$-espace vectoriel, resp. $\mathrm{q}$-espace vectoriel) sur $K$ et $L$ un réseau entier de $V$. Soit
$\gamma:L^{\sharp}\to\mathop{\mathrm{r\acute{e}s}}L$
l'homomorphisme de passage au quotient.}

\medskip
(a) {\em Soit $I$ un sous-module de
$\mathrm{r\acute{e}s}\hspace{1pt}L$ alors les
deux conditions suivantes sont équivalentes~:
\begin{itemize}
\item [(i)] le sous-module $I$ est isotrope~; 
\item [(ii)] le réseau $\gamma^{-1}(I)$ est
entier.
\end{itemize}}

\medskip
(b) {\em L'application $I\mapsto\gamma^{-1}(I)$
est une bijection croissante de l'ensemble des
sous-modules isotropes de 
$\mathrm{r\acute{e}s}\hspace{1pt}L$ sur
l'ensemble des réseaux entiers de $V$
contenant $L$ (et donc contenus dans $L^{\sharp}$).}

\medskip
(c) {\em Soit $I$ un sous-module isotrope de $\mathrm{r\acute{e}s}\hspace{1pt}L$, alors la forme bilinéaire sy\-mé\-trique
$I^{\perp}/I\times I^{\perp}/I\to K/A$ (resp. la forme bilinéaire alternée
$I^{\perp}/I\times I^{\perp}/I\to K/A$ , resp. la forme quadratique $I^{\perp}/I\to K/A$), induite par celle dont $\mathrm{r\acute{e}s}\hspace{1pt}L$ est muni, fait de $I^{\perp}/I$ un $\mathrm{e}$-module (resp. $\mathrm{ae}$-module, resp. $\mathrm{qe}$-module) qui s'identifie au résidu du réseau entier $\gamma^{-1}(I)$.}

\vspace{0,75cm}
\textsc{Foncteurs hyperboliques}

\bigskip
Soit $A$ un anneau commutatif unitaire. Soit $L$ un $A$-module projectif de type fini~; $\mathrm{Hom}_{A}(L,A)$ est encore un $A$-module projectif de type fini que l'on note~$L^{*}$ (et que l'on appelle le {\em dual} de $L$). L'application
$$
L\oplus L^{*}\to A
\hspace{24pt},\hspace{24pt}
(x,\xi)\mapsto\langle x,\xi\rangle
$$
est une forme quadratique non dégénérée qui fait du $A$-module projectif de type fini $L\oplus L^{*}$ un $\mathrm{q}$-module. Ce $\mathrm{q}$-module est noté $\mathrm{H}(L)$ et appelé le\linebreak $\mathrm{q}$-module {\em hyperbolique} sur $L$.

\medskip
La forme bilinéaire symétrique associée est l'application $$
\hspace{24pt}
((x,\xi),(y,\eta))\mapsto \langle x,\eta\rangle+\langle y,\xi\rangle
\hspace{24pt};
$$
$L\oplus L^{*}$ muni de cette forme est appelé le $\mathrm{b}$-module {\em hyperbolique} sur $L$ et est encore noté $\mathrm{H}(L)$.

\medskip
Pareillement, $L\oplus L^{*}$ muni de la forme bilinéaire alternée
$$
((x,\xi),(y,\eta))\mapsto \langle x,\eta\rangle-\langle y,\xi\rangle
$$
est un $\mathrm{a}$-module sur $A$, appelé $\mathrm{a}$-module {\em hyperbolique} sur $L$ et toujours noté~$\mathrm{H}(L)$.

\bigskip
Soit $H$ un $\mathrm{q}$-module (resp. $\mathrm{a}$-module) sur $A$~; rappelons que dans ce contexte un lagrangien de $H$ est un facteur direct $L$ avec $L=L^{\perp}$ et $\mathrm{q}(L)=0$ (resp. $L=L^{\perp}$).

\bigskip
\textbf{Proposition 1.2.} {\em Soient $H$ un $\mathrm{q}$-module (resp. $\mathrm{a}$-module) et $L$ un la\-grangien~; alors l'inclusion de $L$ dans $H$ se prolonge en un isomorphisme de $\mathrm{q}$-modules (resp. $\mathrm{a}$-modules) $\mathrm{H}(L)\simeq H$.}

\bigskip
\textit{Démonstration.} On démontre la ``version quadratique'' de l'énoncé~; nous suivons sans vergogne la démonstration de \cite[Proposition 2.1.5]{BL} qui traite implicitement de la ``version alternée''. Soient $i:L\to H$ l'inclusion de $L$ dans $H$ et $p:H\to L^{*}$ l'homomorphisme (de $A$-modules) composé de l'isomorphisme $H\to H^{*}$ induit par la forme bilinéaire et de l'homomorphisme~$i^{*}$. Puisque $L$ est un lagrangien, la suite de $A$-modules
$$
\begin{CD}
0 @>>> L @>i>> H @>p>> L^{*} @>>> 0
\end{CD}
$$
est exacte. Il s'agit de montrer qu'il existe une section $A$-linéaire $s:L^{*}\to H$ de $p$ vérifiant $\mathrm{q}(s(\xi))=0$ pour tout $\xi$ dans $L^{*}$. Soit $\Sigma$ l'ensemble des sections $A$-linéaires de $p$. L'ensemble $\Sigma$ n'est pas vide puisque $L^{*}$ est projectif~; de plus $\Sigma$ possède une structure
canonique d'espace affine sous $\mathrm{Hom}_{A}(L^{*},L)$ que l'on identifie avec le $A$-module des formes bilinéaires sur $L^{*}$, disons $\mathcal{B}_{L^{*}}$. Soient $\mathcal{Q}_{L^{*}}$ le $A$-module des formes quadratiques sur $L^{*}$ et $\gamma:\Sigma\to\mathcal{Q}_{L^{*}}$ l'application qui associe à une section $s$ la forme quadratique $\xi\mapsto\mathrm{q}(s(\xi))$~; soit $u$ un élément de $\mathrm{Hom}_{A}(L^{*},L)$, on constate que l'on a~:
$$
\gamma\hspace{1pt}(s+u)=
\hspace{1pt}\gamma(s)+\widetilde{\gamma}\hspace{1pt}(u)
\hspace{24pt},
$$
$\widetilde{\gamma}$ désignant l'application $\mathcal{B}_{L^{*}}\to\mathcal{Q}_{L^{*}}$ qui associe à une forme bilinéaire $u$ la forme quadratique $\xi\mapsto u(\xi,\xi)$. Le fait que $\gamma^{-1}(0)$ n'est pas vide résulte maintenant du fait que $\widetilde{\gamma}$ est surjective. Ceci est clair si $L^{*}$ est libre~; le cas général en découle en faisant intervenir un $A$-module $M$ tel que la somme directe $L^{*}\oplus M$ est libre de dimension finie.
\hfill$\square$

\medskip
\textit{Remarque.} La considération du $\mathrm{b}$-espace vectoriel $\langle\begin{bmatrix} 0 & 1 \\ 1 & 1 \end{bmatrix}\rangle$ sur $\mathbb{F}_{2}$ montre que l'énoncé analogue à 1.2 pour les $\mathrm{b}$-modules n'est pas vérifié en général.

\bigskip
On dit qu'un $\mathrm{q}$-module (resp. $\mathrm{a}$-module) est {\em hyperbolique} s'il est isomorphe à un $\mathrm{H}(L)$ pour un certain $A$-module projectif de type fini $L$~; la proposition ci-dessus dit en particulier qu'un $\mathrm{q}$-module est hyperbolique si et seulement s'il possède un lagrangien.

\bigskip
Soient maintenant $A$ un anneau de Dedekind (de corps des fractions $K$) et $I$ un $A$-module de torsion de type fini~; $\mathrm{Hom}_{A}(I,K/A)$ est encore un\linebreak $A$-module de torsion de type fini que l'on note $I^{\vee}$ (et que l'on appelle le {\em dual} de $I$). L'application
$$
I\oplus I^{\vee}\to K/A
\hspace{24pt},\hspace{24pt}
(x,\xi)\mapsto\langle x,\xi\rangle
$$
est une forme quadratique non dégénérée qui fait du $A$-module de torsion de type fini $I\oplus I^{\vee}$ un $\mathrm{qe}$-module. Ce $\mathrm{qe}$-module est noté $\mathrm{H}(I)$ et appelé le $\mathrm{qe}$-module {\em hyperbolique} sur $I$.

\bigskip
Soit $H$ un $\mathrm{qe}$-module sur $A$~; rappelons que dans ce contexte un lagrangien de $H$ est un sous-module $I$ avec $I=I^{\perp}$ et $\mathrm{q}(I)=0$. La démonstration de la proposition suivante est laissée au lecteur (on pourra comparer le point (b) avec le point (b) de \cite[Proposition 2.1.5]{BL}).

\bigskip
\textbf{Proposition-Définition 1.3.} {\em Soit $H$ un $\mathrm{qe}$-module. Soient $I$ et $J$ deux\linebreak lagrangiens de $H$~; on dit que $I$ et $J$ sont {\em transverses} (ou que $J$ est {\em transverse à} $I$) si l'on a $I\cap J=0$.

\medskip
{\em (a)} Soient $I$ et $J$ deux lagrangiens transverses de $H$. Alors la forme d'enlacement de $H$ induit un isomorphisme $J\cong I^{\vee}$ et l'homomorphisme composé
$$
\mathrm{H}(I)=I\oplus I^{\vee}\to I\oplus J\to H
$$
est un isomorphisme de $\mathrm{qe}$-modules (la première flèche est la somme directe de l'identité de $I$ et de l'inverse de l'isomorphisme $J\cong I^{\vee}$).

\medskip
{\em (b)} Soient $I$ un lagrangien de $H$ et $\mathcal{T}_{I}$ l'ensemble (éventuellement vide) des lagrangiens transverses à $I$. Alors $\mathcal{T}_{I}$ possède une structure canonique\linebreak d'espace affine sous le $A$-module $(\Lambda^{2}(I^{\vee}))^{\vee}$ (en clair le $A$-module constitué des applications bilinéaires alternées $I^{\vee}\times I^{\vee}\to K/A$).}

\bigskip
On dit qu'un $\mathrm{qe}$-module est {\em hyperbolique} s'il est isomorphe à un $\mathrm{H}(I)$ pour un certain $A$-module de type fini $I$~; le point (a) de la  proposition ci-dessus dit en particulier qu'un $\mathrm{qe}$-module est hyperbolique si et seulement s'il possède deux lagrangiens transverses.

\vspace{0,75cm}
\textsc{Produits tensoriels de formes}

\bigskip
Soit $A$ un anneau commutatif~; soient $L_{1}$ et $L_{2}$ deux $\mathrm{b}$-modules sur $A$.\linebreak L'homomorphisme de $A$-modules, $L_{1}\otimes_{A}L_{2}\to(L_{1}\otimes_{A}L_{2})^{*}\cong L_{1}^{*}\otimes_{A}L_{2}^{*}$, produit tensoriel des homomorphismes structurels $L_{1}\to L_{1}^{*}$ et $L_{2}\to L_{2}^{*}$ fait de $L_{1}\otimes_{A}L_{2}$ un $\mathrm{b}$-module que l'on appelle le {\em produit tensoriel} des $\mathrm{b}$-modules $L_{1}$ et $L_{2}$. La forme bilinéaire symétrique de $L_{1}\otimes_{A}L_{2}$ est caractérisée par le fait que l'on a $(x_{1}\otimes_{A}x_{2}).(y_{1}\otimes_{A}y_{2})=(x_{1}.y_{1})(x_{2}.y_{2})$ pour tous $x_{1},y_{1}$ dans $L_{1}$ et tous $x_{2},y_{2}$ dans $L_{2}$.
\vfill\eject

\medskip
\textit{Mutatis mutandis} on définit~:

\smallskip
-- le produit tensoriel d'un $\mathrm{b}$-module et d'un $\mathrm{q}$-module qui est un $\mathrm{q}$-module~;

\smallskip
-- le produit tensoriel d'un $\mathrm{b}$-module et d'un $\mathrm{e}$-module qui est un $\mathrm{e}$-module~;

\smallskip
-- le produit tensoriel d'un $\mathrm{b}$-module et d'un $\mathrm{qe}$-module qui est un $\mathrm{qe}$-module~;

\smallskip
-- le produit tensoriel d'un $\mathrm{q}$-module et d'un $\mathrm{e}$-module qui est un $\mathrm{qe}$-module.

\smallskip
Par exemple si $L_{1}$ est un $\mathrm{b}$-module et $L_{2}$ est un $\mathrm{q}$-module alors la forme quadratique de $L_{1}\otimes_{A}L_{2}$ est caractérisée par le fait que l'on a $\mathrm{q}(x_{1}\otimes_{A}x_{2})=(x_{1}.x_{1})\hspace{1pt}\mathrm{q}(x_{2})$ pour tout $x_{1}$ dans $L_{1}$ et tout $x_{2}$ dans $L_{2}$.

\vspace{0,75cm}
\textsc{Discriminant d'un $\mathrm{q}$-module de rang constant pair, déterminant de Dickson-Dieudonné}

\bigskip
Soient $A$ un anneau commutatif et $L$ un $\mathrm{q}$-module sur $A$ de rang constant pair, disons non nul. Soit $\Delta(L)$ le centre de la partie paire, notée $\mathrm{C}^{+}(L)$, de l'algèbre de Clifford de $L$ (voir par exemple \cite[Chap. III]{Ch}).

\medskip
-- La $A$-algèbre commutative $\Delta(L)$ est un ``revêtement double'' de $A$, en clair une $A$-algèbre étale et un $A$-module projectif de rang $2$ \cite[Exposé XII, Proposition 1.5]{sga7} qu'il faut voir comme le {\em discriminant} de $L$. Dans le cas  $L=\mathrm{H}(P)$, avec $P$ un $A$-module projectif de rang constant, ce revêtement est trivial (et même trivialisé) : $\Delta(L)=A\times A$.

\medskip
-- Un automorphisme $\alpha$ du $\mathrm{q}$-module $L$ induit un automorphisme $\Delta(\alpha)$ de la $A$-algèbre $\Delta(L)$. Si l'on identifie le groupe des automorphismes de la\linebreak $A$-algèbre $\Delta(L)$ avec $\mathbb{Z}/2(A)$, alors $\Delta(\alpha)$ s'identifie à un élément de $\mathbb{Z}/2(A)$ que l'on appelle le {\em déterminant de Dickson-Dieudonné} de $\alpha$~; nous le notons~$\mathop{\widetilde{\mathrm{d\acute{e}t}}}\alpha$. Rappelons la définition du foncteur en groupes $A\mapsto\mathbb{Z}/2(A)$ qui apparaît subrepticement ci-dessus~: $\mathbb{Z}/2(A)$ est l'ensemble des éléments $x$ de~$A$ vérifiant $x^{2}=x$ muni de la loi de groupe $(x,y)\mapsto x+y-2xy$.

\medskip
Ces deux points sont une ``globalisation'' sans surprises de résultats bien connus dans le cas où $A$ est un corps. Le cas subtil est celui où $A$ est un corps de caractéristique $2$, voir \cite{Di}\cite[\S 9, Exc. 9]{Bou1}.

\bigskip
Soit $L$ un $\mathrm{b}$-module sur $A$ de rang constant $n$. Le {\em déterminant de $L$} est le $\mathrm{b}$-module $\Lambda^{n}L$ ($\Lambda^{n}L$ est un $A$-module projectif de rang $1$, il est muni de la forme bilinéaire symétrique induite par celle de $L$)~; nous le notons $\mathop{\mathrm{d\acute{e}t}}L$. Dans le cas où $L$ est libre, la classe d'isomorphisme de $\mathop{\mathrm{d\acute{e}t}}L$ s'identifie à un élément  de $A^{\times}/A^{\times\hspace{1pt}2}$. Cet élément est la classe dans $A^{\times}/A^{\times\hspace{1pt}2}$ du déterminant de la matrice de Gram $[e_{i}.e_{j}]$, pour n'importe quelle base $(e_{1},e_{2},\ldots,e_{n})$ de~$L$~;  la classe en question est aussi souvent notée $\mathop{\mathrm{d\acute{e}t}}L$.

\medskip
Soit $L$ un $\mathrm{q}$-module sur $A$ de rang constant pair, disons $2n$ avec $n\geq 1$. La relation entre le discriminant de $L$ et le déterminant du $\mathrm{b}$-module sous-jacent est donnée par l'énoncé 1.4 ci-dessous dont nous devons une démonstration à Pierre Deligne. Pour formuler cet énoncé, il nous faut introduire quelques notations. On note $\mathrm{D}(L)$ est le conoyau de l'unité $\eta: A\to\Delta(L)$. On constate que le $A$-module $\mathrm{D}(L)$ est projectif de rang $1$ et qu'il est  muni d'une forme bilinéaire symétrique non dégé\-né\-rée canonique, disons~$\theta$, que l'on peut définir, par exemple, comme induite par la forme bilinéaire symétrique $\Delta(L)\times\Delta(L)\to A, (x,y)\mapsto\mathrm{tr}_{\Delta(L)/A}((x-\bar{x})y)$ ($\bar{x}$ désignant le ``conjugué'' de~$x$)~; en d'autres termes, $\mathrm{D}(L)$ possède une structure naturelle de $\mathrm{b}$-module sur $A$ de rang $1$. On note enfin $(-1)^{n}\mathop{\mathrm{d\acute{e}t}}L$ le $\mathrm{b}$-module sur $A$ de rang $1$ obtenu en multipliant par $(-1)^{n}$ la forme bilinéaire symétrique dont $\mathop{\mathrm{d\acute{e}t}}L$ est munie. Voici l'énoncé promis~:

\bigskip
\textbf{Proposition 1.4.} {\em Les deux $\mathrm{b}$-modules sur $A$ de rang $1$, $\mathrm{D}(L)$ et $(-1)^{n}\mathop{\mathrm{d\acute{e}t}}L$, sont naturellement isomorphes.}

\vspace{0,75cm}
\textsc{Groupes classiques}

\bigskip
L'objet essentiel de cette dernière rubrique du paragraphe et de fixer les notations et la terminologie relatives aux groupes orthogonaux et symplectiques (et leurs diverses déclinaisons) qui seront utilisées dans ce mémoire.

\bigskip
Soit $A$ un anneau commutatif unitaire.

\bigskip
-- Soit $L$ un $A$-module projectif de type fini. On note $\mathrm{GL}(L)$ le groupe des automorphismes de $L$. Le foncteur $R\mapsto\mathrm{GL}(R\otimes_{A}L)$, défini sur la catégorie des $A$-algèbres commutatives et à valeurs dans la catégorie des groupes, est un $A$-schéma en groupes que l'on note $\mathrm{GL}_{L}$~; on observera que si le rang de $L$ est $1$, alors $\mathrm{GL}(L)$ et $\mathrm{GL}_{L}$ s'identifient respectivement au groupe~$A^{\times}$ et au $A$-schéma en groupes $\mathbb{G}_{\mathrm{m}}$ (si $G$ est un $\mathbb{Z}$-schéma en groupes on note encore~$G$ le $A$-schéma en groupes obtenu par changement de base). Si $L$ est de rang constant, disons~$n$, on note $\mathrm{SL}_{L}$ le noyau de l'homomorphisme ``déterminant'' $\mathrm{d\acute{e}t}:\mathrm{GL}_{L}\to\mathrm{GL}_{\Lambda^{n}L}=\mathbb{G}_{\mathrm{m}}$. On note $\mathrm{PGL}_{L}$ le $A$-schéma en groupes défini comme le foncteur qui associe à une $A$-algèbre commutative $R$ le groupe des automorphismes de la $R$-algèbre $\mathrm{End}_{R}(R\otimes_{A}L)$~; $\mathrm{PGL}_{L}$ peut être aussi vu comme le $A$-schéma en groupes quotient de $\mathrm{GL}_{L}$ par $\mathbb{G}_{\mathrm{m}}$.  Evidemment, dans le cas $A=\mathbb{Z}$ et $L=A^{n}$, la notation  $\mathrm{GL}_{\mathbb{Z}^{n}}$, $\mathrm{SL}_{\mathbb{Z}^{n}}$, $\mathrm{PGL}_{\mathbb{Z}^{n}}$, sera remplacée par la notation $\mathrm{GL}_{n}$, $\mathrm{SL}_{n}$, $\mathrm{PGL}_{n}$.

\bigskip
-- Soit $L$ un $\mathrm{q}$-module (resp. $\mathrm{b}$-module) sur $A$, un endomorphisme $\alpha$ du\linebreak $A$-module sous-jacent est dit {\em orthogonal} si l'on a $\mathrm{q}(\alpha(x))=\mathrm{q}(x)$ pour tout $x$ dans $L$ (resp. $\alpha(x).\alpha(y)=x.y$ pour tous $x$ et $y$ dans $L$)~; les endomorphismes orthogonaux forment un groupe pour la composition (les endomorphismes orthogonaux sont donc en fait des automorphismes) que l'on appelle le {\em groupe orthogonal} de $L$ et que l'on note~$\mathrm{O}(L)$. Le foncteur $R\mapsto\mathrm{O}(R\otimes_{A}L)$, défini sur la catégorie des $A$-algèbres commutatives et à valeurs dans la catégorie des groupes, est un $A$-schéma en groupes que l'on note~$\mathrm{O}_{L}$.

\medskip
En fait les formes quadratiques jouent dans notre mémoire un rôle plus important que les formes bilinéaires symétriques. Une des raisons de cette prééminence est l'énoncé ci-dessous (qui n'est pas  vérifié, en toute généralité, pour les $\mathrm{b}$-modules).

\bigskip
\textbf{Proposition 1.5} {\em Pour tout $\mathrm{q}$-module $L$ sur un anneau commutatif $A$, le $A$-schéma en groupes $\mathrm{O}_{L}$ est lisse sur $A$.}

\medskip
\footnotesize
\textit{Démonstration.} Puisque la propriété que l'on veut vérifier est locale pour la topologie de Zariski, on peut supposer que $L$ est libre, disons $L=A^{n}$ pour un certain entier~$n$, ce que nous faisons ci-après. La forme quadratique s'écrit
$$
\hspace{24pt}
(x_{1},x_{2},\ldots,x_{n})\mapsto\sum_{i,j}q_{i,j}\hspace{2pt}x_{i}x_{j}
\hspace{24pt},
$$
$[q_{i,j}]:=Q$ désignant une matrice $n\times n$, définie à l'addition d'une matrice alternée près (une {\em matrice alternée} est une matrice antisymétrique avec des zéros sur la diagonale, on peut également définir une matrice alternée comme une antisymétrisée), et une matrice $M$, de taille $n\times n$, à coefficients dans une $A$-algèbre $R$, appartient à $\mathrm{O}_{L}(R):=\mathrm{O}(R\otimes_{A}L)$ si et seulement si la matrice ${}^{\mathrm{t}}\hspace{-2pt}MQM-Q$ est alternée. La proposition résulte de ce que les équations fournies par cette caractérisation ($\frac{n(n+1)}{2}$ polynômes en $n^{2}$ indéterminées à coefficients dans $A$) satisfont au critère jacobien de lissité.
\hfill$\square$
\normalsize

\medskip
On suppose maintenant que $L$ est de rang constant pair. L'application qui associe à un automorphisme orthogonal $\alpha$ d'un $\mathrm{q}$-module~$L$ son déterminant de Dickson-Dieudonné $\mathop{\widetilde{\mathrm{d\acute{e}t}}}\alpha$ (voir plus haut) induit un homomorphisme de $A$-schémas en groupes, disons $\widetilde{\mathrm{d\acute{e}t}}:\mathrm{O}_{L}\to\mathbb{Z}/2$. La proposition 1.4 implique que cet homomorphisme relève l'homomorphisme $\mathrm{d\acute{e}t}:\mathrm{O}_{L}\to\mu_{2}$, en clair que le diagramme suivant
$$
\xymatrix{  & \mathbb{Z}/2 \ar[d] \\
\mathrm{O}_{L} \ar[r]^{\mathrm{d\acute{e}t}} \ar[ru]^{\widetilde{\mathrm{d\acute{e}t}}} & \mu_{2}}
$$
est commutatif (rappelons la définition de l'homomorphisme vertical~: soit $R$ une $A$-algèbre commutative, $\mathbb{Z}/2\hspace{1pt}(R)$ est l'ensemble des éléments $x$ de $R$ vérifiant $x^2=x$, muni de la loi de groupe $(x,y)\mapsto x+y-2xy$, et l'homomorphisme  $\mathbb{Z}/2\hspace{1pt}(R)\to\mu_{2}(R)$ envoie $x$ sur $1-2x$). Pour se convaincre de cette implication, contempler le diagramme commutatif suivant
$$
\begin{CD}
0@>>>A@>\eta>>\Delta(L)@>>>\Lambda^{2n}L@>>>0
\\
& & @VV\mathrm{id}V @VV\Delta(\alpha)V @VV\Lambda^{2n}\alpha V
\\
0@>>>A@>\eta>>\Delta(L)@>>>\Lambda^{2n}L@>>>0
\end{CD}
$$
dont les lignes sont exactes.

\medskip
Le $A$-schéma en groupes $\mathrm{SO}_{L}$ est défini comme le noyau de $\widetilde{\mathrm{d\acute{e}t}}$.

\medskip
Soit $L$ un $\mathrm{q}$-module sur $A$. On note $\mathrm{GO}(L)$ le sous-groupe de $\mathrm{GL}(L)\times A^{\times}$ constitué des couples $(\alpha,\nu)$ tels que l'on a $\mathrm{q}(\alpha(x))=\nu\hspace{1pt}\mathrm{q}(x)$ pour tout $x$ dans~$L$. Le $A$-schéma en groupes $\mathrm{GO}_{L}$ est défini comme précédemment. Cette définition fait apparaître $\mathrm{GO}_{L}$ comme un sous-groupe de $\mathrm{GL}_{L}\times\mathbb{G}_{\mathrm{m}}$~; le lecteur se convaincra que la restriction de la projection $\mathrm{GL}_{L}\times\mathbb{G}_{\mathrm{m}}\to\mathrm{GL}_{L}$ à $\mathrm{GO}_{L}$ est une immersion fermée et donc que $\mathrm{GO}_{L}$ peut être aussi vu comme un sous-groupe fermé de $\mathrm{GL}_{L}$. On note $\nu:\mathrm{GO}_{L}\to\mathbb{G}_{\mathrm{m}}$ l'homomorphisme restriction de la projection $\mathrm{GL}_{L}\times\mathbb{G}_{\mathrm{m}}\to\mathbb{G}_{\mathrm{m}}$ à $\mathrm{GO}_{L}$. Le groupe $\mathrm{GO}(L)$, introduit ci-dessus, s'appelle le groupe des {\em similitudes orthogonales} de $L$, l'élément $\nu$ de~$A^{\times}$ s'appelle le {\em facteur de similitude} de $\alpha$. On note $\mathrm{PGO}_{L}$ le $A$-schéma en groupes quotient $\mathrm{GO}_{L}/\mathbb{G}_{\mathrm{m}}$.

\medskip
On suppose à nouveau que $L$ est de rang $2n$. Soit $(\alpha,\nu)$ un élément de $\mathrm{GO}(L)$~; on constate que l'élément $\nu^{-n}\mathop{\mathrm{d\acute{e}t}}\alpha$ de $A^{\times}$, disons $d$, vérifie $d^{2}=1$. On note $\mathrm{d}:\mathrm{GO}(L)\to\mu_{2}(A)$ l'homomorphisme de groupes $(\alpha,\nu)\mapsto d$. On note encore $\mathrm{d}:\mathrm{GO}_{L}\to\mu_{2}$ l'homomorphisme de $A$-schémas en groupes associé~; on obervera que $d$ prolonge l'homomorphisme $\mathrm{d\acute{e}t}:\mathrm{O}_{L}\to\mu_{2}$. La proposition 1.4 implique, comme plus haut, que $\mathrm{d}$ se relève en un homomorphisme\linebreak $\widetilde{\mathrm{d}}:\mathrm{GO}_{L}\to\mathbb{Z}/2$ qui prolonge l'homomorphisme $\widetilde{\mathrm{d\acute{e}t}}:\mathrm{O}_{L}\to\mathbb{Z}/2$. Précisons un peu. On rappelle que l'on identifie le groupe des automorphismes de $\Delta(L)$ avec $\mathbb{Z}/2\hspace{1pt}(A)$~;  $\widetilde{\mathrm{d}}(\alpha,\nu)$ correspond, \textit{via} cette identification, à la composition
$$
\begin{CD}
\hspace{24pt}\Delta(L)@>\Delta(\alpha)>>\Delta(\nu L)@>[\nu]>>\Delta(L)\hspace{24pt},
\end{CD}
$$
$\nu L$ désignant le $\mathrm{q}$-module obtenu en multipliant la forme quadratique de $L$ par $\nu$ et $[\nu]$ désignant l'isomorphisme induit par l'isomorphisme $[\nu]:\mathrm{C}^{+}(\nu L)\to\mathrm{C}^{+}(L)$ introduit dans \cite[Exposé XII, \S1.3]{sga7}. Le $A$-schéma en groupes $\mathrm{GSO}_{L}$ est défini comme le noyau de $\widetilde{\mathrm{d}}$. On constate que $\widetilde{\mathrm{d}}$ induit par passage au quotient un homomorphisme $\mathrm{PGO}_{L}\to\mathbb{Z}/2$~; le $A$-schéma en groupes $\mathrm{PGSO}_{L}$ est défini comme le noyau de cet homomorphisme induit. Alternativement, $\mathrm{PGSO}_{L}$ peut être défini comme le quotient $\mathrm{GSO}_{L}/\mathbb{G}_{\mathrm{m}}$.

\bigskip
-- Soit enfin $L$ un $\mathrm{a}$-module sur $A$. On note $\mathrm{Sp}(L)$ le groupe des automorphismes de cet $\mathrm{a}$-module. Les $A$-schémas en groupes $\mathrm{Sp}_{L}$, $\mathrm{GSp}_{L}$, $\mathrm{PGSp}_{L}$ et l'homomorphisme $\nu:\mathrm{GSp}_{L}\to\mathbb{G}_{\mathrm{m}}$ sont définis \textit{mutatis mutandis}. Si $L$ est de rang constant, disons $2n$,  la théorie du Pfaffien entraîne que l'homomorphisme $\mathrm{d\acute{e}t}:\mathrm{GSp}_{L}\to\mathbb{G}_{\mathrm{m}}$ (induit par l'homomorphisme $\mathrm{d\acute{e}t}:\mathrm{GL}_{L}\to\mathbb{G}_{\mathrm{m}}$) coïncide avec $\nu^{n}$ (relation qui implique en particulier que $\mathrm{Sp}_{L}$ est un sous-groupe de $\mathrm{SL}_{L}$). Soit $n\geq 1$ un entier~; dans le cas où $A=\mathbb{Z}$ et où $L$ est le $\mathrm{a}$-module hyper\-bolique $\mathrm{H}(\mathbb{Z}^{n})$, la notation $\mathrm{Sp}_{\mathrm{H}(\mathbb{Z}^{n})}$, $\mathrm{GSp}_{\mathrm{H}(\mathbb{Z}^{n})}$, $\mathrm{PGSp}_{\mathrm{H}(\mathbb{Z}^{n})}$, sera remplacée par la notation $\mathrm{Sp}_{2n}$, $\mathrm{GSp}_{2n}$, $\mathrm{PGSp}_{2n}$. On rappelle que $\mathrm{Sp}_{2}$ s'identifie à $\mathrm{SL}_{2}$.

\bigskip
Tous les $A$-schémas en groupes introduits ci-dessus sont affines et de présen\-tation finie sur $A$~;  pour abréger, de tels $A$-schémas en groupes seront appelés des {\em $A$-groupes}.

\bigskip
Une dernière remarque pour clore cette rubrique sur les groupes classiques~: soit $\mathrm{P}G=G/\mathbb{G}_{\mathrm{m}}$ l'un des $A$-schémas en groupes ``projectifs'' que nous venons de définir~; au chapitre IV nous aurons seulement à considérer le groupe $\mathrm{P}G(A)$ pour des anneaux $A$ avec $\mathrm{Pic}(A)=0$ si bien que le monomorphisme canonique $G(A)/A^{\times}\to\mathrm{P}G(A)$ sera un isomorphisme.

\section{Sur les $\mathrm{b}$-modules et $\mathrm{q}$-modules sur $\mathbb{Z}$}

L'objet de ce paragraphe est de rappeler quelques points, très classiques, voir par exemple \cite[Chap. V]{serre}\cite[Chap. II]{MH}, de la théorie des $\mathrm{b}$-modules et $\mathrm{q}$-modules sur $\mathbb{Z}$.

\medskip
Pour organiser ces rappels, numérotés 1, 2 et 3 ci-après, il est commode de disposer du concept d'anneau de Witt, voir par exemple \cite[Chap. I, \S7]{MH}. Rappelons-en la définition. Soit $A$ un anneau commutatif, l'ensemble des classes d'iso\-mor\-phisme de $\mathrm{b}$-modules sur $A$, disons $\mathcal{B}(A)$, est un monoïde commutatif pour la somme orthogonale, on note $\mathrm{W}(A)$ le monoïde quotient $\mathcal{B}(A)/\mathcal{N}(A)$, $\mathcal{N}(A)$ désignant le sous-monoïde engendré par les classes d'isomorphismes des $\mathrm{b}$-modules neutres (un $\mathrm{b}$-module est dit {\em neutre} s'il possède un {\em lagrangien},  c'est-à-dire un facteur direct qui est son propre orthogonal)~; le monoïde $\mathrm{W}(A)$ est un groupe et le produit tensoriel des $\mathrm{b}$-modules en fait un anneau commutatif. Le groupe abélien $\mathrm{WQ}(A)$ est défini \textit{mutatis mutandis} en termes de $\mathrm{q}$-modules sur $A$, voir par exemple \cite[App. 1]{MH} (rappelons que dans ce cas les $\mathrm{q}$-modules neutres sont en fait hyperboliques, voir Proposition 1.2)~;  $\mathrm{WQ}(A)$ est naturellement un $\mathrm{W}(A)$-module.

\bigskip
Venons-en maintenant au cas $A=\mathbb{Z}$.

\bigskip
1) Le premier des points auxquels nous faisions allusion plus haut, est la détermination de $\mathrm{W}(\mathbb{Z})$. Le résultat est le suivant~: l'homomorphisme cano\-nique $\mathrm{W}(\mathbb{Z})\to\mathrm{W}(\mathbb{R})$ est un isomorphisme. Ceci peut être reformulé de deux façons~:

\smallskip
-- L'homomorphisme ``signature'', disons $\tau:\mathrm{W}(\mathbb{Z})\to\mathbb{Z}$, est un isomorphisme. Précisons ce que l'on entend ici par {\em signature}. Soit $E$ un $\mathrm{b}$-espace vectoriel sur $\mathbb{R}$, un tel $E$ est isomorphe à un $\mathrm{b}$-espace vectoriel de la forme
$$
\langle +1,+1,\ldots,+1,-1,-1,\ldots,-1\rangle
$$
et la signature de $E$, notée $\tau(E)$ ci-après, est la différence $n_{+}-n_{-}$ entre le nombre de $+1$ et le nombre de $-1$ qui apparaissent ci-dessus. Il est clair que l'homomorphisme $\tau:\mathrm{W}(\mathbb{R})\to\mathbb{Z}$ est un isomorphisme.

\smallskip
-- L'homomorphisme ``unité'', disons $\eta:\mathbb{Z}\to\mathrm{W}(\mathbb{Z})$, est un isomorphisme. Pour une très jolie démonstration de ce résultat qui n'utilise pas le théorème de Hasse-Minkowski, voir \cite[Chap. IV, \S2]{MH}.

\bigskip
\textbf{Scholie 2.1.} {\em Soient $L_{1}$ et $L_{2}$ deux $\mathrm{b}$-modules sur $\mathbb{Z}$. Les deux conditions suivantes sont équivalentes~:
\begin{itemize}
\item [(i)] les deux $\mathrm{b}$-espaces vectoriels sur $\mathbb{Q}$, $\mathbb{Q}\otimes_{\mathbb{Z}}L_{1}$ et $\mathbb{Q}\otimes_{\mathbb{Z}}L_{2}$, sont iso\-morphes~;
\item [(ii)] les deux $\mathrm{b}$-espaces vectoriels sur $\mathbb{R}$, $\mathbb{R}\otimes_{\mathbb{Z}}L_{1}$ et $\mathbb{R}\otimes_{\mathbb{Z}}L_{2}$, sont isomorphes.
\end{itemize}}

\bigskip
2) Le deuxième point concerne la théorie des ``vecteurs de Wu''. Soit $L$ un $\mathrm{b}$-module sur $\mathbb{Z}$~; comme $\mathbb{F}_{2}\otimes_{\mathbb{Z}}L$ est un $\mathrm{b}$-espace vectoriel sur $\mathbb{F}_{2}$, il existe un élément $u$ de $L$, bien défini modulo~$2L$, tel que l'on a
$$
x.x\equiv u.x\pmod{2}
$$
pour tout $x$ dans $L$. Nous appelons $u$ un {\em vecteur de Wu} (on rencontre aussi la terminologie ``vecteur caractéristique'', la terminologie ``vecteur de Wu'' fait référence aux classes définies dans la cohomologie modulo $2$ des variétés compactes par Wen-Tsün Wu \cite{Wu}). On constate que la réduction modulo~$8$ de l'entier $u.u$ est indépendante du choix de $u$ et que la cor\-respondance $L\mapsto u.u$ induit un homomorphisme d'anneaux commutatifs unitaires, disons $\sigma:\mathrm{W}(\mathbb{Z})\to\nolinebreak\mathbb{Z}/8$. On observera que l'on peut ci-dessus remplacer $\mathbb{Z}$ par $\mathbb{Z}_{2}$ et que l'homomor\-phisme $\sigma$ se factorise à travers $\mathrm{W}(\mathbb{Z}_{2})$. Enfin, la réduction modulo $2$ de $\sigma$ se factorise par un homomorphisme $\mathrm{W}(\mathbb{F}_{2})\to\mathbb{Z}/2$ qui coïncide avec l'isomorphisme ``dimension modulo $2$''.

\medskip
Compte tenu de ce que nous avons rappelé plus haut sur $\mathrm{W}(\mathbb{Z})$, $\sigma$ s'identifie à la réduction modulo $8$ de $\mathbb{Z}$ dans $\mathbb{Z}/8$.

\bigskip
\textbf{Scholie 2.2.} {\em {\em (a)} Soient $L$ un $\mathrm{b}$-module sur $\mathbb{Z}$ et $u$ un vecteur de Wu de $L$. Alors on a la congruence
$$
\hspace{24pt}
\tau(L)\equiv u.u\pmod{8}
\hspace{24pt}.
$$

\medskip
{\em (b)} La signature d'un $\mathrm{q}$-module sur $\mathbb{Z}$ est divisible par $8$.}
\vfill\eject

\bigskip
3) Le dernier point est plus technique. Soit $L$ un $\mathrm{b}$-module sur $\mathbb{Z}$ {\em impair}, c'est-à-dire non pair. Soit $M$ le sous-module d'indice $2$ de $L$ constitué des vecteurs $x$ vérifiant $x.x\equiv 0\bmod{2}$. L'application $M\to\mathbb{Z},x\mapsto\frac{x.x}{2}$ fait de $M$ un $\widetilde{\mathrm{q}}$-module dont on détermine le résidu ci-après.

\medskip
On considère la suite exacte
$$
\hspace{24pt}
0\to L/M\to M^{\sharp}/M\to M^{\sharp}/L\to 0
\hspace{24pt}.
$$
Soient $u$ un vecteur de Wu de $L$ et $v$ un élément de $L$ avec $v.v\equiv 1\bmod{2}$ (ou encore $u.v\equiv 1\bmod{2}$); les quotients $L/M$ et $M^{\sharp}/L$ sont des  groupes cycliques d'ordre $2$ respectivement engendrés par les classes de $v$ et $\frac{u}{2}$ ($L/M$ est en fait un lagrangien du $\mathrm{e}$-module sous-jacent au $\mathrm{qe}$-module $\mathop{\mathrm{r\acute{e}s}}M$ et $M^{\sharp}/L$ est canoniquement isomorphe au dual de ce lagrangien).  La suite exacte ci-dessus est scindable si et seulement si $u$ appartient à $M$ c'est-à-dire si l'on a $u.u\equiv 0\bmod{2}$ ou encore, compte tenu de ce que nous avons rappelé dans le point 2, si la dimension de $L$ est paire. On distingue donc deux cas, suivant la parité de cette dimension~:

\medskip
-- Dans le cas $\dim L\equiv 1\bmod{2}$, $\mathop{\mathrm{r\acute{e}s}}M$ est isomorphe à $\mathbb{Z}/4$ et engendré par la classe de $\frac{u}{2}$~; on a l'égalité dans $\mathbb{Q}/\mathbb{Z}$
$$
\mathrm{q}\hspace{1pt}(x\hspace{1pt}\frac{u}{2})=\frac{\tau(L)}{8}\hspace{3pt}x^{2}
$$
pour tout $x$ dans $\mathbb{Z}$.

\medskip
-- Dans le cas $\dim L\equiv 0\bmod{2}$, $\mathop{\mathrm{r\acute{e}s}}M$ est isomorphe à $\mathbb{Z}/2\oplus\mathbb{Z}/2$ et engendré par les classes de $\frac{u}{2}$ et $v$ ou encore les classes de $\frac{u}{2}$ et $\frac{u}{2}-v$~; on a l'égalité dans $\mathbb{Q}/\mathbb{Z}$
$$
\mathrm{q}\hspace{1pt}(x\hspace{1pt}\frac{u}{2}+y\hspace{1pt}(\frac{u}{2}-v))=\frac{\tau(L)}{8}\hspace{3pt}(x^{2}+y^{2})+(\frac{\tau(L)}{4}+\frac{1}{2})\hspace{3pt}xy
$$
pour tous $x$ et $y$ dans $\mathbb{Z}$.

\medskip
On voit par inspection que l'équation $\mathrm{q}(\iota)=0$ avec $\iota\in\mathop{\mathrm{r\acute{e}s}}M-\{0\}$ n'a aucune solution pour $\tau(L)\not\equiv 0\bmod{8}$ et que pour $\tau(L)\equiv 0\bmod{8}$ elle en a exactement $2$ à savoir $\frac{u}{2}$ et $\frac{u}{2}-v$. On a donc au bout du compte obtenu en particulier l'énoncé suivant~:

\bigskip
\textbf{Scholie 2.3.} {\em Soient $L$ un $\mathrm{b}$-module sur $\mathbb{Z}$ impair et $M$ le sous-module d'indice $2$ de $L$ constitué des vecteurs $x$ vérifiant $x.x\equiv 0\bmod{2}$ ($M$ est donc un $\widetilde{\mathrm{q}}$-module sur $\mathbb{Z}$). Les conditions suivantes sont équivalentes~:
\begin{itemize}
\item [(i)] le $\mathrm{qe}$-module $\mathop{\mathrm{r\acute{e}s}}M$ est isomorphe à $\mathrm{H}(\mathbb{Z}/2)$~;
\item [(ii)] la signature de $L$ est divisible par $8$.
\end{itemize}}
\vfill\eject

\vspace{0,75cm}
\textsc{Partition en deux classes des vecteurs de Wu}

\medskip
Soient $L$ un $\mathrm{b}$-module sur $\mathbb{Z}$ impair et $M$ le sous-module d'indice $2$ de $L$ constitué des vecteurs $x$ vérifiant $x.x\equiv 0\bmod{2}$.

\medskip
Notons $\mathrm{Wu}(L)$ l'ensemble des vecteurs de Wu de $L$. L'action de $L$ sur $\mathrm{Wu}(L)$ définie par $(u,x)\mapsto u+2x$ ($u$ dans $\mathrm{Wu}(L)$ et $x$ dans $L$) est libre et transitive. L'action, induite par restriction, de $M$ sur $\mathrm{Wu}(L)$ a exactement deux orbites (celles de $u$ et $u-2v$, $u$ désignant un vecteur de Wu arbitraire de $L$ et $v$ un élément de $L$ avec $v.v\equiv 1\bmod{2}$)~; nous dirons que deux vecteurs de Wu dans la même orbite sont {\em équivalents}. Soient $u_{1}$ et $u_{2}$ deux vecteurs de Wu non équivalents~; on peut paraphraser la discussion qui précède le scholie 2.3 de la façon suivante~:

\smallskip
-- Les classes de $\frac{1}{2}\hspace{1pt}u_{1}$ et $\frac{1}{2}\hspace{1pt}u_{2}$ engendrent le groupe abélien $\mathop{\mathrm{r\acute{e}s}}M$.

\smallskip
-- Si la dimension de $L$ est impaire ces classes sont d'ordre $4$ et opposées.

\smallskip
-- Si la dimension de $L$ est paire ces classes sont d'ordre $2$ et constituent une base du $\mathbb{Z}/2$-espace vectoriel $\mathop{\mathrm{r\acute{e}s}}M$.

\smallskip
-- On a les égalités suivantes dans $\mathbb{Q}/\mathbb{Z}$~:
$$
\hspace{24pt}
\mathrm{q}(\frac{u_{1}}{2})=\frac{\tau(L)}{8}
\hspace{12pt},\hspace{12pt}
\mathrm{q}(\frac{u_{2}}{2})=\frac{\tau(L)}{8}
\hspace{12pt},\hspace{12pt}
\mathrm{q}(\frac{u_{1}}{2}+\frac{u_{2}}{2})=\frac{1+\dim L}{2}
\hspace{24pt};
$$
on observera que ces égalités déterminent dans les deux cas la forme quadratique d'enlacement $\mathrm{q}:\mathop{\mathrm{r\acute{e}s}}M\to\mathbb{Q}/\mathbb{Z}$.

\bigskip
Voici une illustration de ce qui précède. Soit $n$ un entier naturel~; on considère la forme bilinéaire symétrique ``euclidienne''
$$
\hspace{24pt}
\mathbb{Z}^{n}\times\mathbb{Z}^{n}\to\mathbb{Z}
\hspace{12pt},\hspace{12pt}
((x_{1},x_{2},\ldots,x_{n}),(y_{1},y_{2},\ldots,y_{n}))\mapsto\sum_{i=1}^{n}x_{i}y_{i}
\hspace{24pt}.
$$
Muni de cette forme, $\mathbb{Z}^{n}$ est un $\mathrm{b}$-module sur $\mathbb{Z}$ impair que l'on note $\mathrm{I}_{n}$. On note $\mathrm{D}_{n}$ le sous-module d'indice $2$ de $\mathrm{I}_{n}$ constitué des vecteurs $x$ vérifiant $x.x\equiv 0\bmod{2}$ ou encore $\sum_{i=1}^{n}x_{i}\equiv 0\bmod{2}$ (cette notation évoque sciemment la théorie des systèmes de racines, nous reviendrons sur ce thème au paragraphe 3). Soit $(\varepsilon_{1},\varepsilon_{2},\ldots,\varepsilon_{n})$ la base canonique de  $\mathrm{I}_{n}$ ; on observe que les vecteurs $u_{1}:=\varepsilon_{1}+\varepsilon_{2}+\ldots+\varepsilon_{n}$ et $u_{2}:=-\varepsilon_{1}+\varepsilon_{2}+\ldots+\varepsilon_{n}$ sont des vecteurs de Wu pour $\mathrm{I}_{n}$ et que ces vecteurs sont non équivalents. Si $n$ est divisible par $8$ alors le $\mathrm{qe}$-module $\mathop{\mathrm{r\acute{e}s}}\mathrm{D}_{n}$ peut être explicité de la façon suivante~:

\smallskip
-- comme $\mathbb{Z}$-module, $\mathop{\mathrm{r\acute{e}s}}\mathrm{D}_{n}$ est un $\mathbb{Z}/2$-espace vectoriel de dimension $2$ dont les classes des vecteurs $\iota_{1}:=\frac{1}{2}\hspace{1pt}u_{1}$ et $\iota_{2}:=\frac{1}{2}\hspace{1pt}u_{2}$ constituent une base~;

\smallskip
-- sa forme quadratique d'enlacement est déterminée par $\mathrm{q}(\iota_{1})=0$, $\mathrm{q}(\iota_{2})=0$ et $\iota_{1}.\iota_{2}=\frac{1}{2}$.

\bigskip
Les points (b) et (c) de la proposition 1.1 montrent que le réseau de $\mathbb{Q}\otimes_{\mathbb{Z}}\mathrm{I}_{n}$ engendré par $\mathrm{D}_{n}$ et $\frac{1}{2}(\varepsilon_{1}+\varepsilon_{2}+\ldots+\varepsilon_{n})$ est un $\mathrm{q}$-module sur $\mathbb{Z}$ que l'on note $\mathrm{E}_{8}$ pour $n=8$ et $\mathrm{E}_{n}$ ou $\mathrm{D}_{n}^{+}$ pour $n\geq 16$ (et $n\equiv 0\bmod{8}$).

\bigskip
\textbf{Scholie 2.4.} {\em L'homomorphisme composé
$$
\begin{CD}
\mathrm{WQ}(\mathbb{Z})@>\mathrm{oubli}>>\mathrm{W}(\mathbb{Z})@>\tau>>\mathbb{Z}
\end{CD}
$$
induit un isomorphisme (de $\mathrm{W}(\mathbb{Z})$-modules) de $\mathrm{WQ}(\mathbb{Z})$ sur l'idéal $8\hspace{1pt}\mathbb{Z}$~; le groupe $\mathrm{WQ}(\mathbb{Z})$ est infini cyclique engendré par la classe de $\mathrm{E}_{8}$.}

\bigskip
\textit{Démonstration.} Il suffit d'observer que l'oubli $\mathrm{WQ}(\mathbb{Z})\to\mathrm{W}(\mathbb{Z})$ est injectif.

\vspace{0,75cm}
\textsc{Genre d'un $\mathrm{q}$-module sur $\mathbb{Z}$}

\bigskip
Cet intertitre fait référence au point (b) de l'énoncé ci-dessous.

\bigskip
\textbf{Scholie 2.5.} {\em Soient $L$ un $\mathrm{q}$-module sur $\mathbb{Z}$ et $p$ un nombre premier.

\medskip
{\em (a)} Le $\mathrm{q}$-espace vectoriel $\mathbb{F}_{p}\otimes_{\mathbb{Z}}L$ est hyperbolique.

\medskip
{\em (b)} Le $\mathrm{q}$-module $\mathbb{Z}_{p}\otimes_{\mathbb{Z}}L$ est hyperbolique.}

\bigskip
\textit{Démonstration.} Le premier point implique le second~; en effet, deux $\mathrm{q}$-modules $L_{1}$ et $L_{2}$ sur $\mathbb{Z}_{p}$ sont isomorphes si et seulement si les $\mathrm{q}$-espaces vectoriels $\mathbb{F}_{p}\otimes_{\mathbb{Z}_{p}}L_{1}$ et $\mathbb{F}_{p}\otimes_{\mathbb{Z}_{p}}L_{2}$ le sont. Venons-en au point (a). Il suffit de montrer que l'homomorphisme naturel $\mathrm{WQ}(\mathbb{Z})\to\mathrm{WQ}(\mathbb{F}_{p})$ est trivial. Le cas $p$ impair, facile, est laissé au lecteur. Pour $p=2$, on peut invoquer les arguments suivants~:

\smallskip
-- L'invariant de Arf, disons $\mathrm{Arf}:\mathrm{WQ}(\mathbb{F}_{2})\to\mathrm{H}^{1}_{\mathrm{\acute{e}}t}(\mathbb{F}_{2};\mathbb{Z}/2)\cong\mathbb{Z}/2$, est un isomorphisme.

\smallskip
-- Le groupe 
$\mathrm{H}^{1}_{\mathrm{\acute{e}}t}(\mathbb{Z};\mathbb{Z}/2)$ est trivial.

\smallskip
L'argument ci-dessus peut être remplacé par l'argument ci-dessous, plus prosaïque~:

\smallskip
-- L'homomorphisme $\mathrm{WQ}(\mathbb{Z})\to\mathrm{WQ}(\mathbb{F}_{2})$ se factorise à travers $\mathrm{WQ}(\mathbb{Z}_{2})$ et le déterminant d'un $\mathrm{q}$-module $L$ sur $\mathbb{Z}_{2}$ est égal à la classe dans $\mathbb{Z}_{2}^{\times}/\mathbb{Z}_{2}^{\times\hspace{1pt}2}$ de l'élément $(-1)^{\frac{\dim L}{2}}\hspace{2pt}(-3)^{\mathrm{Arf}(\mathbb{F}_{2}\otimes_{\mathbb{Z}_{2}}L)}$ (pour se convaincre de cette égalité observer par exemple qu'un tel $L$ se décompose en une somme orthogonale de $\mathrm{q}$-modules de dimension $2$, chacun muni d'une base $(e,f)$ avec $e.f=1$).
\hfill$\square$

\bigskip
On peut alternativement déduire l'énoncé 2.5 de l'énoncé ci-dessous dont la démonstration est laissée au lecteur~:

\bigskip
\textbf{Proposition 2.6.} {\em Soit $p$ un nombre premier~; soient $L_{1}$ et $L_{2}$ deux $\mathrm{q}$-modules sur $\mathbb{Z}_{p}$. Les conditions suivantes sont équivalentes~:
\begin{itemize}
\item [(i)] $L_{1}$ et $L_{2}$ sont isomorphes~;
\item [(ii)] les deux $\mathrm{q}$-espaces vectoriels sur $\mathbb{F}_{p}$, $\mathbb{F}_{p}\otimes_{\mathbb{Z}}L_{1}$ et $\mathbb{F}_{p}\otimes_{\mathbb{Z}}L_{2}$, sont isomorphes~;
\item [(iii)] $L_{1}$ et $L_{2}$ ont même dimension et même déterminant~;
\item [(iv)] les deux $\mathrm{q}$-espaces vectoriels sur $\mathbb{Q}_{p}$, $\mathbb{Q}_{p}\otimes_{\mathbb{Z}}L_{1}$ et $\mathbb{Q}_{p}\otimes_{\mathbb{Z}}L_{2}$, sont isomorphes~;
\end{itemize}}

\bigskip
En vue de futures références, nous ajoutons aux rappels précédents, les deux énoncés suivants qui sont des raffinements de l'énoncé 2.1.

\bigskip
\textbf{Théorème 2.7.} {\em Soient $L_{1}$ et $L_{2}$ deux $\mathrm{q}$-modules sur $\mathbb{Z}$. On suppose que le $\mathrm{q}$-espace vectoriel sur $\mathbb{R}$, $\mathbb{R}\otimes_{\mathbb{Z}}L_{1}$ est indéfini, alors les deux conditions suivantes sont équivalentes~:
\begin{itemize}
\item [(i)] $L_{1}$ et $L_{2}$ sont isomorphes~;
\item [(ii)] les deux $\mathrm{q}$-espaces vectoriels sur $\mathbb{R}$, $\mathbb{R}\otimes_{\mathbb{Z}}L_{1}$ et $\mathbb{R}\otimes_{\mathbb{Z}}L_{2}$, sont isomorphes.
\end{itemize}}

\bigskip
\textbf{Théorème 2.8.} {\em Soient $L_{1}$ et $L_{2}$ deux $\mathrm{q}$-modules sur $\mathbb{Z}$ et $p$ un nombre premier. Les deux conditions suivantes sont équivalentes~:
\begin{itemize}
\item [(i)] les deux $\mathrm{q}$-espaces vectoriels sur $\mathbb{Z}[\frac{1}{p}]$, $\mathbb{Z}[\frac{1}{p}]\otimes_{\mathbb{Z}}L_{1}$ et $\mathbb{Z}[\frac{1}{p}]\otimes_{\mathbb{Z}}L_{2}$, sont isomorphes~;
\item [(ii)] les deux $\mathrm{q}$-espaces vectoriels sur $\mathbb{R}$, $\mathbb{R}\otimes_{\mathbb{Z}}L_{1}$ et $\mathbb{R}\otimes_{\mathbb{Z}}L_{2}$, sont isomorphes.
\end{itemize}}

\bigskip
Les deux énoncés peuvent se démontrer à l'aide du théorème d'approximation forte pour les groupes $\mathrm{Spin}$ \cite{platonovrap}. On peut aussi démontrer 2.7 de la façon suivante~: commencer par observer que 2.1 implique qu'un $\mathrm{b}$-module sur $\mathbb{Z}$ indéfini représente zéro et achever comme le fait Serre dans \cite[Ch. V, \S 3]{serre}.

\vspace{0,75cm}
\textsc{Terminologie classique dans le cas défini positif}

\medskip
Soit $V$ un espace euclidien, en clair un $\mathbb{R}$-espace vectoriel de dimension finie muni d'une forme bilinéaire symétrique définie positive (appelée produit scalaire).  Dans ce contexte, un {\em réseau} de $V$ est un sous-groupe discret cocompact. Soit $L$ un tel réseau, on dit que $L$ est {\em entier} si le produit scalaire $x.y$ est entier pour tous $x$ et $y$ dans $L$~; muni de la forme bilinéaire symétrique induite par le produit scalaire, $L$ est un $\widetilde{\mathrm{b}}$-module sur $\mathbb{Z}$, défini positif (c'est-à-dire que le $\widetilde{\mathrm{b}}$-espace vectoriel $\mathbb{R}\otimes_{\mathbb{Z}}L$ est défini positif). Réciproquement un $\widetilde{\mathrm{b}}$-module sur $\mathbb{Z}$ défini positif $L$ est un réseau entier dans l'espace euclidien $\mathbb{R}\otimes_{\mathbb{Z}}L$. Un réseau entier $L$ de $V$ est appelé {\em réseau unimodulaire} s'il est de covolume~$1$ (l'adjectif entier est donc implicite dans ce cas), en d'autres termes si le $\widetilde{\mathrm{b}}$-module $L$ est en fait un $\mathrm{b}$-module. Un réseau entier $L$ de $V$ est dit {\em pair} si $x.x$ est pair pour tout $x$ dans $L$, en d'autres termes si le $\widetilde{\mathrm{b}}$-module $L$ est en fait un $\widetilde{\mathrm{q}}$-module. Un $\mathrm{q}$-module sur $\mathbb{Z}$ défini positif $L$ est donc un réseau unimodulaire pair dans l'espace euclidien $\mathbb{R}\otimes_{\mathbb{Z}}L$.

\medskip
Dans ce mémoire nous emploierons plus souvent la terminologie classique, réseau unimodulaire pair (resp. réseau unimodulaire, réseau entier pair, réseau entier), que la terminologie $\mathrm{q}$-module (resp. $\mathrm{b}$-module, $\widetilde{\mathrm{q}}$-module, $\widetilde{\mathrm{b}}$-module) sur $\mathbb{Z}$ défini positif (en vérité, la terminologie ``$\mathrm{q}$-module'', ``$\mathrm{b}$-module'', \ldots, n'est guère employée que par le second auteur~!).

\section{Systèmes de racines et réseaux unimodulaires pairs}

Les réseaux unimodulaires pairs qui apparaissent en dimension $8$, $16$ et $24$, à l'exception du réseau de Leech, sont tous construits à partir de certains systèmes de racines par un procédé que nous allons décrire.

\bigskip
En fait les systèmes de racines en question sont certaines  sommes directes de systèmes de racines de type $\mathbf{A}_{l}$, $\mathbf{D}_{l}$, $\mathbf{E}_{6}$, $\mathbf{E}_{7}$ et $\mathbf{E}_{8}$~; nous dirons que de telles sommes directes sont des sytèmes de racines de {\em type ADE} (la terminologie anglo-saxonne, justifiée par la considération du graphe de Dynkin, est {\em simply laced}). Il est clair que les systèmes de racines irréductibles de type ADE sont caractérisés, parmi tous les systèmes de racines irréductibles, par la propriété que toutes les racines sont de même longueur. Le lecteur vérifiera sans peine  que l'on peut prendre comme définition des systèmes de racines de type ADE, la définition \textit{ad hoc} suivante, variante de \cite[Ch. VI, \S1, Déf. 1]{bourbaki}~:

\bigskip
\textbf{Définition 3.1.} {\em Soient $V$ un espace euclidien et $R$ un sous-ensemble de $V$ constitué de vecteurs $\alpha$ vérifiant $\alpha.\alpha=2$. Nous dirons que $R$ est un {\em système de racines de type ADE} dans $V$ si les conditions suivantes sont vérifiées~:
\begin{itemize}
\item [{\em (I)}] Le sous-ensemble $R$ est fini et engendre $V$.
\item [{\em (II)}] Pour tout $\alpha$ dans $R$ la réflexion orthogonale $x\mapsto x-(\alpha.x)\hspace{2pt}\alpha$ de $V$ (notée $\mathrm{s}_{\alpha}$) laisse stable $R$.
\item [{\em (III)}] Pour tous $\alpha$ et $\beta$ dans $R$ le produit scalaire $\alpha.\beta$ est entier.
\end{itemize}}

\bigskip
Le groupe de Weyl $\mathrm{W}(R)$ engendré par les $\mathrm{s}_{\alpha}$ apparaît ici comme un sous-groupe du groupe orthogonal de $V$~; il en est de même du groupe (fini) $\mathrm{A}(R)$ constitué des automorphismes du $\mathbb{R}$-espace vectoriel $V$ qui laissent stable $R$ (ceci résulte par exemple de \cite[Ch. VI, \S1, Prop. 3 et Prop. 7]{bourbaki}). Si l'on identifie $V$ et $V^{*}$ \textit{via} le produit scalaire alors l'application structurelle $\alpha\mapsto\alpha^{\vee}$ est l'identité, en particulier les deux systèmes de racines $R$ et $R^{\vee}$ coïncident.

\bigskip
Soit $R\subset V$ un système de racines de type ADE~; le réseau de $V$ engendré par $R$, noté $\mathrm{Q}(R)$, est un réseau entier pair et $R$ s'identifie au sous-ensemble de $\mathrm{Q}(R)$ constitué des éléments $\alpha$ vérifiant $\alpha.\alpha=2$ (cette dernière propriété n'est pas, \textit{a priori}, une conséquence immédiate de 3.1, on peut se convaincre de ce qu'elle est satisfaite en observant qu'elle l'est pour les systèmes de racines $\mathbf{A}_{l}$, $\mathbf{D}_{l}$, $\mathbf{E}_{6}$, $\mathbf{E}_{7}$ et $\mathbf{E}_{8}$). On observera que $\mathrm{A}(R)$ s'identifie au groupe des automorphismes de $\mathrm{Q}(R)$, comme $\widetilde{\mathrm{q}}$-module sur $\mathbb{Z}$. La notation $\mathrm{Q}(R)$ est celle de Bourbaki \cite{bourbaki} qui appelle {\em poids radiciels} les éléments de ce réseau. Réciproquement, soient $L$ un réseau entier (en d'autres termes un $\widetilde{\mathrm{b}}$-module sur $\mathbb{Z}$, défini positif), $\mathrm{R}(L)$ le sous-ensemble (fini) de $L$ constitué des éléments $\alpha$ vérifiant $\alpha.\alpha=2$, que l'on appelle les {\em racines} de $L$, et $\mathrm{V}(L)$ le sous-espace de $\mathbb{R}\otimes_{\mathbb{Z}}L$ engendré par $\mathrm{R}(L)$, alors $\mathrm{R}(L)$ est un système de racines dans $\mathrm{V}(L)$ de type ADE (prendre $\Lambda=\{2\}$ au début du \no{4} de \cite[Ch. VI, \S4]{bourbaki}). En conclusion, la classification des systèmes de racines de type ADE coïncide avec celle des réseaux entiers (pairs) engendrés par leurs racines.

\bigskip
Le procédé qui permet d'obtenir un réseau unimodulaire pair à partir de certains systèmes de racines de type $\mathrm{ADE}$ est simplement un cas particulier du procédé général fourni par la proposition 1.1. Soit $R$ un système de racines de type $\mathrm{ADE}$. Supposons que le $\mathrm{qe}$-module $\mathop{\mathrm{r\acute{e}s}}\mathrm{Q}(R)$ possède un lagrangien~$I$,  alors l'image réciproque de $I$ dans le réseau $\mathrm{Q}(R)^{\sharp}$, par l'application canoni\-que $\mathrm{Q}(R)^{\sharp}\to\mathrm{Q}(R)^{\sharp}/\mathrm{Q}(R)=:\mathop{\mathrm{r\acute{e}s}}\mathrm{Q}(R)$, est un $\mathrm{q}$-module (voir les points (b) et (c) de la proposition 1.1) sur $\mathbb{Z}$ défini positif, en d'autres termes un réseau unimodulaire pair. Le réseau $\mathrm{Q}(R)^{\sharp}$ est le réseau des {\em poids} du système de racines $R$~; il est noté $\mathrm{P}(R)$ dans \cite{bourbaki}.

\bigskip
\textit{Exemple}

\medskip
Soit $n\geq 1$ un entier~; on munit le $\mathbb{R}$-espace vectoriel $\mathbb{R}^{n}$ de sa structure euclidienne canonique et on note $(\varepsilon_{1},\varepsilon_{2},\ldots,\varepsilon_{n})$ sa base canonique.

\medskip
Considérons à nouveau le réseau entier pair $\mathrm{D}_{n}\subset\mathrm{I}_{n}:=\mathbb{Z}^{n}\subset\mathbb{R}^{n}$ introduit au paragraphe 2. Le système de racines $\mathbf{D}_{n}$, $n\geq 3$, est défini par l'égalité $\mathbf{D}_{n}:=\mathrm{R}(\mathrm{D}_{n})$~; $\mathrm{R}(\mathrm{D}_{n})$ engendre $\mathrm{D}_{n}$ pour $n\geq 2$. Rappelons ce que l'on a vu au paragraphe 2~:

\smallskip
-- Le $\mathrm{qe}$-module $\mathop{\mathrm{r\acute{e}s}}\mathrm{D}_{n}$ possède un lagrangien si et seulement si $n$ est divisible par $8$, ce que l'on suppose ci-après.

\smallskip
-- On a dans ce cas un isomorphisme de $\mathrm{qe}$-modules $\mathop{\mathrm{r\acute{e}s}}\mathrm{D}_{n}\cong\mathrm{H}(\mathbb{Z}/2)$ et les deux lagrangiens de $\mathop{\mathrm{r\acute{e}s}}\mathrm{D}_{n}$ sont respectivement engendrés par les classes des vecteurs $\frac{1}{2}(\varepsilon_{1}+\varepsilon_{2}+\ldots+\varepsilon_{n})$ et  $\frac{1}{2}(-\varepsilon_{1}+\varepsilon_{2}+\ldots+\varepsilon_{n})$~; on observera incidemment que ces deux vecteurs sont échangés par l'automorphisme involutif $(x_{1},x_{2},\ldots x_{n})\mapsto(-x_{1},x_{2},\ldots x_{n})$ de~$\mathrm{D}_{n}$.

\medskip
On a noté $\mathrm{E}_{n}$ le réseau engendré par $\mathrm{D}_{n}$ et $\frac{1}{2}(\varepsilon_{1}+\varepsilon_{2}+\ldots+\varepsilon_{n})$~;  $\mathrm{E}_{n}$~est l'exemple le plus simple de réseau unimodulaire pair obtenu par le procédé décrit ci-dessus (pour des exemples plus sophistiqués, voir la classification des réseaux unimodulaires pairs de dimension 24 évoquée ci-après).

\medskip
Le système de racines $\mathbf{E}_{8}$ est défini par l'égalité $\mathbf{E}_{8}:=\mathrm{R}(\mathrm{E}_{8})$~; $\mathrm{R}(\mathrm{E}_{8})$ engendre $\mathrm{E}_{8}$, en d'autres termes on a l'égalité $\mathrm{E}_{8}:=\mathrm{Q}(\mathbf{E}_{8})$.

\medskip
Pour $n\geq 16$ (et $n\equiv 0\bmod{8}$) on a $\mathrm{R}(\mathrm{E}_{n})=\mathrm{R}(\mathrm{D}_{n})=\mathbf{D}_{n}$~; le sous-groupe engendré par $\mathrm{R}(\mathrm{E}_{n})$ est d'indice $2$ dans $\mathrm{E}_{n}$.

\medskip
Le réseau unimodulaire pair $\mathrm{E}_{n}$ est aussi noté $\mathrm{D}_{n}^{+}$ (au moins pour $n\geq 16$)~; cette notation est justifiée par le fait que le groupe des automorphismes de~$\mathrm{D}_{n}$ opère transitivement sur l'ensemble des deux lagrangiens de $\mathop{\mathrm{r\acute{e}s}}\mathrm{D}_{n}$.

\vspace{0,75cm}
\textsc{Classification des réseaux unimodulaires pairs en dimension $8$, $16$ et $24$}

\bigskip
Cette classification est due à Louis J. Mordell en dimension $8$, à Ernst Witt en dimension $16$ et à Hans-Volker Niemeier en dimension $24$ \cite{Ni}.  Ci-après, nous évoquons la stratégie ingénieuse développée par Boris Venkov \cite{venkov} pour retrouver la classification de Niemeier. Sa stratégie fonctionne également en dimension 8 et 16. En effet, l'idée de départ de Venkov est de considérer des séries thêta ``à coefficients harmoniques de degré $2$'' et d'observer que ces séries sont identiquement nulles car toute forme modulaire parabolique pour  $\mathrm{SL}_{2}(\mathbb{Z})$ de poids $14=\frac{24}{2}+2$ est nulle, les arguments qu'il utilise ensuite ne font plus intervenir de façon essentielle que la dimension des réseaux est~$24$~; or toute forme modulaire parabolique pour  $\mathrm{SL}_{2}(\mathbb{Z})$ de poids $\frac{n}{2}+2$ est également nulle pour $n=8,16$.

\bigskip
En considérant le ``coefficient de $e^{2\imath\pi\tau}$'' dans les séries thêta évoquées ci-dessus, Venkov obtient l'identité suivante~:

\bigskip
\textbf{Proposition 3.2.} {\em Soit $L$ un réseau unimodulaire pair de dimension $n=8,16,24$. Alors on a l'identité
$$
\sum_{\alpha\in\mathrm{R}(L)}(\alpha.x)^{2}
\hspace{4pt}=\hspace{4pt}
\frac{2\hspace{2pt}\vert\mathrm{R}(L)\vert}{n}
\hspace{4pt}x.x
$$
pour tout $x$ dans l'espace euclidien $\mathbb{R}\otimes_{\mathbb{Z}}L$  (la notation $\vert-\vert$ désignant le cardinal d'un ensemble fini).}

\bigskip
Venkov déduit ensuite de cette identité l'énoncé ci-dessous~:

\bigskip
\textbf{Proposition-Définition 3.3.} {\em Soit $L$ un réseau unimodulaire pair de\linebreak dimension $n=8,16,24$. Si l'ensemble $\mathrm{R}(L)$ des racines de $L$ (rappelons qu'il s'agit là des éléments $\alpha$ de $L$ vérifiant $\alpha.\alpha=2$) est non vide alors il vérifie les propriétés suivantes~:

\medskip
{\em (a)} L'ensemble $\mathrm{R}(L)$ est un système de racines (de type ADE) de rang $n$ (dans $\mathbb{R}\otimes_{\mathbb{Z}}L$), en d'autres termes $\mathrm{R}(L)$ engendre le $\mathbb{R}$-espace vectoriel $\mathbb{R}\otimes_{\mathbb{Z}}L$.

\medskip
{\em (b)} Toutes les composantes irréductibles du système de racines $\mathrm{R}(L)$ ont le même nombre de Coxeter, que l'on appelle le nombre de Coxeter de $L$ et que l'on note $\mathrm{h}(L)$~; un tel système de racines sera dit {\em équicoxeter}.

\medskip
{\em (c)} On a $\vert\mathrm{R}(L)\vert=n\hspace{2pt}\mathrm{h}(L)$.}

\bigskip
\textit{Remarques}

\smallskip
-- Venkov démontre simultanément (b) et (c) à l'aide de 3.2 et \cite[Ch. VI, \S1, Prop. 32]{bourbaki}. On peut également démontrer (c) à partir de (b) en invoquant la relation $\vert R\vert=n\hspace{1pt}h$ qui lie nombre de racines, rang et nombre de Coxeter, pour tout système de racines irréductible réduit (voir \cite[3.18]{Hu}\cite[Ch. VI, \S1, Exc. 20]{bourbaki}).

\smallskip
-- Puisque le nombre de Coxeter de la somme directe de deux systèmes de racines est le p.p.c.m. de leurs nombres de Coxeter (se rappeler que le nombre de Coxeter d'un système de racines est défini comme l'ordre d'une transformation de Coxeter), $\mathrm{h}(L)$ est aussi le nombre de Coxeter du système de racines $\mathrm{R}(L)$.

\bigskip
\textbf{Scholie 3.4.} {\em Soit $L$ un réseau unimodulaire pair de dimension $24$ avec\linebreak $\mathrm{R}(L)\not=\emptyset$. Alors on a l'identité
$$
\sum_{\alpha\in\mathrm{R}(L)}(\alpha.x)^{2}
\hspace{4pt}=\hspace{4pt}
2\hspace{2pt}\mathrm{h}(L)\hspace{4pt}x.x
$$
pour tout $x$ dans l'espace euclidien $\mathbb{R}\otimes_{\mathbb{Z}}L$.}

\bigskip
\textbf{Corollaire 3.5.} {\em Tout réseau unimodulaire pair de dimension $8$ est isomorphe à $\mathrm{E}_{8}$.}

\bigskip
\textit{Démonstration.} Soit $L$ un tel réseau. La série thêta de $L$, modulaire de poids~$4$ pour $\mathrm{SL}_{2}(\mathbb{Z})$, est nécessairement égale à la série d'Eisenstein normalisée~$\mathbb{E}_{4}$ (la surutilisation des E est responsable de cette notation inha\-bituelle). On en déduit $\vert\mathrm{R}(L)\vert=240$ et $\mathrm{h}(L)=30$. Cette dernière égalité force par in\-spection $\mathrm{R}(L)\simeq\mathbf{E}_{8}$. Comme le réseau $\mathrm{E}_{8}=\mathrm{Q}(\mathbf{E}_{8})$ est unimo\-dulaire, on a bien\linebreak $L\simeq\mathrm{E}_{8}$.
\hfill$\square$

\bigskip
\textbf{Corollaire 3.6.} {\em Tout réseau unimodulaire pair de dimension $16$ est isomorphe à $\mathrm{D}_{16}^{+}$ ou à $\mathrm{E}_{8}\oplus\mathrm{E}_{8}$ (et ces deux réseaux ne sont pas isomorphes).}

\bigskip
\textit{Démonstration.} Soit $L$ un tel réseau. La série thêta de $L$ est nécessairement égale à la série d'Eisenstein normalisée~$\mathbb{E}_{8}=\mathbb{E}_{4}^{2}$. On en déduit $\vert\mathrm{R}(L)\vert=480$ et $\mathrm{h}(L)=30$. Cette dernière égalité force par inspection $\mathrm{R}(L)\simeq\mathbf{D}_{16}$ ou $\mathrm{R}(L)\simeq\mathbf{E}_{8}\coprod\mathbf{E}_{8}$ (nous notons $\coprod$ la {\em somme directe des systèmes de racines}).

\smallskip
Dans le cas $\mathrm{R}(L)\simeq\mathbf{E}_{8}\coprod\mathbf{E}_{8}$ on conclut comme précédemment que l'on a $L\simeq\mathrm{E}_{8}\oplus\mathrm{E}_{8}$.

\smallskip
Dans le cas $\mathrm{R}(L)\simeq\mathbf{D}_{16}$, $\mathrm{D}_{16}$ apparaît, à isomorphisme près, comme un sous-réseau de $L$~; on peut donc supposer que l'on a dans $\mathbb{Q}\otimes_{\mathbb{Z}}L$ les inclusions $\mathrm{D}_{16}\subset L\subset\mathrm{D}_{16}^{\sharp}$. Alors $L/\mathrm{D}_{16}$ est un lagrangien du $\mathrm{qe}$-module $\mathop{\mathrm{r\acute{e}s}}\mathrm{D}_{16}$ et l'on voit que l'on a $L\simeq\mathrm{D}_{16}^{+}$.
\hfill$\square$

\bigskip
Venons-en maintenant au choses sérieuses, c'est-à-dire à la détermination des classes d'isomorphisme de réseaux unimodulaires pairs de dimension $24$.

\medskip
Soit $L$ un réseau unimodulaire pair de dimension $24$ avec $\mathrm{R}(L)\not=\emptyset$. Les propriétés (a) et (b) de la proposition 3.3 nous disent que $\mathrm{R}(L)$ est un système de racines de type ADE, de rang $24$ et équicoxeter. Venkov commence par expliciter (par inspection) la liste des classes d'isomorphisme de tels systèmes de racines. Cette liste comporte $23$ éléments $\mathbf{R}_{1},\mathbf{R}_{2},\ldots,\mathbf{R}_{23}$. Pour la liste complète nous renvoyons à la deuxième colonne de \cite[Ch. 16,  Table 16.1]{conwaysloane} (Table \ref{listeniemeier} de notre mémoire)~; le lecteur notera que nous indexons les termes de cette colonne par les entiers $1,2,\ldots,23$, plutôt que par les lettres grecques $\alpha,\beta,\ldots,\psi$, comme le font Conway et Sloane. Voici quelques exemples~:

\medskip
$\mathbf{R}_{1}=\mathbf{D}_{24}$\hspace{6pt},\hspace{6pt}$\mathbf{R}_{2}=\mathbf{D}_{16}\coprod\mathbf{E}_{8}$\hspace{6pt},\hspace{6pt}$\mathbf{R}_{3}=\mathbf{E}_{8}\coprod\mathbf{E}_{8}\coprod\mathbf{E}_{8}$\hspace{6pt},\hspace{6pt}$\mathbf{R}_{4}=\mathbf{A}_{24}$\hspace{6pt},\hspace{6pt}$\mathbf{R}_{5}=\mathbf{D}_{12}\coprod\mathbf{D}_{12}$\hspace{6pt},\hspace{6pt}$\mathbf{R}_{6}=\mathbf{A}_{17}\coprod\mathbf{E}_{7}$\hspace{6pt},\hspace{6pt}$\mathbf{R}_{7}=\mathbf{D}_{10}\coprod\mathbf{E}_{7}\coprod\mathbf{E}_{7}$\hspace{6pt},\hspace{6pt}$\mathbf{R}_{23}=\mathbf{A}_{1}\coprod\mathbf{A}_{1}\linebreak\coprod\ldots\coprod\mathbf{A}_{1}$ ($24$ composantes irréductibles toutes égales à $\mathbf{A}_{1}$).

\medskip
Avant de poursuivre notre exposition des arguments de Venkov il nous faut faire quelques observations et rappels numérotés 1, 2 et 3 ci-aprés.

\bigskip
1) Soit $M$ un $\widetilde{\mathrm{q}}$-module sur $\mathbb{Z}$. Si $M$ est défini positif alors le quotient $\mathop{\mathrm{r\acute{e}s}}M:=M^{\sharp}/M$ possède, en plus de sa structure de $\mathrm{qe}$-module, une structure supplé\-mentaire que nous nous proposons de décrire. Soient $\xi$ un élément de $\mathop{\mathrm{r\acute{e}s}}M$ et $\gamma:M^{\sharp}\to\mathop{\mathrm{r\acute{e}s}}M$ la surjection canonique, on définit une application $\mathrm{qm}:\mathop{\mathrm{r\acute{e}s}}M\to\mathbb{Q}\cap[0,\infty[$ en posant
$$
\hspace{24pt}
\mathrm{qm}(\xi)=\inf_{x\in\gamma^{-1}(\xi)}\mathrm{q}(x)
\hspace{24pt}.
$$
Il est clair que cette application fait commuter le diagramme suivant
$$
\xymatrix{  & \mathbb{Q}\cap[0,\infty[ \ar[d] \\
\mathop{\mathrm{r\acute{e}s}}M \ar[r]^{\mathrm{q}} \ar[ru]^{\mathrm{qm}} & \mathbb{Q}/\mathbb{Z}}
$$
dans lequel la flèche verticale est la restriction de la réduction modulo $\mathbb{Z}$. Nous dirons que $\mathop{\mathrm{r\acute{e}s}}M$ muni de cette structure supplémentaire est un {\em $\mathrm{qe}$-module de Venkov}.

\medskip
\textit{Exemple.} Soit $n>0$ un entier divisible par $8$~; on prend $M=\mathrm{D}_{n}$. On a vu au paragraphe $2$ que $\mathop{\mathrm{r\acute{e}s}}\mathrm{D}_{n}$ est isomorphe à $\mathbb{Z}/2\oplus\mathbb{Z}/2$ muni de la forme quadratique d'enlacement définie par $\mathrm{q}(\bar{0},\bar{0})=0$, $\mathrm{q}(\bar{1},\bar{0})=0$, $\mathrm{q}(\bar{0},\bar{1})=0$ et $\mathrm{q}(\bar{1},\bar{1})=\frac{1}{2}$~; l'application $\mathrm{qm}$ est donnée quant à elle par $\mathrm{qm}(\bar{0},\bar{0})=0$, $\mathrm{qm}(\bar{1},\bar{0})=\frac{n}{8}$, $\mathrm{qm}(\bar{0},\bar{1})=\frac{n}{8}$ et $\mathrm{qm}(\bar{1},\bar{1})=\frac{1}{2}$.

\medskip
La proposition suivante est évidente~:

\bigskip
\textbf{Proposition 3.7.} {\em Soit $M$ un réseau entier pair. Soient $I$ un sous-module de $\mathop{\mathrm{r\acute{e}s}}M$ avec $I\subset I^{\perp}$ et $L$ le réseau entier associé. Alors les deux conditions suivantes sont équivalentes~:
\begin{itemize}
\item [(i)] on a $\mathrm{R}(L)=\mathrm{R}(M)$~;
\item [(ii)] on a $\mathrm{qm}(\xi)>1$ pour tout $\xi$ dans $I-\{0\}$.
\end{itemize}}

\bigskip
\textit{Exemple.} On reprend l'exemple précédent. La proposition montre que l'on a $\mathrm{R}(\mathrm{E}_{n})=\mathrm{R}(\mathrm{D}_{n})$ pour $n\geq 16$ et $\mathrm{R}(\mathrm{E}_{8})\supsetneqq\mathrm{R}(\mathrm{D}_{8})$.

\bigskip
2) Soit $R\subset V$ un système de racines de type ADE~; on rappelle maintenant comment déterminer l'application $\mathrm{qm}:\mathop{\mathrm{r\acute{e}s}}\mathrm{Q}(R)\to\mathbb{Q}\cap [0,\infty[$. Il est clair que l'on peut supposer $R$ irréductible, ce que l'on fait ci-dessous.

\medskip
On fixe une chambre $C$ de $R$~; on note respectivement $(\alpha_{1},\alpha_{2},\ldots,\alpha_{l})$ et $\widetilde{\alpha}$, la base de $R$ et la plus grande racine correspondantes. On rappelle que $C$ est le sous-ensemble de $V$ constitué des éléments $\xi$ vérifiant $\alpha_{i}.\xi>0$ pour $i=1,2,\ldots,l$. On rappelle également que l'on a $\widetilde{\alpha}=n_{1}\alpha_{1}+n_{2}\alpha_{2}+\dots+n_{l}\alpha_{l}$ avec $n_{i}\in\mathbb{N}-\{0\}$ pour $i=1,2,\ldots,l$~; on note $\mathrm{J}$ le sous-ensemble de $\{1,2,\ldots,l\}$ constitué des indices $i$ pour lesquels on a $n_{i}=1$. On note $\varpi_{1},\varpi_{2},\ldots,\varpi_{l}$ les {\em poids fondamentaux}, c'est-à-dire les éléments de $V$ définis par $\alpha_{i}.\varpi_{j}=\delta_{i,j}$ (symbole de Kronecker)~;  $(\varpi_{1},\varpi_{2},\ldots,\varpi_{l})$ est une base du $\mathbb{Z}$-module $\mathrm{Q}(R)^{\sharp}$ et $C$ est le cône ouvert de $V$ engendré par cette base.

\bigskip
\textbf{Proposition 3.8.} {\em Soit $R$ un système de racines irréductible de type ADE, muni d'une chambre $C$. On note $\gamma:\mathrm{Q}(R)^{\sharp}\to\mathop{\mathrm{r\acute{e}s}}\mathrm{Q}(R)$ l'application canonique et $\Pi(R)$ le sous-ensemble $\{0\}\coprod \{\varpi_{j};j\in\mathrm{J}\}$ de~$\mathrm{Q}(R)^{\sharp}$.

\medskip
{\em (a)} La restriction de $\gamma$ à $\Pi(R)$ est une bijection.

\medskip
{\em (b)} Pour tout $\xi$ dans $\Pi(R)$, on a $\mathrm{qm}(\gamma(\xi))=\mathrm{q}(\xi)$.}

\bigskip
Les deux points de la proposition ci-dessus permettent la détermination de l'application $\mathrm{qm}:\mathop{\mathrm{r\acute{e}s}}\mathrm{Q}(R)\to\mathbb{Q}\cap[0,\infty[$. Une référence pour le point (a) est \cite[Ch. VI, \S2, Cor. de Prop. 6]{bourbaki}~; le point (b) est implicite dans \cite{venkov}. Pour le confort du lecteur nous donnons ci-après, en petits caractères, une démonstration de ces deux points fondée sur la proposition 3.9 ci-dessous (pour laquelle une référence est \cite[Ch. V, \S3, Th. 2]{bourbaki}) qui est un résultat fondamental (c'est le cas de le dire~!) concernant l'action du groupe de Weyl affine sur $V$.

\medskip
Avant d'énoncer cette proposition, rappelons la définition de l'{\em alcove}, notée $Alc$ dans ce mémoire, associée à la chambre $C$~:
$$
\hspace{24pt}
Alc
\hspace{4pt}:=\hspace{4pt}
\{\hspace{2pt}\xi\hspace{2pt};\hspace{2pt}\xi\in V,\hspace{4pt}\alpha_{i}.\xi>0\hspace{4pt}\text{pour}\hspace{4pt}i=1,2,\ldots,l
\hspace{4pt}\text{et}\hspace{4pt}\widetilde{\alpha}.\xi<1\hspace{2pt}\}
\hspace{24pt};
$$
la fermeture $\overline{Alc}$ de $Alc$ dans $V$ est bien sûr définie par les inégalités larges $\alpha_{i}.\xi\geq 0$ et $\widetilde{\alpha}.\xi\leq 1$.

\bigskip
\textbf{Proposition 3.9.} {\em Toute orbite de l'action canonique du produit semi-direct $\mathrm{W}(R)\ltimes\mathrm{Q}(R)$ (le groupe de Weyl affine) sur $V$ rencontre $\overline{Alc}$ en un et un seul point.}

\footnotesize
\bigskip
\textit{Démonstration du point (a) de la proposition 3.8.}

\medskip
Par définition $\Pi(R)$ est contenu dans $\overline{Alc}$ puisque pour tout $j$ dans $J$ (qui peut être vide) les produits scalaires $\alpha_{i}.\xi$, $i=1,2,\ldots,l$, et $\widetilde{\alpha}.\xi$ sont $0$ ou $1$. En fait, il n'est pas difficile de se convaincre de ce que l'on a $\Pi(R)=\mathrm{Q}(R)^{\sharp}\cap\overline{Alc}$~: Soit $\xi$ un élément de $\mathrm{Q}(R)^{\sharp}$~; on écrit $\xi=\sum_{i=1}^{l}(\alpha_{i}.\xi)\hspace{1pt}\varpi_{i}$, si $\xi$ est dans $\overline{C}$ alors les entiers $\alpha_{i}.\xi$ sont positifs ou nuls, si $\xi$ est non nul et que l'on a en outre $\widetilde{\alpha}.\xi\leq 1$ alors $\xi$ est forcément l'un des $\varpi_{j}$ avec $j$ dans $J$.

\medskip
Compte tenu de ce qui précède, le point (a) de la proposition 3.8 est conséquence de la proposition 3.9~:

\smallskip
-- La partie ``unicité'' de cette proposition montre que la restriction de $\gamma$ à $\Pi(R)$ est injective.

\smallskip
-- La partie ``existence'' montre qu'elle est aussi surjective. En effet, soit $\xi$ un élément de $\mathrm{Q}(R)^{\sharp}$ alors il existe $\eta$ dans $\overline{Alc}$, $w$ dans $\mathrm{W}(R)$ et $x$ dans $\mathrm{Q}(R)$, tels que l'on a $\xi=w\eta+x$~; on en déduit $\eta\in\Pi(R)$, puis $\gamma(\xi)=\gamma(w\eta)$ et enfin $\gamma(\xi)=\gamma(\eta)$ car l'action de $\mathrm{W}(R)$ sur $\mathop{\mathrm{r\acute{e}s}}\mathrm{Q}(R)$, induite par celle de $\mathrm{W}(R)$ sur $\mathrm{Q}(R)$, est triviale \cite[Ch. VI, \S1, Prop. 27]{bourbaki}.

\bigskip
\textit{Démonstration du point (b) de la proposition 3.8.}

\medskip
Ce point est conséquence de l'implication $(v)\Rightarrow(i)$ de la proposition ci-dessous.
\vfill\eject

\bigskip
\textbf{Proposition 3.10.} {\em Soit $R\subset V$ un système de racines irréductible de type~ADE. Soit $\xi$ un élément de $V$. Les conditions suivantes sont équivalentes~:
\begin{itemize}
\item [(i)] On a $\mathrm{q}(\xi)\leq\mathrm{q}(\xi+x)$ pour tout $x$ dans $\mathrm{Q}(R)$.
\item [(ii)] On a $\mathrm{q}(\xi)\leq\mathrm{q}(\xi+\alpha)$ pour tout $\alpha$ dans $R$.
\item [(iii)] On a $\alpha.\xi\leq1$ pour tout $\alpha$ dans $R$.
\item [(iv)] On a $\vert\alpha.\xi\vert\leq1$ pour tout $\alpha$ dans $R$.
\item [(v)] Il existe un élément $w$ de $\mathrm{W}(R)$ tel que $w\xi$ appartient à $\overline{Alc}$ (on suppose ici $R$ muni d'une chambre $C$).
\end{itemize}
}

\bigskip
\textit{Démonstration.} Les implications $(i)\Rightarrow(ii)\Rightarrow(iii)\iff(iv)$ sont triviales (pour $(ii)\Rightarrow(iii)\iff(iv)$, observer que si $\alpha$ est une racine alors il en est de même pour $-\alpha$).

\medskip
\textit{Démonstration de $(iv)\iff(v)$.} On pose
$$
\hspace{24pt}
\Psi
\hspace{4pt}:=\hspace{4pt}
\{\hspace{2pt}\xi\hspace{2pt};\hspace{2pt}\xi\in V,\hspace{4pt}\alpha.\xi\leq 1\hspace{4pt}\text{pour tout}\hspace{4pt}\alpha\hspace{4pt}\text{dans}\hspace{4pt}R\hspace{2pt}\}
\hspace{24pt};
$$
il s'agit de montrer que l'on a $\Psi=\bigcup_{w\in\mathrm{W}(R)}w\hspace{2pt}\overline{Alc}$.

\smallskip
Soit $\xi$ un élément de $V$, alors il existe un élément $w$ de $\mathrm{W}(R)$ avec $w^{-1}\xi\in\overline{C}$. Si l'on suppose en outre $\xi\in\Psi$ alors on a $(w\widetilde{\alpha}).\xi=\widetilde{\alpha}.(w\xi)\leq 1$ et donc $w^{-1}\xi\in\overline{Alc}$~; ceci montre l'inclusion $\Psi\subset\bigcup_{w\in\mathrm{W}(R)}w\hspace{2pt}\overline{Alc}$.

\smallskip
Soit $\alpha$ une racine positive (pour la chambre $C$)~; comme $\alpha$ s'écrit $\sum_{i=1}^{l}\nu_{i}\alpha_{i}$ avec $0\leq\nu_{i}\leq n_{i}$, on a $0\leq\alpha.\xi\leq\widetilde{\alpha}.\xi$ pour tout $\xi$ dans $\overline{C}$. L'inclusion $\overline{Alc}\subset\Psi$ en résulte~; on en déduit $\bigcup_{w\in\mathrm{W}(R)}w\hspace{2pt}\overline{Alc}\subset\Psi$ puisque $\Psi$ est stable sous l'action de $\mathrm{W}(R)$. 
\hfill$\square$

\medskip
\textit{Démonstration de $(iii)\Rightarrow(i)$.} Soit $\Gamma$ un réseau arbitraire de l'espace euclidien~$V$. On pose
$$
\Phi_{\Gamma}
\hspace{4pt}:=\hspace{4pt}
\{\hspace{2pt}\xi\hspace{2pt};\hspace{2pt}\xi\in V,\hspace{4pt}\mathrm{q}(\xi)\leq\mathrm{q}(\xi-x)\hspace{4pt}\text{pour tout}\hspace{4pt}x\hspace{4pt}\text{dans}\hspace{4pt}\Gamma\hspace{2pt}\}
\hspace{24pt},
$$
en clair, $\Phi_{\Gamma}$ est le sous-ensembe de $V$ constitué des points dont la distance à l'origine est inférieure ou égale à la distance à tout point de $\Gamma$. L'étude de ces sous-ensembles est bien sûr fort classique ! Par définition $\Phi_{\Gamma}$ est l'intersection des demi-espaces $x.\xi\leq\mathrm{q}(x)$, $x$~parcourant $\Gamma$~; il est facile de se convaincre de ce que $\Phi_{\Gamma}$ est compact et intersection d'une sous-famille finie de ces demi-espaces. Le sous-ensemble $\Phi_{\Gamma}$ vérifie les trois propriétés suivantes~:
\begin{itemize}
\item[(1)] $\Phi_{\Gamma}$ est la fermeture dans $V$ de son intérieur~;
\item[(2)] les translatés $\Phi_{\Gamma}+x$, $x$ parcourant $\Gamma$, recouvrent $V$~; 
\item[(3)] les translatés $\mathring{\Phi}_{\Gamma}+x$, $x$ parcourant $\Gamma$, sont deux à deux disjoints (la notation $\mathring{\Phi}_{\Gamma}$ désigne l'intérieur de $\Phi_{\Gamma}$).
\end{itemize}

\medskip
\textbf{Lemme 3.11.} {\em Soit $\Phi$ (resp. $\Phi'$) un sous-ensemble de $V$ vérifiant les propriétés (1), (2) et (3) (resp. (1) et (3)). Si l'on a $\Phi\subset\Phi'$ alors on a $\Phi=\Phi'$.}

\medskip
\textit{Démonstration.} Supposons $\Phi'\not\subset\Phi$. Dans ce cas on a aussi $\mathring{\Phi}'\not\subset\Phi$ puisque $\Phi$ est fermé et que $\Phi'$ est la fermeture de son intérieur. Soit $\xi$ un élément de $V$ avec $\xi\in\mathring{\Phi}'$ et $\xi\not\in\Phi$. Puisque $\Phi$ vérifie (2) il existe $x$ dans $\Gamma$ tel que l'on $\xi+x\in\Phi$ et \textit{a fortiori} $\xi+x\in\Phi'$. On a donc $\mathring{\Phi}'\cap(\Phi'+x)\not=\emptyset$~; par un argument de topologie générale analogue au précédent on a aussi $\mathring{\Phi}'\cap(\mathring{\Phi}'+x)\not=\emptyset$. Puisque $\Phi$ vérifie (3), on a $x=0$. Contradiction.
\hfill$\square$

\medskip
Le lemme ci-dessus fournit une démonstration de l'implication $(iii)\Rightarrow(i)$. On pose $\Phi:=\Phi_{\mathrm{Q}(R)}$~; $\Phi$ est donc le sous-ensemble de $V$ constitué des $\xi$ qui vérifient la condition (i). Il s'agit de montrer que l'on a $\Phi=\Psi$. Il est clair que $\Psi$ vérifie la propriété (1)~; pour conclure il suffit de montrer qu'il vérifie aussi la propriété (3). Soit $\xi$ un élément de $\mathrm{Q}(R)$ avec $\mathring{\Psi}\cap(\mathring{\Psi}+x)\not=\emptyset$~; soit $\xi$ un élément de cette intersection. Les inégalités $-1<\alpha.\xi<1$ et $-1<\alpha.(\xi+x)<1$ impliquent $-2<\alpha.x<2$~; comme $\alpha.x$ est entier, on a aussi $-1\leq\alpha.x\leq 1$, c'est-à-dire $x\in\Psi$. La proposition ci-dessous dit que $x$ est nul.

\medskip
\textbf{Proposition 3.12.} {\em On a $\Psi\cap\mathrm{Q}(R)=\{0\}$.}

\medskip
\textit{Démonstration.} Compte tenu de l'égalité $\Psi=\bigcup_{w\in\mathrm{W}(R)}w\hspace{2pt}\overline{Alc}$ il sufit de montrer que l'on a $\overline{Alc}\cap\mathrm{Q}(R)=\{0\}$. Cette égalité est à nouveau impliquée par la proposition 3.9.
\hfill$\square\square$

\normalsize
\bigskip
\textit{Remarque terminologique.} Les éléments $\xi$ de $\mathrm{Q}(R)^{\sharp}$ vérifiant $\mathrm{q}(\xi)\leq\mathrm{q}(\xi+x)$ pour tout $x$ dans $Q$ sont appelés {\em poids minuscules}~; cependant on réserve souvent l'appellation ``poids minuscules'' aux poids fondamentaux $\varpi_{j}$, $j\in J$, considérés plus haut (voir par exemple \cite{stembridge}).

\bigskip
3) Soit $R$ un système de racines de type ADE. L'action de $\mathrm{A}(R)$ sur $\mathrm{Q}(R)$ induit une action de $\mathrm{A}(R)$ sur le $\mathrm{qe}$-module de Venkov $\mathop{\mathrm{r\acute{e}s}}\mathrm{Q}(R)$. La restriction à $\mathrm{W}(R)$ de l'action de $\mathrm{A}(R)$ sur $\mathop{\mathrm{r\acute{e}s}}\mathrm{Q}(R)$ est triviale (voir plus haut)~; comme $\mathrm{W}(R)$ est distingué dans $\mathrm{A}(R)$ \cite[Ch. VI, \S1, Prop. 16]{bourbaki}, on dispose d'une action canonique (en fait fidèle \cite[Ch. VI, \S4, Exc. 7]{bourbaki}) du groupe quotient $\mathrm{G}(R):=\mathrm{A}(R)/\mathrm{W}(R)$ sur $\mathop{\mathrm{r\acute{e}s}}\mathrm{Q}(R)$.

\bigskip
On revient maintenant à la classification des réseaux unimodulaires pairs de dimension $24$ dont l'ensemble des racines est non vide. Pour achever cette classification, Venkov vérifie au cas par cas l'énoncé (miraculeux~!) suivant~:

\bigskip
\textbf{Proposition 3.13.} {\em Soit $R$ un système de racines de type ADE, de rang $24$ et équicoxeter.

\medskip
{\em (a)} Le $\mathrm{qe}$-module de Venkov $\mathop{\mathrm{r\acute{e}s}}\mathrm{Q}(R)$ possède un lagrangien $I$ avec $\mathrm{qm}(\xi)>1$ pour tout~$\xi$ dans $I-\{0\}$.

\medskip
{\em (b)} Un tel lagrangien est unique à l'action de $\mathrm{G}(R)$ près.}

\bigskip
\textbf{Corollaire 3.14.} {\em L'application $L\mapsto\mathrm{R}(L)$ induit une bijection, de l'ensemble des classes d'isomorphisme de réseaux unimodulaires pairs de dimension $24$ avec $\mathrm{R}(L)\not=\emptyset$, sur l'ensemble des classes d'isomorphisme de systèmes de racines de type ADE de rang $24$ et équicoxeter.}

\bigskip
\textbf{Scholie 3.15.} {\em Soit $L$ un réseau unimodulaire pair de dimension $24$ avec $\mathrm{R}(L)\not=\emptyset$.

\medskip
{\em (a)} L'action du groupe de Weyl $\mathrm{W}(\mathrm{R}(L))$  sur $\mathrm{R}(L)$ se prolonge en une action (orthogonale) sur $L$ si bien que $\mathrm{W}(\mathrm{R}(L))$ s'identifie à un sous-groupe du groupe orthogonal $\mathrm{O}(L)$.

\medskip
{\em (b)} Le groupe $\mathrm{W}(\mathrm{R}(L))$ est distingué dans $\mathrm{O}(L)$ et le quotient $\mathrm{O}(L)/\mathrm{W}(\mathrm{R}(L))$ est canoniquement isomorphe au sous-groupe de $\mathrm{G}(R)$ qui stabilise le lagrangien $L/\mathrm{Q}(\mathrm{R}(L))$ du $\mathrm{qe}$-module $\mathop{\mathrm{r\acute{e}s}}\mathrm{Q}(R)$ (lagrangien qui est l'un de ceux considérés dans la proposition 3.13).}

\bigskip
(Le point (a) du scholie ci-dessus résulte de ce que l'action canonique de $\mathrm{W}(\mathrm{R}(L))$ sur $\mathop{\mathrm{r\acute{e}s}}\mathrm{Q}(R)$ est triviale.)

\bigskip
\textsc{Exemples}

\bigskip
1) $\mathrm{R}(L)\cong\mathbf{D}_{24}$, $\mathrm{R}(L)\cong\mathbf{D}_{16}\coprod\mathbf{E}_{8}$, $\mathrm{R}(L)\cong\mathbf{E}_{8}\coprod\mathbf{E}_{8}\coprod\mathbf{E}_{8}$

\medskip
Il résulte de ce que l'on a déjà vu que l'on a dans ces trois cas respectivement $L\cong\mathrm{D}_{24}^{+}$, $L\cong\mathrm{D}_{16}^{+}\oplus\mathrm{E}_{8}$, $L\cong\mathrm{E}_{8}\oplus\mathrm{E}_{8}\oplus\mathrm{E}_{8}$.

\bigskip
2) $\mathrm{R}(L)\cong\mathbf{A}_{24}$

\medskip
Soit $n\geq 1$ un entier~; on munit le $\mathbb{R}$-espace vectoriel $\mathbb{R}^{n+1}$ de sa structure euclidienne canonique et on note $(\varepsilon_{1},\varepsilon_{2},\ldots,\varepsilon_{n+1})$ sa base canonique. On note $\mathrm{A}_{n}$ le sous-module de $\mathrm{I}_{n+1}:=\mathbb{Z}^{n+1}\subset\mathbb{R}^{n+1}$ constitué des $(n+1)$-uples $(x_{1},x_{2},\ldots,x_{n+1})$ avec $x_{1}+x_{2}+\ldots+x_{n+1}=0$~; $\mathrm{A}_{n}$ est un réseau entier pair (dans l'hyperplan de $\mathbb{R}^{n+1}$ constitué des $(n+1)$-uples $(\xi_{1},\xi_{2},\ldots,\xi_{n+1})$ avec $\xi_{1}+\xi_{2}+\ldots+\xi_{n+1}=0$). Le système de racines $\mathbf{A}_{n}$ est défini par l'égalité $\mathbf{A}_{n}:=\mathrm{R}(\mathrm{A}_{n})$ (si bien que l'on a aussi l'égalité $\mathrm{A}_{n}=\mathrm{Q}(\mathbf{A}_{n})$).

\medskip
Le groupe abélien sous-jacent au $\mathrm{qe}$-module de Venkov $\mathop{\mathrm{r\acute{e}s}}\mathrm{A}_{n}$ est cyclique d'ordre $n+1$ engendré par la classe de la projection orthogonale de $\varepsilon_{1}$ sur l'hyperplan $\sum_{i=1}^{n+1}\xi_{i}=0$, disons $\varpi$. La forme quadratique d'enlacement $\mathrm{q}:\mathbb{Z}/(n+1)\to\mathbb{Q}/\mathbb{Z}$, définie par transport de structure, est donnée par
$$
\hspace{24pt}
\mathrm{q}(\bar{k})=k^{2}\hspace{1pt}\mathrm{q}(\varpi)=\frac{n\hspace{1pt}k^{2}}{2(n+1)}
\hspace{24pt}.
$$
L'application $\mathrm{qm}:\mathbb{Z}/(n+1)\to\mathbb{Q}\cap[0,\infty[$ est donnée quant à elle par
$$
\mathrm{qm}(\bar{k})=\frac{k\hspace{1pt}(n+1-k)}{2(n+1)}
\hspace{24pt}\text{pour}\hspace{6pt}0\leq k\leq n
$$
(dans le cas de $\mathbf{A}_{n}$, tous les poids fondamentaux sont minuscules).

\medskip
Il est clair que $\mathop{\mathrm{r\acute{e}s}}\mathrm{A}_{n}$ possède un lagrangien, au sens bilinéaire, si et seulement si l'entier $n+1$ est un carré, $n+1=r^{2}$, à savoir le sous-module engendré par $r\hspace{1pt}\varpi$, et que dans ce cas ce lagrangien est unique.

\medskip
On suppose $n+1=r^{2}$ et on note $I$ le sous-module engendré par $r\hspace{1pt}\varpi$~; on constate que $I$ est un lagrangien au sens quadratique si et seulement si $n$ est pair, c'est-à-dire $r$ impair, et que dans ce cas $I$ vérifie la condition (ii) de la proposition 3.7 si et seulement on a l'inégalité $r\geq 5$.
\vfill\eject

\medskip
On note $\mathrm{A}_{24}^{+}$ le réseau de $\mathbb{Q}\otimes_{\mathbb{Z}}\mathrm{A}_{24}$ engendré par $\mathrm{A}_{24}$ et $5\hspace{1pt}\varpi$. La discussion précédente montre que $\mathrm{A}_{24}^{+}$ est un réseau unimodulaire pair de dimension $24$ avec $\mathrm{R}(\mathrm{A}_{24}^{+})\cong\mathbf{A}_{24}$ et que cette propriété caractérise, à isomorphisme près, $\mathrm{A}_{24}^{+}$ parmi les réseaux unimodulaires pairs de dimension $24$.

\bigskip
3) $\mathrm{R}(L)\cong\mathbf{D}_{12}\coprod\mathbf{D}_{12}$

\medskip
Le $\mathrm{qe}$-module de Venkov $\mathop{\mathrm{r\acute{e}s}}\mathrm{D}_{12}$ est isomorphe à $\mathbb{Z}/2\oplus\mathbb{Z}/2$ muni des applications $\mathrm{q}$ et $\mathrm{qm}$ définies respectivement par~:

\smallskip
$\mathrm{q}(\bar{0},\bar{0})=0$, $\mathrm{q}(\bar{1},\bar{0})=\frac{1}{2}$, $\mathrm{q}(\bar{0},\bar{1})=\frac{1}{2}$, $\mathrm{q}(\bar{1},\bar{1})=\frac{1}{2}$~;

\smallskip
$\mathrm{qm}(\bar{0},\bar{0})=0$, $\mathrm{qm}(\bar{1},\bar{0})=\frac{3}{2}$, $\mathrm{qm}(\bar{0},\bar{1})=\frac{3}{2}$ et $\mathrm{qm}(\bar{1},\bar{1})=\frac{1}{2}$.

\smallskip
Le $\mathrm{qe}$-module de Venkov $\mathop{\mathrm{r\acute{e}s}}(\mathrm{D}_{12}\oplus\mathrm{D}_{12})$ est isomorphe à $\mathop{\mathrm{r\acute{e}s}}\mathrm{D}_{12}\oplus\mathop{\mathrm{r\acute{e}s}}\mathrm{D}_{12}$. On constate que les lagrangiens de $\mathop{\mathrm{r\acute{e}s}}(\mathrm{D}_{12}\oplus\mathrm{D}_{12})$ sont les graphes des permutations $\phi$ de $\mathop{\mathrm{r\acute{e}s}}\mathrm{D}_{12}$ qui préservent $0$ (une telle permutation est linéaire et préserve la forme quadratique d'enlacement). Le graphe de $\phi$ vérifie la condition (ii) de la proposition 3.7 si et seulement si l'on a $\phi(\mathrm{qm}^{-1}(\frac{1}{2}))\not=\mathrm{qm}^{-1}(\frac{1}{2})$~;  de tels $\phi$ sont au nombre de $4$ et l'on constate que le groupe $\mathrm{G}(\mathbf{D}_{12}\coprod\mathbf{D}_{12})$ (qui est isomorphe à $\mathfrak{S}_{2}\wr\mathrm{G}(\mathbf{D}_{12})=\mathfrak{S}_{2}\wr\mathfrak{S}_{2}$)  agit bien transitivement sur l'ensemble des $4$ graphes correspondants.

\bigskip
4) $\mathrm{R}(L)\cong\mathbf{A}_{17}\coprod\mathbf{E}_{7}$

\medskip
Le $\mathrm{qe}$-module de Venkov $\mathop{\mathrm{r\acute{e}s}}(\mathrm{A}_{17}\oplus\mathrm{E}_{7})$ est isomorphe à $\mathop{\mathrm{r\acute{e}s}}\mathrm{A}_{17}\oplus\mathop{\mathrm{r\acute{e}s}}\mathrm{E}_{7}$. La structure du $\mathrm{qe}$-module de Venkov $\mathop{\mathrm{r\acute{e}s}}\mathrm{A}_{n}$, pour tout $n$, a déja été explicitée dans le deuxième exemple ci-dessus. Le $\mathrm{qe}$-module $\mathop{\mathrm{r\acute{e}s}}\mathrm{E}_{7}$ est quant à lui  isomorphe à $\mathbb{Z}/2$ muni de la forme quadratique d'enlacement définie par $\mathrm{q}(\bar{1})=-\frac{1}{4}$ (voir le point (d) de B.2.2)~; on vérifie que l'on a $\mathrm{qm}(\bar{1})=\frac{3}{4}$. On constate que l'unique lagrangien du $\mathrm{qe}$-module $\mathop{\mathrm{r\acute{e}s}}(\mathrm{A}_{17}\oplus\mathrm{E}_{7})$ satisfait bien la condition (ii) de la proposition 3.7.

\bigskip
5) $\mathrm{R}(L)\cong\mathbf{D}_{10}\coprod\mathbf{E}_{7}\coprod\mathbf{E}_{7}$

\medskip
La structure du $\mathrm{qe}$-module de Venkov $\mathop{\mathrm{r\acute{e}s}}(\mathrm{D}_{10}\oplus\mathrm{E}_{7}\oplus\mathrm{E}_{7})$ est déterminée par celle de $\mathop{\mathrm{r\acute{e}s}}\mathrm{E}_{7}$ que nous venons d'expliciter et celle de $\mathop{\mathrm{r\acute{e}s}}\mathrm{D}_{10}$. Celle-ci est la suivante~: $\mathop{\mathrm{r\acute{e}s}}\mathrm{D}_{10}$ est isomorphe à $\mathbb{Z}/2\oplus\mathbb{Z}/2$ muni de l'application $\mathrm{qm}$ définie par $\mathrm{qm}(\bar{0},\bar{0})=0$, $\mathrm{qm}(\bar{1},\bar{0})=\frac{5}{4}$, $\mathrm{qm}(\bar{0},\bar{1})=\frac{5}{4}$ et $\mathrm{qm}(\bar{1},\bar{1})=\frac{1}{2}$.

\medskip
Le $\mathrm{qe}$-module $\mathop{\mathrm{r\acute{e}s}}(\mathrm{D}_{10}\oplus\mathrm{E}_{7}\oplus\mathrm{E}_{7})$ possède deux lagrangiens, disons $I_{1}$ et $I_{2}$, à savoir les graphes des deux isomorphismes de $\mathrm{qe}$-modules de $\mathop{\mathrm{r\acute{e}s}}\mathrm{D}_{10}$ sur $\langle -1\rangle\otimes\mathop{\mathrm{r\acute{e}s}}(\mathrm{E}_{7}\oplus\mathrm{E}_{7})$. On constate que $I_{1}$ et $I_{2}$ satisfont tous deux la condition (ii) de la proposition 3.7.

\medskip
Le groupe $G:=\mathrm{G}(\mathbf{D}_{10}\coprod\mathbf{E}_{7}\coprod\mathbf{E}_{7})$ s'identife à $\mathbb{Z}/2\times\mathbb{Z}/2$ car il est isomorphe au produit des groupes $\mathrm{G}(\mathbf{D}_{10})$ et $\mathrm{G}(\mathbf{E}_{7}\coprod\mathbf{E}_{7})$ qui sont tous deux  cycliques d'ordre $2$. On constate que $G$ agit bien transitivement sur l'ensemble $\{I_{1},I_{2}\}$ ; on observe également que le $\mathbb{Z}/2$ diagonal agit trivialement (cette observation sera utilisée à la fin du paragraphe 2 de l'appendice B).

\bigskip
6) $\mathrm{R}(L)\cong 24\hspace{1pt}\mathbf{A}_{1}$ (cette notation désigne la somme directe de $24$ copies du système de racines $\mathbf{A}_{1}$)

\medskip
Le $\mathrm{qe}$-module de Venkov $\mathop{\mathrm{r\acute{e}s}}\mathrm{A}_{1}$ est isomorphe à $\mathbb{Z}/2$ muni de l'application $\mathrm{qm}$ définie par $\mathrm{qm}(x)=\frac{\lambda(x)}{4}$, $\lambda:\mathbb{Z}/2\to\mathbb{N}$ désignant l'application définie par $\lambda(\bar{0})=0$ et $\lambda(\bar{1})=1$~; il en résulte que le $\mathrm{qe}$-module de Venkov $\mathop{\mathrm{r\acute{e}s}}(\mathrm{A}_{1}^{\oplus\hspace{1pt}24})$ est isomorphe à $(\mathbb{Z}/2)^{24}$ muni de l'application $\mathrm{qm}$ définie par
$$
\hspace{24pt}
\mathrm{qm}(x_{1},x_{2},\ldots,x_{24})
\hspace{4pt}=\hspace{4pt}
\frac{1}{4}\hspace{2pt}\sum_{i=1}^{24}\hspace{1pt}\lambda(x_{i})
\hspace{4pt}=:\hspace{4pt}
\frac{1}{4}\hspace{2pt}\mathrm{wt}(x_{1},x_{2},\ldots,x_{24})
\hspace{24pt}.
$$
Le groupe $\mathrm{G}(24\hspace{1pt}\mathbf{A}_{1})$ est isomorphe à $\mathfrak{S}_{24}\wr\mathrm{G}(\mathbf{A}_{1})=\mathfrak{S}_{24}$ et son action sur $\mathop{\mathrm{r\acute{e}s}}(\mathrm{A}_{1}^{\oplus\hspace{1pt}24})$ s'identifie à l'action évidente de $\mathfrak{S}_{24}$ sur $(\mathbb{Z}/2)^{24}$.

\medskip
Un lagrangien de $\mathop{\mathrm{r\acute{e}s}}(\mathrm{A}_{1}^{\oplus\hspace{1pt}24})$ s'identifie à un sous-espace vectoriel $I$ de $(\mathbb{Z}/2)^{24}$ avec $\dim I=12$ et $\mathrm{wt}(x)\equiv 0\bmod{4}$ pour tout $x$ dans $I$~; un tel $I$ est ce qu'on appelle un {\em code binaire auto-dual pair} de longueur $24$. A isomorphisme près (en d'autres termes, modulo l'action de $\mathfrak{S}_{24}$), il existe $9$ codes binaires auto-duaux pairs de longueur $24$ \cite{PS}.

\medskip
Le lagrangien $I$ vérifie en outre la condition (ii) de la proposition 3.7 si et seulement si l'on a $\mathrm{wt}(x)\geq 8$ pour tout $x$ dans $I-\{0\}$. Là encore, il existe un code binaire auto-dual pair de longueur $24$ vérifiant cette propriété, et à isomorphisme près un seul, à savoir le code de Golay \cite{Pl}.

\bigskip
L'énoncé suivant permet de compléter l'énoncé 3.14~:

\bigskip
\textbf{Théorème 3.16.} {\em Tout réseau unimodulaire pair dont l'ensemble des racines est vide est isomorphe au réseau de Leech.}

\bigskip
On obtient~:

\bigskip
\textbf{Théorème 3.17.} {\em L'application $L\mapsto\mathrm{R}(L)$ induit une bijection, de l'ensemble des classes d'isomorphisme de réseaux unimodulaires pairs de dimension $24$, sur l'ensemble à $24$ éléments constitué, des classes d'isomorphisme de systèmes de racines de type ADE de rang $24$ et équicoxeter, et de l'ensemble vide.}

\bigskip
Le théorème 3.16 est dû à Conway (voir \cite{Co}). Venkov en propose dans \cite{venkov} une démonstration qui utilise la notion de voisin de Kneser dont traite le prochain chapitre.


\chapter{Voisins à la Kneser}

\section{Variations sur la notion de voisins à la Kneser}

La notion de réseaux unimodulaires $2$-voisins a été introduite par Martin Kneser dans \cite{kneser16}. Dans l'article en question, Kneser décrit, à l'aide de cette notion, un algorithme pour classifier les réseaux unimodulaires (l'exhaustivité est essentiellement due au théorème II.2.8) et met en oeuvre cet algorithme pour expliciter la liste des classes d'isomorphisme de réseaux unimodulaires de dimension inférieure ou égale à $15$ (voir aussi \cite[\S 106F]{O'M}). 

\bigskip
Commençons nos variations par deux observations très générales (nous faisons intervenir ci-dessous un anneau de Dedekind $R$ arbitraire, mais les deux applications que nous avons en tête sont $R=\mathbb{Z}$ et $R=\mathbb{Z}_{p}$)~:

\bigskip
\textbf{Proposition 1.1.} {\em Soient $R$ un anneau de Dedekind et $K$ son corps des fractions. Soit $V$ un $\mathrm{q}$-espace vectoriel sur $K$,  de dimension finie~; soient $L_{1}$ et $L_{2}$ deux réseaux entiers auto-duaux de $V$ (les réseaux $L_{1}$ et $L_{2}$ sont donc en particulier deux $\mathrm{q}$-modules et le réseau $L_{1}\cap L_{2}$ un $\widetilde{ \mathrm{q}}$-module, sur $R$).

\medskip
On pose $I_{1}=L_{1}/(L_{1}\cap L_{2})$ et $I_{2}=L_{2}/(L_{1}\cap L_{2})$.

\medskip
{\em (a)} Les sous-modules $I_{1}$ et $I_{2}$ sont deux lagrangiens transverses du $\mathrm{qe}$-module $\mathop{\mathrm{r\acute{e}s}}(L_{1}\cap L_{2})$. La forme d'enlacement de ce $\mathrm{qe}$-module induit un isomorphisme de $I_{2}$ sur $I_{1}^{\vee}$, disons $\iota$, et l'homomorphisme composé
$$
\begin{CD}
\mathrm{H}(I_{1})=I_{1}\oplus I_{1}^{\vee}
@>\mathrm{id}\oplus\iota^{-1}>>
I_{1}\oplus I_{2}@>>>
\mathop{\mathrm{r\acute{e}s}}(L_{1}\cap L_{2})
\end{CD}
$$
est un isomorphisme de $\mathrm{qe}$-modules sur $R$ (on rappelle que la notation $\mathrm{H}(I_{1})$ désigne le $\mathrm{qe}$-module hyperbolique sur le $R$-module de torsion de type fini $I_{1}$).

\medskip
{\em (b)} Soit $r$ le nombre minimum de générateurs du $R$-module $I_{1}$~; on a l'inégalité
$$
\hspace{24pt}
2\hspace{2pt}r
\hspace{4pt}\leq\hspace{4pt}
\dim_{K}V
\hspace{24pt}.
$$}

\medskip
\textit{Démonstration.} La première partie du point (a) est évidente~; la seconde est la reprise du point (a) de II.1.3. Passons au point (b). Soit $I$ un $R$-module de torsion de type fini, on note $\mathrm{r}(I)$ le nombre minimum de générateurs de cet $R$-module~; $\mathrm{r}(I)$ peut être vu, par exemple, comme le plus grand des entiers $k$ tels que la puissance extérieure $\Lambda^{k}I$ est non nulle. Comme le dual $I^{\vee}$ est isomorphe (non canoniquement) à $I$ (pour s'en convaincre, observer que $R$ est ``localement principal'') on a $\mathrm{r}(I)=\mathrm{r}(I^{\vee})$~; on a donc, d'après le point~(a), l'égalité $\mathrm{r}(\mathrm{r\acute{e}s}(L_{1}\cap L_{2}))=2\hspace{1pt}\mathrm{r}(I_{1})$. Comme $\mathrm{r\acute{e}s}(L_{1}\cap L_{2})$ est un quotient du réseau ${(L_{1}\cap L_{2})}^{\sharp}$ on a $\mathrm{r}(\mathrm{r\acute{e}s}(L_{1}\cap L_{2}))\leq\dim_{K}V$.
\hfill$\square$

\bigskip
On spécialise au cas $R=\mathbb{Z}$.

\bigskip
\textbf{Scholie-Définition 1.2.} {\em Soit $V$ un $\mathrm{q}$-espace vectoriel sur $\mathbb{Q}$~; soient $L_{1}$ et~$L_{2}$ deux réseaux entiers auto-duaux de $V$ (les réseaux $L_{1}$ et $L_{2}$ sont donc en particulier deux $\mathrm{q}$-modules et le réseau $L_{1}\cap L_{2}$ un $\widetilde{ \mathrm{q}}$-module, sur $\mathbb{Z}$).

\medskip
Soit $A$ un groupe abélien fini~; les conditions suivantes sont équivalentes~:
\begin{itemize}
\item [(i)] le quotient $L_{1}/(L_{1}\cap L_{2})$ est isomorphe à $A$~;
\item [(ii)] le quotient $L_{2}/(L_{1}\cap L_{2})$ est isomorphe à $A$.
\end{itemize}

\medskip
Si ces conditions sont vérifiées alors on dit que $L_{1}$ et $L_{2}$ sont {\em $A$-voisins} (ou que $L_{2}$ est un {\em $A$-voisin} de $L_{1}$).}

\vspace{0,75cm}
\textsc{$d$-voisins, point de vue asymétrique}

\bigskip
La spécificité de la notion de $A$-voisins pour $A$ cyclique tient à l'énoncé ci-dessous, qui peut être vu comme un corollaire du point (b) de II.1.3.~:

\bigskip
\textbf{Proposition 1.3.} {\em Soit $A$ un groupe cyclique fini~; alors $A^{\vee}$ est l'unique lagrangien du $\mathrm{qe}$-module $\mathrm{H}(A)$ transverse (au sens de II.1.3) au lagrangien~$A$.}

\bigskip
On fixe maintenant un $\mathrm{q}$-module $L$ sur $\mathbb{Z}$, un entier $d\geq 2$, et on analyse l'ensemble des $\mathbb{Z}/d$-voisins de $L$ dans $\mathbb{Q}\otimes_{\mathbb{Z}}L$ (qui est un $\mathrm{q}$-espace vectoriel sur $\mathbb{Q}$). On allège la terminologie $\mathbb{Z}/d$-voisin en {\em $d$-voisin}.

\medskip
Dans ce contexte, un $d$-voisin de $L$ est un réseau entier $L'$ dans $\mathbb{Q}\otimes_{\mathbb{Z}}L$ avec $L'^{\sharp}=L'$ et $L/(L\cap L')\simeq\mathbb{Z}/d$. On pose $M=L\cap L'$. D'après ce qui précède~:

\smallskip
-- Le réseau $dL'$ est contenu dans $M$.

\smallskip
-- L'image de l'homomorphisme composé $dL'\subset M\subset L\to L/dL$ est une droite isotrope, disons $c$, de $L/dL$ muni de sa structure de $\mathrm{q}$-module sur $\mathbb{Z}/d$. Précisons ce que l'on entend ici par {\em droite isotrope} de $L/dL$~: $c$ est un sous-module libre de dimension $1$ du $\mathbb{Z}/d$-module $L/dL$ (nécessairement facteur direct car $\mathbb{Z}/d$ est un $\mathbb{Z}/d$-module injectif) tel que la restriction de la forme quadratique $\mathrm{q}:L/dL\to\mathbb{Z}/d$ à $c$ est nulle.

\smallskip
-- Le réseau $M$ est l'image inverse par  l'homomorphisme $L\to L/dL$ de $c^{\hspace{1pt}\perp}$, $c^{\perp}$ désignant le sous-module du $\mathbb{Z}/d$-module $L/dL$ orthogonal de la droite $c$.

\smallskip
-- Le réseau $L'$ est l'image inverse par l'homomorphisme $M^{\sharp}\to\mathop{\mathrm{r\acute{e}s}}M$ de l'unique lagrangien transverse au lagrangien $L/M$.

\smallskip
Ce qui précède montre que l'application $L'\mapsto c$ est injective, la proposition ci-dessous montre qu'elle est aussi surjective.

\bigskip
\textbf{Proposition 1.4.} {\em Soit $c$ une droite isotrope de $L/dL$~; soit $M$ le sous-module de $L$ image inverse de $c^{\hspace{1pt}\perp}$ par  l'homomorphisme $L\to L/dL$ . Alors le $\mathrm{qe}$-module $\mathop{\mathrm{r\acute{e}s}}M$ est isomorphe à $\mathrm{H}(\mathbb{Z}/d)$, le quotient $L/M$ en est un lagrangien et l'image inverse par l'homomorphisme $M^{\sharp}\to\mathop{\mathrm{r\acute{e}s}}M$ de l'unique lagrangien qui lui est transverse est un $d$-voisin $L'$ de $L$ avec $L\cap L'=M$.}

\bigskip
Avant de vérifier cet énoncé nous introduisons une terminologie et une notation que nous utiliserons dans la suite du chapitre.

\medskip
Soit $L$ un $\mathbb{Z}$-module libre de dimension finie~; on dit qu'un élément $u$ de~$L$ est indivisible ou {\em primitif} s'il est non nul et si le quotient $L/\mathbb{Z}u$ est sans torsion. Tout élément non nul $u$ de $L$ s'écrit de façon unique $\mathrm{c}(u)\hspace{1pt}v$ avec $v$ primitif et $\mathrm{c}(u)$ dans $\mathbb{N}-\{0\}$. Soit $d\geq 2$ un entier, nous dirons qu'un élément $u$ de~$L$ est {\em $d$-primitif}, s'il est non nul et si $d$ est premier avec $\mathrm{c}(u)$. En d'autres termes, $u$ est $d$-primitif si le sous-module de $L/dL$ engendré par la classe de $u$ est un $\mathbb{Z}/d$-module libre de dimension $1$. Comme nous l'avons déjà observé, un tel sous-module est nécessairement facteur direct si bien qu'un élément $d$-primitif $u$ de $L$ définit un élément de $\mathrm{P}_{L}(\mathbb{Z}/d)$, la notation $\mathrm{P}_{L}$ désignant le schéma dont les $R$-points, pour tout anneau commutatif $R$, sont les facteurs directs de rang~$1$ du $R$-module libre $R\otimes_{\mathbb{Z}}L$. Cet élément sera noté $[u]$.

\bigskip
\textit{Démonstration de 1.4.} Il est clair qu'il suffit de vérifier que $\mathop{\mathrm{r\acute{e}s}}M$ est isomorphe à $\mathrm{H}(\mathbb{Z}/d)$. Soit $u$ un élément $d$-primitif de $L$ dont la classe modulo $d$ engendre la droite $c$~; puisque cette droite est isotrope on a $\mathrm{q}(u)\equiv 0\bmod{d}$. Soit  $v$ un élément de $L$ avec $u.v\equiv 1\bmod{d}$. On constate que $v$ et $\frac{u}{d}$ appartiennent à $M^{\sharp}$ et que $\mathop{\mathrm{r\acute{e}s}}\mathrm{M}$ est un $\mathbb{Z}/d$-module libre de dimension $2$ dont une base est constituée des classes de ces deux éléments. On pose $w=\frac{u}{d}-\frac{\mathrm{q}(u)}{d}\hspace{1pt}v$~; on constate également que l'on a dans $\mathbb{Q}/\mathbb{Z}$ les égalités $\mathrm{q}(v)=0$, $\mathrm{q}(w)=0$ et $v.w=\frac{1}{d}$.
\hfill$\square$

\bigskip
On voit donc au bout du compte que le réseau $L'$ est entièrement déterminé par la donnée de $c$. On met en place ci-après une notation qui souligne cette dépendance.

\medskip
On note $\mathrm{C}_{L}(\mathbb{Z}/d)$ l'ensemble des droites isotropes de $L/dL$. La justification de cette notation est la suivante~: on note $\mathrm{C}_{L}\subset\mathrm{P}_{L}$ la quadrique définie par la forme quadratique dont $L$ est muni. On observera incidemment que la non dégénérescence de cette forme quadratique fait que $\mathrm{C}_{L}$ est lisse sur $\mathbb{Z}$ (version projective du critère jacobien de lissité). Il est clair que les $\mathbb{Z}/d$-points de $\mathrm{C}_{L}$ sont les droites isotropes de $L/dL$. Soit $c$ un élément de $\mathrm{C}_{L}(\mathbb{Z}/d)$, les réseaux $M$ et $L'$ introduits dans la proposition 1.4 seront respectivement notés $\mathrm{M}_{d}(L;c)$ et $\mathrm{vois}_{d}(L;c)$. On note enfin $\mathrm{Vois}_{d}(L)$ l'ensemble des $d$-voisins de $L$ dans $\mathbb{Q}\otimes_{\mathbb{Z}}L$. On a tout fait pour avoir l'énoncé suivant~:

\bigskip
\textbf{Proposition 1.5.} {\em L'application
$$
\mathrm{C}_{L}(\mathbb{Z}/d)\to\mathrm{Vois}_{d}(L)
\hspace{24pt},\hspace{24pt}
c\mapsto\mathrm{vois}_{d}(L;c)
$$
est une bijection.}

\bigskip
Soit $u$ un élément $d$-primitif de $L$ avec $\mathrm{q}(u)\equiv 0\bmod{d}$~; les réseaux $\mathrm{M}_{d}(L;[u])$ et $\mathrm{vois}_{d}(L;[u])$ seront aussi notés $\mathrm{M}_{d}(L;u)$ et $\mathrm{vois}_{d}(L;u)$. Pour mémoire, nous explicitons ci-après l'algorithme $u\leadsto\mathrm{vois}_{d}(L;u)$ fourni par la démonstration de 1.4~:

\medskip
Soit $v$ un élément de $L$ avec $u.v\equiv 1\bmod{d}$ alors $\mathrm{vois}_{d}(L;u)$ est le réseau de $\mathbb{Q}\otimes_{\mathbb{Z}}L$ engendré par
$$
\hspace{24pt}
\mathrm{M}_{d}(L;u):=\{\hspace{2pt}x\hspace{2pt};\hspace{2pt}x\in L\hspace{2pt},\hspace{2pt}u.x\equiv 0\bmod{d}\hspace{2pt}\}
\hspace{24pt}\text{et}\hspace{24pt}\frac{u-\mathrm{q}(u)v}{d}
\hspace{24pt}.
$$
Posons $\widetilde{u}=u-\mathrm{q}(u)v$~; on observera que l'on a $\widetilde{u}\equiv u\bmod{d}$, en d'autres termes $[\widetilde{u}]=[u]$ dans $\mathrm{C}_{L}(\mathbb{Z}/d)$, et $\mathrm{q}(\widetilde{u})\equiv 0\bmod{d^{2}}$. Cette observation conduit de façon évidente à une présentation alternative de l'algorithme~: $u$ étant donné, on détermine un élément $\widetilde{u}$ de $L$ avec $\widetilde{u}\equiv u\bmod{d}$ et $\mathrm{q}(\widetilde{u})\equiv 0\bmod{d^{2}}$, $\mathrm{vois}_{d}(L;u)$ est alors le réseau de $\mathbb{Q}\otimes_{\mathbb{Z}}L$ engendré par $\mathrm{M}_{p}(L;u)$ et $\frac{\widetilde{u}}{d}$.

\vspace{0,75cm}
\textsc{$d$-voisins, point de vue abstrait}

\bigskip
Soient $L_{1}$ et $L_{2}$ deux $\mathrm{q}$-modules sur $\mathbb{Z}$, disons de même dimension $n$. Sans surprises, on dit que $L_{2}$ est $d$-voisin de $L_{1}$ si $L_{2}$ est isomorphe, comme\linebreak $\mathrm{q}$-module, à un $d$-voisin de $L_{1}$ dans $\mathbb{Q}\otimes_{\mathbb{Z}}L_{1}$. Compte tenu de ce qui précède, si $L_{2}$ est $d$-voisin de $L_{1}$ alors $L_{1}$ est $d$-voisin de $L_{2}$~; on dit donc aussi que $L_{1}$ et $L_{2}$ sont $d$-voisins. Pour éviter la confusion, il sera parfois commode de parler de voisinage ``abstrait'' pour la notion que nous venons d'introduire et de voisinage ``concret'' pour celle introduite dans la définition 1.2. Si le $\mathrm{q}$-espace vectoriel sur $\mathbb{R}\otimes_{\mathbb{Z}}L_{1}$ est indéfini alors $L_{2}$ est $d$-voisin de $L_{1}$ si et seulement si les deux $\mathrm{q}$-espaces vectoriels sur $\mathbb{R}\otimes_{\mathbb{Z}}L_{1}$ et $\mathbb{R}\otimes_{\mathbb{Z}}L_{2}$ sont isomorphes (Théorème II.2.7)~: la relation de $d$-voisinage abstrait ne présente pas grand intérêt dans ce cas~! On suppose donc ci-après que  $L_{1}$ et $L_{2}$ sont définis positifs (ce qui force $n$ à être divisible par $8$). Comme convenu, on abandonne la terminologie ``$\mathrm{q}$-module sur $\mathbb{Z}$ défini positif'' pour la terminologie ``réseau unimodulaire pair''.

\medskip
Soient $L_{1}$ et $L_{2}$ deux réseaux unimodulaires pairs. On note $\widetilde{\mathrm{Vois}}_{d}(L_{1},L_{2})$ l'ensemble des isomorphismes de $\mathrm{q}$-espaces vectoriels $\phi:\mathbb{Q}\otimes_{\mathbb{Z}}L_{2}\to\mathbb{Q}\otimes_{\mathbb{Z}}L_{1}$ avec $L_{1}/(L_{1}\cap\phi(L_{2}))$ cyclique d'ordre $d$. On observera que $\widetilde{\mathrm{Vois}}_{d}(L_{1},L_{2})$ est un ensemble fini~; par définition cet ensemble est non vide si et seulement si $L_{1}$ et $L_{2}$ sont $d$-voisins. On observera également que $\widetilde{\mathrm{Vois}}_{d}(L_{1},L_{2})$ est muni d'une action à gauche libre du groupe orthogonal $\mathrm{O}(L_{1})$ et d'une action à droite libre du groupe orthogonal $\mathrm{O}(L_{2})$, actions qui commutent. On note $\mathrm{Vois}_{d}(L_{1},L_{2})$ le sous-ensemble de $\mathrm{Vois}_{d}(L_{1})$ constitué des $d$-voisins de $L_{1}$ dans $\mathbb{Q}\otimes_{\mathbb{Z}}L_{1}$ qui sont isomorphes, comme $\mathrm{q}$-module, à $L_{2}$~; $\mathrm{Vois}_{d}(L_{1},L_{2})$ est canoniquement muni d'une action à gauche de $\mathrm{O}(L_{1})$. Par définition encore, l'application $\widetilde{\mathrm{Vois}}_{d}(L_{1},L_{2})\to\mathrm{Vois}_{d}(L_{1},L_{2})\hspace{2pt},\hspace{2pt}\phi\mapsto\phi(L_{2})$ induit une bijection $\mathrm{O}(L_{1})$-équivariante
$$
\hspace{24pt}
\widetilde{\mathrm{Vois}}_{d}(L_{1},L_{2})/\mathrm{O}(L_{2})
\hspace{4pt}\cong\hspace{4pt}
\mathrm{Vois}_{d}(L_{1},L_{2})
\hspace{24pt}.
$$
Le cardinal de l'ensemble $\mathrm{Vois}_{d}(L_{1},L_{2})$ est noté $\mathrm{N}_{d}(L_{1},L_{2})$.

\medskip
Attention~: l'introduction de cette notation n'est pas anodine, l'étude de ces cardinaux, en dimension 16 et 24, pour $d$ premier, est le thème principal de notre mémoire (voir chapitre I)~!

\medskip
On note $[-]$ la classe d'isomorphisme d'un réseau unimodulaire pair~; il est clair que $\mathrm{N}_{d}(L_{1},L_{2})$ ne dépend que de $[L_{1}]$ et $[L_{2}]$, pour cette raison cet entier sera aussi noté $\mathrm{N}_{d}([L_{1}],[L_{2}])$

\bigskip
\textbf{Lemme 1.6.} {\em L'application $\phi\mapsto\phi^{-1}$ induit une bijection de $\widetilde{\mathrm{Vois}}_{d}(L_{1},L_{2})$ sur $\widetilde{\mathrm{Vois}}_{d}(L_{2},L_{1})$.}

\bigskip
\textit{Démonstration.} Conséquence du ``point de vue symétrique'' sur les $d$-voisins, en clair de la proposition 1.1.
\hfill$\square$

\bigskip
\textbf{Scholie 1.7.} {\em On a la relation
$$
\frac{1}{\vert\mathrm{O}(L_{1})\vert}
\hspace{4pt}\mathrm{N}_{d}(L_{1},L_{2})
\hspace{4pt}=\hspace{4pt}
\frac{1}{\vert\mathrm{O}(L_{2})\vert}
\hspace{4pt}\mathrm{N}_{d}(L_{2},L_{1})
$$
(la notation $\vert-\vert$ désigne le cardinal d'un ensemble fini).}

\bigskip
\textit{Démonstration.} Observer que l'on a $\vert\widetilde{\mathrm{Vois}}_{d}(L_{1},L_{2})\vert=\mathrm{N}_{d}(L_{1},L_{2})\hspace{1pt}\vert\mathrm{O}(L_{2})\vert$.
\hfill$\square$

\bigskip
L'énoncé 1.7 peut être précisé en l'énoncé 1.8 ci-dessous~; pour une illustration de 1.8 voir 3.1 et 3.2.

\medskip
Ce qui précède montre que l'ensemble quotient $\mathrm{O}(L_{1})\backslash\widetilde{\mathrm{Vois}}_{d}(L_{1},L_{2})/\mathrm{O}(L_{2})$ est canoniquement en bijection avec l'ensemble quotient $\mathrm{O}(L_{1})\backslash\mathrm{Vois}_{d}(L_{1},L_{2})$ et l'ensemble quotient $\mathrm{O}(L_{2})\backslash\mathrm{Vois}_{d}(L_{2},L_{1})$ si bien que l'on dispose d'une bijection canonique  entre ces deux derniers ensembles quotients.

\bigskip
\textbf{Scholie 1.8.} {\em Soient $\Omega_{1}$ une $\mathrm{O}(L_{1})$-orbite dans $\mathrm{Vois}_{d}(L_{1},L_{2})$ et $\Omega_{2}$ une $\mathrm{O}(L_{2})$-orbite dans $\mathrm{Vois}_{d}(L_{2},L_{1})$ qui se correspondent par la bijection canonique
$$
\hspace{24pt}
\mathrm{O}(L_{1})\backslash\mathrm{Vois}_{d}(L_{1},L_{2})
\hspace{4pt}\cong\hspace{4pt}
\mathrm{O}(L_{2})\backslash\mathrm{Vois}_{d}(L_{2},L_{1})
\hspace{24pt};
$$
alors on a la relation
$$
\hspace{24pt}
\frac{\vert\Omega_{1}\vert}{\vert\mathrm{O}(L_{1})\vert}
\hspace{4pt}=\hspace{4pt}
\frac{\vert\Omega_{2}\vert}{\vert\mathrm{O}(L_{2})\vert}
\hspace{24pt}.
$$}

\bigskip
On peut se convaincre plus directement de l'égalité ci-dessus, en revenant à la notion concrète de $d$-voisinage.

\bigskip
\textbf{Proposition 1.9.} {\em Soit $V$ un $\mathrm{q}$-espace vectoriel sur $\mathbb{Q}$, défini positif~; on suppose que $V$ contient deux réseaux unimodulaires pairs $L_{1}$ et $L_{2}$ et que ces réseaux sont $d$-voisins dans  $V$. Soient $\Omega_{1}$ la $\mathrm{O}(L_{1})$-orbite de $L_{2}$ et $\Omega_{2}$ la $\mathrm{O}(L_{2})$-orbite de $L_{1}$~; alors on a la relation
$$
\frac{\vert\Omega_{1}\vert}{\vert\mathrm{O}(L_{1})\vert}
\hspace{4pt}=\hspace{4pt}
\frac{\vert\Omega_{2}\vert}{\vert\mathrm{O}(L_{2})\vert}
\hspace{4pt}=\hspace{4pt}
\frac{1}{\vert\mathrm{O}(L_{1})\cap\mathrm{O}(L_{2})\vert}
$$
($\mathrm{O}(L_{1})$ et $\mathrm{O}(L_{2})$ s'identifient à des sous-groupes de $\mathrm{O}(V)$ et l'intersection qui apparaît ci-dessus est effectuée dans $\mathrm{O}(V)$).}

\bigskip
\textit{Démonstration.} Observer que les stabilisateurs de $L_{2}$ pour l'action de $\mathrm{O}(L_{1})$ et de $L_{1}$ pour l'action de $\mathrm{O}(L_{2})$ s'identifie tous deux à $\mathrm{O}(L_{1})\cap\mathrm{O}(L_{2})$.
\hfill$\square$

\vspace{0,75cm}
\textsc{$2$-voisins, le point de vue de Borcherds} \cite[Chap. 17]{conwaysloane}

\bigskip
On commence par deux remarques concernant la proposition 1.9. On en reprend les notations et on considère le réseau $L_{1}\cap L_{2}$.

\medskip
-- Ce réseau est un $\widetilde{\mathrm{q}}$-module dont le résidu est muni d'un couple de lagrangiens, tous deux cycliques d'ordre $d$ et transverses l'un à l'autre, à savoir $\omega:=(L_{1}/(L_{1}\cap L_{2}),L_{2}/(L_{1}\cap L_{2}))$. On observe que le groupe $\mathrm{O}(L_{1})\cap\mathrm{O}(L_{2})$ s'identifie au sous-groupe de $\mathrm{O}(L_{1}\cap L_{2})$, disons $\mathrm{O}(L_{1}\cap L_{2};\omega)$, constitué des éléments qui préservent $\omega$. Les égalités de 1.9 peuvent donc être réécrites sous la forme suivante~:
$$
\hspace{24pt}
\vert\Omega_{1}\vert
\hspace{4pt}=\hspace{4pt}
\frac{\vert\mathrm{O}(L_{1})\vert}{\vert\mathrm{O}(L_{1}\cap L_{2};\omega)\vert}
\hspace{24pt},\hspace{24pt}
\vert\Omega_{2}\vert
\hspace{4pt}=\hspace{4pt}
\frac{\vert\mathrm{O}(L_{2})\vert}{\vert\mathrm{O}(L_{1}\cap L_{2};\omega)\vert}
\hspace{24pt}.
$$

\medskip
-- Dans le cas où $d$ est premier, la paire de lagrangiens sous-jacente au couple de lagrangiens $\omega$, est uniquement déterminée en fonction de $L_{1}\cap L_{2}$. Ceci implique que $\mathrm{O}(L_{1}\cap L_{2};\omega)$ est au plus d'indice $2$ dans $\mathrm{O}(L_{1}\cap L_{2})$. En particulier, on a $\mathrm{O}(L_{1}\cap L_{2};\omega)=\mathrm{O}(L_{1}\cap L_{2})$ si $L_{1}$ et $L_{2}$ ne sont pas isomorphes.

\bigskip
Ces remarques étant faites, nous pourrions continuer à traiter de $d$-voisins avec $d\geq 2$ quelconque, mais pour simplifier l'exposition, nous traitons  ci-dessous seulement de $p$-voisins avec $p$ premier~; en fait, le cas que nous avons en vue est $p=2$.

\medskip
Soit $n>0$ un entier divisible par $8$~; on rappelle que l'on note $\mathrm{X}_{n}$ l'ensemble fini des classes d'isomorphisme des réseaux unimodulaires pairs de dimension~$n$. Soit $p$ un nombre premier~; on introduit trois autres ensembles finis.

\smallskip
-- On note $\mathrm{Y}_{n}(p)$ l'ensemble des classes d'isomorphisme des couples $(L_{1},L_{2})$, $L_{1}$ désignant un réseau unimodulaire pair de dimension $n$ et $L_{2}$ un $p$-voisin de $L_{1}$ dans $\mathbb{Q}\otimes_{\mathbb{Z}}L_{1}$.

\smallskip
-- On note $\mathrm{B}_{n}(p)$ l'ensemble des classes d'isomorphisme des $\widetilde{\mathrm{q}}$-modules $M$ sur~$\mathbb{Z}$, avec $\dim M=n$, $\mathbb{R}\otimes_{\mathbb{Z}}M>0$ et $\mathop{\mathrm{r\acute{e}s}}M\simeq\mathrm{H}(\mathbb{Z}/p)$.

\smallskip
-- On note $\widetilde{\mathrm{B}}_{n}(p)$ l'ensemble des classes d'isomorphisme de couples  $(M;\omega)$ avec $M$ comme précédemment et $\omega$ une bijection de l'ensemble des lagrangiens de $\mathop{\mathrm{r\acute{e}s}}M$ sur l'ensemble $\{1,2\}$. Par définition $\widetilde{\mathrm{B}}_{n}(p)$ est muni d'une action à gauche du groupe symétrique $\mathfrak{S}_{2}$ et l'ensemble quotient $\mathfrak{S}_{2}\backslash\widetilde{\mathrm{B}}_{n}(p)$ s'identifie à $\mathrm{B}_{n}(p)$.

\bigskip
Nous avons tout fait pour que les ensembles $\mathrm{Y}_{n}(p)$ et $\widetilde{\mathrm{B}}_{n}(p)$ soient canoniquement en bijection.  Soit $(M;\omega)$ comme ci-dessus, on note $\mathrm{d}_{i}(M;\omega)$, $i=1,2$, l'image réciproque par la surjection $M^{\sharp}\to\mathop{\mathrm{r\acute{e}s}}M$ de $\omega^{-1}(i)$~; $\mathrm{d}_{1}(M;\omega)$ et $\mathrm{d}_{2}(M;\omega)$ sont deux réseaux unimodulaires pairs (de dimension $n$) $p$-voisins dans $\mathbb{Q}\otimes_{\mathbb{Z}}M$. En passant aux classes d'isomorphisme on obtient deux applications que l'on note encore $\mathrm{d}_{1}$ et $\mathrm{d}_{2}$ de $\widetilde{\mathrm{B}}_{n}(p)$ dans $\mathrm{X}_{n}$.

\medskip
Là encore, nous avons tout fait pour avoir l'énoncé suivant.
\vfill\eject

\bigskip
\textbf{Proposition 1.10.} {\em Soit $p$ un nombre premier~; soient $x_{1}$ et $x_{2}$ deux éléments de $\mathrm{X}_{n}$. On a~:
$$
\mathrm{N}_{p}(x_{1},x_{2})
\hspace{4pt}=\hspace{4pt}
\sum_{\beta\hspace{2pt}\in\hspace{2pt}\mathrm{d}_{1}^{-1}(x_{1})\hspace{1pt}\cap\hspace{1pt}\mathrm{d}_{2}^{-1}(x_{2})}\hspace{8pt}
\frac{\vert\mathrm{O}(x_{1})\vert}{\vert\mathrm{O}(\beta)\vert}
$$
(nous laissons au lecteur le soin de décoder les notations $\vert\mathrm{O}(x_{1})\vert$ et $\vert\mathrm{O}(\beta)\vert$).}

\bigskip
La proposition ci-dessus a un avatar un peu plus concret pour $p=2$ parce que la donnée d'un $\widetilde{\mathrm{q}}$-module $M$ sur~$\mathbb{Z}$, avec $\dim M=n$, $\mathbb{R}\otimes_{\mathbb{Z}}M>0$ et $\mathop{\mathrm{r\acute{e}s}}M\simeq\mathrm{H}(\mathbb{Z}/2)$, est équivalente à celle d'un réseau unimodulaire impair $L$ avec $\dim M=n$. Expliquons pourquoi (en suivant Borcherds).

\smallskip
-- A un $M$ comme ci-dessus on associe le réseau $L$ image réciproque par la surjection $M^{\sharp}\to\mathop{\mathrm{r\acute{e}s}}M$ de la droite du $\mathbb{Z}/2$-espace vectoriel $\mathop{\mathrm{r\acute{e}s}}M$ qui n'est pas isotrope au sens quadratique~; cette droite est isotrope au sens bilinéaire si bien que $L$ est un réseau unimodulaire impair.

\smallskip
-- A un réseau unimodulaire impair $L$ de dimension $n$ on associe le sous-module $M$ constitué des éléments $x$ avec $x.x\equiv 0 \bmod{2}$ (voir Scholie II.2.3).

\medskip
De plus l'ensemble des deux lagrangiens de $\mathop{\mathrm{r\acute{e}s}}M$ est naturellement en bijection avec l'ensemble des deux classes de vecteurs de Wu de $L$ (voir la discussion qui suit le scholie II.2.3).

\bigskip
On est donc amené à introduire les notations suivantes.

\smallskip
-- On note $\mathrm{B}_{n}$ l'ensemble fini des classes d'isomorphisme de réseaux unimodulaires $L$ impairs de dimension $n$.

\smallskip
-- On note $\widetilde{\mathrm{B}}_{n}$ l'ensemble fini des classes d'isomorphisme des couples $(L;\omega)$, $L$ désignant un réseau unimodulaire impair de dimension $n$ et $\omega$ une bijection de l'ensemble des deux classes de vecteurs de Wu de $L$ sur $\{1,2\}$. Par définition, $\widetilde{\mathrm{B}}_{n}$ est muni d'une action à gauche du groupe symétrique $\mathfrak{S}_{2}$ et l'ensemble quotient $\mathfrak{S}_{2}\backslash\widetilde{\mathrm{B}}_{n}$ s'identifie à $\mathrm{B}_{n}$.

\smallskip
-- Soit $L_{i}$, $i=1,2$, le réseau unimodulaire pair de $\mathbb{Q}\otimes_{\mathbb{Z}}L$, engendré par le sous-module de $L$ constitué des éléments $x$ avec $x.x\equiv 0 \bmod{2}$ et $\frac{1}{2}\hspace{1pt}\omega^{-1}(i)$. On note encore $\mathrm{d}_{i}:\widetilde{\mathrm{B}}_{n}\to\mathrm{X}_{n}$, $i=1,2$, l'application qui associe à la classe d'isomorphisme de $(L;\omega)$, la classe d'isomorphisme de $L_{i}$.

\medskip
Par construction, les réseaux unimodulaires $L$ et $L_{i}$, $i=1,2$, sont $2$-voisins ``au sens bilinéaire''~: $L\cap L_{i}$ est d'indice $2$ dans $L$ et $L_{i}$. Le lecteur vérifiera que $L_{1}$  et $L_{2}$ sont caractérisés, parmi les réseaux unimodulaires de $\mathbb{Q}\otimes_{\mathbb{Z}} L$ qui sont $2$-voisins de $L$, par le fait qu'ils sont pairs. Nous dirons que $L_{1}$ et $L_{2}$ sont les {\em $2$-voisins pairs} de $L$.

\medskip
Voici enfin  l'avatar annoncé de la proposition 1.10 pour $p=2$. Son énoncé semble identique~; les seules différences sont en fait les suivantes~: $\mathrm{d}_{1}$ et $\mathrm{d}_{2}$ désigne maintenant les applications de $\widetilde{\mathrm{B}}_{n}$ dans $\mathrm{X}_{n}$ introduites ci-dessus, $\beta$ appartient à $\widetilde{\mathrm{B}}_{n}$ et si $\beta$ est représenté par un réseau unimodulaire impair $L$ muni d'une bijection $\omega$ de l'ensemble de ses deux classes de vecteurs de Wu sur $\{1,2\}$ alors $\mathrm{O}(\beta)$ est le sous-groupe de $\mathrm{O}(L)$ qui préserve $\omega$. 

\bigskip
\textbf{Proposition 1.11.} {\em Soient $x_{1}$ et $x_{2}$ deux éléments de $\mathrm{X}_{n}$. On a~:
$$
\hspace{24pt}
\mathrm{N}_{2}(x_{1},x_{2})
\hspace{4pt}=\hspace{4pt}
\sum_{\beta\hspace{2pt}\in\hspace{2pt}\mathrm{d}_{1}^{-1}(x_{1})\hspace{1pt}\cap\hspace{1pt}\mathrm{d}_{2}^{-1}(x_{2})}\hspace{8pt}
\frac{\vert\mathrm{O}(x_{1})\vert}{\vert\mathrm{O}(\beta)\vert}
\hspace{24pt}.
$$}

\vspace{0,75cm}
\textsc{Les graphes des $p$-voisins}

\bigskip
Soient $n>0$ un nombre entier divisible par $8$ et $p$ un nombre premier~; le {\em graphe des $p$-voisins} $\mathrm{K}_{n}(p)$ est le graphe dont l'ensemble des sommets est l'ensemble $\mathrm{X}_{n}$ des classes de réseaux unimodulaires pairs de dimension~$n$ et dont l'ensemble des arêtes est l'ensemble des parties à deux éléments $\{[L_{1}],[L_{2}]\}$ de $\mathrm{X}_{n}$ avec $L_{1}$ et $L_{2}$ $p$-voisins (on rappelle que l'on note $[L]$ la classe d'isomorphisme d'un réseau unimodulaire pair $L$).

\bigskip
\textbf{Théorème 1.12 } (M. Kneser)\textbf{.} {\em Pour tout $n$ et tout $p$ le graphe $\mathrm{K}_{n}(p)$ est connexe.}

\bigskip
\textit{Démonstration.} Soient $L$ et $M$ deux réseaux unimodulaires pairs de même dimension~; il faut montrer qu'il existe une suite finie de réseaux unimodulaires pairs
$$
L=L_{0}\hspace{1pt},\hspace{1pt}
L_{1}\hspace{1pt},\hspace{1pt}
L_{2}\hspace{1pt},\hspace{1pt}
\ldots\hspace{1pt},\hspace{1pt}
L_{m-1}\hspace{1pt},\hspace{1pt}
L_{m}=M
$$
avec $L_{k}$ et $L_{k+1}$ (abstraitement) $p$-voisins pour $0\leq k\leq m-1$.

\medskip
Le théorème II.2.8 montre qu'il existe un isomorphisme de $\mathrm{q}$-modules $\phi:\mathbb{Z}[\frac{1}{p}]\otimes_{\mathbb{Z}}L\to\mathbb{Z}[\frac{1}{p}]\otimes_{\mathbb{Z}}M$. Quitte à remplacer $M$ par $\phi^{-1}(M)$, on peut donc supposer $M\subset\mathbb{Z}[\frac{1}{p}]\otimes_{\mathbb{Z}}L\subset\mathbb{Q}\otimes_{\mathbb{Z}}L$~; on pose $V=\mathbb{Q}\otimes_{\mathbb{Z}}L=\mathbb{Q}\otimes_{\mathbb{Z}}M$. On rappelle maintenant ce que dit le point (a) de 1.1. On considère $L$ et~$M$ comme des réseaux de $V$ et on pose $N=L\cap M$. On a $N^{\sharp}=L+M$ et $\mathop{\mathrm{r\acute{e}s}}N=L/N\oplus M/N$ (comme groupe abélien). On pose $I=L/N$ et $J=M/N$~; l'accouplement $I\times J\to\mathbb{Q}/\mathbb{Z}$, induit par la forme d'enlacement de $\mathop{\mathrm{r\acute{e}s}}N$, est non dégénéré. On en déduit que $J$ et $\mathop{\mathrm{r\acute{e}s}}N$ sont canoniquement isomorphes, respectivement au dual
de Pontryagin $I^{\vee}$ et au $\mathrm{qe}$-module d'enlacement hyperbolique $\mathrm{H}(I)$. Puisque l'on a $M\subset\mathbb{Z}[\frac{1}{p}]\otimes_{\mathbb{Z}}L$, le groupe abélien fini $I$ est un $p$-groupe~; soit
$$
I=I_{0}\supset I_{1}\supset I_{2}
\supset\ldots\supset I_{m-1}\supset I_{m}=0
$$
une suite finie décroissante de sous-groupes de
$I$ avec $I_{k}/I_{k+1}\simeq\mathbb{Z}/p$ pour
$0\leq k\leq m-1$. Soit
$$
0=J_{0}\subset J_{1}\subset I_{2}
\subset\ldots\subset J_{m-1}
\subset J_{m}=J
$$
la suite ``orthogonale'' de sous-groupes de $J$. On pose $K_{k}=I_{k}\oplus J_{k}$~; on constate que $K_{k}$ est un lagrangien de $\mathop{\mathrm{r\acute{e}s}}N$. On note $L_{k}$ l'image inverse de $K_{k}$ par l'homomorphisme canonique $N^{\sharp}\to
\mathrm{r\acute{e}s}
\hspace{1pt}N$. Par
construction
$$
L=L_{0}\hspace{1pt},\hspace{1pt}
L_{1}\hspace{1pt},\hspace{1pt}
L_{2}\hspace{1pt},\hspace{1pt}
\ldots\hspace{1pt},\hspace{1pt}
L_{m-1}\hspace{1pt},\hspace{1pt}
L_{m}=M
$$
est une suite de réseaux unimodulaires pairs avec  $L_{k}$ et $L_{k+1}$ (concrètement) $p$-voisins pour $0\leq k\leq m-1$.
\hfill$\square$

\vspace{0,75cm}
\textsc{Miscellanées}

\bigskip
On rassemble sous cet intertitre quatre énoncés plus ou moins  techniques (1.13, 1.14, 1.16 et 1.17), concernant la notion de $d$-voisinage, que nous invoquerons par la suite.

\bigskip
Le moins technique des quatre est l'énoncé 1.13 ci-dessous~; sa démonstration, qui est une illustration de la proposition II.1.1, est laissée au lecteur.

\bigskip
\textbf{Proposition 1.13.} {\em Soient $L$ un $\mathrm{q}$-module sur $\mathbb{Z}$ et $d\geq 2$ un nombre entier~; soit $u$ un élément $d$-primitif de $L$ avec $\mathrm{q}(u)\equiv 0\bmod{d^{2}}$. On suppose $d=d_{1}d_{2}$ avec $d_{1}\geq 2$ et $d_{2}\geq 2$. Alors $\frac{u}{d_{1}}$ est un élément $d_{2}$-primitif de $\mathrm{vois}_{d_{1}}(L;u)$ et l'on a, dans $\mathbb{Q}\otimes_{\mathbb{Z}}L$, l'égalité suivante
$$
\hspace{24pt}
\mathrm{vois}_{d}(L;u)
\hspace{4pt}=\hspace{4pt}
\mathrm{vois}_{d_{2}}(\mathrm{vois}_{d_{1}}(L;u);\frac{u}{d_{1}})
\hspace{24pt}.
$$}

\bigskip
\textbf{Proposition 1.14.} {\em Soient $L$ un $\mathrm{q}$-module sur $\mathbb{Z}$ et $d\geq 2$ un nombre entier~; soit $u$ un élément $d$-primitif de $L$ avec $\mathrm{q}(u)=d$. On note $\mathrm{s}_{u}$ la symétrie orthogonale de $\mathbb{Q}\otimes_{\mathbb{Z}}L$ par rapport à l'hyperplan $u^{\perp}$.

\medskip
{\em (a)} On a dans $\mathbb{Q}\otimes_{\mathbb{Z}}L$ l'égalité
$$
\hspace{24pt}
\mathrm{vois}_{d}(L;u)
\hspace{4pt}=\hspace{4pt}
\mathrm{s}_{u}(L)
\hspace{24pt};
$$
en particulier $\mathrm{vois}_{d}(L;u)$ est isomorphe à $L$.
\vfill\eject

\medskip
{\em (b)} On suppose $d=d_{1}d_{2}$ avec $d_{1}\geq 2$ et $d_{2}\geq 2$. Alors $u$ est $d_{i}$-primitif pour $i=1,2$ et l'on a dans $\mathbb{Q}\otimes_{\mathbb{Z}}L$ l'égalité
$$
\hspace{24pt}
\mathrm{vois}_{d_{2}}(L;u)
\hspace{4pt}=\hspace{4pt}
\mathrm{s}_{u}(\mathrm{vois}_{d_{1}}(L;u))
\hspace{24pt};
$$
en particulier $\mathrm{vois}_{d_{1}}(L;u)$ et $\mathrm{vois}_{d_{2}}(L;u)$ sont isomorphes.}

\bigskip
\textit{Démonstration.} Explicitons $\mathrm{s}_{u}$~:
$$
\mathrm{s}_{u}(x)
\hspace{4pt}=\hspace{4pt}
 x-\frac{u.x}{d}\hspace{2pt}u
$$
pour tout $x$ dans $\mathbb{Q}\otimes_{\mathbb{Z}}L$. Cette expression montre que $\mathrm{s}_{u}$ induit un automorphisme de $\mathrm{M}_{d}(L;u)$ comme $\widetilde{\mathrm{q}}$-module sur $\mathbb{Z}$ et donc un automorphisme de $\mathop{\mathrm{r\acute{e}s}}\mathrm{M}_{d}(L;u)$ comme $\mathrm{qe}$-module sur $\mathbb{Z}$. Soient $I$ et $J$ les deux lagrangiens transverses de $\mathop{\mathrm{r\acute{e}s}}\mathrm{M}_{d}(L;u)$ correspondant respectivement, \textit{via} le point (b) de II.1.1, aux réseaux $L$ et $\mathrm{vois}_{d}(L;u)$. Soit $v$ un élément de $L$ avec $u.v\equiv 1\bmod{d}$~; on rappelle que $I$ et $J$ sont respectivement engendrés par les classes dans $\mathop{\mathrm{r\acute{e}s}}\mathrm{M}_{d}(L;u)$ de $v$ et $\frac{1}{d}\hspace{1pt}(u-\mathrm{q}(u)\hspace{1pt}v)$. On constate que l'on a $\mathrm{s}_{u}(I)=J$~; cette égalité implique le point (a) de la proposition.

\medskip
On constate que les réseaux $\mathrm{vois}_{d_{1}}(L;u)$ et $\mathrm{vois}_{d_{2}}(L;u)$ correspondent, \textit{via} la bijection du point (b) de II.1.1, respectivement aux lagrangiens $d_{1}I\oplus d_{2}J$ et $d_{2}I\oplus d_{1}J$ de $\mathop{\mathrm{r\acute{e}s}}\mathrm{M}_{d}(L;u)$. D'après ce qui précède, on a $\mathrm{s}_{u}(d_{1}I\oplus d_{2}J)=d_{2}I\oplus d_{1}J$~; cette égalité implique le point (b) de la proposition.
\hfill$\square$

\medskip
\textit{Remarque.} Si l'on convient que $L$ est le seul $1$-voisin de $L$ dans $\mathbb{Q}\otimes_{\mathbb{Z}}L$ alors le point (a) de la proposition 1.14 apparaît comme un cas particulier de son point (b).

\bigskip
\textbf{Proposition 1.15.} {\em Soit $L$ un $\mathrm{b}$-module sur $\mathbb{Z}$. On suppose qu'il existe un élément $e$ de $L$ avec $e.e=1$ (et donc que $L$ est impair). Alors la symétrie orthogonale de $L$ par rapport à l'hyperplan $e^{\perp}$ échange les deux classes de vecteurs de Wu de $L$.}

\bigskip
\textit{Démonstration.} Soient $\mathrm{s}_{e}$ la symétrie en question et $u$ un vecteur de Wu de~$L$~; l'égalité $\mathrm{s}_{e}(u)=u-2(e.u)e$ et la congruence $e.u\equiv e.e\bmod{2}$ montrent que les vecteurs de Wu $u$ et $\mathrm{s}_{e}(u)$ sont non équivalents.
\hfill$\square$

\bigskip
\textbf{Corollaire 1.16.} {\em Soit $L$ un réseau unimodulaire impair de dimension divisible par $8$. On suppose qu'il existe un élément $e$ de~$L$ avec $e.e=1$ . Alors la symétrie orthogonale de $\mathbb{Q}\otimes_{\mathbb{Z}}L$ par rapport à l'hyperplan $e^{\perp}$ échange les deux réseaux unimodulaires pairs $2$-voisins de $L$.}

\bigskip
La proposition suivante, dont la vérification est immédiate, montre que le point (a) de la proposition 1.14, dans le cas particulier $d=2$, et le corollaire 1.16 sont intimement reliés.

\bigskip
\textbf{Proposition 1.17.} {\em Soient $L$ un réseau unimodulaire pair et $u$ un élément de~$L$ avec $\mathrm{q}(u)=2$ (égalité qui implique que $u$ est $2$-primitif). Soit $B$ le réseau impair de $\mathbb{Q}\otimes_{\mathbb{Z}}L$ dont les $2$-voisins pairs sont $L$ et $\mathrm{vois}_{2}(L;u)$ (voir la fin de la discussion intitulée ``$2$-voisins, le point de vue de Borcherds''). On note $e$ l'élément $\frac{u}{2}$ de $\mathbb{Q}\otimes_{\mathbb{Z}}L$. Alors~:

\smallskip
-- on a $e.e=1$~;

\smallskip
-- le réseau $B$ est engendré dans $\mathbb{Q}\otimes_{\mathbb{Z}}L$ par $\mathrm{M}_{2}(L;u)$ et $e$~;

\smallskip
- les symétries orthogonales $\mathrm{s}_{u}$ et $\mathrm{s}_{e}$ de $\mathbb{Q}\otimes_{\mathbb{Z}}L$ coïncident.}

\section{Sur les opérateurs de Hecke associés à la notion de voisinage}

Soit $n>0$ un entier divisible par $8$~; on rappelle que l'on note $\mathrm{X}_{n}$ l'ensemble fini des classes d'isomorphisme de réseaux unimodulaires de dimension $n$.

\medskip
On note $\mathbb{Z}[\mathrm{X}_{n}]$ le $\mathbb{Z}$-module libre engendré par l'ensemble $\mathrm{X}_{n}$.

\medskip
Soient $A$ un groupe abélien fini et $L$ un réseau unimodulaire pair de dimension~$n$~; on note $\mathrm{Vois}_{A}(L)$ l'ensemble fini constitué des $A$-voisins de $L$ dans $\mathbb{Q}\otimes_{\mathbb{Z}}L$ (la notation $\mathrm{Vois}_{d}(-)$ introduite dans le paragraphe précédent est donc une abréviation de $\mathrm{Vois}_{\mathbb{Z}/d}(-)$).

\medskip
{\em L'opérateur de Hecke} $\mathop{\mathrm{T}_{A}}$ est l'endomorphisme de~$\mathbb{Z}[\mathrm{X}_{n}]$ défini par
$$
\hspace{24pt}
\mathop{\mathrm{T}_{A}}\hspace{1pt}[L]
\hspace{4pt}:=\hspace{4pt}
\sum_{L'\in\mathrm{Vois}_{A}(L)}[L']
\hspace{24pt},
$$
pour tout réseau unimodulaire pair $L$ de dimension $n$.

\medskip
\textit{Remarques}

\smallskip
-- Soit $\mathrm{r}(A)$ le nombre minimum de générateurs du groupe abélien $A$. Le point (b) de 1.1 montre que $\mathop{\mathrm{T}_{A}}$ est nul si l'on a $2\hspace{1pt}\mathrm{r}(A)>n$. On pourrait donc supposer $2\hspace{1pt}\mathrm{r}(A)\leq n$ dans la définition ci-dessus.

\smallskip
-- Soient $A$ et $B$ deux groupes abéliens finis~; il n'est pas trop difficile de se convaincre de ce que les opérateurs de Hecke $\mathop{\mathrm{T}_{A}}$ et $\mathop{\mathrm{T}_{B}}$ commutent si les cardinaux de $A$ et $B$ sont premiers entre eux. En fait, $\mathop{\mathrm{T}_{A}}$ et $\mathop{\mathrm{T}_{B}}$ commutent pour tous $A$ et $B$~; ceci sera démontré et largement généralisé au chapitre IV.

\bigskip
Soit $d\geq 2$ un entier~; on abrégera la notation $\mathop{\mathrm{T}_{\mathbb{Z}/d}}$ en $\mathop{\mathrm{T}_{d}}$. Par définition même des entiers $\mathrm{N}_{d}(x,y)$ on a~:
$$
\mathop{\mathrm{T}_{d}}\hspace{1pt}x
\hspace{4pt}=\hspace{4pt}
\sum_{y\in\mathrm{X}_{n}}
\mathrm{N}_{d}(x,y)\hspace{2pt}y
$$
pour tout $x$ dans $\mathrm{X}_{n}$. En d'autres termes, si l'on considère $\mathop{\mathrm{T}_{d}}$ comme une $(\mathrm{X}_{n},\mathrm{X}_{n})$-matrice, alors son coefficient d'indice $(y,x)$ est $\mathrm{N}_{d}(x,y)$.

\medskip
La proposition 1.5 montre que l'on a aussi
$$
\mathop{\mathrm{T}_{d}}\hspace{1pt}[L]
\hspace{4pt}=\hspace{4pt}
\sum_{c\in\mathrm{C}_{L}(\mathbb{Z}/d)}
[\mathrm{Vois}_{d}(L;c)]
$$
pour tout réseau unimodulaire pair $L$.

\bigskip
Nous dégageons l'énoncé 2.2 ci-dessous en vue de futures références. Cet énoncé est impliqué par 2.1 et 1.5~; l'énoncé 2.1 est quant à lui essentiellement conséquence de II.2.5.

\bigskip
\textbf{Proposition-Définition 2.1.} {\em Soient $L$ un réseau unimodulaire pair et $d\geq 2$ un entier. Alors le $\mathrm{q}$-module $\mathbb{Z}/d\otimes_{\mathbb{Z}}L$ est hyperbolique. En particulier, le cardinal de la quadrique $\mathrm{C}_{L}(\mathbb{Z}/d)$ ne dépend que la dimension de $L$, disons~$n$, on le note $\mathrm{c}_{n}(d)$. On a
$$
\mathrm{c}_{n}(p)
\hspace{4pt}=\hspace{4pt}
\sum_{m=0}^{n-2}\hspace{2pt}p^{m}\hspace{4pt}+\hspace{4pt}p^{\frac{n}{2}-1}
$$
pour tout nombre premier $p$.}

\bigskip
\textit{Remarque (suite).} Le calcul de $\mathrm{c}_{n}(d)$ pour tout $d$, se déduit facilement du calcul pour $d$ premier et du fait que les quadriques $\mathrm{C}_{L}$ sont lisses sur $\mathbb{Z}$.

\bigskip
\textbf{Proposition 2.2.} {\em Soit  $d\geq 2$ un entier. On a
$$
\sum_{y\in\mathrm{X}_{n}}
\mathrm{N}_{d}(x,y)
\hspace{4pt}=\hspace{4pt}
\mathrm{c}_{n}(d)
$$
pour tout $x$ dans $\mathrm{X}_{n}$.}

\bigskip
\textit{Remarques (suite et fin)}

\medskip
-- Le scholie 1.7 peut être reformulé de la façon suivante~:

\bigskip
\textbf{Proposition 2.3.} {\em Soit $d\geq 2$ un entier. L'endomorphisme $\mathrm{T}_{d}$ de $\mathbb{Z}[\mathrm{X}_{n}]$ est auto-adjoint pour le produit scalaire, disons $(-\vert-)$, défini par
$$
(x\vert y)
\hspace{4pt}=\hspace{4pt}
\vert\mathrm{O}(x)\vert\hspace{4pt}\delta_{x,y}
$$
pour $x$ et $y$ dans $\mathrm{X}_{n}$ (la notation $\delta_{x,y}$ désignant le symbole de Kronecker).}

\bigskip
Là encore, cet énoncé sera largement généralisé au chapitre IV.

\bigskip
-- Soit $\epsilon:\mathbb{Z}[\mathrm{X}_{n}]\to\mathbb{Z}$ l'homomorphisme de $\mathbb{Z}$-modules défini par $\epsilon(x)=1$ pour tout $x$ dans $\mathrm{X}_{n}$~; la proposition 2.2 dit que l'on a $\epsilon\circ\mathrm{T}_{d}=\mathrm{c}_{n}(d)\hspace{1pt}\epsilon$, en d'autres termes, que $\epsilon$ est vecteur propre de l'endomorphisme $\mathrm{T}_{d}^{*}$ de $(\mathbb{Z}[\mathrm{X}_{n}])^{*}$, pour la valeur propre $\mathrm{c}_{n}(d)$ (la notation $\mathrm{T}_{d}^{*}$ est une notation alternative pour la notation ${}^{\mathrm{t}}\mathrm{T}_{d}$ utilisée par ailleurs). Cette observation et la proposition 2.3 conduisent à l'énoncé suivant~:

\bigskip
\textbf{Proposition 2.4.} {\em Soit $d\geq 2$ un entier. L'élément 
$$
\sum_{x\in\mathrm{X}_{n}}\frac{1}{\vert\mathrm{O}(x)\vert}\hspace{4pt}x
$$
de $\mathbb{Q}\otimes_{\mathbb{Z}}\mathbb{Z}[\mathrm{X}_{n}]$ est vecteur propre de $\mathrm{T}_{d}$ pour la valeur propre $\mathrm{c}_{n}(d)$.}

\bigskip
-- La proposition 2.3 implique que $\mathrm{T}_{d}$ est diagonalisable, au moins après extension des scalaires à $\mathbb{R}$. En fait, pour $n=8,16,24$, les valeurs propres de $\mathrm{T}_{d}$ sont entières. C'est trivial pour $n=8$ puisque l'on a $\vert\mathrm{X}_{8}\vert=1$~; cela l'est presque autant pour $n=16$ puisque l'on a $\vert\mathrm{X}_{16}\vert=2$ et que l'on connaît déjà une valeur propre entière, à savoir $\mathrm{c}_{16}(d)$. Le cas $n=24$ requiert davantage d'efforts. Nous expliquerons en 3.3 comment Gabriele Nebe et Boris Venkov ont déterminé $\mathrm{T}_{2}$, à partir de travaux de Borcherds~; on constate (merci \texttt{PARI}) que les racines du polynôme caractéristique de $\mathrm{T}_{2}$ sont entières et simples. Puisque $\mathrm{T}_{d}$ commute avec $\mathrm{T}_{2}$ pour tout $d$, les vecteurs propres de $\mathrm{T}_{2}$ sont aussi propres pour $\mathrm{T}_{d}$ et les valeurs propres de $\mathrm{T}_{d}$ sont entières.

\medskip
L'étude des propriétés arithmétiques de ces valeurs propres, pour $d$ premier, est la motivation du présent travail.

\section{Exemples}

\vspace{0,2cm}
\textbf{3.1.} Détermination de $\mathrm{T}_{2}$ pour $n=16$

\medskip
Il est bien connu que l'homomorphisme canonique $\mathrm{O}(\mathrm{E}_{8})\to\mathrm{O}(\mathbb{F}_{2}\otimes_{\mathbb{Z}}\mathrm{E}_{8})$ induit un isomorphisme $\mathrm{O}(\mathrm{E}_{8})/\{\pm{1}\}\cong\mathrm{O}(\mathbb{F}_{2}\otimes_{\mathbb{Z}}\mathrm{E}_{8})$ (voir par exemple \cite[Ch. VI, \S4, Exc. 1]{bourbaki}). Il en résulte que l'action de $\mathrm{O}(\mathrm{E}_{8})$ partitionne\linebreak $\mathbb{F}_{2}\otimes_{\mathbb{Z}}\mathrm{E}_{8}-\{0\}$ en deux orbites à savoir $\mathrm{q}^{-1}(0)$ et $\mathrm{q}^{-1}(1)$~; ces deux orbites ont respectivement $135$ et $120$ éléments.
\vfill\eject

\medskip
On considère maintenant les réseaux $\mathrm{E}_{8}\oplus\mathrm{E}_{8}$ et $\mathrm{E}_{16}$ que l'on plonge de façon standard dans $\mathbb{Q}^{16}$~; on note $(\varepsilon_{1},\varepsilon_{2},\ldots,\varepsilon_{16})$ la base canonique de $\mathbb{Q}^{16}$.

\medskip
Il est clair que le groupe $\mathrm{O}(\mathrm{E}_{8}\oplus\mathrm{E}_{8})$ est canoniquement isomorphe au produit en couronne $\mathfrak{S}_{2}\wr\mathrm{O}(\mathrm{E}_{8})$. D'après ce qui précède, l'action de $\mathrm{O}(\mathrm{E}_{8}\oplus\mathrm{E}_{8})$ partitionne la quadrique $\mathrm{C}_{\mathrm{E}_{8}\oplus\mathrm{E}_{8}}(\mathbb{F}_{2})$ en trois orbites~:

\smallskip
-- l'orbite du point $[2\hspace{1pt}\varepsilon_{1}]$ (le vecteur $2\hspace{1pt}\varepsilon_{1}$ appartient au réseau $\mathrm{E}_{8}\oplus\mathrm{E}_{8}$ et vérifie $\mathrm{q}(2\hspace{1pt}\varepsilon_{1})=2$, la notation $[2\hspace{1pt}\varepsilon_{1}]$ désigne sa classe dans $\mathrm{C}_{\mathrm{E}_{8}\oplus\mathrm{E}_{8}}(\mathbb{F}_{2})$), cette orbite a $2\times135=270$ éléments~;

\smallskip
-- l'orbite du point $[\varepsilon_{1}+\varepsilon_{2}+\varepsilon_{9}+\varepsilon_{10}]$ (on observera que l'on a à nouveau $\mathrm{q}(\varepsilon_{1}+\varepsilon_{2}+\varepsilon_{9}+\varepsilon_{10})=2)$, cette orbite a $120^{\hspace{1pt}2}=14400$ éléments~;

\smallskip
-- l'orbite du point $[2\hspace{1pt}\varepsilon_{1}+2\hspace{1pt}\varepsilon_{9}]$ (on observera que l'on a cette fois $\mathrm{q}(2\hspace{1pt}\varepsilon_{1}+2\hspace{1pt}\varepsilon_{9})=4)$, cette orbite a $135^{\hspace{1pt}2}=18225$ éléments.

\medskip
D'après le point (a) de 1.14, le réseau $\mathrm{vois}_{2}(\mathrm{E}_{8}\oplus\mathrm{E}_{8};c)$ est isomorphe (comme $\mathrm{q}$-module) à $\mathrm{E}_{8}\oplus\mathrm{E}_{8}$ pour tout $c$ dans l'une des deux premières orbites.

\medskip
Comme le graphe des $2$-voisins est connexe (Théorème 1.12), on a forcément $\mathrm{N}_{2}(\mathrm{E}_{8}\oplus\mathrm{E}_{8},\mathrm{E}_{16})=18225$. Compte tenu de 1.7 et 2.2, cette égalité détermine l'opérateur de Hecke $\mathrm{T}_{2}:\mathbb{Z}[\mathrm{X}_{16}]\to\mathbb{Z}[\mathrm{X}_{16}]$~; sa matrice dans la base $(\mathrm{E}_{16},\mathrm{E}_{8}\oplus\nolinebreak\mathrm{E}_{8})$, disons encore $\mathrm{T}_{2}$, est la suivante~:
$$
\hspace{24pt}
\mathrm{T}_{2}=
\begin{bmatrix}
20025 & 18225  \\ 12870 & 14670
\end{bmatrix}
\hspace{24pt}.
$$
En fait, il est facile de vérifier que l'on a $\mathrm{vois}_{2}(\mathrm{E}_{8}\oplus\mathrm{E}_{8};[2\hspace{1pt}\varepsilon_{1}+2\hspace{1pt}\varepsilon_{9}])=\mathrm{E}_{16}$ (ce qui implique $\mathrm{N}_{2}(\mathrm{E}_{8}\oplus\mathrm{E}_{8},\mathrm{E}_{16})=18225$). En effet, on constate que le réseau  $\mathrm{M}_{2}(\mathrm{E}_{8}\oplus\mathrm{E}_{8};[2\hspace{1pt}\varepsilon_{1}+2\hspace{1pt}\varepsilon_{9}])$ est engendré par $\mathrm{D}_{8}\oplus\mathrm{D}_{8}$ et $\frac{1}{2}\sum_{i=1}^{16}\varepsilon_{i}$ [le premier (resp. second) $\mathrm{D}_{8}$ est l'orthogonal modulo $2$ dans le premier (resp. second) $\mathrm{E}_{8}$ du vecteur $2\hspace{1pt}\varepsilon_{1}$ (resp. $2\hspace{1pt}\varepsilon_{9}$)]. Comme l'on a $\mathrm{q}(2\hspace{1pt}\varepsilon_{1}+2\hspace{1pt}\varepsilon_{9})=4$, le réseau $\mathrm{vois}_{2}(\mathrm{E}_{8}\oplus\mathrm{E}_{8},[2\hspace{1pt}\varepsilon_{1}+2\hspace{1pt}\varepsilon_{9}])$ est engendré par $\mathrm{M}_{2}(\mathrm{E}_{8}\oplus\mathrm{E}_{8},[2\hspace{1pt}\varepsilon_{1}+2\hspace{1pt}\varepsilon_{9}])$ et $\varepsilon_{1}+\varepsilon_{9}$ (se rappeler l'algorithme $u\leadsto\mathrm{vois}_{d}(L;u)$). Or le réseau engendré par $\mathrm{D}_{8}\oplus\mathrm{D}_{8}$ et  $\varepsilon_{1}+\varepsilon_{9}$ est $\mathrm{D}_{16}$ si bien que le réseau $\mathrm{vois}_{2}(\mathrm{E}_{8}\oplus\mathrm{E}_{8},[2\hspace{1pt}\varepsilon_{1}+2\hspace{1pt}\varepsilon_{9}])$ coïncide avec le réseau engendré par $\mathrm{D}_{16}$ et $\frac{1}{2}\sum_{i=1}^{16}\varepsilon_{i}$, c'est-à-dire $\mathrm{E}_{16}$.

\bigskip
\textit{Variante.} Pour illustrer 1.8, et accessoirement pour nous rassurer, nous consi\-dé\-rons maintenant les $2$-voisins de~$\mathrm{E}_{16}$.

\medskip
Le groupe $\mathrm{O}(\mathrm{E}_{16})$ s'identifie au sous-groupe de $\mathrm{O}(\mathrm{I}_{16})$ constitué des éléments qui préservent la classe (au sens de la discussion qui suit le scholie II.2.3) du vecteur de Wu  $\sum_{i=1}^{16}\varepsilon_{i}$. On a donc un isomorphisme canonique $\mathrm{O}(\mathrm{E}_{16})\cong\mathfrak{S}_{16}\ltimes(\{\pm 1\}^{16})^{0}$, la notation $(\{\pm 1\}^{16})^{0}$ désignant le sous-groupe de $\{\pm 1\}^{16}$ constitué des $16$-uples $(\eta_{1},\eta_{2},\ldots,\eta_{16})$ avec $\eta_{1}\eta_{2}\ldots\eta_{16}=1$.

\medskip
On constate que l'action de $\mathrm{O}(\mathrm{E}_{16})$ partitionne la quadrique $\mathrm{C}_{\mathrm{E}_{16}}(\mathbb{F}_{2})$ en quatre orbites~:

\smallskip
-- l'orbite du point $[2\hspace{1pt}\varepsilon_{1}]$ (le vecteur $2\hspace{1pt}\varepsilon_{1}$ appartient au réseau $\mathrm{E}_{16}$ et vérifie $\mathrm{q}(2\hspace{1pt}\varepsilon_{1})=2$), cette orbite a un seul élément~;

\smallskip
-- l'orbite du point $[\sum_{i=1}^{4}\varepsilon_{i}]$ (on observera que l'on a $\mathrm{q}(\sum_{i=1}^{4}\varepsilon_{i})=2)$, cette orbite a $2{16\choose 4}=3640$ éléments~;

\smallskip
-- l'orbite du point $[\sum_{i=1}^{8}\varepsilon_{i}]$ (on observera que l'on a $\mathrm{q}(\sum_{i=1}^{8}\varepsilon_{i})=4)$, cette orbite a ${16\choose 8}=12870$ éléments.

\smallskip
-- l'orbite du point $[\frac{1}{2}\sum_{i=1}^{16}\varepsilon_{i}]$ (on observera que l'on a $\mathrm{q}(\frac{1}{2}\sum_{i=1}^{16}\varepsilon_{i})=2)$, cette orbite a $2^{14}=16384$ éléments.

\medskip
Toujours d'après le point (a) de 1.14, le réseau $\mathrm{vois}_{2}(\mathrm{E}_{16};c)$ est isomorphe (comme $\mathrm{q}$-module) à $\mathrm{E}_{16}$ pour tout $c$ qui n'est pas dans la troisième orbite.

\medskip
On conclut comme précedemment~: comme le graphe des $2$-voisins est connexe, on a forcément $\mathrm{N}_{2}(\mathrm{E}_{16},\mathrm{E}_{8}\oplus\mathrm{E}_{8})=12870$.

\bigskip
\textbf{3.2.} Détermination de $\mathrm{T}_{3}$ pour $n=16$

\medskip
On plonge toujours $\mathrm{E}_{16}$ de façon standard dans $\mathbb{Q}^{16}$. On constate que l'action de $\mathrm{O}(\mathrm{E}_{16})$ partitionne la quadrique $\mathrm{C}_{\mathrm{E}_{16}}(\mathbb{F}_{3})$ en cinq orbites, à savoir les orbites des classes des vecteurs suivants de $\mathrm{E}_{16}$ (qui appartiennent en fait à $\mathrm{D}_{16}$)~:

\smallskip
-- $u_{1}=2\hspace{1pt}\varepsilon_{1}+\varepsilon_{2}+\varepsilon_{3}$~;

\smallskip
-- $u_{2}=\varepsilon_{1}+\varepsilon_{2}+\ldots+\varepsilon_{6}$~;

\smallskip
-- $u_{3}=2\hspace{1pt}\varepsilon_{1}+\varepsilon_{2}+\varepsilon_{3}+\ldots+\varepsilon_{9}$~;

\smallskip
-- $u_{4}=\varepsilon_{1}+\varepsilon_{2}+\ldots+\varepsilon_{12}$~;

\smallskip
-- $u_{5}=2\hspace{1pt}\varepsilon_{1}+\varepsilon_{2}+\varepsilon_{3}+\ldots+\varepsilon_{15}$.

\medskip
On constate que l'on a $(\mathrm{q}(u_{i}))_{i=1,2,\ldots,5}=(3,3,6,6,9)$. D'après 1.14 (a), le réseau $\mathrm{vois}_{3}(\mathrm{E}_{16},[u_{i}])$ est isomorphe (comme $\mathrm{q}$-module) à $\mathrm{E}_{16}$ pour $i=1,2$. On constate que le cardinal de l'orbite de $[u_{i}]$ est ${{16}\choose{3\hspace{0,5pt}i}}\hspace{1pt}2^{\hspace{1pt}3\hspace{0,5pt}i-1}$. Ce cardinal n'est pas divisible par $286$ pour $i=4,5$~; d'après 1.8 le réseau $\mathrm{vois}_{3}(\mathrm{E}_{16};[u_{i}])$ est encore isomorphe (comme $\mathrm{q}$-module) à $\mathrm{E}_{16}$ pour $i=4,5$. En effet, on a
$$
\frac{\vert\mathrm{O}(\mathrm{E}_{8}\oplus\mathrm{E}_{8})\vert}{\vert\mathrm{O}(\mathrm{E}_{16})\vert}
\hspace{4pt}=\hspace{4pt}
\frac{405}{286}
$$
avec $405$ et $286$ premiers entre eux~; le scholie 1.8 montre que si $\mathrm{vois}_{3}(\mathrm{E}_{16};[u_{i}])$ est isomorphe à $\mathrm{E}_{8}\oplus\mathrm{E}_{8}$ alors le cardinal de l'orbite de $[u_{i}]$ est divisible par~$286$.

\medskip
Puisque le graphe des $3$-voisins est connexe, le réseau $\mathrm{vois}_{3}(\mathrm{E}_{16};[u_{3}])$ doit être isomorphe (comme $\mathrm{q}$-module) à $\mathrm{E}_{8}\oplus\mathrm{E}_{8}$ (ce que le logiciel \texttt{PARI} confirme)  et on a $\mathrm{N}_{3}(\mathrm{E}_{16},\mathrm{E}_{8}\oplus\mathrm{E}_{8})={{16}\choose{9}}\hspace{1pt}2^{\hspace{1pt}8}=2928640$. Compte tenu de 1.7 et 2.2, cette égalité détermine l'opérateur de Hecke $\mathrm{T}_{3}:\mathbb{Z}[\mathrm{X}_{16}]\to\mathbb{Z}[\mathrm{X}_{16}]$~; sa matrice dans la base $(\mathrm{E}_{16},\mathrm{E}_{8}\oplus\mathrm{E}_{8})$, disons encore $\mathrm{T}_{3}$, est la suivante~:
$$
\hspace{24pt}
\mathrm{T}_{3}=
\begin{bmatrix}
4248000 & 4147200 \\  2928640 & 3029440
\end{bmatrix}
\hspace{24pt}.
$$

\bigskip
\textbf{3.3.} Détermination de $\mathrm{T}_{2}$ pour $n=24$ (d'après Nebe-Venkov \cite{nebevenkov})

\bigskip
Soit $n>0$ un entier divisible par $8$.

\medskip
Nous avons expliqué au paragraphe 1, en suivant Borcherds, pourquoi l'ensemble $\mathrm{Y}_{n}(2)$ des classes d'isomorphisme des couples $(L_{1},L_{2})$, $L_{1}$ désignant un réseau unimodulaire pair de dimension~$n$ et $L_{2}$ un $2$-voisin de $L_{1}$ dans $\mathbb{Q}\otimes_{\mathbb{Z}}L_{1}$, est naturellement en bijection avec l'ensemble $\widetilde{\mathrm{B}}_{n}$ des classes d'isomorphismes des couples $(L;\omega)$, $L$ désignant un réseau unimodulaire impair de dimension $n$ et $\omega$ une bijection de l'ensemble des deux classes de vecteurs de Wu de $L$ sur $\{1,2\}$. Rappelons que nous avons noté $\mathrm{B}_{n}$ l'ensemble des classes d'isomorphisme de réseaux unimodulaires $L$ impairs de dimension~$n$~; nous notons en outre $\mathrm{B}_{n}^{1}$ le sous-ensemble de $\mathrm{B}_{n}$ constitué des $[L]$ telles que $L$ représente $1$ (en clair, tels qu'il existe $e$ dans $L$ avec $e.e=1$) et $\mathrm{B}_{n}^{2}$ le complémentaire $\mathrm{B}_{n}-\mathrm{B}_{n}^{1}$. Dans \cite{borcherdsthese}, Borcherds utilise la bijection $\mathrm{Y}_{24}(2)\cong\widetilde{\mathrm{B}}_{24}$ pour déterminer $\mathrm{B}_{24}$. Il explicite dans \cite[Chap. 17]{conwaysloane} la liste des $156$ éléments $b$ de $\mathrm{B}_{24}^{2}$ et donne, pour chacun de ces $b$, suffisamment d'informations pour déterminer $\vert\mathrm{O}(b)\vert$. La raison pour laquelle il se limite à $\mathrm{B}_{24}^{2}$ est qu'un réseau $L$ qui représente $1$ est isomorphe à une somme orthogonale $\mathrm{I}_{1}\oplus L'$ et que l'énumération des réseaux unimodulaires, ne représentant pas $1$, de dimension strictement inférieure à $23$, avait déjà été faite (voir \cite[Chap. 16, Table~16.7]{conwaysloane}, $\mathrm{B}_{24}^{1}$ a $117$ éléments). On observera, avec Nebe et Venkov, que si un réseau unimodulaire impair  $L$, de dimension $n$, représente $1$, alors les deux réseaux unimodulaires pairs $2$-voisins de $L$ sont isomorphes d'après 1.16. Soit $\mathrm{B}_{n}^{2,0}$ le sous-ensemble de $\mathrm{B}_{n}^{2}$ constitué des classes d'isomorphisme des réseaux unimodulaires impairs $L$, de dimension $n$, tels que les deux réseaux unimodulaires pairs $2$-voisins de $L$, disons $L_{1}$ et~$L_{2}$, sont non isomorphes~; soit $\mathrm{e}$ l'application de~$\mathrm{B}_{n}^{2,0}$ dans l'ensemble des paires d'éléments de $\mathrm{X}_{n}$, qui associe à $[L]$ la paire $\{[L_{1}],[L_{2}]\}$. La proposition 1.11 se spécialise de la façon suivante~:

\bigskip
\textbf{Proposition 3.3.1.} {\em Soient $x_{1}$ et $x_{2}$ deux éléments distincts de $\mathrm{X}_{n}$. On a~:
$$
\hspace{24pt}
\mathrm{N}_{2}(x_{1},x_{2})
\hspace{4pt}=\hspace{4pt}
\sum_{b\hspace{2pt}\in\hspace{2pt}\mathrm{e}^{-1}(\{x_{1},x_{2}\})}\hspace{8pt}
\frac{\vert\mathrm{O}(x_{1})\vert}{\vert\mathrm{O}(b)\vert}
\hspace{24pt}.
$$}

\bigskip
Nebe et Venkov déterminent $\mathrm{T}_{2}$ à l'aide de l'énoncé ci-dessus et de la table de Borcherds, en tenant compte de 2.2 ; attention, nos conventions font que la matrice $(24,24)$ de \cite[page 59]{nebevenkov} est la transposée de notre $\mathrm{T}_{2}$.

\bigskip
\textit{Remarques}

\bigskip
-- La proposition 1.11 donne pour $n=8$ la relation
$$
\hspace{24pt}
\frac{\vert\mathrm{O}(\mathrm{E}_{8})\vert}{\vert\mathrm{O}(\mathrm{I}_{8})\vert}
\hspace{4pt}=\hspace{4pt}
\frac{\mathrm{c}_{8}(2)}{2}
\hspace{24pt}.
$$
Plus généralement, la proposition 1.11 conduit, pour tout $n$ divisible par $8$, à la relation suivante entre formules de masses
$$
\sum_{b\in\mathrm{B}_{n}}\frac{1}{\vert\mathrm{O}(b)\vert}
\hspace{4pt}=\hspace{4pt}
\frac{\mathrm{c}_{n}(2)}{2}
\sum_{x\in\mathrm{X}_{n}}\frac{1}{\vert\mathrm{O}(x)\vert}
\leqno{(*)}
$$
(voir \cite[Chap. 16, \S 2]{conwaysloane}). Expliquons pourquoi. La proposition 1.11 dit que l'on a
$$
\mathrm{N}_{2}(x,y)\hspace{4pt}\frac{1}{\vert\mathrm{O}(x)\vert}
\hspace{4pt}=\hspace{4pt}
\sum_{\mathrm{d}_{1}(\beta)=x\hspace{2pt},\hspace{2pt}\mathrm{d}_{2}(\beta)=y}\hspace{8pt}
\frac{1}{\vert\mathrm{O}(\beta)\vert}
$$
pour tous $x$ et $y$ dans $\mathrm{X}_{n}$~; en sommant en $y$ puis en $x$ on obtient
$$
\hspace{24pt}
\sum_{\beta\in\widetilde{\mathrm{B}}_{n}}\frac{1}{\vert\mathrm{O}(\beta)\vert}
\hspace{4pt}=\hspace{4pt}
\mathrm{c}_{n}(2)
\sum_{x\in\mathrm{X}_{n}}\frac{1}{\vert\mathrm{O}(x)\vert}
\hspace{24pt}.
\leqno{(**)}
$$
Soit $\mathrm{p}:\widetilde{\mathrm{B}}_{n}\to\mathrm{B}_{n}$ l'application évidente~; l'égalité
$$
\sum_{\beta\in\mathrm{p}^{-1}(b)}\frac{1}{\vert\mathrm{O}(\beta)\vert}
\hspace{4pt}=\hspace{4pt}
\frac{2}{\vert\mathrm{O}(b)\vert}
$$
montre que les relations $(*)$ et $(**)$ sont équivalentes.

\bigskip
-- L'analyse que nous avons faite en 3.1 des $2$-voisinages entre réseaux unimodulaires pairs de dimension $16$ conduit à l'énoncé ci-dessous que nous utiliserons dans l'appendice A. Comme en 3.1, on identifie de façon standard les réseaux $\mathrm{E}_{16}$, $\mathrm{E}_{8}\oplus\mathrm{E}_{8}$ et $\mathrm{D}_{8}\oplus\mathrm{D}_{8}$ à des réseaux de $\mathbb{Q}^{16}$ dont on note $(\epsilon_{1},\epsilon_{2},\ldots,\epsilon_{16})$ la base canonique.

\bigskip
\textbf{Scholie-Définition 3.3.2.} {\em Le réseau de $\mathbb{Q}^{16}$ engendré par
$$
\mathrm{D}_{8}\oplus\mathrm{D}_{8}
\hspace{24pt},\hspace{24pt}
\frac{1}{2}\sum_{i=1}^{16}\epsilon_{i}
\hspace{24pt},\hspace{24pt}
-\epsilon_{1}-\epsilon_{9}
+\frac{1}{2}\sum_{i=1}^{8}\epsilon_{i}
$$
est un réseau unimodulaire impair que l'on note $\mathrm{Bor}_{16}$~; ses $2$-voisins pairs sont $\mathrm{E}_{16}$ et $\mathrm{E}_{8}\oplus\mathrm{E}_{8}$. Le réseau $\mathrm{Bor}_{16}$ est, à isomorphisme près, le seul réseau unimodulaire impair de dimension $16$ qui ne représente pas $1$. Le système de racines $\mathrm{R}(\mathrm{Bor}_{16})$ est isomorphe à $\mathbf{D}_{8}\coprod\mathbf{D}_{8}$.}

\bigskip
\textit{Démonstration.}  Soit $B$ un réseau unimodulaire impair de dimension $16$ dont les $2$-voisins pairs sont isomorphes~; l'analyse évoquée ci-dessus et la proposition 1.17 montrent que $B$ représente $1$. Soit maintenant $B$ un réseau unimodulaire impair de dimension $16$ qui ne représente pas 1~; d'après ce qui précède les $2$-voisins pairs  de $B$ sont isomorphes, l'un à $\mathrm{E}_{16}$, l'autre à $\mathrm{E}_{8}\oplus\mathrm{E}_{8}$. Soit $M$ le sous-module de $B$ constitué des éléments $x$ avec $x.x$ pair~; d'après l'analyse en question $M$ est isomorphe (comme $\widetilde{\mathrm{q}}$-module) à la fois, au sous-module de $\mathrm{E}_{16}$ orthogonal modulo $2$ de l'élément $\sum_{i=1}^{8}\epsilon_{i}$ et au sous-module de $\mathrm{E}_{8}\oplus\mathrm{E}_{8}$ orthogonal modulo $2$ de l'élément $2\epsilon_{1}+2\epsilon_{9}$. On constate que ces deux orthogonaux modulo $2$, vus comme des réseaux de $\mathbb{Q}^{16}$, coïncident avec le réseau, disons $\mathrm{M}_{16}$, engendré par $\mathrm{D}_{8}\oplus\mathrm{D}_{8}$ et $\frac{1}{2}\sum_{i=1}^{16}\epsilon_{i}$~; en effet $\mathrm{M}_{16}$ est d'indice $2$ à la fois dans $\mathrm{E}_{16}$ et $\mathrm{E}_{8}\oplus\mathrm{E}_{8}$ et l'on a $(\sum_{i=1}^{8}\epsilon_{i}).x\equiv 0\bmod{2}$ (resp. $(2\epsilon_{1}+2\epsilon_{9}).x\equiv 0\bmod{2}$) pour tout $x$ dans $\mathrm{M}_{16}$. Incidemment, ceci montre que l'on a $\mathrm{M}_{16}=\mathrm{E}_{16}\cap(\mathrm{E}_{8}\oplus\mathrm{E}_{8})$. Soit $\xi$ l'élément $-\epsilon_{1}-\epsilon_{9}+\frac{1}{2}\sum_{i=1}^{8}\epsilon_{i}$ de $\mathbb{Q}^{16}$~; on constate que $\xi$ appartient à $\mathrm{M}_{16}^{\sharp}$ et que l'on a $\xi.\xi=3$, il en résulte que $\xi$ engendre la droite ``non quadratiquement isotrope'' de $\mathop{\mathrm{r\acute{e}s}}\mathrm{M}_{16}$, si bien que $\mathrm{Bor}_{16}$ est le réseau unimodulaire impair correspondant à $\mathrm{M}_{16}$ par la théorie de Borcherds. D'après cette théorie $M\simeq\mathrm{M}_{16}$ implique $B\simeq\mathrm{Bor}_{16}$ ($\mathrm{Bor}_{16}$ ne représente pas $1$ car  ses $2$-voisins pairs ne sont pas isomorphes). La dernière assertion du scholie se démontre à l'aide de II.3.7 (observer que l'on a $\mathrm{R}(\mathrm{M}_{16})=\mathrm{R}(\mathrm{E}_{16})\cap\mathrm{R}(\mathrm{E}_{8}\oplus\mathrm{E}_{8})\simeq\mathbf{D}_{8}\coprod\mathbf{D}_{8}$ et que l'image par la fonction $\mathrm{qm}$ de la classe de $\xi$ dans $\mathop{\mathrm{r\acute{e}s}}\mathrm{M}_{16}$ est forcément $\frac{3}{2}$).
\hfill$\square$

\section{Sur les $d$-voisinages entre un réseau de Niemeier avec racines et le réseau de Leech}

La raison d'être de ce paragraphe est la suivante~:
 
\smallskip
Soit $d\geq 2$ un entier~; la détermination des entiers $\mathrm{N}_{d}(L,\mathrm{Leech})$, pour tout réseau de Niemeier avec racines $L$, entraîne celle de l'opérateur de Hecke $\mathrm{T}_{d}$ (dans l'expression $\mathrm{N}_{d}(L,\mathrm{Leech})$, ``Leech'' est une abréviation pour ``réseau de Leech'', celle-ci sera fréquemment utilisée dans la suite du mémoire).

\medskip
Expliquons pourquoi. On a vu en 3.3, que l'on connaît explicitement, grâce à Nebe-Venkov, l'opé\-rateur de Hecke $\mathrm{T}_{2}:\mathbb{Z}[\mathrm{X}_{24}]\to\mathbb{Z}[\mathrm{X}_{24}]$~; on peut vérifier (merci \texttt{PARI}) que les éléments $\mathrm{T}_{2}^{k}\hspace{1pt}[\mathrm{Leech}]$, $0\leq k\leq 23$, sont linéairement indépendants. Comme les opérateurs de Hecke $\mathrm{T}_{d}$ et $\mathrm{T}_{2}$ commutent, la déter\-mination de $\mathrm{T}_{d}\hspace{1pt}[\mathrm{Leech}]$ entraîne celle de $\mathrm{T}_{d}$. Compte tenu de 1.7 et 2.2, la détermination de $\mathrm{T}_{d}\hspace{1pt}[\mathrm{Leech}]$ est quant à elle équivalente à celle des entiers $\mathrm{N}_{d}(x,[\mathrm{Leech}])$, $x$ décrivant l'ensemble $\mathrm{X}_{24}-\{[\mathrm{Leech}]\}$.

\vspace{0,75cm}
\textbf{4.1.} Conditions nécessaires pour qu'un réseau de Niemeier avec racines  puisse avoir un $d$-voisin sans racines

\bigskip
\textbf{Proposition 4.1.1.} {\em Soient $L$ un réseau de Niemeier avec racines et $d\geq 2$ un entier~; soit $\mathrm{h}(L)$ le nombre de Coxeter de $L$ (voir Proposition-Définition II.3.3). Si $L$ a un $d$-voisin sans racines alors on a l'inégalité
$$
\hspace{24pt}
d
\hspace{4pt}\geq\hspace{4pt}
\mathrm{h}(L)
\hspace{24pt}.
$$}

\bigskip
\textit{Démonstration.} On suppose qu'il existe un élément $c$ de $\mathrm{C}_{L}(\mathbb{Z}/d)$ tel que l'on a $\mathrm{R}(\mathrm{vois}_{d}(L;c))=\emptyset$. On a \textit{a fortiori} $\mathrm{R}(\mathrm{M}_{d}(L;c))=\emptyset$ (on rappelle que $\mathrm{M}_{d}(L;c)$ est l'intersection dans $\mathbb{Q}\otimes_{\mathbb{Z}}L$ des réseaux $L$ et $\mathrm{vois}_{d}(L;c)$) ou encore $\mathrm{R}(L)\cap\mathrm{M}_{d}(L;c)=\emptyset$. Soit $u$ un élément de $L$ qui représente $c$, la condition $\mathrm{R}(\mathrm{vois}_{d}(L;c))=\emptyset$ implique donc
$$
\hspace{24pt}
\alpha.u
\hspace{4pt}\not\equiv\hspace{4pt}
0 \pmod{d} \hspace{6pt}\text{pour tout $\alpha$ dans $\mathrm{R}(L)$}
\hspace{24pt}.
$$
On constate que l'on a nécessairement l'inégalité $d\geq\mathrm{h}(L)$ en prenant dans la proposition ci-dessous, pour $R$ une composante irréductible de $\mathrm{R}(L)$, et pour $f$ la forme linéaire $x\mapsto x.u$.
\hfill$\square$

\bigskip
 \textbf{Proposition 4.1.2.} {\em Soient $V$ un $\mathbb{R}$-espace vectoriel de dimension finie et $R\subset V$ un système de racines irréductible et réduit~; soit $h$ le nombre de Coxeter de $R$. Soit $f:V\to\mathbb{R}$ une forme linéaire dont la restriction à $R$ prend des valeurs entières~; soit $d\geq 2$ un entier avec $d<h$. Alors il existe une racine $\alpha$ dans $R$ telle que l'on a
 $$
 \hspace{24pt}
 f(\alpha)\equiv 0 \pmod{d}
  \hspace{24pt}.
 $$}
 
 \textit{Démonstration.} On fixe une chambre $C$ du système de racines $R$~; on note respectivement $\{\alpha_{1},\alpha_{2},\ldots,\alpha_{l}\}$ et $\widetilde{\alpha}$, la base de $R$ et la plus grande racine correspondantes. On rappelle que l'on a $\widetilde{\alpha}=n_{1}\alpha_{1}+n_{2}\alpha_{2}+\dots+n_{l}\alpha_{l}$ avec $n_{i}\in\mathbb{N}-\{0\}$ pour $i=1,2,\ldots,l$ et
$$
 n_{1}+n_{2}+\dots+n_{l}=h-1\leqno{(\mathrm{hm})}
$$
\cite[Ch. VI, \S1, Prop. 31]{bourbaki} ($\mathrm{hm}$ est pour ``hauteur maximale''). On note enfin $Alc$ l'alcôve de $V^{*}$ contenue dans $C$ à laquelle $0$ est adhérent \cite[Ch. VI, \S2, Prop. 4]{bourbaki}~; $Alc$ (resp. $\overline{Alc}$) est donc l'ouvert (resp. le fermé) de $V^{*}$ constitué des éléments $\phi$ vérifiant les inégalités $\langle\alpha_{i},\phi\rangle>0$ (resp. $\langle\alpha_{i},\phi\rangle\geq 0$) pour  $i=1,2,\ldots,l$ et $\langle\widetilde{\alpha},\phi\rangle<1$ (resp. $\langle\widetilde{\alpha},\phi\rangle\leq 1$).
 
\medskip
Soit $\phi$ un élément de $V^{*}$. Puisque $\overline{Alc}$ est un domaine fondamental pour l'action du groupe de Weyl affine sur $V^{*}$ (voir par exemple \cite[Ch. VI, \S2, \no 1 et \no 2]{bourbaki}), il existe un élément $w$ du groupe de Weyl de $R$ et un élément $\theta$ du réseau $\mathrm{Q}(R^{\vee})$ de $V^{*}$ tels que l'on a les inégalités
$$
\hspace{24pt}
\langle\hspace{1pt}w\alpha_{i}\hspace{1pt},\hspace{1pt}\phi-\theta\hspace{1pt}\rangle\geq 0
\hspace{8pt}\text{pour}\hspace{4pt}i=1,2,\ldots,l
\hspace{24pt}\text{et}\hspace{24pt}
\langle\hspace{1pt}w\widetilde{\alpha}\hspace{1pt},\hspace{1pt}\phi-\theta\hspace{1pt}\rangle\leq 1
\hspace{24pt}.
$$

\medskip
On obtient une preuve de la proposition en prenant $\phi=\frac{1}{d}f$~: il existe $w$ et $\theta$ comme ci-dessus tels que l'on a les inégalités
$$
\hspace{24pt}
f(w\alpha_{i})-d\hspace{1pt}\theta(w\alpha_{i})\geq 0
\hspace{24pt}\text{et}\hspace{24pt}
f(w\widetilde{\alpha})-d\hspace{1pt}\theta(w\widetilde{\alpha})\leq d
\hspace{24pt}.
$$

\smallskip
On observe que les $w\alpha_{i}$ et $w\widetilde{\alpha}$ sont des racines et que les  $f(w\alpha_{i})$, $\theta(w\alpha_{i})$, $f(w\widetilde{\alpha})$ et $\theta(w\widetilde{\alpha})$ sont des entiers~; on pose $x_{i}=f(w\alpha_{i})-d\hspace{1pt}\theta(w\alpha_{i})$ et\linebreak $y=d-(f(w\widetilde{\alpha})-d\hspace{1pt}\theta(w\widetilde{\alpha}))$. On a $x_{i}\geq 0$, $y\geq 0$ et
$$
\hspace{24pt}
 n_{1}\hspace{1pt}x_{1}+n_{2}\hspace{1pt}x_{2}+\dots+n_{l}\hspace{1pt}x_{l}+y=d
 \hspace{24pt}.
 $$
 Compte tenu de l'égalité $(\mathrm{hm})$,  l'un des entiers $x_{1}$, $x_{2}$, \ldots, $x_{l}$, $y$ doit être nul, ce qui démontre bien la proposition.
 \hfill$\square$
 
 \bigskip
\textit{Remarque.} En fait, nous n'utilisons la proposition 4.1.2 que pour des systèmes de racines irréductibles de type ADE, c'est à dire pour $R=\mathbf{A}_{n}$ ($n\geq 1$), $R=\mathbf{D}_{n}$ ($n\geq 3$), $R=\mathbf{E}_{6}$, $R=\mathbf{E}_{7}$, et $R=\mathbf{E}_{8}$. Il existe des démonstrations très simples de l'énoncé 4.1.2 dans les deux premiers cas. Nous traitons ci-dessous le second~; le traitement du premier est similaire (et d'ailleurs plus facile). 

\medskip
On munit $\mathbb{R}^{n}$ de sa structure euclidienne canonique~; on rappelle que l'on~a $\mathbf{D}_{n}=\mathrm{R}(\mathrm{D}_{n})$, $\mathrm{D}_{n}$ désignant le sous-module de $\mathbb{Z}^{n}$ constitué des $n$-uples $(x_{1},x_{2},\ldots,x_{n})$ avec $\sum_{i=1}^{n}x_{i}$ pair. Une forme linéaire $f:\mathbb{R}^{n}\to\mathbb{R}$ qui prend des valeurs entières sur $\mathrm{D}_{n}$ induit une application linéaire $\mathbb{Z}^{n}\to\frac{1}{2}\mathbb{Z}$. On pose $\lambda=0$ (resp. $\lambda=1$)  si l'homomorphisme composé $\mathrm{D}_{n}\overset{f}{\to}\mathbb{Z}\to\mathbb{Z}/2$ est trivial (resp. non trivial)~; on pose $\nu_{i}=f(\epsilon_{i})-\frac{\lambda}{2}$, $\varepsilon_{1},\varepsilon_{2},\ldots,\varepsilon_{n}$ désignant la base canonique de $\mathbb{R}^{n}$ (les $\nu_{i}$ sont donc des entiers). On constate que si l'on a $f(\alpha)\not\equiv 0\pmod{d}$ pour tout $\alpha$ dans $\mathbf{D}_{n}$ alors l'application de $\{1,2,\ldots,n\}$ dans $\mathbb{Z}/d$ qui associe à $i$ la classe modulo $d$ de $\nu_{i}$ induit une injection de $\{1,2,\ldots,n\}$ dans l'ensemble quotient de $\mathbb{Z}/d$ par l'involution $t\mapsto -t-\lambda$ (se rappeler que $\mathbf{D}_{n}$ est constitué des éléments $\pm\varepsilon_{i}\pm\varepsilon_{j}$, $i\not=j$). Or le cardinal du quotient en question est majoré par $\frac{d}{2}+1$, on a donc $n\leq\frac{d}{2}+1$ soit encore $d\geq 2n-2$.
\hfill$\square$

\bigskip
\textit{Notations relatives aux réseaux de Niemeier avec racines}
 
 \medskip
 Comme tous les énoncés de la suite de ce paragraphe 4  concernent les réseaux de Niemeier avec racines, nous récapitulons et complétons ci-dessous la liste des notations que nous utilisons pour ces réseaux.
 
  \medskip
 Soit $L$ un réseau de Niemeier avec racines~:
 
 \smallskip
 --\hspace{8pt}$R=\mathrm{R}(L)$ désigne l'ensemble des racines de $L$~;
 
 \smallskip
 --\hspace{8pt}on pose $V=\mathbb{R}\otimes_{\mathbb{Z}}L$, $V$ est un espace euclidien de dimension $24$ et $R\subset V$ est un système de racines de type ADE, equicoxeter, de rang $24$~;

 \smallskip
 --\hspace{8pt}$W=\mathrm{W}(L)$ désigne le groupe de Weyl du système de racines $R$, on rappelle que $W$ s'identifie à un sous-groupe du groupe orthogonal $\mathrm{O}(L)$ (point (a) du scholie II.3.15)~;
 
 \smallskip
 --\hspace{8pt}$Q=\mathrm{Q}(R)\subset L$ est le réseau de $V$ engendré par $R$~;

\smallskip
 --\hspace{8pt}$Q^{\sharp}=(\mathrm{Q}(R))^{\sharp}\supset L$ est le réseau dual de $Q$~;
 
 \smallskip
 --\hspace{8pt}$f=\mathrm{f}(L)$ désigne l'indice de $Q$ dans $Q^{\sharp}$ (que Bourbaki appelle {\em indice de connexion} du système de racines $R$)~;
 
 \smallskip
 --\hspace{8pt}$g=\mathrm{g}(L)$ désigne l'indice de $Q$ dans $L$ (compte tenu du fait que $L/Q$ est un lagrangien du $\mathrm{qe}$-module $Q^{\sharp}/Q$, on a $f=g^{2}$).
 
\medskip
On note $R_{1},R_{2},\ldots,R_{c}$ les composantes irréductibles de $R$~:

\smallskip
--\hspace{8pt}le nombre de ces composantes irréductibles est donc $c=\mathrm{c}(L)$~;

\smallskip
--\hspace{8pt}toutes ces composantes irréductibles ont le même nombre de Coxeter, à savoir $h=\mathrm{h}(L)$.

\medskip
On choisit une chambre $C$ de $R$~:

\smallskip
--\hspace{8pt}$B\subset R$ désigne la base de $R$ correspondant à ce choix~;

\smallskip
--\hspace{8pt}$R_{+}\subset R$ est l'ensemble des racines positives pour la relation d'ordre sur~$V$ définie par $C$~;

\smallskip
--\hspace{8pt}$H\subset R_{+}\subset R$ désigne l'ensemble des éléments maximaux de $R$ pour l'ordre en question (l'ensemble $H$ a $c$ éléments, plus précisément on a $H\cap R_{k}=\{\widetilde{\alpha}_{k}\}$, $\widetilde{\alpha}_{k}$ désignant  la plus grande racine de $R_{k}$).

\medskip
On note $Alc$ l'alcôve de $V$ contenue dans $C$ à laquelle $0$ est adhérent~; $Alc$ (resp. $\overline{Alc}$) est l'ouvert (resp. le fermé) de $V$ constitué des éléments $x$ vérifiant les inégalités $\alpha.x>0$ (resp. $\alpha.x\geq 0$) pour $\alpha\in B$ et $\alpha.x<1$ (resp. $\alpha.x\leq 1$) pour $\alpha\in H$.

\medskip
On pose enfin $\Pi: =Q^{\sharp}\cap\overline{Alc}$. Le sous-ensemble $\Pi$ de $Q^{\sharp}$ s'identifie au produit d'ensembles $\prod_{k=1}^{c}\Pi(R_{i})$ (la notation $\Pi(-)$ est introduite dans la proposition II.3.8, l'égalité $\Pi(S) =\mathrm{Q}(S)^{\sharp}\cap\overline{Alc}$ pour $S$ un sytème de racines irréductible de type ADE est établie dans la démonstration de II.3.8 (a)).

\bigskip
La démonstration de la proposition 4.1.2 conduit au scholie suivant~:
 
 \bigskip
 \textbf{Scholie 4.1.3.} {\em Soit $L$ un réseau de Niemeier avec racines. Soient $\xi$ un élément de $Q^{\sharp}$ et $d\geq 1$ un entier. Alors il existe un élément $w$ de $W$ et un élément $x$ de $Q$ tels que l'élément $\eta:=w\hspace{1pt}\xi+d\hspace{1pt}x$ de $Q^{\sharp}$ appartient à $d\hspace{2pt}\overline{Alc}$, en clair, tels que l'on a les inégalités  $\alpha.\eta\geq 0$ pour $\alpha\in B$ et $\alpha.\eta\leq d$ pour $\alpha\in H$. De plus si l'on a $\alpha.\xi\not\equiv 0\bmod{d}$ pour tout $\alpha$ dans $R$ alors le couple $(w,x)$ qui apparaît ci-dessus est uniquement déterminé en fonction de $\xi$.}

\bigskip
Soit $x$ un élément de $V$~; on dit que $x$ est {\em régulier} si l'on a $\alpha.x\not=0$ pour tout~$\alpha$ dans $R$ (en d'autres termes, si $x$ est dans une chambre). Soit~$d\geq 1$ un entier~; nous dirons que $x$ de $V$ est {\em $d$-régulier} si l'on a $\alpha.x\not\in d\mathbb{Z}$ pour tout $\alpha$ dans~$R$ (en d'autres termes, si $\frac{1}{d}\hspace{1pt}x$ est dans une alcôve). Soit~$d\geq 2$ un entier~; nous dirons qu'un élément de $\mathrm{P}_{L}(\mathbb{Z}/d)$ est {\em régulier} s'il est représenté par un élément $d$-régulier $u$ de $L$, en clair vérifiant $\alpha.u\not\equiv 0\bmod{d}$ pour tout~$\alpha$ dans $R$ (cette condition est indépendante du choix de $u$)~; le sous-ensemble de $\mathrm{P}_{L}(\mathbb{Z}/d)$ constitué de ces éléments sera noté $\mathrm{P}^{\mathrm{r\acute{e}g}}_{L}(\mathbb{Z}/d)$. On posera enfin $\mathrm{C}^{\mathrm{r\acute{e}g}}_{L}(\mathbb{Z}/d):=\mathrm{C}_{L}(\mathbb{Z}/d)\cap\mathrm{P}^{\mathrm{r\acute{e}g}}_{L}(\mathbb{Z}/d)$. La démonstration que nous avons donnée de la proposition 4.1.1 revient en fait à vérifier l'énoncé plus précis suivant (le point (a) étant évident)~:

\bigskip
\textbf{Scholie 4.1.4.} {\em Soient $L$ un réseau de Niemeier avec racines et $d\geq 2$ un entier.

\medskip
{\em (a)} Soit $c$ un élément de $\mathrm{C}_{L}(\mathbb{Z}/d)$~; si le réseau $\mathrm{vois}_{d}(L;c)$ est sans racines alors $c$ appartient à $\mathrm{C}^{\mathrm{r\acute{e}g}}_{L}(\mathbb{Z}/d)$.

\medskip
{\em (b)} Si l'ensemble $\mathrm{P}^{\mathrm{r\acute{e}g}}_{L}(\mathbb{Z}/d)$ est non vide alors on a l'inégalité
$d\geq\mathrm{h}(L)$.}

\vspace{0,75cm}
\textbf{4.2.} Sur les $h$-voisinages et $h+1$-voisinages entre un réseau de Niemeier avec racines de nombre de Coxeter $h$ et le réseau de Leech
 
 \bigskip
Soit $L$ un réseau de Niemeier avec racines, de nombre de Coxeter $h$. On a vu dans le sous-paragraphe 4.1 qu'une condition nécessaire pour que $L$ ait pour $d$-voisin le réseau de Leech est l'inégalité $d\geq h$. On montre en particulier ci-après que cette inégalité est optimale~; ceci est intimement relié aux \textit{holy constructions} du réseau de Leech dues à Conway et Sloane \cite{conwaysloane23leech}. On détermine également l'entier $\mathrm{N}_{d}(L,\mathrm{Leech})$ pour $d=h,h+1$~; ce calcul sera exploité en \ref{calcfinvpgenre2}.
 
\bigskip
Nous commençons par rappeler ce qu'est un vecteur de Weyl d'un réseau de Niemeier et par rassembler quelques-unes des propriétés de ces vecteurs que nous serons amenés à utiliser.

 \vspace{0,75cm}
 \textsc{Vecteurs de Weyl}
 
 \medskip
Soient $L$ un réseau de Niemeier avec racines et $C\subset V$ une chambre du système de racines $R$. Soit $\rho$ la demi-somme des racines positives (pour la relation d'ordre sur~$V$ définie par $C$)~:
$$
2\hspace{1pt}\rho
\hspace{4pt}=\hspace{4pt}
\sum_{\alpha\in R_{+}}\alpha
$$
(cette égalité montre que $\rho$ appartient à $\frac{1}{2}\hspace{1pt}L$, la proposition 4.2.1 ci-dessous dit que $\rho$ appartient en fait à $L$). On dit que $\rho$ est un {\em vecteur de Weyl} du système de racines $R$ ou du réseau $L$. Soit $\alpha$ une racine de $R$ alors $\alpha$ appartient à $B$ si et seulement si l'on a $\rho.\alpha=1$ (voir \cite[Ch. VI, §1, Prop. 29]{bourbaki})~; cette observation montre que l'application $C\mapsto\rho$ est bijective. Il en résulte que l'action de $W$ sur l'ensemble des vecteurs de Weyl est simplement transitive.

\bigskip
\textbf{Proposition 4.2.1} (Borcherds)\textbf{.} {\em Soient $L$ un réseau de Niemeier avec\linebreak racines et $\rho$ un vecteur de Weyl de $L$. Alors $\rho$ appartient à $L$.}

\bigskip
\textit{Démonstration.} Avant de reproduire l'argument donné par Borcherds dans \cite{borcherdsthese}\cite{borcherdsleech}, nous dégageons quelques énoncés qui nous seront utiles par la suite.

\bigskip
\textbf{Proposition 4.2.2.} {\em Soit $L$ un réseau de Niemeier avec racines.

\medskip
{\em (a)} On a 
$$
h\hspace{2pt}x.y
\hspace{4pt}=\hspace{4pt}
\sum_{\alpha\in R_{+}}(\alpha.x)(\alpha.y)
$$
pour tous $x$ et $y$ dans $V$.

\medskip
{\em (b)}  Soit $\rho$ un vecteur de Weyl de $L$, on a 
$$
h\hspace{1pt}\mathrm{q}(x)-\rho.x
\hspace{4pt}=\hspace{4pt}
\sum_{\alpha\in R_{+}}\frac{(\alpha.x)^{2}-\alpha.x}{2}
$$
pour tout $x$ dans $V$.}
\vfill\eject

\bigskip
\textit{Démonstration.} Le point (a) est un avatar de II.3.4. Le point (b) (dû à Borcherds) est conséquence du point (a) et de la définition de $\rho$. Le point (a) implique l'énoncé suivant~:

\bigskip
\textbf{Corollaire 4.2.3.} {\em Soit $L$ un réseau de Niemeier avec racines. Alors le\linebreak quotient $Q^{\sharp}/Q$ est annulé par $h$.}

\bigskip
En effet, le point (a) de 4.2.2 montre que si $\xi$ et $\eta$ sont deux éléments de $Q^{\sharp}$ alors $h\hspace{1pt}\xi.\eta$ est entier  et donc que $h\hspace{1pt}\xi$ appartient à $Q$.
\hfill$\square$

\bigskip
Soit $\rho$ un vecteur de Weyl de $L$. Comme l'action canonique de $W$ sur $Q^{\sharp}/Q$ est triviale l'image de $\rho$ dans $Q^{\sharp}/Q$ est indépendante du choix de ce vecteur de Weyl. Voici une autre explication de ce phénomène. On fait intervenir la structure de $\mathrm{qe}$-module de $Q^{\sharp}/Q=:\mathop{\mathrm{r\acute{e}s}}Q$. D'après 4.2.3, l'application $Q^{\sharp}/Q\to\mathbb{Q}/\mathbb{Z}\hspace{2pt},\hspace{2pt}\xi\mapsto h\hspace{1pt}\mathrm{q}(\xi)$ est linéaire (et à valeurs dans $\frac{1}{2}\mathbb{Z}/\mathbb{Z}$)~; il existe donc un élément $\sigma$ de $Q^{\sharp}/Q$, uniquement déterminé (et annulé par $2$), tel que l'on a $h\hspace{1pt}\mathrm{q}(\xi)=\sigma.\xi$ pour tout~$\xi$ dans $Q^{\sharp}/Q$. Le point (b) de 4.2.2 montre que $\sigma$ est la classe de~$\rho$. En effet ce point montre que $h\hspace{1pt}\mathrm{q}(\xi)-\rho.\xi$ est entier pour tout $\xi$ dans $Q^{\sharp}$. Nous avons obtenu~:

\medskip
\textbf{Proposition 4.2.4.} {\em Soient $L$ un réseau de Niemeier avec racines et $\rho$ un vecteur de Weyl de $L$. Alors l'image de $\rho$ dans $Q^{\sharp}/Q$, disons $\overline{\rho}$, est caractérisée par la propriété
$$
h\hspace{1pt}\mathrm{q}(\xi)
\hspace{4pt}=\hspace{4pt}
\overline{\rho}.\xi
$$
pour tout $\xi$ dans $Q^{\sharp}/Q$.}

\bigskip
Exprimé avec le formalisme introduit ci-dessus, l'argument de Borcherds pour démontrer 4.2.1 est le suivant~: Soit $I$ un lagrangien du $\mathrm{qe}$-module $Q^{\sharp}/Q$, $\mathrm{q}(I)=0$ implique $\overline{\rho}\in I^{\perp}=I$~; en prenant $I=L/Q$ on obtient $\rho\in L$.
\hfill$\square$

\bigskip
\textbf{Proposition 4.2.5} (Venkov)\textbf{.} {\em Soient $L$ un réseau de Niemeier avec racines et $\rho$ un vecteur de Weyl de $L$. Alors on a
$$
\hspace{24pt}
\mathrm{q}(\rho)
\hspace{4pt}=\hspace{4pt}
h(h+1)
\hspace{24pt}.
$$}

\bigskip
\textit{Démonstration.} Il suffit par exemple de constater que l'on a $\rho.\rho=\frac{n}{12}\hspace{1pt}h(h+1)$ pour tout système de racines de type ADE irréductible de rang $n$ et d'invoquer le point (b) de la proposition-définition II.3.3.
\hfill$\square$

\bigskip
\textbf{Proposition 4.2.6.} {\em Soit $L$ un réseau de Niemeier avec racines. Pour tout élément $x$ de $V$, on a l'inégalité
$$
\hspace{24pt}
\inf_{\alpha\in R}\hspace{4pt}(\alpha.x)^{2}
\hspace{4pt}\leq\hspace{4pt}
\frac{\mathrm{q}(x)}{h(h+1)}
\hspace{24pt}.
$$
De plus, on a l'égalité si et seulement s'il existe un vecteur de Weyl $\rho$ de $L$ et un nombre réel $\lambda\geq 0$ tels que l'on a $x=\lambda\rho$.}

\bigskip
\textit{Démonstration.} Soit $C\subset V$ une chambre de $R$ telle que $x$ appartienne à~$\overline{C}$. Soient respectivement $B$ et $\rho$ la base de $R$ et le vecteur de Weyl associés à~$C$~; soit $\{\varpi_{\alpha}\}_{\alpha\in B}$ la base duale de $B$ (par rapport au produit scalaire). On considère l'égalité
$$
\hspace{24pt}
x.x=\sum_{(\alpha,\beta)\in B\times B}(\alpha.x)(\beta.x)\hspace{4pt}\varpi_{\alpha}.\varpi_{\beta}
\hspace{24pt};
$$
en observant que l'on a $\alpha.x\geq 0$ et $\beta.x\geq 0$ (par définition de $C$), $\varpi_{\alpha}.\varpi_{\beta}\geq 0$ (voir par exemple \cite[Ch. VI, \S1, Th. 2, Rem. 2]{bourbaki}) et $\rho.\rho=\sum_{(\alpha,\beta)}\varpi_{\alpha}.\varpi_{\beta}$ (spécialisation de l'égalité ci-dessus), on constate que l'on a l'inégalité
$$
\hspace{24pt}
x.x\hspace{4pt}\geq\hspace{4pt}
\rho.\rho\hspace{2pt}\inf_{\alpha\in B}(\alpha.x)^{2}
\hspace{24pt}.
\leqno{(*)}
$$
Posons $\lambda=\inf_{\alpha\in B}\alpha.x$~; l'inégalité ci-dessus peut-être affinée en
$$
\hspace{24pt}
x.x\hspace{4pt}\geq\hspace{4pt}
\rho.\rho\hspace{2pt}\inf_{\alpha\in B}(\alpha.x)^{2}
+\sum_{\alpha\in B}\hspace{4pt}(\hspace{1pt}(\alpha.x)^{2}-\lambda^{2}\hspace{1pt})\hspace{4pt}\varpi_{\alpha}.\varpi_{\alpha}
\hspace{24pt}.
$$
Ceci montre que si l'on a égalité dans l'inégalité (*) alors on a $\alpha.x=\lambda$ pour tout $\alpha$ dans $B$ et donc $x=\lambda\rho$.
\hfill$\square$

\bigskip
\textbf{Scholie 4.2.7} {\em Soit $L$ un réseau de Niemeier avec racines. Pour tout élément régulier $\xi$ de $\mathrm{Q}(R)^{\sharp}$, en clair vérifiant $\alpha.\xi\not=0$ pour tout $\alpha$ dans $R$, on a l'inégalité
$$
\hspace{24pt}
\mathrm{q}(\xi)
\hspace{4pt}\geq\hspace{4pt}
h(h+1)
\hspace{24pt}.
$$
De plus on a l'égalité si et seulement $\xi$ est un vecteur de Weyl de $L$.}

\bigskip
Le point (a) de la proposition ci-dessous dit en particulier que le point (b) de 4.1.4 est ``optimal''.

\bigskip
\textbf{Proposition 4.2.8}\textbf{.} {\em Soient $L$ un réseau de Niemeier avec racines et $\rho$ un vecteur de Weyl de $L$.

\medskip
{\em (a)} Le vecteur de Weyl $\rho$ est un élément primitif ({\em a fortiori} $h$-primitif) et $h$-régulier de $L$.

\medskip
{\em (b)} Soit $\xi$ un élément de $Q^{\sharp}$~; les deux conditions suivantes sont équivalentes~:
\begin{itemize}
\smallskip
\item [(i)] l'élément $\xi$ est $h$-régulier~;
\smallskip
\item [(ii)] il existe un élément $w$ de $W$ et un élément $x$ de $Q$ tels que l'on a\linebreak $\xi=w\rho+hx$.
\end{itemize}

\smallskip
De plus si ces conditions sont vérifiées alors le couple $(w,x)$ qui apparaît dans la condition (ii) est uniquement déterminé en fonction de $\xi$.}

\bigskip
\textit{Démonstration du point (a).} L'égalité $\alpha.\rho=1$ pour $\alpha\in B$ montre que $\rho$ est primitif. Avant de démontrer que $\rho$ est $h$-régulier, rappelons la définition et quelques propriétés de la fonction {\em hauteur}, disons $\mathrm{H}:R_{+}\to\mathbb{N}-\{0\}$. Soit $\beta$ un élément de $R_{+}$~; $\beta$ s'écrit $\sum_{\alpha\in B}n_{\alpha}\hspace{1pt}\alpha$ avec $n_{\alpha}$ dans $\mathbb{N}$ \cite[Ch. VI, \S1, Th. 3]{bourbaki} et l'on pose $\mathrm{H}(\beta):=\sum_{\alpha\in B}n_{\alpha}$. La fonction $\mathrm{H}$ vérifie les propriétés suivantes~:

\smallskip
--\hspace{8pt}$\mathrm{H}(\beta)=\beta.\rho$~;

\smallskip
--\hspace{8pt}$\mathrm{H}(\beta)\geq 1$ et $\mathrm{H}(\beta)=1\iff\beta\in B$~;

\smallskip
--\hspace{8pt}$\mathrm{H}(\beta)\leq h-1$ et $\mathrm{H}(\beta)=h-1\iff\beta\in H$.

\smallskip
La dernière propriété ci-dessus résulte de \cite[Ch. VI, \S1, Prop. 31]{bourbaki} (référence que nous avons déjà invoquée lors de la démonstration de 4.1.2) et de la définition même du sous-ensembe $H\subset R_{+}$. L'égalité $R=R_{+}\coprod-R_{+}$ et les inégalités $1\leq\beta.\rho\leq h-1$ pour tout $\beta$ dans~$R_{+}$ montre que l'on a $\alpha.\rho\not\equiv 0\bmod{h}$ pour tout $\alpha$ dans $R$.
\hfill$\square$

\medskip
\textit{Démonstration du point (b).} L'implication $(ii)\Rightarrow(i)$ résulte du fait que $\rho$ est $h$-régulier. On démontre $(i)\Rightarrow(ii)$. Compte tenu de 4.1.3, on peut supposer $\xi\in h\hspace{1pt}\overline{Alc}$, c'est-à-dire~:

\smallskip
--\hspace{8pt}$\alpha.\xi\geq 0$ pour tout $\alpha$ dans $B$~;

\smallskip
--\hspace{8pt}$\widetilde{\alpha}.\xi\leq h$ pour tout $\widetilde{\alpha}$ dans $H$. 

\smallskip
La première inégalité montre que si $\xi$ est $h$-régulier alors on a $\alpha.\xi\geq 1$ pour tout $\alpha$ dans $B$ (par définition $\beta.\xi$ est dans $\mathbb{Z}$ pour tout $\beta$ dans $R$)~; soit encore $\alpha.(\xi-\rho)\geq 0$ pour tout $\alpha$ dans $B$. Pareillement, la seconde montre que l'on a $\widetilde{\alpha}.\xi\leq h-1$ pour tout $\widetilde{\alpha}$ dans $H$~; soit encore $\widetilde{\alpha}.(\xi-\rho)\leq 0$ pour tout $\widetilde{\alpha}$ dans~$H$. Or un élément $\eta$ de $V$ qui vérifie $\alpha.\eta\geq 0$ pour tout $\alpha$ dans $B$ et $\widetilde{\alpha}.\eta\leq 0$ pour tout $\widetilde{\alpha}$ dans $H$ est nul (par exemple, parce qu'il appartient à  $\epsilon\hspace{1pt}\overline{Alc}$ pour tout $\epsilon>0$). On a donc $\xi=\rho$. La dernière partie du point (b), concernant l'unicité du couple $(w,x)$ qui apparaît dans la condition (ii), est un sous-produit de la démonstration que nous venons de faire.
\hfill$\square$

\bigskip
Nous en arrivons au dernier énoncé concernant les vecteurs de Weyl des réseaux de Niemeier que nous voulons dégager~; les points (b) et (c), sont encore dûs à Borcherds.

\bigskip
\textbf{Proposition 4.2.9.}  {\em Soient $L$ un réseau de Niemeier avec racines, $\rho$ un vecteur de Weyl de $L$ et  $\xi$ un élément de $Q^{\sharp}$.

\medskip
{\em (a)} L'élément $\rho-h\hspace{1pt}\xi$ de $Q^{\sharp}$ appartient à $L$ et est $h$-régulier.

\medskip
{\em (b)} On a l'inégalité
$$
\hspace{24pt}
\mathrm{q}(\rho-h\hspace{1pt}\xi)
\hspace{4pt}\geq\hspace{4pt}
\mathrm{q}(\rho)
\hspace{4pt}=\hspace{4pt}
h(h+1)
\hspace{24pt}.
$$

\medskip
{\em (c)} Les conditions suivantes sont équivalentes~:
\begin{itemize}
\smallskip
\item [(i)] on a l'égalité dans l'inégalité ci-dessus~;
\smallskip
\item [(ii)] $\xi$ appartient à $\Pi$~;
\smallskip
\item [(iii)] $\rho-h\hspace{1pt}\xi$ est un vecteur de Weyl de $L$.
\end{itemize}}

\bigskip
\textit{Démonstration.} L'élément $\rho-h\hspace{1pt}\xi$ appartient à $L$ d'après 4.2.1 et 4.2.3. Il est $h$-régulier d'après le point (a) de 4.2.8. Il est \textit{a fortiori} régulier, si bien que l'inégalité du point (b) peut être vue comme une conséquence du scholie 4.2.7. Cependant l'argument de Borcherds \cite{borcherdsthese}\cite{borcherdsleech}, qui utilise le point (b) de 4.2.2, est plus efficace pour traiter des cas d'égalité. En effet, on a $\mathrm{q}(\rho-h\hspace{1pt}\xi)-\mathrm{q}(\rho)=h\hspace{1pt}(h\hspace{1pt}\mathrm{q}(\xi)-\rho.\xi)$. Puisque l'on a $t^2-t\geq 0$ pour tout $t$ dans $\mathbb{Z}$, le second membre de l'égalité du point en question, avec $x=\xi$, est positif ou nul et il est nul si et seulement si l'on a $\alpha.\xi\in\{0,1\}$ pour tout $\alpha$ dans $R_{+}$~; or cette dernière propriété caractérise les éléments de $\Pi$ (voir la démonstration du point (a) de la proposition II.3.8). Ceci démontre l'équivalence $(i)\iff(ii)$ du point (c). L'équivalence $(i)\iff(iii)$ résulte quant à elle de 4.2.7 (cas d'égalité).
\hfill$\square$

\vspace{0,75cm}
\textsc{Holy constructions}
 
\bigskip
Nous en arrivons maintenant à l'énoncé principal du sous-paragraphe 4.2. Cet énoncé mérite le nom de théorème à cause du point (c) qui est implicite dans \cite{conwaysloane23leech}, au moins en ce qui concerne le réseau $\mathrm{vois}_{h}(L;\rho)$.

\bigskip
\textbf{Théorème 4.2.10.} {\em Soient $L$ un réseau de Niemeier avec racines et $\rho$ un vecteur de Weyl de $L$. On note $\mathrm{s}_{\rho}$ la réflexion orthogonale de $\mathbb{Q}\otimes_{\mathbb{Z}}L$ par rapport à l'hyperplan orthogonal à~$\rho$.

\medskip
{\em (a)} La classe de $\rho$ dans l'espace projectif $\mathrm{P}_{L}(\mathbb{Z}/h)$ (resp. $\mathrm{P}_{L}(\mathbb{Z}/(h+1))$) appartient à la quadrique $\mathrm{C}_{L}(\mathbb{Z}/h)$ (resp. $\mathrm{C}_{L}(\mathbb{Z}/(h+1)$).

\medskip
{\em (b)} Les réseaux $\mathrm{vois}_{h}(L;\rho)$ et $\mathrm{vois}_{h+1}(L;\rho)$ sont échangés par la réflexion $\mathrm{s}_{\rho}$~:
$$
\hspace{24pt}
\mathrm{vois}_{h+1}(L;\rho)
\hspace{4pt}=\hspace{4pt}
\mathrm{s}_{\rho}(\mathrm{vois}_{h}(L;\rho))
\hspace{24pt}.
$$

\medskip
{\em (c)} Les réseaux $\mathrm{vois}_{h}(L;\rho)$ et $\mathrm{vois}_{h+1}(L;\rho)$ sont sans racines~:
$$
\hspace{24pt}
\mathrm{vois}_{h}(L;\rho)
\hspace{4pt}\simeq\hspace{4pt}
\mathrm{Leech}
\hspace{24pt},\hspace{24pt}
\mathrm{vois}_{h+1}(L;\rho)
\hspace{4pt}\simeq\hspace{4pt}
\mathrm{Leech}
\hspace{24pt}.
$$

\medskip
{\em (d)} On a~:
$$
\hspace{24pt}
\mathrm{N}_{h}(L,\mathrm{Leech})
\hspace{4pt}=\hspace{4pt}
\frac{\vert W\vert}{\varphi(h)\hspace{1pt}g}
\hspace{24pt},\hspace{24pt}
\mathrm{N}_{h+1}(L,\mathrm{Leech})
\hspace{4pt}=\hspace{4pt}
\frac{\vert W\vert}{\varphi(h+1)}
\hspace{24pt}.
$$
(la notation $\varphi(-)$ désigne ci-dessus l'indicateur d'Euler d'un entier naturel et la notation $\vert-\vert$ le cardinal d'un ensemble fini, on rappelle que $g$ désigne l'indice de $Q$ dans $L$ ou encore la racine carrée de l'indice de connexion $f$ de~$R$).}

\bigskip
\textit{Démonstration des points (a) et (b).} Le point (a) est conséquence de 4.2.5~; le point (b) est un cas particulier du point (b) de 1.14.

\bigskip
\textit{Démonstration du point (c).} Compte tenu du point (b) il suffit de montrer que $\mathrm{vois}_{h}(L;\rho)$ ou $\mathrm{vois}_{h+1}(L;\rho)$ est sans racines. Nous proposons trois démonstrations.

\bigskip
1) La première est très prosaïque. Le logiciel \texttt{PARI} permet de calculer sans problèmes le {\em minimum} d'un réseau entier $\Lambda$, disons $\mathrm{m}(\Lambda)$,  à savoir l'entier $\inf_{x\in\Lambda-\{0\}}x.x$, pour $\dim\Lambda=24$ et donc de vérifier que l'on a bien $\mathrm{m}(\mathrm{vois}_{h}(L;\rho))=4$ pour les 23 réseaux de Niemeier avec racines.

\bigskip
2) La deuxième démonstration consiste à identifier le réseau $\mathrm{vois}_{h}(L;\rho)$ avec la construction du réseau de Leech  donnée par Conway et Sloane dans \cite{conwaysloane23leech} où, comme nous l'avons déjà dit, cette identification est implicite. La construction en question, que Conway et Sloane dénomment {\em holy construction}, est rappelée ci-dessous.

\medskip
Conway et Sloane associent d'abord à la donnée $(L;\rho)$ deux sous-ensembles finis de $L$~:

\smallskip
-- Le premier, disons $F$, est la réunion disjointe $B\coprod-H$~; on rappelle que les éléments de $H$ sont les plus grandes racines, ${\widetilde{\alpha}}_{1},{\widetilde{\alpha}}_{2},\ldots, {\widetilde{\alpha}}_{c}$,  des composantes irréductibles, $R_{1},R_{2},\ldots,R_{c}$, de $R$. Pour allèger la notation on note $\alpha_{1},\alpha_{2},\ldots,\alpha_{24}$ les éléments de $B$ et on pose $\alpha_{24+i}=-{\widetilde{\alpha}}_{i}$~; les éléments de $F$ (qui sont tous des racines de~$R$) correspondent aux sommets du {\em graphe de Dynkin complété} de~$R$.

\smallskip
-- Le second, disons $G$, est constitué des poids minuscules de $R$ (voir la remarque qui suit la proposition II.3.12) qui sont dans $L\cap\overline{C}$. On a donc $G=L\cap\Pi$ et l'application canonique de $G$ dans $L/Q$ est une bijection~; on note $\mu_{0},\mu_{1},\ldots,\mu_{g-1}$ les éléments de $G$ avec $\mu_{0}=0$.

\smallskip
Conway et Sloane considèrent ensuite le réseau de $\mathbb{Q}\otimes_{\mathbb{Z}}L$, disons $\mathrm{HC}(L;\rho)$, constitué des éléments de la forme
$$
\hspace{24pt}
\sum_{i=1}^{24+c}\hspace{2pt}m_{i}\hspace{1pt}\alpha_{i}
\hspace{4pt}+\hspace{4pt}
\sum_{j=0}^{g-1}\hspace{3pt}n_{j}\hspace{1pt}(\frac{\rho}{h}-\mu_{j})
\hspace{24pt},
$$
les $m_{i}$ et $n_{j}$ désignant des entiers vérifiant $\sum_{i}m_{i}+\sum_{j}n_{j}=0$.

\medskip
La proposition suivante montre que le réseau $\mathrm{HC}(L;\rho)$ peut être décrit en termes de voisinage à la Kneser.

\bigskip
\textbf{Proposition 4.2.11.} {\em Soient $L$ un réseau de Niemeier avec racines et $\rho$ un vecteur de Weyl de $L$. Le réseau $\mathrm{HC}(L;\rho)$ de Conway et Sloane coïncide avec le réseau $\mathrm{vois}_{h}(L;\rho)$.}

\bigskip
\footnotesize
\textit{Démonstration.} Il s'agit plutôt d'une vérification d'où les petits caractères. On a vu (voir la discussion qui suit 1.5) que le réseau $\mathrm{vois}_{h}(L;\rho)$ est le sous-$\mathbb{Z}$-module de $\mathbb{Q}\otimes_{\mathbb{Z}}L$ engendré par $M$ et $\frac{\widetilde{\rho}}{h}$, $M$ désignant le noyau de l'homomorphisme de $L$ dans $\mathbb{Z}/h$ qui associe à un élément $x$ de $L$ la réduction $\bmod{\hspace{2pt}h}$ de l'entier $\rho.x$ et $\widetilde{\rho}$ désignant un élément de $L$ avec $\widetilde{\rho}\equiv\rho\bmod{h}$ et $\mathrm{q}(\widetilde{\rho})\equiv 0\bmod{h^{2}}$. On constate que l'on peut prendre $\widetilde{\rho}=\rho-h\alpha_{1}$~; en effet on a $\mathrm{q}(\rho-h\alpha_{1})=h(h+1)-h+h^{2}=2h^{2}$.

\medskip
Ce rappel étant fait, on observe que les $\alpha_{1}-\alpha_{i}$ apartiennent à $M$~; en effet on a $\rho.\alpha_{i}=1$ pour $i\leq 24$ et $\rho.\alpha_{i}=1-h$ pour $i>24$. On observe également que les $\mu_{j}$ appartiennent à~$M$~; en effet le point (c) de 4.2.9 montre que l'on a $\mathrm{q}(\rho-h\mu_{j})=\mathrm{q}(\rho)$, égalité qui équivaut à $\rho.\mu_{j}=h\hspace{1pt}\mathrm{q}(\mu_{j})$.

\medskip
Soit maintenant $x$ un élément de $\mathbb{Q}\otimes_{\mathbb{Z}}L$ de la forme considérée par Conway et Sloane. Les observations ci-dessus montrent que $x$ s'écrit aussi $y+(\sum_{j}n_{j})\frac{\widetilde{\rho}}{h}$ avec $y$ dans $M$, en d'autres termes que l'on a dans $\mathbb{Q}\otimes_{\mathbb{Z}}L$ l'inclusion $\mathrm{HC}(L;\rho)\subset\mathrm{vois}_{h}(L;\rho)$. L'inclusion $\mathrm{vois}_{h}(L;\rho)\subset\mathrm{HC}(L;\rho)$ résulte du fait que $L$ est engendré par $\mathrm{Q}(R)$ et les $\mu_{j}$. Détaillons un peu. On déduit du fait en question que $M$ est engendré par les $\mu_{j}$, les $\alpha_{1}-\alpha_{i}$ et~$h\alpha_{1}$. On se convainc de ce que les $\mu_{j}$ appartiennent à $\mathrm{HC}(L;\rho)$ en contemplant l'égalité $\mu_{j}=(\frac{\rho}{h}-\mu_{0})-(\frac{\rho}{h}-\mu_{j})$. Il est clair qu'il en est de même pour les $\alpha_{1}-\alpha_{i}$. On suppose $\alpha_{1}\in R_{1}$~; on se convainc de ce que $h\alpha_{1}$ appartient à $\mathrm{HC}(L;\rho)$ en observant que l'on a $\widetilde{\alpha}_{1}=\sum_{\beta\in B\cap R_{1}}m_{\beta}\hspace{1pt}\beta$, les $m_{\beta}$ étant des entiers avec  $\sum_{\beta\in B\cap R_{1}}m_{\beta}=h-1$. Enfin, on se convainc de ce que $\frac{\widetilde{\rho}}{h}$ appartient à $\mathrm{HC}(L;\rho)$ en écrivant $\frac{\widetilde{\rho}}{h}=-\alpha_{1}+(\frac{\rho}{h}-\mu_{0})$.
\hfill$\square$
\normalsize

\bigskip
Dans \cite{conwaysloane23leech} Conway et Sloane disent avoir vérifié au cas par cas le fait que $\mathrm{HC}(L;\rho)$ est sans racines... si bien que la démonstration du point (c) grâce à \cite{conwaysloane23leech}  paraît très semblable à celle que nous avons proposée en premier. Il faut cependant signaler que Borcherds \cite{borcherdsthese}\cite{borcherdsleech} a découvert, peu après la parution de \cite{conwaysloane23leech}, une preuve uniforme en termes de réseaux\linebreak ``lorentziens''.

\bigskip
3) La troisième démonstration du point (c) que nous proposons fait un usage systématique de la théorie des voisins à la Kneser. Il s'agit d'une démonstration par l'absurde~: on suppose $\mathrm{R}(\mathrm{vois}_{h}(L;\rho))\not=\emptyset$ et on produit une contradiction.

\medskip
On commence par dégager un énoncé \textit{ad hoc} concernant la théorie générale des $d$-voisins~:

\bigskip
\textbf{Lemme 4.2.12.} {\em Soit $L$ un $\mathrm{q}$-module sur $\mathbb{Z}$~; soient $d\geq 2$ un entier et $u$ un élément $d$-primitif de $L$ avec $\mathrm{q}(u)\equiv 0\bmod{d}$. Soit $x$ un élément de $\mathrm{vois}_{d}(L;u)$, alors $\mathrm{q}(u)\hspace{1pt}x$ appartient à $\mathrm{M}_{d}(L;u)$ et l'on a dans $L$ la congruence suivante~:
$$
\mathrm{q}(u)\hspace{1pt}x
\hspace{4pt}\equiv\hspace{4pt}
(u.x)\hspace{2pt}u
\hspace{8pt}
\bmod{d}
$$
(observer que $u$ appartient à $\mathrm{M}_{d}(L;u)\subset\mathrm{vois}_{d}(L;u)$ si bien que $u.x$ est entier et que $(u.x)\hspace{2pt}u$ appartient également à $\mathrm{M}_{d}(L;u)$).}

\bigskip
\textit{Démonstration.} On pose $M=\mathrm{M}_{d}(L;u)$ et $L'=\mathrm{vois}_{d}(L;u)$. L'élément $\mathrm{q}(u)\hspace{1pt}x$ appartient à $M$ puisque le quotient $L'/M$ est cyclique d'ordre $d$ et que $\mathrm{q}(u)$ est divisible par $d$. Soit $v$ un élément de $L$ avec $u.v\equiv 1\bmod{d}$, puisque $L'$ est engendré par $M$ et $\frac{u-\mathrm{q}(u)v}{d}$ il suffit de vérifier la congruence de l'énoncé pour $x\in M$ et pour $x=\frac{u-\mathrm{q}(u)v}{d}$. Dans le premier cas les deux membres sont divisibles par $d$~; la vérification dans le second cas est immédiate.
\hfill$\square$

\bigskip
\textbf{Scholie 4.2.13.} {\em Soient $L$ un réseau de Niemeier avec racines et $\rho$ un vecteur de Weyl de $L$. Soit $x$ un élément de $\mathrm{vois}_{h}(L;\rho)$, alors $\rho.x$ est entier, $h\hspace{1pt}x$ et $(h+1)\hspace{1pt}\mathrm{s}_{\rho}(x)$ appartiennent respectivement à $\mathrm{M}_{h}(L;\rho)$ et $\mathrm{M}_{h+1}(L;\rho)$, et l'on~a dans $L$ les congruences suivantes~:
$$
\hspace{18pt}
h\hspace{1pt}x
\hspace{4pt}\equiv\hspace{4pt}
(\rho.x)\hspace{2pt}\rho
\hspace{8pt}
\bmod{h}
\hspace{18pt},\hspace{18pt}
(h+1)\hspace{1pt}\mathrm{s}_{\rho}(x)
\hspace{4pt}\equiv\hspace{4pt}
(\rho.x)\hspace{2pt}\rho
\hspace{8pt}
\bmod{(h+1)}
\hspace{18pt}.
$$}

\medskip
\textit{Démonstration.} Elle résulte du lemme 4.2.12 et des égalités $\mathrm{q}(\rho)=h(h+1)$, $\mathrm{vois}_{h+1}(L;\rho)=\mathrm{s}_{\rho}(\mathrm{vois}_{h}(L;\rho))$ et $\rho.\mathrm{s}_{\rho}(x)=-\rho.x$.
\hfill$\square$

\medskip
On analyse maintenant les contraintes que la condition $\mathrm{R}(\mathrm{vois}_{h}(L;\rho))\not=\emptyset$ impose.

\bigskip
\textbf{Proposition 4.2.14.} {\em Soient $L$ un réseau de Niemeier avec racines, $\rho$ un vecteur de Weyl de $L$ et $\alpha'$ une racine du réseau $\mathrm{vois}_{h}(L;\rho)$.

\medskip
{\em (a.1)} L'entier $\rho.\alpha'$ n'est pas divisible par $h$.

\medskip
{\em (a.2)} Les entiers $\rho.\alpha'$ et $h$ ne sont pas premiers entre eux.

\medskip
{\em (b.1)} L'entier $\rho.\alpha'$ n'est pas divisible par $h+1$.

\medskip
{\em (b.2)} Les entiers $\rho.\alpha'$ et $h+1$ ne sont pas premiers entre eux.}
\vfill\eject

\bigskip
\textit{Démonstration.} On pose $L'=\mathrm{vois}_{h}(L;\rho)$ et $M=\mathrm{M}_{h}(L;\rho)$~; on rappelle que l'on a $M=L\cap L'$.

\medskip
\textit{Démonstration du point (a.1).} Le scholie 4.2.13 montre que si $\rho.\alpha'$ est divisible par $h$ alors $\alpha'$ appartient à $L$, ce qui est interdit par $\mathrm{R}(M)=\emptyset$.
\hfill$\square$

\medskip
\textit{Démonstration du point (a.2).} Le scholie 4.2.13 montre que si l'entier $\rho.\alpha'$ est premier avec $h$ alors $h\alpha'$ est un élément $h$-primitif de $L$ et la classe de cet élément dans la quadrique $\mathrm{C}_{L}(\mathbb{Z}/h)$ est égale à celle de $\rho$. Ceci est interdit par l'égalité $\mathrm{q}(h\alpha')=h^{2}$ et 4.2.7.
\hfill$\square$

\medskip
\textit{Démonstration des points (b.1) et (b.2).} On observe tout d'abord que si $\alpha'$ est une racine de $L'$ alors $\mathrm{s}_{\rho}(\alpha')$ est une racine de $\mathrm{s}_{\rho}(L')=\mathrm{vois}_{h+1}(L;\rho)$ et que l'on a $\rho.\mathrm{s}_{\rho}(\alpha')=-\rho.\alpha'$. On procède ensuite comme précédemment. Pour la démonstration du point (b.2), on utilise les énoncés 4.2.15, 4.2.16 et 4.2.17 ci-après~; 4.2.15 est le pendant de 4.2.8, 4.2.16 est semblable au point (b) de 4.2.9 et 4.2.17 à 4.2.7.

\bigskip
\textbf{Proposition 4.2.15.} {\em Soient $L$ un réseau de Niemeier avec racines et $\rho$ un vecteur de Weyl de $L$.

\medskip
{\em (a)} Le vecteur de Weyl $\rho$ est $h+1$-régulier.

\medskip
{\em (b)} Soit $\xi$ un élément de $Q^{\sharp}$~; les deux conditions suivantes sont équivalentes~:
\begin{itemize}
\smallskip
\item [(i)] l'élément $\xi$ est $h+1$-régulier~;
\smallskip
\item [(ii)] il existe un élément $w$ de $W$, un élément $\varpi$ de $\Pi$ et un élément $x$ de $Q$ tels que l'on a $\xi=w(\rho+\varpi)+(h+1)x$.
\end{itemize}
\smallskip
De plus si ces conditions sont vérifiées alors le triplet $(w,\varpi,x)$ qui apparaît dans la condition (ii) est uniquement déterminé en fonction de $\xi$.

\medskip
{\em (c)} Soit $u$ un élément de $L$~; les deux conditions suivantes sont équivalentes~:
\begin{itemize}
\smallskip
\item [(i)] l'élément $u$ est $h+1$-régulier~;
\smallskip
\item [(ii)] il existe un élément $w$ de $W$ et un élément $x$ de $L$ tels que l'on a $u=w\rho+(h+1)x$.
\end{itemize}

\smallskip
De plus si ces conditions sont vérifiées alors le couple $(w,x)$ qui apparaît dans la condition (ii) est uniquement déterminé en fonction de $u$.}

\bigskip
\textit{Démonstration des points (a) et (b).} Ce sont des variantes de celles des points (a) et (b) de 4.2.8.  Démontrons par exemple l'implication $(i)\Rightarrow(ii)$ du point~(b). Compte tenu de 4.1.3, on peut supposer $\xi\in (h+1)\hspace{1pt}\overline{Alc}$ (à nouveau la dernière partie du point (b), concernant l'unicité du triplet $(w,\varpi,x)$, sera un sous-produit de notre démonstration). Si $\xi$ est $h+1$-régulier, alors on a cette fois les inégalités $\alpha.(\xi-\rho)\geq 0$ pour tout $\alpha$ dans $B$ et  $\widetilde{\alpha}.(\xi-\rho)\leq 1$ pour tout $\widetilde{\alpha}$ dans~$H$. On a donc $\xi-\rho\in Q^{\sharp}\cap\overline{Alc}=:\Pi$.
\hfill$\square$

\medskip
\textit{Démonstration du point (c).} L'implication $(ii)\Rightarrow(i)$ du point (c) résulte du fait que $\rho$ est $h+1$-régulier. L'implication $(i)\Rightarrow(ii)$ du point (c) résulte de l'implication $(i)\Rightarrow(ii)$ du point (b). En effet si $u$ est $h+1$-régulier, alors il s'écrit de façon unique $u=w_{0}(\rho+\varpi)+(h+1)x_{0}$ avec $(w_{0},\varpi,x_{0})\in W\times\Pi\times Q$. On observe que puisque $u$ est dans $L$ alors il en est de même pour $\varpi$~: $\varpi\in\Pi\cap L$. On écrit $u=w_{0}(\rho-h\varpi)+(h+1)(w_{0}\varpi+x_{0})$. D'après 4.2.9~(c), $\rho-h\varpi$ est un vecteur de Weyl de $L$, il existe donc un élément $w_{1}$ de $W$, uniquement déterminé en fonction de $\varpi$, tel que l'on a $\rho-h\varpi=w_{1}\rho$. En posant $w=w_{0}w_{1}$ et $x=w_{0}\varpi+x_{0}$ on a bien $u=w\rho+(h+1)x$. L'unicité du couple $(w,x)$ résulte de ce que l'application canonique $\Pi\cap L\to L/Q$ est une bijection.
\hfill$\square$

\bigskip
\textbf{Proposition 4.2.16.} {\em Soit $x$ un élément non nul de $L$, alors on a l'inégalité
$$
\hspace{24pt}
\mathrm{q}(\rho-(h+1)\hspace{1pt}x)
\hspace{4pt}\geq\hspace{4pt}
(h+1)(h+2)
\hspace{24pt}.
$$}

\medskip
\textit{Démonstration.} On adapte l'argument de Borcherds utilisé dans la démonstra\-tion de 4.2.9. Soit $x$ un élément de $L$, on vérifie que l'on a cette fois l'égalité
$$
\hspace{24pt}
\mathrm{q}(\rho-(h+1)\hspace{1pt}x)
\hspace{4pt}=\hspace{4pt}
(h+1)\hspace{2pt}(\hspace{2pt}h+\hspace{2pt}\mathrm{q}(x)+\sum_{\alpha\in R_{+}}\frac{(\alpha.x)^{2}-\alpha.x}{2}\hspace{2pt})
\hspace{24pt}.
$$
Pour estimer le second membre de cette égalité, on distingue deux cas~:

\medskip
1) Il existe une racine positive $\alpha$ telle que l'on a $\alpha.x\not\in\{0,1\}$. On constate que l'on a dans ce cas $(\alpha.x)^{2}-\alpha.x\geq 2$ et donc $\mathrm{q}(x)+\sum_{\alpha\in R_{+}}\frac{(\alpha.x)^{2}-\alpha.x}{2}\geq 2$.

\medskip
2) On a $\alpha.x\in\{0,1\}$ pour tout $\alpha$ dans $R_{+}$. Si $x$ est non nul alors on a $\mathrm{q}(x)\geq 2$~; en effet $x$ ne peut être une racine puisque $x$ ou $-x$ serait une racine positive, disons $\beta$, avec $\beta.x\not\in\{0,1\}$.
\hfill$\square$

\bigskip
\textbf{Corollaire 4.2.17.} {\em Soit $u$ un élément de $L$. Si $u$ est $(h+1)$-primitif et si l'on a $\mathrm{q}(u)={(h+1)}^{2}$ alors $u$ n'est pas $(h+1)$-régulier.}

\bigskip
\textit{Démonstration.} D'après le point (c) de 4.2.15, si un tel $u$ est $(h+1)$-régulier alors il existe un élément $w$ du groupe de Weyl de $R$ et un élément $x$ de $L$ tel que l'on a $u=w\rho+(h+1)\hspace{1pt}x$. Les  propositions 4.2.5 et 4.2.16 montrent que l'on a $\mathrm{q}(u)=h(h+1)$ ou $\mathrm{q}(u)\geq(h+1)(h+2)$.
\hfill$\square$

\bigskip
La proposition 4.2.14 montre que le réseau $\mathrm{vois}_{h}(L;\rho)$ est sans racines si $h$ ou $h+1$ est premier. C'est le cas pour $19$ des réseaux de Niemeier sans racines. Les $4$ qui font de la résistance correspondent à $h=25,14,9,8$. Pour en venir à bout, on raffine l'argument précédent.
\vfill\eject

\bigskip
\textbf{Proposition 4.2.18.} {\em Soit $L$ un réseau de Niemeier avec racines. On suppose que le réseau $\mathrm{vois}_{h}(L;\rho)$ possède aussi des racines et on note $h'$ son nombre de Coxeter. Alors il existe un entier $\nu$ vérifiant les trois conditions suivantes~: 
\begin{itemize}
\smallskip
\item [{\em (1)}] $\nu>0$;
\smallskip
\item [{\em (2)}] $\nu^{2}\leq\frac{h(h+1)}{h'(h'+1)}$~;
\smallskip
\item [{\em (3)}] $\mathrm{p.g.c.d.}(\nu,h)\not=1$ et $\mathrm{p.g.c.d.}(\nu,h+1)\not=1$.
\end{itemize}}

\bigskip
\textit{Démonstration.} On pose $L'=\mathrm{vois}_{h}(L;\rho)$. On applique la proposition 4.2.6 au réseau $L'$ en prenant pour $x$ l'élément $\rho$ : il existe une racine $\alpha'$ de $L'$ telle que l'on a $(\rho.\alpha')^{2}\leq\frac{h(h+1)}{h'(h'+1)}$. On pose $\nu=\vert\rho.\alpha'\vert$. Ce qui précède montre que $\nu$ est un nombre entier qui vérifie les trois conditions ci-dessus.
\hfill$\square$

\bigskip
On note $\mathrm{S}(h,h')$ le sous-ensemble de $\mathbb{Z}$ constitué des entiers $\nu$ vérifiant les trois conditions de la proposition 4.2.18~; il est clair que l'on a $\mathrm{S}(h,h_{1}')\subset\mathrm{S}(h,h_{2}')$ pour $h_{1}'\geq h_{2}'$. On constate que $\mathrm{S}(h,2)$ est vide pour $h\not=25$. Le réseau $\mathrm{vois}_{h}(L;\rho)$ est donc sans racines pour $h\not=25$.

\medskip
Il reste à faire un dernier effort pour $h=25$. On a $\mathrm{S}(25,2)=\{10\}$ et $\mathrm{S}(25,3)=\emptyset$. La seconde égalité force $h'=2$ et donc $\mathrm{R}(L')=24\hspace{1pt}\mathbf{A}_{1}$. On contemple alors l'égalité (Scholie II.3.4)
$$
\hspace{24pt}
\sum_{\beta\in R_{+}'}\hspace{2pt}(\rho.\beta)^{2}
\hspace{4pt}=\hspace{4pt}
2600
\hspace{24pt},
$$
$R_{+}'$ désignant l'ensemble à $24$ éléments constitué des racines positives de $L'$, pour un certain choix de chambre. Cette égalité montre que l'on ne peut avoir $\vert\rho.\beta\vert=10$ pour tout $\beta$ dans $R_{+}'$~; la proposition 4.2.14 montre qu'il existe $\beta_{1}$ dans $R_{+}'$  avec $\vert\rho.\beta_{1}\vert\geq 20$. Il en résulte qu'il existe $\beta_{2}$ dans $R_{+}'$  avec $\vert\rho.\beta_{2}\vert\leq 9$ (observer que l'on a $23\times 10^{2}>2600-20^{2}$). Contradiction avec la proposition 4.2.14.

\medskip
Cette contradiction achève notre troisième démonstration du point (c) du théorème 4.2.10.
\hfill$\square\square$

\bigskip
\textit{Démonstration du point (d) du théorème 4.2.10} 

\medskip
La proposition 4.2.8 montre que la classe de $\rho$ dans  $\mathrm{P}_{L}(\mathbb{Z}/h)$ appartient à $\mathrm{P}^{\mathrm{r\acute{e}g}}_{L}(\mathbb{Z}/h)$ et que l'action de $W$ sur cet ensemble est transitive. Comme la classe de $\rho$ appartient à $\mathrm{C}^{\mathrm{r\acute{e}g}}_{L}(\mathbb{Z}/h)$, on constate que l'on a $\mathrm{C}^{\mathrm{r\acute{e}g}}_{L}(\mathbb{Z}/h)=\mathrm{P}^{\mathrm{r\acute{e}g}}_{L}(\mathbb{Z}/h)$. Pareillement, la proposition 4.2.15 montre que la classe de $\rho$ dans $\mathrm{P}_{L}(\mathbb{Z}/(h+1))$ appartient à $\mathrm{P}^{\mathrm{r\acute{e}g}}_{L}(\mathbb{Z}/(h+1))$, que l'action de $W$ sur cet ensemble est transitive (4.2.15 (c)) et que l'on a $\mathrm{C}^{\mathrm{r\acute{e}g}}_{L}(\mathbb{Z}/(h+1))=\mathrm{P}^{\mathrm{r\acute{e}g}}_{L}(\mathbb{Z}/(h+1))$. Compte tenu du fait qu'une condition nécéssaire pour que $\mathrm{vois}_{d}(L;c)$ soit isomorphe au réseau de Leech est que $c$ appartienne à $\mathrm{C}^{\mathrm{r\acute{e}g}}_{L}(\mathbb{Z}/d)$ (4.1.4 (a)), la démonstration du point (d) de 4.2.10 consiste à vérifier que le cardinal du stabilisateur, pour l'action de $W$, de la classe de $\rho$ dans $\mathrm{P}_{L}(\mathbb{Z}/h)$ (resp. $\mathrm{P}_{L}(\mathbb{Z}/(h+1))$), est $\phi(h)g$ (resp. $\phi(h+1)$). Ceci résulte des énoncés 4.2.19 et 4.2.20 ci-dessous.

\bigskip
\textbf{Proposition 4.2.19.} {\em Soient $L$ un réseau de Niemeier avec racines et $\rho$ un vecteur de Weyl de $L$. Alors le stabilisateur, pour l'action de $W$, de la classe de $\rho$ dans $\mathrm{P}_{L}(\mathbb{Z}/h)$, est une extension canonique de $(\mathbb{Z}/h)^{\times}$ par $L/Q$.}

\bigskip
\textit{Démonstration.}  Soient $S\subset W$ le stabilisateur en question et $w$ un élément de $S$. On a par définition $w\rho=\lambda\rho+hx$ avec $\lambda$ dans $\mathbb{Z}$ premier à $h$ et $x$ dans~$L$. On constate que la classe $\overline{\lambda}$ de $\lambda$ dans $(\mathbb{Z}/h)^{\times}$ ne dépend que de $w$ et que l'application $w\mapsto\overline{\lambda}$ est un homomorphisme de groupes, disons $\pi:S\to(\mathbb{Z}/h)^{\times}$. L'implication $(i)\Rightarrow(ii)$ de 4.2.8 (b) montre que $\pi$ est surjectif puisque $\lambda\rho$ est $h$-régulier pour tout $\lambda$ dans $\mathbb{Z}$ premier à $h$. On considère maintenant le sous-groupe $\ker\pi$. Soit $w$ dans $\ker\pi$~; on a cette fois $w\rho=\rho+hx$ avec $x$ dans $L$. On constate que l'application qui associe à $w$ la classe de $x$ dans $L/Q$ est un homomorphisme de groupes, disons $\iota:\ker\pi\to L/Q$ (on utilise ici que l'action de $W$ sur $Q^{\sharp}/Q$ est triviale). Le même argument que précédemment montre que $\iota$ est surjectif. La propriété d'unicité qui apparaît dans 4.2.8 (b) montre que $\iota$ est injectif.
\hfill$\square$

\bigskip
\textit{Remarques}

\medskip
1) Le $\mathbb{Z}$-module $L/Q$ est un $\mathbb{Z}/h$-module d'après 4.2.3, si bien que l'on dispose d'une action naturelle de $(\mathbb{Z}/h)^{\times}$ sur $L/Q$~; cette action coïncide avec celle définie par l'extension de 4.2.19. Si $\rho$ appartient à $Q$ (ce qui n'est pas toujours le cas, voir 4.2.4), alors la démonstration de 4.2.19 montre implicitement que $S$ est canoniquement isomorphe au produit semi-direct $L/Q\rtimes(\mathbb{Z}/h)^{\times}$. En fait, on peut vérifier que l'extension de 4.2.19 est toujours triviale.

\medskip
2) Par construction l'homomorphisme $\iota$ qui apparaît dans la démonstration de 4.2.19 se factorise par une application ensembliste de $\ker\pi$ dans $L$. Ceci montre que l'on dispose d'une section ensembliste (ne dépendant que du choix de $\rho$) de l'homomorphisme $L\to L/Q$. L'image de cette section est $\Pi\cap L$ et cette seconde remarque est intimement reliée à la remarque qui suit 4.2.9.

\bigskip
\textbf{Proposition 4.2.20.} {\em Soient $L$ un réseau de Niemeier avec racines et $\rho$ un vecteur de Weyl de $L$. Alors le stabilisateur, pour l'action de $W$, de la classe de $\rho$ dans $\mathrm{P}_{L}(\mathbb{Z}/(h+1))$, est canoniquement isomorphe à $(\mathbb{Z}/(h+1))^{\times}$.}

\bigskip
\textit{Démonstration.} Elle est analogue à la précédente, la proposition 4.2.15 remplaçant la proposition 4.2.8. Elle est en réalité plus simple. Cette fois l'homomorphisme $\pi:S\to(\mathbb{Z}/(h+1))^{\times}$ est un isomorphisme. La raison de cette simplification est la suivante : dans le point (b) de 4.2.8 l'élément $x$ appartient à $Q$, dans le point (c) de 4.2.15 il appartient à $L$.
\hfill$\square\square\square$

\vspace{0,75cm}
\textbf{4.3.} Sur les stabilisateurs de l'action de $W$ sur $\mathrm{P}^{\mathrm{r\acute{e}g}}_{L}(\mathbb{Z}/d)$, pour $L$ un réseau de Niemeier avec racines

\bigskip
On majore dans ce sous-paragraphe la taille des stabilisateurs en question. Notre motivation est l'énoncé 4.3.3 dont on verra l'utilité en IX.\ref{calcfinvpgenre2-1}.

\medskip
Soient $L$ un réseau de Niemeier avec racines et $d\geq 2$ un entier.

\medskip
On suppose ci-après que $d$ est premier avec l'indice $g$ de $Q$ dans $L$. Dans ce cas l'homomorphisme canonique $Q/dQ\to L/dL$ est un isomorphisme. On introduit, \textit{mutatis mutandis}, les notations $\mathrm{P}_{Q}$, $\mathrm{C}_{Q}$ (cette fois ce schéma est seulement lisse sur $\mathbb{Z}[\frac{1}{g}]$), $\mathrm{P}_{Q}(\mathbb{Z}/d)$, $\mathrm{P}_{Q}(\mathbb{Z}/d)$, $\mathrm{P}^{\mathrm{r\acute{e}g}}_{Q}(\mathbb{Z}/d)$ et $\mathrm{C}^{\mathrm{r\acute{e}g}}_{Q}(\mathbb{Z}/d)$, ainsi que la terminologie correspondante~; il est clair que les bijections canoniques $\mathrm{P}_{Q}(\mathbb{Z}/d)\cong\mathrm{P}_{L}(\mathbb{Z}/d)$, $\mathrm{C}_{Q}(\mathbb{Z}/d)\cong\mathrm{C}_{L}(\mathbb{Z}/d)$, $\mathrm{P}^{\mathrm{r\acute{e}g}}_{Q}(\mathbb{Z}/d)\cong\mathrm{P}^{\mathrm{r\acute{e}g}}_{L}(\mathbb{Z}/d)$ et $\mathrm{C}^{\mathrm{r\acute{e}g}}_{Q}(\mathbb{Z}/d)\cong\mathrm{C}^{\mathrm{r\acute{e}g}}_{L}(\mathbb{Z}/d)$ sont $W$-équivariantes.

\bigskip
\textbf{Proposition 4.3.1.} {\em Soient $L$ un réseau de Niemeier avec racines et $d\geq 2$ un entier premier avec l'indice de $Q$ dans $L$~; soit $S$ le stabilisateur, pour l'action de $W$, d'un élément de $\mathrm{P}^{\mathrm{r\acute{e}g}}_{L}(\mathbb{Z}/d)$.

\medskip
{\em (a)} Le groupe $S$ est canoniquement isomorphe à un sous-groupe de $(\mathbb{Z}/d)^{\times}$.

\medskip
{\em (b)} Si $d$ est premier, alors l'action de $S$ sur $R$ (induite par celle de $W$) est libre.}

\bigskip
\textit{Démonstration du point (a).} Compte tenu de ce que l'on a dit plus haut, on peut remplacer $\mathrm{P}_{L}(\mathbb{Z}/d)$ par $\mathrm{P}_{Q}(\mathbb{Z}/d)$. Soient $u$ un élément $d$-primitif de $Q$, $S\subset W$ le stabilisateur de la classe de $u$ dans $\mathrm{P}_{Q}(\mathbb{Z}/d)$ et $w$ un élément de~$S$. On procède comme dans la démonstration de 4.2.19. On a par définition $wu=\lambda u+dx$ avec $\lambda$ dans $\mathbb{Z}$ premier à $d$ et $x$ dans $Q$, la classe $\overline{\lambda}$ de $\lambda$ dans $(\mathbb{Z}/d)^{\times}$ ne dépend que de $w$ et  l'application $\pi:S\to(\mathbb{Z}/d)^{\times}\hspace{2pt},\hspace{2pt}w\mapsto\overline{\lambda}$ est un homomorphisme de groupes. On vérifie que $\pi$ est injectif si $u$ est $d$-régulier à l'aide du scholie 4.1.3.
\hfill$\square$

\medskip
\textit{Démonstration du point (b).} L'égalité $wu=\lambda u+dx$ implique la congruence $\alpha.(wu)\equiv\lambda(\alpha.u)\bmod{d}$, pour tout $\alpha$ dans $R$, ou encore $(w^{-1}\alpha).u\equiv\lambda(\alpha.u)\linebreak\bmod{d}$. On constate que si l'on a $w\alpha=\alpha$ (et donc $w^{-1}\alpha=\alpha$) alors on a $(\lambda-1)(\alpha.u)\equiv 0\bmod{d}$, soit encore $(\alpha.u)(\pi(w)-1)=0$ dans $\mathbb{Z}/d$. Si $u$ est $d$-régulier  et si $d$ est premier alors on obtient $\pi(w)=1$. Comme $\pi$ est injectif  quand $u$ est $d$-régulier, on a bien l'implication $w\alpha=\alpha\Rightarrow w=\mathrm{id}$.
\hfill$\square$

\bigskip
\textit{Remarque.} La proposition 4.2.20 est une illustration de l'énoncé 4.3.1 (a)~; la proposition 4.2.19 montre que l'hypothèse faite sur $d$ n'est pas superflue.

\bigskip
\textbf{Corollaire-Définition 4.3.2.} {\em Soient $L$ un réseau de Niemeier avec racines et $p$ un nombre premier~; on note $\mathrm{D}_{p}(L)$ le p.g.c.d. des entiers $p-1$, $24\hspace{1pt}h$ et $\vert W\vert$. Si $p$ ne divise pas l'indice de $Q$ dans $L$, alors le stabilisateur, pour l'action de $W$, d'un élément de $\mathrm{P}^{\mathrm{r\acute{e}g}}_{L}(\mathbb{F}_{p})$ est canoniquement isomorphe à un sous-groupe du groupe $\mu_{\mathrm{D}_{p}(L)}(\mathbb{F}_{p})$ (qui est cyclique d'ordre $\mathrm{D}_{p}(L)$).}

\bigskip
\textit{Démonstration.} Soit $S\subset W$ l'un de ces stabilisateurs. D'après 4.3.1 (a), $S$~s'identifie à un sous-groupe de $\mathbb{F}_{p}^{\times}$~; d'après 4.3.1 (b) le cardinal de $S$ divise le cardinal de $R$ à savoir $24\hspace{1pt}h$.
\hfill$\square$

\bigskip
\textit{Remarque.} On constate par inspection que $24\hspace{1pt}\mathrm{h}(L)$ divise $\vert\mathrm{W}(L)\vert$ sauf pour $\mathrm{R}(L)=24\hspace{1pt}\mathbf{A}_{1}$ auquel cas le p.g.c.d. des entiers  $24\hspace{1pt}\mathrm{h}(L)$ et $\vert\mathrm{W}(L)\vert$ est $16$ ; on a donc $\mathrm{p.g.c.d.}\hspace{1pt}(p-1,24\hspace{1pt}\mathrm{h}(L),\vert\mathrm{W}(L)\vert)=\mathrm{p.g.c.d.}\hspace{1pt}(p-1,24\hspace{1pt}\mathrm{h}(L))$ en dehors de ce cas.

\bigskip
\textbf{Scholie-Définition 4.3.3.} {\em Soient $L$ un réseau de Niemeier avec racines et $p$ un nombre premier~; on note $\mathrm{pas}\hspace{1pt}(L;p)$ l'entier défini par
$$
\hspace{24pt}
\mathrm{pas}\hspace{1pt}(L;p)
\hspace{4pt}:=\hspace{4pt}
\frac{\vert\mathrm{W}(L)\vert}{\mathrm{p.g.c.d.}\hspace{1pt}(p-1,24\hspace{1pt}\mathrm{h}(L),\vert\mathrm{W}(L)\vert)}
\hspace{24pt}.
$$
Si $p$ ne divise pas l'indice de $Q$ dans $L$ alors $\mathrm{N}_{p}(L,\mathrm{Leech})$ est divisible par $\mathrm{pas}\hspace{1pt}(L;p)$~; on note dans ce cas $\mathrm{n}_{p}(L)$ l'entier défini par
$$
\hspace{24pt}
\mathrm{N}_{p}(L,\mathrm{Leech})
\hspace{4pt}=\hspace{4pt}
\mathrm{n}_{p}(L)\hspace{4pt}
\mathrm{pas}\hspace{1pt}(L;p)
\hspace{24pt}.
$$}

\bigskip
\textit{Démonstration.} L'entier $\mathrm{N}_{p}(L,\mathrm{Leech})$ est la somme des cardinaux des $W$-orbites des points $c$ de $\mathrm{C}_{L}(\mathbb{F}_{p})$ avec $\mathrm{vois}_{p}(L;c)\simeq\mathrm{Leech}$, or ces points appartiennent à $\mathrm{P}^{\mathrm{r\acute{e}g}}_{L}(\mathbb{F}_{p})$ d'après 4.1.4 (a).
\hfill$\square$

\bigskip
\textit{Remarque.} L'entier $\mathrm{pas}\hspace{1pt}(L;p)$ est le produit de l'entier $\mathrm{pas}_{1}(L)$ (indépendant de~$p$) et de l'entier $\mathrm{pas}_{2}\hspace{2pt}(L;p)$ respectivement définis par
$$
\hspace{4pt}
\mathrm{pas}_{1}(L)
=
\frac{\vert\mathrm{W}(L)\vert}{\mathrm{p.g.c.d.}\hspace{1pt}(24\hspace{1pt}\mathrm{h}(L),\vert\mathrm{W}(L)\vert)}
\hspace{4pt},\hspace{4pt}
\mathrm{pas}_{2}(L;p)
=
\frac{24\hspace{1pt}\mathrm{h}(L)}{\mathrm{p.g.c.d.}\hspace{1pt}(p-1,24\hspace{1pt}\mathrm{h}(L))}
\hspace{4pt}.
$$

\bigskip
\textsc{Exemples}

\medskip
On illustre ce qui précède en considérant le réseau de Niemeier  $\mathrm{A}_{24}^{+}$ associé au système de racines $\mathbf{A}_{24}$ (voir le deuxième des exemples qui suivent II.3.15) et les nombres premiers $29$ et $31$~; cette illustration est intéressée~: le calcul des entiers $\mathrm{N}_{29}(\mathrm{A}_{24}^{+},\mathrm{Leech})$ et $\mathrm{N}_{31}(\mathrm{A}_{24}^{+},\mathrm{Leech})$ sera utilisé en IX.\ref{calcfinvpgenre2-1}.

\bigskip
On rappelle que l'on a par construction $\mathrm{Q}(\mathbf{A}_{24})=\mathrm{A}_{24}$, $\mathrm{A}_{24}$ désignant le sous-$\mathbb{Z}$-module de $\mathbb{Z}^{25}$ constitué des $25$-uples $(x_{1},x_{2},\ldots,x_{25})$ avec $\sum_{i}x_{i}=0$, muni de la forme bilinéaire paire induite par la structure euclidienne de $\mathbb{R}^{25}$. Le groupe de Weyl $W$ s'identifie au groupe symétrique  $\mathfrak{S}_{25}$, son action sur $\mathrm{A}_{24}$ étant l'action évidente. Il en résulte que le $\mathbb{F}_{p}$-espace vectoriel $\mathbb{F}_{p}\otimes_{\mathbb{Z}}\mathrm{A}_{24}$ s'id\-entifie au sous-espace vectoriel de ${\mathbb{F}_{p}}^{\hspace{-2pt}25}$ constitué des éléments $(x_{1},x_{2},\ldots,x_{25})$ avec  $\sum_{i}x_{i}=0$, l'action induite de $\mathfrak{S}_{25}$ étant toujours l'action évidente.

\medskip
On note respectivement $\widetilde{\mathrm{C}}_{\mathrm{A}_{24}}(\mathbb{F}_{p})$ et $\widetilde{\mathrm{C}}^{\mathrm{r\acute{e}g}}_{\mathrm{A}_{24}}(\mathbb{F}_{p})$ les images inverses de $\mathrm{C}_{\mathrm{A}_{24}}(\mathbb{F}_{p})$ et  $\mathrm{C}^{\mathrm{r\acute{e}g}}_{\mathrm{A}_{24}}(\mathbb{F}_{p})$ dans $\mathbb{F}_{p}\otimes_{\mathbb{Z}}\mathrm{A}_{24}-\{0\}$~; pour $p\not=2$, $\widetilde{\mathrm{C}}_{\mathrm{A}_{24}}(\mathbb{F}_{p})$ s'identifie, comme $(\mathfrak{S}_{25}\times\mathbb{F}_{p}^{\times})$-ensemble, au sous-ensemble de ${\mathbb{F}_{p}}^{\hspace{-2pt}25}-\{(0,0,\ldots,0)\}$ constitué des éléments $(x_{1},x_{2},\ldots,x_{25})$ vérifiant
$$
x_{1}+x_{2}+\ldots+x_{25}=0
\hspace{24pt}\text{et}\hspace{24pt}
x_{1}^{2}+x_{2}^{2}+\ldots+x_{25}^{2}=0
\hspace{24pt}.
\leqno{(*)\hspace{8pt}
}
$$
Par définition le $25$-uple $(x_{1},x_{2},\ldots,x_{25})$ est une application ensembliste de $\{1,2,\ldots,25\}$ dans $\mathbb{F}_{p}$~; on la note $x$. On constate que $x$ appartient à $\widetilde{\mathrm{C}}_{\mathrm{A}_{24}}^{\mathrm{r\acute{e}g}}(\mathbb{F}_{p})$ si et seulement si l'application $x$ est injective. En effet, si l'on voit le système de racines $\mathbf{A}_{24}$, de la façon habituelle, comme un sous-ensemble de l'espace euclidien $\mathbb{R}^{25}$ muni de sa base cano\-nique $\{\varepsilon_{1},\varepsilon_{2},\ldots,\varepsilon_{25}\}$ \cite[Planche~I]{bourbaki}, alors les racines sont les $\varepsilon_{i}-\varepsilon_{j}$, $i\not=j$. On remarquera incidemment que cette observation permet de se convaincre, à peu de frais, que si l'ensemble $\mathrm{P}^{\mathrm{r\acute{e}g}}_{\mathrm{A}_{24}^{+}}(\mathbb{F}_{p})$ est non vide alors on a nécessairement $p\geq 25$, au moins pour $p\not=5$).

\medskip
On en vient maintenant aux deux exemples que nous avions en vue.

\bigskip
1) Sur les réseaux $\mathrm{vois}_{29}(\mathrm{A}_{24}^{+};c)$ pour $c$ dans $\mathrm{C}^{\mathrm{r\acute{e}g}}_{\mathrm{A}_{24}^{+}}(\mathbb{F}_{29})$

\medskip
Soit $x:\{1,2,\ldots,25\}\to\mathbb{F}_{29}$ une application injective~;  on note $\{y_{1},y_{2},y_{3},y_{4}\}$ la partie à $4$ éléments de $\mathbb{F}_{29}$ complémentaire de l'image de $x$. Comme l'on a $\sum_{z\in\mathbb{F}_{29}}z=0$ et $\sum_{z\in\mathbb{F}_{29}}z^{2}=0$, l'application $x$ vérifie (*) si et seulement cette partie à $4$ éléments vérifie
$$
y_{1}+y_{2}+y_{3}+y_{4}=0
\hspace{24pt}\text{et}\hspace{24pt}
y_{1}^{2}+y_{2}^{2}+y_{3}^2+y_{4}^{2}=0
\hspace{24pt}.
$$
On note $\mathcal{C}$ (resp. $\widetilde{\mathcal{C}}$) le sous-schéma de $\mathbf{P}^{3}$ (resp. $\mathbf{A}^{4}-\{0\}$), disons sur $\mathbb{F}_{29}$, défini par les équations ci-dessus~; il est clair que $\mathcal{C}$ est isomorphe à $\mathbf{P}^{1}$. On~note enfin $\mathcal{C}^{\mathrm{r\acute{e}g}}$ (resp. $\widetilde{\mathcal{C}}^{\mathrm{r\acute{e}g}}$) l'ouvert de $\mathcal{C}$ (resp. $\widetilde{\mathcal{C}}$) défini par $y_{i}\not=y_{j}$ pour $i\not=j$~; comme $-2$ n'est pas un carré dans $\mathbb{F}_{29}$ on a en fait $\mathcal{C}(\mathbb{F}_{29})=\mathcal{C}^{\mathrm{r\acute{e}g}}(\mathbb{F}_{29})$ et $\widetilde{\mathcal{C}}(\mathbb{F}_{29})=\widetilde{\mathcal{C}}^{\mathrm{r\acute{e}g}}(\mathbb{F}_{29})$.

\medskip
Nous avons introduit le formalisme ci-dessus pour en arriver aux énoncés suivants~:

\smallskip
-- L'ensemble $\mathfrak{S}_{25}\backslash\widetilde{\mathrm{C}}^{\mathrm{r\acute{e}g}}_{\mathrm{A}_{24}}(\mathbb{F}_{29})$ quotient de l'action de $\mathfrak{S}_{25}$ sur $\widetilde{\mathrm{C}}^{\mathrm{r\acute{e}g}}_{\mathrm{A}_{24}}(\mathbb{F}_{29})$ et l'ensemble $\mathfrak{S}_{4}\backslash\widetilde{\mathcal{C}}^{\mathrm{r\acute{e}g}}(\mathbb{F}_{29})$ quotient de l'action évidente de $\mathfrak{S}_{4}$ sur $\widetilde{\mathcal{C}}^{\mathrm{r\acute{e}g}}(\mathbb{F}_{29})$ sont canoniquement en bijection. De plus ces deux quotients sont munis d'actions naturelles de $\mathbb{F}_{29}^{\times}$ et la bijection que l'on vient d'évoquer est équivariante.

\smallskip
-- On dispose d'une bijection canonique, disons
$$
\begin{CD}
\hspace{24pt}
\kappa:\mathfrak{S}_{4}\backslash\mathcal{C}^{\mathrm{r\acute{e}g}}(\mathbb{F}_{29})\
@>\cong>>
\mathfrak{S}_{25}\backslash\mathrm{C}^{\mathrm{r\acute{e}g}}_{\mathrm{A}_{24}}(\mathbb{F}_{29})
\hspace{24pt}.
\end{CD}
$$
De plus, pour toute $\mathfrak{S}_{4}$-orbite $\mathcal{O}$ de $\mathcal{C}^{\mathrm{r\acute{e}g}}(\mathbb{F}_{29})$, les stabilisateurs de $\mathcal{O}$ et $\kappa(\mathcal{O})$, qui s'identifient tous deux à des sous-groupes de $\mathbb{F}_{29}^{\times}$, sont canoniquement  isomorphes.

\medskip
On analyse enfin l'action de $\mathfrak{S}_{4}$ sur l'ensemble $\mathcal{C}(\mathbb{F}_{29})=\mathcal{C}^{\mathrm{r\acute{e}g}}(\mathbb{F}_{29})$~; pour cela il peut être utile d'observer que le fait que $-3$ n'est pas un carré dans $\mathbb{F}_{29}$ implique que tout élément $(y_{1},y_{2},y_{3},y_{4})$ de $\widetilde{\mathcal{C}}(\mathbb{F}_{29})$ vérifie $y_{i}\not=0$ pour tout~$i$.

\medskip
On vérifie que l'action du groupe $\mathfrak{S}_{4}$ sur l'ensemble à $30$ éléments $\mathcal{C}(\mathbb{F}_{29})$ a exactement deux orbites~:

\smallskip
-- l'orbite, disons $\mathcal{O}_{1}$, de la classe du point $(1,12,-1,-12)$ de $\mathbb{F}_{29}^{4}$ dont le stabilisateur est isomorphe à $\mu_{4}(\mathbb{F}_{29})$ (observer que $\{1,12,-1,-12\}\subset\mathbb{F}_{29}^{\times}$ est le sous-groupe $\mu_{4}(\mathbb{F}_{29})$),

\smallskip
--  l'orbite, disons $\mathcal{O}_{2}$, de la classe du point $(1,4,6,-11)$ qui est libre.

\medskip
Il en résulte que l'action du groupe $\mathfrak{S}_{25}$ sur l'ensemble $\mathrm{C}^{\mathrm{r\acute{e}g}}_{\mathrm{A}_{24}}(\mathbb{F}_{29})$ a exactement deux orbites à savoir $\Omega_{1}=\kappa(\mathcal{O}_{1})$ dont le stabilisateur est isomorphe à $\mu_{4}(\mathbb{F}_{29})$ et $\Omega_{2}=\kappa(\mathcal{O}_{2})$ qui est libre~; on observera que ceci confirme 4.3.2 puisque l'on a $\mathrm{D}_{29}(\mathrm{A}_{24}^{+})=4$.

\bigskip
\texttt{PARI} nous dit que les réseaux $\mathrm{vois}_{29}(\mathrm{A}_{24}^{+};\Omega_{1})$ et $\mathrm{vois}_{29}(\mathrm{A}_{24}^{+};\Omega_{2})$ (l'abus de notation est véniel) sont tous deux isomorphes au réseau de Leech. On obtient donc finalement~:
$$
\mathrm{N}_{29}(\mathrm{A}_{24}^{+},\mathrm{Leech})
\hspace{4pt}=\hspace{4pt}
\frac{5}{4}\hspace{4pt}
\vert\mathrm{W}(\mathbf{A}_{24})\vert
\hspace{4pt}=\hspace{4pt}
19389012554163732480000000
$$
(soit encore, avec la notation introduite en 4.3.3, $\mathrm{n}_{29}(\mathrm{A}_{24}^{+})=5$).

\medskip
En fait, on peut éviter ci-dessus le recours à \texttt{PARI}, en invoquant la proposition \textit{ad hoc} suivante~:

\bigskip
\textbf{Proposition 4.3.4.} {\em Soient $L$ un réseau de Niemeier avec racines, $\rho$ un vecteur de Weyl et $\alpha$ une racine de $L$~; on note $d$ l'entier $2h+1-\rho.\alpha$. Alors~:
\begin{itemize}
\smallskip
\item [{\em (1)}] On a $d\geq h+2$.
\smallskip
\item [{\em (2)}] On a $\mathrm{q}(\rho-h\hspace{1pt}\alpha)=h\hspace{1pt}d$.
\smallskip
\item [{\em (3)}] Il existe $\beta$ appartenant à $B$ (en clair, la base de $\mathrm{R}(L)$ déterminée par $\rho$) avec $\alpha.\beta=0$.
\smallskip
\item [{\em (4)}] L'élément $\rho-h\hspace{1pt}\alpha$ de $L$ est primitif.
\smallskip
\item [{\em (5)}] Le réseau $\mathrm{vois}_{d}(L;\rho-h\hspace{1pt}\alpha)$ (qui est bien défini d'après (2) et (4)) est isomorphe au réseau $\mathrm{vois}_{h}(L;\rho)$ (qui est isomorphe au réseau de Leech d'après 4.2.10 (c)).
\end{itemize}}

\bigskip
\textit{Démonstration.} La propriété (1) résulte de l'inégalité $\vert\rho.\alpha\vert\leq h-1$~;  (2) est immédiate. La propriété (3) est claire si le système de $R$ n'est pas irréductible~; on peut s'en convaincre par inspection dans le cas contraire ($R=\mathbf{A}_{24},\mathbf{D}_{24}$). La propriété (3) implique la propriété (4)~: observer que l'on a $(\rho-h\hspace{1pt}\alpha).\beta=1$. Enfin la propriété (5) est conséquence du point (b) de 1.14.
\hfill$\square$

\bigskip
Pour appliquer cette proposition au cas qui nous intéresse, à savoir $L=\mathrm{A}_{24}^{+}$ et $d=29$, nous devons choisir $\rho$ et $\alpha$ avec $\rho.\alpha=22$. Suivant Bourbaki, on prend $\rho=\sum_{i=1}^{25}(13-i)\hspace{1pt}\varepsilon_{i}$~; il y a alors trois choix possibles pour $\alpha$~: $\alpha_{i}=\varepsilon_{i}-\varepsilon_{i+22}$, $i=1,2,3$. Soit $c_{i}$ la classe de $\rho-25\hspace{1pt}\alpha_{i}$ dans $\mathrm{C}^{\mathrm{r\acute{e}g}}_{\mathrm{A}_{24}}(\mathbb{F}_{29})$ ($c_{i}$ est nécessairement $29$-régulière puisque l'on a $\mathrm{vois}_{29}(\mathrm{A}_{24}^{+};\rho)\simeq\mathrm{Leech}$)~; on constate que $c_{2}$ est dans l'orbite $\Omega_{1}$ et $c_{1}$ et $c_{3}$ dans l'orbite $\Omega_{2}$ (soit $w_{0}\in\mathfrak{S}_{25}$ la permutation $i\mapsto 26-i$, on observera que $w_{0}$ fixe $c_{2}$ et échange $c_{1}$ et $c_{3}$).

\bigskip
2) Sur les réseaux $\mathrm{vois}_{31}(\mathrm{A}_{24}^{+};c)$ pour $c$ dans $\mathrm{C}^{\mathrm{r\acute{e}g}}_{\mathrm{A}_{24}^{+}}(\mathbb{F}_{31})$

\medskip
On peut à nouveau appliquer la méthode précédente pour déterminer l'entier $\mathrm{N}_{31}(\mathrm{A}_{24}^{+},\mathrm{Leech})$. Cette fois $\mathcal{C}$ est la quadrique projective définie par les équations
$$
\hspace{12pt}
y_{1}+y_{2}+y_{3}+y_{4}+y_{5}+y_{6}=0
\hspace{12pt}\text{et}\hspace{12pt}
y_{1}^{2}+y_{2}^{2}+y_{3}^2+y_{4}^{2}+y_{3}^2+y_{4}^{2}=0
\hspace{12pt}.
$$
Cependant les cardinaux des ensembles $\mathcal{C}(\mathbb{F}_{31})$ et $\mathcal{C}^{\mathrm{r\acute{e}g}}(\mathbb{F}_{31})$ sont respectivement $30784$ et 18864 si bien que le volume des calculs est plus important. Nous donnons seulement ci-dessous le résultat de ces calculs, en épargnant les détails au lecteur.

\medskip
Les stabilisateurs non triviaux de l'action de $\mathfrak{S}_{6}$ sur $\mathcal{C}^{\mathrm{r\acute{e}g}}(\mathbb{F}_{31})$ sont les sous-groupes $\mu_{6}(\mathbb{F}_{31})$, $\mu_{3}(\mathbb{F}_{31})$ et $\mu_{2}(\mathbb{F}_{31})$ (en accord avec l'égalité $\mathrm{D}_{31}(\mathrm{A}_{24}^{+})=30$).

\smallskip
-- Il existe une seule orbite avec stabilisateur $\mu_{6}(\mathbb{F}_{31})$~; le voisin de $\mathrm{A}_{24}^{+}$ qui lui est associé est $\mathrm{Leech}$.

\smallskip
-- Il existe une seule orbite avec stabilisateur $\mu_{5}(\mathbb{F}_{31})$~; le voisin associé est à nouveau $\mathrm{Leech}$.

\smallskip
-- Il existe exactement $4$ orbites avec stabilisateur $\mu_{3}(\mathbb{F}_{31})$~; une conduit à $\mathrm{Leech}$, deux au réseau de Niemeier dont le système de racines est $24\mathbf{A}_{1}$ et une au réseau de Niemeier dont le système de racines est $12\mathbf{A}_{2}$.

\smallskip
-- Il existe une seule orbite avec stabilisateur $\mu_{2}(\mathbb{F}_{31})$~; le voisin associé est $\mathrm{Leech}$.

\smallskip
-- Il existe exactement $24$ orbites libres~; $8$ conduisent à $\mathrm{Leech}$, $15$ au réseau de Niemeier dont le système de racines est $24\mathbf{A}_{1}$ et une seule au réseau de Niemeier dont le système de racines est $12\mathbf{A}_{2}$.

\smallskip
(On observera au passage que l'inventaire que nous venons de faire montre en particulier que  4.1.4 (a) n'admet pas de réciproque.)

\medskip
On en déduit~:
$$
\mathrm{N}_{31}(\mathrm{A}_{24}^{+},\mathrm{Leech})
\hspace{4pt}=\hspace{4pt}
\frac{46}{5}\hspace{4pt}
\vert\mathrm{W}(\mathbf{A}_{24})\vert
\hspace{4pt}=\hspace{4pt}
142703132398645071052800000
$$
(soit encore $\mathrm{n}_{31}(\mathrm{A}_{24}^{+})=276$).

\vspace{0,75cm}
\textbf{4.4.} Complément~: Sur les $2$-voisins, associés à un vecteur de Weyl, d'un réseau de Niemeier avec racines

\bigskip
Ce sous-paragraphe s'inspire de \cite{borcherdsthese}, en particulier les commentaires par lesquels il s'achève sont à  comparer aux arguments que donnent Borcherds dans cet article pour prouver l'existence \textit{a priori} d'un réseau unimodulaire pair de dimension $24$ sans racines.

\bigskip
Soient $L$ un réseau de Niemeier avec racines et $\rho$ l'un de ses vecteurs de Weyl. L'égalité $\mathrm{q}(\rho)=h(h+1)$ (Proposition 4.2.5) implique la congruence $\mathrm{q}(\rho)\equiv 0\bmod{2}$, si bien que l'on peut considérer le réseau $\mathrm{vois}_{2}(L;\rho)$ ($\rho$ est primitif donc \textit{a fortiori} $2$-primitif). Nous étudions ci-après ce $2$-voisin de $L$.

\bigskip
On étend la définition du nombre de Coxeter d'un réseau de Niemeier avec racines (Proposition-Définition II.3.3), à tous les réseaux de Niemeier, en convenant que le nombre de Coxeter du réseau de Leech est $0$.

\bigskip
\textbf{Proposition 4.4.1.} {\em Soient $L$ un réseau de Niemeier avec racines et $\rho$ un vecteur de Weyl de $L$. On a l'inégalité
$$
\hspace{24pt}
\mathrm{h}(\mathrm{vois}_{2}(L;\rho))
\hspace{4pt}\leq\hspace{4pt}
\frac{\mathrm{h}(L)+1}{2}
\hspace{24pt}.
$$}

\bigskip
\textit{Démonstration.} On pose $h=\mathrm{h}(L)$ et $h'=\mathrm{h}(\mathrm{vois}_{2}(L;\rho))$.

\medskip
On suppose tout d'abord $h$ pair. Soit $\widetilde{\rho}$ un élément de $L$ avec $\widetilde{\rho}\equiv\rho\bmod{h}$ et $\mathrm{q}(\widetilde{\rho})\equiv 0\bmod{h^{2}}$. On a vu, lors de l'étude de l'algorithme des voisins, qu'un tel élément existe~; ici on peut prendre $\widetilde{\rho}=\rho-h\alpha$, $\alpha$ désignant une racine de la base de $\mathrm{R}(L)$ définie par $\rho$. La proposition 1.13 dit que $\frac{\widetilde{\rho}}{2}$ est un élément $\frac{h}{2}$-primitif de $\mathrm{vois}_{2}(L;\rho)$ et que l'on a dans $\mathbb{Q}\otimes_{\mathbb{Z}}L$ l'égalité
$$
\hspace{24pt}
\mathrm{vois}_{h}(L;\rho)
\hspace{4pt}=\hspace{4pt}
\mathrm{vois}_{\frac{h}{2}}(\mathrm{vois}_{2}(L;\rho);\frac{\widetilde{\rho}}{2})
\hspace{24pt}.
$$
Cette égalité montre que le réseau $\mathrm{vois}_{2}(L;\rho)$ possède un $\frac{h}{2}$-voisin sans racines puisque le point (c) de 4.2.10 dit que  le réseau $\mathrm{vois}_{h}(L;\rho)$ est sans racines. On en déduit $\frac{h}{2}\geq h'$ grâce à 4.1.1.

\medskip
Le cas $h$ impair est similaire. On contemple cette fois l'égalité
$$
\mathrm{vois}_{h+1}(L;\rho)
\hspace{4pt}=\hspace{4pt}
\mathrm{vois}_{\frac{h+1}{2}}(\mathrm{vois}_{2}(L;\rho);\frac{\widetilde{\rho}}{2})
$$
avec par exemple $\widetilde{\rho}=\rho+(h+1)\alpha$.
\hfill$\square$

\bigskip
On se propose maintenant d'affiner l'énoncé 4.4.1.

\bigskip
\textbf{Proposition 4.4.2.} {\em Soient $L$ un réseau de Niemeier avec racines et $\rho$ un vecteur de Weyl de $L$. On a l'égalité
$$
\hspace{24pt}
\mathrm{h}(\mathrm{vois}_{2}(L;\rho))
\hspace{4pt}=\hspace{4pt}
\frac{\mathrm{h}(L)}{2}-\frac{\iota(\mathrm{R}(L))}{8}+2
\hspace{24pt},
$$
la notation $\iota(\mathrm{R}(L))$ désignant ci-dessus le nombre d'exposants impairs du groupe de Weyl de $\mathrm{R}(L)$ {\em \cite[Ch. V, \S6, Déf. 2]{bourbaki}.}}

\bigskip
Avant d'expliquer la démonstration de cette proposition, donnons quelques informations sur le nombre d'exposants impairs, disons $\iota(R)$, du groupe de Weyl d'un système de racines $R$.

\smallskip
(1) L'invariant $\iota$ est additif en le sens suivant~: $\iota(R_{1}\coprod R_{2})=\iota(R_{1})+\iota(R_{2})$.

\smallskip
(2) La valeur de $\iota$ sur les systèmes de racines irréductibles de type ADE est la suivante~: $\iota(\mathbf{A}_{l})=[\frac{l+1}{2}]$, $\iota(\mathbf{D}_{l})=2[\frac{l}{2}]$, $\iota(\mathbf{E}_{6})=4$, $\iota(\mathbf{E}_{7})=7$, $\iota(\mathbf{E}_{8})=8$.

\smallskip
(3) Soit $R$ un système de racines de rang $l$ et de nombre de Coxeter $h$. Comme l'ensemble des exposants de $R$ est stable sous l'involution $m\mapsto h-m$ et que son cardinal est $l$ \cite[Ch. V, \S6, \no2]{bourbaki}, on constate que si $h$ est impair alors on  a $l=2\hspace{1pt}\iota(R)$.

\medskip
Soit $R$ un système de racines de type ADE de rang $l$ et de nombre de Coxeter~$h$. Les points (1) et (2) montrent que l'on a l'inégalité $l\geq 2\iota(R)$ et que l'on a l'égalité si et seulement si toutes les composantes irréductibles de $R$ sont de type $\mathbf{A}_{d}$ avec $d$ pair. Cette dernière condition équivaut à $h$ impair, ce qui est en accord avec le point (3).

\medskip
Ce qui précède montre que la proposition 4.4.2 implique bien la proposition 4.4.1 et que l'on a l'égalité dans l'inégalité de 4.4.1 si et seulement si $\mathrm{R}(L)$ est isomorphe à la somme directe de $\frac{24}{d}$ copies de $\mathbf{A}_{d}$ avec $d$ un diviseur pair de $24$ (se rappeler que $\mathrm{R}(L)$ est équicoxeter). On observera incidemment que les trois conditions suivantes sont équivalentes~:

\smallskip
-- $\mathrm{h}(L)$ est impair~;

\smallskip
-- $\mathrm{R}(L)$ est isomorphe à la somme directe de $\frac{24}{d}$ copies de $\mathbf{A}_{d}$ avec $d$ un diviseur pair de $24$~;

\smallskip
-- on a l'égalité $\mathrm{h}(\mathrm{vois}_{2}(L;\rho))=
\frac{\mathrm{h}(L)+1}{2}$.

\bigskip
La démonstration de la proposition 4.4.2 est basée sur les deux propositions 4.4.3 et 4.4.4 ci-après.

\bigskip
\textbf{Proposition 4.4.3.} {\em Soient $L$ un réseau de Niemeier avec racines et $\rho$ un vecteur de Weyl de $L$. Alors le nombre de racines $\alpha$ de $L$ avec $\rho.\alpha$ pair est  $12\hspace{1pt}\mathrm{h}(L)-\iota(\mathrm{R}(L))$.}

\bigskip
\textit{Démonstration.} Soient $R$ un système de racines, $C$ une chambre de $R$, $B\subset R$ la base et $R_{+}\subset R$ le sous-ensemble des racines positives, définis par $C$~;  la démonstration de 4.4.3 est conséquence de la relation, dégagée par Bertram Kostant dans \cite{kostantsl2}, entre la fonction hauteur $\mathrm{R}_{+}\to\mathbb{N}-{0}$ et les exposants du groupe de Weyl de $R$. Nous rappelons cette théorie ci-dessous.

\medskip
La fonction {\em hauteur}, disons $\mathrm{H}:R_{+}\to\mathbb{N}-\{0\}$, associe à une racine positive la somme de ses coordonnées dans la base $B$~; nous notons $\mathrm{Exp}(R)$ l'ensemble des exposants de $R$.

\medskip
Soient $A$ un groupe abélien et $f:\mathbb{N}-\{0\}\to A$ une application (ensembliste). Soit $F:\mathbb{N}-\{0\}\to A$ la ``primitive'' de $f$, c'est-à-dire l'application définie par
$$
\hspace{24pt}
F(m)
\hspace{4pt}=\hspace{4pt}
\sum_{k=1}^{m}\hspace{4pt}f(k)
\hspace{24pt}.
$$
Alors on a dans $A$ l'égalité suivante~:
$$
\hspace{24pt}
\sum_{\alpha\in R_{+}}\hspace{4pt}f(\mathrm{H}(\alpha))
\hspace{4pt}=\hspace{4pt}
\sum_{m\in\mathrm{Exp}(R)}\hspace{4pt}F(m)
\hspace{24pt}.
\leqno{(\mathrm{Ko})}
$$
Cette égalité est conséquence de la relation de Kostant évoquée plus haut. Expliquons pourquoi. Soient $i$ un élément de $\mathbb{N}-\{0\}$ et $\delta_{(i)}:\mathbb{N}-\{0\}\to\nolinebreak\mathbb{Z}$\linebreak la ``masse de Dirac'' correspondante. L'égalité ci-dessus dit, dans le cas particulier $f=\delta_{(i)}$, que le cardinal de $\mathrm{H}^{-1}(i)$ est égal au cardinal du sous-ensemble de $\mathrm{Exp}(R)$ constitué des $m$ avec $m\geq i$~; ceci est le résultat de Kostant. Le cas général en résulte par linéarité.

\medskip
On revient maintenant à la démonstration de 4.4.3. Soit $\nu$ (resp. $\nu_{+}$) le nombre de racines (resp. de racines positives, pour la chambre associée à~$\rho$) $\alpha$ de $L$ avec $\rho.\alpha$ pair~; il est clair que l'on a $\nu=2\hspace{1pt}\nu_{+}$. Dans le contexte de 4.4.3, on a $\mathrm{H}(\alpha)=\rho.\alpha$~; en prenant pour $f$, dans l'égalité (Ko), la fonction $\mathbb{N}-\{0\}\to\mathbb{Z},k\mapsto (-1)^{k}$, on obtient ~:
$$
\hspace{24pt}
\vert\mathrm{R}_{+}(L;\rho)\vert-2\hspace{2pt}\nu_{+}
\hspace{4pt}=\hspace{4pt}
\iota(\mathrm{R}(L))
\hspace{24pt},
$$
$\mathrm{R}_{+}(L;\rho)$ désignant le sous-ensemble constitué des racines positives pour la chambre associée à $\rho$, soit encore
$$
\hspace{24pt}
\nu
\hspace{4pt}=\hspace{4pt}
\frac{\vert\mathrm{R}(L)\vert}{2}
-\iota(\mathrm{R}(L))
\hspace{24pt}.
$$
Or on a $\vert\mathrm{R}(L)\vert=24\hspace{1pt}\mathrm{h}(L)$ d'après II.3.3 (c).
\hfill$\square$

\bigskip
Soient $\Lambda$ un réseau entier et $k\geq 0$ un entier, on note $\mathrm{r}_{k}(\Lambda)$ le nombre de représentations de $k$ par $\Lambda$, c'est-à-dire le nombre d'éléments $x$ de $\Lambda$ avec $x.x=k$.

\bigskip
\textbf{Proposition 4.4.4} (Borcherds)\textbf{.} {\em Soit $B$ un réseau unimodulaire impair de dimension $24$ ; soient $L_{1}$ et $L_{2}$ les deux $2$-voisins pairs de $B$. Alors on~a~:
$$
\hspace{24pt}
\mathrm{r}_{2}(L_{1})+ \mathrm{r}_{2}(L_{2})
\hspace{4pt}=\hspace{4pt}
3\hspace{2pt}\mathrm{r}_{2}(B)
\hspace{2pt}-\hspace{2pt}
24\hspace{2pt}\mathrm{r}_{1}(B)
\hspace{2pt}+\hspace{2pt}
48
\hspace{24pt}.
$$}

\bigskip
\textit{Démonstration (esquisse).} Soient $n>0$ un entier divisible par $8$ et $B$ un réseau unimodulaire de dimension $n$~; on considère la série thêta
$$
\vartheta_{B}(\tau)
\hspace{4pt}=\hspace{4pt}
\sum_{x\in B} e^{\imath\pi\tau x.x}
$$
($\tau$ dans le demi-plan de Poincaré). La fonction $\vartheta_{B}$ est une forme modulaire de poids $\frac{n}{2}$ pour le sous-groupe, disons $\Gamma'$, du groupe $\Gamma:=\mathrm{SL}_{2}(\mathbb{Z})/\{\pm\mathrm{I}\}$, engendré par les transformations $\tau\mapsto\tau+2$ et $\tau\mapsto\frac{-1}{\tau}$~; elle est en outre modulaire pour $\Gamma$ si $B$ est pair. On note respectivement $\mathrm{M}_{\frac{n}{2}}(\Gamma)$ et $\mathrm{M}_{\frac{n}{2}}(\Gamma')$ les $\mathbb{C}$-espaces vectoriels constitués des formes modulaires de poids $\frac{n}{2}$ pour les groupes $\Gamma$ et $\Gamma'$~; comme $\Gamma'$ est d'indice fini (à savoir $3$) dans $\Gamma$ on dispose d'un homomorphisme de transfert, disons $\mathrm{tr}:\mathrm{M}_{\frac{n}{2}}(\Gamma')\to\mathrm{M}_{\frac{n}{2}}(\Gamma)$. La démonstration de l'énoncé suivant est laissée en exercice au lecteur.

\bigskip
\textbf{Lemme 4.4.5.} {\em Soient $n>0$ un entier divisible par $8$ et $B$ un réseau unimodulaire impair de dimension $n$ ; soient $L_{1}$ et $L_{2}$ les deux $2$-voisins pairs de~$B$. Alors on~a~:
$$
\hspace{24pt}
\vartheta_{L_{1}}+\vartheta_{L_{2}}
\hspace{4pt}=\hspace{4pt}
\mathrm{tr}(\vartheta_{B})
\hspace{24pt}.
$$}

On vérifie que $\mathcal{B}:=(\hspace{1pt}\mathbb{E}_{4}^{3}\hspace{1pt},\hspace{1pt}\mathbb{E}_{4}^{2}\hspace{1pt}\vartheta_{\mathrm{I}_{8}}\hspace{1pt},\hspace{1pt}\mathbb{E}_{4}\hspace{1pt}\vartheta_{\mathrm{I}_{16}}\hspace{1pt},\hspace{1pt}\Delta\hspace{1pt})$ est une base de $\mathrm{M}_{12}(\Gamma')$ (on rappelle que la notation $\mathbb{E}_{4}$ désigne la série d'Eiseinstein normalisée, modulaire de poids~$4$ pour $\Gamma$, que l'on a $\vartheta_{\mathrm{E}_{8}}=\mathbb{E}_{4}$, que $\Delta$ est l'unique forme parabolique normalisée de poids $12$ pour $\Gamma$ et que $(\hspace{1pt}\mathbb{E}_{4}^{3}\hspace{1pt},\hspace{1pt}\Delta\hspace{1pt})$ est une base de $\mathrm{M}_{12}(\Gamma)$).

\medskip
Puisque $\mathbb{E}_{4}$ et $\Delta$ sont modulaires pour $\Gamma$ on a $\mathrm{tr}(\mathbb{E}_{4}^{3})=3\hspace{1pt}\mathbb{E}_{4}^{3}$ et $\mathrm{tr}(\Delta)=\nolinebreak3\hspace{1pt}\Delta$~; d'autre part le lemme 4.4.5 implique $\mathrm{tr}(\vartheta_{\mathrm{I}_{8}})=2\hspace{1pt}\mathbb{E}_{4}$ et $\mathrm{tr}(\vartheta_{\mathrm{I}_{16}})=2\hspace{1pt}\mathbb{E}_{4}^{2}$. Comme le transfert est $\mathrm{M}(\Gamma)$-linéaire, $\mathrm{M}(\Gamma)$ désignant la $\mathbb{C}$-algèbre graduée des formes modulaires pour $\Gamma$, on constate au bout du compte que l'image de la base $\mathcal{B}$, par l'homomorphisme de transfert, est $(\hspace{1pt}3\hspace{1pt}\mathbb{E}_{4}^{3}\hspace{1pt},\hspace{1pt}2\hspace{1pt}\mathbb{E}_{4}^{3}\hspace{1pt},\hspace{1pt}2\hspace{1pt}\mathbb{E}_{4}^{3}\hspace{1pt},\hspace{1pt}3\hspace{1pt}\Delta\hspace{1pt})$.

\medskip
Soit maintenant $B$ un réseau unimodulaire impair de dimension $24$. Soient $(c_{0},c_{1},c_{2},c_{3})$ les coordonnées de $\vartheta_{B}$ dans la base $\mathcal{B}$~:
$$
\hspace{24pt}
\vartheta_{B}=c_{0}\hspace{1pt}\mathbb{E}_{4}^{3}\hspace{1pt}+\hspace{1pt}c_{1}\hspace{1pt}\mathbb{E}_{4}^{2}\hspace{1pt}\vartheta_{\mathrm{I}_{8}}\hspace{1pt}+\hspace{1pt}c_{2}\hspace{1pt}\mathbb{E}_{4}\hspace{1pt}\vartheta_{\mathrm{I}_{16}}\hspace{1pt}+\hspace{1pt}c_{3}\hspace{1pt}\Delta
\hspace{24pt}.
$$
Puisque le terme constant du développement de Fourier de $\vartheta_{B}$ est $1$, on a $c_{0}+c_{1}+c_{2}=1$. Comme $B$ est impair, le lemme 4.4.5 dit que le terme constant du développement de Fourier de $\mathrm{tr}(\vartheta_{B})$ est $2$~; on a donc $3\hspace{1pt}c_{0}+2\hspace{1pt}c_{1}+2\hspace{1pt}c_{2}=2$. On en déduit $c_{0}=0$.

\medskip
On note $\mathrm{M}_{12}^{0}(\Gamma')$ le sous-espace vectoriel de $\mathrm{M}_{12}(\Gamma')$ engendré par $\mathbb{E}_{4}^{2}\hspace{1pt}\vartheta_{\mathrm{I}_{8}}$, $\mathbb{E}_{4}\hspace{1pt}\vartheta_{\mathrm{I}_{16}}$ et $\Delta$. Soient $f$ un élément de $\mathrm{M}_{12}^{0}(\Gamma')$ et
$$
f\hspace{4pt}=\hspace{4pt}
\mathrm{r}_{0}(f)
\hspace{4pt}+\hspace{4pt}
\mathrm{r}_{1}(f)\hspace{2pt}e^{\imath\pi\tau}
\hspace{4pt}+\hspace{4pt}
\mathrm{r}_{2}(f)\hspace{2pt}e^{2\imath\pi\tau}
\hspace{4pt}+\hspace{4pt}\ldots
$$
le début de son développement  de Fourier~; on constate que l'application linéaire $\mathrm{M}_{12}^{0}(\Gamma')\to\mathbb{C}^{3}\hspace{2pt},\hspace{2pt}f\mapsto(\mathrm{r}_{0}(f),\mathrm{r}_{1}(f),\mathrm{r}_{2}(f))$ est un isomorphisme. Il en résulte que le coefficient de $e^{2\imath\pi\tau}$ dans le développement  de Fourier de $\mathrm{tr}(f)$ s'exprime linéairement en fonction de $(\mathrm{r}_{0}(f),\mathrm{r}_{1}(f),\mathrm{r}_{2}(f))$. En résolvant un système linéaire on trouve que ce coefficient est $48\hspace{1pt}\mathrm{r}_{0}(f)-24\hspace{1pt}\mathrm{r}_{1}(f)+3\hspace{1pt}\mathrm{r}_{2}(f)$. Compte tenu de 4.4.5, on obtient 4.4.4 en prenant $f=\vartheta_{B}$.
\hfill$\square$

\bigskip
\textit{Démonstration de 4.4.2 à l'aide de 4.4.3 et 4.4.4}

\medskip
Soit $B$ le réseau unimodulaire impair de dimension $24$ dont les deux $2$-voisins pairs sont $L$ et $\mathrm{vois}_{2}(L;\rho)$. Comme l'on a $\mathrm{r}_{2}(\Lambda)=24\hspace{1pt}\mathrm{h}(\Lambda)$ pour tout réseau unimodulaire pair $\Lambda$ de dimension $24$ (point (c) de II.3.3 pour $\mathrm{r}_{2}(\Lambda)\not=0$ et convention pour $\mathrm{r}_{2}(\Lambda)=0$) la proposition 4.4.4 donne
$$
\hspace{24pt}
24\hspace{2pt}\mathrm{h}(L)+24\hspace{2pt}\mathrm{h}(\mathrm{vois}_{2}(L;\rho))
\hspace{4pt}=\hspace{4pt}
3\hspace{2pt}\mathrm{r}_{2}(B)
\hspace{2pt}-\hspace{2pt}
24\hspace{2pt}\mathrm{r}_{1}(B)
\hspace{2pt}+\hspace{2pt}
48
\hspace{24pt}.
\leqno{(1)}
$$
Par construction le sous-module de $B$ constitué des éléments $x$ avec $x.x$ pair est le réseau $\mathrm{M}_{2}(L;\rho)$ ; on a donc $\mathrm{r}_{2}(B)=\mathrm{r}_{2}(\mathrm{M}_{2}(L;\rho))$. Par construction encore, $\mathrm{r}_{2}(\mathrm{M}_{2}(L;\rho))$ est le nombre de racines $\alpha$ de $L$ avec $\rho.\alpha$ pair~; on a donc
$$
\mathrm{r}_{2}(\mathrm{M}_{2}(L;\rho))
\hspace{4pt}=\hspace{4pt}
12\hspace{2pt}\mathrm{h}(L)-\iota(\mathrm{R}(L))
\leqno{(2)}
$$
d'après 4.4.3. Les égalités (1) et (2) entraînent
$$
\hspace{24pt}
\mathrm{h}(\mathrm{vois}_{2}(L;\rho))
\hspace{4pt}=\hspace{4pt}
\frac{\mathrm{h}(L)}{2}-\frac{\iota(\mathrm{R}(L))}{8}+2-\mathrm{r}_{1}(B)
\hspace{24pt}.
\leqno{(3)}
$$
Il reste à montrer $\mathrm{r}_{1}(B)=0$. On procède par l'absurde. Si l'on a $\mathrm{r}_{1}(B)\not=0$ alors on a $\mathrm{r}_{1}(B)\geq 2$ et l'égalité (3) implique l'inégalité $\mathrm{h}(\mathrm{vois}_{2}(L;\rho))
<\nolinebreak\mathrm{h}(L)$ qui montre que les réseaux $L$ et $\mathrm{vois}_{2}(L;\rho)$ ne sont pas isomorphes. Or le corollaire 1.16 montre que si l'on a $\mathrm{r}_{1}(B)\not=0$ alors les réseaux $L$ et $\mathrm{vois}_{2}(L;\rho)$ sont isomorphes.
\hfill$\square\square$

\vspace{0,75cm}
\textit{Commentaires}

\medskip
On note $\mathrm{XR}_{24}$ le sous-ensemble de $\mathrm{X}_{24}$ constitué des classes d'isomorphisme de réseaux unimodulaires pairs de dimension $24$ avec racines. On note $\mathrm{YR}_{24}$ l'ensemble des classes d'isomorphisme de systèmes de racines équicoxeter de rang $24$~; on note $\mathrm{Y}_{24}$ la réunion disjointe de $\mathrm{YR}_{24}$ et du singleton $\{\emptyset\}$. D'après les points (a) et (b) de II.3.3, l'application $L\mapsto\mathrm{R}(L)$ induit des applications $\mathrm{XR}_{24}\to\mathrm{YR}_{24}$ et $\mathrm{X}_{24}\to\mathrm{Y}_{24}$, la seconde prolongeant la première, que l'on note encore $\mathrm{R}$. On oublie ci-dessous que l'on sait que l'application $\mathrm{R}:\mathrm{X}_{24}\to\mathrm{Y}_{24}$ est bijective.

\medskip
Soit $L$ un réseau unimodulaire pair de dimension $24$ et $\rho$ un vecteur de Weyl de $L$. Comme le groupe $\mathrm{W}(\mathrm{R}(L))$ permute transitivement les vecteurs de Weyl de $L$, l'application $L\mapsto\mathrm{vois}_{2}(L;\rho)$ induit une application $\mathrm{XR}_{24}\to\mathrm{X}_{24}$ que l'on note $\varphi$.

\medskip
Soit $R$ un élément de $\mathrm{YR}_{24}$, on pose
$$
\hspace{24pt}
\mathrm{h}'(R)
\hspace{4pt}=\hspace{4pt}
\frac{\mathrm{h}(R)}{2}-\frac{\iota(R)}{8}+2
\hspace{24pt};
$$
on vérifie que $\mathrm{h}'(R)$ appartient à $\mathbb{N}$. La proposition 4.4.2 dit que l'on a l'égalité suivante~:
$$
\hspace{24pt}
\mathrm{h}(\mathrm{R}(\varphi([L]))
\hspace{4pt}=\hspace{4pt}
\mathrm{h}'(\mathrm{R}([L]))
\hspace{24pt}.
\leqno{(1)}
$$
Cette égalité suffit à déterminer $\mathrm{R}(\varphi([L]))$ si $\mathrm{h}'(\mathrm{R}(L))$ est différent de $12$, $10$ ou $6$. En effet, les fibres $\mathrm{h}^{-1}(k)$ de l'application $\mathrm{h}:\mathrm{X}_{24}\to\mathbb{N}$ ont $0$ ou $1$ élément sauf pour $k=12,10,6$, auxquels cas on a $\mathrm{h}^{-1}(12)=\{\mathbf{A}_{11}\mathbf{D}_{7}\mathbf{E}_{6},4\mathbf{E}_{6}\}$, $\mathrm{h}^{-1}(10)=\{2\mathbf{A}_{9}\hspace{1pt}\mathbf{D}_{6},4\mathbf{D}_{6}\}$ et $\mathrm{h}^{-1}(6)=\{4\mathbf{A}_{5}\hspace{1pt}\mathbf{D}_{4},6\mathbf{D}_{4}\}$. On vérifie par inspection que $\mathrm{h}'^{-1}(12)$ est vide et que l'on a $\mathrm{h}'^{-1}(10)=\{2\mathbf{D}_{12}\}$ et $\mathrm{h}'^{-1}(6)=\{3\mathbf{D}_{8},\mathbf{A}_{11}\mathbf{D}_{7}\mathbf{E}_{6},4\mathbf{E}_{6}\}$. Dans les cas $R=2\mathbf{D}_{12},3\mathbf{D}_{8},\mathbf{A}_{11}\mathbf{D}_{7}\mathbf{E}_{6},4\mathbf{E}_{6}$, on détermine $\mathrm{R}(\varphi([L]))$ à l'aide de la condition (2) ci-après.

\medskip
Soit $R$ un système de racines de type  ADE muni d'une chambre $C$, ou ce qui revient au même d'un vecteur de Weyl $\rho$. On note $R/2$ le sous-système de racines de $R$ constitué des racines de hauteur paire, pour la fonction hauteur définie par $C$. Le système de racines $R/2$ est encore de type ADE~; on observera que $R/2$ est canoniquement muni d'une chambre~: les racines positives pour cette chambre sont celles qui le sont pour $C$. Il est clair que la classe d'isomorphisme de $R/2$  est indépendante du choix de $C$. Au niveau des classes d'isomorphisme, l'application $R\mapsto R/2$ est déterminée par les deux propriétés ci-dessous~:

\smallskip
-- On a l'égalité $(R_{1}\coprod R_{2})/2=R_{1}/2\hspace{1pt}\coprod\hspace{1pt}R_{2}/2$.

\smallskip
--  Pour $R$ irréductible le système de racines $R/2$ est le suivant~: $\mathbf{A}_{2m}/2=\mathbf{A}_{m}\coprod\mathbf{A}_{m-1}$, $\mathbf{A}_{2m+1}/2=\mathbf{A}_{m}\coprod\mathbf{A}_{m}$, $\mathbf{D}_{2m}/2=\mathbf{D}_{m}\coprod\mathbf{D}_{m}$, $\mathbf{D}_{2m+1}/2=\linebreak\mathbf{D}_{m+1}\coprod\mathbf{D}_{m}$, $\mathbf{E}_{6}/2=\mathbf{A}_{5}\coprod\mathbf{A}_{1}$, $\mathbf{E}_{7}/2=\mathbf{A}_{7}$, $\mathbf{E}_{8}/2=\mathbf{D}_{8}$ (avec les conventions $\mathbf{A}_{0}=\emptyset$, $\mathbf{D}_{2}=\mathbf{A}_{1}\coprod\mathbf{A}_{1}$ et $\mathbf{D}_{3}=\mathbf{A}_{3}$).

\medskip
On a par définition $\mathrm{R}(L)/2=\mathrm{R}(\mathrm{M}_{2}(L;\rho))$ et donc
$$
\mathrm{R}([L])/2
\hspace{4pt}\subset\hspace{4pt}
\mathrm{R}(\varphi([L]))
\leqno{(2)}
$$
pour tout $[L]$ dans $\mathrm{XR}_{24}$. Cette inclusion permet de déterminer $\mathrm{R}(\varphi([L]))$ pour $\mathrm{h}'(\mathrm{R}([L]))\in\{10,6\}$. En effet, on a~:

\smallskip
--\hspace{8pt}$(2\mathbf{D}_{12})/2=4\mathbf{D}_{6}\not\subset 2\mathbf{A}_{9}\hspace{1pt}\mathbf{D}_{6}$~;

\smallskip
--\hspace{8pt}$(3\mathbf{D}_{8})/2=6\mathbf{D}_{4}\not\subset 4\mathbf{A}_{5}\hspace{1pt}\mathbf{D}_{4}$~;

\smallskip
--\hspace{8pt}$(\mathbf{A}_{11}\mathbf{D}_{7}\mathbf{E}_{6})/2=\mathbf{A}_{1}\hspace{1pt}\mathbf{A}_{3}\hspace{1pt}\mathbf{A}_{4}\hspace{1pt}3\mathbf{A}_{5}\not\subset 6\mathbf{D}_{4}$~;

\smallskip
--\hspace{8pt}$(4\mathbf{E}_{6})/2=4\mathbf{A}_{1}\hspace{1pt}4\mathbf{A}_{5}\not\subset 6\mathbf{D}_{4}$.

\bigskip
La discussion précédente conduit aux énoncés suivants~:

\bigskip
\textbf{Proposition-Définition 4.4.6.} {\em Il existe une unique application
$$
\hspace{24pt}
\psi:\mathrm{YR}_{24}\to\mathrm{Y}_{24}
\hspace{24pt},
$$
telle que l'on a $\mathrm{h}(\psi(R))=\mathrm{h}'(R)$ et $R/2\subset\psi(R)$.}
\vfill\eject

\bigskip
\footnotesize
\textit{Remarque.} Soit $R$ un élément de $\mathrm{YR}_{24}$, les conditions suivantes sont équivalentes~:

\smallskip
-- $\psi(R)=R/2$~;

\smallskip
-- $\iota(R)=24$~;

\smallskip
-- les composantes irréductibles de $R$ sont de type $\mathbf{A}_{1}$, $\mathbf{D}_{l}$ avec $l$ pair, $\mathbf{E}_{7}$ ou $\mathbf{E}_{8}$.
\normalsize

\medskip
On prolonge $\psi$ en une application $\psi:\mathrm{Y}_{24}\to\mathrm{Y}_{24}$ en posant $\psi(\emptyset)=\emptyset$~; pareillement, on prolonge $\varphi$ en une application $\varphi:\mathrm{X}_{24}\to\mathrm{X}_{24}$ en posant $\varphi([L])=[L]$ si $L$ est sans racines.

\bigskip
\textbf{Proposition 4.4.7.} {\em Soient $R$ un élément de $Y_{24}$ et $k$ un entier naturel. On a $\psi^{k}(R)=\emptyset$ pour $\mathrm{h}(R)<2^{k}+1$.}

\medskip
(Observer que l'on a $\mathrm{h}(\psi(R))-1\leq\frac{1}{2}\hspace{2pt}(\mathrm{h}(R)-1)$.)

\bigskip
\textbf{Proposition 4.4.8.} {\em Le diagramme
$$
\begin{CD}
\mathrm{X}_{24}@>\varphi>>\mathrm{X}_{24} \\
@V\mathrm{R}VV @V\mathrm{R}VV \\
\mathrm{Y}_{24}@>\varphi>>\mathrm{Y}_{24} \\
\end{CD}
$$
est commutatif.}

\bigskip
\textbf{Scholie 4.4.9.} {\em Soient $L$ un réseau unimodulaire pair de dimension $24$ et $k$ un entier naturel. Alors $\varphi^{k}([L])$ est sans racines pour $\mathrm{h}(L)<2^{k}+1$.}

\bigskip
Le graphe orienté de la figure III.1 explicite l'application $\psi$~: ses sommets sont les éléments de $\mathrm{Y}_{24}$ et ses arêtes les couples $(x,y)$ de $\mathrm{Y}_{24}\times\mathrm{Y}_{24}-\mathrm{diagonale}$, avec $x\in\mathrm{YR}_{24}$ et $y=\psi(x)$.


\parindent=0.5cm

\chapter{Formes automorphes et op\'erateurs de
Hecke}\label{fautetophecke}\label{chap4}

\section{R\'eseaux et ensembles de classes des $\Z$-groupes} \label{ensclass}

On note {\rm P} l'ensemble des nombres premiers.  Soient $\widehat{\Z}=
\prod_{p \in {\rm P}} \Z_p$ et $\AAA_f=\Q \otimes \widehat{\Z}$ l'anneau des
ad\`eles finis de $\Q$.  Fixons $G$ un $\Z$-groupe, c'est-\`a-dire un
sch\'ema en groupes affine de type fini sur $\Z$.  Le groupe $G(\AAA_f)$
s'identifie canoniquement au sous-groupe de $\prod_{p \in {\rm P}} G(\Q_p)$
dont les \'el\'ements $(g_p)$ v\'erifient $g_p \in G(\Z_p)$
pour {\it presque tout} $p$, c'est-\`a-dire pour tout $p \in {\rm P}$ sauf \'eventuellement
un nombre fini d'entre eux.  Les groupes $G(\Q)$ et $G(\widehat{\Z})$ se plongent
naturellement dans $G(\AAA_f)$ et satisfont $G(\widehat{\Z})=\prod_{p \in
{\rm P}} G(\Z_p)\, \, \, \, \, {\rm et} \, \, \, \, G(\Z)=G(\Q) \cap
G(\widehat{\Z})$.  Le $G(\AAA_f)$-ensemble
$$\mathcal{R}(G)=G(\AAA_f)/G(\widehat{\Z})$$ jouera un r\^ole important dans
ce chapitre.  Nous le notons $\mathcal{R}$, comme ``r\'eseau'', car il
s'identifie en g\'en\'eral \`a l'ensemble des r\'eseaux d'un certain type
d'un $\Q$-espace vectoriel.  \ps\ps

Un r\'esultat classique de Borel~\cite[\S
5]{borelfini} assure que l'{\it ensemble des classes de $G$} :
$${\rm Cl}(G)=G(\Q)\backslash G(\AAA_f) / G(\widehat{\Z})=G(\Q)\backslash
\mathcal{R}(G)$$ est fini.  Son cardinal $\mathrm{h}(G)=|{\rm Cl}(G)|$ est
appel\'e {\it nombre de classes} de $G$.  Dans ce qui suit nous d\'ecrivons
$\mathcal{R}(G)$ et ${\rm Cl}(G)$ dans les quelques cas standards qui vont nous
int\'eresser (voir par exemple~\cite[\S 2]{borelfini}).  \ps\ps

\subsection{Groupes lin\'eaires}\label{ensclassgln}

Commen\c{c}ons par le cas de $\GL_n$. Si $V$ est un espace vectoriel
de dimension finie $n$ sur le
corps des fractions d'un anneau principal $A$, nous d\'esignerons par
$\mathcal{R}_A(V)$ l'ensemble des r\'eseaux de $V$ relativement \`a $A$, c'est-\`a-dire des sous-$A$-modules libres de rang $n$ de $V$
(\S II.1). Il est muni d'une action transitive de $\GL(V)$, le stabilisateur de $L \in
\mathcal{R}_A(V)$ \'etant $\GL(L)$. \ps\ps

Soit $V$ un $\Q$-espace vectoriel de dimension $n$. Si $p$ est premier,
et si l'on pose $V_p=V \otimes \Q_p$, on dispose d'une application naturelle $\mathcal{R}_\Z(V) \rightarrow \mathcal{R}_{\Z_p}(V_p)$, $M \mapsto M_p : =  M \otimes \Z_p$. 
Fixons $L \in \mathcal{R}_\Z(V)$ et posons $G=\GL_L$. Il est \'el\'ementaire de
v\'erifier, suivant
Eichler~\cite[\S 13]{eichlerqf}, que l'application
\begin{equation}\label{plongementdeichler} \mathcal{R}_\Z(V) \rightarrow \prod_{p \in
{\rm P}} \mathcal{R}_{\Z_p}(V_p), \, \, \, \, \, M \mapsto (M _p),\end{equation}
est une injection de $\mathcal{R}_\Z(V)$ sur le sous-ensemble $\prod_{p \in {\rm P}}' \mathcal{R}_{\Z_p}(V_p) \subset \prod_{p \in
{\rm P}} \mathcal{R}_{\Z_p}(V_p)$ constitu\'e des 
familles $(M_p)$ telles que $M_p=L_p$ pour
presque tout $p$ (ce sous-ensemble ne d\'ependant d'ailleurs pas du choix de $L$). 
L'action naturelle de $G(\AAA_f)$ sur $\prod_{p \in {\rm
P}} \mathcal{R}_{\Z_p}(V_p)$ pr\'eserve $\prod_{p \in {\rm P}}' \mathcal{R}_{\Z_p}(V_p)$, et elle est transitive sur ce dernier. Ainsi, si l'on identifie $\mathcal{R}_\Z(V)$ \`a  $\prod_{p \in {\rm
P}}' \mathcal{R}_{\Z_p}(V_p)$ au moyen de l'application \eqref{plongementdeichler}, ce que l'on fera syst\'ematiquement par la suite, on obtient par transport de structure une action transitive de $G(\AAA_f)$ sur $\mathcal{R}_\Z(V)$ \'etendant l'op\'eration \'evidente de $G(\Q)=\GL(V)$.  Le stabilisateur du r\'eseau $L$ \'etant $G(\widehat{\Z})$, on en d\'eduit un isomorphisme de
$G(\AAA_f)$-ensembles $$\mathcal{R}(G) \isomo \mathcal{R}_\Z(V)$$ Comme $G(\Q)$ agit aussi
transitivement sur $\mathcal{R}_\Z(V)$, il vient en particulier que $$\mathrm{h}({\rm GL}_n)=1.$$ Quand
$G={\rm PGL}_L$ (resp.  $G={\rm SL}_L$), $\mathcal{R}(G)$ s'interpr\`ete \'egalement
comme le quotient de $\mathcal{R}_\Z(V)$ par $\Q^\times$ pour l'action par homoth\'eties
(resp.  comme le sous-ensemble de $\mathcal{R}_\Z(V)$ constitu\'e des $M$ poss\'edant
une $\Z$-base de d\'eterminant $1$ relativement \`a une $\Z$-base de $L$). 
On a encore $\mathrm{h}({\rm PGL}_n)=\mathrm{h}({\rm SL}_n)=1$.\ps\ps

\subsection{Groupes orthogonaux et
symplectiques}\label{ensclassbil}\label{hsimilitudes}

Supposons de plus que le $\Q$-espace vectoriel $V$ est muni d'une forme
bilin\'eaire non d\'eg\'en\'er\'ee $\varphi$, suppos\'ee sym\'etrique ou
altern\'ee. Soit $L \in \mathcal{R}_\Z(V)$. Rappelons que son {\it r\'eseau dual} est
le r\'eseau $L^\sharp \in \mathcal{R}_\Z(V)$ d\'efini par (\S II.1) $$L^\sharp = \{ v \in V,
\varphi(v,x) \in \Z, \forall x \in L\}.$$  Nous dirons que $L$
est {\it homodual}, pour ``homoth\'etique \`a son dual'', s'il existe $\lambda \in \Q^\times$ tel que
$L^\sharp = \lambda L$; il existe alors un unique tel $\lambda$ qui soit $ >0$, que l'on
notera $\lambda_L$.  Le r\'eseau $L$ est dit {\it autodual} si $L^\sharp=L$. Si
$L$ est homodual, et si $\varphi$ est sym\'etrique (resp. 
altern\'ee), la forme bilin\'eaire $\lambda_L \varphi$ munit $L$ d'une
structure de ${\rm b}$-module (resp.  ${\rm a}$-module) sur $\Z$ au sens du
\S II.1.  On dira alors que $L$ est {\it pair} si $\lambda_L \varphi(x,x) \in 2 \Z$ pour tout $x \in
L$.  C'est automatique si $\varphi$ est altern\'ee, et si $\varphi$ est
sym\'etrique cela permet de voir $L$ comme un ${\rm q}$-module sur $\Z$ en
posant ${\rm q}(x)=\lambda_L \, \, \frac{\varphi(x,x)}{2}$ pour $x \in L$.
Notons $$\mathcal{R}_\Z^{\rm a}(V) \subset \mathcal{R}_\Z^{\rm h}(V)$$ les sous-ensembles de $\mathcal{R}_\Z(V)$
constitu\'es des r\'eseaux autoduaux (resp.  homoduaux)
{\underline{pairs}}. \ps\ps

Posons $n=\dim V$. Fixons $L \in \mathcal{R}_{\Z}^{\rm a}(V)$.  L'existence d'un tel $L$
entra\^ine, par r\'eduction modulo $2$, la congruence $n \equiv 0 \bmod 2$. Consid\'erons le sous-$\Z$-groupe $G \subset \GL_L$ d\'efini par
: $$ G = \left\{ \begin{array}{ll} {\rm Sp}_L & {\rm si}\,\, \varphi \,\,
{\mathrm{est\, \, altern\acute{e}e}}, \\ {\rm O}_L & {\rm
sinon.}\end{array}\right.$$ Nous noterons \'egalement $
\widetilde{G}$ le $\Z$-groupe de similitudes correspondant, de sorte que $G \subset \widetilde{G}  \subset \GL_L$,
ainsi que $\mathrm{P}\widetilde{G}$ le $\Z$-groupe des similitudes projectives, quotient de
$\widetilde{G}$ par son sous-$\Z$-groupe central isomorphe \`a
$\mathbb{G}_m$ constitu\'e des homoth\'eties~(\S II.1).

\begin{lemme}\label{lemmeortho}  La restriction \`a $\widetilde{G}(\AAA_f)$ {\rm (}resp. $G(\AAA_f)${\rm )} de l'action de ${\rm GL}_L(\AAA_f)$ sur $\mathcal{R}_\Z(V)$ pr\'eserve
$\mathcal{R}_{\Z}^{\rm h}(V)$ {\rm (}resp. $\mathcal{R}_\Z^{\rm a}(V)${\rm )}. 
\end{lemme}

Avant de proc\'eder \`a la d\'emonstration, introduisons les analogues locaux des d\'efinitions pr\'ec\'edentes. Soit $p$ premier. Si $M \in \mathcal{R}_{\Z_p}(V_p)$, il y a encore un sens
\`a consid\'erer le r\'eseau dual $M^\sharp \in \mathcal{R}_{\Z_p}(V_p)$ (relativement \`a $\Z_p$,
voir~\S II.1). On note $\mathcal{R}_{\Z_p}^{\rm h}(V_p) \subset
\mathcal{R}_{\Z_p}(V_p)$ le sous-ensemble des r\'eseaux $M$ tels qu'il existe $\lambda \in
\Q_p^\times$ tel que $M^\sharp = \lambda M$ et $\lambda
\varphi(x,x) \in 2 \Z_p$ pour tout $x \in M$. 
Notons \'egalement $\mathcal{R}_{\Z_p}^{\rm a}(V_p) \subset         
\mathcal{R}_{\Z_p}^{\rm h}(V_p)$ le sous ensemble des $M$ tels que $M^\sharp = M$. Si $M \in \mathcal{R}_{\Z_p}^{\rm h}(V_p)$, il existe un unique $\lambda_M \in p^\Z$ tel que $M^\sharp = \lambda_M M$. Si $\varphi$ est sym\'etrique (resp. altern\'ee), la forme quadratique $x \mapsto \lambda_M \frac{\varphi(x,x)}{2}$ (resp. la forme altern\'ee $\lambda_M \varphi$) munit alors $M$ d'une structure de ${\rm q}$-module (resp. ${\rm a}$-module) sur $\Z_p$. 

\begin{pf}  Soit $M \in \mathcal{R}_\Z(V)$. Commen\c{c}ons par observer que $M$ est dans $\mathcal{R}_\Z^{\rm h}(V)$  si, et seulement si, $M_p$ est dans
$\mathcal{R}_{\Z_p}^{\rm h}(V_p)$  pour
tout $p$ premier, auquel cas on a de plus $\lambda_M = \prod_p \lambda_{M_p}$ (bien entendu, $\lambda_{M_p}$ vaut $1$ pour presque tout $p$). Cela d\'ecoule en effet de l'identit\'e $\AAA_f^\times=\Q^\times \cdot \widehat{\Z}^\times$ (i.e. $\mathrm{h}(\mathbb{G}_m)=1$) et de la relation imm\'ediate $(N^\sharp)_p = (N_p)^\sharp$, valable pour tout $p$ premier et tout $N \in \mathcal{R}_\Z(V)$. En particulier, $M \in \mathcal{R}_\Z^{\rm a}(V)$ si, et seulement si,  $M_p \in \mathcal{R}_{\Z_p}^{\rm a}(V_p)$ pour tout $p$. \ps
Pour conclure la d\'emonstration, il suffit d'observer que si $g \in \widetilde{G}(\Q_p)$ a pour facteur de similitude $\nu(g)$ (\S II.1), et si $M \in \mathcal{R}_{\Z_p}(V_p)$, on a la relation $g(M)^\sharp = \nu(g)^{-1} g(M^\sharp)$.
\end{pf}

Observons que les homoth\'eties $\Q^\times$ agissent sur $\mathcal{R}_\Z(V)$ en
pr\'eservant $\mathcal{R}_\Z^{\rm h}(V)$.  D'apr\`es le lemme ci-dessus, l'ensemble
quotient $$\underline{\mathcal{R}}_\Z^{\rm h}(V) := \Q^\times \backslash
\mathcal{R}^{\rm h}_\Z(V)$$ est donc muni d'une action de $\mathrm{P}\widetilde{G}(\AAA_f)$ \'etendant l'action
\'evidente de $\mathrm{P}\widetilde{G}(\Q)$. Nous noterons 
$\underline{M}$ la classe d'homoth\'etie d'un $M \in \mathcal{R}_\Z(V)$. Nous disposons au final d'un diagramme commutatif :

$$\xymatrix{ \mathcal{R}(G) \ar@{^{(}->}[r]  \ar@{->}^{\rm \omega_1}[d] &
\mathcal{R}(\widetilde{G}) \ar@{->>}[r] \ar@{->}^{\rm \omega_2}[d] &
\mathcal{R}(\mathrm{P}\widetilde{G})\ar@{->}^{\rm \omega_3}[d] \\
\mathcal{R}_\Z^{\rm a}(V) \ar@{->>}[d] \ar@{^{(}->}[r] & \mathcal{R}_\Z^{\rm
h}(V) \ar@{->>}[d] \ar@{->>}[r] & \underline{\mathcal{R}}_\Z^{\rm h}(V)
\ar@{->>}[d] \\ G(\Q) \backslash \mathcal{R}_\Z^{\rm a}(V) \ar@{->}^{\rm
\xi_1}[r] & \widetilde{G}(\Q) \backslash \mathcal{R}_\Z^{\rm h}(V)
\ar@{->>}^{\xi_2}[r] & \mathrm{P}\widetilde{G}(\Q) \backslash
\underline{\mathcal{R}}_\Z^{\rm h}(V)} $$

Les $\omega_i$, pour $i=1,2,3$, sont respectivement les applications ``orbite'' de $L$, $L$, et $\underline{L}$, sous les actions de $G(\AAA_f)$, $\widetilde{G}(\AAA_f)$ et $\mathrm{P}\widetilde{G}(\AAA_f)$. Toutes les autres fl\`eches d\'esignent les applications canoniques. 

\begin{prop}\label{rortho}\label{diagrammesimilitude} Les applications $\omega_i$ et $\xi_j$ sont bijectives.
En particulier, l'action de $G(\AAA_f)$ sur $\mathcal{R}_{\Z}^{\rm a}(V)$ est transitive, l'orbite de $L$ d\'efinissant un isomorphisme de 
$G(\AAA_f)$-ensembles $\mathcal{R}(G) \isomo \mathcal{R}_{\Z}^{\rm a}(V)$.
\end{prop}

\begin{pf} L'injectivit\'e des $\omega_i$ est \'evidente. Commen\c{c}ons par v\'erifier la derni\`ere assertion, qui n'est
autre que la surjectivit\'e de $\omega_1$. Si $\varphi$ est sym\'etrique, le
scholie~II.2.5 affirme que pour tout $M \in \mathcal{R}_\Z^{\rm
a}(V)$, le ${\rm q}$-module $M_p$ sur $\Z_p$ est hyperbolique. Il est en
particulier isomorphe \`a $L_p$, ce qui conclut car toute isom\'etrie $L_p \rightarrow M_p$ est n\'ecessairement induite par un \'el\'ement de ${\rm O}(V_p)=G(\Q_p)$. Supposons donc $\varphi$ altern\'ee. Il est bien connu que si
$A$ est un anneau principal, il existe \`a \'equivalence pr\`es une et une
seule forme bilin\'eaire altern\'ee non d\'eg\'en\'er\'ee sur le $A$-module $A^n$ ($n$ pair).
On conclut en consid\'erant le cas $A=\Z_p$.  \ps

La surjectivit\'e de $\omega_3$ (resp.  de $\omega_2$) est cons\'equence de celle de $\omega_2$ (resp.  de celles de
$\omega_1$ et $\xi_1$). Montrons la surjectivit\'e de
$\xi_1$. Si $M \in \mathcal{R}_\Z^{\rm h}(V)$, et $g \in \widetilde{G}(\Q)$ a pour facteur de similitude $\nu(g)$, on a
$\lambda_{g(M)} = \pm \, \nu(g)^{-1}\, \lambda_M$. Il suffit donc de voir 
que $\nu(\widetilde{G}(\Q))$ contient l'ensemble $\Q_{>0}$ des rationnels
strictement positifs. C'est \'evident dans le cas altern\'e, et plus g\'en\'eralement si $V$ est hyperbolique. 
Dans le cas sym\'etrique, il s'agit de voir que si $\lambda \in \Q_{>0}$ alors $V$ et $V \otimes \langle \lambda \rangle$ (obtenu en multipliant la forme quadratique sur $V$ par $\lambda$) sont isomorphes en tant que ${\rm q}$-modules sur $\Q$.  Mais ils le sont sur $\Q_p$ pour tout
premier $p$, les $V \otimes \Q_p$ \'etant hyperboliques d'apr\`es le
scholie~II.2.5.  Ils le sont sur $\R$ car $\lambda >0$. 
On conclut par le th\'eor\`eme de Hasse-Minkowski. \ps

L'application $\xi_2$ est bijective car ${\rm P}\widetilde{G}(\Q)=\widetilde{G}(\Q)/\Q^\times$. V\'erifions enfin l'injectivit\'e de $\xi_1$. On peut supposer $\varphi$ sym\'etrique, car l'argument du premier paragraphe montre que $h(G)=1$ si $\varphi$ est altern\'ee. Supposons donc qu'il existe $M \in \mathcal{R}_\Z^{\rm a}(V)$ et $g \in \widetilde{G}(\Q)$ tels que $g(M)=L$. On en d\'eduit $\nu(g) = \pm 1$. Si $\nu(g)=1$ alors $g \in G(\Q)$ : on a termin\'e. Sinon, $M$ est isom\'etrique au ${\rm q}$-module $L \otimes \langle -1 \rangle$, d'espace sous-jacent $L$ mais de forme quadratique oppos\'ee. Cela entra\^ine que $V \otimes \R$ est hyperbolique, puis que $L$ et $M$ sont isomorphes d'apr\`es le th\'eor\`eme II.2.7.
\end{pf}

\begin{cor}\label{classnbpgg} $\mathrm{h}(G)=\mathrm{h}(\widetilde{G})=\mathrm{h}(\mathrm{P}\widetilde{G})$. \end{cor}

Quand $\varphi$ est altern\'ee, la classification sus-mentionn\'ee des formes
altern\'ees non d\'eg\'e\-n\'er\'ees appliqu\'ee \`a l'anneau $\Z$ entra\^ine
\footnote{Les assertions $\mathrm{h}({\rm SL}_n)=\mathrm{h}({\rm Sp}_{2g})=1$
rappel\'ees ci-dessus sont \'egalement des cas tr\`es particuliers du
th\'eor\`eme d'approximation forte de Kneser (voir \cite{kneserboulder} \cite[Thm. 
7.12]{platonovrap}).  Il affirme que $\mathrm{h}(G)=1$ d\`es que le $\C$-groupe
$G_\C$ est semi-simple, simplement connexe, et que le groupe topologique
$G(\R)$ n'a pas de sous-groupe distingu\'e compact connexe non trivial.} 
$\mathrm{h}(G)=1$, puis $\mathrm{h}({\rm Sp}_{2g})=\mathrm{h}({\rm GSp}_{2g})=\mathrm{h}({\rm
PGSp}_{2g})=1$ pour tout $g\geq 1$.  \ps\ps

Supposons $\varphi$  sym\'etrique.  Si le ${\rm q}$-module $L \otimes \R$
n'est pas d\'efini, il r\'esulte de m\^eme du th\'eor\`eme~II.2.7
que $\mathrm{h}({\rm O}_L)=1$.  La situation est bien diff\'erente si $L \otimes \R$
est d\'efini positif, ce que l'on suppose d\'esormais. 
Rappelons que $L$ peut alors \^etre vu comme un r\'eseau unimodulaire pair
de l'espace euclidien $V \otimes \R$ de dimension $n$.  En particulier $n \equiv 0
\bmod 8$.  Dans ce cas, $\mathcal{R}_{\Z}^{\rm a}(V)$ est par d\'efinition l'ensemble
des r\'eseaux unimodulaires pairs de $V \otimes \R$ qui sont inclus dans $L
\otimes \Q$.  Rappelons que ${\rm X}_n$ d\'esigne l'ensemble des classes
d'isom\'etrie de r\'eseaux unimodulaires pairs de l'espace euclidien $V
\otimes \R$.  D'apr\`es le scholie~II.2.1, l'inclusion
naturelle ${\rm O}(V) \backslash \mathcal{R}_\Z^{\rm a}(V) \rightarrow {\rm X}_n$ est bijective,
et induit donc un isomorphisme ${\rm Cl}({\rm O}_L) \isomo {\rm X}_n$. En
particulier, si ${\rm O}_n$ d\'esigne le $\Z$-groupe orthogonal du r\'eseau $L={\rm E}_n$ (\S I.3), on obtient l'\'egalit\'e 
$$\mathrm{h}({\rm O}_n)=|{\rm X}_n|,$$ qui montre que ${\rm h}({\rm O}_n)$ est un nombre tout \`a fait
int\'eressant.  \ps\ps

\subsection{${\rm SO}_L$ versus ${\rm O}_L$}\label{sovso}

On poursuit l'\'etude du paragraphe pr\'ec\'edent en supposant
$\varphi$ sym\'etrique, de sorte que $G,\widetilde{G}$ et $\mathrm{P}\widetilde{G}$ sont
respectivement ${\rm O}_L$, ${\rm GO}_L$ et ${\rm PGO}_L$. On s'int\'eresse
\`a leurs sous-$\Z$-groupes respectifs ${\rm SO}_L$, ${\rm GSO}_L$ et ${\rm
PGSO}_L$ (\S II.1). Les groupes ${\rm SO}_L(\AAA_f)$, ${\rm GSO}_L(\AAA_f)$ et ${\rm PGSO}_L(\AAA_f)$  agissent respectivement sur $\mathcal{R}_\Z^{\rm a}(V)$,
$\mathcal{R}_\Z^{\rm h}(V)$ et $\underline{\mathcal{R}}_\Z^{\rm h}(V)$ (Proposition~\ref{lemmeortho}). Consid\'erons le diagramme commutatif suivant,
\'etendant celui du \S \ref{hsimilitudes}. 
{\scriptsize 
$$ \shorthandoff{!;} \xymatrix@C-1cm@R-.3cm{ & \mathcal{R}({\rm SO}_L) \ar@{->}[rr] \ar@{->}^(0.20){\rm \widetilde{\omega}_1}[dd]|!{[dl];[dr]}\hole
\ar@{->}_{\rm \mu_1}[ld] & & \mathcal{R}({\rm GSO}_L) \ar@{->}[rr] \ar@{->}^(0.2){\rm
\widetilde{\omega}_2} [dd] |!{[dl];[dr]}\hole \ar@{->}_{\mu_2}[ld] & & \mathcal{R}({\rm
PGSO}_L)\ar@{->}^(0.2){\rm \widetilde{\omega}_3}[dd]\ar@{->}_{\mu_3}[ld] \\
\mathcal{R}({\rm O}_L) \ar@{->}[rr] \ar@{->}^(0.2){\rm \omega_1}[dd] & & \mathcal{R}({\rm GO}_L)
\ar@{->}[rr] \ar@{->}^(0.2){\rm \omega_2}[dd] & & \mathcal{R}({\rm
PGO}_L)\ar@{->}^(0.2){\rm \omega_3}[dd] & \\ & \mathcal{R}_\Z^{\rm a}(V)
\ar@{^{(}->}[rr] |!{[ur];[dr]}\hole \ar@{->>}[dd] |!{[dl];[dr]}\hole
\ar@{->}^{\rm Id}[ld] & & \mathcal{R}_\Z^{\rm h}(V) \ar@{->>} [rr] |!{[ur];[dr]}\hole
\ar@{->>}[dd] |!{[dl];[dr]}\hole \ar@{->}^{\rm Id}[ld] & &
\underline{\mathcal{R}}_\Z^{\rm h}(V)\ar@{->>}[dd]\ar@{->}^{\rm Id}[ld] \\ \mathcal{R}_\Z^{\rm
a}(V) \ar@{->}[rr] \ar@{->>}[dd] & & \mathcal{R}_\Z^{\rm h}(V) \ar@{->}[rr]
\ar@{->>}[dd] & & \underline{\mathcal{R}}_\Z^{\rm h}(V) \ar@{->>}[dd] & \\ & {\rm
SO}_L(\Q)\backslash \mathcal{R}_\Z^{\rm a}(V) \ar@{->}_(0.3){\widetilde{\xi}_1}[rr]
|!{[ur];[dr]}\hole \ar@{->>}[ld] & & {\rm GSO}_L(\Q)\backslash \mathcal{R}_\Z^{\rm
h}(V) \ar@{->}_(0.3){\widetilde{\xi}_2}[rr] |!{[ur];[dr]}\hole
\ar@{->>}[ld] & & {\rm PGSO}_L(\Q)\backslash \underline{\mathcal{R}}_\Z^{\rm h}(V)
\ar@{->>}[ld] \\ {\rm O}_L(\Q)\backslash \mathcal{R}_\Z^{\rm a}(V)
\ar@{->}_{\xi_1}[rr] & & {\rm GO}_L(\Q)\backslash \mathcal{R}_\Z^{\rm h}(V)
\ar@{->}_{\xi_2}[rr] & & {\rm PGO}_L(\Q)\backslash \underline{\mathcal{R}}_\Z^{\rm
h}(V) & } $$}

Les applications verticales $\widetilde{\omega}_i$ sont encore les
applications ``orbites'' de $L$ (resp. $L$, resp. $\underline{L}$) et les autres
fl\`eches sont les applications canoniques.

\begin{prop}\label{rso=ro} Les applications $\widetilde{\omega}_i$, $\mu_i$ et $\widetilde{\xi}_j$ sont bijectives. 
En particulier, l'action de ${\rm SO}_L(\AAA_f)$ sur $\mathcal{R}_{\Z}^{\rm a}(V)$ est
encore transitive, l'orbite de $L$ d\'efinissant un isomorphisme $\mathcal{R}({\rm
SO}_L) \isomo \mathcal{R}_{\Z}^{\rm a}(V)$.
\end{prop}

\begin{pf} Nous avons d\'ej\`a vu que l'action naturelle de ${\rm O}_L(\AAA_f)$ sur $\mathcal{R}_{\Z}^{\rm a}(V)$ est transitive (Proposition \ref{rortho}).
Il en va de m\^eme de la restriction de cette action \`a son sous-groupe ${\rm SO}_L(\AAA_f)$ car le groupe orthogonal d'un {\rm q}-module
hyperbolique non trivial sur $\Z_p$ contient toujours un \'el\'ement de
d\'eterminant $-1$.  Cette m\^eme raison montre que les $\mu_i$ sont bijectifs, car  
${\rm O}_L(\Z_p)/{\rm SO}_L(\Z_p) \rightarrow {\rm GO}_L(\Q_p)/{\rm GSO}_L(\Q_p)$ est bijective pour tout premier $p$ (\S II.1). Les $\omega_i$ \'etant bijectifs, on en d\'eduit la bijectivit\'e des $\widetilde{\omega}_i$. \ps
La bijectivit\'e de $\widetilde{\xi}_2$ est \'evidente. La surjectivit\'e de $\widetilde{\xi}_1$ se d\'eduit de celle de $\xi_1$ et de ce que $-1 \in \DET({\rm O}(V))$. Enfin, l'injectivit\'e de $\widetilde{\xi}_1$ se d\'emontre de mani\`ere similaire \`a celle de $\xi_1$ (Proposition \ref{rortho}),  en utilisant que $-1 \in \DET({\rm O}({\rm H}(\Z^{n/2})))$. \end{pf}

\begin{cor} Si $L$ est un ${\rm q}$-module sur $\Z$, on a $\mathrm{h}({\rm SO}_L)=\mathrm{h}({\rm GSO}_L)=\mathrm{h}({\rm PGSO}_L)$. Si de plus $L \otimes \R$ est ind\'efini, ces entiers sont \'egaux \`a $1$.
\end{cor}

\begin{pf} La premi\`ere assertion est cons\'equence de la bijectivit\'e des applications $\xi_{i}$ (Proposition \ref{rso=ro}). Quand $L \otimes \R$ n'est pas
d\'efini, nous avons d\'ej\`a expliqu\'e au \S \ref{ensclassbil} l'\'egalit\'e $\mathrm{h}({\rm O}_L)=1$. Il reste \`a voir qu'il existe $s \in {\rm O}(L)$ tel que $\DET \, s = -1$. L'hypoth\`ese sur $L$ et le th\'eor\`eme II.2.7 montrent qu'il existe un ${\rm q}$-module sur $\Z$, disons $L'$, tel que $L \simeq L' \oplus {\rm H}(\Z)$ (somme orthogonale). On conclut car ${\rm H}(\Z)$ contient un automorphisme de d\'eterminant $-1$. \end{pf}

Supposons finalement $L$ d\'efini positif. On dispose alors comme
pr\'ec\'edemment d'une bijection canonique 
${\rm Cl}({\rm SO}_L) \isomo \widetilde{{\rm X}}_n$, o\`u $\widetilde{{\rm X}}_n$ d\'esigne l'ensemble
des classes d'isom\'etries {\it directes} de r\'eseaux unimodulaires pairs
de $V \otimes \R$ (autrement dit, l'ensemble des orbites de l'action de ${\rm SO}(V \otimes \R)$ sur ces derniers). La classe d'isom\'etrie d'un r\'eseau unimodulaire pair
$M \subset V \otimes \R$ admet exactement $1$ ou $2$ ant\'ec\'edents par la projection canonique
$$\widetilde{{\rm X}}_n \rightarrow {\rm X}_n$$ selon que ${\rm O}(M)$ poss\`ede un
\'el\'ement de d\'eterminant $-1$ ou non. C'est par exemple le cas si $M$ poss\`ede au
moins une racine, c'est-\`a-dire un $\alpha \in M$ tel que $\alpha \cdot
\alpha =2$, car la sym\'etrie orthogonale associ\'ee est dans ${\rm O}(M)$
(\S II.3).  En revanche, si $M$
est le r\'eseau de Leech alors ${\rm O}(M)={\rm SO}(M)$ d'apr\`es Conway~\cite{conwayCo}. Les rappels du~\S II.3 entra\^inent le corollaire suivant. Pour $n\equiv 0 \bmod 8$, on pose ${\rm SO}_n={\rm SO}_{{\rm E}_n}$. 

\begin{cor}\label{hson} On a $\mathrm{h}({\rm SO}_8)=1$, $\mathrm{h}({\rm SO}_{16})=2$ et $\mathrm{h}({\rm
SO}_{24})=25$. \end{cor}

\subsection{Groupes orthogonaux en dimension impaire}\label{soimph}

On se place encore dans le cadre du~\S\ref{ensclassbil}, en supposant  $\varphi$  sym\'etrique. On consid\`ere plut\^ot l'ensemble
$$\mathcal{R}_\Z^{\rm b}(V) \subset \mathcal{R}_\Z(V)$$ constitu\'e des $L
\in \mathcal{R}_\Z(V)$ tels que $\varphi(x,x) \in 2\Z$ pour tout $x \in L$ et tels que $L^\sharp/L \simeq \Z/2\Z$. Cette derni\`ere condition \'equivaut \`a demander que  $\varphi_{|L \times L}$ 
est de d\'eterminant $\pm 2$. Nous renvoyons \`a l'appendice ${\rm B}$ pour une \'etude de ces r\'eseaux. \ps\ps

Fixons $L \in
\mathcal{R}_\Z^{\rm b}(V)$, ce qui n\'ecessite que la dimension $n$ de $V$ soit impaire. Alors ${\rm
SO}_L(\AAA_f)$ agit transitivement sur $\mathcal{R}_\Z^{\rm b}(V)$ d'apr\`es la proposition {\rm B}.2.5, le stabilisateur de
$L$ \'etant ${\rm SO}_L(\widehat{\Z})$. Si $L \otimes \R$ n'est pas d\'efini, le nombre de classes de ${\rm
SO}_L$ est $1$ : c'est un r\'esultat classique qu'il ne serait pas difficile de d\'eduire de la proposition {\rm B}.2.5 (iii) et du th\'eor\`eme II.2.7. La situation est plus int\'eressante si 
$L \otimes \R$ est d\'efini, disons positif pour fixer les id\'ees, ce que l'on suppose d\'esormais. \ps\ps

Dans ce cas, on a la congruence $n \equiv \pm 1 \bmod 8$ et ${\rm Cl}({\rm SO}_L)$ s'identifie \`a l'ensemble des classes d'isom\'etrie de r\'eseaux pairs et de d\'eterminant $2$ dans $\R^n$ (\S ${\rm B}$.2). Il n'y a pas lieu ici de distinguer isom\'etries directes et indirectes car $x \mapsto -x$ est dans ${\rm O}(M)$ et de d\'eterminant $-1$ pour tout $M \in \mathcal{R}_\Z^{\rm b}(V)$. Si $n \equiv 1 \bmod 8$, on pose ${\rm L}_n = {\rm E}_{n-1} \oplus {\rm A}_1$. Si $n \equiv -1 \bmod 8$, on note ${\rm L}_n$ l'orthogonal d'une racine quelconque de ${\rm E}_{n+1}$ ; ces racines \'etant permut\'ees transitivement par le groupe orthogonal de ${\rm E}_{n+1}$, la classe 
d'isom\'etrie d'un tel r\'eseau ne d\'epend d'aucun choix. Si $n \equiv \pm 1 \bmod 8$ le r\'eseau ${\rm L}_n$ est donc pair de d\'eterminant $2$ (\S {\rm B}.2), 
et l'on pose ${\rm SO}_n = {\rm SO}_{{\rm L}_n}$ (\S {\rm B}.1). \ps\ps

Les valeurs connues de ${\mathrm{h}}({\rm SO}_n)$ avec $n$ impair sont r\'esum\'ees dans le corollaire qui suit (voir aussi \cite{conwaysloane}). Les cas $n\leq 23$ sont trait\'es en exemple dans l'appendice ${\rm B}$.2 et le cas $n=25$ est d\^u \`a Borcherds \cite[Table -2]{borcherdsthese}.

\begin{cor}\label{hsonimpair} On a $\mathrm{h}({\rm SO}_1)=\mathrm{h}({\rm SO}_7)=\mathrm{h}({\rm SO}_9)=1$, $\mathrm{h}({\rm SO}_{15})=2$, $\mathrm{h}({\rm SO}_{17})=4$, $\mathrm{h}({\rm SO}_{23})=32$ et $\mathrm{h}({\rm SO}_{25})=121$. 
\end{cor}

\section{Correspondances de Hecke}\label{corrhecke}

\subsection{Formalisme g\'en\'eral}\label{corrheckegeneral}

Soit $\Gamma$ un groupe (abstrait) et soit $X$ un $\Gamma$-ensemble transitif.  L'anneau des {\it correspondances (ou op\'erateurs) de Hecke} de
$X$ est l'anneau $$\mathrm{H}(X) = {\rm End}_{\Z[\Gamma]} (\Z[X]).$$
\`A chaque $T \in {\rm End}_\Z(\Z[X])$ est associ\'e une matrice \`a coefficients entiers
$(T_{x,y})_{(x,y) \in X \times X}$, qui le d\'etermine uniquement, d\'efinie par la formule $$\forall y \in X, \, \, \, \, T(y)=
\sum_{x \in X} T_{x,y} \, \, x.$$
Par d\'efinition, un tel \'el\'ement $T$ est dans ${\rm H}(X)$ si et seulement si la fonction $X \times X \rightarrow \Z, (x,y) \mapsto T_{x,y},$ est constante sur les orbites du groupe $\Gamma$ op\'erant diagonalement sur $X \times X$. La fonction $\Gamma \backslash (X \times X) \rightarrow \Z$ qui en r\'esulte est alors \`a support fini, par finitude de $\{ x \in X,
T_{x,y} \neq 0\}$ pour $y \in X$ et par transitivit\'e de $X$. Nous disposons donc d'une application
\begin{equation}\label{repcanhecke}\mathrm{H}(X) \rightarrow \Hom_{\rm
sf}(\Gamma \backslash (X\times X),\Z), \hspace{1 cm}T \mapsto ((x,y) \mapsto
T_{x,y}),\end{equation}
o\`u $\Hom_{\rm sf}(Y,\Z)$ d\'esigne le groupe ab\'elien des fonctions \`a
support fini de l'ensemble $Y$ dans $\Z$.  
\ps\ps

Si $x \in X$ on note $\Gamma_x \subset
\Gamma$ le stabilisateur de $x$.  On suppose que les propri\'et\'es
\'equivalentes suivantes sont satisfaites : \ps\ps\begin{itemize}

\item[(i)] Pour toute orbite $\Omega$ de $\Gamma$ dans $X \times X$, et pour tout $x \in X$,
l'intersection $\Omega \cap (X \times \{x\})$ est finie. \ps\ps

\item[(ii)] Pour tout $x \in X$, les orbites de $\Gamma_x$ sur $X$ sont
finies.  Autrement dit, pour tout $x,y \in X$, l'intersection $\Gamma_x \cap \Gamma_y$ est
d'indice fini dans $\Gamma_x$.  \ps\ps

\end{itemize}

\noindent Ces conditions assurent que l'application~\eqref{repcanhecke} est bijective.  En particulier, $\mathrm{H}(X)$ est un
$\Z$-module libre ayant pour base naturelle les fonctions caract\'eristiques
des orbites de $\Gamma$ sur $X \times X$.  \ps\ps

Fixons $x \in X$. La transitivit\'e de $X$ assure que l'application $\Gamma \rightarrow X \times X, \gamma \mapsto (\gamma(x),x),$ induit
des bijections \begin{equation}\label{bijmathecke}\Gamma_x\backslash \Gamma/\Gamma_x \lisomo \Gamma_x
\backslash ( X \times \{x\}) \lisomo \Gamma \backslash (X \times X).\end{equation} Cela 
identifie en particulier $\mathrm{H}(X)$ \`a $\Hom_{\rm
sf}(\Gamma_x\backslash \Gamma/\Gamma_x,\Z)$.  Par transport de structure, ce
dernier h\'erite de $\mathrm{H}(X)$ une structure d'anneau : on retombe sur
les pr\'esentations plus standards des anneaux de Hecke, comme par exemple dans \cite{satake}, \cite[\S
3]{shimura}, \cite{cartier}, \cite[Kap.  IV]{freitag} ou~\cite{grossatake}. On prendra
garde que suivant les r\'ef\'erences, la structure d'anneau consid\'er\'ee
sur $\Hom_{\rm sf}(\Gamma_x\backslash \Gamma/\Gamma_x,\Z)$ (d\'efinie en g\'en\'eral par un produit de convolution explicite) peut l\'eg\`erement diff\'erer
de la n\^otre : c'est notamment le cas des articles de Cartier et Gross, auxquels nous r\'ef\`ererons au~\S\ref{parametrisationsatake}, pour lesquels l'anneau ${\rm H}(X)$ est pr\'ecis\'ement l'anneau oppos\'e du n\^otre.

\ps\ps

Observons que la seconde formulation de la condition (ii) \'etant
sym\'etrique en $x,y$, la condition (i) \'equivaut encore
\`a demander que pour toute orbite $\Omega$ de $\Gamma$ dans $X \times X$,
et pour tout $x \in X$, l'intersection $\Omega \cap (\{x\} \times X)$ est
finie.
Ainsi, si $T \in \mathrm{H}(X)$, il existe un unique $T^{\rm t} \in \mathrm{H}(X)$ tel que
$T^{\rm t}_{x,y}=T_{y,x}$ pour tout $x,y \in X$. L'endomorphisme
$$T \mapsto T^{\rm t}$$ de $\mathrm{H}(X)$ est une anti-involution, i.e. satisfait $(ST)^{\rm t} = T^{\rm t}
S^{\rm t}$ et \linebreak $(T^{\rm t})^{\rm t}=T$ pour tout $S,T \in \mathrm{H}(X)$, qui correspond
simplement \`a la transpos\'ee sur les matrices associ\'ees. Observons que cette
anti-involution est l'identit\'e si, et seulement si, les $\Gamma$-orbites de $X
\times X$ sont invariantes par $(x,y) \mapsto (y,x)$, auquel cas $\mathrm{H}(X)$
est commutative : c'est un cas particulier du crit\`ere de Gelfand. \ps\ps

\subsection{Le foncteur des $\Gamma$-modules vers les $\mathrm{H}(X)^{\rm
opp}$-modules}\label{foncteurmx}

Soit $X$ un $\Gamma$-ensemble transitif satisfaisant les conditions (i) et (ii) du \S \ref{corrheckegeneral}. L'anneau $\mathrm{H}(X)$ intervient de la
mani\`ere suivante dans la th\'eorie des repr\'esentations de $\Gamma$. Si
$M$ est un $\Z[\Gamma]$-module, alors le groupe ab\'elien
$$M_X={\rm Hom}_{\Z[\Gamma]}(\Z[X],M)$$
h\'erite d'une action \underline{\`a droite} de $\mathrm{H}(X)$ par composition \`a la source. Il est \'evident que $M \mapsto M_X$ est un foncteur des $\Gamma$-modules (\`a gauche) vers les $\mathrm{H}(X)$-modules \`a droite. \ps\ps

Si $M$ est un $\Z[\Gamma]$-module, et si $x \in X$, l'application $\varphi \mapsto \varphi(x)$ identifie $M_X$ au sous-groupe des invariants 
$M^{\Gamma_x} \subset M$, ce qui munit \'egalement ce sous-groupe d'une structure de
$\mathrm{H}(X)$-module. Supposons que $T \in {\rm H}(X)$ a pour matrice la fonction caract\'eristique de la double classe $\Gamma_x \gamma \Gamma_x$
via l'identification $\Gamma_x \backslash \Gamma/ \Gamma_x \isomo\Gamma \backslash (X \times X)$ choisie au~\S\ref{corrheckegeneral}. On a la formule classique 
\begin{equation}\label{actionconcrethecke} T(m)=\sum_i \gamma_i(m)\, , \hspace{0.5 cm}  \forall m \in M^{{\Gamma}_{x}}, \end{equation} 
pour toute d\'ecomposition $\Gamma_x \gamma \Gamma_x=\coprod_i \gamma_i \Gamma_x$ (r\'eunion finie).\ps\ps

Dans ce contexte, l'anti-involution $T \mapsto T^{\rm t}$ d\'efinie
au~\S\ref{corrheckegeneral} prend le sens suivant. Soient $M$, $M'$ deux $\Z[\Gamma]$-modules, $N$ 
un groupe ab\'elien, et $(-|-): M \times M' \rightarrow N$ une application bilin\'eaire telle que $(\gamma
m | \gamma m')=(m | m')$ pour tout $\gamma \in \Gamma$ et tout $(m,m') \in M \times M'$.
Observons que si $(\varphi,\varphi') \in M_X \times M'_X$, alors $$(\varphi,\varphi'):=(\varphi(x)
| \varphi'(x))$$ ne d\'epend pas du choix de $x \in X$, et d\'efinit une application bilin\'eaire $M_X \times M'_X \rightarrow N$. Si l'on identifie $M_X$ \`a $M^{\Gamma_x}$ comme plus haut, cette application bilin\'eaire n'est autre que la restriction de $(-|-)$ \`a $M^{\Gamma_x} \times {M'}^{\Gamma_x}$. \ps\ps

On dira que $X$ est {\it sym\'etrique} s'il v\'erifie, en plus des conditions (i) et (ii) du \S \ref{corrheckegeneral},  les propri\'et\'es \'equivalentes suivantes\footnote{Cette propri\'et\'e n'est pas automatique si $X$ est infini. On pourra par exemple consid\'erer le groupe $\Gamma = \Q \rtimes \Q^\times$ des transformations affines de $\Q$, et le $\Gamma$-ensemble $X$ constitu\'e des parties de $\Q$ de la forme $a\Z+b$ avec $a \in \Q^\times$ et $b \in \Q$.} :
\ps\ps

\begin{itemize}

\item[(iii)] Pour toute orbite $\Omega$ de $\Gamma$ dans $X \times X$, et pour tout $x \in X$, alors $|\Omega \cap (X \times \{x\})|=|\Omega \cap (\{x\} \times X)|$.\ps\ps

\item[(iv)]  Pour tout $x,y \in X$, $\Gamma_x \cap \Gamma_y$ a m\^eme indice dans $\Gamma_x$ et $\Gamma_y$.\ps\ps

\end{itemize}

\begin{lemme}\label{lemmeadjonction} Supposons $X$ sym\'etrique. 
Si $T \in \mathrm{H}(X)$, et si $(\varphi,\varphi') \in M_X \times M'_X$, alors 
$(T(\varphi),\varphi')=(\varphi,T^{\rm t}(\varphi'))$.
\end{lemme}

\begin{pf} Soit $\psi : X \times X \rightarrow N$ une application constante
sur chaque $\Gamma$-orbite dans $X \times X$, et nulle hors d'un nombre fini
d'entre elles. La propri\'et\'e de sym\'etrie
de $X$ entra\^ine, pour tout $x \in X$, la relation $\sum_{y \in X} \psi(y,x)=\sum_{y
\in X} \psi(x,y)$. On applique ceci \`a la fonction 
$(x,y) \mapsto  T_{x,y} \cdot (\varphi(x)|\varphi'(y))$.
\end{pf}

\begin{remarque} {\rm \label{transposeHXmod} Supposons que $V$ est un ${\rm H}(X)$-module \`a droite. L'application ${\rm H}(X) \times V \rightarrow V, \, \, (T,v) \mapsto T^{{\rm t}} \, v,$ d\'efinit une structure de ${\rm H}(X)$-module (\`a gauche) sur $V$, que l'on notera $V^{{\rm t}}$. }
\end{remarque}

\subsection{L'anneau de Hecke d'un $\Z$-groupe}\label{annheckeg}

Soit maintenant $G$ un $\Z$-groupe. Nous allons appliquer les d\'efinitions
pr\'ec\'edentes \`a $\Gamma=G(\AAA_f)$ et $X=\mathcal{R}(G)$.  {\it L'anneau de Hecke de $G$} est l'anneau $$\mathrm{H}(G):=\mathrm{H}(\mathcal{R}(G)).$$
 \ps\ps

Rappelons que pour tout premier $p$, le groupe $G(\Q_p)$ h\'erite de
$\Q_p$ une structure de groupe topologique localement compact (s\'epar\'e, r\'eunion d\'enombrable de compacts), le sous-groupe $G(\Z_p)$ \'etant \`a la fois compact et ouvert.
Le groupe $G(\AAA_f)$ est \'egalement un groupe topologique localement compact
pour la topologie ayant pour base de voisinages de l'unit\'e les ouverts de
la forme $\prod_{p \in {\rm P}} U_p$ o\`u chaque  $U_p$ pour $p$ premier est un voisinage ouvert de l'unit\'e dans $G(\Q_p)$, tel que $U_p=G(\Z_p)$
pour presque tout $p$.  En particulier, $G(\widehat{\Z})$ est un sous-groupe
ouvert compact de $G(\AAA_f)$.  Cela entra\^ine que $\mathcal{R}(G)$ satisfait la
propri\'et\'e (ii) du~\S\ref{corrheckegeneral}, ainsi que les
$G(\Q_p)$-ensembles $$\mathcal{R}_p(G):=G(\Q_p)/G(\Z_p).$$ 
Le $G(\AAA_f)$-ensemble $\mathcal{R}(G)$, ainsi que les
$\mathcal{R}_p(G)$, sont sym\'etriques au sens du~\S\ref{foncteurmx} si $G(\AAA_f)$ est unimodulaire, ce qui est
notamment le cas si la composante neutre de $G(\C)$ est r\'eductive \cite[\S 5.5]{borelfini}. 
 \ps\ps

Si $p$ est premier, on d\'efinit \'egalement $\mathrm{H}_p(G)$ comme l'anneau de Hecke du
$G(\Q_p)$-ensemble $\mathcal{R}_p(G)$. Observons que le $G(\AAA_f)$-ensemble $\mathcal{R}(G)$ s'identifie canoniquement au sous-ensemble de $\prod_{p \in {\rm P}} \mathcal{R}_p(G)$ constitu\'e des $(x_p)$ tels que $x_p = G(\Z_p)$ pour presque tout $p$. 
Nous avons d\'ej\`a vu une manifestation de ce fait dans le plongement d'Eichler~\eqref{plongementdeichler}.
On dispose en particulier, pour chaque premier $p$, d'un homomorphisme
d'anneaux injectif canonique $$\mathrm{H}_p(G) \rightarrow \mathrm{H}(G)$$ associant
\`a un $T \in \mathrm{H}_p(G)$ l'endomorphisme de $\Z[\mathcal{R}(G)]$ envoyant
$y=(y_\ell)_{\ell \in {\rm P}}$ sur $\sum_x T_{x_p,y_p} \,  x$, o\`u
$x$ parcourt les \'el\'ements de $\mathcal{R}(G)$ tels que $x_\ell=y_\ell$ dans
$\mathcal{R}_\ell(G)$ pour tout $\ell \neq p$.  On \'ecrira simplement $$\mathrm{H}_p(G) \subset
\mathrm{H}(G).$$  Si $p \neq q$, $S \in \mathrm{H}_p(G)$ et $T \in
\mathrm{H}_q(G)$ alors $TS=ST$. \ps\ps

Si l'on se donne pour tout premier $p$ une $G(\Q_p)$-orbite $\Omega_p  \subset \mathcal{R}_p(G) \times \mathcal{R}_p(G)$, et si de plus $\Omega_p$ est l'orbite de $G(\Z_p) \times G(\Z_p)$ pour presque tout $p$, alors le sous-ensemble des \'el\'ements $(\omega_p)$ de $\prod_p \Omega_p$ tels que $\omega_p = G(\Z_p) \times G(\Z_p)$ pour presque tout $p$ s'identifie naturellement \`a une $G(\AAA_f)$-orbite dans $\mathcal{R}(G)  \times \mathcal{R}(G)$. R\'eciproquement, toute  $G(\AAA_f)$-orbite $\Omega \subset \mathcal{R}(G) \times \mathcal{R}(G)$ est de cette forme, et ce pour une unique famille $(\Omega_p)$, la $G(\Q_p)$-orbite $\Omega_p$ \'etant l'image de $\Omega$ par la projection canonique $\mathcal{R}(G) \times \mathcal{R}(G) \rightarrow \mathcal{R}_p(G) \times \mathcal{R}_p(G)$. Il resulte de ces observations et de la bijectivit\'e de l'application \eqref{repcanhecke} que $\mathrm{H}(G)$ est isomorphe au produit tensoriel de ses sous-anneaux  $\mathrm{H}_p(G)$ : $$\bigotimes_{p
\in {\rm P}} \mathrm{H}_p(G)\isomo \mathrm{H}(G).$$ La compr\'ehension de $\mathrm{H}(G)$ se ram\`ene donc
enti\`erement \`a celle des $\mathrm{H}_p(G)$. \ps\ps

L'anneau ${\rm H}_p(G)$ ne d\'epend que du $\Z_p$-groupe $G_{\Z_p} = G \times_\Z \Z_p$.  Lorsque $G_{\Z_p}$ est r\'eductif, les r\'esultats g\'en\'eraux de Satake et
Bruhat-Tits entra\^inent que $\mathrm{H}_p(G)$ est commutative : nous y reviendrons
au~\S\ref{parametrisationsatake}. Il en va donc de m\^eme de $\mathrm{H}(G)$ si $G$ est r\'eductif sur $\Z$.  
Cette propri\'et\'e est toutefois \'el\'ementaire dans les cas les plus classiques, que nous allons
rappeler ci-dessous. \ps\ps

\subsection{Quelques anneaux de Hecke classiques}\label{annheckeclass}\label{defperestroika}

Supposons d'abord  $G={\rm PGL}_n$. Nous avons vu que $\mathcal{R}(G)$ s'identifie
\`a $$\underline{\mathcal{R}}_\Z(V):=\Q^\times \backslash \mathcal{R}_\Z(V)$$ o\`u $V=\Q^n$. Rappelons que
$\underline{M} \in \underline{\mathcal{R}}_\Z(V)$ d\'esigne la classe
d'homoth\'etie d'un r\'eseau $M \in \mathcal{R}_\Z(V)$. \ps\ps

Si $M, N \in \mathcal{R}_\Z(V)$, il existe un plus petit entier $d\geq 1$ tel que
$dN \subset M$. La classe d'isomorphisme du groupe ab\'elien $M/dN$ ne
d\'epend que de la $G(\AAA_f)$-orbite de $(\underline{N},\underline{M})$ dans $\underline{\mathcal{R}}_\Z(V) \times \underline{\mathcal{R}}_\Z(V)$.
La th\'eorie des diviseurs \'el\'ementaires montre alors que l'application
ainsi d\'efinie $$G(\AAA_f) \backslash (\underline{\mathcal{R}}_\Z(V) \times \underline{\mathcal{R}}_\Z(V)) \rightarrow {\rm AF},$$
${\rm AF}$ d\'esignant l'ensemble des classes d'isomorphisme de groupes ab\'eliens
finis, est une injection dont l'image est constitu\'ee des groupes engendr\'es par $n-1$
\'el\'ements. Si $A$ est un tel groupe, l'op\'erateur de Hecke
associ\'e ${\rm T}_A \in \mathrm{H}(G)$ v\'erifie par d\'efinition $${\rm T}_A(\underline{M})=\sum_N \underline{N},$$
o\`u $N$ parcourt les sous-groupes de $M$ tels que $M/N \simeq A$. Ces op\'erateurs ${\rm T}_A$,  avec $A$ parcourant les groupes ab\'eliens finis engendr\'es par $n-1$ \'el\'ements, forment donc une $\Z$-base de $\mathrm{H}(G)$. Il est clair que ${\rm T}_{A \times B}={\rm T}_A {\rm T}_B$ si $|A|$ et $|B|$
sont premiers entre eux, et que ${\rm T}_A \in \mathrm{H}_p(G)$ si, et seulement si, $A$
est un $p$-groupe. \ps\ps

Si $n=2$, on v\'erifie ais\'ement que ${\rm T}_A^{\rm t}={\rm T}_A$ pour tout $A$, en particulier
$\mathrm{H}(G)$ est commutative (la notation ${\rm T}^{\rm t}$ est d\'efinie au \S \ref{corrheckegeneral}).  Le premier point ne vaut plus pour $n>2$, mais $\mathrm{H}(G)$ reste commutative.  On peut le voir
simplement en munissant $V$ d'une forme bilin\'eaire sym\'etrique non
d\'eg\'en\'er\'ee.  L'application $\underline{M} \mapsto
\underline{M^\sharp}$ est une involution de
$\underline{\mathcal{R}}_\Z(V)$. Elle induit une involution lin\'eaire de $\Z[ \,\underline{\mathcal{R}}_\Z(V)\,]$, puis par conjugaison une involution 
$\iota$ de $\mathrm{H}(G)$, qui n'est autre que
$(T_{\underline{N},\underline{M}}) \mapsto
(T_{\underline{N^\sharp},\underline{M^\sharp}})$ sur les matrices
associ\'ees.  Mais si $N \subset M$ alors $N^\sharp/M^\sharp$ est en
dualit\'e parfaite avec $M/N$, et donc $\iota$ co\"incide avec
l'anti-involution canonique de $\mathrm{H}(G)$ : $\iota(T)=T^{\rm t}$ pour
tout $T \in {\rm H}(G)$ (voir aussi~\cite[\S
3]{shimura}).\ps\ps

Discutons maintenant le cas particuli\`erement important pour
ce m\'emoire des $\Z$-groupes orthogonaux et symplectiques \cite{satake} \cite{freitag} \cite{andrianov}. On reprend les notations du~\S\ref{ensclassbil}, en particulier $V=L\otimes \Q$ est de dimension $n$ paire, $\varphi$ est une forme 
bilin\'eaire sur $V$ qui est sym\'etrique (resp. altern\'ee), pour laquelle $L$ est autodual pair, et $G \subset \GL_L$ est le groupe ${\rm O}_L$ (resp. ${\rm Sp}_L$). \ps\ps

Dans ce cas nous avons vu que $\mathcal{R}(G)$ s'identifie au
$G(\AAA_f)$-ensemble $\mathcal{R}_{\Z}^{\rm a}(V)$ des r\'eseaux autoduaux de $V$ (Proposition~\ref{rortho}). Si
$(N,M) \in \mathcal{R}_{\Z}^{\rm a}(V) \times \mathcal{R}_{\Z}^{\rm a}(V)$, la classe d'isomorphisme du
groupe ab\'elien $M/(N \cap M)$ ne d\'epend que de la $G(\AAA_f)$-orbite de
$(N,M)$. On a donc d\'efini une application naturelle
\begin{equation}\label{appheckeorth}G(\AAA_f) \backslash (\mathcal{R}_{\Z}^{\rm a}(V) \times
\mathcal{R}_{\Z}^{\rm a}(V)) \rightarrow {\rm AF}, \, \, \,
(N,M) \mapsto M/(N \cap M).\end{equation}

\begin{prop}\label{baseheckeo} L'application~\eqref{appheckeorth} est une injection dont l'image
est constitu\'ee des groupes engendr\'es par $n/2$ \'el\'ements.
\end{prop}

Cette proposition est bien connue, nous en redonnerons une d\'emonstration \`a la fin de ce
paragraphe pour le confort du lecteur. Soit $A$ un groupe ab\'elien fini engendr\'e par au plus $n/2$ \'el\'ements.
Il lui est donc associ\'e un op\'erateur de Hecke $${\rm T}_A \in \mathrm{H}(G)$$  
d\'efini par ${\rm T}_A(M) = \sum_N N$, la somme portant sur les $N$ tels
que $M/(N\cap M) \simeq A$, soit encore sur les $A$-voisins de $M$ au sens
de la d\'efinition~III.1.2 dans le cas quadratique. Ces op\'erateurs ${\rm T}_A$ forment donc une $\Z$-base de $\mathrm{H}(G)$. On a bien s\^ur encore 
${\rm T}_{A \times B}={\rm T}_A {\rm T}_B$ si $|A|$ et $|B|$
sont premiers entre eux, et ${\rm T}_A \in \mathrm{H}_p(G)$ si, et seulement si, $A$ est un
$p$-groupe. Du point de vue du chapitre~III, un op\'erateur tout particuli\`erement important pour
nous est l'op\'erateur ${\rm T}_{\Z/d\Z}$ pour $d\geq 1$, que l'on note
\'egalement simplement ${\rm T}_d$. 

\begin{prop}\label{heckeocomm} Soit $A$ un groupe ab\'elien fini engendr\'e par $n/2$ \'el\'ements.  Alors
${\rm T}_A^{\rm t}={\rm T}_{A^\vee}={\rm T}_A$. En particulier, l'anneau ${\rm H}(G)$ est commutatif. 
\end{prop}

\begin{pf} En effet, la premi\`ere assertion est cons\'equence du
scholie-d\'efinition III.1.2 quand $\varphi$ est sym\'etrique, et de mani\`ere similaire dans le cas altern\'e. La seconde se d\'eduit de la premi\`ere d'apr\`es la fin du paragraphe~\ref{corrheckegeneral}. 
Voir \'egalement \cite[Ch. III]{satake}, \cite[Kap. IV]{freitag} et le~\S V.\ref{exemplehecke}.\end{pf}

Discutons enfin du groupe des similitudes projectives ${\rm P}\widetilde{G}$. Soient $p$ un nombre premier et $\mathcal{R}^{\rm h}_{\Z_p}(V_p)$ l'ensemble des r\'eseaux homoduaux pairs de
$V_p$, introduit apr\`es l'\'enonc\'e du lemme~\ref{lemmeortho}. Si $\varphi$ est sym\'etrique (resp. altern\'ee), on rappelle qu'un r\'eseau $M \in \mathcal{R}_{\Z_p}(V_p)$ est homodual si et seulement s'il existe $\lambda_M \in p^\Z$, n\'ecessairement unique, tel que $x \mapsto \lambda_M\,\frac{\varphi(x,x)}{2}$ (resp. $\lambda_M \varphi$) munisse $M$ d'une structure de ${\rm q}$-module (resp. ${\rm a}$-module) sur $\Z_p$.  Comme le ${\rm q}$-module $V_p$ est hyperbolique d'apr\`es le scholie~II.2.5, il en va de m\^eme de $M \in \mathcal{R}_{\Z_p}^{\rm h}(V_p)$ en tant que ${\rm q}$-module sur $\Z_p$ d'apr\`es la proposition II.1.2. Cela montre que l'application $g \mapsto g (L)$ induit des isomorphismes  $\mathcal{R}_p(\widetilde{G})\isomo \mathcal{R}^{\rm h}_{\Z_p}(V_p)$ et $\mathcal{R}_p(G) \isomo \mathcal{R}^{\rm a}_{\Z_p}(V_p)$. En particulier, l'ensemble $\underline{\mathcal{R}}^{\rm h}_{\Z_p}(V_p) : = \Q_p^\times \backslash \mathcal{R}^{\rm h}_{\Z_p}(V_p)$ s'identifie naturellement \`a $\mathcal{R}_p({\rm P}\widetilde{G})$. \ps\ps

Soit $M \in \mathcal{R}^{\rm h}_{\Z_p}(V_p)$. On note  ${\rm v}_M \in \Z$ l'unique \'el\'ement tel que $\lambda_M = p^{-{\rm v}_M}$. Si $g \in \widetilde{G}(\Q_p)$, alors ${\rm v}_{g(M)}={\rm v}_M+v$ o\`u $v$ est  la valuation $p$-adique de $\nu(g)$.
Soit $(\underline{N},\underline{M})$ un couple d'\'el\'ements de
$\underline{\mathcal{R}}^{\rm
h}_{\Z_p}(V_p)$. Quitte \`a changer de repr\'esentant
$N$, on peut supposer  ${\rm v}_M-{\rm v}_N \in \{0,1\}$. Le couple $(M/N\cap M,{\rm v}_M-{\rm v}_N)$ ne d\'epend alors que de la
${\rm P}\widetilde{G}(\Q_p)$-orbite de $(\underline{N},\underline{M})$, ce qui d\'efinit une application 
\begin{equation}\label{appheckego} {\rm P}\widetilde{G}(\Q_p) \backslash
(\underline{\mathcal{R}}_{\Z_p}^{\rm h}(V_p) \times
\underline{\mathcal{R}}_{\Z_p}^{\rm h}(V_p)) \rightarrow {\rm AF} \times
\{0,1\}.\end{equation}

\begin{prop}\label{baseheckego} L'application~\eqref{appheckego} est une
injection d'image l'ensemble des couples $(A,-)$ o\`u $A$ est un $p$-groupe ab\'elien engendr\'e par $n/2$ \'el\'ements.
\end{prop}

Nous reportons au~\S V.\ref{exemplehecke} la d\'emonstration de cette
proposition.  Soit $(A,i) \in {\rm AF} \times \{0,1\}$ o\`u $A$ est un
$p$-groupe engendr\'e par au plus $n/2$ \'el\'ements.  Nous dirons que
$\underline{N} \in \underline{\mathcal{R}}_{\Z_p}^{\rm h}(V_p)$ est un $A$-voisin de
type $i$ de $\underline{M} \in \underline{\mathcal{R}}_{\Z_p}^{\rm h}(V_p)$, si l'image
de $(\underline{N},\underline{M})$ par l'application \eqref{appheckego} est
$(A,i)$.  L'op\'erateur de Hecke correspondant est not\'e $${\rm T}_{(A,i)}
\in {\rm H}_p({\rm P}\widetilde{G});$$ ils forment une $\Z$-base de $ {\rm
H}_p({\rm P}\widetilde{G})$. Si $M^\sharp=M$, on constate que
$\underline{N}$ est un $A$-voisin de type $0$ de $\underline{M}$ si, et
seulement si, $\underline{N}$ poss\`ede un repr\'esentant autodual, alors
unique, et si ce dernier est un $A$-voisin de $M$ au sens pr\'ec\'edent.  La
notion de $A$-voisin de type $1$ de $\underline{M}$ est en revanche
``nouvelle''.  L'exemple qui suit sera particuli\`erement important dans ce
m\'emoire.  \ps\ps

Soient $M,N \in \mathcal{R}_{\Z}^{\rm h}(V)$. Suivant Koch et Venkov dans le cas
quadratique \cite{kochvenkov}, nous dirons que $N$ est une {\it perestro\"ika} de $M$
relativement \`a $p$ si $$pM \subsetneq N \subsetneq M.$$ 
Il est ais\'e de v\'erifier que $N$ est une perestro\"ika de $M$
relativement \`a $p$ si, et seulement si, ${\rm v}_M-{\rm v}_{p^{-1}N}=1$ et $\underline{N}$ est un $0$-voisin de $\underline{M}$ de type $1$. De plus, la proposition suivante est imm\'ediate. 

\begin{prop} Soient $M \in \mathcal{R}_{\Z}^{\rm h}(V)$ et $p$ un nombre premier. L'application $N \mapsto N/pM$
d\'efinit une bijection de l'ensemble des perestro\"ikas de $M$ relativement
\`a $p$ sur l'ensemble des lagrangiens de $M \otimes \F_p$. 
\end{prop}

L'op\'erateur des perestro\"ikas relativement \`a $p$ est
l'op\'erateur $${\rm K}_p:={\rm T}_{(0,1)} \in \mathrm{H}_p(\mathrm{P}\widetilde{G}).$$ 
Si $(N,M) \in \mathcal{R}_\Z^{\rm h}(V)$, $N$ est une
perestro\"ika de $M$ relativement \`a $p$ si, et seulement si, $pM$ est une
perestro\"ika de $N$ relativement \`a $p$.  En particulier, ${\rm K}_p^{\rm t}={\rm K}_p$. 
En fait, on a $T^{\rm t}=T$ pour tout $T \in {\rm H}({\rm P}\widetilde{G})$ comme
on le verra au \S V.\ref{exemplehecke}. \ps\ps

Terminons ce paragraphe, comme pr\'evu, par une d\'emonstration de la
proposition~\ref{baseheckeo}.\ps\ps


\begin{pf} On se place dans le cas quadratique, i.e. $\varphi$ sym\'etrique et ${\rm q}(x)=\frac{\varphi(x,x)}{2}$, auquel cas $L$ est un ${\rm q}$-module sur $\Z$. La d\'emonstration dans le cas altern\'e est similaire (et m\^eme plus simple).\ps

Il s'agit de montrer que si $U$ est un ${\rm q}$-module hyperbolique sur $\Q_p$, et si $(L_1,L_2)$ et $(L'_1,L'_2)$ sont deux couples de r\'eseaux autoduaux de 
$U$ tels que $L_1/(L_1 \cap L_2) \simeq L'_1/(L'_1\cap L'_2)$ alors il existe $\alpha \in {\rm O}(U)$ tel que $\alpha(L_i)=L'_i$ pour $i=1,2$. On raisonne par r\'ecurrence sur $\dim(U)$. \ps

Les cas $U=0$ et $L_1=L_2$ sont triviaux. On suppose donc que $L_1 \neq L_2$. L'id\'eal annulateur du quotient $L_{1}/(L_{1}\cap L_{2})$ est donc de la forme $p^{\nu}\mathbb{Z}_{p}$ avec $\nu\geq 1$. Nous affirmons qu'il existe un \'el\'ement $e_{1}$ de $L_{1}$ et un \'el\'ement $e_{2}$ de $L_{2}$ tels que
$$
\hspace{24pt}
\mathrm{q}(e_{1})=0
\hspace{24pt},\hspace{24pt}
\mathrm{q}(e_{2})=0
\hspace{24pt},\hspace{24pt}
e_{1}.e_{2}=p^{-\nu}
\hspace{24pt}.
$$
En effet, il est tout d'abord facile de se convaincre de ce qu'il existe un \'el\'ement $\epsilon_{1}$ de $L_{1}$ et un \'el\'ement $\epsilon_{2}$ de $L_{2}$ avec
$\epsilon_{1}.\epsilon_{2}=p^{-\nu}$. Le lemme de Hensel montre ensuite qu'il existe une matrice
$$
\hspace{24pt}
P=
\begin{bmatrix}
a_{1,1} & a_{1,2} \\ a_{2,1} & a_{2,2}
\end{bmatrix}
\hspace{24pt},
$$
\`a coefficients dans $\mathbb{Z}_{p}$, avec $P\equiv\mathrm{I}\bmod{p}^{\nu}$ telle que l'on a~:
$$
\hspace{24pt}
{}^{\mathrm{t}}\hspace{-2pt}P
\begin{bmatrix}
2\mathrm{q}(\epsilon_{1}) & p^{-\nu} \\
p^{-\nu} & 2\mathrm{q}(\epsilon_{2})
\end{bmatrix}
P
\hspace{4pt}=\hspace{4pt}
\begin{bmatrix}
0 & p^{-\nu} \\ p^{-\nu} & 0
\end{bmatrix}
\hspace{24pt}.
$$
On prend $e_{1}=a_{1,1}\epsilon_{1}+a_{2,1}\epsilon_{2}$ et $e_{2}=a_{1,2}\epsilon_{1}+a_{2,2}\epsilon_{2}\in L_{2}$ (la congruence $P\equiv\mathrm{I}\bmod{p}^{\nu}$ implique $e_{1}\in L_{1}$ et $e_{2}\in L_{2}$). Ce qui termine la d\'emonstration de l'affirmation. \ps
\bigskip
Achevons maintenant la r\'ecurrence. On note respectivement $H$, $H_{1}$ $H_{2}$, le sous-espace vectoriel de $U$ engendr\'e par $e_{1}$ et $e_{2}$, le sous-module de $L_{1}$ engendr\'e par $e_{1}$ et $p^{\nu}e_{2}$, le sous-module de $L_{2}$ engendr\'e par $p^{\nu}e_{1}$ et $e_{2}$. On munit $H$, $H_{1}$ et $H_{2}$ des formes quadratiques induites par celle de $U$. Par construction, $H\approx\mathrm{H}(\mathbb{Q}_{p})$ et $H_i \approx \mathrm{H}(\Z_p)$ pour $i=1,2$. On note respectivement $W$, $M_{1}$ et $M_{2}$, l'orthogonal de $H$ dans $U$, l'orthogonal de $H_{1}$ dans $L_{1}$ et l'orthogonal de $H_{2}$ dans $L_{2}$. On a des d\'ecompositions en somme orthogonale
$$
\hspace{24pt}
U=H\oplus W
\hspace{24pt},\hspace{24pt}
L_{1}=H_{1}\oplus M_{1}
\hspace{24pt},\hspace{24pt}
L_{2}=H_{2}\oplus M_{2}
$$
et des isomorphismes
$$
\hspace{8pt}
L_{1}/(L_{1}\cap L_{2})
\cong
H_{1}/(H_{1}\cap H_{2})
\hspace{2pt}\oplus\hspace{2pt}
M_{1}/(M_{1}\cap M_{2})
\hspace{8pt},\hspace{8pt}
H_{1}/(H_{1}\cap H_{2})\cong\mathbb{Z}_{p}/p^{\nu}\mathbb{Z}_{p}
\hspace{8pt}.
$$
On remplace le couple $(L_{1},L_{2})$ par le couple $(L'_{1},L'_{2})$ et on introduit pareillement les ${\rm q}$-espaces vectoriels $H'$, $W'$ et les ${\rm q}$-modules $H'_{1}$, $H'_{2}$, $M'_{1}$, $M'_{2}$. On obtient l'automorphisme cherch\'e $\alpha:U\to U$ comme somme orthogonale d'isomorphismes ad\'equats de $\mathrm{q}$-espaces vectoriels $H\to H'$ et $W\to W'$, l'existence du second \'etant assur\'ee par l'hypoth\`ese de r\'ecurrence.
\end{pf}

\subsection{$\mathrm{H}({\rm SO}_L)$ versus $\mathrm{H}({\rm O}_L)$}\label{hsovso}

Soit $L$ un ${\rm q}$-module sur $\Z$. Discutons bri\`evement du lien entre $\mathrm{H}({\rm SO}_L)$ et $\mathrm{H}({\rm O}_L)$.
Les cas de ${\rm PGSO}_L$ et ${\rm PGO}_L$ pourraient se traiter de
mani\`ere similaire. \ps\ps

D'apr\`es la proposition~\ref{rso=ro}, l'inclusion ${\rm SO}_L \rightarrow {\rm O}_L$ induit une
bijection ${\rm SO}_L(\AAA_f)$-\'equivariante 
$\mathcal{R}({\rm SO}_L) \isomo \mathcal{R}({\rm O}_L)$. Il en r\'esulte que $\mathrm{H}({\rm O}_L)$ s'identifie canoniquement \`a un 
sous-anneau de $\mathrm{H}(\rm{SO}_L)$ : ce sont les
sous-anneaux de ${\rm End}_\Z(\Z[\mathcal{R}_{\Z}^{\rm a}(V)])$ constitu\'es
respectivement des endomorphismes ${\rm O}_L(\AAA_f)$-\'equivariants et ${\rm
SO}_L(\AAA_f)$-\'equivariants. Le groupe quotient $${\rm
O}_L(\AAA_f)/{\rm SO}_L(\AAA_f) \simeq (\Z/2\Z)^{\rm P}$$ agit naturellement
par conjugaison sur $\mathrm{H}({\rm SO}_L)$, avec pour anneau des invariants $\mathrm{H}({\rm
O}_L)$. Cette action respecte la d\'ecomposition de $\mathrm{H}(G)$ en produit tensoriel des $\mathrm{H}_p(G)$
sur les $p \in {\rm P}$, et identifie \'egalement $\mathrm{H}_p({\rm O}_L)$ avec $\mathrm{H}_p({\rm SO}_L)^{\Z/2\Z}$. \ps \ps

Donnons un exemple d'un \'el\'ement de $\mathrm{H}_p({\rm SO}_L)$ qui n'est pas dans
$\mathrm{H}_p({\rm O}_L)$.  Soit $A=(\Z/p\Z)^{n/2}$ o\`u $n$ est le rang de $L$.  Consid\'erons l'ensemble $\Omega$ des
couples $(N,M)$ d'\'el\'ements de $\mathcal{R}_{\Z}^{\rm a}(V)$ tels que $N$ est un $A$-voisin de $M$. La proposition~\ref{baseheckeo} assure que $\Omega$ est une
${\rm O}_L(\Q_p)$-orbite. En revanche, il est r\'eunion disjointe de deux
orbites sous l'action de ${\rm SO}_L(\Q_p)$.  Pour le voir, on commence par
v\'erifier, par des consid\'erations similaires \`a celles du \S III.1, que l'application $$N \mapsto (M\cap N)/pM$$ induit une surjection (non bijective en g\'en\'eral) entre les $A$-voisins de $M$ et les
lagrangiens du ${\rm q}$-module hyperbolique $M \otimes \F_p$.  Mais il est
bien connu que pour tout corps $k$, et tout $k$-module hyperbolique $V$, il
y a exactement deux orbites de lagrangiens de $V$ sous l'action de ${\rm
SO}(V)$ (et une seule sous ${\rm O}(V)$, par le th\'eor\`eme de Witt). 
Par lissit\'e de ${\rm SO}_M$ sur $\Z_p$, chacune de ces deux orbites d\'efinit donc une ${\rm SO}(M)$-orbite de
$A$-voisins de $M$, et par cons\'equent deux op\'erateurs de Hecke distincts ${\rm T}^\pm_A \in
\mathrm{H}({\rm SO}_L)$ dont la somme vaut ${\rm T}_A$, qui sont \'echang\'es sous
l'action de ${\rm O}_L(\Q_p)/{\rm SO}_L(\Q_p)=\Z/2\Z$. \ps\ps

\subsection{Isog\'enies}\label{parisogenies}

Nous discutons maintenant des isog\'enies entre
$\Gamma$-ensembles transitifs, en exposant une 
variante des consid\'erations dans~\cite[II \S 7]{satake}. \ps\ps 

Soient $X$ un $\Gamma$-ensemble et $X'$ un $\Gamma'$-ensemble. Rappelons qu'un morphisme $X \rightarrow X'$ est un couple $(f,g)$ o\`u $g : X \rightarrow X'$ est une application et $f : \Gamma \rightarrow \Gamma'$ est un morphisme de groupes, tels que $g(\gamma x)=f(\gamma) g(x)$ pour tout $x \in X$ et tout $\gamma \in \Gamma$. Dans ce qui suit, il sera commode de convenir qu'un ensemble transitif est non vide. \ps\ps

\begin{lemme}\label{lemmestabS} Soient $X$ un $\Gamma$-ensemble transitif, $X'$ un $\Gamma'$-ensemble, et $(f,g)$ un morphisme $X \rightarrow X'$ tel que $f(\Gamma)$ est distingu\'e dans $\Gamma'$. Soit $S$ le stabilisateur de $g(X)$ dans $\Gamma'$, c'est-\`a-dire $S=\{\gamma \in \Gamma', \gamma g(X) \subset g(X)\}$. \begin{itemize}\ps\ps
\item[(i)] Pour tout $x \in g(X)$ on a $S=f(\Gamma) \Gamma'_x$. \ps \ps
\item[(ii)] $S=\{ \gamma \in \Gamma', \gamma g(X) \cap g(X) \neq \emptyset\}$.
\end{itemize}
\end{lemme}

\begin{pf} Soit $x \in g(X)$. Le sous-groupe $f(\Gamma)$ \'etant distingu\'e dans $\Gamma'$, le sous-ensemble $E_x:=f(\Gamma) \Gamma'_x \subset \Gamma'$ est un sous-groupe. La transitivit\'e de $X$ montre alors que : \ps\ps

\noindent --   $E_x$ ne d\'epend pas du choix de $x \in g(X)$, \ps\ps

\noindent --   $E_x$ est  l'ensemble des $\gamma \in \Gamma'$ tels que $\gamma(x) \in g(X)$. \ps\ps

\noindent On a donc $S=\bigcap_{x \in g(X)} E_x = \bigcup_{x \in g(X)} E_x =\{ \gamma \in \Gamma', \gamma g(X) \cap g(X) \neq \emptyset\}$. \end{pf}

Soient $X$ un $\Gamma$-ensemble transitif, $X'$ un $\Gamma'$-ensemble, et $(f,g)$ un morphisme $X \rightarrow X'$. On suppose comme dans le lemme ci-dessus que $f(\Gamma)$ est distingu\'e dans $\Gamma'$ et de plus que l'application $g$ est injective.\footnote{Nous renvoyons \`a
l'article de Satake pour une variante sans l'hypoth\`ese
d'injectivit\'e de $g$. Le lecteur ne perdrait pas grand chose ici \`a supposer que $\Gamma \subset \Gamma'$, $X \subset X'$, et que $f$ et $g$ sont les inclusions correspondantes.}  Soit $S$ le stabilisateur de $g(X)$ dans $\Gamma'$. L'application $(s,x) \mapsto g^{-1}(s(g(x))$, qui a un sens par injectivit\'e
de $g$, d\'efinit une action de $S$ sur $X$, dont la restriction \`a $f : \Gamma \rightarrow S$ 
est le $\Gamma$-ensemble $X$. Elle induit donc une action de $S/f(\Gamma)$ sur ${\rm H}(X)$ par automorphismes d'anneaux, 
dont nous noterons ${\rm H}(X)^{\rm inv} \subset {\rm H}(X)$ le
sous-anneau des invariants, qui est \'egalement ${\rm End}_{\Z[S]}(\Z[X])$. 

\begin{defprop}\label{defpropiso} Soit $u=(f,g) : X \rightarrow X'$ un morphisme entre le $\Gamma$-ensemble transitif $X$ et le $\Gamma'$-ensemble transitif $X'$.
On suppose que $f(\Gamma)$ est distingu\'e dans $\Gamma'$ et que $g$ est injectif.  \ps\ps \begin{itemize}
\item[(i)] Si $T \in {\rm H}(X)^{\rm inv}$ il existe un unique $T' \in \mathrm{H}(X')$ qui s'annule sur $(X' - 
g(X)) \times g(X)$ et tel que $T'_{g(x),g(y)}=T_{x,y}$ pour tout $x,y \in X$. \ps\ps
\item[(ii)] L'application ainsi d\'efinie ${\rm H}(u) : \mathrm{H}(X)^{\rm inv} \rightarrow \mathrm{H}(X'), T \mapsto T',$
est un homomorphisme injectif d'anneaux.\ps\ps
\end{itemize}
\end{defprop}

\begin{pf} L'assertion d'unicit\'e dans le (i) d\'ecoule de l'injectivit\'e de $g$ et de la transitivit\'e de $X'$. L'assertion (ii) se d\'eduit 
imm\'ediatement du (i). Il ne reste donc qu'\`a justifier l'existence de $T'$ dans le (i). Mais le point (ii) du lemme \ref{lemmestabS} montre que l'injection $g : X \rightarrow X'$ induit une bijection ${\rm Ind}_S^{\Gamma'} X \isomo X'$, et donc un isomorphisme
$\Z[\Gamma']\otimes_{\Z[S]}\Z[X] \isomo \Z[X']$. Ainsi, toute application lin\'eaire $S$-\'equivariante $T : \Z[X] \rightarrow \Z[X]$, compos\'ee avec $g: \Z[X] \rightarrow \Z[X']$, s'\'etend de mani\`ere unique en une application $\Gamma'$-\'equivariante $T' : \Z[X'] \rightarrow \Z[X']$; elle a les propri\'et\'es de l'\'enonc\'e.
\end{pf}

Dans tous les exemples que nous aurons \`a consid\'erer, il se trouve que le groupe $S$ pr\'eservera chaque $\Gamma$-orbite de $X \times X$, de sorte que l'on aura $\mathrm{H}(X)^{\rm inv}=\mathrm{H}(X)$. Un cas particuli\`erement simple est celui o\`u $\Gamma'=\Gamma$, $X'=X$, et $f$ et $g$ sont bijectifs. Dans ce cas $S=f(\Gamma)$ et ${\rm H}(u)$ est par 
d\'efinition l'automorphisme de $\mathrm{H}(X)$ qui vaut matriciellement $(T_{x,y}) \mapsto (T_{g^{-1}x,g^{-1}y})$.\ps\ps

Pla\c{c}ons-nous dans les hypoth\`eses de la proposition-d\'efinition~\ref{defpropiso}. Si $M$ est un $\Gamma'$-module, notons $M_{|\Gamma}$ le $\Gamma$-module obtenu par restriction
de $M$ par $f : \Gamma \rightarrow \Gamma'$. On dispose alors d'une application injective canonique 
$$M_{X'} \rightarrow (M_{|\Gamma})_{X}, \, \, \, \, \, \varphi \mapsto
\varphi_{|X}:=\varphi \circ g.$$
Le lemme suivant est imm\'ediat.

\begin{lemme}\label{compatibilitesatake} On reprend les hypoth\`eses de la proposition-d\'efinition~\ref{defpropiso}. Soient $M$ un $\Gamma'$-module,  $T \in \mathrm{H}(X)^{\rm inv}$ et $\varphi \in M_{X'}$, alors
$T(\varphi_{|X})={\rm H}(u)(T)(\varphi)$.
\end{lemme}
\ps\ps

\begin{example}\label{exampleiso} {\rm En guise d'exemple, repla\c{c}ons-nous dans le contexte des groupes de
similitudes (\S\ref{hsimilitudes}) et consid\'erons le $\Z$-morphisme
naturel $\mu : G \rightarrow \mathrm{P}\widetilde{G}$.  Le paragraphe
pr\'ec\'edent s'applique et d\'efinit un morphisme d'anneaux $${\rm H}(\mu) : \mathrm{H}(G) \rightarrow \mathrm{H}(\mathrm{P}\widetilde{G})$$ 
tel que ${\rm H}(\mu)({\rm T}_A)={\rm T}_{(A,0)}$ pour tout groupe ab\'elien fini $A$ engendr\'e par au plus $n/2$
\'el\'ements.}
\end{example}
\ps \ps

\noindent En effet, on consid\`ere $\Gamma=G(\AAA_f)$, $X=\mathcal{R}(G)$, $\Gamma'=\mathrm{P}\widetilde{G}(\AAA_f)$, $X'=\mathcal{R}(\mathrm{P}\widetilde{G})$, et on prend pour $f$ et $g$ les
applications naturelles d\'eduites de $\mu$.  Le groupe $\Gamma$ est distingu\'e dans
$\widetilde{G}(\AAA_f)$, ainsi donc que $f(\Gamma)$ dans $\Gamma'$.  De
plus, $g$ s'identifie \`a l'injection naturelle $\mathcal{R}_\Z^{\rm a}(V) \rightarrow
\underline{\mathcal{R}}_\Z^{\rm h}(V)$, $M \mapsto \underline{M}$, d'apr\`es la
proposition~\ref{diagrammesimilitude}. Le groupe $S$ est le sous-groupe des \'el\'ements $g \in \widetilde{G}(\AAA_f)$ tels que $\nu(g)$ est de la forme $a^2b$ avec $a \in \AAA_f^\times$ et $b \in \widehat{\Z}^\times$. 
Il agit trivialement sur $\Gamma
\backslash (\mathcal{R}_\Z^{\rm a}(V) \times \mathcal{R}_\Z^{\rm a}(V))$.  En effet, si $N,M \in \mathcal{R}_\Z^{\rm a}(V)$, $g \in
\widetilde{G}(\AAA_f)$, et si $p$ est un nombre premier, $g$
induit un isomorphisme $M_p/(N_p\cap M_p) \simeq g(M)_p/(g(N)_p \cap g(M)_p)$, ce qui conclut par la
proposition~\ref{baseheckeo}. L'assertion sur ${\rm T}_A$ se d\'eduit de la discussion cons\'ecutive \`a la
proposition~\ref{baseheckego}.

\section{Formes automorphes d'un $\Z$-groupe}\label{paysageaut}
L'anneau des ad\`eles de $\Q$ est l'anneau $\AAA= \R
\times \AAA_f$. Soit $G$ un $\Z$-groupe. Le groupe $G(\R)$ est un groupe de
Lie de mani\`ere naturelle, et le groupe $$G(\AAA)=G(\R) \times G(\AAA_f)$$
est encore localement compact et s\'epar\'e pour la topologie produit, la topologie sur
$G(\AAA_f)$ ayant d\'ej\`a \'et\'e rappel\'ee au~\S\ref{annheckeg}. 
Le groupe $G(\Q)$ se plonge naturellement diagonalement dans $G(\AAA)$ : c'est
un sous-groupe ferm\'e discret (voir~\cite[Ch. II \S 3]{GGPS} pour les rudiments sur ces constructions).

\subsection{Formes automorphes de carr\'e int\'egrable}\label{fautcarreint}\label{fautpar}

Rappelons quelques r\'esultats classiques d\^us \`a Borel et Harish-Chandra,
pour lesquels nous renvoyons \`a~\cite[\S 5]{borelfini}.  Supposons que la
composante neutre de $G(\C)$ est semi-simple \cite{humphreyslag} \cite{borelgp}. Le groupe localement compact $G(\AAA)$
est alors unimodulaire.  L'espace homog\`ene $$G(\Q)
\backslash G(\AAA)$$ h\'erite d'apr\`es Weil d'une mesure de Radon $\mu$
positive (non nulle) invariante par translations \`a droite sous
$G(\AAA)$~\cite[Chap. II]{weil} \cite[Chap. 2]{rudin}.  Il est de mesure finie.  \ps\ps

L'espace des {\it formes automorphes de carr\'e int\'egrable pour $G$} est le sous-espace
$$\mathcal{A}^2(G) \subset {\rm L}^2(G(\Q) \backslash G(\AAA),\mu)$$ des
\'el\'ements qui sont invariants par translations \`a droite sous
$G(\widehat{\Z})$ \cite[Ch. 3]{GGPS} \cite[\S 4]{boreljacquet}.  C'est un espace de Hilbert pour le produit scalaire hermitien $$\langle \,f,f' \,\rangle_{\rm Pe} =\int \overline{f} f' \, \, \, {\rm d}\mu,$$
appel\'e aussi {\it produit de Petersson}. Alternativement, $\mathcal{A}^2(G)$ peut \^etre vu comme l'espace des fonctions de carr\'e int\'egrable sur $G(\Q) \backslash G(\AAA)/G(\widehat{\Z})$, muni de la mesure de Radon image de $\mu$ par l'application (propre) canonique $G(\Q) \backslash G(\AAA) \rightarrow G(\Q) \backslash G(\AAA)/G(\widehat{\Z})$. L'espace $\mathcal{A}^2(G)$ est muni de deux structures additionnelles importantes que nous d\'ecrivons maintenant.  \ps\ps

D'une part, l'espace $\mathcal{A}^2(G)$ \'etant l'espace des $G(\widehat{\Z})$-invariants du
$G(\AAA_f)$-module ${\rm L}^2(G(\Q) \backslash G(\AAA),\mu)$ pour les translations
\`a droite, est muni d'une action \`a droite de l'anneau
de Hecke $\mathrm{H}(G)$ (\S\ref{foncteurmx},\,\S\ref{annheckeg}).  Cette action est une
{\it $\star$-action} pour le produit de Petersson. Nous entendons par l\`a que
l'adjoint de $T \in \mathrm{H}(G)$ est l'op\'erateur $T^{\rm t}$ d\'efini
au~\S\ref{corrheckegeneral} : si $f,f' \in
\mathcal{A}^2(G)$ et $T \in \mathrm{H}(G)$,
\begin{equation}\label{staraction}\langle \,T(f), f' \,\rangle_{\rm Pe} = \langle\, f,
T^{\rm t} (f') \,\rangle_{\rm Pe}.\end{equation}
C'est en effet cons\'equence du lemme~\ref{lemmeadjonction} et de l'unimodularit\'e de $G(\AAA_f)$.\ps\ps

D'autre part, $\mathcal{A}^2(G)$ est stable sous l'action de $G(\R)$ par
translations \`a droite, et cette action commute \`a celle de $\mathrm{H}(G)$.  Elle
fait de $\mathcal{A}^2(G)$ une repr\'esentation unitaire du groupe de Lie
$G(\R)$ (nous renvoyons \`a~\cite{knapp} comme r\'ef\'erence g\'en\'erale sur les repr\'esentations unitaires).  Une description plus classique de
cette repr\'esentation s'obtient en \'ecrivant \begin{equation}
\label{fnbclassesG} G(\AAA_f)=\coprod_{i=1}^{\mathrm{h}(G)} G(\Q) g_i
G(\widehat{\Z})\end{equation} pour certains \'el\'ements $g_i  \in G(\AAA_f)$,  par finitude du nombre de classes de $G$.
Pour tout $i$, $G(\Q) g_i G(\widehat{\Z})$ est un ouvert de $G(\AAA_f)$ et 
le groupe ``de congruences'' $$\Gamma_i= G(\Q) \cap g_i
G(\widehat{\Z}) g_i^{-1}$$ est un sous-groupe discret de $G(\R)$
commensurable \`a $G(\Z)$. L'application $f \mapsto (f_{|G(\R) \times
g_i})_i$ induit un isomorphisme $G(\R)$-\'equivariant \begin{equation}
\label{fautreel}\mathcal{A}^2(G) \isomo \prod_{i=1}^{\mathrm{h}(G)} {\rm L}^2(\Gamma_i
\backslash G(\R)), \end{equation} chaque $\Gamma_i\backslash G(\R)$
h\'eritant naturellement d'une mesure de Radon $>0$, invariante \`a droite par $G(\R)$, de masse finie, uniquement d\'etermin\'ee par $\mu$.  Cette
repr\'esentations de $G(\R)$ comporte en g\'en\'eral une partie ``discr\`ete''
notoirement difficile \`a d\'ecrire, ainsi qu'une partie ``continue'' dont
l'\'etude a \'et\'e ramen\'ee par Langlands \`a des parties discr\`etes pour
des groupes $G'$ annexes \cite{langlandseis}.  

\subsection{L'ensemble $\Pi_{\rm disc}(G)$}\label{fautdiscgen}

Nous ne nous int\'eresserons ici qu'\`a la
partie discr\`ete de ${\mathcal A}^2(G)$, c'est-\`a-dire au sous-espace $$\mathcal{A}_{\rm
disc}(G) \subset \mathcal{A}^2(G)$$ d\'efini comme \'etant l'adh\'erence de
la somme des sous-$G(\R)$-repr\'esentations ferm\'ees et topologiquement irr\'eductibles de
$\mathcal{A}^2(G)$.  C'est une repr\'esentation de $G(\R)$ qui est somme
orthogonale d'irr\'eductibles par construction\footnote{Il est utile \`a ce stade de rappeler la version suivante du lemme de Schur. Soient $U$ et $V$ des espaces de Hilbert munis de repr\'esentations unitaires d'un groupe $\Gamma$. On suppose que $U$ est topologiquement irr\'eductible et que $u : U \rightarrow V$ est une application lin\'eaire continue $\Gamma$-\'equivariante non nulle. Alors l'adjoint $u^\ast : V \rightarrow U$ (qui est $\Gamma$-\'equivariant) satisfait $u^\ast \circ u = \lambda {\rm Id}_U$ pour un certain $\lambda \in \R^\times$. En effet, $u^\ast \circ u \in {\rm End}(U)$ est hermitien non nul et commute \`a $\Gamma$, son spectre est donc r\'eduit \`a un point $\{ \lambda \}$ d'apr\`es le th\'eor\`eme spectral. Il vient que $V$ est somme orthogonale de ${\rm Im}(u)$ (qui est ferm\'e) et ${\rm Ker}(u^\ast)$. }, chacune n'intervenant
qu'avec une multiplicit\'e finie d'apr\`es un r\'esultat fondamental d\^u
\`a Harish-Chandra (voir l'introduction de~\cite{harishchandra} ainsi que le
th\'eor\`eme 1 du chapitre 1 loc.  cit., voir
\'egalement~\cite{boreljacquet}).  Autrement dit, si $U$ est une
repr\'esentation irr\'eductible unitaire de $G(\R)$, l'espace
$$\mathcal{A}_U(G):={\rm Hom}_{G(\R)}(U,\mathcal{A}_{\rm disc}(G))={\rm
Hom}_{G(\R)}(U,\mathcal{A}^2(G))$$ est de dimension finie sur $\C$.  On a bien
entendu un isomorphisme canonique \begin{equation}\label{decompreel}
\widehat{\bigoplus_{U \in {\rm Irr}(G(\R))}} U
\otimes \mathcal{A}_U(G) \lisomo \mathcal{A}_{\rm disc}(G) ,\end{equation} ${\rm Irr}(H)$ d\'esignant
l'ensemble des classes d'isomorphisme de repr\'esentations unitaires topologiquement
irr\'eductibles du groupe localement compact $H$.  \ps\ps

La structure de $\mathrm{H}(G)$-module \`a droite de $\mathcal{A}^2(G)$ induit naturellement une structure de $\mathrm{H}(G)$-module \`a droite 
sur $\mathcal{A}_U(G)$. Ce dernier h\'erite \'egalement d'un produit scalaire hermitien pour lequel l'action de ${\rm H}(G)$ est encore une $\star$-action. Par exemple, si $e
\in U$ est non nul fix\'e, et si $\varphi,\varphi' \in \mathcal{A}_U(G)$, on
peut poser $\langle \varphi,\varphi'\rangle= \langle \varphi(e), \varphi'(e)
\rangle_{\rm Pe}$.  Or il est bien connu qu'une sous-$\C$-alg\`ebre de
${\rm M}_n(\C)$ stable par $M \mapsto {}^{\rm t}\!\overline{M}$ est
semi-simple : si $X$ est dans son radical de Jacobson, la
matrice hermitienne $X\,{}^{\rm t}\!\overline{X}$ est nilpotente, donc nulle, ce qui
entra\^ine que $X$ est nulle.  En particulier, $\mathcal{A}_U(G)$ est
semi-simple vue comme repr\'esentation de la $\C$-alg\`ebre
$\mathrm{H}(G)^{\rm opp}\otimes \C$.  \ps\ps

Appelons  {\it repr\'esentation de $(G(\R),{\rm H}(G))$} la
donn\'ee d'un espace de Hilbert muni d'une repr\'esentation unitaire
de $G(\R)$ et d'une structure de module \`a droite sur ${\rm H}(G)$ tels que
l'action de tout \'el\'ement de $G(\R)$ commute avec celle de tout
\'el\'ement de ${\rm H}(G)$. Ces repr\'esentations forment une cat\'egorie $\C$-lin\'eaire de mani\`ere naturelle : un morphisme $E \rightarrow F$ est une application $\C$-lin\'eaire continue $E \rightarrow F$ commutant aux actions de $G(\R)$ et ${\rm H}(G)$. Si $U$ est une repr\'esentation unitaire de $G(\R)$, et si $V$ est un ${\rm H}(G)^{\rm opp}\otimes \C$-module suppos\'e de dimension finie comme $\C$-espace vectoriel, alors $U \otimes V$ est une repr\'esentation de $(G(\R),{\rm H}(G))$ de mani\`ere naturelle (le produit tensoriel \'etant pris au dessus de $\C$). Nous d\'esignerons par $\Pi(G)$ l'ensemble des 
classes d'isomorphisme de repr\'esentations de $(G(\R),\mathrm{H}(G))$ de cette forme, telles que de plus $U$ est topologiquement irr\'eductible et $V$ est simple. La restriction \`a $G(\R)$ d'une telle repr\'esentation unitaire $\pi$ est isomorphe \`a $U^{\dim V}$, de sorte que la classe d'isomorphisme de la repr\'esentation unitaire $U$ est enti\`erement d\'etermin\'ee par la repr\'esentation unitaire de $G(\R)$ sous-jacente \`a $\pi$ : on la note $\pi_\infty$. De m\^eme, le ${\rm H}(G)^{\rm opp}\otimes \C$-module sous-jacent \`a $\pi$ est semi-simple et $V$-isotypique, de sorte que la classe d'isomorphisme du ${\rm H}(G)^{\rm opp}\otimes \C$-module $V$ est uniquement d\'etermin\'ee par celle de $\pi$ : on la note $\pi_f$. En particulier, on a $\pi \simeq \pi_\infty \otimes \pi_f$ pour tout $\pi \in \Pi(G)$. Enfin, le lemme de Schur entra\^ine que tout $\pi \in \Pi(G)$ est topologiquement irr\'eductible en tant que repr\'esentation de $(G(\R),\mathrm{H}(G))$. \ps\ps

D'apr\`es la discussion ci-dessus, si $U \in {\rm Irr}(G(\R))$ alors l'espace $U \otimes \mathcal{A}_U(G)$ est une repr\'esentation de $(G(\R),{\rm H}(G))$ de mani\`ere naturelle, ainsi bien entendu que $\mathcal{A}_{\rm disc}(G)$, l'isomorphisme \eqref{decompreel} commutant trivialement aux actions de $G(\R)$ et ${\rm H}(G)$. Il en r\'esulte que l'on dispose d'une d\'ecomposition 
comme somme hilbertienne d'\'el\'ements de $\Pi(G)$ raffinant la d\'ecomposition~\eqref{decompreel}
\begin{equation}\label{decompdisc} \mathcal{A}_{\rm disc}(G) = \underset{\pi \in
\Pi(G)}{\widehat{\bigoplus}}m(\pi) \, \pi \end{equation} o\`u $m(\pi) \geq 0$ est un entier
appel\'e {\it multiplicit\'e de $\pi$}. Par d\'efinition, si $\pi \in \Pi(G)$ et si $U \simeq \pi_\infty$, alors ${\rm m}(\pi)$ est la multiplicit\'e de $\pi_f$ dans le ${\rm H}(G)^{\rm opp}\otimes \C$-module $\mathcal{A}_U(G)$, qui est semi-simple et de dimension finie. On note $$\Pi_{\rm disc}(G) \subset
\Pi(G)$$ le sous-ensemble des $\pi$ telles que $m(\pi) \neq 0$. \ps\ps

Les \'el\'ements de $\Pi_{\rm disc}(G)$ seront appel\'es {\it repr\'esentations automorphes discr\`etes}\footnote{Le lecteur prendra garde que la d\'efinition utilis\'ee ici ne d\'epend pas que de $G_\Q$ mais bien de $G$ comme $\Z$-groupe. Dans la litt\'erature, nos repr\'esentations automorphes discr\`etes de $G$ sont plus commun\'ement nomm\'ees ``repr\'esentations automorphes discr\`etes de $G(\AAA)$, sph\'eriques (ou non ramifi\'ees) relativement \`a $G(\widehat{\Z})$''. La perte apparente de g\'en\'eralit\'e dans notre pr\'esentation est toutefois illusoire \`a ce stade, \'etant donn\'e que tout sous-groupe compact ouvert de $G(\AAA_f)$ est de la forme $G'(\widehat{\Z})$ pour un $\Z$-groupe $G'$ bien choisi et tel que $G'_\Q \simeq G_\Q$.} de $G$. Le seul exemple vraiment \'evident de telle repr\'esentation est la {\it repr\'esentation triviale}, not\'ee $1_G$,
r\'ealis\'ee comme sous-espace (de dimension $1$) des fonctions constantes de $\mathcal{A}^2(G)$ (noter que $\mu$ est de masse
finie). L'action de $G(\R)$ dans $1_G$ est bien entendu l'action triviale, celle de ${\rm H}(G)$ \'etant la multiplication par le ``degr\'e'' (voir l'exemple VI.\ref{degretriv}).
En g\'en\'eral, l'ensemble $\Pi_{\rm disc}(G)$ est infini d\'enombrable, contrairement \`a $\Pi(G)$. Nous donnerons quelques exemples concrets dans
les chapitres suivants. \ps\ps

Un \'el\'ement $F \in \mathcal{A}_U(G)$ sera dit {\it
propre} s'il est non nul et s'il engendre un ${\rm H}(G)^{\rm opp} \otimes
\C$-module irr\'eductible. Quand ${\rm H}(G)$ est commutative, il est \'equivalent de demander que $F
\neq 0$ soit un vecteur propre de tous les op\'erateurs de Hecke
dans ${\rm H}(G)$. Si $F$ est propre, et si $V \subset \mathcal{A}_U(G)$ d\'esigne le ${\rm H}(G)^{\rm opp}\otimes \C$-module engendr\'e par $F$, l'image de $U \otimes V$ dans $\mathcal{A}_{{\rm disc}}(G)$ par l'application canonique \eqref{decompreel} est une sous-repr\'esentation topologiquement irr\'eductible de $(G(\R),{\rm H}(G))$ not\'ee $\pi_F$ ; c'est la  {\it repr\'esentation (automorphe, discr\`ete) engendr\'ee par $F$}. On note souvent \'egalement encore $\pi_F$ sa classe d'isomorphisme, qui est un \'el\'ement de $\Pi_{{\rm disc}}(G)$. \ps\ps

Enfin, suivant Gelfand, Graev et
Piatetski-Shapiro  \cite[Ch.  3 \S 7]{GGPS}, on dispose toujours du
sous-espace $\mathcal{A}_{\rm cusp}(G) \subset \mathcal{A}^2(G)$ constitu\'e des {\it formes paraboliques} (la d\'efinition d'une forme parabolique est rappel\'ee ci-dessous). Il s'agit d'un sous-espace ferm\'e, stable par les actions de $G(\R)$ et ${\rm H}(G)$. Gelfand, Graev et Piatetski-Shapiro d\'emontrent l'inclusion 
\begin{equation} \label{inclusioncuspdisc} \mathcal{A}_{\rm
cusp}(G) \subset \mathcal{A}_{\rm disc}(G)\end{equation}
\noindent (voir aussi \cite[Theorem 16.2]{borelsl2}).  On note 
$$\Pi_{\rm cusp}(G) \subset \Pi_{\rm disc}(G)$$ l'ensemble des $\pi \in
\Pi(G)$ qui interviennent dans le sous-espace $\mathcal{A}_{\rm cusp}(G)$. \ps\ps

Quand $G_{\Q}$ n'admet pas de sous-$\Q$-groupe parabolique strict, ce qui \'equivaut \`a dire que $G(\Q)$ n'a pas d'\'el\'ement unipotent non trivial, on a l'\'egalit\'e \'evidente $\mathcal{A}_{\rm cusp}(G)=\mathcal{A}^2(G)$. Dans ce cas\footnote{En fait, un r\'esultat fameux de Godement montre que sous cette m\^eme hypoth\`ese sur $G$ le groupe $G(\Q)$ est cocompact dans $G(\AAA)$, ce qui entra\^ine plus directement l'\'egalit\'e $\mathcal{A}_{\rm disc}(G)={\mathcal{A}^2}(G)$ dans ce cas particulier (voir par exemple \cite[Lemma 16.1]{borelsl2}).}, l'inclusion \eqref{inclusioncuspdisc} entra\^ine $\mathcal{A}_{\rm
disc}(G)={\mathcal{A}^2}(G)$.  \ps\ps \ps

Rappelons la d\'efinition d'une forme parabolique. Soit $P \subset G_\Q$ un sous-$\Q$-groupe {\it parabolique  strict}, c'est-\`a-dire tel que $P(\C)$ soit connexe, contienne un sous-groupe de Borel de la composante neutre de $G(\C)$, et ne soit pas \'egal \`a cette composante toute enti\`ere \cite{humphreyslag}\cite{borelgp}. Si $N$ d\'esigne le radical unipotent de $P$, alors le groupe localement compact $N(\AAA)$ est unimodulaire, et son sous-groupe $N(\Q)$ est discret et cocompact. On note ${\rm d}n$ une mesure de Radon ${\rm N}(\AAA)$-invariante $>0$ sur $N(\Q)\backslash N(\AAA)$. Soient $f : G(\Q) \backslash G(\AAA) \rightarrow \C$ une fonction bor\'elienne de carr\'e int\'egrable et $g \in G(\AAA)$. La fonction $n \mapsto f(ng)$, $N(\Q) \backslash N(\AAA) \rightarrow \C$, est alors bor\'elienne, de carr\'e int\'egrable pour presque tout $g \in G(\AAA)$. On dit que $f$ est parabolique si pour tout sous-$\Q$-groupe parabolique strict $P$ de $G_\Q$, on a $\int_{N(\Q)\backslash N(\AAA)} f(ng) \, \, {\rm d}n=0$ pour presque tout $g \in N(\AAA)\backslash G(\AAA)$. On v\'erifie que le sous-ensemble de ${\rm L}^2(G(\Q)\backslash G(\AAA),\mu)$ constitu\'e des classes de fonctions paraboliques est un sous-espace vectoriel ferm\'e (voir par exemple \cite[Prop. 8.2]{borelsl2}). Il est trivialement stable par les translations \`a droite par les \'el\'ements de $G(\AAA)$.

\section{Formes automorphes pour ${\rm O}_n$}\label{cascompact}

\subsection{Formes automorphes des $\Z$-groupes $G$ tels que $G(\R)$ est compact.}\label{fautgrcompact}

On se replace dans le cadre du~\S\ref{fautcarreint}. Supposons que le $\Z$-groupe $G$ a la propri\'et\'e que $G(\R)$ est compact.
Dans ce cas, les groupes
$\Gamma_i=G(\Q) \cap g_i G(\widehat{\Z})g_i^{-1}$ de la formule~\eqref{fnbclassesG} sont des sous-groupes
finis de $G(\R)$, car discrets dans un compact.  La description \eqref{fautreel} montre que $G(\Q)\backslash G(\AAA)$ est compact. De plus, l'\'egalit\'e $\mathcal{A}_{\rm disc}(G)=\mathcal{A}^2(G)$ se d\'eduit directement du th\'eor\`eme de Peter-Weyl.  
Nous allons donner une autre description des $\mathrm{H}(G)$-modules $\mathcal{A}_U(G)$. \ps\ps

Si $U$ est un $\Z[G(\Q)]$-module, d\'esignons par ${\rm M}_U(G)$
l'espace des fonctions $$F : \mathcal{R}(G) \longrightarrow U$$ telles que
$F(\gamma x )=\gamma \cdot F(x)$ pour tous $\gamma \in G(\Q)$ et $x \in
\mathcal{R}(G)$. Observons que ${\rm M}_U(G)$ s'identifie canoniquement \`a 
${\rm Hom}_{\Z[G(\Q)]}(\Z[\mathcal{R}(G)],U)$, ce qui lui conf\`ere une action \`a droite de l'anneau $\mathrm{H}(G)$. Mieux, $U \mapsto
{\rm M}_U(G)$ d\'efinit un foncteur des 
$G(\Q)$-modules vers les $\mathrm{H}(G)^{\rm opp}$-modules. Sa structure additive est tr\`es simple, 
car $F \mapsto
(F(g_i))$ induit un isomorphisme
\begin{equation} \label{mwproduit} {\rm M}_U(G)\longrightarrow \prod_{i=1}^{\mathrm{h}(G)} U^{\Gamma_i}.\end{equation}  
En particulier ${\rm M}_{U\oplus V}(G) \simeq {\rm M}_U(G) \oplus {\rm
M}_V(G)$. Observons d'ailleurs que la construction ci-dessus 
a un sens pour tous les $\Z$-groupes $G$. \ps\ps

Supposons maintenant que $U$ est une repr\'esentation complexe, continue, de dimension finie,  de $G(\R)$, et notons $U^\ast$ son dual.
Si $F \in {\rm M}_U(G)$ et $u \in U^\ast$, on note $\varphi_F(u)$ la fonction $(h, x) \mapsto \langle u, 
h^{-1} F(x) \rangle$, $G(\R) \times \mathcal{R}(G) \rightarrow \C$. C'est une fonction continue ; elle est donc dans $\mathcal{A}^2(G)$ car  $G(\Q) \backslash (G(\R) \times \mathcal{R}(G))$ est compact d'apr\`es \eqref{fautreel}. La relation \'evidente $\varphi_F(g u) = g \cdot (\varphi_F(u))$, valable pour $u \in U^\ast$ et $g \in G(\R)$, montre que la fonction $\varphi_F,  u \mapsto \varphi_F(u),$ est un \'el\'ement de $\mathcal{A}_{U^\ast}(G)$. La d\'emonstration de lemme suivant est imm\'ediate, et laiss\'ee au lecteur.

\begin{lemme}\label{fautcomp} Supposons que $U$ est une repr\'esentation irr\'eductible de $G(\R)$. Alors  $F \mapsto \varphi_F$ est un isomorphisme
$\mathrm{H}(G)$-\'equivariant ${\rm M}_U(G) \isomo \mathcal{A}_{U^\ast}(G)$.\end{lemme}

Depuis l'article de Gross~\cite{grossalg}, les \'el\'ements de ${\rm M}_U(G)$ sont
parfois appel\'es {\it formes modulaires alg\'ebriques} de poids $U$ du $\Z$-groupe $G$. Par exemple, si $U=\C$ est la
repr\'esentation triviale, le $\mathrm{H}(G)^{\rm opp}$-module ${\rm M}_\C(G)$ s'identifie canoniquement \`a l'espace
des fonctions ${\rm Cl}(G) \rightarrow \C$, soit encore au dual du
$\mathrm{H}(G)$-module $\C[{\rm Cl}(G)]$. \ps\ps

Terminons ces g\'en\'eralit\'es par une assertion de compatibilit\'e \`a certains morphismes de $\Z$-goupes. 
Soit $\mu : G \rightarrow G'$ un morphisme de $\Z$-groupes. Il induit de mani\`ere \'evidente un morphisme du $G(\AAA_{f})$-ensemble $\mathcal{R}(G)$ vers le $G'(\AAA_{f})$-ensemble $\mathcal{R}(G')$, not\'e $(f_{\mu},g_{\mu})$, au sens du \S \ref{parisogenies}. On suppose que $f_{\mu}(G(\AAA_f))$ et distingu\'e dans $G'(\AAA_f)$, que $g_\mu$ est injective, et que de plus l'action du groupe $S$ d\'efini {\it loc. cit.} sur $\mathcal{R}(G)$ est triviale. C'est par exemple trivialement le cas si $\mu$ est un isomorphisme.  On dispose alors d'un homomorphisme injectif d'anneaux ${\rm H}(\mu) : \mathrm{H}(G) \rightarrow \mathrm{H}(G')$ d\'efini {\it loc. cit}. Soit $U'$ un $G'(\Q)$-module, et soit $U$ sa restriction \`a $G(\Q)$. Le lemme suivant paraphrase le lemme~\ref{compatibilitesatake}.

\begin{lemme}\label{compatibilitesatake2} Le morphisme $\mu^\ast : {\rm M}_{U'}(G') \longrightarrow {\rm M}_{U}(G)$, $\varphi \mapsto (x \mapsto \varphi(g_\mu(x)))$, 
satisfait $T \circ \mu^\ast = \mu^\ast \circ {\rm H}(\mu)(T)$ pour tout $T \in \mathrm{H}(G')$.
\end{lemme}  \ps\ps

\subsection{Cas des groupes ${\rm O}_n$ et ${\rm SO}_n$.}\label{fauton}
\label{fautson}Sp\'ecifions maintenant cette construction au $\Z$-groupe
orthogonal ${\rm O}_n$ du r\'eseau unimodulaire pair ${\rm E}_n \subset
\R^n$, pour $n \equiv 0 \bmod 8$ (\S II.3, choisir un autre r\'eseau conduirait \`a
une th\'eorie \'equivalente \`a celle qui suit). \ps\ps

Dans ce cas, nous avons vu au~\S\ref{ensclassbil} que le ${\rm
O}_n(\AAA_f)$-ensemble $\mathcal{R}({\rm O}_n)$ s'identifie canoniquement
\`a l'ensemble des r\'eseaux unimodulaires pairs de $\R^n$ inclus dans ${\rm
E}_n \otimes \Q$, puis que ${\rm Cl}({\rm O}_n)\isomo {\rm X}_n$.  En
particulier, $${\rm M}_\C({\rm O}_n)=\C[{\rm X}_n]^\ast.$$ L'action \`a
droite de $\mathrm{H}({\rm O}_n)$ sur ${\rm M}_\C({\rm O}_n)$ d\'efinit par
transpos\'ee une action \`a gauche de $\mathrm{H}({\rm O}_n)$ sur $\C[{\rm
X}_n]$.  L'op\'erateur ${\rm T}_{\Z/d\Z} \in \mathrm{H}({\rm O}_n)$ d\'efini
au~\S\ref{annheckeclass}, vu comme endomorphisme de $\C[{\rm X}_n]$, est
notamment l'op\'erateur ${\rm T}_d$ du~\S III.2.  La description de la structure
du $\mathrm{H}({\rm O}_n)^{\rm opp}$-module ${\rm M}_\C({\rm O}_n)$ quand
$n\leq 24$ est donc le th\`eme principal de ce m\'emoire. \ps\ps

L'anneau $\mathrm{H}({\rm O}_n)$ est commutatif d'apr\`es la
proposition~\S\ref{heckeocomm}.  Fixons $U$ une repr\'esentation (complexe,
continue, de dimension finie) de ${\rm O}_n(\R)$.  D'apr\`es le
lemme~\ref{fautcomp} et les r\'esultats g\'en\'eraux rappel\'es
au~\S\ref{fautcarreint}, l'action de $\mathrm{H}({\rm O}_n)$ est donc
co-diagonalisable sur chaque ${\rm M}_U({\rm O}_n)$.  Les valeurs propres de
ces op\'erateurs ont des significations arithm\'etiques importantes. Nous
verrons au~\S VI.\ref{correpgalOn} qu'elles sont en fait reli\'ees, de
mani\`ere assez surprenante {\it a priori}, aux repr\'esentations du groupe de
Galois absolu de $\Q$. La droite des fonctions constantes dans ${\rm
M}_\C({\rm O}_n)$ est par exemple trivialement stable par ${\rm T}_A$ pour tout
$A$, la valeur propre de ${\rm T}_p$ sur cette droite \'etant bien s\^ur ${\rm c}_n(p)$
(Prop.-D\'ef. III.2.1). Nous donnerons des exemples nettement plus int\'eressants
aux chapitres suivants.  \ps\ps


\begin{remarque}\label{lnonaf}{\rm  Soit $\mathcal{L}_n$ l'ensemble de tous les r\'eseaux unimodulaires pairs de $\R^n$ d\'ej\`a consid\'er\'e dans l'introduction (\S \ref{chapintro}). Il contient ${\mathcal R}({\rm O}_{n})$ et l'action naturelle de ${\rm O}_n(\R)$ sur $\mathcal{L}_{n}$ prolonge l'action naturelle de ${\rm O}_n(\Q)$ sur ${\mathcal R}({\rm O}_n)$. L'application ${\rm O}_n(\R) \times {\rm O}_n(\AAA_f) \rightarrow \mathcal{L}_n$, $(g_\infty,g_f) \mapsto g_\infty^{-1}(g_f({\rm E}_n)),$ se factorise donc en une application 
$${\rm O}_n(\Q) \backslash {\rm O}_n(\AAA) / {\rm O}_n(\widehat{\Z}) \rightarrow \mathcal{L}_n.$$
C'est une bijection : la surjectivit\'e se d\'eduit du scholie II.2.1 et l'injectivit\'e est imm\'ediate. }
\end{remarque}

Discutons maintenant du cas de ${\rm SO}_n$. D'apr\`es la
proposition~\ref{rso=ro} et le~\S\ref{hsovso}, l'inclusion ${\rm SO}_n
\rightarrow {\rm O}_n$ induit une bijection $\mathcal{R}({\rm SO}_n) \isomo
\mathcal{R}({\rm O}_n)$ et $\mathrm{H}({\rm O}_n)$ est naturellement un
sous-anneau de $\mathrm{H}({\rm SO}_n)$.  Soit $U$ un ${\rm
SO}_n(\Q)$-module, et soit $$U'= {\rm Ind}_{{\rm SO}_n(\Q)}^{{\rm O}_n(\Q)}
U.$$  La propri\'et\'e universelle des modules induits fournit un
isomorphisme canonique ${\rm ind} : {\rm Hom}_{\Z[{\rm SO}_n(\Q)]}(\Z[\mathcal{R}({\rm
O}_n)]_{|{\rm SO}_n(\Q)},U) \isomo {\rm Hom}_{\Z[{\rm
O}_n(\Q)]}(\Z[\mathcal{R}({\rm O}_n)],U'),$ qui s'\'ecrit encore $${\rm ind} : {\rm
M}_U({\rm SO}_n) \isomo {\rm M}_{U'}({\rm O}_n).$$  
Cet isomorphisme est trivalement
$\mathrm{H}({\rm O}_n)$-\'equivariant, de sorte que l'\'etude des
$\mathrm{H}({\rm O}_n)$-modules ${\rm M}_U({\rm SO}_n)$ se ram\`ene \`a
celle des ${\rm M}_W({\rm O}_n)$ o\`u $W$ est un ${\rm
O}_n(\Q)$-module. Rajoutons que si $U$ est la restriction \`a ${\rm
SO}_n(\Q)$ d'un ${\rm SO}_n(\R)$-module $V$, et si $V'$ d\'esigne l'induite
de $V$ \`a ${\rm O}_n(\R)$, alors $V'_{|{\rm O}_n(\Q)} = {\rm Ind}_{{\rm
SO}_n(\Q)}^{{\rm O}_n(\Q)} U$. \ps\ps
Supposons enfin que $W$ est un ${\rm O}_n(\Q)$-module et que $W'$ d\'esigne
sa restriction \`a ${\rm SO}_n(\Q)$. Le groupe ${\rm O}_n(\Q)$ agit naturellement sur ${\rm
M}_{W'}({\rm SO}_n)$, par
$(\gamma,f) \mapsto (x \mapsto \gamma(f(\gamma^{-1}(x))))$, le sous-groupe ${\rm SO}_n(\Q)$
agissant trivialement. Soit $s \in {\rm End}({\rm M}_{W'}({\rm
SO}_n))$ l'op\'erateur induit par l'\'el\'ement non-trivial du quotient ${\rm
O}_n(\Q)/{\rm SO}_n(\Q)\simeq \Z/2\Z$. La restriction
des fonctions par l'application bijective $\mathcal{R}({\rm SO}_n) \rightarrow
\mathcal{R}({\rm O}_n)$
d\'efinit alors une injection ${\rm H}({\rm O}_n)$-\'equivariante 
$${\rm res} : {\rm M}_W({\rm O}_n) \rightarrow {\rm M}_{W'}({\rm SO}_n)$$
dont l'image est ${\rm M}_{W'}({\rm SO}_n)^{s={\rm id}}$.\ps\ps

\begin{example}\label{exdim1det} {\rm En guise d'exemple, ${\rm ind}$ induit une d\'ecomposition canonique $${\rm M}_\C({\rm SO}_n) \simeq {\rm M}_\C({\rm O}_n) \oplus {\rm M}_{\DET}({\rm O}_n)$$ o\`u
$\DET$ est la repr\'esentation de dimension $1$ donn\'ee par le
d\'eterminant. Regardant \'egalement $\C$ comme restriction \`a ${\rm
SO}_n(\R)$ de la repr\'esentation triviale de ${\rm O}_n(\R)$, cela munit ${\rm M}_\C({\rm
SO}_n)$ d'une sym\'etrie $s$ qui pr\'eserve la d\'ecomposition ci-dessus,
admettant pour points fixes ${\rm M}_\C({\rm O}_n)$. } \end{example}

Nous renvoyons \`a~\cite[\S 2]{chrenard2} pour une \'etude
des espaces ${\rm M}_U({\rm SO}_8)$ en fonction de la
repr\'esentation $U$, notamment de leur dimension, ainsi qu'au \S \ref{tableexempleso8} pour des exemples.

\subsection{Un produit scalaire hermitien invariant}\label{prodherminv}

Consid\'erons le cas d'un $\Z$-groupe $G$ g\'en\'eral tel que $G(\R)$ est compact. Soit $U$ une repr\'esentation complexe, continue, de dimension finie, de $G(\R)$. L'isomorphisme ${\rm M}_U(G) \isomo \mathcal{A}_{U^\ast}(G)$ munit par transport de structure ${\rm M}_U(G)$ d'un produit scalaire hermitien naturel, pour lequel l'action de $\mathrm{H}(G)$ est une $\star$-action, d'apr\`es le~\S\ref{fautcarreint}, et qu'il ne reste qu'\`a expliciter.
Fixons pour cela un produit scalaire hermitien $G(\R)$-invariant sur $U$, disons  $\langle-, -\rangle_U$. Choisissons \'egalement des \'el\'ements $g_{i} \in G(\AAA_{f})$ v\'erifiant l'\'egalit\'e \eqref{fnbclassesG} ; on rappelle que $\Gamma_{i}=G(\Q) \cap g_{i} G(\widehat{\Z}) g_{i}^{-1}$ est un groupe fini.
\ps\ps

\begin{prop}\label{hermitien} Si $F,F' \in {\rm M}_U(G)$, $\left(F,F'\right)=\sum_{i=1}^{\mathrm{h}(G)}
\frac{1}{|\Gamma_i|}\langle F(g_i),F'(g_i)\rangle_U$ d\'efinit un
produit scalaire hermitien sur ${\rm M}_U(G)$ ind\'ependant du choix des $g_{i}$, et pour lequel l'action de $\mathrm{H}(G)$ est
une $\star$-action.
\end{prop}

Par manque de r\'ef\'erence ad\'equate nous donnons une d\'emonstration.\ps\ps

\begin{pf} Fixons $e \in U^\ast$ non nul. D'apr\`es l'isomorphisme~\eqref{fautcomp} et le~\S \ref{paysageaut}
$$\left(F,F'\right) : = \int_{G(\Q) \backslash G(\AAA)}  \overline{\varphi_F(e)} \varphi_{F'}(e) \,{\rm d}m$$
est un produit scalaire hermitien sur ${\rm M}_U(G)$ pour lequel l'action de $\mathrm{H}(G)$ est une $\star$-action. Nous allons v\'erifier qu'il est proportionnel \`a celui de l'\'enonc\'e. \ps

Soient $\Omega_i \subset G(\AAA)$ l'ouvert compact $g_i (G(\R) \times G(\widehat{\Z}))$, $\pi : G(\AAA) \rightarrow G(\Q) \backslash G(\AAA)$ la projection canonique, et $\overline{\Omega_{i}}=\pi(\Omega_{i})$. Par d\'efinition, $G(\Q) \backslash G(\AAA)$ est r\'eunion disjointe (finie) des $\overline{\Omega_{i}}$. V\'erifions d'abord qu'il existe une mesure de Haar $m$ sur $G(\AAA)$ telle que pour toute fonction continue $\psi$ sur (le compact) $G(\Q) \backslash G(\AAA)$, on ait la formule \begin{equation}\label{formintgqga} \int_{G(\Q) \backslash G(\AAA)} \psi \, {\rm d}\mu = \sum_{i=1}^{{\rm h}(G)} \frac{1}{|\Gamma_{i}|} \int_{\Omega_{i}} \psi \circ \pi \, \, {\rm d}m.\end{equation}
En effet, rappelons que si $f$ est continue \`a support compact sur $G(\AAA)$, alors $\tilde{f}(g):= \sum_{\gamma \in G(\Q)} f(\gamma g)$ est continue \`a support compact sur $G(\Q)\backslash G(\AAA)$. De plus, par la propri\'et\'e caract\'eristique de la mesure quotient $\mu$ il existe une unique mesure de Haar $m$ sur $G(\AAA)$ telle que pour toute fonction continue \`a support compact $f$ sur $G(\AAA)$  alors 
$\int_{G(\AAA)} f \,{\rm d}m= \int_{G(\Q)\backslash G(\AAA)} \tilde{f} \, {\rm d}\mu$ (voir \cite[Chap. II]{weil}). \ps

Si $g \in G(\AAA)$ posons $n_{i}(g)=|G(\Q)g \cap \Omega_{i}|$. On a clairement $n_{i}(\gamma g k)=n_{i}(g)$ pour tout $\gamma \in G(\Q)$ et tout $k \in 1 \times G(\widehat{\Z})$. Par d\'efinition, on a \'egalement $n_{i}(g_{j})= |\Gamma_{i}| \delta_{i,j}$, la notation $\delta_{{i,j}}$ d\'esignant le symbole de Kronecker. Soit $\psi$ une fonction continue sur $G(\Q) \backslash G(\AAA)$. La fonction $f_{i} = {\rm 1}_{\Omega_{i}} \times \psi \circ \pi$, $G(\AAA) \rightarrow \C$, est continue \`a support dans $\Omega_{i}$, et satisfait $\tilde{f_{i}}(g)=\psi(\pi(g))\, n_{i}(g)$ pour tout $g \in G(\AAA)$ (nous notons ${\rm 1}_{A}$ la fonction carcat\'eristique de l'ensemble $A$). Autrement dit, on a $\psi \times {\rm 1}_{\overline{\Omega_{i}}} = \frac{1}{ |\Gamma_{i}|} \tilde{f_{i}}$. Cela d\'emontre la formule \eqref{formintgqga}.\ps

Appliquons cette formule \`a la fonction $\psi = \overline{\varphi_F(e)} \, \varphi_{F'}(e)$. Observons que si $U=\C$, de sorte que $\psi$ est constante \'egale \`a $|e(1)|^{2}\overline{F(g_{i})} F'(g_{i})$ sur $\Omega_{i}$, la proposition r\'esulte de ce que $m(\Omega_{i})=m(G(\R) \times G(\Z))$ est ind\'ependant de $i$. En g\'en\'eral, introduisons ${\rm d} g$ la mesure de Haar sur $G(\R)$ de masse totale $1$, et $m_f$ la mesure de Haar sur $G(\AAA_f)$ telle que ${\rm d} m =  {\rm d} g \times {\rm d} m_f$. L'invariance de $\psi$ \`a droite par $1 \times G(\widehat{\Z})$ et le th\'eor\`eme de Fubini entra\^inent $$\int_{\Omega_{i}} \psi \circ \pi \, \, {\rm d}m = m_{f}(G(\widehat{\Z})) \int_{G(\R)} \overline{\langle e, g^{-1}F(g_{i}) \rangle }\langle e, g^{-1}F'(g_{i}) \rangle \, {\rm d} g.$$
Soit $E \in U$ tel que $ \langle E, x \rangle_{U} = \langle e, x \rangle$ pour tout $x \in U$. Les relations d'orthogonalit\'e des coefficients matriciels pour les repr\'esentations irr\'eductibles des groupes compacts montrent $$\int_{G(\R)} \overline{\langle e, g^{-1}F(g_{i}) \rangle } \, \langle e, g^{-1}F'(g_{i}) \rangle \, {\rm d} m_{\infty} = \frac{1}{{\rm dim} \,U} \langle E, E \rangle_{U} \langle F(g_{i}), F'(g_{i}) \rangle_{U},$$ ce qui conclut la d\'emonstration de la proposition. \end{pf}

Supposons par exemple $G={\rm O}_n$ et $U=\C$. Si $L_i \in \mathcal{R}_\Z^{\rm a}({\rm E}_n\otimes \Q)$ d\'esigne le r\'eseau
$g_i(L)$, on a $\Gamma_i={\rm O}(L_i) \subset {\rm O}_n(\Q)$.  La relation ${\rm T}_A={\rm T}_A^{\rm t}$ de la proposition~\ref{heckeocomm} et la proposition pr\'ec\'edente s'\'ecrivent alors
$${\rm N}_A(L,M)|{\rm O}(M)|={\rm N}_A(M,L)|{\rm O}(L)|,$$
o\`u ${\rm N}_A(L,M)$ d\'esigne le nombre des $A$-voisins de $L$ 
isom\'etriques \`a $M$ (avec $L,M \in \mathcal{R}({\rm O}_n)$). C'est la g\'en\'eralisation 
que nous avions promise de la proposition III.2.3. \ps\ps

\begin{cor} \label{corhermitien} La forme bilin\'eaire 
sur ${\rm M}_{U^{\ast}}(G) \times {\rm M}_{U}(G)$ d\'efinie par $$\left(F,F'\right)=\sum_i \frac{1}{|\Gamma_i|} \langle F(g_i), F'(g_i) \rangle$$
est ind\'ependante du choix des $g_{i}$, et non-d\'eg\'en\'er\'ee. Elle satisfait l'identit\'e $\left(T(F),F'\right)=(F,T^{\rm t}(F'))$ pour tous $T \in \mathrm{H}(G)$, $F \in {\rm M}_{U^{\ast}}(G)$ et  $F' \in {\rm M}_{U}(G)$. En particulier, elle d\'efinit  un isomorphisme canonique entre le $\mathrm{H}(G)$-module ${\rm M}_{U^*}(G)^*$ et le ${\rm H}(G)$-module  ${\rm M}_U(G)^{{\rm t}}$ (voir la remarque~\ref{transposeHXmod}).
\end{cor}

\begin{pf} Si $V$ est un $\C$-espace vectoriel, on note $\overline{V}$ le $\C$-espace vectoriel conjugu\'e (\`a savoir le groupe ab\'elien $V$ muni de l'action de $\C$ d\'efinie par $(\lambda,v) \mapsto \overline{\lambda}v$, $\C \times V \rightarrow V$). Si $U$ est comme dans l'\'enonc\'e, $\overline{U}$ est une repr\'esentation de $G(\R)$ de mani\`ere naturelle, et l'application $v \mapsto (u \mapsto \langle v, u \rangle_{U})$ induit un isomorphisme de repr\'esentations $\overline{U} \isomo U^{\ast}$. On dispose donc d'un isomorphisme naturel ${\rm M}_{U^{\ast}}(G) \isomo {\rm M}_{\overline{U}}(G) = \overline{{\rm M}_{U}(G)}$. Via cet isomorphisme, la forme bilin\'eaire de l'\'enonc\'e co\"incide avec la forme $\overline{{\rm M}_{U}(G)} \times {\rm M}_{U}(G) \rightarrow \C$, 
$(F,F') \mapsto \sum_{i} \frac{1}{|\Gamma_i|} \langle F(g_i), F'(g_i) \rangle_{U}$, qui n'est autre que la forme hermitienne sur ${\rm M}_{U}(G)$ donn\'ee par la proposition \ref{hermitien}. On en d\'eduit les deux premi\`eres assertions, la derni\`ere \'etant alors \'evidente.
\end{pf}

Terminons par une derni\`ere observation. Si $L \in \mathcal{R}(G)$ et $u \in U$, l'application $ F \mapsto \langle F(L),u\rangle$ est une forme lin\'eaire sur ${\rm M}_{U^{\ast}}(G)$, not\'ee ${\rm ev}_{{L,u}}$. On dispose d'une unique application lin\'eaire $$\Z[\mathcal{R}(G)] \otimes U \longrightarrow {\rm M}_{U^{\ast}}(G)^{\ast}$$
envoyant $[L] \otimes u$ sur ${\rm ev}_{{L,u}}$ pour tout $L \in \mathcal{R}(G)$ et tout $u \in U$. Le $\C$-espace vectoriel $\Z[\mathcal{R}(G)] \otimes U$ est muni d'une action diagonale de $G(\Q)$, et l'application ci-dessus est constante sur les orbites de cette action. Elle se factorise donc en une application lin\'eaire
\begin{equation} \label{coinvmu} (\Z[\mathcal{R}(G)] \otimes U)_{G(\Q)} \rightarrow {\rm M}_{U^{\ast}}(G)^{\ast},
\end{equation}
la notation $V_{\Gamma}$ d\'esignant les co-invariants du $\Gamma$-module $V$. C'est un isomorphisme : cela se d\'eduit en effet simplement de la finitude de $G(\Q)\backslash \mathcal{R}(G)$ et de l'isomorphisme naturel $ (U^{\ast})^{\Gamma} \isomo (U_{\Gamma})^{\ast} $, valable pour tout sous-groupe fini $\Gamma$ de $G(\R)$. L'isomorphisme \eqref{coinvmu} commute trivialement aux actions naturelles (\`a gauche) de ${\rm H}(G)$. Si on le compose avec l'isomorphisme 
${\rm M}_{U^*}(G)^* \longrightarrow {\rm M}_U(G)^{{\rm t}}$ donn\'e par le corollaire pr\'ec\'edent, 
on en d\'eduit un isomorphisme canonique de ${\rm H}(G)$-modules
\begin{equation} \label{coinvmudual} (\Z[\mathcal{R}(G)] \otimes U)_{G(\Q)} \isomo {\rm M}_{U}(G)^{{\rm t}}. \end{equation}
Il envoit (la classe de) l'\'el\'ement $[L] \otimes u$ sur un \'el\'ement de ${\rm M}_{U}(G)$ que nous noterons $[L,u]$. Concr\`etement, $[L,u]$ est l'unique fonction $F \in {\rm M}_{U}(G)$ nulle hors de $G(\Q) \cdot L$ telle que 
$F(L) =\sum_{\gamma \in \Gamma} \gamma(u)$, o\`u $\Gamma=G(\Q)_{L}$ est le stabilisateur de $L$ dans $G(\Q)$.
L'isomorphisme \eqref{coinvmudual} jouera un (petit) r\^ole dans les consid\'erations de s\'eries th\^eta au \S V.\ref{sertharm} et au chapitre VII.

\section{Formes modulaires de Siegel}\label{formesiegel} 

Commen\c{c}ons par quelques rappels sur les formes modulaires de Siegel
(voir \cite{semcartan1}, \cite{semcartan2}, \cite{freitag},
\cite{andrianov}).  Nous suivrons de pr\`es l'exposition de van der Geer
\cite{vandergeer}, \`a laquelle nous renvoyons notamment pour un historique du
sujet.  \newcommand{\Sp}{{\mathrm{Sp}}}
\newcommand{\GSp}{{\mathrm{GSp}}}\newcommand{\PGSp}{{\mathrm{PGSp}}}

\subsection{Point de vue classique}\label{fsiegelclass}

Soit $g\geq 1$ un entier. Si $R$ est un anneau nous noterons
${\rm Mat}_g(R)$ l'ensemble des matrices de taille $g \times g$ \`a
coefficients dans $R$, et ${\rm Sym}_g(R) \subset {\rm Mat}_g(R)$ le
sous-ensemble des matrices sym\'etriques. On note $1_g$ la matrice identit\'e dans ${\rm Mat}_g(R)$ et on d\'esigne par ${\rm J}_{2g} \in {\rm Mat}_{2g}(R)$ l'\'el\'ement $${\rm J}_{2g} = \begin{pmatrix} 0 & 1_g \\ -1_g & 0 \end{pmatrix}.$$

Le demi-espace de Siegel de genre $g$ est l'ouvert $$\mathbb{H}_g \subset
{\rm Sym}_g(\C)$$ des matrices dont la partie imaginaire est d\'efinie
positive.  Le $\Z$-groupe ${\rm GSp}_{2g}$ sera vu comme le sous-sch\'ema en
groupes de $\GL_{2g}$ constitu\'e des $\gamma$ tels que $\gamma \,{\rm
J}_{2g} \,{}^{\rm t}\gamma = \nu(\gamma) {\rm J}_{2g}$; le morphisme $\nu : {\rm GSp}_{2g} \rightarrow \mathbb{G}_{m}$ d\'esignant le facteur de similitude.
Ses \'el\'ements sont de la forme
$$\gamma = \begin{pmatrix} a_{\gamma} & b_{\gamma} \\ c_{\gamma} & d_{\gamma} \end{pmatrix}$$ avec $a_{\gamma},b_{\gamma},c_{\gamma},d_{\gamma} \in {\rm
Mat}_g$ satisfaisant les relations $a_{\gamma} \,{}^{\rm t}b_{\gamma} = b_{\gamma} \,{}^{\rm t}a_{\gamma}$, $c_{\gamma} \, {}^{\rm t}d_{\gamma} = d_{\gamma} \,{}^{\rm t}c_{\gamma}$ et
$a_{\gamma}\,{}^{\rm t}d_{\gamma}  - b_{\gamma}\, {}^{\rm t}c_{\gamma} = \nu(\gamma) 1_g$.  \ps\ps

Soit $\GSp_{2g}(\R)^{+}$ le sous-groupe de $\GSp_{2g}(\R)$ constitu\'e des \'el\'ements de facteur de similitude $>0$. Si $\gamma \in \GSp_{2g}(\R)^{+}$ et $\tau \in \mathbb{H}_g$, on v\'erifie que
l'\'el\'ement ${\rm j}(\gamma,\tau):=c_{\gamma}\tau+d_{\gamma}$ est dans $ \GL_g(\C)$ et que
$$(\gamma,\tau) \mapsto \gamma\tau=(a_{\gamma}\tau+b_{\gamma})(c_{\gamma}\tau+d_{\gamma})^{-1}$$ d\'efinit une
action transitive de $\GSp_{2g}(\R)^{+}$ sur $\mathbb{H}_g$ par transformations
biholomorphes.  De plus, on v\'erifie la relation de $1$-cocyle ${\rm
j}(\gamma\gamma',\tau)={\rm j}(\gamma,\gamma'\tau){\rm j}(\gamma',\tau)$
pour tous $\gamma,\gamma'\in {\rm GSp}_{2g}(\R)^{+}$ et tout $\tau \in {\mathbb
H}_g$.  \ps\ps

Soit $W$ un $\C$-espace vectoriel de dimension finie muni d'une
$\C$-repr\'esenta\-tion $\rho : \GL_g \rightarrow \GL_W$.  Une {\it
forme modulaire de Siegel de poids $W$ et de genre $g\geq 1$} est 
une fonction holomorphe $f : \mathbb{H}_g \rightarrow W$ telle que $$
f(\gamma\tau)=\rho(\,{\rm j}(\gamma,\tau)\,) \cdot f(\tau), \, \, \forall
\tau \in \mathbb{H}_g, \, \, \, \forall \gamma \in \Sp_{2g}(\Z).$$ Lorsque
$g=1$ on rajoute l'hypoth\`ese que $f$ est born\'ee sur $\{\tau \in
\mathbb{H}_1, {\rm Im}(\tau) >1\}$.  Elles forment un
$\C$-espace vectoriel not\'e $${\rm M}_W({\rm
Sp}_{2g}(\Z)),$$
dont la dimension est finie d'apr\`es Siegel. \ps\ps

Lorsque $(\rho,W)= (\DET^k,\C)$ pour $k \in \Z$, on parle
de {\it formes de Siegel classiques}, ou scalaires, de poids $k$, et par
opposition de formes \`a valeurs vectorielles sinon.  Dans ce cas, on note
aussi ${\rm M}_k({\rm Sp}_{2g}(\Z))$ l'espace ${\rm M}_W({\rm Sp}_{2g}(\Z))$. 
Quand $g=1$, on retrouve comme cas particulier les formes modulaires
usuelles pour le groupe ${\rm SL}_2(\Z)$ qui sont par exemple trait\'ees en
d\'etail dans le cours d'arithm\'etique de Serre~\cite{serre}. Observons enfin que la pr\'esence de l'\'el\'ement $- 1_{{\rm 2g}} \in {\rm Sp}_{{2g}}(\Z)$, et la relation ${\rm j}(-1_{2g},\tau) = -1_g$, entra\^inent ${\rm M}_{W}({\rm Sp}_{2g}(\Z))=0$ si $\rho(- {\rm 1}_{g}) = - {\rm id}_W$. \ps\ps

Terminons cette section par une reformulation de la notion de forme modulaire de Siegel. Supposons que la repr\'esentation $(\rho,W)$ est irr\'eductible, ou plus g\'en\'eralement qu'il existe un \'el\'ement ${\rm m}_{W} \in \Z$, n\'ecessairement unique, tel que $\rho(z {\rm 1}_{g})=z^{{\rm m}_{W}} \,{\rm id}_W$ pour tout $z \in \C^\times$. Si $f$ est une application $\mathbb{H}_g \rightarrow W$, on pose $$f_{{|}_{W}}\gamma : \mathbb{H}_g \rightarrow W, \, \,  \, \, \, \tau \mapsto \,\nu(\gamma)^{{\rm m}_{W}/2}\, \,\rho({\rm j}(\gamma,\tau))^{{-1}} \,\,f(\gamma \tau).$$ L'application
$(\gamma,f) \mapsto f_{{|}_{W}}\gamma$ d\'efinit une action \`a droite du groupe ${\rm GSp}_{2g}(\R)^{+}$ sur
l'espace des fonctions holomorphes $\mathbb{H}_g \rightarrow W$ ; par construction, cette action est  triviale sur le sous-groupe des homoth\'eties de rapport $>0$ dans ${\rm GSp}_{2g}(\R)^{+}$.  Une {\it
forme modulaire de Siegel de poids $W$ et de genre $g\geq 2$} est par d\'efinition un
\'el\'ement ${\rm Sp}_{2g}(\Z)$-invariant pour cette action. 

\subsection{D\'eveloppements de Fourier et formes paraboliques}\label{devfouriersiegel}

Si $n \in {\rm Sym}_g(\C)$, on pose $$q^n=e^{2i\pi\, {\rm tr}(n\,
\tau)}=\prod_{1\leq i,j \leq g} e^{2 i\pi \,n_{i,j}\tau_{i,j}},$$ c'est une
fonction homolorphe sur $\mathbb{H}_g$.  Si $n$ est {\it demi-entier}, c'est-\`a-dire si $n \in \frac{1}{2}{\rm Sym}_g(\Z)$ et si $n_{i,i} \in \Z$ pour
tout $i=1,\dots, g$, alors $q^n$ est invariante par translations par
${\rm Sym}_g(\Z)$. On montre que chaque $f \in {\rm M}_W({\rm Sp}_{2g}(\Z))$ admet 
un d\'eveloppement de Fourier, normalement convergent sur tout compact de $\mathbb{H}_{g}$, de la forme $$f=\sum_{n\geq 0} a_n q^n$$
o\`u $n \in \frac{1}{2}{\rm Sym}_g(\Z)$ parcourt les \'el\'ements positifs
(au sens des matrices sym\'etriques r\'eelles) demi-entiers, et o\`u les
$a_n$ sont dans $W$ \cite[\S 4]{vandergeer}.  Pour $g\geq 2$, l'op\'erateur de Siegel est
l'application $$\Phi_g : {\rm M}_W(\Sp_{2g}(\Z)) \longrightarrow {\rm
M}_{W'}({\rm Sp}_{2g-2}(\Z))$$ d\'efinie par $\Phi_g(\sum_n a_n q^n)=\sum_{n'}
a_{n'} q^{n'}$, o\`u l'on voit ${\rm Sym}_{g-1}(-)$ comme sous-ensemble de
${\rm Sym}_g(-)$ de derni\`ere ligne et colonne nulles, et $W'=W_{|\GL_{g-1} \times 1}$ \cite[\S 5]{vandergeer}.  Le
sous-espace des formes paraboliques est $${\rm S}_W({\rm Sp}_{2g}(\Z)):={\rm
Ker}(\Phi_g) \subset {\rm M}_W({\rm Sp}_{2g}(\Z)).$$ Une forme de Siegel est
donc parabolique si son d\'eveloppement de Fourier $\sum_n a_n q^n$
satisfait $a_n=0$ pour tout $n$ tel que $\DET(n)=0$. Quand $(W,\rho) =
(\C,\DET^k)$ pour $k \in \Z$, on note ${\rm S}_k({\rm Sp}_{2g}(\Z))$ pour
${\rm S}_W({\rm Sp}_{2g}(\Z))$.  \ps\ps

\subsection{Relation entre ${\rm S}_W({\rm Sp}_{2g}(\Z))$ et
$\mathcal{A}^2(\PGSp_{2g})$}\label{fautsiegel}

Nous allons maintenant rappeler le lien classique entre ${\rm S}_W({\rm
Sp}_{2g}(\Z))$ et l'espace $\mathcal{A}_{\rm cusp}(\PGSp_{2g})$.  Une
r\'ef\'erence r\'ecente agr\'eable sur ce point est
l'article~\cite{asgarischmidt} auquel nous renverrons d\`es qu'il sera
possible de formuler l'\'enonc\'e principal (voir aussi \cite[\S 5]{taibisiegel}).\ps\ps

On pose $G={\rm PGSp}_{2g}$. Le facteur de similitude $\nu : {\rm GSp}_{2g}
\rightarrow \mathbb{G}_m$ induit un homomorphisme $\nu_{\infty} : G(\R) \rightarrow
\R^\times/\R^\times_{>0}$ dont nous noterons $G(\R)^+$ le noyau.  Le
morphisme canonique ${\rm Sp}_{2g}(\R) \rightarrow G(\R)$ induit un
isomorphisme $${\rm Sp}_{2g}(\R)/\{\pm 1\} \isomo G(\R)^+.$$ 
On pose \'egalement $G(A)^+=G(A) \cap G(\R)^+$ si $A$ est un sous-anneau de $\R$.  \ps\ps

D'apr\`es le \S \ref{hsimilitudes} on a $\mathrm{h}(G)=1$. 
Comme $\nu_{{\infty}}(G(\Z))=\{\pm 1\}$ on en d\'eduit 
\begin{equation}\label{decompsymp} G(\AAA)=G(\Q)(G(\R)^+ \times G(\widehat{\Z})) \end{equation} 
et suivant~\eqref{fautreel} que la restriction $f \mapsto f_{|G^+(\R) \times 1}$ induit un isomorphisme $G(\R)^+$-\'equivariant 
\begin{equation}\label{approxfortesp} \mathcal{A}^2({\rm PGSp}_{2g}) \isomo
{\rm L}^2(G(\Z)^+\backslash G(\R)^+).\end{equation}
L'action de $\GSp_{2g}(\R)^{+}$ sur $\mathbb{H}_g$ rappel\'ee au~\S\ref{fsiegelclass} se factorise en une action de $G(\R)^+$. 
Cette derni\`ere est fid\`ele, transitive, et ses
stabilisateurs sont les sous-groupes compacts maximaux de $G(\R)^+$. Si $K$ d\'esigne le
stabilisateur dans ${\rm Sp}_{2g}(\R)$ de l'\'el\'ement 
$i {\rm 1}_g\in \mathbb{H}_g$, et si $K^+$ d\'esigne son image dans $G(\R)^+$,
on a donc une identification naturelle  $$G(\R)^+/K^+ \isomo
\mathbb{H}_g.$$

\noindent Soit $(\rho,W)$ une $\C$-repr\'esentation de $\GL_g$ comme au~\S\ref{fsiegelclass}, suppos\'ee maintenant irr\'eductible et telle que ${\rm m}_{W} \equiv 0 \bmod 2$. Fixons $w \in W^\ast$ et $f \in {\rm S}_W({\rm Sp}_{2g}(\Z))$, nous allons leur associer une fonction $\varphi_{w,f} \in \mathcal{A}^2(G)$. Observons la fonction $\varphi: G(\R)^+ \longrightarrow \C$
d\'efinie par $$\varphi(\gamma)=\langle w,  (f_{|_{W}} \gamma)( i {\rm
1}_g)\,\rangle.$$
Par construction, $\varphi$ est continue et invariante \`a gauche par $G(\Z)^+$, c'est donc la restriction \`a $G(\R)^+ \times 1$ d'une unique fonction continue $\varphi': G(\Q)\backslash G(\AAA) \rightarrow \C$, invariante par translations \`a droite par $G(\widehat{\Z})$, d'apr\`es la formule~\eqref{decompsymp}. On pose alors $$\varphi_{w,f}:=\varphi'.$$ 
D'apr\`es \cite[Lemma 5]{asgarischmidt}, on a $\varphi_{w,f} \in \mathcal{A}_{\rm cusp}(G)$. \ps\ps

Avant d'\'enoncer la proposition finale, il nous faut encore d\'efinir la notion d'\'el\'ement
holomorphe de $\mathcal{A}^2(G)$.  Soient $\mathfrak{g}$ et $\mathfrak{k}$
les alg\`ebres de Lie respectives de $G(\R)^+$ et $K$, et soit
$\mathfrak{g}=\mathfrak{k}\oplus \mathfrak{p}$ la d\'ecomposition de Cartan
associ\'ee.  Soit $d : \mathfrak{g} \rightarrow {\rm T}_{i1_g}$ la
diff\'erentielle en l'identit\'e de l'application $G(\R)^+ \rightarrow
\mathbb{H}_g$, $h \mapsto h(i1_g)$.  Elle induit un isomorphisme
$\R$-lin\'eaire $$\mathfrak{p} \isomo {\rm T}_{i1_g}={\rm Sym}_g(\C).$$ La
structure de $\C$-espace vectoriel de ${\rm Sym}_g(\C)$ munit donc
$\mathfrak{p}$ d'une structure de $\C$-espace vectoriel qui d\'ecompose
$\mathfrak{p}\otimes_\R \C$ en $\mathfrak{p}^+\oplus \mathfrak{p}^-$ de
sorte que $d$ induise un isomorphisme $\C$-lin\'eaire $\mathfrak{p}^+ \isomo
{\rm T}_{i1_g}$.  Un \'el\'ement $f \in \mathcal{A}^2(G)$ sera dit {\it
holomorphe} s'il est continu et si pour tout $g \in G(\AAA)$ la fonction
$G(\R) \rightarrow \C, h \mapsto f(gh)$, est infiniment diff\'erentiable et
annul\'ee par $\mathfrak{p}^-$.

\begin{prop}\label{caracsiegel} L'application $(w,f) \mapsto \varphi_{w,f}$
d\'efinit une injection $\C[K]$-lin\'eaire $$W^\ast \otimes {\rm
S}_W({\rm Sp}_{2g}(\Z)) \longrightarrow \mathcal{A}_{\rm cusp}(\PGSp_{2g})$$
dont l'image est l'ensemble des $f \in \mathcal{A}_{\rm cusp}(\PGSp_{2g})$
qui sont holomorphes et $W^\ast$-isotypiques sous l'action de $K$.  \end{prop}

Pr\'ecisons cet \'enonc\'e. L'application $h \mapsto {\rm j}(h,i 1_g)$ est
un morphisme de groupes $K \rightarrow \GL_g(\C)$, qui r\'ealise $\GL_g(\C)$
comme complexification du groupe unitaire compact $K$.  Cela permet
notamment de voir $W$ comme une repr\'esentation de $K$ par restriction, qui
est irr\'eductible car $W$ l'est comme repr\'esentation de $\GL_g$.  Nous
renvoyons \`a \cite[Thm.  1, \S 4.5]{asgarischmidt} pour une
d\'emonstration de cette proposition, 
\`a l'assertion de surjectivit\'e pr\`es, qui est v\'erifi\'ee dans \cite[\S 5.2]{taibisiegel}.  \ps\ps

\subsection{Action des op\'erateurs de Hecke}\label{heckesiegel}

Il r\'esulte de la proposition~\ref{caracsiegel} que l'image de l'application de cet \'enonc\'e est stable sous
l'action $\mathrm{H}({\rm PGSp}_{2g})$, de sorte que l'espace ${\rm S}_W({\rm Sp}_{2g}(\Z))$ h\'erite
de $\mathcal{A}^2({\rm PGSp}_{2g})$ une action de $\mathrm{H}({\rm
PGSp}_{2g})$. \`A des constantes de normalisation pr\`es \'eventuellement introduites suivant les auteurs pour 
des raisons d'int\'egralit\'e, cette action co\"incide avec l'action d\'efinie traditionnellement
sur ${\rm S}_W({\rm Sp}_{2g}(\Z))$, et m\^eme sur ${\rm M}_W({\rm
Sp}_{2g}(\Z))$, que nous rappelons ci-dessous (voir aussi~\cite[Kap. 
IV]{freitag}, \cite[\S 16]{vandergeer} et~\cite[\S 4.3]{asgarischmidt}). Sans rentrer dans les d\'etails, mentionnons qu'elle est particuli\`erement naturelle lorsque l'on voit ${\rm
Sp}_{2g}(\Z)\backslash \mathbb{H}_g$ comme l'espace des vari\'et\'es
ab\'eliennes complexes de dimension $g$ munies d'une polarisation principale\footnote{Une polarisation principale sur un r\'eseau $L \subset \C^g$ est la donn\'ee d'une forme bilin\'eaire altern\'ee non-d\'eg\'en\'er\'ee $\eta : L \times L \rightarrow \Z$ dont l'extension des scalaires $\eta_\R$ \`a $L \otimes \R = \C^g$ satisfait $\eta_\R(ix,iy)=\eta_\R(x,y)$ pour tout $x,y \in \C^g$, et telle que la forme hermitienne associ\'ee $(x,y) \mapsto \eta_\R(ix,y)+i\eta_\R(x,y)$ sur $\C^g$ est d\'efinie positive. La th\'eorie de Riemann permet d'identifier naturellement ${\rm Sp}_{2g}(\Z)\backslash \mathbb{H}_g$ \`a l'ensemble des $\GL_g(\C)$-orbites de couples $(L,\eta)$ o\`u $L \subset \C^g$ est un r\'eseau et $\eta$ est une polarisation principale sur $L$.}~\cite[\S 10]{vandergeer}. \ps\ps

Soient $(W,\rho)$ une $\C$-repr\'esentation irr\'eductible de $\GL_g$, $p$ un nombre premier, et $G$ le $\Z$-groupe ${\rm PGSp}_{2g}$. L'application naturelle $$a: G(\Z[\frac{1}{p}])^{+} / G(\Z)^{+} \rightarrow G(\Q_{p})/G(\Z_{p})$$ 
est bijective puisque ${\rm h}(G)=1$ (Corollaire~\ref{classnbpgg}) et $\nu_{\infty}(G(\Z))=\{\pm 1\}$ (\S \ref{fautsiegel}). Elle induit donc de mani\`ere \'evidente un homomorphisme injectif entre l'anneau ${\rm H}_p(G)$ et l'anneau de Hecke du $G(\Z[\frac{1}{p}])^{+}$-ensemble $G(\Z[\frac{1}{p}])^{+}/G(\Z)^{+}$. Cet homomorphisme est un isomorphisme ; en effet, cela suit de l'isomorphisme \eqref{repcanhecke} et de ce que $a$ induit \'egalement une bijection 
\begin{equation}\label{bijdbclasssp}G(\Z)^{+}\backslash G(\Z[\frac{1}{p}])^{+} / G(\Z)^{+} \rightarrow G(\Z_{p}) \backslash G(\Q_{p})/G(\Z_{p}),\end{equation} comme le montre la th\'eorie des diviseurs \'el\'ementaires (Propositions \ref{baseheckeo} et \ref{baseheckego}, voir \'egalement le \S VI.\ref{deuxbases}). \ps\ps

Supposons que l'\'el\'ement $T \in {\rm H}_{p}(G)$ ait pour matrice la fonction caract\'eristique de la double classe $G(\Z_{p})\gamma G(\Z_{p})$ avec $\gamma \in G(\Z[\frac{1}{p}])^{+}$, au sens de l'identification \eqref{bijmathecke}. Si l'on \'ecrit $$G(\Z)^{+} \,\gamma \,G(\Z)^{+} = \coprod_{i} \gamma_{i} \,G(\Z)^{+},$$  on constate imm\'ediatement \`a l'aide de la formule \eqref{actionconcrethecke} que le diagramme suivant est commutatif 
\begin{equation}\label{diaghecke} \xymatrix{ {\rm S}_{W}({\rm Sp}_{{2g}}) \ar@{->}[rr] \ar@{->}_{f \mapsto \sum_{i} f_{|_{W}}\gamma_{i}^{-1} }[d] & & {\rm Hom}_{K}(W,{\cal A}_{{\rm cusp}}({\rm PGSp}_{2g})) \ar@{->}^{T}[d] \\
{\rm S}_{W}({\rm Sp}_{{2g}}) \ar@{->}[rr] & &  {\rm Hom}_{K}(W,{\cal A}_{{\rm cusp}}({\rm PGSp}_{2g})) }\end{equation}
les applications verticales \'etant celles d\'efinies par la proposition~\ref{caracsiegel} (voir \cite[Lemma 9]{asgarischmidt} pour le d\'etail de l'argument). \'Etant donn\'e l'\'egalit\'e $T=T^{\rm t}$ pour tout $T \in {\rm H}(G)$, il ne sera d'ailleurs pas n\'ecessaire de se souvenir de l'inversion des $\gamma_{i}$ dans la formule ci-dessus.\ps\ps

La formule \eqref{diaghecke} permet de d\'eterminer le lien exact entre les op\'erateurs de Hecke consid\'er\'es ici et diverses d\'efinitions que l'on retrouve dans la litt\'erature. Nous nous contenterons dans ce qui suit de donner le dictionnaire avec les d\'efinitions de
Serre~\cite[Ch.  VII \S 2, \S 5]{serre} dans le cas $g=1$. Des cas particuliers en genre $g=2$ seront abord\'es au chapitre IX. \ps\ps

Soit $k\geq 0$ un entier pair. Serre d\'efinit dans~\cite[Ch. VII \S
5.3]{serre} pour tout entier $n\geq 1$ un endomorphisme de ${\rm M}_k({\rm
SL}_2(\Z))$ qu'il note ${\rm T}(n)$, et dont il d\'etermine l'effet sur les $q$-d\'eveloppements.  On dispose \'egalement d'un autre
endomorphisme donn\'e par l'op\'eration d\'efinie ci-dessus de l'op\'erateur
${\rm T}_A\in \mathrm{H}({\rm PGL}_2)$ d\'efini au~\S\ref{annheckeclass}, o\`u $A$ est un groupe cyclique. 
La traduction est alors :
\begin{equation}\label{tradserregl2} n^{-\frac{k-1}{2}}{\rm T}(n) =
n^{-\frac{1}{2}}\sum_{d^2|n}
{\rm T}_{\Z/\frac{n}{d^2}\Z}.\end{equation} Cela vient notamment du fait que la correspondance
${\rm T}(n)$ chez Serre associe \`a un r\'eseau l'ensemble de ses sous-groupes
d'indice $n$, plut\^ot que ceux de quotient $\Z/n\Z$.\ps\ps


\subsection{$\mathcal{A}_{\rm disc}({\rm Sp}_{2g})$ se d\'eduit de
$\mathcal{A}_{\rm disc}({\rm PGSp}_{2g})$}\label{spvspgsp}

Observons que le morphisme ${\rm Sp}_{2g}(\AAA) \rightarrow {\rm
PGSp}_{2g}(\AAA)$ induit, par restriction des fonctions, un isomorphisme
$${\rm Res} : \mathcal{A}^2({\rm PGSp}_{2g}) \isomo \mathcal{A}^2({\rm
Sp}_{2g}).$$ Cela se d\'eduit en effet de la formule~\eqref{fautreel},
sachant que $$\mathrm{h}({\rm Sp}_{2g})=\mathrm{h}({\rm PGSp}_{2g})=1$$ et que
l'homomorphisme naturel ${\rm Sp}_{2g}(\R) \rightarrow {\rm PGSp}_{2g}(\R)$
induit un hom\'eomorphisme ${\rm Sp}_{2g}(\Z) \backslash {\rm Sp}_{2g}(\R) \isomo
{\rm PGSp}_{2g}(\Z)\backslash {\rm PGSp}_{2g}(\R)$.  \ps\ps

On rappelle que nous avons d\'efini au~\S\ref{defperestroika} un morphisme d'anneaux injectif $\mathrm{H}({\rm Sp}_{2g}) \rightarrow \mathrm{H}({\rm PGSp}_{2g})$, que l'on verra d\'esormais par un l\'eger abus de langage comme une inclusion.
La source et le but du morphisme ${\rm Res}$ sont donc tous deux des $\mathrm{H}({\rm Sp}_{2g})$-modules.

\begin{prop} \label{compsppgsp} L'application ${\rm Res}$ commute aux actions de ${\rm Sp}_{2g}(\R)$ et de $\mathrm{H}({\rm Sp}_{2g})$. Elle envoie $\mathcal{A}^2_{\rm disc}({\rm PGSp}_{2g})$ sur 
$\mathcal{A}_{\rm disc}^2({\rm Sp}_{2g})$. 
\end{prop}

\begin{pf} La premi\`ere assertion est \'evidente et la seconde d\'ecoule du
lemme~\ref{compatibilitesatake}.  La derni\`ere est cons\'equence de la premi\`ere et de
ce que l'image de ${\rm Sp}_{2g}(\R)$ dans ${\rm PGSp}_{2g}(\R)$ est
d'indice fini (\'egal \`a $2$).  \end{pf}

\parindent=0.5cm

\chapter{S\'eries th\^eta et r\'eseaux unimodulaires pairs}\label{seriestheta}\label{Chap5}

\section{S\'eries th\^eta de Siegel}\label{thetasiegel}

Soient $L \subset \R^n$ un r\'eseau unimodulaire pair et  $g\geq 1$ un
entier.  Si $v$ est un $g$-uple d'\'el\'ements de $L$, c'est-\`a-dire $v=(v_{i}) \in L^{g}$, nous notons
$v.v:=(v_i \cdot v_j)_{{i,j}} \in {\mathrm{Sym}}_g(\Z)$ la matrice de
Gram associ\'ee; elle est positive et $\frac{v.v}{2}$
est demi-enti\`ere au sens du \S IV.\ref{devfouriersiegel}.  La {\it s\'erie
th\^eta de Siegel de genre $g$ de $L$} est la fonction holomorphe sur
$\mathbb{H}_g$ : $$\vartheta_g(L) = \sum_{v \in L^g} q^{\frac{v.v}{2}}.$$
Elle ne d\'epend que de la classe d'isom\'etrie de $L$.  Son d\'eveloppement
de Fourier s'\'ecrit $\sum_{n \geq 0} a_n q^n$ o\`u $a_n$ est le nombre de
$g$-uples d'\'el\'ements de $L$ de matrice de Gram $n$.  Quand $g=1$, on
retrouve bien entendu les s\'eries th\^eta classiques, trait\'ees
dans~\cite{serre}, $a_n$ \'etant simplement le nombre des $x \in L$ tels que
$x \cdot x = 2n$.  Siegel a d\'emontr\'e  $$\vartheta_g(L) \in {\rm
M}_{\frac{n}{2}}({\rm Sp}_{2g}(\Z))$$ \cite[Kap.  I \S 0]{freitag}.  Soit ${\rm X}_{n}$ l'ensemble des classes d'isom\'etrie de r\'eseaux unimodulaires pairs dans $\R^{n}$ d\'ej\`a introduit dans l'introduction. Il est
avantageux de lin\'eariser la construction ci-dessus en consid\'erant
l'application lin\'eaire $$\vartheta_g : \C[{\rm X}_n] \longrightarrow {\rm
M}_{\frac{n}{2}}({\rm Sp}_{2g}(\Z)), \, \, \, \, \, \, [L] \mapsto \vartheta_g(L).$$ De plus, on note
$\vartheta_0 : \C[{\rm X}_n] \rightarrow \C$ l'application lin\'eaire
envoyant $[L]$ sur $1$ pour tout $L$.  Rappelons que l'espace $\C[{\rm
X}_n]$ s'identifie canoniquement au dual de l'espace de formes modulaires alg\'ebriques ${\rm M}_\C({\rm O}_n)$ (\S IV.\ref{fauton}).  Il est donc muni d'une
action (\`a gauche) de l'anneau ${\rm H}({\rm O}_n)$, par transposition.  \ps\ps

Dans le but de comprendre le $\mathrm{H}({\rm O}_n)^{{\rm opp}}$-module ${\rm M}_\C({\rm O}_n)$, un
fait important est que l'application $\vartheta_g$ entrelace l'action de
$\mathrm{H}({\rm O}_n)$ sur $\C[{\rm X}_n]$ et celle de $\mathrm{H}({\rm Sp}_{2g}) \subset \mathrm{H}({\rm
PGSp}_{2g})$~(\S IV.\ref{fautsiegel}), et ce selon des recettes bien pr\'ecises. 
Ces relations ont \'et\'e d\'ecouvertes par Eichler dans des cas
particuliers quand $g=1$ \cite[Satz
21.3]{eichlerqf} \cite{rallis1} et sont appel\'ees depuis {\it relations
de commutation d'Eichler}. Elles entra\^inent en particulier que le noyau de $\vartheta_g$ est
stable par l'action de ${\rm H}({\rm O}_n)$.  Le cas du genre $g$ quelconque
a \'et\'e \'etudi\'e sous divers aspects par de nombreux auteurs depuis
Eichler, dont Rallis~\cite{rallis2}, Freitag~\cite[Ch.  IV \S 5]{freitag},
Andrianov~\cite[Ch.  V]{andrianov}.  Nous nous contenterons pour l'instant
d'\'enoncer le cas particulier suivant, explicit\'e par Walling
dans~\cite{walling}.  Nous renvoyons le lecteur au \S\ref{appendiceC} pour une d\'emonstration 
dans le cas du genre $g=1$, ce qui constitue un exercice sans grande
difficult\'e.  \ps\ps

On rappelle que l'op\'erateur de Hecke ${\rm T}_p \in \mathrm{H}({\rm O}_n)$ est l'op\'erateur
${\rm T}_{\Z/p\Z}$ des $p$-voisins~(\S IV.\ref{annheckeclass},\S IV.\ref{fautson}).  On d\'efinit
un op\'erateur de Hecke $${\rm S}_p \in \mathrm{H}({\rm Sp}_{2g})$$ en consid\'erant les
couples $(M,N) \in \mathcal{R}_\Z^{\rm a}(\Q^{2g})$ tels que soit $M\cap N$
est d'indice $p$ dans $M$ et $N$, soit $M=N$~(\S IV.\ref{annheckeclass}).
Autrement dit, c'est l'op\'erateur ${\rm T}_{\Z/p\Z}+1$ dans les notations {\it loc. cit.}

\begin{prop}\label{commeichlernaif} On suppose $1 \leq g \leq \frac{n}{2}$. Pour tout nombre premier $p$, le diagramme 

$$\xymatrix{ \C[{\rm X}_{n}] \ar@{->}^{\vartheta_{g}}[rrr] \ar@{->}_{{\rm T}_{p}}[d] &  & & {\rm M}_{\frac{n}{2}}({\rm Sp}_{2g}(\Z))   \ar@{->}^{\, \, p^{\frac{n}{2}-1-g} \, {\rm S}_p\,+ \, p^g\,\frac{p^{n-2g-1}-1}{p-1} }[d]\\ 
\C[{\rm X}_{n}] \ar@{->}^{\vartheta_{g}}[rrr] & & & {\rm M}_{\frac{n}{2}}({\rm Sp}_{2g}(\Z)) }$$
est commutatif.
En particulier, ${\rm Ker}\, \vartheta_g$ est stable par ${\rm T}_p$. 
\end{prop}

\begin{pf} C'est un cas particulier de~\cite[Thm. 2.1]{walling}, observant
que notre op\'erateur ${\rm S}_p$ co\"incide avec l'op\'erateur not\'e
${\rm T}_1(p^2)+1$ par Walling.
\end{pf}


Si l'on remplace ${\rm T}_p$ par un op\'erateur de Hecke ${\rm T}_A \in
\mathrm{H}({\rm O}_n)$
plus g\'en\'eral, les relations analogues donn\'ees par exemple par
Andrianov et Freitag prennent une forme d'apparence assez absconse. 
Elles sont toutefois particuli\`erement
limpides (notamment celle ci-dessus) dans la pr\'esentation de
Rallis, qui requiert n\'eanmoins les
points de vue de Satake et Langlands sur les op\'erateurs de Hecke. Nous
reviendrons sur cela au~\S VI.\ref{releichlerbisrallis}. \ps\ps

Observons que si $\Phi_g$ est l'op\'erateur de Siegel rappel\'e au~\S IV.\ref{devfouriersiegel}, on a la relation \'evidente $\Phi_g(\vartheta_g(L))=\vartheta_{g-1}(L)$ pour $g\geq 2$. Cette relation s'\'etend \`a $g=1$ si l'on pose $\Phi_1(\sum_{n\geq 0} a_n q^n)=a_0$. Il en d\'ecoule que la suite des sous-espaces $${\rm Ker}\, \vartheta_g \subset \C[{\rm X}_n],$$ 
pour $g\geq 0$, est d\'ecroissante. De plus, si $g\geq 1$ alors $\vartheta_g$ induit une injection $\C$-lin\'eaire
$$  {\rm Ker} \, \vartheta_{g-1} /{\rm Ker}\, \vartheta_{g} \longrightarrow {\rm S}_{\frac{n}{2}} ({\rm Sp}_{2g}(\Z)).$$

Il est imm\'ediat que l'on a ${\rm Ker} \,\vartheta_n=\{0\}$ et que ${\rm Ker}\,
\vartheta_0$ est de codimension $1$ dans $\C[{\rm X}_n]$. De plus, le vecteur 
$\sum_{L \in {\rm X}_n} \frac{1}{|{\rm O}(L)|}\, [L]$
est propre sous l'action de ${\rm H}({\rm O}_n)$ de valeurs propres
explicites (Proposition~III.2.4 et \S IV.\ref{prodherminv}), la droite qu'il engendre
\'etant un suppl\'ementaire de ${\rm Ker}\, \vartheta_0$. La question de
d\'eterminer toute la filtration ${\rm Ker} \, \vartheta_g$, ainsi que la
structure des $\mathrm{H}({\rm O}_n)$-modules ${\rm Ker}\, \vartheta_{g-1} /{\rm
Ker}\, \vartheta_g$ quand $g\geq 1$, est en revanche tout \`a fait non triviale d\`es que
$n>8$.  C'est \'evidemment un probl\`eme plus fin que celui de comprendre le
$\mathrm{H}({\rm O}_n)$-module ${\rm M}_\C({\rm O}_n)$, le but de ce m\'emoire... Quand
$n=16$ et $n=24$, cette filtration a \'et\'e \'etudi\'ee en d\'etail
par plusieurs auteurs, dont nous allons rappeler
les contributions ci-dessous. Cela conduira \`a une
d\'emonstration directe du th\'eor\`eme I.\ref{thmintro16} (cas $n=16$) ainsi qu'\`a point de d\'epart pour notre d\'emonstration du th\'eor\`eme I.\ref{thmintro24} (cas $n=24$). \ps\ps

\section{S\'eries th\^eta de ${\rm E_8}\oplus {\rm E}_8$ et ${\rm
E}_{16}$} \label{casdim16}

On rappelle ${\rm X}_{16}=\{{\rm E}_8\oplus {\rm E}_8, {\rm E}_{16}\}$
(Witt).  L'espace ${\rm M}_8({\rm SL}_2(\Z))$ \'etant de dimension $1$, on a l'identit\'e
bien connue
\begin{equation}\label{egaltheta1}\vartheta_1({\rm E}_8 \oplus {\rm
E}_8)=\vartheta_1({\rm E}_{16}).\end{equation} En particulier, l'\'el\'ement
$[{\rm E}_8\oplus {\rm E}_8]-[{\rm E}_{16}]$ engendre ${\rm Ker}\,
\vartheta_1={\rm Ker}\, \vartheta_0$, et c'est un vecteur propre des ${\rm T}_p$. 
Un fait absolument remarquable, conjectur\'e par Witt dans~\cite{witt}, est
que l'identit\'e~\eqref{egaltheta1} persiste jusqu'en genre $3$ :
\begin{equation} \vartheta_g({\rm E}_8 \oplus {\rm E}_8)=\vartheta_g({\rm
E}_{16}) \, \, \, {\rm si}\, \, \, \, g=1,2,3.\end{equation} Cela a \'et\'e
d\'emontr\'e par Witt {\it loc.  cit.} pour $g=2$, et beaucoup plus tard par
Igusa et Kneser, ind\'ependamment, pour $g=3$ \cite{igusa1}\cite{kneser}. 
Igusa d\'emontre ${\rm S}_8({\rm Sp}_{2g}(\Z))=0$ pour $g \leq 3$. Nous
renvoyons \`a l'appendice~\ref{appendicekneser} pour un exposition de la
d\'emonstration remarquable de Kneser, qui est tr\`es diff\'erente. Au final, ${\rm Ker } \, \vartheta_3 = {\rm Ker} \, \vartheta_0$, et si
$$F:=\vartheta_4({\rm E}_8 \oplus {\rm E}_8)- \vartheta_4({\rm E}_{16}),$$
alors $F \in {\rm S}_8({\rm Sp_8}(\Z))$. Il est bien connu que l'on a $F
\neq 0$ ; pr\'ecisons ce fait. \ps\ps

\begin{prop}\label{Fnonnulle} Soit ${\rm c}_Q$ le coefficient de Fourier de $F$
correspondant \`a une matrice de Gram d'une $\Z$-base d'un r\'eseau pair $Q$ de rang
$4$ {\rm (}il ne d\'epend pas du choix d'une telle base{\rm )}. On a 
$$\frac{{\rm c}_{{\rm D}_4}}{|{\rm O}({\rm D}_4)|} = 4480 \hspace{1 cm}{\rm et}
\hspace{1 cm} \frac{{\rm c}_{{\rm A}_4}}{|{\rm O}({\rm A}_4)|} = -21504.$$
En particulier, ${\rm c}_{{\rm D}_4}=-{\rm c}_{{\rm A}_4}$. 
\end{prop}

\begin{pf} En effet, un examen des syst\`emes de racines
${\bf D}_n$ et ${\bf A}_m$ montre que le nombre de sous-r\'eseaux de ${\rm
D}_n$ isom\'etriques \`a ${\rm D}_4$ (resp.  ${\rm A}_m$) est
${{n}\choose{4}}$ (resp.  $2^m{{n}\choose{m+1}}$).  D'autre part, il se trouve
que si $R={\bf D}_4$ ou ${\bf A}_4$, il existe exactement une orbite de
sous-r\'eseaux de ${\rm E}_8$ isom\'etriques \`a ${\rm Q}(R)$ sous l'action
de ${\rm O}({\rm E}_8)$, son cardinal \'etant $\frac{|{\rm O}({\rm E}_8)|}{|{\rm W}(R)| \cdot |{\rm A}(R)|}$. Indiquons bri\`evement comment justifier
cette derni\`ere affirmation, dont la proposition d\'ecoule alors par un simple
calcul. On traite simultan\'ement les deux cas $R={\bf D}_4$ ou ${\bf A}_4$, et on
pose $Q={\rm Q}(R)$, $E=\,${\rm r\'es}\,$Q$ et $\Gamma={\rm A}(R)/{\rm W}(R)$ (la notation {\rm r\'es} a \'et\'e d\'efinie au \S II.1). \ps\ps

On commence par observer qu'un r\'eseau
euclidien pair $L \subset \R^4$ tel que {\rm r\'es}~$L \simeq E$ est n\'ecessairement isomorphe \`a 
$Q$. En effet, nous laissons au lecteur le soin de v\'erifier que cela se d\'eduit du fait bien connu suivant : tout r\'eseau entier
et unimodulaire de $\R^d$ pour $d \in \{4,5\}$ est isom\'etrique au r\'eseau
carr\'e ${\rm I}_d$ (\S II.2). La proposition 2.2 de l'appendice B montre que si $L \subset {\rm E}_8$ est isom\'etrique \`a $Q$, son
orthogonal $L^\perp$, qui est pair de rang $4$, admet un r\'esidu isomorphe
\`a $-E$. Mais dans les deux cas $E \simeq -E$, donc $L^\perp$ est isom\'etrique \`a $Q$. 
On consid\`ere enfin les lagrangiens du r\'esidu $E \oplus E$ de $Q \oplus
Q$; par la proposition II.1.1, ils donnent chacun naissance \`a un ${\rm q}$-module sur $\Z$ de rang $8$
contenant $Q \oplus Q$, n\'ecessairement isomorphe \`a ${\rm E}_8$.
Le ${\rm qe}$-module $E$ \'etant anisotrope, les lagrangiens de $E \oplus -E$ 
sont les graphes des automorphismes de $E$, et on constate dans les deux
cas qu'ils sont permut\'es simplement transitivement par $1 \times
\Gamma$, car l'homomorphisme naturel $\Gamma \rightarrow {\rm Aut}(E)$
est bijectif. Il suit de ces observations que l'ensemble des sous-r\'eseaux
de ${\rm E}_8$ isom\'etriques \`a $Q$ forment une seule orbite sous l'action
de ${\rm O}({\rm E_8})$, de stabilisateur isomorphe \`a ${\rm A}(R) \times {\rm W}(R)$.\end{pf}

\ps\ps
Ainsi, ${\rm Ker}\, \vartheta_4 =0$. Cela termine la description de la filtration de $\C[{\rm
X}_{16}]$ et nous ram\`ene \`a comprendre l'action de ${\rm H}({\rm Sp}_8)$
sur la forme propre $F \in {\rm S}_8({\rm Sp_8}(\Z))$. Cette forme est particuli\`erement int\'eressante.  En effet, Igusa a d\'emontr\'e
dans~\cite{igusa2} qu'elle est proportionnelle \`a la fameuse forme de
Schottky. Poor et Yuen en ont obtenu dans~\cite{py} une seconde
d\'emonstration, en v\'erifiant
\begin{equation} \label{dims8sp8} \dim {\rm S}_8({\rm Sp_8}(\Z)) = 1.\end{equation} Nous renvoyons \`a \cite{dukei1} pour une seconde d\'emonstration de cette \'egalit\'e, ainsi qu'au Th\'eor\`eme IX.\ref{thmsiegelsmallweight} pour une
troisi\`eme !\ps\ps

Soit $\tau(n)$ la fonction de Ramanjuan, d\'efinie par $\Delta=q\prod_{n\geq
1}(1-q^n)^{24}=\sum_{n\geq 1} \tau(n)q^n$.  Un calcul \'el\'ementaire montre
que le th\'eor\`eme~\ref{thmintro16} de l'introduction est
cons\'equence imm\'ediate du (i) du th\'eor\`eme suivant, l'apparition des
termes $286$ et $405$ provenant de la relation $\frac{|{\rm O}({\rm
E}_{16})|}{{|{\rm O}({\rm E}_8 \oplus {\rm E}_8)|}}=\frac{286}{405}$. Rappelons que l'op\'erateur de Hecke ${\rm S}_{p} \in {\rm H}_{p}({\rm Sp}_{{2g}})$ a \'et\'e introduit au \S \ref{thetasiegel}.

\begin{thm}\label{vpnontrivdim16} Soit $p$ un nombre premier. \ps \ps
\begin{itemize} \item[(i)] La valeur propre de ${\rm T}_p$ sur $[{\rm E}_8 \oplus
{\rm E}_8]-[{\rm E}_{16}]$ est $$p^4 \, \frac{p^7-1}{p-1} + p^7 + \tau(p)
\,\frac{p^4-1}{p-1}.$$\ps\ps

 \item[(ii)] La valeur propre de ${\rm S}_p$ sur la droite
${\rm S}_8({\rm Sp}_8(\Z))$ est $$p^4 + \tau(p) \, p^{-3} \,
\frac{p^{4}-1}{p-1}.$$ \end{itemize} \end{thm}

\noindent Nous avons vu ci-dessus que $\vartheta_4$ induit un isomorphisme
${\rm Ker} \, \vartheta_0 \isomo {\rm S}_8({\rm Sp}_8(\Z))$.  Les assertions (i)
et (ii) sont alors \'equivalentes d'apr\`es la
proposition~\ref{commeichlernaif} (relations de commutation d'Eichler). 
Il se trouve que l'assertion (ii) est une application imm\'ediate des travaux
d'Ikeda~\cite{ikeda1} (preuve de la conjecture de Duke-Imamo\u{g}lu~\cite{bkuss}).  En
effet, si $k$ et $g$ sont des entiers pairs tels que $k \equiv g \bmod 4$,
Ikeda construit {\it loc.  cit.} une application lin\'eaire injective $${\rm
I}_g : {\rm S}_k({\rm SL}_2(\Z)) \longrightarrow {\rm S}_{\frac{k+g}{2}}({\rm
Sp}_{2g}(\Z))$$ ayant la propri\'et\'e de compatibilit\'e suivante aux
op\'erateurs de Hecke.  Si $f=q+\sum_{n\geq 2}a_n q^n \in {\rm S}_k({\rm
SL}_2(\Z))$ est propre pour $\mathrm{H}({\rm PGL}_2)$, alors pour tout nombre premier
$p$, ${\rm I}_g(f)$ est propre pour ${\rm S}_p$ de valeur propre
\begin{equation}\label{ikedalift} p^g(1 + a_p \,\,p^{-\frac{k+g}{2}+1}
\,\frac{p^{g}-1}{p-1}). \end{equation}
\noindent (Nous renvoyons le n\'eophyte au~\S \ref{retourchap4} pour une explication du passage de l'\'enonc\'e d'Ikeda \`a
celui-l\`a). Dans la litt\'erature, la forme ${\rm I}_g(f)$ est souvent appel\'ee {\it rel\`evement d'Ikeda} en
genre $g$ de $f$.  Si $k=12$, $g=4$, et si $f$ est la forme modulaire $\Delta \in {\rm S}_{12}({\rm
SL}_2(\Z))$, on constate que son rel\`evement d'Ikeda ${\rm I}_4(\Delta)$ est un
\'el\'ement (non nul) de ${\rm S}_8({\rm Sp}_8(\Z))$, ayant la valeur propre de
l'\'enonc\'e (ii), ce qui conclut.  $\square$

\ps\ps

Mentionnons qu'il avait d\'ej\`a \'et\'e observ\'e par Breulmann et
Kuss~\cite[\S 3]{bkuss} que lorsque $p=2$, la valeur propre de l'op\'erateur
de Hecke ${\rm S}_p$ sur la forme de Schottky est bien donn\'ee par la formule du
(ii), conform\'ement \`a la conjecture de Duke-Imamo\u{g}lu. Leur m\'ethode
consiste \`a r\'ealiser $F$ comme s\'erie th\^eta \`a coefficients
harmoniques construite \`a partir de ${\rm E}_8$ (ce que nous ferons d'ailleurs \'egalement un peu plus loin!). Observons qu'une v\'erification
similaire, qui para\^it plus \'econome, est fournie par le calcul de ${\rm T}_2$
effectu\'e au~\S III.3.1. \ps\ps

L'apparition de la forme $\Delta$ dans l'argument
ci-dessus est tr\`es indirecte, cons\'equence du r\'esultat profond
d'Ikeda. Nous allons donner au \S\ref{i4dpartrialite} une autre d\'emonstration du (ii)
que nous avons d\'ecouverte et qui est ind\'ependante du travail d'Ikeda.
Elle repose notamment sur la trialit\'e pour
le groupe ${\rm PGSO}_{{\rm E}_8}$ sur $\Z$. Nous obtiendrons \'egalement plus tard une troisi\`eme
d\'emonstration de l'assertion (ii) du th\'eor\`eme, nettement plus sophistiqu\'ee, reposant sur la th\'eorie d'Arthur : c'est le cas particulier (particuli\`erement simple!) o\`u $k=8$ dans l'\'enonc\'e du th\'eor\`eme IX.\ref{thmsiegelsmallweight}. Enfin, nous expliquerons aussi comment l'assertion (i) du th\'eor\`eme
se d\'eduit trivialement d'une conjecture tr\`es g\'en\'erale concernant la th\'eorie d'Arthur, qui sera expos\'ee au chapitre VIII : voir les exemples du \S VIII.\ref{formexpliciteson0}.

\section{S\'eries th\^eta des r\'eseaux de Niemeier}\label{thetaniemeier}

Consid\'erons maintenant le cas $n=24$. D'apr\`es
Erokhin \cite{erokhin}, les s\'eries th\^eta de genre $12$ des $24$ r\'eseaux de
Niemeier sont lin\'eairement ind\'ependantes, i.e. ${\rm Ker}\, \vartheta_{12}=0$.
Cela ne vaut pas en genre $11$. En effet, ainsi que l'ont observ\'e Borcherds, Freitag et
Weissauer~\cite{bfw} par une construction ing\'enieuse, ${\rm Ker}\, \vartheta_{11}$ est de
dimension $1$. C'est un analogue spectaculaire en dimension $24$ de la d\'ecouverte de Witt
\'etudi\'ee au paragraphe pr\'ec\'edent. \ps\ps

Une \'etude plus fine du $\mathrm{H}({\rm O}_{24})$-module filtr\'e $\C[{\rm
X}_{24}]$ a \'et\'e initi\'ee par Nebe et Venkov dans leur d\'elicieux
article~\cite{nebevenkov}.  Leur point de d\'epart est le calcul de ${\rm
T}_2$ rappel\'e au~\S III.3.3, qu'ils d\'eduisent de la th\`ese de
Borcherds~\cite{borcherdsthese}.  Ils observent que ${\rm T}_2$ poss\`ede
$24$ valeurs propres enti\`eres distinctes, et disposent pour chacune
d'entre elle d'un vecteur propre explicite, qui est n\'ecessairement propre pour
l'action de tout $\mathrm{H}({\rm O}_n)$.  Il sera commode, suivant ces
auteurs, de num\'eroter ces vecteurs propres ${\rm v}_i$, $i=1,\dots,24$, de
sorte que les valeurs propres $\lambda_i$ de ${\rm T}_2$ associ\'ees soient
rang\'ees par ordre d\'ecroissant (voir la table~\ref{tablevp2}).  La
d\'etermination de la filtration de $\C[{\rm X}_{24}]$ \'equivaut alors \`a
celle du {\it degr\'e} de chaque ${\rm v}_i$, c'est-\`a-dire du plus petit entier
${\rm g}_i \geq 0$ tel que $\vartheta_{{\rm g}_i}({\rm v}_i) \neq 0$.  Nebe et Venkov y
parviennent pour $22$ des $24$ valeurs propres, et proposent une conjecture
pour les deux cas restants, que nous avons soulign\'es dans le tableau
ci-dessous.

\begin{table}[h]
\centering
{\footnotesize
\begin{tabular}{|c||c|c|c|c|c|c|c|c|}
\hline $\lambda$ & 8390655 & 4192830 & 2098332 & 1049832 & 533160 & 519120 & 268560 & 244800 \\
\hline degr\'e       & 0 & 1 & 2 & 3 & 4 & 4 & 5 & 5 \\
\hline $\lambda$ & 145152 & 126000 & 99792 & 91152 & 89640 & 69552 & 51552 & 45792 \\
\hline degr\'e & 6 & 6 & 6 & 7 & 8 & 7 &
8 & 7 \\
\hline $\lambda$ & 35640 & 21600 & 17280 & 5040 & -7920 & -16128 & -48528 & -98280 \\
\hline degr\'e & 8 & 8 & \underline{9} & 9 & \underline{10} & 10 & 11 & 12 \\
\hline
\end{tabular}
\caption{La filtration de $\C[{\rm X}_{24}]$, d'apr\`es Nebe et Venkov.}
\label{tablevp2}
}
\end{table}

Commentons un peu cette table. La valeur propre triviale
$$(2^{12}-1)(2^{11}+1)=8390655=\lambda_1$$ est bien entendu associ\'ee \`a l'unique vecteur propre
degr\'e $0$.  De plus, il est facile de v\'erifier que $\vartheta_1$ induit
une surjection ${\rm ker} \,\vartheta_0 \rightarrow {\rm S}_{12}({\rm
SL}_2(\Z))$, de sorte qu'un, et un seul, des vecteurs propres de ${\rm T}_2$ est
de degr\'e $1$. D'apr\`es le th\'eor\`eme \ref{appendiceCthm}, c'est celui de valeur propre
$$\tau(2)^2-2^{11}+2\,(2^{21}-1)=4192830=\lambda_2.$$ 
Les autres valeurs propres sont nettement plus subtiles \`a
comprendre.  Par exemple, Nebe et Venkov mentionnent \`a la fin de leur
article, que pour $i=3, 5, 11, 13$ et $24$ (auquel cas ${\rm g}_i$ vaut
respectivement $2, 4, 6, 8, 12$), alors $\vartheta_{{\rm g}_i}({\rm v}_i)$ est
proportionnel au rel\`evement d'Ikeda $${\rm I}_{{\rm g}_i}(\Delta_{23-{\rm g}_i})$$
o\`u $\Delta_w$ d\'esigne un g\'en\'erateur de ${\rm S}_{w+1}({\rm SL}_2(\Z))$
quand $w \in \{11, 15, 17, 19, 21\}$ (voir le corollaire \S VI.\ref{cor24ikedabocherer} pour une justification de ce fait). L'action de $\mathrm{H}({\rm O}_{24})$ sur ces ${\rm v}_i$
est donc explicitement connue (modulo les
relations d'Eichler) en fonction des coefficients de Fourier des $5$ formes
modulaires $\Delta_w$ ci-dessus. Quand $i=24$, cela avait d\'ej\`a
\'et\'e remarqu\'e par Borcherds, Freitag et Weissauer, imm\'ediatement
apr\`es l'annonce de~\cite{ikeda1} (voir
\'egalement~\cite{bkuss}). Par exemple,
$$\tau(2)(2^{12}-1)=-98280=\lambda_{24} ;$$ plus
g\'en\'eralement la valeur propre de ${\rm T}_p$ sur ${\rm v}_{24}$ est
$\tau(p)\frac{p^{12}-1}{p-1}$. \ps\ps

Un pas suppl\'ementaire spectaculaire a alors \'et\'e obtenu \`a nouveau par
Ikeda dans son article~\cite{ikeda2}, comme cons\'equence des r\'esultats
sus-cit\'es de Nebe-Venkov, et de sa r\'esolution partielle d'une conjecture
de Miyawaki. Ikeda parvient \`a exprimer l'action de $\mathrm{H}({\rm O}_{24})$ sur
tous les ${\rm v}_i$, toujours en terme des $\Delta_w$ ci-dessus, except\'e pour
quatre d'entre elles, \`a savoir les suivantes :

\begin{table}[h]
\centering
\begin{tabular}{|c||c|c|c|c|}
\hline $\lambda$ & 126000 & 51552 & 17280 &  -7920  \\
\hline degr\'e       &  6 & 8 & \underline{9} & \underline{10} \\
\hline
\end{tabular}
\caption{Les $4$ vecteurs propres myst\'erieux}
\end{table}

Nous reviendrons sur l'\'enonc\'e exact d\'emontr\'e par Ikeda au \S \ref{so24etnv}, dans lequel
nous expliquerons \'egalement les valeurs popres manquantes ci-dessus. \ps\ps

Terminons ce paragraphe par une discussion de ${\rm M}_{\DET}({\rm
O}_{24})$. Nous avons d\'ej\`a observ\'e au paragraphe IV.\ref{sovso} que
l'application naturelle $\widetilde{{\rm X}}_n \rightarrow {\rm X}_n$ est
bijective pour $n<24$, mais que pour $n=24$ la classe du r\'eseau de Leech
(et seulement celle-ci) poss\`ede deux ant\'ec\'edents, que nous noterons
${\rm Leech}^{\pm}$. Il suit que ${\rm M}_{\DET}({\rm O}_{24})$ est de
dimension $1$, constitu\'e des fonctions nulles sur les $23$ r\'eseaux de
Niemeier avec racines et prenant des valeurs oppos\'ees sur ${\rm Leech}^+$ et
${\rm Leech}^-$. 

\begin{prop}\label{leechmoins} La valeur propre de ${\rm T}_p$ sur la droite ${\rm
M}_{\DET}({\rm O}_{24})$ est $$\tau(p)\, \frac{p^{12}-1}{p-1}.$$ \end{prop}

Nous d\'emontrerons ce r\'esultat au~\S \ref{demoleechmoins}, et plus pr\'ecis\'ement
que l'anneau ${\rm H}({\rm O}_{24})$ agit de la m\^eme mani\`ere sur les
droites $\C {\rm v}_{24} \subset {\rm M}_\C({\rm O}_{24})$ et ${\rm M}_{\DET}({\rm O}_{24})$. Cela r\'epondra notamment \`a une question
de Schulze-Pillot~\cite[Remark, \S 1]{pillot}. \ps\ps


Cette proposition admet une traduction frappante en terme de r\'eseaux. Disons qu'un $p$-voisin $M$ d'un r\'eseau unimodulaire pair $L$ de $\R^n$ est {\it propre} s'il existe $g \in {\rm SO}(\R^n)$ tel que $g(M)=L$. On note ${\rm N}_p^+(L,M)$ le nombre des $p$-voisins propres de $L$ isom\'etriques \`a $M$ et on pose ${\rm N}_p^-(L,M)={\rm N}_p(L,M)-{\rm N}_p^+(L,M)$.  Si $L$ poss\`ede une isom\'etrie de d\'eterminant $-1$, tous ses $p$-voisins sont propres. C'est bien entendu le cas de tous les r\'eseaux de Niemeier avec racines. Le cas du r\'eseau de Leech est en revanche plus int\'eressant, comme le montre le corollaire suivant imm\'ediat de la proposition \ref{leechmoins}.

\begin{cor}\label{voispropreleech}  Pour tout nombre premier $p$, on a la relation $${\rm N}_p^+({\rm Leech},{\rm Leech})-{\rm N}^-_p({\rm Leech},{\rm Leech})=\tau(p) \frac{p^{12}-1}{p-1}.$$
\end{cor}

Une cons\'equence amusante de ce corollaire est que la fameuse conjecture de Lehmer est \'equivalente \`a ``${\rm N}_p^+({\rm Leech},{\rm Leech}) \neq {\rm N}^-_p({\rm Leech},{\rm Leech})$ pour tout premier $p$'' !

\section{Une construction alternative de ${\rm I}_4(\Delta)$ par trialit\'e} \label{i4dpartrialite}

Comme promis au~\S\ref{casdim16}, nous allons donner ci-dessous une seconde
d\'emonstration du th\'eor\`eme~\ref{vpnontrivdim16} (et donc du
th\'eor\`eme I.\ref{thmintro16}), ind\'ependante du th\'eor\`eme d'Ikeda~\cite{ikeda1}. 
Cela nous permettra au passage de donner des exemples non triviaux de formes
automorphes pour ${\rm O}_8$ et d'illustrer les techniques du chapitre IV. 
\ps\ps

\subsection{S\'eries th\^eta harmoniques}\label{sertharm}

On se place dans l'espace euclidien $V=\R^n$. Soit $1 \leq g \leq \frac{n}{2}$ un
entier. L'espace vectoriel $V^g=V \otimes \R^g$ est muni d'une repr\'esentation $\R$-lin\'eaire naturelle de
${\rm O}(V) \times \GL_g(\R)$. Pour tout entier $d\geq 0$, consid\'erons l'espace 
${\rm H}_{d,g}(V)$ des polyn\^omes $P : V^g \rightarrow \C$ tels que : \ps \ps

(i) $P$ est harmonique relativement au laplacien euclidien de $V^g$, \ps\ps

(ii) $P \circ h = \DET(h)^d P$ pour tout $h \in \GL_g(\R)$. \ps \ps Cet
espace est stable sous l'action de ${\rm O}(V)$.  On en construit des
\'el\'ements de la mani\`ere suivante.  Soit $I \subset V \otimes \C$ un
sous-espace isotrope de dimension $g$, et soit $e_1,\dots,e_g$ une $\C$-base de
$I$.  Il n'est pas difficile de v\'erifier que
\begin{equation}\label{vecteurextremal} (v_1,\dots,v_g) \mapsto \DET [ e_i \cdot v_j]_{1\leq i,j \leq g}^d\end{equation} est
un \'el\'ement de ${\rm H}_{d,g}(V)$. Il engendre une droite ne d\'ependant
que de $I$.  Ces droites sont permut\'ees transitivement par ${\rm O}(V)$ et
engendrent ${\rm H}_{d,g}(V)$, qui est une repr\'esentation irr\'eductible
de ${\rm O}(V)$ d'apr\`es~\cite[(0.9),(5.7)]{kw}.\ps\ps

Si $L \subset V$ est un r\'eseau unimodulaire pair, et si $P \in {\rm H}_{d,g}(V)$, on pose
$$\vartheta_g(L,P)= \sum_{v \in L^g} P(v) q^{\frac{v.v}{2}}.$$ 
L'\'equation fonctionnelle de la fonction $\vartheta$ de Jacobi permet de
d\'emontrer  $\vartheta_g(L,P) \in {\rm M}_{\frac{n}{2}+d}({\rm Sp}_{2g}(\Z))$ 
\cite[Kap. III \S 3]{freitag}. 
Observons que 
\begin{equation}\label{invarianceotheta}
\vartheta_g(L,P)=\vartheta_g(\gamma(L),\gamma(P)), \hspace{1 cm} \forall \gamma \in {\rm O}(V).\end{equation}
En particulier, $\vartheta_g(L,P)=0$ pour tout $P$ si ${\rm H}_{d,g}(V)^{{\rm O}(L)}=0$.
Commen\c{c}ons par donner un exemple, sans doute bien connu, o\`u cet espace d'invariants n'est pas nul
pour $n=8$.

\begin{lemme}\label{invariantscox} Soient $R \subset V$ un syst\`eme de
racines, $W \subset {\rm O}(V)$ son groupe de Weyl et $W^+=W \cap {\rm
SO}(V)$.  Alors $$H_W(t):=\sum_{d\geq 0} \, (\dim {\rm H}_{d,1}(V)^W)\, t^d=
(1-t^2)\prod_i (1-t^{m_i+1})^{-1},$$ o\`u les $m_i$ sont les exposants de
$W$, et $$\sum_{d\geq 0} \,(\dim {\rm H}_{d,1}(V)^{W^+})\, t^d=
(1+t^{|R|/2)})H_W(t).$$ \end{lemme}

\begin{pf} Soient $A$ et $B$ les s\'eries g\'en\'eratrices respectives des
suites $\dim {\rm Pol}_d(V)^{W}$ et $\dim ({\rm Pol}_d(V) \otimes \DET)^{W}$
(``anti-invariants''), pour $d\geq 0$.  D'apr\`es \cite[Ch.  V \S 6]{bourbaki}, on a
$B=t^{|R|/2}A$ et $A=\prod_i (1-t^{m_i+1})^{-1}$.  Soit ${\rm Pol}_d(V)$
l'espace des polyn\^omes $V \rightarrow \C$ qui sont homog\`enes de degr\'e
$d$.  Notons ${\bf \Delta}=\sum_{i=1}^{n} \frac{\partial^2}{\partial x_i^2}$ le laplacien standard de $\R^{n}$. On dispose
pour tout $d \in \Z$ d'une suite exacte ${\rm O}(V)$-\'equivariante $$ 0
\longrightarrow {\rm H}_{d,1}(V) \longrightarrow {\rm Pol}_d(V)
\overset{{\bf \Delta}}{\longrightarrow} {\rm Pol}_{d-2}(V) \longrightarrow
0,$$ la surjectivit\'e de ${\bf \Delta}$ \'etant un r\'esultat classique (voir
par exemple~\cite[\S 5.2.3]{GW}).  On en d\'eduit (i) puis (ii).  \end{pf}

Consid\'erons par exemple le r\'eseau ${\rm E}_8 \subset \R^8$, dont le
syst\`eme de racines ${\rm R}({\rm E}_8)$ est pr\'ecis\'ement de type ${\bf
E}_8$.  Ses exposants sont les $8$ entiers $1 \leq m \leq 30$ premiers \`a
$30$.  Comme ${\rm W}({\bf E}_8)={\rm O}({\rm E}_8)$, on
trouve {\small \begin{equation}\label{invariants8} \sum_{d\geq 0} (\dim {\rm
H}_{d,1}(V)^{{\rm O}({\rm E}_8)}) \,t^d= 1\,+\, t^8\, + \,t^{12}\,
+\,t^{14}\,+\,t^{16}\,+\,t^{18}\,+\,2\, t^{20}\, +\, \cdots.\end{equation}}
\noindent Le plus petit invariant est donc pour $d=8$.

\begin{prop}\label{polynomeA} Le polyn\^ome 
$A(x)=- 30 \,( x \cdot x )^4 +  \sum_{\alpha \in {\rm R}({\rm E}_8)} (\alpha \cdot x)^8 $ 
est dans ${\rm H}_{8,1}(V)^{{\rm O}(E_8)}$. Il satisfait $A(\alpha) = 144$ pour toute racine $\alpha \in {\rm R}({\rm E}_8)$. En particulier, on a l'\'egalit\'e
$\vartheta_1({\rm E}_8,A) = 240 \cdot 144\,\, \Delta$.
\end{prop}

\begin{pf} L'invariance de $A$ par ${\rm O}({\rm E}_8)$ est \'evidente. V\'erifions que $A$ est harmonique. Avec les notations de la
d\'emonstration du lemme~\ref{invariantscox}, le polyn\^ome ${\bf \Delta} A
\in {\rm Pol}_6(V)$ est invariant.  Pour tout $d \geq 2$, on a la d\'ecomposition
$${\rm Pol}_d(V)=(x \cdot x) \, {\rm Pol}_{d-2}(V)\oplus {\rm H}_{d,1}(V).$$
La s\'erie de Poincar\'e~\eqref{invariants8} montre donc que ${\bf \Delta} A
$ est proportionnel \`a $ (x \cdot x)^3$. Il ne reste qu'\`a voir que
${\bf \Delta}A$ s'annule sur les racines de ${\rm E}_8$.  C'est une
v\'erification simple laiss\'ee au lecteur, utilisant d'une part les
formules $${\bf \Delta} (\alpha \cdot x)^k = k(k-1)(\alpha \cdot \alpha)
(\alpha \cdot x)^{k-2},\, \, \, \, \, {\bf \Delta} (x\cdot
x)^k=2k(\dim(V)+2k-2) (x\cdot x)^{k-1},$$ et d'autre part le fait que dans
${\rm E}_8$, il y a exactement $114$ racines non orthogonales \`a une racine
$\alpha_0$ donn\'ee d'apr\`es \cite[Chap.  VI, \S 1.11, Prop. 
32]{bourbaki}, leur produit scalaire avec $\alpha_0$ \'etant $\pm 1$ pour
$112$ d'entre elles.  Cette m\^eme propri\'et\'e entra\^ine aussi
$A(\alpha)=144$ pour tout $\alpha \in {\rm R}({\rm E}_{8})$, puis la derni\`ere assertion.  \end{pf}

Il sera utile d'exprimer l'identit\'e ci-dessus en terme de formes
automorphes pour ${\rm O}_8={\rm O}_{{\rm E}_8}$ (\S IV.\ref{fauton}). Nous allons pour cela lin\'eariser la d\'efinition des s\'eries th\'eta harmoniques, \`a la mani\`ere du \S \ref{thetasiegel}. Rappelons que pour tout entier $n \equiv 0 \bmod 8$, le ${\rm O}_{n}(\AAA_{f})$-ensemble ${\mathcal{R}}({\rm O}_{n})$ s'identifie naturellement \`a celui des r\'eseaux unimodulaires pairs de $\R^{n}$ inclus dans ${\rm E}_{n} \otimes \Q$ (\S IV.\ref{ensclassbil}). \ps\ps

Si $L \in \mathcal{R}({\rm O}_{n})$, l'application $P \mapsto \vartheta_{g}(L,P)$, ${\rm H}_{{d,g}}(V) \rightarrow {\rm M}_{{\frac{n}{2}+d}}({\rm O}_{n})$, est $\C$-lin\'eaire. De plus, $\vartheta_{g}(\gamma L, \gamma P)= \vartheta_{g}(L,P)$ pour tout $\gamma \in {\rm O}_{n}(\Q)$ (formule \eqref{invarianceotheta}). On dispose donc d'une unique application lin\'eaire $$(\Z[\mathcal{R}({\rm O}_{n})] \otimes {\rm H}_{{d,g}}(V))_{{\rm O}_{n}(\Q)} \rightarrow {\rm M}_{{\frac{n}{2}+d}}({\rm Sp}_{{2g}}(\Z))$$
envoyant la classe de $[L] \otimes P$ sur $\vartheta_{g}(L,P)$ pour tout $L \in \mathcal{R}({\rm O}_n)$ et tout $P \in {\rm H}_{d,g}(V)$. D'apr\`es l'isomorphisme IV.\eqref{coinvmu}, le ${\rm H}({\rm O}_{n})$-module de gauche s'identifie canoniquement \`a ${\rm M}_{{\rm H}_{{d,g}}(V)^{\ast}}({\rm O}_{n})^{\ast}$. Cette double dualit\'e n'\'etant pas tr\`es heureuse, il est agr\'eable de la supprimer en utilisant l'isomorphisme canonique entre les ${\rm H}({\rm O}_{n})$-modules ${\rm M}_{{\rm H}_{{d,g}}(V)^{\ast}}({\rm O}_{n})^{\ast}$ et ${\rm M}_{{\rm H}_{{d,g}}}(V)({\rm O}_{n})^{\rm t}$ (Corollaire IV.\ref{corhermitien}). Rappelons que le ${\rm t}$ en exposant dans ce dernier espace signifie que l'action de ${\rm H}({\rm O}_{n})$ est tordue par $T \mapsto T^{{\rm t}}$, au sens de la remarque IV.\ref{transposeHXmod}. Comme $T^{{\rm t}}=T$ pour tout $T \in {\rm H}({\rm O}_{n})$ (Proposition~IV.\ref{heckeocomm}), cette torsion est sans effet ici, et sera all\`egrement ignor\'ee par la suite. \ps\ps

Rappelons enfin que nous notons $[L,P]$ l'\'el\'ement de ${\rm M}_{{\rm H}_{{d,g}}(V)}({\rm O}_{n})$ image de $[L] \otimes P$ par l'isomorphisme canonique compos\'e $(\Z[\mathcal{R}({\rm O}_{n})] \otimes {\rm H}_{{d,g}}(V))_{{\rm O}_{n}(\Q)} \isomo {\rm M}_{{\rm H}_{{d,g}}(V)}({\rm O}_{n})$ (voir la fin du \S IV.\ref{prodherminv}).\ps\ps

\begin{defprop}\label{defpropserietheta} Il existe une unique application lin\'eaire $$\vartheta_{d,g}: {\rm M}_{{\rm
H}_{d,g}(V)}({\rm O}_n) \rightarrow {\rm M}_{\frac{n}{2}+d}({\rm Sp}_{2g}(\Z))$$
envoyant $[L,P]$ sur $\vartheta_g(L,P)$ pour tout $L \in \mathcal{R}({\rm O}_n)$ et tout $P \in {\rm H}_{d,g}(V)$.  Si $d>0$ alors ${\rm Im}(\vartheta_{d,g}) \subset {\rm S}_{\frac{n}{2}+d}({\rm Sp}_{2g}(\Z))$.
 \end{defprop}
 
 \begin{pf} L'existence et l'unicit\'e d'une telle application suivent de la discussion pr\'ec\'edente. Concr\^etement, si $F \in  {\rm M}_{{\rm
H}_{d,g}(V)}({\rm O}_n)$ on a $$\vartheta_{d,g}(F) = \sum_i \frac{1}{|{\rm O}(L_i)|}\, \vartheta_g(L_i,P_i)$$ o\`u les $L_i$ forment un syst\`eme de repr\'esentants des ${\rm O}_n(\Q)$-orbites de ${\mathcal R}({\rm O}_n)$, et o\`u l'on a pos\'e $P_i = F(L_i)$. \end{pf}
 
Quand $d=0$, on a un isomorphisme ${\rm H}_{d,g}(V) \simeq \C$ (repr\'esentation triviale), et l'application $\vartheta_{d,g}$ n'est autre que la compos\'ee de l'isomorphisme ${\rm M}_{\C}({\rm O}_{n}) \isomo \C[{\rm X}_{n}]$ donn\'e par le corollaire IV.\ref{corhermitien} et de l'application $\vartheta_g$ du~\S\ref{thetasiegel}. \ps\ps

\begin{cor} \label{dimm8o8} L'application $\vartheta_{8,1} : {\rm M}_{{\rm
H}_{8,1}(\R^8)}({\rm O}_8) \rightarrow {\rm S}_{{12}}({\rm SL}_{2}(\Z))$ est un isomorphisme entre espaces de dimension $1$. 
\end{cor}

\begin{pf} Il est bien connu que $\dim {\rm S}_{{12}}({\rm SL}_{2}(\Z))=1$. Rappelons, d'apr\`es Mordell, que l'on a ${\rm
X}_8=\{{\rm E}_8\}$, de sorte que pour tout entier $d\geq 0$ on a $\dim {\rm M}_{{\rm H}_{d,1}(\R^8)}({\rm O}_8) =
\dim {\rm H}_{d,1}(\R^8)^{{\rm O}({\rm E}_8)}$. Un examen de la s\'erie~\eqref{invariants8} montre que $\dim  {\rm M}_{{\rm
H}_{8,1}(\R^8)}({\rm O}_8) = 1$. Il suffit donc de voir que $\vartheta_{8,1}$ est non nulle, mais cela r\'esulte de la proposition \ref{polynomeA}.
\end{pf}

On dispose encore de relations de commutation
d'Eichler  pour l'application $\vartheta_{d,g} : {\rm M}_{{\rm H}_{d,g}(V)}({\rm O}_n) \rightarrow {\rm M}_{\frac{n}{2}+d}({\rm
Sp}_{2g}(\Z))$ (voir le \S VI.\ref{releichlerbisrallis} pour une justification). En particulier, pour tout nombre premier $p$ on a 
\begin{equation}\label{eichlerthetaharmgen} \vartheta_{d,g} \circ {\rm T}_p = (p^{\frac{n}{2}-1-g} \, {\rm S}_p\,+ \, p^g\,\frac{p^{n-2g-1}-1}{p-1} )\circ \vartheta_{d,g}. 
\end{equation}

\noindent Nous renvoyons au \S \ref{appendiceC} pour une d\'emonstration de cette formule pour $g=1$. Combin\'ee au corollaire \ref{dimm8o8}, on obtient le corollaire suivant. \ps\ps

\begin{cor}\label{cortpm888} La valeur propre de ${\rm T}_p$ sur la droite
${\rm M}_{{\rm H}_{8,1}(\R^8)}({\rm O}_8)$ est $$p^{-8} (\tau(p)^2-p^{11})+
\,p\, \frac{p^5-1}{p-1}.$$ \end{cor}

\subsection{L'op\'erateur de Hecke des
perestro\"ikas}\label{utilperestroika}

\noindent Rappelons que nous avons d\'efini dans l'exemple~IV.\ref{exampleiso} une injection
naturelle $${\rm H}(\mu) : \mathrm{H}({\rm O}_n) \rightarrow \mathrm{H}({\rm
PGO}_n)$$ 
associ\'ee au $\Z$-morphisme canonique $\mu : {\rm O}_n \rightarrow {\rm
PGO}_n$. Nous noterons d\'esormais ${\rm H}(\mu)$ comme une inclusion $\mathrm{H}({\rm O}_n) \subset \mathrm{H}({\rm
PGO}_n)$ pour ne pas
surcharger les notations.  Soit $W$ une repr\'esentation de ${\rm
PGO}_n(\R)$, et soit $W'$ sa restriction \`a ${\rm O}_n(\R)$.  D'apr\`es la
proposition~\S IV.\ref{compatibilitesatake2}, la restriction des fonctions
suivant $\mathcal{R}({\rm O}_n) \rightarrow \mathcal{R}({\rm PGO}_n)$
d\'efinit une application $\mathrm{H}({\rm O}_n)$-\'equivariante $$\mu^\ast
: {\rm M}_{W}({\rm PGO}_n) \rightarrow {\rm M}_{W'}({\rm O}_n).$$

\begin{lemme}\label{inflationmonGO} L'application $\mu^\ast$ est un isomorphisme.
\end{lemme}

\begin{pf} En effet, l'application ${\rm O}_n(\Q)
\backslash \mathcal{R}({\rm O}_n) \rightarrow {\rm PGO}_n(\Q) \backslash \mathcal{R}({\rm
PGO}_n)$ est bijective d'apr\`es la proposition~IV.\ref{diagrammesimilitude} (c'est
$\xi_2 \circ \xi_1$). On conclut car si $M \subset \R^n$ est un r\'eseau
euclidien, ou plus g\'en\'eralement un $\widetilde{{\rm b}}$-module sur $\Z$ d\'efini positif, alors ${\rm O}(M)={\rm
GO}(M)$. 
\end{pf}

Une mani\`ere de voir ce lemme est de dire que l'action de $\mathrm{H}({\rm O}_n)$ sur
${\rm M}_{W'}({\rm O}_n)$ s'augmente en une action de l'anneau plus gros $\mathrm{H}({\rm
PGO}_n)$. Nous allons appliquer ceci aux ${\rm H}_{d,1}(\R^n)$. Ces espaces sont munis d'une repr\'esentation naturelle de
${\rm GO}(\R^n)$, sur laquelle l'homoth\'etie de rapport $\lambda$ agit par
le scalaire $\lambda^{-d}$. En particulier, si $d$ est pair alors $${\rm H}_{d,1}(\R^n) \otimes
\nu^{d/2}$$ 
se factorise en une repr\'esentation de ${\rm PGO}_n(\R)$ dont la
restriction \`a ${\rm O}_n(\R)$ est simplement ${\rm H}_{d,1}(\R^n)$. On rappelle que l'op\'erateur ${\rm K}_p \in \mathrm{H}_p({\rm PGO}_n)$ des perestro\"ikas 
relativement \`a $p$ a \'et\'e d\'efini au~\S IV.\ref{defperestroika} (et ${\rm K}_p={\rm K}_p^{\rm t}$).

\begin{lemme}\label{commeichlergenus1} Soient $d \equiv 0 \bmod
2$, $W={\rm H}_{d,1}(\R^n) \otimes \nu^{d/2}$, $p$ un nombre premier, et $\ell_{2r}(p)=\prod_{i=0}^{r-1}(1+p^{i})$  {\rm (}{\rm i.e.} le nombre de lagrangiens de ${\rm
H}(\F_p^{r})$ {\rm )}. Le diagramme suivant est commutatif :
$$\xymatrix{ {\rm M}_{W}({\rm PGO}_{n}) \ar@{->}^{\vartheta_{d,1} \circ \mu^{\ast}}[rrr] \ar@{->}_{{\rm K}_{p}}[d] &  & & {\rm M}_{\frac{n}{2}+d}({\rm SL}_{2}(\Z))   \ar@{->}^{\ell_{n-2}(p)\,
p^{-d/2} \,{\rm T}(p) }[d]\\ 
{\rm M}_{W}({\rm PGO}_{n})  \ar@{->}^{\vartheta_{d,1} \circ \mu^{\ast}}[rrr] & & & {\rm M}_{\frac{n}{2}+d}({\rm SL}_{2}(\Z)) }$$
\end{lemme}

\begin{pf} C'est une variante harmonique du r\'esultat d'Eichler 
\cite[Satz 21.3]{eichlerqf} (voir aussi \cite{rallis1}). Nous redonnons
l'argument. \ps

\noindent On fixe un r\'eseau unimodulaire $M$ dans ${\rm E}_n \otimes \Q$ (ou, plus g\'en\'eralement, un r\'eseau homodual de ${\rm E}_n \otimes \Q$), ainsi qu'un polyn\^ome $P$ dans l'espace ${\rm H}_{d,1}(\R^n)\otimes \nu^{d/2}$.
Soit ${\rm Per}_p(M)$ l'ensemble des perestro\"ikas de $M$ relativement \`a $p$. \ps

\noindent Si $N \in {\rm Per}_p(M)$,
alors $\mu^\ast([\underline{N},P])=[\gamma(N),p^{-d/2}\,\gamma (P) ]$ o\`u $\gamma \in {\rm GO}(\Q)$ est
un \'el\'ement quelconque de facteur de similitude $p^{{-1}}$. Comme
le r\'eseau $p^{-1/2}N \subset \R^n$ est isom\'etrique \`a $\gamma(N)$, la
relation~\eqref{invarianceotheta} entra\^ine 
$$\vartheta_{d,1}\circ \mu^\ast([N,P])=\vartheta_{d,1}(\frac{N}{\sqrt{p}},P).$$
Si $m\geq 0$, le m-i\`eme coefficient de Fourier de $\vartheta_{d,1} \circ \mu^{\ast} \circ {\rm K}_{p} ([M,P])$ 
vaut donc
\begin{equation} \label{thetakpcoeffm} \sum_{(N,v)}P(\frac{v}{\sqrt{p}}) \end{equation}
la somme portant sur tous les couples $(N,v)$ avec $N \in {\rm Per}_p(M)$ et
$v \in N$ tels que $v \cdot v=2mp$. \ps

\noindent Supposons $m$ premier \`a $p$. Observons qu'un \'el\'ement $v \in M$ tel que $v \cdot v = 2mp$ est non nul modulo $pM$ et
isotrope. Il appartient donc \`a exactement $\ell_{n-2}(p)$
perestro\"ikas de $M$ (i.e. le nombre de lagrangiens de $M \otimes \F_p$
contenant une droite isotrope donn\'ee). La somme \eqref{thetakpcoeffm} vaut donc simplement
$$p^{-d/2} \ell_{n-2}(p) a_{mp},$$ o\`u $a_r$ d\'esigne le r-i\`eme coefficient de
Fourier de $\vartheta_{d,1}(M,P)$. Pour conclure, on pourrait traiter le cas $m$ multiple de $p$ de mani\`ere similaire, ou simplement invoquer le lemme \ref{appendiceClemme}.
\end{pf}

\begin{cor}\label{corperesdelta} Pour tout nombre premier $p$, la valeur propre de ${\rm K}_p$ sur la droite ${\rm M}_{{\rm H}_{8,1}(\R^8) \otimes
\nu^{4}}({\rm PGO}_8)$
est $2 \,p^{-4}\, \frac{p^4-1}{p-1} \tau(p)$. 
\end{cor}

\begin{pf}Cela suit en effet de l'identit\'e $\ell_6(p)=2(1+p)(1+p^2)=2\frac{p^4-1}{p-1}$, des lemmes \ref{inflationmonGO} et \ref{commeichlergenus1}, et du corollaire \ref{dimm8o8}.\end{pf}

\subsection{Passage de ${\rm PGO}_n$ \`a ${\rm PGSO}_n$}\label{pgpso}

L'application $\mu_3$ de la proposition~IV.\ref{rso=ro} \'etant bijective et ${\rm
PGSO}_n(\AAA_f)$-\'equi\-variante, les arguments donn\'es
aux~\S IV.\ref{hsovso} et~\S IV.\ref{fauton} concernant la
comparaison entre ${\rm SO}_n$ et ${\rm O}_n$ s'\'etendent {\it mutatis mutandis} 
au cas de ${\rm PGSO}_n$ et ${\rm PGO}_n$. On dispose notamment d'une action
du groupe $${\rm PGO}_n(\AAA_f)/{\rm PGSO}_n(\AAA_f) \simeq (\Z/2\Z)^{\rm P}$$
sur $\mathrm{H}({\rm PGSO}_n)$ qui pr\'eserve chaque $\mathrm{H}_p({\rm PGSO}_n)$, et 
dont les invariants sont exactement $\mathrm{H}({\rm PGO}_n)$. En
particulier, $\mathrm{H}({\rm PGO}_n) \subset \mathrm{H}({\rm PGSO}_n)$.
Si $U$ est une repr\'esentation de ${\rm PGSO}_n(\R)$, d'induite ${\rm Ind}\, U$ \`a ${\rm
PGO}_n(\R)$, on dispose d'un isomorphisme canonique $\mathrm{H}({\rm
PGO}_n)$-\'equivariant $${\rm ind} : {\rm M}_{U}({\rm PGSO}_n) \rightarrow
{\rm M}_{{\rm Ind} \,U}({\rm PGO}_n).$$ \indent Enfin, si $W$ est une
repr\'esentation de ${\rm PGO}_n(\R)$, et si $W'$ d\'esigne sa restriction
\`a ${\rm PGSO}_n(\R)$, l'action par conjugaison de ${\rm PGO}_{n}(\Q)/{\rm PGSO}_{n}(\Q) \simeq \Z/2\Z$ sur ${\rm M}_{W'}({\rm PGSO}_n)$ muni ce dernier espace d'une sym\'etrie naturelle que l'on notera $s$. La restriction des fonctions selon la bijection $\mathcal{R}({\rm PGSO}_{n}) \rightarrow \mathcal{R}({\rm PGO}_{n})$ d\'efinit une injection ${\rm
H}({\rm PGO}_n)$-\'equivariante $${\rm res} : {\rm M}_{W}({\rm PGO}_n) \longrightarrow {\rm M}_{W'}({\rm
PGSO}_n)$$
d'image ${\rm M}_{W'}({\rm PGSO}_n)^{s={\rm id}}$. \ps\ps

Soit $W$ la repr\'esentation ${\rm H}_{d,1}(\R^n)\otimes \nu^{d/2}$ de ${\rm
PGO}_n(\R)$ d\'efinie au~\S IV.\ref{utilperestroika} ($d \equiv 0 \bmod 2$). Si $n>2$, sa restriction $W'$ \`a
${\rm PGSO}_n(\R)$ est irr\'eductible \cite[\S 5.2]{GW}. Pour des raisons g\'en\'erales, on a enfin ${\rm Ind} \, W' \simeq W \oplus W \otimes \epsilon$, o\`u $\varepsilon$ est le caract\`ere d'ordre $2$ de ${\rm PGO}_n(\R)$ de
noyau ${\rm PGSO}_n(\R)$. 

\begin{lemme}\label{lemmedesc8} Soient $W={\rm H}_{8,1}(\R^8) \otimes \nu^{4}$ et $W'=W_{|{\rm PGSO}_8(\R)}$. L'application de restriction ${\rm res} : {\rm M}_{W}({\rm
PGO}_8) \rightarrow {\rm M}_{W'}({\rm PGSO}_8)$ est
bijective.  \end{lemme}
 
\begin{pf} Il s'agit de montrer que l'espace ${\rm M}_{W'}({\rm PGSO}_8)^{s=-{\rm
id}}$, qui est naturellement isomorphe \`a ${\rm M}_{W \otimes \varepsilon}({\rm PGO}_8)$,
est nul. D'apr\`es le lemme~IV.\ref{inflationmonGO}, il
est \'equivalent de montrer que ${\rm M}_{{\rm H}_{8,1}(\R^8) \otimes \varepsilon}({\rm O}_8)=0$.
Cela d\'ecoule de ce que ${\rm X}_8=\{{\rm E}_8\}$, ${\rm O}({\rm E}_8)=W({\bf
E}_8)$, et du lemme~\ref{invariantscox}, qui assure que les $W({\bf
E}_8)^+$-invariants de ${\rm H}_{8,1}(\R^8)$ sont en fait $W({\bf E}_8)$-invariants (et de dimension $1$).
\end{pf}

Nous avons vu au \S IV.\ref{hsovso} que l'\'el\'ement ${\rm T}_{(\Z/p\Z)^{{n/2}}}$ de ${\rm H}_{p}({\rm O}_{n})$ se d\'ecompose de mani\`ere naturelle comme somme de deux \'el\'ements ${\rm T}_{(\Z/p\Z)^{{n/2}}}^{\pm}$ de ${\rm H}_{p}({\rm SO}_{n})$. Un ph\'enom\`ene similaire se produit avec l'op\'erateur de Hecke ${\rm K}_p \in {\rm H}_{p}({\rm PGO}_n)$ des perestro\"ikas, qui se raffine lorsqu'on le voit dans ${\rm H}_{p}({\rm PGSO}_n)$.  \ps\ps

En effet, observons que
$M \in \mathcal{R}_\Z^{\rm a}({\rm E}_n \otimes \Q)$ \'etant donn\'e, l'ensemble $\Omega$ des $\underline{N} \in
\underline{\mathcal{R}}_\Z^{\rm h}({\rm E}_n \otimes \Q)$ tels que $N$ est une perestro\"ika de $M$ relativement \`a
$p$ est constitu\'e d'exactement deux orbites sous l'action de ${\rm
GSO}(M)$, et d'une seule sous ${\rm GO}(M)$. En effet, cela d\'ecoule de
la lissit\'e de ${\rm GO}_M$ et ${\rm GSO}_M$ sur $\Z$, et de ce que l'action de ${\rm GO}(M
\otimes \F_p)$ (resp. ${\rm GSO}(M \otimes \F_p)$) sur l'ensemble des lagrangiens de $M \otimes \F_p
\simeq {\rm H}(\F_p^{n/2})$ est transitive d'apr\`es le th\'eor\`eme de Witt (resp. admet deux orbites). Ces deux orbites donnent donc
naissance \`a deux op\'erateurs de Hecke ${\rm K}_p^\pm \in \mathrm{H}_{p}({\rm PGSO}_n)$ tels que $${\rm K}_p
= {\rm K}_p^+ + {\rm K}_p^-$$
et qui sont \'echang\'es sous l'action de ${\rm PGO}_n(\Q_p)/{\rm PGSO}_n(\Q_p) \simeq \Z/2\Z$. 

\begin{cor}\label{coractionKpm} Soit $W=({\rm H}_{8,1}(\R^8) \otimes
\nu^{4})_{|{\rm PGSO}_8(\R)}$. Pour tout nombre premier $p$, la valeur propre de ${\rm K}_p^{\pm }$ sur la droite
${\rm M}_W({\rm PGSO}_8)$ est $p^{-4}\frac{p^4-1}{p-1}\,\tau(p)$. 
\end{cor}

\begin{pf} Soit $s_0 \in {\rm PGO}({\rm E}_8)$ l'image d'une r\'eflexion par
rapport \`a une racine.  La conjugaison par cet \'el\'ement d\'efinit un
$\Z$-automorphisme de ${\rm PGSO}_8$, et induit en particulier un
isomorphisme ${\rm H}(s_0)$ de ${\rm H}({\rm PGSO}_8)$. D'autre part, la d\'emonstration du lemme \ref{lemmedesc8} r\'ev\`ele que la sym\'etrie $s$ de ${\rm M}_W({\rm PGSO}_8)$,
\'egalement induite par $s_0$, est l'identit\'e. Mais le 
lemme IV.\ref{compatibilitesatake2} assure que 
$$T \circ s_{0} = s_{0} \circ {\rm H}(s_0)(T) , \hspace{1
cm}\forall\,\, T \in  \mathrm{H}({\rm PGSO_8}).$$
Or on a ${\rm H}(s_0)({\rm K}_p^+)={\rm K}_p^{-}$ car l'image de $s_0$ dans ${\rm PGO}_8(\Q_p)/{\rm PGSO}_8(\Z_p)$ est non triviale pour tout nombre premier $p$. Il suit que ${\rm
K}_p^+$ et ${\rm K}_p^-$ ont m\^eme valeur propre sur la droite ${\rm M}_W({\rm
PGSO}_8)$, \`a savoir la moiti\'e de celle de ${\rm K}_p$. 
On conclut par le corollaire~\ref{corperesdelta} et le lemme~\ref{lemmedesc8}. \end{pf}

\subsection{Trialit\'e pour ${\rm PGSO}_8$}\label{trialitepgso8}

L'\'etape suivante repose sur la trialit\'e pour le $\Z$-groupe ${\rm PGSO}_8$. L'existence de la trialit\'e dans ce contexte est bri\`evement discut\'ee par Gross dans~\cite[\S 4]{grossinv}, qui consid\`ere plus pr\'ecis\'ement le cas du rev\^etement ${\rm Spin}_8$ de $G$ (on peut alors redescendre au groupe $G$ par quotient). 
\'Etant donn\'e son importance ici, il ne semble pas inutile de pr\'eciser un peu cette
construction. \ps\ps

Nous suivrons l'approche de~\cite{blijspringer2} dans le cas des corps, qui s'\'etend \`a tout anneau par les r\'esultats de~\cite{kps}. On rappelle que si $A$ est un anneau commutatif, une $A$-alg\`ebre d'octonions\footnote{Dans~\cite{kps}, les auteurs parlent aussi d'{\it alg\`ebre de Cayley}, ou encore d'{\it alg\`ebre de composition de rang $8$}, cf. p. 51 et 56 {\it loc. cit}. Soulignons que l'on ne r\'eclame pas l'associativit\'e de $\star$, qui n'est de fait jamais satisfaite.} $C$ est la donn\'ee d'un ${\rm q}$-module sur $A$ de rang $8$ muni d'une structure de $A$-alg\`ebre unitaire $(x,y) \mapsto x \star y$, tels que ${\rm q}(x \star y) = {\rm q}(x) {\rm q}(y)$ pour tout $x,y \in C$.
Le point de d\'epart est d'enrichir le ${\rm q}$-module ${\rm E}_8$ d'une structure de $\Z$-alg\`ebre d'octonions, de laquelle d\'ependra la construction d'une trialit\'e sur ${\rm PGSO}_8$. Ainsi que l'ont observ\'e van der Blij et Springer~\cite[(4.5)]{blijspringer1}, il existe une telle structure sur ${\rm E}_8$, et m\^eme une seule modulo ${\rm O}({\rm E}_8)$, \`a savoir l'anneau des 
octonions de Coxeter~\cite{coxeter}.  \ps\ps

Soit $C$ une $A$-alg\`ebre d'octonions. Consid\'erons la propri\'et\'e suivante portant sur $C$ et $\gamma \in {\rm GSO}(C)$, not\'ee $\mathcal{P}(C,\gamma)$ ({\it principe de trialit\'e de Cartan}):
\begin{equation}\label{proppp}  \exists \,\,\gamma',\gamma'' \in {\rm GSO}(C), \, \, \gamma'(x \star y) = \gamma(x) \star \gamma''(y), \,\,\, \forall x,y \in C. \end{equation}
La proposition~\cite[Prop. 4.5]{kps} affirme que : \begin{itemize} \ps\ps
\item[(i)] Si $\mathcal{P}(C,\gamma)$ est satisfaite, alors le couple $(\gamma',\gamma'')$ est unique modulo action diagonale de $A^\times$. \ps\ps
\item[(ii)] Si ${\rm Pic}(A)=0$ alors $\mathcal{P}(C,\gamma)$ est satisfaite pour tout $\gamma \in {\rm GSO}(C)$ (quand $A$ est un corps, c'est~\cite[Thm. 1]{blijspringer2}).\ps\ps
\item[(iii)] Si $\gamma \in {\rm GSO}(C)$, alors il existe une partition de l'unit\'e $1=\sum_i f_i$ dans $A$, telle que pour tout $i$, $\mathcal{P}(\gamma \otimes A_{f_i} ,C \otimes A_{f_i})$ est satisfaite.
\end{itemize}
\ps\ps
On rappelle que ${\rm PGSO}_C$ est par d\'efinition le quotient de ${\rm GSO}_C$ par son sous-$A$-groupe central $\mathbb{G}_m$ constitu\'e des homoth\'eties (\S II.1). Les propri\'et\'es (i) et (iii) justifient 
alors imm\'ediatement la d\'efinition suivante. 

\begin{defprop}\label{deftrialite} {\rm (Trialit\'e)} Soit $C$ une $A$-alg\`ebre d'octonions et soit $\pi : {\rm GSO}_C \rightarrow {\rm PGSO}_C$ le morphisme naturel.  Il existe un unique automorphisme $\tau$ du $A$-groupe ${\rm PGSO}_C$ ayant la propri\'et\'e suivante : pour
toute $A$-alg\`ebre commutative $B$, et tout $\gamma \in {\rm GSO}_C(B)$ tel
que $\mathcal{P}(C \otimes B,\gamma)$ est satisfaite, alors
$\tau(\pi(\gamma))=\pi(\gamma'')$. 
\end{defprop}

D'apr\`es~\cite[Prop.  4.6]{kps} et~\cite[\S 1
Cor.  2]{blijspringer2}, cet automorphisme $\tau$ satisfait $\tau^3=1$. On l'appelle {\it trialit\'e du $A$-groupe ${\rm PGSO}_{C}$} (qui ne d\'epend que du ${\rm q}$-module sur $A$ sous-jacent \`a $C$) {\it associ\'ee \`a la structure d'octonions $C$}. Il existe de nombreux points de vue sur la trialit\'e dans la
litt\'erature.  Une propri\'et\'e g\'eom\'etrique fascinante, d\'ecouverte par E.  Study et d\'evelopp\'ee par E. Cartan, est
la suivante.\ps\ps


\begin{lemme}\label{trialitestudy} Soit $C$ une alg\`ebre d'octonions sur le corps $k$ dont le ${\rm q}$-module sous-jacent est hyperbolique.  Soient $Q_1$, $Q_2$ et $Q_3$, les classes de conjugaison des sous-groupes de ${\rm PGSO}(C)$ stabilisant respectivement une droite isotrope de $C$ ou l'un des deux types de lagrangiens de $C$.
La trialit\'e de ${\rm PGSO}_{C}$ permute transitivement les trois classes $Q_i$ dans ${\rm PGSO}(C)$. \ps\ps
\end{lemme}

\begin{pf}  C'est~\cite[Thm. 8]{blijspringer2}. Pr\'ecis\'ement, si $a \mapsto \overline{a}:=-a+ (a \cdot 1) $ d\'esigne l'involution canonique de $C$, alors $\tau$ envoie 
le stabilisateur de la doite isotrope $k \overline{a} \subset C$ sur le stabilisateur du lagrangien $a \star C$, et ce dernier sur le stabilisateur du lagrangien $C \star a$ (qui est de type oppos\'e). \end{pf}

Soit $\tau \in {\rm Aut}( {\rm PGSO}_8)$ la trialit\'e associ\'ee \`a une structure fix\'ee d'octonions de Coxeter sur ${\rm E}_8$.  En tant qu'automorphisme du $\Z$-groupe ${\rm PGSO}_8$, elle agit naturellement sur
$\mathcal{R}({\rm PGSO}_8)$, les $\mathcal{R}_p({\rm PGSO}_8)$, et sur $\mathrm{H}({\rm PGSO}_8)$ en pr\'eservant les $\mathrm{H}_p({\rm PGSO}_8)$ (Lemme IV.\ref{compatibilitesatake2}). Les inclusions naturelles $\mathrm{H}({\rm O}_8) \subset \mathrm{H}({\rm PGO}_8) \subset \mathrm{H}({\rm PGSO}_8)$ permettent de voir l'op\'erateur de Hecke ${\rm T}_p$ comme un \'el\'ement de $\mathrm{H}({\rm PGSO}_8)$. 

\begin{cor}\label{lemmetrialitep} Pour tout nombre premier $p$, l'application $T \mapsto {\rm H}(\tau)(T)$ induit un $3$-cycle du sous-ensemble $\{{\rm T}_p, {\rm K}_p^+, {\rm K}_p^-\}  \subset \mathrm{H}_p({\rm PGSO}_8)$.  \end{cor} 

\begin{pf} Soient $G={\rm PGSO}_8$, $M={\rm E}_8$ et $V_0={\rm E}_8 \otimes \Q$.  L'\'el\'ement
$\underline{M} \in \underline{\mathcal{R}}_\Z^{\rm h}(V_0) = \mathcal{R}(G)$ a bien s\^ur
pour stabilisateur $G(\widehat{\Z})$ sous l'action de $G(\AAA_f)$, il est
donc pr\'eserv\'e par $\tau$.  Soient $\Omega_1,\Omega_2$ et $\Omega_3$ les
$G(\widehat{\Z})$-orbites dans $\underline{\mathcal{R}}_\Z^{\rm h}(V_0)$
constitu\'ees respectivement des classes d'isom\'etrie de $p$-voisins de
$M$, ou des deux types de perestro\"ikas de $M$ relativement \`a $p$.  Ces
orbites se factorisent manifestement en des $G(\F_p)$-orbites, les classes
de conjugaison de stabilisateurs de $G(\F_p)$-correspondantes \'etant
exactement les classes $Q_i$ du lemme pr\'ec\'edent appliqu\'e \`a $C={\rm
E}_8 \otimes \F_p$, ce qui conclut.  \end{pf}

Si $d\geq 0$, on a d\'ej\`a dit que la repr\'esentation ${\rm H}_{d,r}(\R^{2r})$ est irr\'eductible en
tant que repr\'esentation de ${\rm GO}(\R^{2r})$. Cependant, il suit de~\cite[Thm. 6.13]{kw} que sa restriction \`a ${\rm
GSO}(\R^{2r})$ se d\'ecompose en somme directe de deux repr\'esentations irr\'eductibles non
isomorphes $${\rm H}_ {d,r}(\R^{2r}) = {\rm H}_ {d,r}(\R^{2r})^+ \oplus
{\rm H}_{d,r}(\R^{2r})^-$$ 
que nous ne chercherons pas \`a distinguer. Concr\`etement, si $e_1,\dots,e_r$ est une base d'un lagrangien $I \subset \C^{2r}$, la fonction $(v_1,\dots,v_r) \mapsto \DET [ e_i \cdot v_j]^d$ est dans 
${\rm H}_ {d,r}(\R^{2r})^{\pm}$ le signe $\pm$ \'etant uniquement d\'etermin\'e par le type du lagrangien $I$. Observons que si $d=0$, auquel cas
${\rm H}_{d,r}(\R^{2r}) \simeq \Lambda^r (\R^{2r})^{\ast}$, ce ph\'enom\`ene est bien connu ! \ps\ps

Soient $\Gamma$ un groupe, $U$ un $\Gamma$-module et $\sigma \in {\rm Aut}(\Gamma)$. On note $U^\sigma$ le $\Gamma$-module obtenu par restriction de $U$ via $\sigma : \Gamma \rightarrow \Gamma$. 

\begin{cor}\label{lemmetrialiteinfini} Soient $d \equiv 0 \bmod 2$. L'application $U \mapsto U^{\tau}$ induit un $3$-cycle 
de l'ensemble constitu\'e des classes d'isomorphisme des trois repr\'esentations de ${\rm PGSO}_8(\R)$ $${\rm H}_{d,1}(\R^8) \otimes \nu^{d/2} \, \, \, {\rm et} \, \, \, {\rm
H}_{d/2,4}(\R^8)^{\pm} \otimes \nu^d.$$  
\end{cor} 
\begin{pf} Observons que les espaces ${\rm H}_{d,g}(\R^8)$ sont munis d'actions
naturelles de ${\rm GO}_8(\C)$ \'etendant les op\'erations
de ${\rm GO}_8(\R)$ consid\'er\'ees pr\'e\-c\'edemment. En particulier, les trois
repr\'esentations de l'\'enonc\'e se factorisent en des repr\'esentations de ${\rm
PGSO}_8(\C)$, qui sont bien s\^ur irr\'eductibles. Rappelons que $V=\R^{8}$.\ps\ps

Soit $D \subset V \otimes \C$ une droite isotrope. Le stabilisateur $S_D \subset {\rm PGSO}_8(\C)$ de la droite $D$ est un sous-groupe parabolique isomorphe \`a ${\rm GSO}(\C^6) \ltimes
\C^6$. Son action naturelle sur $V^{\otimes 2} \otimes
\nu^{-1}$ pr\'eserve la droite $D^{\otimes 2}$, sur laquelle il agit donc
par multiplication par un caract\`ere que nous noterons $\eta_{D}$. Soit $\ell(D) \subset {\rm H}_{d,1}(\R^8) \otimes 
\nu^{d/2}$ la droite de polyn\^omes harmoniques associ\'es \`a $D$ (formule~\eqref{vecteurextremal}). On constate que 
$\ell(D)$ est propre de caract\`ere $\eta_{D}^{d/2}$ sous
l'action de $S_D$. \ps\ps 

De m\^eme, soit $I \subset V\otimes \C$ un lagrangien. Le stabilisateur $S_I
\subset {\rm PGSO}_8(\C)$ de $I$ est un
sous-groupe parabolique isomorphe \`a ${\rm GL}(\C^4)/\{\pm 1\} \ltimes {\rm
Sym}^2(\C^4)$.  Son action naturelle sur $\Lambda^4 V \otimes \nu^{-2}$
pr\'eserve la droite $\Lambda^4 I$, sur laquelle il agit donc par un
caract\`ere que nous noterons $\eta_{I}$. La doite $\ell(I) \subset {\rm H}^{\pm}_{d/2,4}(\R^8) \otimes
\nu^{d}$ de polyn\^omes harmoniques associ\'es \`a $I$,  le signe $\pm$ d\'ependant du type du lagrangien $I$,
est manifestement propre de caract\`ere $\eta_{I}^{d/2}$ sous l'action de $S_I$.\ps\ps

Rappelons que d'apr\`es la th\'eorie du plus haut poids de Cartan-Weyl, \'etant donn\'e un 
sous-groupe parabolique $S$ du groupe semi-simple $G={\rm PGSO}_8(\C)$, et un caract\`ere
polyn\^omial $\eta : S \rightarrow \C^\times$, il existe \`a isomorphisme
pr\`es au plus une repr\'esentation polyn\^omiale irr\'eductible de $G$ dont la restriction
\`a $S$ contienne le caract\`ere $\eta$. De plus, si une telle
repr\'esentation existe pour le couple $(S,\eta)$ avec $\eta 
\neq 1$, elle n'existe pas pour le couple $(S,\eta^{-1})$ (propri\'et\'e de
dominance). \ps\ps

Les observations ci-dessus caract\'erisent donc uniquement les
trois repr\'esentations de l'\'enonc\'e. Pour conclure, on observe que 
$\tau$ permute les trois types de
sous-groupes paraboliques consid\'er\'es ci-dessus d'apr\`es le
lemme~\ref{trialitestudy}, et ce avec leur caract\`ere distingu\'e $S
\rightarrow \C^\times$ not\'e $\eta_{\ast}$ ci-dessus. En effet, cette derni\`ere propri\'et\'e est automatique car si
$S$ est un tel parabolique, on observe que le groupe des caract\`eres
polyn\^omiaux $S \rightarrow \C^\times$ est isomorphe \`a $\Z$, le caract\`ere
$\eta$ \'etant l'unique g\'en\'erateur dominant d'apr\`es le cas $d=2$. \end{pf}

Soit $d$ un entier pair. Quitte \`a \'echanger les signes $\pm$, on peut supposer que $({\rm H}_{d,1}(\R^{8}) \otimes \nu^{{d/2}})^{\tau} \simeq {\rm H}_{d/2,4}(\R^{8})^{+} \otimes \nu^{d}$ d'apr\`es le corollaire ci-dessus. D'apr\`es
le lemme IV.\ref{compatibilitesatake2}, l'automorphisme $\tau$ induit un isomorphisme 
$$\tau^\ast : {\rm M}_{{\rm H}_{d/2,4}(\R^{8})^{+} \otimes \nu^{d}}({\rm PGSO}_8) \isomo {\rm M}_{{\rm H}_{d,1}(\R^{8}) \otimes \nu^{{d/2}}}({\rm PGSO}_8)$$
tel que $T \circ \tau^\ast = \tau^\ast \circ {\rm H}(\tau)(T)$ pour tout $T
\in {\rm H}({\rm PGSO}_8)$. 

\begin{cor}\label{cortrialdim} Soit $d$ un entier pair. On a une suite d'isomorphismes :
$$\xymatrix{ {\rm M}_{{\rm H}_{d/2,4}(\R^{8})^{+} \otimes \nu^{d}}({\rm PGSO}_8) \ar@{->}_{\tau^{\ast}}^{{\wr}}[d]  \ar@{->}^{\, \,  {\rm ind}}_{{\sim}}[r] & {\rm M}_{{\rm H}_{d/2,4}(\R^{8}) \otimes \nu^{d}}({\rm PGO}_{8}) \ar@{->}^{\hspace{.5cm} \mu^{\ast}}_{\hspace{.5cm} {\sim}}[r] & {\rm M}_{{\rm H}_{d/2,4}(\R^{8})}({\rm O}_{8}) \\
{\rm M}_{{\rm H}_{d,1}(\R^{8}) \otimes \nu^{{d/2}}}({\rm PGSO}_8)\ar@{->}^{\hspace{.8cm} \mu^{\ast}}_{\hspace{.8cm} \sim}[r] & {\rm M}_{{\rm H}_{d,1}(\R^{8})}({\rm SO}_{8})  & }
$$
\end{cor}

\begin{pf} Tous ces isomorphismes ont d\'ej\`a \'et\'e d\'ecrits, \`a l'exception de celui de la ligne du bas, not\'e $\mu^{\ast}$. C'est le morphisme d\'efinit par le restriction des fonctions par la bijection ${\mathcal R}({\rm SO}_{n}) \rightarrow {\mathcal R}({\rm PGSO}_{n})$, qui est un isomorphisme pour des raisons identiques \`a celles \'evoqu\'ees dans la d\'emonstration du lemme \ref{inflationmonGO}.
\end{pf}

\begin{cor}\label{derniercor8} La valeur propre de ${\rm T}_p$ sur la droite
${\rm M}_{{\rm H}_{4,4}(\R^8)}({\rm O}_8)$ est $$p^{-4}\, \frac{p^4-1}{p-1}\,
\tau(p).$$
\end{cor}

\begin{pf} On applique le corollaire \ref{cortrialdim} \`a $d=8$. Le corollaire \ref{coractionKpm} et le lemme \ref{lemmetrialitep} montrent que la valeur propre de l'op\'erateur de Hecke ${\rm T}_{p}$ sur la droite ${\rm M}_{{\rm H}_{4,4}(\R^{8})^{+} \otimes \nu^{8}}({\rm PGSO}_8)$ est celle de l'\'enonc\'e. On conclut car $\mu^{\ast}$ et ${\rm ind}$ sont ${\rm H}({\rm O}_{8})$-\'equivariants (Lemme \ref{inflationmonGO}, \S \ref{pgpso}).
\end{pf}

\subsection{Une derni\`ere s\'erie th\^eta et fin de la d\'emonstration}\label{constthetag4}

Consid\'erons pour finir l'application $$\vartheta_{4,4} :
{\rm M}_{{\rm H}_{4,4}(\R^8)}({\rm
O}_8) \longrightarrow {\rm S}_8({\rm Sp}_8(\Z)).$$

\begin{prop}\label{thetanonnul44} $\vartheta_{4,4}$ est un isomorphisme. \end{prop}

\begin{pf} En effet, les deux espaces \'etant de dimension $1$ (formule \eqref{dims8sp8}, Corollaire \ref{derniercor8}), il suffit de voir que
cette application est non nulle, ce qui a d\'ej\`a \'et\'e v\'erifi\'e par
Breulmann et Kuss dans~\cite{bkuss}. Expliquons rapidement comment
on peut proc\'eder. \ps

Soient $e=(e_1,\dots,e_4)$ un quadruplet d'\'el\'ements de ${\rm E}_8
\otimes \C$ engendrant un lagrangien, et $P_e(v_1,\dots,v_4)=\DET [ e_i
\cdot v_j]_{1\leq i,j \leq 4}$; pour tout entier $d\geq 0$ on a $P_e^d \in
{\rm H}_{d,4}(\R^8)$.  Soient $Q \subset {\rm E}_8$ un sous-r\'eseau de rang
$4$, ainsi que $v_1,\dots,v_4$ une $\Z$-base de $Q$.  La relation
$P_e(\gamma(v_1),\dots,\gamma(v_4)))=\det(\gamma)P_e(v_1,\dots,v_4)$ pour
tout $\gamma \in {\rm GL}(Q)$ montre que $P_e(v_1,\dots,v_4)^d$ ne d\'epend
pas du choix des $v_i$ lorsque $d$ est pair; il y a donc un sens \`a
le noter $P_e(Q)^d$.  En particulier, si $d$ est pair le coefficient
de Fourier de la s\'erie th\^eta $\vartheta_{d,4}({\rm E}_8,P_e^{d})$
correspondant \`a la matrice de Gram d'une $\Z$-base de $Q$ est
$${\rm c}_Q(P_e^d)=|{\rm O}(Q)|\sum_M P_e(M)^d,$$ la somme portant sur l'ensemble
des sous-r\'eseaux $M \subset {\rm E}_8$ isom\'etriques \`a $Q$.  Nous
donnons quelques valeurs num\'eriques dans la table~\ref{thetagenre4}. \ps\ps

\begin{table}[htp]
\caption{Quelques valeurs de $\frac{{\rm c}_Q(P_e^d)}{|{\rm O}(Q)|}$ o\`u $e=(\varepsilon_{2j-1}+i\varepsilon_{2j})_{1\leq j \leq 4}$}
{\small \renewcommand{\arraystretch}{1.4} \medskip
\begin{center}
\begin{tabular}{|c||c|c|c|c|c|c|}
\hline  $Q \hspace{3 mm}{\backslash} \hspace{3 mm} d$ &  $0$ & $2$ & $4$ & $6$ & $8$ & $10$ \\

\hline ${\rm D}_4$ & $3150$ & $0$ & $4800$ & $-4800$ & $43200$ & $-81600$ \\

\hline ${\rm A}_4$ & $24192$ & $0$ &  $-23040$ & $-46080$ & $-69120$ & $-92160$ \\

\hline ${\rm A}_1 \oplus {\rm A}_3$ & $151200$ & $0$ & $115200$ & $1267200$ & $6566400$ & $7718400$ \\

\hline  ${{\rm A}_2}^{2}$ &  $67200$ & $0$ & $115200$ & $-1382400$ & $4492800$ & $-43084800$ \\

\hline ${{\rm A}_1}^2\oplus {\rm A}_2$ &  $302400$ & $0$ & $-691200$ &  $2073600$ & $85017600$ & $214963200$\\ 

\hline ${{\rm A}_1}^4$ &  $122850$ & $0$ & $576000$ & $-6796800$ & $191808000$ & $-343641600$ \\
\hline \hline
{\scriptsize $\dim {\rm M}_{{\rm H}_{d,4}(\R^8)}({\rm O}_8)$} & $1$ & $0$ & $1$ & $1$ & $1$ & $2$ \\
\hline
\end{tabular}
\end{center}}
\label{thetagenre4}
\end{table}

Dans cette table, $Q \simeq {\rm Q}(R)$ o\`u $R$ est un syst\`eme de racines (de
type ADE) de rang $4$ (\S I.3) et $(\varepsilon_{j})_{1 \leq j \leq 8}$ d\'esigne la base canonique de $\R^{8}$.  Il n'est pas difficile d'\'enumerer \`a l'aide de
l'ordinateur tous les sous-r\'eseaux de ${\rm E}_8$ isom\'etriques \`a $Q$. Par exemple, si $\Phi$ d\'esigne un syst\`eme positif de ${\rm R}(E_8)$, et
si $<$ d\'esigne un ordre total arbitraire fix\'e de $\Phi$, les
sous-r\'eseaux de ${\rm E}_8$ isom\'etriques \`a ${\rm D}_4$ sont en
bijection avec les quadruplets $(r_1,r_2,r_3,r_4)$ d'\'el\'ements de $\Phi$
tels que $r_1<r_2<r_3$, et tels que les \'el\'ements $r_1,r_2$ et $r_3$
soient deux \`a deux orthogonaux et de produit scalaire $-1$ avec $r_4$. Nous renvoyons \`a la feuille de calculs \cite{clcalc} pour une impl\'ementation de cet algorithme sous \texttt{PARI} et pour une justification de la table \ref{thetagenre4}. \ps \ps
La proposition d\'ecoule du fait que ${\rm c}_{{\rm D}_4}(P_e^4) \neq 0$. 
Constatons \'egalement, pour se rassurer, que l'on retrouve bien
l'\'egalit\'e ${\rm c}_{{\rm D}_4}(P_e^4)=-{\rm c}_{{\rm A}_4}(P_e^4)$, ce qui est en
accord avec la proposition \ref{Fnonnulle}.  \end{pf}

\begin{remark}{\rm La derni\`ere ligne de la table \ref{thetagenre4} d\'ecoule de
l'isomorphisme ${\rm M}_{{\rm H}_{d,4}(\R^8)}({\rm O}_8) \simeq {\rm
M}_{{\rm H}_{2d,1}(\R^8)}({\rm SO}_8)$ ( $ \simeq  {\rm H}_{2d,1}(\R^8)^{{\rm W}({\rm E}_{8})^{+}}$) donn\'e par le corollaire \ref{derniercor8}, ainsi que du lemme~\ref{invariantscox}. La nullit\'e de ces espaces quand $d=2$ explique celle de la colonne $d=2$ de la table. La table montre donc que $\vartheta_{d,4}$ est injectif quand $d$ est pair et $\leq 8$. En faisant varier la base lagrangienne $e$, il est facile de
v\'erifier que $\vartheta_{10,4}$ est \'egalement injective.} \end{remark}

Cela conclut la d\'emonstration du (ii) du th\'eor\`eme~\ref{vpnontrivdim16},
en vertu du corollaire~\ref{derniercor8} et des relations de commutation
d'Eichler (\ref{eichlerthetaharmgen}). La suite d'isomorphismes suivante r\'esume assez bien notre
d\'emonstration. On a pos\'e $W={\rm H}_{{8,1}}(\R^{8})$, $U={\rm H}_{{4,4}}(\R^{8})$ et $U^{+}={\rm H}_{{4,4}}(\R^{8})^{+}$.

$$\xymatrix{{\rm M}_{W \otimes \nu^{4}}({\rm PGSO}_8) \ar@{<-}^{\wr}_{\tau^\ast}[d] & {\rm M}_{W \otimes \nu^{4}}({\rm PGO}_8)
\ar@{->}^{\hspace{.3cm}\sim}_{\hspace{.3cm}{\rm res}}[l] \ar@{->}_{\hspace{.4cm} \sim}^{\hspace{.4cm} \mu^\ast}[r] & {\rm M}_{W}({\rm O}_8)
\ar@{->}_{\hspace{-.3cm}\sim}^{\hspace{-.3cm}\vartheta_{8,1}}[r] & {\rm S}_{12}({\rm SL}_2(\Z)) \\
{\rm M}_{U^{+} \otimes \nu^{8}}({\rm PGSO}_8) \ar@{->}^{\sim}_{{\rm ind}}[r] & {\rm M}_{U \otimes \nu^{8}}({\rm
PGO}_8) \ar@{->}^{\hspace{.4cm} \sim}_{\hspace{.4cm} \mu^\ast}[r] & {\rm M}_U({\rm O}_8)
\ar@{->}^{\hspace{-.3cm}\sim}_{\hspace{-.3cm}\vartheta_{4,4}}[r] & {\rm S}_8({\rm Sp}_8(\Z))
}$$

\section{Appendice : un exemple simple de relation d'Eichler}\label{appendiceC}

On se propose de d\'emontrer la formule~IV.\eqref{eichlerthetaharmgen} pour $g=1$. Soient $L$ un r\'eseau unimodulaire pair de rang $r$ et $P : L \otimes \R \rightarrow \C$ un polyn\^ome harmonique homog\`ene de
degr\'e $d$. On rappelle que la s\'erie th\^eta associ\'ee $\vartheta(L,P)=\sum_{v \in L} P(v) q^{\frac{v \cdot v}{2}}$ 
est un \'el\'ement de ${\rm M}_{d+r/2}(\SL_2(\Z))$. \ps\ps

\begin{thm} \label{appendiceCthm} Soient $L$ un r\'eseau unimodulaire pair de rang $r$, $P : L \otimes \R \rightarrow \C$ un polyn\^ome harmonique homog\`ene de
degr\'e $d$, et $p$ un nombre premier. On a la relation
$$\sum_{L'} \vartheta(L',P) = \left( p\,\frac{p^{r-3}-1}{p-1}\, +\, p^{-d}\,{\rm T}(p^2) \right) \vartheta(L,P),$$
o\`u $L' \subset \frac{1}{p} L$ parcourt les $p$-voisins de $L$. \end{thm}

\begin{pf}
On notera ${\rm a}_n(g)$ le $n$-i\`eme coefficient de Fourier de la forme modulaire $g \in {\rm M}_k(\SL_2(\Z))$. On rappelle la relation  \cite[p. 164]{serre} $${\rm a}_n({\rm T}(p^2)(g))=\sum_{d|(p^2,n)} d^{k-1} {\rm a}_{\frac{np^2}{d^2}}(g).$$ On pose $f=\sum_{L'} \vartheta(L',P)$ o\`u $L' \subset \frac{1}{p} L$ parcourt les $p$-voisins de $L$. \ps \medskip

Fixons $n \geq 1$ un entier, pour l'instant premier \`a $p$, et posons ${\rm q}(x)=\frac{x \cdot x }{2}$ pour tout $x \in L \otimes \R$.
Consid\'erons l'ensemble $X$ des couples $(L',w)$ o\`u $L'$ est un $p$-voisin de $L$ et $w$ est un \'el\'ement de $L'$ tel que ${\rm q}(w)=n$. Soit aussi $Y=\{v \in L, {\rm q}(v)=np^2\}$. On dispose d'une application \'evidente $$\pi : X \rightarrow Y, \, \, \, \, (L',w) \mapsto pw.$$
Nous allons voir que $\pi$ est surjective et examiner ses fibres. 
Soit $v \in Y$, il y a deux cas : \begin{itemize} \ps\ps \medskip

\item[(1)] $v$ est dans $pL$. Dans ce cas, $w=v/p$ est dans $L$ et il y a donc
autant de $p$-voisins $L'$ de $L$ contenant $w$ que de sous-r\'eseaux
d'indice $p$ de $L$ contenant $w$ et qui sont l'orthogonal d'une droite isotrope
modulo $p$ (les ``$M$'' du \S III.1). Comme on a $({\rm q}(w),p)=1$, $w$ est non isotrope mod $pL$ : on trouve 
$|\pi^{-1}(\{v\})|=\frac{p^{r-2}-1}{p-1}$. \ps \ps

\item[(2)] $v$ n'est pas dans $pL$. Dans ce cas, il engendre une droite isotrope dans $L/pL$. Si $M$ est son orthogonal modulo $p$, c'est-\`a-dire l'ensemble des $x \in L$ tels que $x.v \in p\Z$, alors $L'=M+\Z \frac{v}{p}$ est un $p$-voisin de $L$
car $p^2$ divise ${\rm q}(v)$. Mieux, c'est l'unique $p$-voisin de $L$ contenant $w=v/p$. 
En effet, si $K$ est un tel $p$-voisin, et si $N=L \cap K$, alors $v$ est dans $N$ et pour tout $x$ dans $K$ on a $x.\frac{v}{p} \in \Z$. En particulier, $N \subset M$, puis $N=M$ et $K=L'$. Ainsi, $|\pi^{-1}(\{v\})|=1$.\ps\end{itemize}
\ps \medskip

\noindent Supposons d'abord $P=1$ (et donc $h=0$). Cette analyse montre  
$${\rm a}_n(f)=|X|=\frac{p^{r-2}-1}{p-1}{\rm a}_n(\vartheta(L,1)) +
({\rm a}_{np^2}(\vartheta(L,1))-{\rm a}_n(\vartheta(L,1))),$$
qui est la formule cherch\'ee du moins pour les coefficients
d'indice premier \`a $p$. Quand $P$ est quelconque, ${\rm a}_n(f)$ est la somme des
$p^{-d}P(v)$ pour $(L',v/p)$ parcourant $X$, et donc
$${\rm a}_n(f)=\frac{p^{r-2}-1}{p-1}{\rm a}_n(\vartheta(L,P)) +
(p^{-d}{\rm a}_{np^2}(\vartheta(L,P))- {\rm a}_n(\vartheta(L,P))),$$

\bigskip

\noindent On conclut par le lemme suivant :

\bigskip

\begin{lemme}\label{appendiceClemme} Soient $g \in {\rm M}_k(\SL_2(\Z))$, $k>0$, et $p$ premier. On suppose 
que ${\rm a}_n(g)=0$ pour tout $n$ premier \`a $p$. Alors $g=0$.\end{lemme}

\begin{pf} En effet, la fonction holomorphe $g(\tau)$ est alors invariante par $\tau \mapsto \tau + \frac{1}{p}$, et on conclut car le sous-groupe de ${\rm SL}_2(\R)$ engendr\'e par cette translation et ${\rm SL}_2(\Z)$ n'est pas discret. \end{pf}

On peut aussi  v\'erifier la formule pour tous les coefficients, en introduisant encore $\pi : X \rightarrow Y$ comme pr\'ec\'edemment. On constate que le d\'ecompte n'est pas modifi\'e pour le cas (2), mais uniquement pour le (1). \ps \medskip

{\it Premier sous-cas} : $\frac{n}{p} \in \Z - p\Z$. Soit $v$ dans $Y$ de la forme $pw$ avec $w \in L$. Remarquons que $w \notin pL$ car $p^2$ ne divise pas ${\rm q}(w)=n$, par contre $w$ est isotrope dans $L/pL$. Mais si $x \in L/pL$ est un vecteur isotrope, il existe $$1 + p \,{\rm c}_{r-2}(p) = \frac{p^{r-2}-1}{p-1} +  p^{r/2-1}$$
droites isotropes de $L/pL$ orthogonales \`a $x$, c'est donc aussi $|\pi^{-1}(\{v\})|$. Ici, ${\rm c}_i(p)=\frac{p^{i-1}-1}{p-1}+p^{i/2-1}$ est le cardinal de la quadrique hyperbolique de rang $i$ sur $\Z/p\Z$. On en d\'eduit
$${\rm a}_n(f)=(\frac{p^{r-2}-1}{p-1}+p^{r/2-1}){\rm a}_n(\vartheta(L,P)) + (p^{-d}{\rm a}_{np^2}(\vartheta(L,P))-{\rm a}_n(\vartheta(L,P))),$$
ce qui conclut. \ps \medskip

{\it Second sous-cas} : $p^2$ divise $n$.  Soit $v$ dans $Y$ de la forme $pw$ avec $w$ dans $L$. Alors $w$ est isotrope dans $L/pL$. Si il est nul, il est dans tous les $p$-voisins de $L$ : $|\pi^{-1}(\{v\})|={\rm c}_r(p)$. Sinon, on a comme ci-dessus $|\pi^{-1}(\{v\})|=\frac{p^{r-2}-1}{p-1} +  p^{r/2-1}$. Mais $${\rm c}_r(p)-\frac{p^{r-2}-1}{p-1} -  p^{r/2-1}=p^{r-2}$$
On en d\'eduit l'identit\'e 
$${\rm a}_n(f)=p^{d+r-2}{\rm a}_{n/p^2}(\vartheta(L,P))+(\frac{p^{r-2}-1}{p-1}+p^{r/2-1}){\rm a}_n(\vartheta(L,P)) $$ $$+ (p^{-d}{\rm a}_{np^2}(\vartheta(L,P))-{\rm a}_n(\vartheta(L,P))).$$
Ce qui conclut. \end{pf}


\chapter{Param\'etrisation de Langlands}\label{chap5}

\section{Rappels et notations sur les $k$-groupes r\'eductifs}\label{rapredchap6}

Soit $k$ un corps alg\'ebriquement clos. Nous renvoyons aux trait\'es de
Springer~\cite{springerlivre} et Borel~\cite{borelgp} pour la th\'eorie des
$k$-groupes r\'eductifs.  Notre convention est qu'un tel $k$-groupe est
connexe.  Si $k$ est de caract\'eristique nulle, on rappelle qu'un
$k$-groupe connexe $G$ est dit r\'eductif si la cat\'egorie de ses
$k$-repr\'esentations de dimension finie est semi-simple.  \ps\ps 

Si $A$ est un anneau quelconque, un $A$-groupe $G$ 
est dit r\'eductif s'il est lisse sur $A$ et si pour tout homomorphisme de $A$ 
vers un corps alg\'ebriquement clos $k$, le groupe $G \times_A k$ est r\'eductif
: voir \cite{sga} et \cite{conradsga}. Les $A$-groupes classiques \'etudi\'es
au~\S~II. 1 sont donc r\'eductifs \cite[\S 23]{borelgp} \cite[App. C]{conradsga}, 
except\'e le groupe orthogonal en dimension paire qui n'est pas connexe, ainsi que les groupes de
similitudes et similitudes projectives associ\'es pour la m\^eme raison. Un $A$-groupe est dit semi-simple s'il est r\'eductif et si son centre est fini sur
$A$. Le centre (sch\'ematique) d'un $A$-groupe r\'eductif $G$ sera not\'e ${\rm Z}(G)$.
Une {\it isog\'enie centrale} $G \rightarrow G'$ entre deux $A$-groupes
r\'eductifs est un morphisme fini et plat de $A$-groupes qui est surjectif et
de noyau (au sens sch\'ematique) inclus dans ${\rm Z}(G)$. Plus g\'en\'eralement, un morphisme de $A$-groupes $G \rightarrow G'$ 
sera dit {\it central} si le morphisme induit $G \times {\rm Z}(G') \rightarrow G'$ est plat, surjectif et de noyau inclus dans ${\rm Z}(G) \times {\rm Z}(G')$. 

\subsection{La donn\'ee radicielle bas\'ee d'un $k$-groupe r\'eductif}
\label{gpdepl}

Soit $k$ un corps alg\'ebriquement clos. La th\'eorie des syst\`emes de
racines des $k$-groupes r\'eductifs, convenablement formul\'ee, produit une
\'equivalence de cat\'egorie canonique $$\Psi : \mathcal{C}_k \isomo
\mathcal{D}$$ entre la cat\'egorie $\mathcal{C}_k$ des $k$-groupes
r\'eductifs ``\`a automorphismes int\'erieurs pr\`es'' et la cat\'egorie
$\mathcal{D}$ des {\it donn\'ees radicielles bas\'ees}.  Cette
classification est due \`a Chevalley dans le cas des groupes semi-simples et
au s\'eminaire Demazure-Grothendieck \cite[Exp.  XXI]{sga} en g\'en\'eral, o\`u elle est m\^eme
\'etudi\'ee sur un anneau $k$ quelconque.  Nous nous contenterons
ci-dessous de pr\'eciser son \'enonc\'e, en renvoyant \`a
\cite{springerlivre} pour un expos\'e d\'etaill\'e, \'egalement r\'esum\'e
dans~\cite{springercorvallis} et~\cite[II, Ch. 1]{jantzen}, ainsi qu'\`a Kottwitz~\cite[\S 1]{kottwitz}
pour la formulation intrins\`eque adopt\'ee ici.  \ps\ps 

Si $G$ est un $k$-groupe r\'eductif, on note ${\rm Int}(G)$ le groupe
des automorphismes {\it int\'erieurs} de $G$, i.e.  de la forme ${\rm int}_g\, : x \mapsto gxg^{-1}$, pour $g \in G(k)$. Si $G$ et $G'$ sont deux $k$-groupes
r\'eductifs on note \'egalement ${\rm Hom}_c(G,G')$ l'ensemble des morphismes centraux 
de $G$ vers $G'$. Il est muni d'une action \'evidente de ${\rm Int}(G')$. 
On observe qu'il y a un sens \`a consid\'erer la cat\'egorie $\mathcal{C}_k$ dont
les objets sont les $k$-groupes r\'eductifs et ayant pour morphismes $G
\rightarrow G'$ l'ensemble quotient ${\rm Hom}_c(G,G')/{\rm Int}(G')$, la 
composition des morphismes se d\'eduisant de celle des morphismes centraux par
passage au quotient. \ps \medskip

\noindent Rappelons qu'une {\it donn\'ee radicielle bas\'ee} est la donn\'ee : \begin{itemize} \ps \medskip

\item[-] de deux groupes ab\'eliens libres de rang fini $X$ et $X^\vee$ munis d'un accouplement parfait $\langle-,-\rangle : X \times X^\vee
\rightarrow \Z$, \ps\ps  
\item[-] de parties finies $\Phi \subset X$ et $\Phi^\vee \subset X^\vee$
munies d'une bijection $\Phi \rightarrow \Phi^\vee$ not\'ee $\alpha \mapsto \alpha^\vee$, \ps\ps 
\item[-] de sous-ensembles $\Delta \subset \Phi$ et $\Delta^\vee \subset \Phi^\vee$ tels que $\Delta^\vee=\{\alpha^\vee, \alpha \in
\Delta\}$,  \par \medskip
\end{itemize} \noindent  soumis aux conditions suivantes : \par  \medskip
\begin{itemize}
\item[-]  pour tout $\alpha \in \Phi$, on a $\langle \alpha, \alpha^\vee \rangle =2$, \ps\ps 

\item[-]  si $s_\alpha \in {\rm End}(X)$ d\'esigne la reflexion $x \mapsto x -
\langle x, \alpha^\vee \rangle \alpha$, et si $s_{\alpha^\vee} \in {\rm
End}(X^\vee)$ est d\'efinie de mani\`ere similaire en \'echangeant $\alpha$ et $\alpha^\vee$, alors pour tout $\alpha \in \Phi$ on a $s_\alpha(\Phi)=\Phi$ et $s_{\alpha^\vee}(\Phi^\vee)=\Phi^\vee$. \ps \medskip
\end{itemize}
\noindent Il r\'esulte de ces axiomes que si ${\mathrm Q}(\Phi) \subset X$ d\'esigne
le groupe ab\'elien engendr\'e par les \'el\'ements de $\Phi$, alors $\Phi$
est un syst\`eme de racines dans ${\rm Q}(\Phi)\otimes \Q$ au sens de
Bourbaki\, \cite[Ch. VI]{bourbaki}. On supposera enfin que \begin{itemize} \ps \smallskip

\item[-] $\Phi$ est r\'eduit\footnote{Soulignons que 
cette hypoth\`ese ne fait pas partie des axiomes dans les r\'ef\'erences sus-cit\'ees, elle nous \'evitera quelques circonvolutions.} et $\Delta$ en est une
base. \ps \medskip
\end{itemize}

Un morphisme $\psi_1 \rightarrow \psi_2$ entre deux donn\'ees radicielles bas\'ees
$\psi_i=(X_i,\Phi_i,\\ \Delta_i,X_i^\vee,\Phi_i^\vee,\Delta_i^\vee)$ est la
donn\'ee d'une application lin\'eaire $X_2 \rightarrow X_1$ induisant une
bijection $\Phi_2 \rightarrow
\Phi_1$ envoyant $\Delta_2$ sur $\Delta_1$, et dont la transpos\'ee $X_1^\vee \rightarrow X_2^\vee$ induit
\'egalement une bijection $\Phi_1^\vee \rightarrow \Phi_2^\vee$ envoyant
$\Delta_1^\vee$ sur $\Delta_2^\vee$.  Cela d\'efinit la cat\'egorie $\mathcal{D}$. Une {\it isog\'enie} $\psi_1 \rightarrow \psi_2$ est un morphisme comme ci-dessus induisant un isomorphisme $X_2 \otimes \Q \rightarrow X_1 \otimes \Q$.
Il ne nous reste qu'\`a rappeler la d\'efinition du foncteur $\Psi$. \ps\ps

Si $G$ est un $k$-groupe r\'eductif, la donn\'ee radicielle
bas\'ee $\Psi(G)$ qui lui est associ\'ee s'obtient ainsi. On choisit un tore
maximal $T$ de $G$ ainsi qu'un sous-groupe de Borel $B$ contenant $T$. On note $${\mathrm X}^*(T)={\rm
Hom}(T,{\mathbb G}_m), \, \, \, {\mathrm X}_*(T)={\rm Hom}({\mathbb
G}_m,T),$$ les groupes ab\'eliens libres de rang fini respectifs
constitu\'es des caract\`eres et des co-caract\`eres du tore $T$. Ils sont
munis d'une dualit\'e parfaite \'evidente $\langle-,-\rangle : {\mathrm
X}^*(T) \times {\mathrm X}_*(T) \longrightarrow \Hom({\mathbb G}_m,{\mathbb
G}_m)=\Z$. On pose alors
$$\Psi(G,T,B)=(X^\ast(T),\Phi(G,T),\Delta(G,T,B),X_\ast(T),\Phi^\vee(G,T),\Delta^\vee(G,T,B)),$$ o\`u 
$\Phi(G,T)$ (resp. $\Phi^\vee(G,T)$) d\'esigne l'ensemble des
racines (resp. coracines\footnote{Soit $\alpha \in \Phi(G,T)$, soit $T_\alpha \subset T$ la composante neutre du noyau de $\alpha : T \rightarrow \mathbb{G}_m$,  et soit $Z_\alpha$ le groupe d\'eriv\'e du centralisateur dans $G$ de $T_\alpha$. C'est un $k$-groupe isomorphe \`a ${\rm
SL}_2$ ou ${\rm PGL}_2$. On rappelle que la coracine 
$\alpha^\vee \in {\rm X}_*(T)$ est l'unique cocaract\`ere d'image dans
$Z_\alpha$ tel que $\langle \alpha,\alpha^\vee
\rangle = 2$. }) de $G$ relativement \`a $T$, et o\`u
$\Delta(G,T,B)$ est la base de $\Phi(G,T)$ associ\'ee au syst\`eme positif 
de $\Phi(G,T)$ intervenant dans ${\rm Lie}(B)$. C'est une donn\'ee radicielle bas\'ee. \par \smallskip Si l'on change le couple $(T,B)$ en $(T',B')$, il existe un \'el\'ement $g \in G(k)$ unique modulo $T(k)$ tel que $g T g^{-1} = T'$ et $g B g^{-1} = B'$. L'automorphisme int\'erieur $\, \, {\rm int}_g \, \, $ induit un isomorphisme $\Psi(G,T,B) \isomo \Psi(G,T',B')$ dans $\mathcal{D}$ qui est ind\'ependant du choix de $g$. \ps \smallskip Suivant Kottwitz~\cite[\S 1]{kottwitz},
on d\'efinit $\Psi(G)$ comme \'etant la limite inductive (ou projective!) des $\Psi(G,T,B)$, index\'ee
par les couples $(T,B)$, avec pour morphismes de transition les isomorphismes induits par des \'el\'ements de ${\rm
Int}(G)$. La construction $G \mapsto \Psi(G)$ est fonctorielle en les morphismes centraux et fait correspondre \`a une isog\'enie centrale une isog\'enie de donn\'ees radicielles. 
En particulier, le groupe ${\rm Aut}(G)$ des automorphismes du $k$-groupe $G$ agit sur ${\rm Aut}_{\mathcal{D}}(\Psi(G))$, le sous-groupe ${\rm Int}(G)$ agissant trivialement. \ps\ps 

Jusqu'ici $k$ \'etait un corps alg\'ebriquement clos. Mentionnons que la d\'efinition du foncteur $\Psi$ s'\'etend
verbatim au cas d'un anneau quelconque $k$ si l'on se restreint \`a la sous-cat\'egorie
de $\mathcal{C}_k$ constitu\'ee des $k$-groupes r\'eductifs {\it d\'eploy\'es} \cite[XXII
Prop. 1.14]{sga}, \cite[II Ch. 1]{jantzen}. Si $k$ est un anneau int\`egre tel que ${\rm
Pic}(k)=0$ (par exemple un corps ou un anneau principal), ce sont les $k$-groupes r\'eductifs poss\'edant un tore maximal
d\'eploy\'e, i.e. isomorphe
\`a une puissance de $\mathbb{G}_m$. En particulier, il y a un sens \`a parler de la donn\'ee radicielle
bas\'ee $\Psi(G)$ d'un tel $k$-groupe. Elle s'identifie canoniquement dans
$\mathcal{D}$ \`a celle de $G \times_k K$ pour tout homomorphisme de $k$
vers un corps alg\'ebriquement clos $K$. \ps
\medskip
\noindent {\sc Vocabulaire}\ps\ps 

Lorsqu'un $k$-groupe r\'eductif $G$ a une donn\'ee radicielle bas\'ee bien
d\'efinie, on parlera librement du syst\`eme de
racines de $G$, des racines simples/positives de $G$, du groupe de Weyl de
$G$, etc... pour d\'esigner les objets analogues d\'eduits de $\Psi(G)$. Par
exemple, le groupe de Weyl de $G$ est le sous-groupe $W
\subset {\rm Aut}(X)$ engendr\'e par l'ensemble des $s_\alpha$ avec
$\alpha \in \Delta$, o\`u $\Psi(G)=(X,\Phi,\Delta,X^\vee,\Phi^\vee,\Delta^\vee)$. Un \'el\'ement de $X$ est appel\'e {\it poids} de $G$, et le 
groupe ab\'elien $X$ le {\it r\'eseau des poids} de $G$; de m\^eme $X^\vee$ est le {\it r\'eseau des co-poids} de $G$.

\subsection{Dual de Langlands}\label{dualitedorad}

Si $\psi=(X,\Phi,\Delta,X^\vee,\Phi^\vee,\Delta^\vee)$ est une donn\'ee radicielle
bas\'ee, alors 
$$\psi^\vee=(X^\vee,\Phi^\vee,\Delta^\vee,X,\Phi,\Delta)$$ en est encore une de
mani\`ere \'evidente, appel\'ee {\it donn\'ee duale} de $\psi$. L'association
$\psi \mapsto \psi^\vee$ d\'efinit un endofoncteur contravariant
involutif de $\mathcal{D}$. Lorsque $k$ est alg\'ebriquement clos, il
induit via l'\'equivalence de cat\'egories $\Psi$ une involution de $\mathcal{C}_k$ : c'est le
point de d\'epart de la notion de dual de Langlands, \`a ceci pr\`es que l'on fait intervenir le corps de complexes. \ps\ps 
Pr\'ecis\'ement, si $G$
est un $k$-groupe r\'eductif d\'eploy\'e, un {\it groupe dual} de $G$ au sens de Langlands
est la donn\'ee d'un $\C$-groupe r\'eductif $\widehat{G}$ et d'un
isomorphisme $\Psi(\widehat{G}) \isomo \Psi(G)^\vee$ dans $\mathcal{D}$. Le $\C$-groupe $\widehat{G}$ est alors uniquement d\'etermin\'e par $G$ \`a isomorphismes int\'erieurs pr\`es. Par abus de language, on l'appelle 
le {\it dual de Langlands de $G$}. 

\subsection{Exemples} \label{exdorad}Nous laissons au lecteur le soin de v\'erifier que 
$\widehat{{\rm GL}}_n(\C) \simeq \GL_n(\C)$ et $\widehat{{\rm PGL}}_n(
\C) \simeq \SL_n(\C)$. En revanche, il nous sera utile de d\'etailler les cas (tr\`es
classiques!) des groupes orthogonaux et
symplectiques. \ps\ps 

Nous utiliserons \`a plusieurs reprises la construction suivante. Soient $\psi=(X,\Phi,\Delta,X^\vee,\Phi^\vee,\Delta^\vee)$ une donn\'ee
radicielle bas\'ee et $Y \subset X \otimes \Q$ un sous-groupe de type
fini contenant $\Phi$. On
suppose $$\Phi^\vee
\subset Y^\sharp:=\{ x \in X^\vee \otimes \Q, \, \, \, \langle y, x \rangle \subset
\Z\, \, \, \, \, \forall y \in Y\}.$$
\noindent L'orthogonal $Y^\bot$ de $Y$ dans $X^\vee \otimes \Q$ est alors
d'intersection nulle avec ${\rm Q}(\Phi^\vee)$, et si  $\pi : {\rm
Q}(\Phi^\vee) \rightarrow Y^\sharp/Y^\perp$ d\'esigne l'application canonique, on constate que
$$\psi'=(Y,\Phi,\Delta,Y^\sharp/Y^\bot,\pi(\Phi^\vee),\pi(\Delta^\vee))$$ est une
donn\'ee radicielle bas\'ee de mani\`ere \'evidente. 
Si l'on a une inclusion $Y \subset X$ (resp.  $X \subset Y$), elle induit un
morphisme $\psi \rightarrow \psi'$ (resp.  une isog\'enie $\psi' \rightarrow
\psi$).  Cette construction appliqu\'ee \`a $\psi^\vee$ fournit par ailleurs
une construction similaire dans laquelle caract\`eres et co-caract\`eres
sont \'echang\'es.  \ps \bigskip

\noindent Dans ce qui suit $k$ est un anneau quelconque.\ps \medskip

\noindent {\sc Groupe sp\'ecial orthogonal pair et ses variantes} \medskip

Soient $r\geq 2$ un entier,
$U=k^r$ et
$V={\rm H}(U)=U\oplus U^*$ le $\rm{q}$-module hyperbolique sur $U$ (\S~II.1).  Le
$k$-groupe $\widetilde{G}=\rm{GSO}_V$ est r\'eductif et d\'eploy\'e.  \ps\ps 

Si
$(e_i)_{i=1}^{r}$ est une $k$-base de $U$, et si $e_i^* \in U^*$ d\'esigne
sa base duale, le $k$-sous-groupe $\widetilde{T}$ de $\widetilde{G}$
pr\'eservant chacune des droites $k e_i$ et $k e_j^*$ est un tore maximal
d\'eploy\'e de $\widetilde{G}$.  Le $k$-sous-groupe de $\widetilde{G}$
pr\'eservant le drapeau complet de $U$ associ\'e \`a $\{e_1\}$,
$\{e_1,e_2\}$, $\dots$, est un sous-groupe de Borel contenant
$\widetilde{T}$.  \ps\ps 

Soit $\varepsilon_i \in \mathrm{X}^*(\widetilde{T})$ le
caract\`ere de $\widetilde{T}$ agissant sur $k e_i$, soit $\nu :
\widetilde{G} \rightarrow {\mathbb{G}}_m$ le facteur de similitude, et
soit $\varepsilon_0$ la restriction de $\nu$ \`a $\widetilde{T}$.  On constate que
$\widetilde{T}$ agit
sur $k e_j^\ast$ par multiplication par le caract\`ere
$-\varepsilon_j+\varepsilon_0$.  Les $\varepsilon_i$, $i=0,\dots,r$ forment
une $\Z$-base de $\mathrm{X}^*(\widetilde{T})$.  \ps\ps 

L'ensemble
$\Phi(\widetilde{G},\widetilde{T})$ est
constitu\'e des $\pm (\varepsilon_i - \varepsilon_j)$ et $\pm
(\varepsilon_i+\varepsilon_j-\varepsilon_0)$ pour $1\leq i < j \leq r$.  De
plus, $\Delta(\widetilde{G},\widetilde{T},\widetilde{B})$ est r\'eunion des
$\varepsilon_i-\varepsilon_{i+1}$, $i=1,\dots,r-1$ et de
$\varepsilon_{r-1}+\varepsilon_r-\varepsilon_0$.  Soit $\varepsilon_i^* \in
\mathrm{X}_*(\widetilde{T})$ la $\Z$-base duale de la base de ${\rm
X}^\ast(\widetilde{T})$
constitu\'ee des $\varepsilon_i$, $i=0,\dots,r$.  On constate que pour
$1\leq i<j$, alors
$(\varepsilon_i-\varepsilon_j)^\vee=\varepsilon_i^\ast-\varepsilon_j^\ast$
et
$(\varepsilon_i+\varepsilon_j-\varepsilon_0)^\vee=\varepsilon_i^\ast+\varepsilon_j^\ast$. \ps\ps 

Soit $s \in {\rm O}_V(k)$ l'\'el\'ement fixant $e_i$ et $e_i^\ast$ pour
$i<r$, et \'echangeant $e_r$ et $e_r^\ast$.  La conjugaison par $s$
induit un automorphisme de $\widetilde{G}$ pr\'eservant $\widetilde{T}$ et
$\widetilde{B}$. 
Soit $\Psi(s)$ l'automorphisme induit de $\Psi(\widetilde{G})$ : il fixe
$\varepsilon_i$ pour $i=0,\cdots,r-1$ et envoie $\varepsilon_r$ sur
$\varepsilon_0-\varepsilon_r$. Si $r \neq 4$, c'est l'unique involution non triviale du
``diagramme de Dynkin'' de $\Delta(\widetilde{G},\widetilde{T},\widetilde{B})$.\ps\ps 

Consid\'erons maintenant le $k$-groupe $G={\rm SO}_V$. Sa donn\'ee
radicielle bas\'ee relative \`a $T:=\widetilde{T} \cap G$ et
$B=\widetilde{B} \cap G$ se d\'eduit de
celle de $\widetilde{G}$ par les recettes rappell\'ees plus haut, 
en consid\'erant le sous-groupe de cocaract\`eres
${\rm X}_\ast(T) = \varepsilon_0^\bot =\oplus_{i=1}^r \Z \varepsilon_i^\ast \subset
{\rm X}_\ast(\widetilde{T})$ et le groupe de caract\`eres
${\rm X}^\ast(T)={\rm X}^\ast(\widetilde{T})/\Z\varepsilon_0$. 
Autrement dit, ``on impose $\varepsilon_0=0$ dans $\Psi(\widetilde{G},\widetilde{T},\widetilde{B})$''. \ps\ps 

Soit  $\underline{\varepsilon_i}$ l'image de $\varepsilon_i$ dans ${\rm
X}^\ast(T)$, de sorte que ${\rm X}^\ast(T) = \oplus_{i=1}^r \Z \underline{\varepsilon_i}$.  On constate que
l'application lin\'eaire ${\rm X}_\ast(T) \rightarrow {\rm X}^\ast(T)$ envoyant
$\varepsilon_i^\ast$ sur $\underline{\varepsilon_i}$ induit un isomorphisme $\Psi({\rm
SO}_V) \isomo \Psi({\rm SO}_V)^\vee$; en particulier
$$\widehat{\rm{SO}}_V(\C) \simeq \rm{SO}_{2r}(\C)$$ 
\noindent o\`u ${\rm SO}_{2r}$ d\'esigne le $\C$-groupe sp\'ecial orthogonal du ${\rm
q}$-module standard sur $\C^{2r}$. \ps\ps 

De m\^eme, la donn\'ee radicielle de ${\rm P}\widetilde{G}={\rm PGSO}_V$
relative aux images respectives ${\rm P}\widetilde{T}$ et ${\rm
P}\widetilde{B}$ de $\widetilde{T}$ et $\widetilde{B}$ dans
${\rm P}\widetilde{G}$ s'obtient en consid\'erant le sous-groupe de
caract\`eres ${\rm X}^\ast({\rm P}\widetilde{T})= \zeta^\bot \subset
{\rm X}^\ast(\widetilde{T})$, o\`u $\zeta$ d\'esigne le co-caract\`ere
central $\varepsilon_0^\ast + \sum_{i=1}^r \varepsilon_i^\ast$, et le groupe de co-caract\`eres ${\rm X}_\ast({\rm
P}\widetilde{T})={\rm X}_\ast(\widetilde{T})/\Z \zeta$. Autrement dit, 
on impose $-2\varepsilon_0^\ast=\sum_{i=1}^r
\varepsilon_i^\ast$ dans la donn\'ee de $\widetilde{G}$.  Le groupe
$\widehat{{\rm PGSO}_V}(\C)$ est isomorphe au groupe des spineurs ${\rm
Spin}_{2r}(\C)$ du ${\rm q}$-module standard sur $\C^{2r}$.  \ps\ps 

\bigskip

\noindent {\sc Groupe sp\'ecial orthogonal impair} \par \medskip Soient $r\geq 1$ un entier, $U=k^r$
et soit $V$ le $k$-module ${\rm H}(U) \oplus k$ muni de la forme quadratique
somme orthogonale du ${\rm q}$-module ${\rm H}(U)$ et de $x \mapsto x^2$. 
Le $k$-groupe $G=\rm{SO}_V$ est alors semi-simple et d\'eploy\'e (\S {\rm B}.1). \ps\ps 

On
d\'efinit un tore maximal d\'eploy\'e $T$ \`a partir d'une $k$-base $(e_i)$
de $U$, un sous-groupe de Borel $B$ contenant $T$, ainsi qu'une $\Z$-base
$\varepsilon_i$ de ${\rm X}^*(T)$, comme pr\'ec\'edemment.  On constate
cette fois-ci que $\Phi(G,T)$ est r\'eunion des $\pm \varepsilon_i \pm
\varepsilon_j$ pour $1\leq i<j \leq r$ et des $\pm \varepsilon_i$ pour
$i=1,\dots,r$.  De plus, $\Delta(G,T,B)$ est r\'eunion des
$\varepsilon_i-\varepsilon_{i+1}$ pour $i<r$ et de $\varepsilon_r$.  \ps\ps 

Les groupes de similitudes et similitudes projectives associ\'es \`a $V$
sont peu diff\'erents de $G$ dans ce cadre, nous ne les consid\`ererons pas.  En revanche, le groupe des
spineurs de $V$ jouera un r\^ole.  On peut le d\'efinir suivant Chevalley au
moyen de l'alg\`ebre de Clifford de $V$.  Sur un corps alg\'ebriquement
clos, on peut se contenter de d\'ecrire sa donn\'ee radicielle bas\'ee :
c'est la donn\'ee associ\'ee au sous-groupe $Y={\rm X}^\ast(T)\,+\,\Z\,
\frac{1}{2}(\sum_{i=1}^r \varepsilon_i ) \subset {\rm X}^\ast(T) \otimes
\Q$.\ps
\bigskip
\noindent {\sc Groupe symplectique et ses variantes} \par \medskip Enfin, les $k$-groupes de la
s\'erie symplectique sont aussi r\'eductifs d\'eploy\'es.  Consid\'erons
d'abord le $k$-groupe $\widetilde{G}={\rm GSp}_{2g}$ des similitudes
symplectiques de la forme altern\'ee hyperbolique sur $U=k^g$.  \ps\ps 

On d\'efinit
$\widetilde{T}$, $\widetilde{B}$, les $\varepsilon_i$ et les
$\varepsilon_i^\ast$, pour $i=0,\cdots,g$, comme dans le cas orthogonal pair. 
L'ensemble $\Phi(\widetilde{G},\widetilde{T})$ est cette fois-ci constitu\'e
des $\pm (\varepsilon_i - \varepsilon_j)$ pour $1\leq i<j \leq g$, et des
$\pm (\varepsilon_i+\varepsilon_j-\varepsilon_0)$ pour $1\leq i \leq j \leq
g$.  De plus, $\Delta(\widetilde{G},\widetilde{T},\widetilde{B})$ est
r\'eunion des $\varepsilon_i-\varepsilon_{i+1}$ pour $1\leq i <g$ et de
$2\varepsilon_g-\varepsilon_0$.  Enfin,
$(\varepsilon_i-\varepsilon_j)^\vee=\varepsilon_i^\ast-\varepsilon_j^\ast$,
$(\varepsilon_i+\varepsilon_j-\varepsilon_0)^\vee=\varepsilon_i^\ast+\varepsilon_j^\ast$
pour $i<j$, et $(2\varepsilon_i -\varepsilon_0)^\vee=\varepsilon_i^\ast$. 
\ps\ps 

Les donn\'ees radicielles des $k$-groupes $G={\rm Sp}_{2g}$ et ${\rm
P}\widetilde{G}={\rm PGSp}_{2g}$ se d\'eduisent verbatim de celle de $\widetilde{G}$
comme dans le cas orthogonal pair. Au final, on constate $\widehat{{\rm Sp}_{2g}}(\C) \simeq {\rm
SO}_{2g+1}(\C)$ et $\widehat{{\rm PGSp}_{2g}}(\C) \simeq {\rm
Spin}_{2g+1}(\C)$.

\subsection{Repr\'esentations des groupes r\'eductifs d\'eploy\'es en
caract\'eristique nulle}
\label{repalg} Soit $k$ un corps alg\'ebriquement clos de caract\'eristique
nulle, soit $G$ un $k$-groupe r\'eductif, et soit
$\Psi(G)=(X,\Phi,\Delta,X^\vee,\Phi^\vee,\Delta^\vee)$ sa donn\'ee
radicielle bas\'ee.  Soit $$X_+=\{\lambda \in X, \langle
\lambda,\alpha^\vee\rangle \geq 0 \, \, \, \forall \alpha \in \Delta\}$$ le
sous-mono\"ide additif de $X$ constitu\'e des {\it poids dominants} de $G$. 
C'est un domaine fondamental de $X$ sous l'action du groupe de Weyl $W$ de
$G$.  \ps\ps 

On munit $X$ d'un ordre partiel pour la relation dite de {\it
dominance} : $\lambda \leq \mu \Leftrightarrow$ $\mu-\lambda$ est une somme
finie d'\'el\'ements de $\Delta$ \cite{stembridge}.  Une propri\'et\'e
notable de cette relation est que si $\lambda, \mu \in X_+$ sont tels que $\lambda < \mu$, il existe une racine
$\alpha \in \Phi$ positive relativement \`a $\Delta$ telle que $\mu - \alpha
\in X_+$ et $\lambda \leq \mu -\alpha$ \,\cite[Cor.  2.7]{stembridge}.  Par
exemple, un \'el\'ement $\lambda \in X_+$ est minimal si et seulement si
$\lambda -\alpha \notin X_+$ pour toute racine positive $\alpha \in \Phi$. 
\ps\ps 

Une $k$-repr\'esentation de $G$ est la donn\'ee d'un $k$-espace vectoriel $V$ de dimension finie 
et d'un morphisme de $k$-groupes $G \rightarrow
\GL_V$. Elles forment une cat\'egorie ab\'elienne de mani\`ere \'evidente,
qui est semi-simple car $G$ est r\'eductif. Le produit tensoriel des repr\'esentations d\'efinit une structure d'anneau commutatif ${\rm Rep}(G)$ sur le groupe de Grothendieck de cette cat\'egorie, que l'on notera ${\rm Rep}(G)$.
L'association $G \mapsto {\rm Rep}(G)$ d\'efinit de mani\`ere naturelle un foncteur de la cat\'egorie
$\mathcal{C}_k$ dans les anneaux commutatifs (\S\ref{gpdepl}). \ps\ps 

Si $\lambda \in X_+$, la th\'eorie du plus haut poids de Cartan-Weyl d\'emontre qu'il existe
une $k$-repr\'esentation
irr\'eductible $V_\lambda$ de $G$, unique \`a isomorphisme pr\`es, dont c'est le plus haut poids.
De plus, toute irr\'eductible s'obtient ainsi. Rappelons bri\`evement de quoi
il s'agit. Soit $T$ un tore maximal de $G$ et
soit $B$ un sous-groupe de Borel contenant $T$, de sorte que $\Psi(G)$ s'identifie
canoniquement \`a $\Psi(G,T,B)$. L'action de $T$ sur toute $k$-repr\'esentation
$V$ de $G$ est diagonalisable, et l'on note ${\rm Poids}(V) \subset X$ 
l'ensemble des caract\`eres de $T$ intervenant dans $V$. Il est stable pour
l'action de $W$. Si $V$ est irr\'eductible, on d\'emontre que l'espace des invariants
$V^{B(k)}$ est de dimension $1$ et que $T$ y agit par un \'el\'ement de $X_+$ : c'est le {\it plus haut poids} de $V$. 
Le plus haut poids $\lambda$ de $V$ a alors la propri\'et\'e suivante : pour tout $\mu \in {\rm Poids}(V)$ alors $\mu
\leq \lambda$. De plus, on a $${\rm Poids}(V) \cap X_+ = \{ \mu
\in X_+, \mu \leq \lambda\}$$ 
(voir par exemple \cite[\S 13.2 et \S 21.3]{humphreys}). \ps\ps

\section{Param\'etrisation de Satake}  \label{parametrisationsatake}

\subsection{L'isomorphisme de Satake}\label{isomsatake} Soit $G$ un $\Z_p$-groupe. \`A la mani\`ere du~\S\ref{chap4}.\ref{annheckeg}, notons $\mathcal{R}_p(G)$ le $G(\Q_p)$-ensemble $G(\Q_p)/G(\Z_p)$ et ${\rm H}_p(G)$ l'anneau de Hecke de $\mathcal{R}_p(G)$~(\S\ref{corrhecke}). \ps\ps 

On suppose que $G$ est r\'eductif d\'eploy\'e (\S\ref{gpdepl}). Ainsi que l'a observ\'e
Gross \cite[Prop. 1.1]{grossinv}, cette derni\`ere hypoth\`ese est satisfaite si $G$ provient par
extension des scalaires \`a $\Z_p$ d'un $\Z$-groupe r\'eductif, ce qui sera
toujours le cas dans nos applications. Soit $\widehat{G}$ le dual de Langlands de $G$, c'est-\`a-dire un
$\C$-groupe r\'eductif $\widehat{G}$ muni d'un isomorphisme
$\Psi(\widehat{G}) \isomo \Psi(G)^\vee$ (\S\ref{dualitedorad}). Son anneau de Grothendieck ${\rm Rep}(\widehat{G})$
est alors canoniquement d\'efini (\S\ref{repalg}).  L'isomorphisme de
Satake~\cite{satake}, revisit\'e par Langlands~\cite[\S 2]{langlandsyale},
est un isomorphisme d'anneaux canonique\footnote{Au sens strict, il faudrait remplacer ${\rm H}_p(G)$ par l'anneau oppos\'e ${\rm H}_p(G)^{\rm opp}$ dans l'\'enonc\'e de Satake, mais comme le r\'esultat entra\^ine la commutativit\'e de ${\rm H}_p(G)$ nous ne nous embarrasserons pas de cette d\'ecoration.}
$${\rm Sat} : {\rm H}_p(G) \otimes \Z[ \,p^{-\frac{1}{2}}\, ] \isomo {\rm Rep}(\widehat{G})\otimes  \Z[ \,p^{-\frac{1}{2}}\,]$$

Nous renvoyons \`a l'article de Satake~\cite{satake}, ainsi qu'aux
articles de survol de Cartier~\cite[\S IV]{cartier} et de
Gross~\cite{grossatake}, pour des d\'etails sur sa d\'efinition et ses propri\'et\'es
g\'en\'erales, que nous ne rediscutons que bri\`evement ci-dessous. La construction originale de Satake suppose satisfaites certaines propri\'et\'es axiomatiques du couple de groupes
$(G(\Z_p),G(\Q_p))$, qu'il v\'erifie pour les groupes classiques,
et qui ont \'et\'e d\'emontr\'ees par Tits en g\'en\'eral~\cite{titscorvallis}. Le point de vue utilis\'e ici, consistant \`a faire
appara\^itre la ``struture enti\`ere'' ${\rm Rep}(\widehat{G})$ plut\^ot que les fonctions centrales sur
$\widehat{G}$, est emprunt\'e \`a l'article de Gross sus-cit\'e. Ainsi que l'observe
Gross, on peut remplacer $\Z[\,p^{-\frac{1}{2}}\,]$ par $\Z[\,p^{-1}\,]$ dans l'isomorphisme de Satake lorsque la
demi-somme des racines positives de $G$ est un poids de $G$.\ps\bigskip

{\sc D\'efinition de l'homomorphisme de Satake :} Soient $T$
un $\Z_p$-tore maximal d\'eploy\'e de $G$, $B$ un $\Z_p$-groupe de Borel de $G$
contenant $T$ et $N$ le radical unipotent de $B$.  Si $V$
est un $G(\Q_p)$-module, le groupe ab\'elien $V_N$ des co-invariants de $V$
sous l'action de $N(\Q_p)$ est muni d'une structure de $T(\Q_p)$-module car
$T(\Q_p)$ normalise $N(\Q_p)$; cela d\'efinit un foncteur des
$G(\Q_p)$-modules vers les $T(\Q_p)$-modules appel\'e {\it foncteur de
Jacquet}.  Les d\'ecompositions ensemblistes $G(\Q_p)=B(\Q_p)G(\Z_p)$ et
$B(\Q_p) = T(\Q_p) \times N(\Q_p)$ assurent que l'inclusion \'evidente
$\mathcal{R}_p(T) \rightarrow \mathcal{R}_p(G)$ induit une bijection 
(``projection horocyclique'')
$$\mathcal{R}_p(T) \isomo N(\Q_p)\backslash \mathcal{R}_p(G).$$
 Il en r\'esulte
que si l'on pose $V=\Z[\mathcal{R}_p(G)]$ alors $V_N$ s'identifie canoniquement \`a
$\Z[\mathcal{R}_p(T)]$, d'o\`u l'on tire un homomorphisme d'anneaux $${\rm s}_1
: {\rm H}_p(G) \rightarrow {\rm H}_p(T).$$ Soit ${\rm X}_\ast(T)={\rm
Hom}(\mathbb{G}_m,T)$ le groupe des co-caract\`eres de $T$.  L'application
naturelle ${\rm X}_\ast(T) \rightarrow T(\Q_p)$, $\lambda \mapsto
\lambda(p)$, induit une bijection ${\rm X}_\ast(T) \isomo \mathcal{R}_p(T)$
puis un isomorphisme d'anneaux $\eta: \Z[{\rm X}_\ast(T)] \isomo {\rm H}_p(T)$ par
commutativit\'e de $T$. On consid\`ere alors l'homomorphisme $$s_2 : {\rm H}_p(G) \rightarrow \Z[{\rm X}_\ast(T)][\,p^{-\frac{1}{2}}\,]$$
d\'efini comme le compos\'e de $(\eta^{-1} \circ s_1) \otimes \Z[\, p^{-\frac{1}{2}}]$ par l'automorphisme de
$\Z[{\rm X}_\ast(T)][\,p^{-\frac{1}{2}}\,]$ envoyant un co-caract\`ere
$\lambda$ sur $p^{-\langle \lambda, \rho \rangle}\, \lambda $, o\`u $\rho$
est la demi-somme des racines positives de $G$ relativement \`a $(T,B)$.  Satake d\'emontre que $s_2$
est un isomorphisme sur l'anneau des invariants $\Z[{\rm
X}_\ast(T)]^W[\,p^{-\frac{1}{2}}]$, o\`u $W$ d\'esigne le groupe de Weyl de
$G$.  Mais d'apr\`es la d\'efinition du groupe dual $\widehat{G}$, l'anneau $\Z[{\rm X}_\ast(T)]^W$ s'identifie
canoniquement \`a $\Z[{\rm X}^\ast(\widehat{T})]^W$, soit encore \`a ${\rm Rep}(\widehat{G})$ d'apr\`es Chevalley : cela d\'efinit l'isomorphisme
${\rm Sat}$. $\square$ \ps \bigskip


\ps\ps 

D\'esignons par $\widehat{G}(\C)_{\rm ss}$ l'ensemble (bien d\'efini!) des
classes de conjugaison d'\'el\'ements semi-simples de $\widehat{G}(\C)$. 
Soit $c \in \widehat{G}(\C)_{\rm ss}$.  L'application $V \mapsto {\rm trace}(c
\, |\, V)$, associant \`a une $\C$-repr\'esentation de dimension finie $V$
de $\widehat{G}$ la trace de $c$ dans $V$, s'\'etend en un homomorphisme
d'anneaux ${\rm tr}(c): {\rm Rep}(\widehat{G}) \rightarrow \C$.  D'apr\`es un
r\'esultat classique d\^u \`a Chevalley, l'application ainsi d\'efinie
$${\rm tr} : \widehat{G}(\C)_{\rm ss} \rightarrow {\rm Hom}_{\rm anneaux}({\rm
Rep}(\widehat{G}),\C)$$ est une bijection.  Le scholie suivant, l'un des
points de d\'epart des travaux de Langlands, en r\'esulte imm\'ediatement.\ps\ps

\begin{scholie}\label{corsatake} L'application $c \mapsto {\rm tr}(c) \circ {\rm Sat}$ d\'efinit une bijection
$$\widehat{G}(\C)_{\rm ss} \isomo \Hom_{\rm anneaux}({\rm H}_p(G),\C).$$ 
\end{scholie}

Mentionnons enfin que l'involution $T \mapsto T^{\rm t}$ de ${\rm H}_p(G)$
correspond par l'homomorphisme de Satake \`a l'involution de ${\rm
Rep}(\widehat{G})$ induite par la dualit\'e sur les repr\'esentations, ou
encore \`a l'inversion sur $\widehat{G}(\C)_{\rm ss}$. \ps\ps

\begin{example}\label{degretriv} {\rm Repla\c{c}ons-nous tout
d'abord dans le cadre g\'en\'eral du~\S IV.\ref{corrheckegeneral}, o\`u $X$
d\'esigne un $\Gamma$-ensemble transitif quelconque.  Regardons $\Z$ comme
$\Gamma$-module pour l'action triviale.  Le ${\rm H}(X)^{\rm opp}$-module
$\Z_X$ est libre de rang $1$ sur $\Z$, et d\'efinit donc un morphisme
d'anneaux appel\'e {\it degr\'e} $${\rm deg} : {\rm H}(X) \rightarrow \Z$$ qui n'est autre que
$\deg(h)=\sum_{x \in X} h_{x,y}$, o\`u $y \in X$ est un \'el\'ement
quelconque.  Si $X=\mathcal{R}_p(G)$, on peut se demander quel est
l'\'el\'ement $s \in \widehat{G}(\C)_{\rm ss}$ correspondant \`a l'homomorphisme
$\deg$ par le scholie ci-dessus. Comme le foncteur de Jacquet du
$G(\Q_p)$-module trivial $\Z$ est le $T(\Q_p)$-module trivial $\Z$, il
d\'ecoule\footnote{Soient $V$ un $G(\Q_p)$-module et $\pi
: V^{G(\Z_p)} \rightarrow V_N$ la projection canonique. L'anneau ${\rm
H}(G)$ op\`ere sur $V^{G(\Z_p)}$ (\S IV.\ref{foncteurmx}). On observera que
par construction de $s_2$, on a $\pi \circ T = s_2(T) \circ \pi$ pour tout $T
\in {\rm H}_p(G)$. L'assertion s'en d\'eduit en consid\'erant
$V=\Z$ et en se rappelant du d\'ecalage par $\rho$ intervenant dans la
d\'efinition de l'homomorphisme de Satake.}
de la d\'efinition de l'homomorphisme de Satake rappel\'ee ci-dessus que $s$
est la classe de conjugaison de
$\rho(p)\,=(2\rho)(p^{\,\frac{1}{2}})$, o\`u $2 \rho$ est vu comme un co-caract\`ere de $\widehat{G}$. }
\end{example}

\ps \ps
\noindent {\sc Isog\'enies}
\ps\ps

Soient $G$ et $G'$ deux $\Z_p$-groupes r\'eductifs d\'eploy\'es et soit $f : G \rightarrow G'$ un morphisme central.  Ce morphisme induit d'une part un homomorphisme
d'anneaux ${\rm Rep}(f) : {\rm Rep}(\widehat{G}) \rightarrow {\rm
Rep}(\widehat{G'})$ par l'\'equivalence $\Psi$ et la dualit\'e sur les
donn\'ees radicielles.  D'autre part, Satake a d\'efini dans~\cite[\S
7]{satake} un homomorphisme d'anneaux canonique ${\rm H}(f) : {\rm H}_p(G)
\longrightarrow {\rm H}_p(G')$. \ps\ps

Lorsque $f$ est une isog\'enie centrale, 
${\rm H}(f)$ co\"incide avec l'homomorphisme ${\rm H}_p(G) \rightarrow {\rm H}_p(G')$ associ\'e par
la proposition-d\'efinition~\S\ref{defpropiso} au morphisme \'evident 
$\mathcal{R}_p(G) \rightarrow \mathcal{R}_p(G')$ d\'efini par
$f$. En effet, v\'erifions que ce dernier satisfait les hypoth\`eses du~\S\ref{parisogenies}.
D'une part, un argument galoisien direct assure que $f(G(\Q_p))$ contient le
groupe d\'eriv\'e de $G'(\Q_p)$. De plus, la d\'ecomposition de Cartan-Tits
(\S\ref{deuxbases}) montre que $G(\Z_p)$ est un sous-groupe
compact maximal de $G$, \'egal \`a $f^{-1}(G'(\Z_p))$, d'o\`u
l'injectivit\'e de $\mathcal{R}_p(G) \rightarrow \mathcal{R}_p(G')$. Mieux,
cette d\'ecomposition entra\^ine l'injectivit\'e de $G(\Z_p)\backslash G(\Q_p)/G(\Z_p)
\rightarrow G'(\Z_p)\backslash G'(\Q_p)/G'(\Z_p)$, et donc que l'action sur
${\rm H}_p(G)$ de $G'(\Z_p)$, puis du groupe $S$ {\it loc. cit.}, est
triviale. \ps\ps 

Le second th\'eor\`eme de
Satake~\cite[\S 7 Thm. 4]{satake} est la commutativit\'e du diagramme
\begin{equation}\label{diagcomsat}
\xymatrix{ {\rm H}_p(G) \ar@{->}[d]_{{\rm Sat}} \ar@{->}[rr]^{{\rm H}(f)} & & {\rm H}_p(G')
\ar@{->}[d]^{{\rm Sat}} \\
 {\rm Rep}(\widehat{G}) \otimes \Z[\,p^{-\frac{1}{2}}\,] \ar@{->}[rr]_{{\rm Rep}(f)} & & {\rm Rep}(\widehat{G}' ) \otimes
\Z[\,p^{-\frac{1}{2}}\,]}
\end{equation}
auquel on r\'ef\`erera comme la ``compatibilit\'e de l'isomorphisme de Satake aux isog\'enies''. \ps\ps 

\ps\ps\noindent
{\sc Exemple du groupe sp\'ecial orthogonal pair}
\ps\ps

Donnons un exemple d'application de la discussion pr\'ec\'edente dans le cas d'un
automorphisme de $G$.  Soit $r\geq 1$ un entier, soit $V$ le ${\rm
q}$-module hyperbolique sur $\Z_p^r$, et soit $G={\rm SO}_V$, de sorte que
$\widehat{G}$ est le $\C$-groupe ${\rm SO}_{2r}$ (\S\ref{exdorad}). Le groupe ${\rm O}(V)$
agit par $\Z_p$-automorphismes sur $G$ (par conjugaison), ainsi donc sur
$\Psi(G)$. L'homomorphisme induit ${\rm O}(V)/{\rm SO}(V) \rightarrow {\rm
Aut}_{\mathcal{D}}(\Psi(G))$ est bijectif, et 
l'\'el\'ement non-trivial est induit par l'\'el\'ement $\Psi(s)$ d\'efini loc.cit. 
Ce groupe agit \'egalement sur ${\rm Rep}(\widehat{G})$ par fonctorialit\'e,
et cette action co\"incide avec l'action naturelle de ${\rm O}_{2r}(\C)/{\rm
SO}_{2r}(\C)$ pour la m\^eme raison. D\'esignons par ${\rm H}_p({\rm O}_V)$ et
${\rm H}_p({\rm SO}_V)$ les anneaux de Hecke respectifs des $\Z_p$-groupes ${\rm
O}_V$ et ${\rm SO}_V$.  Nous avons d\'efini au~\S\ref{hsovso} un
homomorphisme canonique ${\rm H}_p({\rm O}_V) \rightarrow {\rm H}_p({\rm SO}_V)$ identifiant
${\rm H}_p({\rm O}_V)$ \`a l'anneau des invariants ${\rm H}_p({\rm SO}_V)^{{\rm O}(V)}$.  On
en d\'eduit par composition avec l'isomorphisme de Satake de ${\rm SO}_V$ un
isomorphisme canonique  \par 
\begin{equation}\label{descsato} {\rm H}_p({\rm O}_V) \otimes \Z[\, p^{-\frac{1}{2}}\,]
\isomo ({\rm Rep}({\rm SO}_{2r}(\C))\otimes \Z[\, p^{-\frac{1}{2}}\,]
)^{{\rm O}_{2r}(\C)}.\end{equation}\par \smallskip

\begin{scholie}\label{corsatakeo} L'isomorphisme de Satake de ${\rm SO}_V$ induit une bijection entre ${\rm Hom}_{\rm anneaux}({\rm H}_p({\rm O}_V),\C)$ et l'ensemble des classes de ${\rm
O}_{2r}(\C)$-conjugaison d'\'el\'ements semi-simples de ${\rm
SO}_{2r}(\C)$. 
\end{scholie}
\ps\ps

\subsection{Les deux bases naturelles de l'anneau de Hecke de
$G$}\label{deuxbases}

Soit $G$ un $\Z_p$-groupe r\'eductif d\'eploy\'e, de dual de Langlands
$\widehat{G}$. \'Ecrivons
$\Psi(\widehat{G})=(X,\Phi,\Delta,X^\vee,\Phi^\vee,\Delta^\vee)$ et d\'esignons par $X_+ \subset X$ l'ensemble
ordonn\'e des poids dominants de $\widehat{G}$
comme au~\S\ref{repalg}. Nous allons rappeler les deux $\Z$-bases naturelles
de ${\rm H}_p(G)$ et ${\rm Rep}(\widehat{G})$ index\'ees par $X_+$, et indiquer quelques liens
entre ces derni\`eres dont nous aurons besoin, suivant Gross~\cite{grossatake}. \ps\ps 

Une cons\'equence de la r\'eductivit\'e de $G$ sur $\Z_p$ est l'existence d'une
{\it d\'ecomposition de Cartan}, due \`a Tits dans cette g\'en\'eralit\'e mais classique dans nos
exemples (th\'eorie des ``diviseurs \'el\'ementaires''). Soit $T$ un tore maximal d\'eploy\'e de $G$ et $B$ un sous-groupe
de Borel de $G$ contenant $T$, ce qui identifie canoniquement
$\Psi(\widehat{G})$ \`a $\Psi(G,T,B)^\vee$, et en particulier $X$ au groupe de cocaract\`eres ${\rm
X}_\ast(T)={\rm Hom}(\mathbb{G}_m,T)$. La d\'ecomposition en question s'\'ecrit 
$$G(\Q_p)=\coprod_{\lambda \in X_+}G(\Z_p)\lambda(p)G(\Z_p).$$
Si $\lambda \in X$ on d\'esignera par ${\rm c}_\lambda \in {\rm H}_p(G)$
la fonction caract\'eristique de la double classe
$G(\Z_p)\lambda(p)G(\Z_p)$, ou selon le point de vue, des couples $(x,y)$
dans $G(\Q_p)/G(\Z_p)$ tels
que $y^{-1}x \in G(\Z_p)\lambda(p)G(\Z_p)$ (\S IV.\ref{corrheckegeneral}).
L'\'el\'ement ${\rm c}_\lambda \in {\rm H}_p(G)$ ne d\'epend pas du choix de $(T,B)$. Il est clair que ${\rm c}_{\lambda}^{\rm t}={\rm c}_{-\lambda}$ (\S\ref{corrheckegeneral}) et
${\rm c}_{w(\lambda)}={\rm c}_{\lambda}$ pour tout $\lambda \in X$ et $w
\in W$ (le groupe de Weyl de $G$). D'apr\`es la d\'ecomposition de
Cartan-Tits, les ${\rm c}_\lambda$ avec $\lambda \in X_+$ forment une $\Z$-base de ${\rm H}_p(G)$.  Si $\lambda,\mu \in X_+$, on a 
\begin{equation}\label{propclambda} 
{\rm c}_\lambda \cdot {\rm c}_\mu = {\rm c}_{\lambda+\mu} + \sum_{\nu < \lambda + \mu} n_{\lambda,\mu,\nu} {\rm c}_\nu
\end{equation}
pour certains entiers $n_{\lambda,\mu,\nu}$ \cite[(2.9)]{grossatake}. L'anneau ${\rm H}_p(G)$ admet donc une filtration \'evidente ind\'ex\'ee par le
mono\"ide ordonn\'e $X_+$, d'anneau gradu\'e associ\'e $\Z[X_+]$. En particulier, si l'on d\'esigne par 
$\Omega \subset X_+$ une famille g\'en\'eratrice de $X_+$, 
l'homomorphisme d'anneaux $\Z[\{x_\omega\}_{\omega \in \Omega}] \rightarrow {\rm H}_p(G)$ envoyant $x_\omega$ sur ${\rm c}_\omega$ est
surjectif. Si $X^\vee={\rm Q}(\Phi^\vee)$, auquel cas $X_+
\simeq \N^r$, et si $\Omega$ est la base de $X_+$ (co-poids
fondamentaux), cet homomorphisme est un isomorphisme. \ps\ps 

De m\^eme, les classes $[V_\lambda] \in {\rm Rep}(\widehat{G})$ des
repr\'esentations irr\'eductibles $V_\lambda$ pour $\lambda \in X_+$
fournissent une $\Z$-base de ${\rm Rep}(\widehat{G})$ d'apr\`es la th\'eorie 
du plus haut poids rappel\'ee au \S\ref{repalg}.  Bien qu'elles soient toutes
deux index\'ees par $X_+$, le lien entre les ${\rm Sat}({\rm c}_\lambda)$ et
les $[V_\lambda]$ est non trivial.  Nous renvoyons \`a l'article de Gross
{\it loc.  cit.} pour une discussion d\'etaill\'ee de cette question, dans
laquelle des travaux de Lustig~\cite{lusztig} jouent un r\^ole essentiel. \ps\ps 

Soient $\widehat{T} \subset \widehat{G}$ un tore maximal et $\widehat{B} \subset \widehat{G}$ un sous-groupe de Borel contenant $\widehat{T}$, de sorte que $\Psi(\widehat{G})$ s'identifie \`a $\Psi(\widehat{G},\widehat{T},\widehat{B})$. Si $V$ est une $\C$-repr\'esentation de $\widehat{G}$ et si $\mu \in X$ on note $V(\mu)
\subset V$ l'espace propre pour le caract\`ere $\mu$ sous l'action de $\widehat{T}$.\ps\ps 

\begin{prop}\label{propsatake} {\rm (Gross)} Soient $G$ un $\Z_p$-groupe semi-simple
d\'eploy\'e et $X_+$ l'ensemble ordonn\'e des poids dominants de $\widehat{G}$.  Soit $\lambda \in
X_+$.\ps\ps 
\begin{itemize} \item[(i)] Si $\lambda$ est
un \'el\'ement minimal alors $p^{\langle \lambda,\rho
\rangle} [V_\lambda] = {\rm Sat}({\rm c}_\lambda)$, o\`u $2 \rho$ est la somme des racines positives de $G$.  \ps\ps  
\item[(ii)] Si $\mu \in X_+$ et $\dim(V_\lambda(\mu))=1$,
alors $d_\lambda(\mu)=1$. \ps\ps  
\item[(iii)] {\rm (Lusztig)} Si $V_\lambda={\rm Lie}(G)$ est la repr\'esentation adjointe, alors
$d_\lambda(0)=\sum_i p^{m_i-1}$ o\`u les $m_i$ sont les exposants du groupe
de Weyl de $G$. 
\end{itemize} \end{prop}

En effet, ainsi que l'explique Gross~\cite[\S 3]{grossatake}, on
dispose pour tout $\lambda \in X_+$ d'une identit\'e de la forme
\begin{equation} \label{idsatakegen} p^{\langle \lambda,\rho \rangle} [V_\lambda] = {\rm Sat}({\rm c}_\lambda) +
\sum_{\{\mu \in {\rm P}^+, \,\mu < \lambda\} } d_\lambda(\mu)\,{\rm
Sat}({\rm c}_\mu),\end{equation} pour certains entiers $d_\lambda(\mu)$ d\'ependant de $p$. 
Le (i) en r\'esulte. Suivant Lusztig et S. Kato, Gross donne \'egalement une formule
explicite (bien que difficilement praticable) pour $d_\lambda(\mu)$ sous l'hypoth\`ese que $G$ est
adjoint, i.e. de centre
trivial. Il en d\'eduit (ii) et (iii) dans le cas adjoint (formules (4.5) et (4.6) \S 4
{\it loc.  cit.}). Pour achever la d\'emonstration de la proposition, il ne reste qu'\`a 
expliquer comment nous ramener \`a ce cas pour un $G$ semi-simple g\'en\'eral. Le lemme suivant est cons\'equence imm\'ediate des d\'efinitions
(voir~\cite[(7.4)]{satake}). \ps\ps

\begin{lemme} Soit $f : G \rightarrow G'$ un morphisme central entre
$\Z_p$-groupes r\'eductifs d\'eploy\'es, $X$ (resp. $X'$) le r\'eseau des poids de $\widehat{G}$ (resp. $\widehat{G'}$). Pour tout poids dominant $\mu \in X$,
 on a ${\rm H}(f)({\rm c}_\mu)={\rm c}_{\mu'}$ et ${\rm
Rep}(f)([V_{\mu}])=[V_{\mu'}]$, o\`u $\mu'$ est l'image de $\mu$ par l'application $\Psi(f)^\vee : X \rightarrow X'$.
\end{lemme}

Dans les notations de ce lemme, et si de plus $f$ est une isog\'enie
centrale de sorte que l'application ${\rm H}(f)$ soit injective,
l'ind\'ependance lin\'eaire des ${\rm c}_{\mu'}$ dans ${\rm H}_p(G')$ entra\^ine
donc que $d_{\lambda}(\mu)=d_{\lambda'}(\mu')$ pour tout $\lambda,\mu \in
X_+$.  On conclut la d\'emonstration de la proposition en consid\'erant
l'isog\'enie canonique $G \rightarrow G/{\rm Z}(G)$.  Cet argument montre
\'egalement que la formule de Kato et Lusztig susmentionn\'ee vaut pour tout
$\Z_p$-groupe semi-simple $G$, par r\'eduction au cas adjoint.  $\square$\ps\ps


\subsection{Cas des groupes classiques : un formulaire} \label{TPSO}\label{exemplehecke} 
\par \medskip
\noindent {\sc Groupe orthogonal pair et ses variantes} \par \medskip Soient $r\geq 2$ un entier et
$L$ le ${\rm q}$-module hyperbolique sur $\Z_p^r$. On dispose d'un carr\'e
commutatif d'injections naturelles (\S IV.\ref{hsovso}, Exemple~IV.\ref{exampleiso}, \S\ref{isomsatake})
$$\xymatrix{ {\rm H}_p({\rm SO}_L) \ar@{^{(}->}[r] & {\rm H}_p({\rm PGSO}_L) \\
\mathrm{H}_p({\rm O}_L) \ar@{^{(}->}[r]\ar@{^{(}->}[u]& {\rm H}_p({\rm PGO}_L)
\ar@{^{(}->}[u]}$$

L'injection du haut commute aux actions naturelles du groupe \`a deux
\'el\'ements ${\rm O}_L(\Z_p)/{\rm SO}_L(\Z_p)$, celle du bas \'etant alors
l'injection qui s'en d\'eduit sur les invariants.  Pour all\'eger les
notations ces injections seront vues comme des inclusions.  Il sera utile de
d\'ecrire d'abord ${\rm H}_p({\rm PGSO}_L)$.  \ps\ps

Nous reprenons les notations du~\S\ref{exdorad} relatives au $\Z_p$-groupes
${\rm GSO}_L$, ${\rm SO}_L$ et ${\rm PGSO}_L$.  Soit $\lambda \in {\rm
X}_\ast({\rm P}\widetilde{T})$.  Il admet un unique ant\'ec\'edent par
l'application canonique ${\rm X}_\ast(\widetilde{T}) \rightarrow {\rm
X}_\ast({\rm P}\widetilde{T})$, disons $$\widetilde{\lambda}= \sum_{i=0}^r
m_i \varepsilon_i^\ast,$$ tel que $m_0= \langle \varepsilon_0,
\widetilde{\lambda} \rangle \in \{0,1\}$.  Le r\'eseau homo-dual
$M=\widetilde{\lambda}(p)L$ a pour $\Z_p$-base l'ensemble des $p^{m_i}e_i$
et des $p^{m_0-m_i}e_i^\ast$ pour $1 \leq i \leq r$.  Il v\'erifie donc
$$M/M\cap L \simeq \prod_{i=1}^r (\Z/p^{d_i}\Z),$$ o\`u $d_i={\rm
Max}(m_i-m_0,-m_i)=|m_i-\frac{m_0}{2}|-\frac{m_0}{2}$ si $i \in \{ 1,\cdots
r\}$, et $M^\sharp = p^{-m_0} M$.  Notons $A_\lambda \in {\rm AF}$ la classe
d'isomorphisme du groupe ab\'elien ci-dessus et posons $v_\lambda = \langle
\varepsilon_0,\widetilde{\lambda}\rangle=v_\lambda=m_0$.  L'application
$$\eta : \lambda \mapsto (A_\lambda,v_\lambda)$$ induit manifestement une
surjection de ${\rm X}_\ast({\rm P}\widetilde{T})$ sur l'ensemble des
couples $(A,v)$ o\`u $A$ est un $p$-groupe ab\'elien fini engendr\'e par $r$
\'el\'ements, et $v \in \{0,1\}$.  Il n'est pas difficile de v\'erifier que
$\eta$ est constante sur les orbites du sous-groupe de ${\rm Aut}({\rm
X}_\ast({\rm P}\widetilde{T}))$ engendr\'e par $W$ et l'automorphisme
$\tau:=\Psi(s)$ introduit au~\S\ref{exdorad}.  De plus, le co-poids
$\lambda$ est dominant si et seulement si $m_1 \geq m_2 \geq \cdots \geq
m_{r-1} \geq d_r$.  Deux co-poids dominants $\lambda,\lambda'$ ont donc
m\^eme image par $\eta$ si, et seulement si, $\lambda' \in
\{\lambda,\tau(\lambda)\}$.  Si l'on compare ces consid\'erations avec
celles du \S~IV.\ref{annheckeclass}, on en d\'eduit le :

\begin{scholie} Si $\lambda$ est un co-poids de ${\rm PGSO}_L$, on a l'\'egalit\'e 
${\rm T}_{(A_\lambda,v_\lambda)}=\sum_{\mu \in \{\lambda,\tau(\lambda)\}} {\rm c}_\mu$.
De plus, on a $T^{\rm t}=T$ pour tout $T \in {\rm H}_p({\rm PGO}_L)$.
\end{scholie}

L'injection naturelle ${\rm X}_\ast(T) \rightarrow {\rm X}_\ast({\rm P}\widetilde{T})$
identifie les co-poids de ${\rm SO}_L$ avec ceux de ${\rm PGSO}_L$ tels
que $v_\lambda=0$. Si $\lambda$ est un co-poids de ${\rm SO}_L$, on en d\'eduit
l'\'egalit\'e ${\rm T}_{A_\lambda}=\sum_{\mu \in \{\lambda,\tau(\lambda)\}}
{\rm c}_\mu$ dans ${\rm H}_p({\rm O}_L)$. \ps\ps  

\'Etant donn\'e que ${\rm H}_p({\rm PGO}_L)$ s'identifie aux invariants de
${\rm H}_p({\rm PGSO}_L)$ sous l'action de conjugaison de ${\rm O}_L(\Q_p)/{\rm
SO}_L(\Q_p)$, la proposition~IV.\ref{baseheckego} est cons\'equence du fait que les ${\rm c}_\lambda$, avec $\lambda$
dominant, forment une $\Z$-base de ${\rm H}_p({\rm PGSO}_L)$. Si $1 \leq i \leq r$,
notons $\lambda_i \in {\rm  X}_\ast({\rm P}\widetilde{T})$ l'image de $\varepsilon_1^\ast + \cdots +
\varepsilon_i^\ast \in {\rm X}_\ast(\widetilde{T})$, 
et notons \'egalement $\lambda_{r+1}\in {\rm
X}_\ast({\rm P}\widetilde{T})$ l'image de $-\varepsilon_0^\ast$. En
particulier, on a $2\lambda_{r+1} = \lambda_r$. D'apr\`es le scholie, on a les
relation ${\rm c}_{\lambda_i}={\rm T}_{(\Z/p\Z)^i}$ dans ${\rm H}_p({\rm O}_L)$ pour $i<r$,
${\rm c}_{\lambda_r}+{\rm c}_{\tau(\lambda_r)}={\rm T}_{(\Z/p\Z)^r}$, et 
${\rm c}_{\lambda_{r+1}}+{\rm c}_{\tau(\lambda_{r+1})}={\rm K}_p$ dans ${\rm H}_p({\rm PGO}_L)$.
L'\'enonc\'e suivant est bien connu \cite{satake} \cite[\S 4]{rallis1}.\ps\ps

\begin{cor}\label{corgenhecke}\begin{itemize}\item[(i)] L'homomorphisme $\Z[X_1,\cdots,X_r] \rightarrow {\rm H}_p({\rm PGO}_L)$ 
envoyant $X_i$ sur ${\rm T}_{(\Z/p\Z)^i}$ pour $1\leq i \leq r-1$, et $X_r$ sur ${\rm K}_p$, est un isomorphisme d'anneaux.\ps\ps
\item[(ii)] L'homomorphisme $\Z[Y_1,\cdots,Y_r] \rightarrow {\rm H}_p({\rm O}_L)$ envoyant $Y_i$ sur 
${\rm T}_{(\Z/p\Z)^i}$ pour $1\leq i \leq r$, est un isomorphisme d'anneaux.\ps\ps
\end{itemize}
\end{cor} 

\begin{pf} Le groupe ${\rm PGSO}_L$ \'etant adjoint, la discussion du~\S\ref{deuxbases}
entra\^ine que ${\rm H}_p({\rm PGSO}_L)$ est l'anneau de polyn\^omes sur les
${\rm c}_{\omega}$, o\`u $\omega$ parcourt les co-poids fondamentaux de ${\rm PGSO}_L$.  Ce sont les \'el\'ements
$\lambda_{r+1}, \tau(\lambda_{r+1})$ et les $\lambda_i$ pour
$i=1,\cdots,r-2$ d'apr\`es \cite[Planche IV]{bourbaki}.  Les $r-2$ derniers
sont invariants par $\tau$, et les deux premiers sont \'echang\'es. 
Rappelons que si $A$ est un anneau commutatif, le sous-anneau de $A[U,V]$
constitu\'e des $P(U,V)$ tels que $P(U,V)=P(V,U)$ est $A[UV,U+V]$. Ainsi, ${\rm H}_p({\rm PGO}_L)$ est l'anneau de polyn\^omes sur les
${\rm c}_{\lambda_i}$, $1\leq i<r-1$, ${\rm K}_p$ et
${\rm c}_{\lambda_{r+1}}{\rm c}_{\tau(\lambda_{r+1})}$. Mais les seuls co-poids
dominants de ${\rm PGSO}_L$ qui sont strictement inf\'erieurs \`a
$\lambda_{r+1}+\tau(\lambda_{r+1})=\lambda_{r-1}$ sont les
$\lambda_i$ avec $0 \leq i<r-1$ et $i \equiv r-1 \bmod 2$, avec la
convention $\lambda_0=0$. Il existe donc des entiers $a_j \in \Z$ et une identit\'e de la forme $$
{\rm c}_{\lambda_{r+1}}{\rm c}_{\tau(\lambda_{r+1})} = {\rm c}_{\lambda_{r-1}} +
\sum_{0 \leq j < r-1} a_j {\rm c}_{\lambda_j}.$$ 

Cela d\'emontre le (i). Le (ii) se d\'emontre par des arguments similaires. Le mono\"ide des
co-poids dominants de ${\rm SO}_L$ est engendr\'e par les $\lambda_i$ pour
$i \leq r-1$, $\lambda_r$ et $\tau(\lambda_r)$. Le sous-anneau 
${\rm H}_p({\rm SO}_L) \subset {\rm H}_p({\rm PGSO}_L)$ est donc engendr\'e par
les ${\rm c}_\lambda$ o\`u $\lambda$ parcourt cette liste. Mais si 
$S=\{\lambda_{r+1},\tau(\lambda_{r+1})\}$ et si
$s,t \in S$, on observe comme plus haut que ${\rm c}_s{\rm c}_t-{\rm c}_{s+t}$ est une combinaison lin\'eaire \`a
coefficients entiers des ${\rm c}_{\lambda_i}$ pour $0 \leq i \leq r-2$. Cela
entra\^ine que ${\rm H}_p({\rm SO}_L)$ est \'egalement engendr\'e par l'anneau
$\Z[{\rm c}_{\lambda_1},\cdots,{\rm c}_{\lambda_{r-2}}]$ et les trois \'el\'ements
${\rm c}_s{\rm c}_t$ o\`u $s,t \in S$. On conclut en observant que si $A$ est un anneau commutatif, le sous-anneau de
$A[U^2,V^2,UV] \subset A[U,V]$ constitu\'e des polyn\^omes sym\'etriques en
$U$ et $V$ est $A[UV,(U+V)^2]$. \ps\ps 
\end{pf}

D'apr\`es le scholie ci-dessus, pour $m\geq 0$ l'op\'erateur ${\rm T}_{p^m} \in {\rm H}_p({\rm O}_L)$ des
$p^m$-voisins co\"incide avec ${\rm c}_{m \lambda_1}$. L'op\'erateur
${\rm T}_{(\Z/p\Z)^2}$ interviendra \`a plusieurs reprises par la suite, et sera
\'egalement not\'e ${\rm T}_{p,p}$. 

\begin{example}\label{calcultp2} {\rm Dans ${\rm H}_p({\rm O}_L)$, on a la relation 
$$({\rm T}_p)^2 = {\rm T}_{p^2} + (p+1)\,\, {\rm T}_{p,p} +
\frac{(p^r-1) (p^{r-1}+1)}{p-1}.$$}
\end{example}

\begin{pf}  Les co-poids dominants de ${\rm SO}_L$ strictement inf\'erieurs \`a $2
\lambda_1$
sont $\lambda_2$ et $0$, d'o\`u l'existence de
$a,b \in \Z$ tels que $({\rm T}_p^2)={\rm T}_{p^2}+ a\, {\rm T}_{p,p} + b$ (formule~\eqref{propclambda}).
Comme $L$ est un $p$-voisin de chacun de ses $p$-voisins (III. D\'ef. 1.2), l'entier $b$
est simplement le nombre $b={\rm c}_{2r}(p)$ de $p$-voisins d'un r\'eseau autodual
(III. D\'ef. 2.1). Calculons $a$ \`a l'aide de l'homomorphisme {\it degr\'e}, $\deg : {\rm
H}({\rm SO}_L) \rightarrow \Z$, introduit dans l'exemple~\ref{degretriv}. \ps Le degr\'e de ${\rm T}_A$ est le nombre de $A$-voisins de $L$, en particulier 
$\deg({\rm T}_p)={\rm c}_{2r}(p)$ et $\deg({\rm T}_{p^2})=p^{2r-2}{\rm c}_{2r}(p)$ d'apr\`es la
proposition III. 1.4. Dans le m\^eme esprit que le \S III. 1, on
v\'erifierait ais\'ement que si $1\leq i\leq r$, le nombre de
$(\Z/p\Z)^i$-voisins de $L$ est le produit du
nombre de sous-espaces isotropes de rang $i$ de $L \otimes \F_p$ par le
nombre de lagrangiens du ${\rm q}$-module hyperbolique sur $(\Z/p\Z)^i$ qui sont transverses \`a $(\Z/p\Z)^i$
(i.e. $p^{\frac{i(i-1)}{2}}$ d'apr\`es la proposition-d\'efinition I. 1.3 (b)).  Pour $i=2$, on en trouve donc $
\frac{{\rm c}_{2r}(p) {\rm c}_{2r-2}(p)}{(p+1)}\cdot p$.  Un petit calcul conduit \`a
$a=p+1$.
\end{pf}

\begin{remark} {\rm Il serait int\'eressant de savoir si les ${\rm T}_{p^i}$, $i=1,\cdots,r$,
engendrent
la $\Q$-alg\`ebre ${\rm H}_p({\rm O}_L) \otimes \Q$.} \end{remark} Terminons ce formulaire sur les groupes orthogonaux pairs par certaines
propri\'et\'es de l'isomorphisme de Satake. La demi-somme des racines positives
de ${\rm GSO}_L$ est $\rho = (r-1) \varepsilon_1 + (r-2) \varepsilon_2 + \cdots +
\varepsilon_{r-1} - \frac{r(r+1)}{4} \varepsilon_0$. 
Le seul co-poids dominant minimal de ${\rm SO}_L$ est $\lambda_1$,
et ${\rm PGSO}_L$ en admet deux autres qui sont $\lambda_{r+1}$ et
$\tau(\lambda_{r+1})$. Le premier est le poids dominant
de la repr\'esentation standard $V_{\rm St}$ (de dimension $2r$) de $\widehat{{\rm SO}_L}(\C)={\rm SO}_{2r}(\C)$,
et le deux autres sont les poids dominants des deux repr\'esentations {\rm Spin} 
$V^{\pm}_{\rm Spin}$ de $\widehat{{\rm PGSO}_L}(\C) = {\rm Spin}_{2r}(\C)$. La
proposition~\ref{propsatake} (i) et le scholie entra\^inent les identit\'es : 
\begin{equation}\label{formulesattp} p^{r-1}[V_{\rm St}]={\rm Sat}({\rm T}_p)\, \, \, \, \, {\rm et}\, \, \, \, \,
p^{\frac{r(r-1)}{4}}([V_{\rm Spin}^+]+[V_{\rm Spin}^-])={\rm Sat}({\rm K}_p).\end{equation}

Consid\'erons maintenant la repr\'esentation $\Lambda^2 V_{\rm St}$, qui
n'est autre que la repr\'esentation adjointe de ${\rm SO}_{2r}(\C)$. Son
plus haut poids est $\lambda_2$; l'unique poids
dominant inf\'erieur \`a ce dernier est le poids $0$. Les points (i) et
(iii) de la proposition~\ref{propsatake} entra\^inent donc
\begin{equation}\label{formulesattpp} p^{2r-3}[\Lambda^2 V_{\rm St}] = {\rm Sat}({\rm T}_{p,p}) +
p^{r-2}+\sum_{i=0}^{r-2} p^{2i}.\end{equation}
(Nous pourrions \'egalement invoquer l'exemple~\ref{degretriv} \`a la place
du (iii) de la proposition~\ref{propsatake}.)


\par \medskip
\noindent {\sc Groupe symplectique et ses variantes} \par \medskip Soient $g\geq 1$ un entier et
$L$ le ${\rm a}$-module hyperbolique sur $\Z_p^g$. Nous reprenons les
notations du~\S\ref{exdorad} relatives au $\Z_p$-groupes ${\rm GSp}_{2g}$,
${\rm Sp}_{2g}$ et ${\rm PGSp}_{2g}$. On constate de m\^eme que plus haut 
que si $A$ est un $p$-groupe ab\'elien fini engendr\'e par $g$ \'el\'ements,
disons $A \simeq \prod_{i=1}^g \Z/p^{m_i}\Z$ avec $m_1 \geq \cdots \geq m_g
\geq 0$, alors 
\begin{equation}\label{tradsatsymp}
{\rm T}_A = {\rm c}_{\sum_{i=1}^g m_i \varepsilon_i^\ast}, \, \, \,{\rm et}\, \, \, 
{\rm T}_{(A,1)}={\rm c}_{\varepsilon_0^\ast+\sum_{i=1}^g  (m_i+1) \varepsilon_i^\ast}.\end{equation}
D'apr\`es Shimura \cite{shimurasp}, l'anneau ${\rm H}_p({\rm PGSp}_{2g})$ est l'anneau de polyn\^omes sur ${\rm K}_p$ et les ${\rm
T}_{(\Z/p\Z)^i}$ pour $i<g$, le sous-anneau ${\rm H}_p({\rm Sp}_{2g})$ \'etant
engendr\'e par les ${\rm T}_{(\Z/p\Z)^i}$ pour $i \leq g$ (la situation est
en fait plus simple que celle du corollaire~\ref{corgenhecke}, les
mono\"ides des co-poids dominants de ${\rm Sp}_{2g}$ et ${\rm PGSp}_{2g}$ \'etant libres). 
\ps\ps

La demi-somme des racines positives $\rho$ de ${\rm GSp}_{2g}$ est $-
\frac{g(g+1)}{4}\, \varepsilon_0^\ast + g \, \varepsilon_1^\ast + \break (g-1)
\, \varepsilon_2^\ast + \cdots + \varepsilon_g^\ast$. Soit $V_{\rm
St}$ la repr\'esentation standard (de dimension $2g+1$) de
$\widehat{{\rm Sp}_{2g}}(\C)={\rm SO}_{2g+1}(\C)$, et soit $V_{\rm Spin}$
la repr\'esentation Spin de ${\rm Spin}_{2g+1}(\C)$. Le plus haut poids de
$V_{\rm St}$ est $\varepsilon_1^\ast$, qui a pour unique poids dominant
strictement inf\'erieur $0$. Le plus haut poids de $V_{\rm Spin}$ est
$-\varepsilon_0^\ast$ qui est minimal. La proposition~\ref{propsatake}
entra\^ine donc
\begin{equation}\label{satakesp}p^g [V_{\rm St}] = {\rm Sat}({\rm T}_p)+1, \, \, \, {\rm et}\, \, \,
p^{\frac{g(g+1)}{4}} [V_{\rm Spin}] = {\rm Sat}({\rm K}_p).\end{equation}

\par \medskip
\noindent {\sc Groupe sp\'ecial orthogonal impair}\par \medskip Nous n'insisterons que sur les diff\'erences avec les autres cas, qui sont minimes. Soient $r\geq 1$ un entier, 
$L$ le $\Z_p$-module $\Z_p^r \oplus (\Z_p^r)^\ast \oplus \Z_p$ muni de la forme quadratique
somme orthogonale du ${\rm q}$-module hyperbolique sur $\Z_p^r$ et de $x \mapsto x^2$, $V=L \otimes \Q_p$ et $G$ le $\Z_p$-groupe ${\rm SO}_L$ (\S ${\rm B}$.1). C'est un exercice de v\'erifier que l'application $g \mapsto g(L)$ identifie $G(\Q_p)/G(\Z_p)$ au sous-ensemble $\mathcal{R}_{\Z_p}^{\rm b}(V) \subset \mathcal{R}_{\Z_p}(V)$ des r\'eseaux $M \subset V$ tels que $M^\sharp = M$ si $p>2$, ou tels que $M$ est un sous-r\'eseau d'indice $2$ de $M^\sharp$ si $p=2$ (\S IV.\ref{annheckeclass}). \ps\ps 

Soit $V_{\rm St}$ la repr\'esentation standard de $\widehat{{\rm
SO}_L}(\C)={\rm Sp}_{2r}(\C)$ (de dimension $2r$); son plus haut poids est
$\varepsilon_1^\ast$ qui est minimal (nous reprenons les notations du~\S \ref{exdorad} relatives au $\Z_p$-groupe $G$). On constate que l'op\'erateur de Hecke ${\rm c}_{\varepsilon_1^\ast}$ est associ\'e aux couples $(N,M) \in \mathcal{R}_{\Z_p}^{\rm b}(V)^2$ tels que $M \cap N$ est d'indice $p$ dans $M$ : c'est l'op\'erateur ${\rm T}_p$ des $p$-voisins au sens de l'appendice ${\rm B}$.3. La demi-somme des racines positives de
${\rm SO}_L$ \'etant $\frac{2r-1}{2} \, \varepsilon_1 + \frac{2r-3}{2}
\, \varepsilon_2 + \cdots + \frac{1}{2}\, \varepsilon_r$, on a donc
$$p^{\frac{2r-1}{2}}[V_{\rm St}]={\rm T}_p.$$

\section{L'isomorphisme de Harish-Chandra} \label{secisohc}

\subsection{Le centre de l'alg\`ebre enveloppante d'un $\C$-groupe
r\'eductif} Soient $G$ un $\C$-groupe r\'eductif, $\mathfrak{g}$ la
$\C$-alg\`ebre de Lie de $G$, ${\rm U}(\mathfrak{g})$ son alg\`ebre
enveloppante, et ${\rm Z}({\rm U}(\mathfrak{g}))$ le centre de ${\rm
U}(\mathfrak{g})$~\cite[Ch.  2]{dixmier}.  Soit $V$ un ${\rm
U}(\mathfrak{g})$-module. On dit que $V$ admet un caract\`ere central s'il
existe un homomorphisme de $\C$-alg\`ebres $${\rm c}_V : {\rm Z}({\rm
U}(\mathfrak{g})) \rightarrow \C$$ tel que $z \cdot v = {\rm c}_V(z) v$ pour
tout $v \in V$ et tout $z \in{\rm Z}({\rm U}(\mathfrak{g}))$; on appelle alors
${\rm c}_V$ le {\it caract\`ere central} de $V$. D'apr\`es \cite[Prop. 
2.6.8]{dixmier}, tout ${\rm U}(\mathfrak{g})$-module simple admet un
caract\`ere central. Suivant Harish-Chandra et Langlands, nous allons rappeler dans ce
paragraphe comment voir ces caract\`eres centraux comme des classes de
conjugaison semi-simples dans la $\C$-alg\`ebre de Lie
$\widehat{\mathfrak{g}}$ du $\C$-groupe r\'eductif dual $\widehat{G}$ de
$G$.  \ps\ps 

Soit ${\rm Pol}(\widehat{\mathfrak{g}})={\rm
Sym}({\,\widehat{\mathfrak{g}}\,}^\ast)$ la $\C$-alg\`ebre des fonctions
polynomiales sur $\widehat{\mathfrak{g}}$.  Elle est munie d'une action
naturelle de $\widehat{G}(\C)$ issue de l'action adjointe sur
$\widehat{\mathfrak{g}}$, dont nous notons ${\rm
Pol}(\widehat{\mathfrak{g}})^{\widehat{G}}$ l'alg\`ebre des invariants. 
L'isomorphisme de Harish-Chandra est un isomorphisme canonique $${\rm HC}
:{\rm Z}({\rm U}(\mathfrak{g}))\overset{\sim}{\longrightarrow} {\rm
Pol}(\widehat{\mathfrak{g}})^{\widehat{G}}$$ \cite[Thm.  7.4.5 \&
7.3.5]{dixmier} \cite[\S 2]{langlandsyale}. Soit
$\widehat{\mathfrak{g}}_{\rm ss}$ l'ensemble des classes de conjugaison
d'\'el\'ements semi-simples de $\widehat{\mathfrak{g}}$.  Chaque telle classe
$X \in \widehat{\mathfrak{g}}_{\rm ss}$ d\'efinit par \'evaluation un
homomorphisme de $\C$-alg\`ebres ${\rm
Pol}(\widehat{\mathfrak{g}})^{\widehat{G}} \rightarrow \C, \, P \mapsto
P(X)$.  Un r\'esultat classique de Chevalley assure que l'application
$\widehat{\mathfrak{g}}_{\rm ss} \rightarrow {\rm Hom}_{\C-{\rm alg}}({\rm
Pol}(\widehat{\mathfrak{g}})^{\widehat{G}},\C)$ ainsi d\'efinie est
bijective.  \ps\ps

\begin{scholie}\label{corhc} L'isomorphisme de Harish-Chandra
induit une bijection canonique ${\rm Hom}_{\C-{\rm alg}}({\rm
Z}({\rm U}(\mathfrak{g})),\C) \isomo \widehat{\mathfrak{g}}_{\rm ss}$.  \end{scholie}\ps\ps 

Si $X$ est le r\'eseau des poids de $G$, les \'el\'ements de $X \otimes \C$ peuvent \^etre vus comme des \'el\'ements de $\widehat{\mathfrak{g}_{\rm ss}}$. En effet, soient $\widehat{T}$ un tore maximal de $\widehat{G}$ et $\widehat{B} \subset \widehat{G}$ un sous-groupe de Borel contenant $\widehat{T}$. La donn\'ee $\Psi(G)^\vee$ s'identifie \`a $\Psi(\widehat{G},\widehat{T},\widehat{B})$; en particulier $X$ s'identifie \`a ${\rm X}_\ast(\widehat{T})$. L'application exponentielle d\'efinit un isomorphisme naturel entre ${\rm X} \otimes \C $ et l'alg\`ebre de Lie complexe $\widehat{\mathfrak{t}}$ de $\widehat{T}$. Si $W$ est le groupe de Weyl de $G$, on en d\'eduit un bijection canonique
\begin{equation}\label{poidscarinf} (X \otimes \C)/W \overset{\sim}{\longrightarrow} \widehat{\mathfrak{g}}_{\rm ss}.\end{equation}

\begin{example}\label{exemplecentcar} {\rm Soient $\lambda \in X_+$ un poids
dominant de $G$ et $V_\lambda$ la $\C$-re\-pr\'esentation irr\'eductible de
$G$ de plus haut poids $\lambda$ (\S\ref{repalg}).  Cette repr\'esentation
munit $V_\lambda$ d'une structure de ${\rm U}(\mathfrak{g})$-module.  Ce
module est simple et son caract\`ere central correspond \`a la classe de
conjugaison de $\lambda + \rho$, o\`u $\rho$ est la demi-somme des racines
positives de $G$.  \par Plus g\'en\'eralement, fixons une paire $T \subset
B$ dans $G$ identifiant $\Psi(G)$ \`a $\Psi(G,T,B)$. Soient $\mathfrak{t} \subset
\mathfrak{b}$ leurs $\C$-alg\`ebres de Lie respectives, et $V$ un ${\rm U}(\mathfrak{g})$-module engendr\'e par un \'el\'ement $e \in V$ tel que $\mathfrak{b} \, e \subset
\C e$ ({\it module de plus haut poids}).  Soit $\lambda \in
\mathfrak{t}^\ast$ la forme lin\'eaire d\'efinie par $h e = \lambda(h) e $
pour tout $h \in \mathfrak{t}$ (on peut \'egalement la voir de mani\`ere duale comme un
\'el\'ement de $\widehat{\mathfrak{t}}$).  Alors $V$ admet un caract\`ere
central d'apr\`es~\cite[Prop.  7.1.8]{dixmier}, et il r\'esulte assez
directement de la d\'efinition de l'homomorphisme de Harish-Chandra que la
classe de conjugaison qui lui correspond est celle de $\lambda + \rho$ \cite[\S
7.4.6]{dixmier}. 
 } \end{example}

\subsection{Caract\`ere infinit\'esimal d'une repr\'esentation unitaire}
\label{carinf}

Soit $G$ un $\R$-groupe r\'eductif. On applique les consid\'erations et
notations du paragraphe pr\'ec\'edent au $\C$-groupe $G_\C:=G \times_\R \C$.
Nous renvoyons \`a \cite{knapp} et \cite{wallach} pour un expos\'e g\'en\'eral
de la th\'eorie des repr\'esentations unitaires des groupes de Lie r\'eductifs. \ps\ps 

Soit $V$ un espace de Hilbert muni d'une repr\'esentation unitaire du groupe
de Lie $G(\R)$.  Soit $V^\infty \subset V$ le sous-espace des vecteurs $\mathcal{C}^\infty$, c'est-\`a-dire des $v
\in V$ tels que l'application $g \mapsto g v, G(\R) \rightarrow V$, soit de
classe $\mathcal{C}^\infty$; il est dense dans $V$ (G\r{a}rding) et stable par $G(\R)$.  Si la
repr\'esentation unitaire $V$ est irr\'eductible alors le ${\rm U}(\mathfrak{g})$-module $V^\infty$ admet
un caract\`ere central~\cite[\S 1.6.5]{wallach}, appel\'e {\it caract\`ere infinit\'esimal}
de $V$, nous le notons ${\rm inf}_V$. Ainsi que l'a d\'emontr\'e Harish-Chandra, c'est un invariant assez fin de la repr\'esentation $V$ : \`a 
isomorphisme pr\`es il n'y a qu'un nombre fini (\'eventuellement nul) de repr\'esentations 
unitaires irr\'eductibles de $G$ ayant un caract\`ere infinit\'esimal donn\'e 
(c'est un r\'esultat difficile, voir \cite[Cor. 10.37]{knapp}). L'isomorphisme d'Harish-Chandra
permet de voir ${\rm inf}_V$ comme une classe de conjugaison semi-simple dans
$\widehat{\mathfrak{g}}$.  Nous allons donner deux exemples.  \ps\ps 

Supposons d'abord que $G(\R)$ est un groupe compact, auquel cas il est
n\'ecessairement connexe d'apr\`es Chevalley \cite[V.  24.6 (c)
(ii)]{borelgp}.  Toute $\C$-repr\'esentation $V$ de $G_\C$ d\'efinit par
restriction une repr\'esentation $V_{|G(\R)}$ de $G(\R)$.  Le foncteur $V
\mapsto V_{|G(\R)}$ est une \'equivalence de cat\'egories entre
$\C$-repr\'esentations de $G_\C$ et repr\'esentations complexes, continues,
de dimension finie, de $G(\R)$.  En particulier, toute repr\'esentation
irr\'eductible de $G(\R)$ est isomorphe \`a $(V_\lambda)_{|G(\R)}$ pour un
unique poids dominant $\lambda$ de $G_\C$; on la notera en g\'en\'eral
$V_\lambda$ pour all\`eger les notations.  D'apr\`es
l'exemple~\ref{exemplecentcar}, son caract\`ere infinit\'esimal est la
classe de conjugaison de $\lambda+\rho$ dans $\widehat{\mathfrak{g}}$. Il d\'etermine en particulier uniquement $V_\lambda$.  \ps\ps 

Supposons maintenant que $G$ est le $\R$-groupe ${\rm Sp}_{2g}$. Reprenons
quelques notations du~\S\ref{fautsiegel} \`a ceci pr\`es que nous notons ici
$\mathfrak{g}$ l'alg\`ebre de Lie {\it complexifi\'ee} de $G(\R)$.  On
choisit pour sous-groupe compact maximal $K \subset G(\R)$ le fixateur de $i
1_g$ dans $\mathbb{H}_g$, d'alg\`ebre de Lie $\mathfrak{k}$.  C'est un
groupe unitaire \`a $g$ variables : l'homomorphisme $k \mapsto {\rm j}(k,i 1_g)$,
$K \rightarrow \GL_g(\C)$, identifie $\GL_g(\C)$ au complexifi\'e de $K$,
puis (par diff\'erentiation) $\mathfrak{k}_\C$ \`a l'alg\`ebre de Lie $\mathfrak{gl}_g(\C)$.  La d\'ecomposition de Cartan complexifi\'ee
s'\'ecrit $$\mathfrak{g}=\mathfrak{k}_\C \oplus \mathfrak{p}^+ \oplus
\mathfrak{p}^-,$$ o\`u $\mathfrak{p}^\pm$ sont des sous-alg\`ebres de Lie
{\it ab\'eliennes} stables par ${\rm ad}(K)$.  Le point clef est que
l'action adjointe sur $\mathfrak{p}$ de l'\'el\'ement
$\frac{1}{\sqrt{2}}\left(\begin{array}{cc}1 & 1 \\ -1 & 1\end{array}\right)$
du centre de $K$ (qui est isomorphe \`a ${\rm U}(1)$) induit la structure complexe
naturelle du $\R$-espace vectoriel $\mathfrak{p} \simeq {\rm Sym}_g(\C)$. 
\ps\ps 

Soit $T$ un tore maximal de $\GL_g$, $B$ un sous-groupe de Borel contenant
$T$, d'alg\`ebres de Lie complexes respectives $\mathfrak{t} \subset
\mathfrak{b} \subset \mathfrak{gl}_g(\C)=\mathfrak{k}_\C$.  Les propri\'et\'es sus-mentionn\'ees
de $\mathfrak{p}^-$ assurent que $\mathfrak{t}$ est une sous-alg\`ebre de
Cartan de $\mathfrak{g}$ et que $\mathfrak{b} \oplus \mathfrak{p}^-$ en est
une sous-alg\`ebre de Borel.

\begin{prop}\label{infcarholo} {\rm (Harish-Chandra)} Soient $V$ une
repr\'esentation unitaire de ${\rm Sp}_{2g}(\R)$, $e$ un \'el\'ement de $V^\infty$ et $U$ une
$\C$-repr\'esentation irr\'eductible de $\GL_g$ tels que : \ps\ps 
\begin{itemize} \item[(i)] $\mathfrak{p}^- e = 0$,\ps\ps  \item[(ii)] la
repr\'esentation de $K$ engendr\'ee par $e$ est isomorphe \`a $U_{|K}$. 
\end{itemize} \ps\ps  \noindent Alors : \ps\ps \begin{itemize} \item[(a)] Le ${\rm U}(\mathfrak{g})$-module ${\rm
U}(\mathfrak{g})e \subset V^\infty$ admet un caract\`ere central.  Sa classe de
conjugaison semi-simple associ\'ee est $\lambda + \rho$ o\`u $\lambda \in
\mathfrak{t}^\ast$ est le plus haut poids de $U$ relativement \`a $B$, et
$\rho$ est la demi-somme des racines de $\mathfrak{t}$ dans $\mathfrak{b}
\oplus \mathfrak{p}^-$.  \ps\ps  \item[(b)] La sous-repr\'esentation ferm\'ee
$V' \subset V$ engendr\'ee par $e$ sous l'action de ${\rm Sp}_{2g}(\R)$ est irr\'eductible.
De plus, si $f \in (V')^\infty$ satisfait (i) et (ii) alors $f \in \C[K].e$.
\end{itemize} \ps\ps 

\noindent \`A isomorphisme pr\`es, il existe au plus une repr\'esentation unitaire irr\'eductible
du groupe ${\rm Sp}_{2g}(\R)$ poss\'edant un vecteur $e$ qui est $\mathcal{C}^\infty$ et satisfait les
propri\'et\'es (i) et (ii).

\end{prop}

Ce r\'esultat est bien connu des sp\'ecialistes de la th\'eorie des
repr\'esentations unitaires des groupes de Lie, nous en donnons une
d\'emonstrations pour le confort du lecteur.\ps\ps 

\begin{pf}  D'apr\`es le (ii), l'espace $E=\C[K].e \subset
V^\infty$ est une repr\'esentation de $K$ isomorphe \`a $U_{|K}$.  Il est
stable par $\mathfrak{k}_\C$ et annul\'e par $\mathfrak{p}^-$ puisque ${\rm
ad}(K) \mathfrak{p}^- \subset \mathfrak{p}^-$.  Il est donc \'egalement
stable par la sous-alg\`ebre parabolique $\mathfrak{q} = \mathfrak{k}_\C
\oplus \mathfrak{p^-}$ de $\mathfrak{g}$.  Soit $Y$ le ${\rm
U}(\mathfrak{g})$-module induit ${\rm U}(\mathfrak{g}) \otimes_{{\rm
U}(\mathfrak{q})} E$.  L'inclusion de $E$ dans $V^\infty$ induit donc un
morphisme ${\rm U}(\mathfrak{g})$-\'equivariant $$u : Y \longrightarrow
V^\infty,$$ dont on notera l'image $X=u(Y)$.  Comme $K$ est connexe, $E$ est
irr\'eductible comme ${\rm U}(\mathfrak{k}_\C)$-module.  Quitte \`a
remplacer $e$ par un \'el\'ement bien choisi de $E$ on peut supposer que
$\mathfrak{b} e \subset \C e$ et que $h e = \lambda(h) e$ pour tout $h \in
\mathfrak{t}$ (th\'eorie du plus haut poids de Cartan-Weyl).  L'\'el\'ement
$1 \otimes e$ engendre donc $Y$ et satisfait $(\mathfrak{b} \oplus
\mathfrak{p}^-) ( 1 \otimes e) \subset \C (1 \otimes e)$ d'apr\`es (i), on en d\'eduit
que $Y$, et donc $X={\rm U}(\mathfrak{g})e$, admet un caract\`ere
infinit\'esimal satisfaisant l'assertion (a), d'apr\`es le second paragraphe
de l'exemple~\ref{exemplecentcar}.  \ps\ps  Observons \'egalement que $X$ est
stable par $K$.  De plus, l'action adjointe de $K$ sur ${\rm
U}(\mathfrak{g})$, ainsi que son action naturelle sur $E$, d\'efinissent une
structure de $K$-module sur $Y$ telle que $u$ soit $K$-\'equivariant.  Ces
structures font de $Y$ et $X$ des $(\mathfrak{g},K)$-modules, que nous
noterons $Y'$ et $X'$, et de $u$ un morphisme de $(\mathfrak{g},K)$-modules
(nous renvoyons \`a \cite[\S 3.3]{wallach} pour ces notions).  Nous avons
d\'ej\`a vu que $Y$ est un ${\rm U}(\mathfrak{g})$-module de plus haut
poids.  D'apr\`es~\cite[Prop.  7.1.8]{dixmier}, cela entra\^ine d'une part
que $Y'$ admet un unique quotient simple, et \'egalement que $Y'$ est
admissible (cela signifie que chaque repr\'esentation irr\'eductible de $K$ appara\^it avec une multiplicit\'e finie, et se d\'eduit de ce que chaque poids de $Y$ a une multiplicit\'e finie, voir
{\it loc.  cit.}).  Mais $X'$ admet un produit hermitien invariant (il est
unitaire au sens de~\cite[\S 9.3.3]{wallach}), car la repr\'esentation $V$
est unitaire par hypoth\`ese.  \'Etant admissible comme quotient de $Y'$, $X'$ est
irr\'eductible, c'est donc l'unique quotient irr\'eductible de $Y'$. 
Comme $X'$ admet un caract\`ere central, un r\'esultat de Harish-Chandra assure 
que tous ses vecteurs sont en fait analytiques \cite[\S 1.6, \S 3.4.9]{wallach},
\cite[VIII 8.7]{knapp}, et donc que l'adh\'erence $\overline{X}$ de $X$ dans $V$ est
stable par $G(\R)$ \cite[\S 1.6.6]{wallach}.  Elle admet $X'$ pour
$(\mathfrak{g},K)$-module : c'est donc l'unique repr\'esentation
irr\'eductible unitaire de $G(\R)$ de $(\mathfrak{g},K)$-module $X'$
\cite[\S 3.4.11]{wallach}.  Cela d\'emontre le premier point du (b).  Les
deux assertions restantes suivent du fait d\'ej\`a vu que $X'$ est l'unique
quotient irr\'eductible de $Y'$. \end{pf}

 Soit $(X,\Phi,\Delta,X^\vee,\Phi^\vee,\Delta^\vee)$ la donn\'ee radicielle
bas\'ee associ\'ee \`a $(\GL_g,T,B)$.  On \'ecrit de mani\`ere standard
$X=\oplus_{i=1}^g \Z \varepsilon_i, \, \Phi=\{ \pm (\varepsilon_i -
\varepsilon_j), 1 \leq i < j \leq g\}$, $\Delta=\{
\varepsilon_i-\varepsilon_{i+1}, 1\leq i <g\}$, $X^\vee= \oplus_{i=1}^g \Z \varepsilon_i^\ast$ et
$(\varepsilon_i-\varepsilon_j)^\vee=\varepsilon_i^\ast-\varepsilon_j^\ast$
pour $i<j$.  Les poids dominants de $\GL_g$ sont donc les $\lambda \in X$
tels que $\lambda=\sum_{i=1}^g m_i \varepsilon_i$ avec $m_1 \geq m_2 \geq
\cdots \geq m_g$.  \ps\ps  

\begin{cor}\label{infcarsiegel} Soit $W$ la $\C$-repr\'esentation
irr\'eductible de $\GL_g$ de plus haut poids $\sum_{i=1}^g m_i
\varepsilon_i$. Supposons qu'il existe une repr\'esentation unitaire
irr\'eductible $\pi'_W$ de ${\rm Sp}_{2g}(\R)$ satisfaisant les conditions de la
proposition~\ref{infcarholo} pour $U=W^\ast$. La classe de conjugaison semi-simple de
$\mathfrak{so}_{2g+1}(\C)$ qui correspond \`a ${\rm inf}_{\pi'_W}$ a pour valeurs propres les $2g+1$ entiers 
$$\pm (m_i-i), \, \, \, \, i=1,\cdots,g,\, \, \, \, {\rm et}\, \, \, \, 0.$$
\end{cor}

\begin{pf} En effet, un calcul sans difficult\'es montre que la
repr\'esentation adjointe de $K$ sur $\mathfrak{p}^-$ est isomorphe \`a la
restriction par l'homomorphisme ${\rm j}$ de la repr\'esentation $(g,X) \,\mapsto \,g\,
X \,{}^{\rm t}\!g$ de $\GL_g(\C)$ sur ${\rm Sym}_g(\C)$ (carr\'e sym\'etrique de
la repr\'esentation standard).  L'ensemble de ses poids est donc
$$\{\varepsilon_i+\varepsilon_j, 1\leq i \leq j \leq g\}.$$ Cette
description montre que la base correspondant \`a $\mathfrak{b}\oplus
\mathfrak{p}^-$ du syst\`eme de racines de $G_\C$ associ\'e \`a $T$ n'est
autre que la base standard introduite au~\S\ref{exdorad}.  En particulier,
l'\'el\'ement $\rho$ de la proposition~\ref{infcarholo} (a) n'est autre que
$g \varepsilon_1 + (g-1) \varepsilon_2 + \cdots + \varepsilon_g$.  Le poids 
dominant $\lambda$ de $W^\ast$ relativement \`a $B$ est $\sum_{i=1}^g -m_{g+1-i} \varepsilon_i$, 
et donc $$\lambda + \rho = \sum_{i=1}^g
(i-m_i)\varepsilon_{g+1-i}.$$ Les poids de $\widehat{{\rm Sp}_{2g}}(\C)={\rm
SO}_{2g+1}(\C)$ dans sa repr\'esentation standard sur $\C^{2g+1}$ \'etant
$0$ et les $\pm \varepsilon_i^\ast$ d'apr\`es le~\S\ref{exdorad}, cela conclut.\end{pf}

Soit $W$ une $\C$-repr\'esentation irr\'eductible de $\GL_g$ dans laquelle
$-1_g$ agit trivialement.  Soit $f \in {\rm S}_W({\rm Sp}_{2g}(\Z))$ une
forme modulaire de Siegel parabolique non nulle et de poids $W$.  La
proposition ci-dessus s'applique \`a $U=\mathcal{A}_{\rm cusp}({\rm
PGSp}_{2g})$, tout \'el\'ement $e$ dans l'image de $W^\ast \otimes f$
(Prop.~IV.\ref{caracsiegel}) et $U=W^\ast$.  Elle d\'emontre que si $w \in
W^\ast$ est non nul, la fonction $\varphi_{w,f} \in \mathcal{A}_{\rm
cusp}({\rm PGSp}_{2g})$ d\'efinie au~\S IV.\ref{fautsiegel} engendre topologiquement, sous l'action de ${\rm
Sp}_{2g}(\R)$, une sous-repr\'esentation irr\'eductible de $\mathcal{A}_{\rm
cusp}({\rm PGSp}_{2g})$, n\'ecessairement isomorphe \`a la repr\'esentation
$\pi'_W$ de l'\'enonc\'e pr\'ec\'edent.  Cela d\'emontre l'existence de
$\pi'_W$ d\`es que ${\rm S}_W({\rm Sp}_{2g}(\Z)) \neq 0$.  En fait, 
Harish-Chandra a d\'emontr\'e l'existence de $\pi'_W$ pour toute $W$
dont le plus haut poids satisfait $m_g>g$ (c'est la {\it s\'erie discr\`ete
holomorphe}, voir~\cite[Ch.  VI, \S 4 Thm.  6.6]{knapp}). 
Si cette hypoth\`ese sur $W$ est satisfaite nous dirons que $W$ est {\it
positive}; c'est le seul cas qui nous int\'eressera dans ce m\'emoire. 
Observons que si $W$ est positive, les $2g+1$ entiers de l'\'enonc\'e
ci-dessus sont distincts. \ps\ps

Supposons que $W$ est positive et que $-1_g$ agit trivialement dans $W$
(i.e. $\sum_i m_i  \equiv 0 \bmod 2$), de sorte que $\pi'_W$ se factorise
par ${\rm Sp}_{2g}(\R)/\{\pm 1_{2g}\}$. Il n'est pas difficile de v\'erifier
que $\pi'_W$ n'est pas isomorphe \`a sa
conjugu\'ee ext\'erieure\footnote{Cela vient de ce que cette conjugu\'ee
ext\'erieure admet un vecteur $\mathcal{C}^\infty$ annul\'e par
$\mathfrak{p}^+$ et engendrant $W$ sous l'action de $K$ ({\it plus bas poids}).
Son $(\mathfrak{g},K)$-module s'\'etudie de mani\`ere tout-\`a-fait
similaire \`a celui de $\pi'_W$ : il ne peut \^etre
isomorphe \`a celui de $\pi'_W$ que s'il est de dimension finie, i.e.  si 
$\pi'_W$ (et donc $W$) est triviale. Cela ne se produit pas car la
repr\'esentation triviale de ${\rm Sp}_{2g}(\R)$ n'appara\^it pas dans 
$\mathcal{A}_{\rm cusp}({\rm PGSp}_{2g})$.} par un \'el\'ement de ${\rm
PGSp}_{2g}(\R) \backslash {\rm Sp}_{2g}(\R)$. Autrement dit, la
repr\'esentation unitaire de ${\rm PSGp}_{2g}(\R)$
$$\pi_W = {\rm Ind}_{{\rm Sp}_{2g}(\R)}^{{\rm PGSp}_{2g}(\R)}
\pi'_W,$$
induite d'un sous-groupe d'indice $2$, est irr\'eductible. Bien entendu, $\pi_W$ et $\pi'_W$ ont m\^eme caract\`ere
infinit\'esimal, car ${\rm Sp}_{2g}(\R)$ et ${\rm PGSp}_{2g}(\R)$ ont m\^eme
alg\`ebre de Lie.  \ps\ps

Fixons $v_W \in \pi_W$ un vecteur $\mathcal{C}^\infty$ non nul, annul\'e par $\mathfrak{p}^-$,
engendrant $W^\ast$ sous l'action de $K$, et propre pour l'action de $\mathfrak{b}$. Un tel vecteur est
unique \`a multiplication par un \'el\'ement de $\C^\ast$ pr\`es d'apr\`es
la proposition~\ref{infcarholo}. On fixe de m\^eme $e_W \in
W^\ast$ non nul et de plus haut poids relativement \`a $B$. 

\begin{cor}\label{corapiw} Supposons $W$ positive. Si $F \in {\rm S}_W(\Sp_{2g}(\Z))$, il
existe un unique $u_F \in {\mathcal{A}}_{\pi_W}({\rm PGSp}_{2g})$ tel que
$u_F(v_W)=\varphi_{e_W,F}$. L'application $F \mapsto u_F$ d\'efinit un
isomorphisme ${\rm H}^{\rm opp}({\rm PGSp}_{2g})$-\'equivariant 
$${\rm S}_W({\rm Sp}_{2g}(\Z)) \isomo {\mathcal{A}}_{\pi_W}({\rm
PGSp}_{2g}).$$
\end{cor}

\begin{pf} Les propositions IV.\ref{caracsiegel} et \ref{infcarholo}
montrent que l'application de l'\'enonc\'e induit un isomorphisme ${\rm
H}^{\rm opp}({\rm PGSp}_{2g})$-\'equivariant entre ${\rm S}_W({\rm
Sp}_{2g}(\Z))$ et le sous-espace ${\rm Hom}_{G(\R)}(\pi_W, \mathcal{A}_{\rm
cusp}({\rm PGSp}_{2g})) \subset {\mathcal{A}}_{\pi_W}({\rm PGSp}_{2g})$.  On
conclut par le fait g\'en\'eral suivant : si $G$ est un $\Z$-groupe tel que
$G_\Q$ est semi-simple, et si $U$ est une s\'erie discr\`ete de $G(\R)$,
alors l'inclusion ${\rm Hom}_{G(\R)}(U,\mathcal{A}_{\rm cusp}(G)) \subset
\mathcal{A}_U(G)$ est une \'egalit\'e \cite[Thm.  4.3]{wallachart}.
\end{pf}

Ainsi, si $F \in {\rm S}_W(\Sp_{2g}(\Z))$ est non nulle et vecteur propre
de tous les op\'erateurs de Hecke dans ${\rm H}(\PGSp_{2g})$ il y a un
sens \`a parler de la repr\'esentation $\pi_F \in \Pi_{\rm
disc}(\PGSp_{2g})$ engendr\'ee par $F$, suivant la d\'efinition g\'en\'erale
du~\S\ref{fautdiscgen}. Elle satisfait
$(\pi_F)_\infty=\pi_W$. \ps\ps 

Le discours pr\'ec\'edent peut \'egalement \^etre tenu
pour ${\rm Sp}_{2g}$ plut\^ot que ${\rm PGSp}_{2g}$ et d\'emontre l'existence d'un isomorphisme ${\rm H}^{\rm opp}({\rm Sp}_{2g})$-\'equivariant entre ${\rm S}_W({\rm Sp}_{2g}(\Z))$ et ${\mathcal{A}}_{\pi'_W}({\rm Sp}_{2g})$; son contenu n'en est que moins fin, d'apr\`es la proposition~\ref{compsppgsp}. 

\ps \ps
\noindent {\sc Isomorphismes exceptionnels en genre $1$ et $2$ }
\ps\ps

Dans ce qui suit, on suppose $W$ positive, de plus haut poids $\sum_i m_i \varepsilon_i$, et
telle que $-1_g$ agit trivialement dans $W$.\ps\ps 

Supposons $g=1$. Dans ce cas, $W$ est la repr\'esentation $\DET^k$ avec
$k=m_1 > 1$ et $k \equiv 0 \bmod 2$.  L'isomorphisme $\mathfrak{sl}_2(\C)
\simeq {\mathfrak{so}}_3(\C)$ (carr\'e sym\'etrique) permet de
voir le caract\`ere infinit\'esimal de $\pi_W$ comme la classe de
conjugaison semi-simple dans $\mathfrak{sl}_2(\C)$ de valeurs propres $\pm
\frac{k-1}{2}$.  En fait, la classification bien connue du dual
unitaire de ${\rm SL}_2(\R)$ (Bargmann \cite{bargmann}) montre qu'\`a isomorphisme pr\`es,
l'unique repr\'esentation unitaire irr\'eductible de $\PGL_2(\R)$ ayant un
caract\`ere infinit\'esimal de valeurs propres $\pm \frac{k-1}{2}$, avec
$k>3$ entier pair, est la repr\'esentation $\pi_{\DET^k}$. Quand $k=2$, il
faut lui rajouter les deux repr\'esentations de dimension $1$.  \ps\ps 

Supposons $g=2$. Dans ce cas, conform\'ement aux notations
de~\cite{vandergeer},
$W$ est la repr\'esentation ${\rm Sym}^j(\C^2) \otimes \DET^k$ avec
$j=m_1-m_2$ et $k=m_2$, de plus $k>2$ et $j \equiv 0 \bmod 2$.  L'isomorphisme
exceptionnel $\mathfrak{sp}_4(\C) \simeq \mathfrak{so}_5(\C)$ permet 
de voir le caract\`ere infinit\'esimal de $\pi_W$ comme la classe
de conjugaison semi-simple dans $\mathfrak{sp}_4(\C)$ de valeurs propres
pour son action sur $\C^4$ : $$\pm \frac{w_1}{2},\,\,\,\, \pm
\frac{w_2}{2}$$ o\`u $w_1=m_1+m_2-3=2k+j-3$ et $w_2=m_1-m_2+1=j+1$.

\section{La conjecture d'Arthur-Langlands}

\subsection{Param\'etrisation de Langlands de $\Pi(G)$ pour $G$ semi-simple sur
$\Z$}\label{parlanpig} \ps\ps 

Soit $H$ un $\C$-groupe, de composante neutre $H^0$ et d'alg\`ebre de Lie complexe $\mathfrak{h}$. 
Notons $H(\C)_{\rm ss}$ (resp. $\mathfrak{h}_{\rm ss}$) l'ensemble des classes de
$H(\C)$-conjugaison d'\'el\'ements semi-simples de
$H^0(\C)$ (resp. de $\mathfrak{h}$), et consid\'erons l'ensemble
$$\mathcal{X}(H)$$
des familles $(c_v)_{v \in {\rm P} \cup \{\infty\}}$, o\`u $c_\infty \in
\mathfrak{h}_{\rm ss}$ et $c_p \in H(\C)_{\rm ss}$ pour
tout $p \in {\rm P}$. Dans la discussion qui va suivre, $H$ sera connexe
(et m\^eme semi-simple), mais des exemples non connexes associ\'es aux groupes
orthogonaux pairs appara\^itront par la suite. Tout morphisme de
$\C$-groupes $r : H \rightarrow H'$ d\'efinit une application encore not\'ee
$r : \mathcal{X}(H) \rightarrow \mathcal{X}(H')$, envoyant $(c_v)$ sur $(r(c_v))$. \ps\ps 

Soit $G$ un $\Z$-groupe semi-simple. Comme nous l'avons d\'ej\`a dit, pour tout nombre premier $p$ le $\Z_p$-groupe $G_{\Z_p}$ est r\'eductif et d\'eploy\'e \cite[Prop. 1.1]{grossinv}; il admet donc une donn\'ee radicielle bas\'ee $\Psi(G_{\Z_p})$. De plus, si $\overline{\Q}$ (resp. $\overline{\Q_p}$) est une cl\^oture alg\'ebrique de $\Q$ (resp. $\Q_p$), et si $\overline{\Q} \rightarrow \overline{\Q_p}$ et $\overline{\Q} \rightarrow \C$ sont deux plongements, les isomorphismes associ\'es de donn\'ees radicielles bas\'ees 
$$\Psi(G_{\Z_p}) \isomo \Psi(G_{\overline{\Q_p}})  \overset{\sim}{\leftarrow} \Psi(G_{\overline{\Q}})  \isomo \Psi(G_{\C})$$
ne d\'ependent d'aucun de ces choix de plongements effectu\'es\footnote{L'argument de Gross est le suivant. C'est un fait g\'en\'eral que l'action naturelle de ${\rm Gal}(\overline{\Q}/\Q)$ sur $\Psi(G_{\overline{\Q}})$ se factorise en une action fid\`ele du groupe de Galois d'un corps de nombres $K$ galoisien sur $\Q$. Le fait que $G$ est r\'eductif sur $\Z_p$ entra\^ine que $K$ est non ramifi\'e en $p$, et donc $K=\Q$ d'apr\`es un r\'esultat fameux de Minkowski. Cela entra\^ine que $G$ est d\'eploy\'e sur $\Z_p$ et le reste des assertions ci-dessus.}. Le dual de Langlands de $G_{\overline{\Q}}$ est donc canoniquement le dual de Langlands des $G_{\Z_p}$ pour tout $p$ et de $G_\C$ : nous le noterons $\widehat{G}$. \ps\ps 

\ps

Suivant Langlands~\cite{langlandsyale}, on dispose d'une application canonique
$${\rm c} : \Pi(G) \longrightarrow \mathcal{X}(\widehat{G}), \pi \mapsto
({\rm c}_v(\pi)),$$
d\'efinie comme suit. Soit $\pi=\pi_\infty \otimes \pi_f \in \Pi(G)$. On note
${\rm c}_\infty(\pi)$ le caract\`ere infinit\'esimal de $\pi_\infty$ (\S\ref{carinf}). L'isomorphisme de Satake entra\^ine que ${\rm H}(G)=\bigotimes_p {\rm H}_p(G)$ est commutative, de sorte que
$\pi_f$ est de dimension $1$ et peut \^etre vue comme un homomorphisme
d'anneaux de ${\rm H}(G)^{\rm opp}={\rm H}(G)$ dans $\C$, ou ce qui revient au m\^eme
comme une collection de morphismes d'anneaux $$\pi_p :
{\rm H}_p(G) \rightarrow \C,$$ o\`u $\pi_p$ est la restriction de $\pi_f$ \`a ${\rm H}_p(G)$ au sens du
\S IV.\ref{annheckeg}. Le scholie~\ref{corsatake} associe donc \`a chaque $\pi_p$ un unique
\'el\'ement ${\rm c}_p(\pi) \in \widehat{G}(\C)_{\rm ss}$. Par d\'efinition,
${\rm c}(\pi)$ d\'etermine $\pi_f$ et le caract\`ere infinit\'esimal de
$\pi_\infty$, l'application ${\rm c}$ est donc \`a fibres finies (Harish-Chandra~\S\ref{carinf}). \ps\ps 

\begin{example}\label{exempletriv} {\rm (La repr\'esentation triviale) Soit $\pi=1_G \in
\Pi_{\rm disc}(G)$ la repr\'esentation triviale de $G$
(\S IV.\ref{fautdiscgen}).  D'apr\`es l'exemple~\ref{exemplecentcar}, $2
{\rm c}_\infty(\pi)$ est la classe de conjugaison du co-poids $2\rho$ de
$\widehat{G}$.  De m\^eme, d'apr\`es l'exemple~\ref{degretriv}, la classe de
conjugaison ${\rm c}_p(\pi)$ est celle de $\rho(p)=(2\rho)(p^{\frac{1}{2}})$.  }
\end{example}

\subsection{Quelques exemples} \label{exparlan}\ps\ps 
 Consid\'erons
d'abord le $\Z$-groupe ${\rm PGL}_2$, de groupe dual ${\rm SL}_2$. 
Soit $k>0$ un entier pair, $F=\sum_{n\geq 1} a_n q^n \in {\rm S}_k({\rm SL}_2(\Z))$ une forme modulaire
propre pour tous les op\'erateurs de Hecke de ${\rm H}(\PGL_2)$ telle que
$a_1=1$ (elles forment une base de ${\rm S}_k({\rm SL}_2(\Z))$ \cite[VII \S 5.4]{serre}).
Soit $\pi \in \Pi_{\rm cusp}(\PGL_2)$ la repr\'esentation engendr\'ee par
$F$ (voir \S IV.\ref{fautdiscgen}, \S V.\ref{carinf}). On a d\'ej\`a d\'etermin\'e 
${\rm c}_\infty(\pi)$, c'est-\`a-dire ${\rm inf}_{\pi_W}$, en fonction de $k$ {\it loc. cit.} Un $\Z$-isomorphisme $\PGL_2 \simeq
\PGSp_2$ induit un isomorphisme ${\rm H}_p(\PGSp_2) \isomo {\rm H}_p(\PGL_2)$
envoyant ${\rm K}_p$ sur ${\rm T}_{\Z/p\Z}$, les relations~\eqref{satakesp}
et~\eqref{tradserregl2}, ainsi que \cite[VII Thm. 7]{serre}, montrent donc que pour tout $p$ premier 
$$p^{\frac{k-1}{2}}\, {\rm Trace}({\rm c}_p(\pi)) = a_p.$$

Consid\'erons maintenant le $\Z$-groupe ${\rm PGSp}_4$, de groupe dual le
$\C$-groupe ${\rm Sp}_4$ (qui est aussi ${\rm Spin}_5$). Soit $W$ la repr\'esentation ${\rm Sym}^j(\C^2) \otimes
\DET^k$ de $\GL_2(\C)$ avec $k\geq 3$, soit $F \in {\rm S}_W({\rm
Sp}_4(\Z))$ une forme propre (non nulle), et soit $\pi \in \Pi_{\rm
cusp}({\rm PGSp}_4)$ la repr\'esentation engendr\'ee par $F$. On a d\'ej\`a
d\'etermin\'e ${\rm c}_\infty(\pi)$ (i.e. ${\rm inf}_{\pi_W}$) en fonction de $j$ et $k$ au~\S\ref{carinf}. Si $p$
est premier, l'\'el\'ement ${\rm c}_p(\pi) \in {\rm Sp}_4(\C)_{\rm ss}$ est uniquement caract\'eris\'e par sa
trace et celle de son carr\'e altern\'e dans la repr\'esentation
tautologique $V_{\rm Spin} \simeq \C^4$ de ${\rm Sp}_4(\C)$. Si ${\rm K}_p(F)= a_p F$ et ${\rm
T}_p(F)=b_p F$, les relations~\eqref{satakesp} montrent que 
$$p^{3/2}\, {\rm Trace}({\rm c}_p(\pi) \,|\, V_{\rm Spin}) = a_p \, \, \, \, \, 
{\rm et}\, \, \, \, p^2 \,{\rm
Trace}({\rm c}_p(\pi)\,|\, \Lambda^2 V_{\rm Spin}) = b_p+p^2+1.$$

Pour un $g\geq 1$ g\'en\'eral, consid\'erons une forme propre $F \in {\rm
S}_W({\rm Sp}_{2g}(\Z))$ telle que ${\rm T}_p(F)=b_p F$.  Si $\pi \in
\Pi_{\rm cusp}({\rm Sp}_{2g})$ d\'esigne la repr\'esentation engendr\'ee par
$F$, alors ${\rm c}_\infty(\pi) \subset {\rm so}_{2g+1}(\C)_{\rm ss}$ est donn\'e
en fonction de $W$ par la proposition~\ref{infcarsiegel}.  Si $V_{\rm St} \simeq
\C^{2g+1}$ d\'esigne la repr\'esentation tautologique de ${\rm
SO}_{2g+1}(\C)$, alors pour tout premier $p$, $$p^g\, {\rm Trace}({\rm c}_p(\pi)\,
|\, V_{\rm
st})=b_p+1.$$  \ps\ps 

Soit maintenant $n \equiv 0 \bmod 8$ et $G={\rm SO}_n$ le groupe sp\'ecial
orthogonal de ${\rm E}_n$ (\S IV.\ref{fautson}), de sorte que $\widehat{G}(\C)={\rm SO}_n(\C)$
Soit $W$ la repr\'esentation irr\'eductible de plus haut poids
$\sum_{i=1}^{n/2} m_i \varepsilon_i$ avec $m_1 \geq \cdots m_{n/2-1} \geq
|m_{n/2}|$ dans les notations du \S\ref{exdorad}.  Soit $F \in {\rm M}_W({\rm
SO}_n)$ une forme propre et $\pi \in \Pi_{\rm disc}({\rm SO}_n)$ la
repr\'esentation engendr\'ee.  Par d\'efinition, on a $\pi_\infty \simeq
W^\ast$, mais $W^\ast \simeq W$ car $n \equiv 0 \bmod 4$, de sorte que les
$n$ valeurs propres de ${\rm c}_\infty(\pi) \in \mathfrak{so}_n(\C)_{\rm ss}$ sont
$$\pm ( m_i + \frac{n}{2} - i ), \, \, \, i = 1, \cdots, \frac{n}{2}$$
d'apr\`es le~\S\ref{carinf}.  Soit $p$ un nombre premier.  Supposons que
${\rm T}_p(F)=\lambda_p F$, ${\rm T}_{p^2}(F) =\lambda_{p^2} F$ et ${\rm
T}_{p,p}(F)=\lambda_{p,p} F$
(\S\ref{exemplehecke}). Les
relations~\eqref{formulesattp},\eqref{formulesattpp}, ainsi que
l'exemple~\ref{calcultp2}, s'\'ecrivent
$$p^{\frac{n}{2}-1}\, {\rm Trace}( {\rm c}_p(\pi) \,|\, V_{\rm St}) = \lambda_p,$$
$$p^{n-3}\,{\rm Trace}( {\rm c}_p(\pi) \,|\, \Lambda^2 V_{\rm St} )= \lambda_{p,p} +
p^{\frac{n}{2}-2}+\frac{p^{n-2}-1}{p^2-1},$$
et $(p+1) \lambda_{p,p}=\lambda_p^2 -\lambda_{p^2} -
\frac{(p^{\frac{n}{2}}-1)(p^{\frac{n}{2}-1}+1)}{p-1}$. \ps\ps


\subsection{La conjecture d'Arthur-Langlands} \label{parconjarthlan}

Soient $G$ un $\Z$-groupe semi-simple et $r : \widehat{G} \rightarrow {\rm SL}_n$ une $\C$-repr\'esentation. Cette repr\'esentation induit une application
$\mathcal{X}(\widehat{G}) \rightarrow \mathcal{X}(\SL_n), (c_v) \mapsto
(r(c_v))$ que nous notons encore $r$.  Si $\pi \in \Pi(G)$, on peut lui associer l'\'el\'ement $$\psi(\pi,r):=r({\rm c}(\pi)) \in \mathcal{X}(\SL_n).$$  
Cet \'el\'ement est appel\'e {\it param\`etre de Langlands du couple $(\pi,r)$}. Si $\pi \in \Pi_{\rm disc}(G)$, les conjectures de Langlands \cite{langlandspb}, pr\'ecis\'ees par Arthur \cite{arthurunipotent}, affirment que $\psi(\pi,r)$ s'exprime en terme des $\Pi_{\rm cusp}(\PGL_m)$ pour $m\geq 1$. Avant de rappeler comment, nous devons introduire quelques notations.  \ps \medskip

-- Notons ${\rm St}_m$ la $\C$-repr\'esentation tautologique de $\SL_m$ sur $\C^m$. Si $a$ et $b$ sont des entiers, la somme directe et le produit tensoriel des repr\'esentations ${\rm St}_a$ et ${\rm St}_b$ d\'efinissent des $\C$-repr\'esentations de ${\rm SL}_a \times {\rm SL}_b$ de dimensions respectives $a+b$ et $ab$, ainsi donc que des applications naturelles $$\mathcal{X}(\SL_a) \times \mathcal{X}(\SL_b) \rightarrow \mathcal{X}(\SL_{a+b})\, \, \, \, {\rm et}\, \, \, \, \mathcal{X}(\SL_a) \times \mathcal{X}(\SL_b) \rightarrow \mathcal{X}(\SL_{ab}).$$
Nous les noterons respectivement $(c,c') \mapsto c \oplus c'$ et $(c,c') \mapsto c \otimes c'$. Ces op\'erations sont commutatives, associatives et distributives en un sens \'evident. \ps\ps 

-- Suivant Arthur~\cite{arthurunipotent}, consid\'erons l'\'el\'ement $e \in \mathcal{X}({\rm SL}_2)$ d\'efini par 
$$e_p = \left[ \begin{array}{cc} p^{-\frac{1}{2}} & 0 \\ 0 & p^{\frac{1}{2}}
\end{array}\right] \, \, \, \forall p \in {\rm P}, \, \, \,\, \, {\rm et}\, \,\, \, \, 
e_\infty=\left[\begin{array}{cc} -\frac{1}{2} & 0 \\ 0 & \frac{1}{2}
\end{array}\right].$$ 
Il donne naissance pour tout entier $d\geq 1$ \`a l'\'el\'ement ${\rm Sym}^{d-1}(e) \in \mathcal{X}(\SL_d)$, o\`u ${\rm Sym}^{d-1}$ d\'esigne la repr\'esentation ${\rm Sym}^{d-1} {\rm St}_2$ de $\SL_2$. Nous le noterons $[d]$; par exemple $[2]=e$. Ces \'el\'ements joueront par la suite un r\^ole particuli\`erement important. Observons d'ores et d\'ej\`a que $[d]={\rm c}(1_{\PGL_d})$,  d'apr\`es l'exemple~\ref{exempletriv}. On pose plus g\'en\'eralement, pour $m, d\geq 1$ entiers et $c \in \mathcal{X}(\SL_m)$,
$$c [d] := c \otimes[d].$$

-- Si $\pi \in \Pi_{\rm cusp}(\PGL_m)$, l'\'el\'ement ${\rm c}(\pi) \in \mathcal{X}(\SL_m)$ sera simplement not\'e $\pi$. Cet abus de langage sera en g\'en\'eral innocent, car ${\rm c}(\pi)$ d\'etermine $\pi$ d'apr\`es le {\it th\'eor\`eme de multiplicit\'e $1$ forte} de Piatetski-Shapiro, Jacquet
et Shalika \cite{jasha1}. (On prendra garde que l'injectivit\'e de l'application de param\'etrisation ${\rm c}$ est tr\`es sp\'ecifique des
$\Z$-groupes ${\rm PGL}_m$).  

 \ps \ps

Ainsi, si $n_1,\dots,n_k$, $d_1,\dots,d_k$ sont des entiers $\geq 1$, si $\pi_i \in \Pi_{\rm cusp}(\PGL_{n_i})$ pour tout $i=1,\dots,k$, et si $n=\sum_{i=1}^k n_i d_i$, on dispose d'un \'el\'ement bien d\'efini 
$$\pi_1[d_1]\oplus \pi_2[d_2] \oplus \cdots \oplus \pi_k[d_k] \,\,\,\in \,\,\,\mathcal{X}(\SL_n).$$
Il ne d\'epend que du multi-ensemble $\{(\pi_i,d_i), i=1,\dots,k\}$. Le sous-ensemble de $\mathcal{X}(\SL_n)$ form\'e des \'el\'ements qui sont de cette forme, pour un quadruple $(k,(n_i),(d_i),(\pi_i))$ quelconque tel que $n=\sum_{i=1}^k n_i d_i$, sera not\'e $$\mathcal{X}_{\rm AL}(\SL_n).$$ On dispose d'un r\'esultat d'unicit\'e remarquable, d\^u \`a Jacquet-Shalika \cite{jasha2} (voir aussi~\cite{Lgl}) : 

\begin{prop}\label{jacquetshalika} Supposons que $\oplus_{i=1}^k \pi_i[d_i] = \oplus_{j=1}^{l} \pi'_{j}[d'_{j}]$  dans $\mathcal{X}(\SL_n)$.  Alors $k=l$ et il existe $\sigma \in \got{S}_k$ tel que pour tout $1 \leq i \leq k$, $(\pi'_i,d'_i)=(\pi_{\sigma(i)},d_{\sigma(i)})$. 
\end{prop}

Le cas particulier des conjectures d'Arthur et Langlands que nous souhaitons mettre en \'evidence est le suivant. 

\begin{conj}\label{conjarthlan} {\rm (Langlands~\cite{langlandspb},
Arthur~\cite{arthurunipotent})} Soient $G$ un $\Z$-groupe semi-simple et $r :
\widehat{G} \rightarrow {\rm SL}_n$ une $\C$-repr\'esentation. Si $\pi \in \Pi_{\rm disc}(G)$ alors $\psi(\pi,r) \in \mathcal{X}_{\rm AL}(\SL_n)$. 
\end{conj} \ps\ps 

Autrement dit, pour tout $\pi \in \Pi_{\rm disc}(G)$, il existe un entier $k\geq 1$, des entiers $n_1,\dots,n_k$, $d_1,\dots,d_k$ et des repr\'esentations $\pi_i \in \Pi_{\rm cusp}(\PGL_{n_i})$ pour tout $1 \leq i \leq k$, tels que $\psi(\pi,r)=\oplus_{i=1}^k \pi_i[d_i]$ (et cette \'ecriture est unique \`a permutation des facteurs pr\`es d'apr\`es la proposition~\ref{jacquetshalika}). Concr\`etement, cela dit que pour
tout $v \in {\rm P} \cup \{\infty\}$, les $n$ valeurs propres de la classe de conjugaison
semi-simple $\rho({\rm c}_v(\pi))$ sont les $\lambda \,p^\mu$ pour $v \in {\rm P}$ (resp. $\lambda+\mu$ pour $v=\infty$) o\`u $i$ parcourt $\{1,\dots,k\}$, et o\`u pour tout $i \in \{1,\dots,k\}$ : \begin{itemize} \ps\ps 
\item[-] $\lambda$ parcourt les valeurs propres de ${\rm
St}_{r_i}({\rm c}_v(\pi_i))$ compt\'ees avec multiplicit\'es, \ps\ps 
\item[-] $\mu$ parcourt $\frac{1-d_i}{2}, \,\frac{3-d_i}{2}, \,\cdots, \,\frac{d_i-3}{2},\,\frac{d_i-1}{2}$. 
\end{itemize}
 \ps\ps
\medskip

Le groupe $G$ et $\rho$ \'etant donn\'es, Langlands et Arthur proposent
\'egalement une description conjecturale de l'image de $\pi \mapsto
\psi(\pi,r)$ qui est nettement plus difficile \`a formuler en
g\'en\'eral. Disons simplement qu'elle comporte deux ingr\'edients : le
premier est un certain groupe hypoth\'etique appel\'e {\it groupe de
Langlands de $\Q$} et le second est la {\it formule de multiplicit\'e
d'Arthur-Langlands} (voir \cite{Lgl},\cite{arthurconjlan}, ainsi que
\cite{chrenard2} pour une discussion d'une variante ``sur $\Z$'' du groupe de
Langlands de $\Q$). Dans les exemples qui suivent, nous ne discutons que de cas, nettement plus
simples, o\`u $G$ est soit $\PGL_n$ soit un groupe classique, et o\`u $r$ est sa
repr\'esentation ``tautologique''. \ps\ps

\subsection{Quelques exemples}\label{exconjal}

\ps\ps
\noindent {\sc Le cas de $\PGL_n$}\label{pargelbartjacquet}
\ps\ps

Si $\pi \in \Pi_{\rm cusp}(\PGL_n)$, on a tautologiquement ${\rm c}(\pi)=\psi(\pi,{\rm
St}_n)=\pi$.  Si $d$ est un diviseur de $n$ et si $\varpi \in \Pi_{\rm
cusp}(\PGL_{n/d})$, Speh a construit \`a l'aide de r\'esidus de s\'eries
d'Eisenstein un $\pi \in \Pi_{\rm disc}(\PGL_n)$ tel que $\psi(\pi,{\rm
St}_n)=\varpi[d]$.  La conjecture de Jacquet, d\'emontr\'ee par Moeglin et
Waldspurger~\cite{moeglinwaldspurger}, affirme que tout $\pi \in \Pi_{\rm
disc}(\PGL_{n})$ est de cette forme.  Cela d\'emontre la
conjecture~\ref{conjarthlan} quand $G={\rm PGL}_n$ et $r={\rm St}_n$.  \ps\ps 

Un autre cas fameux et connu de la conjecture~\ref{conjarthlan} concerne la repr\'esentation
carr\'e sym\'etrique ${\rm Sym}^2 : \widehat{\PGL_2}={\rm SL}_2 \rightarrow
\SL_3$. Plus
pr\'ecis\'ement, Gelbart et Jacquet ont d\'emontr\'e que si $\pi \in \Pi_{\rm cusp}(\PGL_2)$,
il existe une unique $\pi' \in \Pi_{\rm cusp}(\PGL_3)$ telle que
$\psi(\pi,{\rm Sym}^2)=\pi'$ \cite{GJ}. On note par abus $\pi'={\rm Sym}^2\, \pi$. Par
exemple, si $\pi$ est engendr\'ee par une forme propre $F \in {\rm S}_k({\rm
SL}_2(\Z))$ comme au~\S\ref{parlanpig}, alors ${\rm c}_\infty({\rm Sym}^2\, \pi)$ admet pour valeurs propres
$0$ et $\pm (k-1)$, et ${\rm c}_p({\rm Sym}^2\, \pi) \in \SL_3(\C)_{\rm ss}$ admet pour polyn\^ome
caract\'eristique $$X^3- (p^{1-k}a_p^2-1) (X^2-X)+1.$$ 
\ps\ps
\noindent {\sc Groupes classiques}
\ps\ps

La {\it classification d'Arthur} \cite{arthur}, sur laquelle nous
reviendrons au chapitre suivant, d\'emontre la conjecture~\ref{conjarthlan} 
quand $G$ est un {\it groupe classique} (au sens l\'eg\`erement restrictif pr\'ecis\'e
ci-dessous) et $r$ la {\it repr\'esentation
standard} de $\widehat{G}$; Arthur d\'ecrit aussi l'image de $\pi \mapsto
\psi(\pi,r)$.  Contentons-nous ici de pr\'eciser les terminologies en italique.  \ps\ps 

Notons ${\rm Class}_\C$ l'ensemble constitu\'e des $\C$-groupes ${\rm
Sp}_{2g}$ pour un entier $g\geq 1$ et des $\C$-groupes ${\rm SO}_m$ pour un
entier $m \neq 2$. Par exemple, ${\rm PGSp}_{2g}$ n'est pas isomorphe \`a un groupe
appartenant \`a ${\rm Class}_\C$ si $g>2$. Le lecteur
observera que les $\C$-groupes des deux familles sus-mentionn\'ees sont deux
\`a deux non isomorphes. Tout $\C$-groupe dans ${\rm Class}_\C$ poss\`ede une
$\C$-rep\'esentation distingu\'ee tautologique, sur $\C^{2g}$ pour ${\rm
Sp}_{2g}$, sur $\C^m$ pour ${\rm SO}_m$, appel\'ee repr\'esentation {\it standard},
et not\'ee ${\rm St}$. Elle est irr\'eductible et fid\`ele, et de dimension minimale
pour ces propri\'et\'es. \ps\ps 

Les $\Z$-groupes semi-simples $G$ auxquels s'appliquent les travaux d'Arthur
sus-cit\'es sont ceux tels que $G_\C$ est isomorphe \`a un \'el\'ement de ${\rm
Class}_\C$, auquel cas $\widehat{G}$ l'est aussi (\S\ref{exdorad}). C'est le cas du $\Z$-groupe
${\rm Sp}_{2g}$ et du $\Z$-groupe ${\rm SO}_L$ o\`u $L$ est un ${\rm
q}$-module sur $\Z$ de dimension $>2$ (\S II.1). C'est \'egalement le cas des $\Z$-groupes sp\'eciaux
orthogonaux de rang impair \'etudi\'es dans l'appendice ${\rm B}$ (voir la proposition 1.7 de cet appendice). On peut d\'emontrer qu'il n'y en a pas
d'autre \`a isomorphismes pr\`es \cite{grossinv}.  \ps\ps 

\ps\ps
\noindent {\sc Cas des $\Z$-groupes ${\rm O}_n$ et ${\rm PGO}_n$}
\ps\ps

Consid\'erons le $\Z$-groupe orthogonal $G'={\rm O}_n$ du ${\rm q}$-module
${\rm E}_n$ (\S IV.\ref{fautson}).  N'\'etant pas connexe, on ne peut lui
appliquer verbatim les consid\'erations du~\S\ref{parlanpig}. N\'eanmoins,
nous avons d\'ej\`a observ\'e au corollaire \ref{corsatakeo} que pour tout nombre
premier $p$, l'isomorphisme de Satake de $G={\rm SO}_n$ induit une bijection
$${\rm Hom}_{\rm anneaux}({\rm H}_p({\rm O}_n),\C) \isomo {\rm
O}_n(\C)_{\rm ss};$$ on prendra garde que le terme de droite d\'esigne
l'ensemble des classes de ${\rm O}_n(\C)$-conjugaison d'\'el\'ements
semi-simples de ${\rm SO}_n(\C)$ (\S\ref{parlanpig}).  De m\^eme, il est facile de v\'erifier que
si $V'$ est une repr\'esentation irr\'eductible de $G'(\R)$, sa restriction
$V$ \`a $G(\R)$ est soit irr\'eductible, soit somme de deux
repr\'esentations non isomorphes mais conjugu\'ees ext\'erieures sous
l'action de ${\rm O}_n(\R)$; la ${\rm O}_n(\C)$-orbite de l'\'el\'ement de
$(\mathfrak{so}_n)_{\rm ss}$ associ\'e au caract\`ere infinit\'esimal de
chacun des constituants de $V$ est donc ind\'ependante du constituant
choisi, et nous dirons par abus que c'est le caract\`ere infinit\'esimal de
$V'$. Nous avons donc d\'efini une
application de param\'etrisation $$\Pi({\rm O}_n) \longrightarrow
\mathcal{X}({\rm O}_n(\C))$$ que nous noterons encore ${\rm c} : \pi \mapsto
({\rm c}_v(\pi))$. Il est donc naturel de poser
$$\mathcal{X}(\widehat{{\rm O}_n}):=\mathcal{X}({\rm O}_n(\C)).$$

Suivant le~\S\ref{parlanpig}, toute $\C$-repr\'esentation $r : {\rm O}_n(\C) \rightarrow
\GL_n(\C)$ induit une application 
$$r : \mathcal{X}(\widehat{{\rm O}_n}) \rightarrow \mathcal{X}({\rm GL}_n)$$
de sorte qu'il y a encore un sens \`a poser $\psi(\pi,r)=r({\rm c}(\pi))$ quand $r$ est comme ci-dessus
et $\pi \in \Pi({\rm O}_n)$.
\ps\ps 

On dispose pour ${\rm O}_n$ d'une conjecture analogue \`a la
conjecture~\ref{conjarthlan}; elle se d\'eduit en fait de cette conjecture
appliqu\'ee au $\Z$-groupe ${\rm SO}_n$, expliquons comment. Soient 
$\pi' \in \Pi_{\rm disc}({\rm O}_n)$, $V'=(\pi')_\infty^\ast$ et
$V=V'_{|{\rm SO}_n(\R)}$. D'apr\`es le \S IV.\ref{fauton}, on dispose d'une injection naturelle ${\rm H}({\rm O}_n)$-\'equivariante
$${\rm res}: {\rm M}_{V'}({\rm O}_n) \longrightarrow  {\rm M}_V({\rm SO}_n).$$
Soit $\pi \in \Pi_{\rm disc}({\rm SO}_n)$ la repr\'esentation engendr\'ee par
l'une quelconque des formes propres appartenant au ${\rm H}({\rm
SO}_n)$-module engendr\'e par ${\rm res}(F)$, o\`u $F \in {\rm M}_{V'}({\rm O}_n)$
d\'esigne une forme propre quelconque engendrant $\pi'$ (\S IV.\ref{fautdiscgen}).
 Les repr\'esentations $\pi \in \Pi_{\rm disc}({\rm
SO}_n)$ obtenues ainsi seront appel\'ees {\it les constituants de la restriction de $\pi'$
\`a ${\rm SO}_n$}; elles forment un ensemble fini non vide. Si $\pi$ est un tel constituant, ${\rm c}(\pi')$ est par
d\'efinition l'image de ${\rm c}(\pi)$ par l'homomorphisme canonique $\mathcal{X}(\widehat{{\rm
SO}_n}) \rightarrow \mathcal{X}(\widehat{{\rm O}_n})$. La proposition qui suit est donc \'evidente. \ps\ps 

\begin{prop}\label{paramarthlangon} Soient $\pi' \in \Pi_{\rm disc}({\rm O}_n)$, $\pi \in \Pi_{\rm
disc}({\rm SO}_n)$ un constituant de la restriction de $\pi'$ \`a ${\rm
SO}_n$, $r' : {\rm O}_n(\C) \rightarrow \GL_m(\C)$ une
$\C$-repr\'esentation, et $r$ la restriction $r'$ \`a ${\rm SO}_n(\C)$. Alors $\psi(\pi',r')=\psi(\pi,r)$. En particulier, la conjecture d'Arthur-Langlands est satisfaite pour $(\pi,r)$ si et seulement si elle l'est pour $(\pi',r)$. 
\end{prop}

Cette proposition est particuli\`erement utile dans le cas de la
repr\'esentation tautologique ${\rm St}$ de ${\rm O}_n(\C)$ sur $\C^n$, dont
la restriction \`a ${\rm SO}_n(\C)$ est sa repr\'esentation standard. \ps\ps 

La discussion de ce paragraphe admet enfin un analogue naturel au cas du
$\Z$-groupe ${\rm PGO}_n$, pour lequel on dispose d'une application de
param\'etrisation $${\rm c} : \Pi({\rm PGO}_n) \rightarrow \mathcal{X}(\widehat{{\rm
PGO}_n}):=\mathcal{X}({\rm Pin}_n)$$
o\`u ${\rm Pin}_n$ d\'esigne le $\C$-groupe ${\rm Pin}$ du ${\rm q}$-module
standard $V_n$ de dimension $n$ sur $\C$. Suivant~\cite{ABS}, on rappelle
qu'il s'agit du sous-groupe des
\'el\'ements $x$ de l'alg\`ebre de Clifford ${\rm C}(V_n)$ de $V_n$ tels que $x
x^{t} =1$ et $\alpha(x)V_nx^{-1} \subset V_n$, $x \mapsto \alpha(x)$ et $x \mapsto x^t$ d\'esignant
respectivement l'involution et l'anti-involution canonique de ${\rm
C}(V_n)$. Sa composante neutre, d\'efinie part $\alpha={\rm id}$, est le
$\C$-groupe ${\rm Spin}_n$; elle est d'indice $2$. Tout \'el\'ement $e \in V_n$ 
tel que ${\rm q}(e)=1$ (i.e.
$e^2=1$) appartient \`a ${\rm Pin}(V_n)={\rm Pin}_n(\C)$ et d\'efinit une section du
morphisme canonique ${\rm Pin}_n \rightarrow \Z/2\Z$. Enfin, on dispose
d'un morphisme surjectif naturel ${\rm Pin}_n \rightarrow
{\rm O}_{V_n}, x \mapsto (v \mapsto \alpha(x)vx^{-1})$. Son noyau est $\pm 1$
et l'image de tout
\'el\'ement $e \in V_n$ tel que ${\rm q}(e)=1$ est la sym\'etrie orthogonale
par rapport \`a $e$. \ps
\medskip
\ps\ps
\noindent {\sc Le cas de la repr\'esentation triviale}
\ps\ps

Supposons maintenant que $G$ est un $\Z$-groupe semi-simple quelconque et consid\'erons la repr\'esentation triviale $1 \in \Pi_{\rm disc}(G)$. 
Soit $\mu : {\rm SL}_2 \rightarrow \widehat{G}$ un $\C$-morphisme principal au sens de Kostant. Une formulation \'equivalente de la description de ${\rm c}(1)$ donn\'ee dans l'exemple \ref{exempletriv} est 
$$\mu(e)  = {\rm c}(1)$$
(cette observation est explicitement donn\'ee dans \cite[\S 7]{grossatake}, mais elle est sans doute plus ancienne). En particulier, si $r : \widehat{G} \rightarrow {\rm SL}_n$ d\'esigne un $\C$-morphisme arbitraire, et s'il on d\'ecompose la repr\'esentation $r \circ \mu$ de ${\rm SL}_2$ sous la forme $\oplus_{i=1}^k {\rm Sym}^{d_i-1}\, \C^2$ avec $d_1,\dots,d_k$ des entiers $\geq 1$, on a donc 
$$\psi(1,r) = \oplus_{i=1}^k [d_i].$$
Ceci est bien en accord avec la conjecture d'Arthur-Langlands (en l'occurence, avec les conjectures d'Arthur \cite{arthurunipotent}), qui vaut donc pour le couple $(1,r)$ quelque soit la repr\'esentation $r$. \ps\ps 

Supposons par exemple $G={\rm SO}_n$ avec $n \equiv 0 \bmod 8$. On a alors ${\rm St} \circ \mu = {\rm Sym}^{n-2} \,\C^2 \oplus  1$, autrement dit $\psi(1,{\rm St})=[n-1]\oplus [1]$. Cette relation vaut \'egalement si $G={\rm O}_n$, car la repr\'esentation triviale de ${\rm SO}_n$ est \'evidemment la restriction \`a ${\rm SO}_n$ de la repr\'esentation triviale de ${\rm O}_n$, elles ont donc m\^eme param\`etre standard d'apr\`es la proposition \ref{paramarthlangon}.

 \ps\ps
\subsection{Liens avec les fonctions ${\rm L}$}\label{lienfonctionl}

Soient $G$ un $\Z$-groupe semi-simple, $r : 
\widehat{G} \rightarrow \GL_n$ une $\C$-repr\'esentation et $\pi \in \Pi_{\rm disc}(G)$.
D'apr\`es Langlands~\cite[\S 3]{langlandsyale}, le produit eul\'erien 
$${\rm L}(s,\pi,r)=\prod_{p \in {\rm P}} \DET(1-p^{-s} r({\rm c}_p(\pi)))^{-1}$$
est absolument convergent pour ${\rm Re}(s)$ assez grand (voir aussi~\cite[\S 2.5]{shahidi}).\ps\ps 

Si $\pi \in \Pi_{\rm cusp}(\PGL_n)$ on pose ${\rm L}(s,\pi)={\rm
L}(s,\pi,{\rm St}_n)$ (on rappelle que ${\rm St}_n$ d\'esigne la repr\'esentation
tautologique de ${\rm SL}_n$, \S \ref{parconjarthlan}).  On sait depuis Godement et Jacquet que ${\rm
L}(s,\pi)$ admet un prolongement holomorphe \`a $\C$ tout entier, \`a moins
que $n=1$ auquel cas $\pi=1_{\PGL_1}$ et cette fonction ${\rm L}(s,\pi)$
n'est autre que la fonction $\zeta(s)$ de Riemann.  D'apr\`es Jacquet et
Shalika \cite{jasha1}, le produit eul\'erien ${\rm L}(s,\pi)$ est m\^eme
absolument convergent pour ${\rm Re}(s)>1$. Si $\psi(\pi,r)=\oplus_{i=1}^k \pi_i[d_i]$ il d\'ecoule des d\'efinitions que l'on a $${\rm L}(s,\pi,r)=
\prod_{i\in I} \prod_{j=0}^{d_i-1} {\rm
L}(s+j+\frac{1-d_i}{2},\pi_i),$$ et donc ${\rm
L}(s,\pi,r)$ admet \'egalement un prolongement m\'eromorphe \`a $\C$ tout
entier, dont tous les p\^oles sont expliqu\'es par les occurences de la
repr\'esentation triviale dans $\psi(\pi,r)$.\ps\ps 

Le lecteur prendra garde que les normalisations employ\'ees ici font de
$s=\frac{1}{2}$ le centre naturel des \'equations fonctionnelles en jeu. Supposons
par exemple que $\pi \in \Pi_{\rm cusp}({\rm PGL}_2)$ d\'esigne la repr\'esentation associ\'ee
\`a une forme propre normalis\'ee $F=\sum_{n\geq 0} a_n q^n$, de poids $k$, comme au
\S\ref{exparlan}. On a alors ${\rm L}(s-\frac{k-1}{2},\pi)=\sum_{n\geq 1}
\frac{a_n}{n^s}$ pour ${\rm Re}(s)>1$, ainsi que la relation de Hecke
\begin{equation}\label{relheckegl2} (2\pi)^{-s} \Gamma(s) {\rm L}(s-\frac{k-1}{2},\pi) = \int_0^\infty F(it) t^s
\frac{dt}{t}, \, \, \, \, \forall s \in \C. \end{equation}

\subsection{La conjecture de Ramanujan g\'en\'eralis\'ee}\label{ramanujan}

Soit $G$ un $\Z$-groupe semi-simple et soit $\pi \in \Pi(G)$.  On dit que
$\pi$ satisfait la conjecture de Ramanujan, ou qu'elle est {\it temp\'er\'ee}, si\footnote{On ajoute en g\'en\'eral une condition sur
$\pi_\infty$,
conjecturalement automatique, que nous omettons ici.} pour tout $p \in {\rm P}$, les valeurs propres de
${\rm c}_p(\pi)$ dans une (et donc toute) $\C$-repr\'esentation fid\`ele de
$\widehat{G}$ sont toutes de module $1$.  \ps\ps 

La {\it conjecture de Ramanujan
g\'en\'eralis\'ee} affirme que si $\pi \in \Pi_{\rm cusp}(\PGL_n)$ alors
$\pi$ est temp\'er\'ee. L'exemple typique de repr\'esentation qui n'est pas
temp\'er\'ee est la repr\'esentation triviale $1 \in \Pi_{\rm
disc}(\PGL_2)$, les valeurs propres de ${\rm c}_p(1)=e_p$ \'etant $p^{\pm 1/2}$. Plus g\'en\'eralement, la repr\'esentation triviale $1_G$ de $G$ n'est pas temp\'er\'ee si $G \neq 1$ (exemple~\ref{exempletriv}). La
conjecture~\ref{conjarthlan} exprime donc notamment le d\'efaut de la conjecture
de Ramanujan en g\'en\'eral. \ps\ps 

La conjecture de Ramanujan g\'en\'eralis\'ee est toujours ouverte, m\^eme
pour $G=\PGL_2$.  Par des travaux de nombreux auteurs, elle est cependant
connue pour la classe importante des repr\'esentations $\pi \in \Pi_{\rm
cusp}(\PGL_n)$ dites {\it alg\'ebriques r\'eguli\`eres polaris\'ees}
(voir \S VIII.\ref{repgal}).  En g\'en\'eral, on dispose toutefois de l'estim\'ee de
Jacquet-Shalika~\cite{jasha1} : pour tout $\pi\in \Pi_{\rm cusp}(\PGL_n)$,
tout $p \in {\rm P}$ et toute valeur propre $\lambda$ de ${\rm
c}_p(\pi)$, alors $p^{-1/2} < |\lambda| < p^{1/2}$.  \ps\ps

\chapter{Quelques cas de la conjecture d'Arthur-Langlands} \label{contenuchap5}\label{chap7}

\section{Retour sur les relations d'Eichler}\label{retoureichler}

Dans toute cette partie, $g$ et $n$ sont des entiers $\geq 1$ fix\'es tels que $n \equiv 0 \bmod 8$. On
consid\`ere les $\Z$-groupes ${\rm Sp}_{2g}$ et ${\rm O}_n$ (rappelons que ce dernier est le $\Z$-groupe orthogonal du r\'eseau ${\rm E}_n$).\ps\ps

\subsection{Le point de vue de Rallis} \label{releichlerbisrallis}
La s\'erie th\'eta de Jacobi permet de construire une application $\C$-lin\'eaire
naturelle, unique \`a multiplication par un scalaire non nul pr\`es,
$$\vartheta : {\rm M}_{U}({\rm O}_n) \longrightarrow {\rm M}_V({\rm
Sp}_{2g})$$ pour certains couples $(U,V)$ o\`u $U$ est une
$\C$-repr\'esentation irr\'eductible de ${\rm O}_n(\C)$ et $V$ une
$\C$-repr\'esentation irr\'eductible de ${\rm GL}_{g}(\C)$
\cite{kw} \cite{freitagharmts}. Les couples
$(U,V)$ permis seront dit {\it compatibles} et seront d\'ecrits plus bas. Deux cas
particuliers de cette construction ont d'ailleurs d\'ej\`a jou\'e un r\^ole
dans ce m\'emoire : le couple (1,\, $\DET^{n/2}$) au \S \ref{thetasiegel} et le
couple (${\rm H}_{d,g}$,\,$\DET^{n/2+d}$) au \S \ref{i4dpartrialite} pour $2g
\leq n$. Une propri\'et\'e importante de l'application $\vartheta$, d\'ej\`a
discut\'ee dans divers cas {\it loc. cit.}, est
qu'elle entrelace certains op\'erateurs de Hecke de ${\rm O}_n$ et ${\rm
Sp}_{2g}$ (``relations de commutations d'Eichler''). Le but de ce paragraphe
est de rappeler le point de vue de Rallis~\cite{rallis2} sur ces formules. \ps\ps 

Posons $\mathfrak{z}_{{\rm Sp}_{2g}}={\rm Z}({\rm U}({\mathfrak{sp}}_{2g}(\C)))$ et
$\mathfrak{z}_{{\rm O}_n}={\rm Z}({\rm U}({\mathfrak{so}}_{n}(\C)))^{{\rm
O}_n(\C)}$. Rallis construit dans~\cite{rallis2} un morphisme surjectif de $\C$-alg\`ebres $${\rm Ral} : \left\{ \begin{array}{ll}
{\rm H}({\rm O}_n) \otimes \mathfrak{z}_{{\rm O}_n}\rightarrow {\rm H}({\rm Sp}_{2g}) \otimes \mathfrak{z}_{{\rm Sp}_{2g}} & {\rm si\,}\,\,
n>2g,\\ {\rm H}({\rm Sp}_{2g})\otimes \mathfrak{z}_{{\rm Sp}_{2g}}  \rightarrow {\rm H}({\rm O}_n)\otimes \mathfrak{z}_{{\rm O}_n}, & {\rm
si\,}\,\, n\leq 2g.  \end{array} \right.$$ 
respectant les sous-anneaux ${\rm H}_p(\ast) \otimes \C$ ($p$ premier) et $\mathfrak{z}_{\ast}$ de part et d'autre. Il a les propri\'et\'es suivantes : \begin{itemize}\ps\ps 

\item[(i)] Si $(U,V)$ est compatible alors :  ${\rm Inf}_V \circ {\rm Ral}_{|\mathfrak{z}_{{\rm O}_n}} = {\rm Inf}_U$ si $n>2g$, et ${\rm Inf}_U \circ {\rm Ral}_{|\mathfrak{z}_{{\rm Sp}_{2g}}} = {\rm Inf}_V$ sinon. \ps\ps 

\item[(ii)] (La relation de commutation d'Eichler-Rallis) $\vartheta \circ T = {\rm Ral}(T) \circ \vartheta$ pour tout $T \in {\rm H}({\rm O}_n)$ si $n>2g$, et $\vartheta
\circ {\rm Ral}(T) = T \circ \vartheta$ pour tout $T \in {\rm H}({\rm
Sp}_{2g})$ sinon. \ps\ps 
\end{itemize}

\noindent Rajoutons que la condition n\'ecessaire d'admissibilit\'e d'un couple
$(U,V)$ donn\'ee au (i), portant sur leurs caract\`eres infinit\'esimaux,
est suffisante si de plus la restriction de $U$ \`a ${\rm SO}_n(\C)$ n'est
pas irr\'eductible.  Si cette restriction est irr\'eductible, un et un seul
des couples $(U,V)$ et $(U \otimes \DET,V)$ est admissible : voir~\cite{kw}
pour une description d\'etaill\'ee.\ps\ps 

Rallis donne enfin une interpr\'etation du morphisme ${\rm Ral}$ en terme
des isomorphismes de Satake et Harish-Chandra de ${\rm O}_n$ et ${\rm
Sp}_{2g}$, que nous allons \'egalement rappeler.  Si $a \geq 1$, notons ${\rm O}_a$ (resp.  ${\rm SO}_a$) le $\C$-groupe orthogonal (resp. 
sp\'ecial orthogonal) complexe standard \`a $a$ variables. Cette notation est \`a
priori en conflit avec celle des $\Z$-groupes ${\rm O}_n$ et ${\rm
SO}_n$, uniquement d\'efinis pour $n \equiv 0 \bmod 8$, mais c'est sans
incidence puisque lorsque les deux notations co\"incident elles d\'esignent
sur $\C$ le m\^eme objet. Le groupe
${\rm O}_a(\C)$ agit par conjugaison sur $\mathcal{X}({\rm SO}_a)$, et cette
action est non triviale si $a$ est pair. Si $a < b$ sont des entiers $\geq
1$ tels que $a \not\equiv b \bmod 2$, il existe un
$\C$-morphisme $$\rho_{a,b} : {\rm O}_a
\times {\rm SL}_2 \longrightarrow {\rm O}_b,$$ uniquement d\'etermin\'e
modulo conjugaison \`a l'arriv\'ee par ${\rm O}_b(\C)$, tel que la
repr\'esentation ${\rm St} \circ \rho_{a,b}$ soit isomorphe \`a 
la somme directe de la repr\'esentation standard du facteur
${\rm O}_a$ et de la repr\'esentation ${\rm Sym}^{b-a-1} {\rm St}_2$ de
${\rm SL}_2$. Les duaux de Langlands respectifs de ${\rm SO}_{n}$ et ${\rm
Sp}_{2g}$ sont les $\C$-groupes ${\rm SO}_n$ et ${\rm SO}_{2g+1}$. La
discussion de~\cite[\S 6]{rallis2} se traduit en l'\'enonc\'e suivant (on rappelle que $e \in \mathcal{X}({\rm SL}_2)$ d\'esigne l'\'el\'ement d'Arthur, d\'efini au \S VI.\ref{parconjarthlan}):\ps\ps
 
\begin{prop} \label{rallissatakeretour} 
\begin{itemize} \item[(i)] Si
$n>2g$, l'application $\mathcal{X}(\widehat{{\rm Sp}_{2g}}) \rightarrow
\mathcal{X}(\widehat{{\rm O}}_n)$ induite par ${\rm Ral}$ est $x \mapsto
\rho_{2g+1,n}(x,e)$.  \ps\ps  \item[(ii)] Si $n \leq 2g$, l'application
$\mathcal{X}(\widehat{{\rm O}_{n}}) \rightarrow
\mathcal{X}(\widehat{{\rm Sp}}_{2g})$ induite par ${\rm Ral}$ est $x \mapsto
\rho_{n,2g+1}(x,e)$.\end{itemize} \end{prop}

Ce r\'esultat, combin\'e \`a la diagonalisabilit\'e de ${\rm H}({\rm O}_n)$
sur les ${\rm M}_U({\rm O}_n)$, ram\`ene l'\'etude des relations de
commutations d'Eichler \`a des propri\'et\'es des isomorphismes de Satake de
${\rm SO}_n$ et ${\rm Sp}_{2g}$.  Par exemple, on observe que la
proposition V.\ref{commeichlernaif} et la relation V.\eqref{eichlerthetaharmgen}
se d\'eduisent imm\'ediatement des formules VI.\eqref{formulesattp}
et VI.\eqref{satakesp}.  \ps\ps 

De m\^eme, la partie archim\'edienne de la
proposition~\ref{rallissatakeretour} \'eclaircit substantiellement la num\'erologie concernant
les caract\`eres infinit\'esimaux des couples compatibles $(U,V)$.
Illustrons ceci sur le couple (${\rm H}_{d,g}$,\,$\DET^{n/2+d}$). Pour
$g<\frac{n}{2}$ le plus haut poids de ${\rm H}_{d,g}(\R^n)$ est manifestement $d \, \sum_{i=1}^g
\varepsilon_i$, et celui de ${\rm H}_{d,\frac{n}{2}}(\R^n)^\pm$ est $d (\pm
\varepsilon_{\frac{n}{2}} + \sum_{i=1}^{\frac{n}{2}-1} \varepsilon_i)$ (voir \S
IV.\ref{trialitepgso8}, \S\ref{exparlan}), ce qui donne un caract\`ere infinit\'esimal dans
$\mathfrak{so}_n(\C)$ avec pour valeurs propres les $\pm (d+\frac{n}{2}-i)$,
pour $i=1,\dots,g$, et $\pm (\frac{n}{2}-i)$ pour $i=g+1,\dots,\frac{n}{2}$ si
$g<\frac{n}{2}$. D'un autre c\^ot\'e, la repr\'esentation $\pi'_{\det^k}$
(\S VI.\ref{carinf}) a 
un caract\`ere infinit\'esimal dans
$\mathfrak{so}_{2g+1}(\C)$ avec pour valeurs propres $0$ et les $\pm (k-i)$ pour
$i=1,\dots,g$, ce qui est bien compatible avec la
proposition~\ref{rallissatakeretour}. \ps\ps

Disons que les repr\'esentations $\pi_{\rm O} \in \Pi_{\rm disc}({\rm O}_n)$ et $\pi_{\rm Sp}  \in
\Pi_{\rm disc}({\rm Sp}_{2g})$ sont $\vartheta$-correspondantes s'il existe
un couple compatible $(U,V)$ et des formes propres $F \in {\rm M}_U({\rm
O}_n)$ et $G \in {\rm S}_V({\rm Sp}_{2g})$, engendrant respectivement
$\pi_{\rm O}$ et $\pi_{\rm Sp}$, telles que $\vartheta(F)=G$.  


\begin{cor}\label{corparamrallis} Supposons que $\pi_{\rm O}  \in \Pi_{\rm disc}({\rm O}_n)$ et $\pi_{\rm Sp} \in
\Pi_{\rm disc}({\rm Sp}_{2g})$ sont $\vartheta$-correspondantes.  Alors $\psi(\pi_{\rm O},{\rm
St})=\psi(\pi_{\rm Sp},{\rm St}) \oplus [n-2g-1]$ si $n>2g$, et 
$\psi(\pi_{\rm Sp},{\rm St})=\psi(\pi_{\rm O},{\rm St}) \oplus [2g+1-n]$ sinon.\end{cor}\ps\ps

\subsection{Un raffinement : passage aux groupes {\rm
Spin}}\label{rallisspin}

Dans ce paragraphe, nous discutons d'un raffinement de la proposition~\ref{rallissatakeretour} 
probablement bien connu des sp\'ecialistes, mais pour lequel nous n'avons
pas trouv\'e de r\'ef\'erence dans la litt\'erature (il est cependant
implicite dans \cite{rallis1} dans le cas $g=1$). 
 Supposons le couple $(U,V)$ compatible. Observons que si ${\rm M}_U({\rm O}_n) \neq 0$
alors $-1 \in {\rm O}_n(\Z)$ agit trivialement dans $U$, ce que nous
supposons d\`es \`a pr\'esent. Cela entra\^ine que $U$ se factorise en une
repr\'esentation $U'$ de ${\rm PGO}_n(\C)$, et le lemme IV.\ref{inflationmonGO} fournit
un isomorphisme naturel ${\rm H}({\rm O}_n)$-\'equivariant 
$${\rm M}_{U'}({\rm PGO}_n) \isomo {\rm M}_U({\rm O}_n).$$
Autrement dit, ${\rm M}_U({\rm O}_n)$ est naturellement muni d'une action de
l'anneau de Hecke plus gros ${\rm H}({\rm PGO}_n)$. De m\^eme, ${\rm M}_V({\rm Sp}_{2g}(\Z))$
est muni d'une action de ${\rm H}({\rm PGSp}_{2g})$, et on peut se demander comment
l'\'enonc\'e de Rallis s'\'etend \`a ces op\'erateurs. \ps\ps 

Pr\'ecis\'ement, soit $F \in {\rm M}_U({\rm PGO}_n)$ propre pour ${\rm H}({\rm PGO}_n)$ et telle
que $\vartheta(F) \neq 0$. Les relations de Rallis
assurent que $\vartheta(F)$ est propre pour ${\rm H}({\rm Sp}_{2g})$.
D'apr\`es le \S VI.\ref{exemplehecke}, pour tout nombre premier $p$ l'anneau ${\rm H}_p({\rm PGSp}_{2g})$ (resp. ${\rm
H}_p({\rm PGO}_n)$) est engendr\'e par ${\rm H}({\rm Sp}_{2g})$ (resp. ${\rm
H}({\rm O}_n)$) et l'op\'erateur des perestro\"ikas ${\rm K}_p$
correspondant. Or on dispose d'une relation d'Eichler suppl\'ementaire, 
en fait la plus simple de toutes \cite[Thm. 4.5]{freitagtheta}, qui prend la
forme suivante quand $n>2g$ : 
\begin{equation}\label{releichlpere} \vartheta \circ {\rm K}_p =
p^{\frac{g(\frac{n}{2}-g-1)}{2}}\left[ \prod_{i=0}^{\frac{n}{2}-g+1}(p^i+1) \right] {\rm K}_p \circ
\vartheta.\end{equation}
Le cas de cette formule dont nous aurons r\'eellement besoin dans l'application au
th\'eor\`eme~\ref{arthurtrialite} est celui de la proposition
IV.\ref{commeichlergenus1} (\`a comparer avec la
formule~IV.\ref{tradserregl2}). Cela montre que $\vartheta(F)$ est propre pour ${\rm H}({\rm PGSp}_{2g})$. Nous dirons que 
$\pi_{\rm PGO} \in \Pi_{\rm disc}({\rm PGO}_n)$ et $\pi_{\rm PGSp}  \in
\Pi_{\rm disc}({\rm PGSp}_{2g})$ sont $\vartheta$-correspondantes s'il existe
un couple compatible $(U,V)$ et des formes propres $F \in {\rm M}_U({\rm
PGO}_n)$ et $G \in {\rm S}_V({\rm Sp}_{2g}(\Z))$, engendrant respectivement
$\pi_{\rm PGO}$ et $\pi_{\rm PGSp}$, telles que $\vartheta(F)=G$. \ps\ps 
Si $n>2g$, le $\C$-morphisme $\rho_{2g+1,n} : {\rm SO}_{2g+1} \times {\rm SL}_2 \rightarrow {\rm SO}_n$ se rel\`eve en un
$\C$-morphisme $\widetilde{\rho}_{2g+1,n} : {\rm Spin}_{2g+1} \times {\rm SL}_2 \rightarrow {\rm
Spin}_n$. De m\^eme, si $n \leq 2g$ le $\C$-morphisme $\rho_{n,2g+1}$ se
rel\`eve en un morphisme $\widetilde{\rho}_{n,2g+1} : {\rm Pin}_n \times {\rm SL}_2 \rightarrow {\rm
Pin}_{2g+1}$ (voir le \S VI.\ref{exconjal} pour des rappels concernant ${\rm Pin}$). \ps\ps

\begin{prop}\label{generalrallis} Supposons que les repr\'esentations $\pi_{\rm PGO} \in \Pi_{\rm
disc}({\rm PGO}_n)$ et $\pi_{\rm PGSp}  \in \Pi_{\rm disc}({\rm PGSp}_{2g})$
sont $\vartheta$-correspondantes. \begin{itemize} \ps \ps
\item[(i)] Si $n>2g$, alors ${\rm c}(\pi_{\rm
PGO})$ est l'image de $\widetilde{\rho}_{2g+1,n}({\rm c}(\pi_{\rm PGSp}), e)$
par l'application naturelle $\mathcal{X}({\rm
Spin}_n) \rightarrow \mathcal{X}({\rm Pin}_n)$.\ps\ps 
\item[(ii)] Si $n \leq 2g$, alors ${\rm c}(\pi_{\rm
PGSp})=\widetilde{\rho}_{n,2g+1}({\rm c}(\pi_{\rm PGO}), e)$.
\end{itemize}
\end{prop} 

\begin{pf} Supposons $n>2g$. L'\'egalit\'e \`a d\'emontrer vaut apr\`es projection dans
$\mathcal{X}(\widehat{{\rm O}_n})$ d'apr\`es Rallis Prop.~\ref{rallissatakeretour} (i).
D'apr\`es une observation d\'ej\`a donn\'ee ci-dessus, il ne reste qu'\`a la v\'erifier
apr\`es application de l'op\'erateur de Hecke ${\rm K}_p \in {\rm H}({\rm PGO}_n)$, vu par l'isomorphisme de Satake comme une fonction sur ${\rm Spin}_n(\C)_{\rm ss}$ invariante sous l'action de ${\rm Pin}_n(\C)$. D'apr\`es la formule VI.\eqref{formulesattp}, $${\rm Sat}({\rm K}_p)=p^{\frac{\frac{n}{2}(\frac{n}{2}-1)}{4}}([V_{\rm Spin}^+]+[V_{\rm Spin}^-])$$
o\`u $V_{\rm Spin^{\pm}}$ sont les deux repr\'esentations Spin de ${\rm
Spin}_n(\C)$ (conjugu\'ees l'une de l'autre sous ${\rm Pin}_n(\C)$). Mais il est bien connu que la restriction de chacune d'elle \`a ${\rm Spin}_{2g+1}(\C) \times {\rm Spin}_{n-2g-1}(\C)$ est
isomorphe au produit tensoriel des repr\'esentations Spin de chacun des deux facteurs.
Soit $r_a : {\rm SL}_2 \rightarrow {\rm SL}_{2^a}$ la $\C$-repr\'esentation obtenue en relevant la
repr\'esentation irr\'eductible de dimension impaire $2a+1$ dans ${\rm
Spin}_{2a+1}$ puis en composant par la repr\'esentation Spin de ce dernier.
C'est un exercice de v\'erifier que $${\rm
trace}(\,r_a(e_p)\,)\,=\,\prod_{i=1}^a(p^{-\frac{i}{2}}+p^{\frac{i}{2}}).$$
On conclut par un calcul imm\'ediat \`a partir des
formules VI.\eqref{formulesattp}, VI.\eqref{satakesp} et \eqref{releichlpere}. Le
cas $n \leq 2g$ est similaire. La relation d'Eichler montr\'ee par Freitag \cite[Thm.
4.5]{freitagtheta} est 
\begin{equation}\label{releichlpere2} {\rm K}_p \circ \vartheta =
p^{-\frac{g(\frac{n}{2}-g-1)}{2}}\left[ \prod_{i=1}^{g-\frac{n}{2}}(p^{-i}+1)
\right] \vartheta \circ {\rm K}_p.\end{equation}
On conclut comme pr\'ec\'edemment en utilisant que la restriction \`a ${\rm
Spin}_n(\C) \times {\rm Spin}_{2g+1-n}(\C)$ de la repr\'esentation Spin de
${\rm Spin}_{2g+1}(\C)$ est le produit tensoriel de la repr\'esentation
$V_{\rm Spin}^+\oplus V_{\rm Spin}^-$ de ${\rm Spin}_n(\C)$, et de la
repr\'esentation Spin de ${\rm Spin}_{2g+1-n}(\C)$.
\end{pf} \ps\ps

\begin{remark}\label{pinspin} {\rm Observons pour terminer que si $n>2g$, un
\'el\'ement $(c_v)$ dans l'image de la compos\'ee $\mathcal{X}({\rm
Spin}_{2g+1}) \times \mathcal{X}({\rm SL}_2) \rightarrow \mathcal{X}({\rm
Spin}_n) \rightarrow \mathcal{X}({\rm Pin}_n)$ a la priori\'et\'e que pour
tout $v$, la classe de ${\rm Pin}_n(\C)$-conjugaison $c_v$ est en fait une
simple classe de ${\rm Spin}_n(\C)$-conjugaison. Une mani\`ere de le voir,
par exemple pour $v$ premier, est d'observer que si $\gamma \in {\rm
Spin}(V_n)$ a la propri\'et\'e que son image $\gamma'$ dans ${\rm SO}(V_n)$
admet la valeur propre $1$, alors il existe $e \in V_n$ tel que ${\rm
q}(e)=1$ tel que $e \gamma = \gamma e$ (voir le \S VI.\ref{exconjal} pour ces
notations).  En effet, il suffit de choisir un quelconque $e$ dans l'espace
$V_n^{\gamma'=1}$ (qui est non-d\'eg\'en\'er\'e et non nul) tel que ${\rm
q}(e)=1$.  On observe que l'on a $\alpha(\gamma)e\gamma^{-1}=\gamma'(e)=e$ et
$\alpha(\gamma)=\gamma$, puis $\gamma e = e \gamma$.} \end{remark}

\section{$\Pi_{\rm disc}({\rm O}_8)$ et trialit\'e}\label{retourchap4}

La premi\`ere partie du r\'esultat suivant est due \`a
Waldspurger~\cite{waldspurger79}; la proposition V.\ref{polynomeA} en est
une v\'erification \'el\'ementaire dans le cas particulier $k=12$. 
Sa seconde partie est une forme \`a la fois plus pr\'ecise et plus conceptuelle de l'id\'ee principale du
\S \ref{i4dpartrialite}.  On rappelle que nous avons introduit
la repr\'esentation irr\'eductible ${\rm H}_{d,g}(\R^n)$ de ${\rm O}_n(\R)$
au~\S V.\ref{sertharm}.\ps\ps

\begin{thm}\label{arthurtrialite} Soit $\pi \in \Pi_{\rm cusp}(\PGL_2)$ la
repr\'esentation engendr\'ee par une forme propre de ${\rm S}_k({\rm
SL}_2(\Z))$, $k$ \'etant un entier pair $\geq 12$.  \begin{itemize} \ps\ps 
\item[(i)] Il existe $\pi' \in \Pi_{\rm disc}({\rm O}_8)$ telle que
$\pi'_\infty \simeq {\rm H}_{k-4,1}(\R^8)$ et $$\psi(\pi',{\rm St})= {\rm
Sym}^2 \pi \oplus [5].$$ \ps\ps  \item[(ii)] Il existe $\pi'' \in \Pi_{\rm
disc}({\rm O}_8)$ telle que $\pi''_\infty \simeq {\rm H}_{k/2-2,4}(\R^8)$ et
$$\psi(\pi'',{\rm St})=\pi[4].$$\end{itemize} \end{thm}

\begin{pf} Soit $U={\rm H}_{k-4,1}(\R^8)$. D'apr\`es~\cite[Thm. 1]{waldspurger79}, on a
$$\vartheta_{k-4,1}({\rm M}_U({\rm O}_8))={\rm
M}_k({\rm SL}_2(\Z)).$$ D'apr\`es les relations de commutation d'Eichler on
peut donc trouver une forme propre $F \in {\rm M}_U({\rm O}_8)$ dont l'image $G=\vartheta_{k-4,1}(F)$ engendre
$\pi$.  Soit $\pi' \in \Pi_{\rm cusp}(\SL_2)$ la repr\'esentation
engendr\'ee par $G$. Consid\'erons l'isog\'enie $i: {\rm SL}_2(\C)=\widehat{{\rm
PGL}_2}(\C) \rightarrow \widehat{{\rm SL}_2}(\C) = {\rm SO}_3(\C)$. La
proposition~\ref{compsppgsp} et la compatibilit\'e de l'isomorphisme de
Satake aux isog\'enies assurent que ${\rm c}(\pi')=i({\rm c}(\pi))$. Mais
${\rm St} \circ i$ n'est autre que la repr\'esentation ${\rm Sym}^2 {\rm St}_2$ de ${\rm
SL}_2(\C)$, de sorte que ${\rm St}({\rm c}(\pi'))={\rm Sym}^2 {\rm c}(\pi)$. Le (i) se d\'eduit
alors du corollaire~\ref{corparamrallis}.  \ps\ps  

V\'erifions le (ii). Quitte \`a modifier $F$, on peut supposer que $F_0=F$ est
propre pour ${\rm H}({\rm PGO}_8)$, comme au \S\ref{rallisspin}. Soit $\pi_0 \in \Pi_{\rm disc}({\rm
PGO}_8)$ la repr\'esentation qu'elle engendre et soit $U'=U \otimes
\nu^{k/2-2}$. Soit $F_1$ l'image de $F_0$ par l'application naturelle  ${\rm res} : {\rm M}_{U'}({\rm PGO}_8)  \rightarrow {\rm M}_{U'}({\rm
PGSO}_8)$ (\S IV.\ref{pgpso}). Cette application \'etant injective ${\rm
H}({\rm PGSO}_8)$-\'equivariante d'apr\`es
{\it loc. cit.}, la forme $F_1$ est non nulle et si $\pi_1 \in
\Pi_{\rm disc}({\rm PGSO}_8)$ d\'esigne la repr\'esentation engendr\'ee,
alors ${\rm c}(\pi_0) \in \mathcal{X}({\rm Pin}_8)$ est l'image de ${\rm c}(\pi_1)
\in \mathcal{X}({\rm Spin}_8)$ par l'homomorphisme naturel. Cette derni\`ere
propri\'et\'e d\'etermine d'ailleurs uniquement ${\rm c}(\pi_1)$ par la
remarque~\ref{pinspin}, de sorte que la proposition~\ref{generalrallis}
s'\'ecrit 
\begin{equation} \label{genrallis2} {\rm c}(\pi_1)= \widetilde{\rho}_{3,8}({\rm c}(\pi),e) \end{equation}
\par
\ps
Consid\'erons maintenant, suivant le~\S V.\ref{trialitepgso8}, un automorphisme de trialit\'e $\tau$ de ${\rm PGSO}_8$ d\'efini
au moyen d'une structure d'octonions de Coxeter sur ${\rm E}_8$. En
particulier, $(U')^\tau$ est isomorphe \`a la repr\'esentation $V={\rm
H}_{k/2-2,4}^{\pm} \otimes \nu^{k-4}$ d'apr\`es le lemme V.\ref{lemmetrialiteinfini}. Soit
${\rm tri} : {\rm M}_{U'}({\rm PGSO}_8) \isomo {\rm M}_{V}({\rm PGSO}_8)$ l'isomorphisme
not\'e $(\tau^\ast)^{-1}$ {\it loc. cit.}. On pose enfin 
$$F_2={\rm tri}(F_1), \, \, \, F_3={\rm ind}\, (F_2)
\in {\rm  M}_{{\rm Ind}(V)}({\rm PGO}_8)\, \, \, {\rm  et}\, \, \,  F_4 = \mu^\ast(F_3) \in {\rm M}_{\rm Ind V}({\rm O}_8)$$
(\S IV.\ref{fauton}, \S V.\ref{pgpso}).
Ces fonctions sont non nulles et propres, et engendrent donc des
repr\'esentations automorphes $\pi_2$, $\pi_3$ et $\pi_4$ des $\Z$-groupes respectifs ${\rm PGSO}_8$, ${\rm
PGO}_8$
et ${\rm O}_8$. \ps

La compatibilit\'e de l'isomorphisme de Satake aux
isog\'enies assure que ${\rm c}(\pi_2)$ est l'image par $\tau^{\pm 1}$ de
${\rm c}(\pi_1)$, puis que ${\rm c}(\pi_4)$ est l'image de ${\rm c}(\pi_2)$ par l'homomorphisme naturel
$\eta : {\rm Spin}_8(\C) \rightarrow {\rm SO}_8(\C)$. Mais il est bien connu que la repr\'esentation
${\rm St} \circ \eta \circ \tau^{\pm 1}$ n'est autre que la repr\'esentation $V_{\rm Spin^{\pm}}$ de
${\rm Spin}_8$. On a d\'ej\`a dit que la restiction de $V_{\rm Spin^{\pm}}$ \`a ${\rm Spin}_3 \times {\rm Spin}_5$
est le produit tensoriel des repr\'esentations {\rm Spin} de ${\rm Spin}_3 \simeq
{\rm SL}_2$ (de dimension $2$) et ${\rm Spin}_5 \simeq {\rm Sp}_4$ (de
dimension $4$); en particulier elle ne d\'epend pas du signe $\pm$. Or la repr\'esentation ${\rm Sym}^4
\, {\rm St}_2$ de ${\rm SL}_2$, vue dans ${\rm
SO}_5$, puis relev\'ee dans ${\rm Spin}_5 \simeq {\rm Sp}_4$ et compos\'ee avec la
repr\'esentation standard de ${\rm Sp}_4$, est la repr\'esentation ${\rm
Sym}^3\, {\rm St}_2$. Nous avons donc d\'emontr\'e la suite d'\'egalit\'es
$$\pi[4]=\psi(\pi_1,V_{\rm Spin}^{\pm})=\psi(\pi_2,{\rm St} \circ \eta)=\psi(\pi_4,{\rm
St}),$$
et la repr\'esentation $\pi^{''}=\pi_4$ convient, ce qui conclut. \end{pf}
\ps\ps
Observons que le th\'eor\`eme de Gelbart-Jacquet (\S VI.\ref{pargelbartjacquet}) entra\^ine que le couple $(\pi',{\rm St})$ satisfait la conjecture d'Arthur-Langlands. Il est de plus manifeste que le couple $(\pi'',{\rm St})$ satisfait la conjecture d'Arthur-Langlands.\ps\ps

Donnons une seconde formulation du r\'esultat pr\'ec\'edent. On rappelle que
l'homomorphisme ${\rm SO}_8 \rightarrow {\rm PGSO}_8$ d\'etermine par
dualit\'e de Langlands un
$\C$-morphisme $\eta : \widehat{{\rm PGSO}_8} \rightarrow \widehat{{\rm
SO}_8}$; les $3$ repr\'esentations irr\'eductibles de dimension $8$ de $\widehat{{\rm
PGSO}_8}(\C)
\simeq {\rm Spin}_8(\C)$ sont donc ${\rm St} \circ \eta$, $V_{\rm Spin}^+$
et $V_{\rm Spin}^-$. \ps\ps

\begin{thm}\label{arthurtrialitev2} Soit $\pi \in \Pi_{\rm cusp}(\PGL_2)$ la
repr\'esentation engendr\'ee par une forme propre de ${\rm S}_k({\rm    
SL}_2(\Z))$, $k$ \'etant un entier pair $\geq 12$. Il existe $\pi' \in \Pi_{\rm
disc}({\rm PGSO}_8)$ telle que 
$$\psi(\pi',{\rm St} \circ \eta)=\pi[4], \, \, \, \psi(\pi',V_{\rm Spin}^+)=\pi[4] \,
\, \, {\rm et}\, \, \, \psi(\pi',V_{\rm Spin}^-)={\rm Sym}^2\pi \oplus [5].$$
\end{thm}

\begin{pf} La repr\'esentation $\pi_2$ de la d\'emonstration du
th\'eor\`eme~\ref{arthurtrialite} satisfait $\psi(\pi_2,{\rm St} \circ
\eta)=\pi[4]$ et $$\{\psi(\pi_2,V_{\rm Spin}^+),\psi(\pi_2,V_{\rm
Spin}^-)\}=\{\pi[4],{\rm Sym}^2\pi \oplus [5]\}.$$ Soit $\pi'_2$ la
repr\'esentation de ${\rm PGSO}_8$ engendr\'ee par l'image par la sym\'etrie $s$
d'une forme pour ${\rm PGSO}_8$ engendrant $\pi_2$ (\S V.\ref{pgpso}). Alors l'une
des deux repr\'esentations $\pi_2$ et $\pi'_2$ convient pour $\pi'$. 
\end{pf}

Le principe de la d\'emonstration du th\'eor\`eme~\ref{arthurtrialite} est de port\'ee plus large,
et s'applique notamment aux s\'eries th\^eta de genre sup\'erieur. Il permet
de produire des repr\'esentations de ${\rm O}_8$ de param\`etres
de Langlands standards int\'eressants, fonctions de ceux des \'el\'ements de $\Pi_{\rm
cusp}({\rm PGSp}_{2g})$ pour $1\leq g\leq 3$. \ps\ps

\begin{thm}\label{trialitegenus2} Supposons que $\pi \in \Pi_{\rm cusp}({\rm PGSp}_{2g})$ admette un $\theta$-corres\-pondant dans $\Pi_{\rm disc}({\rm PGO}_8)$.\ps\ps 
\begin{itemize}
\item[(i)] Supposons $g=2$. Soient $V_4$ et $V_5$ les repr\'esentations irr\'eductibles de ${\rm Sp}_4(\C)=\widehat{{\rm PGSp}_4}(\C)$ de dimensions respectives $4$ et $5$, i.e. $V_4$ est la repr\'esentation standard et $\Lambda^2 V_4 \simeq  V_5 \oplus 1$.  Alors il existe : \ps\ps 

-- $\pi' \in \Pi_{\rm disc}({\rm SO}_8)$ telle que $\psi(\pi',{\rm St})= \psi(\pi,V_5)\oplus [3]$, \ps\ps  
-- $\pi'' \in \Pi_{\rm disc}({\rm SO}_8)$ telle que $\psi(\pi'',{\rm St})=\psi(\pi,V_4)[2]$. \ps\ps 

\item[(ii)] Supposons $g=3$. Soient $V_{\rm Spin}$ la repr\'esentation Spin de ${\rm Spin}_7(\C)=\widehat{{\rm PGSp}_6}(\C)$ et $V_7$ sa repr\'esentation naturelle de dimension $7$. Alors il existe : \ps\ps 
-- $\pi' \in \Pi_{\rm disc}({\rm SO}_8)$ telle que $\psi(\pi',{\rm St})= \psi(\pi,V_7)\oplus [1]$, \ps\ps  
-- $\pi'' \in \Pi_{\rm disc}({\rm SO}_8)$ telle que $\psi(\pi'',{\rm St})=\psi(\pi,V_{\rm Spin})$. \ps\ps 

\end{itemize}
\end{thm}

\begin{pf} L'existence de $\pi'$ dans les deux cas est classique et d\'ecoule du corollaire~\ref{corparamrallis}. En ce qui concerne celle de $\pi''$,  sa d\'emonstration est tr\`es similaire \`a celle de la repr\'esentation du m\^eme nom dans le th\'eor\`eme~\ref{arthurtrialite}, 
elle est donc laiss\'ee en exercice pour le lecteur. Par exemple dans le cas (i), on montre d'abord l'existence de $\pi_0 \in \Pi_{\rm disc}({\rm PGSO}_8)$ telle que $\psi(\pi_0,V_{\rm Spin}^\pm)=\psi(\pi,V_4)[2]$ et $\psi(\pi_0,{\rm St} \circ \eta)=\psi(\pi,V_5)\oplus [3]$; l'application de la trialit\'e \`a $\pi_0$ conduit alors \`a la repr\'esentation $\pi''$.  
\end{pf}
\ps\ps
\begin{remark}{\rm (Travaux de B\"ocherer) \label{rembocherer} Disons un mot
sur l'hypoth\`ese du th\'eor\`eme et la question reli\'ee de la
surjectivit\'e de l'application $\vartheta$ en g\'en\'eral, qui est un probl\`eme
classique remontant \`a Eichler ({\it probl\`eme de la base d'Eichler}).  On dispose du
r\'esultat marquant suivant d\^u \`a B\"ocherer
\cite{bocherertheta2}\cite{bochererbowdoin}, g\'en\'eralisant le travail de
Waldspurger sus-cit\'e pour $g=1$ : pour $d>0$ l'application $$\vartheta_{d,g} : {\rm M}_{{\rm
H}_{d,g}(\R^n)}({\rm O}_n) \rightarrow {\rm S}_{\frac{n}{2}+d}({\rm
Sp}_{2g}(\Z))$$ est surjective d\`es que $n>4g$ (voir
\'egalement~\cite{bocherertheta1} pour le cas $d=0$ ainsi que la remarque
VIII.\ref{rembochamelioration}). De mani\`ere plus
pr\'ecise, B\"ocherer donne une condition n\'ecessaire et suffisante pour qu'un forme propre $F \in {\rm S}_{\frac{n}{2}+d}({\rm
Sp}_{2g}(\Z))$ soit dans l'image de $\vartheta_{d,g}$, quand $n\geq 2g$. Elle porte
sur la fonction ${\rm L}(s,\pi,{\rm St})$, o\`u $\pi \in \Pi_{\rm cusp}({\rm
Sp}_{2g})$ est engendr\'ee par $F$, dont on sait qu'elle admet un
prolongement m\'eromorphe \`a
$\C$ tout entier (voir \S \ref{complboc}).  Si $n > 2g$ (resp. 
$n = 2g$), il montre que $F$ est dans l'image de $\vartheta_{d,g}$ si et
seulement si ${\rm L}(s,\pi,{\rm St})$ ne s'annule pas en $s=\frac{n}{2}-g$ (resp. 
si et seulement si ${\rm L}(s,\pi,{\rm St})$ admet un p\^ole simple en
$s=1$) : voir~\cite[Thm. $4_1$, Thm. 5]{bochererbowdoin}. Cette condition est automatiquement satisfaite si $n>4g$.
B\"ocherer a \'egalement \'etudi\'e la question de
l'injectivit\'e de $\vartheta_g=\vartheta_{0,g}$, pour laquelle il obtient des crit\`eres du
m\^eme type~\cite{bochererkern}.} \end{remark}
\ps\ps
\noindent Revenons \`a l'\'enonc\'e du th\'eor\`eme~\ref{arthurtrialitev2}.\ps\ps

\begin{cor}\label{cortrikeda} Supposons que $\pi \in \Pi_{\rm cusp}(\PGL_2)$ et $\pi' \in
\Pi_{\rm disc}({\rm PGSO}_8)$ satisfont les hypoth\`eses et conclusions du
th\'eor\`eme~\ref{arthurtrialitev2}. Supposons de plus que $\pi'$ poss\`ede
un $\vartheta$-correspondant $\pi'' \in \Pi_{\rm cusp}(\PGSp_8)$,
c'est-\`a-dire qu'il existe $F \in {\rm M}_{{\rm H}_{k/2-2,4}^\pm(\R^8)}({\rm
PGSO}_8) \simeq {\rm M}_{{\rm H}_{k/2-2,4}(\R^8)}({\rm PGO}_8)$ engendrant
$\pi'$ et telle que $\vartheta(F)$ est un \'el\'ement non nul de ${\rm S}_{k/2+2}({\rm
Sp}_8(\Z))$, $\pi''$ d\'esignant alors la repr\'esentation engendr\'ee par
$\vartheta(F)$. Alors 
$$\psi(\pi'', V_{\rm St})=\pi[4] \oplus [1] \, \, \, \, \, {\rm et}\, \,
\, \, \, \psi(\pi'',{\rm V}_{\rm Spin})=\pi[4]\oplus {\rm Sym}^2 \pi \oplus [5].$$
\end{cor}

\begin{pf} Cela d\'ecoule imm\'ediatement du th\'eor\`eme~\ref{arthurtrialitev2} et des
relations d'Eichler-Rallis raffin\'ees (Proposition~\ref{generalrallis}
dans le cas $g=4=n/2$), car la restriction \`a ${\rm Spin}_8 \rightarrow {\rm Spin}_9$ de la
repr\'esentation Spin de ${\rm Spin}_9(\C)$ est la repr\'esentation $V_{\rm
Spin}^+ \oplus V_{\rm Spin}^-$ de ${\rm Spin}_8(\C)$.
\end{pf}
\ps\ps
Lorsque $\pi$ est engendr\'ee par $\Delta \in {\rm S}_{12}({\rm SL}_2(\Z))$,
nous avons v\'erifi\'e dans la proposition~V.\ref{thetanonnul44} que
l'hypoth\`ese sur $\pi'$ est satisfaite (cela devrait pouvoir \'egalement se d\'eduire d'une variante harmonique de~\cite{bochererkern}, voir~\cite[\S XI]{bochererbowdoin}). On rappelle que ${\rm S}_8({\rm
Sp}_8(\Z))$ est de dimension $1$ engendr\'e par la forme de Schottky $J$
(\S \ref{casdim16}). Nous noterons $\Delta_{11} \in
\Pi_{\rm cusp}({\rm PGL_2})$ la repr\'esentation engendr\'ee par la forme
modulaire $\Delta$. 
\ps\ps
\begin{cor}\label{corpij} \begin{itemize} \item[(i)] Si $\pi_J \in \Pi_{\rm cusp}({\rm PGSp}_8)$ d\'esigne la repr\'esentation
engendr\'ee par la forme de Schottky, alors $\psi(\pi_J,V_{\rm St})=\Delta_{11}[4]\oplus
[1]$ et
$$\psi(\pi_J,V_{\rm Spin})=\Delta_{11}[4]
\oplus {\rm Sym}^2 \Delta_{11} \oplus [5].$$\ps\ps 
\item[(ii)] Soit $\pi \in \Pi_{\rm disc}({\rm PGSO}_{16})$ l'unique repr\'esentation non-triviale
telle que $\pi_\infty=\C$. Alors $\psi(\pi,V_{\rm St})=\Delta_{11}[4] \oplus [7] \oplus [1]$
et
$$\hspace{-9mm}{\small \psi(\pi,V_{\rm Spin}^{\pm})=\,\psi(\pi_J,V_{\rm Spin}) \,\oplus \, \psi(\pi_J,V_{\rm Spin})[7]}$$
\end{itemize}
\end{cor}

\begin{pf} L'assertion (i) r\'esulte du corollaire~\ref{cortrikeda} et de la
discussion pr\'ec\'edant l'\'enonc\'e.  L'assertion (ii) s'en d\'eduit,
\'etant donn\'e la relation $J=\vartheta_4({\rm E}_8 \oplus {\rm
E}_8)-\vartheta_4({\rm E}_{16})$ (\S~IV.\ref{casdim16}) et la
proposition~\ref{generalrallis} (ii). On remarquera que si le $\C$-morphisme $g : {\rm SL}_2 \rightarrow
{\rm SO}_7$ est tel que ${\rm St} \circ g \simeq {\rm Sym}^6 {\rm St}_2$, et si $f  : {\rm SL}_2 \rightarrow {\rm Spin}_7$ rel\`eve $g$, 
alors la restriction \`a $f$ de la repr\'esentation Spin de
${\rm Spin}_7$ est isomorphe \`a ${\rm Sym}^6 {\rm St}_2 \oplus 1$ (voir par exemple~\cite[\S 7]{grossminuscule}).  
\end{pf}

\begin{cor}\label{perdim16} Pour tout nombre premier $p$, le nombre des perestro\"ikas de ${\rm E}_8
\oplus {\rm E}_8$ relativement \`a $p$ qui sont isomorphes \`a ${\rm E}_{16}$ est 
$$ \frac{405}{691}\left( \prod_{i=0}^3(p^i+1) \right)
(p^{11}+p^7+p^6+p^5+p^4+1+\tau(p))
(p^{11}+1 - \tau(p)).$$
\end{cor}

\section{Quelques cons\'equences des travaux d'Ikeda et B\"ocherer}\label{ikedabocherer}

Comme nous l'avions d\'ej\`a expliqu\'e au \S \ref{casdim16}, le premier
point du corollaire \ref{corpij} est \'egalement cons\'equence du th\'eor\`eme suivant d\^u
\`a Ikeda~\cite{ikeda1} (d\'emonstration de la {\it conjecture de Duke-Imamo\u{g}lu}). Il \'etend un r\'esultat de
Andrianov, Maass et Zagier dans le cas du genre $g=2$ (d\'emonstration de la {\it conjecture de Saito-Kurokawa}, voir~\cite{kurokawa}\cite{zagiersk}\cite{eichlerzagier}).

\begin{thm}\label{thmikeda2} {\rm \cite{ikeda1}} Soit $\pi \in \Pi_{\rm
cusp}({\rm PGL}_2)$ la repr\'esentation engendr\'ee par une forme propre de
poids $k$ pour ${\rm SL}_2(\Z)$, et soit $g \geq 1$ un entier tel que $k
\equiv g \bmod 4$, alors il existe une repr\'esentation $\pi' \in \Pi_{\rm
cusp}({\rm Sp}_{2g})$, engendr\'ee par une forme modulaire de Siegel
scalaire de poids $\frac{k+g}{2}$ pour ${\rm Sp}_{2g}(\Z)$, telle que
$\psi(\pi',{\rm St}) = \pi[g] \oplus [1]$.  \end{thm}

Supposons $\pi$ et $\pi'$ comme dans l'\'enonc\'e ci-dessus. \'Etant donn\'e que l'on conna\^it la fonction ${\rm L}(s,\pi',{\rm St})$, les r\'esultats sus-cit\'es de B\"ocherer (Remarque \ref{rembocherer}, \cite{bochererbowdoin}) donnent une condition n\'ecessaire et suffisante pour que $\pi'$ admette un $\vartheta$-correspondant. 

\begin{thm}\label{corikedabocherer} Soient $k, g$ et $n$ des entiers pairs
non nuls tels que $k \equiv g \bmod 4$, $n \equiv 0 \bmod 8$, et $2g \leq n
\leq k+g$. Soient $\pi \in \Pi_{\rm cusp}({\rm PGL}_2)$ la repr\'esentation engendr\'ee par une forme
propre de poids $k$ pour ${\rm SL}_2(\Z)$ et $\pi' \in \Pi_{\rm cusp}({\rm
Sp}_{2g})$ satisfaisant les conclusions du th\'eor\`eme d'Ikeda relativement
\`a $\pi$.  \begin{itemize}\ps\ps 

\item[(i)] La repr\'esentation $\pi'$ admet un $\vartheta$-correspondant
$\pi'' \in \Pi_{\rm disc}({\rm O}_n)$ tel que $\pi''_\infty \simeq {\rm
H}_{\frac{k+g-n}{2},g}(\R^n)$ si, et seulement si, $n>3g$ ou ${\rm L}(\frac{1}{2},\pi)
\neq 0$.  \ps\ps  \item[(ii)] Supposons $n>3g$ ou ${\rm L}(\frac{1}{2},\pi) \neq
0$ {\rm (}auquel cas $k \equiv g \equiv 0 \bmod 4${\rm )}.  Si $n>2g$ {\rm
(}resp.  $n=2g${\rm )}, il existe $\pi'' \in \Pi_{\rm cusp}({\rm O}_n)$
telle que $\psi(\pi'',{\rm St})=\pi[g]\oplus [n-2g-1] \oplus [1]$ {\rm
(}resp.  $\psi(\pi'',{\rm St})=\pi[g]${\rm )} et $\pi''_\infty \simeq {\rm
H}_{\frac{k+g-n}{2},g}(\R^n)$.  \ps\ps  \end{itemize} \end{thm}

\begin{pf} Soit $\mathcal{D} = \frac{1}{2}\Z - \Z$ (les demi-entiers non
entiers).  Comme on a $\psi(\pi',{\rm St})=\pi[g]\oplus [1]$, la fonction ${\rm
L}(s,\pi',{\rm St})$ est le produit de la fonction $\zeta(s)$ de Riemann par
les fonctions ${\rm L}(s+j,\pi)$ o\`u $j$ parcourt les \'el\'ements de
$\mathcal{D}$ tels que $|j|\leq \frac{g-1}{2}$ (noter que $g \equiv 0 \bmod
2$). On consid\`ere le couple compatible $({\rm
H}_{\frac{k+g-n}{2},g}(\R^n),\DET^{\frac{k+g}{2}})$, qui est bien d\'efini car $2g \leq n \leq k+g$. Si $n>2g$ (resp.  $n=2g$), B\"ocherer montre que $\pi'$ admet un
$\vartheta$-correspondant comme dans l'\'enonc\'e si, et seulement si, ${\rm
L}(\frac{n}{2}-g,\pi',{\rm St}) \neq 0$ (resp.  ${\rm L}(s,\pi',{\rm St})$
admet un p\^ole simple en $s=1$).  Comme $\zeta(s)$ admet un p\^ole simple
en $s=1$ et ne s'annule pas si ${\rm Re}(s)>1$, il est \'equivalent de
demander que ${\rm L}(s,\pi) \neq 0$ pour tout $s \in \mathcal{D}$ tel que
$|s - \delta-\frac{n}{2}+g|\leq \frac{g-1}{2}$, o\`u $\delta=1$ si $n=2g$ et
$0$ sinon.  \ps\ps

Si $s \in \mathcal{D}$ et $s \neq \frac{1}{2}$ alors ${\rm L}(s,\pi) \neq
0$.  En effet, pour ${\rm Re}(s)>1$ c'est une cons\'equence de la convergence
absolue du produit eul\'erien d\'efinissant ${\rm L}(s,\pi)$ (par exemple,
d'apr\`es Deligne, ou Rankin-Selberg). En g\'en\'eral, on
rappelle que d'apr\`es Hecke la fonction
$\xi(s,\pi)=(2\pi)^{-s-\frac{k-1}{2}} \Gamma(s+\frac{k-1}{2}){\rm L}(s,\pi)$
est une fonction enti\`ere de $s$ satisfaisant (voir \S\ref{lienfonctionl})
\begin{equation} \label{eqfpi2} \xi(1-s,\pi)=i^k \xi(s,\pi).  \end{equation}
On conclut car la fonction $\Gamma(s)$ n'admet pas de p\^ole en les
\'el\'ements de $\mathcal{D}$.  Pour conclure le (i), on observe que
$|\frac{1}{2}-\delta-\frac{n}{2}+g|\leq \frac{g-1}{2}$ si, et seulement si,
$n \leq 3g$.  \ps\ps
V\'erifions enfin l'assertion (ii). Sous les hypoth\`eses de l'\'enonc\'e, on dispose d'un $\vartheta$-correspondant $\pi'' \in \Pi_{\rm disc}({\rm O}_n)$ de la repr\'esentation $\pi'$, d'apr\`es le point (i). Appliquons le corollaire \ref{corparamrallis}. Sous l'hypoth\`ese $n>2g$, il entra\^ine $\psi(\pi'',{\rm St}) = \psi(\pi',{\rm St}) \oplus [n-2g-1]$, ce que l'on voulait. Dans le cas d'\'egalit\'e $n=2g$, il s'\'ecrit $$\psi(\pi'',{\rm St})\, \oplus\,[1]\,=\, \psi(\pi',{\rm St})\,=\,\pi[g] \,\oplus \,[1],$$ ce qui \'equivaut manifestement \`a $\psi(\pi'',{\rm St}) = \pi[g]$. On observe enfin que l'hypoth\`ese ${\rm L}(\frac{1}{2},\pi) \neq 0$ entraine $k \equiv 0 \bmod
4$ d'apr\`es l'\'equation fonctionnelle~\eqref{eqfpi2}.  \end{pf}

\begin{remark}\label{nonvanish}{\rm  Supposons que $\pi \in \Pi_{\rm
cusp}(\PGL_2)$ est engendr\'ee par une forme propre $F$ de poids $k \equiv 0
\bmod 4$ pour ${\rm SL}_2(\Z)$.  Si $12 \leq k \leq 20$ alors $F=\Delta
\,\vartheta_1({\rm E}_8)^{\frac{k-12}{4}}$; comme cette derni\`ere ne prend
que des valeurs strictement positives sur l'axe imaginaire, il vient que
$\Gamma(s+\frac{k-1}{2}){\rm L}(s,\pi)>0$ pour tout $s \in \R$ (formule
VI.\eqref{relheckegl2}).  En particulier ${\rm L}(\frac{1}{2},\pi) \neq 0$. 
Il semblerait que l'on ne connaisse aucun exemple pour lequel ${\rm
L}(\frac{1}{2},\pi)=0$ : voir \`a ce sujet~\cite{conreyfarmer}, qui le
v\'erifient pour tout $k \leq 500$ (d'ailleurs en lien avec le
th\'eor\`eme~\ref{corikedabocherer}, qu'ils n'explicitent pas toutefois).}
\end{remark}

Il est int\'eressant de confronter le th\'eor\`eme~\ref{corikedabocherer}
aux r\'esultats pr\'ec\'edents.  Tout d'abord, si on l'applique \`a $n=16$,
$g=4<n/3$ et $k=12$ (de sorte que $d=0$), on d\'eduit ais\'ement du (i) que
$\vartheta_4({\rm E}_8\oplus {\rm E}_8) - \vartheta_4({\rm E}_{16}) \in {\rm
M}_8({\rm Sp}_8(\Z))$ est une forme parabolique non nulle de param\`etre
standard $\Delta_{11}[4] \oplus [1]$.  Cet argument ``massue'' red\'emontre donc
\`a la fois la conjecture de Witt (Chap.  IV \eqref{egaltheta1}) et
l'assertion concernant $\psi(\pi,V_{\rm St})$ dans le (ii) du
corollaire~\ref{corpij} (et donc le th\'eor\`eme IV.\ref{vpnontrivdim16} !). 
De m\^eme, le th\'eor\`eme~\ref{corikedabocherer} s'applique pour
$n=24=k+g$, ce qui se produit pour les $5$ couples $$(k,g) \in \{(12,12),
(16,8), (18,6), (20,4), (22,2)\},$$ en prenant pour $\pi$ la
repr\'esentation engendr\'ee par l'unique forme normalis\'ee de poids $k$
pour ${\rm SL}_2(\Z)$ quand $k \leq 22$.  Nous noterons $\Delta_{k-1} \in
\Pi_{\rm cusp}(\PGL_2)$ cette repr\'esentation.  L'hypoth\`ese du th\'eor\`eme est satisfaite dans
les trois derniers cas car $g < 24/3=8$, ainsi que pour les deux premiers
car ${\rm L}(\frac{1}{2},\Delta_{k-1}) \neq 0$ si $k=12$ ou $16$.

\begin{cor}\label{cor24ikedabocherer} Pour tout $k \in \{12,16,18,20,22\}$, il existe $\pi \in \Pi_{\rm disc}({\rm O}_{24})$ telle que $\pi_\infty=\C$ et $\psi(\pi,{\rm St})=\Delta_{k-1}[24-k] \oplus [2k-25] \oplus [1]$ si $k>12$, 
$\psi(\pi,{\rm St})=\Delta_{11}[12]$ si $k=12$.

\end{cor}

En revanche, on observe que le cas $n=2g=8$ conduit \`a une version
affablie, et paradoxalement plus ch\`ere, du
th\'eor\`eme~\ref{arthurtrialite} (ii), puisqu'il ne le d\'emontre que sous
l'hypoth\`ese suppl\'ementaire ${\rm L}(\frac{1}{2},\pi) \neq 0$ (et en
particulier $k \equiv 0 \bmod 4)$.  Le cas $k \equiv 2 \bmod 4$ du
th\'eor\`eme~\ref{arthurtrialite} (ii) appara\^it donc particuli\`erement
int\'eressant de ce point de vue.  Ce ph\'enom\`ene un peu troublant, ainsi
que la num\'erologie d'apparence un peu particuli\`ere de l'\'enonc\'e du
th\'eor\`eme~\ref{corikedabocherer}, sera tr\`es nettement \'eclair\'e
lorsque nous expliquerons les r\'esultats d'Arthur au chapitre suivant (voir
notamment les \S VIII.\ref{formexpliciteson0} et \S \ref{comptheta}).  \ps\ps 

Terminons ce paragraphe par un dernier exemple montrant qu'en g\'en\'eral,
l'inclusion $\vartheta_g(\C[{\rm X}_n]) \subset {\rm M}_{\frac{n}{2}}({\rm
Sp}_{2g}(\Z))$ est stricte.

\begin{cor}\label{cex32} L'application $\vartheta_{14} : \C[{\rm X}_{32}] \rightarrow {\rm M}_{16}({\rm Sp}_{28}(\Z))$ n'est pas surjective. \end{cor}

\begin{pf} En effet, consid\'erons le cas $n=32$, $k=18$ et $g=14$ (et donc
encore $d=0$).  Le th\'eor\`eme d'Ikeda assure l'existence d'une forme
propre $F$ dans ${\rm S}_{16}({\rm Sp}_{28}(\Z))$ engendrant une
repr\'esentation de param\`etre standard $\Delta_{17}[14] \oplus [1]$.  On a
$n<3g$ et ${\rm L}(\frac{1}{2},\Delta_{17})=0$ car $18 \equiv 2 \bmod 4$, de
sorte que le th\'eor\`eme~\ref{corikedabocherer} (i) affirme que $F \notin
{\rm Im}(\vartheta_{14})$.  \end{pf}

Cet exemple semble avoir \'echapp\'e \`a Nebe et Venkov \cite[\S
2.2]{nebevenkov}, qui pr\'esentent la question de l'\'egalit\'e entre
$\vartheta_g(\C[{\rm X}_n])$ et ${\rm M}_{\frac{n}{2}}({\rm Sp}_{2g}(\Z))$,
pour tout $n \equiv 0 \bmod 8$ et tout $g\geq 1$, comme un probl\`eme
ouvert.  Comme nous le verrons au \S VIII.\ref{formexpliciteson0}, la
th\'eorie d'Arthur sugg\`ere en fait qu'il n'existe pas de $\pi \in \Pi_{\rm
disc}({\rm O}_{32})$ tel que $\psi(\pi,{\rm St})=\Delta_{17}[14] \oplus [3]
\oplus [1]$.

\section{Une table des premiers \'el\'ements de $\Pi_{\rm disc}({\rm
SO}_8)$}\label{tableexempleso8}

D'apr\`es~\cite[Ch. 2]{chrenard2}, on dispose d'une formule pour la
dimension de ${\rm M}_{U_\lambda}({\rm SO}_8)$ en fonction du plus haut
poids $\lambda=\sum_{i=1}^4 m_i
\varepsilon_i$ de la repr\'esentation $U_\lambda$ (\S VI.\ref{exparlan}).  Pour
des petites valeurs de $\lambda$, par exemple d\`es que $m_1(\lambda):=m_1$
est $\leq 9$, on constate par inspection que ces dimensions sont $\leq 1$,
et m\^eme presque toujours nulles~\cite[App.  C.  Table 2]{chrenard2}. 
Quand cette dimension vaut $1$, il existe donc une unique repr\'esentation
$\pi \in \Pi_{\rm disc}({\rm SO}_8)$ telle que $\pi_\infty \simeq
U_{\lambda}$, notons-la $\pi(\lambda)$.  \ps\ps 

Les consid\'erations de ce chapitre permettent de
d\'emontrer l'existence d'un certain nombre d'\'el\'ements de $\Pi_{\rm
disc}({\rm O}_8)$ ou $\Pi_{\rm disc}({\rm SO}_8)$.  On peut se demander si
ces \'el\'ements suffisent \`a expliquer tous les $\pi(\lambda)$ ci-dessus. 
La r\'eponse \`a cette question est donn\'ee par la table~\ref{tableSO8},
qui indique la liste de tous les $\psi(\pi(\lambda),{\rm St})$ quand
$m_1(\lambda) \leq 8$.  Pour des raisons de num\'erologie, il est plus
significatif de faire appara\^itre dans cette table l'\'el\'ement
$\lambda+\rho$, caract\`ere infinit\'esimal de $U_\lambda$
(\S VI.\ref{exparlan}), que nous coderons par le quadruple $z(\lambda)=(2
m_1+6, 2 m_2+4, 2 m_3+2, 2 |m_4|)$, si $\lambda =
\sum_{i=1}^4 m_i \varepsilon_i$. \ps\ps 

\begin{table}[htp]
\caption{Param\`etres standards des $\pi(\lambda)$ quand $m_1(\lambda) \leq 8$.}
{\scriptsize \renewcommand{\arraystretch}{1.8} \medskip
\begin{center}
\begin{tabular}{|c|c|c|c|}
\hline  $z(\lambda)$ & $\psi(\pi(\lambda),{\rm St})$ & $z(\lambda)$ & $\psi(\pi(\lambda),{\rm St})$ \\

\hline $(6,4,2,0)$ & $[7] \oplus [1]$ & $(22,16,14,0)$ & {\color{blue}${\rm Sym}^2 \Delta_{11} \oplus \Delta_{15}[2] \oplus [1]$} \\

\hline $(14,12,10,8)$ & $\Delta_{11}[4]$ &  $(22,18,16,0)$ & {\color{blue}${\rm Sym}^2 \Delta_{11} \oplus \Delta_{17}[2] \oplus [1]$} \\

\hline $(18,16,2,0)$ & $\Delta_{17}[2] \oplus [3] \oplus [1]$ &
$(22,20,2,0)$ & $\Delta_{21}[2] \oplus [3] \oplus [1]$ \\

\hline  $(18,16,14,12)$ & $\Delta_{15}[4]$ & $(22,20,6,4)$ & $\psi_{4,10}[2]$ \\

\hline $(20,18,8,6)$ & $\psi_{6,8}[2]$ & $(22,20,10,8)$ & $\psi_{8,8}[2]$ \\ 

\hline $(20,18,16,14)$ & $\Delta_{17}[4]$ & $(22,20,14,12)$ & $\psi_{12,6}[2]$  \\

\hline $(22,4,2,0)$ & ${\rm Sym}^2 \Delta_{11} \oplus [5]$ & $(22,20,18,0)$ & ${\rm Sym}^2 \Delta_{11} \oplus \Delta_{19}[2] \oplus [1]$ \\

\hline $(22,12,10,0)$ & {\color{blue}${\rm Sym}^2 \Delta_{11} \oplus \Delta_{11}[2] \oplus [1]$} & $(22,20,18,16)$ & $\Delta_{19}[4]$ \\
\hline
\end{tabular}
\end{center}}
\label{tableSO8}
\end{table}

Discutons un peu cette table. La repr\'esentation de param\`etre $[7] \oplus
[1]$ est bien entendu la repr\'esentation triviale.  Rappelons que la
notation $\Delta_w \in \Pi_{\rm cusp}(\PGL_2)$ a \'et\'e introduite au
\S\ref{ikedabocherer}.  Les quatres \'el\'ements $\psi_{j,k} \in \mathcal{X}(\SL_4)$ seront expliqu\'es au \S IX.\ref{remtripsijk}.  \ps \medskip
\begin{itemize}

\item[(i)] L'existence des repr\'esentations de param\`etre $\Delta_w[4]$
d\'ecoule du th\'eor\`eme \ref{arthurtrialite}.  Dans le cas du quadruple $(14,12,10,8)$,
c'est la repr\'esentation utilis\'ee dans la d\'emonstration du
th\'eor\`eme~\ref{casdim16}.  \ps \ps

\item[(ii)]  L'existence d'une repr\'esentation de $\Pi_{\rm disc}({\rm
O}_8)$ de param\`etre $\Delta_w[2] \oplus [3] \oplus [1]$ pour $w \equiv 1
\bmod 4$ se d\'eduit du th\'eor\`eme~\ref{corikedabocherer} : c'est le cas $n=8$, $k=w+1$ et $g=2$, qui s'applique car ${\rm L}(\frac{1}{2},\Delta_w) \neq 0$ pour $w=17, 21$ (Exemple~\ref{nonvanish}). Dans ces cas particuliers, un r\^ole important est jou\'e par  les formes propres dans ${\rm
S}_{\frac{w+3}{2}}({\rm Sp}_4(\Z))$ de param\`etre standard $\Delta_w[2]
\oplus [1]$ (les deux ``premi\`eres'' formes de Saito-Kurokawa, cas $g=2$ du
th\'eor\`eme~\ref{thmikeda2}). Lorsque $w=17, 21$, la surjectivit\'e de $\vartheta_{\frac{w-5}{2},2} : {\rm
M}_{{\rm H}_{\frac{w-5}{2},2}(\R^8)}({\rm O}_8) \rightarrow {\rm
S}_{\frac{w+3}{2}}({\rm Sp}_{4}(\Z))$ se v\'erifie par un simple calcul de
coefficient de s\'eries th\^eta \'etant donn\'e que
l'on sait depuis Igusa~\cite{igusagenus2} que ${\rm S}_{\frac{w+3}{2}}({\rm Sp}_{4}(\Z))$ est
de dimension $1$. Ce calcul sera justifi\'e dans la proposition IX.\ref{calcimagetheta}.  \ps\ps

\item[(iii)] Le cas $z(\lambda)= (22,20,18,0)$ admet \'egalement une
histoire int\'eressante, puisqu'il a \'et\'e \'etudi\'e par Miyawaki dans
\cite{miyawaki}.  Il a montr\'e que $\pi(\lambda)$ admet un
$\theta$-correspondant dans ${\rm S}_{12}({\rm Sp}_6(\Z))$ qu'il a
conjectur\'e de param\`etre standard ${\rm Sym}^2 \Delta_{11} \oplus
\Delta_{19}[2]$, ce qu'a par la suite d\'emontr\'e Ikeda~\cite{ikeda2}.  Le
cas des autres param\`etres, d'apparence pourtant similaire, de la forme ${\rm Sym}^2
\Delta_{11} \oplus \Delta_w[2] \oplus [1]$ avec $w \in \{11,15,17\}$ est
plus subtile.  Quand $w=11,15$ nous pourrions les justifier comme pour
$w=19$ si nous disposions d'un analogue de la construction d'Ikeda pour des
formes non scalaires, car pour les deux valeurs de $\lambda$ pertinentes on
pourrait sans doute v\'erifier que $\pi(\lambda)$ admet un $\vartheta$-correspondant pour
${\rm Sp}_8$ (voir le \S \ref{comptheta}). Ces param\`etres (indiqu\'es en bleu dans la table) sont pr\'edits par la th\'eorie d'Arthur, comme nous le verrons au
chapitre VIII. \ps\ps 
\end{itemize}

Nous renvoyons au chapitre VIII.\ref{formexpliciteson0} pour une confirmation directe mais conditionnelle de
l'int\'egralit\'e de cette table utilisant la th\'eorie d'Arthur, 
ainsi qu'\`a~\cite[Ch.  2
\& 7]{chrenard2} pour des tables nettement plus fournies. Bien qu'il serait sans doute possible de le faire par des construction de s\'eries th\^eta et les m\'ethodes du \S \ref{fsiegelpoids12}, nous ne justifierons pas inconditionnellement les trois param\`etres en bleu de la table \ref{tableSO8}. \ps\ps

\section{${\rm M}_{\DET}({\rm
O}_{24})$}\label{demoleechmoins}

\noindent Dans ce paragraphe nous d\'emontrons la proposition V.\ref{leechmoins}. 

\begin{prop}\label{leechmoins2} Le param\`etre standard de l'unique 
repr\'esentation $\pi$ dans $\Pi_{\rm disc}({\rm O}_{24})$ telle que
$\pi_\infty \simeq \DET$ est $\Delta_{11}[12]$. 
\end{prop}

\begin{pf} D'apr\`es le corollaire~\ref{cor24ikedabocherer}, il existe $\pi \in \Pi_{\rm disc}({\rm O}_{24})$ telle que $\pi_\infty=\C$ et $\psi(\pi,{\rm
St})=\Delta_{11}[12]$. On a $$\tau(2)(2^{12}-1)=2^{11} \,{\rm trace}({\rm c}_2(\pi), V_{\rm
St})$$ qui n'est autre que la valeur propre $\lambda_{24}$ dans les notations du~\S \ref{thetaniemeier}. Soit $f_{24} \in {\rm M}_\C({\rm O}_{24})$ la fonction image de
l'\'el\'ement $v_{24}$
par l'isomorphisme de ${\rm H}({\rm O}_{24})$-modules ${\rm M}_\C({\rm O}_{24}) \simeq \C[{\rm X}_{24}]$
d\'efini au corollaire IV.\ref{corhermitien}. \ps\ps

On rappelle que ${\rm M}_\C({\rm
SO}_{24})$ est muni d'une action de ${\rm H}({\rm SO}_{24})$, ainsi que
d'une involution $s$ (le ``changement d'orientation'', voir la fin du \S
IV.\ref{fauton}) dont la d\'ecomposition en espaces propres s'\'ecrit $${\rm
M}_{\C}({\rm SO}_{24}) = {\rm M}_{\C}({\rm O}_{24}) \oplus {\rm
M}_{\DET}({\rm O}_{24}).$$ Cette d\'ecomposition n'est
pas respect\'ee \`a priori par ${\rm H}({\rm SO}_{24})$, mais seulement par
${\rm H}({\rm O}_{24})$, la relation dont on dispose \'etant : $T \circ s =
s \circ {\rm H}(s)(T)$ pour tout $T \in {\rm H}({\rm SO}_{24})$. \ps\ps

Soit $V
\subset {\rm M}_{\C}({\rm SO}_{24})$ le ${\rm H}({\rm SO}_{24})[s]$-module
engendr\'e par $f_{24}$ et soit $f'_{24} \in V$ une forme propre pour ${\rm
H}({\rm SO}_{24})$. Elle engendre naturellement une repr\'esentation $\pi'
\in \Pi_{\rm disc}({\rm SO}_{24})$ de param\`etre standard $\Delta_{11}[12]$. En particulier, la
classe de conjugaison ${\rm c}_2(\pi') \subset {\rm SO}_{24}(\C)$ n'admet
pas la valeur propre $\pm 1$, car les valeurs propres de ${\rm
c}_2(\Delta_{11})$
sont non r\'eelles (ce sont les racines de
$x^2+\frac{24}{2^{11/2}}x+1$). Cette classe de conjugaison sous ${\rm
SO}_{24}(\C)$ n'est donc pas stable par l'action de conjugaison de ${\rm
O}_{24}(\C)$. La compatibilit\'e de l'isomorphisme de Satake aux
isomorphismes montre donc que si $f_{24}''=s(f_{24}') \in V$, et si $\pi'' \in
\Pi_{\rm disc}({\rm SO}_{24})$ est engendr\'ee par $f_{24}''$, alors on a $${\rm
c}_2(\pi')\neq {\rm c}_2(\pi'').$$ En particulier, $f_{24}'$ et $f_{24}''$ ne
sont pas proportionnelles et ont m\^eme valeurs propres que $f_{24}$ sous
l'action de ${\rm H}({\rm O}_{24})$. On conclut car l'\'el\'ement non nul
$f_{24}'-f_{24}''=(1-s)f_{24}'$ engendre la droite ${\rm M}_{\DET}({\rm O}_{24})$.
\end{pf}

\newcommand{\pcb}{\Pi_{\rm cusp}^\bot}

\chapter{La classification d'Arthur pour les $\Z$-groupes classiques}\label{classarthur}\label{introchaparthur}

\section{Param\`etres standards pour les groupes
classiques}\label{enonceparamst} Le but de ce
chapitre est d'expliquer la description de $\Pi_{\rm disc}(G)$ issue des travaux
d'Arthur~\cite{arthur} lorsque $G$ est un $\Z$-groupe semi-simple {\it classique}. {\it Nous
entendrons dor\'enavant par l\`a que $G$ est de la forme ${\rm Sp}_{2g}$ pour $g \geq 1$, ou ${\rm SO}_L$, o\`u $L$
est soit un ${\rm q}$-module sur $\Z$ de dimension $\neq 2$, soit un ${\rm q}$-${\rm i}$-module sur $\Z$ au sens de l'appendice {\rm B}  (\S V.\ref{exconjal})}.  Si $r$ est un entier $\geq 1$,
on pose $${\rm SO}_{r,r}={\rm SO}_{{\rm H}(\Z^r)} \, \, \, {\rm et}\, \, \,
{\rm SO}_{r+1,r}={\rm SO}_{{\rm H}(\Z^r)\oplus {\rm A}_1}.$$ Les
$\Z$-groupes classiques semi-simples qui sont des groupes de Chevalley\footnote{Nous utiliserons dans ce chapitre le terme de {\it groupe de Chevalley} comme synonyme de {\it $\Z$-groupe semi-simple d\'eploy\'e}.} sont
donc ${\rm Sp}_{2g}$ pour $g\geq 1$, ${\rm SO}_{r,r}$ pour $r\geq 2$ et
${\rm SO}_{r+1,r}$ pour $r\geq 1$.  Il sera commode de consid\'erer le
$\Z$-groupe trivial comme un groupe classique de Chevalley, que l'on notera
\'egalement ${\rm SO}_{1,0}$.  De plus, un r\^ole important sera jou\'e par
les $\Z$-groupes ${\rm SO}_n$, d\'efinis pour tout entier $n \geq 1$ tel que
$n \equiv -1,0,1 \bmod 8$ par $${\rm SO}_n={\rm SO}_{{\rm L}_n}$$ o\`u :
${\rm L}_n={\rm E}_n$ si $n\equiv 0 \bmod 8$, ${\rm L}_n={\rm E}_{n-1}
\oplus {\rm A}_1$ si $n \equiv 1 \bmod 8$, et o\`u ${\rm L}_n$ d\'esigne
l'orthogonal d'une racine\footnote{Le choix de cette racine est fix\'e une
fois pour toutes, et ne jouera aucun r\^ole dans la suite.  Comme toutes ces
racines sont permut\'ees transitivement par ${\rm W}({\bf D}_{n+1}) \subset
{\rm O}({\rm E}_{n+1})$, la classe d'isomorphisme du $\Z$-groupe ${\rm
SO}_n$ ne d\'epend par exemple que de $n$.  } de ${\rm E}_{n+1}$ si $n
\equiv -1 \bmod 8$ (\S IV.\ref{ensclassbil} \& \S IV.\ref{soimph}).  \ps\ps

Rappelons que si $G$ est un $\Z$-groupe
classique semi-simple, le $\C$-groupe $\widehat{G}$ est un groupe classique
complexe, et poss\`ede une repr\'esentation irr\'eductible distingu\'ee qui
est sa repr\'esentation standard (\S V.\ref{parlanpig}, \S V.\ref{exconjal})
$${\rm St} : \widehat{G} \rightarrow {\rm SL}_n(\C).$$

\begin{thmv}\label{arthurst} {\rm \cite[Thm. 1.5.2]{arthur}} Si $G$ est un $\Z$-groupe classique de Chevalley
et si $\pi
\in \Pi_{\rm disc}(G)$ alors $\psi(\pi,{\rm St}) \in \mathcal{X}_{\rm AL}(\SL_n)$. 

\end{thmv}

Ce r\'esultat, qui est un cas particulier de la conjecture g\'en\'erale
V.\ref{conjarthlan}, est aussi un cas tr\`es particulier du th\'eor\`eme
sus-cit\'e  (cas des repr\'esentations ``non ramifi\'ees en tous les nombres
premiers''). Il repose sur un \'edifice assez formidable de r\'esultats
difficiles, notamment de multiples variantes de la formule des traces
d'Arthur-Selberg (Arthur), la d\'ecomposition spectrale des espaces des formes automorphes
(Langlands), la th\'eorie de l'endoscopie (Langlands, Shelstad, Kottwitz),
et la d\'emonstration du fameux {\it lemme fondamental} (Waldspurger \cite{waltran}, Ng\^o \cite{ngo}, Laumon,
Chaudouard \cite{chaulau1, chaulau2}). Ainsi que l'explique Arthur dans son livre, les r\'esultats
de \cite{arthur} d\'ependent d'une variante ``avec torsion'' de ses
travaux sur la ``stabilisation'' de la formule des traces. La formule requise a r\'ecemment \'et\'e \'etablie par Moeglin et Waldspurger, dans une longue s\'erie d'articles \cite{stablewal, moewal2} (voir aussi~\cite{wallab}). Arthur mentionne \'egalement une autre hypoth\`ese concernant une extension des travaux de Shelstad concernant l'endoscopie tordue pour les groupes de Lie r\'eels, et qui a fait l'objet depuis de travaux de Shelstad \cite{shelstad1, shelstad2, shelstad3} et Mezo \cite{mezo1, mezo2}. Les \'enonc\'es de ce m\'emoire d\'ependant de ces travaux tous r\'ecents, par l'interm\'ediaire des \'enonc\'es de \cite{arthur}, seront affubl\'es d'une \'etoile verte ${\color{green}\star}$ en exposant.
\ps\ps 


Le travail d'Arthur~\cite{arthur} porte comme nous l'avons dit sur les
$\Z$-groupes classiques de Chevalley.  Cependant, ses travaux ant\'erieurs
permet\-tent d'en d\'eduire \'egalement une classification de $\Pi_{\rm
disc}(G)$ lorsque $G$ est un $\Z$-groupe classique quelconque.  Cette
classification est annonc\'ee dans le chapitre 9 de \cite{arthur}, mais la
r\'edaction compl\`ete n'est pas encore disponible, c'est pourquoi nous l'\'enoncerons sous forme de conjecture.\footnote{Mentionnons sur
ce point qu'aucune des difficult\'es \'evoqu\'ees par Arthur dans son
chapitre 9 ne semble s'appliquer \`a la situation qui nous int\'eresse, qui
ne concerne que des ``formes int\'erieures pures'' des groupes de Chevalley \cite{kaletha}.}
Comme nous le verrons, ces r\'esultats ne seront pas n\'ecessaires \`a
l'\'etablissement de nos r\'esultats principaux.  En revanche, ils
\'eclairent substantiellement les questions qui nous pr\'eoccupent, de sorte
qu'il serait dommage de les laisser sous silence. Par exemple, ils conduisent 
comme nous le verrons \`a une description directe et tr\`es pr\'ecise de
$\Pi_{\rm disc}({\rm SO}_n)$, le th\`eme principal de ce m\'emoire.

\begin{conj}\label{arthurstint} {\rm \cite[Ch. 9]{arthur}} L'\'enonc\'e du
th\'eor\`eme~\ref{arthurst} vaut encore si $G$ est un
$\Z$-groupe classique quelconque.
\end{conj}

\medskip

\pn
{\sc Notations} : Tout comme au \S V.\ref{exconjal}, nous notons ${\rm Class}_\C$ l'ensemble
constitu\'e des $\C$-groupes ${\rm Sp}_{2g}$ ($g\geq 1$) et ${\rm SO}_m$
($m\geq 0$, $m \neq 2$). Un groupe $H \in {\rm Class}_\C$ est uniquement d\'etermin\'e
par le couple associ\'e $({\rm
n}_H,{\rm w}_H) \in \N \times \{0, 1\}$ d\'efini de la mani\`ere suivante : \ps\ps

-- ${\rm n}_H$ est la dimension de la repr\'esentation standard de $H$, \ps\ps

-- ${\rm w}_H=0$ si, et seulement si, $H \simeq {\rm SO}_m$ pour un entier $m \geq 0$.\ps\ps

\section{Repr\'esentations autoduales de $\PGL_n$}\label{repselfdual}

La description d'Arthur, et plus g\'en\'eralement la conjecture
d'Arthur-Langlands, fait intervenir des repr\'esentations dans $\Pi_{\rm
cusp}(\PGL_m)$ pour tout $m\geq 1$, que nous n'avons rencontr\'ees jusqu'\`a
pr\'esent essentiellement que dans le cas $m=2$ et plut\^ot du point de vue
de l'identit\'e ${\rm GSp}_2 = {\rm GL}_2$.  Nous revenons maintenant sur
celles-l\`a.

\subsection{Dualit\'e dans $\Pi_{\rm disc}(\PGL_n)$}\label{dualitepgln} Soient $n\geq 1$ un entier et $\mathcal{R}_n$ l'ensemble des sous-groupes
discrets de rang $n$ de $\R^n$. Il est muni d'une action naturelle
transitive de ${\rm GL}_n(\R)$ et l'orbite du r\'eseau $\Z^n$ s'identifie \`a ${\rm GL}_n(\R)/{\rm
GL}_n(\Z)\isomo \mathcal{R}_n$. Le sous-groupe ${\rm GL}_n(\Q)$ pr\'eserve 
le sous-espace $\mathcal{R}_\Z(\Q^n) \subset \mathcal{R}_n$ des r\'eseaux de
$\Q^n$; cette action est transitive et s'\'etend naturellement en une action
de ${\rm GL}_n(\AAA_f)$ (\S IV.\ref{ensclassgln}). L'application ${\rm GL}_n(\R) \times {\rm
GL}_n(\AAA_f) \rightarrow \mathcal{R}_n$, $(g_\infty,g_f) \mapsto
g_\infty^{-1}(g_f(\Z^n))$ a donc un sens; elle induit une bijection ${\rm GL}_n(\Q)
\backslash {\rm GL}_n(\AAA)/{\rm GL}_n(\widehat{\Z}) \isomo \mathcal{R}_n$,
puis un isomorphisme naturel 
$$\mathcal{A}^2({\rm PGL}_n) \isomo {\rm L}^2(\underline{\mathcal{R}}_n),$$
o\`u $\underline{\mathcal{R}}_n$ d\'esigne le quotient de $\mathcal{R}_n$ par le
groupe $\R^\times$ des homoth\'eties, muni d'une mesure ${\rm GL}_n(\R)$ invariante
non nulle (\S IV.\ref{fautcarreint}). Les actions naturelles de $\PGL_n(\R)$ et de ${\rm H}(\PGL_n)$ 
sur ${\rm L}^2(\underline{\mathcal{R}}_n)$ d\'eduites par transport
de structure sont alors les actions
\'evidentes. En particulier, si $f \in {\rm L}^2(\underline{\mathcal{R}}_n)$ est
continue, si ${\rm T}_A \in {\rm H}({\rm PGL}_n)$ est l'op\'erateur d\'efini au \S
IV.\ref{annheckeclass}, et si $L \in
\mathcal{R}_n$, alors ${\rm T}_A(f)(\underline{L})=\sum_M f(\underline{M})$ o\`u $M$ parcourt les sous-groupes de $L$ tels que
$L/M \simeq A$, la notation $\underline{N}$ d\'esignant la classe d'homoth\'etie
de $N \in \mathcal{R}_n$. \ps\ps 

Le $\Z$-groupe $\GL_n$ admet pour automorphisme $g \mapsto {}^t\!g^{-1}$,
qui agit donc \'egalement sur $\mathcal{A}^2({\rm PGL}_n)$ par une
involution not\'ee $\theta$, qui pr\'eserve les sous-espaces
$\mathcal{A}_{\rm disc}({\rm PGL}_n)$ et $\mathcal{A}_{\rm cusp}({\rm
PGL}_n)$. Concr\`etement, munissons le $\Z$-module $\Z^n$ de la forme
bilin\'eaire sym\'etrique non d\'eg\'en\'er\'ee standard, et notons
$L^\sharp$ le r\'eseau dual de $L \in \mathcal{R}_n$ pour cette forme.
L'involution en question de ${\rm L}^2(\underline{\mathcal{R}}_n)$, encore not\'ee $\theta$, est 
simplement d\'efinie par
$\theta(f)(\underline{L})=f(\underline{L^\sharp})$.  Elle
satisfait donc $\theta(T(f))= \iota(T)(\theta(f))$ pour tout
$T \in {\rm H}({\rm PGL}_n)$ et tout $f \in {\rm
L}^2(\underline{\mathcal{R}}_n)$, o\`u $\iota$ est l'automorphisme involutif de ${\rm
H}({\rm PGL}_n)$ d\'efini au
\S IV.\ref{annheckeclass}. \ps\ps

Si $\pi= \pi_\infty \otimes \pi_f \in \Pi(\PGL_n)$, on note $\pi^\vee \in \Pi(\PGL_n)$
l'\'el\'ement d\'efini de la mani\`ere suivante. D'une part, $(\pi^\vee)_\infty$ est la
repr\'esentation ayant m\^eme espace que $\pi_\infty$ mais compos\'ee par
l'automorphisme $g \mapsto {}^{\rm t}\!g^{-1}$ de $\PGL_n(\R)$. D'autre
part, si l'on voit $(\pi^\vee)_f$ et $\pi_f$ comme des homomorphismes ${\rm
H}({\rm PGL}_n) \rightarrow \C$, alors $(\pi^\vee)_f = \pi_f \circ \iota$. \ps\ps 

Si $c$ est la classe de conjugaison d'un \'el\'ement semisimple $g$ de ${\rm SL}_n(\C)$, nous noterons $c^{-1}$ (resp. $\overline{c}$) la classe de conjugaison de $g^{-1}$ (resp. du conjugu\'e complexe de $g$). De plus, si $c$ est la classe de conjugaison d'un \'el\'ement semisimple $X$ de ${\got{sl}}_n(\C)$, nous noterons $-c$ la classe de conjugaison de $-X$.
La proposition suivante est bien connue.

\begin{prop} \label{propdesparamdual} Soit $\pi \in \Pi(\PGL_n)$. \begin{itemize}\ps\ps
\item[(i)] Si $\pi \in \Pi_{\rm cusp}(\PGL_n)$ {\rm (}resp. $\pi \in \Pi_{\rm disc}(\PGL_n)${\rm )}, il en va de m\^eme de $\pi^\vee$.\ps\ps
\item[(ii)] Pour tout premier $p$, on a l'\'egalit\'e ${\rm c}_p(\pi^\vee) = {\rm c}_p(\pi)^{-1}$. Si l'on a de plus $\pi \in \Pi_{\rm disc}(\PGL_n)$, on a \'egalement ${\rm c}_p(\pi^\vee) = \overline{{\rm c}_p(\pi)}$.\ps\ps
\item[(iii)] La repr\'esentation $(\pi^\vee)_\infty$ est le dual de la repr\'esentation unitaire $\pi_\infty$, et l'on a l'\'egalit\'e ${\rm c}_\infty(\pi^\vee) = -{\rm c}_\infty(\pi)$.\ps\ps
\end{itemize}
\end{prop}

\begin{pf} Les paragraphes qui pr\'ec\'edent l'\'enonc\'e justifient l'assertion (i). Soit $p$ un nombre premier. Pour tout \'el\'ement $T \in {\rm H}(\PGL_n)$, on a la relation $\iota(T)=T^{\rm t}$ d'apr\`es le \S IV.\ref{annheckeclass}. La discussion suivant le
scholie V.\ref{corsatake} montre donc l'\'egalit\'e ${\rm c}_p(\pi^\vee) = {\rm c}_p(\pi)^{-1}$. D'autre part, l'action de ${\rm H}(\PGL_n)$ sur ${\mathcal A}_{\rm cusp}({\rm PGL}_n)$ \'etant une $\star$-action pour le produit de Petersson sur ce dernier d'apr\`es le \S {\rm IV}.\ref{fautcarreint}, on en d\'eduit \'egalement que les \'el\'ements $\pi_p$ et $\pi_p^\vee$, vus comme morphismes d'anneaux ${\rm H}_p({\rm PGL}_n) \rightarrow \C$, sont conjugu\'es complexes l'un de l'autre. Du point de vue de l'isomorphisme de Satake cela s'\'ecrit ${\rm c}_p(\pi^\vee)=\overline{{\rm c}_p(\pi)}$.  Le (iii) est cons\'equence de la th\'eorie des caract\`eres d'Harish-Chandra et du fait que tout \'el\'ement de ${\rm GL}_n(\R)$ est conjugu\'e \`a sa transpos\'ee. Il n'est pas difficile de v\'erifier sur la d\'efinition de l'isomorphisme d'Harish-Chandra que l'on a ${\rm c}_\infty(\pi^\vee) = -{\rm c}_\infty(\pi)$.
\end{pf}

\begin{definition} Soit $\pi \in \Pi(\PGL_n)$. La repr\'esentation $\pi^\vee$ est appel\'ee repr\'esentation duale de $\pi$. On dit que $\pi$ est autoduale si $\pi^\vee \simeq \pi$.  On note $\Pi_{\rm cusp}^\bot(\PGL_n) \subset \Pi_{\rm
cusp}(\PGL_n)$ le sous-ensemble des repr\'esentations autoduales. 
\end{definition}

Notons que d'apr\`es le th\'eor\`eme de multiplicit\'e $1$ de Jacquet-Shalika, pour que $\pi \in
\Pi_{\rm cusp}(\PGL_n)$ soit autoduale il suffit que ${\rm
c}_p(\pi)={\rm c}_p(\pi)^{-1}$ pour tout premier $p$ (voire m\^eme, pour tout premier $p$ sauf un nombre fini). \ps \ps

La repr\'esentation triviale de ${\rm PGL}_1$ est bien s\^ur autoduale. De
plus, l'automorphisme $g \mapsto {}^t\!g^{-1}$ \'etant int\'erieur quand
$n=2$, l'inclusion $\Pi_{\rm cusp}^\bot(\PGL_2) \subset
\Pi_{\rm cusp}(\PGL_2)$ est une \'egalit\'e.  Cela ne vaut pas d\`es que
$n>2$. L'int\'er\^et principal des repr\'esentations autoduales pour nos
pr\'eoccupations vient du th\'eor\`eme suivant, pr\'ecisant l'\'enonc\'e du th\'eor\`eme \ref{arthurst}.

\begin{thmv}\label{arthurstbis} {\rm \cite[Thm. 1.5.2]{arthur}} Soient $G$ un $\Z$-groupe classique de Chevalley et $\pi
\in \Pi_{\rm disc}(G)$. Alors $\psi(\pi,{\rm St})$ est de la forme $\oplus_{i=1}^k \pi_i[d_i]$ avec $\pi_i \in \Pi_{\rm cusp}^\bot(\PGL_{n_i})$ pour tout $i=1,\dots,k$ et $\sum_{i=1}^k n_i d_i = {\rm n}_{\widehat{G}}$. De plus, cette \'ecriture est unique et les couples $(\pi_i,d_i)$, $i=1,\dots,k$, sont deux \`a deux distincts.
\end{thmv}

Un raffinement similaire de la conjecture \ref{arthurstint} est \'egalement attendu. Remarquons que tous les param\`etres de Langlands mis en \'evidence aux chapitres V et VII ont bien la
propri\'et\'e ci-dessus. En particulier, on constate que la repr\'esentation de Gelbart-Jacquet ${\rm Sym}^2 \pi \in \Pi_{\rm
cusp}(\PGL_3)$, o\`u $\pi \in \Pi_{\rm cusp}(\PGL_2)$, est autoduale. Une justification partielle de ces propri\'et\'es d'autodualit\'e est donn\'ee par la proposition \'el\'ementaire suivante.

\begin{prop}\label{propselfdquasi} Soient $G$ un $\Z$-groupe classique et $\pi \in
\Pi(G)$. Supposons que $\psi(\pi,{\rm St})=\oplus_{i=1}^k \pi_i[d_i]$ o\`u $\pi_i \in \Pi_{\rm cusp}(\PGL_{n_i})$ pour tout $i=1,\dots,k$. Alors pour tout $i$ il existe $j$ tel que $\pi_j =\pi_i^\vee$ et $d_j = d_i$. \end{prop}

\begin{pf} En effet, la repr\'esentation ${\rm St}$ de $\widehat{G}$ \'etant
autoduale, on a ${\rm St}({\rm c}_p(\pi))={\rm St}({\rm c}_p(\pi))^{-1}$ pour tout premier $p$ (resp. ${\rm St}({\rm c}_\infty(\pi))=-{\rm St}({\rm c}_\infty(\pi))$) et donc 
\'egalement $\psi(\pi,{\rm St})=\oplus_{i=1}^k \pi_i^\vee[d_i]$. On
conclut d'apr\`es la proposition V.\ref{jacquetshalika}
(Th\'eor\`eme de Jacquet-Shalika).
\end{pf}

\subsection{Repr\'esentations alg\'ebriques r\'eguli\`eres}\label{paralgreg}

Soit $\pi \in \Pi_{\rm cusp}(\PGL_n)$. Les {\it poids} de $\pi$ sont les valeurs propres de la classe de conjugaison semi-simple 
${\rm c}_\infty(\pi) \subset {\rm M}_n(\C)$. Leur ensemble sera not\'e ${\rm Poids}(\pi) \subset \C$.

\begin{definition}\label{defalg} Soit $\pi \in \Pi_{\rm cusp}(\PGL_n)$. On
dit que $\pi$ est {\it alg\'ebrique}\footnote{Le lecteur prendra garde
qu'il existe plusieurs notions de repr\'esentations automorphes
alg\'ebriques dans la litt\'erature.  La d\'efinition ci-dessus,
qui est essentiellement celle consid\'er\'ee par exemple dans~\cite[\S 18.2]{borelcorvallis}, mais qui
n'est pas celle utilis\'ee par Clozel dans~\cite{clozel}, est r\'eminescente
de la notion de caract\`ere de Hecke de type ${\rm A}_0$ au sens de Weil
(voir~\cite{buzzardgee} pour des clarifications entre les diverses notions).} si
${\rm{Poids}}(\pi) \subset \frac{1}{2}\Z$ et si $\forall w,w' \in {\rm{Poids}}(\pi), \, \, w-w' \in \Z$.  \ps Si $\pi \in \Pi_{\rm cusp}(\PGL_n)$ est
alg\'ebrique, son {\rm poids motivique} est le plus grand $w \in \Z$ tel que
$-\frac{w}{2} \in {\rm{Poids}}(\pi)$, on le note ${\rm w}(\pi)$.  En
particulier, ${\rm Poids}(\pi) \subset \frac{{\rm w}(\pi)}{2}+\Z$. 
\end{definition}

Bien que les repr\'esentations alg\'ebriques forment une partie infime de
$\Pi_{\rm cusp}(\PGL_n)$, ce seront les seules \`a jouer un r\^ole dans ce
travail. Un indicateur \`a cela est donn\'e par la proposition 
suivante. \ps\ps 


\begin{prop}\label{toutestalg} Soient $G$ un $\Z$-groupe semi-simple, $\pi \in \Pi(G)$ et $r : \widehat{G} \rightarrow {\rm SL}_n$ une
$\C$-repr\'esentation.  On suppose que : \begin{itemize} \ps\ps  
\item[(i)] $\pi_\infty$ a m\^eme caract\`ere infinit\'esimal qu'une $\C$-repr\'esentation irr\'eductible de dimension
finie de $G_\C$. \ps\ps 
\item[(ii)] $\psi(\pi,r)=\oplus_{i=1}^k \pi_i[d_i]$ avec $\pi_i \in \Pi_{\rm cusp}(\PGL_{n_i})$ pour $i=1,\dots,k$ (\S V.\ref{parconjarthlan}). \ps\ps 
\end{itemize}
Alors $\pi_i$ est alg\'ebrique pour tout $i$. De plus, la classe de ${\rm w}(\pi_i)+d_i-1$
dans $\Z/2\Z$ ne d\'epend que de $r$ (et non de l'entier $i$, ou m\^eme de $\pi$). 
 \end{prop}

\begin{pf}  Soit $\mu$ le plus haut poids de $r$ (un
co-poids de $G_\C$). Le caract\`ere infinit\'esimal de $\pi_\infty$ est de
la forme $\lambda+\rho$ o\`u $\lambda$ est un poids dominant de $G_\C$ et
$\rho$ la demi-somme des racines positives (\S V.\ref{carinf}). Les valeurs propres de $r({\rm c}_\infty(\pi))$ sont
par d\'efinition de la forme $\langle \lambda + \rho, \mu' \rangle$, o\`u
$\mu' $ est un poids de $r$.  Mais $2 \rho$ est un poids de $G_\C$ et
$\langle \mu - \mu', \rho \rangle \in \Z$ si $\mu'\leq \mu$, ces valeurs propres
sont donc toutes dans $\frac{1}{2}\Z$ et diff\`erent deux \`a deux d'un
\'el\'ement de $\Z$.  Cette propri\'et\'e est h\'erit\'ee par
les poids des $\pi_i$. \end{pf}

Cette proposition s'applique notamment pour tout $\pi \in \Pi(G)$ si $G(\R)$ est compact (\S V.\ref{carinf}).
Elle s'applique \'egalement si $G={\rm Sp}_{2g}$ et si 
$\pi$ est engendr\'ee par une forme propre dans ${\rm S}_W({\rm
Sp}_{2g}(\Z))$ avec $W$ positif, d'apr\`es le corollaire~\S V.\ref{infcarsiegel}
(ainsi que la discussion qui suit) et le \S V.\ref{exdorad}. \ps\ps 

Soient $G$ un $\Z$-groupe classique et ${\rm St} : \widehat{G} \rightarrow {\rm SL}_n$ la repr\'esentation standard de $\widehat{G}$. Nous aurons besoin de pr\'eciser l'analyse pr\'ec\'edente dans ce contexte.  Soit ${\rm Irr}(G_\C)$ l'ensemble des classes d'isomorphisme de $\C$-repr\'esentations irr\'eductibles de dimension finie de $G_\C$.
Pour chaque $V \in {\rm Irr}(G_\C)$ observons la classe de conjugaison semi-simple ${\rm St}({\rm Inf}_V) \subset {\rm M}_n(\C)$ (on rappelle que ${\rm Inf}_V \in \widehat{\mathfrak{g}}_{\rm ss}$ d\'esigne le caract\`ere infinit\'esimal de $V$). On constate par inspection des donn\'ees radicielles (\S VI.\ref{exdorad}) qu'il y a trois cas  bien distincts : \begin{itemize} \ps\ps  

\item[I.] Si $\widehat{G}={\rm SO}_n(\C)$ avec $n=2g+1$ impair, alors $V
\mapsto {\rm St}({\rm inf}_V)$ induit une bijection entre ${\rm
Irr}(\widehat{G})$ et l'ensemble des classes de conjugaison semi-simples $X
\subset {\rm M}_n(\C)$ telles que $-X=X$ et dont toutes les valeurs propres
sont distinctes et dans $\Z$.  \ps\ps 

\item[II.] Si $\widehat{G}={\rm Sp}_n(\C)$ (et donc $n$ est pair), alors $V
\mapsto {\rm St}({\rm Inf}_V)$ induit une bijection entre ${\rm
Irr}(\widehat{G})$ et l'ensemble des classes de conjugaison semi-simples $X
\subset {\rm M}_n(\C)$ telles que $-X=X$ et dont toutes les valeurs propres
sont distinctes et dans $\frac{1}{2}\Z-\Z$.  \ps\ps 

\item[III.] Si $\widehat{G}={\rm SO}_n(\C)$ avec $n$ pair, alors $V \mapsto
{\rm St}({\rm Inf}_V)$ induit une surjection entre ${\rm Irr}(\widehat{G})$
et l'ensemble des classes de conjugaison semi-simples $X \subset {\rm
M}_n(\C)$ telles que $-X=X$, dont toutes les valeurs propres sont dans $\Z$
et distinctes, sauf peut-\^etre la valeur propre $0$, dont la multiplicit\'e
est $\leq 2$.  De plus, ${\rm St}({\rm Inf}_V)={\rm St}({\rm Inf}_{V'})$ si
et seulement si $V$ et $V'$ sont conjugu\'ees l'une de l'autre sous l'action
ext\'erieure de ${\rm O}_n(\C)$ (ce qui entra\^ine $V=V'$ si, et seulement
si, $0$ est valeur propre de ${\rm Inf}_V={\rm Inf}_{V'}$).  \ps\ps 
\end{itemize}

\noindent Dans tous les cas, les valeurs propres de ${\rm St}({\rm Inf}_V)$ sont donc dans $\frac{{\rm w}_{\widehat{G}}}{2} + \Z$.

\begin{definition}\label{definitionreguliere} On dit que $\pi \in \Pi_{\rm cusp}(\PGL_n)$ est {\rm r\'eguli\`ere} si $|{\rm Poids}(\pi)|=n$.
\end{definition}

\begin{prop}\label{parstalg} Soient $G$ un $\Z$-groupe classique et ${\rm St} : \widehat{G} \rightarrow {\rm SL}_n$ la repr\'esentation standard de $\widehat{G}$. On suppose que $\psi=\oplus_{i=1}^r \pi_i[d_i] \in \mathcal{X}_{\rm AL}({\rm SL}_n)$, o\`u $\pi_i \in
\Pi_{\rm cusp}(\PGL_{n_i})$ pour tout $i=1,\dots,k$, et que $\psi_\infty={\rm St}({\rm Inf}_V)$ o\`u $V \in {\rm Irr}(G_\C)$.  \ps
Soit $i \in \{1,\dots,k\}$. Alors $\pi_i$ est alg\'ebrique, de poids motivique ${\rm w}(\pi_i) \equiv  d_i - 1 +  {\rm w}_{\widehat{G}} \bmod 2$, et autoduale. De plus, $\pi_i$ est r\'eguli\`ere, \`a moins que
l'on ne soit dans le cas exceptionnel suivant : \ps\ps  \begin{itemize}
\item[(a)] $\widehat{G}(\C) \simeq {\rm SO}_{2r}(\C)$ et $\psi_\infty$ admet $0$ pour valeur propre double,
\ps\ps  \item[(b)] $d_i=1$, $n_i \equiv 0 \bmod 2$ et
$|{\rm Poids}(\pi_i)|=n_i-1$,\ps\ps  \item[(c)] et pour tout $j \neq i$,
$\pi_j$ est r\'eguli\`ere et $n_j \equiv 0 \bmod 2$.  \ps\ps  \end{itemize} \end{prop}
	
\begin{pf} Hormis l'assertion concernant autodualit\'e de $\pi_i$, la proposition d\'ecoule imm\'ediatement de l'analyse des cas I, II et III ci-dessus et de la proposition~\ref{toutestalg}. \ps
V\'erifions l'autodualit\'e de $\pi_i$. D'apr\`es la proposition~\ref{propselfdquasi}, il existe $j$ tel que $\pi_j=\pi_i^\vee$ et $d_j=d_i$. D'apr\`es la proposition \ref{propevidentealg} ci-dessous, $\pi_j=\pi_i^\vee$ entra\^ine  ${\rm c}_\infty(\pi_j)={\rm c}_\infty(\pi_i)$.  Compte tenu de l'hypoth\`ese sur $\psi_\infty$, et de l'analyse des cas I, II et III, cela entra\^ine $j=i$ ou $n_i=n_j=1$. Dans ce dernier cas, on a n\'ecessairement $\pi_j=\pi_i=1$, et donc dans tous les cas nous avons bien $\pi_i^\vee=\pi_i$. \end{pf}	
	
Dans la d\'emonstration ci-dessus nous avons fait appel \`a la proposition suivante (\'evidente si $\pi$ est suppos\'ee autoduale).

\begin{prop} \label{propevidentealg} Si $\pi \in \Pi_{\rm cusp}(\PGL_n)$ est alg\'ebrique alors $\pi_\infty$ est isomorphe \`a son dual. En particulier, on a ${\rm c}_\infty(\pi) = {\rm c}_\infty(\pi^\vee)$ et  : \ps\ps  
\begin{itemize} 
\item[(i)]  $w \mapsto -w$ est une bijection de ${\rm Poids}(\pi)$ {\rm (}pr\'eservant les multiplicit\'es comme valeurs propres de ${\rm c}_\infty(\pi)${\rm )}. \ps\ps 
\item[(ii)] si $n \equiv 1 \bmod 2$, alors $0 \in {\rm Poids}(\pi)$ et ${\rm w}(\pi) \equiv 0 \bmod 2$.\ps\ps 
\end{itemize}
\end{prop}	
	
Pour l'expliquer, il sera n\'ecessaire d'\'etudier plus en d\'etail les composantes archim\'ediennes des \'el\'ements de $\Pi_{\rm cusp}(\PGL_n)$. Une autre motivation \`a cela est que les repr\'esentations alg\'ebriques r\'eguli\`eres autoduales, ainsi que celles intervenant dans le cas exceptionnel de l'\'enonc\'e de la proposition \ref{parstalg}, satisfont quelques propri\'et\'es suppl\'ementaires cach\'ees, qu'il nous faudra pr\'eciser. Cette analyse sera \'egalement n\'ecessaire pour appliquer les \'enonc\'es d'Arthur. En effet, 
ils feront intervenir des {\it facteurs
$\varepsilon$ de paires de repr\'esentations} alg\'ebriques, qui se lisent aussi
(comme les facteurs $\Gamma$) sur leurs composantes
archim\'ediennes. Ce travail a d\'ej\`a \'et\'e men\'e dans~\cite[\S 3.11]{chrenard2}, dont nous allons rappeler
quelques r\'esultats dans les paragraphes suivants. 


\subsection{Repr\'esentations de $\GL_n(\R)$}\label{apparitionWR} Soit ${\rm W}_\R$ le groupe
de Weil du corps $\R$ \cite[\S 1]{tate}.  C'est un groupe topologique, extension non-triviale
de ${\rm Gal}(\C/\R)=\Z/2\Z$ par $\C^\times$, pour l'action naturelle de conjugaison. Il est engendr\'e par
son sous-groupe ouvert $\C^\times$ ainsi qu'un \'el\'ement $j$, avec pour relations $j^2=-1
$ et $j z j^{-1} = \overline{z}$ pour tout $z \in \C^\times$. \ps\ps 

Suivant Langlands~\cite{langlandsreel}, les repr\'esentations continues et semi-simples ${\rm W}_\R \rightarrow
\GL_n(\C)$ joueront un r\^ole important. Rappelons la forme des
irr\'eductibles, qui sont de dimension $1$ ou $2$.  Soit 
 $$\eta : {\rm W}_\R \rightarrow
\R^\times$$ l'unique morphisme de groupes tel que $\eta(j)=-1$ et
$\eta(z)=z\overline{z}$ pour tout $z \in \C$; il induit un isomorphisme ${\rm
W}_\R^{\rm ab} \isomo \R^\times$.
Les morphismes continus ${\rm
W}_\R \rightarrow \C^\times$ sont donc les $|\eta|^s$ et
$\epsilon_{\C/\R}|\eta|^s$, o\`u $s \in \C$ et
$\epsilon_{\C/\R}=\eta/|\eta|$. Si $w$ est un entier $\geq 0$, 
consid\'erons la repr\'esentation induite
$${\rm I}_w ={\rm Ind}_{\C^\times}^{{\rm W}_\R}
\, \, \, \left(z \mapsto (\frac{z}{|z|})^w \right).$$ Elle est irr\'eductible si et seulement si $w \neq 0$, et de
plus ${\rm I}_0 \simeq 1 \oplus \epsilon_{\C/\R}$.  Les irr\'eductibles de
dimension $2$ de ${\rm W}_\R$ sont les ${\rm I}_w \otimes |\eta|^s$ o\`u
$w \neq 0$ et $s \in \C$.\ps\ps

Notons $\Phi(\GL_n(\R))$ l'ensemble des classes d'isomorphisme de
repr\'esentations continues semi-simples ${\rm W}_\R \rightarrow \GL_n(\C)$. 
La {\it param\'etrisation de Langlands} associe \`a chaque repr\'esentation
irr\'eductible unitaire $U$ de $\GL_n(\R)$ (et plus g\'en\'eralement \`a
tout $(\mathfrak{g},K)$-module irr\'eductible de ce dernier) un \'el\'ement ${\rm
L}(U) \in \Phi(\GL_n(\R))$ qui d\'etermine $U$ \`a isomorphisme pr\`es
\cite{langlandsreel} \cite{knappmotives}.  Bien que l'application $U \mapsto {\rm
L}(U)$ puisse \^etre enti\`erement explicit\'ee~\cite{knappmotives}, cela ne nous
sera pas r\'eellement utile : nous nous contenterons des assertions de
compatibilit\'es
suivantes, valables pour toute repr\'esentation unitaire irr\'eductible
$U$ de $\GL_n(\R)$.
\ps\ps  \begin{itemize}\ps\ps

\item[(i)] (dualit\'e) ${\rm L}(U^\ast) \simeq {\rm L}(U)^\ast$. \ps\ps 


\item[(ii)] (caract\`ere central) $\DET \, {\rm L}(U) = \chi_U \circ \eta$ o\`u
$\chi_U : \R^\times \rightarrow \C^\times$ d\'esigne le caract\`ere central de $U$. \ps\ps 

\item[(iii)] (caract\`ere infinit\'esimal) \'Ecrivons ${\rm
L}(U)_{|\mathbb{C}^\times} \simeq \oplus_{i=1}^n \chi_i$, o\`u les $\chi_i$ pour
$1\leq i \leq n$ sont des caract\`eres $\C^\times
\rightarrow \C^\times$. Pour tout $i$, il existe donc un unique couple $(\lambda_i,\mu_i) \in \C^2$, avec  $\lambda_i-\mu_i \in
\Z$, tel que\footnote{Suivant Langlands, il est suggestif de noter
$z^\lambda\overline{z}^\mu$ l'\'el\'ement $(\frac{z}{|z|})^{\lambda-\mu}|z|^{\lambda+\mu}$
quand $z \in \C^\times$, et $\lambda,\mu \in \C$ sont tels que $\lambda- \mu
\in \Z$.} $\chi_i(z)=(\frac{z}{|z|})^{\lambda_i-\mu_i}|z|^{\lambda_i+\mu_i}$.
Alors ${\rm Inf}_U$ est la classe de conjugaison semi-simple de ${\rm
M}_n(\C)$ ayant pour valeurs propres les $\lambda_i$, $i=1,\dots,n$.\ps\ps 
\end{itemize}

La param\'etrisation ci-dessus s'applique notamment aux repr\'esentations
irr\'eductibles unitaires de ${\rm PGL}_n(\R)$, vues comme 
repr\'esentations de $\GL_n(\R)$ de caract\`ere central trivial. \ps\ps 

L'assertion (i) de la proposition suivante est le {\it lemme de puret\'e} de
Clozel~\cite[Lemma 4.9]{clozel}, elle entra\^ine la proposition \ref{propevidentealg}.

\begin{prop}\label{contraintesalg} Soit $\pi \in \Pi_{\rm cusp}(\PGL_n)$ alg\'ebrique de poids
les $\frac{\omega_i}{2}$ (compt\'es avec multiplicit\'es), o\`u $\omega_1 \geq \dots \geq \omega_n$. Soient $E$
et $F$ les sous-ensembles de $\{1,\dots,n\}$ d\'efinis par $E=\{i, \omega_i >0\}$ et 
$F=\{i, \omega_i=0\}$. \ps\ps

\begin{itemize}
\item[(i)] Il existe un unique $(m_j) \in \{0,1\}^F$ tel que
$${\rm L}(\pi_\infty) \simeq \bigoplus_{i \in E} {\rm I}_{\omega_i} \oplus
\bigoplus_{j \in F} \epsilon_{\C/\R}^{m_j},$$
en particulier $|F| \equiv n \bmod 2$ et $\pi_\infty$ est isomorphe \`a son dual.
\ps\ps
\item[(ii)] Si ${\rm w}(\pi) \equiv 0 \bmod 2$ alors $\sum_{j \in F} m_j \equiv
|E| \bmod 2$. \ps\ps 
\item[(iii)] Si $\pi$ est r\'eguli\`ere et ${\rm w}(\pi) \equiv n \equiv 0
\bmod 2$, alors $n \equiv 0 \bmod 4$.\ps\ps 
\item[(iv)] Si $|{\rm Poids}(\pi)|=n-1$ alors $n \equiv 0 \bmod 2$ et
$F=\{\frac{n}{2},\frac{n}{2}+1\}$. De plus, $n \equiv 0 \bmod
4$ si et seulement si ${\rm L}(\pi_\infty) \simeq \oplus_{i=1}^{n/2} {\rm I}_{w_i}$. \ps\ps 

\end{itemize}

\end{prop}

\begin{pf} (voir \cite[\S 3.11]{chrenard2}) Rappelons l'argument du lemme de puret\'e de Clozel. D'apr\`es Jacquet et Shalika, si ${\rm I}_w \otimes |\eta|^s$ (resp. $|\eta|^s$ ou $\varepsilon_{\C/\R} |\eta|^s$) est une sous-repr\'esentation de ${\rm L}(\pi_\infty)$, alors $-\frac{1}{2} < {\rm Re} \,s < \frac{1}{2}$. D'autre part, $\pm \frac{w}{2}+s$ (resp. $s$) est un poids de $\pi$, par compatibilit\'e de la param\'etrisation de Langlands au
caract\`ere infinit\'esimal.. L'hypoth\`ese ${\rm Poids}(\pi) \subset \frac{1}{2}\Z$ entra\^ine donc $s=0$. Ainsi, ${\rm L}(\pi_\infty)$ est somme directe de repr\'esentations de la forme ${\rm I}_w$, $1$ et $\epsilon_{\C/\R}$. La premi\`ere assertion du (i) est alors cons\'equence de la compatibilit\'e de la param\'etrisation de Langlands au
caract\`ere infinit\'esimal. La congruence $|F| \equiv n \bmod 2$ s'en d\'eduit. L'autodualit\'e de $\pi_\infty$ est cons\'equence de celle de ${\rm L}(\pi_\infty)$ (comme repr\'esentation de ${\rm W}_\R$) et de la compatibilit\'e de  la param\'etrisation de Langlands au dual. \ps\ps 
Pour le (ii), on remarque que $\DET({\rm
L}(\pi_\infty))=1$ (compatibilit\'e au caract\`ere central), et on conclut car $\DET ( {\rm I}_w ) = \epsilon_{\C/\R}^{w+1}$ et 
$\omega_i \equiv {\rm w}(\pi) \bmod 2$ pour tout $i$. Le (iii) d\'ecoule du (ii) car
si $\pi$ est r\'eguli\`ere et $n \equiv 0 \bmod 2$, le (i) montre que
$F=\emptyset$ et $|E|=n/2$. Si $|{\rm Poids}(\pi)|=n-1$, le (i) montre que $0$ est l'unique poids double et $n \equiv 0 \bmod 2$, donc ${\rm w}(\pi) \equiv 0 \bmod 2$, $|F|=2$ et $|E|=n/2-1$,  puis (ii) entra\^ine (iv).
\end{pf}

Cet \'enonc\'e et la proposition~\ref{parstalg} montrent que si $G$ un $\Z$-groupe classique, et si $\pi \in \Pi_{\rm
disc}(G)$ est telle que ${\rm Inf}_{\pi_\infty}={\rm Inf}_V$ pour un certain $V \in {\rm Irr}(\widehat{G})$,  alors $\psi(\pi,{\rm St})$
satisfait certaines
contraintes de nature combinatoire, que l'on r\'esume dans
l'\'enonc\'e suivant \cite[Lemma 3.23]{chrenard2}.
\begin{cor}\label{cormodulo4} Soient $G$ un $\Z$-groupe classique et ${\rm St} : \widehat{G} \rightarrow {\rm SL}_n$ la repr\'esentation standard de $\widehat{G}$. On suppose que $\psi=\oplus_{i=1}^r \pi_i[d_i] \in \mathcal{X}_{\rm AL}({\rm SL}_n)$, o\`u $\pi_i \in
\Pi_{\rm cusp}(\PGL_{n_i})$ pour tout $i=1,\dots,k$, et que $\psi_\infty={\rm St}({\rm Inf}_V)$ o\`u $V \in {\rm Irr}(G_\C)$.
\begin{itemize}\ps  \ps
\item[(i)] Si $\widehat{G}(\C) \simeq {\rm SO}_{2g+1}(\C)$, alors il
existe un unique $1 \leq i_0 \leq k$ tel que $n_{i_0} d_{i_0} \equiv 1 \bmod
2$. De plus, $n_i d_i \equiv 0 \bmod 4$ pour tout $i \neq i_0$. \ps\ps 
\item[(ii)] Si $\widehat{G}(\C) \simeq {\rm Sp}_{2g}(\C)$, alors $n_i
d_i \equiv 0 \bmod 2$ pour tout $i$. \ps\ps 
\item[(iii)] Si $\widehat{G}(\C) \simeq {\rm SO}_n(\C)$ avec $n \equiv 0 \bmod 4$, alors $n_i 
d_i \equiv 0 \bmod 4$ pour tout $i$, except\'e dans le cas exceptionnel
suivant : $0$ est une valeur propre double de ${\rm St}({\rm
c}_\infty(\pi))$ et il existe exactement $2$ entiers $i$, disons $i_1$ et
$i_2$, tels que $n_id_i \not \equiv 0 \bmod 4$. Ils satisfont
$n_{i_1} d_{i_1} n_{i_2} d_{i_2} \equiv 3 \bmod 4$.
\end{itemize}
\end{cor}

\begin{pf} Supposons $\widehat{G}(\C) \simeq {\rm SO}_{2g+1}(\C)$. Comme on a $\sum_{i=1}^k n_id_i=2g+1$, il existe au moins un entier $i_0$ tel que $n_{i_0}d_{i_0} \equiv 1 \bmod 2$. Pour un tel entier, ${\rm c}_\infty(\pi_{i_0})$ admet la valeur propre $0$ car $\pi_{i_0}=\pi_{i_0}^\vee$ et $n_{i_0}$ est impair. Comme $d_{i_0}$ est \'egalement impair, $[d_{i_0}]_\infty$ admet \'egalement la valeur propre $0$, ainsi donc que $(\pi_{i_0}[d_{i_0}])_\infty$. La premi\`ere partie du (i) s'en d\'eduit car $0$ est une valeur propre simple de $\psi_\infty$ d'apr\`es le cas I de l'analyse du~\S\ref{paralgreg}. Pour la seconde, on observe que pour tout $i=1,\dots,k$ on a ${\rm w}(\pi_i) \equiv d_i -1 \bmod 2$, ce qui conclut quand ${\rm w}(\pi_i)$ est impair, le cas \'ech\'eant r\'esultant du (iii) de la proposition~\ref{contraintesalg}. \par 

Si $\widehat{G}(\C) \simeq {\rm Sp}_{2g}(\C)$, la relation ${\rm w}(\pi_i) + d_i -1 \equiv 1 \bmod 2$ pour tout $i=1,\dots,k$ montre que si $d_i$ est impair alors ${\rm w}(\pi_i)$ l'est aussi, et donc $n_i$ est pair (Proposition \ref{propevidentealg}). Cela d\'emontre le (ii). \par

Supposons enfin $\widehat{G}(\C) \simeq {\rm SO}_n(\C)$ avec $n \equiv 0 \bmod 4$. En particulier, ${\rm w}(\pi_i) \equiv d_i-1 \bmod 2$ pour tout $i=1,\cdots,k$. Si ${\rm w}(\pi_i) \equiv 1 \bmod 2$ (et donc $n_i \equiv 0 \bmod 2$), alors $d_i$ est pair et donc $n_i d_i \equiv 0 \bmod 4$.  Si ${\rm w}(\pi_i) \equiv 0 \bmod 2$ et $n_i \equiv 0 \bmod 2$, la proposition~\ref{contraintesalg} (iii) assure que $n_i \equiv 0 \bmod 4$, sauf peut-\^etre si $\pi_i$ n'est pas r\'eguli\`ere. Cette \'eventualit\'e, si elle se produit, se produit pour un unique $i_0$ d'apr\`es la proposition~\ref{parstalg}, et dans ce cas $n_i \equiv 0 \bmod 2$ pour tout $i \neq i_0$, de sorte qu'au final $n_i d_i \equiv 0 \bmod 4$ pour tout $i \neq i_0$. Le r\'esultat est cons\'equence dans ce cas de l'identit\'e $n = \sum_{i=1}^k n_i d_i \equiv n_{i_0} d_{i_0} \bmod 4$ et de l'hypoth\`ese $n \equiv 0 \bmod 4$. On peut donc \'ecarter l'\'eventualit\'e ci-dessus et supposer que $\pi_i$ est r\'eguli\`ere pour tout $i$, en particulier $n_i d_i \equiv 0 \bmod 4$ pour tout $i$ tel que $n_i$ est pair. Soit $J \subset \{1,\dots,k\}$ l'ensemble des $i$ tels que $n_i$ est impair (auquel cas $d_i$ est impair et ${\rm w}(\pi_i)$ est pair), on peut supposer $J$ non vide. D'apr\`es l'argument donn\'e au (i), cela entra\^ine que $0$ est valeur propre double de $\psi_\infty$ et que $|J|\leq 2$. On conclut car $n \equiv \sum_{j \in J} n_j d_j \bmod 4$. 
\end{pf}

\subsection{Conjecture de Ramanujan et repr\'esentations
galoisiennes}

\label{repgal}

Un cas particulier des conjectures de Langlands, dans l'esprit de la fameuse
conjecture de Shimura-Taniyama-Weil, est que l'ensemble des fonctions ${\rm L}$ de la forme ${\rm L}(s+\frac{{\rm w}(\pi)}{2}+m,\pi)$, o\`u
$m \in \Z$ et $\pi$ parcourt les repr\'esentations alg\'ebriques des
$\Pi_{\rm cusp}({\rm PGL}_n)$ pour $n\geq 1$, devrait co\"incider exactement
avec celui des fonctions ${\rm L}$ des motifs sur $\Q$ ayant bonne
r\'eduction partout (et disons ``\`a coefficients dans $\overline{\Q}$'' et ``simples'')
\cite{Lgl} \cite{motives}.  La
conjecture de Ramanujan pour un $\pi$ alg\'ebrique (\S V.\ref{ramanujan}) serait alors
cons\'equence de l'existence du motif ${\rm M}(\pi)$ associ\'e et des
conjectures de Weil, d\'emontr\'ees par Deligne.  Gr\^ace aux travaux de
nombreux math\'ematiciens (dont Eichler-Shimura, Deligne, Langlands,
Kottwitz, Clozel, Harris-Taylor, Waldspurger, Ng\^o, Laumon, 
Clozel-Harris-Labesse, Shin, Chenevier-Harris) on dispose
d\'esormais d'une construction affaiblie de ${\rm M}(\pi)$ pour les $\pi$ 
alg\'ebriques r\'eguli\`eres et polaris\'ees qui est toutefois suffisante 
pour d\'emontrer le th\'eor\`eme suivant.  Si $\pi \in \Pi_{\rm
cusp}(\PGL_n)$ est une telle repr\'esentation, il est connu que le
polyn\^ome caract\'eristique\footnote{Dans la d\'efinition de ${\rm
P}_p(\pi)$, il est sous-entendu que ${\rm c}_p(\pi)\,p^{w(\pi)/2}$ d\'esigne la classe de
conjugaison semisimple de ${\rm GL}_n(\C)$ obtenue en faisant le produit de
la classe ${\rm c}_p(\pi)$, vue dans ${\rm GL}_n(\C) \supset {\rm
SL}_n(\C)$, par le scalaire
$p^{w(\pi)/2} \in \C^\ast \subset \GL_n(\C)$.}
$${\rm P}_p(\pi) = \DET(t -  {\rm c}_p(\pi)\,p^{w(\pi)/2}) \, \, \in \, \, \C[t]$$ est \`a coefficients dans le sous-corps
$\overline{\Q} \subset \C$ des nombres alg\'ebriques.  Le th\'eor\`eme
suivant est d\'emontr\'e dans \cite{chl},\cite{shin},\cite{chharris} et
\cite{clozelram} (voir aussi \cite{caraiani}).

\begin{thm}\label{existencegalois} Soit $\pi \in \Pi_{\rm cusp}^\bot(\PGL_n)$ alg\'ebrique et
r\'eguli\`ere.  \begin{itemize} \ps\ps  \item[(i)] $\pi$
satisfait la conjecture de Ramanujan.  \ps\ps  \item[(ii)] Soient $\ell \in {\rm
P}$, $\overline{\Q}_\ell$ une cl\^oture alg\'ebrique de $\Q_\ell$, et $\iota : \overline{\Q} \rightarrow
\overline{\Q}_\ell$ un plongement.
Il existe une repr\'esentation continue $\rho_{\pi,\iota} : {\rm Gal}(\overline{\Q}/\Q)
\longrightarrow \GL_n(\overline{\Q}_\ell)$, unique \`a isomorphisme pr\`es, qui est semi-simple, non ramifi\'ee  
hors de $\ell$, et telle que  $$\DET(t - \rho_{\pi,\iota}({\rm
Frob}_p))=\iota({\rm P}_p(\pi))$$
pour tout nombre premier $p \neq \ell$. 
\end{itemize} 
\end{thm}

Dans cet \'enonc\'e, ${\rm Frob}_p$ d\'esigne une classe de conjugaison
d'\'el\'ements de Frobenius arithm\'etiques en $p$. On sait de plus
d'apr\`es {\it loc. cit.} que la restriction de la repr\'esentation $\rho_{\pi,\iota}$ \`a ${\rm
Gal}(\overline{\Q}_\ell/\Q_\ell)$ est cristalline au sens de Fontaine, de
poids de Hodge-Tate les $\lambda + \frac{{\rm w}(\pi)}{2}$ o\`u $\lambda \in   
{\rm Poids}(\pi)$. On conjecture que $\rho_{\pi,\lambda}$ est irr\'eductible mais ce n'est connu que si $n\leq 5$ \cite{calegarigee}.
Observons que l'autodualit\'e de $\pi$ et le th\'eor\`eme de Cebotarev 
entra\^inent l'isomorphisme
\begin{equation}\label{rhopilautodual} \rho_{\pi,\iota}^\ast \simeq \rho_{\pi,\iota} \otimes
\omega_\ell^{-{\rm w}(\pi)},\end{equation}
o\`u $\omega_\ell : {\rm Gal}(\overline{\Q}/\Q) \rightarrow \Z_\ell^\times$ d\'esigne le caract\`ere cyclomotique $\ell$-adique. On sait d\'emontrer que si ${\rm w}(\pi) \equiv 1 \bmod 2$ (resp. ${\rm w}(\pi) \equiv 0 \bmod 2$) alors il existe un accouplement ${\rm Gal}(\overline{\Q}/\Q)$-\'equivariant non-d\'eg\'en\'er\'e $\rho_{\pi,\iota} \otimes \rho_{\pi,\iota} \rightarrow \omega_\ell^{{\rm w}(\pi)}$ qui est altern\'e (resp. sym\'etrique) \cite{bchsign}.

\begin{remark} \label{remrepgal} {\rm On s'attend \`a ce que le {\it (ii)}
du th\'eor\`eme soit valable
sans supposer $\pi$ r\'eguli\`ere ou autoduale. Des travaux r\'ecents de
Harris-Lan-Taylor-Thorne et Sholze montrent que l'on peut s'affranchir de 
l'hypoth\`ese d'autodualit\'e (mais ne d\'emontrent pas {\it (i)} pour ces
$\pi$).
Mentionnons enfin que si $\pi \in \Pi_{\rm cusp}^\bot(\PGL_n)$ est alg\'ebrique
et satisfait $|{\rm Poids}(\pi)| =n-1$ et $n \equiv 0 \bmod 4$, on sait \'egalement d\'emontrer
{\it (ii)}~\cite{goldring}, mais pas {\it (i)}. }
\end{remark}

\begin{cor}\label{corgaloisrepo} Soient $G$ un $\Z$-groupe classique et $\pi
\in \Pi_{\rm disc}(G)$ telle que $\pi_\infty$ a m\^eme caract\`ere
infinit\'esimal qu'une repr\'esentation irr\'eductible de dimension finie de
$G(\C)$.  On suppose la conjecture d'Arthur-Langlands satisfaite pour
$(\pi,{\rm St})$.  Si $\widehat{G}(\C) \simeq {\rm SO}_m(\C)$, on suppose de
plus $m \not \equiv 2 \bmod 4$.  \ps\ps  Soient $\ell \in {\rm P}$, $\overline{\Q}_\ell$ une cl\^oture alg\'ebrique de $\Q_\ell$, et $\iota :
\overline{\Q} \rightarrow \overline{\Q}_\ell$ un plongement.  Il existe une
unique repr\'esentation continue, semi-simple, non ramifi\'ee hors de $\ell$
$$\rho_{\pi,\iota} : {\rm Gal}(\overline{\Q}/\Q) \longrightarrow \GL_{{\rm
n}_{\widehat{G}}}(\overline{\Q}_\ell)$$ telle que $\forall p \in {\rm P}
- \{\ell\}$, $\DET(t - \rho_{\pi,\iota}({\rm
Frob}_p))=\iota(\DET(t  -  \, \, {\rm St}({\rm c}_p(\pi))\, p^{\frac{{\rm w}_{\widehat{G}}}{2}}))$. 
\end{cor}

Le fait que $\det(t  -  \, \, {\rm St}({\rm c}_p(\pi))\, p^{\frac{{\rm w}_{\widehat{G}}}{2}}) \in \overline{\Q}[t]$ fait
partie de l'assertion (et peut se v\'erifier directement simplement dans les cas qui nous int\'eressent). 

\begin{pf} \'Ecrivons $\psi(\pi,{\rm St})=\oplus_{i=1}^k \pi_i[d_i]$. 
Le th\'eor\`eme~\ref{existencegalois} et la remarque ci-dessus s'appliquent
aux repr\'esentations automorphes $\pi_i$, d'apr\`es la proposition \ref{parstalg} et le corollaire~\ref{cormodulo4}
(iii). Il suffit alors de poser 
$$\rho_{\pi,\iota} = \bigoplus_{i=1}^k \rho_{\pi_i,\iota} \otimes
(\oplus_{j=0}^{d_i-1} \omega_\ell^j) \otimes 
\omega_\ell^{\frac{{\rm w}_{\widehat{G}}-{\rm w}(\pi_i)+1-d_i}{2}}.$$
L'assertion d'unicit\'e d\'ecoule du th\'eor\`eme de Cebotarev.
\end{pf}

Sp\'ecifions ce r\'esultat dans le cas de ${\rm O}_n$, en utilisant la formule
VI.\eqref{formulesattp}. 

\begin{cor}\label{correpgalOn} Soient $n\equiv 0 \bmod 8$ et $F \in {\rm M}_U({\rm O}_n)$
une forme propre pour ${\rm H}({\rm O}_n)$, disons telle que ${\rm T}_p(F)= \lambda_p F$ 
pour tout premier $p$. On suppose v\'erif\'ee la conjecture
d'Arthur-Langlands pour le couple $(\pi,{\rm St})$ o\`u $\pi \in \Pi_{\rm
disc}({\rm O}_n)$ est la repr\'esentation engendr\'ee par $F$ (\S VI.\ref{exconjal}). \ps
Soient $\ell \in {\rm P}$, $\overline{\Q}_\ell$ une cl\^oture alg\'ebrique
de $\Q_\ell$, et $\iota : \overline{\Q} \rightarrow \overline{\Q}_\ell$ un plongement. Il existe une unique repr\'esentation continue, semi-simple, non ramifi\'ee  
hors de $\ell$, $\rho_{F,\iota}: {\rm Gal}(\overline{\Q}/\Q) \longrightarrow {\rm
GL}_n(\overline{\Q}_\ell)$, 
telle que $\forall p \in {\rm P} \backslash \{\ell\}$, ${\rm trace} \, \rho_{F,\iota}({\rm Frob}_p) =
\iota(\lambda_p)$. 
\end{cor}

\subsection{Fonctions ${\rm L}$ de paires de repr\'esentations alg\'ebriques}\label{parfaceps}

Soient $\pi \in \Pi_{\rm cusp}(\PGL_n)$ et $\pi' \in \Pi_{\rm cusp}(\PGL_{n'})$, la fonction ${\rm L}$ du couple $\{\pi,\pi'\}$ est d\'efinie par le produit eul\'erien
$${\rm L}(s,\pi \times \pi' ) = \prod_{p \in {\rm P}} \DET ({\rm I}_{nn'} - p^{-s} \, \, {\rm c}_p(\pi) \otimes {\rm c}_p(\pi'))^{-1}.$$
C'est un cas particulier de la construction de Langlands rappel\'ee au \S V.\ref{lienfonctionl} o\`u ${\rm G} = {\rm PGL}_n \times {\rm PGL}_{n'}$ et o\`u $r$ est le produit tensoriel des repr\'esentations standards de ${\rm SL}_n$ et ${\rm SL}_{n'}$. On retrouve ${\rm L}(s,\pi)$ lorsque $\pi' = 1$ est la repr\'esentation triviale de $\PGL_1$. Cette fonction ${\rm L}$ a \'et\'e \'etudi\'ee par Rankin et Selberg quand $n=n'=2$ et pour tout $n,n'$ par Jacquet, Piatetski-Shapiro et Shalika. Ces auteurs d\'emontrent que le produit eul\'erien ci-dessus est absolument convergent quand ${\rm Re} \,s > 1$ et qu'il admet un prolongement m\'eromorphe \`a $\C$ tout entier. De plus, ce prolongement analytique est une fonction enti\`ere de $s$, sauf si $\pi' \simeq \pi^\vee$ auquel cas $s=1$ en est l'unique p\^ole (et il est simple). Enfin, si ${\rm L}_\infty(s,\pi \times \pi')$ est un produit convenable de {\it facteurs $\Gamma$}, et si $\xi(s,\pi \times \pi')\,=\,{\rm L}_\infty(s,\pi \times \pi')\,{\rm L}(s,\pi\times \pi')$, on dispose d'une \'equation fonctionnelle de la forme $$\xi(s,\pi \times \pi')\, =\, \varepsilon(\pi \times \pi')\, \xi (1-s,\pi^\vee \times (\pi')^\vee),$$
o\`u $\varepsilon(\pi \times \pi') \in \C^\times$. Nous renvoyons au cours de Cogdell~\cite[\S 9]{cogdell} pour un expos\'e synth\'etique de ces r\'esultats. Ajoutons que si $\pi$ et $\pi'$ sont autoduales, la relation $\xi(s,\pi \times \pi') \,=\, \varepsilon(\pi \times \pi')\, \xi(1-s,\pi \times \pi')$ entra\^ine que $\varepsilon(\pi \times \pi') = \pm 1$ est un simple signe. Il nous sera utile de rappeler la recette pr\'ecise pour $\varepsilon(\pi \times \pi')$ et ${\rm L}_\infty(s,\pi \times \pi')$. 
Ils ne d\'ependent tous deux que des composantes archim\'ediennes de $\pi$ et $\pi'$. Pour simplifier, nous restreindrons cette discussion au cas o\`u $\pi$ et $\pi'$ sont des repr\'esentations alg\'ebriques, le seul cas dont nous aurons besoin. \ps\ps 
Soit ${\rm Rep}_{\rm alg}({\rm W}_\R)$ l'ensemble des classes d'isomorphisme de repr\'esentations continues et semi-simples de ${\rm W}_\R$, sur des $\C$-espaces vectoriels de dimension finie, qui sont triviales sur le sous-groupe $\R_{>0} \subset \C^\times$ de ${\rm W}_\R$. Les \'el\'ements de ${\rm Rep }_{\rm alg}({\rm W}_\R)$ sont exactement les sommes directes de repr\'esentations de la forme $1, \epsilon_{\C/\R}$ o\`u ${\rm I}_w$ pour $w>0$ (\S\ref{apparitionWR}). Suivant Weil, il existe une unique mani\`ere d'associer \`a tout $\rho \in {\rm Rep}_{\rm alg}({\rm W}_\R)$ une racine $4$\`eme de l'unit\'e $\varepsilon(\rho) \in \{1,i,-1,-i\}$ et une fonction m\'eromorphe $\Gamma(s,\rho)$ de la variable complexe $s$, de sorte que pour tout $\rho, \rho' \in {\rm Rep}_{\rm alg}({\rm W}_\R)$ on ait $$\varepsilon(\rho \oplus \rho')= \varepsilon(\rho)\varepsilon(\rho'), \, \, \Gamma(s,\rho\oplus\rho')=\Gamma(s,\rho)\Gamma(s,\rho'),$$ ainsi que : \begin{itemize} \ps\ps 
\item[(i)] $\varepsilon({\rm I}_w)=i^{w+1}$ et $\Gamma(s,{\rm
I}_w)=\Gamma_\C(s+\frac{w}{2})$ pour tout $w\geq0$,\ps\ps
\item[(ii)] $\varepsilon(1)=1$ et $\Gamma(s,1)=\Gamma_\R(s)$. \ps\ps
\end{itemize}\ps\ps 

\noindent On rappelle que $\Gamma(s)= \int_{0}^{\infty} e^{-t} \,t^s \,\frac{dt}{t}$
si ${\rm Re}\, s >0$, et qu'il est coutume de poser
$$\Gamma_\R(s)=\pi^{-\frac{s}{2}}\Gamma(\frac{s}{2}) \, \, \,{\rm et}\, \, \, \Gamma_\C(s)=2(2\pi)^{-s}\Gamma(s),$$ de sorte que
$\Gamma_\C(s)=\Gamma_\R(s)\Gamma_\R(s+1)$ (formule de duplication). Noter
que du cas $w=0$ nous d\'eduisons 
$\varepsilon(\epsilon_{\C/\R})=i$ et $\Gamma(s,\epsilon_{\C/\R})=
\Gamma_\R(s+1)$. \ps \ps

\begin{prop}\label{epsgammaalg} Soient $\pi \in \Pi_{\rm cusp}(\PGL_n)$ et $\pi' \in \Pi_{\rm cusp}(\PGL_{n'})$ alg\'ebriques. Posons $\rho={\rm L}(\pi_\infty)\otimes {\rm L}(\pi'_\infty)$. On a
$$\varepsilon(\pi \times \pi') = \varepsilon(\rho)\, \, {\rm et}\, \,  {\rm L}_\infty(s,\pi \times \pi')=\Gamma(s,\rho).$$
\end{prop}

\begin{pf} L'\'enonc\'e portant sur sur ${\rm L}_\infty(s,\pi \times \pi')$ a bien un sens car ${\rm L}(\pi_\infty)$ et ${\rm L}(\pi'_\infty)$ sont dans ${\rm Rep }_{\rm alg}({\rm W}_\R)$ d'apr\`es la proposition \ref{contraintesalg} (i). L'assertion sur ${\rm L}_\infty(s,\pi \times \pi')=\Gamma(s,\rho)$ vaut par d\'efinition \cite[Ch. 9]{cogdell}. Une inspection des formules dans \cite[\S 3]{tate} montre que pour tout $\rho \in {\rm Rep}_{\rm alg}({\rm W}_\R)$ le nombre $\varepsilon(\rho)$ d\'efini ci-dessus est exactement celui not\'e $\varepsilon(\rho,\psi,dx)$ {\it loc. cit.} o\`u $dx$ est la mesure de
Lebesgue sur $\R$ et $\psi : \R \rightarrow \C^\times$ le caract\`ere $x \mapsto e^{2i\pi x}$. Comme $\pi$ et $\pi'$ sont ``non ramifi\'ees \`a toutes les places finies'' dans la terminologie usuelle, ce facteur $\varepsilon(\rho)$ co\"incide donc avec $\varepsilon(\pi \times \pi')$ \cite[Ch. 9]{cogdell}.  
\end{pf}

Il d\'ecoule de ces formules que $\varepsilon(\pi\times \pi')$ est une fonction explicite des poids de $\pi$ et $\pi'$. Il est utile \`a ce stade d'observer que ${\rm I}_w \otimes \epsilon_{\C/\R} \simeq {\rm I}_w$ et 
$${\rm I}_w \otimes {\rm I}_{w'} \simeq {\rm I}_{w+w'} \oplus {\rm I}_{|w-w'|}$$
pour tous entiers $w, w' \geq 0$. En particulier, $\varepsilon({\rm I}_w \otimes {\rm I}_{w'})=(-1)^{1+{\rm Max}(w,w')}$. 


\section{La formule de multiplicit\'e d'Arthur} \label{paramarthur}

\subsection{L'alternative symplectique-orthogonal d'Arthur}\label{altsoarth}

Rappelons que si $H$ est un $\C$-groupe classique, nous notons ${\rm n}_H$ la dimension de sa
repr\'esentation standard.

\begin{thmv}\label{seed} {\rm (Arthur)} Soit $\pi \in \pcb({\rm PGL}_n)$. Il existe un $\Z$-groupe
classique de Chevalley $G^\pi$, unique \`a isomorphisme pr\`es, ayant les propri\'et\'es suivantes :\ps\ps 
\begin{itemize}
\item[(i)] ${\rm n}_{\widehat{G^\pi}}=n$, \ps\ps 
\item[(ii)] il existe $\pi' \in \Pi_{\rm disc}(G^\pi)$ telle que ${\rm c}(\pi)=\psi(\pi',{\rm St})$. \ps\ps
\end{itemize}
\end{thmv}

C'est un cas particulier de~\cite[Thm.1.4.1 \& Thm. 1.5.2]{arthur} (voir aussi la {\it m\'ethode de descente} de Ginzburg, Rallis et Soudry \cite{grs} pour un \'enonc\'e affaibli). Par d\'efinition, le groupe $G^\pi$ satisfait ${\rm n}_{\widehat{G^\pi}}=n$. Lorsque $n$ est impair, la seule possibilit\'e est donc
$G^\pi \simeq {\rm Sp}_{n-1}$, mais lorsque $n$ est pair $G^\pi$ est
isomorphe \`a ${\rm SO}_{\frac{n}{2},\frac{n}{2}}$ ou ${\rm
SO}_{\frac{n}{2}+1,\frac{n}{2}}$ (exclusivement). Si $n=2$ alors $G^\pi \simeq {\rm SO}_{2,1} \simeq {\rm PGL}_2$, \'etant donn\'e que ${\rm SO}_{1,1}
\simeq \mathbb{G}_m$ n'est pas semi-simple. Lorsque $n=1$, de sorte que $\pi$
est la repr\'esentation triviale de $\PGL_1$, on a enfin $G^\pi={\rm
SO}_{1,0}$ (le
$\Z$-groupe trivial). \ps\ps 

La repr\'esentation $\pi \in \Pi_{\rm cusp}^\bot({\rm PGL}_n)$ est dite
orthogonale si $\widehat{G^\pi}(\C) \simeq {\rm SO}_m(\C)$ pour un certain
entier $m\geq 1$ (ou encore si ${\rm w}_{\widehat{G^\pi}}=0$ dans les
notations du~\S\ref{introchaparthur}), et symplectique sinon.  \ps\ps 

\begin{propv}\label{altso} Soit $\pi \in \Pi_{\rm cusp}^\bot(\PGL_n)$
alg\'ebrique.  On suppose que $\pi$ poss\`ede au moins un poids qui soit
valeur propre simple de ${\rm c}_\infty(\pi)$.  Alors $\pi$ est
symplectique si, et seulement si, ${\rm w}(\pi) \equiv 1 \bmod 2$. 
\end{propv}

\begin{pf} C'est une variante de~\cite[Cor. 3.8]{chrenard2}. D'apr\`es Arthur \cite[Thm. 11.4.2]{arthur}, la repr\'esentation ${\rm L}(\pi_\infty)$ de ${\rm W}_\R$ sur $\C^n$ pr\'eserve une forme bilin\'eaire non-d\'eg\'en\'er\'ee $b$, qui est altern\'ee si $\pi$ est symplectique et sym\'etrique sinon. L'hypoth\`ese faite sur $\pi$ entra\^ine qu'au moins une des repr\'esentations $1, \varepsilon_{\C/\R},$ o\`u ${\rm I}_w$ avec $w>0$, intervient dans ${\rm L}(\pi_\infty)$ avec multiplicit\'e $1$ (Proposition~\ref{contraintesalg} (i)) ; notons $E \subset \C^n$ le sous-espace correspondant. Chacune de ces repr\'esentations \'etant irr\'eductible et autoduale, la restriction de $b$ \`a $E$ est non-d\'eg\'en\'er\'ee. \'Etant donn\'e la relation $\DET \,{\rm I}_w = \varepsilon_{\C/\R}^{w+1}$, on constate que $b$ est altern\'ee si, et seulement si, on a $E \simeq {\rm I}_w$ avec $w \equiv 1 \bmod 2$. 
\end{pf} 

Les r\'esultats d'Arthur admettent \'egalement des cons\'equences concernant
les fonctions ${\rm L}$ de paires de repr\'esentations autoduales : voir
\cite[Thm.  1.5.3]{arthur}.  En particulier, si $\pi \in \Pi_{\rm
cusp}^\bot(\PGL_n)$ et $\pi' \in \Pi_{\rm cusp}^\bot(\GL_m)$ sont soit
toutes deux symplectiques, soit toutes deux orthogonales, alors
$\varepsilon(\pi\times \pi')=1$ (c'est un Th\'eor\`eme
{\color{green}\!\!${}^\star$}).  Dans le cas o\`u $\pi'=1$, on en d\'eduit que
$$\varepsilon(\pi):=\varepsilon(\pi \times 1)$$ vaut $1$ si $\pi$ est
orthogonale.  Quand $\pi$ est alg\'ebrique autoduale et orthogonale, cela
donne une relation non-triviale sur ses poids : voir~\cite[Prop. 
1.8]{chrenard2}.

\subsection{La formule de multiplicit\'e: hypoth\`eses g\'en\'erales}\label{amfhypgen}

Soient $G$ un $\Z$-groupe classique et $n={\rm n}_{\widehat{G}}$ la dimension de la repr\'esentation standard ${\rm St}$ de $\widehat{G}$. Fixons un entier $k\geq 1$, ainsi que pour tout $i=1,\dots,k$ un couple $(\pi_i,d_i)$ o\`u $d_i\geq 1$ est un entier et $\pi_i \in \Pi_{\rm cusp}^\bot(\PGL_{n_i})$. On suppose que $n=\sum_{i=1}^k n_i d_i$ et on consid\`ere l'\'el\'ement
$$\psi = \oplus_{i=1}^k \pi_i[d_i]$$
de $\mathcal{X}(\SL_n)$. \ps\ps 

Soit  $U$ une repr\'esentation irr\'eductible unitaire de $G(\R)$. La formule de multiplicit\'e d'Arthur, conjectur\'ee en toute g\'en\'eralit\'e dans \cite{arthurunipotent}, et d\'emontr\'ee dans \cite{arthur} lorsque $G$ est un groupe classique de Chevalley, donne une condition n\'ecessaire est suffisante pour qu'il existe un $\pi \in \Pi_{\rm disc}(G)$ telle que $\pi_\infty \simeq U$ et $\psi(\pi,{\rm St})=\psi$. Elle s'exprime sous la forme d'une formule d'orthogonalit\'e entre deux caract\`eres sur un $2$-groupe ab\'elien fini \'el\'ementaire ${\rm C}_\psi$, que nous allons expliciter dans les paragraphes qui suivent. Le premier de ces caract\`eres, not\'e $\varepsilon_\psi$ et d\'ecrit 
au~\S\ref{parepsilonpsi},  est ind\'ependant de $U$. Il est introduit en grande g\'en\'eralit\'e par Arthur dans \cite{arthurunipotent} et tient compte des signes $\varepsilon(\pi_i \times \pi_j)$ (\S\ref{parfaceps}) selon une combinatoire bien pr\'ecise. Le second de ces caract\`eres, qui est le plus fin des deux, a son origine qui remonte aux travaux de Shelstad \cite{shelstadlind} (voir aussi \cite{knappzuckerman,labesselanglands,langlandsshelstad,AJ,arthurunipotent,arthur,abv,shelstad0,shelstad3}). Il ne d\'epend essentiellement que de $U$ et d'un certain morphisme ${\rm SL}_2(\C) \times {\rm W}_\R \rightarrow \widehat{G}$ associ\'e \`a $\psi$;  il sera d\'ecrit aux~\S\ref{parpaqarch} et~\S\ref{paramaj}. \ps\ps

Les travaux d'Arthur \cite{arthur} sont tr\`es g\'en\'eraux, et nous ne les appliquerons que dans des cas tr\`es particuliers, pour lesquels les \'enonc\'es sont sensiblement simplifi\'es. Nous ferons les deux hypoth\`eses suivantes :  \begin{itemize}\ps\ps

\item[(H1)] $n \not \equiv 2 \bmod 4$, \ps\ps  

\item[(H2)] $\psi_\infty={\rm St}({\rm Inf}_V)$ o\`u $V \in {\rm Irr}(G_\C)$ (\S\ref{paralgreg}). \ps\ps 
\end{itemize} 

La premi\`ere hypoth\`ese n'est une contrainte que si $G$ est un groupe sp\'ecial orthogonal pair. Dans ce cas, $G(\R)$ a pour signature $(p,q)$ avec $p \equiv q \bmod 8$ d'apr\`es le Scholie II.2.2 (b), de sorte que l'hypoth\`ese s'\'ecrit \'egalement $p \equiv q \equiv 0  \bmod 2$ (elle est bien s\^ur satisfaite si $G={\rm SO}_n$ avec $n \equiv 0 \bmod 8$). \ps\ps 


La seconde hypoth\`ese, portant sur $\psi_\infty$, a \'et\'e explicit\'ee au \S\ref{paralgreg} (cas I, II ou III), o\`u nous en avons donn\'e quelques cons\'equences de nature combinatoire portant sur les $\pi_i$. En particulier :\begin{itemize}\ps\ps

\item[(a)] pour tout $i=1,\dots,k$, $\pi_i$ est autoduale alg\'ebrique (et m\^eme r\'eguli\`ere, hormis dans un cas exceptionnel),\ps\ps 
\item[(b)] pour tout $i=1,\dots,k$, ${\rm w}(\pi_i) + d_i -1 \equiv {\rm w}_{\widehat{G}} \bmod 2$, \ps\ps 
\item[(c)] pour tous $i \neq j$, si $(n_i,d_i)=(n_j,d_j)$ alors $\pi_j \not \simeq \pi_i$. \ps\ps 
\end{itemize}
\ps\ps

En effet, (a) et (b) d\'ecoulent de la proposition~\ref{parstalg}. Le point (c), qui n'est non-trivial sous (H2) que si $n_i=n_j=1$ et $d_i=d_j=1$, d\'ecoule de (H1) et du corollaire \ref{cormodulo4} (iii). \ps\ps 

\subsection{Le groupe ${\rm C}_\psi$ et le caract\`ere $\varepsilon_\psi$}\label{parepsilonpsi}

On conserve les hypoth\`eses et notations du paragraphe pr\'ec\'edent. D'apr\`es le (b) ci-dessus et la proposition~\ref{altso}, observons qu'il existe un $\C$-morphisme 
$$ \nu : {\rm SL_2} \times \prod_{i=1}^k \widehat{G^{\pi_i}} \longrightarrow \widehat{G}$$
tel que la $\C$-repr\'esentation ${\rm St} \circ \nu$, d'espace sous-jacent $V \simeq \C^n$, se d\'ecompose comme somme directe $$V= \oplus_{i=1}^k V_i,$$
o\`u $V_i$ est isomorphe au produit tensoriel de la repr\'esentation ${\rm Sym}^{d_i-1} {\rm St}_2$ de ${\rm SL}_2$ par la repr\'esentation standard de $\widehat{G^{\pi_i}}$ (les autres facteurs $\widehat{G^{\pi_j}}$ pour $j \neq i$ agissant trivialement). 
Un tel morphisme $\nu$ n'est pas unique : il l'est seulement modulo composition \`a l'arriv\'ee par un automorphisme du $\C$-groupe $\widehat{G}$. Il sera fix\'e une fois pour toutes dans ce qui suit; nous discuterons en temps voulu de la d\'ependance en ce choix dans la formule finale. \ps\ps 

Soit ${\rm C}_\nu$ le centralisateur dans $\widehat{G}(\C)$ de l'image de $\nu$. La repr\'esentation ${\rm St}$ l'identifie au sous-groupe de ${\rm SL}(V)$ constitu\'e des \'el\'ements $g$ pr\'eservant chaque $V_i$ et tels que $g_{|V_i}= \epsilon_i \, {\rm Id}_{V_i}$, o\`u $(\epsilon_i) \in \{\pm 1\}^k$. Comme on a $\dim(V_i)=n_id_i$, le groupe ${\rm C}_\nu$ est donc dans une suite exacte naturelle 
$$ 1 \longrightarrow {\rm C}_\nu \overset{\rm St}{\longrightarrow} \{\pm 1\}^k \overset{\delta}{\longrightarrow}\{ \pm 1\},$$
o\`u $\delta(\epsilon_i)=\prod_{i=1}^k \epsilon_i^{n_id_i}$. Cette description abstraite de ${\rm C}_\nu$ est manifestement ind\'ependante du choix de $\nu$, c'est pourquoi nous le noterons simplement ${\rm C}_\psi$. \ps\ps 

Le centre ${\rm Z}_{\widehat{G}}$ de $\widehat{G}(\C)$ est un sous-groupe de ${\rm C}_\psi$. Notons $I \subset \{1,\dots,k\}$ le sous-ensemble constitu\'e des entiers $i$ tels que $n_i d_i \equiv 0 \bmod 2$, et pour tout $i \in I$ notons $$s_i \in {\rm C}_\psi$$ l'\'el\'ement agissant par $-1$ sur $V_i$ et $1$ sur $V_j$ pour $j \neq i$. D'apr\`es (H1) et (H2), le corollaire~\ref{cormodulo4} s'applique. Il entra\^ine $|I|\geq k-1$ ainsi que le lemme suivant. \ps\ps

\begin{lemme}\label{structcnu} ${\rm C}_\psi$ est engendr\'e par ${\rm Z}_{\widehat{G}}$ et les $s_i$ pour $i \in I$. \end{lemme}\ps\ps

Arthur d\'efinit ensuite \cite[p. 47]{arthur} un homomorphisme
$\varepsilon_\psi : {\rm C}_\psi \longrightarrow \{\pm 1\}$ trivial sur
${\rm Z}_{\widehat{G}}$.  Pour le d\'ecrire, il suffit de donner sa valeur
sur l'\'el\'ement $s_i$ pour $i \in I$.  Arthur consid\`ere pour cela la
restriction \`a $\nu$ de la repr\'esentation adjointe de $\widehat{G}$ sur
${\rm Lie} \,\widehat{G}$, c'est donc une repr\'esentation du produit ${\rm
C}_\nu \times \SL_2 \times (\prod_{i=1}^k \widehat{G^{\pi_i}})$.  L'entier
$i \in I$ \'etant fix\'e, c'est un exercice de v\'erifier que le sous-espace
de ${\rm Lie} \,\widehat{G}$ sur lequel $s_i$ agit par $-1$ est isomorphe
\`a $\bigoplus_{j \neq i} V_j \otimes V_i$ comme repr\'esentation de $\SL_2
\times (\prod_{i=1}^k \widehat{G^{\pi_i}})$.  Mais si pour $d\geq 1$ on note
$r_d$ la repr\'esentation ${\rm Sym}^{d-1} {\rm St}_2$ de ${\rm SL}_2$
(o\`u ${\rm St}_2$ d\'esigne la repr\'esentation standard, \S
VI.\ref{parconjarthlan}), et
si $a \geq b$ sont des entiers $\geq 1$, il est bien connu que $$r_a \otimes
r_b \simeq \oplus_{i=1}^b r_{a-b+2i-1};$$ en particulier, $r_a \otimes r_b$
admet ${\rm Min}(a,b)$ facteurs irr\'eductibles pour tout $a,b\geq 1$.  La
recette pour $\varepsilon_\psi$ d\'ecrite par Arthur {\it loc.  cit.} prend
donc la forme explicite suivante, o\`u l'on a incorpor\'e le
r\'esultat$^{\color{green}\ast}$ d'Arthur affirmant que $\varepsilon(\pi
\times \pi')=1$ si $\pi$ et $\pi'$ sont simultan\'ement symplectiques ou
orthogonales.\ps\ps

\begin{defprop} Il existe un unique homomorphisme $\varepsilon_\psi : {\rm C}_\psi \rightarrow \{\pm 1\}$ trivial sur ${\rm Z}_{\widehat{G}}$ et tel que $\forall i \in I, \, \varepsilon_\psi(s_i) = \prod_{j \neq i} \varepsilon(\pi_i \times \pi_j)^{{\rm Min}(d_i,d_j)}$. 
\end{defprop}

Pr\'ecisons que le produit ci-dessus est pris sur tous les $j=1,\dots,k$ distincts de $i$. D'apr\`es le r\'esultat$^{\color{green}\ast}$ d'Arthur sus-cit\'e, on peut m\^eme se restreindre aux entiers $j=1,\dots,k$ tels que ${\rm w}(\pi_j) \not \equiv {\rm w}(\pi_i) \bmod 2$. Pour justifier directement l'existence de $\varepsilon_\psi$, le lecteur observera que les $s_i, i \in I,$ sont lin\'eairement ind\'ependants sur $\F_2$ dans ${\rm C}_\psi$, et que s'ils engendrent un sous-groupe rencontrant non trivialement ${\rm Z}_{\widehat{G}}$,  alors $|I|=k$ et ${\rm Z}_{\widehat{G}}$ est engendr\'e par $\prod_{i=1}^k s_i$.\ps\ps

\subsection{Le cas des groupes de Chevalley}\label{parpaqarch}

On maintient les notations et hypoth\`eses des \S \ref{amfhypgen} et \S
\ref{parepsilonpsi}.  Arthur consid\`ere un morphisme de groupes
$$\nu_\infty : {\rm SL}_2(\C) \times {\rm W}_\R \rightarrow
\widehat{G}(\C)$$ d\'efini de la mani\`ere suivante.  \ps\ps

Pour tout $i=1,\dots,k$, un argument similaire \`a celui donn\'e dans la
d\'emonstration de la proposition~\ref{altso}, bas\'e sur la
proposition~\ref{contraintesalg}, assure qu'il existe un morphisme de
groupes $\mu_i : {\rm W}_\R \rightarrow \widehat{G^{\pi_i}}(\C)$ dont la
compos\'ee avec la repr\'esentation standard de $\widehat{G^{\pi_i}}$ est
isomorphe \`a ${\rm L}((\pi_i)_\infty)$.  Cette propri\'et\'e d\'etermine
uniquement $\mu_i$ modulo composition \`a l'arriv\'ee par ${\rm
Aut}(\widehat{G^{\pi_i}})$, mais il sera commode de fixer arbitrairement un
tel $\mu_i$.  Le morphisme $\nu_\infty$ est par d\'efinition la compos\'ee
du morphisme diagonal $(g,w) \mapsto (g,\prod_{i=1}^k \mu_i(w))$ par le
morphisme $\nu$.  La ${\rm Aut}(\widehat{G})$-orbite de $\nu_\infty$ dans
l'ensemble ${\rm Hom}({\rm SL}_2(\C) \times {\rm W}_\R, \widehat{G}(\C))$
sera not\'ee $\psi_\R$; il ne faudra pas la confondre avec la classe de
conjugaison semi-simple $\psi_\infty$, qui renferme une information nettement
moins fine.  L'orbite $\psi_\R$ ne d\'epend que de $\psi$ (et non du choix
de $\nu$ ou des $\mu_i$), et m\^eme mieux, que de l'ensemble des $k$ couples
$((\pi_i)_\infty,d_i)$ pour $i=1,\dots k$.  \ps\ps

Soit ${\rm C}_{\nu_\infty}$ le centralisateur de l'image de $\nu_\infty$ dans $\widehat{G}$. Il est \'evident que 
$${\rm C}_\nu \subset {\rm C}_{\nu_\infty}.$$
Il est ais\'e de d\'ecrire ${\rm C}_{\nu_\infty}$ \`a la mani\`ere de la description pr\'ec\'edente de ${\rm C}_\nu$. \ps\ps

\begin{lemme}\label{lemmepsir} La repr\'esentation ${\rm St} \circ \nu_\infty$ de ${\rm SL}_2(\C) \times {\rm W}_\R$ est semi-simple, sans multiplicit\'es, et tous ses constituants irr\'eductibles sont autoduaux. En particulier, ${\rm C}_{\nu_{\infty}}$ est un $2$-groupe ab\'elien fini \'el\'ementaire.
\end{lemme}

\begin{pf}  C'est \cite[Lemma 3.15]{chrenard2}. La seconde assertion d\'ecoule de la premi\`ere. Le seul point de cette derni\`ere qui ne d\'ecoule pas imm\'ediatement des d\'efinitions et de {\rm (}H2{\rm)} est l'assertion de multiplicit\'e $1$, qui est non triviale quand $\psi_\infty$ admet $0$ pour valeur propre double (et donc ${\rm w}_{\widehat{G}}=0$). Dans ce cas, ${\rm St} \circ \nu_\infty$ pourrait alors contenir soit $1$, soit $\epsilon_{\C/\R}$ avec multiplicit\'e $2$ (avec action triviale du facteur ${\rm SL}_2(\C)$ dans les deux cas). D'apr\`es {\rm (}H1{\rm )}, le corollaire \ref{cormodulo4} (iii) et la proposition \ref{contraintesalg} (iii), ces caract\`eres ne peuvent toutefois intervenir avec multiplicit\'e $2$ dans un m\^eme ${\rm L}((\pi_i)_\infty)$. Ils apparaissent donc (n\'ecessairement avec multiplicit\'e $1$) dans ${\rm L}((\pi_i)_\infty)$ et ${\rm L}((\pi_j)_{\infty})$ avec $i\neq j$, avec de plus $d_i=d_j=1$ et $n_i$ et $n_j$ sont impairs et non congrus modulo $4$ d'apr\`es {\it loc. cit}. Mais comme ${\rm w}(\pi_i)={\rm w}(\pi_j)=0$, cela contredit la proposition~\ref{contraintesalg} (ii).  
\end{pf} 

Afin de poursuivre l'analyse de la formule d'Arthur, nous supposerons dans un
premier temps que $G={\rm Sp}_{2g}$ ou ${\rm SO}_{r+1,r}$.  Notons
$\Pi_{\rm unit}(H)$ l'ensemble des classes d'isomorphisme de
repr\'esentations unitaires irr\'eductibles du groupe de Lie r\'eel $H$. 
Arthur~\cite[Thm.  1.5.1]{arthur} associe \`a $\nu_\infty$ un ensemble
fini\footnote{Pr\'ecisons que ce que nous notons ici $\psi_\R$, ${\rm
C}_{\nu_\infty}$ et $\Pi(\psi_\R)$ est not\'e respectivement $\psi$,
$\mathcal{S}_\psi$ et $\widetilde{\Pi}_{\psi}$ dans l'\'enonc\'e d'Arthur;
de plus il ne donne pas de nom \`a $\iota$ et si $u \in \Pi(\nu_\infty)$ il
note $x \mapsto \langle x, u \rangle$ le caract\`ere $\chi_u$.  L'image de
$\iota$, un sous-ensemble fini de $\Pi_{\rm unit}(G(\R))$, est commun\'ement
appel\'e {\it paquet d'Arthur associ\'e \`a $\psi_\R$}.  } $\Pi(\nu_\infty)$
muni de deux applications $$\Pi_{\rm unit}(G(\R))
\overset{\iota}{\longleftarrow} \Pi(\nu_\infty) \overset{u \mapsto
\chi_u}{\longrightarrow} {\rm Hom}_{\rm groupes}({\rm
C}_{\nu_\infty},\C^\times).$$ L'ensemble $\Pi(\nu_\infty)$ et $\iota$, ainsi
en fait que $\chi$ une fois clarifi\'ee la d\'ependence de ${\rm
C}_{\nu_\infty}$ en le choix de $\nu_\infty$, ne d\'ependent que de la
classe de ${\rm Aut}(\widehat{G})$-conjugaison $\psi_\R$ de $\nu_\infty$,
c'est pourquoi nous noterons aussi $\Pi(\psi_\R)$ pour $\Pi(\nu_\infty)$. 
Arthur d\'emontre une propri\'et\'e qui caract\'erise enti\`erement le
triplet $(\Pi(\psi_\R),\iota,\chi)$ \cite[Thm.  2.2.1]{arthur}, sans
toutefois expliciter concr\`etement ce triplet. Nous reviendrons sur ce point au
paragraphe suivant.  Comme l'explique Arthur \cite[p.  42]{arthur}, on
s'attend \`a ce que $\iota$ soit injective, de sorte que $\Pi(\psi_\R)$
serait d\'efini comme une partie de $\Pi_{\rm unit}(G(\R))$ et $\iota$
serait simplement ignor\'ee.\ps\ps

\begin{thmv}\label{amfchev}{\rm  (Formule de multiplicit\'e d'Arthur
\cite[Thm.  1.5.2]{arthur})} Supposons que $G={\rm Sp}_{2g}$ ou ${\rm
SO}_{r+1,r}$.  Fixons $\psi \in \mathcal{X}_{\rm AL}({\rm SL}_{{\rm
n}_{\widehat{G}}})$ satisfaisant {\rm(}H2{\rm )} ainsi que $U \in \Pi_{\rm
unit}(G(\R))$.  \ps\ps

Soit $\pi \in \Pi(G)$ l'unique repr\'esentation telle que $\pi_\infty \simeq
U$ et ${\rm St}({\rm c}_p(\pi))=\psi_p$ pour tout premier $p$.  Alors $\pi
\in \Pi_{\rm disc}(G)$ si, et seulement si, il existe $u \in \Pi(\psi_\R)$
tel que $$U=\iota(u) \, \, \, {\rm et} \, \, \, {\chi_u}_{|{\rm C}_\nu} =
\varepsilon_\psi.$$ Plus pr\'ecis\'ement, la multiplicit\'e ${\rm m}(\pi)$
de $\pi$ dans ${\mathcal{A}}_{\rm disc}(G)$ {\rm (\S IV.\ref{fautdiscgen})}
est exactement le nombre des \'el\'ements $u \in \Pi(\psi_\R)$ ayant la
propri\'et\'e ci-dessus.  \end{thmv}

\ps
D\'ecryptons cet \'enonc\'e.  Observons d'abord que l'existence et
l'unicit\'e de $\pi$ viennent de ce que si $G={\rm Sp}_{2g}$ ou ${\rm SO}_{r+1,r}$, l'application ${\rm St} :
\widehat{G}(\C)_{\rm ss} \rightarrow {\rm SL}_{{\rm n}_{\widehat{G}}}(\C)_{\rm ss}$ est injective
et d'image l'ensemble des classes \'egales \`a leur inverse.  En
particulier, la repr\'esentation $\pi$ de l'\'enonc\'e est l'unique
candidate possible telle que $\pi_\infty \simeq U$ et $\psi(\pi,{\rm
St})=\psi$.  Le th\'eor\`eme affirme d'abord que si $U$ n'est pas dans
$\iota(\Pi(\psi_\R))$ alors ${\rm m}(\pi)=0$.  Supposons donc que $U \in
\iota(\Pi(\psi_\R))$, et supposons \'egalement pour simplifier que l'on
sache que $\iota^{-1}(U)$ est un singleton $\{u\}$.  Le th\'eor\`eme affirme
alors que ${\rm m}(\pi) \neq 0$ si et seulement si ${\chi_u}_{|{\rm C}_\psi}
= \varepsilon_\psi$, auquel cas $m(\pi)=1$.  \ps\ps

Pour des raisons d'exhaustivit\'e, d\'ecrivons maintenant le cas restant du groupe $G={\rm SO}_{r,r}$, o\`u $r \equiv 0 \mod 2$ par l'hypoth\`ese {\rm (}H1{\rm )}. Dans ce cas, $\widehat{G}(\C)={\rm SO}_{2r}(\C)$. L'image de l'application ${\rm St} : \widehat{G}(\C)_{\rm ss} \rightarrow {\rm SL}_{2r}(\C)_{\rm ss}$ est encore l'ensemble des classes \'egales \`a leur inverse, mais elle n'est plus injective : deux \'el\'ements semi-simples de ${\rm SO}_{2r}(\C)$ conjugu\'es dans ${\rm SL}_{2r}(\C)$ ne sont pas toujours conjugu\'es dans ${\rm SO}_{2r}(\C)$, mais seulement dans ${\rm O}_{2r}(\C)$. Les fibres non vides de l'application ci-dessus sont donc exactement les orbites de l'action naturelle de $${\rm O}_{2r}(\C)/{\rm SO}_{2r}(\C)  = {\rm Out}(\widehat{G}) = \Z/2\Z$$ sur  $\widehat{G}(\C)_{\rm ss}$. Un ph\'enom\`ene analogue se produit pour ${\rm St} : \widehat{\mathfrak{g}}_{\rm ss} \rightarrow (\mathfrak{sl}_{2r})_{\rm ss}$. L'action de ${\rm O}_{r,r}(\R)$ par conjugaison sur $G(\R)$ d\'efinit \'egalement une action de $\Z/2\Z={\rm O}_{r,r}(\R)/{\rm SO}_{r,r}(\R)$ sur $\Pi_{\rm unit}(G(\R))$, dont l'ensemble des orbites est not\'e par Arthur $\widetilde{\Pi}_{\rm unit}(G(\R))$. Le th\'eor\`eme \cite[Thm. 1.5.1]{arthur} d'Arthur associe alors \`a $\nu_\infty$ un triplet $(\Pi(\psi_\R),\iota,\chi)$ comme plus haut, \`a la seule diff\'erence pr\`es que $\iota$ est cette-fois ci une application $\Pi(\psi_\R) \rightarrow \widetilde{\Pi}_{\rm unit}(G(\R))$. \ps

\begin{thmv}\label{amfchev2}{\rm  (Formule de multiplicit\'e d'Arthur pour ${\rm SO}_{r,r}$ \cite[Thm. 1.5.2]{arthur})} Supposons que $G={\rm SO}_{r,r}$ avec $r \equiv 0 \bmod 2$.  Fixons $\psi \in \mathcal{X}_{\rm AL}({\rm SL}_{2r})$ satisfaisant {\rm(}H2{\rm )} ainsi que $U \in \widetilde{\Pi}_{\rm unit}(G(\R))$. Soit $E \subset \Pi(G)$ l'ensemble des $\pi \in \Pi(G)$ telles que $\psi(\pi,{\rm St})=\psi$ et  $\pi_\infty \in U$. Soit $F \subset \Pi(\psi_\R)$ l'ensemble des \'el\'ements $u$ tels que $U \in \iota(u)$ et ${\chi_u}_{|{\rm C}_\nu} = \varepsilon_\psi$. Alors
$$\sum_{\pi \in E} {\rm m}(\pi)=m_\psi |F|,$$
o\`u $m_\psi=1$ \`a moins que $\psi=\oplus_{i=1}^k \pi_i [d_i]$ avec $d_i \equiv 0 \bmod 2$ pour tout $i$, auquel cas $m_\psi=2$. 

\end{thmv}

(Arthur s'attend encore \`a l'injectivit\'e de $\iota$, et donc \`a ce que $|F| \in \{0,1\}$). En vue d'appliquer les th\'eor\`emes \ref{amfchev} ou \ref{amfchev2}, il est \'evidemment crucial d'en savoir plus sur le triplet $(\Pi(\psi_\R),\iota,\chi)$. L'objectif des paragraphes qui suivent est de rappeler des choses connues \`a ce sujet. \ps

\section{S\'eries discr\`etes}\label{paramdisc} 

\subsection{S\'eries discr\`etes, d'apr\`es Harish-Chandra.}
\label{paramdiscA} Les r\'esultats de ce paragraphe sont d\^us \`a
Harish-Chandra.  Il sera commode de renvoyer \`a l'ouvrage de Knapp
\cite{knapp} pour les d\'emonstrations.  \ps\ps

Soit $H$ un $\R$-groupe semi-simple.  On rappelle que $\pi
\in \Pi_{\rm unit}(H(\R))$ est une {\it s\'erie discr\`ete} si elle
intervient comme sous-repr\'esentation de la repr\'esentation r\'eguli\`ere
${\rm L}^2(H(\R))$ \cite[Ch. IX \S 3]{knapp}.  Harish-Chandra a d\'emontr\'e que $H(\R)$ admet des
s\'eries discr\`etes si, et seulement si, $H$ poss\`ede un {\it tore maximal anisotrope}, c'est-\`a-dire un 
$\R$-sous-groupe $T \subset H$ tel que $T_\C \subset H_\C$ est un tore maximal et $T(\R)$ est
compact \cite[Thm. 12.20]{knapp}.  C'est le cas de ${\rm Sp}_{2g}$, ainsi que du groupe sp\'ecial
orthogonal r\'eel de signature $(p,q)$ si et seulement si $pq \equiv 0 \bmod
2$ ou $p+q \equiv 1 \bmod 2$, de sorte que la condition (H1) du \S \ref{amfhypgen} est en fait
\'equivalente \`a demander que $G(\R)$ admet des s\'eries discr\`etes.  \ps\ps

On suppose d\'esormais que $H(\R)$ admet des s\'eries discr\`etes. Les
tores maximaux anisotropes de $H$ forment une seule
orbite sous l'action de conjugaison de $H(\R)$ ; fixons l'un d'entre eux que nous notons $T$.  Introduisons
quelques objets associ\'es au couple $(H,T)$.  Les caract\`eres et cocaract\`eres de $T_\C$ \'etant tous d\'efinis sur $\R$, on notera
simplement ${\rm X}_\ast(T)$ (resp. ${\rm X}^\ast(T)$) pour ${\rm X}_\ast(T_\C)$ (resp. ${\rm X}^\ast(T_\C)$). Soient $$\Phi=\Phi(H_\C,T_\C) 
\subset {\rm X}^\ast(T)\, \, \, {\rm et}\, \, \,
W={\rm W}(H_\C,T_\C)$$ le groupe de Weyl de $\Phi$. Posons \'egalement
$$W_{\rm r} = {\rm W}(H,T) \overset{{\mathrm{d\acute{e}f}}.}{=}{\rm
N}_{H(\R)}(T(\R))/T(\R),$$ o\`u ${\rm N}_{H(\R)}(T(\R))$ est le
normalisateur de $T(\R)$ dans $H(\R)$ ({\it groupe de Weyl r\'eel} de
$(H,T)$).  On dispose enfin d'un unique sous-groupe compact maximal $K$ de
$H(\R)$ contenant $T(\R)$ (et donc ${\rm N}_{H(\R)}(T(\R))$) et l'on note
$\Phi_{\rm c} \subset \Phi$ le syst\`eme de racines (dites {\it compactes})
de $K$ relativement \`a $T(\R)$.  On dispose d'inclusions naturelles ${\rm
W}(\Phi_{\rm c}) \subset W_{\rm r} \subset W$, qui sont strictes en
g\'en\'eral.  \ps\ps 

Si $V \in {\rm Irr}(H_\C)$, notons $\Pi_V \subset \Pi_{\rm unit}(H(\R))$ l'ensemble des s\'eries discr\`etes ayant m\^eme caract\`ere infinit\'esimal que $V$. Harish-Chandra a d\'emontr\'e que toute s\'erie discr\`ete de $H(\R)$ appartient \`a $\Pi_V$, pour un unique $V$. Si $\Delta$ est une base du syst\`eme de racines $\Phi$, il a d\'efini une repr\'esentation $\pi_{\Delta,V} \in \Pi_V$, uniquement d\'etermin\'ee par les valeurs prises par son caract\`ere $\Theta_{\Delta,V}$ sur l'ensemble $T(\R)_{\rm reg} \subset T(\R)$ des \'el\'ements $t \in T(\R)$ tels que $\alpha(t) \neq 1$ pour tout $\alpha \in \Phi$. Pr\'ecis\'ement, en notant $t^\mu$ pour $\mu(t)$ si $t \in T(\R)$ et $\mu \in {\rm X}^\ast(T)$,
$$\Theta_{\Delta,V}(t) =  (-1)^{\frac{1}{2}\dim H(\R)/K} \frac{ \sum_{w \in W_{\rm r}}  \varepsilon(w) t^{w(\lambda+\rho)-\rho}}{\prod_{\alpha \in \Phi^+}(1-t^{-\alpha})},\,\, \, \, \, \forall t \in T(\R)_{\rm reg},$$
o\`u $\lambda \in {\rm X}^\ast(T)$ d\'esigne le plus haut poids de $V$ relativement au sous-groupe de Borel $T_\C \subset B \subset H_\C$ d\'efini par $\Delta$, $\Phi^+ \subset \Phi$ est le syst\`eme positif d\'efini par $\Delta$, et $2\rho$ la somme des \'el\'ements de $\Phi^+$ \cite[Thm. 9.20 \& 12.7]{knapp}. Il d\'emontre que tout \'el\'ement de $\Pi_V$ est de la forme $\pi_{\Delta,V}$, et que $\pi_{\Delta,V} \simeq \pi_{\Delta',V}$ si, et seulement si, $\Delta'=w(\Delta)$ avec $w \in W_{\rm r}$ \cite[Thm. 12.21]{knapp}. En particulier, $|\Pi_V|=|W_{\rm r}\backslash W|$ est ind\'ependant de $V$.\ps\ps 

Supposons qu'il existe un sous-groupe $S \subset K$ isomorphe \`a
$\mathbb{S}^1$ et dont l'action adjointe sur ${\rm Lie} \, H(\R)$ admet
exactement ${\rm Lie}(K)$ pour points fixes (en particulier, $S$ est dans le
centre de la composante neutre de $K$).  Cela se produit par exemple quand
$H$ est ${\rm Sp}_{2g}$ ou le groupe sp\'ecial orthogonal de signature
$(p,2)$ avec $p\geq 1$.  Dans ce cas, $H(\R)$ poss\`ede des s\'eries
discr\`etes {\it holomorphes} \cite[Chap.  VI]{knapp}, et nous en avons
d\'ej\`a rencontr\'e des exemples au \S VI.\ref{carinf}.  Dans la
classification d'Harish-Chandra ce sont exactement les repr\'esentations
$\pi_{\Delta,V}$ obtenues lorsque $\Delta$ est la base d'un syst\`eme
positif de la forme $\{\alpha \in \Phi, \varphi(\alpha)>0\}$ o\`u $\varphi :
{\rm X}^\ast(T) \rightarrow \R$ est une forme lin\'eaire telle que $$\forall \, \alpha \in \Phi_c, \, \, \forall \, \beta \in \Phi-\Phi_c, \, \, \, \, \,
0<|\varphi(\alpha)|< |\varphi(\beta)|$$ (les racines compactes ``sont plus
petites'' que les racines non compactes), d'apr\`es la remarque suivant \cite[Thm.  9.20]{knapp}.\ps

\subsection{Param\'etrisation canonique de Shelstad, cas des groupes d\'eploy\'es}\label{paramdiscB} Il sera n\'ecessaire de rappeler une seconde param\'etrisation des \'el\'ements de $\Pi_V$, $V \in {\rm
Irr}(H_\C)$ \'etant fix\'e, qui comme nous l'avons d\'ej\`a dit
au~\S\ref{amfhypgen} est due \`a Shelstad, intervenant dans l'\'enonc\'e de la formule de multiplicit\'e d'Arthur. Elle doit son origine aux identit\'es existantes entre les caract\`eres des s\'eries discr\`etes de $H$ et ceux des s\'eries discr\`etes d'une collection de $\R$-groupes associ\'es, appel\'es {\it endoscopiques} par Langlands \cite{shelstadlind}. Une exposition d\'etaill\'ee de ces identit\'es d\'epasserait tr\`es nettement le cadre de ce m\'emoire (et la comp\'etence des auteurs), et nous ne nous y aventurerons pas. Nous suivrons l'expos\'e synth\'etique de Shelstad \cite{shelstadnotesbanff2011}, qui d\'ecrit notamment les normalisations pr\'ecises utilis\'ees par Arthur, tout en empruntant parfois le point de vue \'eclairant d'Adams \cite{adams}.  \ps\ps

Faisons d\'esormais l'hypoth\`ese suppl\'ementaire que $H$ est d\'eploy\'e sur $\R$ (semi-simple, et tel que $H(\R)$ poss\`ede des s\'eries discr\`etes).  
Soit $T$ un tore maximal anisotrope de $H$ ; on reprend les notations du \S\ref{paramdiscA} pour les objets associ\'es. La
param\'etrisation de Shelstad n'est compl\`etement canonique que si $H$ est
adjoint. En g\'en\'eral, elle d\'ependra du choix d'une $W_{\rm r}$-orbite de
l'ensemble $\mathcal{B}(T)$ des bases $\Delta$ du syst\`eme de racines $\Phi$
telles que $\Delta \cap \Phi_{\rm c} = \emptyset$. Cet ensemble est non
vide et muni d'une action naturelle, simplement transitive, du groupe de
Weyl r\'eel $W_{\rm r}^{\rm ad}$ de $(H/{\rm
Z}(H))(\R)$ relativement \`a l'image de $T(\R)$, qui v\'erifie $W_{\rm r} \subset W_{\rm r}^{\rm ad} \subset W$.  Le ${\rm W}(H,T)$-ensemble $\mathcal{B}(T)$ d\'epend \'evidemment du choix de $T$, mais le fait que l'ensemble de ces $T$ ne forme qu'une classe de $H(\R)$-conjugaison assure que l'ensemble quotient ${\rm W}(H,T)\backslash \mathcal{B}(T)$ ne d\'epend que de $H$, il sera not\'e $\mathcal{B}(H)$. Des r\'esultats de Kostant
et Vogan montrent que le choix d'une $W_{\rm r}$-orbite $O \subset
\mathcal{B}(T)$, ce que l'on notera aussi $O \in \mathcal{B}(H)$, est \'equivalent \`a celui d'une {\it classe d'\'equivalence de donn\'ees de Whittaker} $D$ pour $H(\R)$, une notion que nous n'introduirons pas ici, mais qui est
pr\'ecis\'ement la donn\'ee de r\'ef\'erence chez Arthur : l'unique s\'erie
discr\`ete dans $\Pi_V$ qui est g\'en\'erique pour $D$ est $\pi_{\Delta,V}$
pour $\Delta \in O$. 
 \ps\ps

\begin{definition}\label{defparshel} Soient $T$ un tore maximal anisotrope de $H$, $\Delta \in \mathcal{B}(T)$ et $\rho^\vee$ la
demi-somme des co-racines de $T$ positives relativement \`a $\Delta$.  La
param\'etrisation de Shelstad de $\Pi_V$ relative \`a $(T,\Delta)$ est l'application $$\kappa_\Delta : \Pi_V
\longrightarrow {\rm X}_\ast(T) \otimes \Z/2\Z, \pi \mapsto \kappa_\Delta(\pi),$$
d\'efinie par $\kappa_\Delta(\pi_{w^{-1} \Delta,V}) \equiv w \rho^\vee -\rho^\vee
\bmod 2{\rm X}_\ast(T)$ pour tout $w \in W$.  \end{definition}

Soulignons que $\rho^\vee \in \frac{1}{2}{\rm X}_\ast(T)$ n'est pas dans
${\rm X}_\ast(T)$ en g\'en\'eral.  En revanche, le terme $w \rho^\vee
-\rho^\vee$ est bien dans ${\rm X}_\ast(T)$ : il appartient m\^eme au
sous-groupe engendr\'e par les co-racines. Nous verrons plus loin que $\pi \mapsto \kappa_\Delta(\pi)$
est injective. \ps\ps 

Expliquons pourquoi cette d\'efinition co\"incide avec celle donn\'ee par
Shelstad \cite[\S 8]{shelstadnotesbanff2011}.  Suivant Langlands et
Shelstad, consid\'erons le premier groupe de cohomologie $H^1(\R,T)$ de
$T(\C)$ vu comme $\Z[{\rm Gal}(\C/\R)]$-module.  Comme on a ${\rm
Gal}(\C/\R)=\Z/2\Z=\langle 1, \sigma\rangle$, c'est simplement le quotient
du groupe ab\'elien $$Z^1(\R,T)=\{t \in T(\C), t\sigma(t)=1\}$$ par le
sous-groupe des \'el\'ements de la forme $t\sigma(t)^{-1}$ avec $t \in
T(\C)$.  Rappelons d'abord que $H^1(\R,T)$ s'identifie canoniquement avec le
but de l'application $\kappa_\Delta$ (dualit\'e de Tate-Nakayama).  En effet,
soit $T_2=\{t \in T(\R), t^2=1\}=T(\R) \cap Z^1(\R,T)$.  La compos\'ee des
applications naturelles \begin{equation}\label{isomgamma} T_2
\hookrightarrow Z^1(\R,T) \rightarrow H^1(\R,T)\end{equation} est
manifestement un isomorphisme, le tore $T$ \'etant $\R$-isomorphe \`a un
produit fini de copies de $\mathbb{S}^1$.  Si $\mu \in {\rm X}_\ast(T)
\otimes \C$, notons $e^\mu$ l'unique \'el\'ement $z \in T(\C)$ tel que
$\lambda(z)=e^{ \langle \lambda, \mu \rangle}$ pour tout $\lambda \in {\rm
X}^\ast(T)$.  L'application $\mu \mapsto e^{i\pi\mu}$ induit donc un
isomorphisme naturel ${\rm X}_\ast(T) \otimes \Z/2\Z \isomo T_2$.  En
mettant bout-\`a-bout les isomorphismes ci-dessus, on obtient l'isomorphisme
naturel annonc\'e $H^1(\R,T) \isomo {\rm X}_\ast(T) \otimes \Z/2\Z$.  \ps\ps 
 
On dispose \'egalement d'une action naturelle de $W$ sur $H^1(\R,T)$ induite par l'inclusion $T \subset H$,
d\'efinie de la mani\`ere suivante \cite{borovoi14}.  Si $N$ est le normalisateur de
$T(\C)$ dans $H(\C)$, alors $n\sigma(n)^{-1} \in T(\C)$ pour tout $n \in N$,
et donc $(n,t) \mapsto n \star t := n t \sigma(n)^{-1}$ d\'efinit une action de $N$ sur
$Z^1(\R,T)$. Elle induit par passage au quotient une action de $W=N/T(\C)$
sur $H^1(\R,T)$, que nous noterons encore $(w,x) \mapsto w \star x$.  Ceci \'etant
dit, si $w \in W$ alors l'\'el\'ement $\kappa_\Delta(\pi_{w^{-1}\Delta,V}) \in {\rm
X}_\ast(T) \otimes \Z/2\Z$ consid\'er\'e par Shelstad \cite[p. 
15]{shelstadnotesbanff2011} est par d\'efinition l'image de $w \star 1 \in H^1(\R,T)$, $1 \in H^1(\R,T)$ d\'esignant l'\'el\'ement
neutre, par l'isomorphisme
$H^1(\R,T) \isomo {\rm X}_\ast(T) \otimes \Z/2\Z$ rappel\'e au paragraphe
pr\'ec\'edent.  \ps\ps 

Prenons garde que l'isomorphisme $\gamma : T_2 \rightarrow H^1(\R,T)$,
d\'efini par la formule \eqref{isomgamma}, n'entrelace pas l'action
\'evidente par conjugaison de $W$ sur $T_2$ et celle sur $H^1(\R,T)$
d\'efinie ci-dessus. L'identit\'e $g t \sigma(g)^{-1} = g t g^{-1} g
\sigma(g)^{-1}$ montre toutefois que c'est vrai pour l'action du sous-groupe
$W_{\rm r}$. Pour s'\'eviter des confusions, nous noterons $(w,t) \mapsto w(t)$ l'action usuelle de $W$
sur $T(\C)$ et $(w,t) \mapsto w \cdot t = \gamma^{-1}(w \star \gamma(t))$
l'action ``tordue'' sur $T_2$. La relation pr\'ecise entre les deux est donn\'ee par le lemme suivant (ii),
essentiellement d\^u \`a Langlands \cite[\S 3]{shelstadannens}\cite[Lemma
5.1]{kottwitzannarbor}\cite[Lemma 7.8]{adams}.  Fixons $\Delta \in \mathcal{B}(T)$ et notons $\rho^\vee \in \frac{1}{2}{\rm X}_\ast(T)$ la demi-somme des coracines de $(H_\C,T_\C)$ positives relativement \`a $\Delta$. Suivant \cite{adams}, posons $$t_{\rm b}=e^{i\pi
\rho^\vee} \in T(\C).$$ On a $t_{\rm b}^4=1$ donc $t_{\rm b} \in T(\R)$.\ps\ps

\begin{lemme}\label{lemmelankott}\begin{itemize} \item[(i)] Le centralisateur de $t_{\rm b}$ dans $H(\R)$ est $K$ et $t_{\rm b}^2 \in {\rm Z}(H)$. \ps\ps 
\item[(ii)] Pour tout $t \in T_2$ et $w \in W$ on a $w \cdot t =  w(t t_{\rm b}) t_{\rm b}^{-1}$.
\end{itemize}
\end{lemme}
\ps
\noindent
La relation $w \cdot 1 = e^{i\pi (w \rho^\vee- \rho^\vee)}$ s'en d\'eduit et conclut. Observons que le (i) entra\^ine que l'\'el\'ement $t_b$ ne d\'epend que de la $W_{\rm r}$-orbite de $\Delta$ dans $\mathcal{B}(T)$. 
\ps
\begin{pf} Si $\alpha \in \Phi$ alors $\alpha(t_{\rm b}) = (-1)^{\langle
\rho^\vee,\alpha \rangle}$.  Soit $s : \Phi \rightarrow \Z/2\Z$
l'application telle que $s^{-1}(0)=\Phi_c$.  La d\'ecomposition de Cartan de
${\rm Lie} \, H(\R)$ relativement \`a $K$ montre que
$s(\alpha+\beta)=s(\alpha)+s(\beta)$ chaque fois que $\alpha,\beta$ et
$\alpha+\beta$ sont dans $\Phi$.  Par r\'ecurrence sur la ``hauteur''
$|\langle \rho^\vee, \alpha \rangle|$, on a  $s(\alpha) \equiv \langle
\rho^\vee, \alpha \rangle \bmod 2$ pour tout $\alpha \in \Phi$.  Cela montre
que l'automorphisme int\'erieur de $H$ d\'efini par $t_{\rm b}$ est
l'involution de Cartan de $H(\R)$ relativement \`a $K$, puis le (i).  \ps\ps 
Comme on a $t_{\rm b}^2 \in {\rm Z}(H)$, la fonction $f (w) = w(t_{\rm b})t_{\rm
b}^{-1}$ d\'efinit un $1$-cocyle de $W$ \`a valeurs dans $T_2$.  De m\^eme, 
$g(w)=w \cdot 1$ est aussi un $1$-cocyle \`a valeurs dans $T_2$.  Comme $w
\cdot t = w(t) \,w \cdot 1$, il suffit de
voir que $f$ et $g$ co\"incident sur les ${\rm s}_\alpha$ avec $\alpha \in \Delta$,
soit encore que ${\rm s}_\alpha \cdot 1 = f(s_\alpha) = e^{i\pi
\alpha^\vee}$. C'est exactement
le calcul de Langlands~\cite[Prop.  2.1]{shelstadannens}.  Cela d\'emontre
le (ii).  \end{pf}

\noindent Observons enfin que si $\Delta \in \mathcal{B}(T)$, $w \in W_{\rm r}$ et $\pi \in \Pi_V$, on a 
\begin{equation}\label{chgtbasekappa} \kappa_{w\Delta}(\pi) = w \kappa_\Delta(\pi),\end{equation}
\`a cause de l'identit\'e imm\'ediate $(ww') \cdot 1=w \cdot (w' \cdot 1)$ pour tout $w' \in W$. En particulier, la $W_{\rm r}$-orbite de $\kappa_\Delta(\pi)$ ne d\'epend que de celle de $\Delta$ dans $\mathcal{B}(T)$.\ps\ps

\subsection{Interpretation duale et lien avec les paquets
d'Arthur}\label{paramdiscC}On conserve les notations et hypoth\`eses des
\S\ref{paramdiscA} et \S\ref{paramdiscB}.  Avant de donner un exemple,
donnons l'interpr\'etation duale utile de la param\'etrisation de Shelstad. 
Soit $\widehat{H}$ le $\C$-groupe dual de $H_\C$ (rappelons que $H$ est
d\'eploy\'e sur $\R$).  Suivant Langlands~\cite{langlandsreel}, il existe
une bijection naturelle entre ${\rm Irr}(H_\C)$ et l'ensemble des
$\widehat{H}$-classes de conjugaison de {\it param\`etres discrets} $\varphi
: {\rm W}_\R \rightarrow \widehat{H}$.  Rappelons qu'un tel objet est par
d\'efinition un morphisme de groupes continu tel que le sous-groupe
$\varphi({\rm W}_\R) \subset \widehat{H}$ est constitu\'e d'\'el\'ements
semi-simples et poss\`ede un centralisateur fini dans $\widehat{H}$, not\'e
${\rm C}_\varphi$.  Cette bijection est caract\'eris\'ee par le fait que le
caract\`ere infinit\'esimal d'un tel $\varphi$ co\"incide avec celui de la
repr\'esentation dans ${\rm Irr}(H_\C)$ qui lui correspond.\ps\ps 

Expliquons cette assertion. Le sous-groupe $\varphi(\C^\times)
\subset \widehat{H}$ est commutatif, connexe, et constitu\'e d'\'el\'ements
semi-simples.  Son adh\'erence Zariski dans $\widehat{H}$ est donc un tore,
et son centralisateur dans $\widehat{H}$, not\'e $S$, un sous-groupe de Levi
(d'un sous-groupe parabolique) de $\widehat{H}$.  Par d\'efinition, ${\rm C}_\varphi \subset S$ est le
sous-groupe fix\'e par la conjugaison par $\varphi(j)$.  Comme il est fini,
cela force $S$ \`a \^etre un tore maximal et $\varphi(j)$ \`a agir par
l'inversion $s \mapsto s^{-1}$ (en particulier, $-1$ est un \'el\'ement du
groupe de Weyl de $\widehat{H}$, ce qui est bien le cas sous notre
hypoth\`ese sur $H$).  Soient $\lambda_\phi,\mu_\phi \in {\rm X}_\ast(S)
\otimes \C$ les uniques \'el\'ements tels que $\lambda_\phi-\mu_\phi \in
{\rm X}_\ast(S)$ et $\xi(\varphi(z))=z^{\langle \xi,\lambda_\phi \rangle}
\overline{z}^{\langle \xi, \mu_\phi\rangle}$ pour tout $z \in \C^\times$ et
$\xi \in {\rm X}^\ast(S)$ (voir la note de bas de page \S\ref{apparitionWR})
Le caract\`ere infinit\'esimal de $\varphi$ est par d\'efinition la classe
de $\widehat{H}$-conjugaison de $\lambda_\phi$, vue dans ${\rm Lie}
\,\widehat{H}$ (formule VI.\eqref{poidscarinf}).  De plus,
$\mu_\phi=-\lambda_\phi$.\ps \ps

Fixons maintenant un param\`etre discret $\varphi : {\rm W}_\R \rightarrow
\widehat{H}$.  Observons que si $\varphi'$ est dans la classe de
$\widehat{H}$-conjugaison de $\varphi$, et si $h \in \widehat{H}$ est tel
que $\varphi'=h \varphi h^{-1}$, l'isomorphisme ${\rm C}_\varphi \rightarrow
{\rm C}_{\varphi'}$ induit par conjugaison par $h$ est ind\'ependant du
choix de $h$, de sorte que ${\rm C}_\varphi$ est un groupe ab\'elien
canonique.  Soit $S$ le tore maximal de $\widehat{H}$ contenant
$\varphi(\C^\times)$ et $B \subset \widehat{H}$ l'unique sous-groupe de
Borel contenant $S$ tel que l'\'el\'ement $\lambda_\phi \in \frac{1}{2} {\rm
X}_\ast(S)$ d\'efini ci-dessus soit dominant par
rapport \`a $B$.  La donn\'ee de $\Delta \in \mathcal{B}(T)$, qui d\'efinit un unique
sous-groupe de Borel $T_\C \subset Q \subset H_\C$, et celle de $\widehat{H}$,
d\'efinissent un isomorphisme unique $\Psi(H_\C,T,Q)^\vee \simeq
\Psi(\widehat{H},S,B)$, et en particulier un isomorphisme privil\'egi\'e
$\widehat{T} \rightarrow S$, soit encore un isomorphisme $$i_\Delta : {\rm X}_\ast(T)
\isomo {\rm X}^\ast(S).$$ Ceci \'etant dit, ${\rm C}_\varphi=S_2=\{s \in S,
s^2=1\}$ et l'application naturelle $\beta : {\rm X}_\ast(S) \otimes \Z/2\Z \rightarrow {\rm Hom}(S_2,\C^\times)$ est un isomorphisme de groupes.
Au final, l'application $\kappa_\Delta$ de Shelstad (d\'efinition~\ref{defparshel}) induit une application naturelle
\begin{equation}\label{chiintdual} \chi_O : \Pi_V
\longrightarrow {\rm Hom}_{\rm groupes}({\rm C}_\varphi,\C^\times),
\end{equation}
d\'efinie par $\chi_O := \beta \circ (i_\Delta \otimes \Z/2\Z) \circ \kappa_\Delta$, o\`u $O \in \mathcal{B}(T)$ d\'esigne la $W_{\rm r}$-orbite de $\Delta$. On observe en effet que si $w \in W$ alors $i_{w\Delta} = i_\Delta \circ w^{-1}$, de sorte que l'application $\chi_O$ ne d\'epend bien que $O$, d'apr\`es la relation \eqref{chgtbasekappa}. Son image est constitu\'ee d'homomorphismes triviaux sur ${\rm Z}(\widehat{H})$ car $w \rho^\vee - \rho^\vee$ est une somme de racines de $(\widehat{H},S)$ pour tout $w \in W$. Elle est injective mais non surjective en g\'en\'eral.\ps\ps 

Le lien avec le \S\ref{parpaqarch} est que lorsque l'homomorphisme
$\nu_\infty$ d\'efini {\it loc.  cit.} est trivial sur le facteur ${\rm
SL}_2(\C)$, ce qui \'equivaut \`a demander $d_i=1$ pour tout $i=1,\dots,k$,
alors c'est un param\`etre discret ${\rm W}_\R \rightarrow \widehat{G}$
(Lemme~\ref{lemmepsir}) ayant m\^eme caract\`ere infinit\'esimal que la
repr\'esentation $V \in {\rm Irr}(G_\C)$ fix\'ee par l'hypoth\`ese (H2)
\S\ref{amfhypgen}.  Il sera commonde de traiter \`a part le cas de ${\rm
SO}_{r,r}$. \ps\ps

(a) Supposons que $G={\rm Sp}_{2g}$ ou ${\rm SO}_{r+1,r}$. Alors 
d'apr\`es~\cite{shelstadnotesbanff2011, shelstad1} et \cite{mezo2},
l'ensemble $\Pi(\psi_\R)$ consid\'er\'e par Arthur est $\Pi_V$,
l'application $\iota$ est l'inclusion \'evidente $\Pi_V \subset \Pi_{\rm
unit}(G(\R))$, et l'application $\chi$ {\it loc.  cit.} est l'application $\chi_O$
d\'efinie ci-dessus.  Quand $G={\rm SO}_{r+1,r}$
(adjoint) il n'y a qu'un choix possible de $O \in
\mathcal{B}(G_\R)$, de sorte
que tout est canoniquement d\'efini.  Lorsque $G={\rm Sp}_{2g}$, il y en a
exactement deux (voir ci-dessous), et il faut bien s\^ur s'en tenir au choix
fait par Arthur \cite[p.  55]{arthur}.  \'Etant donn\'e que l'application
naturelle ${\rm PGSp}_{2g}(\Z) \rightarrow \pi_0({\rm PGSp}_{2g}(\R))=\Z/2\Z$ est surjective, ce choix ne jouera en fait aucun r\^ole
dans nos applications.  \ps\ps 

(b) Supposons enfin $G={\rm SO}_{r,r}$ avec $r \equiv 0 \bmod 2$. Le groupe \`a deux \'el\'ements ${\rm O}_{r,r}(\R)/{\rm SO}_{r,r}(\R)= \{1,\theta\}$ agit naturellement sur $\Pi_{\rm unit}(G(\R))$ et ${\rm Irr}(G_\C)$. On observe que $\theta$ induit une bijection $\Pi_V \isomo \Pi_{\theta(V)}$. 

\begin{lemme}\label{lemmethetacar}Soient $H={\rm SO}_{r,r}$ avec $r \equiv 0\bmod 2$, $T$ un tore maximal anisotrope de $H$, $\Phi=\Phi(H_\C,T_\C)$, $\Delta$ une base de $\Phi$, et $V \in {\rm Irr}(H_\C)$. Alors $\theta(\pi_{\Delta,V}) \simeq \pi_{\Delta,\theta(V)}$. 
\end{lemme}

\begin{pf} Il est \'equivalent de se donner $T$ ou une d\'ecomposition de ${\rm H}(\R^r)$ comme somme orthogonale $\oplus_{i \in I} P_i$ de plans $P_i$ suppos\'es d\'efinis (postifs ou n\'egatifs). Ceci \'etant fait, il est \'equivalent de se donner $\Delta$, un sous-groupe de Borel de $H_\C$ contenant $T_\C$, ou encore un ordre total sur l'ensemble $I$ avec pour chaque $i \in I$ un choix de l'une des deux droites isotropes de $P_i \otimes \C$ (\S VI.\ref{exdorad}). Soit $i_0$ le plus grand \'el\'ement de $I$. Soit $s \in {\rm O}({\rm H}(\R^r))$ l'unique \'element qui agit par l'identit\'e sur $P_i$ pour $i<i_0$, et qui \'echange les deux droites isotropes de $P_{i_0}$. C'est un repr\'esentant de $\theta$ pr\'eservant $T$ ainsi que la base $\Delta \subset \Phi$. 
On conclut car la propri\'et\'e caract\'eristique de $\pi_{\Delta,V}$ montre que $\theta(\pi_{\Delta,V}) = \pi_{s(\Delta),\theta(V)}$ et on a vu que $s(\Delta)=\Delta$.
\end{pf}

En particulier, si $\theta(V) \simeq V$ on constate que $\theta$ agit par
l'identit\'e sur $\Pi_V$ : tout $\pi \in \Pi_V$ se prolonge \`a ${\rm
O}({\rm H}(\R^r))$.  Retournons au $\Z$-groupe $G={\rm SO}_{r,r}$. 
D'apr\`es les r\'esultats sus-cit\'es de Shelstad et Mezo, l'ensemble
$\Pi(\psi_\R)$ consid\'er\'e par Arthur est l'image $\widetilde{\Pi}_V$ de $\Pi_V$ dans
$\widetilde{\Pi}_{\rm unit}(G(\R))$, et l'application $\iota$ est
l'inclusion \'evidente.  Si $\pi \in \Pi_V$, le lemme \ref{lemmethetacar}
assure que $\pi$ et $\theta(\pi)$ ont m\^eme caract\`ere de Shelstad, ce qui
fournit par passage au quotient une application bien d\'efinie
$\widetilde{\Pi}_V \longrightarrow {\rm Hom}_{\rm
groupes}({\rm C}_{\nu_\infty},\C^\times) $ : c'est l'application
consid\'er\'ee par d\'efinition chez Arthur.  Pour \^etre tout-\`a-fait
exact, comme dans le cas $G={\rm PGSp}_{2g}$ l'ensemble $\mathcal{B}(G_\R)$
contient deux \'el\'ements, et il faut choisir celui correspondant \`a la donn\'ee de Whittaker fix\'ee par Arthur, mais
l\`a encore ce choix
ne jouera pas de r\^ole dans nos applications.

\subsection{Exemple : la s\'erie discr\`ete holomorphe de ${\rm
Sp}_{2g}(\R)$}\label{paramdiscD} Consid\'erons l'exemple du groupe $H={\rm
Sp}_{2g}$, de $\C$-groupe dual $\widehat{H}={\rm SO}_{2g+1}$.  Posons
$E={\rm H}(\R^g)$, muni de sa forme altern\'ee hyperbolique ${\rm a}$, de
sorte que $H={\rm Sp}_E$.  Rappelons que si $I \in {\rm Sp}(E)$ est tel que
$I^2=-{\rm id}_{E}$, il munit $E$ d'une structure complexe ainsi que d'une
forme hermitienne non d\'eg\'en\'er\'ee ${\rm h}$ pour cette structure, de
forme bilin\'eaire associ\'ee $$(e,f) \mapsto {\rm a}(I e,f)+i {\rm
a}(e,f).$$ Le centralisateur de $I$ dans ${\rm Sp}_E$ est alors le
$\R$-groupe unitaire ${\rm U}_{\rm h}$.  Choisissons l'\'el\'ement $I$ de
sorte que ${\rm h}$ soit d\'efinie positive, auquel cas ${\rm U}_{\rm
h}(\R)$ est un sous-groupe compact maximal de ${\rm Sp}_E$.  Il existe une
unique classe de conjugaison de tels \'el\'ements sous ${\rm Sp}(E)$.  Par
exemple, on peut prendre $I={\rm J}_{2g}$ dans les notations du \S
IV.\ref{fsiegelclass}, ${\rm U}_{{\rm h}}(\R)$ \'etant alors le groupe $K$
du \S IV.\ref{fautsiegel} (de complexifi\'e ${\rm j}(-,i {\rm 1}_g) : K
\rightarrow \GL_g(\C)$ {\it loc.  cit.}).  Le choix d'une d\'ecomposition de
l'espace hermitien $(E,{\rm h})$ en somme orthogonale de droites d\'etermine
un tore maximal anisotrope $T$ tel que $T(\R) \subset K$.  Dans les
notations du \S \ref{paramdiscA} (et conform\'ement au \S VI.\ref{carinf}),
on a ${\rm X}^\ast(T)=\oplus_{i=1}^g \Z \varepsilon_i$, $$\Phi=\{\pm 2
\varepsilon_i, i=1,\dots,g\} \cup \{\pm \varepsilon_i \pm \varepsilon_j,
1\leq i < j \leq g\}$$ et $\Phi_{\rm c}=\{\pm (\varepsilon_i-\varepsilon_j),
1 \leq i < j \leq g\}$.  En particulier, $W_{\rm r}={\rm W}(\Phi_c)$ n'est
autre que le groupe sym\'etrique $\got{S}_g$ agissant sur ${\rm X}^\ast(T)
\simeq \Z^g$ de mani\`ere habituelle et $W_{\rm r}^{\rm ad}=W_{\rm r} \times
\{\pm {\rm id}\}$.  L'ensemble de bases $\mathcal{B}(T)$ est constitu\'e de
deux $W_{\rm r}$-orbites, \'echang\'ees entre elles par $x \mapsto -x$,
l'une d'entre \'etant par exemple
$$\Delta=\{2\varepsilon_g,-\varepsilon_g-\varepsilon_{g-1},\varepsilon_{g-1}+\varepsilon_{g-2},-\varepsilon_{g-2}-\varepsilon_{g-3},\dots,(-1)^{g-1}(\varepsilon_2+\varepsilon_1)\}.$$
Si $\varepsilon_i^\ast \in {\rm X}_\ast(T)$ est la base duale de
$(\varepsilon_i)$, observons que $\Delta$ est la base du syst\`eme positif
$\{ \alpha \in \Phi, f(\alpha)>0\}$ o\`u $f = \varepsilon_g^\ast - 2
\,\varepsilon_{g-1}^\ast + 3 \,\varepsilon_{g-2}^\ast + \cdots + (-1)^{g-1}
\,(g-1)\, \varepsilon_1^\ast$.  Un petit calcul montre alors que $\rho^\vee
= \sum_{i=1}^g (-1)^{i-g}\, \frac{2g-2i+1}{2}\, \varepsilon_i^\ast$, de
sorte que $$ \rho^\vee \equiv \frac{1}{2} \sum_{i=1}^g \varepsilon_i^\ast
\bmod 2 {\rm X}_\ast(T).$$ Observons que cet \'el\'ement est bien invariant
par $W_{\rm r}$.  Soit $A \subset W$ le sous-groupe constitu\'e des
\'el\'ements $a \in W$ tels que $a(\varepsilon_i)= \pm \varepsilon_i$; il
est manifestement isomorphe \`a $\{\pm 1 \}^g$.  On constate que tout
\'el\'ement de ${\rm X}_\ast(T)$ est congru modulo $2$ \`a un \'el\'ement de
la forme $a \rho^\vee-\rho^\vee$ pour un unique $a \in A$: l'action de $A$
sur $H^1(\R,T)$ est donc simplement transitive, et $\kappa_{\Delta}$ (resp. 
$\chi_O$) est bijective.  Observons enfin que si l'on remplace $O$ par la
${\rm W}_{\rm r}$-orbite $-O$, ce qui revient \`a changer $w$ en $-w$, on a
$$\kappa_\Delta \equiv \kappa_{-\Delta} + 2 \rho^\vee \bmod 2{\rm
X}_\ast(T).$$ \ps\ps 

Consid\'erons maintenant la s\'erie discr\`ete holomorphe dans $\Pi_V$ 
(\S VI.\ref{carinf},\S\ref{paramdiscA}). On v\'erifie ais\'ement qu'il
existe exactement deux ${\rm W}_r$-orbites de bases de $\Phi$ dont le
syst\`eme positif associ\'e est tel que les racines compactes sont plus petites que les racines non-compactes, \`a
savoir celles des bases $\pm \Delta'$ o\`u 
$$\Delta'=\{2\varepsilon_g\} \cup \{\varepsilon_{i+1}-\varepsilon_i, i=1,\dots,g-1\},$$
(consid\'erer par exemple la forme lin\'eaire $\sum_{i=1}^g (g+i)\, \varepsilon_i^\ast$).
Suivant les conventions choisies, qu'il ne sera pas n\'ec\'essaire de
pr\'eciser, l'une de ces bases donne lieu aux repr\'esentations
not\'ees $\pi'_W$ au \S VI.\ref{carinf}, et son oppos\'ee \`a sa conjugu\'ee
ext\'erieure sous ${\rm PGSp}_{2g}(\R)$ (dans la litt\'erature on trouve
parfois le nom de s\'eries discr\`etes {\it holomorphes} et {\it
anti-holomorphes} pour distinguer ces deux types). On constate alors que $\Delta'=w^{-1}\Delta$ 
o\`u $w \in W$ est l'\'el\'ement envoyant $\varepsilon_i$ sur $(-1)^{g-i}\varepsilon_i$. Ainsi, 
\begin{equation}\label{formulecarhol} \kappa_\Delta(\pi_{\Delta',V}) \equiv w\rho^\vee - \rho^\vee \equiv \sum_{i \not \equiv g \bmod 2} \varepsilon_i^\ast \bmod 2,\end{equation}
conform\'ement au calcul fait dans~\cite[Lemma 9.1]{chrenard2}. 
On obtient de m\^eme $\kappa_\Delta(\pi_{-\Delta',V}) \equiv \kappa_\Delta(\pi_{\Delta',V})+2\rho^\vee \bmod 2{\rm X}_\ast(T)$. \ps\ps 

\subsection{Formes pures des groupes d\'eploy\'es}\label{paramdiscE} Dans ce
paragraphe, nous rappelons comment la param\'etrisation du \S\ref{paramdiscB}
s'\'etend \`a toutes les formes pures du $\R$-groupe d\'eploy\'e $H$,
suivant Vogan, Kottwitz, Arthur, Shelstad et Adams \cite{abv}\cite[\S
1]{arthurkgroup}.  Notre expos\'e est largement inspir\'e de la
pr\'esentation agr\'eable d'Adams \cite{adams} ainsi que des notes de
Shelstad \cite{shelstadbanff08,shelstadnotesbanff2011}. Nous renvoyons \`a
\cite[Ch. III \S 1]{serrecg}, \cite[\S 2]{abv} et \cite{borovoi14} pour des g\'en\'eralit\'es
sur les formes des groupes r\'eels. \ps\ps 

Consid\'erons d'abord un $\R$-groupe $G$ quelconque. L'ensemble $G(\C)$ est
muni d'une action de ${\rm Gal}(\C/\R)=\{1, \sigma \}$.  On
consid\`ere $$Z^1(\R,G)=\{x \in G(\C), x\sigma(x)=1\}.$$  Le groupe $G(\C)$
op\`ere sur $Z^1(\R,G)$ par $(g,x) \mapsto g x\sigma(g)^{-1}$, l'ensemble
quotient \'etant l'ensemble de cohomologie usuel $H^1(\R,G)$.  \`A tout
\'el\'ement $x \in Z^1(\R,G)$ est associ\'e une involution $$\sigma_x={\rm
int}_x \circ \sigma$$ de $G(\C)$.  C'est l'involution galoisienne
d'une unique structure r\'eelle du $\C$-groupe $G_\C$, dont nous noterons
$G_x$ le $\R$-groupe associ\'e.  En particulier, $$G_x(\R)=\{g \in G(\C),
\sigma(g)=x^{-1}gx\}.$$ Une telle forme r\'eelle de $G$ est dite {\it pure}. 
Le lemme suivant est \'evident.

\begin{lemme}\label{lemmestabreel} Le stabilisateur dans $G(\C)$ de $x \in
Z^1(\R,G)$ est $G_x(\R)$.  \end{lemme}

Si $x$ et $x' \in Z^1(\R,G)$ ont m\^eme classe dans $H^1(\R,G)$, et si $h
\in G(\C)$ est tel que $hx\sigma(h)^{-1}=x'$, alors $\sigma_{x'} \circ {\rm
int}_h = {\rm int}_h \circ \sigma_x$, de sorte que $${\rm int}_h : G_x
\rightarrow G_{x'}$$ est un isomorphisme d\'efini sur $\R$.  Un point
important est qu'il ne d\'epend pas du $h$ choisi ci-dessus, du moins {\it
modulo les ${\rm int}_g$ avec $g \in G_{x}(\R)$}, d'apr\`es le lemme
ci-dessus.  Autrement dit, les $\R$-groupes $G_x$ et $G_{x'}$ sont
naturellement isomorphes, et m\^eme canoniquement modulo ``automorphismes
int\'erieurs'', si $x$ et $x'$ sont \'equivalents.  Tout en s'imposant une
certaine prudence, il y a donc un sens \`a parler du $\R$-groupe $G_c$
d\'efini ``\`a automorphismes int\'erieurs pr\'es'' comme \'etant le
$\R$-groupe $G_x$ pour n'importe quel $x \in Z^1(\R,G)$ dans la classe $c
\in H^1(\R,G)$.  En particulier, il y a un sens \'evident \`a parler de
$\Pi_{\rm unit}(G_c(\R))$ pour $c \in H^1(\R,G)$.  \ps\ps

\begin{example}\label{exemplefpure} {\rm Donnons quelques exemples classiques
\cite[Ch. III \S 1.2]{serrecg}. Si $E$ est un $\R$-espace
vectoriel de dimension finie alors ${\rm Z}^1(\R,{\rm GL}_E)$ est l'ensemble
des involutions semi-lin\'eaires de $E \otimes_\R \C$.  D'apr\`es le
th\'eor\`eme de Hilbert 90, $x \mapsto E_x=(E \otimes_\R \C)^{x={\rm Id}}$
identifie ${\rm Z}^1(\R,{\rm GL}_E)$ \`a l'ensemble des structures r\'eelles du $\C$-espace vectoriel
$E \otimes_\R \C$, et $({\rm GL}_E)_x={\rm GL}_{E_x}$ pour tout $x \in {\rm
Z}^1(\R,{\rm GL}_E)$.  On v\'erifie imm\'ediatement que si $E$ est un ${\rm
q}$-module sur $\R$, alors ${\rm Z}^1(\R,{\rm O}_E) \subset {\rm Z}^1(\R,{\rm
GL}_E)$ s'identifie aux structures r\'eelles $E_x \subset E \otimes_\R
\C$ telles que ${\rm q}(E_x) \subset \R$.  Si $x \in {\rm Z}^1(\R,{\rm
O}_E)$ on dispose donc d'un ${\rm q}$-module $E_x$ pour la forme ${\rm q}_{|E_x}$ et
l'on a $({\rm O}_{E})_x={\rm O}_{E_x}$. L'application $x \mapsto E_x$ induit alors
une bijection entre $H^1(\R,{\rm O}_{E})$ et l'ensemble \`a $\dim(E)+1$
\'el\'ements des classes d'isomorphisme de ${\rm q}$-modules sur $\R$ de rang $\dim \,E$
(la signature). 
On v\'erifie enfin que $x \in {\rm Z}^1(\R,{\rm SO}_E)$ si, et seulement si, $E$ et $E_x$ ont m\^eme
discriminant, auquel cas $({\rm SO}_{E})_x={\rm SO}_{E_x}$, puis que $x
\mapsto E_x$ induit une bijection entre $H^1(\R,{\rm SO}_{E})$ et l'ensemble 
des classes d'isomorphisme de ${\rm q}$-modules sur $\R$ de m\^eme rang et discriminant que $E$. On dispose
d'un discours analogue dans le cas altern\'e, qui est seulement plus simple
car il n'y a qu'une forme altern\'ee non d\'eg\'en\'er\'ee de chaque
dimension paire ($H^1(\R,{\rm Sp}_{2g})=1$). }
\end{example} 

Revenons maintenant \`a notre $\R$-groupe d\'eploy\'e $H$, muni d'un tore maximal anisotrope $T$. Si $c \in
H^1(\R,H)$ nous allons rappeler la param\'etrisation de Shelstad de l'ensemble $\Pi^c_V$ des s\'eries discr\`etes de $H_c(\R)$ ayant m\^eme caract\`ere infinit\'esimal que $V \in {\rm Irr}(H_\C)={\rm Irr}((H_c)_\C)$.  L'inclusion $T \subset H$ induit une injection naturelle $Z^1(\R,T) \rightarrow Z^1(\R,H)$ ainsi qu'une application
\begin{equation}\label{isomh1}W \backslash H^1(\R,T) \longrightarrow H^1(\R,H),\end{equation}
l'action de $W$ sur $H^1(\R,T)$ \'etant celle rappel\'ee au \S\ref{paramdiscB}. Cette application est bijective d'apr\`es Shelstad (voir aussi~\cite{borovoi14}). En particulier toute forme pure de $H$ est isomorphe comme $\R$-groupe \`a $H_t$ pour $t \in T_2$.\ps\ps 
 Les formes $H_t$ avec $t \in T_2$, et plus g\'en\'eralement $t \in Z^1(\R,T)$, ont de sympatique qu'elles partagent toutes $T$ pour tore maximal anisotrope, car $\sigma_t$ co\"incide avec $\sigma$ sur $T$. Comme de plus $H_\C=(H_t)_\C$, le syst\`eme de racines $\Phi$ de $(H_\C,T_\C)$ est canoniquement celui de $((H_t)_\C,T_\C)$, et il en va de m\^eme de son groupe de Weyl $W$ (le sous-groupe de Weyl r\'eel ${\rm W}(H_t,T) \subset W$ d\'epend en revanche de $t$ bien entendu).
Si $\Delta$ est une base de $\Phi$ et $t \in T_2$, il y a donc un sens \`a consid\'erer la s\'erie discr\`ete $\pi_{\Delta,V,t}$ de $H_t(\R)$ associ\'ee par Harish-Chandra \`a la base $\Delta$ de $\Phi$, dans les notations du \S\ref{paramdiscA}. On rappelle que l'ensemble $\mathcal{B}(T)$ d\'efini au \S\ref{paramdiscB} est relatif au couple $(H,T)$.

\begin{definition}\label{defparshel2}Soient $c \in H^1(\R,H)$ et $V \in {\rm
Irr}(H_\C)$.  Soient $T$ un tore maximal anisotrope de $H$, $\Delta \in \mathcal{B}(T)$ et $\rho^\vee$ la demi-somme des
co-racines de $T$ positives relativement \`a $\Delta$.  L'application de
param\'etrisation de Shelstad relativement \`a $\Delta$ est l'unique application $$\kappa_\Delta^c : \Pi_V^c
\longrightarrow {\rm X}_\ast(T) \otimes \Z/2\Z, \pi \mapsto \kappa_\Delta^c(\pi),$$
telle que pour tout $t \in T_2$ dans la classe $c$, et tout $w \in W$, on
ait $\kappa_\Delta^c(\pi_{w^{-1}\Delta,V,t}) \equiv w(\mu+\rho^\vee) - \rho^\vee
\bmod 2 {\rm X}_\ast(T)$, o\`u $\mu \in {\rm X}_\ast(T)$ est tel que
$t=e^{i\pi \mu}$.  \end{definition} 

Cach\'e derri\`ere cette d\'efinition se
trouve le fait suivant : soient $w \in W$ et $t \in T_2$, soit $n \in H(\C)$
un repr\'esentant de $w$ tel que $n \star t = w \cdot t$, de sorte que 
${\rm int}_n : H_t \rightarrow H_{w \cdot t}$ d\'efinisse un $\R$-isomorphisme, 
alors la restriction de la repr\'esentation $\pi_{w\Delta,V,w \cdot t}$ de $H_{w\cdot t}(\R)$ \`a ${\rm int}_n$ est
isomorphe \`a la repr\'esentation $\pi_{\Delta,V,t}$ de $H_t(\R)$.  Lorsque
$c=1$ alors $\kappa_\Delta^1=\kappa_\Delta$, et l'on retrouve bien entendu la
d\'efinition~\ref{defparshel}.  On v\'erifierait comme au \S
\ref{paramdiscB} que $\kappa_{w \Delta}^c = w \circ \kappa_{\Delta}^c$ pour tout $w \in {\rm W}(H,T)$. Le lemme ci-dessous (i) assure de plus qu'elle est injective, de sorte que
l'application $$\coprod_{c \in H^1(\R,G)} \kappa^c_\Delta : \coprod_{c \in H^1(\R,H)}
\Pi_V^c \longrightarrow {\rm X}_\ast(T) \otimes \Z/2\Z$$ est bijective.  \ps\ps 

La concordance de cette pr\'esentation avec la d\'efinition donn\'ee par
Shelstad {\it loc.  cit.}, qui associe \`a $\pi_{w^{-1}\Delta,V,t}$
la classe de $w \cdot t$ dans $H^1(\R,T)$, d\'ecoule encore imm\'ediatement
du lemme~\ref{lemmelankott}.  C'est essentiellement le point de vue donn\'e
par Adams dans \cite{adams}, \`a ceci pr\`es que son point de d\'epart est
un $\R$-groupe \`a points r\'eels compact, plut\^ot que d\'eploy\'e (et
qu'il consid\`ere les formes int\'erieures g\'en\'erales).  La d\'efinition
ci-dessus admet une interpr\'etation duale identique \`a celle mentionn\'ee
au \S\ref{paramdiscC} en terme des param\`etres discrets de $H$, et conduit \`a une application canonique
$$\chi_O^c = \beta \circ (i_\Delta \otimes \Z/2\Z) \circ \kappa_\Delta^c, \, \, \, \, \Pi_V^c \rightarrow  {\rm Hom}_{\rm groupes}({\rm C}_\varphi,\C^\times),$$
qui ne d\'epend que de la $W_{\rm r}$-orbite $O$ de $\Delta$ dans $\mathcal{B}(T)$ (soulignons encore que ce choix est relatif \`a $H$, et non \`a $H_c$). Mentionnons que
dans son expos\'e, et suivant Arthur et Kottwitz, Shelstad se limite aux $c
\in H^1(\R,H)$ qui sont dans l'image de l'application naturelle
$H^1(\R,H_{\rm sc}) \rightarrow H^1(\R,H)$, $H_{\rm sc} \rightarrow H$
d\'esignant le rev\^etement simplement connexe de $H$.  La r\'eunion
disjointe des $H_c$ index\'ee par les telles classes $c$ forme alors un
$K$-{\it groupe} au sens d'Arthur.  Cela revient \`a se restreindre aux
\'el\'ements $t \in T_2$ de la forme $e^{i\pi \mu}$ o\`u $\mu \in {\rm
X}_\ast(T)$ est une somme de coracines de $(H_\C,T_\C)$, soit encore aux
\'el\'ements de ${\rm X}^\ast(\widehat{T}) \otimes \Z/2\Z$ triviaux sur le
centre de $\widehat{H}$. \ps\ps 

\begin{lemme}\label{stabwreel} Soient $T$ un tore maximal anisotrope de $H$ et $t \in T_2$.
\begin{itemize} \ps\ps 
\item[(i)] ${\rm W}(H_t,T)=\{w \in W, w \cdot t = t\}$. \ps\ps 
\item[(ii)] ${\rm Int}\, tt_{\rm b}$ est une involution de Cartan de $H_t$, o\`u $t_{\rm b} \in T_2$ est l'\'el\'ement d\'efini au \S\ref{paramdiscB} et associ\'e \`a une $W_{\rm r}$-orbite dans $\mathcal{B}(T)$.  \ps\ps 
\end{itemize} 

\end{lemme}

\begin{pf} On constate que $w \cdot t = t$ si, et seulement si, il existe un
repr\'esentant $n \in H(\C)$ de $w$ qui fixe $t \in Z^1(\R,H)$, i.e.  dans
$H_t(\R)$ d'apr\`es le lemme~\ref{lemmestabreel}.  Le (ii) r\'esulte du cas
particulier $t=1$, qui est le lemme~\ref{lemmelankott} (i).  \end{pf}

Si $c \in H^1(\R,H)$ et si $V \in {\rm Irr}(H_\C)$ alors $V$ peut \^etre vue
par restriction comme une repr\'esentation irr\'eductible de dimension finie
de $H_c(\R)$ bien d\'efinie \`a isomorphisme pr\`es.  Lorsque $H_c(\R)$ est
compact, c'est l'unique \'el\'ement du singleton $\Pi_V^c$
(\S\ref{paramdiscA}).  Le point (ii) du lemme ci-dessus montre que cela se
produit si, et seulement si, $t_{\rm b} \in T_2$ et $c$ est la classe d'un
$t \in t_b {\rm Z}(H)$.  On rappelle que ${\rm Z}(H) \subset T_2$ car $-1
\in W$ (cela vient de ce que $H$ est \`a la fois d\'eploy\'e sur $\R$ et
poss\`ede un tore maximal compact).  

\begin{cor}\label{corcompact} Soient $c \in H^1(\R,H)$ telle que
$H_c(\R)$ est compact, $V \in {\rm Irr}(H_\C)$, $\Delta \in \mathcal{B}(T)$ et $\rho^\vee \in \frac{1}{2}{\rm X}_\ast(T)$ associ\'e \`a $\Delta$. Alors $\rho^\vee
\in {\rm X}_\ast(T)$ et $\kappa^c_\Delta(V) \equiv \rho^\vee + \nu \bmod 2 {\rm
X}_\ast(T)$, o\`u $e^{i\pi \nu} \in {\rm Z}(H)$.\end{cor}

Consid\'erons l'exemple plus int\'eressant de la s\'erie discr\`ete
holomorphe du groupe sp\'ecial orthogonal de signature $(m,2)$ pour $m\geq
1$ impair (qui n'est d\'eploy\'e que si $m \leq 3$).  Partons du $\R$-groupe
d\'eploy\'e $H={\rm SO}_{r+1,r}$.  \'Ecrivons ${\rm H}(\R^r) \oplus \R$
comme somme orthogonale d'une droite $D$ et de plans $P_i$, $i=1,\dots r$,
chacun suppos\'e d\'efini et de signe que nous pr\'eciserons plus tard. 
Cette d\'ecomposition d\'efinit un unique tore anisotrope $T$ de $H$ pr\'eservant
$D$ et chacun des $P_i$. Choisissons de mani\`ere
arbitraire, pour $i=1,\dots, r$, l'une des deux droites isotropes $\ell_i$ de
${\rm P}_i \otimes \C$, et notons $\varepsilon_i \in {\rm X}^\ast(T_\C)$ le
caract\`ere de $T$ sur $\ell_i$. La suite $\ell_1,\dots,\ell_r$ d\'efinit
comme au \S VI.\ref{exdorad} un unique sous-groupe de Borel de $H_\C$
contenant $T_\C$, correspondant \`a la base standard $$\Delta =\{
\varepsilon_i, i=1,\dots,r\} \cup \{\varepsilon_i \pm \varepsilon_j, 1\leq
i< j \leq r\}$$ de $\Phi=\Phi(H_\C,T_\C)$.  Supposons maintenant que l'on a
choisi les $P_i$ de signe $(-1)^{i-1}$, et $D$ de signe $(-1)^r$, ce qui est
loisible.  On constate qu'aucun \'el\'ement de $\Delta$ n'est compact, de
sorte que $\Delta \in \mathcal{B}(T)$ (il n'est pas difficile de voir que
$\mathcal{B}(T)$ ne forme qu'une ${\rm W}_r$-orbite).  De plus, le \S
VI.\ref{exdorad} montre que la demi-somme $\rho^\vee$ des coracines de
$(H_\C,T_\C)$ positives relativement \`a $\Delta$ est $$\rho^\vee =
\sum_{i=1}^r (r-i+1) \varepsilon_i^\ast \in {\rm X}_\ast(T_\C),$$ o\`u
$\varepsilon_i^\ast \in {\rm X}_\ast(T_\C)$ d\'esigne la base duale de
$(\varepsilon_i)$, de sorte que $\varepsilon_i(t_{\rm b})=(-1)^{r-i+1}$. \ps\ps 

Ceci \'etant dit int\'eressons-nous aux formes r\'eelles de $H$.  D'apr\`es
l'exemple \ref{exemplefpure}, il existe une unique classe $c \in H^1(\R,H)$
telle que $H_c$ soit isomorphe au groupe sp\'ecial orthogonal de signature
$(2,2r-1)$.  Concr\`etement, si $t \in T_2$ agit par multiplication par $s_j=\pm 1$ sur le plan $P_j$, et par $1$ sur $D$, la forme r\'eelle du
${\rm q}$-module ${\rm H}(\R^r) \oplus \R$ associ\'ee \`a $t$ est la somme
directe de $D$, des $P_j$ tels que $s_j=1$, et des $i P_j$ tels que
$s_j=-1$.  En particulier, $H_t$ est le groupe sp\'ecial orthogonal de
signature $(2 a,b)$ o\`u $a$ est le nombre des $1 \leq j \leq r$ tels que
$s_j (-1)^{j-1} \neq (-1)^r$ et $2a+b=2r+1$.  Autrement
dit, $H_t \simeq H'$ si, et seulement si, il existe un entier $1\leq s\leq
r$ tel que $t = t_{\rm b} e^{i\pi \varepsilon_s^\ast}$. Ces \'el\'ements de $T_2$ ne
formant qu'une $W$-orbite pour l'action tordue, ils appartiennent bien \`a une
m\^eme classe $c \in {\rm H}^1(\R,H)$. Identifions par exemple $H'$ \`a
$H_{t'}$ o\`u $$t'=t_{\rm b}e^{i\pi \varepsilon_1^\ast}$$ (tout
$\R$-isomorphisme de $H'$ \'etant de la forme ${\rm Int}\, h$ avec $h \in
H'(\R)$, le choix de l'identification n'importe pas).  L'ensemble
$\Phi'_{\rm c} \subset \Phi$ des racines compactes de $(H_{t'},T)$ est
manifestement $$\Phi'_{\rm c}=\{\pm \varepsilon_i, i=2,\dots,r\} \cup \{ \pm
\varepsilon_i \pm \varepsilon_j, 2 \leq i<j\leq r\},$$
d'apr\`es le lemme \ref{stabwreel}. En consid\'erant la
forme lin\'eaire $(2r-2) \, \varepsilon_1^\ast+\sum_{i=2}^r (r-i+1) \,
\varepsilon_i^\ast$ on constate que le syst\`eme positif de la base $\Delta$
plus haut est tel que tout \'el\'ement de $\Phi'_{\rm c}$ est plus petit que
tout \'el\'ement de $\Phi-\Phi'_c$ (\S \ref{paramdiscA}).  C'est m\^eme
l'unique telle base modulo action du groupe de Weyl r\'eel de $(H_{t'},T)$. 
Il existe donc une unique s\'erie discr\`ete holomorphe $\pi_{{\rm
hol},V}$ de $H'$ ayant m\^eme caract\`ere infinit\'esimal que $V \in {\rm
Irr}(H_\C)$, et l'on a \begin{equation}\label{carholsoimpair}\kappa_\Delta^c(\pi_{\rm
hol,V}) \equiv (r-1) (\varepsilon_1^\ast + \varepsilon_2^\ast) +
\sum_{i=3}^r (r-i+1) \varepsilon_i^\ast \bmod 2{\rm
X}_\ast(T_\C).\end{equation} En effet, c'est la d\'efinition~\ref{defparshel2}
appliqu\'ee \`a $w=1$ et $t=t'=e^{i\pi (\rho^\vee+\varepsilon_1^*)}$.

\subsection{Paquets d'Adams-Johnson}\label{paramaj}

Soient $H$ un $\R$-groupe semi-simple d\'eploy\'e poss\'edant des s\'eries
discr\`etes et $T$ un tore maximal anisotrope de $H$.  On note encore
$\Phi=\Phi(H_\C,T_\C) \subset {\rm X}^\ast(T)$ son syst\`eme de racines et
$W$ le groupe de Weyl de $\Phi$.\ps\ps 

Dans les paragraphes qui suivent nous allons rappeler bri\`evement certains
ensembles, ou {\it paquets}, de repr\'esentations irr\'eductibles unitaires
de $H(\R)$ qui ont \'et\'e d\'efinis par Adams et Johnson dans \cite{AJ} et
qui jouent un r\^ole important dans la th\'eorie d'Arthur \cite[\S
5]{arthurunipotent} (voir aussi \cite{abv} pour un contexte tr\`es
g\'en\'eral).  Le point de d\'epart est la donn\'ee d'un {\it param\`etre
d'Adams-Johnson} $$\varphi : {\rm SL}_2(\C) \times {\rm W}_\R \rightarrow
\widehat{H}(\C)$$ qui est un morphisme de groupes satisfaisant certaines
propri\'et\'es, qui seront d'abord discut\'ees informellement et
pr\'ecis\'ees plus bas.  Adams et Johnson lui associent un sous-ensemble
fini $$\Pi_{\rm AJ}(\varphi) \subset \Pi_{\rm unit}(H(\R)),$$ ne d\'ependant
que de la classe de $\widehat{H}(\C)$-conjugaison de $\varphi$.  \ps\ps 

Le param\`etre $\varphi$ d\'etermine d'abord une repr\'esentation $V_\varphi
\in {\rm Irr}(H_\C)$.  Par exemple, les param\`etres d'Adams-Johnson
triviaux sur le facteur ${\rm SL}_2(\C)$ s'identifient exactement avec les
param\`etres discrets de Langlands rappel\'es au \S\ref{paramdiscC}, et pour
un tel $\varphi$ alors par d\'efinition $\Pi_{\rm
AJ}(\varphi)=\Pi_{V_\varphi}$.  En g\'en\'eral, $\Pi_{\rm AJ}(\varphi)$ est
constitu\'e de repr\'esentations ayant m\^eme caract\`ere infinit\'esimal
que $V_\varphi$; mieux, elles poss\`edent de la
$(\mathfrak{h},K)$-cohomologie \`a coefficients dans $V_\varphi^\ast$
\cite{vz}.  Concr\`etement, \`a chaque base $\Delta \subset \Phi$ le
param\`etre $\varphi$ fait correspondre un sous-groupe parabolique
$P_{\Delta,\varphi} \subset H_\C$ contenant $T_\C$.  Soit $L_{\Delta,\varphi}
\subset H$ le sous-groupe de Levi de $P_{\Delta,\varphi}$ contenant $T$,
il est n\'ecessairement d\'efini sur $\R$ car $T(\R)$ est compact (bien
s\^ur, aucun des paraboliques propres de $H_\C$ contenant $T$ n'est d\'efini
sur $\R$).  Enfin, $\varphi$ d\'etermine une repr\'esentation $\rho$ de
$L_{\Delta,\varphi}(\R)$ dimension $1$. Sa description exacte par Adams et
Johnson est assez d\'elicate, du moins quand $L_{\Delta,\varphi}(\R)$ n'est
pas connexe\footnote{Lorsque $L_{\Delta,\varphi}(\R)$ est connexe, le caract\`ere $\rho$ est d\'etermin\'e par sa diff\'erentielle en l'identit\'e, elle-m\^eme caract\'eris\'ee par la propri\'et\'e que la repr\'esentation $\pi_{\Delta,\varphi}$ d\'efinie ci-apr\`es doit avoir le m\^eme caract\`ere infinit\'esimal que $V$.}, mais dans une large mesure sa compr\'ehension ne sera
heureusement pas n\'ecessaire \`a notre discussion.  La donn\'ee de
$P_{\Delta,\varphi}$ et de $\rho$ permet alors de d\'efinir une
repr\'esentation $$\pi_{\Delta,\varphi} \in \Pi_{\rm unit}(H(\R))$$ par
induction cohomologique en un degr\'e appropri\'e \cite{voganbook,vogan}. 
L'ensemble de ces repr\'esentations, lorsque $\Delta$ parcourt les bases de
$\Phi$, est par d\'efinition le paquet $\Pi_{\rm AJ}(\varphi)$.  Soit
$W_{\rm r}={\rm W}(H,T)\subset W$ (\S\ref{paramdiscA}).  Une base
$\Delta \subset \Phi$ \'etant fix\'ee, et si $L=L_{\Delta,\varphi}$,
l'application $W \rightarrow \Pi_{\rm AJ}(\varphi)$, $w \mapsto \pi_{w \Delta,\varphi}$, induit alors une bijection (d\'ependante de $\Delta$) $$W_{\rm r} \backslash W /
{\rm W}(L_\C,T_\C) \isomo \Pi_{\rm AJ}(\varphi).$$ \ps\ps 

Pr\'ecisons maintenant les axiomes  {\rm (AJ1) \& (AJ2)} d\'efinissant les
param\`etres d'Adams-Johnson, suivant \cite{AJ}\cite[p. 
195]{kottwitzannarbor}\cite[App.  A]{chrenard2}.  Soit $\varphi : {\rm
SL}_2(\C) \times {\rm W}_\R \rightarrow \widehat{H}(\C)$ un morphisme de
groupes suppos\'e continu, alg\'ebrique sur le facteur ${\rm SL}_2(\C)$, et
tel que $\varphi(1 \times {\rm W}_\R)$ est constitu\'e d'\'el\'ements
semi-simples.  Consid\'erons l'homomorphisme $\widetilde{\varphi} : {\rm
W}_\R \rightarrow \widehat{H}(\C)$ obtenu en composant $\varphi$ par le
morphisme d'Arthur $$ {\rm W}_\R \rightarrow {\rm SL}_2(\C) \times {\rm
W}_\R, \, g \mapsto \left[ \begin{array}{cc} |\eta(g)|^{1/2} & 0 \\ 0 &
|\eta(g)|^{-1/2} \end{array}\right] \times g,$$ o\`u $\eta: {\rm W}_\R
\rightarrow \R^\times$ est le caract\`ere rappel\'e au \S\ref{apparitionWR}. 
Le sous-groupe $\widetilde{\varphi}(\C^\ast) \subset \widehat{H}(\C)$ est
connexe et constitu\'e d'\'el\'ements semi-simples, on peut donc l'inclure
dans un tore maximal
 $$S \subset \widehat{H}.$$ 
Il existe alors un unique couple $\lambda,\mu \in {\rm X}_\ast(S) \otimes \C$ tels que $\lambda-\mu \in {\rm X}_\ast(S)$ et $\xi(\widetilde{\varphi}(z))=z^{\langle \xi,\lambda \rangle} \overline{z}^{\langle \xi,\mu \rangle}$ pour tout $\xi \in {\rm X}^\ast(S)$ et tout $z \in \C^\ast$ (voir la note de bas de page \S \ref{apparitionWR}). \ps\ps

{\rm (AJ1)} -- {\it  La classe de $\widehat{H}(\C)$-conjugaison de
$\lambda$, vue dans ${\rm Lie} \, \widehat{H}$, est le caract\`ere
infinit\'esimal d'une repr\'esentation de dimension finie $V_\varphi \in
{\rm Irr}(H_\C)$.} \ps\ps

Cela entra\^ine en particulier que $S$ est l'unique tore maximal de
$\widehat{H}$ contenant $\widetilde{\varphi}(\C^\ast)$.  Cela munit
\'egalement $\widehat{H}$ d'un unique sous-groupe de Borel $B$ contenant
$S$ tel que $\lambda$ soit dominant relativement \`a $B$.  On consid\`ere
ensuite le centralisateur $M \subset \widehat{H}(\C)$ du sous-groupe connexe
commutatif constitu\'e d'\'el\'ements semi-simples $\varphi(1 \times
\C^\times)$ ; c'est donc un sous-groupe de Levi (d'un parabolique) de
$\widehat{H}$.  Il contient $S$.  \ps\ps 

La donn\'ee d'une base $\Delta$ de $\Phi$ permet d'identifier la donn\'ee
radicielle bas\'ee $({\rm X}^\ast(T),\Phi,\Delta,\cdots)^\vee$ \`a
$\Psi(\widehat{H},S,B)$, et fournit en particulier un isomorphisme
privil\'egi\'e $$i_\Delta : {\rm X}_\ast(T) \isomo {\rm X}^\ast(S)$$
envoyant $\Delta^\vee$ sur la base de $\Phi(\widehat{H},S)$ associ\'ee \`a $B$.  
Soit $L_{\Delta,\varphi} \subset H_\C$ l'unique sous-groupe de Levi (de
sous-groupes paraboliques) contenant $T$
tel que $i_\Delta(\Phi^\vee(L_{\Delta,\varphi},T_\C))=\Phi(M,S)$ 
(en particulier, $M \simeq \widehat{L_{\Delta,\varphi}}$). Soit
$P_{\Delta,\Phi} \subset H_\C$ l'unique sous-groupe parabolique de
sous-groupe de Levi $L_{\Delta,\varphi}$ contenant 
le sous-groupe de Borel de $H_\C$ contenant $T$ et
associ\'e \`a $\Delta$, c'est le sous-groupe mentionn\'e dans la description
informelle plus haut. L'axiome restant sert \`a la d\'efinition du
caract\`ere $\chi$ (voir \cite{AJ} et la r\'edaction de Ta\"ibi \cite[\S 4.2.2]{taibisiegel} pour plus de
pr\'ecisions \`a ce sujet). 
\ps\ps 

{\rm (AJ2)} \,\,--\,\, {\it L'homomorphisme ${\rm SL}_2(\C) \rightarrow M$
induit par $\varphi$ est principal, c'est-\`a-dire qu'il induit un
$\mathfrak{sl}_2$-triplet de ${\rm Lie}\, M$ {\rm principal} au sens de
Kostant {\rm{\cite{kostantsl2}}}.  De plus, le centralisateur ${\rm
C}_\varphi$ de ${\rm Im}\, \varphi$ dans $\widehat{H}(\C)$ est fini.} \ps\ps

\noindent La premi\`ere hypoth\`ese entra\^ine que le centralisateur de
$\varphi({\rm SL}_2(\C)\times \C^\times)$ dans $\widehat{H}(\C)$ est le
centre ${\rm Z}(M)$ de $M$.  Le groupe ${\rm C}_\varphi$ est donc le
sous-groupe de ${\rm Z}(M)$ fix\'e par la conjugaison par $\varphi(1 \times
j)$.  La seconde hypoth\`ese assure que $${\rm C}_\varphi = {\rm
Z}(M)_2\overset{{\mathrm{d\acute{e}f}}.}{=}\{z \in {\rm Z}(M), z^2=1\}.$$
\noindent (voir par exemple~\cite[Lemma A.1]{chrenard2}) Comme l'a
remarqu\'e Ta\"ibi \cite[\S 4.2.2]{taibisiegel}, sous la premi\`ere hypoth\`ese de
(AJ2), alors la seconde est aussi \'equivalente \`a demander que $\varphi$
est trivial sur $1 \times \R_{>0}$ (c'est \'evidemment n\'ecessaire, car $\R^\times$ est le centre de ${\rm W}_\R$, mais
\'egalement suffisant).  Quand $H$ est un groupe classique, il est facile 
d\'eterminer tous ses param\`etres d'Adams-Johnson : voir l'exemple \ref{exsp2gaj} et le lemme \ref{ajclassique}.

\ps\ps
\noindent {\sc Cas des formes r\'eelles pures}
\ps\ps
Consid\'erons pour clore ce paragraphe le cas g\'en\'eral des formes r\'eelles pures de
$H$.  Soit $c \in H^1(\R,H)$.  La construction  {\it
loc.  cit.} d'Adams et Johnson, qui n'est pas sp\'ecifique aux groupes d\'eploy\'es, associe \'egalement \`a tout param\`etre d'Adams-Johnson $\varphi$ de $H$, un
ensemble de repr\'esentations $$\Pi_{\rm AJ}^c(\varphi) \subset \Pi_{\rm
unit}(H_c(\R)).$$ Si $t \in T_2$ est dans la classe $c$, et si $\Delta$ est
une base de $\Phi$, on dispose encore du sous-groupe parabolique
$T_\C \subset P_{\Delta,\varphi} \subset ({H_t})_\C=H_\C$. Son sous-groupe
de Levi $T \subset L_{\Delta,\varphi,t} \subset H_t$, d\'efini sur $\R$, est la forme r\'eelle pure de $L_{\Delta,\varphi}$ associ\'ee \`a $\sigma_t$ (\S\ref{paramdiscE}). 
Adams et Johnson lui associent encore un caract\`ere $\rho : L_{\Delta,\varphi,t}(\R) \rightarrow \C^\ast$, et d\'efinissent une repr\'esentation
$$\pi_{\Delta,\varphi,t} \in \Pi_{\rm unit}(H_t(\R)) \isomo  \Pi_{\rm
unit}(H_c(\R))$$ par induction chomologique de $\rho$ \`a $H_t(\R)$,
relativement \`a $P_{\Delta,\varphi}$, en un degr\'e ad\'equat.  Ils posent
$\Pi_{\rm AJ}^c(\varphi)=\{\pi_{\Delta,\varphi,t}\}$, o\`u $\Delta$ parcourt
les bases de $\Phi$ et $t \in T_2$ les repr\'esentants de $c$. Comme au
\S\ref{paramdiscE}, les
variables $\Delta$ et $t$ sont redontantes et reli\'ees par la relation
$\pi_{w \Delta,\varphi, w \cdot t} \simeq \pi_{\Delta,\varphi,t}$ valable pour toute
base $\Delta \subset \Phi$,  tout $t \in T_2$ et tout $w \in W$. 
Si $t \in T_2$ et une base $\Delta \subset \Phi$ sont
fix\'es, et si $L=L_{\Delta,\varphi,t} \subset H_t$ est le sous-groupe de Levi associ\'e,
l'application $W \rightarrow \Pi_{\rm unit}(H(\R))$, $w \mapsto \pi_{w
 \Delta,\varphi,t}$, induit cette fois-ci une bijection $${\rm W}(H_t,T) \backslash W
/ {\rm W}(L_\C,T_\C) \isomo \Pi_{\rm AJ}^c(\varphi).$$ Le cas particulier
$c=1$ donne un autre point de vue sur $\Pi_{\rm
AJ}^1(\varphi)=\Pi_{\rm AJ}(\varphi)$.

\subsection{Exemple : param\`etres d'Adams-Johnson de ${\rm Sp}_{2g}$ }\label{exsp2gaj}

Consid\'erons le $\R$-groupe $H={\rm
Sp}_{2g}$.  Notons $V = \C^{2g+1}$ l'espace sous-jacent \`a la
repr\'esentation standard de $\widehat{H}={\rm SO}_{2g+1}$, muni de la forme
quadratique ${\rm q}$ standard.  Soit $\varphi : {\rm SL}_2(\C) \times {\rm
W}_\R \rightarrow {\rm SO}(V)$ un homomorphisme continu, alg\'ebrique sur le
facteur ${\rm SL}_2(\C)$, et trivial sur $1 \times \R_{>0}$.  Comme ${\rm
W}_\R/\R_{>0}$ est compact, le groupe ${\rm SL}_2(\C) \times {\rm W}_\R$
agit de mani\`ere semi-simple sur $V$, qui se d\'ecompose donc en somme
orthogonale $$V = \bigoplus_{j \in J} V_j$$ de sous-espaces irr\'eductibles
$V_j$.  En tant que repr\'esentation de ${\rm SL}_2(\C) \times {\rm W}_\R$
triviale sur $1 \times \R_{>0}$, $V_j$ est n\'ecessairement autoduale.
Observons aussi que si ${\rm q}_{|V_j}$ est d\'eg\'en\'er\'ee, alors ${\rm
q}(V_j)=0$ et il existe $j' \neq j$ tel que $V_{j'} \simeq V_j^\ast \simeq
V_j$. \ps\ps 

Supposons que $\varphi$ satisfait (AJ1), nous allons analyser
$\varphi$ et constater notamment que (AJ2) est automatiquement
satisfaite. D'apr\`es le cas I du \S\ref{paralgreg}, les $V_j$ sont
deux-\`a-deux non isomorphes (donc non-d\'eg\'en\'er\'es), et un seul d'entre eux est de dimension
impaire, disons $V_0$.  Il est alors \'evident que ${\rm C}_\varphi$ est le
sous-groupe fini de ${\rm SL}(V)$ constitu\'e des \'el\'ements $g$ tels que
$g(V_j) \subset V_j$ pour tout $j$ et $g_{|V_j}=\pm {\rm id}_{V_j}$.  De
plus, $V_0$ est une repr\'esentation irr\'eductible de ${\rm SL}_2(\C)$,
avec ${\rm W}_\R$ agissant par multiplication par un caract\`ere $\chi_0$,
et pour $j\neq 0$ on peut \'ecrire $$V_j \simeq Q_j \otimes R_j, \, \, \, Q_j \simeq {\rm
Sym}^{d_j-1} {\rm St}_2, \, \, \, R_j \simeq {\rm I}_{r_j},$$ avec $r_j>0$ et $d_j
+ r_j \equiv 1 \bmod 2$.  On peut munir $Q_j$ et $R_j$ de formes
bilin\'eaires non-d\'eg\'en\'er\'ees, pr\'eserv\'ees respectivement par
${\rm SL}_2(\C)$ et ${\rm W}_\R$, et dont le produit tensoriel est la forme
bilin\'eaire sur $V_j$ associ\'ee \`a ${\rm q}$.  La restriction de $R_j$
\`a $\C^\ast \subset {\rm W}_\R$ est somme directe des deux droites stables
et isotropes, disons $\ell_j$ et $\ell'_j$, l'\'el\'ement $z \in \C^\ast$
agissant par multiplication par $(\frac{z}{|z|})^{r_j}$ sur $\ell_j$.  Ainsi,
$M \subset {\rm SO}_V$ est le sous-groupe $${\rm SO}_{V_0} \times \prod_{j
\neq 0} M_j$$ o\`u $M_j \subset {\rm SO}_{V_j}$ est le stabilisateur des
lagrangiens transverses $Q_j \otimes \ell_j$ et $Q_j \otimes \ell'_j$.  En particulier, $M_j \simeq \GL_{d_j}$, $M$ est bien le sous-groupe de
Levi d'un parabolique de ${\rm SO}_V$, et ${\rm C}_\varphi={\rm Z}(M)_2$. 
Rappelons qu'un $\C$-morphisme ${\rm SL}_2 \rightarrow L$ avec $L$
classique (resp.  $\GL_d$) est principal si, et seulement si, la
repr\'esentation de ${\rm SL}_2$ compos\'ee de $f$ et de la repr\'esentation
standard (resp.  tautologique) de $L$ est irr\'eductible, ou somme de la
repr\'esentation triviale et d'une irr\'eductible non triviale si $L(\C)
\simeq {\rm SO}_{2r}(\C)$.  Cela montre que (AJ2) est satisfaite.  
Une analyse similaire \`a celle faite jusqu'\`a pr\'esent montre plus g\'en\'eralement le lemme
suivant.

\begin{lemme}\label{ajclassique} Supposons que $H_\C \in {\rm Class}_\C$ et
soit ${\rm St} : \widehat{H}(\C) \rightarrow {\rm SL}_n(\C)$ la
repr\'esentation standard.  Soit $\varphi : {\rm SL}_2(\C) \times {\rm W}_\R
\rightarrow \widehat{H}(\C)$ un morphisme continu, alg\'ebrique sur le
facteur ${\rm SL}_2(\C)$, et trivial sur $1 \times \R_{>0}$.  Alors
$\varphi$ est un param\`etre d'Adams-Johnson si, et seulement si, il
satisfait {\rm (AJ1)} et si la repr\'esentation ${\rm St} \circ \varphi :
{\rm W}_\R \rightarrow {\rm SL}_n(\C)$ est sans multiplicit\'es.  \end{lemme}

Poursuivons l'\'etude de l'exemple pr\'ec\'edent : d\'ecrivons le tore $S$ et le sous-groupe de Borel $B$ de ${\rm SO}_V$ associ\'es \`a $\varphi$. Convenons que $Q_0:=V_0$.  Observons que
la restriction de $Q_j$ au tore diagonal de ${\rm SL}_2(\C)$ est somme directe
de $d_j$ droites stables canoniques $\ell_{j,n}$, o\`u $n=\frac{d_j-1}{2},
\cdots, \frac{1-d_j}{2}$, l'\'el\'ement $\left[\begin{array}{cc}|z|^{1/2} & 0 \\ 0 & |z|^{-1/2}\end{array}\right]$
agissant sur $\ell_{j,n}$ par multiplication par $|z|^n$.  On a d'une part une
d\'ecomposition orthogonale $$V=V^+ \oplus V^- \oplus \ell_{0,0},$$ o\`u $V^+$
(resp.  $V^-$) est le lagrangien somme directe des $Q_j \otimes \ell_j$
(resp.  $Q_j \otimes \ell'_j$) et des $\ell_{0,n}$ avec $n>0$ (resp. 
$n<0$). D'autre part, l'espace $V^+$ est lui-m\^eme somme directe des $g$ droites 
$\ell_{0,n}$ (avec $n>0$) et $\ell_{j,m} \otimes \ell_j$ ($j \neq 0$, $m$ arbitraire). La condition (AJ1)
determine une unique mani\`ere d'ordonner ces $g$ droites isotropes, disons
$\C e_1, \dots, \C e_g \subset V^+$, sorte que pour tout $z \in \C^\times$ l'\'el\'ement $\widetilde{\varphi}(z)$ agisse
par multiplication par $z^{w_i} \overline{z}^{w'_i}$ sur $e_i$, o\`u $(w_i,w'_i)$
est un couple d'entiers tel que $$w_1>w_2>\cdots>w_g>0$$ (ils s'expriment
bien s\^ur simplement en terme des $r_j$ et $d_j$, mais il ne sera pas n\'ecessaire
d'expliciter comment). Le tore $S$ est le stabilisateur dans ${\rm SO}_V$ 
des droites $\C e_i$, $i=1,\dots,g$, et de $V^-$. Le sous-groupe de Borel
$S \subset B \subset {\rm SO}_{2g+1}$ est le stabilisateur du drapeau
associ\'e aux $e_i$ comme au \S VI.\ref{exdorad}, il rend l'\'el\'ement 
$\lambda_\varphi$ dominant par d\'ecroissance et positivit\'e des $w_i$. 
\ps\ps 

D\'ecrivons enfin les $L_{\Delta,\varphi}$. Consid\'erons le ${\rm a}$-module hyperbolique $E={\rm H}(\R^g)$, $I \in {\rm Sp}(E)$ de carr\'e $-{\rm id}_E$ et ${\rm h}$ la forme hermitienne d\'efinie positive sur $E$ associ\'ee comme dans le \S \ref{paramdiscD}. Choisissons une d\'ecomposition de l'espace
hermitien $(E,{\rm h})$ en somme orthogonale de $\R$-plans $I$-stables, et notons $T \subset {\rm U}_{\rm h}$
le tore associ\'e, constitu\'e des \'el\'ements de ${\rm Sp}_E$ stabilisant chacun de ces plans ; c'est un tore maximal anisotrope de $H$. Si $P \subset E$ est un $\R$-plan $I$-stable, on dispose d'une d\'ecomposition canonique $$P \otimes_\R \C = P^+ \oplus P^-$$
en somme de droites propres de $I$ pour les valeurs propres respectives de $+i$ et $-i$. Le choix d'une base $\Delta$ de $\Phi(H_\C,T_\C)$ \'equivaut \`a celui d'une
num\'erotation $P_1,\dots,P_g$ des $\R$-plans d\'efinissant $T$, ainsi que pour tout $i=1,\dots,g$ d'un signe $s_i \in \{+,-\}$ : le sous-groupe de Borel contenant $T$ associ\'e \`a une telle donn\'ee est le stabilisateur du drapeau form\'e sur $$P_1^{s_1}, P_2^{s_2}, \dots, P_g^{s_g}$$ (voir le \S VI.\ref{exdorad}). Un tel choix \'etant
fait, d\'ecrivons le $\R$-groupe $L_{\Delta,\varphi}$ associ\'e. Notons
$\eta_i \in {\rm X}^\ast(T_\C)$ le caract\`ere de $T_\C$ sur $P_i^{s_i}$, ainsi que $(\eta_i^\ast)$ la base duale des $(\eta_i)$ dans ${\rm X}_\ast(T_\C)$. D'apr\`es le \S VI.\ref{exdorad}, l'isomorphisme $i_\Delta : {\rm X}_\ast(T_\C) \simeq {\rm X}^\ast(S)$ envoie $\eta_i^\ast$ sur le caract\`ere de $S$ sur la droite $\C e_i$, et par d\'efinition 
$i_\Delta(\Phi^\vee(L_{\Delta,\varphi},T_\C))=\Phi(M,S)$. On en d\'eduit imm\'ediatement le $\C$-groupe $L_{\Delta,\varphi}$. Concr\`etement, il existe une unique d\'ecomposition $E = \oplus_{j \in J} E_j$, o\`u $E_j$ est la somme directe
des $\R$-plans $P_i$ pour les indices $i$ tels que $\C e_i \subset V_j$. En particulier, $\dim_\R E_j = \dim_\C V_j$ si $j \neq 0$ et
$\dim_\R E_0 = \dim_\C V_0 - 1$. De plus, si $j \neq 0$ il existe une d\'ecomposition en somme de lagrangiens transverses $$E_j \otimes \C = E_j^+ \oplus E_j^-$$ o\`u $E_j^+$ (resp. $E_j^-$) est la somme directe des $P_i^{s_i}$ (resp. $P_i^{-s_i}$) appartenant \`a $E_j \otimes \C$. Par d\'efinitions, le $\C$-groupe $L_{\Delta,\varphi}$ est le sous-groupe de ${\rm Sp}_{E \otimes \C}$
$${\rm Sp}_{E_0 \otimes \C} \times \prod_{j \neq 0} L_j,$$
o\`u $L_j \simeq {\rm GL}_{E_i^+}$ est le stabilisateur dans ${\rm Sp}_{E_i \otimes \C}$ des sous-espaces $E_j^+$ et $E_j^-$. Il ne reste qu'\`a expliciter sa structure r\'eelle. Soit $I' \in T(\R) \subset {\rm Sp}(E)$ l'\'el\'ement de carr\'e $-{\rm id}_{E}$ co\"incidant avec $s_i I$ sur $P_i$. Soient ${\rm h}'$ la forme hermitienne sur $E$ associ\'ee \`a $I'$, d\'efinie par ${\rm h}'(u)={\rm a}(I'u,u)$, et ${\rm h}'_j : E_j \rightarrow \R$ la restriction de ${\rm h}'$ \`a $E_j$; elle est de signature $(p_j,q_j)$ o\`u $q_j$ est le nombre des indices $i$ tels que $P_i^{-} \subset E_j^+$. L'\'el\'ement $I'$ induit un \'el\'ement central de $L_j$ pour tout $j \neq 0$ et donc $L_j = {\rm U}_{{\rm h}'_j}$. Au final, on a l'isomorphisme de $\R$-groupes 
$$L_{\Delta,\varphi} \simeq {\rm Sp}_{E_0} \times \prod_{j \neq 0} {\rm U}_{{\rm h}'_j}.$$
En particulier, $L_{\Delta,\varphi}(\R)$ est connexe.

\subsection{Param\'etrisation duale de $\Pi^c_{\rm
AJ}(\varphi)$}\label{parajs}

Soient $\varphi$ un param\`etre d'Adams-Johnson du $\R$-groupe d\'eploy\'e
$H$ et $c \in H^1(\R,H)$ une forme r\'eelle pure de $H$.  Suivant \cite[\S
3]{AJ} et \cite[\S 5]{arthurunipotent}, le paquet $\Pi_{\rm AJ}^c(\varphi)$
est encore muni d'une application de param\'etrisation naturelle $$\Pi_{\rm
AJ}^c(\varphi) \longrightarrow {\rm Hom}_{\rm groupes}({\rm
C}_\varphi,\C^\times), \pi \mapsto \chi_{O,\varphi}^c(\pi),$$ induite par
celle de Shelstad, que nous allons maintenant rappeler (voir aussi \cite[\S
8]{kottwitzannarbor} \cite[App.  A]{chrenard2} \cite[\S
4.2.2]{taibisiegel}).  Elle ne d\'epend, comme au \S\ref{paramdiscB}, que du
choix d'un $O \in \mathcal{B}(H)$, et que l'on fixe d\'esormais.  \ps\ps 

Soient $M \supset S \subset B \subset \widehat{H}$ associ\'es \`a $\varphi$ comme au \S\ref{paramaj}. En particulier, ${\rm C}_\varphi={\rm Z}(M)_2$. Fixons $\Delta \in O$, ce qui d\'efinit $\rho^\vee \in \frac{1}{2}{\rm X}_\ast(T)$ (\S\ref{paramdiscB}), un isomorphisme privil\'egi\'e $i_\Delta : {\rm X}_\ast(T) \isomo {\rm X}^\ast(S)$ (\S\ref{paramdiscC}), ainsi qu'un sous-groupe de Levi $L=L_{\Delta,\varphi} \subset H$ (\S\ref{paramaj}). 
L'inclusion $T \rightarrow L$ induit une bijection $${\rm W}(L_\C,T_\C) \backslash
H^1(\R,T) \isomo H^1(\R,L)$$ d'apr\`es Shelstad (voir aussi~\cite{borovoi14}). Cette bijection et le lemme \ref{stabwreel} (i) montrent que l'application $$f_\Delta : \Pi_{\rm AJ}^c(\varphi) \rightarrow H^1(\R,L),$$ envoyant $\pi_{w^{-1}\Delta,\varphi,t}$, pour $w \in W$ et $t \in T_2$ dans la classe $c$, sur l'image de l'\'el\'ement $w \star t \in H^1(\R,T)$ dans $H^1(\R,L)$, est bien d\'efinie. De plus, elle est injective, d'image \'egale \`a la fibre de l'application naturelle $H^1(\R,L) \rightarrow H^1(\R,H)$ au dessus de la classe $c$.  On dispose \'egalement d'un diagramme commutatif
$$\xymatrix{ H^1(\R,L) \ar@{->}^{\hspace{-1.5cm} g_\Delta}[rr] & & {\rm Hom}_{\rm groupes}({\rm Z}(M)_2,\C^\times) \\ 
H^1(\R,T)  \ar@{->}^{h_\Delta} [rr] \ar@{->}^{\rm can}[u] &  & {\rm X}^\ast(S)\otimes \Z/2\Z  \ar@{->}_{\rm can}[u]}$$
Les applications verticales sont les applications \'evidentes (celle de droite \'etant notamment induite par l'inclusion ${\rm Z}(M)_2 \subset S$), l'application $h_\Delta$ est la compos\'ee de $i_\Delta \otimes \Z/2\Z$ par l'isomorphisme canonique $H^1(\R,T) \isomo {\rm X}_\ast(T) \otimes \Z/2\Z$ rappel\'e au \S\ref{paramdiscC}, et l'application $g$ (un cas particulier des constructions g\'en\'erales de Kottwitz \cite{kottwitz}) est l'unique application faisant commuter ce diagramme. Concr\`etement, si $x \in H^1(\R,L)$, et si $t \in T_2$ est dans la classe $x$, que l'on \'ecrit $t=e^{i \pi \mu}$ avec $\mu \in {\rm X}_\ast(T)$, alors $g(x)$ est la restriction \`a ${\rm Z}(M)_2$ du caract\`ere $i_\Delta(\mu) \in {\rm X}^\ast(S)$ (elle ne d\'epend ni du choix de $t$ dans la classe $x$, ni bien s\^ur de celui de $\mu$).
Au final, on pose $$\chi_{O,\varphi}^c(\pi) \overset{{\mathrm{d\acute{e}f}}.}{=} g_\Delta \circ f_\Delta (\pi)$$ 
(il ne d\'epend que $O$ et non de $\Delta \in O$). Concr\`etement, pour tout $w \in W$ et $t=e^{i\pi \mu} \in T_2$ dans la classe $c$, le caract\`ere $\chi_{O,\varphi}^c(\pi_{w^{-1} \Delta,\varphi,t})$ est la restriction \`a ${\rm Z}(M)_2$ de $i_\Delta ( w(\mu+\rho^\vee)-\rho^\vee)$, o\`u $\rho^\vee$ est relatif \`a $\Delta$.\ps\ps

On prendra garde que comme l'ont remarqu\'e Adams et Johnson, $\pi \mapsto
\chi_{O,\varphi}^c(\pi)$ n'est pas injective en g\'en\'erale.  Terminons par quelques
observations simples mais importantes concernant le comportement des
s\'eries discr\`etes, suivant \cite[\S 8]{kottwitzannarbor} et \cite[App. A]{chrenard2}. Tout d'abord, nous avons d\'ej\`a dit que lorsque $\varphi$ est trivial sur
le facteur ${\rm SL}_2$, alors la construction d'Adams et Johnson redonne
$\Pi_{{\rm AJ}}^c(\varphi)=\Pi_{V_\varphi}$.  Plus pr\'ecis\'ement, dans ce
cas on a $M=S$, $L_{\Delta,\varphi}=T$ pour toute base $\Delta$ de $\Phi$, et la repr\'esentation $\pi_{\Delta,\varphi,t}$
co\"incide avec $\pi_{\Delta,V_\varphi}$. Il est alors \'evident que les deux
param\'etrisations $\chi_{O,\varphi}^c$ et $\chi_O^c$ co\"incident. Correctement
formul\'ee, cette propri\'et\'e s'\'etend aux s\'eries discr\`etes de $H_c(\R)$ intervenant dans
$\Pi^c_{\rm AJ}(\varphi)$ pour tout $\varphi$, comme l'a observ\'e Kottwitz
\cite[p. 196]{kottwitzannarbor}, rappelons comment. Pour cela, fixons un param\`etre d'Adams-Johnson $\varphi$ de $H$. Il lui est associ\'e une suite d'inclusions
\begin{equation} {\rm C}_\varphi={\rm Z}(M)_2 \subset S \subset
B \subset  \widehat{H} \end{equation}
d\'efinie au \S\ref{paramaj}. Suivant Kottwitz, observons qu'il existe un
param\`etre discret $\varphi_{\rm disc} : {\rm W}_\R \rightarrow
\widehat{H}(\C)$, unique \`a conjugaison pr\`es par $S(\C)$, ayant disons m\^eme caract\`ere infinit\'esimal que $V_\varphi$,
tel que $\varphi_{\rm disc}(\C^\times) \subset S(\C)$, et de caract\`ere
infinit\'esimal dominant relativement \`a $B$ (\S\ref{paramdiscC}). En
particulier, on dispose d'une inclusion canonique 
\begin{equation} \label{inclcomp} {\rm C}_\varphi={\rm Z}(M)_2 \rightarrow {\rm C}_{\varphi_{\rm
disc}}=S_2.
\end{equation}

\begin{prop}\label{propaj} Soient $\varphi$ un param\`etre d'Adams-Johnson de $H$ et $c \in
H^1(\R,H)$. Soient $T \subset H$ un tore maximal anisotrope, $O$ une $W_{\rm r}$-orbite dans $\mathcal{B}(T)$, $t_{\rm b} \in T_2$ l'\'el\'ement associ\'e \`a $O$, 
$t \in T_2$ dans la classe $c$, et $\Delta$ une base de $\Phi(H_\C,T_\C)$. \ps\ps 
\begin{itemize}
\item[(i)] La s\'erie discr\`ete $\pi_{\Delta,V_\varphi,t}$ appartient \`a $\Pi_{\rm AJ}^c(\varphi)$ si, et seulement si, $t t_{\rm
b} \in {\rm Z}(L_{\Delta,\varphi})$, auquel cas $\pi_{\Delta,V_\varphi,t} \simeq
\pi_{\Delta,\varphi,t}$.\ps\ps 
\item[(ii)] Si $\pi \in \Pi_{\rm AJ}^c(\varphi)$ est une s\'erie discr\`ete, alors $\chi_{O,\varphi}^c(\pi)$ est la
restriction de $\chi^c_O(\pi)$ \`a ${\rm C}_\varphi$ via l'homomorphisme canonique \eqref{inclcomp}.
\end{itemize}
\end{prop} 

\begin{pf} \cite[Lemma A.3 \& A.5]{chrenard2} et \cite[p.
196]{kottwitzannarbor}.
\end{pf}

\begin{example} \label{exemplecompactfinal} {\rm Soit $c \in H^1(\R,H)$
telle que $H_c(\R)$ est compact.  Alors $\Pi_{\rm AJ}^c(\varphi) =
\{V_\varphi\}$ pour tout param\`etre d'Adams-Johnson $\varphi$.  De plus,
d'apr\`es le lemme~\ref{corcompact} et la proposition~\ref{propaj}, $\delta^c_{V_\varphi}$ est la
restriction \`a ${\rm C}_\varphi = S_2 \rightarrow S$ de
$\rho^\vee + \mu \in {\rm X}^\ast(S)$.  Ici, $S$ d\'esigne le tore maximal de
$\widehat{H}$ associ\'e \`a $\varphi$, $\rho^\vee$ d\'esigne la demi-somme
des racines positives de $(\widehat{H},S)$ relativement \`a l'unique
sous-groupe de Borel rendant le caract\`ere infinit\'esimal de $\varphi$
dominant (\S\ref{paramaj}), et l'image de $\mu$ dans ${\rm X}^\ast(S) \otimes
\Z/2\Z$ est d\'efinie {\it loc.  cit}.  De mani\`ere alternative, elle
appartient au noyau de l'application naturelle $ {\rm X}^\ast(S) \otimes \Z/2\Z \longrightarrow {\rm X}^\ast(S_{\rm sc}) \otimes \Z/2\Z$, o\`u le tore ${\rm S}_{\rm sc}$ d\'esigne
l'image inverse de $S$ dans le rev\^etement universel de $\widehat{H}$
\cite[Lemma A.6]{chrenard2}.  }\end{example}

\begin{example} \label{exempleholaj} {\rm Soient $H={\rm Sp}_{2g}(\R)$, $V
\in {\rm Irr}(H_\C)$ et $\pi$ la s\'erie discr\`ete holomorphe ou
anti-holomorphe de $H(\R)$ ayant m\^eme caract\`ere infinit\'esimal que $V$
(\S \ref{paramdiscD}).  Soient $T \subset H$ un tore maximal anisotrope et
$K \subset H(\R)$ le sous-groupe compact maximal contenant $T(\R)$, choisis
par exemple comme au \S \ref{paramdiscD}.  Nous avons vu dans ce paragraphe
que $\pi \simeq \pi_{\pm \Delta',V_\varphi,1}$ o\`u $\Delta'$ est la base
donn\'ee {\it loc.  cit.}.  Soit $\varphi : {\rm SL}_2(\C) \times {\rm W}_\R
\rightarrow \widehat{H}(\C)$ un param\`etre d'Adams-Johnson tel que
$V_\varphi \simeq V$.  La proposition \ref{propaj} montre donc que $\pi \in
\Pi_{\rm AJ}(\varphi)$ si, et seulement si, $t_{\rm b} \in {\rm Z}(L_{\pm
\Delta',\varphi})$.  D'apr\`es le lemme \ref{lemmelankott} (i), il est
encore \'equivalent \`a demander que $L_{\pm \Delta',\varphi}(\R) \subset
K$.  Mais $L_{\pm \Delta',\varphi}$ est d\'ecrit par la recette du
\S\ref{exsp2gaj} ; sachant que la base $\Delta'$ a par construction la
propri\'et\'e que tous les $s_i$ sont de m\^eme signe, on constate que
$L_{\pm \Delta',\varphi}(\R) \subset K$ si, et seulement si, $\dim E_0 =
d_0-1 = 0$ dans les notations {\it loc.  cit.} Au final, $\pi \in \Pi_{\rm
AJ}(\varphi)$ si, et seulement si, le seul constituant de dimension impaire
de la repr\'esentation semi-simple ${\rm St} \circ \varphi$ est de dimension
$1$ (on a red\'emontr\'e le \cite[Lemma 9.4]{chrenard2}).  } \end{example}

\subsection{Lien avec les paquets d'Arthur}\label{parconjajenonce}

La conjecture suivante fait partie du folklore \cite[\S 5]{arthurunipotent}
\cite[p.  43]{arthur}.

\begin{conj}\label{conjaj} Soient $H$ le $\R$-groupe d\'eploy\'e ${\rm Sp}_{2g}$ ou ${\rm
SO}_{r+1,r}$, ainsi que $\varphi$ une classe de $\widehat{H}(\C)$-conjugaison de
param\`etres d'Adams-Johnson de $H$. Si $(\Pi(\varphi),\iota,\chi)$
d\'esigne le triplet qui lui est associ\'e par Arthur {\rm \cite[Thm. 1.5.1]{arthur}}, suivant les notations
du \S\ref{parpaqarch}, alors : \ps\ps
\begin{itemize}
\item[(a)] $\iota : \Pi(\varphi) \rightarrow \Pi_{\rm unit}(H(\R))$ est une
injection d'image $\Pi_{\rm AJ}(\varphi)$. \ps\ps 
\item[(b)] le caract\`ere $\chi \circ \iota^{-1} : \Pi_{\rm AJ}(\varphi)
\rightarrow {\rm Hom}_{\rm groupes}({\rm C}_{\varphi},\C^\times)$, bien d\'efini d'apr\`es le (a), 
co\"incide avec le caract\`ere $\pi \mapsto \chi_{O,\varphi}^1(\pi)$ d\'efini au \S\ref{parajs}.\end{itemize}\ps\ps 
\noindent Dans le (b) ci-dessus,  $O \in \mathcal{B}(H)$ correspond \`a la donn\'ee de Whittaker choisie par Arthur {\rm \cite[p. 55]{arthur}}. 
\end{conj}

\begin{remarque}  \label{etudiantmoeglin}{\rm  Comme nous l'avons d\'ej\`a dit au \S \ref{paramdiscC}, cette conjecture est
connue si $\varphi$ est trivial sur le facteur ${\rm SL}_2(\C)$ par les
travaux de Shelstad et Mezo. En g\'en\'eral, des progr\`es concernant cette conjecture ont  r\'ecemment \'et\'e obtenus par Colette Moeglin 
et Nicolas Arancibia (th\`ese en pr\'eparation). Plus pr\'ecis\'ement, Arancibia annonce la d\'emonstration de la conjecture \ref{conjaj} dans le cas particulier o\`u chaque constituant irr\'eductible de la
repr\'esentation ${\rm St} \circ \varphi$, que l'on peut \'ecrire sous la
forme $U \otimes V$ avec $U$ {\rm (}resp.  $V${\rm )} repr\'esentation
irr\'eductible de ${\rm SL}_2(\C)$ {\rm (}resp.  ${\rm W}_\R${\rm )},
satisfait $\dim \,V \,=1$ ou $\dim \,U \,\leq 4$. }
\end{remarque}

%

\begin{remarque} {\rm La conjecture \ref{conjaj} (et le th\'eor\`eme ci-dessus) admettent une variante pour le
$\R$-groupe $H={\rm SO}_{r,r}$ avec $r \equiv 0 \bmod 2$ dans laquelle
$\Pi_{\rm unit}(H(\R))$ est remplac\'e par $\widetilde{\Pi}_{\rm
unit}(H(\R))$ (\S \ref{parpaqarch}), et $\Pi_{\rm AJ}(\varphi)$ par son
image dans $\widetilde{\Pi}_{\rm unit}(H(\R))$ (on dispose d'un analogue
pour $\Pi_{\rm AJ}(\varphi)$ du lemme \ref{lemmethetacar}). 
}
\end{remarque}

Ainsi que l'annonce Arthur \cite[Chap. 9]{arthur}, sa description de
$\Pi_{\rm disc}(G)$ et la formule de multiplicit\'e \'enonc\'ee au \S
\ref{parpaqarch} pour les $\Z$-groupes $G$ de Chevalley admet un analogue aux
autres $\Z$-groupes classiques, \`a savoir les $\Z$-groupes sp\'eciaux
orthogonaux ${\rm SO}_L$ introduits au \S \ref{enonceparamst}.  Fixons donc un
tel $L$ et posons $G={\rm SO}_L$ ; on suppose de plus que $G(\R)$ a des s\'eries discr\`etes,
{\it i.e.} que $L \otimes \R$
est de signature paire s'il est de rang pair (le cas important pour la suite est $G={\rm
SO}_n$). Notons \'egalement $G^\star$ le groupe (special orthogonal) de Chevalley tel que $G_\C
\simeq G^\star_\C$. \ps\ps 

Dans la g\'en\'eralit\'e consid\'er\'ee par Arthur,
l'\'enonc\'e m\^eme de sa formule d\'epend d'une donn\'ee auxiliaire, \`a
savoir d'une ``r\'ealisation'' de $G$ comme forme int\'erieure de $G^\star$.
Le cas particulier du $\Z$-groupe $G={\rm SO}_L$ est particuli\`erement agr\'eable
car $G$ peut \^etre construit comme forme int\'erieure pure
``sur $\Z$'' de $G^\star$ ; disons simplement que cela vient de ce que
$\langle \pm 1 \rangle \otimes L$ est localement isomorphe pour la topologie \'etale sur ${\rm
Spec}(\Z)$ \`a ${\rm H}(\Z^r)$ si le rang de $L$ est pair, et \`a ${\rm H}(\Z^r) \oplus {\rm A}_1$
sinon (\S
II.2, Appendice \ref{appendixb}). Cette propri\'et\'e semble sensiblement
simplificatrice pour les questions de normalisation des facteurs de
transferts \'evoqu\'ees par Arthur {\it loc.  cit.} et \'etudi\'ees par
Kaletha \cite{kaletha}. \ps\ps 
Plus concr\`etement, appelons {\it r\'ealisation r\'eelle} de $G$ la donn\'ee d'un couple
$\xi=(c,f)$ o\`u $c \in {\rm H}^1(\R,G^\star_\R)$ et $f : G_\R \isomo
(G^\star_\R)_c$ est un isomorphisme. Il existe toujours des r\'ealisations
r\'eelles de $G$ : on l'a
vu au cours de l'exemple trait\'e \`a la fin du \S \ref{paramdiscE} quand $L$ est de rang impair,
le cas du rang pair \'etant analogue (il utilise que la signature de $L \otimes
\R$ est paire). Une telle r\'ealisation \'etant donn\'ee, il y a un sens \`a
d\'efinir pour tout param\`etre d'Adams-Johnson $\varphi$ de $G_\R$, ainsi
que $O \in \mathcal{B}(G_\R^\star)$, un couple $$(\Pi^\xi_{\rm
AJ}(\varphi),\chi_{O,\varphi}^\xi),$$
avec $\Pi^\xi_{\rm AJ}(\varphi) \subset \Pi_{\rm unit}(G(\R))$ et
$\chi_{O,\varphi}^\xi : \Pi_{\rm AJ}(\varphi) \rightarrow {\rm
Hom}_{\rm groupes}({\rm C}_\varphi,\C^\times)$, en transportant simplement
le couple $(\Pi_{\rm
AJ}^c(\varphi),\chi_{O,\varphi}^c)$ par la bijection $\Pi_{\rm unit}(G(\R)) \isomo \Pi_{\rm unit}((G^\star_\R)_c(\R))$
induite par $f$. \ps\ps 

Un second ph\'enom\`ene int\'eressant est que le groupe de similitudes projectives ${\rm P}\widetilde{G}(\Z)$ rencontre
toutes les composantes connexes de ${\rm P}\widetilde{G}(\R)$ (pour le voir, utiliser par exemple le th\'eor\`eme II.2.7). Comme cela a
d\'ej\`a \'et\'e observ\'e dans \cite{chrenard2}, cela sugg\`ere que 
la formule finale d'Arthur doit est compl\`etement canonique. Tout comme au
\S\ref{parpaqarch}, elle fait intervenir certains ensembles de
repr\'esentations unitaires de $G(\R)$ (d\'ependant d'un choix de r\'ealisation
r\'eelle de
$G$) dont l'existence, ainsi qu'une caract\'erisation, est annonc\'ee par
Arthur. On s'attend encore \`a ce que ceux qui contiennent des s\'eries discr\`etes
de $G(\R)$ soient exactement les paquets d'Adams-Johnson. Afin, de ne pas
multiplier les \'enonc\'es, et comme nous avons d\'ej\`a trait\'e en d\'etail le cas des groupes de Chevalley, nous ne formulerons pr\'ecis\'ement que la conjecture
finale attendue, \'etant entendu qu'il s'agit de la concat\'enation de deux
\'enonc\'es. Si $U,U' \in \Pi_{\rm unit}(G(\R))$, nous noterons $U \sim
U'$ si $U \simeq U'$, ou si $G(\R)$ est un groupe sp\'ecial orthogonal pair
et $U$ et $U'$ sont conjugu\'es ext\'erieurs l'un de l'autre par le groupe
orthogonal pair r\'eel correspondant ; si $X \subset \Pi_{\rm unit}(G(\R))$
on note enfin $\widetilde{X}=\{U \in \Pi_{\rm unit}(G(\R))\, |\,  \exists U' \in
X, U \sim U'\}$.

\begin{conj}\label{conjaj2} Soient $G$ un $\Z$-groupe classique, ${\rm St} :
\widehat{G} \rightarrow {\rm SL}_n$ la repr\'esentation standard, $\psi \in
\mathcal{X}_{\rm AL}({\rm SL}_n)$ et $U \in \Pi_{\rm unit}(G(\R))$ une
s\'erie discr\`ete telle que ${\rm St}({\rm Inf}_U)=\psi_\infty$. \ps\ps 
	Soient $\xi$ une r\'ealisation r\'eelle de $G$, $E \subset \Pi(G)$ l'ensemble des repr\'esentations $\pi$ telles
que $\pi_\infty \sim U$ et $\psi(\pi,{\rm St})=\psi$ (c'est un singleton si
$G$ n'est pas un groupe sp\'ecial orthogonal pair), ainsi que $\nu_\infty :
{\rm SL}_2(\C) \times {\rm W}_\R \longrightarrow \widehat{G}$ un param\`etre
d'Adams-Johnson associ\'e \`a $\psi$ par la recette d\'ecrite au
\S\ref{parpaqarch}. Alors $E \cap \Pi_{\rm
disc}(G) = \emptyset$, sauf si $U \in \Pi_{\rm AJ}^\xi(\nu_\infty)$ et 
\begin{equation}\label{criteramfgen} (\chi^\xi_{O,\nu_\infty}(U))_{|{\rm
C}_\psi}=\varepsilon_\psi, \end{equation} auquel cas 
$\sum_{\pi \in E} {\rm m}(\pi) = m_\psi$ o\`u l'entier $m_\psi$ vaut $1$, sauf si $G$ est sp\'ecial orthogonal pair et $\psi = \oplus_i
\pi_i[d_i]$ avec $d_i \equiv 0 \bmod 2$ pour tout $i$, auquel cas
$m_\psi=2$.
\end{conj}

L'assertion de canonicit\'e mentionn\'ee plus haut signifie pr\'ecis\'ement
que bien que le caract\`ere $\chi_{O,\nu_\infty}^\xi$ d\'epende des choix de
$O$ et de $\xi$, sa restriction \`a ${\rm C}_\psi \subset
{\rm C}_{\nu_\infty}$ n'en d\'epend pas. En effet, c'est \'evident si
$G$ est sp\'ecial orthogonal impair, car tout $\R$-automorphisme de $G_\R$
est int\'erieur et $|\mathcal{B}(G_\R)|=1$. Si $G \simeq G^\ast \simeq {\rm Sp}_{2g}$, le
groupe ${\rm Aut}_\R(G_\R)/{\rm Int}(G(\R)) = \Z/2\Z$ agit simplement
transitivement sur $|\mathcal{B}(G_\R)|=2$, de sorte qu'il y a exactement
deux choix \`a consid\'erer ; comme cela avait d\'ej\`a \'et\'e observ\'e
dans \cite[Lemma 9.5 \& 9.6]{chrenard2}, le crit\`ere \eqref{criteramfgen} est
en fait le m\^eme dans les deux cas (cela sera manifeste dans la d\'emonstration du
th\'eor\`eme~\ref{amfexplsp}). La situation est similaire si 
$G$ est sp\'ecial orthogonal pair, pour lequel il faut de plus tenir compte de l'automorphisme ext\'erieur provenant du groupe orthogonal correspondant (il est d'ailleurs d\'efini sur $\Z$ d\`es que $L$ poss\`ede une racine, i.e. un \'el\'ement $\alpha$ tel que ${\rm q}(\alpha)=1$).

\section{Formules de multiplicit\'e explicites : un formulaire}\label{formulairesAMF}

Si l'on confronte le th\'eor\`eme g\'en\'eral d'Arthur Thm. \ref{amfchev}
(ou la conjecture~\ref{conjaj2}) aux
consid\'erations et exemples du paragraphe \ref{paramdisc}, on obtient des formes
explicites de la formule de multiplicit\'e d'Arthur. Dans cette partie, nous
nous proposons de les d\'ecrire, \`a la mani\`ere de \cite[\S 3.29]{chrenard2}, dans les cas
particuliers d'importance pour ce m\'emoire, \`a savoir quand $G={\rm Sp}_{2g}$ et
pour une composante archim\'edienne $\pi_\infty$ dans la s\'erie discr\`ete holomorphe,
ou quand $G={\rm SO}_n$ pour $n \equiv -1,0,1 \bmod 8$. 

\subsection{Formule explicite pour ${\rm Sp}_{2g}$}\label{formexplicitesp}

Dans tout ce paragraphe, $g$ d\'esigne un entier $\geq 1$. Soit
$$\psi=\oplus_{i=1}^k \pi_i[d_i] \in  \mathcal{X}_{\rm AL}({\rm
SL}_{2g+1}),$$ 
 o\`u $k\geq 1$ est un entier et $\pi_i \in \Pi_{\rm
cusp}^\bot(\PGL_{n_i})$ et $d_i \geq 1$ pour tout $i=1,\dots,k$. On suppose que $\psi_\infty$
satisfait la condition {\rm (H2)} relativement \`a ${\rm Sp}_{2g}$
(\S\ref{amfhypgen}),  ce qui signifie que $\psi_\infty$ admet pour valeurs propres $2g+1$
entiers distincts
$$w_1>\cdots>w_g >0 >-w_g>\cdots>-w_1.$$
(\S \ref{paralgreg} cas I.). D'apr\`es le
lemme~\ref{cormodulo4} (i), il existe un unique entier $i_0 \in
\{1,\dots,k\}$ tel que $n_{i_0}d_{i_0}$ est impair. Quitte \`a
r\'eindexer les $\pi_i$, on peut supposer que $i_0=k$ sans perte de
g\'en\'eralit\'e. Observons que d'apr\`es ce m\^eme lemme, on a \'egalement
$n_i \equiv 0 \bmod 2$ et $n_i d_i \equiv 0 \bmod 4$ pour tout $i\neq k$, et
$n_k d_k \equiv 2g+1 \bmod 4$. \ps\ps 

Consid\'erons l'homomorphisme de groupes multiplicatifs $\chi : \{\pm 1\}^{k-1} \rightarrow \{\pm 1\}$
d\'efini de la mani\`ere suivante. Fixons $1 \leq i \leq k-1$ et notons $s_i \in \{\pm 1\}^{k-1}$
l'\'el\'ement d\'efini par $(s_i)_j=-1 \Leftrightarrow j=i$. Il y a deux cas
:\ps\ps
\begin{itemize}

\item[(i)] Si $d_i \equiv 0 \bmod 2$, on pose
$\chi(s_i)=(-1)^{\frac{n_id_i}{4}}$. \ps\ps 

\item[(ii)] Si $d_i \equiv 1 \bmod 2$, on pose $\chi(s_i)=(-1)^{|K_i|}$ o\`u
$K_i$ est l'ensemble des indices $1
\leq j \leq g$ qui sont impairs et tels que $w_j \in {\rm Poids}(\pi_i)$. 
\end{itemize}

\begin{thmv}\label{amfexplsp} Soit $\psi = \oplus_{i=1}^k \pi_i[d_i] \in
\mathcal{X}_{\rm AL}({\rm SL}_{2g+1})$ tel que $\psi_\infty$ admet $2g+1$
valeurs propres entieres et distinctes, et pour lequel $\pi_k \in \Pi_{\rm
cusp}(\PGL_{n_k})$ est telle que $d_k n_k \equiv 1 \bmod 2$.  Soit $\pi \in
\Pi({\rm Sp}_{2g})$ l'unique repr\'esentation telle que $\psi(\pi,{\rm
St})=\psi$ et $\pi_\infty$ est s\'erie discr\`ete holomorphe.  Supposons la
conjecture~\ref{conjaj} satisfaite pour ${\rm Sp}_{2g}$ et le morphisme
$\nu_\infty$ associ\'e \`a $\psi$ d\'efini au \S \ref{parpaqarch} (par
exemple, $d_i=1$ pour tout $i=1,\dots,k$).  \ps\ps 

	Alors $\pi \in \Pi_{\rm disc}({\rm Sp}_{2g})$ si, et seulement si, les
deux conditions 
suivantes sont satisfaites :\ps \ps
\begin{itemize}
\item[(a)] $d_k=1$, \ps\ps 
\item[(b)] pour tout $i=1,\dots,k-1$, on a $\chi(s_i)=\underset{1\leq j
\leq k, j \neq i}{\prod} \varepsilon(\pi_i\times\pi_j)^{{\rm
Min}(d_i,d_j)}$,\ps\ps 
\end{itemize}
\noindent o\`u $\chi$ d\'esigne le caract\`ere de $\{\pm 1\}^{k-1}$ associ\'e \`a $\psi$
qui est d\'efini ci-dessus. 
Enfin, si ces conditions sont satisfaites alors ${\rm m}(\pi)=1$.

\end{thmv}

\begin{pf} Reprenons quelques constructions du \S \ref{parpaqarch}.
On choisit notamment $\nu : {\rm SL}_2(\C) \times \prod_{i=1}^k
\widehat{G^{\pi_i}}(\C) \rightarrow {\rm SO}_{2g+1}(\C)$ associ\'e \`a
$\psi$ comme dans ce paragraphe, ainsi que pour chaque $i=1,\dots,k$, un homomorphisme $\mu_i : {\rm
W}_\R \rightarrow \widehat{G^{\pi_i}}(\C)$ tel que ${\rm St} \circ \mu_i
\simeq {\rm L}((\pi_i)_\infty)$. Le groupe ${\rm C}_\nu$ d\'efini {\it loc.
cit.} s'identifie naturellement avec $\{\pm 1\}^{k-1}$, avec leurs
\'el\'ements distingu\'es $s_1,\dots,s_{k-1}$. On dispose de plus d'un homomorphisme
$$\nu_\infty : {\rm SL}_2(\C) \times {\rm W}_\R \longrightarrow {\rm
SO}_{2g+1}(\C)$$ d\'eduit de $\nu$ et des $\mu_i$.
La condition {\rm (H2)} et les lemmes \ref{ajclassique} et \ref{lemmepsir}
assurent que c'est un param\`etre d'Adams-Johnson de ${\rm Sp}_{2g}$ (et m\^eme un
param\`etre de Langlands discret si $d_i=1$ pour tout $i$). Cela donne un
sens \`a l'hypoth\`ese, affirmant que la conjecture~\ref{conjaj} est satisfaite
pour ${\rm Sp}_{2g}(\R)$ et $\nu_\infty$. \ps\ps 

D'apr\`es
l'exemple~\ref{exempleholaj}, la s\'erie discr\`ete holomorphe $\pi_{\rm
hol}$ de ${\rm
Sp}_{2g}(\R)$ de caract\`ere infinit\'esimal $z$ tel que ${\rm
St}(z)=\psi_\infty$ est dans $\Pi(\nu_\infty)$ si, et seulement si, ${\rm St}
\circ \nu_\infty$ ne contient pas de repr\'esentation de la forme ${\rm
Sym}^{d-1}{\rm St}_2 \otimes \chi$ avec $d>1$ et $\chi \in
\{1,\epsilon_{\C/\R}\}$. Il est \'equivalent de demander que $d_k=1$. En effet, si $i=1,\dots,k$ 
on constate que ${\rm c}_\infty(\pi_i)$ admet la valeur propre $0$ si, et seulement si,
$n_i$ est impair, {\rm i.e.} $i=k$ (\S\ref{paralgreg}). D'apr\`es le
th\'eor\`eme~\ref{amfchev} et la conjecture~\ref{conjaj} pour le couple $({\rm
Sp}_{2g},\nu_\infty)$, il ne reste donc qu'\`a d\'emontrer que 
$$\chi_{O,\nu_\infty}^1(\pi_{\rm hol})_{|{\rm C}_\nu} = \chi.$$

Nous allons pour cela sp\'ecifier les constructions du \S\ref{paramaj}
relatives \`a $\nu_\infty$, \`a
la mani\`ere de l'\'etude faite dans l'exemple \ref{exsp2gaj}. On
consid\`ere l'homomorphisme $\widetilde{\nu_\infty} : {\rm W}_\R \rightarrow {\rm
SO}_{2g+1}(\C)$ d\'eduit en composant $\nu_\infty$ et le morphisme d'Arthur
comme au \S \ref{paramaj}. L'analyse du \S \ref{exsp2gaj} montre qu'il
existe un unique couple $(V^+,V^-)$ de lagrangiens transverses du ${\rm q}$-module
$V=\C^{2g+1}$ tels que : \begin{itemize} \ps\ps 

\item[--] $V^+$ et $V^-$ sont stables par $\widetilde{\nu_\infty}(\C^\times)$, \ps\ps 
\item[--] $V^+$ poss\`ede une $\C$-base $e_1,\dots,e_g$ v\'erifiant pour tout
$j=1,\dots,g$ et tout $z \in \C^\times$, la relation
$\widetilde{\nu_\infty}(z)(e_j)=z^{w_j}\overline{z}^{w'_j}e_j$. \ps\ps 
\end{itemize}

Dans cette relation $w_1,\dots,w_g$ sont les entiers d\'eduits de
$\psi_\infty$ d\'efinis avant
l'\'enonc\'e du th\'eor\`eme, et les $w'_i$ sont \'egalement des entiers uniquement
d\'etermin\'es (qu'il ne sera pas n\'ecessaire de pr\'eciser). La donn\'ee
de cette $\C$-base lagrangienne $(e_j)_{1 \leq j \leq g}$ d\'etermine
comme au \S \ref{exsp2gaj} un unique tore maximal $S \subset {\rm
SO}_{2g+1}$ ainsi qu'unique sous-groupe de Borel de $B$ contenant $S$.
\'Ecrivons comme au \S VI.\ref{exdorad} $${\rm X}^\ast(S)=\oplus_{j=1}^g \Z
\varepsilon_j,$$ $\varepsilon_j$ d\'esignant le caract\`ere de $S$ sur $\C
e_j$. Soit $\varepsilon_j^\ast \in {\rm X}_\ast(S)$ la base duale de
$(\varepsilon_j)$. Soit $\lambda \in \frac{1}{2}{\rm X}_\ast(S)$ est l'\'el\'ement associ\'e \`a
$\widetilde{\nu_\infty}$ intervenant dans la d\'efinition {\rm (AJ1)}
du~\S\ref{paramaj}. On a manifestement 
$$\lambda = \sum_{j=1}^g w_j \varepsilon_j^\ast,$$
de sorte que $\lambda$ est dominant relativement \`a $B$ d'apr\`es
les in\'egalit\'es $w_1>w_2>\dots>w_g>0$. Nous avons
explicit\'e jusqu'ici la suite d'inclusions 
$${\rm C}_\nu \subset {\rm C}_{\nu_\infty} \subset S \subset B \subset {\rm
SO}_{2g+1}$$
associ\'ee \`a $\nu$ et $\nu_\infty$. La proposition~\ref{propaj} (ii) ainsi que l'exemple \S\ref{paramdiscD}
entra\^inent que le caract\`ere $\chi^1_{O,\nu_\infty}(\pi_{\rm hol}) : {\rm
C}_{\nu_\infty} \rightarrow \C^\times$ est la restriction \`a ${\rm
C}_{\nu_\infty}$ de l'un des deux \'el\'ements suivants de ${\rm
X}^\ast(S)$ :
$$ \chi_{0}=\sum_{j \equiv 0 \bmod 2} \varepsilon_j \, \, \, \, \, 
{\rm ou}\, \, \,  \, \, \chi_1=\sum_{j \equiv 1 \bmod 2} \varepsilon_j,$$
les deux sommes portant sur tous les $j \in \{1,\dots,g\}$ de parit\'e
indiqu\'ee. \ps\ps 

V\'erifions enfin que la restriction \`a ${\rm C}_\nu$ de n'importe lequel
de ces deux caract\`eres $\chi_u$, $u \in \{0,1\}$, co\"incide avec le caract\`ere $\chi$. 
Fixons $1 \leq i \leq k-1$. Soit $J_i$ le sous-ensemble de $\{1,\dots,g\}$
constitu\'e des entiers $j$ tels que $w_j$ est valeur propre de
$(\pi_i[d_i])_\infty$. Observons tout d'abord que par construction,
l'image de l'\'el\'ement $s_i \in {\rm C}_\nu=\{\pm 1\}^{k-1}$ par
l'inclusion naturelle ${\rm C}_\nu \subset S$ est d\'etermin\'ee par la
relation suivante, satisfaite pour tout $j \in \{1,\dots,g\}$ :
$$\varepsilon_j(s_i)=-1 \, \, \Leftrightarrow
\, \, j \in J_i.$$
\'Ecrivons $J_i = J_i^{\rm 0} \coprod J_i^1$ o\`u
$J_i^u=\{j \in J_i, j \equiv u \bmod 2\}$. Par d\'efinition, 
$$\chi_u(s_i)=(-1)^{|J_i^u|}.$$
Observons que si $d_i$ est pair, auquel cas $\pi_i$ est
symplectique et $n_i$ est \'egalement pair, alors $J_i$ est r\'eunion
disjointe de $\frac{|J_i|}{2}=\frac{n_id_i}{2}$ paires d'entiers
cons\'ecutifs, de sorte que 
$$\chi_0(s_i)=\chi_1(s_i)=(-1)^{\frac{n_id_i}{4}}=\chi(s_i).$$
Supposons enfin $d_i$ impair, de sorte que $\pi_i$ est orthogonale et $n_i \equiv
0 \bmod 4$ d'apr\`es le corollaire~\ref{cormodulo4}. Soit $P_i \subset J_i$ le
sous-ensemble des $j \in J_i$ tels que $w_j \in {\rm Poids}(\pi_i)$. Alors
$J_i$ est la r\'eunion disjointe de $P_i$ et de ses translat\'es par $$\pm 1,
\pm 2, \dots, \pm \frac{d_i-1}{2}.$$ Mais si $1 \leq d \leq
\frac{d_i-1}{2}$, on a n\'ecessairement $w_i+d=w_{i-d}$ et $w_i-d=w_{i+d}$.
Les indices $i-d$ et $i+d$ \'etant congrus modulo $2$, on en d\'eduit 
$$\chi_u(s_i)=(-1)^{|P_i^u|}$$
o\`u $P_i^u=\{j \in P_i, j \equiv u \bmod 2\}$. On a donc
$\chi_1(s_i)=\chi(s_i)$, et aussi $\chi_0(s_i)=\chi_1(s_i)$ car
$|P_i|=\frac{n_i}{2} \equiv 0 \bmod 2$. 
\end{pf}


\begin{example}\label{verifikeda} {\rm En guise d'exemple, comparons
l'\'enonc\'e du th\'eor\`eme \ref{amfexplsp} au th\'eor\`eme
VII.\ref{thmikeda2} d'Ikeda.  Soit $k>0$ un entier pair et soit $\pi \in
\Pi_{\rm cusp}(\PGL_2)$ la repr\'esentation engendr\'ee par une forme propre
dans ${\rm S}_k({\rm SL}_2(\Z))$.  Consid\'erons le param\`etre $\psi =
\pi[g] \oplus [1]$.  Comme ${\rm Poids}(\pi)=\{\pm \frac{k-1}{2}\}$, on
constate que $\psi$ satisfait la condition {\rm (H2)} si, et seulement si,
$g \equiv 0 \bmod 2$ et $k>g$, auquel cas les valeurs propres de
$\psi_\infty$ sont
$$\frac{k-1}{2},\frac{k-3}{2},\cdots,\frac{k-g}{2},0,-\frac{k-g}{2},\cdots,-\frac{k-3}{2},-\frac{k-1}{2}.$$
On a bien s\^ur $k=2$, $d_2=1$, et ${\rm C}_\psi=\{\pm 1\}$ est engendr\'e
par l'\'el\'ement $s_1$.  On a \'egalement
$$\varepsilon_\psi(s_1)=\varepsilon(\pi \times
1)=\varepsilon(\pi)=i^k=(-1)^{k/2}.$$ De plus, on constate que $g$ \'etant
pair, $\chi(s_1)=(-1)^{g/2}$.  Sous la conjecture~\ref{conjaj}, la condition
n\'ecessaire et suffisante d'existence de $\pi' \in \Pi_{\rm cusp}({\rm
Sp}_{2g})$ telle que $\pi'_\infty$ est s\'erie discr\`ete holomorphe et
$\psi(\pi',{\rm St})=\psi$ s'\'ecrit donc $$(-1)^{k/2}=(-1)^{g/2}$$ soit $k \equiv g \bmod
4$.  C'est bien la condition intervenant dans l'\'enonc\'e d'Ikeda. Le r\'esultat d'Ikeda est en fait plus fort, d'abord
car il est inconditionnel, mais aussi car il n'a pas besoin de supposer
$k>g$ (et il serait int\'eressant d'\'etudier \'egalement les paquets d'Arthur
correspondant \`a ce cas plus g\'en\'eral). Mentionnons \'egalement que dans 
son suppl\'ement \cite{ikeda3} \`a \cite{ikeda2}, Ikeda d\'emontre que
si $k \equiv g \bmod 4$ (resp. $k \not \equiv g \bmod 4$) alors ${\rm m}(\pi')=1$
(resp. ${\rm m}(\pi')=0$) : voir \cite[Thm. 7.1, \S 15]{ikeda3}. 
}\end{example}

\ps \medskip

Terminons ce paragraphe par une traduction dans le langage classique de l'assertion de multiplicit\'e $1$ dans le th\'eor\`eme \ref{amfexplsp}.

\begin{corv}\label{mult1classsp} Soit $W$ la $\C$-repr\'esentation de
$\GL_g$ de plus haut poids $\sum_{i=1}^g m_i \varepsilon_i$ avec $m_1> m_2 >
\cdots > m_g > g+1$ (\S \ref{carinf}).  Si $F,G \in {\rm S}_W({\rm
Sp}_{2g}(\Z))$ sont deux formes propres pour ${\rm H}({\rm Sp}_{2g})$, et si
tout \'el\'ement de ${\rm H}({\rm Sp}_{2g})$ admet m\^eme valeur propre sur
$F$ et $G$, alors $F$ et $G$ sont proportionnelles.  Lorsque $g=2$, la
m\^eme assertion vaut en supposant seulement $m_1>m_2>2$.  \end{corv}

\begin{pf} Rappelons que l'on dispose d'un isomorphisme de ${\rm H}^{\rm
opp}({\rm Sp}_{2g})$-modules ${\rm S}_W({\rm Sp}_{2g}(\Z)) \isomo
\mathcal{A}_{\pi'_W}({\rm Sp}_{2g})$ (\S VI.\ref{corapiw} et la remarque qui
suit).  Il suffit donc de voir que si $\pi$ d\'esigne la repr\'esentation
engendr\'ee par une forme propre de ${\rm S}_W({\rm Sp}_{2g}(\Z))$ sous
l'action de ${\rm H}({\rm Sp}_{2g})$, alors ${\rm m}(\pi)=1$.  On rappelle
que $\pi_\infty \simeq \pi'_W$ et ${\rm St}({\rm Inf}_{\pi'_W})$ admet pour
valeurs propres $0$ et les $\pm (m_r -r)$ pour $r=1,\dots,g$ (Proposition \S
VI.\ref{infcarsiegel},): elles sont donc deux-\`a-deux distinctes et non
cons\'ecutives par hypoth\`ese. Si l'on \'ecrit $\psi(\pi,{\rm
St})=\oplus_i \pi_i[d_i]$, ce qui est loisible d'apr\`es le th\'eor\`eme
\ref{arthurst}, on a donc $d_i=1$ pour tout $i$, et on conclut par le
th\'eor\`eme~\ref{amfexplsp}.  Lorsque $g=2$ et $m_1 > m_2=2$, on conclut
encore car l'unique autre possibilit\'e $\psi(\pi,{\rm St})=\pi_1 \oplus
[3]$ avec $\pi_1 \in \Pi_{\rm cusp}(\PGL_2)$ est impossible, puisqu'un tel
$\pi_1$ est n\'ecessairement symplectique (\S \ref{altsoarth} ou proposition
\ref{contraintesalg} (i)).  \end{pf}

On s'attend \`a ce que le corollaire soit valable pour tout $W$ : cela d\'ecoule de la conjecture~\ref{conjaj} d\`es que $m_g>g$ !

\subsection{Formule explicite pour ${\rm SO}_n$ avec $n \equiv \pm 1 \bmod
8$.}\label{formexpliciteson1}

Supposons maintenant que $n$ est un entier $\equiv \pm 1 \bmod 8$ et
consid\'erons le $\Z$-groupe ${\rm SO}_n$.  Soit $$\psi=\oplus_{i=1}^k
\pi_i[d_i] \in \mathcal{X}_{\rm AL}({\rm SL}_{n-1}),$$ 
o\`u $k\geq 1$ est un entier et $\pi_i \in \Pi_{\rm
cusp}^\bot(\PGL_{n_i})$ et $d_i \geq 1$ pour tout $i=1,\dots,k$.  On suppose que
$\psi_\infty$ satisfait la condition {\rm (H2)} relativement \`a ${\rm
SO}_n$ (\S\ref{amfhypgen}), ce qui signifie que $\psi_\infty$ admet pour
valeurs propres $n-1$ demi-entiers (non entiers) distincts
$$w_1>\cdots>w_{\frac{n-1}{2}}
>-w_{\frac{n-1}{2}}>\cdots>-w_1.$$ (\S \ref{paralgreg} cas II.).  D'apr\`es le
lemme~\ref{cormodulo4} (ii), pour tout $i=1,\dots,k$ on a $n_id_i \equiv 0
\bmod 2$.  \ps\ps 

Consid\'erons l'homomorphisme de groupes multiplicatifs $\chi : \{\pm 1\}^k
\rightarrow \{\pm 1\}$ d\'efini de la mani\`ere suivante.  Fixons $1 \leq i \leq k$ et
notons $s_i \in \{\pm 1\}^k$ l'\'el\'ement d\'efini par $(s_i)_j=-1
\Leftrightarrow j=i$. Il y a deux cas :  \ps\ps

\begin{itemize}

\item[(i)] L'entier $d_i$ est pair. Si $n_i$ est pair, on pose
$\chi(s_i)=(-1)^{\frac{n_id_i}{4}}$. Si $n_i$ est impair, on pose $\chi(s_i)=\epsilon_i
\cdot (-1)^{\frac{(n_i-1)d_i}{4}}$, o\`u $\epsilon_i=(-1)^{[\frac{3d_i}{4}]}$
vaut $-1$ si $\frac{d_i}{2} \equiv 1,2 \bmod 4$ et $1$ sinon. (Dans tous les cas, 
$\chi(s_i)=(-1)^{[\frac{3\,n_id_i}{4}]}$.)
\ps\ps 

\item[(ii)] L'entier $d_i$ est impair. On pose $\chi(s_i)=(-1)^{|K_i|}$ o\`u
$K_i$ est l'ensemble des indices $1
\leq j \leq \frac{n-1}{2}$ qui sont $\equiv \frac{n-1}{2} \bmod 2$ et tels que $w_j \in {\rm
Poids}(\pi_i)$.\ps\ps 
\end{itemize}

\begin{thm}\label{amfexplso1} {\rm (Cas $n \equiv \pm 1 \bmod 8$)} Soit
$\psi = \oplus_{i=1}^k \pi_i[d_i] \in \mathcal{X}_{\rm AL}({\rm SL}_{n-1})$
tel que $\psi_\infty$ admet $n-1$ valeurs propres distinctes et dans
$\frac{1}{2}\Z-\Z$.  Soit $\pi \in \Pi({\rm SO}_n)$ l'unique
repr\'esentation telle que $\psi(\pi,{\rm St})=\psi$.  Supposons la
conjecture~\ref{conjaj2} satisfaite pour ${\rm SO}_n$ et le morphisme
$\nu_\infty$ associ\'e \`a $\psi$ d\'efini au \S \ref{parpaqarch}. Alors
$\pi \in \Pi_{\rm disc}({\rm SO}_n)$ si, et seulement si, pour tout
$i=1,\dots,k$, on a $$\chi(s_i)=\underset{1\leq j \leq k, j \neq i}{\prod}
\varepsilon(\pi_i\times\pi_j)^{{\rm Min}(d_i,d_j)},$$
\noindent o\`u $\chi$ d\'esigne le caract\`ere de $\{\pm 1\}^k$ associ\'e
\`a $\psi$ qui est d\'efini ci-dessus.
Enfin, si ces conditions sont satisfaites alors ${\rm m}(\pi)=1$.
\end{thm}

\begin{pf} (voir \cite[\S 3.30.1]{chrenard2}) La d\'emonstration
est similaire \`a celle du
th\'eor\`eme~\ref{amfexplsp}, c'est pourquoi nous n'insisterons que sur les
diff\'erences avec celle-l\`a. Une analyse de $\nu_\infty$ similaire \`a
celle du \S \ref{exsp2gaj}, tenant compte de ce que le groupe dual est celle
fois-ci ${\rm Sp}_{n-1}$, conduit \`a une explicitation de la suite
d'inclusions canoniques 
$${\rm C}_\nu=\{\pm 1\}^k \subset {\rm C}_{\nu_\infty} \subset S \subset B \subset {\rm
Sp}_{n-1}.$$
On invoque alors le corollaire \ref{corcompact} \`a la place de l'exemple
\ref{paramdiscD} (et la conjecture \ref{conjaj2} \`a la place de
la conjecture \ref{conjaj} et du th\'eor\`eme \ref{amfchev}). Comme le centre
de ${\rm SO}_n$ est trivial, ce corollaire affirme que le caract\`ere 
$\chi^\xi_{O,\nu_\infty}(\pi_\infty)$ est la restriction \`a ${\rm
C}_{\nu_\infty}$ de la demi-somme $\rho^\vee$ des racines positives de $T$
relativement \`a $B$, qu'il suffit donc d'expliciter. Posons $n=2r+1$
et consid\'erons pour cela la donn\'ee radicielle bas\'ee standard de ${\rm
Sp}_{2r}$
rappel\'ee au \S VI.\ref{exdorad}. Dans les notations {\it loc. cit.}, on
constate que $$\rho^\vee = \sum_{i=1}^r (r-i+1) \varepsilon_i \equiv \varepsilon_r
+ \varepsilon_{r-2} + \cdots \bmod 2.$$
On conclut comme dans la d\'emonstration du th\'eor\`eme~\ref{amfexplsp} que
la restriction de ce caract\`ere \`a ${\rm C}_\psi$ est le caract\`ere
$\chi$ de l'\'enonc\'e. On prendra garde au cas o\`u l'entier
$i \in \{1,\dots,k\}$ v\'erifie $d_i \equiv 0 \mod 2$ et $n_i \equiv 1 \bmod
2$, car alors $0$ est un poids de $\pi_i$. Dans ce cas, l'ensemble $J_i$
d\'efini dans cette d\'emonstration est r\'eunion disjointe de
$\frac{(n_i-1)d_i}{2}$ paires d'entiers cons\'ecutifs, ainsi que de
l'ensemble $\{r,r-1,\dots,r+1-\frac{d_i}{2}\}$, d'o\`u la modification \`a
apporter \`a
la d\'efinition de $\chi$ dans le cas (i).
\end{pf}

\subsection{Formule explicite pour ${\rm SO}_n$ avec $n \equiv 0 \bmod
8$.}\label{formexpliciteson0}

Consid\'erons enfin le $\Z$-groupe ${\rm SO}_n$ pour $n \equiv 0 \bmod 8$.
Soit $$\psi=\oplus_{i=1}^k \pi_i[d_i] \in \mathcal{X}_{\rm AL}({\rm
SL}_{n}),$$ o\`u $k\geq 1$ est un entier
et $\pi_i \in \Pi_{\rm cusp}^\bot(\PGL_{n_i})$ et $d_i \geq 1$ pour tout
$i=1,\dots,k$.  On suppose que $\psi_\infty$ satisfait la condition {\rm
(H2)} relativement \`a ${\rm SO}_n$ (\S\ref{amfhypgen}), ce qui signifie que
$\psi_\infty$ admet pour valeurs propres $n$ entiers $$w_1>\cdots>w_{n/2}
 \geq -w_{n/2}>\cdots>-w_1$$ (\S \ref{paralgreg} cas III.) 
Notons $I_1 \subset \{1,\cdots,k\}$ le sous-ensemble des indices $i$ tels que
$n_i d_i \equiv 1 \bmod 2$ et posons $I_0=\{1,\dots,k\} - I_1$. D'apr\`es le
lemme~\ref{cormodulo4} (ii), on a soit $I_1 = \emptyset$, soit $|I_1|=2$ 
(ce qui ne peut se produire que si $w_{n/2}=0$). De plus, $n_id_i
\equiv 0 \bmod 4$ si $i \in I_0$. \ps\ps 

Consid\'erons l'homomorphisme de groupes multiplicatifs $\chi : \{\pm
1\}^{I_0} \rightarrow \{\pm 1\}$ d\'efini de la mani\`ere suivante.  Fixons $i \in I_0$ et   
notons $s_i \in \{\pm 1\}^{I_0}$ l'\'el\'ement d\'efini par $(s_i)_j=-1
\Leftrightarrow j=i$.  \ps\ps

\begin{itemize}

\item[(i)] Si $d_i \equiv 0 \bmod 2$, on pose $\chi(s_i)=(-1)^{\frac{n_id_i}{4}}$. \ps\ps 

\item[(ii)] Si $d_i \equiv 1 \bmod 2$, on pose $\chi(s_i)=(-1)^{|K_i|}$ o\`u
$K_i$ est l'ensemble des indices $1 \leq j \leq \frac{n}{2}$ qui sont
impairs et tels que $w_j \in {\rm Poids}(\pi_i)$.\ps\ps 
\end{itemize}

\begin{thm}\label{amfexplso0} {\rm (Cas $n \equiv 0 \bmod 8$)} Soient     
$\psi = \oplus_{i=1}^k \pi_i[d_i] \in \mathcal{X}_{\rm AL}({\rm SL}_n)$
tel que $\psi_\infty$ admet $n$ valeurs propres enti\`eres, au moins $n-1$
d'entre elles \'etant distinctes, ainsi que $\{1,\dots,k\}=I_0 \coprod I_1$
la partition associ\'ee ci-dessus. Soit 
$\Pi \subset \Pi({\rm SO}_n)$ le sous-ensemble des repr\'esentations $\pi$
telles que $\psi(\pi,{\rm St})=\psi$ ; c'est un singleton si $I_1 \neq \emptyset$.
Supposons la conjecture~\ref{conjaj2} satisfaite pour ${\rm SO}_n$ et le morphisme
$\nu_\infty$ associ\'e \`a $\psi$ d\'efini au \S \ref{parpaqarch}. 
Alors $\Pi \cap \Pi_{\rm disc}({\rm SO}_n) \neq \emptyset$ si, et seulement
si, on a 
\begin{equation}\label{critso8n} \chi(s_i)=\underset{1\leq j \leq k, j \neq
i}{\prod} \varepsilon(\pi_i\times\pi_j)^{{\rm Min}(d_i,d_j)}, \, \, \,
\forall i \in I_0,\end{equation}
o\`u $\chi$ d\'esigne le caract\`ere de $\{\pm 1\}^{I_0}$ associ\'e
\`a $\psi$ qui est d\'efini ci-dessus. Enfin, si cette condition est satisfaite alors $\sum_{\pi \in \Pi} {\rm
m}(\pi)=1$ si $I_1 \neq \emptyset$, et $\sum_{\pi \in \Pi} {\rm m}(\pi)=2$
sinon.  \end{thm}

\begin{pf} (voir \cite[\S 3.30.2]{chrenard2}) La d\'emonstration
est similaire \`a celle du th\'eor\`eme~\ref{amfexplso1}. On explicite \`a
la mani\`ere du \S \ref{exsp2gaj} la suite d'inclusions canoniques $${\rm C}_\nu
\subset {\rm C}_{\nu_\infty} \subset S \subset B \subset {\rm
SO}_n.$$
L'homomorphisme $\{\pm 1\}^{I_0} \rightarrow {\rm C}_\nu$ envoyant l'\'el\'ement
$s_i$ pour $i \in I_0$ d\'efini ci-dessus sur l'\'el\'ement du m\^eme nom
d\'efini au \S \ref{parepsilonpsi} induit une surjection $\{\pm 1\}^{I_0}
\rightarrow {\rm C}_\nu/{\rm Z}({\rm SO}_n)$. Dans la donn\'ee radicielle bas\'ee 
standard de $({\rm SO}_n,S,B)$ rappel\'ee au \S VI.\ref{exdorad}, on constate
cette fois-ci que la demi-somme des racines positives vaut 
$$\rho^\vee = \sum_{i=1}^{r} (r-i) \varepsilon_i \equiv \varepsilon_{r-1} +
\varepsilon_{r-3} + \cdots \bmod 2,$$
o\`u $r=\frac{n}{2}$. Comme le centre de ${\rm SO}_n$ est non trivial, le corollaire \ref{corcompact} 
affirme seulement que $\chi^\xi_{O,\nu_\infty}(\pi_\infty)$ est la
restriction \`a ${\rm C}_{\nu_\infty} \subset S$ de $\rho^\vee$ ou
$\rho^\vee+\nu$, o\`u 
$$\nu \equiv  \sum_{i=1}^{r} \varepsilon_i \bmod 2.$$
Mais $\nu_{|{\rm C}_\nu}=1$ car $n_i d_i \equiv 0 \bmod 4$ pour tout $i \in
I_0$, et le reste s'en d\'eduit. 
\end{pf}

\begin{remarque}\label{remmult2} {\rm (Sur la multiplicit\'e $2$)} {\rm Soit $\alpha \in {\rm E}_n$ tel que $\alpha \cdot
\alpha =2$ ainsi que ${\rm s}_\alpha \in {\rm O}({\rm E}_n)$ la sym\'etrie
orthogonale associ\'ee \`a cette racine. La conjugaison ext\'erieure par ${\rm s}_\alpha$
induit une involution, disons not\'ee $S$, de $\Pi({\rm SO}_n)$, dont nous 
avons d\'ej\`a \'etudi\'e certains aspects au \S IV.\ref{fauton}, ainsi que
dans les exemples du \S VI.\ref{isomsatake} et du \S VI.\ref{exconjal}. 
Ces paragraphes montrent que si $\pi \in \Pi({\rm SO}_n)$ alors ${\rm m}(\pi)={\rm
m}(S(\pi))$, et pour tout $v \in {\rm P}\cup \{\infty\}$ les classes de conjugaison 
${\rm c}_v(\pi)$ et ${\rm c}_v(S(\pi))$ sont images l'une de l'autre
par l'action de l'\'el\'ement non-trivial de ${\rm O}_n(\C)/{\rm SO}_n(\C)$.
En particulier, $S$ pr\'eserve $\Pi_{\rm disc}({\rm SO}_n)$ et $\psi(S(\pi),{\rm St})=\psi(\pi,{\rm
St})$. De plus, cela entra\^ine que $S(\pi)$ est isomorphe \`a $\pi$ si, et
seulement si, $\pm 1$ (resp. $0$) est valeur propre de ${\rm c}_p(\pi)$
(resp. de ${\rm c}_\infty(\pi)$) pour tout premier $p$. Pla\c{c}ons-nous maintenant 
dans les hypoth\`eses du th\'eor\`eme ci-dessus, et 
supposons $I_1=\emptyset$ ainsi que \eqref{critso8n} satisfaite. Les
observations pr\'ec\'edentes montrent que l'ensemble $\Pi$ de l'\'enonc\'e est
stable par $S$. En particulier, si $S$ n'admet aucun
point fixe dans $\Pi$ le th\'eor\`eme assure que $\Pi \cap \Pi_{\rm disc}({\rm SO}_n)$ est
constitu\'e deux repr\'esentations \'echang\'ees par $S$, chacune \'etant de
multiplicit\'e $1$. Voici deux cas particuliers o\`u ceci s'applique. \ps\ps 

{\rm I.} $0$ n'est pas valeur propre de $\psi_\infty$ : $S(\pi)_\infty \neq \pi_\infty$ pour tout $\pi \in \Pi$. \ps\ps

{\rm II.} $\pi_i$ est symplectique pour $i=1,\dots k$ : $S(\pi)_p
\neq \pi_p$ pour tout $\pi \in \Pi$ et $p$ premier. En effet, 
les valeurs propres $\lambda$ de ${\rm
c}_p(\pi_i)$ satisfont $p^{-\frac{1}{2}} < |\lambda| < p^{\frac{1}{2}}$
d'apr\`es Jacquet-Shalika (et m\^eme $|\lambda|=1$ par la conjecture de Ramanujan
\S VI.\ref{ramanujan}), alors que les valeurs propres $\mu$ de $[d_i]_p$
satisfont $|\mu|\geq p^{\frac{1}{2}}$ ou $|\mu| \leq p^{-\frac{1}{2}}$ (car
$d_i \equiv 0 \bmod 2$), de
sorte que $|\lambda\mu| \neq 1$ et $\pm 1$ n'est pas valeur propre de
$\psi_p$. \ps

\noindent Dans le cas g\'en\'eral, une combinaison de ces id\'ees montre que
pour qu'il existe une repr\'esentation de multiplicit\'e $>1$ dans un
$\Pi_{\rm disc}({\rm SO}_n)$ il faut et suffit qu'il existe $\pi \in
\Pi_{\rm cusp}(\PGL_{4m})$ autoduale, alg\'ebrique, orthogonale, telle que $|{\rm
Poids}(\pi)|=4m-1$ et ${\rm c}_p(\pi)$ admet la valeur propre $\pm 1$ pour
tout $p$. On pourrait montrer qu'il n'en existe pas quand $m=1$.

}
\end{remarque} 

\ps \medskip
\noindent {\sc Deux crit\`eres}
\ps \medskip

\noindent Donnons deux crit\`eres pour que la relation~\eqref{critso8n} soit satisfaite. 

\begin{critere}\label{critamfsonA} On suppose $\psi$ de la
forme $\left( \bigoplus_{i=1}^{k-2}\pi_i[d_i] \right) \oplus [d_{k-1}]
\oplus [1]$ avec $d_i \equiv 0 \bmod 2$ pour tout $i=1,\dots,k-2$ et $d_{k-1}>1$
impair. Alors
la relation~\eqref{critso8n} est satisfaite si, 
et seulement si, pour tout
$i=1,\dots,k-2$ on a $$(-1)^{\frac{n_id_i}{4}}=\varepsilon(\pi_i)^{1+{\rm
Min}(d_i,d_{k-1})}.$$
\end{critere}

\begin{pf} C'est une application imm\'ediate des formules et du fait que
$\varepsilon(\pi_i \times \pi_j)=1$ si $i,j \in  I_0=\{1,\dots,k-2\}$. Ce dernier point peut soit
\^etre vu comme cas particulier du r\'esultat g\'en\'eral d'Arthur car les
repr\'esentations $\pi_i$ et $\pi_j$ sont symplectiques (\S \ref{altsoarth}), soit se d\'emontrer
directement car $\varepsilon({\rm I}_a \otimes {\rm I}_b)=(-1)^{1+{\rm
Max}(a,b)}=1$ si $a$ et $b$ sont impairs (\S \ref{parfaceps}). 
\end{pf}

\begin{critere}\label{critamfsonB} On suppose $\psi$ de la forme $\left( \bigoplus_{i=1}^{k-2}\pi_i[d_i] \right) \oplus \pi_{k-1}\oplus
[d_k]$ avec : \begin{itemize} \ps\ps 
\item[(i)] $d_i \equiv 0 \bmod 2$ pour tout $i=1,\dots,k-2$, \ps\ps 
\item[(ii)] et $\pi_{k-1} \in \Pi_{\rm cusp}(\PGL_3)$ telle que ${\rm
w}(\pi_{k-1})> \underset{1 \leq i \leq k-2}{{\rm Max}} \, \,{\rm w}(\pi_i)$.\ps\ps 
\end{itemize}
La relation~\ref{critso8n} est satisfaite si, et seulement si, pour tout
$i=1,\dots,k-2$ 
$$(-1)^{\frac{n_i}{2}(\frac{d_i}{2}-1)}=\varepsilon(\pi_i)^{1+{\rm Min}(d_i,d_k)}.$$
\end{critere}

\begin{pf} Rappelons que $\varepsilon({\rm I}_w \otimes {\rm I}_{w'})=(-1)^{1+{\rm Max}(w,w')}$
et ${\rm I}_w \otimes \epsilon_{\C/\R} \simeq {\rm I}_w$. Comme ${\rm w}(\pi_i)$ est
impair si $i<k-1$ et pair si $i=k-1$, l'hypoth\`ese sur ${\rm w}(\pi_{k-1})$ assure que
pour tout $i \leq k-2$,  $$\varepsilon(\pi_i \times
\pi_{k-1})=\varepsilon({\rm L}(\pi_i)_\infty \otimes {\rm I}_{{\rm
w}(\pi_{k-1})})\varepsilon({\rm L}(\pi_i)_\infty \otimes
\epsilon_{\C/\R})=(-1)^{\frac{n_i}{2}}\varepsilon(\pi_i).$$ On conclut car 
$\varepsilon(\pi_i \times \pi_j)=1$ si $i,j \leq k-2$. 
\end{pf}

\ps \medskip
\noindent {\sc Exemples }
\ps \medskip

Admettant la conjecture~\ref{conjaj2}, le th\'eor\`eme~\ref{amfexplso0}
devient un outil formidable pour v\'erifier (et peut-\^etre, dans un futur proche,
red\'emontrer !) les r\'esultats pr\'ec\'edemment
obtenus dans ce m\'emoire. Son application la plus imm\'ediate concerne le
cas o\`u $\psi$ est tel que soit $I_0 = \emptyset$, soit $k=1$, puisqu'alors
la condition~\eqref{critso8n} est trivialement satisfaite. On retrouve par exemple
la conclusion du th\'eor\`eme VII.\ref{arthurtrialite} (ii), ainsi que celle
du (i) si l'on demande simplement que la repr\'esentation $\pi'$ soit dans
$\Pi_{\rm disc}({\rm SO}_8)$ plut\^ot que dans $\Pi_{\rm disc}({\rm
O}_8)$. On retrouve aussi l'assertion concernant $\Delta_{11}[12]$ du corollaire 
VII.\ref{cor24ikedabocherer} (voir le \S IX.\ref{preuve1thm24}). \ps\ps

En guise d'autre exemple, consid\'erons le cas o\`u $$\psi = \pi_1[d_1] \oplus [d_2] \oplus [1]$$
avec $\pi_1 \in \Pi_{\rm cusp}(\PGL_2)$ de poids $\{ \pm \frac{k-1}{2}
\}$, $k$ \'etant un entier pair $>0$. On a $\varepsilon(\pi_1)=(-1)^{k/2}$.
D'apr\`es le crit\`ere \ref{critamfsonA}, la relation \eqref{critso8n} est satisfaite si, et seulement
si :
\ps\ps  
\begin{itemize}
\item[(I)] soit $d_1 < d_2$ et $d_1 \equiv k \bmod 4$, \ps\ps 
\item[(II)] soit $d_1 > d_2$  et $d_1 \equiv 0 \bmod 4$. \ps\ps 
\end{itemize}
\ps\ps 
\noindent Par exemple, si $\psi = \Delta_{11}[4] \oplus [7] \oplus [1]$,
auquel cas $n=16$ et $\psi_\infty={\rm St}({\rm Inf}_1)$, on est dans le
cas (I) ci-dessus. Cela red\'emontre imm\'ediatement l'assertion concernant $\psi(\pi,V_{\rm St})$ du
corollaire VII.\ref{corpij}, ainsi donc que le (i) du th\'eor\`eme V.\ref{vpnontrivdim16}. 
On retrouve aussi les assertions du corollaire~\ref{cor24ikedabocherer} pour
$k>12$ (voir IX.\S \ref{preuve1thm24}) : on est dans le cas (II) pour $k=16$
car $(d_1,d_2)=(8,7)$, et dans le cas (I) pour $k>16$ car $24-k \equiv k \bmod 4$
et $2k-25>24-k$. En ce qui concerne $\psi=\Delta_{17}[14]\oplus [3] \oplus
[1]$, on n'est ni dans le cas (I), ni dans le cas (II), ce qui corrobore
comme promis le corollaire VII.\ref{cex32}. \ps\ps 

Expliquons enfin la table du \S \ref{tableexempleso8}. Compte tenu de
l'analyse ci-dessus, il ne reste qu'\`a comprendre les cases bleues de cette table. 
Mais si $\psi = \Delta_w[2] \oplus {\rm Sym}^2 \Delta_{11} \oplus [1]$ avec $w <22={\rm w}({\rm Sym}^2 \Delta_{11})$,
le crit\`ere \ref{critamfsonB} montre que la relation \eqref{critso8n} est
satisfaite, ce qui conclut.

\section{Compatibilit\'e \`a la correspondance th\^eta}\label{comptheta}

Soient $n\equiv 0 \bmod 8$ et $g\geq 1$ des entiers tels que $n>2g$. Soit
$$\psi_S = \oplus_{i=1}^k \pi_i[d_i] \in \mathcal{X}_{\rm AL}({\rm
SL}_{2g+1}),$$
o\`u $k\geq 1$ et, pour tout $i=1,\dots,k$, $\pi_i \in \Pi_{\rm cusp}(\PGL_{n_i})$. On suppose que 
$$\psi_O := \psi_S \oplus [n-2g-1] \in \mathcal{X}_{\rm AL}({\rm
SL}_{n})$$
satisfait la condition {\rm (H2)} relativement \`a ${\rm SO}_n$. Il est
\'equivalent de demander que les valeurs propres de $(\psi_S)_\infty$ soient
$2g+1$ entiers $$w_1>\cdots>w_g>0>-w_g>\cdots>-w_1$$ avec de plus $w_g \geq
\frac{n}{2}-g$. Il sera commode de choisir l'indexation des $\pi_i$ de sorte que $n_kd_k \equiv 1
\bmod 2$ (\S\ref{formexplicitesp}). On suppose enfin $d_k=1$, c'est d'ailleurs automatique si $n \neq 2g+2$. \ps\ps 

Soit $\pi_O \in \Pi({\rm SO}_n)$ {\rm (resp. $\pi_S \in \Pi({\rm
Sp}_{2g})$)} l'unique repr\'esentation telle que $\psi(\pi_O,{\rm
St})=\psi_O$ {\rm (resp. $\psi(\pi_S,{\rm St})=\psi_S$)}. Notons ${\rm m}(\pi_O)$
et ${\rm m}(\pi_S)$ les multiplicit\'es respectives de $\pi_O$ et $\pi_S$
dans $\Pi_{\rm disc}({\rm SO_n})$ et $\pi_{\rm disc}({\rm Sp}_{2g}))$. Elles valent
chacune $0$ ou $1$ d'apr\`es la formule de multiplicit\'e Arthur, sous la conjecture \ref{conjaj2}. 

\begin{prop}\label{compmomsp} Admettons la conjecture \ref{conjaj2}. 
On a ${\rm m}(\pi_O)={\rm m}(\pi_S)$ si, et seulement si,
pour tout $i=1,\dots,k-1$ tel que $d_i \equiv 0 \bmod 2$ et $d_i \geq n-2g$,
$$\varepsilon(\pi_i)=1.$$
\end{prop}

\begin{pf} C'est imm\'ediat sur les formules explicites donn\'ees dans les
\S \ref{formexplicitesp} et \S \ref{formexpliciteson0}. On observe en effet que
l'injection naturelle 
$${\rm C}_{\psi_S} \rightarrow {\rm C}_{\psi_O} $$
induit un isomorphisme ${\rm C}_{\psi_S} \isomo {\rm C}_{\psi_O}/{\rm
Z}({\rm SO}_n)$. Le groupe ${\rm C}_{\psi_S}$ s'identifie naturellement \`a
$\{\pm 1\}^{k-1}$ (voir la d\'emonstration du th\'eor\`eme \ref{amfexplsp}),
et l'injection ci-dessus l'identifie (avec ses \'el\'ements privil\'egi\'es
$s_i$) au sous-groupe $\{ \pm 1\}^{I_0}$ d\'efini au
\S\ref{formexpliciteson0} (en particulier, ${\rm
I}_0=\{1,\dots,k-1\}$). Via cet identification, on constate que le caract\`ere $\chi$
de $\{\pm 1\}^{k-1}$ d\'efini au \S \ref{formexplicitesp} co\"incide avec le
caract\`ere $\chi$ de $\{\pm 1\}^{I_0}$ d\'efini au \S \ref{formexpliciteson0}.
D'apr\`es les th\'eor\`emes \ref{amfexplsp} et \ref{amfexplso0}, la condition pour que ${\rm m}(\pi_O)={\rm m}(\pi_S)$
\'equivaut donc \`a 
$$\varepsilon(\pi_i)^{{\rm Min}(d_i,n-2g-1)}=1$$
tout $i=1,\dots,k-1$. C'est automatique quand $d_i$ est impair car alors
$\pi_i$ est orthogonale.
\end{pf}

Lorsque ${\rm m}(\pi_S)=1$, le crit\`ere de B\"ocherer (remarque
VII.\ref{rembocherer}) donne une 
condition n\'ecessaire est suffisante, du moins si les entiers $w_i$ sont
cons\'ecutifs, pour que la forme propre de ${\rm S}_{w_1+1}({\rm
Sp}_{2g})$ (bien d\'efinie \`a un scalair pr\`es) engendrant $\pi_S$ admette un $\vartheta$-correspondant $\pi'$ dans
$\Pi_{\rm disc}({\rm O}_n)$ : il faut et suffit que ${\rm L}(s,\pi_S,{\rm
St})$ ne s'annule pas en $s=\frac{n}{2}-g$. Par restriction \`a ${\rm
SO}_n$ (\S IV.\ref{fauton}), l'existence de $\pi'$ entra\^ine que ${\rm m}(\pi_O)$ est non nul
(corollaire
VII.\ref{corparamrallis}). Il est donc important de v\'erifier que le crit\`ere de B\"ocherer
est compatible \`a la proposition ci-dessus, et c'est bien le cas : 

\begin{prop}\label{checkcomptheta} ${\rm L}(s,\pi_{\rm S},{\rm St})$ ne
s'annule pas en $s=\frac{n}{2}-g$ si, et seulement si, pour tout
$i=1,\dots,k-1$ tel que $d_i \geq n-2g$ et $d_i \equiv 0 \bmod 2$, on a
${\rm L}(\frac{1}{2},\pi_i) \neq 0$ (et donc $\varepsilon(\pi_i)=1$).  En
particulier, ${\rm L}(\frac{n}{2}-g,\pi_{\rm S},{\rm St})\neq 0$ d\`es que
$n>3g$.  \end{prop}

\begin{pf} La fonction ${\rm L}(s,\pi_{\rm S},{\rm St})$ est produit des ${\rm L}(s+j,\pi_i)$ pour $i=1,\dots,k$ et $j \in \frac{d_i-1}{2}+\Z$ tel que $|j| \leq \frac{d_i-1}{2}$. On rappelle que si $\pi_i \neq 1$, le produit eul\'erien de ${\rm L}(s,\pi_i)$ est absolument convergent pour ${\rm Re} \, s >1$, et que $\xi(s,\pi_i)=\Gamma(s,{\rm L}((\pi_i)_\infty)){\rm L}(s,\pi_i)$ admet un prolongement holomorphe sur $\C$ tel que $\xi(1-s,\pi_i)=\varepsilon(\pi_i)\xi(s,\pi)$ (\S VI.\ref{lienfonctionl}). De plus, on a ${\rm L}(1,\pi_i) \neq 0$ \cite{jasha0}. Par hypoth\`ese, si $\pi_i=1$ alors 
$i=k$, $d_k=1$, et ${\rm L}(s,\pi_i)$ est la fonction $\zeta$ de Riemann. Observons que si ce cas se produit alors $n-2g-1 \equiv 3 \bmod 4$ ; en particulier s'il on pose $$s_0=\frac{n}{2}-g$$ alors $s_0 \geq 2$ n'est ni un z\'ero, ni un p\^ole, de $\zeta$. Comme $d_k =1$ et $s_0 \leq 1$, on a ${\rm L}(s_0,\pi_{\rm S},{\rm St})=0$ si, et seulement si, il existe $1 \leq i \leq k-1$ et $j \in \frac{d_i-1}{2}+\Z$ avec $|j| \leq \frac{d_i-1}{2}$ tels que  ${\rm L}(s_0+j,\pi_i)=0$. \ps\ps 
Fixons $i<k$. La repr\'esentation $\pi_k$ \'etant la seule parmi les $\pi_s$ \`a poss\'eder le poids $0$, ${\rm L}((\pi_i)_\infty)$ est somme directe des ${\rm I}_w$ o\`u $w$ parcourt les poids $>0$ de $\pi_i$. La fonction $\Gamma(s)$ ne s'annule nulle part sur l'axe
r\'eel et admet pour seuls p\^oles les entiers $\leq 0$. La recette pour $\Gamma(s,{\rm L}((\pi_i)_\infty))$ (\S\ref{parfaceps}) et les propri\'et\'es rappel\'ees ci-dessus de $\xi(s,\pi_i)$ montrent donc que si ${\rm L}(s,\pi_i)=0$, disons pour $s \in \R$, alors soit $0<s<1$, soit $s \leq - w$, o\`u $w$ est le plus petit poids $>0$ de $\pi_i$. Mais $w-\frac{d_i-1}{2} \geq w_g \geq \frac{n}{2}-g = s_0$, de sorte que pour $j\geq \frac{1-d_i}{2}$ on a 
$$s_0+j > - s_0+\frac{1-d_i}{2}  \geq -w.$$
Enfin, $0 < s_0+j <1$ \'equivaut \`a $s_0+j=\frac{1}{2}$ si $j \in
\frac{1}{2}\Z$. On conclut la premi\`ere assertion car $\frac{1}{2}$ est de la
forme $\frac{n}{2}-g+j$ avec $\frac{1-d_i}{2} \leq j \leq \frac{1-d_i}{2}$
et $j \in \frac{d_i-1}{2} + \Z$ si, et seulement si, $d_i \equiv 0 \bmod 2$
et $n-2g \leq d_i$. La seconde provient de la relation \'evidente $d_i \leq g$. \end{pf} 

\begin{remarque}\label{rembochamelioration} {\rm Notons que la condition suffisante $n>3g$ obtenue ci-dessus, l\'egitime en pratique d'apr\`es le th\'eor\`eme$^{\color{green}\star}$ \ref{arthurst} d'Arthur, s'av\`ere meilleure que la condition g\'en\'erale $n>4g$ de B\"ocherer. } \end{remarque}

En th\'eorie, il pourrait exister des param\`etres $\psi_{\rm S}$ tels que
${\rm m}(\pi_O)={\rm m}(\pi_S)=1$ mais tels que $\pi_O$ et $\pi_S$ ne soient
pas $\vartheta$-correspondants. Pour fabriquer un tel exemple, il faudrait
trouver une repr\'esentation alg\'ebrique r\'eguli\`ere autoduale symplectique
$\varpi$ telle que ${\rm L}(\frac{1}{2},\varpi)=0$ mais $\varepsilon(\varpi)=1$. Les
auteurs ne connaissent aucun tel exemple (comparer avec la remarque VII.\ref{nonvanish}).
C'est un fait que nous exploitons d'ailleurs positivement \`a plusieurs reprises dans ce m\'emoire !\ps\ps 

\section{Compatibilit\'e \`a la fonction ${\rm L}$ de
B\"ocherer}\label{complboc}

Soient $g\geq 1$ un entier, $k \in \Z$, $F \in {\rm S}_k({\rm Sp}_{2g}(\Z))$
une forme propre, ainsi que $\pi_F \in \Pi_{\rm cusp}({\rm Sp}_{2g})$ la
repr\'esentation engendr\'ee par $F$. B\"ocherer a d\'emontr\'e dans \cite{boclfunction} que le produit eul\'erien
${\rm L}(s,\pi_F,{\rm St})$ (d\'efini au \S VI.\ref{lienfonctionl}) est absolument convergent si ${\rm Re} \,s > g+1$, et que la fonction 
$$\xi_{\rm B}(s,\pi_F,{\rm
St})\,\,{:=}\,\,\left(\Gamma(s,\epsilon_{\C/\R}^g)\prod_{i=1}^g\Gamma_\C(s+k-i)\right){\rm
L}(s,\pi_F,{\rm St}),$$ 
admet un prolongement
m\'eromorphe \`a $\C$ ainsi qu'une \'equation fonctionelle
$$\xi_{\rm B}(s,\pi_F,{\rm St})=\xi_{\rm B}(1-s,\pi_F,{\rm St})$$ (voir aussi
\cite{langlandseis, andkal, psr}). Rappelons
que l'on a $\Gamma(s,1)=\Gamma_\R(s)$ et
$\Gamma(s,\epsilon_{\C/\R})=\Gamma_\R(s+1)$ dans les notations du \S \ref{parfaceps}. Les p\^oles de $\xi_{\rm B}(s,\pi_F,{\rm St})$ ont \'et\'e \'etudi\'es par Mizumoto \cite[Theorem 1, \S 3 Corollary] {mizumoto}. Il d\'emontre les \'enonc\'es suivants, dans lesquels $\nu$ d\'esigne l'\'el\'ement de $\{0,1\}$ d\'efini par la congruence $\nu \equiv k \bmod 2$ : \ps\ps\begin{itemize}

\item[(a)] si $k\geq g$ alors $\xi_{\rm B}(s,\pi_F,{\rm St})$ admet au plus des p\^oles simples en $s=0$ et $s=1$, et elle est holomorphe ailleurs.  \ps\ps
\item[(b)] si $g>k$ alors $\xi_{\rm B}(s,\pi_F,{\rm St})$ admet au plus des p\^oles simples en $s=g-k+\nu+1,g-k+\nu,k-g-\nu+1,k-g-\nu$, au plus des p\^oles doubles aux entiers $s \in \Z$ tel que $g-k+\nu-1 \geq s \geq 2+k-g-\nu$, et elle est holomorphe ailleurs. \ps \ps \end{itemize}

\noindent Supposons maintenant que l'on ait la relation $$\psi(\pi_F,{\rm St})=\oplus_{i=1}^r \pi_i[d_i],$$ avec $\pi_i \in \Pi_{\rm
cusp}^\bot(\PGL_{n_i})$ et $d_i\geq 1$ pour tout $1 \leq i \leq r$, ce qui est loisible d'apr\`es Arthur (Th\'eor\`eme \ref{arthurstbis}). La th\'eorie des fonctions ${\rm L}$
standards des \'el\'ements de $\Pi_{\rm cusp}(\PGL_m)$, par Godement et Jacquet,
montre que la fonction d\'efinie par $\xi(s,\pi_i)=\Gamma(s,{\rm L}((\pi_i)_\infty)){\rm L}(s,\pi_i)$
(\S \ref{parfaceps}) poss\'ede un prolongement m\'eromorphe \`a $\C$ et
une \'equation fonctionnelle $s \mapsto 1-s$. Cela fournit donc une seconde 
fa\c{c}on naturelle de compl\'eter ${\rm L}(s,\pi_F,{\rm St})$ en posant simplement 
$$\xi_{\rm A}(s,\pi_F,{\rm St})\,\,:=\,\,\prod_{i=1}^r \prod_{j=0}^{d_i-1}
\xi(s+j-\frac{d_i-1}{2},\pi_i).$$
Cette fonction est \'egalement m\'eromorphe sur $\C$ et invariante par $s \mapsto
1-s$ (\`a un signe pr\`es \`a priori, mais qui est en fait \'egal \`a
$1$ car d'apr\`es Arthur \cite[Thm. 1.5.3 (b)]{arthur} on a $\varepsilon(\pi)=1$ pour toute repr\'esentation $\pi \in \Pi_{{\rm cusp}}(\PGL_m)$ autoduale et orthogonale). \ps\ps

Ainsi, $\frac{\xi_{\rm B}(s,\pi_F,{\rm
St})}{\xi_{\rm A}(s,\pi_F,{\rm St})}$ est un quotient ``explicite'' de produits
de facteurs
$\Gamma$. Lorsque $k>g+1$, il est facile de voir sur les recettes
respectives pour ces facteurs que ce quotient vaut $1$, {\rm
i.e.}  $\xi_{\rm B}(s,\pi_F,{\rm
St})=\xi_{\rm A}(s,\pi_F,{\rm St})$. La situation s'av\`ere plus int\'eressante
lorsque $k \leq g+1$, auquel cas la comparaison de ces facteurs, combin\'ee
aux propri\'et\'es rappel\'ees ci-dessus des p\^oles de $\xi_{\rm B}(s,\pi_F,{\rm
St})$ et $\xi_{\rm A}(s,\pi_F,{\rm St})$, admet des cons\'equences non
triviales sur $\psi(\pi_F,{\rm St})$. L'\'enonc\'e de la proposition \ref{cormizudisc} ci-dessous est sugg\'er\'e par le point (a) du
th\'eor\`eme \ref{amfexplsp} (et donc par la conjecture \ref{conjaj}) lorsque $k=g+1$. Il n\'ecessite l'introduction de quelques quantit\'es pr\'eliminaires. \ps\ps

Pour tout entier $a\geq 1$ on pose $$\delta(\pi_F,a)\,=\,{\rm ord}_{s=a} \prod_{\{i, \pi_i \neq
1\}}\prod_{j=0}^{d_i-1} \xi(s+j-\frac{d_i-1}{2},\pi_i).$$ 
\noindent On rappelle que si $\pi \in \Pi_{\rm cusp}({\rm PGL}_m)$ est telle que $\pi
\neq 1$ alors $\xi(s,\pi)$ est une fonction enti\`ere de $s$. De plus, si $s \in
\frac{1}{2}\Z$ satisfait $\xi(s,\pi)=0$ alors d'apr\`es Jacquet-Shalika on a $s= \frac{1}{2}$ et  ${\rm ord}_{s=\frac{1}{2}}\Gamma(s,{\rm L}(\pi_\infty)) = 0$. Pour tout entier $a \geq 1$ on a donc l'\'egalit\'e
\begin{equation}\label{autreexpdeltaa}\delta(\pi_F,a)\,:=\, \sum_{\{i \, \, |\, \, \,  d_i \equiv 0 \bmod 2, \, d_i \geq 2 a\}} {\rm ord}_{s=\frac{1}{2}} {\rm
L}(s,\pi_i).\end{equation}
\noindent 
En particulier, on a les in\'egalit\'es $0 \leq \delta(\pi_F,b) \leq \delta(\pi_F,a)$ pour $b \geq a \geq 1$. Il sera commode \'egalement de poser, pour tout entier $n\geq 0$ :
$$p_n(s)\,=\,\frac{\Gamma_\C(s+n)}{\Gamma_\C(s-n)} \, \, \, {\rm et} \, \, \,
\gamma_n(s)\,=\,\frac{\Gamma(s,\varepsilon_{\C/\R}^n)}
{\Gamma_\R(s)}\prod_{i=1}^n
\frac{\Gamma_\C(s+i)}{\Gamma_\R(s-i)\Gamma_\R(s+i)}.$$
\`A l'aide des formules $\Gamma_\C(s+1)\,=\,\frac{s}{2\pi}\,\Gamma_\C(s)$ et
$\Gamma_\C(s)=\Gamma_\R(s)\Gamma_\R(s+1)$, on v\'erifie que les \'egalit\'es suivantes sont valables pour tout entier $n\geq 0$ :
$$p_n(s)\,=\,(2\pi)^{-2n}\prod_{1-n \leq m \leq n}(s-m)\, \, \, \, {\rm et}\, \, \, \, \gamma_n(s)=\prod_{0 \leq 2m \leq n} p_{n-2m}(s).$$
En particulier, $p_n$ et $\gamma_n$ sont des polyn\^omes en $s$ invariants par $s \mapsto 1-s$.\ps\ps

\begin{propv}\label{cormizudisc} Soit $F \in {\rm S}_k({\rm Sp}_{2g}(\Z))$ une forme
propre avec $k=g$ ou $k=g+1$. Supposons $\psi(\pi_F,{\rm St})=\oplus_{i=1}^r
\pi_i[d_i]$ avec $\pi_r=1$ et $d_r>1$. Alors on a l'in\'egalit\'e
$\delta(\pi_F,\frac{d_r+1}{2}) > 0$.
\end{propv}

\begin{pf} Traitons d'abord le cas $k=g+1$. Comme $\psi(\pi_F,{\rm St})_\infty$ admet $0$ pour valeur propre
simple (ce sont les $2g+1$ entiers $n$ tels que $|n|\leq g$ d'apr\`es le \S VI.\ref{infcarsiegel}), on a $\pi_i \neq 1$ si
$i<r$, ainsi donc que l'\'egalit\'e 
\begin{equation} \label{inegpreuvecomizu}{\rm ord}_{s=\frac{d_r+1}{2}}\,\,\xi_{\rm A}(s,\pi_F,{\rm
St}) \,\,=\,\,
\delta(\pi,\frac{d_r+1}{2})- 1\end{equation} d'apr\`es les rappels ci-dessus. On constate aussi que si $i<r$ alors $\pi_i$ est
r\'eguli\`ere, chacun de ses poids \'etant $>\frac{d_i-1}{2} \geq 0$, on a donc par d\'efinitions les \'egalit\'es
$$\frac{\xi_{\rm B}(s,\pi_F,{\rm St})}{\xi_{\rm A}(s,\pi_F,{\rm St})} = \frac{\Gamma(s,\epsilon_{\C/\R}^g)}{\Gamma_\R(s)}
\prod_{i=1}^{\frac{d_r-1}{2}}
\frac{\Gamma_\C(s+i)}{\Gamma_\R(s+i)\Gamma_\R(s-i)} =
\frac{\Gamma(s,\epsilon_{\C/\R}^g)}{\Gamma(s,\epsilon_{\C/\R}^{\frac{d_r-1}{2}})}
\gamma_{\frac{d_r-1}{2}}(s).$$
Ce terme ne s'annule pas, et il est fini, en $s=\frac{d_r+1}{2}$. En supposant $\delta(\pi,\frac{d_r+1}{2})=0$, on a
donc ${\rm ord}_{s=\frac{d_r+1}{2}} \xi_{\rm B}(s,\pi_F,{\rm St})=-1$. Mais d'apr\`es Mizumoto, 
les seuls p\^oles \'eventuels de $\xi_{\rm B}(s,\pi_F,{\rm St})$ sont en $s=0$ ou $s=1$, ce qui montre
$d_r=1$. \ps
Le cas $k=g$ est similaire. Dans ce cas, $\psi(\pi_F,{\rm St})_\infty$ admet $0$ pour valeur propre triple et les entiers $\pm 1, \dots, \pm (g-1)$ pour valeurs propres simples. En particulier, il est possible que $\pi_i$ admette le poids $0$ pour $i<r$, mais dans ce cas on a $d_i=1$. L'\'egalit\'e \eqref{inegpreuvecomizu} est donc encore satisfaite. De plus, ce m\^eme argument montre 
$\frac{\xi_{\rm B}(s,\pi_F,{\rm St})}{\xi_{\rm A}(s,\pi_F,{\rm St})} = \mu(s) \, \gamma_{\frac{d_r-1}{2}}(s)$ avec $$\mu(s)=\frac{\Gamma_\C(s)\Gamma(s,\epsilon_{\C/\R}^g)}
{\Gamma(s,\epsilon_{\C/\R}^{\frac{d_r-1}{2}})\Gamma_\R(s+e_1)\Gamma_\R(s+e_2)},$$
pour certains \'el\'ements $e_1,e_2 \in \{0, 1\}$ qu'il sera inutile de pr\'eciser (voir la remarque \ref{remgammargammac}). La fonction $\mu(s)$ est finie et non nulle en $s=\frac{d_r+1}{2}$, et on conclut comme dans le cas $k=g+1$ par le r\'esultat de Mizumoto.
\end{pf}

\begin{remarque}\label{remgammargammac}{\rm  Soient $a,b \in \Z$. Observons que si la fonction m\'eromorphe $\Gamma_\R(s)^a\Gamma_\R(s+1)^b$ est invariante par $s \mapsto 1-s$, alors on a\,\,$a=b=0$. En effet, l'ordre d'annulation de cette fonction en $s\,=\,2,1,0$ et $-1$ vaut respectivement $0$, $0$, $-a$ et $-b$. L'invariance par $s \mapsto 1-s$ des fonctions $\xi_{\rm A}(s,\pi_F,{\rm St})$, $\xi_{\rm B}(s,\pi_F,{\rm St})$, ainsi que $\gamma_n(s)$ pour tout $n\geq 0$, permet donc de pr\'eciser la d\'emonstration pr\'ec\'edente. Tout d'abord, on en d\'eduit la congruence $\frac{d_r-1}{2} \equiv g \bmod 2$ pour $k=g+1$, que nous avions d\'ej\`a obtenue autrement au \S \ref{formexplicitesp}. Dans le cas $k=g$, on en d\'eduit l'\'egalit\'e des images dans $\Z/2\Z$ des ensembles $\{e_1,e_2\}$ et $\{\frac{d_r+1}{2},g\}$.  }\end{remarque}

La proposition suivante, qui sera utile au \S \ref{fsiegelpoids12}, sugg\`ere des propri\'et\'es cach\'ees du param\`etre standard de $\pi_F$ dans le cas $g \geq
k$, soit encore des paquets d'Arthur de ${\rm Sp}_{2g}(\R)$ contenant
les repr\'esentations de la forme $\pi'_W$ qui ne sont pas des s\'eries discr\`etes (\S
VI.\ref{carinf}). Il serait int\'eressant de les clarifier. La notation $[x]$ d\'esigne la partie enti\`ere du r\'eel $x$. \ps\ps

\begin{propv}\label{cormizunondisc} Soit $F \in {\rm S}_k({\rm Sp}_{2g}(\Z))$ une forme  
propre avec $g \geq k$. On suppose $\psi(\pi_F,{\rm
St})=\oplus_{i=1}^r \pi_i[d_i]$ et, pour tout $i=1,\dots,r$ tel que $\pi_i \neq 1$, la propri\'et\'e suivante satisfaite :
\begin{equation} \label{propcormizu} \forall \lambda \in {\rm Poids}(\pi_i), \, \, \, \, \, |\lambda|\geq \frac{d_i-1}{2}.\end{equation}
On a alors l'\'egalit\'e \begin{equation} \label{quotxiaxib}\frac{\xi_{\rm
B}(s,\pi_F,{\rm St})}{\xi_{\rm A}(s,\pi_F,{\rm St})} =
\left(\prod_{i=1}^{g-k} \frac{1}{p_i(s)}\right) \, \left( \prod_{i \in
I} \gamma_{\frac{d_i-1}{2}} \right).\end{equation} 
\noindent Si l'on pose $I=\{ 1\leq i \leq r, \pi_i=1\}$ et $I'=\{i \in I, d_i=1\}$, on a alors
\begin{equation} \label{incormizu} \delta(\pi_F,1) + \sum_{i \in I} [\frac{d_i+1}{4}] \geq g-k -2 + |I'|+2|I-I'|.\end{equation} \par

\noindent Si $I$ est non vide, on pose enfin $d \, =\, {\rm max}\{d_i, i \in I\}$ ;  on a l'in\'egalit\'e \,$d\,\geq \,2  |I-I'|+1$. D\'efinissons $\nu \in \{0,1\}$ par la congruence $\nu \equiv k \bmod 2$. Alors on est dans un des cas suivants : \begin{itemize}\ps\ps
\item[(a)] $|I|=3$, \hspace{4pt}$\frac{d+1}{2} \geq g-k + 4 -\delta(\pi_F,1)$,  \hspace{4pt}et l'une des deux in\'egalit\'es suivantes est v\'erifi\'ee :  \hspace{4pt}$\delta(\pi_F,1) > 2-\nu$ \hspace{4pt} ou  \hspace{4pt}$\delta(\pi_F,\frac{d+1}{2})>0$.\ps\ps
\item[(b)] $|I|=2$,  \hspace{4pt}$|I'|=0$,  \hspace{4pt}$\frac{d+1}{2} \geq g - k + 3 - \delta(\pi_F,1)$, \hspace{4pt}et l'une des deux in\'egalit\'es suivantes est v\'erifi\'ee :  \hspace{4pt}$\delta(\pi_F,1) > 1-\nu$ \hspace{4pt} ou  \hspace{4pt}$\delta(\pi_F,\frac{d+1}{2})>0$.\ps\ps
\item[(c)] $|I|=2$,  \hspace{4pt}$|I'|=1$, \hspace{4pt} $\frac{d+1}{4} \geq g - k + 1 - \delta(\pi_F,1)$,  \hspace{4pt}et l'une des deux in\'egalit\'es suivantes est v\'erifi\'ee :  \hspace{4pt}$\delta(\pi_F,1) > \frac{g-k-\nu}{2}$ \hspace{4pt} ou  \hspace{4pt}$\delta(\pi_F,\frac{d+1}{2})>0$.\ps\ps
\item[(d)] $|I|=1$,  \hspace{4pt}$|I'|=0$,  \hspace{4pt}$\frac{d+1}{4} \geq g - k  - \delta(\pi_F,1)$,  \hspace{4pt}et l'une des deux in\'egalit\'es suivantes est v\'erifi\'ee :  \hspace{4pt}$\delta(\pi_F,1) > \frac{g-k-\nu-2}{2}$ \hspace{4pt} ou  \hspace{4pt}$\delta(\pi_F,\frac{d+1}{2})>0$.\ps\ps
\item[(e)] $I'=I$ \hspace{4pt}  et  \hspace{4pt} $\delta(\pi_F,1) \geq g-k-2 + |I'|$.
\ps \end{itemize}

\end{propv}

\begin{pf} Par d\'efinition, on a la relation $${\rm ord}_{s=1} \,\,\xi_{\rm
A}(s,\pi_F,{\rm St})\,=\,-|I'|- 2 |I-I'| + \delta(\pi_F,1).$$ On rappelle que
$\psi(\pi_F,{\rm St})_\infty$ admet pour valeurs propres les $2g+1$ entiers
$0$ et $\pm (k-j)$, $j=1,\dots,g$. En particulier, $0$ en est valeur propre de
multiplicit\'e $3$ et les $\pi_i$ sont alg\'ebriques. De plus, si l'on a $\pi_i \neq 1$ et si $0$ est un poids de $\pi_i$, alors on a n\'ecessairement
$d_i=1$ d'apr\`es l'hypoth\`ese \eqref{propcormizu}. Ainsi, l'hypoth\`ese \eqref{propcormizu} et la recette pour les
${\rm L}(s,(\pi_i)_\infty)$ (\S \ref{contraintesalg}, \S \ref{parfaceps}) entra\^inent qu'il existe $e_1,e_2,e_3
\in \{0,1\}$ tels que si l'on pose
$\mu(s)=\frac{\Gamma_\C(s)\Gamma(s,\epsilon_{\C/\R}^g)}
{\Gamma_\R(s+e_1)\Gamma_\R(s+e_2)\Gamma_\R(s+e_3)}$, on a
$$\frac{\xi_{\rm
B}(s,\pi_F,{\rm St})}{\xi_{\rm A}(s,\pi_F,{\rm St})} =
\left(\prod_{i=1}^{g-k} \frac{1}{p_i(s)}\right) \, \left( \prod_{i \in
I} \gamma_{\frac{d_i-1}{2}} \right) \, \mu(s).$$ 
La remarque \ref{remgammargammac}, et l'invariance des $p_n(s)$ par $s \mapsto 1-s$, montrent $\mu(s)=1$, d'o\`u la formule \eqref{quotxiaxib}.\ps

 Comme on a ${\rm ord}_{s=1} \,p_n(s) =
1$ si $n\geq 1$, et $p_0=1$, on a aussi ${\rm ord}_{s=1} \gamma_n(s) =
[\frac{n+1}{2}]$. On a donc $$ {\rm ord}_{s=1} \frac{\xi_{\rm B}(s,\pi_F,{\rm
St})}{\xi_{\rm A}(s,\pi_F,{\rm St})} \,= \,- (g-k) \,+\, \sum_{i \in I}\,\,
[\frac{d_i+1}{4}]$$ 
d'apr\`es la relation \eqref{quotxiaxib}. L'in\'egalit\'e \eqref{incormizu} de l'\'enonc\'e d\'ecoule alors de
l'in\'egalit\'e ${\rm ord}_{s=1} \xi_{\rm B}(s,\pi_F,{\rm St})
\geq -2$ de Mizumoto.  \ps\ps Supposons $|I| \neq 0$ et v\'erifions la
derni\`ere assertion.  D'apr\`es le th\'eor\`eme \ref{arthurstbis} d'Arthur, les $d_i$ avec
$i \in I$ sont deux \`a deux distincts. En particulier, on a $|I'|\leq 1$. Comme ${\rm St}\, {\rm c}_\infty(\pi_F)$ admet $\pm 1$ pour valeur propre de multiplicit\'e $\leq 2$, on a aussi $|I-I'|\leq 2$ (et donc $|I|\leq 3$). Les $d_i$ \'etant \'evidemment impairs, on a enfin $d \geq 2|I-I'|+1$. \ps

Supposons d'abord $|I|=3$, auquel cas on a $|I'|=1$.
L'in\'egalit\'e  \eqref{incormizu}, associ\'ee \`a l'in\'egalit\'e \'evidente $\sum_{i \in I} [\frac{d_i+1}{4}] \leq \sum_{i \in I-I'} \frac{d_i+1}{4} < \frac{d+1}{2}$, entra\^inent
$$\frac{d+1}{2} \geq g-k + 4 - \delta(\pi_F,1)\, =\, (g-k+2+\nu)+(2-\nu-\delta(\pi_F,1)).$$
Supposant l'in\'egalit\'e $\delta(\pi_F,1) \leq 2-\nu$, on en d\'eduit $\frac{d+1}{2}
\geq g-k+2+\nu$. On a donc ${\rm ord}_{s=\frac{d+1}{2}} \, \, \xi_{\rm B}(s,\pi_F,{\rm
St}) \geq 0$
d'apr\`es Mizumoto. D'autre part, on a aussi ${\rm ord}_{s=\frac{d+1}{2}} \, \xi_{\rm A}(s,\pi_F,{\rm
St}) \,=\, -1 + \delta(\pi_F,\frac{d+1}{2})$, et $${\rm ord}_{s=\frac{d+1}{2}} \,\,\frac{\xi_{\rm B}(s,\pi_F,{\rm St})}{\xi_{\rm A}(s,\pi_F,{\rm St})} \leq 0$$ d'apr\`es la formule \eqref{quotxiaxib}, on obtient l'in\'egalit\'e attendue $\delta(\pi_F,\frac{d+1}{2}) \geq 1$. Cela d\'emontre le cas (a). On proc\`ede de m\^eme dans les cas (b), (c), (d) et (e). 
\end{pf}

\chapter{D\'emonstration des th\'eor\`emes principaux}\label{chap9}

\section{Les formes modulaires de genre $2$ de Tsushima}\label{prelimsp4}

Si $j \geq 0$ et
$k \in \Z$ sont des entiers, notons ${\rm S}_{j,k}$ l'espace
${\rm S}_W({\rm Sp}_4(\Z))$ o\`u $W$ est la repr\'esentation ${\rm Sym}^j \, \C^2 \otimes \DET^k$ de $\GL_2(\C)$ (\S \ref{formesiegel}, \S
VI.\ref{carinf}).  Il est nul si $j \equiv 1 \bmod 2$, car alors ${\rm -1}_2$
agit par $-{\rm id}$ sur $W$, ou si $k \leq 0$ (Freitag \cite[Prop. 
4.6]{vandergeer}), c'est pourquoi nous supposerons toujours $j \equiv 0 \bmod 2$ et $k>0$. 
\ps\ps

\subsection{La formule de dimension de Tsushima} \label{dimtsushima}

Une formule explicite pour $\dim \, {\rm S}_{j,k}$ a \'et\'e
d\'etermin\'ee par R.  Tsushima quand $k\geq 5$ \cite{tsushima}, \'etendant
un r\'esultat d'Igusa concernant les formes scalaires (cas $j=0$, $k \in \Z$
quelconque \cite{igusagenus2}).  Lorsque $j + 2 k -3 \leq 21$, qui
s'av\`erera \^etre le cas d'int\'er\^et pour ce m\'emoire, la formule de
Tsushima montre que ${\rm S}_{j,k}=0$, sauf pour $6$ valeurs
$(j,k)$ donn\'ees dans la table suivante et pour lesquelles $\dim \,
{\rm S}_{j,k}=1$. Comme nous le verrons plus loin, $\dim
{\rm S}_{j,k}$ est \'egalement nul quand $k \leq 4$ et
$j+2k-3 \leq 21$ : voir la remarque \ref{formtsushimataibi}.  La ligne $(w,v)$ sera expliqu\'ee au \S \ref{parstdtsushima}.

\begin{table}[h]
\renewcommand{\arraystretch}{1.5}
\begin{center}
\begin{tabular}{|c||c|c|c|c|c|c|}
\hline $(j,k)$ & $(0,10)$ & $(6,8)$ & $(0,12)$ & $(4,10)$ & $(8,8)$ & $(12,6)$\cr
\hline $(w,v)$ & $(17,1)$ & $(19,7)$ & $(21,1)$ & $(21,5)$ & $(21,9)$ & $(21,13)$\cr
\hline
\end{tabular}
\end{center}

\caption{\small Les couples $(j,k)$ tels que $\dim {\rm S}_{j,k} \neq 0$, pour {\scriptsize $j+2k-3 \leq 21$} et {\scriptsize $k \geq 5$}, d'apr\`es Tsushima.}
\label{dimsjknonnul}
\end{table}

Pour chacun des $6$ couples $(j,k)$ ci-dessus, notons $F_{j,k}$ un
g\'en\'erateur de ${\rm S}_{j,k}({\rm Sp}_4(\Z))$. \'Etant donn\'e le grand
r\^ole jou\'e par ces formes de Siegel dans la suite, expliquons comment
montrer directement leur existence, par une construction de s\'erie th\^eta
b\^atie sur le r\'eseau ${\rm E}_8$. \ps\ps 

Fixons pour cela $j\geq 0$ pair et $k\geq
4$, des entiers pour l'instant quelconque. Il existe un unique isomorphisme de $\C$-alg\`ebres 
$\C[X,Y] \isomo {\rm Sym} \,\C^2$ envoyant respectivement $X$ et $Y$ sur les \'el\'ements $(1,0)$ et $(0,1)$ de $\C^2$. Par transport de structure, cette isomorphisme munit $\C[X,Y]$ d'une repr\'esentation de $\GL_2(\C)$. Le sous-espace $\C[X,Y]_j \subset \C[X,Y]$ des polyn\^omes homog\`enes de degr\'e $j$ est une sous-repr\'esentation isomorphe \`a ${\rm Sym}^j \,\C^2$. \ps\ps 

Soit $I \subset {\rm E}_8 \otimes \C$ un sous-espace isotrope 
de dimension $2$, ainsi que $u,v$ et $w$ trois \'el\'ements de $I$.
Consid\'erons l'application ${\rm E}_8^2 \rightarrow \C[X,Y]_j$
d\'efinie par 
$${\rm P}_{j,k,u,v,w}(x,y) = \DET^{k-4} \left[ \begin{array}{cc} x \cdot u & x
\cdot v \\ y \cdot u & y \cdot v \end{array} \right] ( (x \cdot w) X + (y
\cdot w) Y)^j.$$
L'\'equation fonctionnelle de la fonction $\vartheta$ de Jacobi (\`a deux
variables) permet de d\'emontrer que la fonction 
$$\vartheta_2({\rm E}_8,{\rm P}_{j,k,u,v,w})=\sum_{(x,y) \in {\rm E}_8 \times {\rm E}_8 }
{\rm P}_{j,k,u,v,w}(x,y)
\, \, q^{\frac{1}{2}\left[ \begin{array}{cc} x \cdot x & x \cdot y \\  y\cdot x &
y \cdot y \end{array}\right]}$$
est une forme modulaire de Siegel pour ${\rm Sp}_4(\Z)$ \`a coefficients dans
la repr\'esentation ${\rm Sym}^j \otimes \DET^k$ \cite[\S 2]{freitagharmts}. Elle est
manifestement parabolique si $k>4$, et si $k=4$ le coefficient en
$X^j$ de son image par l'op\'erateur $\Phi_1$ de Siegel est la s\'erie
th\^eta du polyn\^ome harmonique $x \mapsto (x \cdot w)^j$ sur ${\rm E}_8$ (\S V.\ref{sertharm}), un \'el\'ement de ${\rm
M}_{j+4}({\rm SL}_2(\Z))$ qui est parabolique si $j>0$. \ps\ps 

\bigskip

\noindent {\sc Application num\'erique} : Il n'est pas difficile, \`a l'aide de l'ordinateur, de d\'eterminer les
coefficients de Fourier de $f_{j,k,u,v,w}=\vartheta_2({\rm E}_8,{\rm P}_{j,k,u,v,w})$ en
des matrices de Gram de petit discriminant : nous renvoyons \`a la feuille de calculs \cite{clcalc} pour une justification des affirmations qui vont suivre. D\'ecrivons le r\'esultat du calcul en discriminant $\leq
12$, obtenu en prenant $u=(2,i,i,i,i,0,0,0)$, $v=(0,0,0,i,-i,i,i,2)$, $w=u+v$, et en \'enum\'erant b\^etement tous les
couples $(x,y) \in {\rm E}_8^2$ de matrice de Gram l'une des $7$
matrices apparaissant dans la table~\ref{tablecoeffdjk}. On constate que pour chacun des $6$ couples $(j,k)$ en question, 
tous les coefficients calcul\'es sont non nuls, hormis l'un d'eux pour $(j,k)=(6,8)$. La table~\ref{tablecoeffdjk} indique pr\'ecis\'ement les coefficients de Fourier de $\frac{1}{\lambda_{j,k}} f_{j,k,u,v,u+v}$, o\`u $\lambda_{j,k} \in \Z-\{0\}$ est une constante qui n'a pas de signification particuli\`ere (tout comme notre choix de $u, v, w$), et que nous n'expliciterons pas ici. En guise de v\'erification, mentionnons que lorsque $(j,k)=(0,10)$ nos calculs sont compatibles avec la table IV de
\cite{ressal}.
Comme on a ${\rm
S}_{14}({\rm SL}_2(\Z))=0$, cela red\'emontre ${\rm S}_{j,k} \neq 0$ dans tous les cas. $\square$ \ps \bigskip

Le couple $(j,k)$ \'etant fix\'e, on peut \'egalement v\'erifier que si l'on fait varier les param\`etres $u,v$ et $w$ dans le calcul ci-dessus (ou m\^eme par appel \`a un calcul formel), alors le quadruplet de
coefficients calcul\'es (des polyn\^omes~!) n'est
modifi\'e que par un scalaire, comme il se doit car $\dim {\rm S}_{j,k}=1$.  Cela
peut \'egalement se d\'emontrer autrement, de la mani\`ere suivante.  \ps\ps 

Posons $W_{j,k}={\rm
Sym}^j \C^2 \otimes \det^k$ et notons $U_{j,k}$ la
repr\'esentation naturelle de ${\rm O}_8(\C)$ sur l'espace des polyn\^omes
${\rm E}_8 \otimes \C^2 \rightarrow W_{j,k-4}$ qui sont
$\GL_2(\C)$-\'equivariant et {\it pluriharmoniques} \cite{kw}
\cite{freitagharmts}; la fonction ${\rm P}_{j,k,u,v,w}$ est un \'el\'ement
typique de $U_{j,k}$. Ces r\'ef\'erences assurent que si $k\geq 4$, le couple $(U_{j,k},W_{j,k})$ est compatible au sens
du \S VII.\ref{releichlerbisrallis}. Pr\'ecis\'ement, on dispose d'une application lin\'eaire 
\begin{equation}\label{mapthetajk}\vartheta : {\rm M}_{U_{j,k}}({\rm O}_8) \longrightarrow {\rm
M}_{W_{j,k}}({\rm Sp}_4(\Z)).\end{equation}
envoyant l'\'el\'ement $[{\rm E}_8,{\rm P}_{j,k,u,v,w}]$, d\'efini \`a la fin paragraphe IV.\ref{prodherminv}, sur la s\'erie th\^eta $\vartheta_2({\rm E}_8,{\rm P}_{j,k,u,v,w})$, pour tout triplet d'\'el\'ements $u,v,w$ appartenant \`a un m\^eme sous-espace isotrope de rang $2$ de ${\rm E}_8 \otimes \C$. \ps\ps 

On v\'erifie ais\'ement que la repr\'esentation $U_{j,k}$, qui est
irr\'eductible restreinte \`a ${\rm SO}_8(\C)$ d'apr\`es \cite{kw}, admet un
plus haut poids de la forme $(j+k-4) \varepsilon_1 + (k-4) \varepsilon_2$ dans les
notations du \S VI.\ref{exparlan}. Or les tables\footnote{Voir
\url{http://gaetan.chenevier.perso.math.cnrs.fr/table/dim_SO8_dom.txt}} de \cite{chrenard2} montrent que pour les $6$ couples
$(j,k)$ d'int\'er\^et, on a $\dim {\rm M}_{U_{j,k}}({\rm O}_8) = 1$. Nous en d\'eduisons, comme promis, que l'espace $\vartheta({\rm M}_{U_{j,k}}({\rm O}_8))$ est de dimension $1$ pour ces couples. \'Etant donn\'e que l'on a aussi $\dim {\rm S}_{j,k}=1$, on obtient la proposition suivante.

\begin{prop}\label{calcimagetheta} Si $(j,k)$ est l'un des $6$ couples de la
table \ref{dimsjknonnul}, alors l'application \eqref{mapthetajk} induit un isomorphisme
${\rm M}_{U_{j,k}}({\rm O}_8) \isomo {\rm S}_{j,k}$ entre espaces de dimension $1$.
\end{prop}

\subsection{Param\`etres standards des $6$ premi\`eres formes de genre $2$} \label{parstdtsushima}

Soit $F \in {\rm S}_{j,k}$ une forme propre pour l'action de ${\rm H}({\rm PGSp}_4)$. Notons $\pi_F \in \Pi_{\rm cusp}(\PGSp_4)$ la repr\'esentation engendr\'ee par $F$ (Corollaire VI.\ref{corapiw}).  Observons que le groupe de Chevalley ${\rm PGSp}_4$ est un $\Z$-groupe classique, \'etant isomorphe au $\Z$-groupe ${\rm SO}_{3,2}$ ; son dual de Langlands est le $\C$-groupe ${\rm Sp}_4$, il est muni de sa repr\'esentation standard de dimension $4$. D'apr\`es la fin du  \S VI.\ref{carinf}, la classe de conjugaison semi-simple ${\rm St}({\rm c}_\infty(\pi_F)) \subset {\rm M}_4(\C)$ a pour valeurs propres $\pm \frac{w}{2}$ et $\pm \frac{v}{2}$, o\`u 
$$(w,v)=(j+2k-3,j+1),$$
ce qui explique la seconde ligne de la table \ref{dimsjknonnul}. Observons que l'application $(j,k) \mapsto (w,v)$ est une bijection entre l'ensemble des couples $(j,k)$ avec $j \geq 0$ pair et $k\geq 3$, et l'ensemble des couples $(w,v)$ avec $w,v$ impairs tels que $w>v>0$. \ps\ps 

Pour chacun des $6$ couples $(j,k)$ de la table \ref{dimsjknonnul}, l'action de ${\rm H}({\rm PGSp}_{4})$ sur ${\rm S}_{j,k}$ est alors trivialement scalaire, {\it i.e.} $F_{j,k}$ est propre, et nous allons nous int\'eresser au param\`etre
$$\psi_{j,k} = \psi(\pi_{F_{j,k}},{\rm St}) \in \mathcal{X}({\rm SL}_4).$$
Le cas  de la forme scalaire $F_{0,10}$ a une histoire fameuse puisque c'est la premi\`ere forme de Saito-Kurokawa, associ\'ee \`a la forme modulaire de poids $18$ pour ${\rm SL}_2(\Z)$ (\cite{kurokawa}, \S \ref{ikedabocherer}). Du fait que l'on consid\`ere $\pi_{F_{0,10}}$ comme repr\'esentation de ${\rm PGSp}_4$, plut\^ot que ${\rm Sp}_4$, on a la relation (voir \cite{eichlerzagier},\cite{zagiersk})
$$\psi_{0,10} = \Delta_{17} \oplus [2],$$
ce qui est manifestement compatible avec l'\'egalit\'e $(w,v)=(17,1)$ (la notation $\Delta_w$ est introduite au \S\ref{ikedabocherer}). Le cas de la forme $F_{0,12}$ est similaire, et l'on a $\psi_{0,12}=\Delta_{21} \oplus [2]$ d'apr\`es Andrianov, Maass et Zagier.  Comme l'avait \'egalement devin\'e Kurokawa, et expliqu\'e Arthur \cite{arthursp4}, la situation est tr\`es diff\'erente pour les $4$ autres repr\'esentations.

\begin{propv}\label{amfpgsp4} On suppose $j>0$ pair et $k \geq 3$.
\begin{itemize} \ps \ps
 \item[(i)] {\rm (Multiplicit\'e 1)} Si $F,G \in  {\rm S}_{j,k}$ sont deux
formes propres pour l'action de ${\rm H}({\rm Sp}_4)$, et si tout
\'el\'ement de ${\rm H}({\rm Sp}_4)$ a m\^eme valeur propre sur $F$ et
$G$, alors $F$ et $G$ sont proportionnelles.  
\ps\ps  \item[(ii)] Si $F \in {\rm
S}_{j,k}$ est propre pour l'action de ${\rm H}({\rm PGSp}_4)$ alors $\psi(\pi_F,{\rm
St})=\pi$ avec $\pi \in \Pi_{\rm cusp}^\bot(\PGL_4)$. 
\ps\ps  \item[(iii)] L'application $F \mapsto \psi(\pi_F,{\rm St})$ induit une
bijection entre l'ensemble des droites propres de ${\rm S}_{j,k}({\rm
Sp}_4(\Z))$ sous l'action de ${\rm H}({\rm PGSp}_4)$ et l'ensemble des $\pi
\in \Pi_{\rm cusp}^\bot(\PGL_4)$ telles que ${\rm Poids}(\pi)=\{\pm
\frac{j+2k-3}{2}, \pm \frac{j+1}{2}\}$.  \end{itemize} \end{propv}

\begin{pf} Le (i) est le cas particulier $g=2$ du corollaire VIII.\ref{mult1classsp} : dans les notations {\it loc. cit.} on a $(m_1,m_2)=(j+k,k)$ de sorte que $m_1>m_2$. \ps

Posons $(w,v)=(j+2k-3,j+1)$. Soit $F \in {\rm S}_{j,k}({\rm
Sp}_4(\Z))$ une forme propre pour ${\rm H}({\rm PGSp}_4)$.  On applique le
th\'eor\`eme~\ref{arthurst} au $\Z$-groupe classique ${\rm PGSp}_4 \simeq
{\rm SO}_{3,2}$ et \`a sa repr\'esentation $\pi_F \in \Pi_{\rm disc}({\rm
PGSp}_4)$. \'Etant donn\'e qu'un $\pi \in \Pi_{\rm cusp}^\bot(\PGL_2)$ est
symplectique (Proposition \ref{proppigl2}), et que $\pm \frac{1}{2}$ n'est pas valeur propre de ${\rm
St}({\rm c}_\infty(\pi_F))$ car $j>0$, il n'y a que deux possibilit\'es pour $\psi(\pi_F,{\rm St})$ 
(corollaire VIII.\ref{cormodulo4} (ii)) : \ps\ps  
\begin{itemize}
\item[(a)] soit $\psi(\pi_F,{\rm St})=\pi_1 \in \Pi_{\rm
cusp}^\bot(\PGL_4)$, \ps\ps 
\item[(b)] soit $\psi(\pi_F,{\rm St}) = \pi_1 \oplus \pi_2$ avec $\pi_1,
\pi_2 \in \Pi_{\rm cusp}(\PGL_2)$ tels que ${\rm w}(\pi_1)=w$ et ${\rm
w}(\pi_2)=v$.
\end{itemize}
\ps
\noindent Pour d\'emontrer le (ii), il faut donc montrer que le cas (b) ne se produit
pas.  Observons que pour les quatres couples $(j,k)$ de la table
\ref{dimsjknonnul} on a $v \in \{5, 7, 9, 13\}$, de sorte que cela d\'ecoule
directement de la proposition \ref{proppigl2} et de ce que $\dim {\rm
S}_{v+1}({\rm SL}_2(\Z))=0$ pour ces valeurs de $v$. Pour un couple $(j,k)$ g\'en\'eral, c'est
plut\^ot une cons\'equence de la formule de multiplicit\'e d'Arthur pour
${\rm SO}_{3,2}$ (Thm.  VIII.\ref{amfchev}).  \ps\ps 
Supposons en effet que $\psi=\pi_1 \oplus \pi_2 \in {\mathcal{X}}_{\rm AL}({\rm SL}_4)$ o\`u $\pi_1$ et
$\pi_2$ sont comme dans le (b) ci-dessus. Consid\'erons des homomorphismes
$\nu$ et $\nu_\infty$ associ\'es \`a $\psi$ comme au \S VIII.\ref{parpaqarch}. Par
d\'efinition, $\nu : {\rm SL}_2 \times ({\rm SL}_2 \times {\rm SL}_2) \rightarrow {\rm
Sp}_4$ est trivial sur le premier facteur ${\rm SL}_2$ et ${\rm St} \circ
\nu$ est la somme directe des repr\'esentations tautologiques de chacun des
deux autres facteurs ${\rm SL}_2$, de sorte que l'inclusion 
$\{\pm 1 \}^2 = {\rm C}_\nu \hookrightarrow {\rm C}_{\nu_\infty}$ est
une \'egalit\'e. On en d\'eduit que
$\varepsilon_\psi=1$ (car ``$d_i=1$ pour tout $i$'') et que $\Pi(\nu_\infty)$ 
est l'ensemble des deux s\'eries discr\`etes de ${\rm SO}_{3,2}(\R)$ de caract\`ere
infinit\'esimal $\psi_\infty$ d'apr\`es le \S VIII.\ref{paramdiscC}. L'une est holomorphe, disons $\pi_{\rm hol}$, et l'autre g\'en\'erique (la
notion est canonique ici car ${\rm SO}_{3,2}$ est adjoint), de sorte que le caract\`ere
de Shelstad $\chi_{\pi_{\rm hol}}$ est le caract\`ere non trivial de ${\rm
C}_{\nu_\infty}$ qui est trivial sur le centre de $\widehat{{\rm
SO}_{3,2}}={\rm Sp}_4$, \`a savoir le sous-groupe $\{\pm1\}$ diagonal dans ${\rm
C}_{\nu_\infty}$. Une autre mani\`ere de le voir consiste \`a appliquer
simplement la formule VIII.\eqref{carholsoimpair} dans le cas particulier $r=2$. Ainsi, la
restriction la restriction de $\chi_{\pi_{\rm hol}}$ \`a ${\rm C}_\nu={\rm
C}_{\nu_\infty}$ est non triviale, et la formule de multiplicit\'e
d'Arthur affirme que l'unique $\pi \in \Pi({\rm PGSp}_4)$ telle que $\pi
\simeq \pi_{\rm hol}$ et $\psi(\pi,{\rm St})=\psi$ est de multiplicit\'e
nulle (Thm.  VIII.\ref{amfchev}). Comme ${\rm m}(\pi_F)>0$, on est dans le cas (a) ci-dessus, ce qui
prouve le (ii). La m\^eme formule de multiplicit\'e d'Arthur affirme alors
que ${\rm m}(\pi_F)=1$ et plus g\'en\'eralement que pout tout
$\pi \in \Pi_{\rm cusp}(\PGL_4)$ de poids $\{\pm
\frac{w}{2},\pm\frac{v}{2}\}$ il existe un (unique) $\pi \in \Pi_{\rm disc}({\rm
PGSp}_4)$ tel que $\psi(\pi,{\rm St})=\pi$ et $\pi_\infty \simeq \pi_{\rm
hol}$, et qu'il v\'erifie ${\rm m}(\pi)=1$ (c'est le cas o\`u il n'y a rien \`a v\'erifier car
${\rm C}_\psi={\rm Z}(\widehat{G})$). L'assertion (iii) suit alors du corollaire VI.\ref{corapiw}.\end{pf}
\ps
\noindent Nous avons fait appel \`a l'\'enonc\'e tr\`es classique suivant.\ps

\begin{prop} \label{proppigl2} Soient $k \geq 2$ un entier pair et $\mathcal{F}_k \subset {\rm
S}_k({\rm SL}_2(\Z)$ l'ensemble des formes modulaires propres pour ${\rm
H}({\rm PGL}_2)$ et normalis\'ees (i.e. de coefficient de Fourier en $q$ \'egal
\`a $1$). L'application associant \`a $F \in \mathcal{F}_k$ la
repr\'esentation engendr\'ee $\pi_F \in \Pi_{\rm cusp}(\PGL_2)$ induit une
bijection entre $\mathcal{F}_k$ et l'ensemble des $\pi \in \Pi_{\rm
cusp}(\PGL_2)$ telles que ${\rm Poids}(\pi)=\{\pm \frac{k-1}{2}\}$. 
\end{prop}

\begin{pf} On rappelle qu'une forme propre $F \in {\rm S}_k({\rm SL}_2(\Z))$
a toujours son coefficient en $q$ non nul, et que si elle est normalis\'ee
elle est uniquement d\'etermin\'ee par ses valeurs propres sous ${\rm
H}({\rm PGL}_2)$\cite[Chap. VII Thm. 7]{serre}, en particulier $\mathcal{F}_k$ est une base de l'espace
vectoriel ${\rm S}_k({\rm SL}_2(\Z))$. Soit $U_k$ la s\'erie discr\`ete de ${\rm PGL}_2(\R)$ telle que ${\rm Inf}_{U_k} \subset {\rm M}_2(\C)$ a pour valeurs propres
$\pm \frac{k-1}{2}$. Un cas particulier bien connu de la proposition IV.\ref{corapiw} est que l'on dispose d'un isomorphisme ${\rm H}({\rm PGL_2})$-\'equivariant entre ${\rm S}_k({\rm SL}_2(\Z))$ et 
$\mathcal{A}_{U_k}({\rm PGL}_2)={\rm Hom}_{{\rm PGL}_2(\R)}(U_k,\mathcal{A}_{\rm cusp}({\rm PGL}_2))$ \cite[Ch. I \S 4]{GGPS}. Cela montre que l'application de l'\'enonc\'e est bien d\'efinie et injective. \ps\ps 
Une justification sophistiqu\'ee de la surjectivit\'e consiste \`a invoquer la proposition~\ref{contraintesalg} (i). On peut aussi utiliser la classification de Bargmann \cite{bargmann} du dual unitaire de ${\rm SL}_2(\R)$. Elle montre que si $U$ est une repr\'esentation irr\'eductible unitaire de ${\rm PGL}_2(\R)$ telle que ${\rm Inf}_{U} \subset {\rm M}_2(\C)$ a pour valeurs propres
$\pm \frac{k-1}{2}$, alors soit $U \simeq U_k$, soit $\dim U=1$ et $k=2$. En effet, les repr\'esentations de la s\'erie principale ont un caract\`ere infinit\'esimal dont les valeurs propres sont de la forme $\pm i s$ avec $s \in \R$ (cas ``temp\'er\'e'') ou de la forme $\pm s$ avec $s \in ]-\frac{1}{2},\frac{1}{2}[$ (``s\'erie compl\'ementaire'') ; de plus, le caract\`ere infinit\'esimal de la ``limite de s\'erie discr\`ete'' est $0$. Pour \'eliminer le cas $\dim U=1$, on observe que la fomule IV.\eqref{decompsymp} entra\^ine que les seuls \'el\'ements de $\mathcal{A}^2(\PGL_2)$ invariants par ${\rm PGL}_2(\R)^+$ sont les fonctions constantes, qui ne sont pas cuspidales, ce qui conclut.
\end{pf}
\ps
\noindent La proposition \ref{amfpgsp4} et la table de Tsushima justifient la
d\'efinition-proposition suivante.\ps

\begin{defpropv}\label{definitiondeltawv} Si $(w,v) \in \{(19,7),(21,5),(21,9),(21,13)\}$, il existe une unique repr\'esentation dans $\Pi_{\rm cusp}^\bot(\PGL_4)$ qui est alg\'ebrique 
de poids l'ensemble $\{\pm \frac{w}{2},\pm \frac{v}{2}\}$. On la note $\Delta_{w,v}$. 
\end{defpropv}

\noindent Ainsi, pour les $4$ derni\`eres colonnes de la table \ref{dimsjknonnul} on a
la relation $\psi_{j,k}=\Delta_{w,v}$ et $\psi_{j,k}$ ne s'exprime pas \`a
l'aide de formes de genre $1$. \ps

\subsection{Quelques valeurs propres des op\'erateurs de Hecke}\label{definitiontaujk}

Soient $(j,k) \in \{(6,8), (4,10), (8,8), (12,6)\}$, $(w,v)=(j+2k-3,j+1)$, $p$ un nombre premier et $n\geq 1$ un entier. On pose 
\begin{equation}\label{eqdeftaujk} \tau_{j,k}(p^n)
\, \, \, \,   = \, \, \, \, p^{\frac{n w}{2}}\, {\rm trace}\, {\rm St}({\rm c}_p(\pi_{F_{j,k}})^n)\, \, \, \,   = \, \, 
\, \, p^{\frac{n w}{2}}\, 
{\rm trace} \, (\, {\rm c}_p(\Delta_{w,v})^n\, ) .\end{equation} \smallskip
La classe de conjugaison $ {\rm St}({\rm c}_p(\pi_{F_{j,k}})) \subset {\rm SL}_4(\C)$ est \'egale \`a son inverse, elle a donc pour polyn\^ome caract\'eristique 
\begin{equation}\label{polcarsp4} t^4 - \tau_{j,k}(p) \,t^3 +
\frac{\tau_{j,k}(p)^2-\tau_{j,k}(p^2)}{2} \,t^2 - \tau_{j,k}(p) p^{j+2k-3}
\,t + p^{2j+4k-6}.\end{equation}
En particulier, le nombre complexe $\tau_{j,k}(p^n)$ est un polyn\^ome \`a coefficients entiers en $\tau_{j,k}(p)$ et $\frac{\tau_{j,k}(p)^2-\tau_{j,k}(p^2)}{2}$. \ps \ps

La proposition suivante \'etait connue de Shimura \cite{shimurasp}.  Comme l'avait d\'ej\`a expliqu\'e Gross dans \cite[\S 6]{grossatake}, c'est \'egalement une cons\'equence imm\'ediate de la relation VI.\eqref{satakesp} (voir aussi le \S VI.\ref{exparlan}). \ps\ps

\begin{prop}\label{coeffsatsp4shimura} Soient $j,k$ comme ci-dessus et $p$ un nombre premier. \ps \smallskip \begin{itemize}
\item[(a)] Le nombre complexe $\tau_{j,k}(p)$ est la valeur propre de l'op\'erateur de Hecke $p^{\frac{j+2k-6}{2}}{\rm K}_p$ agissant sur la droite ${\rm S}_{j,k}$. \ps \smallskip
\item[(b)] Le nombre complexe $\frac{\tau_{j,k}(p)^2-\tau_{j,k}(p^2)}{2}$ est la valeur propre de l'op\'erateur de Hecke $p^{j+2k-5}({\rm T}_p+1)+p^{j+2k-3}$ agissant sur la droite ${\rm S}_{j,k}$.  \ps \smallskip
\end{itemize}
\end{prop}

D'apr\`es la relation IV.\ref{diaghecke}, l'op\'erateur $p^{\frac{j+2k-6}{2}}\,{\rm K}_p$ co\"incide avec celui not\'e ${\rm T}(p)$ chez van der Geer \cite[\S 16]{vandergeer}, du moins lorsque ce dernier est d\'efini en incluant la normalisation entre parenth\`eses dans la d\'efinition 16.5 {\it loc. cit.} L'op\'erateur not\'e ici ${\rm T}_p$ est parfois not\'e ${\rm T}_1(p^2)$ dans la litt\'erature, au facteur pr\`es d'une puissance de $p$ d\'ependant des auteurs. \ps \ps

Le probl\`eme de d\'eterminer les valeurs propres des op\'erateurs de Hecke
agissant sur les espaces ${\rm S}_{j,k}$ est sensiblement plus difficile en
pratique que son analogue en genre $1$. Une raison est la difficult\'e \`a
d\'eterminer les coefficients de Fourier des formes de genre $2$, notamment
en des matrices de Gram de grand d\'eterminant. De plus, la relation entre
coefficients de Fourier et valeurs propres, \'etudi\'ee dans ce contexte par
Andrianov \cite{andrianov2} pour les formes scalaires et \'etendue aux
formes vectorielles par Arakawa \cite{arakawa}, est plus subtile qu'en genre
$1$. Dans ce qui suit, nous rappelons cette relation dans le cas des op\'erateurs de Hecke ${\rm K}_{p}$ et ${\rm T}_{p}$. Nous en d\'eduirons \`a la fois la proposition suivante, ainsi que le calcul de quelques valeurs des entiers $\tau_{j,k}(p)$. \ps \medskip

\begin{prop} \label{intcongtaujk} Soient $(j,k) \in \{(6,8), (4,10), (8,8), (12,6)\}$, $p$ un nombre premier et $n\geq 1$ un entier. On a $\tau_{j,k}(p^n) \in \Z$ ainsi que la congruence $$\tau_{j,k}(p^2) \, \, \equiv \, \, \tau_{j,k}(p)^2 \, \, \bmod \, \, 2 \,p^{k-2}.$$ 
\end{prop}

\bigskip

\noindent Posons $\Gamma = {\rm Sp}_{4}(\Z)$ et consid\'erons les \'el\'ements suivants de $\mathrm{GSp}_{4}(\Z[\frac{1}{p}])^+$ 
$${\small
\hspace{24pt}
\gamma
\hspace{4pt}:=\hspace{4pt}
\begin{bmatrix}
1 & 0 & 0 & 0 \\ 0 & 1 & 0 & 0 \\ 0 & 0 & p & 0 \\ 0 & 0 & 0 & p
\end{bmatrix},\hspace{24pt} 
\, \, \, 
\gamma'
\hspace{4pt}:=\hspace{4pt}
\begin{bmatrix}
1 & 0 & 0 & 0 \\ 0 & p & 0 & 0 \\ 0 & 0 & p^{2} & 0 \\ 0 & 0 & 0 & p
\end{bmatrix}.}
$$

\bigskip
\begin{lemme}\label{lemmdbclandrianov} \begin{itemize} \item[(a)] L'op\'erateur de Hecke ${\rm K}_p \in {\rm H}_p({\rm PGSp}_4)$ est de degr\'e $(1+p)(1+p^2)$. Sa matrice est la fonction caract\'eristique de l'image de $\Gamma \gamma^{-1} \Gamma$ dans ${\rm PGSp}_4(\Z[\frac{1}{p}])^+$, au sens des identifications {\rm IV}.\eqref{bijmathecke} et {\rm IV}.\eqref{bijdbclasssp}. La double-classe $\Gamma \gamma \Gamma$ est r\'eunion disjointe des classes \`a droite $\Gamma \gamma_{i}$, o\`u $\gamma_{i}$ parcourt la liste des \'el\'ements de la forme suivante:$${\small
\begin{bmatrix}
p & 0 & 0 & 0 \\ 0 & p & 0 & 0 \\ 0 & 0 & 1 & 0 \\ 0 & 0 & 0 & 1
\end{bmatrix},\hspace{15pt} 
\begin{bmatrix}
1 & 0 & a & b \\ 0 & 1 & b & c \\ 0 & 0 & p & 0 \\ 0 & 0 & 0 & p
\end{bmatrix},  \hspace{15pt} 
\begin{bmatrix}
p & 0 & 0 & 0 \\ -d & 1 & 0 & e \\ 0 & 0 & 1 & d \\ 0 & 0 & 0 & p
\end{bmatrix}\, \, \,  {\rm ou}\, \, \,  
\begin{bmatrix}1 & 0 & f & 0 \\ 0 & p & 0 & 0 \\ 0 & 0 & p & 0 \\ 0 & 0 & 0 & 1
\end{bmatrix},}
$$
avec $a,b,c,d,e,f$ des entiers compris entre $0$ et $p-1$~;\par\medskip
\item[(b)]  De m\^eme, ${\rm T}_p \in {\rm H}_p({\rm PGSp}_4)$ est de degr\'e $p\, \frac{p^{4}-1}{p-1}$ et a pour matrice la fonction caract\'eristique de l'image de $\Gamma {\gamma'}^{-1} \Gamma$ dans ${\rm PGSp}_4(\Z[\frac{1}{p}])^+$. La double-classe $\Gamma {\gamma'} \Gamma$ est la r\'eunion disjointe des classes \`a droite $\Gamma \gamma'_{i}$, o\`u $\gamma'_{i}$ parcourt la liste des \'el\'ements de la forme suivante : $${\small
\begin{bmatrix}
p & 0 & 0 & 0 \\ 0 & p^{2} & 0 & 0 \\ 0 & 0 & p & 0 \\ 0 & 0 & 0 & 1
\end{bmatrix}, \hspace{15pt}
\begin{bmatrix}
p^{2} & 0 & 0 & 0 \\ -ap & p & 0 & 0 \\ 0 & 0 & 1 & a \\ 0 & 0 & 0 & p
\end{bmatrix}, \hspace{15pt} 
\begin{bmatrix}
p & 0 & b & c \\ 0 & p & c & d \\ 0 & 0 & p & 0 \\ 0 & 0 & 0 & p
\end{bmatrix},}$$
$${\small \begin{bmatrix} p & 0 & 0 & pe \\ -f & 1 & e & ef+g \\ 0 & 0 & p & pf \\ 0 & 0 & 0 & p^{2}
\end{bmatrix}\, \, \,  {\rm ou}\, \, \,  
\begin{bmatrix} 1 & 0 & h & i \\ 0 & p & pi & 0 \\ 0 & 0 & p^{2} & 0 \\ 0 & 0 & 0 & p
\end{bmatrix},}
$$
avec $a,b,c,d,e,f,i$ des entiers compris entre $0$ et $p-1$, tels que $c^2 \equiv bd\bmod p$ et $(b,c,d) \neq (0,0,0)$, et avec $g$ et $h$ des entiers compris entre $0$ et $p^{2}-1$.\par\medskip
\end{itemize}
\end{lemme}

\begin{pf} Le fait que les matrices de ${\rm K}_p$ et ${\rm T}_p$ sont les fonctions caract\'eristiques respectives des images des doubles-classes $\Gamma \gamma^{-1} \Gamma$ et $\Gamma {\gamma'}^{-1} \Gamma$ est la formule VI.\eqref{tradsatsymp}. \ps\ps
Le degr\'e de ${\rm K}_p \in {\rm H}_p({\rm PGSp}_{2g})$ est le nombre des lagrangiens du ${\rm a}$-module hyperbolique ${\rm H}(\F_p^g)$, \`a savoir $\prod_{i=1}^g(1+p^i)$. De m\^eme, le degr\'e de ${\rm T}_p$ est le nombre de droites (isotropes) de ${\rm H}(\F_p^g)$, multipli\'e par le nombre des droites (isotropes) de ${\rm H}(\F_p)$ transverses \`a une droite donn\'ee, ce qui fait $p\frac{p^g-1}{p-1}$. \ps\ps
Les assertions relatives aux d\'ecompositions des doubles-classes sont dues \`a Andrianov \cite{andrianov2,andrianov}. Justifions-les bri\`evement, en suivant les notes de Buzzard \cite{buzzardsiegel}. Un \'el\'ement de ${\rm GL}_{4}$, disons donn\'e par blocs de taille $2 \times 2$, et suppos\'e de la forme $$\left[\begin{array}{cc} a & b \\ 0 & d \end{array}\right],$$ est dans ${\rm GSp}_{4}$ pour le facteur de similitude $\nu$ si, et seulement si, $a\,{}^{\rm t}b = b\,{}^{\rm t}a$ et $a\,{}^{\rm t}d = \nu 1_{g}$ (\S IV.\ref{fsiegelclass}). Cela montre que chacun des \'el\'ements de l'\'enonc\'e est dans ${\rm GSp}_{4}(\Z[\frac{1}{p}])$, avec pour facteur de similitude $p$ dans le cas (a) et $p^{2}$ dans le cas (b). \ps\ps
Soient $h \in {\rm GSp}_{4}(\Z[1/p])\cap {\rm M}_{4}(\Z)$ et $\overline{h} \in {\rm M}_4(\Z/p)$ la r\'eduction modulo $p$ de $h$. La th\'eorie des ``diviseurs \'el\'ementaires symplectiques'' montre que $h$ est dans $\Gamma \gamma' \Gamma$ (resp. $ \Gamma \gamma \Gamma$) si, et seulement si, $\nu(h)=p^2$ et le rang de $\overline{h}$ est $1$ (resp. $\nu(h)=p$). Cela montre $\gamma_i \in\Gamma \gamma  \Gamma$ et $\gamma'_i \in \Gamma\gamma' \Gamma$ pour tout $i$. \ps \ps
Enfin, on v\'erifie que $\gamma_{i}\gamma_{j}^{-1}\in \Gamma$ (resp. $\gamma'_{i}{\gamma'}_{j}^{-1}\in \Gamma$) force $i=j$ ; il est commode pour cela d'observer que tous les \'el\'ements ci-dessus sont dans un m\^eme sous-groupe de Borel de ${\rm GSp}_{4}$, la ``projection sur la diagonale'' \'etant un homomorphisme. On conclut car dans chacun des deux cas le cardinal de la liste est le degr\'e de l'op\'erateur de Hecke.
\end{pf}

\vspace{1cm}

Soit  $j\geq 0$ un entier. Notons $\rho_j$ la repr\'esentation naturelle de $\GL_2(\C)$ sur l'espace $W_j:={\rm Sym}^j \C^2$. On rappelle que pour  $w \in W_{j}$ et $n \in {\rm M}_2(\C)$, la notation $w \, {\rm q}^n$ d\'esigne la fonction $\mathbb{H}_2 \rightarrow W_{j}, \, \, \tau \mapsto e^{2 i \pi\,  {\rm tr} ( n \tau)} \,w$ (\S IV.\ref{devfouriersiegel}).  \ps\ps 

Soient $k \in \Z$ et $F$ une forme modulaire de Siegel pour ${\rm Sp}_4(\Z)$ \`a coefficients dans la repr\'esentation $W_{j,k}:=W_j \otimes \DET^k$ de $\GL_2(\C)$. Par d\'efinition, cette repr\'esentation a pour espace sous-jacent $W_j$, et $\GL_2(\C)$ y agit par $g \mapsto  \rho_{j}(g) \DET(g)^k$. On rappelle que la forme $F$ admet un d\'eveloppement de Fourier, que nous \'ecrirons ici
$$F = \sum_{n \in \mathcal{N}}\, \,  {\rm a}(n ; F) \, \, {\rm q}^n,$$ 
o\`u $\mathcal{N} \subset \frac{1}{2}{\rm M}_2(\Z)$ est le sous-ensemble des matrices sym\'etriques, positives, et de coefficients diagonaux dans $\Z$,  et o\`u ${\rm a}(n ; F) \in W_j$ pour tout $n \in \mathcal{N}$ (\S IV.\ref{devfouriersiegel}). Il sera commode de poser ${\rm a}(n, F)=0$ si $n \in {\rm M}_2(\Q) - \mathcal{N}$. \ps\ps 

On rappelle que l'on a d\'efini au \S IV.\ref{fsiegelclass} une action \`a droite de ${\rm GSp}_4(\R)^+$ sur l'espace des fonctions $\mathbb{H}_2 \rightarrow W_{j,k}$, not\'ee $(f,\gamma) \mapsto f_{|_{W_{j,k}}} \gamma$. Soient $w \in W_j$, $n \in {\rm M}_2(\C)$, et $\gamma=\left[\begin{array}{cc} a & b \\ 0 & d \end{array}\right]$ dans ${\rm GSp}_{4}(\R)^+$ de facteur de similitude $\nu$, on a
\begin{equation}\label{fourierheckecalc} {w \, \, {\rm q}^n\, \, }_{|W_{j,k}} \gamma \, \, = \, \, \, \nu^{-\frac{j+2k}{2}} \, \cdot \, \DET(a)^k \, \cdot \, e^{\frac{2 i \pi}{\nu} {\rm tr} ({}^{\rm t}a n b)}\, \cdot \, \rho_j({}^{\rm t}a) w\, \, \, {\rm q}^{\frac{{}^{\rm t}a n a}{\nu}}\end{equation}
(on rappelle la relation $d^{-1} = \nu^{-1}\, {}^{\rm t} a$). Il sera commode de faire agir le groupe ${\rm GL}_2(\C)$ sur ${\rm M}_2(\C)$ par $(g,s) \mapsto g \cdot s\,:= \,g\,s \,{}^{\rm t}\!g$. 
\bigskip
\begin{cor}\label{corKpTp} Soient $j,k\in \Z$ tels que $j\geq 0$, $F \in S_{j,k}$, $p$ un nombre premier et $n \in \mathcal{N}$. On a :
$$ p^{\frac{j+2k}{2}} {\rm a}(n ; {\rm K}_p \, F) =  \, \, \, p^{j+2k} \,\, {\rm a}(n/p ; F)  \, \, \, + \, \, \, p^{3} \, \, {\rm a}( np ; F) \, \, \,$$ 
$$+ \, \, \,  p^{k+1}\, \, \sum_{d=0}^{p-1} \rho_j( \left[\begin{array}{cc} p & -d \\ 0 & 1 \end{array}\right]) \, \,{\rm a}( \left[\begin{array}{cc} 1 & d \\ 0 & p \end{array}\right]\cdot \frac{n}{p} ; F)\,\,\,$$
$$ +   \, \, \, p^{k+1}\, \, \rho_j ( \left[\begin{array}{cc} 1 & 0 \\ 0 & p \end{array}\right]) \, \, {\rm a} (\,\left[\begin{array}{cc} p & 0 \\ 0 & 1 \end{array}\right] \cdot \frac{n}{p} ; F),$$ 
et
$$ p^{j+2k} {\rm a}(n ; {\rm T}_p \, F) =  \, \, \, p^{j+3k} \,\, \rho_j ( \left[\begin{array}{cc} 1 & 0 \\ 0 & p \end{array}\right]) \, \, {\rm a} (\,\left[\begin{array}{cc} p & 0 \\ 0 & 1 \end{array}\right] \cdot \frac{n}{p^2} ; F) \, \, \, $$ 
$$+\, \, \, p^{j+3k} \,\, \sum_{a=0}^{p-1} \rho_j ( \left[\begin{array}{cc} p & -a \\ 0 & 1 \end{array}\right]) \, \, {\rm a} (\,\left[\begin{array}{cc} 1 & a \\ 0 & p \end{array}\right] \cdot \frac{n}{p^2} ; F) \, \, \, + \, \,\delta(n,p)\, \,  p^{j+2k} \, \, {\rm a}(n ; F) \, \, \, $$
$$+   \, \, \, p^{k+3}\, \, \sum_{f=0}^{p-1}\, \, \rho_j ( \left[\begin{array}{cc} p & -f \\ 0 & 1 \end{array}\right]) \, \, {\rm a} (\,\left[\begin{array}{cc} 1 & f \\ 0 & p \end{array}\right] \cdot n ; F) \,\,\,$$
$$ +   \, \, \, p^{k+3}\, \, \rho_j ( \left[\begin{array}{cc} 1 & 0 \\ 0 & p \end{array}\right]) \, \, {\rm a} (\,\left[\begin{array}{cc} p & 0 \\ 0 & 1 \end{array}\right] \cdot n ; F),$$ \ps
\noindent o\`u $\delta(n,p) \in \Z$ est d\'efini par la formule~\eqref{defdeltanp} ci-dessous  ; $\delta(n,p) \equiv -1 \bmod p$. 
\end{cor}

\begin{pf} D'apr\`es le lemme \ref{lemmdbclandrianov} et le diagramme IV.\ref{diaghecke}, on a $${\rm K}_p F = \sum_i F_{|W_{j,k}} \gamma_i\, \, \, \, {\rm et}\, \, \,\, {\rm T}_p F = \sum_i F_{|W_{j,k}} \gamma'_i.$$ Compte tenu de la convergence uniforme du d\'eveloppement de Fourier de $F$ sur tout compact de $\mathbb{H}_2$, le corollaire est une application directe de la formule \eqref{fourierheckecalc}. On observera que si $a,d \in \GL_2(\C)$, $b \in {\rm M}_2(\C)$, $\nu \in \C^\ast$ et $m, n \in {\rm M}_2(\C)$ sont tels que $n =\frac{ {}^{\rm t}a m a}{\nu}$ et $d^{-1} = \nu^{-1}\, {}^{\rm t} a$, alors $m = d \cdot \frac{n}{\nu}$ et ${\rm tr} ({}^{\rm t} a m b)\, =\, {\rm tr} (b n\, {}^{\rm t}\! d)$.\ps \smallskip

 \noindent En guise d'exemple, d\'eterminons la contribution des $p^2-1$ \'el\'ements $\gamma'_i$ de la forme $$\begin{bmatrix}
p & 0 & b & c \\ 0 & p & c & d \\ 0 & 0 & p & 0 \\ 0 & 0 & 0 & p
\end{bmatrix},$$
avec $b, c, d$ comme dans le lemme \ref{lemmdbclandrianov} (ii), \`a la somme d\'efinissant $p^{j+2k} \,{\rm a}(n, {\rm T}_p F)$. Elle s'\'ecrit $p^{j+2k} \, \, \delta(n,p)\, \, {\rm a}(n ; F)$ o\`u 
\begin{equation}\label{defdeltanp}\delta(n,p) := \sum_{v \in V-\{0\}, \, \, \, \DET(v)=0} e^{\frac{2 i \pi}{p} {\rm tr} (vn)} \end{equation}
et $V \subset {\rm M}_2(\Z/p\Z)$ d\'esigne le sous-espace des matrices sym\'etriques. La forme quadratique $\DET : V \rightarrow \Z/p\Z$ admet $p+1$ droites isotropes. Si $a$ d\'esigne le nombre de ces droites qui sont dans le noyau de la forme lin\'eaire $v \mapsto  {\rm tr} (vn)$, de sorte que $a \in \{0, 1, 2, p+1\}$, alors $\delta(n,p) = (p-1) \cdot a - (p+1-a) = p (a -1) -1$. \end{pf}

\bigskip

\noindent Le mono\"ide $\rho_j({\rm M}_2(\Z) \cap \GL_2(\C))$ pr\'eservant le r\'eseau ${\rm Sym}^j \Z^2 \subset W_j$, on en tire le corollaire suivant.

\begin{cor}\label{scholieintegralitegenre2} Soient $j,k$ des entiers avec $j\geq 0$ et $k\geq 2$. Pour tout nombre premier $p$, les op\'erateurs de Hecke $p^{\frac{j+2k-6}{2}}\, {\rm K}_p$ et $p^{j+k-3}\,( {\rm T}_p + 1)$ pr\'eservent le sous-groupe de ${\rm S}_{j,k}$ constitu\'e des \'el\'ements dont tous les coefficients de Fourier sont \`a valeurs dans le sous-groupe ${\rm Sym}^j \Z^2 \subset W_j$. 
\end{cor}

\bigskip

\noindent {\it D\'emonstration de la proposition~\ref{intcongtaujk}}. Notons ${\rm S}_{j,k}^{\rm ent} \subset {\rm S}_{j,k}$ le sous-groupe d\'efini dans l'\'enonc\'e. Le $\C$-espace vectoriel ${\rm S}_{j,k}$ \'etant de dimension finie, il existe une partie finie $N \subset \mathcal{N}$ telle que l'application lin\'eaire $$F \mapsto ({\rm a}(n ; F))_{n \in N}, \, \, \, \,  {\rm S}_{j,k} \rightarrow W_j^N,$$ soit injective. Elle envoit le $\Z$-module ${\rm S}_{j,k}^{\rm ent}$ dans $({\rm Sym}^j \Z^2)^N$. Cela montre d'une part que le $\Z$-module ${\rm S}_{j,k}^{\rm ent}$ est libre de rang fini, et d'autre part que l'application naturelle 
$\eta : {\rm S}_{j,k}^{\rm ent} \otimes_\Z \C \rightarrow {\rm S}_{j,k}$ est injective, puisque c'est le cas de l'application naturelle ${\rm Sym}^j \Z^2 \otimes \C \rightarrow W_j$. \ps
Il est plausible qu'en toute g\'en\'eralit\'e l'application $\eta$ soit  
bijective, mais nous n'avons pas trouv\'e de r\'ef\'erence sur ce point. D\'emontrons-le lorsque $(j,k)$ est dans la liste de l'\'enonc\'e de la proposition \ref{intcongtaujk}. Dans ce cas, ${\rm S}_{j,k}$ est de dimension $1$, de sorte qu'il suffit de v\'erifier que ${\rm S}_{j,k}^{\rm ent}$ est non nul. Consid\'erons pour cela la forme modulaire $f_{j,k,u,v,w}= \vartheta_2({\rm E}_8, {\rm P}_{j,k,u,v,w})$ dans ${\rm S}_{j,k}$ que nous avons construite  au \S \ref{dimtsushima}.  Elle a ses coefficients dans $\Z[i] [X,Y]$, o\`u $\Z[i]$ d\'esigne l'anneau des entiers de Gauss, puisque le polyn\^ome harmonique ${\rm P}_{j,k,u,v,w}$ envoie ${\rm E}_8^2$ dans $\Z[i][X,Y]$ ; les quelques coefficients de Fourier non nuls de $f_{j,k,u,v,w}$ que nous avons d\'etermin\'es sont m\^eme dans $\Z$ (Table \ref{tablecoeffdjk}). On observe que pour tout $n \in \mathcal{N}$, on a $\overline{{\rm a}(n ; f_{j,k,u,v,w})}={\rm a}(n ; f_{j,k,\overline{u},\overline{v},\overline{w}})$, $z \mapsto \overline{z}$ d\'esignant respectivement la conjugaison complexe sur $\C[X,Y]$ et sur ${\rm E}_8 \otimes \C$. Ainsi, $f_{j,k,u,v,w}+f_{j,k,\overline{u},\overline{v},\overline{w}}$ est un \'el\'ement non nul de ${\rm S}_{j,k}^{\rm ent}$. \ps
Pour terminer la d\'emonstration de la proposition~\ref{intcongtaujk} , il suffit d'appliquer la proposition \ref{coeffsatsp4shimura} et le corollaire \ref{scholieintegralitegenre2}. $\square$ \ps \medskip

\noindent Soient $n = \left[\begin{array}{cc} n_{11} & \frac{n_{12}}{2} \\ \frac{n_{12}}{2} & n_{22} \end{array}\right] \in \mathcal{N}$  et $p$ premier. On constate les identit\'es 
{\small $$ \left[\begin{array}{cc} p & 0 \\ 0 & 1 \end{array}\right] \cdot \frac{n}{p} \, \, = \, \, \left[\begin{array}{cc} p n_{11} & \frac{n_{12}}{2} \\ \frac{n_{12}}{2} &  n_{22}/p \end{array}\right],$$ $$\left[\begin{array}{cc} 1 & d \\ 0 & p \end{array}\right]\cdot \frac{n}{p} \,\,= \, \, 
\left[\begin{array}{cc} \frac{n_{11}+d n_{12} + d^2 n_{22}}{p} & \frac{n_{12}}{2}+ d n_{22} \\  \frac{n_{12}}{2}+ d n_{22}  & p n_{22} \end{array}\right].$$}\par
\noindent Ainsi, si la forme quadratique sur $\Z^2$ d\'efinie par $n$ est anisotrope modulo le nombre premier $p$, alors ni $n/p$, ni l'une des deux matrices ci-dessus n'est dans $\mathcal{N}$. La proposition \ref{corKpTp} admet donc la cons\'equence suivante. \ps\ps 

\begin{scholie}\label{scholiearakawa} Supposons que $F \in {\rm S}_{j,k}$ est vecteur propre de l'op\'erateur ${\rm p}^{\frac{j+2k}{2}-3}\, {\rm K}_p$ de valeur propre $\lambda$. Si $n \in \mathcal{N}$ et $2n$ est une matrice de Gram d'une forme quadratique sur $\Z^2$ qui est anisotrope modulo le nombre premier $p$, on a la relation $\lambda \, \, \, {\rm a}(n ; F) = {\rm a}(pn ; F)$. En particulier, cette relation d\'etermine $\lambda$ si ${\rm a}(n ; F) \neq 0$.
\end{scholie}

\bigskip
\noindent 

Ce scholie s'applique par exemple lorsque
$$2 \,n \,\,= \,\,\left[ \begin{array}{cc} 2 & -1 \\ -1 & 2\end{array}
\right],$$ qui n'est autre que la matrice de Gram standard du r\'eseau ${\rm A}_2$, lorsque ${\rm A}_2 \otimes \F_p$ ne repr\'esente pas $0$, {\it i.e.} $p \equiv -1 \bmod 3$. On d\'eduit donc de la table
\ref{tablecoeffdjk}
le corollaire suivant.

\begin{cor}\label{taujk2} Les entiers $\tau_{6,8}(2)$, $\tau_{4,10}(2)$,
$\tau_{8,8}(2)$ et $\tau_{12,6}(2)$ valent respectivement $0, -1680, 1344$
et $-240$.
\end{cor}

\smallskip

\begin{remarque}\label{scholiearakawa2} {\rm Supposons que $F \in {\rm S}_{j,k}$ est vecteur propre de l'op\'erateur $2^{j+2k-5}({\rm T}_2+1)+2^{j+2k-3}$ de valeur propre $\lambda$ et posons {\footnotesize
$$n = \frac{1}{2} \left[ \begin{array}{cc} 2 & -1 \\ -1 & 2 \end{array}\right]\, \, \, \, \,  {\rm et}\, \, \, \, \,  m = \frac{1}{2} \left[ \begin{array}{cc} 2 & 0 \\ 0 & 6 \end{array}\right].$$}
Nous laissons au lecteur le soin de d\'eduire du corollaire~\ref{corKpTp} la relation 
{\scriptsize $$(\lambda - 2^{j+2\,k-4}) \, \, \, {\rm a}(n ; F) = 2^{k-2} \left( \rho_j(  \left[ \begin{array}{cc} 2 & 0 \\ -1 & 1 \end{array}\right]) + \rho_j(  \left[ \begin{array}{cc} -1 & 1 \\  2 & 0 \end{array}\right]) +  \rho_j(  \left[ \begin{array}{cc} 1 & -1 \\  1 & 1 \end{array}\right]) \right) \, \, {\rm a}(m; F).$$}
\noindent (V\'erifier en particulier $\delta(n,2) = -3$). En utilisant les valeurs de la table \ref{tablecoeffdjk}, cette formule permet de d\'emontrer que les entiers $\tau_{6,8}(4)$, $\tau_{4,10}(4)$,
$\tau_{8,8}(4)$ et $\tau_{12,6}(4)$ valent respectivement  $409600$, $-700160$, $348160$ et $4276480$. }
\end{remarque}



Il se trouve que les valeurs propres des $F_{j,k}$ ont \'et\'e \'etudi\'ees par 
Faber et van der Geer \cite{fabervandergeerIetII}\cite[\S
24]{vandergeer} de mani\`ere compl\`etement diff\'erente, en comptant {\it \`a la Deligne} des courbes de genre $2$
sur des corps fini. Si $q$ d\'esigne une puissance d'un nombre premier, ils ont pu ainsi d\'eterminer $\tau_{j,k}(q)$ pour tout $q \leq 37$
(en admettant toutefois une propri\'et\'e attendue de la cohomologie de certains faisceaux sur l'espace de Siegel de genre $2$ \cite[\S 24]{vandergeer}) ; ils donnent {\it loc.  cit.} un certain nombre de
valeurs, y compris la valeur $\tau_{j,k}(2)$ ci-dessus.  Nous
exposerons un peu plus loin une m\'ethode tr\`es diff\'erente permettant de
d\'eterminer $\tau_{j,k}(q)$, qui nous permettra de d\'emontrer le th\'eor\`eme suivant.\ps\ps 

\begin{thmv} \label{thmtablestjk} Soient $p$ un nombre premier et $(j,k)$ l'un des couples $(6,8)$, $(4, 10)$, $(8,8)$ ou $(12,6)$. \begin{itemize}\ps\ps  \medskip
\item[(i)] Si $p \leq 113$, l'entier $\tau_{j,k}(p)$ est donn\'e par la table  \ref{taujk}. \ps  \medskip
\item[(ii)] Si $p \leq 29$, l'entier $\tau_{{j,k}}(p^{2})$  est donn\'e par la table  \ref{taujkp2}. \ps\ps
\end{itemize}
\end{thmv}
\medskip\smallskip

\noindent Signalons enfin que dans un travail r\'ecent \cite[\S 8]{cvdg}, Cl\'ery et van der Geer ont r\'e-obtenu les valeurs $\tau_{6,8}(q)$ pour $q \leq 49$ par une m\'ethode encore diff\'erente.\ps \bigskip

\subsection{O\`u l'on explique l'apparition des
$\psi_{j,k}$ dans la table \ref{tableSO8}. }\label{remtripsijk}

Fixons l'un des quatres couples
$(j,k)$ de la table \ref{dimsjknonnul} tel que $j>0$.  Soit $U'_{j,k}$ la
repr\'esentation irr\'eductible de ${\rm SO}_8(\C)$ de plus haut poids
$(\frac{j}{2}+k-4) (\varepsilon_1 + \varepsilon_2) +
\frac{j}{2}(\varepsilon_3 + \varepsilon_4)$ (\S VI.\ref{exparlan}), elle se
factorise par ${\rm PGSO}_8(\C)$.  On dispose alors d'isomorphismes naturels
{\scriptsize $$ {\rm M}_{U'_{j,k}}({\rm SO}_8) \isomol {\rm
M}_{U'_{j,k}}({\rm PGSO}_8) \isomo {\rm M}_{U_{j,k}}({\rm PGSO}_8) \isomo
{\rm M}_{U_{j,k}}({\rm SO}_8) \isomol {\rm M}_{U_{j,k}}({\rm O}_8) \isomo
{\rm S}_{j,k}.$$ } \noindent En effet, le premier et le troisi\`eme sont
g\'en\'eraux (variante du lemme V.\ref{inflationmonGO} bas\'ee sur la
proposition IV.\ref{diagrammesimilitude}).  Le dernier isomorphisme est celui de la
proposition~\ref{calcimagetheta}.  L'avant dernier morphisme est injectif
pour des raisons g\'en\'erales (\S IV.\ref{fauton}), bijectif car $\dim {\rm
M}_{U_{j,k}}({\rm SO}_8)=1$ d'apr\`es \cite[Table 2]{chrenard2}.  Celui du milieu est
induit par la trialit\'e.  En effet, par un calcul laiss\'e au lecteur,
bas\'e sur l'action bien connue de la trialit\'e sur le diagramme de Dynkin
de type ${\bf D}_4$, on observe que si une $\C$-repr\'esentation
irr\'eductible de ${\rm PGSO}_8(\C)$ admet pour plus haut poids
$\sum_{i=1}^4 n_i \varepsilon_i$, alors ses deux repr\'esentations triales
ont pour plus haut poids $\sum_{i=1}^4 m_i \varepsilon_i$ o\`u 
$(m_1,m_2,m_3,\pm m_4)$ vaut {\small $$(\frac{n_1+n_2+n_3+n_4}{2},\frac{n_1+n_2-n_3-n_4}{2},
\frac{n_1-n_2+n_3-n_4}{2},|\frac{n_1-n_2-n_3+n_4}{2}|).$$} 
L'apparition des $\psi_{j,k}$ dans la table \ref{tableSO8} est donc cons\'equence de la suite d'isomorphismes ci-dessus et du th\'eor\`eme
VII.\ref{trialitegenus2} (i).  \ps\ps 

Mentionnons enfin que ce paragraphe sugg\`ere une m\'ethode alternative pour d\'eterminer les $\tau_{j,k}(q)$, bas\'ee sur un calcul de valeurs propres d'op\'erateurs de Hecke  pour ${\rm O}_8$. Nous renvoyons au travail \`a venir de M\'egarban\'e \cite{megarbane} \`a ce sujet.

\section{$\Pi_{\rm disc}({\rm SO}_{24})$ et la conjecture de Nebe-Venkov}\label{so24etnv}

\subsection{Une caract\'erisation de la table \ref{table24}} \noindent Consid\'erons le sous-ensemble suivant de $\coprod_{n\geq 1} \Pi_{\rm alg}(\PGL_n)$ :
$$\Pi = \{ 1, \Delta_{11}, \Delta_{15}, \Delta_{17}, \Delta_{19}, \Delta_{21},
{\rm Sym}^2 \Delta_{11}, \Delta_{19,7}, \Delta_{21,5}, \Delta_{21,9},
\Delta_{21,13}\}.$$
\noindent La proposition ci-dessus donne une ``d\'efinition'' directe de la table \ref{table24}. \ps\ps

\begin{propv}\label{propliste24} L'ensemble des $\psi \in \mathcal{X}({\rm SL}_{24})$ tels que : \begin{itemize}\ps\ps

\item[(i)] $\psi_\infty$ admet pour
valeurs propres simples les entiers $\pm 11, \pm 10, \dots, \pm 1$, ainsi
que $0$ pour valeur propre double,\ps\ps

\item[(ii)] $\psi$ est de la forme $\oplus_{i=1}^k \pi_i[d_i]$ avec $\pi_i \in
\Pi$ pour tout $i$,\ps\ps 
\end{itemize}
est exactement celui donn\'e par la table \ref{table24}. Il a $24$ \'el\'ements. \ps\ps
\end{propv}

\begin{pf} C'est un simple exercice de combinatoire que l'on peut traiter de la mani\`ere suivante. Consid\'erons plus g\'en\'eralement, pour tout entier $n \geq 1$, l'ensemble $\Psi_n$ des \'el\'ements $\psi \in \mathcal{X}({\rm SL}_{n})$ verifiant l'assertion (ii) de l'\'enonc\'e et tels que : \ps \ps

-- $\psi_\infty$ admet pour valeurs propres simples les entiers $\pm \frac{n-1}{2}, \pm \frac{n-3}{2}, \dots,\pm 1$ et $0$ si $n$ est impair,\ps\ps
-- $\psi_\infty$ admet pour valeurs propres simples les entiers $\pm \frac{n-2}{2}, \pm \frac{n-4}{2}, \dots,\pm 1$, ainsi que $0$ pour valeur propre double, si $n$ est pair. \ps \ps

La question est de d\'eterminer $\Psi_{24}$. Nous allons plus g\'en\'eralement expliciter $\Psi_n$ pour tout $1 \leq n \leq 24$ par r\'ecurrence sur $n$. Pour $c \in  \mathcal{X}({\rm SL}_{a})$ et $\Psi \subset  \mathcal{X}({\rm SL}_{b})$ il sera commode de noter $c \oplus \Psi$ l'ensemble des \'el\'ements de $\mathcal{X}({\rm SL}_{a+b})$ de la forme $c \oplus \psi $ avec $\psi \in \Psi$. \ps\ps

Soient $1 \leq n\leq 24$ un entier pair et $\psi \in \Psi_n$. \'Ecrivons $\psi = \oplus_{i=1}^k \pi_i[d_i]$ comme dans l'assertion (ii).  L'in\'egalit\'e $n \leq 24$ entra\^ine que pour tout $i$, les valeurs propres de $(\pi_i[d_i])_\infty$ sont $\leq 11$. Fixons un entier $i$ tel que $(\pi_i[d_i])_\infty$ admette la valeur propre $0$. Le lemme \ref{lemmecombigrandpi} montre que l'on est dans l'un des cas suivants :  \ps \ps 
-- $\pi_i=1$, \ps\ps
-- $\pi_i=\Delta_{11}$, $d_i=12$, $n=24$ et donc $\psi = \Delta_{11}[12]$, \ps\ps

-- Si $\pi_i = {\rm Sym}^2 \Delta_{11}$, $d_i=1$ et $n=24$. \ps\ps 

\noindent En particulier, on constate que sous l'hypoth\`ese $n \leq 22$, il existe deux entiers $i$ tels que $\pi_i=1$, un et un seul d'entre eux v\'erifiant de plus $d_i=1$. On a donc l'\'egalit\'e $\Psi_n = [1] \oplus \Psi_{n-1}$ pour $n$ pair $\leq 22$. De plus, cette analyse montre $$\Psi_{24} \,=\, \{\,\,\Delta_{11}[12]\,\,\} \hspace{4pt}\cup \hspace{4pt} [1] \oplus \Psi_{23} \hspace{4pt} \cup \hspace{4pt} {\rm Sym}^2 \Delta_{11} \oplus \Psi_{21}.$$
  \ps\ps

\noindent Il ne reste qu'\`a d\'ecrire $\Psi_n$ pour $n$ impair $\leq 23$. Les \'el\'ements non triviaux de $\Pi$ \'etant de poids motivique $\geq 11$, on a \'evidemment  $\Psi_n = \{ \,\,[n]\,\, \}$ pour tout entier $1\leq n \leq 11$ impair. La seule repr\'esentation dans $\Pi$ de poids motivique $<15$ \'etant $\Delta_{11}$, on a de plus
$$\Psi_{13}\,=\, \{\,\, [13], \, \, \Delta_{11}[2] \,\oplus \,[9]\,\,\}\, \, \, \, {\rm et} \,\,\, \, \Psi_{15} \,=\,\{\,\,[15], \, \, \Delta_{11}[4]\, \oplus\, [7]\,\,\},$$ 
On constate de m\^eme  les assertions suivantes :  \ps \ps 

--\hspace{4pt} $\Psi_{17}$ est r\'eunion de $\{ \,\,[17], \, \, \Delta_{11}[6] \oplus [5]\, \, \}$ et $\Delta_{15}[2] \oplus \Psi_{13}$, \ps \ps

--\hspace{4pt} $\Psi_{19}$ est r\'eunion de  $ \{\,\,[19],\,\, \Delta_{11}[8] \oplus [3], \, \, \Delta_{15}[4] \oplus [11]\,\,\}$ et $\Delta_{17}[2] \oplus \Psi_{15}$, \ps \ps

--\hspace{4pt} $\Psi_{21}$ est r\'eunion des ensembles $\Delta_{17}[4] \oplus \Psi_{13}$, $\Delta_{19}[2] \oplus \Psi_{17}$, et {\small $$\{ \,\,[21], \,\, \Delta_{11}[10] \oplus [1],\,\,\Delta_{15}[6] \oplus [9], \,\, \Delta_{19,7}[2]  \oplus \Delta_{15}[2]\oplus \Delta_{11}[2] \oplus [5]\, \, \},$$} 
\indent --\hspace{4pt} $\Psi_{23}$ est r\'eunion de $\{ \,\,[23], \, \,\Delta_{15}[8] \oplus [7], \, \,  \Delta_{17}[6] \oplus [11], \, \,  {\rm Sym}^2 \Delta_{11} \oplus \Delta_{11}[10] \,\,\}$, $\Delta_{19}[4] \oplus \Psi_{15}$, $\Delta_{21}[2] \oplus \Psi_{19}$, et de l'ensemble :
{\small $$\{\Delta_{21,5}[2] \oplus \Delta_{17}[2] \oplus \Delta_{11}[4] \oplus [3], \, \, \Delta_{21,9}[2] \oplus \Delta_{15}[4] \oplus [7], \, \,  \Delta_{21,13} \oplus \Delta_{17}[2] \oplus [11]\,\,\}.$$} 

Cette analyse montre au final que l'ensemble $\Psi_{24}$ est l'ensemble des $24$ \'el\'ements de la table \ref{table24}. Une mani\`ere de justifier que ces $24$ \'el\'ements sont distincts consiste \`a invoquer la proposition VI.\ref{jacquetshalika} (Th\'eor\`eme de Jacquet-Shalika). On observera que l'intersection de $[1] \oplus \Psi_{23}$ et ${\rm Sym}^2 \Delta_{11} \oplus \Psi_{21}$ est le singleton  $\{\, {\rm Sym}^2 \Delta_{11} \oplus \Delta_{11}[10] \oplus [1] \}$. Plus g\'en\'eralement, la table \ref{cardpsin} indique le cardinal de $\Psi_n$ en fonction de $n \leq 24$. (Une autre m\'ethode pour justifier que les $24$ \'el\'ements de la table \ref{table24} sont distincts, sur laquelle nous reviendrons au \S \ref{enoncesconjnv}, consisterait \`a observer que les composantes en le nombre premier $2$ de ces \'el\'ements ont des traces distinctes. )
\end{pf}

\begin{table}[h!]
\begin{center}
\renewcommand{\arraystretch}{1.5}
 \begin{tabular}{c||c|c|c|c|c|c|c|c|c|c|c|c|c|}
$n$ & $ \leq 12 $ & $13$ & $14$ & $15$ & $16$ & $17$ & $18$ & $19$ & $20$ & $21$ & $22$ & $23$ & $24$  \cr
\hline
$|\Psi_n|$ & $1$ & $2$ & $2$ & $2$ & $2$ & $4$ & $4$ & $5$ & $5$ & $10$ & $10$ & $14$ & $24$ \cr \end{tabular} 
\caption{Cardinal du sous-ensemble $\Psi_n \subset \mathcal{X}({\rm SL}_n)$ introduit dans la d\'emonstration de la proposition \ref{propliste24}.}
\label{cardpsin}
\end{center}
\end{table}
\renewcommand{\arraystretch}{1}
\ps \medskip
\noindent 


\begin{lemme}\label{lemmecombigrandpi} Soient $\pi \in \Pi - \{1\}$, $d$ un entier $\geq 1$, $\psi = \pi[d]$ et $\Lambda \subset \C$ l'ensemble des valeurs propres de $\psi_\infty$. On suppose $\Lambda \subset \Z$ et $|\lambda| \leq 11$ pour tout $\lambda \in \Lambda$. Alors on a $$|\lambda| \geq \frac{d-1}{2}$$ pour tout $\lambda \in {\rm Poids}(\pi)$. De plus, les assertions suivantes sont satisfaites\, : \begin{itemize} \ps\ps
\item[(i)] Si $0 \in \Lambda$, alors on a $\pi = {\rm Sym}^2 \Delta_{11}$ et $d=1$, ou $\pi = \Delta_{11}$ et $d=12$.\ps \ps
\item[(ii)] Si $1 \in \Lambda$, alors on a $\pi = \Delta_{11}$ et $d \in \{10,12\}$. \ps\ps
\end{itemize}
\end{lemme}

\begin{pf} Il s'agit d'une inspection imm\'ediate de la liste $\Pi$.\end{pf}

\subsection{\'Enonc\'es et survol des d\'emonstrations} \label{enoncesconjnv}

Soit $\psi$ l'un des $24$ \'el\'ements list\'es dans la table \ref{table24}. D'apr\`es la proposition \ref{propliste24} et les exemples du \S VI.\ref{exparlan}, on a $\psi_\infty = {\rm St}({\rm inf}_V)$, o\`u $V$ est la repr\'esentation triviale de ${\rm SO}_{24}(\R)$. L'ensemble $\Pi$ \'etant constitu\'e de repr\'esentations autoduales, l'\'enonc\'e suivant (qui est aussi le th\'eor\`eme I.\ref{thmintro24}) n'est pas absurde !\ps

\begin{thmv}\label{thm24} Les $24$ repr\'esentations $\pi \in \Pi_{\rm
disc}({\rm O}_{24})$ telles que $\pi_\infty = 1$ ont pour param\`etres
standards $\psi(\pi,{\rm St})$ les $24$ \'el\'ements de la table \ref{table24}.  \end{thmv}
\ps

Soulignons que dans son travail~\cite{ikeda2} (reposant notamment sur \cite{bfw}, \cite{ikeda1} et \cite{nebevenkov}), Ikeda \'etait parvenu \`a
identifier $20$ de ces $24$ param\`etres, \`a savoir tous ceux de cette
liste ne ``contenant pas'' l'une des $4$ repr\'esentations $\Delta_{w,v}$.  \ps \ps

\ps\ps

Afin d'en dire un peu plus, rappelons quelques notations introduites au \S \ref{thetaniemeier}. On
note $\lambda_i$, $i=1,\dots,24$, les $24$ valeurs propres distinctes de
l'op\'erateur ${\rm T}_2$ agissant sur $\C[{\rm X}_{24}]$, rang\'ees par
ordre d\'ecroissant \`a la mani\`ere de Nebe et Venkov (Table
\eqref{tablevp2}).  Fixons ${\rm v}_i \in \Z[{\rm X}_{24}]$ un vecteur propre de
${\rm T}_2$, et donc de ${\rm H}({\rm O}_{24})$, associ\'e \`a $\lambda_i$.  On
note $\pi_i \in \Pi_{\rm disc}({\rm O}_{24})$ la repr\'esentation
engendr\'ee par ${\rm v}_i$ et $$\psi_i=\psi(\pi_i,{\rm St})$$ son param\`etre
standard.  Le th\'eor\`eme~\ref{thm24} affirme que ces $24$ param\`etres $\psi_i$ sont ceux de la table \ref{table24}. 
\'Etant donn\'e que les $24$ valeurs propres de ${\rm T}_2$ sont
diff\'erentes, la relation $$\lambda_i \, \, = \, \, 2^{11}\, {\rm trace} \,
\, (\psi_i)_2 \,=\, 2^{11}\, {\rm trace}\, {\rm St}\,{\rm c}_2(\pi_i)$$
caract\'erise uniquement l'application $i \mapsto \psi_i$. Cela fournit par
ailleurs une premi\`ere v\'erification du th\'eor\`eme~\ref{thm24} car on
peut montrer que les $4$ valeurs $\tau_{j,k}(2)$ 
(Corollaire \ref{taujk2}, \cite[\S 24, \S 27]{vandergeer}), ainsi que les coefficients $\tau_k(2)$ en $q^2$ des formes modulaires normalis\'ees
pour ${\rm SL}_2(\Z)$ de poids $k\leq 22$, sont compatibles avec le calcul des $\lambda_i$ par Nebe et
Venkov. \ps\ps 

En guise d'exemple consid\'erons le param\`etre $$\psi = \Delta_{21,5}[2] \oplus \Delta_{17}[2] \oplus
\Delta_{11}[4] \oplus [1] \oplus [3].$$ On a $\tau_{4,10}(2)=-1680$,
$\tau_{12}(2)=-24$ et $\tau_{18}(2)=-528$, de sorte que 
$$2^{11}\, {\rm trace}({\rm St}(\psi_2)) = 
(1+2) \cdot (-1680)+2^2\cdot (1+2)\cdot (-528)$$ $$+2^4\cdot (1+2+2^2+2^3)\cdot
(-24)+2^{11}+2^{10}\cdot(1+2+2^2)=-7920.$$
On retrouve bien la valeur propre $\lambda_{21}$ de Nebe-Venkov, {\it i.e.
$\psi=\psi_{21}$}. \ps\ps 
On note enfin $g_i$ le degr\'e de ${\rm v}_i$, d\'efini au \S \ref{thetaniemeier} : la valeur propre $\lambda_i$ ayant multiplicit\'e $1$, c'est le
plus petit entier $g$ tel que $\pi_i$ admette un $\vartheta$-correspondant
dans $\Pi_{\rm cusp}({\rm Sp}_{2g})$ au sens du \S VII.\ref{releichlerbisrallis}. Par convention on a $g_1=0$. Comme nous l'avons d\'ej\`a expliqu\'e au \S \ref{thetaniemeier}, les $g_i$ ont \'et\'e d\'etermin\'es par Nebe
et Venkov pour $i \neq 19,21$ dans \cite{nebevenkov} ;  ils ont de plus conjectur\'e
$g_{19}=9$ et $g_{21}=10$. 

\begin{thmv}\label{thmnv} \begin{itemize} \item[(i)] Les $g_i$ sont donn\'es par la table~\ref{table24bis}. En particulier, la conjecture de Nebe-Venkov {\rm \cite[\S
3.3]{nebevenkov}} est vraie.  \ps\ps 
\item[(ii)] Pour tout $i\leq 23$, $g_i$ est le plus petit entier $m\geq 0$ tel que $\psi_i$ soit de la forme 
$[23-2 m] \oplus \psi'_i$  avec $\psi'_i \in \mathcal{X}_{\rm AL}({\rm SL}_{2 m+1})$. Enfin, on a $g_{24}=12$. 
\end{itemize}
\end{thmv}

On prendra garde que dans notre table \ref{table24bis}, les
degr\'es $g_i$ sont rang\'es par ordre croissant, mais pas tout \`a fait les indices
$i$. Nous allons d'abord d\'emontrer au \S \ref{thhm24implthmnv} que le th\'eor\`eme \ref{thm24} entra\^ine
le th\'eor\`eme \ref{thmnv}. Nous donnerons ensuite trois d\'emonstrations
du th\'eor\`eme \ref{thm24}, les deux premi\`eres \'etant conditionnelles :  \ps\ps

-- La premi\`ere, sans doute la plus naturelle, est une application directe de la fomule
de multiplicit\'e d'Arthur pour ${\rm SO}_{24}$. Son inconv\'enient \'evident est
qu'elle d\'epend de l'\'etablissement de la formule de multiplicit\'e d'Arthur pour les $\Z$-groupes
${\rm SO}_n$, non encore r\'edig\'e, ainsi que de l'analogue \`a ces
groupes de la conjecture VIII.\ref{conjaj}, que nous avons englob\'es en une seule
conjecture VIII.\ref{conjaj2} au chapitre pr\'ec\'edent. Cette d\'emonstration conditionnelle est expliqu\'ee au \S \ref{preuve1thm24}. \ps\ps

-- Nous donnerons ensuite au \S \ref{preuve2thm24} une seconde d\'emonstration conditionnelle, n'utilisant cette fois-ci que
la th\'eorie d'Arthur pour les groupes de Chevalley, la conjecture VIII.\ref{conjaj}, et des arguments de ``correspondance $\vartheta$''. Dans cette seconde approche, nous d\'emontrerons en fait en m\^eme temps les deux th\'eor\`emes ci-dessus. \ps\ps

-- Nous donnerons enfin une derni\`ere d\'emonstration du th\'eor\`eme \ref{thm24}, cette fois-ci inconditionnelle, au \S \ref{demothmE}. Cette d\'emonstration, assez diff\'erente dans l'esprit et d\'ej\`a discut\'ee dans l'introduction, ne reposera plus du tout sur la formule de multiplicit\'e d'Arthur. Elle fournira une justification plus profonde de l'\'enonc\'e du th\'eor\`eme \ref{thm24}.  \ps\ps 

\begin{remarque} \label{remt3} {\rm \begin{itemize} \item[(i)] Curieusement,
l'op\'erateur de Hecke ${\rm T}_3$ sur $\C[{\rm X}_{24}]$ admet la valeur
propre $1827360$ avec multiplicit\'e $2$.  Cela traduit l'\'egalit\'e un peu
miraculeuse des traces des composantes en $p=3$ des param\`etres $\psi_{19}=
{\rm Sym}^2 \Delta_{11} \oplus \Delta_{19,7}[2] \oplus \Delta_{15}[2] \oplus
\Delta_{11}[2] \oplus [5]$ et $\psi_{21}=\Delta_{21,5}[2] \oplus
\Delta_{17}[2] \oplus \Delta_{11}[4] \oplus [1] \oplus [3]$, comme on le
v\'erifie \`a l'aide la table \ref{taujk}.  \ps\ps  \item[(ii)] L'op\'erateur de
Hecke ${\rm T}_2$ agissant sur $\C[{\rm X}_{32}]$ admet des valeurs propres
non enti\`eres.  En effet, soit $\Delta_{23}$ l'une des deux formes propres
normalis\'ees de poids $24$ pour ${\rm SL}_2(\Z)$.  Il est bien connu que
ses coefficients de Fourier sont dans $\Q(\sqrt{144169})$, par exemple le
second vaut $540 \pm 12 \sqrt{144169}$.  On conclut car d'apr\`es Ikeda et
B\"ocherer (\S \ref{ikedabocherer}), il existe $\pi \in \Pi_{\rm disc}({\rm
O}_{32})$ telle que $\psi(\pi,{\rm St})=\Delta_{23}[8]\oplus [15] \oplus
[1]$.  \end{itemize} } \end{remarque}

\subsection{Le th\'eor\`eme~\ref{thm24} entra\^ine le
th\'eor\`eme~\ref{thmnv}.}\label{thhm24implthmnv} \ps\ps

\noindent Nous allons d\'eduire de la liste des $\psi_i$ dans la
table~\ref{table24bis} la valeur des $g_i$. \ps\ps

\begin{lemme}\label{lemmegenusfaible} On a les in\'egalit\'es $g_{23} \leq 11$ et $g_i \leq 10$ pour $i \leq 22$,  ainsi que l'\'egalit\'e $g_{24}=12$.
\end{lemme}

\begin{pf} En effet, ainsi qu'on l'a rappel\'e au \S
\ref{thetaniemeier}, Erokhin a montr\'e $g_i \leq 12$ pour tout $i$ dans \cite{erokhin}. Ce r\'esultat a \'et\'e r\'e-obtenu par Borcherds,
Freitag et Weissauer dans \cite{bfw}, qui ont v\'erifi\'e de plus l'in\'egalit\'e $g_i \leq 11$ pour tout $i$ sauf
pour exactement l'un d'entre eux (c'est l'assertion que ${\rm Ker}\, \vartheta_{11}$ est
de dimension $1$). Ils calculent pour cela explicitement les coefficients des s\'eries
th\^eta des r\'eseaux de Niemeier en les matrices de Gram des r\'eseaux de
la forme ${\rm Q}(R)$ o\`u $R$ est un syst\`eme de racines irr\'eductible
de type {\rm ADE} et de rang $\leq 12$ \cite[Table p. 146]{bfw}. Comme l'ont observ\'e Nebe et
Venkov \cite[\S 3.1, Lemma 3.3]{nebevenkov}, ces calculs montrent plus exactement les in\'egalit\'es $g_{23}
\leq 11$, $g_i \leq 10$ pour tout $i \leq 22$, et $g_{24}=12$ (au moins l'un des $g_i$
doit \^etre \'egal \`a $12$ d'apr\`es le r\'esultat de \cite{bfw}
sus-mentionn\'e).
\end{pf}

\ps\ps

\begin{pf} (Le th\'eor\`eme~\ref{thm24} entra\^ine le
th\'eor\`eme~\ref{thmnv}) Soit $1 \leq i \leq 23$ et soit $\psi_i' \in \mathcal{X}({\rm SL}_{2g_i+1})$ le param\`etre standard du
$\vartheta$-correspondant dans $\Pi_{\rm cusp}({\rm Sp}_{2g_i})$ de
$\pi_i$. D'apr\`es le lemme \ref{lemmegenusfaible} on a $g_i \leq 11$, de sorte que d'apr\`es Rallis (corollaire VII.\ref{corparamrallis}) on a l'identit\'e
$$\psi_i = \psi'_i \oplus [23-2 g_i].$$
D'apr\`es le th\'eor\`eme VIII.\ref{arthurst} d'Arthur, on a $\psi'_i \in
\mathcal{X}_{\rm AL}({\rm SL}_{2 g_i+1})$. L'unicit\'e des param\`etres
d'Arthur-Langlands (Jacquet-Shalika, Proposition
V.\ref{jacquetshalika}) montre donc que $g_i$ a la propri\'et\'e que
$[23-2 g_i]$ est un ``constituant'' de $\psi_i$ dans la table \ref{table24bis}. \ps\ps 

Si $i\leq 22$, le lemme \ref{lemmegenusfaible} entra\^ine $g_i\leq 10$, soit encore $23-2 g_i\geq 3$. Dans ce cas, $23-2 g_i$ est
l'unique entier $m_i>1$ tel que $[m_i]$ est un constituant de $\psi_i$, l'autre constituant \'eventuel de $\psi_i$ de la forme $[d]$ \'etant $[1]$. Cela
d\'etermine trivialement $g_i$ pour $i\leq 22$ par inspection directe de la
table \ref{table24bis}. L'in\'egalit\'e $g_{23}\leq 11$ et l'identit\'e $\psi_{23} = {\rm Sym}^2
\Delta_{11} \oplus \Delta_{11}[10] \oplus [1]$ montrent de m\^eme que l'on a n\'ecessairement $g_{23}=11$.  Cela termine la d\'emonstration. \end{pf}

\ps\ps
Observons que l'approche utilis\'ee ici ne repose pas sur les calculs fins
des \S 3.3 et \S 3.4 de \cite{nebevenkov}, mais ``seulement'' de la table \cite[p. 146]{bfw} de Borcherds-Freitag-Weissauer et du calcul de ${\rm T}_2$ par Nebe et Venkov.

\subsection{Premi\`ere d\'emonstration, conditionnelle, du th\'eor\`eme
\ref{thm24}}\label{preuve1thm24} Admettons la conjecture VIII.\ref{conjaj2} et appliquons le th\'eor\`eme
VIII.\ref{amfexplso0} \`a chaque param\`etre $\psi$ de la table
\ref{table24} (on est dans ses hypoth\`eses avec $\psi_\infty={\rm St}({\rm
Inf}_1)$).  Nous affirmons que l'\'egalit\'e VIII.\eqref{critso8n} est
toujours satisfaite, ce qui est un ph\'enom\`ene en soit assez miraculeux. 
C'est bien entendu quelque chose que l'on peut simplement v\'erifier
b\^etement dans chacun des $24$ cas. On peut aussi faire les remarques
suivantes : \ps\ps

\begin{itemize} 

\item[(a)] Si $\psi$ ne contient pas ${\rm Sym}^2 \Delta_{11}$, et $\psi
\neq \Delta_{11}[12]$, alors $\psi$
satisfait l'hypoth\`ese du crit\`ere VIII.\ref{critamfsonA}, c'est-\`a-dire qu'il
est de la forme $$(\oplus_{i=1}^{k-2} \pi_i[d_i])\oplus [d_{k-1}] \oplus [1]$$ avec
$\pi_i$ est symplectique pour tout $i\leq k-2$. Une rapide inspection
montre que l'on a toujours soit $d_i<d_{k-1}$ et
$\varepsilon(\pi_i)=(-1)^{\frac{n_id_i}{4}}$, soit $d_i>d_{k-1}$ et $d_i
\equiv 0 \bmod 4$ (ce qui ne se produit que pour $\psi_{13},\psi_{21}$ et $\psi_{22}$).
Concr\`etement, on observe que le facteur epsilon $\varepsilon(\Delta_{j,k})=(-1)^k$ vaut $1$ pour les $4$ couples
$(j,k)$ d'int\'er\^et, et que chaque fois qu'un constituant de la forme $\Delta_w[d]$ appara\^it dans $\psi$ alors on a $d
\equiv w+1 \bmod 4$. On conclut donc par le crit\`ere
VIII.\ref{critamfsonA}. \ps\ps

\item[(b)] Si $\psi$ contient ${\rm Sym}^2 \Delta_{11}$,  alors $\psi$ est de
la forme $$(\oplus_{i=1}^{k-2} \pi_i[d_i])\oplus {\rm Sym}^2 \Delta_{11}
\oplus [1]$$ avec $\pi_i$ symplectique pour tout $i\leq k-2$; il satisfait
donc l'hypoth\`ese du crit\`ere VIII.\ref{critamfsonB}. On conclut encore
par ce crit\`ere en observant que l'on a toujours soit $d_i<d_{k-1}$ et
$\varepsilon(\pi_i)(-1)^{\frac{n_i}{2}}=(-1)^{\frac{n_id_i}{4}}$, soit $d_i>d_{k-1}$ et
$d_i \equiv 2 \bmod 4$ (ce qui ne se produit que pour $\psi_{20}$ et
$\psi_{23}$). \ps\ps 
\end{itemize}
\ps \medskip
\noindent Ainsi, pour tout $i \leq 23$ l'unique repr\'esentation $\pi'_i \in
\Pi({\rm SO}_{24})$ telle que $\psi(\pi'_i,{\rm St})=\psi_i$ satisfait ${\rm
m}(\pi'_i)=1$. Enfin, si l'on a $\psi=\psi_{20}=\Delta_{11}[12]$ l'\'egalit\'e VIII.\eqref{critso8n} est
trivialement satisfaite, et la remarque VIII.\ref{remmult2} assure qu'il existe
exactement deux repr\'esentations $\pi'_{24},\pi^{''}_{24}={\rm S}(\pi'_{24}) \in
\Pi_{\rm disc}({\rm SO}_{24})$ de param\`etre standard $\psi_{24}$, et
que de plus ${\rm m}(\pi'_{24})={\rm m}(\pi^{''}_{24})=1$. \ps\ps

Cette discussion nous fournit $25$ \'el\'ements distincts de $\Pi_{\rm
disc}({\rm SO}_{24})$, chacun \'etant de multiplicit\'e $1$.
Comme ${\rm h}({\rm SO}_{24})=|\widetilde{{\rm X}}_{24}|=25$ (Corollaire
IV.\ref{hson}), ce sont exactement les $\pi \in \Pi_{\rm
disc}({\rm SO}_{24})$ telles que $\pi_\infty =1$. La d\'ecomposition ${\rm
H}({\rm O}_{24})$-\'equivariante (\S IV.\ref{fauton})
$${\rm M}_\C({\rm SO}_{24}) = {\rm M}_\C({\rm O}_{24}) \oplus {\rm
M}_{\DET}({\rm O}_{24}),$$
combin\'ee au fait que ${\rm T}_2$ a $24$ valeurs propres distinctes sur
${\rm M}_\C({\rm O}_{24})$, montre que les $24$ repr\'esentations $\pi \in
\Pi_{\rm disc}({\rm O}_{24})$ telles que $\pi_\infty =1$ ont exactement pour
param\`etres standards les $\psi_i$ de la table~\ref{table24}, et aussi que
l'unique repr\'esentation $\pi \in \Pi_{\rm disc}({\rm O}_{24})$ telle que
$\pi_\infty=\DET$ a pour param\`etre $\Delta_{11}[12]$ (autrement dit, on
retrouve \'egalement la proposition VII.\ref{leechmoins2} !).

\subsection{Deuxi\`eme d\'emonstration du
th\'eor\`eme \ref{thm24}, modulo la conjecture VIII.\ref{conjaj}}\label{preuve2thm24}

Donnons maintenant une d\'emonstration ``moins conditionnelle'' du th\'ero\`eme \ref{thm24}, qui
se passe de la conjecture VIII.\ref{conjaj2}, et repose seulement sur la conjecture VIII.\ref{conjaj}. 
Commen\c{c}ons par deux
observations. \ps\ps

\noindent {\it Observation 1.} Soit $1 \leq g < 12$ et consid\'erons  l'application du \S \ref{thetasiegel} $$\vartheta_g : \C[{\rm X}_{24}] \longrightarrow
{\rm M}_{12}({\rm Sp}_{2g}(\Z)).$$ Supposons que $F
\in {\rm S}_{12}({\rm Sp}_{2g})$ est propre pour ${\rm H}({\rm Sp}_{2g})$ et
notons $\pi_F \in \Pi_{\rm disc}({\rm Sp}_{2g})$ la repr\'esentation
engendr\'ee.  D'apr\`es le th\'eor\`eme VIII.\ref{arthurst} d'Arthur, on
peut \'ecrire $$\psi(\pi_F,{\rm St})=\oplus_{j=1}^k \pi_j[d_j] \in
\mathcal{X}_{\rm AL}({\rm SL}_{2g+1}).$$ D'apr\`es B\"ocherer, $F$ est dans
l'image de $\vartheta_g$ si, et seulement si, on a 
${\rm L}(12-g,\pi_F,{\rm St})
\neq 0$ (Remarque VII.\ref{rembocherer}).  
D'apr\`es la proposition VIII.\ref{checkcomptheta}, on a
${\rm L}(12-g,\pi_F,{\rm St}) \neq 0$ si, et
seulement si, 
\begin{equation}\label{critbfinpreuve} \forall \,\,1\, \leq\, j\, \leq k \, \, {\rm tel}\, \,  {\rm que}\, \,
d_j \equiv 0 \bmod 2
\, \, {\rm et}\, \,  d_j > 23 -2 g, \, \, \, {\rm L}(\frac{1}{2},\pi_j) \neq 0. 
\end{equation}
Si cette condition est
satisfaite, la repr\'esentation $\pi_F$ admet donc un
$\vartheta$-correspondant $\pi'_F$ dans $\Pi_{\rm disc}({\rm O}_{24})$ tel que
$(\pi_F')_\infty =1$, et $$\psi(\pi'_F,{\rm St})=\psi(\pi_F,{\rm St}) \oplus
[23- 2 g]$$
d'apr\`es Rallis (Corollaire VII.\ref{corparamrallis}). \ps\ps

\noindent {\it Observation 2.} Regardons la table~\ref{table24bis}. On constate que pour 
tout $2 \leq i\leq 23$,  il existe un unique $\psi'_i \in \mathcal{X}_{\rm AL}({\rm SL}_{2 g_i +1 })$
tel que l'on ait $$\psi_i = \psi'_i \oplus [23 - 2 g_i].$$  Clairement,
$(\psi'_i)_\infty$ admet pour valeurs propres les $2 g_i+1$ entiers $\pm 11, \pm
10, \dots, \\ \pm (12-g_i)$ ainsi que $0$. 
Soit $\varpi_i \in \Pi({\rm
Sp}_{2g_i})$ l'unique repr\'esentation telle que $\psi(\varpi_i,{\rm
St})=\psi'_i$ et telle que $(\varpi_i)_\infty \simeq \pi'_{\DET^{12}}$ (la
s\'erie discr\`ete holomorphe introduite au \S VI.\ref{carinf}). Comme nous
l'avons d\'ej\`a expliqu\'e \`a plusieurs reprises, la multiplicit\'e ${\rm m}(\varpi_i)$ est non nulle si, et seulement si, il existe une forme propre $F_i \in {\rm S}_{12}({\rm
Sp}_{2g_i}(\Z))$ telle que $\pi_{F_i} \simeq \varpi_i$ (Corollaire
VI.\ref{corapiw}). \ps\ps

Ces deux observations sugg\`erent la strat\'egie optimiste suivante pour d\'emontrer le
th\'eor\`eme~\ref{thm24}. \begin{itemize} \ps\ps  

\item[1.] Montrer ${\rm m}(\varpi_i) \neq 0$ pour tout $2\leq i\leq 23$.\ps\ps 
\item[2.] V\'erifier, par le crit\`ere de B\"ocherer \eqref{critbfinpreuve}, 
pour tout $1 \leq i \leq 23$, que s'il existe une forme propre $F_i \in {\rm S}_{12}({\rm Sp}_{2g_i}(\Z))$
telle que $\pi_{F_i} \simeq \varpi_i$, alors $F_i \in {\rm
Im}(\vartheta_{g_i})$.
\ps\ps 
\end{itemize}

\noindent En effet, ceci \'etant fait, nous en d\'eduisons pour tout $2 \leq i \leq
23$ l'existence d'une repr\'esentation dans $\Pi_{\rm disc}({\rm
O}_{24})$ de composante archim\'edienne triviale et dont le param\`etre
standard est $\psi_i$, \`a savoir un $\vartheta$-correspondant de $\varpi_i$.
L'existence d'une repr\'esentation dans $\Pi_{\rm disc}({\rm SO}_{24})$ de param\`etre standard $\psi_1 = [23] \oplus [1]$ est \'evidente : la repr\'esentation triviale convient (Exemples VI.\ref{exconjal}). Enfin, l'existence d'une repr\'esentation dans $\Pi_{\rm disc}({\rm SO}_{24})$ de param\`etre standard $\psi_{24} = \Delta_{11}[12]$ a d\'ej\`a \'et\'e d\'emontr\'ee dans le corollaire VII.\ref{cor24ikedabocherer} (travaux d'Ikeda et B\"ocherer). \ps\ps

\noindent {\it V\'erification du point 2.} On constate par inspection des
$\psi'_i$ qu'il n'y a rien \`a v\'erifier, car aucun entier $j$ ne v\'erifie
$d_j \equiv 0 \bmod 2$ et $d_j >23-2g_i$ (Crit\`ere \eqref{critbfinpreuve}),
\`a moins que $i \in \{13, 20,
21, 22, 23\}$ auquel cas le crit\`ere s'\'ecrit simplement ${\rm
L}(\frac{1}{2},\Delta_{15}) \neq 0$ pour $i=13$, ${\rm
L}(\frac{1}{2},\Delta_{11})
\neq 0$ pour $i \geq 20$. On conclut car ces deux valeurs de fonctions ${\rm
L}$ sont bien non nulles d'apr\`es la remarque VII.\ref{nonvanish}. \ps\ps

\noindent {\it V\'erification conjecturale du point 1}. La valeur de ${\rm m}(\varpi_i)$ est bien entendu d\'etermin\'ee par le th\'eor\`eme VIII.\ref{amfexplsp}. Pour appliquer
ce th\'eor\`eme, il faut savoir d'une part que le morphisme $\nu_\infty$
satisfait la conjecture \ref{conjaj}, et d'autre part v\'erifier les
condition (a) et (b) de son \'enonc\'e. La condition (a) est clairement
toujours satisfaite, par inspection des $\psi'_i$. En ce qui concerne la condition (b), nous affirmons
qu'elle est \'egalement toujours satisfaite. C'est un miracle du m\^eme
acabit que celui rencontr\'e au \S \ref{preuve1thm24}, et que nous pourrions 
v\'erifier de la m\^eme mani\`ere (ou au cas par cas!). Ce n'est en fait pas
n\'ecessaire, car d'apr\`es la proposition VIII.\ref{compmomsp} (plus exactement
d'apr\`es la d\'emonstration de cette proposition), cette v\'erification se d\'eduit formellement de celle
effectu\'ee au \S \ref{preuve1thm24} si l'on v\'erifie que pour tout
constituant de $\psi'_i$ de la forme $\pi[d]$ avec $d \equiv 0 \bmod 2$ et
$d > 23-2g_i$, on a $\varepsilon(\pi)=1$. Or une inspection imm\'ediate de
la table \ref{table24bis} montre qu'un tel constituant existe uniquement pour
les param\`etres d'indice $i \in \{13,20,21,22,23\}$ et que dans tous les cas $\pi =
\Delta_{11}$ ou $\Delta_{15}$, de sorte que l'on a bien $\varepsilon(\pi)=1$.
Le fait que l'on retrouve exactement les exceptions que nous avons d\^u
consid\'erer dans la v\'erification du point 2 n'est bien s\^ur pas un hasard,
comme nous l'avons expliqu\'e au \S \ref{comptheta}. Ainsi, si pour un $i$
donn\'e nous savons que la conjecture \ref{conjaj} vaut pour le
$\nu_\infty$ associ\'e \`a $\psi'_i$, on en d\'eduit ${\rm m}(\varpi_i)=1$. $\square$ \ps

\begin{remarque}\label{preuvemodmoeglin} {\rm Hormis le cas $i=2$ o\`u $\psi'_i={\rm Sym}^2 \Delta_{11}$ (et o\`u il est
\'evident que ${\rm m}(\varpi_i)=1$ !), le crit\`ere ``$d_j=1$ pour tout $j$''
ne s'applique malheureusement jamais. En revanche, on observe que si l'on disposait du cas particulier de la conjecture VIII.\ref{conjaj} 
annonc\'e par Arancibia d\'ecrit dans la remarque VIII.\ref{etudiantmoeglin}, \`a savoir du cas o\`u ``$\pi_j \neq 1 \Rightarrow d_j \leq 4$'', on pourrait conclure que ${\rm m}(\varpi_i) =1$ d\`es que  $i \notin \{11,12,13,20,22,23\}$. Comme d'autre part, le travail d\'ej\`a mentionn\'e d'Ikeda \cite{ikeda2} assure
que ${\rm m}(\varpi_i) \neq 0$ d\`es que $i \notin \{10,15,19,21\}$, cela conduirait \`a une d\'emonstration inconditionnelle du th\'eor\`eme.}
\end{remarque}

\section{Repr\'esentations alg\'ebriques de poids motivique $\leq
22$}\label{mrw}

\subsection{Un \'enonc\'e de classification} Le but de cette partie est de d\'emontrer le th\'eor\`eme suivant.

\begin{thm} \label{classpoids22} Soient $n\geq 1$ et $\pi \in \Pi_{\rm cusp}(\PGL_n)$
alg\'ebrique de poids motivique $w\leq 22$.  Alors on est dans
l'un des cas suivants : \ps\ps

\begin{itemize} \item[(i)] $n=1$, $w=0$ et $\pi$ est la repr\'esentation triviale. \ps\ps

\item[(ii)] $n=2$, $w \in \{ 11, 15, 17, 19, 21\}$, et $\pi$ est la repr\'esentation $\Delta_w$, engendr\'ee par l'unique forme modulaire parabolique normalis\'ee de poids $w+1$ pour le groupe ${\rm SL}_2(\Z)$.  \ps \ps

\item[(iii)] $n=3$, $w=22$ et $\pi$ est le carr\'e sym\'etrique de $\Delta_{11}$. \ps \ps

\item[(iv)] $n=4$ et ${\rm Poids}(\pi)$ est de la forme $\{\pm \frac{w}{2}, \pm \frac{v}{2}\}$ avec $$(w,v)\hspace{8pt}=\hspace{8pt}(19,7), \hspace{8pt} (21,5), \hspace{8pt}(21,9) \hspace{15pt}{
\rm ou}\hspace{15pt}(21,13).$$ 
\noindent Dans ce cas, $\pi$ est l'unique repr\'esentation dans $\Pi_{{\rm alg}}(\PGL_{4})$ de poids $\{\pm \frac{w}{2}, \pm \frac{v}{2}\}$ ;  en particulier,  on a $\pi \simeq \pi^{\vee}$. \ps \ps

\item[(v)] $n=4$, $w=22$ et ${\rm Poids}(\pi)=\{\pm 11,\pm v\}$ avec $v=4,5$ ou $6$. Dans ce cas, les repr\'esentations $\pi$ et  $\pi^{\vee}$ ne sont pas isomorphes, et ce sont les seules repr\'esentations dans $\Pi_{{\rm alg}}(\PGL_{4})$ de poids $\{\pm 11, \pm v\}$.\ps\ps
\end{itemize}
\medskip
\noindent De plus, s'il existe des repr\'esentations dans $\Pi_{{\rm alg}}(\PGL_4)$ de poids respectifs $\{\pm \frac{21}{2},\pm\frac{9}{2}\}$ et $\{\pm \frac{21}{2},\pm\frac{13}{2}\}$, alors le cas (v) ne se produit pas. 

\end{thm}

\ps\ps
\noindent Soulignons que l'assertion (iv) du th\'eor\`eme affirme seulement l'unicit\'e, mais non l'existence, d'une repr\'esentation $\pi \in \Pi_{{\rm alg}}(\PGL_{4})$ dont les poids sont de la forme $\{\pm \frac{w}{2}, \pm \frac{v}{2}\}$ avec $(w,v)$ parcourant la liste des $4$ couples de l'\'enonc\'e. Cependant nous avons montr\'e  dans la D\'efinition-Proposition \ref{definitiondeltawv}, \`a partir de la formule de Tsushima et des r\'esultats d'Arthur, qu'il existe en effet une telle repr\'esentation, que nous avons not\'ee $\Delta_{{w,v}}$. La raison pour laquelle nous avons formul\'e ainsi le th\'eor\`eme~\ref{classpoids22} est que sa d\'emonstration, comme nous le verrons, ne n\'ecessite pas de conna\^itre l'existence de ces repr\'esentations, et en particulier ne repose pas sur les travaux d'Arthur. Nous en d\'eduisons donc le th\'eor\`eme suivant\footnote{Mentionnons qu'\`a notre connaissance, personne n'a encore d\'emontr\'e l'existence d'une repr\'esentation $\pi \in \Pi_{{\rm alg}}(\PGL_{n})$, $n\geq 1$, qui ne soit pas autoduale.} (th\'eor\`eme \ref{introw22} de l'introduction) :  
\begin{thmv} \label{classpoids22bis} Soient $n\geq 1$ et $\pi \in \Pi_{\rm cusp}(\PGL_n)$
alg\'ebrique de poids motivique $\leq 22$. Alors $\pi$ appartient \`a la liste des $11$ repr\'esentations suivantes :
$$1, \hspace{4pt}\Delta_{11}, \hspace{4pt}\Delta_{15},\hspace{4pt} \Delta_{17},\hspace{4pt}\Delta_{19},\hspace{4pt}\Delta_{19,7},\hspace{4pt}\Delta_{21},\hspace{4pt}\Delta_{21,5},\hspace{4pt}\Delta_{21,9},\hspace{7pt}\Delta_{21,13}\hspace{7pt}{\rm et}\hspace{7pt}{\rm Sym}^2 \Delta_{11}.$$
\end{thmv}

\ps \medskip

Comme nous l'avons d\'ej\`a expliqu\'e dans l'introduction, notre d\'emonstration du th\'eor\`eme \ref{classpoids22} repose sur un analogue dans
le cadre des fonctions ${\rm L}$ automorphes des {\it formules explicites} de Riemann et Weil en th\'eorie
des nombres premiers. Nous renvoyons aux expos\'es de Poitou \cite{poitou,poitou2} \`a ce sujet.  Cet analogue a
\'et\'e d\'evelopp\'e par Mestre \cite{mestre}, puis appliqu\'e par Fermigier aux fonctions
${\rm L}(s,\pi)$ pour $\pi \in \Pi_{\rm cusp}(\PGL_n)$ \cite{fermigier}. Ce dernier en a
d\'eduit {\it loc. cit.} l'annulation de la cohomologie ``cuspidale'' du groupe ${\rm SL}_n(\Z)$ \`a coefficients rationnels
quand $n<24$. Ce r\'esultat a \'et\'e par la suite \'etendu \`a $n<27$ par Miller \cite{miller}, inspir\'e par
des travaux de Rudnick et Sarnak \cite{rudsar}, en consid\'erant plut\^ot la fonction ${\rm L}(s,\pi \times \pi^\vee)$, qui 
pr\'esente plusieurs avantages de convergence et positivit\'e\footnote{Dans un contexte proche, l'avantage \`a consid\'erer ce type de fonctions ${\rm L}$ avait d\'ej\`a \'et\'e observ\'e par Serre \cite[p. 15]{poitou2}.}. Ce sera \'egalement notre point de d\'epart. Plus g\'en\'eralement, nous allons examiner en d\'etail les in\'egalit\'es donn\'ees par la formule explicite appliqu\'ee \`a la fonction ${\rm L}(s,\pi\times \pi')$ pour toute paire de repr\'esentations $\{\pi,\pi'\}$ avec $\pi \in \Pi_{{\rm alg}}(\PGL_{n})$ et $\pi' \in \Pi_{{\rm alg}}(\PGL_{n'})$. \ps

\subsection{La formule explicite pour les fonctions ${\rm L}$ de paires} \label{preliminairesfexpl}

Les formules explicites d\'ependent du choix d'une ``fonction test''. Suivant l'analyse de Poitou et Weil \cite[p. 6]{poitou}, nous entendrons par l\`a toute fonction paire $F :\R \rightarrow \R$ satisfaisant les axiomes suivants, dans lesquels $F_\epsilon$, pour $\epsilon$ r\'eel $>0$, d\'esigne la fonction $\R \rightarrow \R, \,\,x \mapsto F(x)e^{(\frac{1}{2}+\epsilon)x}$  : \ps \ps
 \begin{itemize} 
\item[(T1)] il existe $\epsilon>0$ tel que $F_{\epsilon}$ soit sommable sur $\R_{>0}$; \ps \ps
\item[(T2)] il existe $\epsilon>0$ tel que $F_{\epsilon}$ est \`a variation born\'ee sur $\R$, \'egale en chaque point \`a la moyenne de ses limites \`a droite et \`a gauche; \ps \ps
\item[(T3)] la fonction $\frac{1}{x}(F(x) -F(0))$ est \`a variation born\'ee sur $\R-\{0\}$. \ps\ps
\end{itemize}\ps
\noindent 
En particulier, toute fonction paire de classe $\mathcal{C}^2$ sur $\R$ et \`a support compact est une fonction test. Ce sera le cas dans nos applications, o\`u nous prendrons pour $F$ une modification simple de la fonction d'Odlyzko. Dans un premier temps, il sera cependant plus clair de consid\'erer des fonctions tests arbitraires. \ps\ps

Si $F$ est une fonction test, et si $\epsilon>0$ est tel que $F_\epsilon$ est sommable sur $\R$, l'int\'egrale d\'ependant du param\`etre complexe $s$ d\'efinie par
\begin{equation}\label{formulegrandphi} \Phi_F(s) = \int_\R F(x) e^{(s-\frac{1}{2})x } {\rm d}x\end{equation}
est absolument convergente dans la r\'egion $- \epsilon < {\rm Re}\, \, s < 1+ \epsilon$. La fonction $\Phi_F(s)$ est en particulier bien d\'efinie et holomorphe dans cette r\'egion. La parit\'e de $F$ entra\^ine l'\'egalit\'e $\Phi_F(s)=\Phi_F(1-s)$. De plus, la relation $\Phi_F(\overline{s})=\overline{\Phi_F(s)}$ montre $\Phi_F(s) \in \R$ pour tout r\'eel $s$ tel que $0\leq s \leq 1$. \ps\ps

\ps\ps

On pose $$\Pi_{\rm alg}=\coprod_{n\geq 1} \Pi_{\rm alg}(\PGL_n).$$ 
Soient $\pi, \pi' \in \Pi_{\rm alg}$. Ainsi que nous l'avons rappel\'e au \S {\rm VIII}.\ref{parfaceps}, on dispose suivant Jacquet et Shalika  \cite[Theorem 5.3]{jasha1} d'un produit Eul\'erien $${\rm L}(s, \pi \times \pi')=\prod_p \det(1-p^{-s} {\rm c}_p(\pi) \otimes {\rm c}_p(\pi'))^{-1}$$ 
qui est bien d\'efini et absolument convergent pour ${\rm Re}\, s >1$. De plus, la fonction 
$$\xi(s, \pi \times \pi') \,:= \, \Gamma(s,{\rm L}(\pi_\infty) \otimes {\rm L}(\pi'_\infty)) \,\, {\rm L}(s, \pi \times \pi')$$
admet un prolongement m\'eromorphe \`a $\C$ tout entier v\'erifiant l'\'equation fonctionnelle $\xi(s,\pi \times \pi') \,=\, \varepsilon(\pi \times \pi')\, \xi(1-s,\pi^\vee \times {\pi'}^\vee)$ avec $$\varepsilon(\pi \times \pi')=\varepsilon({\rm L}(\pi_\infty) \otimes {\rm L}(\pi'_\infty)).$$ En particulier, tous les z\'eros de $\xi$ sont dans la bande critique $0 \leq {\rm Re}\, s\, \leq 1$ (Shahidi a m\^eme d\'emontr\'e que ces z\'eros sont dans l'int\'erieur de cette bande). Enfin, la fonction $\xi$ est holomorphe sur $\C - \{0,1\}$, et admet un p\^ole en $s=1$ si, et seulement si, $\pi' = \pi^\vee$, auquel cas il s'agit d'un p\^ole simple. \ps

\begin{defprop}\label{prelimformexpl} Soient $\pi,\pi' \in \Pi_{\rm alg}$, $\xi(s)=\xi(s,\pi^\vee \times \pi')$, $F$ une fonction test et $T$ un r\'eel $>0$. La somme finie $$\sum_{\{s \in \C,\,\, |{\rm Im}\, s| < T, \, \, \xi(s)=0\}} \Phi_F(s)\,\,{\rm ord}_{z=s} \,\xi(z)$$ 
est r\'eelle et admet une limite finie quand $T \rightarrow +\infty$ ; nous noterons ${\rm Z}^F(\pi,\pi')$ cette limite. 
\end{defprop}\ps \ps

En effet, cet \'enonc\'e est un cas particulier des r\'esultats de Mestre \cite[\S I]{mestre}, qui g\'en\'eralisent eux-m\^eme de mani\`ere assez directe ceux de Riemann, Weil, et Poitou \cite{poitou}. Supposons $\pi \in \Pi_{\rm alg}(\PGL_n)$ et $\pi' \in \Pi_{\rm alg}(\PGL_{n'})$. Dans les notations de Mestre, on prend $M=M'=nn'$, $c=0$, $L_1(s)={\rm L}(s,\pi^\vee \times \pi')$, $L_2(s)={\rm L}(s,\pi \times (\pi')^\vee)$, $w=\varepsilon(\pi^\vee \times \pi')$,  $\Lambda_1(s)=\xi(s,\pi^\vee \times \pi')$ et $\Lambda_2(s)=\xi(s,\pi \times (\pi')^\vee)$. Par d\'efinition, nous incorporons les facteurs qu'il note $A^s$ et $B^s$, ainsi que ses coefficients $a_i, a_i', b_i, b_i'$, dans nos facteurs archim\'ediens $\Gamma$, et il n'y a pas de contribution du conducteur (qui vaut $1$ ici). Ceci \'etant dit, les hypoth\`eses (i), (ii) et (iii) {\it loc. cit.} d\'ecoulent de l'\'equation fonctionnelle et de la finitude du nombre de p\^oles des $\Lambda_i$, qui ont d\'ej\`a \'et\'e justifi\'es ci-dessus. L'hypoth\`ese (iii), \`a savoir que la fonction $\xi$ diminu\'ee de ses parties singuli\`eres est born\'ee dans toute bande verticale, est un th\'eor\`eme de Gelbart et Shahidi \cite[Cor. 2]{GSha}. Enfin, seule une version affaiblie de sa derni\`ere hypoth\`ese (iv) est satisfaite, \`a savoir  la convergence absolue des produits Eul\'eriens $L_1$ et $L_2$ ainsi que leur non-annulation sur ${\rm Re}\, s>1$, mais c'est tout ce dont il a besoin par la suite : voir \cite[p. 213-214]{mestre} et surtout l'argument donn\'e par Poitou \cite[p. 2-3]{poitou}. \ps\ps

La conclusion de cette discussion est que tous les r\'esultats de \cite[\S I.2]{mestre} s'appliquent. L'assertion de convergence dans la d\'efinition-proposition ci-dessus est notamment justifi\'ee par Mestre page 213. Le fait que la somme finie apparaissant dans la d\'efinition-proposition \ref{prelimformexpl} est r\'eelle vient du fait que la r\'egion $\{ s \in \C, \, |{\rm Im} \,s|<T\}$ est stable par $s \mapsto 1-\overline{s}$ et des \'egalit\'es
$$\xi(1-\overline{s},\pi^\vee \times \pi') \, = \,\varepsilon(\pi^\vee \times \pi') \,\xi(\overline{s},\pi \times \pi'^\vee) \,=\, \varepsilon(\pi^\vee \times \pi')\, \overline{\xi(s,\pi^\vee \times \pi')}.$$

\noindent Mestre \'etablit \'egalement {\it loc. cit.} la formule explicite que nous allons utiliser. C'est le r\'esultat de l'int\'egration de la $1$-forme $\Phi_F(s) \,{\rm d} {\rm log} \,\xi(s)$ sur le bord du rectangle $\{ s \in \C, \, \,  -\epsilon \leq {\rm Re} \, s \, \leq 1+\epsilon, \, \, |{\rm Im}\, s |\leq A\}$, $A$ et $\epsilon$ \'etant des r\'eels $>0$ convenables, puis d'un passage \`a la limite $\epsilon \rightarrow 0$ et $A \rightarrow \infty$. Afin de pouvoir l'\'enoncer sous une forme agr\'eable, introduisons d'abord certaines quantit\'es pr\'eliminaires ``locales''. L'assertion de convergence dans la d\'efinition suivante est justifi\'ee dans \cite[p. 213-214]{mestre} et \cite[p. 2-3]{poitou}. \ps \ps

\begin{defprop}\label{prelimformexplbf} Soient $\pi, \pi' \in \Pi_{\rm alg}$ et $F$ une fonction test. La somme $$\sum_{p,k} \,\, F(k {\rm log}(p))\frac{{\rm log}(p)}{p^{k/2}}\,\,\,\overline{{\rm tr}\, ({\rm c}_p(\pi)^k)} \,\, {\rm tr} \,({\rm c}_p(\pi')^k)\, ,$$ portant sur tous les couples $(p,k)$ avec $p$ un nombre premier et $k$ un entier $\geq 1$, est absolument convergente ;  nous la notons $\widetilde{\cbil}_f^F(\pi,\pi')$. On a les relations \'evidentes $\widetilde{\cbil}_f^F(\pi,\pi')=\overline{\widetilde{\cbil}_f^F( \pi',\pi)}=\widetilde{\cbil}_f^F((\pi')^\vee,\pi^\vee)$. On pose enfin $$\cbil_f^F(\pi,\pi')\, :=\, \Re\,\widetilde{\cbil}_f^F(\pi,\pi').$$
\end{defprop}

\ps \ps \noindent

Soit ${\rm W}_\R^{\rm alg}$ le quotient du groupe de Weil ${\rm W}_\R$ par la composante connexe de son centre, \`a savoir $\R_{>0}$.
Soit ${\rm K}_\infty$ l'anneau de Grothendieck de la cat\'egorie des
repr\'esentations complexes, continues, et de dimension finie, du groupe ${\rm W}_\R^{\rm alg}$. D'apr\`es les rappels du \S VIII.\ref{apparitionWR}, c'est le
groupe ab\'elien libre sur les (classes des) repr\'esentations $$1, \, \, \, \, \, 
\epsilon_{\C/\R}\, \, \, \, \,  {\rm et}\, \, \, \, \,  {\rm I}_{w} \, \, \, {\rm pour} \, \, \, w \in \Z_{>0}.$$ 
Remarquons que tout \'el\'ement de ${\rm K}_\infty$ est \'egal \`a son dual, car c'est le cas des classes des trois repr\'esentations ci-dessus.\ps\ps

L'application $V
\mapsto \Gamma(-,V)$ introduite au \S VIII.\ref{parfaceps} s'\'etend de mani\`ere naturelle en un homomorphisme
$\Gamma : {\rm K}_\infty \rightarrow \mathcal{M}(\C)^\times$, o\`u
$\mathcal{M}(\C)^\times$ d\'esigne le groupe multiplicatif du corps $\mathcal{M}(\C)$ des fonctions
m\'eromorphes sur $\C$. Son image est constitu\'ee de fonctions
n'ayant ni p\^ole, ni z\'ero, dans le demi-plan ${\rm Re}\, s >0$. \ps\ps

\begin{defprop}\label{lemmemestre} Si $F$ est une fonction test, l'application
$$ {\rm J}_F : {\rm K}_\infty \rightarrow \R, \, \, \, \, \, \,  V \mapsto -\frac{1}{2\pi
i}\int_{{\rm
Re}(s)=\frac{1}{2}} \Phi_F (s) \frac{\Gamma'(s,V)}{\Gamma(s,V)}{\rm d} s,$$
est bien d\'efinie et $\Z$-lin\'eaire. De plus, on a les identit\'es : \ps \ps\begin{itemize}\item[(i)] ${\rm J}_F(1)\,=\,\frac{1}{2}{\rm log}(\pi)\,F(0)+\sigma_F(\frac{1}{2},0)$,\ps\ps
\item[(ii)]  ${\rm J}_F(\epsilon_{\C/\R})\,=\,\frac{1}{2}{\rm log}(\pi)\, F(0)\,+\sigma_F(\frac{1}{2},\frac{1}{2})$, \ps \ps
\item[(iii)] ${\rm J}_F({\rm I}_w)\,=\,{\rm log}(2\pi)\,F(0)\,+\sigma_F(1,\frac{w}{2})$ pour $w\geq 0$.\ps \ps
\end{itemize}
o\`u l'on a pos\'e $\sigma_F(a,b)=a\int_0^{+\infty} \left( F(ax) \frac{e^{-(a/2+b)x}}{1-e^{-x}} - F(0)\frac{e^{-x}}{x} \right){\rm d}x$. \ps\ps
Enfin, l'application ${\rm K}_\infty \times {\rm K}_\infty \rightarrow \R, \, (V,W) \mapsto {\rm J}_F(V^*\otimes W)$, est bilin\'eaire  ; nous la noterons $\cbil_\infty^F$. La relation \'evidente $\cbil_\infty^F(V,W)={\rm  J}_F(V\otimes W)$ pour tout $V,W \in {\rm K}_\infty$ montre que $\cbil_\infty^F$ est sym\'etrique.

\end{defprop}

\begin{pf} Notons $\psi \in \mathcal{M}(\C)$ la fonction digamma, d\'efinie par $\psi(s) = \Gamma'(s)/\Gamma(s)$. Posons \'egalement 
$\psi(s,V)=\frac{\Gamma'(s,V)}{\Gamma(s,V)}$ pour $V \in {\rm K}_\infty$. L'application $V \mapsto \psi(-,V)$, ${\rm K}_\infty \rightarrow \mathcal{M}(\C)$, est manifestement $\Z$-lin\'eaire. Suivant la d\'efinition des facteurs $\Gamma(s,V)$ rappel\'ee au \S VIII.\ref{parfaceps}, on a les identit\'es $$\psi(s,1)= - \frac{1}{2} \log\,\pi + \frac{1}{2} \psi(\,\frac{s}{2}\,)\hspace{10pt}{\rm et}\hspace{10pt}\psi(s,{\rm I}_w) = - \log(2\pi) + \psi(s+\frac{w}{2})$$ pour tout $w\geq 0$. La fonction $F$ est sommable sur $\R$ d'apr\`es {\rm (T1)}, \'egale en tout point \`a la moyenne de ses limites \`a droite et \`a gauche d'apr\`es (T2), et $\frac{F(x)-F(0)}{x}$ est born\'ee au voisinage de $0$ d'apr\`es (T3). La formule d'inversion de Fourier est donc valable en $x=0$, autrement dit l'int\'egrale $-\frac{1}{2\pi
i}\int_{{\rm Re}(s)=\frac{1}{2}} \Phi_F (s) {\rm d} s $ est convergente et vaut $-F(0)$. Il reste \`a examiner la convergence d'une int\'egrale de la forme 
$$-\frac{a}{2i\pi} \int_{{\rm
Re}(s)=\frac{1}{2}} \Phi_F (s)\psi(a s+ b){\rm d}s$$
avec $b \in \R_{\geq 0}$ et $a  \in \R_{>0}$. La convergence (simple) de cette int\'egrale est v\'erifi\'ee dans \cite[p. 6-04]{poitou} et \cite[Lemme I.2.1]{mestre}
ainsi que l'\'egalit\'e de sa somme avec l'int\'egrale (\'egalement convergente)  $\sigma_F(a,b)$ de l'\'enonc\'e.  Cela entra\^ine tous les points de la proposition. \end{pf}

Nous disposons maintenant de tous les ingr\'edients n\'ecessaires pour \'enoncer la formule explicite.  Il sera commode d'introduire le groupe ab\'elien libre sur l'ensemble $\Pi_{\rm alg}$ : $${\rm K} = \Z [\Pi_{\rm alg}].$$ 

Soit $F$ une fonction test. Chacune des trois fonctions $\Pi_{\rm alg} \times \Pi_{\rm alg} \rightarrow \R$ envoyant respectivement $(\pi,\pi')$ sur ${\rm Z}^F(\pi,\pi')$, $\cbil_f^F(\pi, \pi')$ et $\delta_{\pi,\pi'}$ (symbole de Kronecker), s'\'etend en une application bilin\'eaire sym\'etrique ${\rm K} \times {\rm K} \rightarrow \R$ associ\'ee \`a $F$. Nous noterons ces trois formes ${\rm Z}^F$, $\cbil_f^F$ et $\delta$ respectivement. De plus, l'application $\Pi_{\rm alg} \rightarrow {\rm K}_\infty$, $\pi \mapsto {\rm L}(\pi_\infty)$, s'\'etend en un homomorphisme ${\rm L} : {\rm K} \rightarrow {\rm K}_\infty$. L'application ${\rm K} \times {\rm K} \rightarrow \R, \, \, (V,W) \mapsto \cbil_\infty^F({\rm L}(V),{\rm L}(W))$, est donc \'egalement bilin\'eaire et sym\'etrique ; nous nous autoriserons l'abus de langage consistant \`a la noter \'egalement $\cbil_\infty^F$. On pose enfin

$$\cbil^F=\cbil_f^F+\cbil_\infty^F\hspace{10pt}: \hspace{10pt} {\rm K} \times {\rm K} \rightarrow \R.$$

\ps\ps

\begin{prop} {\rm (``Formule explicite'')} \label{formexplpipi} Pour toute fonction test $F$, on a l'\'egalit\'e suivante entre formes bilin\'eaires ${\rm K} \times {\rm K} \rightarrow \R$ : 
$$\cbil^F \,+\, \frac{1}{2}\, {\rm Z}^F\, =\, \Phi_F(0)\, \, \,\delta.$$
\end{prop}

\begin{pf} Par bilin\'earit\'e, il suffit de montrer cette \'egalit\'e pour un couple $(\pi,\pi') \in \Pi_{\rm alg} \times \Pi_{\rm alg}$. Dans ce cas, la formule donn\'ee par Mestre \cite[\S I.2]{mestre} s'\'ecrit, compte tenu de la parit\'e de la fonction $F$ :
{\small $$\cbil_\infty^F(\pi,\pi')\,+\,\cbil_\infty^F(\pi',\pi)\,+\,\widetilde{\cbil}_f^F(\pi,\pi')\,+\,\widetilde{\cbil}_f^F(\pi',\pi)\,+\,{\rm Z}^F(\pi,\pi')\,=\, \delta_{\pi,\pi'}\,(\Phi_F(0)+\Phi_F(1)).$$}\par
\noindent On conclut par les identit\'es $\Phi_F(1)=\Phi_F(0)$, $\cbil_\infty^F(\pi,\pi')=\cbil_\infty^F(\pi',\pi)$,  et $2 \,\cbil_f^F(\pi,\pi')\,=\,\widetilde{\cbil}_f^F(\pi,\pi')\,+\,\widetilde{\cbil}_f^F(\pi',\pi)$.
\end{pf}

\ps \ps 

On dira qu'un \'el\'ement $\pi$ de ${\rm K}$ est effectif si c'est une somme finie d'\'el\'ements de $\Pi_{\rm alg}$ (autrement dit, une combinaison lin\'eaire \`a coefficients $\geq 0$).  Si $\pi$ et $\pi'$ sont effectifs, disons $\pi=\sum_i \pi_i$ et $\pi'=\sum_j \pi'_j$ avec $\pi_i,\pi'_j \in \Pi_{\rm alg}$, on notera $\xi(s,\pi \times \pi')$ (resp. ${\rm L}(s,\pi\times \pi')$), le produit des $\xi(s,\pi_i \times \pi'_j)$ (resp. ${\rm L}(s,\pi_i\times \pi'_j)$). Consid\'erons de plus l'assertion suivante portant sur une fonction test $F$ : \ps\ps\begin{itemize}
\item[(T4)] $\Re \, \, \Phi_F(s)$ est $\geq 0$ pour tout $s \in \C$ tel que $ 0 \leq {\rm Re}\, s \leq 1$. \ps\ps
\end{itemize}
\noindent On rappelle que $\Phi_F(s)$ est un nombre r\'eel pour tout r\'eel $s$ dans $[-1,1]$, en particulier il est $\geq 0$ si $F$ satisfait (T4). Dans ce m\'emoire, nous utiliserons uniquement le corollaire suivant de la formule explicite.

\begin{cor}\label{inegexpl} Soient $F$ une fonction test et $\pi,\pi' \in {\rm K}$ effectifs. On suppose que $F$ satisfait la condition (T4). On a l'in\'egalit\'e $$\cbil^F(\pi,\pi')\, \leq \,\Phi_F(0) \,\delta(\pi,\pi')\, - \,\frac{1}{2}\,\Phi_F(\,\frac{1}{2}\,)\,{\rm ord}_{s=\frac{1}{2}} \xi(s,\pi^\vee \times \pi').$$
\end{cor}

\begin{pf} Sous l'hypoth\`ese (T4) sur $F$, on a l'in\'egalit\'e $${\rm Z}^F(\pi,\pi')\, = \, \Re\,{\rm Z}^F(\pi,\pi')\, \geq \Phi_F(\,\frac{1}{2}\,)\,{\rm ord}_{s=\frac{1}{2}} \xi(s,\pi^\vee \times \pi').$$ On conclut par la proposition \ref{formexplpipi}. 
\end{pf}
\ps\ps
\begin{remarque}{\rm Si $V \in {\rm K}_\infty$, on constate que la fonction m\'eromorphe $\Gamma(s,V)$ est finie et non nulle en $s=\frac{1}{2}$. On en d\'eduit, pour tous $\pi,\pi' \in {\rm K}$ effectifs, l'\'egalit\'e ${\rm ord}_{s=\frac{1}{2}} \xi(s,\pi \times \pi')={\rm ord}_{s=\frac{1}{2}} {\rm L}(s,\pi \times \pi')$. Le corollaire ci-dessus vaut donc \'egalement en rempla\c{c}ant  $\xi(s,\pi^\vee \times \pi')$ par ${\rm L}(s,\pi^\vee \times \pi')$. }\end{remarque}

L'ordre d'annulation des fonctions $\xi(s,\pi \times \pi')$ en $s=\frac{1}{2}$ est r\'eput\'e myst\'erieux. On le minore traditionnellement de la mani\`ere suivante. Si $\pi,\pi' \in \Pi_{\rm alg}$ on pose ${\rm e}^\bot(\pi,\pi')=1$ si $\pi$ et $\pi'$ sont toutes deux autoduales et satisfont $\varepsilon(\pi \times \pi')=-1$ ; sinon on pose ${\rm e}^\bot(\pi,\pi')=0$. La fonction ${\rm e}^\bot : \Pi_{\rm alg} \times \Pi_{\rm alg} \rightarrow \Z$ s'\'etend en une forme bilin\'eaire sym\'etrique ${\rm e}^\bot : {\rm K} \times {\rm K} \rightarrow \Z$. Ainsi, pour tous $\pi,\pi' \in {\rm K}$ effectifs on a l'in\'egalit\'e ${\rm ord}_{s=\frac{1}{2}} \xi(s,\pi \times \pi') \geq {\rm e}^\bot(\pi,\pi')$.

\ps\ps
\begin{cor}\label{inegexploo} Soient $F$ une fonction test et $\pi,\pi' \in {\rm K}$ effectifs. On suppose que $F$ satisfait (T4) et $F\geq 0$.
On pose $${\rm C}^F(\pi,\pi')=\Phi_F(0) \,\delta(\pi,\pi')-\, \frac{1}{2}\, \Phi_F(\,\frac{1}{2}\,) \,{\rm e}^\bot(\pi,\pi')\,-\cbil_\infty^F(\pi,\pi').$$ On a les in\'egalit\'es : \ps \ps
\begin{itemize}
\item[(i)] ${\rm C}^F(\pi,\pi) \geq 0$. En particulier, on a $\cbil_\infty^F(\pi,\pi) \leq \Phi_F(0) \,\delta(\pi,\pi)$. \ps \ps
\item[(ii)] ${\rm C}^F(\pi,\pi') + \sqrt{{\rm C}^F(\pi,\pi) {\rm C}^F(\pi',\pi')} \geq 0$. \ps \ps
\end{itemize}
\end{cor}

\begin{pf} Le corollaire \ref{inegexpl} et la discussion pr\'ec\'edente montrent que sous l'hypoth\`ese (T4) sur $F$ on a  $\cbil^F_f(\pi,\pi') \leq {\rm C}^F(\pi,\pi')$ pour tous $\pi,\pi' \in {\rm K}$ effectifs. Observons que l'hypoth\`ese de positivit\'e de la fonction $F$ entra\^ine la positivit\'e de la forme bilin\'eaire $\cbil_f^F$ sur ${\rm K}$. En effet, 
un \'el\'ement $\varpi \in {\rm K} \otimes \R$ peut s'\'ecrire sous la forme $\varpi = \sum_i \lambda_i \pi_i$ o\`u les $\pi_i$ sont dans $\Pi_{\rm alg}$ et les $\lambda_i$ sont dans $\R$. L'identit\'e
$$\cbil^F_f(\varpi,\varpi)\,=\sum_{p,k} \,F(k\, {\rm log} \,p)\,\frac{\,{\rm log}\,p}{p^{\frac{k}{2}}} \,\,\,|\,\,\sum_i \,\,\lambda_i\,\, {\rm tr}\, {\rm c}(\pi_i)^k \,\,|^2$$
entra\^ine donc $\cbil^F_f(\varpi,\varpi) \geq 0$. L'assertion (ii) se d\'eduit de l'in\'egalit\'e de Cauchy-Schwarz appliqu\'ee \`a $\cbil_f^F$ :
$$|\cbil_f^F(\pi,\pi')|\leq \sqrt{\cbil_f^F(\pi,\pi)\,\cbil_f^F(\pi',\pi')}$$
et de l'in\'egalit\'e \'evidente $|\cbil_f^F(\pi,\pi')| \,\geq \,- \cbil_f^F(\pi,\pi')$. (On prendra garde que la forme bilin\'eaire ${\rm C}^F$ est a priori seulement positive sur les \'el\'ements effectifs de ${\rm K}$, d'o\`u la formulation de l'assertion (ii).)
\end{pf}

Dans \cite{miller}, l'in\'egalit\'e (i) du corollaire ci-dessus est exploit\'ee dans le cas particulier $\pi \in \Pi_{\rm alg}$. Pour un $\pi$ effectif quelconque elle admet le corollaire \ref{cortaibi} ci-dessous, d'abord observ\'e par Olivier Ta\"ibi. L'in\'egalit\'e (ii) semble nouvelle, elle rendra de grands services dans les applications.  \ps \noindent 

\begin{definition}\label{deffonctionm}  Soit $V \in {\rm K}_\infty$. On note ${\rm m}(V)$ le nombre de repr\'esentations $\pi \in \Pi_{\rm alg}$ satisfaisant ${\rm L}(\pi_\infty) \simeq V$ ;  d'apr\`es Harish-Chandra on a ${\rm m}(V) < +\infty$ {\rm (\S IV.\ref{fautdiscgen}{\rm )}}. On note \'egalement ${\rm m}^\bot(V)$  le nombre de repr\'esentations $\pi \in \Pi_{\rm alg}$ autoduales et satisfaisant ${\rm L}(\pi_\infty) \simeq V$. On a l'in\'egalit\'e ${\rm m}^\bot(V)  \leq {\rm m}(V)$. 
\end{definition}

Un \'el\'ement $V$ de ${\rm K}_\infty$ sera dit {\rm effectif} si c'est la classe d'une repr\'esentation de dimension finie, continue, et \`a coefficients dans $\C$, de ${\rm W}_\R^{\rm alg}$. Il est \'evident que si ${\rm m}(V) \neq 0$ alors $V$ est effectif.

\begin{cor} \label{cortaibi} {\rm (Ta\"ibi)} Soient $V \in {\rm K}_\infty$ effectif et $F$ une fonction test, suppos\'ee $\geq 0$ et satisfaisant l'hypoth\`ese {\rm (T4)}. On a l'in\'egalit\'e 
$${\rm m}(V)\,\, \cbil_\infty^F(V,V) \leq \Phi_F(0).$$
\end{cor}

\begin{pf} Si ${\rm m}(V)=0$ il n'y a rien \`a d\'emontrer. Supposons donc qu'il existe un entier $r\geq 1$ et des repr\'esentations  $\pi_1,\dots,\pi_r \in \Pi_{\rm alg}$ distinctes telles que ${\rm L}((\pi_i)_\infty) \simeq V$ pour tout $i$. On applique le point (i) du corollaire~\ref{inegexploo} \`a l'\'el\'ement $\pi = \pi_1 + \pi_2+ \dots + \pi_r$ de ${\rm K}$. On a d'une part les \'egalit\'es $$\cbil^F_\infty(\pi,\pi)\,=\,\cbil^F_\infty(rV,rV)\,=\,r^2 \,\cbil^F_\infty(V,V)$$ et d'autre part $\delta(\pi,\pi)=r$. On en d\'eduit l'in\'egalit\'e $r^2\,\,   \cbil^F_\infty(V,V)\, \leq \,r\, \Phi_F(0)$, et donc $r \,\cbil^F_\infty(V,V) \leq \Phi_F(0)$.
\end{pf}
\ps
\noindent 
Observons que sous l'hypoth\`ese suppl\'ementaire $\cbil^F_\infty(V,V)>0$, le corollaire ci-dessus fournit une majoration explicite de ${\rm m}(V)$. En particulier, il red\'emontre le r\'esultat d'Harish-Chandra susmentionn\'e ${\rm m}(V) < \infty$ (\S IV.\ref{fautdiscgen}). \ps \ps

\begin{cor}\label{inexploo2} Soient $V,V' \in {\rm K}_\infty$ effectifs et $F$ une fonction test $\geq 0$ satisfaisant (T4). On suppose $V \neq V'$ et ${\rm m}(V)\,{\rm m}(V') \,\neq \,0$. \ps \ps \begin{itemize} \item[(i)] Si l'on pose ${\rm n}(V,V')=\frac{{\rm m}^\bot(V)\, {\rm m}^\bot(V')}{4\,{\rm m}(V)\,{\rm m}(V')}(1-\varepsilon(V \otimes V'))$, on a l'in\'egalit\'e

{\footnotesize $$\,{\rm n}(V,V') \,\Phi_F(\frac{1}{2}) \,+ \,\cbil_\infty^F(V,V')\, \,\leq \, \sqrt{(\frac{\Phi_F(0)}{{\rm m}(V)}-\cbil_\infty^F(V,V))(\frac{\Phi_F(0)}{{\rm m}(V')}-\cbil_\infty^F(V',V'))},$$}\par
\item[(ii)] Si l'on a en outre ${\rm m}^\bot(V) {\rm m}^\bot(V') \neq 0$, on a \'egalement l'in\'egalit\'e \ps \ps\end{itemize}

{\footnotesize $$\,\frac{1-\varepsilon(V \otimes V')}{4} \,\Phi_F(\frac{1}{2}) \,+ \,\cbil_\infty^F(V,V')\, \,\leq \, \sqrt{(\frac{\Phi_F(0)}{{\rm m}^\bot(V)}-\cbil_\infty^F(V,V))(\frac{\Phi_F(0)}{{\rm m}^\bot(V')}-\cbil_\infty^F(V',V'))},$$}

\end{cor}
\ps
\begin{pf} Soient $r,r'$ des entiers $\geq 1$ et des repr\'esentations distinctes  $\pi_1,\dots,\pi_r,\pi'_1,\dots,\pi'_{r'}\in \Pi_{\rm alg}$ telles que ${\rm L}(\pi_i)=V$ pour tout $i$ et ${\rm L}(\pi'_j)=V'$ pour tout $j$. On pose $\pi=\sum_i \pi_i$ et $\pi'=\sum_j \pi'_j$. On a $\delta(\pi,\pi')=0$ car $V \neq V'$. Notons $s \leq r$ (resp. $s'\leq r'$) le nombre de repr\'esentations autoduales parmi les  $\pi_i$ (resp. $\pi'_j$).
On a  l'\'egalit\'e \'evidente $${\rm e}^\bot(\pi,\pi')\,=\,s\,s'\, \,\frac{1-\varepsilon(V \otimes V')}{2}.$$ On a de plus $${\rm C}^F(\pi,\pi')\,+\, \frac{{\rm e}^\bot(\pi,\pi')}{2} \,\Phi_F(\frac{1}{2})\, = -\,\cbil_\infty^F(rV,r'V')\,=\,- rr'\,\cbil_\infty^F(V,V').$$  De m\^eme, on constate les in\'egalit\'es ${\rm C}^F(\pi,\pi)\,\leq \,r\, \Phi_F(0)\,- \,r^2\, \cbil^F_\infty(V,V)$ et ${\rm C}^F(\pi',\pi')\,\leq \,r'\, \Phi_F(0)\,- \,(r')^2 \,\cbil_\infty^F(V',V')$. En divisant par $rr' \neq 0$ l'in\'egalit\'e donn\'ee par le point (ii) du corollaire~\ref{inegexploo}, on en d\'eduit l'in\'egalit\'e 
{\footnotesize $$\,\frac{ss'}{4rr'}(1-\varepsilon(V\otimes V')) \,\Phi_F(\frac{1}{2}) \,+ \,\cbil_\infty^F(V,V')\, \,\leq \, \sqrt{(\frac{\Phi_F(0)}{r}-\cbil_\infty^F(V,V))(\frac{\Phi_F(0)}{r'}-\cbil_\infty^F(V',V'))}.$$}\par
\noindent Cette in\'egalit\'e est valable pour tout quadruplet d'entiers $(r,r',s,s')$ tels que $1 \leq r \leq {\rm m}(V)$, $1\leq r' \leq {\rm m}(V')$, $s\leq r$ et $s'\leq r'$. L'assertion (i) r\'esulte du cas particulier $(r,r',s,s')=({\rm m}(V),{\rm m}(V'),{\rm m}^\bot(V),{\rm m}^\bot(V'))$, et l'assertion (ii) du cas $(r,r',s,s')=({\rm m}(V)^\bot,{\rm m}^\bot(V'),{\rm m}^\bot(V),{\rm m}^\bot(V'))$. 

\end{pf}

Comme nous le verrons, ce corollaire permet typiquement de montrer que si $V \neq V'$, l'existence de certains \'el\'ements $\pi \in \Pi_{\rm alg}$ tels que ${\rm L}(\pi_\infty)=V$ entra\^ine l'inexistence d'\'el\'ements $\pi \in \Pi_{\rm alg}$ tels que ${\rm L}(\pi'_\infty)=V'$. Il admet divers raffinements sur lesquels nous reviendrons au \S \ref{interludegeom}.  Terminons par un crit\`ere simple bien connu \cite{poitou} permettant de cons\-truire des fonctions tests satisfaisant {\rm (T4)} ; nous redonnons l'argument pour le confort du lecteur. \ps\ps

\begin{lemme}\label{lemeposphi} Soit $g : \R \rightarrow \R$ une fonction paire, sommable et de carr\'e sommable. Sa transform\'ee de Fourier $\widehat{g}$ est bien d\'efinie et \`a valeurs r\'eelles sur $\R$. Supposons $\widehat{g} \geq 0$ et consid\'erons la fonction $F : \R \rightarrow \R$ d\'efinie par $$F(x)=\frac{g(x)}{{\rm cosh}(x/2)}.$$ Alors $\Phi_F(s)$ est bien d\'efinie pour tout $s \in \C$ tel que $0\leq  {\rm Re}\, s\,  \leq 1$ {\rm (}formule~\eqref{formulegrandphi}{\rm )}, et dans cette r\'egion on a \,\,\,\,$\Re\,\, \Phi_F(s)\,  \geq 0$. 
\end{lemme}
\noindent Notre convention pour la transform\'ee de Fourier d'une fonction sommable $g$ est $\widehat{g}(y) = \int_{\R} g(x) \, e^{-2i\pi xy}\, {\rm d} y$, o\`u $y \in \R$.\ps

\begin{pf} Soit $y \in \C$ tel que $|{\rm Im}\, \,  y | < \frac{1}{2}$, on a les \'egalit\'es
$$\Phi_{F}(\,\frac{1}{2}+i y) = 2 \pi \int_{\R} g(2\pi x)\, \frac{ e^{2i\pi x y}}{{\rm cosh} \,\pi x} \,{\rm d} x = \int_{\R} \frac{\widehat{g}(\frac{x}{2 \pi})}{{\rm cosh} \, \pi (x-y)} \,{\rm d} x.$$
En effet, la premi\`ere est triviale pour tout $y \in \C$. La seconde est par exemple une application de la formule de Plancherel, les fonctions $x \mapsto g(2 \pi x )$ et $x \mapsto \frac{e^{2i\pi x y}}{{\rm cosh} \pi x}$ \'etant sommables et de carr\'e sommable sur $\R$, et du fait que la fonction $\frac{1}{{\rm cosh}\, \pi x}$ est \'egale \`a sa transform\'ee de Fourier. On conclut car d'une part on a $\widehat{g}(z) \geq 0$ pour tout $z \in \R$ par hypoth\`ese, et d'autre part on constate l'in\'egalit\'e ${\rm Re} \,\frac{1}{{\rm cosh}\, \pi z} \,> 0$ pour tout $z \in \C$ tel que $|{\rm Im}\, z \,|< \frac{1}{2}$.
\end{pf}

\subsection{La fonction d'Odlyzko}\label{parfonodly}	Suivant Odlyzko \cite[\S 3]{poitou}, notons $g$ le double du carr\'e de convolution de
la fonction $u : \R \rightarrow \R$ d\'efinie par $u(x)={\rm cos}(\pi x )$
si $|x|\leq \frac{1}{2}$, et $u(x)=0$ sinon. Concr\`etement, la fonction $g$ est nulle
hors du segment $[-1,1]$, et pour $|x|\leq 1$ elle est donn\'ee par la formule  $$g(x)=(1-|x|)\, {\rm cos} (\pi
x)\,+\,\frac{1}{\pi}{\rm sin}(\pi |x|).$$   \ps
\noindent On v\'erifie imm\'ediatement que $g$ est une fonction paire,
positive, \`a support compact, de classe $\mathcal{C}^2$, et que l'on a $g(0)=1$. Sa transform\'ee de Fourier, \`a savoir $2 \, \widehat{u}^2$, est manifestement positive, car $u$ est r\'eelle, paire et sommable. En particulier, pour tout $\lambda \in \R_{>0}$ la fonction ${\rm F}_{\lambda} : \R \rightarrow \R$ d\'efinie par $${\rm F}_{\lambda}(x) = g(x/\lambda) /{\rm cosh}(x/2)$$ est une fonction test $\geq 0$ satisfaisant l'hypoth\`ese {\rm (T4)} (Lemme~\ref{lemeposphi}). \ps\ps


Expliquons comment \'evaluer num\'eriquement la forme forme lin\'eaire ${\rm J}_{{\rm F}_{\lambda}}$ sur $ {\rm K}_\infty$. On rappelle que $\psi(z)=\frac{\Gamma'(z)}{\Gamma(z)}$ d\'esigne la fonction digamma. On pose \'egalement, respectivement pour $z \in \C - \N$ et $z \in \C - \{\pm i\pi\}$,
$$\phi(z) = \frac{1}{2} \psi(\,\frac{z+1}{2}\,)-\frac{1}{2} \psi(\,\frac{z}{2}\,) \hspace{10pt}{\rm et}\hspace{10pt}{\rm r}(z)=2 \pi^{2} \frac{e^{-z}}{(z^{2}+\pi^{2})^{2}}.$$
On observera que l'on a la formule $\phi(z)= \sum_{n\geq 0} \frac{(-1)^n}{z+n}$. Nous remercions Henri Cohen de nous avoir mis sur la voie de la proposition suivante.

\begin{prop}\label{calculexplicite} Soit $\lambda$ un nombre r\'eel $>0$. \par \medskip \begin{itemize} \item[(i)] Pour tout entier $w\geq 0$ on a l'\'egalit\'e
$${\rm J}_{{\rm F}_{\lambda}}({\rm I}_w) = \log \, \pi \,- \,\Re \,\psi(b+\frac{i \pi}{2\lambda}) \,+\,\frac{1}{\pi} \Im \,\psi(b+\frac{i \pi}{2\lambda})\, -\, \frac{1}{2\lambda} \Re\, \psi'(b+\frac{i \pi}{2\lambda})  \,+ \,r_1(w,\lambda),$$
avec $b=\frac{1}{2}+\frac{w}{4}$ et $r_1(w,\lambda)= 2\lambda \sum_{n=0}^\infty {\rm r}(2\lambda(b+n))$. \ps\ps\ps
\item[(ii)] De plus, on a l'\'egalit\'e
$${\rm J}_{{\rm F}_{\lambda}}(1-\epsilon_{\C/\R}) = 1 \,+\, \frac{2\pi}{\lambda} \,\Im\, \phi(1+\frac{i\pi}{\lambda}) \,+ \,\frac{2\pi}{\lambda^2} \,\Im \, \phi'(1+\frac{i\pi}{\lambda})\,+\, r_2(\lambda),$$
avec  $r_2(\lambda)= 2 \lambda  \sum_{n=1}^\infty (-1)^{n+1} \, n \,{\rm r}(\lambda n)$. \ps \ps \ps
\item[(iii)] Enfin, on a $\Phi_{{\rm F}_{\lambda}}(0) = \frac{8}{\pi^2} \lambda$ et 
$$\Phi_{{\rm F}_{\lambda}}(\,\frac{1}{2}\,)\, =\, 4 \, \Re\, \phi(\,\frac{1}{2}+\frac{i\pi}{\lambda}) \, -\,  \frac{4}{\pi} \,\Im \, \phi(\,\frac{1}{2}+\frac{i\pi}{\lambda}) \, + \, \frac{4}{\lambda}\, \Re \, \phi'(\,\frac{1}{2}+\frac{i\pi}{\lambda}) \, + \, r_{3}(\lambda),$$ 
avec $r_3(\lambda)=  4 \lambda \sum_{n=0}^\infty (-1)^n\, {\rm r}(\lambda (n+\frac{1}{2}))$.
\end{itemize}
\end{prop}

\ps \medskip \medskip

\begin{pf} Si $\alpha>0$, on pose $h(\alpha) = \int_0^1 \, g(x) \, e^{-\alpha x} {\rm d}x$. En utilisant par exemple la d\'efinition $g = 2\, u \ast u$, on v\'erifie tout d'abord l'identit\'e
$g''(x)+\pi^2 g(x) = 2\pi |\sin \pi x|$ pour $|x|\leq  1$, puis la relation
 
$$
\hspace{24pt}
h(\alpha) \hspace{4pt} =\hspace{4pt}
\frac{\alpha}{\alpha^{2}+\pi^{2}}
+2\hspace{1pt}\pi^{2}\hspace{2pt}\frac{1+e^{-\alpha}}{(\alpha^{2}+\pi^{2})^{2}}
\hspace{24pt}.
$$
\noindent On a de plus $\int_{0}^{\infty}
g(x/\lambda)\hspace{2pt}e^{-\alpha x}\hspace{2pt}\mathrm{d}x = \lambda\, h( \lambda\, \alpha)$. \ps \medskip

Soit $w\geq 0$ un entier ; on pose $b=\frac{1}{2}+\frac{w}{4}$. La proposition~\ref{lemmemestre} appliqu\'ee \`a la fonction $F={\rm F}_{\lambda}$ s'\'ecrit ${\rm J}_\lambda({\rm I}_{w}) \hspace{4pt} = \hspace{4pt}\log(2\pi) \hspace{4pt}+\hspace{4pt} \sigma_{{\rm F}_{\lambda}}(1,\frac{w}{2})$. Des manipulations \'el\'ementaires conduisent \`a la relation
\begin{multline*}
\sigma_{{\rm F}_{\lambda}}(1,\frac{w}{2}) \hspace{4pt} = \hspace{4pt}
\int_{0}^\infty (\hspace{4pt}\frac{2\hspace{1pt}e^{-2b\hspace{1pt}x}}{1-e^{-2x}}-\frac{e^{-x}}{x}\hspace{4pt})
\hspace{2pt}\mathrm{d}x
\hspace{4pt} +  \hspace{4pt}
\int_{0}^{\infty}
(g(\frac{x}{\lambda})-1)
\hspace{4pt}\frac{2\hspace{1pt}e^{-2b\hspace{1pt}x}}{1-e^{-2x}}
\hspace{4pt}\mathrm{d}x
\hspace{24pt}.
\end{multline*}
\'Ecrivons cette somme sous la forme \'evidente $\sigma_{\lambda}(1,\frac{w}{2})=S_1+S_2$. La formule de Gauss $\psi(z) = - \int_0^\infty \left(\frac{e^{-zx}}{1-e^{-x}}\hspace{4pt} - \hspace{4pt}\frac{e^{-x}}{x}\right) {\rm d} x$, et l'identit\'e $\int_{0}^{\infty} \frac{e^{-\alpha x}-e^{{-x}}}{x} {\rm d} x = - \log\, \alpha$ pour $\alpha>0$,  entra\^inent l'\'egalit\'e $S_1 \, \, \, =\, \, \, -\,{\rm log} \, \,2\,\, -\,\, \psi(b)$.  D'autre part, le d\'eveloppement $\frac{e^{-2b x}}{1-e^{-2x}}=\sum_{n\geq 0} \,e^{-2(b+n)x}$ conduit \`a la relation
$$
\hspace{24pt}
S_2
\hspace{4pt}=\hspace{4pt}
2\hspace{1pt}\lambda
\hspace{2pt}\sum_{n=0}^{\infty}
\hspace{4pt}\left(
h(2\lambda\hspace{1pt}(b+n))-\frac{1}{2\lambda\hspace{1pt}(b+n)}\hspace{4pt}\right)
\hspace{24pt}.
$$
On a $h=h_1+h_2+{\rm r}$ avec $h_1(\alpha)\, =\, \frac{\alpha}{\alpha^{2}\,+\,\pi^{2}}$ et $h_2(\alpha)\,=\,\frac{2\, \pi^2}{(\alpha^{2}+\pi^{2})^{2}}$. De plus, si $u$ et $v$ sont des r\'eels $\neq 0$ on a les identit\'es $\frac{u}{u^2+v^2} - \frac{1}{u} =  \Re \, \,(\frac{1}{u+vi}\,-\,\frac{1}{u})$ et
$$\frac{2v}{(u^2+v^2)^2} \, \, = \, \,-\, \, \frac{1}{v^2}\,\, \Im \,\,(\frac{1}{u+vi}\,-\,\frac{1}{u})\, \, - \, \, \frac{1}{v}\, \, \Re\, \, \frac{1}{(u+vi)^2}.$$
\noindent Appliquons-les \`a $v=\frac{\pi}{2\lambda}$ et $u=b+n$ pour tout entier $n\geq 0$ et sommons. La formule $\psi(b)-\psi(b+z)=\sum_{n\geq 0}  \frac{1}{b+z+n}- \frac{1}{b+n}$ avec $z=\frac{i\pi}{2\lambda}$ entra\^ine alors 
{\small $$\sum_{n\geq 0} \,2\lambda \hspace{4pt} h_1(\,2 \lambda (b+n)\,) - \frac{1}{b+n}\, =\, \sum_{n\geq 0} \,\frac{(b+n)}{(b+n)^{2}+(\frac{\pi}{2\lambda})^{2}} -\frac{1}{b+n}\,=\, \psi(b)\,- \Re\, \psi(b+\frac{i\pi}{2\lambda}).$$ }
De m\^eme, compte tenu de l'identit\'e $\psi'(z) = \sum_{n\geq 0} \frac{1}{(z+n)^2}$, on trouve 
$$\sum_{n\geq 0} \,2\lambda \hspace{4pt} h_2(\,2 \lambda (b+n)\, )\,=\, \frac{1}{\pi}\, \Im \, \psi(b+\frac{i\pi}{2\lambda}) - \frac{1}{2\lambda} \psi'(b+\frac{i\pi}{2\lambda}).$$
En mettant toutes les formules bout \`a bout, on obtient l'assertion (i). La d\'emonstration de l'assertion (ii) est similaire. \`A partir de la proposition~\ref{lemmemestre}, on commence par \'etablir l'\'egalit\'e : 
$${\rm J}_{{\rm F}_{\lambda}}(1-\epsilon_{\C/\R}) \hspace{4pt} = 1 + \hspace{4pt} \int_0^\infty \left(\, g(\frac{x}{2\lambda})-1\right) \hspace{4pt}\frac{e^{-x/2}}{(1\hspace{4pt}+\hspace{4pt}e^{-x/2})^2} \hspace{4pt}{\rm d} x\hspace{4pt}.$$
Le d\'eveloppement $\frac{e^{{-x/2}}}{(1+e^{{-x/2}})^2}=\sum_{n\geq 1}\, (-1)^{n+1}\, n\, e^{-nx/2}$ permet d'\'ecrire 
$${\rm J}_{{\rm F}_{\lambda}}(1-\epsilon_{\C/\R}) \hspace{4pt} \,=\, 1 \,+ \,2 \hspace{4pt}\sum_{n\geq 1} \hspace{4pt}(-1)^{n+1}\hspace{4pt} ( \lambda n \, h(\lambda n) \hspace{4pt}-\hspace{4pt}1).$$
De plus, si $u$ et $v$ sont des r\'eels non nuls on a les identit\'es 
$$  \frac{u^2}{u^2+v^2} - 1 \, \, =\, \, v \, \Im\, \, \frac{1}{u+vi}\, \,\hspace{10pt}{\rm et}\hspace{10pt} \frac{2\,u\,v^2}{(u^2+v^2)^2} \, \, =\, \, -v\, \Im \,\, \frac{1}{(u+vi)^2}.$$
En posant $u=n \lambda$ et $v=\frac{\pi}{\lambda}$, et en observant $\phi(z+1)=\sum_{{n\geq 1}} \frac{(-1)^{n+1}}{z+n}$, on a $$2\,\sum_{n\geq 1} \hspace{4pt}(-1)^{n+1}\hspace{4pt} ( \lambda n \, h_{1}(\lambda n) \hspace{4pt}-\hspace{4pt}1)\,=\, \frac{2\pi}{\lambda} \, \Im \,\phi(1+\frac{i\pi}{\lambda}),$$ 
puis $2 \,\sum_{n\geq 1} \hspace{4pt}(-1)^{n+1}\hspace{4pt}  \lambda n \, h_{2}(\lambda n)\,=\,  \frac{2 \pi}{\lambda^{2}} \, \Im \, \phi'(1+\frac{i\pi}{\lambda})$. Cela d\'emontre le (ii). Il ne reste qu'\`a v\'erifier l'assertion (iii). Par d\'efinition de ${\rm F}_{\lambda}$ et $g$, on a les \'egalit\'es $\Phi_{{\rm F}_{\lambda}}(0)=\int_\R g(x/\lambda) {\rm d}\, x= 2 \lambda\, \widehat{u}(0)^2$ ; la valeur de  $\Phi_{{\rm F}_{\lambda}}(0)$ se d\'eduit donc de la relation imm\'ediate $\widehat{u}(0)=\frac{2}{\pi}$. Pour d\'eterminer $\Phi_{{\rm F}_{\lambda}}(\frac{1}{2})$, on proc\`ede comme pour les assertions (i) et (ii) \`a partir des identit\'es imm\'ediates
$$\Phi_{{\rm F}_{\lambda}}(\frac{1}{2})=4 \int_{0}^{\infty} g(x/\lambda) \frac{e^{-x/2}}{1+e^{{-x}}} {\rm d}x= 4 \lambda \sum_{{n\geq 0}} (-1)^{n} \, h(\lambda(n+1/2)).$$

\end{pf}
\ps \medskip

\noindent {\sc Commentaires sur les calculs num\'eriques effectu\'es dans les paragraphes qui vont suivre}

\bigskip

\noindent 1) Les formules de la proposition \ref{calculexplicite}, bien que peu esth\'etiques, sont tr\`es efficaces pour \'evaluer num\'eriquement ${\rm J}_{\mathrm{F}_{\lambda}}$ et $\Phi_{\mathrm{F}_{\lambda}}(\frac{1}{2})$, et ce avec une pr\'ecision arbitraire. Dans nos applications, nous aurons $0 \leq w \leq 46$ et $\lambda = \log N$ avec $2 \leq N \leq 100$.

\bigskip
\noindent  2) John L. Spouge a \'elabor\'e en 1994 un remarquable algorithme pour d\'eterminer les valeurs des fonctions gamma, $\psi$ (digamma) et $\psi'$ (trigamma).  Ces fonctions sont impl\'ement\'ees dans \texttt{PARI}. Cependant, Henri Cohen nous a fait savoir que \texttt{PARI} utilise pour leur calcul la formule d'Euler-MacLaurin, donc les nombres de Bernoulli, et qu'apr\`es le calcul d'une premi\`ere valeur celui des suivantes est acc\'el\'er\'e par le stockage des nombres de Bernoulli.

\bigskip
\noindent  3) Les trois fonctions $r_{1}(w,\lambda)$, $r_{2}(\lambda)$ et $r_{3}(\lambda)$ qui apparaissent dans l'\'enonc\'e \ref{calculexplicite} sont d\'efinies comme des sommes de s\'eries~; on estime ci-dessous ``les restes'' $\sum_{n=N+1}^{\infty}$ de ces s\'eries.

\ps\ps
\noindent  3.1) On a l'in\'egalit\'e
$$
\hspace{24pt}
0
\hspace{4pt}\leq\hspace{4pt}
r_{1}(w,\lambda)-2\lambda\sum_{n=0}^{N}\mathrm{r}(2\lambda(b+n))
\hspace{4pt}\leq\hspace{4pt}
\frac{2\lambda\hspace{1pt}\mathrm{r}(2\lambda(b+N+1))}{1-e^{-2\lambda}}
\hspace{24pt};
$$
cette in\'egalit\'e r\'esulte de ce que l'on a
$$
0
\hspace{4pt}\leq\hspace{4pt}
\mathrm{r}(2\lambda(b+n))
\hspace{4pt}\leq\hspace{4pt}e^{-2\lambda(n-N-1)}\hspace{1pt}\mathrm{r}(2\lambda(b+N+1))
$$
pour tout $n$ avec $n\geq N+1$.

\ps\ps
\noindent  3.2) Comme la fonction d'une variable r\'eelle $x\mapsto x\hspace{1pt}\mathrm{r}(x)$ est \`a valeurs positives et d\'ecroissante pour $x\geq 0.773$, on a l'in\'egalit\'e
$$
\hspace{24pt}
0
\hspace{4pt}\leq\hspace{4pt}
(-1)^{N}\hspace{2pt}(\hspace{2pt}r_{2}(\lambda)-2\lambda\sum_{n=1}^{N}(-1)^{n+1} n \hspace{1pt}\mathrm{r}(\lambda n)\hspace{2pt})
\hspace{4pt}\leq\hspace{4pt}
2\lambda (N+1) \hspace{1pt}\mathrm{r}(\lambda(N+1))
\hspace{20pt},
$$
sous l'hypoth\`ese $\lambda (N+1)\geq 0.773$.

\ps\ps
\noindent  3.3) Comme la fonction d'une variable r\'eelle $x\mapsto\mathrm{r}(x)$ est \`a valeurs positives et d\'ecroissante, on a l'in\'egalit\'e
$$
\hspace{18pt}
0
\hspace{4pt}\leq\hspace{4pt}
(-1)^{N+1}\hspace{2pt}(\hspace{2pt}r_{3}(\lambda)-4\lambda\sum_{n=0}^{N}(-1)^{n}\hspace{1pt}\mathrm{r}(\lambda(n+\frac{1}{2}))\hspace{2pt})
\hspace{4pt}\leq\hspace{4pt}
4\lambda\hspace{1pt}\mathrm{r}(\lambda(N+\frac{3}{2}))
\hspace{18pt}.
$$

\bigskip
\noindent  4) Les calculs effectu\'es \`a l'aide des formules de la proposition \ref{calculexplicite} sont confirm\'es, avec une pr\'ecision arbitraire, par les routines d'int\'egration num\'erique de \texttt{PARI}.

\bigskip
\noindent  5) La fonction $w\mapsto\mathrm{J}_{\mathrm{F}_{\lambda}}(\mathrm{I}_{w})$, $\lambda$ fix\'e, est d\'ecroissante. On s'en convainc en invoquant le point (iii) de la proposition-d\'efinition \ref{lemmemestre} qui montre que sous les hypoth\`eses $F\geq 0$ et $w'\geq w$ la diff\'erence $\mathrm{J}_{F}(\mathrm{I}_{w})-\mathrm{J}_{F}(\mathrm{I}_{w'})$ est l'int\'egrale sur $[0,+\infty[$ d'une fonction positive ou nulle. \ps\ps

\noindent De m\^eme, la fonction $\lambda\mapsto\mathrm{J}_{\mathrm{F}_{\lambda}}(\mathrm{I}_{w})$, $w$ fix\'e, est croissante. En effet, comme la fonction d'Odlyzko est d\'ecroissante sur $[0,+\infty[$, le point (iii) de la propo\-sition-d\'efinition \ref{lemmemestre} montre que sous l'hypoth\`ese $\lambda\leq\lambda'$ la diff\'erence $\mathrm{J}_{\mathrm{F}_{\lambda'}}(\mathrm{I}_{w})-\mathrm{J}_{\mathrm{F}_{\lambda}}(\mathrm{I}_{w})$ est l'int\'egrale sur $[0,+\infty[$ d'une fonction positive ou nulle (on observera que l'on a $\mathrm{F}_{\lambda}(0)=1$ pour tout $\lambda$).

\bigskip

\noindent Ce qui pr\'ec\'ede implique que l'on l'encadrement $
-1.40
\hspace{4pt}\leq\hspace{4pt}
\mathrm{J}_{\mathrm{F}_{\lambda}}(\mathrm{I}_{w})
\hspace{4pt}\leq\hspace{4pt}
2.63
$
pour $\log 2\leq\lambda\leq\log 100$ et $0\leq w\leq 50$. On trace ci-contre le graphe de la fonction $w\mapsto\mathrm{J}_{\mathrm{F}_{\lambda}}(\mathrm{I}_{w})$ pour quelques valeurs du param\`etre $\lambda$.

\begin{figure}[!h]
\hspace{-1 cm}\includegraphics[width=15cm, height=10cm]{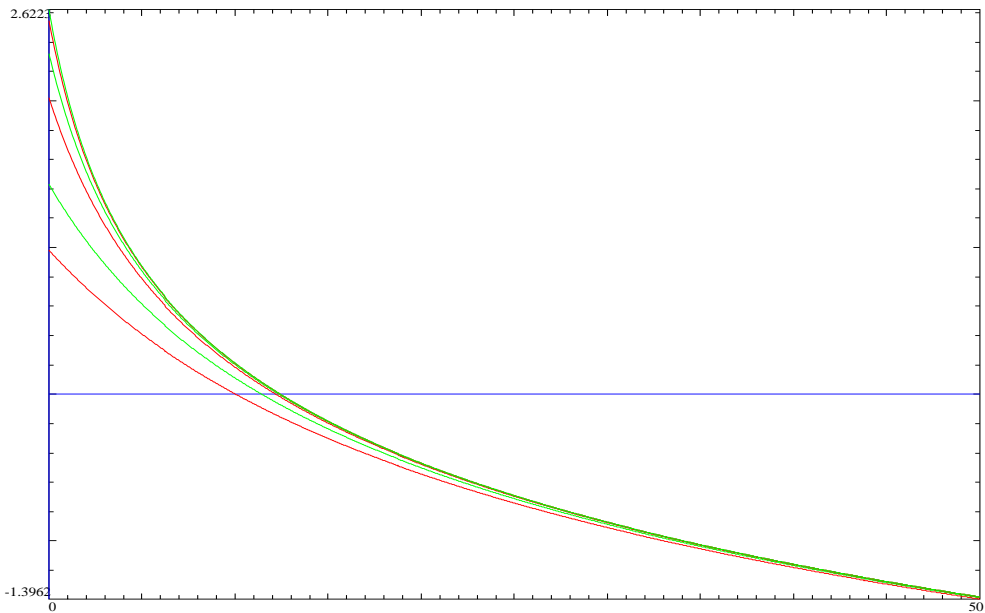} 
\caption{Graphe de la fonction $w \mapsto {\rm J}_{{\rm F}_{\lambda}}({\rm I}_w)$ sur l'intervalle $0 \leq w \leq 50$, pour les param\`etres $\lambda=\log \,2, \, \log\,3,\, \log\,8,\,\log\,20,\log\,60$ et $\,\log\,100$.}
\label{grapheJwl}
\end{figure}

\bigskip
\noindent  6) Dans l'appendice de \cite{mestre} Mestre d\'ecrit un calcul de $\sigma_{\mathrm{F}_{\lambda}}(1,\frac{1}{2})$ (autre\-ment dit $\mathrm{J}_{\mathrm{F}_{\lambda}}(\mathrm{I}_{1})-\log 2\pi)$ pour $\lambda<\pi$ par une m\'ethode totalement diff\'erente de la n\^otre~; la restriction $\lambda<\pi$ tient \`a ce que Mestre utilise implicitement le d\'eveloppement en s\'erie enti\`ere de la fonction holomorphe $\frac{z}{e^{z}-1}$ dans le disque $\vert z\vert<2\pi$.

\bigskip

\subsection{D\'ebut de d\'emonstration du th\'eor\`eme \ref{classpoids22} : cas $w\leq 20$} Pour tout entier $w \geq 0$, on consid\`ere le sous-groupe de ${\rm K}_{\infty}$ d\'efini de la mani\`ere suivante : 
$${\rm K}_{\infty}^{\leq w} \,=\, \left\{ \begin{array}{lll}  \,\left(\, \bigoplus_{1 \leq j \leq w/2} \Z \, {\rm I}_{2j}\right)\, \oplus\, \Z \,1 \,\oplus\, \Z \, \epsilon_{\C/\R}, & & {\rm si}\, w \equiv 0 \bmod 2, \\ \\ \bigoplus_{1 \leq j \leq \frac{w+1}{2}} \Z \, {\rm I}_{2j-1}\,,
&  & {\rm si}\, w \equiv 1 \bmod 2. \end{array} \right.$$
\noindent L'int\'er\^et de cette d\'efinition est que si $\pi \in \Pi_{{\rm alg}}(\PGL_{n})$ est de poids motivique $w$, on a ${\rm L}(\pi_{\infty}) \in {\rm K}_{\infty}^{\leq w}$ (Proposition VIII.\ref{contraintesalg}). On dispose de deux autres contraintes simples sur ${\rm L}(\pi_\infty)$. D'une part, c'est un \'el\'ement effectif de ${\rm K}_\infty$. D'autre part, on a la relation $\DET\, {\rm L}(\pi_\infty)=1$. On rappelle l'identit\'e $\DET \, {\rm I}_v = \epsilon_{\C/\R}^{v+1}$ pour tout $v\geq 0$. \bigskip

Le principe g\'en\'eral de la d\'emonstration va consister \`a montrer que sous une hypoth\`ese convenable sur $\pi_\infty$, il n'existe pas -- ou peu -- de repr\'esentations $\pi \in \Pi_{\rm alg}$ de poids motivique $\leq 22$, et ce \`a l'aide des in\'egalit\'es donn\'ees par la proposition \ref{formexplpipi} appliqu\'ee aux fonctions tests introduites au \S \ref{parfonodly}. Dans toute cette partie, $\lambda$ d\'esignera donc un nombre r\'eel $>0$, et l'on consid\`ere la fonction test associ\'ee ${\rm F}_{\lambda}$ d\'efinie au \S \ref{parfonodly}. Dans la d\'emonstration, nous ferons un usage syst\'ematique de la forme bilin\'eaire $\cbil_{\infty}^{{\rm F}_{\lambda}}$ sur ${\rm K}_\infty$ introduite au \S\ref{preliminairesfexpl}. Nous aurons notamment besoin d'\'evaluer explicitement cette forme bilin\'eaire, ce que nous ferons bien entendu \`a l'aide de l'ordinateur, en utilisant les formules d\'ecrites au \S \ref{parfonodly}. Nous renvoyons le lecteur \`a la feuille de calculs \cite{clcalc} donn\'ee en annexe pour une justification des calculs num\'eriques que nous ferons ci-dessous. Afin d'acqu\'erir une intuition sur la m\'ethode, nous allons commencer par d\'emontrer le cas $n = 2$ du th\'eor\`eme. La repr\'esentation triviale \'etant le seul \'el\'ement de $\Pi(\PGL_1)=\Pi_{\rm cusp}(\PGL_1)$, observons qu'il n'y a rien \`a d\'emontrer pour $n=1$ ! \ps \ps 
\bigskip
\noindent {\sc Cas $n=2$ ou $w \leq 10$}\ps\ps
\medskip
 Supposons donn\'ee une repr\'esentation $\pi \in \Pi_{\rm alg}(\PGL_2)$ de poids motivique $w$. La condition $\DET \,\,{\rm L}(\pi_\infty)\,=\,1$ montre qu'il y a deux cas : soit $w$ est impair et ${\rm L}(\pi_\infty) = {\rm I}_w$, soit $w=0$ et l'on a ${\rm L}(\pi_\infty)=2 \cdot 1$ ou ${\rm L}(\pi_\infty)=2 \cdot \epsilon_{\C/\R}$. \ps
 
 Si $w$ est impair, la proposition \ref{proppigl2} assure que $\pi$ est la repr\'esentation engendr\'ee par une unique forme modulaire propre et normalis\'ee de poids $w+1$ pour le groupe ${\rm SL}_2(\Z)$ ; r\'eciproquement, toute telle forme engendre bien \'evidemment un tel $\pi$. Autrement dit, compte tenu de la d\'efinition \S \ref{deffonctionm} nous avons pour tout $w$ impair $\geq 1$ l'\'egalit\'e :
 $${\rm m}({\rm I}_w) = \dim  {\rm S}_{w+1}({\rm SL}_2(\Z))$$
\noindent  Ainsi, le cas particulier $n=2$ et poids motivique $\neq 0$ du th\'eor\`eme~\ref{classpoids22} est cons\'equence de la description bien connue de ${\rm S}_{w+1}({\rm SL}_2(\Z))$ \cite[Ch. VII]{serre}. De telles descriptions de $\Pi_{\rm alg}(\PGL_n)$ n'existant pas en dimension $n>2$, il y a un int\'er\^et \`a expliquer comment proc\'eder autrement, si possible. Dans ce qui suit, on se propose de d\'emontrer diff\'eremment que ${\rm S}_{w+1}({\rm SL}_2(\Z))$ est nul si $w<11$ ou $w=13$, et de dimension $\leq 1$ si $15 \leq w \leq 21$, \`a l'aide du corollaire \ref{inegexpl}. \ps \ps

La fonction ${\rm F}_{\lambda}$ \'etant continue et \`a support dans $[-\lambda,\lambda]$, observons que pour tous $\pi,\pi' \in \Pi_{\rm alg}$ on a 
\begin{equation}\label{formeBlambda} \widetilde{\cbil}_f^{{\rm F}_{\lambda}}\,= \,\sum_{p^k < e^\lambda} \,\,{\rm F}_\lambda(k {\rm log}\,p)\frac{{\rm log}\,p}{p^{k/2}}\,\,\,\,\overline{{\rm tr}\, ({\rm c}_p(\pi)^k)} \,\, {\rm tr} \,({\rm c}_p(\pi')^k).\end{equation}
En particulier, on constate $\cbil_f^{{\rm F}_{\lambda}}=0$, {\rm i.e.} $\cbil^{{\rm F}_{\lambda}}=\cbil_\infty^{{\rm F}_{\lambda}}$, pour $\lambda \leq  \log\,2$. On consid\`ere maintenant $\pi \in \Pi_{\rm alg}(\PGL_2)$ de poids motivique impair $w$ et $\pi'=1 \in \Pi_{\rm alg}(\PGL_1)$, auquel cas on a 
$\cbil_\infty^{{\rm F}_{\lambda}}(\pi,1)={\rm J}_{{\rm F}_{\lambda}}({\rm I}_w)$. L'in\'egalit\'e du corollaire \ref{inegexpl} entra\^ine donc en particulier $${\rm J}_{{\rm F}_{{\log \,2}}}({\rm I}_w) \leq 0.$$\noindent La table \ref{tablejwlog2} donne l'\'evaluation num\'erique \`a $10^{-2}$ pr\`es de ${\rm J}_{{\rm F}_{{\log \,2}}}({\rm I}_w)$, pour $1 \leq w \leq 21$ impair. Elle contredit l'in\'egalit\'e ci-dessus si $w<11$ : la repr\'esentation $\pi$ n'existe pas. Soulignons que cet argument n'est pas nouveau : c'est par exemple exactement la m\'ethode utilis\'ee par Mestre dans \cite{mestre}. Sa pr\'ecision est assez \'etonnante, puisque l'on sait bien que l'on a ${\rm S}_{12}({\rm SL}_2(\Z)) \neq 0$. Nous verrons plus loin de nombreux autres exemples de la pr\'ecision fascinante des formules explicites.\ps\medskip

\begin{table}[h!]
\renewcommand{\arraystretch}{1.5}
{\scriptsize \begin{tabular}{c||c|c|c|c|c|c|c|c|c|c|c}
$w$ & $1$ & $3$ & $5$ & $7$ & $9$ & $11$ & $13$ & $15$ & $17$ & $19$ & $21$ \cr
\hline
${\rm J}_{{\rm F}_{{\log \,2}}}({\rm I}_w)$ & $0.85$ & $0.61$ & $0.41$ & $0.23$ & $0.07$ & $-0.06$ & $-0.19$ & $-0.30$ & $-0.40$ & $-0.50$ & $-0.59$\cr \end{tabular} }
\caption{Valeurs \`a $10^{-2}$ pr\`es de ${\rm J}_{{\rm F}_{{\log \,2}}}({\rm I}_w)$ pour $1\leq w \leq 21$ impair.}
\label{tablejwlog2}
\end{table}
\renewcommand{\arraystretch}{1}
\ps \medskip
\noindent 
Comme le remarque essentiellement Mestre \cite[Remarque 1 \S III]{mestre}, l'argument pr\'ec\'edent est de port\'ee plus vaste : si $\pi \in \Pi_{\rm alg}(\PGL_n)$ est de poids motivique $\leq 10$, alors $n=1$ et $\pi$ est la repr\'esentation triviale. En effet, si $\pi \neq 1$ on a l'in\'egalit\'e $\cbil_\infty^{{\rm F}_{\log 2}}(\pi,1) \leq 0$ alors que l'on v\'erifie num\'eriquement $\cbil_\infty^{{\rm F}_{\log 2}}(V,1)={\rm J}_{{\rm F}_{\log \,2}}(V)>0.02$ pour $V=1,\epsilon_{\C/\R}$ et ${\rm I}_w$ pour $w<11$. \ps \ps

Expliquons maintenant comment aller plus loin dans le cas $n=2$, \`a commencer par \'eliminer le cas $w=13$. Comme tout \'el\'ement de $\Pi(\PGL_2)$, la repr\'esentation $\pi$ est autoduale ; on a donc ${\rm L}(s,\pi)=\varepsilon(\pi){\rm L}(1-s,\pi)$ avec $\varepsilon(\pi)=\varepsilon({\rm I}_w)=i^{w+1}$. Ce facteur epsilon vaut $-1$ si $w \equiv 1 \bmod 4$. D'apr\`es le corollaire \ref{inegexpl}, on doit donc avoir  $${\rm J}_{{\rm F}_{{\log \,2}}}({\rm I}_w) \leq -\frac{1}{2}\Phi_{{\rm F}_{\log\,2}}(\,\frac{1}{2}\,)$$ pour $w \equiv 1 \bmod 4$. Une \'evaluation num\'erique montre que le nombre $\frac{1}{2}\Phi_{{\rm F}_{\log\,2}}(\,\frac{1}{2}\,)$ vaut $0.28$ \`a $10^{-2}$ pr\`es. La table pr\'ec\'edente montre donc \'egalement que $\pi$ n'existe pas pour $w=13$.  \ps \ps

Montrons maintenant ${\rm m}({\rm I}_w)\leq 1$ pour $w \leq 21$. On a $${\rm m}(w) \,\,\cbil_\infty^{{\rm F}_{\lambda}}({\rm I}_w,{\rm I}_w) \,\,\leq \Phi_{{\rm F}_{\lambda}}(0) \hspace{15pt}{\rm (on\,\, rappelle}\,\, \Phi_{{\rm F}_{\lambda}}(0)=\,\,\frac{8}{\pi^2} \, \lambda\,\,{\rm )}$$ 
pour tout $\lambda >0$ d'apr\`es le corollaire \ref{cortaibi}. D'autre part, on a $$\cbil_\infty^{{\rm F}_{\lambda}}({\rm I}_w,{\rm I}_w)\,=\,{\rm J}_{{\rm F}_{\lambda}}({\rm I}_{2 w})+{\rm J}_{{\rm F}_{\lambda}}({\rm I}_0)$$ d'apr\`es la relation ${\rm I}_u \otimes {\rm I}_v = {\rm I}_{u+v}+{\rm I}_{|u-v|}$.  L'in\'egalit\'e ci-dessus s'av\`ere remarquablement bonne exp\'erimentalement pour des petites valeurs de $w$. Par exemple, lorsque $\lambda = \,{\rm log}\, 8$ on v\'erifie num\'eriquement  ${\rm m}({\rm I}_{11}) \leq 1.17$, ${\rm m}({\rm I}_{15}) \leq 1.48$, ${\rm m}({\rm I}_{17}) \leq 1.66$, ${\rm m}({\rm I}_{19}) \leq 1.86$ et ${\rm m}({\rm I}_{21}) \leq 2.08$.  On en d\'eduit les in\'egalit\'es annonc\'ees ${\rm m}({\rm I}_w) \leq 1$ pour $w\leq 19$, ainsi que ${\rm m}({\rm I}_{21})\leq 2$. Cette derni\`ere \'egalit\'e ne semble toutefois pas am\'eliorable en consid\'erant simplement d'autres valeurs de $\lambda$. \ps \ps

Pour conclure ${\rm m}({\rm I}_{21})\leq 1$ nous allons utiliser le corollaire \ref{inexploo2} (ii) appliqu\'e \`a $V=1$ et $V'={\rm I}_{21}$. On effet, on a ${\rm m}(1)={\rm m}^\bot(1)=1$, ${\rm m}({\rm I}_{21})={\rm m}^\bot({\rm I}_{21})$, et $\varepsilon({\rm I}_{21})=-1$. On en d\'eduit pour tout $\lambda >0$ l'in\'egalit\'e
$$\frac{1}{2} \,\Phi_{{\rm F}_{\lambda}}(\frac{1}{2}) \,+ {\rm J}_{{\rm F}_{\lambda}}({\rm I}_{21})\, \,\leq \, \sqrt{(\Phi_{{\rm F}_{\lambda}}(0)-{\rm J}_{{\rm F}_{\lambda}}(1))\,(\,\frac{\Phi_{{\rm F}_{\lambda}}(0)}{{\rm m}({\rm I}_{21})}-\cbil_\infty^{{\rm F}_{\lambda}}({\rm I}_{21},{\rm I}_{21}))}.$$
Mais pour $\lambda \,=\, \log \,6$ on v\'erifie qu'\`a $10^{-2}$ pr\`es le terme de gauche vaut $0.17$, alors que celui de droite vaut $0.13$  si ${\rm m}({\rm I}_{21})=2$ (et $0.51$ si ${\rm m}({\rm I}_{21})=1$ !). $\square$
\ps \ps
 
\ps\ps \medskip

\noindent {\sc Cas $w$ impair $\leq 19$}
\ps \ps \medskip

Notre but dans ce qui suit est de d\'emontrer que si $\pi \in \Pi_{\rm alg}(\PGL_n)$ est de poids motivique $\leq 19$, et si $n>2$, alors $n=4$ et $\pi$  est l'unique repr\'esentation autoduale telle que ${\rm L}(\pi_\infty)={\rm I}_{19}\oplus {\rm I}_{17}$.  Notre point de d\'epart est le r\'esultat suivant. \ps\ps

\begin{lemme}\label{lemmedefpos} Si $\lambda={\rm log}\,9$, la restriction de $\cbil_\infty^{{\rm F}_{\lambda}}$ 
 \`a ${\rm K}_\infty^{\leq 19}$ est d\'efinie positive.
\end{lemme}

Si $w\geq 0$, nous noterons ${\rm Gram}(w,\lambda)$ la matrice de Gram de la forme bilin\'eaire $\cbil_\infty^{{\rm F}_{\lambda}}$ sur ${\rm K}_{\infty}^{\leq w}$ dans la $\Z$-base naturelle d\'efinissant ${\rm K}_\infty^{\leq w}$, \`a savoir : \ps\ps

\noindent -- dans le cas o\`u $w$ est impair, les ${\rm I}_v$ avec $1\leq v\leq w$ et $v$ impair, \ps\ps

\noindent -- dans le cas o\`u $w$ est pair, les repr\'esentations $1$, $\epsilon_{\C/\R}$ et ${\rm I}_v$ avec $0\leq v \leq w$.\ps\ps

\begin{pf}  Soit $B={\rm Gram}(19,\log 9)$. Les formules de la proposition \ref{calculexplicite} (et le logiciel \texttt{PARI}) permettent de calculer les coefficients de $B$, et de ${\rm Gram}(w,\lambda)$ en g\'en\'eral, avec une pr\'ecision th\'eoriquement arbitraire (dans la feuille de calculs \cite{clcalc} fournie, plus de $20$ chiffres significatifs), et l'on constate que $B$ est bien d\'efinie positive. L'argument que nous d\'egageons ci-dessous montre que nous n'avons besoin en fait que de tr\`es peu de chiffres significatifs ; cet argument sera surtout utile dans la d\'emonstration du lemme \ref{borneVlemm}, o\`u il s'agira de d\'eterminer les vecteurs $v \in {\rm K}_\infty^{\leq 19}$ tels que $\cbil_{\infty}^{{\rm F}_{\lambda}}(v,v)$ est inf\'erieur \`a une certaine constante. \ps\ps


On constate d'abord num\'eriquement que tous les coefficients de $B$ sont de valeur absolue dans l'intervalle $]0.01,3.48[$. Soit $A \in 10^{-4}\,{\rm M}_{10}(\Z)$ la matrice sym\'etrique obtenue en arrondissant l'approximation de $B$ donn\'ee par l'ordinateur \`a l'\'el\'ement le plus proche de $10^{-4}\,\Z$, de sorte que la matrice $A-B$ ait tous ses coefficients $\leq 10^{-4}$ en valeur absolue. Un calcul exact fait par l'ordinateur d\'emontre que $A$ est d\'efinie positive (il suffit d'appliquer le crit\`ere de Sylvester). Soit $||.||$ la norme sur ${\rm M}_n(\R)$ subordonn\'ee \`a la norme $\sup_i |x_i|$ sur $\R^n$, de sorte que $||(m_{i,j})|| = {\rm sup}_i \sum_j |m_{i,j}|$. On constate  $||A^{-1}||\leq 3.23$ et l'on a de plus $||A-B|| \leq 10 \cdot 10^{-4}=10^{-3}$. Cela entra\^ine que le rayon spectral de $A^{-1}(A-B)$ est $\leq 0.00323$, et en particulier que $B$ est d\'efinie positive d'apr\`es le lemme \ref{inegatrivquad}.
\end{pf}

\begin{lemme}\label{inegatrivquad} Soient $V$ un $\R$-espace vectoriel de dimension finie, $b_1,b_2$ deux formes bilin\'eaires sym\'etriques sur $V$, $e=(e_1,\dots,e_n)$ une base de $V$, et $M_i$ la matrice de Gram de $b_i$ dans la base $e$. On suppose $b_2$ d\'efinie positive.
Pour tout $x \in V$ on a l'in\'egalit\'e
$$|b_1(x,x) - b_2(x,x) |\, \leq\, \rho(M_2^{-1}(M_1-M_2)) \,b_2(x,x),$$
o\`u $\rho(M)$ d\'esigne le rayon spectral de la matrice $M$. Si l'on a de plus l'in\'egalit\'e $\rho(M_2^{-1}(M_1-M_2))<1$ alors $b_1$ est d\'efinie positive, et on a pour tout $x \in V$ $$b_2(x,x) \leq (1- \rho(M_2^{-1}(M_1-M_2)))^{-1}\, b_1(x,x).$$
\end{lemme}

\begin{pf} C'est une cons\'equence classique de la diagonalisablilit\'e des endomorphismes autoadjoints d'un espace euclidien.  \end{pf}

 
 \begin{lemme} \label{borneVlemm} Soit $V \in {\rm K}_\infty^{\leq 19}$ effectif et tel que $\cbil_\infty^{{\rm F}_{\log\, 9}}(V,V) \leq \frac{8 \, {\rm log}\, 9}{\pi^{2}}$, alors :\ps\ps \begin{itemize}
\item[(i)] soit $V={\rm I}_w$ avec $9 \leq w \leq 19$,\ps\ps
\item[(ii)] soit $V={\rm I}_{19}+{\rm I}_v$ avec $5 \leq v \leq 13$.
\end{itemize}
\end{lemme}
\ps\ps

 \begin{pf} Le lemme \ref{lemmedefpos} assure qu'il n'existe qu'un
nombre fini d'\'el\'ements non nuls $V$ dans le ``r\'eseau'' ${\rm K}_{\infty}^{\leq
19}$ tels que $\cbil_\infty^{{\rm F}_{\log\, 9}}(V,V) \leq \frac{8 \, {\rm log}\, 9}{\pi^{2}}$.  Il ne reste qu'\`a les \'enum\'erer, ce que l'on va faire en utilisant l'algorithme de Fincke et Pohst \cite{FP} impl\'ement\'e dans \texttt{PARI} (commande \texttt{qfminim}). Pour ne pas \`a avoir \`a justifier les erreurs d'arrondis sur l'algorithme sus-cit\'e, il est commode de r\'eutiliser l'approximation $A \in 10^{-4} \,{\rm M}_{10}(\Z)$ de ${\rm Gram}(19,\log 9)$ introduite dans la d\'emonstration du lemme \ref{lemmedefpos}. Soit ${\rm q}_A : \Z^{10} \rightarrow \frac{1}{10^4}\Z$ la forme quadratique d\'efinie positive $x \mapsto {}^{\rm t} x A x$. Si $V \in {\rm K}_\infty$ v\'erifie $\cbil_\infty^{{\rm F}_{\log\, 9}}(V,V) \leq \frac{8 \, {\rm log}\, 9}{\pi^{2}}$, alors ses coordonn\'ees $(x_1,x_3,\dots,x_{19})$ dans la base ${\rm I}_1,\dots,{\rm I}_{19}$ v\'erifient $$10^4\, {\rm q}_A(x_1,\dots,x_{19}) \leq 10^4\, (1-0.00323)^{-1}\,\cdot \frac{8 \, {\rm log}\, 9}{\pi^{2}} < 17868$$ d'apr\`es le lemme \ref{inegatrivquad}. L'algorithme de Fincke et Pohst assure qu'il existe exactement $24$ couples $\pm (x_i) \in \Z^{10}$ v\'erifiant cette in\'egalit\'e, et nous en donne la liste. On ne retient bien entendu parmi ceux-ci que les \'el\'ements $x=(x_i)$ appartenant \`a $\N^{10}$ (correspondant \`a des $V$ effectifs), il ne reste que $11$ tels \'el\'ements, list\'es dans l'\'enonc\'e.
\end{pf}

\begin{remarque}{\rm  Pr\'ecisons que $9$ n'est pas le plus petit entier $m\geq 2$ tel que ${\rm Gram}(19,{\rm log} \, m )$ soit d\'efinie positive (et encore moins le seul!). Par exemple  tout entier $5 \leq m \leq 100$ conviendrait. N\'eanmoins, ce choix ${\rm log}(9)$, obtenu par t\^atonnement, a l'avantage de minimiser la taille de la liste obtenue dans l'\'enonc\'e du lemme~\ref{borneVlemm}. }
\end{remarque}

\ps\ps

Supposons maintenant que $\pi \in \Pi_{\rm alg}(\PGL_n)$ est de poids motivique impair $w \leq 19$, et posons $V={\rm L}(\pi_\infty) \in {\rm K}_\infty^{\leq 19}$. Le corollaire \ref{inegexploo} (i)  montre que l'on a $\cbil_\infty^{{\rm F}_{\log\, 9}}(V,V) \leq \frac{8 \, {\rm log}\, 9}{\pi^{2}}$, de sorte que $V$ est dans la liste donn\'ee dans le lemme ci-dessus. On peut supposer $n=\dim V > 2$, car le cas $n=2$ a d\'ej\`a \'et\'e trait\'e plus haut.  On en d\'eduit $n=4$, $w=19$ et $V\,=\, {\rm I}_{19}+{\rm I}_v$ avec $5 \leq v \leq 13$ et $v$ impair. On veut montrer $v=7$ et l'unicit\'e de la repr\'esentation $\pi$. D\'emontrons d'abord que $\pi$ est unique si elle existe. 

\begin{lemme} \label{bornemlemm} Pour tout $v \in \{5, 7, 9, 11, 13\}$ on a ${\rm m}({\rm I}_{19}+{\rm I}_v)\leq 1$.
\end{lemme}

\begin{pf} On applique le corollaire~\ref{cortaibi}. Il suffit de voir que pour $V={\rm I}_{19}+{\rm I}_v$, avec $v$ comme dans l'\'enonc\'e, on a 
$$\cbil_\infty^{{\rm F}_{\log 9}}(V,V) > \frac{1}{2}\, \frac{8 \, {\rm log}\, 9}{\pi^{2}}.$$
Mais un calcul num\'erique montre qu'\`a $10^{-2}$ pr\`es on a $\frac{1}{2} \,\,\frac{8 \, {\rm log}\, 9}{\pi^{2}} \simeq 0.89$ alors que le terme de gauche de l'in\'egalit\'e ci-dessus vaut $1.65$, $1.47$, $1.42$, $1.49$ ou $1.70$ lorsque $v$ vaut respectivement $5, 7, 9, 11$ ou $13$.
\end{pf}

\begin{lemme}\label{m1autoduale} Soit $V \in {\rm K}_\infty$. Supposons qu'il existe une unique repr\'esentation $\pi \in \Pi_{\rm alg}$ telle que ${\rm L}(\pi_\infty)=V$. Alors $\pi$ est autoduale. Plus g\'en\'eralement, on a ${\rm m}(V) \equiv {\rm m}^\bot(V) \bmod 2$.
\end{lemme}

\begin{pf} En effet, si $\pi \in \Pi_{\rm alg}$ on a \'egalement $\pi^\vee \in \Pi_{\rm alg}$, ainsi que les \'egalit\'es ${\rm L}((\pi^\vee)_\infty)={\rm L}(\pi_\infty)^\ast = {\rm L}(\pi_\infty)$.
\end{pf}

Pour achever la d\'emonstration du th\'eor\`eme \ref{classpoids22} dans le cas $w$ impair $\leq 19$, il reste \`a d\'emontrer qu'il n'existe pas de repr\'esentation autoduale $\pi \in \Pi_{\rm alg}$ telle que ${\rm L}(\pi_\infty)={\rm I}_{19}+{\rm I}_v$ avec $v \in \{5,9,11,13\}$. Une premi\`ere mani\`ere de proc\'eder serait d'utiliser la proposition \ref{amfpgsp4} (iii) ainsi que la table \ref{dimsjknonnul} (cette table contient les renseignements requis car $19-v\geq 6$ dans tous les cas). Ce n'est toutefois pas n\'ecessaire, car nous allons voir que cela se d\'eduit du corollaire \ref{inexploo2} et de l'existence des repr\'esentations $1, \Delta_{11}$ et $\Delta_{15}$. Pour usage futur, explicitons le crit\`ere suivant. \ps\ps \ps
\ps

\begin{scholie}\label{scholieinexpl} Soient $V,V' \in {\rm K}_\infty$ et $\lambda$ un r\'eel $>0$. On pose {\scriptsize $${\rm t}(V,V', \lambda) \, = \,\sqrt{|(\Phi_{{\rm F}_{\lambda}}(0)-\cbil_\infty^{{\rm F}_{\lambda}}(V,V))(\Phi_{{\rm F}_{\lambda}}(0)-\cbil_\infty^{{\rm F}_{\lambda}}(V',V'))|}\,+\,\frac{\varepsilon(V \otimes V')-1}{4} \,\Phi_{{\rm F}_{\lambda}}(\frac{1}{2}) \,- \,\cbil_\infty^{{\rm F}_{\lambda}}(V,V')\,$$}\noindent Supposons que $V$ et $V'$ sont distincts, effectifs, et v\'erifient ${\rm m}^\bot(V)\geq 1$ et ${\rm m}^\bot(V')\geq 1$. Alors on a ${\rm t}(V,V',\lambda)\geq 0$. En particulier, si ${\rm m}(V)={\rm m}(V')=1$ on a ${\rm t}(V,V',\lambda)\geq 0$.
\end{scholie}

\begin{pf}  La premi\`ere assertion est une cons\'equence imm\'ediate du point (ii) du corollaire \ref{inexploo2}. La seconde d\'ecoule de la premi\`ere et du lemme \ref{m1autoduale}.
\end{pf}
Nous savons que l'on a ${\rm m}(V)=1$ pour $V=1$ ou $V={\rm L}((\Delta_w)_\infty)$ pour $w\leq 21$. Observons la table \ref{tableauelim1}. \begin{table}[h!]
\renewcommand{\arraystretch}{1.5}
\begin{center}
{\small{\begin{tabular}{|c||c|c|c|c|c|c|}\hline
\backslashbox{$v$}{$\pi'$} & $1$ & $\Delta_{11}$ & $\Delta_{15}$ & $\Delta_{17}$ & $\Delta_{19}$ & $\Delta_{21}$ \cr
\hline
$13$ & $0.141$ &$0.074$ & ${\color{blue}-0.006}$& $0.166$& $0.088$& $0.990$ \cr \hline
$11$ & $0.697$ &${\color{blue}-0.492}$ &$0.396$ &$0.498$ &$0.376$ & $1.251$  \cr \hline
$9$ & ${\color{blue}-0.094}$& ${\color{blue}-0.046}$& $0.636$& $0.689$& $0.536$& $1.388$ \cr \hline
$7$ & $0.308$ & $0.223$ & $0.762$ & $0.778$ & $0.597$ & $1.430$   \cr \hline
$5$ & ${\color{blue}-0.660}$ & $0.359$ & $0.771$ & $0.751$ & $0.546$ & $1.357$   \cr \hline
\end{tabular}}}
\end{center}
\caption{Valeur \`a $10^{-3}$ pr\`es de ${\rm t}({\rm I}_{19}+{\rm I}_v,{\rm L}(\pi'_\infty), \log 5)$.}
\label{tableauelim1}
\end{table}
\renewcommand{\arraystretch}{1} On constate que si $v \neq 7$, il existe toujours une repr\'esentation $\pi' \in \{1, \Delta_{11},\Delta_{15}\}$ telle que ${\rm t}({\rm I}_{19}+{\rm I}_v,{\rm L}(\pi'_\infty), \log 5)<0$. Par le scholie \ref{scholieinexpl}, cela d\'emontre ${\rm m}({\rm I}_{19}+{\rm I}_v) \neq 1$, et conclut donc la d\'emonstration. $\square$ \ps
\ps

\begin{remarque}{\rm Le param\`etre $\lambda=\log 5$ a \'et\'e choisi par t\^atonnement. Une variation de ce param\`etre montre que l'existence d'une repr\'esentation dans le cas $v=7$ para\^it assez miraculeuse de ce point de vue. La propri\'et\'e qui semble importante est que les \'ecarts successifs entre les $4$ poids, \`a savoir $\frac{19-v}{2}, v$ et $\frac{19-v}{2}$, sont presque \'egaux pour $v=7$. }
\end{remarque}

\ps \ps \medskip
\noindent {\sc Cas $w$ pair $\leq 20$}
\ps \ps 

On proc\`ede de mani\`ere strictement similaire au cas $w$ impair $\leq 19$, c'est pourquoi nous donnerons moins de d\'etails. 

\begin{lemme}\label{lemmedefpos20} La restriction de $\cbil_\infty^{{\rm F}_{\log 9}}$ 
 \`a ${\rm K}_\infty^{\leq 20}$ est d\'efinie positive. Soit $V \in {\rm K}_\infty^{\leq 20}$ effectif, de d\'eterminant $1$, et v\'erifiant $\cbil_\infty^{{\rm F}_{\log\, 9}}(V,V) \leq \Phi_{{\rm F}_{\log 9}} (0)$, alors :\ps\ps \begin{itemize}
\item[(i)] soit $V=1$, \ps\ps
\item[(ii)] soit $V={\rm I}_w+\epsilon_{\C/\R}$ avec $w=18$ ou $20$, \ps\ps
\item[(iii)] soit $V={\rm I}_w+{\rm I}_v$ avec soit $w=18$ et $8 \leq v \leq 10$, soit  $w=20$ et $4\leq v\leq 14$,\ps\ps
\item[(iv)] soit $V={\rm I}_{20}+{\rm I}_v+1$ avec $10 \leq v \leq 16$.\ps\ps
\end{itemize}
De plus, on a ${\rm m}(V) \leq 1$.
\end{lemme}

\begin{pf} Les deux premi\`eres assertions se d\'eduisent d'une \'etude de la matrice ${\rm Gram}(20,\log 9)$ analogue \`a celle conduisant aux d\'emonstrations des lemmes \ref{lemmedefpos} et \ref{borneVlemm}  : voir la feuille de calculs en annexe \cite{clcalc}.  Pour les $15$ \'el\'ements $V$ de l'\'enonc\'e, on v\'erifie que l'on a ${\rm m}(V) \leq 1.6$ \`a l'aide du corollaire \ref{cortaibi} appliqu\'e \`a $\lambda = \log 9$, d'o\`u la derni\`ere assertion.
\end{pf}

Il en r\'esulte que si $\pi \in \Pi_{\rm alg}$ est de poids motivique pair $\leq 20$, alors $V={\rm L}(\pi_\infty)$ est dans la liste ci-dessus, ${\rm m}(V) \leq 1$ et  $\pi$ est autoduale. On \'elimine les possibilit\'es $$V={\rm I}_{18} + {\rm I}_{10}, \hspace{10pt}{\rm I}_{20} + {\rm I}_{10}+1 \hspace{10pt} {\rm ou}  \hspace{10pt}{\rm I}_{20}+{\rm I}_{16}+1$$ \`a l'aide du crit\`ere $\cbil_\infty^{{\rm F}_{\lambda}}(V,V) \leq \Phi_{{\rm F}_{\lambda}}(0)$ en prenant respectivement $\lambda = \log\,10$, $\log\, 16$ et  $\log\, 16$. Mentionnons que si $V={\rm I}_{18} + {\rm I}_{10}$ et $\lambda = \log 10$, la  quantit\'e $\Phi_{{\rm F}_{\lambda}}(0) - \cbil_\infty^{{\rm F}_{\lambda}}(V,V)$ vaut $-0.00012$ \`a $10^{-5}$ pr\`es, ce qui passe assez juste! \ps \ps

 
 Il reste \`a montrer que $\pi$ n'existe pas si $V \neq 1$. Par chance, nous y parvenons \`a l'aide du scholie \ref{scholieinexpl}, appliqu\'e dans les cas particuliers indiqu\'es dans la table  \ref{tableauelim2}. Cela termine la d\'emonstration du th\'eor\`eme \ref{classpoids22} lorsque le poids motivique $w$ est pair $\leq 20$.

\begin{table}[!h]

\renewcommand{\arraystretch}{1.5}
\begin{center}{\footnotesize{\begin{tabular}{|c|c||c|c||c|c|} \hline
$V$ & $\pi'$ &$V$ & $\pi'$ & $V$ & $\pi'$  \cr \hline
${\rm I}_{18}+\epsilon_{\C/\R}$ & $1$ &  ${\rm I}_{20}+{\rm I}_{12}$ & $\Delta_{15}$ & ${\rm I}_{20}+{\rm I}_{4}$ & $1$ \cr \hline
${\rm I}_{20}+\epsilon_{\C/\R}$ & $1$ & ${\rm I}_{20}+{\rm I}_{10}$ & $1$ &  ${\rm I}_{20}+{\rm I}_{14}+1$ & $\Delta_{15}$ \cr \hline
${\rm I}_{18}+{\rm I}_8$ & $\Delta_{11}$ & ${\rm I}_{20}+{\rm I}_{8}$ & $1$ &  ${\rm I}_{20}+{\rm I}_{12}+1$ & $1$ \cr \hline
${\rm I}_{20}+{\rm I}_{14}$ & $\Delta_{15}$ & ${\rm I}_{20}+{\rm I}_{6}$ & $\Delta_{11}$ & & \cr \hline
\end{tabular}}}
\end{center}
\caption{Des couples $(V,\pi')$ v\'erifiant ${\rm t}(V,{\rm L}(\pi'_\infty),\log 5)<-0.05$.}
\label{tableauelim2}

\end{table}
\renewcommand{\arraystretch}{1}

\subsection{Interlude : un crit\`ere g\'eom\'etrique}\label{interludegeom}

\begin{lemme} \label{lemmegeoelim}Soient $E$ un espace euclidien, $m\geq 1$ un entier,  $x_0,x_1,\dots,x_m$ des \'el\'ements de $E$, et $C_0,C_1,\dots,C_m$ des nombres r\'eels. On suppose que pour tout $i=0,\dots,m$ on a l'in\'egalit\'e $ x_0 \cdot x_i \leq C_i$. Alors :\ps \ps
\begin{itemize}
\item[(i)] On a $C_0\geq 0$ et, pour tout $i=1,\dots,m$, $C_i+\sqrt{C_0 \, (x_i \cdot x_i)} \geq 0$. \ps\ps
\end{itemize}

\noindent Soient $G=( x_i \cdot x_j)_{1\leq i,j, \leq m} \in {\rm M}_m(\R)$ la matrice de Gram des vecteurs $x_1,\dots,x_m$ et $C$ le vecteur colonne $(C_i)_{1\leq i \leq m} \in \R^m$. On suppose $\DET(G) \neq 0$, i.e. que les vecteurs $x_1,\dots,x_m$ sont lin\'eairement ind\'ependants. Alors l'une des deux assertions suivantes est v\'erifi\'ee : \ps\ps \begin{itemize}
\item[(ii)] au moins une des coordonn\'ees du vecteur $G^{-1} C$ est $> 0$. \ps\ps
\item[(ii)'] on a l'in\'egalit\'e (entre nombres r\'eels) ${}^{\rm t}C \,G^{-1} \,C \,\leq \,C_0$.\ps\ps
\end{itemize}
\end{lemme}

\begin{pf} L'in\'egalit\'e $C_0\geq 0$ est \'evidente. De plus, pour $i=1,\dots,m$ l'in\'egalit\'e de Cauchy-Schwarz donne $$C_i \geq x _0\cdot x_i \geq - |x_0 \cdot x_i| \geq - \sqrt{ (x_0\cdot x_0) ( x_i\cdot x_i)} \geq - \sqrt{C_0\, \,\,x_i \cdot x_i}.$$ 
V\'erifions la seconde assertion. Soit $(x_i^\ast)_{1\leq i \leq m}$ la base duale de $(x_i)_{1\leq i \leq m}$ dans l'espace euclidien $F={\rm Vect}_\R(x_1,\dots,x_m)$. Soit $H=(h_{i,j}) \in {\rm M}_m(\R)$ la matrice d\'efinie par les \'egalit\'es $x_j^\ast=\sum_{i=1}^m h_{i,j} x_i$ pour $j=1,\dots,m$. Par d\'efinition de la base duale, $H$ est \'egalement la matrice de Gram $(x_i^\ast \cdot x_j^\ast)$ et on a la relation $H=G^{-1}$. Ainsi, les coefficients du vecteur $G^{-1} \,C = H \,C$ ne sont autres que les produits scalaires  $x_i^\ast \cdot v$, $i=1,\dots,m$, avec $$v \,:= \,\sum_{j=1}^m C_j \,x_j^\ast.$$  Supposons que l'assertion (ii) n'est pas satisfaite, i.e. $x_i^\ast \cdot v\leq 0$ pour tout $i=1,\dots,m$. Nous allons voir que l'on a $x_0 \cdot x_0 \geq v \cdot v$, ce qui est l'assertion (ii)'. Le vecteur $x_0$ s'\'ecrit de mani\`ere unique sous la forme $x_0=v-\sum_{i=1}^m t_i x_i^\ast +w$ avec $w \in F^\perp$ et $t_i \in \R$ pour $i=1,\dots,m$. Si $i=1,\dots,m$, la condition $x_0 \cdot x_i \leq C_i$ \'equivaut \`a $t_i \geq 0$. Posons $\Vert x \Vert^2=x \cdot x$ pour $x\in E$. On conclut en contemplant l'\'egalit\'e suivante :
$$\Vert x_0\Vert^2 \,= \,\Vert v \Vert^2 \,-\, 2 \,\sum_{i=1}^m \,t_i\, \,x_i^\ast\cdot v\, +\, \Vert w-\sum_{i=1}^m\, t_i \,x_i^\ast \Vert^2.$$
\end{pf}
\ps\ps
\begin{remarque}{\rm L'interpr\'etation g\'eom\'etrique du lemme ci-dessus est la suivante. Par hypoth\`eses, le point $x_0$ est dans l'intersection de la boule $B$ de centre $0$ et de rayon $\sqrt{C_0}$, et du ``poly\`edre'' $P$ intersection des $m$ demi-espaces $x \cdot x_i \leq C_i$ pour $i=1,\dots,m$. L'assertion (i) affirme que chacun de ces demi-espaces rencontre $B$, une condition n\'ecessaire \'evidente! La quantit\'e ${}^{\rm t}C \,G^{-1} \,C$ de l'\'enonc\'e est le carr\'e de la distance \`a l'origine de l'espace affine $\{ x \in V\,|\, \,\,x\cdot x_i = C_i, \, \,  \,i=1\dots m\}$, soit encore $v+F^\perp$ dans les notations de la d\'emonstration. Si (et seulement si) la condition (ii) n'est pas satisfaite, cette distance est \'egalement la distance de $P$ \`a l'origine, d'o\`u le r\'esultat. }
\end{remarque}\ps  \ps

Le corollaire \ref{inegexpl} montre que les param\`etres de Satake d'une repr\'esentation $\pi \in \Pi_{\rm alg}$ telle que ${\rm L}(\pi_\infty)=V$ sont soumis \`a un ensemble de contraintes rentrant dans le cadre du lemme pr\'ec\'edent. Pour mettre ceci en place, il sera commode de noter ${\rm Q} \subset \N$ le sous-ensemble des puissances des nombres premiers ; tout $q \in {\rm Q}$ s'\'ecrit donc de mani\`ere unique sous la forme $q = p^k$ avec $p$ premier et $k$ un entier $\geq 1$. Si $\pi \in \Pi(\PGL_n)$ et si $q \in {\rm Q}$ on lui associe le nombre complexe $$x_q(\pi) = {\rm tr} \, \, {\rm c}_{p}(\pi)^{k} \in \C,$$
o\`u l'on a \'ecrit $q=p^k$ avec $p$ premier et $k\geq 1$. Si $\pi$ est autoduale et dans $\Pi_{\rm alg}$ on a m\^eme $x_q(\pi) \in \R$ d'apr\`es la proposition VIII.\ref{propdesparamdual}. Fixons $\lambda$ un r\'eel $>0$, et posons $${\rm Q}_\lambda=\{q \in {\rm Q},\,\,\, q < e^\lambda\} \hspace{10pt}{\rm et}\hspace{10pt}{\rm E}_\lambda = \prod_{q \in {\rm Q}_\lambda} \C.$$
On munit le $\R$-espace vectoriel sous-jacent \`a ${\rm E}_\lambda$ d'une structure d'espace euclidien au moyen du produit scalaire
 $$(x_q) \cdot (y_q) = \sum_{q \in {\rm Q}_\lambda}\,\, {\rm F}_{\lambda}(\log q) \, \,\frac{\log p}{\sqrt{q}}\, \,\Re\,\overline{x}_q\, y_q,$$ o\`u dans cette somme $p$ d\'esigne le diviseur premier de $q$. Pour tout entier $n\geq 1$ et tout $\pi \in \Pi(\PGL_n)$ on dispose d'un vecteur $$x_\lambda(\pi):=(x_q(\pi))_q \in {\rm E}_\lambda.$$ Par d\'efinition de $\cbil_f^{{\rm F}_{\lambda}}$ (Proposition-D\'efinition~\ref{prelimformexplbf}), on a $$\cbil_f^{{\rm F}_{\lambda}}(\pi,\pi') = x_\lambda(\pi) \cdot x_\lambda(\pi')$$\par
\noindent pour tous $\pi,\pi' \in \Pi_{\rm alg}$. \ps

Notons $\Pi_{\rm alg}^\bot \subset \Pi_{\rm alg}$ le sous-ensemble des repr\'esentations autoduales. Pour ne pas multiplier les \'enonc\'es, et \'etant donn\'e les applications que nous avons en vue, nous restreindrons l'analyse ci-dessous aux \'el\'ements de $\Pi_{\rm alg}^\bot$. 
Fixons $\pi_0 \in \Pi_{\rm alg}^\bot$ et posons $V_0={\rm L}((\pi_0)_\infty)$. La proposition~\ref{inegexpl} entra\^ine le syst\`eme d'in\'egalit\'es suivantes :

\begin{equation}\label{methodedudisque}  \left\{ \begin{array}{ll}  x_\lambda(\pi_0) \cdot x_\lambda(\pi) \leq \,\Phi_{{\rm F}_{\lambda}}(0)\, -\, \cbil_\infty^{{\rm F}_{\lambda}}(V_0,V_0),\\ \\ \hspace{40pt}{\rm et}\hspace{10pt} \forall \pi \in \Pi_{\rm alg}^\bot-\{\pi_0\}, \\ \\
x_\lambda(\pi_0) \cdot x_\lambda(\pi) \,\leq  \,-\, \frac{1-\varepsilon(V_0 \otimes {\rm L}(\pi_\infty))}{4} \, \Phi_\lambda(\frac{1}{2})\,-\, \cbil_\infty^{{\rm F}_{\lambda}}(V_0,{\rm L}(\pi_\infty)).
\end{array}\right.
\end{equation}
\ps\ps
\noindent Nous sommes donc manifestement dans les hypoth\`eses du lemme \ref{lemmegeoelim}, il entra\^ine imm\'ediatement le scholie suivant. 

\begin{scholie}\label{scholiemethodedudisque} Soient $V_0 \in {\rm K}_\infty$, $\lambda$ un r\'eel $>0$, $m\geq 1$ un entier, et $\pi_1,\dots,\pi_m$ des \'el\'ements distincts de $\Pi_{\rm alg}^\bot$. On pose $C_0= \Phi_{{\rm F}_{\lambda}}(0)\, -\, \cbil_\infty^{{\rm F}_{\lambda}}(V_0,V_0)$,  
$$C_i\, =\,-\, \frac{1-\varepsilon(V_0 \otimes {\rm L}((\pi_i)_\infty))}{4} \, \Phi_\lambda(\frac{1}{2})\,-\, \cbil_\infty^{{\rm F}_{\lambda}}(V_0,{\rm L}((\pi_i)_\infty))$$
pour $i=1,\dots,m$, et $C=(C_i) \in \R^m$. On suppose que la matrice  de ${\rm M}_m(\R)$
$$G=(x_\lambda(\pi_i) \cdot x_\lambda(\pi_j))_{1\leq i,j \leq m}$$
est inversible, que les coordonn\'ees du vecteur\, $G^{-1} \,C$ \,sont toutes $<0$, et que le r\'eel $\,C_0\,-\,C^{t} \,G^{-1} \,C\,$ est $ <0$. \ps

Si $\pi \in \Pi_{\rm alg}^\bot$ v\'erifie ${\rm L}(\pi_\infty)=V_0$, alors il existe $1\leq i \leq m$ tel que $\pi = \pi_i$. En particulier, si ${\rm L}((\pi_i)_\infty) \neq V_0$ pour tout $i$ alors ${\rm m}^\bot(V_0)=0$.
\end{scholie}
\ps\ps
\par \noindent Dans la suite, nous allons appliquer ce crit\`ere aux \'el\'ements $\pi_i$ de l'ensemble
$$\mathcal{R}=\{1, \Delta_{11}, \Delta_{15},\Delta_{17},\Delta_{19},\Delta_{21},{\rm Sym}^2 \Delta_{11}\} \subset \Pi_{\rm alg}^\bot.$$
Les vecteurs $x_\lambda(\pi)$, pour $\pi \in \mathcal{R}$ et $\lambda$ raisonnable, sont consid\'er\'es comme connus. Par exemple on a $x_\lambda(1)=(1,1,1,\dots)$ et
$$x_{\lambda}(\Delta_{11}) \hspace{4pt} =\hspace{4pt} ( \tau(2)\,2^{-\frac{11}{2}},\hspace{4pt}  \tau(3)\,{3^{-\frac{11}{2}}}, \hspace{4pt} {(\tau(4)-2^{11})}\,{4^{-\frac{11}{2}}}, \hspace{4pt}\dots)$$
$$\hspace{-.9cm} \simeq \hspace{4pt} (-0.530, \hspace{4pt}0.599, \hspace{4pt}-1.719, \hspace{4pt}\dots).$$
\ps\ps
\noindent 

\noindent En guise d'application, d\'emontrons le r\'esultat suivant, dont nous aurons besoin par la suite.\ps\ps

\begin{lemme}\label{applcgeom} Supposons que $V \in {\rm K}_\infty$ appartienne \`a la liste des $8$ \'el\'ements 
$${\rm I}_{21}+{\rm I}_{17}+{\rm I}_{7},\, \,\,\,{\rm I}_{22}+{\rm I}_{4}, \, \,\,\,{\rm I}_{22}+{\rm I}_{12}, \, \, \,\,{\rm I}_{22}+{\rm I}_{16}+1, \, \,\,\, {\rm I}_{22}+{\rm I}_{12}+1, $$ $${\rm I}_{22}+{\rm I}_{16}+{\rm I}_{10}+\epsilon_{\C/\R},\,\,\,\, {\rm I}_{22}+{\rm I}_{20}+{\rm I}_{10}+\epsilon_{\C/\R},\,\,\, \,{\rm I}_{22}+{\rm I}_{20}+{\rm I}_{14}+{\rm I}_4.$$\ps\ps
\noindent Alors on a ${\rm m}^\bot(V)=0$.
\end{lemme}

\begin{pf} On applique le scholie \ref{scholiemethodedudisque} avec $m=2$ et en prenant pour $V_0$, $\{\pi_1,\pi_2\}$, $\lambda$ chacun des triplets indiqu\'es dans le tableau \ref{tableauelim3}. On v\'erifie au cas par cas que les hypoth\`eses du scholie sont bien satisfaites (voir la feuille de calculs \cite{clcalc}). 
 \end{pf}

\begin{table}
\renewcommand{\arraystretch}{1.5}
\begin{center}{\scriptsize {\begin{tabular}{|c|c|c||c|c|c|} \hline
$V_0$ & $\{\pi_1, \pi_2\}$ & $\lambda$ & $V_0$ & $\{\pi_1, \pi_2\}$ & $\lambda$  \cr \hline
${\rm I}_{21}+{\rm I}_{17}+{\rm I}_{7}$  & $\{\Delta_{15}, \, \Delta_{17}\}$ & $\log\, 14$ & 
${\rm I}_{22}+{\rm I}_{4}$ & $\{1, \,\Delta_{21}\}$ & $\log\, 7$ \cr \hline
${\rm I}_{22}+{\rm I}_{12}$ & $\{\Delta_{11},\,\Delta_{15}\}$ & $\log \, 5$ &
${\rm I}_{22}+{\rm I}_{16}+1$ & $\{1, \,\Delta_{17}\}$ & $\log \,8$ \cr \hline
${\rm I}_{22}+{\rm I}_{12}+1$ & $\{1, \,\Delta_{11}\}$ & $\log \,5$ &
${\rm I}_{22}+{\rm I}_{16}+{\rm I}_{10}+\epsilon_{\C/\R}$ & $\{\Delta_{11},\,\Delta_{15}\}$ & $\log \, 9$ \cr \hline
${\rm I}_{22}+{\rm I}_{20}+{\rm I}_{10}+\epsilon_{\C/\R}$ & $\{\Delta_{19},\, {\rm Sym}^2 \, \Delta_{11}\}$  & $\log \,38$ &
${\rm I}_{22}+{\rm I}_{20}+{\rm I}_{14}+{\rm I}_4$ & $\{1,\, \Delta_{21}\}$ & $\log \,40$\cr \hline
\end{tabular}}}
\end{center}
\caption{\small Des triplets $(V_0,\{\pi_1,\pi_2\},\lambda)$ v\'erifiant les hypoth\`eses du scholie~\ref{methodedudisque}.}
\label{tableauelim3}
\end{table}
\renewcommand{\arraystretch}{1}

\subsection{Fin de la d\'emonstration du th\'eor\`eme \ref{classpoids22} : cas des poids motiviques $21$ et $22$.}

\begin{lemme} \label{elimprinc21} Soit $V \in {\rm K}_\infty^{\rm \leq 21}$ tel que $V-{\rm I}_{21}$ est effectif, non nul, et tel que
$\cbil_\infty^{{\rm F}_{\lambda}}(V,V) \leq \Phi_\lambda(0)$ pour $\lambda = \log 28$. Alors $V$ est l'un des $26$ \'el\'ements suivants :  \ps\ps\begin{itemize}
\item[(i)] ${\rm I}_{21} + {\rm I}_{v}$ avec $17 \geq v \geq 3$,\ps \ps
\item[(ii)] ${\rm I}_{21} + {\rm I}_v + {\rm I}_u$ avec $19 \geq
v \geq 13$, $9 \geq u \geq 3$, et $(v,u) \neq (13,9)$, ou avec $(v,u)=(17,11)$,\ps \ps
\item[(iii)] ${\rm I}_{21} + {\rm I}_{19} + {\rm I}_{13} +
{\rm I}_v$ avec $5 \geq v \geq 3$. \ps \ps
\end{itemize}
De plus, on a ${\rm m}(V) \leq 1$. 
\end{lemme}

\begin{pf} La premi\`ere assertion se d\'eduit d'une \'etude de la matrice (d\'efinie positive!) ${\rm Gram}(21,\log 28)$ analogue \`a celle conduisant aux d\'emonstrations des lemmes \ref{lemmedefpos} et \ref{borneVlemm}  : voir la feuille de calculs en annexe. Pour les $26$ \'el\'ements $V$ de l'\'enonc\'e, on v\'erifie que l'on a ${\rm m}(V) < 1.8$ \`a l'aide du corollaire \ref{cortaibi} appliqu\'e \`a $\lambda = \log 28$, d'o\`u la derni\`ere assertion.
\end{pf}

\begin{lemme} \label{elimprinc21suite} Soit $V \in {\rm K}_\infty$ l'un des 26 \'el\'ements list\'es dans l'\'enonc\'e du lemme~\ref{elimprinc21}.
On suppose ${\rm t}(V,{\rm L}(\pi_\infty), \log 27)>0$ pour tout $\pi \in \mathcal{R}$ {\rm (voir le Scholie \ref{scholieinexpl} et le \S \ref{interludegeom})}. Alors $V$ est l'un des \'el\'ements 
${\rm I}_{21}+{\rm I}_5, \hspace{7pt} {\rm I}_{21}+{\rm I}_9, \hspace{7pt} {\rm I}_{21}+{\rm I}_{13}\hspace{7pt}{\rm ou}\hspace{7pt} {\rm I}_{21}+{\rm I}_{17}+{\rm I}_{7}$.
\end{lemme}

\begin{pf} C'est un simple calcul num\'erique pour lequel nous renvoyons \`a la fiche de calculs. 
\end{pf}

Pour d\'emontrer le cas $w=21$ du th\'eor\`eme \ref{classpoids22}, il ne reste donc qu'\`a montrer que l'on a ${\rm m}^\bot({\rm I}_{21}+{\rm I}_{17}+{\rm I}_{7})=0$. Mais cela a d\'ej\`a \'et\'e d\'emontr\'e dans le lemme \ref{applcgeom}. Cela termine la d\'emonstration du th\'eor\`eme~\ref{classpoids22} dans le cas $w=21$. $\square$ \ps \ps


\begin{lemme} \label{elimprinc22} Soit $V \in {\rm K}_\infty^{\rm \leq 22}$. On suppose $V-{\rm I}_{22}$ effectif et l'in\'egalit\'e
 $\cbil_\infty^{{\rm F}_{\lambda}}(V,V) \leq \Phi_\lambda(0)$ satisfaite pour $\lambda = \log 80$. \begin{itemize} \ps \ps 
 \item[(i)] On a ${\rm m}(V) \leq 1$, sauf si $V$ est parmi les \'el\'ements suivants 
 $${\rm I}_{22}+{\rm I}_{12}, \, \, \, \, \, {\rm I}_{22}+{\rm I}_{10},\, \, \, \, \, {\rm I}_{22}+{\rm I}_{8},$$
auquel cas on a seulement ${\rm m}(V) \leq 2$. \ps \ps

\item[(ii)] Supposons de plus $V\neq {\rm I}_{22}+\epsilon_{\C/\R}$ et l'in\'egalit\'e ${\rm t}(V,{\rm L}(\pi_\infty),\log 77)>0$ satisfaite pour tout $\pi \in \mathcal{R}$. Alors $V$ appartient \`a la liste des $8$ repr\'esentations suivantes : \end{itemize}
{\small
\begin{equation} \label{liste8el}\left\{ \begin{array}{c} {\rm I}_{22}+{\rm I}_{12}, \, \,\,\, {\rm I}_{22}+{\rm I}_{8} ,\, \,\,\,{\rm I}_{22}+{\rm I}_{4}, \, \, \,\,{\rm I}_{22}+{\rm I}_{16}+1, \, \,\,\, {\rm I}_{22}+{\rm I}_{12}+1, \\ \\ {\rm I}_{22}+{\rm I}_{16}+{\rm I}_{10}+\epsilon_{\C/\R},\,\,\,\, {\rm I}_{22}+{\rm I}_{20}+{\rm I}_{10}+\epsilon_{\C/\R},\,\,\, \,{\rm I}_{22}+{\rm I}_{20}+{\rm I}_{14}+{\rm I}_4.
\end{array} \right.
\end{equation}
}

\end{lemme}
\begin{pf} On v\'erifie tout d'abord que la matrice  $B={\rm Gram}(22,\log 80)$ est d\'efinie positive par la m\^eme m\'ethode que celle employ\'ee dans la d\'emonstration du lemme  \ref{lemmedefpos} : voir la feuille de calculs en annexe, dans laquelle on \'etudie la matrice sym\'etrique $10^6 \,A$ obtenue en arrondissant \`a l'entier le plus proche tous les coefficients de la matrice $10^6\, B$. \ps

On proc\`ede ensuite comme dans la d\'emonstration du lemme \ref{borneVlemm}. L'algorithme \texttt{qfminim} de \texttt{PARI} appliqu\'e \`a $10^6 A \in {\rm M}_{12}(\Z)$ renvoie un ensemble de $701$ couples $\pm V$ contenant tous les \'el\'ements $V \in {\rm K}_\infty^{\leq 22}$ v\'erifiant $\cbil_\infty^{{\rm F}_{\lambda}}(V,V) \leq \Phi_\lambda(0)$ pour $\lambda=\log 80$. Si l'on ne retient dans cet ensemble que le sous-ensemble, not\'e $\mathcal{L}$, constitu\'e des $V$ tels que $V-{\rm I}_{22}$ est effectif et v\'erifiant $\DET \,V \,=\,1$, il ne reste ``plus que'' $158$ possibilit\'es pour $V$, autrement dit nous avons $|\mathcal{L}|=158$. \ps

Nous d\'eterminons ensuite le sous-ensemble des $V \in \mathcal{L}$ v\'erifiant de plus $0 \,\leq \,2 \,\cbil_\infty^{{\rm F}_{\lambda}}(V,V) \leq \Phi_\lambda(0)$ pour $\lambda = {\rm log}\, 77$ : on constate qu'il ne reste que les trois \'el\'ements de la forme ${\rm I}_{22}+{\rm I}_v$ avec $v=12,10$ ou $8$. Comme dans chacun de ces cas on v\'erifie d'autre part l'in\'egalit\'e $3\,\cbil_\infty^{{\rm F}_{\lambda}}(V,V) \,>\, \Phi_\lambda(0)$ (toujours pour $\lambda = {\rm log}\, 77$), le corollaire \ref{cortaibi} d\'emontre la premi\`ere assertion. \ps
Pour d\'emontrer la seconde assertion, nous calculons simplement les $|\mathcal{L}| \cdot |\mathcal{R}| = 1106$ quantit\'es ${\rm t}(V,{\rm L}(\pi_\infty),{\rm log} 77)$, avec $V \in \mathcal{L}$ et $\pi \in \mathcal{R}$. Nous renvoyons \`a la feuille de calculs pour la justification des r\'esultats. 
\end{pf}

\begin{pf} (Fin de la d\'emonstration du th\'eor\`eme \ref{classpoids22}) Soit $\pi \in \Pi_{\rm alg}$ de poids motivique $22$. On pose $V={\rm L}(\pi_\infty)$. L'\'el\'ement $V-{\rm I}_{22}$ est effectif et l'on a $\cbil_\infty^{{\rm F}_{\lambda}}(V,V) \leq \Phi_{{\rm F}_{\lambda}}(0)$ pour tout $\lambda>0$ d'apr\`es le corollaire \ref{inegexploo} (i). \ps 

Supposons d'abord $\pi$ autoduale, en particulier ${\rm m}^\bot(V) \geq 1$. Le scholie \ref{scholieinexpl} et le lemme \ref{elimprinc22} (ii) montrent que soit $V={\rm I}_{22}+\epsilon_{\C/\R}$, soit $V$ est dans la liste \eqref{liste8el} ci-dessus. De plus, on constate l'in\'egalit\'e $${\rm t}({\rm I}_{22}+{\rm I}_8,{\rm I}_{11},\log 5)<0$$
de sorte que l'on a \'egalement $V \neq {\rm I}_{22}+{\rm I}_8$. Mais d'apr\`es le lemme~\ref{applcgeom}, pour les $7$ \'el\'ements $W$ restants de la liste  \eqref{liste8el} on a ${\rm m}^\bot(W) =0$. En conclusion, on a $V={\rm I}_{22}+\epsilon_{\C/\R}$. Mais le lemme \ref{elimprinc22} (i) entra\^ine ${\rm m}(V)\leq 1$, et donc ${\rm m}^\bot(V)=1$ et $\pi = {\rm Sym}^2 \Delta_{11}$. \ps

Supposons maintenant que $\pi$ n'est pas autoduale. D'apr\`es le lemme \ref{elimprinc22} (i) on a donc 
${\rm m}(V)=2$ et $V={\rm I}_{22}+{\rm I}_v$ avec $v \in \{8,10,12\}$. En particulier, les deux repr\'esentations $\varpi \in \Pi_{\rm alg}$ v\'erifiant ${\rm L}(\varpi_\infty)=V$ sont $\pi$ et $\pi^\vee$. \ps

Pour conclure, il suffit de d\'emontrer que l'\'egalit\'e ${\rm m}({\rm I}_{22}+{\rm I}_{v})=2$ pour $v=8$ (resp. $10$, $12$) entra\^ine ${\rm m}({\rm I}_{21}+{\rm I}_{u})=0$ pour $u=9$ (resp. $9$, $13$). Soient $(v,u)$ l'un des trois couples $(8,9)$, $(10,9)$ ou $(12,13)$,  $V={\rm I}_{22}+{\rm I}_v$ et $V'= {\rm I}_{21}+{\rm I}_{u}$. D'apr\`es le corollaire \ref{inexploo2} (i), il suffit de v\'erifier qu'il existe $\lambda>0$ tel que l'on ait l'in\'egalit\'e 
$$\,\sqrt{(\frac{\Phi_{{\rm F}_{\lambda}}(0)}{2}-\cbil_\infty^{{\rm F}_{\lambda}}(V,V))(\Phi_{{\rm F}_{\lambda}}(0)-\cbil_\infty^{{\rm F}_{\lambda}}(V',V'))} \,- \,\cbil_\infty^{{\rm F}_{\lambda}}(V,V')<0.$$ 
Mais on v\'erifie que pour $\lambda = \log 22$, et disons \`a $10^{-2}$ pr\`es, le terme de gauche vaut respectivement $-0.14$, $-0.03$ et $-0.23$ quand $(v,u)$ vaut $(8,9)$, $(10,9)$ et $(12,13)$. 
\end{pf}

\subsection{Compl\'ements}\label{complementfexpl}

Le premier compl\'ement concerne l'ordre d'annulation en $s=\frac{1}{2}$ de ${\rm L}(s,\pi)$ lorsque $\pi \in \Pi_{\rm alg}$ est de poids motivique $\leq 22$ (comparer avec la remarque VII.\ref{nonvanish}).

\begin{prop}\label{ann12djk} Soit $\pi \in \Pi_{\rm alg}^\bot$ de poids motivique $\leq 22$. On a 
\begin{equation} {\rm ord}_{s=\frac{1}{2}} \,{\rm L}(s,\pi)\,\, = \,\, \left\{ \begin{array}{ll} 0 & {\rm si}\hspace{10pt} \varepsilon(\pi)=1,\\ \\ 1 & {\rm sinon}. \end{array} \right.\end{equation}
De plus on a $\varepsilon(\pi)=-1$ si, et seulement, si $\pi = \Delta_{17}$ ou $\pi = \Delta_{21}$. 
\end{prop}


\begin{pf} Posons $r\,=\,{\rm ord}_{s=\frac{1}{2}}\, {\rm L}(s,\pi)$, ou ce qui revient au m\^eme $r\,=\,{\rm ord}_{s=\frac{1}{2}} \,\xi(s,\pi)$ d'apr\`es la remarque suivant le corollaire \ref{inegexpl}. L'\'equation fonctionnelle $\xi(s,\pi)\, =\, \varepsilon(\pi) \,\xi(1-s,\pi)$ montre que $r$ est pair si $\varepsilon(\pi)\,=\,1$, et impair sinon. Il suffit donc de montrer $r< 2$. On peut supposer $\pi \neq 1$ car on a $\varepsilon(1)=1$ et $\zeta(\,\frac{1}{2}\,) \neq 0$. Un argument similaire \`a celui donn\'e dans la d\'emonstration du corollaire \ref{inegexploo} (ii), appliqu\'e \`a $\pi$ et $1$, montre que sous l'hypoth\`ese $r\geq 2$ on a l'in\'egalit\'e
\begin{equation}\label{eqzero12} - \Phi_{{\rm F}_{\lambda}} (\, \frac{1}{2}\,) \,-\,{\rm J}_{{\rm F}_{\lambda}}(V) + \sqrt{(\Phi_{{\rm F}_{\lambda}}(0)-{\rm J}_{{\rm F}_{\lambda}}(1))(\Phi_{{\rm F}_{\lambda}}(0)-\cbil_\infty^{{\rm F}_{\lambda}}(V,V))} \geq 0\end{equation}
pour tout $\lambda >0$ : il suffit de minorer ${\rm ord}_{s=1/2} \, \xi(s,\pi)$ par $2$ plut\^ot que par ${\rm e}^\bot(\pi,1)$. Mais lorsque $V$ vaut respectivement $${\rm I}_{11},\hspace{4pt} {\rm I}_{15}, \hspace{4pt} {\rm I}_{17},\hspace{4pt}  {\rm I}_{19},\hspace{4pt} {\rm I}_{21},\hspace{4pt}  {\rm I}_{19}+{\rm I}_7, \hspace{4pt} {\rm I}_{21}+{\rm I}_{5}, \hspace{4pt} {\rm I}_{21}+{\rm I}_9, \hspace{4pt} {\rm I}_{21}+{\rm I}_{13}, \hspace{4pt} {\rm I}_{22}+\epsilon_{\C/\R},$$et disons $\lambda = \log 4$, on v\'erifie que le terme de gauche de \eqref{eqzero12} vaut respectivement $-1.07$, $-0.64$, $-0.49$, $-0.35$, $-0.23$, $-0.79$, $-0.86$, $-0.35$, $-0.05$, $-0.82$ \`a $10^{-2}$ pr\`es. La premi\`ere assertion de la proposition d\'ecoule donc du th\'eor\`eme \ref{classpoids22}. On observe que dans la liste $10$ des \'el\'ements $V$ ci-dessus, on a $\varepsilon(V)=1$ sauf pour $V={\rm I}_{17}$ et ${\rm I}_{21}$ ; cela termine la d\'emonstration.
\end{pf}

Une cons\'equence tr\`es simple, mais surprenante, du th\'eor\`eme \ref{classpoids22} est le fait qu'il n'existe qu'un nombre fini de $\pi \in \Pi_{\rm alg}$ telles que ${\rm w}(\pi) \leq 22$. Un retour sur notre d\'emonstration montre que cette assertion de finitude, qui est en fait notre point de d\'epart, est cons\'equence de la propri\'et\'e suivante : si l'on a  $w\leq 22$, il existe des r\'eels $\lambda>0$ tels que la restriction de la forme bilin\'eaire sym\'etrique $\cbil_\infty^{{\rm F}_{\lambda}}$ \`a ${\rm K}_\infty^{\leq w}$ est d\'efinie positive (Lemmes \ref{lemmedefpos}, \ref{lemmedefpos20}, \ref{elimprinc21} et \ref{elimprinc22}). Il se trouve que cette propri\'et\'e persiste pour $w=23$ (mais cesse d'\^etre vraie pour $w>23$ !).  

\begin{prop}\label{finitude23} Il n'existe qu'un nombre fini de repr\'esentations dans $\Pi_{\rm alg}$ de poids motivique $\leq 23$.
\end{prop}

\begin{pf} Un simple calcul montre en effet que ${\rm Gram}(23, 9.74)$ est d\'efinie positive. 
\end{pf}

Nous reportons \`a un travail ult\'erieur l'\'etude d\'etaill\'ee des repr\'esentations dans $ \Pi_{\rm alg}$ dont le poids motivique est $\geq 23$, qui nous entra\^inerait ici trop loin de nos pr\'eoccupations. Mentionnons toutefois deux travaux en relation avec ces questions. Dans \cite{chrenard2}, les auteurs d\'emontrent en supposant la conjecture VIII.\ref{conjaj2} une formule explicite et calculable pour ${\rm m}^\bot(V)$ lorsque : \ps\ps
-- $V$ est {\it sans multiplicit\'e}, c'est-\`a-dire que ses coefficients dans la base $1, \epsilon_{\C/\R}, \{{\rm I}_w, w>0\}$ sont tous dans $\{0,1\}$. \ps\ps
-- $\dim V \leq 8$, les r\'esultats n'\'etant que partiels pour $ \dim V = 7$.\ps\ps
\noindent Dans un remarquable tour de force \cite{taibisiegel}, Ta\"ibi a ensuite red\'emontr\'e ces formules en supposant uniquement la conjecture VIII.\ref{conjaj}, et surtout, il les a \'etendues au cas plus g\'en\'eral $\dim V \leq 14$. Ses r\'esultats sont m\^eme ind\'ependants de toute conjecture si $V$ a ses poids ``assez \'ecart\'es''.\footnote{Pr\'ecisons cette notion. Supposons $V$ effectif, disons de dimension $n= \dim V$, et notons $\lambda_1,\dots,\lambda_n$ le multi-ensemble de $n$ nombres complexes associ\'es \`a $V_{|\C^\ast}$ comme dans l'assertion de compatibilit\'e de la param\'etrisation de Langlands au caract\`ere infinit\'esimal \S VIII.\ref{apparitionWR} (iii). On peut supposer $\lambda_i \in \frac{1}{2}\Z$ pour tout $i$ et $\lambda_i-\lambda_j \in \Z$ pour tout $i,j$. Si les $\lambda_i$ sont entiers, on dit que $V$ a ses poids assez \'ecart\'es si l'on a l'in\'egalit\'e $|\lambda_i- \lambda_j| \neq 1$ pour tout $1 \leq i,j \leq n$. Si les $\lambda_i$ sont dans $\frac{1}{2}\Z-\Z$, on demande $\lambda_{n/2}>1$ ainsi que pour tout entier $3 \leq  i \leq n/2$ l'in\'egalit\'e $\lambda_{i-2} > \lambda_i + 2$. } \ps\ps

Revenons au cas du poids motivique $23$ et admettons la conjecture \ref{conjaj}. La th\'eorie des formes modulaires montre bien entendu l'\'egalit\'e ${\rm m}({\rm I}_{23})=2$.
Consid\'erons ensuite $V \in {\rm K}_\infty$ sans multiplicit\'e tel que $V-{\rm I}_{23}$ est effectif et v\'erifie ${\rm m}^\bot(V) \geq 1$. La th\'eorie des formes de Siegel de genre $2$, et la formule de Tsushima (voir la remarque \ref{formtsushimataibi}), montrent que si $\dim V =2$ on a ${\rm m}(V)=1$ et $V-{\rm I}_{23} \in \{I_{7},\, \, \, \, {\rm I}_9, \, \, \, {\rm I}_{13}\}$. Les calculs de \cite[Cor. I.1.5]{chrenard2}, bas\'es sur l'\'etude des invariants polynomiaux du groupe orthogonal du r\'eseau ${\rm E}_7$, montrent de plus que si $6 \leq \dim V \leq 8$ alors on a $\dim V = 6$, ${\rm m}^\bot(V)=1$ et $V-{\rm I}_{23}$ parcourt la liste suivante : $${\rm I}_{13}+{\rm I}_5,\hspace{4pt} {\rm I}_{15}+{\rm I}_3,\hspace{4pt} {\rm I}_{15}+{\rm I}_7,\hspace{4pt}{\rm I}_{17}+{\rm I}_5, \hspace{4pt}{\rm I}_{17}+{\rm I}_9, \hspace{4pt} {\rm I}_{19}+{\rm I}_3, \hspace{4pt}{\rm I}_{19}+{\rm I}_{11}.$$ 
Enfin, les r\'esultats sus-cit\'es de Ta\"ibi montrent que si $8<\dim V \leq 14$ alors $\dim V=10$, ${\rm m}^\bot(V)=1$ et $V={\rm I}_{23}+{\rm I}_{21}+{\rm I}_{17}+{\rm I}_{11}+{\rm I}_3$. Cela fait en tout $13$ repr\'esentations dans $\Pi_{\rm alg}$ de poids motivique ${\rm 23}$, la question restante est de savoir s'il y en a d'autres. Mentionnons \'egalement que nous ne connaissons que trois $\pi \in \Pi_{\rm alg}$ de poids motivique $24$, de dimensions respectives $7$, $8$ et $8$. Celui de dimension $7$ v\'erifie ${\rm L}(\pi_\infty)=\epsilon_{\C/\R}\oplus \bigoplus_{i=1}^3 {\rm I}_{8i}$ et il est reli\'e \`a la trialit\'e et au groupe ${\rm G}_2$ \cite[Cor. I.1.10]{chrenard2}.


\begin{remarque}\label{formtsushimataibi} {\rm Un cas particulier tr\`es simple des formules de Ta\"ibi est que la formule donn\'ee par Tsushima \cite{tsushima} pour calculer $\dim {\rm S}_{j,k}$ quand $k\geq 5$ est plus g\'en\'eralement valable pour $k\geq 3$ (\S \ref{prelimsp4}), \`a l'exception du cas $(j,k)=(0,3)$. Une inspection de ses valeurs montre notamment l'annulation $\dim {\rm S}_{j,k}=0$ quand $k=3,4$ et $j+2k-3 \leq 21$. Une seconde d\'emonstration de cette annulation est fournie par le th\'eor\`eme \ref{classpoids22} et la proposition \ref{amfpgsp4}.}
\end{remarque}

\section{D\'emonstration du th\'eor\`eme \ref{thmintro24}}

\subsection{Une nouvelle d\'emonstration du th\'eor\`eme~\ref{thmintro16}}

Commen\c{c}ons par donner une nouvelle d\'emonstration du th\'eor\`eme \ref{thmintro16}, en utilisant une m\'ethode qui se g\'en\'eralisera \`a la dimension $24$. On rappelle que l'on a $|{\rm X}_{16}|=2$ et que l'op\'erateur de Hecke ${\rm T}_2$ admet des valeurs propres distinctes sur $\Z[{\rm X}_{16}]$ (Corollaire II.3.6, \S III.3.1). En particulier, il existe exactement deux repr\'esentations $\pi \in \Pi_{\rm disc}({\rm O}_{16})$ telles que $\pi_\infty$ est la repr\'esentation triviale de ${\rm O}_{16}(\R)$. Nous avons d\'ej\`a expliqu\'e au \S \ref{casdim16} que la proposition suivante entra\^ine le th\'eor\`eme \ref{thmintro16} ; nous avons d'ailleurs d\'ej\`a donn\'e une premi\`ere d\'emonstration de cette proposition au corollaire VII.\ref{corpij} (ii).\ps\ps

\begin{propv} \label{encoredim16} Les param\`etres standards $\psi(\pi,{\rm St})$ des deux repr\'esentations $\pi \in \Pi_{\rm disc}({\rm O}_{16})$ telles que $\pi_\infty$ est la repr\'esentation triviale sont $\,[15]\,\oplus \,[1]$ et $\Delta_{11}[4] \oplus [7] \oplus [1]$.
\end{propv}

\begin{pf} D'apr\`es l'exemple donn\'e \`a la fin du \S VI.\ref{exconjal}, nous savons que la repr\'esentation triviale $1 \in \Pi_{\rm disc}({\rm O}_{16})$ satisfait $$\psi(1,{\rm St})=[15]\oplus [1].$$ Il s'agit donc de montrer que si $\pi$ d\'esigne la repr\'esentation non triviale de $\Pi_{\rm disc}({\rm O}_{16})$ telle que $\pi_\infty$ est triviale, alors $\psi(\pi,{\rm St})=\Delta_{11}[4] \oplus [7] \oplus [1]$. Ajoutons que l'on a $\psi(\pi,{\rm St}) \neq \psi(1,{\rm St})$. Il y a une raison g\'en\'erale \`a cela, mais une mani\`ere de le voir ici est d'utiliser que les deux valeurs propres de ${\rm T}_2$ agissant sur $\Z[{\rm X}_{16}]$, \`a savoir \,\, $2^7 \,{\rm tr}\,\, \psi(\pi,{\rm St})_2$\,\, et\,\, $2^7\,{\rm tr}\, \,\psi(1,{\rm St})_2$\,\, (Formule VI.\eqref{formulesattp}), sont distinctes. \ps

D'apr\`es la proposition V.\ref{Fnonnulle}, la repr\'esentation $\pi$ admet un $\vartheta$-correspondant en genre $1 \leq g\leq 4$. \'Etant donn\'e l'in\'egalit\'e $16 > 2\,g$, le th\'eor\`eme VIII.\ref{arthurst} d'Arthur et le corollaire VII.\ref{corparamrallis} montrent que le couple $(\pi,{\rm St})$ v\'erifie la conjecture d'Arthur-Langlands. Autrement dit, il existe un entier $k\geq 1$, et pour $i=1,\dots,k$ des repr\'esentations $\pi_i \in \Pi_{\rm cusp}(\PGL_{n_i})$ et des entiers $d_i\geq 1$, tels que $$\psi(\pi,{\rm St})\,=\,\oplus_{i=1}^k \pi_i[d_i].$$
L'hypoth\`ese sur $\pi_\infty$ impose que les valeurs propres de ${\rm St}\, {\rm c}_\infty(\pi)$ sont les $14$ entiers $\pm 7, \pm 6,\dots, \pm 1$, ainsi que $0$ avec multiplicit\'e $2$. On en d\'eduit que les $\pi_i$ sont alg\'ebriques (Proposition VIII.\ref{toutestalg}), de poids motivique v\'erifiant $\frac{{\rm w}(\pi_i)+d_i-1}{2} \leq 7$. En particulier, on a ${\rm w}(\pi_i) \leq 14$ pour tout $i$.\ps
D'apr\`es le th\'eor\`eme \ref{classpoids22}, pour tout $i=1,\dots,k$ on a $\pi_i=1$ ou $\pi_i = \Delta_{11}$. De plus, les poids de $\Delta_{11}$ sont $\pm \frac{11}{2}$. En consid\'erant la valeur propre $\pm 7$, qui ne peut qu'  ``appartenir'' \`a un constituant de la forme $\Delta_{11}[4]$ ou $[15]$, on constate que les seules possibilit\'es pour $\psi(\pi,{\rm St})$ sont les deux de l'\'enonc\'e. On peut aussi invoquer \`a ce stade l'\'egalit\'e $\Psi_{16}\,=\,\{\,\,[15]\oplus[1],\,\, \Delta_{11}[4]\oplus[7]\oplus [1]\,\,\}$ qui a \'et\'e v\'erifi\'ee au cours de la d\'emonstration de la proposition \ref{propliste24}. Cela termine la d\'emonstration (et red\'emontre m\^eme la conjecture de Witt $g=4$ !). \end{pf}

\subsection{D\'emonstration du th\'eor\`eme E} \label{demothmE}

\begin{thmv}\label{corretourtab24} Les \'el\'ements $\psi \in \mathcal{X}_{\rm AL}({\rm SL}_{24})$ tels que $\psi_\infty$ admette pour valeurs propres simples les entiers $\pm 11, \pm 10, \dots,\pm 1$, ainsi que $0$ pour valeur propre double, sont exactement les $24$ param\`etres de la table \ref{table24}. 
\end{thmv}

\begin{pf} En effet, soit $\psi \in \mathcal{X}_{\rm AL}({\rm SL}_{24})$ tel que $\psi_\infty$ v\'erifie la propri\'et\'e de l'\'enonc\'e. \'Ecrivons $\psi = \oplus_{i=1}^k \pi_i[d_i]$ avec $\pi_i \in \Pi_{\rm alg}(\PGL_{n_i})$ pour tout entier $1 \leq i \leq k$. Soit $i$ un tel entier. La condition sur $\psi_\infty$ entra\^ine $\frac{{\rm w}(\pi_i) + d_i-1}{2} \leq 11$, puis l'in\'egalit\'e ${\rm w}(\pi_i) \leq 22$. D'apr\`es le th\'eor\`eme \ref{classpoids22}, pour tout entier $i$ la repr\'esentation $\pi_i$ est donc dans l'ensemble $\Pi$ introduit 
juste avant l'\'enonc\'e de la proposition \ref{propliste24}. On conclut par cette proposition. 
\end{pf}

\begin{pf} (D\'emonstration du th\'eor\`eme E) Reprenons les notations du \S \ref{enoncesconjnv} : on dispose notamment de $24$ \'el\'ements $\psi_i \in \mathcal{X}({\rm SO}_{24})$, $i=1,\dots,24$, distincts d'apr\`es Nebe et Venkov, et nous devons d\'emontrer que ce sont les \'el\'ements de la table \ref{table24}. Comme cette table a \'egalement exactement $24$ \'el\'ements, et \'etant donn\'e le th\'eor\`eme \ref{corretourtab24}, il ne reste qu'\`a d\'emontrer que l'on a $\psi_i \in \mathcal{X}_{\rm AL}({\rm SL}_{24})$ pour tout $i$. \ps

Nous savons d\'ej\`a, suivant Ikeda et B\"ocherer, que l'on a $\psi_{24}= \Delta_{11}[12]$ (Corollaire IV.\ref{cor24ikedabocherer}, voir aussi le d\'ebut de la d\'emonstration de la proposition \ref{leechmoins2}). Soit $1 \leq i \leq 23$ et soit $\psi_i' \in \mathcal{X}({\rm SL}_{2g_i+1})$ le param\`etre standard du
$\vartheta$-correspondant dans $\Pi_{\rm cusp}({\rm Sp}_{2g_i})$ de
$\pi_i$ (les notations $\pi_i$ et $g_i$ sont rappel\'ees au \S \ref{enoncesconjnv}). D'apr\`es le lemme \ref{lemmegenusfaible} on a $g_i \leq 11$, de sorte que d'apr\`es Rallis (Corollaire VII.\ref{corparamrallis}) on a l'identit\'e
$$\psi_i = \psi'_i \oplus [23-2 g_i].$$
D'apr\`es le th\'eor\`eme VIII.\ref{arthurst} d'Arthur, on a $\psi'_i \in
\mathcal{X}_{\rm AL}({\rm SL}_{2 g_i+1})$. Cela entra\^ine $\psi_i \in \mathcal{X}_{\rm AL}({\rm SL}_{24})$ et conclut la d\'emonstration. \end{pf}

\begin{remarque}
{\rm
\label{remfinthmE}
Dans la d\'emonstration ci-dessus, il a \'et\'e commode de traiter \`a part le cas du param\`etre $\Delta_{11}[12]$, ce qui a \'et\'e possible gr\^ace \`a \cite{ikeda1} et \cite{bfw}. Un examen de la d\'emonstration que nous donnerons du th\'eor\`eme \ref{mainthmsiegel12} montrerait en fait que l'on peut se passer de ces deux r\'ef\'erences, et utiliser simplement le r\'esultat d'Erokhin \cite{erokhin}.
}
\end{remarque}

\section{Formes de Siegel de poids $\leq 12$} \label{fsiegelpoids12}

Le but de cette partie est d'\'etudier l'espace
${\rm S}_k({\rm Sp}_{2g}(\Z))$, avec $g\geq 1$ arbitraire et $k \leq 12$. Il s'agira en premier lieu de d\'eterminer sa dimension, puis de d\'ecrire, pour toute forme $F \in {\rm S}_k({\rm Sp}_{2g}(\Z))$ propre pour l'action de ${\rm H}({\rm Sp}_{2g})$, le param\`etre standard $\psi(\pi_F,{\rm St})$ de la repr\'esentation $\pi_F \in \Pi_{\rm cusp}({\rm Sp}_{2g})$ engendr\'ee par $F$. Pour faire court, nous dirons simplement que $\psi(\pi_F,{\rm St})$ est {\it le param\`etre standard de la forme propre $F$}. \ps \ps

La question de d\'eterminer la dimension de ${\rm S}_k({\rm Sp}_{2g}(\Z))$ a fait l'objet de travaux de nombreux auteurs, dans divers cas particuliers. Nous renvoyons par exemple aux articles de Poor et Yuen \cite{py, py2}, et \`a celui de Nebe et Venkov \cite{nebevenkov}, pour une discussion des dimensions qui \'etaient connues avant ce travail et l'article r\'ecent \cite{taibisiegel}. Rappelons que le cas $g=1$ est classique, et que dans les cas $g=2$ et $3$ une formule pour $\dim {\rm S}_k({\rm Sp}_{2g}(\Z))$ valable pour tout $k>g$ a \'et\'e respectivement d\'emontr\'ee par Igusa \cite{igusagenus2} et Tsuyumine \cite{tsuyumine}. La situation pour $g>3$ est longtemps rest\'ee tr\`es partielle, au sens o\`u  $\dim {\rm S}_k({\rm Sp}_{2g}(\Z))$ n'avait pu \^etre d\'etermin\'e que pour un petit nombre (fini!) de couples $(g,k)$ avec $g>3$, $k\geq 0$ et $g k \equiv 0 \bmod 2$. Elle a substantiellement \'evolu\'ee r\'ecemment avec l'algorithme de Ta\"ibi \cite{taibisiegel}, qui a aboutit \`a une formule concr\`ete pour $\dim {\rm S}_k({\rm Sp}_{2g}(\Z))$ valable pour tout $k>g$ et tout $g \leq 7$ ; cette formule est toutefois encore conditionnelle \`a la conjecture VIII.\ref{conjaj}, du moins au moment o\`u nous \'ecrivons ces lignes ! \ps\ps

La m\'ethode que nous allons utiliser est largement ind\'ependante de celles des auteurs sus-cit\'es. En particulier, elle n'utilise aucun des calculs mentionn\'es ci-dessus en genre $g>2$, et donne donc des nouvelles d\'emonstrations des cas pr\'ec\'edemment connus en ces genres. Dans l'esprit, elle est proche de la d\'emonstration par Duke et Imamo\u{g}lu \cite{dukei1} de l'annulation ${\rm S}_k({\rm Sp}_{2g}(\Z))=0$ pour tout $g\geq 1$ et tout $k\leq 6$. En effet, notre d\'emonstration et la leur ont de commun l'utilisation des ``formules explicites'' (au sens du \S \ref{mrw}). Duke et Imamo\u{g}lu les appliquent \`a la fonction ${\rm L}$ standard d'une forme de Siegel propre, en se basant sur les travaux de B\"ocherer et Mizumoto rappel\'es au \S \ref{complboc}, alors que nous les avons appliqu\'ees aux fonctions ${\rm L}$ de paires de repr\'esentations automorphes cuspidales de groupes lin\'eaires (Jacquet, Piatetski-Shapiro, Shalika). Notre approche est bien s\^ur permise par le th\'eor\`eme d'Arthur VIII.\ref{arthurst}. Sous cet angle, l'ingr\'edient crucial des d\'emonstrations qui vont suivre devient le th\'eor\`eme \ref{classpoids22}. \ps\ps

Rappelons enfin l'annulation tr\`es utile ${\rm
S}_k({\rm Sp}_{2g}(\Z))=0$ si $k < g/2$, d'apr\`es la th\'eorie des formes
singuli\`eres de Freitag, Resnikoff et Weissauer \cite[\S 14]{vandergeer}. En particulier, si $k \leq 12$ nous pouvons supposer $g \leq 24$, de sorte qu'il n'y a en r\'ealit\'e qu'un nombre fini de couples $(k,g)$ \`a consid\'erer.

\subsection{Formes de poids $12$ et d\'emonstration du th\'eor\`eme \ref{pbeichler24} de l'introduction} Rappelons que l'on dispose pour tout $g \geq 1$ d'une application lin\'eaire $\vartheta_g : \C[{\rm X}_{24}] \longrightarrow {\rm
M}_{12}({\rm Sp}_{2g}(\Z))$, ainsi que de $\vartheta_0 : \C[{\rm X}_{24}] \rightarrow \C$ (\S \ref{thetasiegel}).\ps\ps

\begin{thmv} \label{mainthmsiegel12}\begin{itemize} \item[(i)] La dimension des espaces ${\rm S}_{12}({\rm Sp}_{2g}(\Z))$ avec
$g\leq 12$ est donn\'ee par la table : \ps \medskip

\renewcommand{\arraystretch}{1.5}
\begin{center}
\begin{tabular}{|c||c|c|c|c|c|c|c|c|c|c|c|c|}
\hline $g$ & $1$ & $2$ & $3$ & $4$ & $5$ & $6$ & $7$ & $8$ & $9$ & $10$ & $11$ & $12$ \cr
\hline $\dim\, {\rm S}_{12}({\rm Sp}_{2g}(\Z))$ & $1$ & $1$ & $1$ & $2$ & $2$ &
$3$ & $3$ & $4$ & $2$ & $2$ & $1$ & $1$ \cr
\hline
\end{tabular}
\end{center}
\renewcommand{\arraystretch}{1}
\ps  \medskip \noindent
En particulier, $\bigoplus_{1\leq g\leq 12}\, {\rm S}_{12}({\rm
Sp}_{2g}(\Z))$ est de dimension $23$. \ps\ps

\item[(ii)] Pour tout $1 \leq g \leq 12$, 
l'application $\vartheta_g : \C[{\rm X}_{24}] \rightarrow {\rm M}_{12}({\rm Sp}_{2g}(\Z))$ 
induit un isomorphisme $ {\rm Ker}\, \vartheta_{g-1}/ {\rm Ker}\,
\vartheta_g \isomo {\rm S}_{12}({\rm Sp}_{2g}(\Z))$.\ps\ps

\item[(iii)] Il existe exactement $23$ repr\'esentations dans les $\Pi_{\rm cusp}({\rm
Sp}_{2g}(\Z))$, avec $1 \leq g\leq 12$, qui sont engendr\'ees par une forme de Siegel propre de poids
$12$. Elles ont pour param\`etres standards ceux de la table \ref{tablek=12}. \ps\ps

\end{itemize}

\end{thmv}

\begin{pf} Nous allons commencer par \'etablir un lemme pr\'eliminaire expliquant le contenu de la table \ref{tablek=12}. Rappelons l'ensemble 
\begin{equation}\label{ensgrandpi} \Pi=\{{\rm Sym}^2\Delta_{11},\Delta_{21,13},\Delta_{21,9},\Delta_{21,5},\Delta_{21}, \Delta_{19,7},\Delta_{19}, \Delta_{17},\Delta_{15},\Delta_{11},1\},\end{equation} 
introduit avant l'\'enonc\'e de la proposition \ref{propliste24}. Consid\'erons, pour $1 \leq g \leq 12$, l'ensemble $\Phi_g$ de tous les \'el\'ements $\phi \in \mathcal{X}({\rm SL}_{2g+1})$ tels que : \begin{itemize} \ps\ps
\item[(a)] $\phi_\infty$ admet pour valeurs propres les $2g+1$ entiers $0$ et $\pm (12-j)$ avec $j=1,\dots,g$. \ps\ps
\item[(b)] il existe $r\geq 1$, des entiers $d_1,\dots,d_r \geq 1$, et des repr\'esentations $\pi_1,\dots,\pi_r \in \Pi$, tels que $\phi=\oplus_{i=1}^r \pi_i[d_i]$. \ps\ps
\end{itemize}
\noindent 

\noindent Si $1 \leq g \leq 12$ et $\phi \in \Phi_g$, nous dirons que $\phi$ {\it satisfait la condition} ${\rm (C)}$ si dans sa d\'ecomposition (b) ci-dessus, il n'existe pas d'entier $1\leq i \leq r$ v\'erifiant  $\pi_i=1$ et $d_i>1$, et s'il existe au plus un entier $i$ tel que $\pi_i=1$. La condition ${\rm (C)}$ est toujours satisfaite si $g<11$, car dans ce cas $1$ n'est pas valeur propre de $\phi_\infty$ et $0$ en est valeur propre simple. 

\begin{lemme}\label{caractablek=12} L'ensemble des \'el\'ements de $\coprod_{1 \leq g \leq 12} \Phi_g$ satisfaisant la condition ${\rm (C)}$ est exactement l'ensemble des param\`etres r\'eunis dans la table \ref{tablek=12}.
\end{lemme}

\begin{pf} C'est un exercice de combinatoire du m\^eme acabit que celui effectu\'e dans la d\'emonstration de la proposition \ref{propliste24}. On peut le d\'eduire de cette proposition de la mani\`ere suivante. \ps Soit $\phi \in \Phi_g$ avec $g \leq 11$. D'apr\`es la proposition \ref{propliste24},  $\phi \oplus [23-2g]$ est un \'el\'ement de la table \ref{table24}. La propri\'et\'e $\phi \in {\mathcal X}_{\rm AL}({\rm SL}_{2g+1})$, la condition (C), et le th\'eor\`eme de Jacquet-Shalika (Proposition VI.\ref{jacquetshalika}), d\'eterminent alors uniquement $\phi$. C'est ainsi que nous avons en fait d\'efini les $22$ \'el\'ements de la table \ref{tablek=12} correspondants aux genres $g<12$. \ps
Il ne reste qu'\`a voir que le seul \'el\'ement $\phi \in \Phi_{12}$ satisfaisant {\rm (C)} est $\Delta_{11}[12]\oplus [1]$. Pour cela, \'ecrivons $\phi=\oplus_{i=1}^r \pi_i[d_i]$ avec $\pi_i \in \Pi$ pour tout $i$. Par la condition {\rm (C)}, il existe au plus un entier $i$ tel que $\pi_i=1$. La valeur propre $11$ de $\phi_\infty$ \'etant simple, il existe \'egalement au plus un entier $i$ tel que $\pi_i = {\rm Sym}^2 \Delta_{11}$. Comme $0$ est valeur propre triple de $\phi_\infty$, le lemme \ref{lemmecombigrandpi} (i)  montre qu'il existe un entier $i$ tel que $\pi_i = \Delta_{11}$ et $d_i=12$. La seule possibilit\'e est alors $\phi = \Delta_{11}[12]\oplus [1]$. 
\end{pf}

Ce lemme pr\'eliminaire \'etant \'etabli, consid\'erons $F \in {\rm S}_{12}({\rm Sp}_{2g}(\Z))$ une forme propre pour ${\rm
H}({\rm Sp}_{2g})$ avec $1 \leq g\leq 12$ et $\psi \in \mathcal{X}({\rm SL}_{2g+1})$ son param\`etre standard.
D'apr\`es le corollaire VI.\ref{infcarsiegel}, $\psi_\infty$ admet pour valeurs propres
les $2g+1$ entiers $\pm 11, \pm 10, \dots, \pm (12-g)$ et $0$. Toutes ces
valeurs propres sont simples, hormis la valeur propre $0$ pour $g=12$ qui
est de multiplicit\'e $3$. D'apr\`es le th\'eor\`eme \ref{arthurst} d'Arthur, on peut \'ecrire 
\begin{equation} \label{decarthk12} \psi= \oplus_{i=1}^r \pi_i[d_i] \end{equation}
avec $d_i \geq 1$ et $\pi_i \in \Pi_{\rm cusp}(\PGL_{n_i})$ pour $i=1,\dots,r$.  On observe que les $\pi_i$ sont
alg\'ebriques de poids motivique $ \leq 22$. D'apr\`es le th\'eor\`eme \ref{classpoids22}, on en d\'eduit que pour tout $i$ la repr\'esentation $\pi_i$ est dans l'ensemble $\Pi$. Autrement dit, on a $\psi \in \Phi_g$. 

\begin{lemmev}\label{lemmeclek12} Soient $1\leq g\leq 12$, $F \in {\rm S}_{12}({\rm Sp}_{2g}(\Z))$ une forme propre et $\psi$ le param\`etre standard de $F$. Alors $\psi \in \Phi_g$ et $\psi$ v\'erifie la condition {\rm (C)}. \end{lemmev}

\begin{pf} Nous venons de v\'erifier $\psi \in \Phi_g$. La condition {\rm (C)} \'etant automatique quand $g\leq10$, on peut supposer $g \geq 11$. On \'ecrit $\psi = \oplus_{i=1}^r \pi_i[d_i]$ avec $\pi_i \in \Pi$ pour tout $i$. \ps \ps

{\sc Cas $g=11$.} Si $g=11$ alors $0$ est valeur propre simple de $\psi_\infty$. Nous pouvons donc supposer, quitte \`a r\'eindexer les $\pi_i$, que l'on a\footnote{Le lecteur ayant dig\'er\'e les consid\'erations du chapitre VIII remarquera que cette hypoth\`ese est en contradiction avec la conjecture VIII.\ref{conjaj}, par exemple d'apr\`es le th\'eor\`eme VIII.\ref{amfexplsp} ; nous allons en effet aboutir \`a une contradiction, mais en utilisant plut\^ot les r\'esultats du \S \ref{complboc}.} $\pi_r=1$ et $d_r>1$. Posons $$g'= \frac{23-d_r}{2}.$$
C'est un entier v\'erifiant $0 \leq g' \leq 10$ car $d_r$ est impair $>1$. Ainsi, soit on a $\psi = [23]$, soit $d_{r} <23$ et $\psi$ s'\'ecrit sous la forme $\psi' \oplus [d_r]$ avec $\psi' \oplus [1] \in  \Phi_{g'}$. L'in\'egalit\'e $g' \leq 10$ assure que $\psi' \oplus [1]$ satisfait {\rm (C)} : c'est l'un des $12$ \'el\'ements de la table \ref{tablek=12} contenant $[1]$ et v\'erifiant $g' \leq 10$ (Lemme \ref{caractablek=12}). Il y a donc en tout $13$ possibilit\'es pour $\psi$, et il reste \`a d\'emontrer qu'aucune d'entre elles n'est possible. D'apr\`es la proposition VIII.\ref{cormizudisc} appliqu\'ee \`a la forme $F$ (cas $k=g+1$), il suffit pour cela de v\'erifier que l'on a $\delta(\pi_F,\frac{d_r+1}{2})=0$ dans chacun des cas, dans les notations {\it loc. cit.} (voir la formule \eqref{autreexpdeltaa}). C'est \'evident si $\psi = [23]$. Dans les autres cas, on conclut par le lemme \ref{lemmeverifg1112} ci-dessous et la relation $d_r + 1 = 24 -2g'$.  \ps \ps

{\sc Cas $g=12$}.  On proc\`ede de mani\`ere similaire en genre $g=12$. Supposons en effet que $\psi$ ne satisfait pas la condition ${\rm (C)}$. Nous affirmons que l'un des $\pi_i$ vaut ${\rm Sym}^2 \Delta_{11}$ et que deux des $\pi_i$ valent $1$. En effet, $\psi_\infty$ admet $0$ pour valeur propre triple, et les entiers $\pm 1, \dots, \pm 11$ pour valeurs propres simples. Si $(\pi_i[d_i])_\infty$, pour $i=1,\dots,r$, admet la valeur propre $0$, le lemme \ref{lemmecombigrandpi} (i) montre que l'on a soit $\pi_i={\rm Sym}^2 \Delta_{11}$ et $d_i=1$, soit $\pi_i=\Delta_{11}$ et $d_i=12$, soit $\pi_i=1$. Le second cas est exclu car il entra\^ine $\psi=\Delta_{11}[12] \oplus [1]$, qui v\'erifie ${\rm (C)}$. La valeur propre $11$ de $\psi_\infty$ \'etant simple, cela d\'emontre l'affirmation plus haut. \ps

Ainsi, on a $r\geq 3$ et quitte \`a r\'eindexer les $\pi_i$ on peut supposer que l'on a $\pi_1={\rm Sym}^2 \Delta_{11}$, $\pi_r=\pi_{r-1}=1$, $d_{r-1}=1$ et $d_{r}>1$ (d'apr\`es le th\'eor\`eme d'Arthur VIII.\ref{arthurstbis}, on rappelle qu'on ne peut avoir $d_r=d_{r-1}=1$). En particulier, si l'on pose \`a nouveau $g'=\frac{23-d_{r}}{2}$, on a $$\psi=\psi' \oplus [1] \oplus [d_r]$$ avec $\psi' \in \Phi_{g'}$ contenant ${\rm Sym}^2 \Delta_{11}$. L'in\'egalit\'e $d_r \geq 3$, {\it i.e.} $g' \leq 10$, montre que $\psi'$ v\'erifie la condition {\rm (C)} : il est dans la table \ref{tablek=12} par le lemme \ref{caractablek=12}. On constate par inspection qu'il y a $9$ possibilit\'es pour $\psi'$, et donc pour $\psi$. \ps

On exclut ensuite chacune de ces $9$ possibilit\'es en utilisant la proposition VIII.\ref{cormizudisc} appliqu\'ee \`a $F$ (cas $k=g$). Cette proposition s'applique encore \`a cause du lemme \ref{lemmeverifg1112}, qui contredit l'in\'egalit\'e $\delta(\pi_F,\frac{d_r+1}{2})>0$ (observer la relation $d_r + 1 = 24 - 2g'$). Cela termine la d\'emonstration du lemme \ref{lemmeclek12}.
\end{pf}

\begin{lemme}\label{lemmeverifg1112} Soient $1 \leq g' \leq 10$ et $\phi \in \Phi_{g'}$. \'Ecrivons $\phi = \oplus_{i=1}^s \varpi_i[q_i]$. Alors on a ${\rm L}(\frac{1}{2},\varpi_i )\neq 0$ pour tout $1 \leq i \leq s$ tel que $q_i \geq 24-2 g'$ et $\pi_i \neq 1$.
\end{lemme}
\ps\ps
\begin{pf} Seuls les $\phi$ de la table \ref{tablek=12} contenant un facteur de la forme $\Delta_{17}[d]$ (resp. $\Delta_{21}[2]$) m\'eritent une attention d'apr\`es la proposition \ref{ann12djk}. On conclut car on observe par inspection que l'on a toujours $24-2g'>d$ (resp. $g'\leq 10$).
\end{pf}

\noindent Terminons la d\'emonstration du th\'eor\`eme. Pour cela, v\'erifions tout d'abord l'assertion (ii), \`a savoir que pour tout $g=1,\dots,12$, 
l'injection $${\rm Ker}\, \vartheta_{g-1}/{\rm Ker} \,\vartheta_g \rightarrow {\rm S}_{12}({\rm Sp}_{2g}(\Z)),$$ induite par $\vartheta_g$, est surjective. 
Il s'agit de voir que si $g\leq 12$, toute forme propre $F \in  {\rm S}_{12}({\rm Sp}_{2g}(\Z))$ est dans l'image de $\vartheta_g$. Les deux lemmes ci-dessus assurent que son param\`etre standard $\psi(\pi_F,{\rm St})$ est dans la table \ref{tablek=12}. On conclut en constatant que dans tous les cas le crit\`ere de B\"ocherer s'applique : cela a d'ailleurs d\'ej\`a \'et\'e justifi\'e lors de la v\'erification du point 2 dans le \S \ref{preuve2thm24}. L'assertion (ii) est donc d\'emontr\'ee. \ps

Observons maintenant que si $g\leq 12$ et si $G,H \in {\rm S}_{12}({\rm Sp}_{2g}(\Z))$ sont deux formes propres pour ${\rm H}({\rm Sp}_{2g})$ telles que $\psi(\pi_G,{\rm St})=\psi(\pi_H,{\rm St})$, alors $G$ et $H$ sont proportionnelles. En effet, le paragraphe pr\'ec\'edent assure qu'il existe $G',H' \in \C[{\rm X}_{24}]$ tels que $\vartheta_g(G')=G$ et $\vartheta_g(H')=H$. Les formes $G$ et $H$ \'etant propre, la relation de commutation d'Eichler (Proposition V.\ref{commeichlernaif}) assure que l'on peut supposer que $G'$ et $H'$ sont propres pour ${\rm T}_2$. Mais cette m\^eme relation et l'identit\'e $\psi(\pi_{G},{\rm St})=\psi(\pi_H,{\rm St})$ imposent que $G'$ et $H'$ ont m\^eme valeur propre de ${\rm T}_2$ : elles sont donc proportionnelles d'apr\`es le calcul de Nebe et Venkov, et donc $G$ et $H$ sont proportionnelles. \ps

Soit $\Phi$ l'ensemble des param\`etres $\phi$ de la table \ref{tablek=12} tels qu'il existe une forme propre $G \in {\rm S}_{12}({\rm Sp}_{2g}(\Z))$  satisfaisant $\phi=\psi(\pi_G,{\rm St})$ (l'entier $g$ \'etant bien s\^ur uniquement d\'etermin\'e par $\phi$). Le paragraphe pr\'ec\'edent d\'emontre 
$|\Phi| = \sum_{g=1}^{12} \dim {\rm S}_{12}({\rm Sp}_{2g}(\Z))$. L'assertion (ii) du th\'eor\`eme entra\^ine alors  
$$|\Phi| =  \sum_{g=1}^{12} \dim({\rm Ker}\, \vartheta_{g-1}/{\rm Ker} \,\vartheta_g) = \dim ({\rm Ker} \,\vartheta_0) - \dim({\rm Ker}\, \vartheta_{12}) =23,$$
la derni\`ere \'egalit\'e provenant du r\'esultat d'Erokhin ${\rm Ker} \,\vartheta_{12} =0$ \cite{erokhin}.  
Comme il n'y a que $23$ param\`etres dans la table \ref{tablek=12}, $\Phi$ est l'ensemble de tous les param\`etres de cette table, ce qui d\'emontre les assertions (i) et (iii) du th\'eor\`eme (et justifie la remarque \ref{remfinthmE}).  \end{pf}\ps\ps
\begin{prop}\label{annuls12gnon23} Si $g>12$ et $g \neq 24$ alors on a ${\rm S}_{12}({\rm Sp}_{2g}(\Z)) = 0$. 
\end{prop}
\begin{pf} Les auteurs s'excusent du manque d'\'el\'egance de la d\'emonstration qui va suivre. Une raison \`a cela est que l'on ne dispose pas des outils ad\'equats pour la traiter, \`a savoir dans l'id\'eal d'un analogue du th\'eor\`eme VIII.\ref{amfexplsp} dans le cas o\`u l'\'el\'ement $\psi_\infty$ consid\'er\'e {\it loc. cit.} n'a pas toutes ses valeurs propres distinctes. \`A la place, nous utiliserons de mani\`ere ad hoc le substitut tr\`es partiel fourni par la proposition VIII.\ref{cormizunondisc}.\ps\ps

 Il s'agit de d\'emontrer qu'il n'existe pas de forme $F \in {\rm S}_{12}({\rm Sp}_{2g}(\Z))$ propre pour ${\rm Sp}_{2g}(\Z)$ pour $24> g>12$. Supposons donc qu'il existe une telle forme et notons $\psi$ son param\`etre standard. Les valeurs propres de $\psi_\infty$ sont les $2g+1$ entiers $0$ et $\pm (12-j)$ pour $1\leq j \leq g$. D'apr\`es un argument d\'ej\`a utilis\'e, 
le th\'eor\`eme d'Arthur, l'in\'egalit\'e $g<24$, et le th\'eor\`eme \ref{classpoids22}, montrent que l'on a $\psi = \oplus_{i=1}^r \pi_i[d_i]$ avec $\pi_i \in \Pi$ pour $i=1,\dots,r$ (l'ensemble $\Pi$ \'etant d\'efini par la formule \eqref{ensgrandpi}). Il sera commode de noter $\mathcal{X}_{\Pi}({\rm SL}_n)$ l'ensemble des \'el\'ements de $\mathcal{X}_{\rm AL}({\rm SL}_n)$ qui sont somme directe (au sens du \S VI.\ref{parconjarthlan}) d'\'el\'ements de la forme $\varpi[d]$ avec $\varpi \in \Pi$ et $d\geq 1$ ; on pose aussi $\mathcal{X}_{\Pi} = \coprod_{n\geq 1}   \mathcal{X}_{\rm AL}({\rm SL}_n)$.\ps\ps
 
 La premi\`ere assertion du lemme \ref{lemmecombigrandpi} montre que l'hypoth\`ese \eqref{propcormizu} de la proposition VIII.\ref{cormizunondisc} est satisfaite, ainsi donc que ses conclusions. Conform\'ement \`a l'\'enonc\'e de cette proposition, nous posons $I=\{i, \, \, \pi_i=1\}$, $I'=\{i \in I, d_i=1\}$, ainsi que $d={\rm Max}\{d_i, i \in I\}$ si $I \neq \emptyset$. Nous allons d'abord d\'emontrer les deux propri\'et\'es suivantes : \ps \ps \begin{itemize}

\item[(1)]  $\pi_i[d_i] \neq \Delta_{11}[12]$ pour tout $i=1,\dots, r$.\ps \ps

\item[(2)] $\pi_i\neq {\rm Sym}^2 \Delta_{11}$ pour tout $i=1,\dots,r$. \ps\ps
\end{itemize}

\noindent\ps\ps

\noindent Commen\c{c}ons par une observation. Soit $\pi \in \Pi$ v\'erifiant ${\rm L}(\frac{1}{2},\pi)=0$. D'apr\`es la proposition \ref{ann12djk} on a $\pi  \in \{\Delta_{17},\Delta_{21}\}$ et ${\rm ord}_{s=\frac{1}{2}} \, {\rm L}(s,\pi)=1$. En particulier, si ${\rm L}(\frac{1}{2},\pi_i)=0$ on constate que $\pi_i[d_i]$ est dans l'ensemble 
$$ \Psi_0\,\,=\, \{ \Delta_{17}[2], \, \, \Delta_{17}[4], \, \, \Delta_{17}[6], \, \, \Delta_{21}[2]\, \,\}.$$
Il sera utile d'avoir sous les yeux les valeurs propres de $\phi_\infty$ pour $\phi \in \Psi_0$ : elles sont donn\'ees dans la table \ref{tablepsi0} ci-contre. Chacun de ces \'el\'ements de $\Psi_0$ apparaissant au plus avec multiplicit\'e $1$ dans l'\'ecriture de $\psi$ d'apr\`es Arthur, on constate par exemple $\delta(\pi_F,a) \leq 4-a$ pour tout entier $1 \leq a \leq 4$. 

\begin{table}[h]
\renewcommand{\arraystretch}{1.5}
\begin{center}
\begin{tabular}{|c||c|c|c|c|}
\hline
$\phi$ & $\Delta_{17}[2]$ & $\Delta_{17}[4]$ & $\Delta_{17}[6]$ & $\Delta_{21}[2]$\cr
\hline ${\rm E}$ & $\{9,8\}$ & $\{10,9,8,7\}$ & $\{11,10,9,8,7,6\}$ & $\{11,10\}$\cr
\hline
\end{tabular}
\end{center}
\caption{Ensemble ${\rm E}$ des valeurs propres $\geq 0$ de $\phi_\infty$ pour $\phi  \in \Psi_0$.}
\renewcommand{\arraystretch}{1}
\label{tablepsi0}
\end{table}

 \ps\ps \medskip


\noindent {\it D\'emonstration de l'assertion {\rm (1)}.}  Supposons que l'un des $\pi_i[d_i]$, avec $i=1,\dots,r$, soit \'egal \`a $\Delta_{11}[12]$. On a donc une d\'ecomposition de la forme $\psi \,=\, \Delta_{11}[12] \,\oplus \,\psi'$ avec $\psi' \in \mathcal{X}_\Pi$. Observons que $\psi'_\infty$ a pour valeurs propres (simples) les entiers $\pm (g-12), \dots , \pm 1$ et $0$.  La contemplation de la table \ref{tablepsi0} montre donc : $\delta(\pi_F,1)=0$ si $g-12 <9$, $\delta(\pi_F,1) \leq 1$ si $g<23$, et $\delta(\pi_F,1) \leq 2$ dans tous les cas. \ps\ps

\noindent Si l'on a $I=I'$, l'assertion (e) de la proposition VIII.\ref{cormizunondisc} entra\^ine  $\delta(\pi_F,1) \geq g-14$. On en d\'eduit $g \leq 16$ et $\psi' = [2g-23]$, car $\psi' \in \mathcal{X}_\Pi$, ce qui contredit $I=I'$ et $g>12$. On a donc $I'=\emptyset$, $|I|=1$, et aussi $\frac{d-1}{2} \leq g-12$ car $[d]$ est un constituant de $\psi'$. L'assertion (d) de la proposition VIII.\ref{cormizunondisc} entra\^ine d'autre part $\frac{d+1}{4} \geq g-12- \delta(\pi_F,1) \geq g-14$. On en tire successivement $g \leq 17$, $\delta(\pi_F,1)=0$, $g=13$ et $\psi'=[3]$. La formule VIII.\eqref{quotxiaxib} s'\'ecrit dans ce cas $\xi_{\rm A}(s,\pi_F,{\rm St})=\xi_{\rm B}(s,\pi_F,{\rm St})$. Une consid\'eration de l'ordre d'annulation en $s=1$ de ces fonctions conduit \`a une contradiction :  la non-annulation de ${\rm L}(\frac{1}{2},\Delta_{11})$ montre ${\rm ord}_{s=1} \,\xi_{\rm A}(s,\pi_F,{\rm St}) \,=\, -2$, mais d'apr\`es Mizumoto (\S \ref{complboc}) nous avons aussi ${\rm ord}_{s=1} \, \xi_{\rm B}(s,\pi_F,{\rm St}) \,\geq \,-1$. \ps\ps \medskip

\noindent {\it  D\'emonstration de l'assertion {\rm (2)}.} Supposons maintenant qu'il existe un entier $i$ tel que $\pi_i={\rm Sym}^2 \Delta_{11}$. Dans ce cas, on a $d_i=1$ et il n'existe donc d'apr\`es Arthur qu'un seul tel entier $i$. D'apr\`es l'assertion (1) et le lemme \ref{lemmecombigrandpi} (i), on a $|I|=2$, et donc $d\geq 3$. Si l'on a $I' \neq \emptyset$, la consid\'eration de la valeur propre double $1$ de $\psi_\infty$ et le lemme \ref{lemmecombigrandpi} (ii) montrent que l'on a une d\'ecomposition $\psi \,=\, {\rm Sym}^2 \Delta_{11}\, \oplus \,\Delta_{11}[10] \,\oplus\, [1] \,\oplus\, [d]\, \oplus \,\psi'$ avec $\psi' \in \mathcal{X}_\Pi$. En particulier, on a $\frac{d-1}{2} \leq g-12$ et $\delta(\pi_F,1)\leq 2$. La proposition VIII.\ref{cormizunondisc} (c) entra\^ine $\frac{d+1}{4} \geq g-13$, puis $g \leq 15$. Les valeurs propres de $\psi'_\infty$ sont donc parmi $\pm 3, \pm 2, \pm 1$, ce qui contredit $\psi' \in \mathcal{X}_\Pi$.  Au final, nous avons donc $I'=\emptyset$, $|I|=2$, $d\geq 5$ et $\frac{d+1}{2} \geq  g - 9 - \delta(\pi_F,1)$ d'apr\`es le point (b) {\it loc. cit}.  \ps\ps
 
\noindent Supposons $\delta(\pi_F,\frac{d+1}{2})>0$. La description de $\Psi_0$ montre que cela entra\^ine $d=5$ et que l'on a une d\'ecomposition de la forme $\psi = {\rm Sym}^2 \Delta_{11} \oplus \Delta_{17}[6] \oplus [5] \oplus [3] \oplus \psi'$ avec $\psi' \in \mathcal{X}_\Pi$. En particulier, on a $g=23$ et $\delta(\pi_F,1) \geq 11$, ce qui est absurde. On a donc $\delta(\pi_F,1) \geq 2$ d'apr\`es le second point du (b) de la proposition {\it loc. cit.} . \ps\ps

\noindent Supposons enfin que $\psi$ ``contienne'' un \'el\'ement $\phi \in \Psi_0$ tel que $11$ est valeur propre de $\phi_\infty$. Dans ce cas, on a $g=23$ et $\frac{d-1}{2} \leq 10$. L'in\'egalit\'e $\frac{d+1}{2} \geq  g - 9 - \delta(\pi_F,1)$ donne alors $\delta(\pi_F,1) \geq 3$. Par inspection de la table \ref{tablepsi0} cela force $\frac{d-1}{2} \leq 7$, puis $\delta(\pi_F,1) \geq 6$, ce qui est absurde. Ainsi, la seule possibilit\'e restante est celle d'une d\'ecomposition de la forme  $\psi \,=\, {\rm Sym}^2 \Delta_{11}\, \oplus \,\Delta_{17}[2]\, \oplus \,\Delta_{17}[4]\, \oplus\, [d] \,\oplus\, \psi'$ avec $\psi' \in \mathcal{X}_\Pi$. Dans ce cas, on a les in\'egalit\'es $g-12 \geq 9$ et $\frac{d-1}{2} \leq 7$. Elles entra\^inent $\delta(\pi_F,1) \geq 4$ : une contradiction.

\medskip
\noindent {\it Fin de la d\'emonstration.} La conclusion des assertions (1) et (2), compte tenu du lemme \ref{lemmecombigrandpi} (i), est que l'on a $|I|=3$ et $|I'|=1$. Supposons $\delta(\pi_F,\frac{d+1}{2})>0$. Cela entra\^ine $d=5$ et que l'on a une d\'ecomposition de la forme $\psi \,=\, \Delta_{17}[6]\, \oplus [5]\, \oplus [3]\, \oplus \,[1] \,\oplus \,\psi'$ avec $\psi' \in \mathcal{X}_\Pi$. De plus, on dispose de l'in\'egalit\'e $\delta(\pi_F,1)  \geq g-11$ donn\'ee par la proposition VIII.\ref{cormizunondisc} (a). On en tire $g \leq 14$, ce qui entra\^ine $g=14$, et $\delta(\pi_F,1)=3$ : c'est absurde \'etant donn\'e la description de $\Psi_0$. En conclusion, la seconde assertion de la proposition VIII.\ref{cormizunondisc} (a) entra\^ine $\delta(\pi_F,1) = 3$. Un examen de la table \ref{tablepsi0} montre que l'on a une d\'ecomposition de la forme $\psi = \Delta_{21}[2] \oplus \Delta_{17}[2] \oplus \psi'$ avec $\psi' \in \mathcal{X}_\Pi$ et que l'\'el\'ement $\psi'_\infty$ admet les valeurs propres $\pm 10, \pm 9$ et $\pm 8$. En particulier, on a $g-12 \geq 10$ et $\frac{d-1}{2} \leq 7$. Mais l'assertion (a) {\it loc. cit.} montre ausi $\frac{d+1}{2} \geq g-11 \geq 11$, qui entra\^ine $d \geq 23$, ce qui est absurde.\end{pf}
\ps \ps
\begin{cor}\label{dimm12sp2g} Soit $1 \leq g \leq 23$. L'application $\vartheta_g : \C[{\rm X}_{24}] \rightarrow {\rm M}_{12}({\rm Sp}_{2g}(\Z))$ est surjective, et la dimension de ${\rm M}_{12}({\rm Sp}_{2g}(\Z))$ est donn\'ee par la table : \ps \ps
\renewcommand{\arraystretch}{1.5}
\begin{center}
\begin{tabular}{|c||c|c|c|c|c|c|c|c|c|c|c|c|}
\hline $g$ & $1$ & $2$ & $3$ & $4$ & $5$ & $6$ & $7$ & $8$ & $9$ & $10$ & $11$ & $\geq 12$ \cr
\hline $\dim\,{\rm M}_{12}({\rm Sp}_{2g}(\Z))$ & $2$ & $3$ & $4$ & $6$ & $8$ &  $11$ & $14$ & $18$ & $20$ & $22$ & $23$ & $24$ \cr
\hline
\end{tabular}
\end{center}
\renewcommand{\arraystretch}{1}

\end{cor}

\begin{pf} On proc\`ede par r\'ecurrence sur l'entier $g$. Le r\'esultat est bien connu si $g=1$. Si $g>1$ on rappelle que l'on dispose de l'op\'erateur de Siegel 
$\Phi_g : {\rm M}_{12}(\Sp_{2g}(\Z)) \longrightarrow {\rm M}_{12}({\rm Sp}_{2g-2}(\Z))$ ; il v\'erifie la relation $\Phi_g \circ \vartheta_g = \vartheta_{g-1}$. Par hypoth\`ese de r\'ecurrence, l'application $\Phi_g \circ \vartheta_g$ est donc surjective. Mais l'application $\vartheta_g \, : \, {\rm Ker}\,\, \vartheta_{g-1} \rightarrow {\rm S}_{12}({\rm Sp}_{2g}(\Z))$ est \'egalement surjective si $g \leq 12$ d'apr\`es l'assertion (ii) du th\'eor\`eme \ref{mainthmsiegel12}, ainsi que si $12 < g \leq 23$ d'apr\`es la proposition \ref{annuls12gnon23}. Cela d\'emontre la surjectivit\'e de $\vartheta_g$ et $\Phi_g$. \end{pf}

\ps\ps
\begin{remarque}{\rm Il n'est pas vrai que toute forme de Siegel parabolique de poids $16$ pour ${\rm Sp}_{2g}(\Z)$ est une combinaison lin\'eaire de s\'eries th\^eta d'\'el\'ements de ${\rm X}_{32}$, comme le d\'emontre le corollaire VII.\ref{cex32} (le contre-exemple donn\'e \'etant en genre $g=14$). 
}
\end{remarque}

\subsection{Formes de poids $\leq 11$} 

\begin{thmv}\label{thmsiegelsmallweight} Soient $k, g \in \Z$ tels que $ k \leq 11$ et $g \geq 1$.
Alors on a ${\rm S}_k({\rm Sp}_{2g}(\Z))=0$, \`a moins que l'on ne soit
dans l'un des cas suivants : \ps \ps \medskip

\begin{itemize}
\item[(i)] $k=8$ et $g=4$. Dans ce cas ${\rm S}_8({\rm Sp}_8(\Z))$ est de
dimension $1$, engendr\'e par la forme de Schottky, de param\`etre
standard $\Delta_{11}[4] \oplus [1]$.\ps \ps
\item[(ii)] $k=10$ et $g=2$. Dans ce cas ${\rm S}_{10}({\rm Sp}_4(\Z))$ est de
dimension $1$, engendr\'e par la forme de Saito-Kurorawa $F_{10}$, de param\`etre
standard $\Delta_{17}[2] \oplus [1]$. \ps \ps
\item[(iii)] $k=10$ et $g=4$. Dans ce cas ${\rm S}_{10}({\rm Sp}_8(\Z))$ est
de dimension $1$, engendr\'e par la forme d'Ikeda de param\`etre
standard $\Delta_{15}[4] \oplus [1]$. \ps \ps
\item[(iv)] $k=10$ et $g=6$. Dans ce cas ${\rm S}_{10}({\rm Sp}_{12}(\Z))$ est
de dimension $1$, engendr\'e par une forme de param\`etre standard $\Delta_{17}[2]\oplus
\Delta_{11}[4] \oplus [1]$. \ps \ps
\item[(v)] $k=10$ et $g=8$. Dans ce cas ${\rm S}_{10}({\rm Sp}_{16}(\Z))$
est de dimension $1$, engendr\'e par la forme d'Ikeda de param\`etre
standard $\Delta_{11}[8] \oplus [1]$. \ps \ps
\item[(vi)] $k=11$ et $g=6$. Dans ce cas toute forme propre de ${\rm S}_{11}({\rm
Sp}_{12}(\Z))$ a pour param\`etre standard $\Delta_{17}[4]\oplus \Delta_{11}[2] \oplus [1]$. De plus, si l'on admet la conjecture VIII.\ref{conjaj}, on a ${\rm S}_{11}({\rm
Sp}_{12}(\Z))=0$.
\end{itemize}
\end{thmv}

\ps \ps

\noindent Observons que le cas particulier $g=2$ du th\'eor\`eme ci-dessus se d\'eduit du
travail d'Igusa \cite{igusagenus2}, qui montre que ${\rm S}_k({\rm
Sp}_4(\Z))$ est nul si $k \leq 11$, de dimension $1$ et engendr\'e par $F_{10}$ si $k=10$ : voir le \S
\ref{prelimsp4}. Rappelons aussi que l'annulation ${\rm S}_k({\rm Sp}_{2g}(\Z))=0$ pour tout $g\geq 1$ et tout $k\leq 6$ et due \`a Duke et Imamo\u{g}lu \cite{dukei1}.  \ps

\begin{pf} Soient $k\leq 11$ et $F \in {\rm S}_k({\rm Sp}_{2g}(\Z))$ une 
forme propre pour ${\rm H}({\rm Sp}_{2g})$, de repr\'esentation engendr\'ee
$\pi_F \in \Pi_{\rm disc}({\rm Sp}_{2g})$, et de param\`etre standard $\psi
= \psi(\pi_F,{\rm St})=\oplus_{i=1}^r \pi_i[d_i]$ d'apr\`es Arthur. 
D'apr\`es la th\'eorie des formes singuli\`eres, on a $g \leq
2k$. Rappelons que $\psi(\pi_F,{\rm St})_\infty$ admet pour valeurs propres les
$2g+1$ entiers $0$ et $\pm (k-j), j=1,\dots,g$ d'apr\`es le corollaire
VI.\ref{infcarsiegel}. Les repr\'esentations $\pi_i$ sont donc alg\'ebriques de poids motivique $\leq 20$, sauf peut-\^etre si $g=2k=22$ auquel cas l'une d'entre elles peut \^etre de poids motivique $21$ ou $22$ (les autres \'etant de poids motivique $\leq 20$). Supposons dans un premier temps $(g,k) \neq (22,11)$. Le th\'eor\`eme \ref{classpoids22} assure que
les $\pi_i$ sont parmi l'ensemble $$\Pi=\{1, \Delta_{11},
\Delta_{15}, \Delta_{17}, \Delta_{19}, \Delta_{19,7}\}.$$ 
Il s'ensuit que $\psi$ est une somme, sans multiplicit\'es d'apr\`es Arthur
\cite[Thm. 1.5.2]{arthur}, d'\'el\'ements de la forme $[m]$ avec $m\geq 1$
impair, ou de l'ensemble 
$$\Psi := \{ \Delta_w[d] \,\, |\, \,  d \equiv 0 \bmod 2, w+d-1 \leq 20\}  \coprod \{
\Delta_{19,7}[2]\}.$$

\bigskip
\begin{lemmev}\label{lemmek11} Si $F \in {\rm S}_k({\rm Sp}_{2g}(\Z))$ est une forme propre avec $k \leq 11$ et $g \geq 1$ alors on a $k>g$. De plus, si $\psi = \oplus_{i=1}^k \pi_i[d_i]$ d\'esigne le param\`etre standard de $F$, il existe un unique $i \in \{1,\dots, r\}$ tel que
$\pi_i=1$, et l'on a $d_i=1$. 
\end{lemmev}

\ps\ps Admettons temporairement ce lemme et poursuivons la d\'emonstration. On observera que son \'enonc\'e est \'evident sous l'hypoth\`ese $k>g+1$, qui est un cas d\'ej\`a suffisament int\'eressant. On a donc $k>g$ et $\psi$ est une somme sans multiplicit\'es de
$[1]$ et d'\'el\'ements de $\Psi$. Nous nous sommes ainsi ramen\'es \`a un
assez petit nombre de possibilit\'es \`a \'etudier, ce que nous ferons au cas par cas.  

\ps \medskip
{\sc Cas $k \leq 6$}. Dans ce cas, on a ${\rm w}(\pi_i) \leq
10$ pour tout $i$, et donc $\psi=[1]$ et la forme $F$ n'existe
pas ! 

\ps \medskip
{\sc Cas $k=7$}. Dans ce cas, on a ${\rm w}(\pi_i) \leq 11$, et donc $\pi_i \in \{
\Delta_{11},1\}$, pour tout $i$. La seule possibilit\'e est 
$\psi=\Delta_{11}[2]\oplus [1]$, en particulier $g=2$, un cas que nous avons d\'ej\`a trait\'e
: la forme $F$ n'existe pas non plus car d'apr\`es Igusa on a l'annulation $\dim {\rm S}_7({\rm Sp}_4(\Z))=0$. 

\ps \medskip
{\sc Cas $k=8$}. La seule possibilit\'e est $\psi = \Delta_{11}[4] \oplus [1]$.  En
particulier, $g=8$ et $F \in {\rm S}_8({\rm Sp}_8(\Z))$.  Observons incidemment que
si l'on applique notre raisonnement au cas o\`u $F=J=\vartheta_{4}({\rm
E}_{8} \oplus {\rm E}_8)-\vartheta_{4}({\rm E}_{16})$ est la forme de
Schottky (voir le \S \ref{casdim16}), nous obtenons une nouvelle d\'emonstration du fait que le
param\`etre standard de la repr\'esentation engendr\'ee par $J$ est
$\Delta_{11}[4] \oplus [1]$ (corollaire VII.\ref{corpij} (i)), puisque c'est
l'unique param\`etre possible.  \ps

Pour conclure la d\'emonstration du th\'eor\`eme \ref{thmsiegelsmallweight}
pour $k=8$, il suffit d'invoquer le fait que ${\rm S}_8({\rm Sp}_8(\Z))$ est
de dimension $1$ (et engendr\'e par la forme de Schottky) d'apr\`es un
r\'esultat de Poor et Yuen \cite{py} d\'ej\`a mentionn\'e au \S
\ref{casdim16}.  Donnons un autre argument. Observons que $F$ est dans l'image de
l'application lin\'eaire $\vartheta_4 : \C[{\rm X}_{16}] \longrightarrow
{\rm S}_8({\rm Sp}_8(\Z))$. En effet, le crit\`ere de B\"ocherer s'applique puisque la fonction $${\rm L}(s,\pi_F,{\rm
St})=\zeta(s)\prod_{i=0}^3{\rm L}(s+i-3/2,\Delta_{11})$$ ne s'annule pas en
$s=4$ (\S VII.\ref{rembocherer}). Comme ${\rm S}_k({\rm Sp}_{2g}(\Z))$ est engendr\'e par des formes
propres, cela montre ${\rm S}_8({\rm Sp}_8(\Z)) \subset {\rm Im} \,
\vartheta_4$. Cela conclut car il est \'evident que ${\rm S}_8({\rm Sp}_8(\Z)) \cap {\rm Im}
\, \vartheta_4$ est
engendr\'e par $J$ (\S \ref{casdim16}). Soulignons que cet argument pour d\'emontrer
$\dim {\rm S}_8({\rm Sp}_8(\Z))=1$ n'est pas nouveau : il avait d\'ej\`a
\'et\'e observ\'e par Duke et Imamo\u{g}lu dans \cite{dukei1}, qui ont pu
d\'emontrer que ${\rm L}(s,\pi_F,{\rm St})$ admet n\'ecessairement un p\^ole
simple en $s=1$, sans toutefois pouvoir en d\'eduire directement la forme exacte
\`a priori du param\`etre $\psi$. \ps

Donnons enfin un troisi\`eme argument pour d\'emontrer que ${\rm S}_8({\rm Sp}_8(\Z))$
est de dimension $1$. En effet, c'est une cons\'equence imm\'ediate de
notre analyse ci-dessus et du r\'esultat g\'en\'eral suivant d'Ikeda, pr\'ecisant son propre th\'eor\`eme
VII.\ref{thmikeda2}, 
qui s'av\`erera bien utile dans la suite de cette d\'emonstration (voir
aussi l'exemple VIII.\ref{verifikeda}). 

\begin{lemme}\label{lemmeikeda3} {\rm \cite[Thm. 7.1, \S 15]{ikeda3}} Soient $m$ et $g$ des entiers pairs, ainsi que $\pi \in \Pi_{\rm cusp}({\rm PGL}_2)$ la
repr\'esentation engendr\'ee par une forme propre de poids $m$ pour ${\rm
SL}_2(\Z)$. Il existe une forme $G \in {\rm S}_{\frac{m+g}{2}}({\rm
Sp}_{2g}(\Z))$ propre pour ${\rm H}({\rm Sp}_{2g})$ et telle que $\psi(\pi_G,{\rm St})=\pi[g] \oplus [1]$ si, et
seulement si, $m \equiv g \bmod 4$. De plus, si cette condition est
satisfaite alors la forme $G$ est unique \`a un scalaire pr\`es.
\end{lemme}

\ps \medskip
{\sc Cas $k=9$}. On trouve cette fois-ci trois possibilit\'es pour $\psi$, correspondant
respectivement \`a
des genres $2, 4$ et $6$. Celle de genre $2$ est $\psi = \Delta_{15}[2]\oplus [1]$,
qui ne se produit pas car $\dim {\rm S}_{9}({\rm Sp}_4(\Z))=0$ d'apr\`es Igusa. La seconde est
$$\psi = \Delta_{15}[2] \oplus \Delta_{11}[2] \oplus [1],$$ pour laquelle $g=4$. Consid\'erons 
l'application lin\'eaire  (\S V.\ref{sertharm}) $$\vartheta_{5,4} : 
{\rm M}_{{\rm H}_{5,4}(\R^8)}({\rm
O}_8) \rightarrow {\rm S}_9({\rm Sp}_8(\Z)).$$ Comme on a ${\rm L}(\frac{1}{2},\Delta_w) \neq 0$ pour $w=11,15$, 
le produit $\zeta(s)\prod_{i=0}^1{\rm L}(s+i-\frac{1}{2},\Delta_{11}){\rm
L}(s+i-\frac{1}{2},\Delta_{15})$ admet un p\^ole simple en $s=1$, de sorte
que le crit\`ere de B\"ocherer s'applique et montre que $\vartheta_{5,4}$
est surjective. Pour \'eliminer ce second cas, il suffit donc de montrer 
${\rm M}_{{\rm H}_{5,4}(\R^8)}({\rm O}_8)=0$. Mais cette annulation est cons\'equence des
tables de \cite[\S 2]{chrenard2} ; nous avons d'ailleurs d\'ej\`a rencontr\'e cette
propri\'et\'e au \S \ref{tableexempleso8}. Une mani\`ere plus directe de
l'obtenir et de constater que par trialit\'e (\S V.\ref{trialitepgso8}), on a pour tout entier $d\geq
0$ pair l'\'egalit\'e 
\begin{equation} \label{pfeqtri} {\rm dim} \, {\rm M}_{{\rm
H}_{d,1}(\R^8)}({\rm SO}_8) = {\rm dim} \, {\rm M}_{{\rm H}_{d/2,4}(\R^8)}({\rm
O}_8).
\end{equation}
On conclut alors car ${\rm M}_{{\rm H}_{10,1}(\R^8)}({\rm SO}_8)=0$
d'apr\`es le lemme V.\ref{invariantscox}.\ps\ps

Pour se rassurer, v\'erifions que ce dernier raisonnement est coh\'erent avec la formule du th\'eor\`eme
\ref{amfexplsp}, autrement dit, que la formule de multiplicit\'e d'Arthur sugg\`ere
bien que $ \Delta_{15}[2] \oplus \Delta_{11}[2] \oplus [1]$ n'est pas le
param\`etre standard d'une forme de Siegel. Mais cela vient de ce que
dans les notations de ce th\'eor\`eme, on a $\chi(s_2)=-1$ et $\varepsilon(\Delta_{15} \times \Delta_{11})=\varepsilon(\Delta_{11})=1$. \ps\ps

La derni\`ere possibilit\'e est $\psi = \Delta_{11}[6] \oplus [1]$, pour
laquelle $g=6$. Mais elle ne se produit pas d'apr\`es Ikeda car $6 \not \equiv 12 \bmod 4$
(Lemme \ref{lemmeikeda3}). \ps\ps

\ps \medskip
{\sc Cas $k=10$}. On trouve cette fois-ci quatre possibilit\'es pour $\psi$, correspondant
respectivement \`a des genres $2, 4, 6$ et $8$, \`a savoir 
$$\Delta_{17}[2] \oplus [1], \Delta_{15}[4]\oplus [1], \Delta_{17}[2]\oplus
\Delta_{11}[4] \oplus [1], \Delta_{11}[8]\oplus
[1].$$
Ils sont tous trait\'es par le lemme \ref{lemmeikeda3}, sauf celui de genre
$6$, \`a savoir $\psi=\Delta_{17}[2]\oplus \Delta_{11}[4] \oplus [1]$. Cela d\'emontre les
assertions (ii) \`a (v) du th\'eor\`eme, except\'e le fait que $\dim {\rm S}_{10}({\rm Sp}_{12}(\Z))=1$ dans le (iv). \ps\ps

D\'emontrons maintenant cette assertion. Nous avons vu ci-dessus que toute forme propre de ${\rm S}_{10}({\rm Sp}_{12}(\Z))$ a pour
param\`etre standard $\Delta_{17}[2]\oplus \Delta_{11}[4] \oplus [1]$. Une
inspection de la formule de multiplicit\'e d'Arthur sugg\`ere qu'une telle forme existe et a multiplicit\'e $1$. Proc\'edons autrement. 
Combin\'e au crit\`ere de B\"ocherer, cette propri\'et\'e montre
aussi que l'application
lin\'eaire $$\vartheta_{2,6} : {\rm M}_{{\rm H}_{2,6}(\R^{16})}({\rm O}_{16}) \rightarrow {\rm S}_{10}({\rm
Sp}_{12}(\Z))$$ 
est surjective. Pour conclure, il suffit donc de montrer 
$\dim {\rm M}_{{\rm H}_{2,6}(\R^{16})}({\rm O}_{16}) = 1$ et $\vartheta_{2,6}\neq 0$. Mais l'assertion de dimension d\'ecoule du (i) du corollaire \ref{invariantsdim16}. La v\'erification de la non-nullit\'e de $\vartheta_{2,6}$ sera report\'ee \`a la fin du paragraphe pour ne pas interrompre le fil de cette d\'emonstration. \ps\ps

\ps \medskip
{\sc Cas $k=11$}. Ce cas est d\'ej\`a relativement fastidieux. On trouve cette fois-ci
$8$ param\`etres $\psi$ possibles : $1$ en genre $2$, $2$ en genre $4$, $2$ en
genre $6$, $2$ en genre $8$ et $1$ en genre $10$. \ps\ps

Les deux param\`etres de genre $4$ sont $\Delta_{17}[4]\oplus [1]$ et
$\Delta_{19}[2]\oplus \Delta_{15}[2] \oplus [1]$. Le premier ne se produit
pas d'apr\`es Ikeda (Lemme \ref{lemmeikeda3}). Consid\'erons
l'application $\vartheta_{7,4} : {\rm M}_{{\rm H}_{7,4}(\R^8)}({\rm
O}_8) \rightarrow {\rm S}_{11}({\rm Sp}_8(\Z))$. Le crit\`ere de
B\"ocherer montre alors que $\vartheta_{7,4}$ est surjective. Nous allons
voir qu'elle est nulle, ce qui montrera que le cas
$\psi=\Delta_{19}[2]\oplus \Delta_{15}[2] \oplus [1]$ ne se produit pas non
plus. Observons pour cela que la relation \eqref{pfeqtri} assure $\dim {\rm
M}_{{\rm H}_{7,4}(\R^8)}({\rm O}_8) = 1$, \'etant donn\'e $\dim {\rm M}_{{\rm
H}_{14,1}(\R^8)}({\rm SO}_8)=1$ (Lemme V.\ref{invariantscox}).
D'apr\`es le th\'eor\`eme VII.\ref{arthurtrialite}, on sait qu'il existe une forme propre
$G \in {\rm M}_{{\rm H}_{7,4}(\R^8)}({\rm O}_8)$ engendrant une
repr\'esentation dans $\Pi_{\rm disc}({\rm O}_8)$ de param\`etre standard
$\Delta_{17}[4]$ ; on a donc ${\rm M}_{{\rm H}_{7,4}(\R^8)}({\rm O}_8)=\C G$.
Il suffit alors de prouver $\vartheta_{7,4}(G)=0$. Si cette forme est non
nulle, les relations d'Eichler-Rallis montrent que $\vartheta_{7,4}(G)$ est propre et de
param\`etre standard $\Delta_{17}[4] \oplus [1]$ (Corollaire
VII.\ref{corparamrallis}), ce qui contredit le lemme \ref{lemmeikeda3} (ainsi que
le crit\`ere de B\"ocherer!). \ps\ps

Les deux param\`etres possibles de genre $6$ sont
$\Delta_{15}[6]\oplus[1]$ et $\Delta_{17}[4]\oplus \Delta_{11}[2]\oplus
[1]$. Le premier est encore exclu par Ikeda (Lemme \ref{lemmeikeda3}). Le second ne devrait pas se produire d'apr\`es la
formule de multiplicit\'e d'Arthur. En effet, dans les notations du th\'eor\`eme \ref{amfexplsp}) on a $\chi(s_2)=-1$ et $\varepsilon(\Delta_{11} \times \Delta_{17})^2\varepsilon(\Delta_{11}) = 1$. Cependant, nous ne voyons pas comment l'\'eliminer directement, ce qui explique le point
(vi) de l'\'enonc\'e du th\'eor\`eme. Observons que l'on ne tire rien cette fois-ci de l'application $\vartheta_{3,6} : {\rm M}_{{\rm
H}_{3,6}(\R^{16})}({\rm O}_{16}) \rightarrow {\rm S}_{11}({\rm
Sp}_{12}(\Z))$, car le crit\`ere de B\"ocherer montre
$\vartheta_{3,6}=0$ (en fait on a ${\rm M}_{{\rm
H}_{3,6}(\R^{16})}({\rm O}_{16})=0$ d'apr\`es le lemme \ref{invariantsdim16}). \ps\ps

Les deux param\`etres possibles de genre $8$ sont
$\Delta_{19,7}[2]\oplus \Delta_{15}[2] \oplus \Delta_{11}[2] \oplus [1]$ et 
$\Delta_{19}[2]\oplus \Delta_{11}[6]\oplus [1]$. Le crit\`ere de B\"ocherer
montre alors que l'application
$$\vartheta_{3,8} : {\rm M}_{{\rm H}_{3,8}(\R^{16})}({\rm
O}_{16}) \rightarrow {\rm S}_{11}({\rm Sp}_{16}(\Z))$$
est surjective, comme on le v\'erifie imm\'ediatement \`a l'aide de la proposition \ref{ann12djk}. On en d\'eduit ${\rm S}_{11}({\rm Sp}_{16}(\Z))=0$ car on a  ${\rm M}_{{\rm H}_{3,8}(\R^{16})}({\rm O}_{16})=0$ d'apr\`es le lemme
\ref{invariantsdim16} (ii).\ps\ps

Le dernier param\`etre, de genre $10$, est $\Delta_{11}[10]\oplus [1]$, mais
il n'appara\^it pas d'apr\`es le lemme \ref{lemmeikeda3}, ce qui conclut la d\'emonstration du
th\'eor\`eme. Il ne reste qu'\`a d\'emontrer le lemme \ref{lemmek11}. \end{pf} \ps\ps

\begin{pf} (D\'emonstration du lemme \ref{lemmek11}) Supposons $k \leq g$ et posons $I=\{ i , \pi_i=1\}$. Supposons de plus $(k,g) \neq (11,22)$ pour commencer. L'analyse pr\'ec\'edant l'\'enonc\'e du lemme montre $\psi = \oplus_{i=1}^r \pi_i[d_i]$, o\`u pour tout entier $i$, soit $\pi_i[d_i]$ appartient \`a l'ensemble $\Psi$ introduit {\it loc. cit.}, soit $\pi_i=1$. Le lemme \ref{lemmecombigrandpi} montre $|I|=3$ et que l'hypoth\`ese \eqref{propcormizu} de l'\'enonc\'e de la proposition VIII.\ref{cormizunondisc} est satisfaite. En particulier, les conclusions du point (a) de cette proposition affirment que l'on a soit $\delta(\pi_F,1)> 2-\nu$,
o\`u $\nu \in \{0,1\}$ est tel que $\nu \equiv k \bmod 2$, soit $\delta(\pi_F,3)> 0$. \ps\ps

D'apr\`es la proposition \ref{ann12djk}, la quantit\'e $\delta(\pi_F,a)$ est le nombre de constituants de $\psi$ de la forme
$\Delta_{17}[d']$ avec $d'\geq 2a$. Si un tel constituant est dans $\Psi$, cela force $d \in \{2, 4\}$, ce qui montre $\delta(\pi_F,3)=0$ et $\delta(\pi_F,1) \leq 2$ dans tous les cas. La seule possibilit\'e est donc $\delta(\pi_F,1) = 2$, $\nu=1$ ({\rm i.e.} $k$ impair), et $$\psi\,=\,\Delta_{17}[2]\,\oplus \,\Delta_{17}[4] \,\oplus\,
\phi,$$ o\`u $\phi$ est une somme sans multiplicit\'es d'\'el\'ements de $\Psi \coprod \{ [m], m\geq
1\}$. En particulier, $10$ (resp. $9$) est valeur propre (resp. valeur propre double) de $\psi_\infty$. Comme $k$ est impair cela entra\^ine $k=11$ et $20 \leq g \leq 22$. Comme il n'y a pas de forme de Siegel non nulle de poids et genre tous deux impairs,
et puisqu'on a suppos\'e $(k,g) \neq (11,22)$, la seule possibilit\'e est $g = 20$. Le point (a) de la proposition VIII.\ref{cormizunondisc} entra\^ine alors que $[d]$ est un constituant de $\phi$ avec $d \geq 21$, ce qui est absurde car les valeur propres de $\phi_\infty$ sont de valeur absolue $\leq 7$. \ps\ps

Nous avons donc d\'emontr\'e $k>g$. En particulier, il existe un unique entier
$i$ tel que $\pi_i=1$, et il ne reste qu'\`a montrer $d_i=1$. C'est
\'evident si $k>g+1$, et l'on peut donc supposer que $k=g+1$. La proposition
VIII.\ref{cormizudisc} entra\^ine 
\begin{equation} \label{eqpfk11} \delta(\pi_F,\frac{d_i+1}{2})>0.
\end{equation} 
Mais un argument similaire \`a
celui donn\'e ci-dessus montre $\delta(\pi_F,2)=0$ si $k> g$, et donc $d_i=1$, \`a moins que
l'on ne soit dans le cas particulier $\psi\,=\,\Delta_{17}[4]\, \oplus\, \phi$, $\phi$ \'etant la somme de $[d_i]$ et d'\'el\'ements de
$\Psi$. Dans ce cas, on a $k=11$ et $g=10$, de sorte que $\phi_\infty$ admet
pour valeurs propres $0$ et les entiers $\pm j$ avec $1\leq j \leq 6$. Une
inspection des \'el\'ements de $\Psi$ montre que cela force
$\phi\, =\, \Delta_{11}[2]\,\oplus \,[9]$ ou $\phi=[13]$. En particulier, on a $d_i\geq 9$, et
la formule \eqref{eqpfk11} contredit $\delta(\pi_F,3)=0$.\ps\ps

Il ne reste qu'\`a exclure la possibilit\'e $(k,g)=(11,22)$. Le th\'eor\`eme \ref{classpoids22} entra\^ine dans ce cas que chacun des $\pi_i$ est dans l'ensemble 
\begin{equation} \label{liste22prop} \{1,\Delta_{11},\Delta_{15},\Delta_{17},\Delta_{19},\Delta_{21},\Delta_{19,7},\Delta_{21,5},\Delta_{21,9},\Delta_{21,13},{\rm Sym}^2 \Delta_{11}\}.\end{equation}
Nous allons encore appliquer la proposition VIII.\ref{cormizunondisc}, dont les hypoth\`eses sont satisfaites gr\^ace au lemme \ref{lemmecombigrandpi}. Si l'on pose $I'=\{i \in I, d_i=1\}$, on a donc 
\begin{equation}\label{inegk1122lemme} \sum_{i \in I-I'} [\frac{d_i+1}{4}] \geq 9-\delta(\pi_F,1) + |I'| + 2|I-I'|.\end{equation}
D'apr\`es la proposition \ref{ann12djk}, si pour un entier $1 \leq i \leq r$ on a ${\rm L}(\frac{1}{2},\pi_i)=0$ et $d_i \equiv 0 \bmod 2$, alors on a $\pi_i \in \{\Delta_{17}, \Delta_{21}\}$ et ${\rm ord}_{s=1/2} {\rm L}(s,\pi_i)=1$. \'Etant donn\'e que $\psi_\infty$ admet la valeur propre $11$ (resp. $9$) avec multiplicit\'e $1$ (resp. $2$), cela montre l'in\'egalit\'e $\delta(\pi_F,1) \leq 3$. \ps \ps

Supposons $\delta(\pi_F,1) = 3$. Dans ce cas un examen de la table \ref{tablepsi0} montre que l'on a une d\'ecomposition $\psi \,=\, \Delta_{21}[2]\, \oplus \,\Delta_{17}[4]\, \oplus \, \Delta_{17}[2] \, \oplus \, \phi$ o\`u $\phi$ est une somme sans multiplicit\'es d'\'el\'ements de $\Psi \coprod \{ [m], m\geq
1\}$, et $\phi_\infty$ a ses valeurs propres $\leq 7$ en valeur absolue. Le lemme \ref{lemmecombigrandpi} entra\^ine $|I|=3$ et si l'on pose $d = {\rm Max}\, \{ d_i, i \in I\}$ on a donc $\frac{d-1}{2} \leq 7$, alors que l'in\'egalit\'e  \eqref{inegk1122lemme} entra\^ine $\frac{d+1}{2} > 11$, ce qui est absurde. \ps\ps

On a donc $\delta(\pi_F,1) \leq 2$, de sorte que le membre de droite de l'in\'egalit\'e \eqref{inegk1122lemme} est $\geq 7$. Cela entra\^ine $I \neq I'$. Mieux, sous l'hypoth\`ese $I-I'=\{j\}$ on obtient $\frac{d_j+1}{4} \geq 8$, et donc $d_j \geq 31$ : c'est absurde car on a trivialement $d_i \leq 23$ si $i \in I$. Ainsi l'ensemble $I-I'$ a deux \'el\'ements $\{a,b\}$, avec disons $1<d_a<d_b$. L'in\'egalit\'e \eqref{inegk1122lemme} entra\^ine $\frac{d_b+1}{2}>11+|I'|$, qui s'\'ecrit aussi $d_b>21+2|I'|$. Or on a \'evidemment $d_b\leq 23$ et $d_i \equiv 1 \bmod 2$ pour $i \in I$. Cela impose $d_b=23$ et $|I'|=1$. On en d\'eduit $d_a \geq 23$ en utilisant \`a nouveau \eqref{inegk1122lemme}, ce qui contredit $d_a<d_b$. \end{pf}

\begin{prop}\label{tableauinvso16}  Soit $L \subset \R^{16}$ un r\'eseau unimodulaire pair. On note $V_\lambda$ la repr\'esentation irr\'eductible de ${\rm SO}(\R^{16})$ de plus haut poids $$\lambda = m_1 \varepsilon_1 + m_2 \varepsilon_2 + \cdots + m_8 \varepsilon_8, \, \, \, \, \, {\rm avec}\,\, \, \,  \,   m_1 \geq m_2 \geq \dots \geq m_8 \geq 0$$ (\S VI.\ref{exparlan}),  et $V_\lambda^{{\rm SO}(L)} \subset V_\lambda$ le sous-espace des invariants par le groupe fini ${\rm SO}(L) \subset {\rm SO}(\R^{16})$. On suppose $m_1 \leq 4$. \ps
\begin{itemize}
\item[(i)] Si $L \simeq {\rm E}_{16}$, alors $V_\lambda^{{\rm SO}(L)}=0$, \`a mois que $\lambda$ ne soit \'egal \`a $0$ ou de la forme $4 ( \sum_{i=1}^k \varepsilon_i)$ avec $1 \leq k \leq 8$, auxquels cas $\dim V_\lambda^{{\rm SO}(L)}=1$.\ps
\item[(ii)] Si $L \simeq {\rm E}_8 \oplus {\rm E}_8$, alors les couples $(\lambda, \dim V_\lambda^{{\rm SO}(L)})$ tels que $V_\lambda^{{\rm SO}(L)} \neq 0$ sont donn\'es par la table \ref{tableSO16}.
\end{itemize}
\end{prop}

\begin{pf}  Il s'agit d'un calcul bas\'e sur la formule du caract\`ere de Weyl, \`a la mani\`ere de ceux effectu\'es dans \cite[\S 2]{chrenard2}. Nous remercions Olivier Ta\"ibi de nous avoir fait b\'en\'eficier de son propre algorithme, plus rapide que celui utilis\'e {\it loc. cit.}\,, pour l'\'evaluation finale. Il n\'ecessite, comme travail pr\'eliminaire, une \'enum\'eration des polyn\^omes caract\'eristiques des \'el\'ements de ${\rm SO}(L)$, ainsi que leurs multiplicit\'es. Nous renvoyons \`a la feuille de calculs \cite{clcalc} pour une justification des affirmations qui vont suivre. \ps

Dans le cas du r\'eseau $L={\rm E}_{16}$, on a ${\rm O}(L)={\rm W}({\bf D}_{16}) \simeq \{\pm 1\}^{15} \ltimes \mathfrak{S}_{16}$, et l'\'enum\'eration ne pose pas de difficult\'e \`a l'aide de l'ordinateur (si l'on se restreint aux \'el\'ements de d\'eterminant $1$, on trouve par exemple $823$ polyn\^omes). Dans le cas du r\'eseau $L={\rm E}_8 \oplus {\rm E}_8$, le groupe ${\rm O}(L)$ est le produit semi-direct de $\Z/2\Z$ par ${\rm O}({\rm E}_8)^2$ induit par l'\'echange des deux facteurs. Les polyn\^omes caract\'eristiques des classes de conjugaison d'\'el\'ements de ${\rm O}({\rm E}_8)={\rm W}({\bf E}_8)$, ainsi que les cardinaux de ces classes, ont \'et\'e d\'etermin\'es par Carter \cite[Table 11]{carter}. Cela permet de conclure, en observant que le d\'eterminant d'une matrice par blocs de la forme $\left[\begin{array}{cc} X\, {\rm I}_m & g \\ h & X \,{\rm I}_m \end{array}\right]$, avec $g,h \in {\rm M}_m$ et $X$ une ind\'etermin\'ee, vaut $\det(X^2 {\rm I}_m - gh)$. Apr\`es calcul, on trouve par exemple $1544$ polyn\^omes caract\'eristiques pour ${\rm SO}({\rm E}_8 \oplus {\rm E}_8)$. \ps Mentionnons qu'avec un peu de patience, qui \`a ce stade manque aux auteurs, il devrait \^etre \'egalement possible de d\'emontrer "\`a la main'' la proposition !
\end{pf}

\begin{cor}\label{invariantsdim16} \begin{itemize}
\item[(i)] On a  les \'egalit\'es $$\dim {\rm M}_{{\rm H}_{2,6}(\R^{16})}({\rm SO}_{16}) = \dim {\rm M}_{{\rm H}_{2,8}(\R^{16})}({\rm O}_{16})=1.$$\ps
\item[(ii)] Pour tout entier $1 \leq g \leq 8$, on a ${\rm M}_{{\rm H}_{3,g}(\R^{16})}({\rm O}_{16})=0$. \ps
\end{itemize}
\end{cor}
\begin{pf} On rappelle que si $U$ est une repr\'esentation de ${\rm O}_n(\R)$, dont on note $U'$ la restriction \`a ${\rm SO}_n(\R)$, alors 
${\rm M}_U({\rm O}_n)$ est un sous-espace de ${\rm M}_{U'}({\rm SO}_n)$ (c'est l'application ${\rm res}$ introduite au \S IV.\ref{fautson}). De plus, la formule IV.\eqref{mwproduit} entra\^ine l'\'egalit\'e
$$\dim {\rm M}_{U'}({\rm SO}_n) = \dim (U'^{{\rm SO}({\rm E}_{16})}) + \dim (U'^{{\rm SO}({\rm E}_8 \oplus {\rm E}_8)}).$$ 
Si $g<r$, la repr\'esentation irr\'eductible ${\rm H}_{d,g}(\R^{2r})$ de ${\rm O}(\R^{2r})$ reste irr\'eductible quand on la restreint \`a ${\rm SO}(\R^{2r})$; si $g=r$ elle se d\'ecompose en une somme de repr\'esentations non isomorphes ${\rm H}_{d,g}(\R^{2r})^{\pm}$, conjugu\'ees sous l'action ext\'erieure de ${\rm O}(\R^{2r})$ (\S V.\ref{trialitepgso8}). On rappelle que d'apr\`es la formule V.\eqref{vecteurextremal}, le plus haut poids de ${\rm H}_{d,g}(\R^{2r})$ est $d \, \sum_{i=1}^g \varepsilon_i$ si $g<r$, et que ceux de ${\rm H}_{d,g}(\R^{2g})^{\pm}$ sont $ d(\pm \varepsilon_g + \sum_{i=1}^{g-1} \varepsilon_i)$. Le (ii) d\'ecoule alors de la proposition \ref{tableauinvso16}, car les deux dimensions en question sont nulles pour les poids de la forme $3(\sum_{i=1}^g \varepsilon_i)$, $g\geq 1$. Le (i) se d\'emontre de la m\^eme mani\`ere,  les invariants non nuls provenant seulement de ${\rm E}_8 \oplus {\rm E}_8$, en observant l'isomorphisme $${\rm M}_{{\rm H}_{d,8}(\R^{16})}({\rm O}_{16}) \simeq {\rm M}_{{\rm H}_{d,8}(\R^{16})^{\pm}}({\rm SO}_{16})$$ (c'est l'application ${\rm ind}$ du \S IV.\ref{fautson}). \end{pf}
\ps

Pour terminer ce paragraphe, disons un mot sur la construction par s\'erie th\^eta d'un \'el\'ement de ${\rm S}_{10}({\rm Sp}_{12}(\Z))$. Partons du r\'eseau $L={\rm E}_8 \oplus {\rm E}_8$ de $\R^8 \oplus \R^8$.  Soient $e=(e_1,\dots,e_6)$  un sextuplet d'\'el\'ements de $L \otimes \C$ engendrant un sous-espace isotrope de dimension $6$, ainsi que $P_e(v_1,\dots,v_6)=\DET \left[ e_i \cdot v_j \right]_{1 \leq i,j \leq 6}$, de sorte que $P_e^2 \in {\rm H}_{2,6}(L \otimes \R)$. Soient $Q \subset L$ un r\'eseau de rang $6$ et $v_1,\dots,v_6$ une $\Z$-base de $Q$. Pour des raisons semblables \`a celles \'evoqu\'ees au \S V.\ref{constthetag4}, l'\'el\'ement $P_e(v_1,\dots,v_6)^2$ ne d\'epend pas du choix de la $\Z$-base $v_i$ de $Q$, et peut donc \^etre not\'e $P_e(Q)^2$. Le coefficient de Fourier de $\vartheta_{2,6}({\rm E}_8\oplus {\rm E}_8,P_e^2)$ correspondant \`a la matrice de Gram d'une $\Z$-base de $Q$ est alors donn\'e par la formule 
$${\rm c}_Q(P_e^2) = |{\rm O}(Q)| \sum_{M \subset {\rm E}_8 \oplus {\rm E}_8} P_e(M)^2,$$
la somme portant sur les sous-r\'eseaux $M$ isom\'etriques \`a $Q$. Nous allons appliquer ceci au cas particulier o\`u $Q \simeq {\rm Q}(R)$ avec $R={\bf D}_6$ ou $R={\bf E}_6$. On posera ${\rm E}_6:={\rm Q}({\bf E}_6)$. On observe tout d'abord que si un tel sous-r\'eseau est dans ${\rm E}_8 \oplus {\rm E}_8$, il est n\'ecessairement inclus dans l'un des deux facteurs ${\rm E}_8$. De plus, il n'est pas difficile de d\'emontrer que les sous-r\'eseaux de 
${\rm E}_8$ isom\'etriques \`a ${\rm D}_6$ (resp. ${\rm E}_6$) sont exactement les orthogonaux des sous-r\'eseaux de ${\rm E}_8$ isom\'etriques \`a ${\rm A}_1 \oplus {\rm A}_1$ (resp. ${\rm A}_2$), ils sont donc ais\'es \`a \'enum\'erer \`a l'aide de l'ordinateur (voir la feuille de calculs \cite{clcalc}). \ps\ps

Terminons par une application num\'erique. Notons $\varepsilon_1,\dots,\varepsilon_8$ la base canonique du premier facteur $\R^8$, $\varepsilon_1',\dots,\varepsilon_8'$ celle du second, et prenons $e=(\varepsilon_j+i\varepsilon'_j)_{1\leq j \leq 6}$. La premi\`ere observation ci-dessus assure 
que si $Q \simeq {\rm D}_6, {\rm E}_6$ alors ${\rm c}_Q(P_e^2) = 2  |{\rm O}(Q)| \sum_{M \subset {\rm E}_8 \oplus 0} P_e(M)^2$, la somme portant cette fois-ci sur les sous-r\'eseaux 
du premier facteur ${\rm E}_8$ isom\'etriques \`a $Q$. La seconde observation permet d'\'evaluer cette somme. L'ordinateur nous donne \cite{clcalc}
$$\frac{{\rm c}_{{\rm E}_6}(P_e^2)}{2 |{\rm O}({\rm E}_6)|} = 120 \, \, \, \, \,{\rm et} \, \, \, \, \, \frac{{\rm c}_{{\rm D}_6}(P_e^2)}{2 |{\rm O}({\rm D}_6)|} = 540.$$
En particulier, ces deux coefficients sont non nuls et on constate 
$$\frac{{\rm c}_{{\rm D}_6}(P_e^2)}{{\rm c}_{{\rm E}_6}(P_e^2)} = 2.$$

\noindent Le corollaire qui suit fournit au final un proc\'ed\'e concret de construction de toutes les formes de Siegel de poids $\leq 11$ par s\'eries th\^eta. 

\begin{corv} Les applications  $\vartheta_{d,g}$ induisent des isomorphismes entres espaces de dimension $1$ :
$$\vartheta_{4,4} : {\rm M}_{{\rm H}_{4,4}(\R^8)}({\rm O}_8) \isomo {\rm S}_8({\rm Sp}_8(\Z)), \, \, \, \vartheta_{6,2} : {\rm M}_{{\rm H}_{6,2}(\R^8)}({\rm O}_8) \isomo {\rm S}_{10}({\rm Sp}_4(\Z)),$$
$$\vartheta_{6,4} : {\rm M}_{{\rm H}_{6,4}(\R^8)}({\rm O}_8) \isomo {\rm S}_{10}({\rm Sp}_8(\Z)), \, \, \, \vartheta_{2,6} : {\rm M}_{{\rm H}_{2,6}(\R^{16})}({\rm O}_{16}) \isomo {\rm S}_{10}({\rm Sp}_{12}(\Z)),$$
$$\vartheta_{2,8} : {\rm M}_{{\rm H}_{2,8}(\R^{16})}({\rm O}_{16}) \isomo {\rm S}_{10}({\rm Sp}_{16}(\Z)).$$
\end{corv}

\begin{pf} L'assertion sur $\vartheta_{4,4}$ est la proposition V.\ref{thetanonnul44}, et celle sur $\vartheta_{6,4}$ est de m\^eme contenue dans la table \ref{thetagenre4}. L'assertion sur $\vartheta_{6,2}$ est un cas particulier de la proposition \ref{calcimagetheta}. Celle sur $\vartheta_{2,6}$ r\'esulte de la discussion pr\'ec\'edant l'\'enonc\'e du corollaire. Enfin, l'assertion sur $\vartheta_{2,8}$ se d\'eduit de $\dim {\rm M}_{{\rm H}_{2,8}(\R^{16})}({\rm O}_{16})=1$ (Lemme \ref{invariantsdim16} (i)) et du crit\`ere de B\"ocherer, puisque l'on a ${\rm L}(\frac{1}{2},\Delta_{11}) \neq 0$.
\end{pf}

\section{Vers une nouvelle d\'emonstration de l'\'egalit\'e $|{\rm X}_{24}|=24$}\label{nou2425}

L'int\'er\^et du th\'eor\`eme suivant r\'eside dans le fait que sa d\'emonstration n'utilise aucun calcul appartenant \`a la th\'eorie des r\'eseaux unimodulaires. En particulier, elle n'utilise pas bien entendu la d\'etermination par Niemeier de ${\rm X}_{24}$. En revanche, elle repose sur la th\'eorie d'Arthur et le th\'eor\`eme \ref{classpoids22}. On rappelle que les ensembles ${\rm X}_n$ et $\widetilde{{\rm X}}_n$ ont \'et\'e introduits au \S IV.\ref{ensclassbil}.

\begin{thm}\label{niemeiersansreseau} Admettons les conjectures {\rm VIII.\ref{arthurstint}} et {\rm VIII.\ref{conjaj2}}. Alors on a les \'egalit\'es $$|{\rm X}_{24}|=24 \hspace{1cm} et \hspace{1cm} |\widetilde{{\rm X}}_{24}|=25.$$ 
\end{thm}

\ps \ps

\begin{pf} D\'emontrons d'abord l'\'egalit\'e $|\widetilde{{\rm X}}_{24}|=25$. Soit $\pi$ un \'el\'ement de $\Pi_{\rm disc}({\rm SO}_{24})$. D'apr\`es la conjecture VIII.\ref{arthurstint}, on a  $\psi(\pi,{\rm St}) \in \mathcal{X}_{\rm AL}({\rm SL}_{24})$. Si l'on suppose de plus la repr\'esentation $\pi_\infty$ triviale,  le th\'eor\`eme \ref{corretourtab24} entra\^ine que $\psi(\pi,{\rm St})$ est l'un des param\`etres list\'es dans la table \ref{table24}. La relation $$|\widetilde{\rm X}_{24}|\,\,=\,\,\dim\, {\rm M}_\C\,({\rm SO}_{24})$$ montre qu'il ne reste qu'\`a d\'eterminer, pour chacun des param\`etres $\psi$ de cette table, la somme $m_\psi$ des multiplicit\'es des repr\'esentations $\pi' \in \Pi({\rm SO}_{24})$ telles que $\pi'_\infty = 1$ et $\psi(\pi',{\rm St})=\psi$. D'apr\`es la conjecture {\rm VIII.\ref{conjaj2}}, on a toujours $m_\psi \in \{0,1,2\}$, sa valeur exacte pouvant \^etre d\'etermin\'ee \`a l'aide des formulaires donn\'es au \S \ref{formulairesAMF} ; en l'occurence, nous allons utiliser ici le th\'eor\`eme VIII.\ref{amfexplso0}. Cette d\'etermination des $24$ entiers $m_\psi$, essentiellement faite au ``cas par cas'', a d\'ej\`a \'et\'e men\'ee \`a bien au \S \ref{preuve1thm24}. La conclusion de l'analyse {\it loc. cit.} est la suivante : si $\psi \in \mathcal{X}_{\rm AL}({\rm SL}_{24})$ est l'un des \'el\'ements list\'es dans la table \ref{table24} alors \begin{itemize} \ps\ps
\item[(a)] soit $\psi \neq \Delta_{11}[12]$ et l'on a $m_\psi=1$, \ps
\item[(b)] soit $\psi = \Delta_{11}[12]$ et l'on a $m_\psi = 2$. Dans ce cas, il existe exactement deux repr\'esentations distinctes $\pi',\pi''$ dans $\Pi_{\rm disc}({\rm SO}_{24})$ v\'erifiant $\psi(\pi',{\rm St})=\psi(\pi'',{\rm St})=\psi$, et l'on a ${\rm m}(\pi')={\rm m}(\pi'')=1$ et  $\pi''\,=\,S(\pi')$ au sens de la remarque VIII.\ref{remmult2}. \ps\ps
\end{itemize}
On a bien d\'emontr\'e $|\widetilde{{\rm X}}_{24}|=23+2=25$. Montrons enfin $|{\rm X}_{24}|=24$. On part de l'\'egalit\'e $\dim \,{\rm M}_\C({\rm O}_{24})\, =\, |{\rm X}_{24}|$. On rappelle que le groupe \`a deux \'el\'ements ${\rm O}_{24}(\Q)/{\rm SO}_{24}(\Q)=\langle s \rangle $ agit sur ${\rm M}_\C({\rm SO}_{24})$ avec pour sous-espace fixe ${\rm M}_\C({\rm O}_{24})$ (voir l'exemple IV.\ref{exdim1det}).  Soit $f' \in {\rm M}_{\rm \C}({\rm SO}_{24})$ une fonction propre sous ${\rm H}({\rm SO}_{24})$ engendrant la repr\'esentation $\pi'$ introduite dans le point (b) ci-dessus. L'assertion $S (\pi') \neq \pi'$ montre que la droite $\C f'$ n'est pas stable par l'action de l'\'el\'ement $s$, et entra\^ine que le sous-espace vectoriel $$V := \langle f, s f \rangle \subset {\rm M}_\C({\rm SO}_{24})$$ est de dimension $2$ et v\'erifie $\dim \,( {\rm M}_\C({\rm O}_{24}) \cap V)  = 1$ (ce fait a d'ailleurs d\'ej\`a \'et\'e exploit\'e dans la d\'emonstration de la proposition VII.\ref{leechmoins2}). On en d\'eduit d'une part l'in\'egalit\'e $|{\rm X}_{24}| <25$ et d'autre part l'existence d'une repr\'esentation $\pi \in \Pi_{\rm disc}({\rm O}_{24})$ telle que $\psi(\pi,{\rm St})=\Delta_{11}[12]$. \ps\ps
Pour conclure $|{\rm X}_{24}|=24$, il suffit de montrer que pour chacun des $23$ \'el\'ements $\psi \neq \Delta_{11}[12]$ list\'es dans la table \ref{table24}, il existe une repr\'esentation $\pi \in \Pi_{\rm disc}({\rm O}_{24})$ telle que $\psi(\pi,{\rm St}) =\psi$. C'est \'evident dans le cas particulier $\psi = [23] \oplus [1]$, qui correspond \`a la repr\'esentation triviale (Exemples VI.\ref{exconjal}). Pour les autres, l'argument que nous proposons est assez indirect. Il consiste \`a appliquer verbatim la m\'ethode du \S \ref{preuve2thm24}. Elle se r\'esume ainsi : on observe d'abord que $\psi$ s'\'ecrit de mani\`ere unique sous la forme $\psi' \oplus [23-2g]$ avec $1 \leq g<12$ et $\psi' \in {\mathcal X}_{\rm AL}({\rm SL}_{2g+1})$, puis on utilise la conjecture VIII.\ref{conjaj2} pour v\'erifier (au cas pas cas!) que $\psi'$ est le param\`etre standard d'une forme modulaire de Siegel parabolique, propre, de poids $12$, pour ${\rm Sp}_{2g}(\Z)$, et enfin on v\'erifie que cette forme modulaire est une combinaison lin\'eaire de s\'eries th\^eta d'\'el\'ements de ${\rm X}_{24}$ \`a l'aide du crit\`ere de B\"ocherer. Aucun de ces arguments n'utilise la d\'etermination de ${\rm X}_{24}$. Cela conclut la d\'emonstration.
\end{pf}

\section{Quelques \'el\'ements de $\Pi_{\rm disc}({\rm SO}_{n})$ pour $n=15,17$ et $23$}

\begin{thmv}\label{thmpar151723} Soit $n$ un entier impair $\geq 1$. Les \'el\'ements $\psi \in {\mathcal X}_{\rm AL}({\rm SL}_{n-1})$ tels que la classe de conjugaison $\psi_\infty$ admette pour valeurs propres les $n-1$ demi-entiers $\pm \frac{n-2}{2},\,\, \pm \frac{n-4}{2},\,\, \dots, \,\,\pm \frac{1}{2}$ sont les suivants :  \ps \ps
\begin{itemize}
\item[(i)] l'unique \'el\'ement \,\,$[n-1]$\,\, si $n \leq 11$,\ps \ps
\item[(ii)] $[12]$ {\rm et} $\Delta_{11}\oplus[10]$ si $n=13$, \ps\ps
\item[(iii)] $[14]$ {\rm et} $\Delta_{11}[3]\oplus[8]$ si $n=15$,\ps \ps
\item[(iv)] $[16]$, \, $\Delta_{15}\oplus [14]$, \,
$\Delta_{15}\oplus \Delta_{11}[3] \oplus
[8]$\,  {\rm et} \,$\Delta_{11}[5] \oplus [6]$ si $n=17$, \ps \ps
\item[(v)] les $32$ param\`etres list\'es dans la table \ref{tableso23} si $n=23$. \ps \ps
\end{itemize}
\end{thmv}

\begin{pf} Il s'agit d'une cons\'equence du th\'eor\`eme \ref{classpoids22}, dont la d\'emonstration est similaire \`a celles du th\'eor\`eme \ref{corretourtab24} et de la proposition \ref{propliste24}. Mentionnons que les cas $n=19$ et $21$, bien que non explicit\'es dans l'\'enonc\'e, se d\'eduisent imm\'ediatement \'egalement de cette m\'ethode, et m\^eme de la table \ref{tableso23}.
\end{pf}

\noindent Ce th\'eor\`eme a des cons\'equences sur la classification des $\pi \in \Pi_{\rm disc}({\rm SO}_n)$ v\'erifiant $\pi_\infty = 1$ si $n$ est impair, les valeurs propres de ${\rm St}\,{\rm c}_\infty(\pi) \subset {\got{sl}}_{n-1}(\C)$ \'etant dans ce cas les $n-1$ demi-entiers $\pm \frac{n-2}{2},\,\, \pm \frac{n-4}{2},\,\, \dots, \,\,\pm \frac{1}{2}$. \ps\ps

\begin{thm}\label{thm1517} Admettons la conjecture VIII.\ref{arthurstint}. Les param\`etres
standards $\psi(\pi,{\rm St})$ des repr\'esentations $\pi \in \Pi_{\rm
disc}({\rm SO}_n)$ telles que $\pi_\infty =1$ sont : \ps \medskip
\begin{itemize}
\item[(i)] $[14]$ {\rm et} $\Delta_{11}[3]\oplus[8]$ si $n=15$,\ps \medskip
\item[(ii)] $[16]$, \, $\Delta_{15}\oplus [14]$, \,
$\Delta_{15}\oplus \Delta_{11}[3] \oplus
[8]$\,  {\rm et} \,$\Delta_{11}[5] \oplus [6]$ si $n=17$.\ps \medskip
\end{itemize}
\end{thm}

\begin{pf} Sous l'hypoth\`ese $n=15$ ou $n=17$, on rappelle que la dimension $\dim {\rm M}_\C({\rm SO}_n)=|{\rm X}_n|$ vaut respectivement $2$ et $4$ d'apr\`es le corollaire IV.\ref{hsonimpair}. L'op\'erateur ${\rm T}_2$ de $\C[{\rm X}_n]$ est
d\'etermin\'e dans l'appendice {\rm B} \S 5.  Dans les deux cas ses valeurs propres sont
bien compatibles avec le th\'eor\`eme ci-dessus; la chose qui importera ici est qu'elles soient distinctes. Ainsi, si $n=15$ (resp. $n=17$) il existe exactement $2$ (resp. $4$) \'el\'ements de $\Pi_{\rm disc}({\rm SO}_n)$ v\'erifiant $\pi_\infty=1$, chacun ayant multiplicit\'e $1$. L'entier $n$ \'etant impair, les param\`etres standards $\psi(\pi,{\rm St})$ de ces \'el\'ements sont \'evidemment distincts. D'apr\`es la conjecture VIII.\ref{arthurstint}, ces param\`etres satisfont les hypoth\`eses du th\'eor\`eme \ref{thmpar151723}. Mais d'apr\`es les conclusions de ce th\'eor\`eme, il n'existe que $2$ param\`etres possibles pour $n=15$, et $4$ pour $n=17$, ce qui termine la d\'emonstration.
\end{pf}

\begin{thm}\label{thm23} Admettons la conjecture {\rm VIII.\ref{conjaj2}}. Les param\`etres
standards $\psi(\pi,{\rm St})$ des repr\'esentations $\pi \in \Pi_{\rm
disc}({\rm SO}_{23})$ telles que $\pi_\infty =1$ sont les $32$ \'el\'ements list\'es dans la table \ref{tableso23} si $n=23$.
\end{thm}

\begin{pf} On rappelle l'\'egalit\'e $|{\rm X}_{23}|=32$ (Corollaire IV.\ref{hsonimpair}). Soient $\psi$ l'un des $32$ \'el\'ements de la table \ref{tableso23} et $\pi \in \Pi({\rm SO}_{23})$ l'unique repr\'esentation telle que $\psi(\pi,{\rm St})=\psi$. Il s'agit de d\'emontrer que la multiplicit\'e de $\pi$ est non nulle. D'apr\`es la  conjecture {\rm VIII.\ref{conjaj2}}, il suffit pour cela d'appliquer la formule de multiplicit\'e d'Arthur sous sa forme donn\'ee par le 
Th\'eor\`eme VIII.\ref{amfexplso1}. Il s'agit d'une \'etude au cas par cas, que l'on peut l\'eg\`erement simplifier en d\'egageant des crit\`eres du m\^eme type que ceux donn\'es au \S VIII.\ref{formexpliciteson0} dans le cas $n \equiv 0 \bmod 8$, ce que nous laissons en exercice au lecteur car les raisonnements sont tr\`es similaires \`a ceux \'etudi\'es en d\'etail {\it loc. cit.} On constate au final, \^o miracle, que l'on a bien ${\rm m}(\pi)=1$, quel que soit l'\'el\'ement $\psi$ choisi initialement.
\end{pf}

\begin{remarque} \label{remark25} {\rm \begin{itemize} \ps \ps \item[(i)]  Le th\'eor\`eme \ref{thm23} a \'et\'e
\'etendu \`a la dimension $25$ dans \cite[Thm.  1.14]{chrenard2}.  Les $121$
param\`etres en question font intervenir notamment les $7$ repr\'esentations
dans $\Pi_{\rm alg}(\PGL_6)$ de poids motivique $23$ mentionn\'ees au \S \ref{complementfexpl}. \ps \ps
\item[(ii)] Un argument semblable \`a celui de la d\'emonstration du th\'eor\`eme \ref{niemeiersansreseau} permet de red\'emontrer,  \`a partir du th\'eor\`eme \ref{classpoids22} et conditionnellement aux conjectures {\rm VIII.\ref{arthurstint}} et {\rm VIII.\ref{conjaj2}}, mais sans aucun calcul appartenant \`a la th\'eorie des r\'eseaux euclidiens, les \'egalit\'es $|{\rm X}_{7}|=|{\rm X}_9|=1$, $|{\rm X}_{15}|=2$, $|{\rm X}_{17}|=4$ et $|{\rm X}_{23}|=32$.
\end{itemize}
\ps
 } \end{remarque} \ps\ps
 
 Terminons par une curieuse observation. 
Pour tout entier $n \leq 24$ tel que $n \equiv -1,0,1 \bmod 8$, nous avons
d\'ecrit dans ce m\'emoire le sous-ensemble $$\Phi_n \subset
\mathcal{X}({\rm SL}_{2[\frac{n}{2}]})$$ constitu\'e des \'el\'ements de la
forme $\psi(\pi,{\rm St})$ avec $\pi \in \Pi_{\rm disc}({\rm SO}_n)$ tel
que $\pi_\infty=1$.  Par exemple $\Phi_{24}$ et $\Phi_{25}$ sont respectivement donn\'es par les tables \ref{table24} et \ref{tableso23}. 
Cette description est encore conditionnelle quand $n$ impair, auquel cas
elle s'\'etend m\^eme \`a $n=25$ d'apr\`es la remarque \ref{remark25} (i),
mais l\`a n'est pas la question dans cette discussion, o\`u nous serions
pr\^ets volontiers \`a admettre la conjecture VIII.\ref{conjaj2}. 
Consid\'erons le second probl\`eme, en apparence assez diff\'erent, de d\'eterminer le
sous-ensemble $$\Phi'_n \subset \mathcal{X}_{\rm AL}({\rm
SL}_{2[\frac{n}{2}]})$$ constitu\'e de tous les $\psi$  tels que $\psi_\infty$
est le caract\`ere infinit\'esimal de la repr\'esentation triviale de ${\rm
SO}_n(\R)$ (une condition simple sur ses valeurs propres).  La th\'eorie d'Arthur assure d'abord $\Phi_n  \subset \Phi'_n$ ; elle donne aussi un crit\`ere explicite, la ``formule de
multiplicit\'e d'Arthur'', permettant de d\'eterminer si un \'el\'ement $\psi \in \Phi'_n$
donn\'e est dans $\Phi_n$ : ce sont les formules
VIII.\ref{amfexplso1} et VIII.\ref{amfexplso0}. Mais le th\'eor\`eme \ref{classpoids22}
permet de d\'eterminer $\Phi'_n$ pour tout $n \leq 24$. La propri\'et\'e miraculeuse, v\'erifi\'ee
dans tous les cas, est alors l'\'egalit\'e \begin{equation} \label{equationmagique}
\Phi_n = \Phi'_n, \, \, \, \, \forall\, n \leq 24.\end{equation} Il
est m\^eme concevable que cette \'egalit\'e s'\'etende \`a $n=25$. 
Concr\`etement, cela signifie que pour $n\leq 24$ et tout $\psi \in \Phi'_n$,
la formule de multiplicit\'e d'Arthur appliqu\'ee \`a $\psi$ conduit
toujours \`a une multiplicit\'e non nulle.  Il serait int\'eressant de
trouver une raison plus profonde \`a ce ph\'enom\`ene.\ps\ps

On pourrait esp\'erer que l'\'egalit\'e \eqref{equationmagique} se produise pour tout $n \equiv -1,0,1$ $\bmod \, 8$, du moins si l'on remplace $\Phi'_n$ par son sous-ensemble constitu\'e des $\psi$ tels que $\psi_p = \psi_p^{-1}$ pour tout $p$. Il n'en est rien. En effet, la formule de multiplicit\'e d'Arthur montre que si $n=32$, le param\`etre $\Delta_{17}[14]\oplus [3] \oplus [1]$ de $\Phi'_{32}$ ne doit pas appartenir \`a $\Phi_{32}$. De m\^eme, l'\'el\'ement $\Delta_{17}[13] \oplus [4]$ de $\Phi'_{31}$ ne doit pas appartenir \`a $\Phi_{31}$. Si ${\rm S}_{14}({\rm SL}_2(\Z))$ n'\'etait pas nul, de tels exemples existeraient \'egalement en dimensions $23$ et $24$.

\chapter{Applications}

\parindent=0cm

\section{$24$ repr\'esentations galoisiennes $\ell$-adiques}\label{repgal24}

\newcommand{\rhob}{\overline{\rho}}

Rappelons que l'espace vectoriel $\Q[{\rm X}_{24}]$ admet une $\Q$-base ${\rm v}_1,\dots,{\rm v}_{24}$ constitu\'ee de vecteurs propres
communs \`a tous les \'el\'ements de ${\rm H}({\rm O}_{24})$
(\S\ref{so24etnv}). Chacun de ces vecteurs
${\rm v}_i$ engendre une repr\'esentation automorphe $\pi_i \in \Pi_{\rm disc}({\rm
O}_{24})$ dont le param\`etre standard $\psi(\pi_i,{\rm St})$ est
d\'etermin\'e par le th\'eor\`eme${}^{\color{green} \ast}$ \ref{thm24} (Table
\ref{table24bis}). En particulier, le couple
$(\pi_i,{\rm St})$ satisfait la
conjecture d'Arthur-Langlands, de sorte que le corollaire
VIII.\ref{corgaloisrepo} s'applique,
et lui associe des repr\'esentations $\ell$-adiques
de dimension $24$ du groupe ${\rm Gal}(\overline{\Q}/\Q)$. Dans ce qui suit,
nous allons pr\'eciser l'\'enonc\'e obtenu.
\ps

Si $p$ est un nombre premier et $1 \leq i \leq 24$, nous notons $\lambda_i(p)
\in \Z$ la valeur propre de ${\rm T}_p$ sur le vecteur ${\rm v}_i$.
D'apr\`es la formule VI.\eqref{formulesattp}, on a $\lambda_i(p) \, = \, p^{11}\, {\rm trace}
\, {\rm St} ({\rm c}_p(\pi_i))$. Plus g\'en\'eralement, le lemme suivant montre que le
polyn\^ome $\DET ( t - p^{11}\, {\rm St}({\rm c}_p(\pi_i)) )$ est \`a coefficients entiers,
chacun de ses coefficients pouvant \^etre vu comme la valeur propre d'un op\'erateur de Hecke
bien choisi dans ${\rm H}({\rm O}_{24})$. \ps

\begin{lemme} \label{satakecoefcar} Soient $G$ un $\Z_p$-groupe semi-simple
d\'eploy\'e, $\lambda$ un poids dominant de $\widehat{G}$, $V_\lambda$ la 
repr\'esentation irr\'eductible de $\widehat{G}$ associ\'ee, et $\rho$ la 
demi-somme des racines positives de $G_\C$.  Pour tout entier $m \geq 1$, il
existe un unique \'el\'ement $T \in {\rm H}(G)$ tel que pour tout $c \in
\widehat{G}_{\rm ss}$, on ait $$p^{m \langle \lambda, \rho \rangle}\, {\rm
trace}(c\,  | \, \Lambda^m V_\lambda) = {\rm tr}(c) ( {\rm Sat}(T) ).$$ De  
plus, si $G={\rm SO}_n$ et si la repr\'esentation $V_\lambda$ s'\'etend \`a
${\rm O}_n(\C)$, par exemple si $V={\rm St}$, alors $T$ appartient au sous-anneau ${\rm H}({\rm O}_n)   
\subset {\rm H}({\rm SO}_n)$.  \end{lemme}

\begin{pf} L'existence et unicit\'e d'un \'el\'ement $T \in {\rm H}(G)  
\otimes \Z[p^{-\frac{1}{2}}]$ tel que ${\rm Sat}(T) = p^{m \langle \lambda,
\rho \rangle} [\Lambda^m V_\lambda]$ est cons\'equence imm\'ediate de
l'isomorphisme de Satake (\S VI.\ref{isomsatake}).  Il s'agit donc de voir
que $T \in {\rm H}(G)$.  Rappelons que l'on a $$p^{\langle
\rho, \mu \rangle} [V_{\mu}] \in {\rm Sat}({\rm H}(G))$$ pour tout
poids dominant $\mu$ de $\widehat{G}$, d'apr\`es la formule 
VI.\eqref{idsatakegen}. Le
cas $m=1$ suit imm\'ediatement.  Si $m\geq 1$ est quelconque, on observe que
les constituants irr\'eductibles de $\Lambda^m V_\lambda$, qui sont de la  
forme $V_\mu$ o\`u $\mu$ est un poids dominant de $\widehat{G}$, v\'erifient
$\mu \leq m \lambda$ (\S VI.\ref{repalg}).  En effet, cette in\'egalit\'e
vaut plus g\'en\'eralement pour tous les poids $\mu$ de $V_\lambda^{\otimes
m}$, et $\Lambda^m \,V_\lambda$ est un quotient de ce dernier.  On conclut
la premi\`ere assertion du lemme car pour un tel $\mu$, $\langle
\rho, m \lambda - \mu \rangle$ est un entier $\geq 0$.  \ps Pour v\'erifier
la seconde assertion, concernant $G={\rm SO}_n$, il d'observer que si     
$V_\lambda$ se prolonge \`a ${\rm O}_n(\C) \supset \widehat{G}={\rm
SO}_n(\C)$, alors il en va de m\^eme de $\Lambda^m V_\lambda$.  La relation
VI.\eqref{descsato} montre alors que l'\'el\'ement $T \in {\rm H}({\rm SO}_n)$
d\'efini ci-dessus appartient \`a ${\rm H}({\rm O}_n)$.  \end{pf}

\begin{remarque}\label{remcoeffpolcar}{\rm Il est en g\'en\'eral difficile de d\'eterminer
explicitement l'op\'erateur $T$ donn\'e par le lemme ci-dessus, disons en terme de 
la $\Z$-base des
${\rm c}_\mu$ (\S VI.\ref{deuxbases}), et ce m\^eme dans le cas particulier du groupe 
$G={\rm SO}_n$ et de la repr\'esentation standard $V_\lambda={\rm St}$ de
$\widehat{G}$. Dans ce cas, on a toutefois d\'ej\`a dit que $T={\rm
T}_p$ pour $m=1$, et l'on a de plus $T= p\,{\rm T}_{p,p} +
p^{n/2-1}+\sum_{i=0}^{n/2-2} p^{2i+1}$ 
pour $m=2$ (Formule VI.\eqref{formulesattpp}).}\end{remarque}

Nous avons rappel\'e au \S \ref{repselfdual} l'existence et quelques
propri\'et\'es des repr\'esentations galoisiennes $\rho_{\pi,\iota}$
associ\'ees
\`a une repr\'esentation automorphe alg\'ebrique, autoduale, r\'eguli\`ere $\pi \in  
\Pi_{\rm cusp}(\PGL_n)$ et \`a un plongement $\iota : \overline{\Q}
\rightarrow \overline{\Q}_\ell$, o\`u $\overline{\Q}_\ell$ d\'esigne une
cl\^oture alg\'ebrique de $\Q_\ell$. La repr\'esentation galoisienne $\rho_{\pi,\iota}$ est
continue, semi-simple, non ramifi\'ee hors de $\ell$, et sa classe
d'isomorphisme est uniquement d\'etermin\'ee par la relation \begin{equation}\label{relgaloisladique}
\DET(t -\rho_{\pi,\iota}({\rm Frob}_p)) = \iota( \DET(t - {\rm
c}_p(\pi)p^{\frac{{\rm w}(\pi)}{2}}))\end{equation}
pour tout premier $p \neq \ell$ (l'assertion d'unicit\'e \'etant cons\'equence du th\'eor\`eme de
Cebotarev). Si le polyn\^ome $\DET(t - {\rm c}_p(\pi) p^{\frac{{\rm w}(\pi)}{2}})$ est dans
$\Q[t]$ pour tout premier
$p$, la repr\'esentation galoisienne $\rho_{\pi,\iota}$ d\'epend uniquement de $\ell$,
et non du choix de 
$\iota$ ; elle sera simplement not\'ee $\rho_{\pi,\ell}$. \ps

Les repr\'esentations automorphes $\pi$ d'int\'er\^et ici seront les $\Delta_w$ et les 
$\Delta_{w,v}$, et dans ces cas   
les repr\'esentations galoisiennes $\rho_{\pi,\iota}$ (satisfaisant les conditions (i) et (ii) du
th\'eor\`eme VIII.\ref{existencegalois}) ont \'et\'e respectivement
construites par Deligne \cite{deligne} (g\'en\'eralisant une construction
ant\'erieure d'Eichler, Shimura, Kuga et Sato \cite{kugashimura}) et Weissauer \cite{weissauergalois}
(voir aussi les travaux ant\'erieurs de Chai-Faltings \cite{chaifaltings} et Taylor
\cite{taylor}).\ps

Lorsque $\pi$ est engendr\'ee par une forme modulaire $f = q \,+ \,a_2\, q^2 + ... 
$ dans ${\rm S}_k({\rm SL}_2(\Z))$ qui est propre pour les op\'erateurs de
Hecke, la relation \eqref{relgaloisladique} se r\'eduit \`a la relation bien
connue $\DET(t-\rho_{\pi,\iota}({\rm Frob}_p))= t^2 - \iota(a_p) t + p^{k-1}$.  Ce
polyn\^ome est \`a coefficients entiers si $k \leq 22$ (car la forme $f$
l'est) ; il est par ailleurs bien connu que $\rho_{\pi,\ell}$ peut dans ce
cas \^etre choisie \`a coefficients dans le corps $\Q_\ell$. \ps

Lorsque $\pi$
est l'une des $4$ repr\'esentations $\Delta_{w,v}$ d\'efinie au \S
\ref{prelimsp4}, rappelons que le membre de droite de l'\'egalit\'e
\eqref{relgaloisladique} s'\'ecrit aussi 
$$t^4 - \tau_{j,k}(p) \,t^3 +
\frac{\tau_{j,k}(p)^2-\tau_{j,k}(p^2)}{2} \,t^2 - \tau_{j,k}(p)\, p^{j+2k-3}
\,t \,+\, p^{2j+4k-6},$$ o\`u $(j,k)=(v-1,(w-v)/2+2)$ (Formule IX.\eqref{polcarsp4}).  L\`a
encore, ce polyn\^ome est \`a coefficients rationnels, et m\^eme entiers d'apr\`es la proposition IX.\ref{intcongtaujk}. \ps


\begin{thmv}\label{thmgal24} Soient $i=1,\dots,24$, $\ell$ un nombre premier
et $\overline{\Q}_\ell$ une cl\^oture alg\'ebrique de $\Q_\ell$.  Il existe
une repr\'esentation continue et semi-simple $\rho_{i,\ell} : {\rm
Gal}(\overline{\Q}/\Q) \longrightarrow {\rm GL}_{24}(\overline{\Q}_\ell)$,
unique \`a isomorphisme pr\`es, qui est non ramifi\'ee hors de $\ell$ et
telle que pour tout premier $p \neq \ell$ on ait l'\'egalit\'e dans $\Z[t]$
\begin{equation} \label{carrhoil} \DET ( t \, - \rho_{i,\ell}({\rm Frob}_p))
= \DET ( t - p^{11}\, {{\rm St}( \rm c}_p(\pi_i)) ).\end{equation} \noindent En
particulier, on a $\lambda_i(p) = {\rm trace} \, \rho_{i,\ell}({\rm
Frob}_p)$ pour tout premier $p \neq \ell$.  \end{thmv}

\begin{pf} Soit $i=1,\dots,24$. Le th\'eor\`eme IX.\ref{thm24} affirme que
$\psi(\pi_i,{\rm St})$ est de la forme $\oplus_{j=1}^k \varpi_j[d_j]$, o\`u
les $\varpi_j$ sont parmi les repr\'esentations automorphes $1$, $\Delta_w$,
$\Delta_{w,v}$, et ${\rm Sym}^2 \Delta_{11}$. L'existence de
$\rho_{i,\ell}$ en d\'ecoule en posant 
\begin{equation}\label{decomprhoil} \rho_{i,\ell} =
\bigoplus_{j=1}^k \rho_{\varpi_j,\ell} \otimes (\oplus_{m=0}^{d_j-1}
\omega_\ell^m) \otimes \omega_\ell^{\frac{22-{\rm
w}(\varpi_j)+1-d_j}{2}},\end{equation}
la notation $\rho_{\varpi_j,\ell}$ \'etant introduite quelques lignes plus
haut. L'assertion d'unicit\'e est cons\'equence du th\'eor\`eme de Cebotarev. Bien entendu, lorsque
$\varpi_j$ est la repr\'esentation automorphe triviale dans $\Pi_{\rm cusp}(\PGL_1)$, 
$\rho_{\varpi_j,\ell}$ d\'esigne la repr\'esentation galoisienne triviale (de dimension
$1$). De plus, on peut prendre pour 
$\rho_{{\rm Sym}^2 \Delta_{11},\ell}$ la repr\'esentation ${\rm Sym}^2
\rho_{\Delta_{11},\ell}$.
\end{pf}

\begin{remarque}\label{remcoeffrhoil} {\rm En suivant la construction de
Weissauer \cite{weissauergalois}, il devrait \^etre possible de montrer que les repr\'esentations
$\rho_{\Delta_{w,v},\ell}$, ainsi donc que les $\rho_{i,\ell}$, sont
d\'efinies sur $\Q_\ell$, puisque $\dim {\rm S}_{j,k}=1$ pour les 
$4$ couples $(j,k)$ correspondants (\S IX.\ref{dimtsushima}).} \end{remarque}

Il serait int\'eressant d'\'etudier en d\'etail les images des
repr\'esentations galoisiennes $\rho_{\Delta_{w,v},\ell}$, \`a la mani\`ere des travaux
de Serre et Swinnerton-Dyer sur les repr\'esentations $\rho_{\Delta_w,\ell}$
\cite{swinnertondyer}. Nous nous contenterons au \S \ref{parfinharder} de d\'emontrer quelques congruences
satisfaites par ces repr\'esentations, dans l'esprit de la congruence de
Ramanujan. 

\begin{corv}\label{thmgal24red} Soient $i=1,\dots,24$ et $\ell$ un nombre
premier. Il existe une repr\'esentation continue et semi-simple $\rhob_{i,\ell} : {\rm
Gal}(\overline{\Q}/\Q) \longrightarrow {\rm GL}_{24}(\F_\ell)$,
unique \`a isomorphisme pr\`es, qui est non ramifi\'ee hors de $\ell$ et  
telle que pour tout premier $p \neq \ell$ on ait la congruence
\begin{equation} \label{carrhobil} \DET ( t \, - \rhob_{i,\ell}({\rm Frob}_p))
 \equiv \DET ( t - p^{11}\, {{\rm St}( \rm c}_p(\pi_i)) ) \bmod \ell.\end{equation} \noindent
En
particulier, on a $\lambda_i(p) \equiv {\rm trace} \, \rhob_{i,\ell}({\rm
Frob}_p) \bmod \ell$ pour tout premier $p \neq \ell$.  \end{corv}

L'assertion d'unicit\'e est cons\'equence du th\'eor\`eme de Cebotarev et d'un r\'esultat classique de
Brauer-Nesbitt\footnote{Rappelons l'\'enonc\'e de ce dernier. Soient $G$ un
groupe, $k$ un corps, ainsi que $\rho_1,\rho_2 : G \rightarrow \GL_m(k)$ deux
repr\'esentations semi-simples. Les repr\'esentations $\rho_1$ et $\rho_2$
sont isomorphes si, et seulement si, $\DET (t- \rho_1(g))=\DET(t-\rho_2(g))$
pour tout $g \in G$.}. L'existence de $\rhob_{i,\ell}$ se d\'eduit de celle de $\rho_{i,\ell}$ par
un proc\'ed\'e g\'en\'eral standard, rappel\'e ci-dessous. \ps

Fixons une cl\^oture alg\'ebrique
$\overline{\Q}_\ell$ (resp.  $\overline{\mathbb{F}}_\ell$) de
$\mathbb{Q}_\ell$ (resp.  $\mathbb{F}_\ell$), ainsi qu'un morphisme   
d'anneaux $\mathcal{O} \rightarrow \overline{\mathbb{F}}_\ell$, o\`u $\mathcal{O}
\subset \overline{\Q}_\ell$ d\'esigne la cl\^oture int\'egrale de $\Z_\ell$.
Soient $G$ un groupe profini et $\rho : G \rightarrow 
\GL_n(\overline{\Q}_\ell)$ une repr\'esentation continue.  On d\'emontre que :
\ps
\begin{itemize}
\item[(i)] $\DET(t-\rho(g)) \in \mathcal{O}[t]$ pour tout $g \in G$, \ps
\item[(ii)] il existe une repr\'esentation continue semi-simple $\rhob : G \rightarrow
\GL_n(\overline{\F}_\ell)$, unique \`a isomorphisme pr\`es, telle que pour tout $g \in G$, le polyn\^ome caract\'eristique   
$\DET(t - \rhob(g))$ est l'image de $\DET(t-\rho(g)) \in \mathcal{O}[t]$
dans $\overline{\F}_\ell[t]$. \ps
\item[(iii)] si de plus $\DET(t-\rho(g)) \in \Q_\ell[t]$ pour tout $g \in G$,
alors $\rhob$ peut \^etre choisie \`a coefficients dans $\mathbb{F}_\ell$,
et sa classe d'isomorphisme ne d\'epend pas du choix du morphisme $\mathcal{O} \rightarrow \overline{\mathbb{F}}_\ell$. \ps
\end{itemize} 
\noindent On dit alors que $\rhob$ est ``la'' repr\'esentation r\'esiduelle de $\rho$.
Cette construction (le point (iii) inclus) s'applique aux 
$\rho_{i,\ell}$, ainsi qu'aux repr\'esentations de la forme 
$\rho_{\pi,\ell}$ introduites plus haut ; elle donne lieu \`a des
repr\'esentations r\'esiduelles $\rhob_{i,\ell}$ et $\rhob_{\pi,\ell}$ \`a
coefficients dans $\mathbb{\F}_\ell$. Bien entendu, la relation 
\eqref{decomprhoil} induit une d\'ecomposition
similaire
\begin{equation}\label{decomprhobil}
 \rhob_{i,\ell} \simeq
\bigoplus_{j=1}^k \rhob_{\varpi_j,\ell} \otimes (\oplus_{m=0}^{d_j-1}
\overline{\omega_\ell}^{\, m}) \otimes \overline{\omega_\ell}^{\, \, \frac{22-{\rm
w}(\varpi_j)+1-d_j}{2}}.
\end{equation}\ps

Indiquons bri\`evement comment l'on d\'emontre les points (i), (ii)
et (iii). Une application classique du lemme de Baire assure d'abord que
$\rho(G) \subset \GL_n(F)$ o\`u $F \subset \overline{\Q}_\ell$ est une
extension finie de $\Q_\ell$.  Quitte \`a conjuguer $\rho$ par un
\'el\'ement de $\GL_n(F)$, la compacit\'e de $\rho(G)$ permet de supposer
que $\rho(G) \subset \GL_n(\mathcal{O}_F)$ o\`u $\mathcal{O}_F=\mathcal{O}
\cap F$.  Le (i) suit.  En composant
par le morphisme d'anneaux $\mathcal{O} \rightarrow
\overline{\mathbb{F}}_\ell$, il en r\'esulte une repr\'esentation continue
$G \rightarrow \GL_n(\overline{\mathbb{F}}_\ell)$.  Pour d\'emontrer le
point (ii) il suffit de prendre pour $\rhob$ une semi-simplifi\'ee de cette
derni\`ere repr\'esentation. L'assertion d'unicit\'e, ainsi que la
derni\`ere assertion du (iii), d\'ecoulent du r\'esultat de Brauer-Nesbitt. 
L'obstruction de Schur \'etant triviale sur les corps finis, on en d\'eduit
la premi\`ere assertion du (iii) (voir \cite[Lemma 2]{taylor}).

\begin{remarque}\label{pairingsp4} {\rm Soit $\rho=\rho_{\Delta_{w,v},\ell}$.
Nous avons d\'ej\`a observ\'e au \S VIII.\ref{repgal}, formule VIII.\eqref{rhopilautodual}, que l'on dispose d'un isomorphisme $\rho^\ast
\simeq \rho \otimes \omega_\ell^{-w}$. Mieux, d'apr\`es \cite[Corollary
1.3]{bchsign}, il existe un accouplement non
d\'eg\'en\'er\'e, altern\'e, et ${\rm Gal}(\overline{\Q}/\Q)$-\'equivariant
$\rho \otimes \rho \rightarrow
\omega_\ell^{w}$ (voir aussi  \cite{weissauergalois}). Il n'est pas difficile d'en d\'eduire l'existence d'un accouplement non
d\'eg\'en\'er\'e, altern\'e, et ${\rm
Gal}(\overline{\Q}/\Q)$-\'equivariant, $\rhob
\otimes \rhob \rightarrow
\overline{\omega_\ell}^{\,w}$. }
\end{remarque}

\parindent=0cm

\section{Retour sur les $p$-voisins des réseaux de Niemeier} \label{parfinniemeier}

\medskip
On numérote $\mathrm{L}_{1},\mathrm{L}_{2},\ldots,\mathrm{L}_{24}$ les $24$ réseaux de Niemeier (ou plutôt leurs classes d'isomorphisme) suivant la convention de Conway et Sloane \cite[Chap. 16, Table 16.1]{conwaysloane} (les entiers $1,2,\ldots,24$ remplaçant les lettres grecques $\alpha,\beta,\ldots,\omega$). On a donc $\mathrm{R}(\mathrm{L}_{i})=\mathbf{R}_{i}$ pour $i\leq 23$ (voir II.3) et $\mathrm{L}_{24}$ est le réseau de Leech que nous avons aussi noté $\mathrm{Leech}$ en III.4. 

\bigskip
Soit $p$ un nombre premier~; nous notons encore $\mathrm{T}_{p}$ la matrice de l'opérateur de Hecke $\mathrm{T}_{p}:\mathbb{Z}[\mathrm{X}_{24}]\to\mathbb{Z}[\mathrm{X}_{24}]$  dans la base $(\mathrm{L}_{1},\mathrm{L}_{2},\ldots,\mathrm{L}_{24})$.

\bigskip
Comme nous l'avons déjà dit Nebe et Venkov ont déterminé $\mathrm{T}_{2}$ et constaté que les valeurs propres de cet opérateur sont entières et distinctes. On les note
$$
\hspace{24pt}
\lambda_{1}>\lambda_{2}>\ldots>\lambda_{24}
\hspace{24pt}.
$$
On note $\mathrm{v}_{j}$, $1\leq j\leq 24$, un vecteur propre associé à $\lambda_{j}$ dont les coordonnées, dans la base évoquée plus haut, sont entières et premières entre elles (un tel vecteur propre est déterminé au signe près) ; on note $\mathrm{V}$ la matrice à $24$ lignes et colonnes dont la $j$-ième colonne est la colonne des cordonnées de $\mathrm{v}_{j}$.

\bigskip
On note $\lambda_{j}(p)$ l'entier défini par l'égalité $\mathrm{T}_{p}\hspace{1pt}\mathrm{v}_{j}=\lambda_{j}(p)\hspace{1pt}\mathrm{v}_{j}$ (se rappeler que $\mathrm{T}_{2}$ et $\mathrm{T}_{p}$ commutent)~; on a donc par définition $\lambda_{j}(2)=\lambda_{j}$. 

\bigskip
On pose\hspace{4pt}$\theta_{1}(p)=\tau_{12}(p)(=\tau(p))\hspace{2pt},\hspace{2pt}\theta_{2}(p)=\tau_{16}(p)\hspace{2pt},\hspace{2pt}\theta_{3}(p)=\tau_{18}(p)\hspace{2pt},\hspace{2pt}\theta_{4}(p)=\tau_{20}(p)\hspace{2pt},\hspace{2pt}\theta_{5}(p)=\tau_{22}(p)$\hspace{4pt}; en clair, on note $\theta_{r}(p)$ le $p$-ième coefficient de Fourier de la forme modulaire (pour $\mathrm{SL}_{2}(\mathbb{Z})$) parabolique normalisée res\-pectivement de poids $12,16,18,20,22$ pour $r=1,2,3,4,5$. On pose $\theta_{6}(p)=(\theta_{1}(p))^{2}-p^{11}$ ($\theta_{6}(p)=p^{11}\hspace{1pt}\mathrm{tr}(\mathrm{Sym}^{2}\mathrm{c}_{p}(\Delta_{11}))$).

\bigskip
On pose enfin $\theta_{7}(p)=\tau_{6,8}(p)\hspace{2pt},\hspace{2pt}\theta_{8}(p)=\tau_{8,8}(p)\hspace{2pt},\hspace{2pt}\theta_{9}(p)=\tau_{12,6}(p)\hspace{2pt},\hspace{2pt}\theta_{10}(p)=\tau_{4,10}(p)$.

\bigskip
D'après le théorème \ref{thm24} du chapitre IX et la formule \eqref{formulesattp} du chapitre VI, il existe des polynômes $\mathrm{C}_{j,r}$ de $\mathbb{Z}[X]$, $1\leq j\leq 24$, $0\leq r\leq 10$, uniquement déterminés, tels que l'on a
\begin{equation}\label{ljpformel}
\lambda_{j}(p)
\hspace{4pt}=\hspace{4pt}
\mathrm{C}_{j,0}(p)+\sum_{r=1}^{10}\hspace{4pt}\mathrm{C}_{j,r}(p)\hspace{2pt}\theta_{r}(p)
\end{equation}
pour tout nombre premier $p$.

\medskip
Rappelons la valeur de quelques-uns des polynômes $\mathrm{C}_{j,r}$.

\medskip
On a $\mathrm{C}_{1,0}=\sum_{k=0}^{k=22}X^{k}+X^{11}$ et $\mathrm{C}_{1,r}=0$ pour $r\geq 1$, en clair on a  $\lambda_{1}(p)=\mathrm{c}_{24}(p):=\sum_{k=0}^{k=22}p^{k}+p^{11}$ (voir III.2.4 et III.2.1).

\medskip
On a $\mathrm{C}_{2,0}=\sum_{k=1}^{k=21}X^{k}$, $\mathrm{C}_{2,6}=1$ et $\mathrm{C}_{2,r}=0$ pour $r\not=0,6$.

\medskip
Pour $r\geq 7$ la valeur des polynômes $\mathrm{C}_{j,r}$ est la suivante:

\smallskip
-- $\mathrm{C}_{j,7}=0$ pour $j\not=19$ et $\mathrm{C}_{19,7}=X(X+1)$~;

\smallskip
-- $\mathrm{C}_{j,8}=0$ pour $j\not=15$ et $\mathrm{C}_{15,8}=X+1$~;

\smallskip
-- $\mathrm{C}_{j,9}=0$ pour $j\not=10$ et $\mathrm{C}_{10,9}=X+1$~;

\smallskip
-- $\mathrm{C}_{j,10}=0$ pour $j\not=21$ et $\mathrm{C}_{21,10}=X+1$.

\bigskip
En contemplant la formule
\begin{equation}\label{tpformel}
\mathrm{T}_{p}
\hspace{4pt}=\hspace{4pt}
\mathrm{V}\hspace{4pt}\mathrm{diag}(\lambda_{1}(p),\lambda_{2}(p),\ldots,\lambda_{24}(p))\hspace{4pt}\mathrm{V}^{-1}
\end{equation}
(la notation $\mathrm{diag}(\lambda_{1}(p),\lambda_{2}(p),\ldots,\lambda_{24}(p))$ désignant la matrice diagonale dont les coefficients diagonaux sont les $\lambda_{j}(p)$) on obtient l'énoncé suivant :

\medskip
\begin{thm}\label{npformel} Soient $L$ et $L'$ deux réseaux unimodulaires pairs de dimension $24$. Il existe des polynômes $\mathrm{P}_{r}(L,L';X)$ de $\mathbb{Q}[X]$, $0\leq r\leq 10$, uniquement déterminés en fonction des classes d'isomorphismes de $L$ et $L'$, tels que l'on a
$$
\mathrm{N}_{p}(L,L')
\hspace{4pt}=\hspace{4pt}
\mathrm{P}_{0}(L,L';p)+\sum_{r=1}^{10}\hspace{4pt}\mathrm{P}_{r}(L,L';p)\hspace{2pt}\theta_{r}(p)
$$
pour tout nombre premier $p$.
\end{thm}

\bigskip
\textit{Remarque.} Par définition, on a pour tout $r$ avec $0\leq r\leq 10$, l'égalité suivante de matrices, à $24$ lignes et colonnes, à coefficients dans $\mathbb{Q}[X]$~:
$$
\hspace{24pt}
[\hspace{2pt}\mathrm{P}_{r}(\mathrm{L}_{j},\mathrm{L}_{i};X)\hspace{2pt}]
\hspace{4pt}=\hspace{4pt}
\mathrm{V}\hspace{4pt}\mathrm{diag}(\mathrm{C}_{1,r}(X),\mathrm{C}_{2,r}(X),\ldots,\mathrm{C}_{24,r}(X))\hspace{4pt}\mathrm{V}^{-1}
\hspace{24pt};
$$
comme les colonnes de $\mathrm{V}$ sont deux à deux orthogonales pour le produit scalaire de matrice $\mathrm{diag}(\vert\mathrm{O}(\mathrm{L}_{1})\vert,\vert\mathrm{O}(\mathrm{L}_{2})\vert,\ldots,\vert\mathrm{O}(\mathrm{L}_{24})\vert)$ (Proposition III.2.3) cette égalité montre que l'on a
$$
\frac{1}{\vert\mathrm{O}(L)\vert}
\hspace{4pt}\mathrm{P}_{r}(L,L';X)
\hspace{4pt}=\hspace{4pt}
\frac{1}{\vert\mathrm{O}(L')\vert}
\hspace{4pt}\mathrm{P}_{r}(L',L;X)
$$
pour tout $r$, $L$ et $L'$ (ce qui est bien sûr compatible avec le scholie III.1.7).

\bigskip
On note $\mathrm{Proj}_{1}$ le projecteur orthogonal, pour le produit scalaire introduit en III.2.3, de $\mathbb{Q}[\mathrm{X}_{24}]$ sur la droite engendrée par $\mathrm{v}_{1}$. On note $\mathrm{w}$ le vecteur $\sum_{x\in\mathrm{X}_{24}}\frac{1}{\vert\mathrm{O}(x)\vert}\hspace{1pt}x$ de de $\mathbb{Q}[\mathrm{X}_{24}]$~; il résulte de III.2.4 que le vecteur $\mathrm{v}_{1}$ est colinéaire au vecteur $\mathrm{w}$. Soit $y$ un élément de $\mathrm{X}_{24}$~; les égalités $\mathrm{w}.y=1$ et $\mathrm{w}.\mathrm{w}=\sum_{x\in\mathrm{X}_{24}}\frac{1}{\vert\mathrm{O}(x)\vert}$ entraînent~:
$$
\mathrm{Proj_{1}}\hspace{1pt}(y)
\hspace{4pt}=\hspace{4pt}
\frac{\mathrm{w}}
{\sum_{x\in\mathrm{X}_{24}}\frac{1}{\vert\mathrm{O}(x)\vert}}
\hspace{4pt}=\hspace{4pt}
\sum_{x\in\mathrm{X}_{24}}\mu(x)\hspace{2pt}x
$$
en posant
$$
\hspace{24pt}
\mu(x)
\hspace{4pt}=\hspace{4pt}
\frac{\frac{1}{\vert\mathrm{O}(x)\vert}}
{\sum_{x\in\mathrm{X}_{24}}\frac{1}{\vert\mathrm{O}(x)\vert}}
\hspace{24pt}.
$$
On observera que $\mu(x)$ est le quotient de la masse de $x$ par la masse du genre des réseaux unimodulaires pairs de dimension $24$, masses au sens de Minkowski-Siegel~: $\sum_{x\in\mathrm{X}_{24}}\mu(x)\hspace{1pt}\delta_{x}$ est {\em la mesure de probalité de Minkowski-Siegel} sur l'ensemble $\mathrm{X}_{24}$.

\medskip
\begin{thm}\label{tpapprox}
Lorsque le nombre premier $p$ tend vers l'infini, on a
$$
\mathrm{T}_{p}
\hspace{4pt}=\hspace{4pt}
p^{22}\hspace{4pt}\mathrm{Proj}_{1}
\hspace{4pt}+\hspace{4pt}
\mathrm{O}\hspace{1pt}(p^{21})
$$
(la notation $\mathrm{O}(- )$ qui apparaît ci-dessus est la notation de Landau\ldots et n'a pas de rapport avec les groupes orthogonaux~!). 
\end{thm}

\medskip
\textit{Démonstration.} Compte tenu de III.2.3, on a $\mathrm{T}_{p}=\sum_{j=1}^{24}\lambda_{j}(p)\mathrm{Proj}_{j}$, $\mathrm{Proj}_{j}$ étant le projecteur orthogonal de $\mathbb{Q}[\mathrm{X}_{24}]$ sur la droite engendrée par $\mathrm{v}_{j}$. Or on a $\lambda_{1}(p)=p^{22}+\mathrm{O}(p^{21})$ et les inégalités de Ramanujan pour les $\theta_{r}$  impliquent $\lambda_{j}(p)=\mathrm{O}(p^{21})$ pour $j\geq 2$. Rappelons ce que sont ces inégalités. On a $\vert\tau_{k}(p)\vert\leq 2\hspace{1pt}p^{\frac{k-1}{2}}$
pour $k=12,16,18,20,22$, $\vert\tau(p)^{2}-p^{11}\vert\leq 3\hspace{1pt}p^{11}$, $\vert\tau_{6,8}(p)\vert\leq 4\hspace{1pt}p^{\frac{19}{2}}$ et $\tau_{j,k}(p)\vert\leq 4\hspace{1pt}p^{\frac{21}{2}}$ pour $(j,k)=(8,8),(12,6),(4,10)$. Les 5 premières inégalités sont dues à Deligne \cite{deligne}, les 4 dernières à Weissauer \cite{weissauergalois} (la 6-ième est conséquence de la première~!) ; pour une discussion générale concernant ce type d'inégalité voir \S VII.\ref{repgal}.
\hfill$\square$

\medskip
\begin{scholie}\label{npapprox}
Soient $x$ et $y$ deux éléments de $\mathrm{X}_{24}$~; on a
$$
\frac{\mathrm{N}_{p}(x,y)}{\mathrm{c}_{24}(p)}
\hspace{4pt}=\hspace{4pt}
\mu(y)
\hspace{4pt}+\hspace{4pt}
\mathrm{O}\hspace{1pt}(\frac{1}{p})
$$
lorsque $p$ tend vers l'infini.
\end{scholie}

\medskip
\textit{Commentaires.} Soient $L$ un réseau unimodulaire pair de dimension $24$ et $y$ un élément de $\mathrm{X}_{24}$. Le quotient
$$
\frac{\mathrm{N}_{p}([L],y)}{\mathrm{c}_{24}(p)}
\hspace{4pt}=\hspace{4pt}
\frac{\mathrm{N}_{p}([L],y)}{\vert\mathrm{C}_{L}(\mathbb{F}_{p})\vert}
$$
est la proportion des points $c$ de la quadrique $\mathrm{C}_{L}(\mathbb{F}_{p})$ tels que la classe d'isomorphisme du $p$-voisin de $L$ associé à $c$ (voir II.1.5) est $y$. Le scholie \ref{npapprox} dit que cette proportion tend (en $p^{-1}$) vers $\mu(y)$ quand $p$ tend vers l'infini.

\bigskip
\textit{Remarque.} La mesure de probabilité de Minkowski-Siegel sur $\mathrm{X}_{24}$ est très loin d'être uniforme. On a par exemple $\mu([\mathrm{E}_{24}])\approx 2.42\times10^{-17}$ (c'est le minimum de la fonction $\mu$ sur  $\mathrm{X}_{24}$), $\mu([\mathrm{L}_{21}])\approx 0.426$ (c'est le maximum) et $\mu([\mathrm{L}_{20}])+\mu([\mathrm{L}_{21}])+\mu([\mathrm{L}_{22}])\approx 0.906$.

\bigskip
\textit{Remarque.} On constate en inspectant les égalités \eqref{ljpformel} que les inégalités de Ramanujan donnent $\lambda_{j}(p)=\mathrm{C}_{0,j}(p)+\mathrm{O}(p^{\frac{33}{2}})$ pour tout $j$. On en déduit comme précédemment des estimations de la forme
$$
\hspace{24pt}
\frac{\mathrm{N}_{p}(x,y)}{\mathrm{c}_{24}(p)}
\hspace{4pt}=\hspace{4pt}
\mu(y)\hspace{4pt}(\hspace{4pt}
1+
\sum_{n=1}^{5}\hspace{4pt}\kappa_{n}(x,y)\hspace{2pt}\frac{1}{p^n}
\hspace{4pt}+\hspace{4pt}
\mathrm{O}\hspace{1pt}(p^{-\frac{11}{2}})\hspace{4pt})
\hspace{24pt};
$$
le lien entre le vecteur propre $\mathrm{v}_{2}$ et les séries thêta de genre $1$ des réseaux unimodulaires pair de dimension $24$ permet d'expliciter $\kappa_{1}(x,y)$~:
$$
\hspace{24pt}
\kappa_{1}(x,y)\hspace{4pt}=\hspace{4pt}\frac{37092156523}{34673184000}\hspace{6pt}(\hspace{2pt}\mathrm{h}(x)-\frac{2730}{691}\hspace{2pt})\hspace{4pt}(\hspace{2pt}\mathrm{h}(y)-\frac{2730}{691}\hspace{2pt})
\hspace{24pt}.
$$
\footnotesize
Donnons quelques précisions sur le lien évoqué ci-dessus entre $\mathrm{v}_{2}$ et séries thêta de genre~$1$~; ces précisions sont en fait un développement du deuxième commentaire qui suit la table~\ref{tablevp2}. D'une part le théorème V.\ref{appendiceCthm} montre que le diagramme
$$
\begin{CD}
\mathbb{C}[\mathrm{X}_{24}]@>\vartheta_{1}>>\mathrm{M}_{12}(\mathrm{SL}_{2}(\mathbb{Z})) \\
@V\mathrm{T}_{p}VV @VVp\frac{p^{21}-1}{p-1}+\mathrm{T}(p^{2})V \\
\mathbb{C}[\mathrm{X}_{24}]@>\vartheta_{1}>>\mathrm{M}_{12}(\mathrm{SL}_{2}(\mathbb{Z}))
\end{CD}
$$
est commutatif (relation d'Eichler). D'autre part $\{\mathbb{E}_{12},\Delta\}$ est une base du $\mathbb{C}$-espace vectoriel $\mathrm{M}_{12}(\mathrm{SL}_{2}(\mathbb{Z}))$, propre pour les opérateurs de Hecke~; on a en particulier $\mathrm{T}(p^{2})\hspace{1pt}(\Delta)=\tau(p^{2})\hspace{1pt}\Delta=(\tau(p)^{2}-p^{11})\hspace{1pt}\Delta$. Soient $\mathrm{coord}_{\Delta}:\mathrm{M}_{12}(\mathrm{SL}_{2}(\mathbb{Z}))\to\mathbb{C}$ la forme linéaire coordonnée ``d'indice $\Delta$'' définie par cette base et $\eta:\mathbb{C}[\mathrm{X}_{24}]\to\mathbb{C}$ la forme linéaire composée $\mathrm{coord}_{\Delta}\circ\vartheta_{1}$~; il résulte de ce qui précède que $\eta$ est vecteur propre de l'endomorphisme $\mathrm{T}_{p}^{*}$ de $(\mathbb{C}[\mathrm{X}_{24}])^{*}$, pour la valeur propre $\sum_{k=1}^{21}p^{k}+(\tau(p)^{2}-p^{11})=\lambda_{2}(p)$. On vérifie facilement que l'on a
$$
\eta(x)
\hspace{4pt}=\hspace{4pt}
\vert\mathrm{R}(x)\vert-\frac{65520}{691}
\hspace{4pt}=\hspace{4pt}
24\hspace{4pt}(\hspace{2pt}\mathrm{h}(x)-\frac{2730}{691}\hspace{2pt})
$$
pour tout $x$ dans $\mathrm{X}_{24}$ (voir par exemple \cite[\S 6.6, formule (108)]{serre}). Comme en II.2 (Proposition 2.4), il en résulte que le vecteur
$$
\sum_{x\in\mathrm{X}_{24}}\hspace{2pt}\frac{1}{\vert\mathrm{O}(x)\vert}\hspace{2pt}(\hspace{2pt}\mathrm{h}(x)-\frac{2730}{691}\hspace{2pt})\hspace{4pt}x
$$
est vecteur propre de l'endomorphisme $\mathrm{T}_{p}$ de $\mathbb{C}[\mathrm{X}_{24}]$, pour la valeur propre $\lambda_{2}(p)$.
\normalsize

\vspace{0,75cm}
\textsc{Sur le diamètre du graphe des $p$-voisins en dimension $24$}

\bigskip
La formule du théorème \ref{npformel} montre que si l'on a $\mathrm{N}_{p}(L,L')=0$ alors on a
$$
\hspace{24pt}
\mathrm{P}_{0}(L,L';p)^{\hspace{2pt}2}
\hspace{4pt}=\hspace{4pt}
(\hspace{2pt}\sum_{r=1}^{10}\hspace{4pt}\mathrm{P}_{r}(L,L';p)\hspace{2pt}\theta_{r}(p)\hspace{2pt})^{2}
\hspace{24pt}.
$$
D'après l'inégalité de Schwarz on a
$$
(\hspace{2pt}\sum_{r=1}^{10}\hspace{4pt}\mathrm{P}_{r}(L,L';p)\hspace{2pt}\theta_{r}(p)\hspace{2pt})^{2}
\hspace{4pt}\leq\hspace{4pt}
(\hspace{2pt}\sum_{r=1}^{10}\hspace{4pt}\gamma_{k}\hspace{2pt}\mathrm{P}_{r}(L,L';p)^{2}\hspace{2pt})
\hspace{4pt}
(\hspace{2pt}\sum_{r=1}^{10}\hspace{4pt}\gamma_{k}^{-1}\hspace{2pt}\theta_{r}(p)^{2}\hspace{2pt})
$$
pour tout $10$-uple $(\gamma_{1},\gamma_{2},\ldots,\gamma_{10})$ de nombres réels strictement positifs. En prenant
$$
(\gamma_{1},\gamma_{2},\ldots,\gamma_{10})
=
(4\hspace{1pt}p^{11}, 4\hspace{1pt}p^{15}, 4\hspace{1pt}p^{17}, 4\hspace{1pt}p^{19}, 4\hspace{1pt}p^{21}, 9\hspace{1pt}p^{22}, 16\hspace{1pt}p^{19}, 16\hspace{1pt}p^{21}, 16\hspace{1pt}p^{21}, 16\hspace{1pt}p^{21})
$$
on obtient, compte tenu des inégalités de Ramanujan, l'inégalité
$$
\hspace{24pt}
(\hspace{2pt}\sum_{r=1}^{10}\hspace{4pt}\mathrm{P}_{r}(L,L';p)\hspace{2pt}\theta_{r}(p)\hspace{2pt})^{2}
\hspace{4pt}\leq\hspace{4pt}
10\hspace{4pt}
(\hspace{2pt}\sum_{r=1}^{10}\hspace{4pt}\gamma_{k}\hspace{2pt}\mathrm{P}_{r}(L,L';p)^{2}\hspace{2pt})
\hspace{24pt}.
$$
On pose
$$
(\Gamma_{1}(X),\Gamma_{2}(X),\ldots,\Gamma_{10}(X))
=
(4\hspace{1pt}X^{11}, 4\hspace{1pt}X^{15},\ldots,16\hspace{1pt}X^{21})
$$
et
$$
\hspace{24pt}
\mathrm{Q}(L,L';X)
\hspace{4pt}=\hspace{4pt}
\mathrm{P}_{0}(L,L';X)^{\hspace{2pt}2}
-10\hspace{4pt}(\hspace{2pt}\sum_{r=1}^{10}\hspace{4pt}\Gamma_{k}(X)\hspace{2pt}\mathrm{P}_{r}(L,L';X)^{2}\hspace{2pt})
\hspace{24pt};
$$
on constate que $\mathrm{Q}(L,L';X)$ est un polynôme de $\mathbb{Q}[X]$ dont le monôme de plus haut degré est $\mu(L')^{2}\hspace{2pt}X^{44}$. On observera que la remarque qui suit le théorème \ref{npformel} implique l'égalité $\mu(L)^{2}\hspace{2pt}\mathrm{Q}(L,L';X)=\mu(L')^{2}\hspace{2pt}\mathrm{Q}(L',L;X)$. On note $\rho(L,L')$ la plus grande racine réelle du polynôme $\mathrm{Q}(L,L';X)$ (on pourrait convenir que l'on a $\rho(L,L')=-\infty$ si $\mathrm{Q}(L,L';X)$ n'a pas de racine réelle mais ce polynôme a en fait toujours des racines réelles)~; on note enfin $\mathrm{p}(L,L')$ le plus petit nombre premier strictement plus grand que  $\rho(L,L')$. On a tout fait pour avoir $\mathrm{N}_{p}(L,L')>0$ pour $p\geq\mathrm{p}(L,L')$.

\bigskip
\textit{Exemple.} Le $24$-uple $(\hspace{1pt}\rho(\mathrm{L}_{i},\mathrm{Leech})\hspace{1pt})_{1\leq i\leq 24}$ est approximativement le suivant~:
\begin{multline*}
(\hspace{1pt}46.77\hspace{1pt},\hspace{1pt} 30.11\hspace{1pt},\hspace{1pt} 30.88\hspace{1pt},\hspace{1pt} 23.97\hspace{1pt},\hspace{1pt} 21.71\hspace{1pt},\hspace{1pt} 17.80\hspace{1pt},\hspace{1pt} 17.59\hspace{1pt},\hspace{1pt} 15.63\hspace{1pt},\hspace{1pt} 13.72\hspace{1pt},\hspace{1pt} 12.00\hspace{1pt},\hspace{1pt} 11.27\hspace{1pt},\\ \hspace{1pt} 12.14\hspace{1pt},\hspace{1pt} 9.36\hspace{1pt},\hspace{1pt} 9.58\hspace{1pt},\hspace{1pt} 8.48\hspace{1pt},\hspace{1pt} 7.03\hspace{1pt},\hspace{1pt} 6.19\hspace{1pt},\hspace{1pt} 5.21\hspace{1pt},\hspace{1pt} 5.86\hspace{1pt},\hspace{1pt} 4.12\hspace{1pt},\hspace{1pt} 3.10\hspace{1pt},\hspace{1pt} 2.13\hspace{1pt},\hspace{1pt} 1.37\hspace{1pt},\hspace{1pt} 1.68\hspace{1pt})\hspace{4pt}.
\end{multline*}
On constate par inspection que pour tout réseau de Niemeier $L$ le nombre premier $\mathrm{p}(L,\mathrm{Leech})$ est le plus petit nombre premier supérieur ou égal au nombre de Coxeter $\mathrm{h}(L)$. Cette inspection et la proposition III.4.1.1 montrent que l'énoncé ``$\mathrm{N}_{p}(L,L')>0$ pour $p\geq\mathrm{p}(L,L')$'' est optimal pour $L'=\mathrm{Leech}$.

\bigskip
On constate, également par inspection, que l'on a
$$
\rho(L,L')\hspace{4pt}\leq\hspace{4pt}\rho(\mathrm{E}_{24},\mathrm{Leech})\hspace{4pt}<47
$$
pour tous réseaux de Niemeier $L$ et $L'$. On voit donc que le  graphe des $p$-voisins en dimension $24$ est le graphe complet d'ensemble de sommets $\mathrm{X}_{24}$ pour $p\geq 47$. Comme nous avons calculé les $\tau_{j,k}(p)$ pour $p\leq 43$ (et même $p\leq 113$, voir le prochain paragraphe), nous sommes maintenant en mesure de déterminer le diamètre du graphe des $p$-voisins en dimension $24$~; voici le résultat~:

\medskip
\begin{thm}\label{diametre}
Soit $p$ un nombre premier. Le diamètre du graphe $\mathrm{K}_{24}(p)$ est le suivant~: $5$ pour $p=2$, $4$ pour $p=3$, $3$ pour $p=5$, $2$ pour $7\leq p\leq 43$\linebreak et $1$ pour $p\geq 47$.
\end{thm}
\medskip

\section{Détermination des $\tau_{j,k}(q)$ pour de petites valeurs de $q$}\label{calcfinvpgenre2}

\medskip
Les $\tau_{j,k}(q)$ en question, $q=p^{n}$ avec $p$ premier et $n\geq 1$ entier, sont définis en IX.\ref{definitiontaujk}~; on a vu que la détermination des $\tau_{j,k}(p^{n})$ pour $n>2$ se ramène à celle des $\tau_{j,k}(p^{n})$ pour $n=1,2$. Les valeurs des $\tau_{j,k}(p)$ pour $p\leq 113$ sont rassemblées dans la table \ref{taujk}, celles des $\tau_{j,k}(p^{2})$ pour $p\leq 29$ dans la table~\ref{taujkp2}.

\subsection{Détermination des $\tau_{j,k}(p)$ pour $p\leq 113$}\label{calcfinvpgenre2-1}

\medskip
Les entiers $\theta_{r}(p)$, $r\leq 6$, ne sont pas difficiles à calculer (au moins pour un nombre premier dont la taille n'est pas déraisonnable~!). Par contre, comme nous l'avons déjà dit, les tables pour les $\theta_{r}(p)$, $r\geq 7$, sont assez courtes. Nous nous proposons de montrer que la théorie développée dans ce mémoire et les informations que nous avons collectées en III.4 sur la dernière colonne de $\mathrm{T}_{p}$ (c'est-à-dire sur les nombres de voisins $\mathrm{N}_{p}(L,\mathrm{Leech})$ pour $L$ un réseau unimodulaire pair de dimension $24$ représentant $2$) permettent de déterminer ces $\theta_{r}(p)$ pour $p\leq 113$. Les formules \eqref{ljpformel} et \eqref{tpformel} (ou ce qui revient au même le théorème \ref{npformel}), et la détermination de $\tau_{6,8}(p)$, $\tau_{8,8}(p)$, $\tau_{12,6}(p)$ et $\tau_{4,10}(p)$ pour $p\leq 113$  permettent d'expliciter l'opérateur de Hecke $\mathrm{T}_{p}:\mathbb{Z}[\mathrm{X}_{24}]\to\mathbb{Z}[\mathrm{X}_{24}]$ pour $p\leq 113$.

\bigskip
On écrit $\mathrm{N}_{p}(L,L')=\mathrm{N}^{1}_{p}(L,L')+\mathrm{N}^{2}_{p}(L,L')$ avec
$$
\mathrm{N}^{1}_{p}(L,L')
\hspace{4pt}=\hspace{4pt}
\mathrm{P}_{0}(L,L';p)+\sum_{r=1}^{6}\hspace{4pt}\mathrm{P}_{r}(L,L';p)\hspace{2pt}\theta_{r}(p)
$$
(c'est le terme ``aisément calculable'') et
$$
\mathrm{N}^{2}_{p}(L,L')
\hspace{4pt}=\hspace{4pt}
\sum_{r=7}^{10}\hspace{4pt}\mathrm{P}_{r}(L,L';p)\hspace{2pt}\theta_{r}(p)
$$
(c'est le terme ``mystérieux''). Compte tenu de ce que l'on a rappelé plus haut sur les polynômes $\mathrm{C}_{j,r}$, $j\geq 1,r\geq 7$, qui interviennent dans l'expression des $\lambda_{j}(p)$, on a
\begin{equation}\label{npmysterieux}
\hspace{12pt}
\mathrm{N}^{2}_{p}(L,L')
\hspace{4pt}=\hspace{4pt}
\mathrm{c}_{7}(L,L')\hspace{2pt}p\hspace{1pt}(p+1)\hspace{2pt}\theta_{7}(p)+\sum_{r=8}^{10}\hspace{4pt}\mathrm{c}_{r}(L,L')\hspace{2pt}(p+1)\hspace{2pt}\theta_{r}(p)
\hspace{12pt},
\end{equation}
les $\mathrm{c}_{r}(L,L')$, $r\geq 7$, désignant des nombres rationnels déterminés en fonction des classes d'isomorphisme des réseaux $L$ et $L'$.

\bigskip
On voit donc que si l'on connaît les entiers $\mathrm{N}_{p}(L,L')$ pour quatre couples $(L,L')$ (dont les orbites sous l'action de $\mathfrak{S}_{2}$ par échange des facteurs sont distinctes) alors on peut espérer déterminer les $\theta_{r}(p)$ pour $r\geq 7$ en résolvant un système linéaire.

\bigskip
Ci-après, on prend pour $L'$ le réseau de Leech et pour $L$ les quatre réseaux de Niemeier avec racines dont le nombre de Coxeter est le plus grand~: $\mathrm{L}_{1}:=\mathrm{E}_{24}$ ($h=46$), $\mathrm{L}_{2}:=\mathrm{E}_{16}\oplus\mathrm{E}_{8}$ ($h=30$), $\mathrm{L}_{3}:=\mathrm{E}_{8}\oplus\mathrm{E}_{8}\oplus\mathrm{E}_{8}$ ($h=30$) et $\mathrm{L}_{4}:=\mathrm{A}_{24}^{+}$ ($h=25$).

\bigskip
Les quatre égalités $\mathrm{N}^{2}_{p}(\mathrm{L}_{i},\mathrm{Leech})=\sum_{r=7}^{10}\hspace{4pt}\mathrm{P}_{r}(\mathrm{L}_{i},\mathrm{Leech};p)\hspace{2pt}\theta_{r}(p)$, $i=1,2,3,4$, peuvent encore s'écrire
\begin{equation}\label{systeme}
\hspace{24pt}
\mathrm{A}(p)
\begin{bmatrix}
\theta_{7}(p) \\ \theta_{8}(p) \\  \theta_{9}(p) \\ \theta_{10}(p)
\end{bmatrix}
\hspace{4pt}=\hspace{4pt}\mathrm{B}(p)
\hspace{24pt},
\end{equation}
$\mathrm{A}(p)$ désignant la matrice $4\times 4$ produit de la matrice 
\begin{equation}\label{mata}
\mathrm{a}
\hspace{4pt}:=\hspace{4pt}
\begin{bmatrix}
\frac{-20360704}{297} & \frac{31085824}{1495} & \frac{210852224}{15795} & \frac{182174720}{9963} \\ \\
\frac{-1048320}{2057} & \frac{110194560}{116909} & \frac{901568}{2223} & \frac{-15608320}{77121} \\ \\
\frac{16329600}{22627} & \frac{16092820800}{16717987} & \frac{12615840}{35321} & \frac{27014400}{94259} \\ \\
\frac{994175}{4752} & \frac{-36575}{208} & \frac{37053115}{202176} & \frac{-4447625}{79704}
\end{bmatrix}
\end{equation}
et de la matrice diagonale $\mathrm{diag}(p(p+1),p+1,p+1,p+1)$, et  $\mathrm{B}(p)$ la matrice colonne
$$
\hspace{24pt}
\begin{bmatrix}
\mathrm{N}_{p}(\mathrm{E}_{24},\mathrm{Leech})-\mathrm{N}^{1}_{p}(\mathrm{E}_{24},\mathrm{Leech}) \\
\vspace{-10pt} \\
\mathrm{N}_{p}(\mathrm{E}_{16}\oplus\mathrm{E}_{8},\mathrm{Leech})-\mathrm{N}^{1}_{p}(\mathrm{E}_{16}\oplus\mathrm{E}_{8},\mathrm{Leech}) \\
\vspace{-10pt} \\
\mathrm{N}_{p}(\mathrm{E}_{8}\oplus\mathrm{E}_{8}\oplus\mathrm{E}_{8},\mathrm{Leech})-\mathrm{N}^{1}_{p}(\mathrm{E}_{8}\oplus\mathrm{E}_{8}\oplus\mathrm{E}_{8},\mathrm{Leech}) \\
\vspace{-10pt} \\
\mathrm{N}_{p}(\mathrm{A}^{+}_{24},\mathrm{Leech})
-\mathrm{N}^{1}_{p}(\mathrm{A}^{+}_{24},\mathrm{Leech})
\end{bmatrix}
\hspace{24pt}.
$$
On constate que l'on $\mathop{\mathrm{d\acute{e}t}}\mathrm{a}\not=0$ (merci \texttt{PARI})~; le système linéaire \eqref{systeme} permet donc de déterminer les $\theta_{r}(p)$, $r\geq 7$, si l'on connaît les entiers $\mathrm{N}_{p}(\mathrm{L}_{i},\mathrm{Leech})$ pour $i=1,2,3,4$.
\vfill\eject

\bigskip
\textit{Les cas $p\leq 23$}

\bigskip
Soit $L$ un réseau unimodulaire pair de dimension $24$ avec racines. La proposition III.4.1.1 dit que l'on a $\mathrm{N}_{p}(L,\mathrm{Leech})=0$ pour $p<\mathrm{h}(L)$. On a donc en particulier $\mathrm{N}_{p}(\mathrm{L}_{i},\mathrm{Leech})=0$ pour $i=1,2,3,4$~: le système linéaire \eqref{systeme} permet de déterminer les $\theta_{r}(p)$, $r\geq 7$, pour $p\leq 23$.

\bigskip
\textit{Les cas $p=29$ et $p=31$}

\bigskip
On a calculé $\mathrm{N}_{p}(\mathrm{A}^{+}_{24},\mathrm{Leech})$ en III.4.3 pour  $p=29$ et $p=31$. D'après III.4.1.1, on a $\mathrm{N}_{29}(L,\mathrm{Leech})=0$ pour $L\hspace{1pt}=\hspace{1pt}\mathrm{E}_{24}\hspace{1pt},\hspace{1pt}\mathrm{E}_{16}\oplus\mathrm{E}_{8}\hspace{1pt},\hspace{1pt}\mathrm{E}_{8}\oplus\mathrm{E}_{8}\oplus\mathrm{E}_{8}$ et $\mathrm{N}_{31}(L,\mathrm{Leech})=0$ pour $L\hspace{1pt}=\hspace{1pt}\mathrm{E}_{24}$~; d'autre part  le point (d) de III.4.2.10 donne la valeur de $\mathrm{N}_{31}(L,\mathrm{Leech})$ pour $L\hspace{1pt}=\hspace{1pt}\mathrm{E}_{16}\oplus\mathrm{E}_{8}\hspace{1pt},\hspace{1pt}\mathrm{E}_{8}\oplus\mathrm{E}_{8}\oplus\mathrm{E}_{8}$. Le système linéaire \eqref{systeme} permet encore de déterminer les $\theta_{r}(p)$, $r\geq 7$, pour $p=29$ et $p=31$.

\bigskip
\textit{Les cas $p=3$ et $7\leq p\leq 59$}

\bigskip
Le calcul de $\mathrm{N}_{31}(\mathrm{A}^{+}_{24},\mathrm{Leech})$, évoqué ci-dessus, bien qu'élémentaire est assez acrobatique. Nous proposons ci-après une méthode, loin d'être aussi élémentaire mais autrement plus efficace, pour déterminer les entiers $\mathrm{N}_{p}(L,\mathrm{Leech})$, $L$ désignant un réseau de Niemeier avec racines, lorsque $p$ n'est ``pas trop grand en fonction de $L$'', au moins si $p$ ne divise pas l'indice, noté $\mathrm{g}(L)$ en III.4, du sous-module de $L$ engendré par ses racines. Cette méthode est fondée sur les deux observations suivantes~:

\medskip
-- Soit $p$ un nombre premier qui ne divise pas $\mathrm{g}(L)$. L'énoncé III.4.3.3 dit que $\mathrm{N}_{p}(L,\mathrm{Leech})$ appartient à une progression arithmétique (contenant $0$) de\linebreak raison
$$
\hspace{24pt}
\mathrm{pas}(L;p)
\hspace{4pt}:=\hspace{4pt}
\frac{\vert\mathrm{W}(L)\vert}{\mathrm{p.g.c.d.}(\hspace{1pt}p-1\hspace{1pt},\hspace{1pt}24\hspace{1pt}\mathrm{h}(L)\hspace{1pt},\hspace{1pt}\vert\mathrm{W}(L)\vert\hspace{1pt})}
\hspace{24pt}.
$$

\medskip
-- Soient $L$ et $L'$ deux réseaux unimodulaires pairs de dimension $24$. Les inégalités de Ramanujan pour les $\theta_{r}(p)$, $r\geq 7$, fournissent un encadrement
$$
\mathrm{N}^{\mathrm{inf}}_{p}(L,L')
\hspace{4pt}\leq\hspace{4pt}
\mathrm{N}_{p}(L,L')
\hspace{4pt}\leq\hspace{4pt}
\mathrm{N}^{\mathrm{sup}}_{p}(L,L')
$$
dont la largeur est $2\hspace{1pt}\mathrm{K}(L,L')\hspace{1pt}(p+1)\hspace{1pt}p^{\frac{21}{2}}$ avec $\mathrm{K}(L,L'):=4\hspace{1pt}\sum_{r=7}^{10}\hspace{2pt}\vert\mathrm{c}_{r}(L,L')\vert$ (notations de \eqref{npmysterieux}).

\medskip
On prend pour $L$ un réseau de Niemeier avec racines avec $\mathrm{g}(L)$ non divisible par $p$ et pour $L'$ le réseau de Leech. Si la largeur en question est strictement inférieure à $\mathrm{pas}(L;p)$ alors $\mathrm{N}_{p}(L,\mathrm{Leech})$ est déterminé.

\medskip
Détaillons un peu.

\medskip
On a $\vert\theta_{7}(p)\vert\leq 4p^{\frac{19}{2}}$ et $\vert\theta_{r}(p)\vert\leq 4p^{\frac{21}{2}}$ pour $r=8,9,10$. On en déduit l'inégalité
$$
\hspace{24pt}
\vert\mathrm{N}^{2}_{p}(L,L')\vert
\hspace{4pt}\leq\hspace{4pt}
\mathrm{K}(L,L')\hspace{1pt}(p+1)\hspace{1pt}p^{\frac{21}{2}}
\hspace{24pt}.
$$
On pose
$$
\hspace{24pt}
\mathrm{N}^{\mathrm{inf}}_{p}(L,L')
\hspace{4pt}=\hspace{4pt}
\mathrm{N}^{1}_{p}(L,L')-\mathrm{K}(L,L')\hspace{1pt}(p+1)\hspace{1pt}p^{\frac{21}{2}}
\hspace{24pt},
$$
$$
\hspace{24pt}
\mathrm{N}^{\mathrm{sup}}_{p}(L,L')
\hspace{4pt}=\hspace{4pt}
\mathrm{N}^{1}_{p}(L,L')+\mathrm{K}(L,L')\hspace{1pt}(p+1)\hspace{1pt}p^{\frac{21}{2}}
\hspace{24pt}.
$$
Rappelons que nous avons noté $\mathrm{n}_{p}(L)$ l'entier défini par l'égalité $\mathrm{N}_{p}(L,\mathrm{Leech})\linebreak=\mathrm{n}_{p}(L)\hspace{2pt}\mathrm{pas}(L;p)$. On pose
$$
\hspace{24pt}
\nu^{\mathrm{inf}}_{p}(L)
\hspace{4pt}=\hspace{4pt}
\frac{\mathrm{N}^{\mathrm{inf}}_{p}(L,\mathrm{Leech})}{\mathrm{pas}(L;p)}
\hspace{24pt},\hspace{24pt}
\nu^{\mathrm{sup}}_{p}(L)
\hspace{4pt}=\hspace{4pt}
\frac{\mathrm{N}^{\mathrm{sup}}_{p}(L,\mathrm{Leech})}{\mathrm{pas}(L;p)}
\hspace{24pt};
$$
on a donc l'encadrement $\nu^{\mathrm{inf}}_{p}(L)\leq\mathrm{n}_{p}(L)\leq\nu^{\mathrm{sup}}_{p}(L)$.

\bigskip
\textit{Exemples.} On illustre par quelques exemples l'efficacité de cet encadrement~:

\bigskip
-- On a $\nu^{\mathrm{inf}}_{3}(\mathrm{L}_{23})\approx 0.99953$ et $\nu^{\mathrm{sup}}_{3}(\mathrm{L}_{23})\approx 1.00041$~; on en déduit $\mathrm{n}_{3}(\mathrm{L}_{23})=1$ et $\mathrm{N}_{3}(\mathrm{L}_{23},\mathrm{Leech})=8388608$ en accord avec le point (d) de III.4.2.10.

\bigskip
-- On a $\nu^{\mathrm{inf}}_{31}(\mathrm{A}_{24}^{+})\approx 275.99920$ et $\nu^{\mathrm{sup}}_{31}(\mathrm{A}_{24}^{+})\approx 276.00061$~; on en déduit $\mathrm{n}_{31}(\mathrm{A}_{24}^{+})=276$ et $\mathrm{N}_{31}(\mathrm{A}_{24}^{+},\mathrm{Leech})=142703132398645071052800000$,  en\linebreak accord avec le calcul effectué en III.4.3.

\bigskip
-- On constate que l'on a $\nu^{\mathrm{sup}}_{p}(\mathrm{E}_{8}\oplus\mathrm{E}_{8}\oplus\mathrm{E}_{8})<8\times 10^{-6}$ pour $p\leq 29$~; on a donc $\mathrm{N}_{p}(\mathrm{E}_{8}\oplus\mathrm{E}_{8}\oplus\mathrm{E}_{8},\mathrm{Leech})=0$ pour $p\leq 29$, en accord avec III.4.1.1.

\bigskip
-- On a $\nu^{\mathrm{inf}}_{47}(\mathrm{E}_{24})\approx 0.99992$ et $\nu^{\mathrm{sup}}_{47}(\mathrm{E}_{24})\approx1.00006$~; on en déduit $\mathrm{n}_{47}(\mathrm{E}_{24})=1$ et $\mathrm{N}_{47}(\mathrm{E}_{24},\mathrm{Leech})=113145617964492744063713280000$, en accord avec le point (d) de III.4.2.10.

\bigskip
On pose $\mathrm{n}^{\mathrm{inf}}_{p}(L)=\lceil\nu^{\mathrm{inf}}_{p}(L)\rceil$ et $\mathrm{n}^{\mathrm{sup}}_{p}(L)=\lfloor\nu^{\mathrm{sup}}_{p}(L)\rfloor$. Rappelons la notation~: soit $\nu$ un nombre réel~; $\lceil\nu\rceil$ est le plus petit entier $n$ avec $\nu\leq n$ et $\lfloor\nu\rfloor$ le plus grand entier $n$ avec $n\leq\nu$. On a donc par définition l'encadrement
$$
\hspace{24pt}
\mathrm{n}^{\mathrm{inf}}_{p}(L)
\hspace{4pt}\leq\hspace{4pt}
\mathrm{n}_{p}(L)
\hspace{4pt}\leq\hspace{4pt}
\mathrm{n}^{\mathrm{sup}}_{p}(L)
\hspace{24pt}.
$$
On note $\mathrm{e}_{p}(L)$ l'entier naturel $\mathrm{n}^{\mathrm{sup}}_{p}(L)-\mathrm{n}^{\mathrm{inf}}_{p}(L)$. Si l'on a $\mathrm{e}_{p}(L)=0$ alors $\mathrm{N}_{p}(L,\mathrm{Leech})$ est déterminé~: $\mathrm{N}_{p}(L,\mathrm{Leech})=\mathrm{n}^{\mathrm{inf}}_{p}(L)\hspace{2pt}\mathrm{pas}(L;p)$.

\bigskip
\textit{Exemple.}  On prend $L=\mathrm{E}_{24}$. Comme l'on a $\mathrm{g}(\mathrm{E}_{24})=2$,  il faut supposer $p\geq 3$ pour que la théorie précédente s'applique. \texttt{PARI} nous dit que $\mathrm{e}_{p}(\mathrm{E}_{24})$ est nul pour pour $3\leq p\leq 131$. L'entier $\mathrm{N}_{p}(\mathrm{E}_{24},\mathrm{Leech})$ est déterminé pour ces nombres premiers, on a par exemple (pour le plaisir enfantin d'écrire un très gros nombre entier~!)~:
$$
\hspace{12pt}
\mathrm{N}_{131}(\mathrm{E}_{24},\mathrm{Leech})
\hspace{4pt}=\hspace{4pt}123625448053001992116952381878687498240000
\hspace{12pt}.
$$

\bigskip
On considère maintenant le quadruplet
$$
\hspace{24pt}
\underline{\mathrm{e}}_{p}
\hspace{4pt}:=\hspace{4pt}
(\hspace{1pt}
\mathrm{e}_{p}(\mathrm{E}_{24})
\hspace{1pt},\hspace{1pt}
\mathrm{e}_{p}(\mathrm{E}_{16}\oplus\mathrm{E}_{8})
\hspace{1pt},\hspace{1pt}
\mathrm{e}_{p}(\mathrm{E}_{8}\oplus\mathrm{E}_{8}\oplus\mathrm{E}_{8})
\hspace{1pt},\hspace{1pt}
\mathrm{e}_{p}(\mathrm{A}_{24}^{+})
\hspace{1pt})
\hspace{24pt};
$$
comme l'on a $(\hspace{1pt}
\mathrm{g}(\mathrm{E}_{24})
\hspace{1pt},\hspace{1pt}
\mathrm{g}(\mathrm{E}_{16}\oplus\mathrm{E}_{8})
\hspace{1pt},\hspace{1pt}
\mathrm{g}(\mathrm{E}_{8}\oplus\mathrm{E}_{8}\oplus\mathrm{E}_{8})
\hspace{1pt},\hspace{1pt}
\mathrm{g}(\mathrm{A}_{24}^{+})
\hspace{1pt})=(2,2,1,5)$ on suppose $p\not=2,5$. \texttt{PARI} nous dit que l'on a $\underline{\mathrm{e}}_{p}=(0,0,0,0)$ pour $p=3$ et $7\leq p\leq 59$. Les entiers $\mathrm{N}_{p}(\mathrm{L}_{i},\mathrm{Leech})$ sont déterminés pour ces nombres premiers et $i=1,2,3,4$~; du coup le système linéaire \eqref{systeme} permet de calculer les $\theta_{r}(p)$, $r\geq 7$, pour les nombres premiers en question.

\bigskip
\textit{Les cas $61\leq p\leq 107$}

\bigskip
Pour $61\leq p\leq 107$, on n'a plus $\underline{\mathrm{e}}_{p}=(0,0,0,0)$~; on va montrer cependant que l'on peut encore déterminer les $\theta_{r}(p)$, $r\geq 7$. La méthode est décrite ci-après. On pose $\mathrm{n}_{p}(\mathrm{L}_{i})=\mathrm{n}^{\mathrm{inf}}_{p}(\mathrm{L}_{i})+x_{i}$, $1\leq i\leq 4$~; $(x_{1},x_{2},x_{3},x_{4})$ est donc un quadruplet d'entiers naturels qu'il s'agit de déterminer. Ce quadruplet est soumis aux contraintes suivantes, numérotées (1), (2), (3), (4)~:

\bigskip
(1) On a les inégalités $0\leq x_{i}\leq\mathrm{e}_{p}(\mathrm{L}_{i})$, $1\leq i\leq 4$.

\bigskip
(2) Le quadruplet $(x_{1},x_{2},x_{3},x_{4})$ vérifie des congruences affines modulo $p$  ou modulo des diviseurs de $p+1$. Détaillons un peu. La relation \eqref{systeme} peut être écrite sous la forme suivante~:
\begin{equation}\label{systemeaffine}
\hspace{24pt}
\begin{bmatrix} p\hspace{1pt}(p+1)\hspace{1pt}\theta_{7}(p) \\ (p+1)\hspace{1pt}\theta_{8}(p) \\ (p+1)\hspace{1pt}\theta_{9}(p) \\ (p+1)\hspace{1pt}\theta_{10}(p) 
\end{bmatrix}
\hspace{4pt}=\hspace{4pt}
\mathrm{F}(p)
\begin{bmatrix} x_{1} \\ x_{2} \\ x_{3} \\ x_{4}
\end{bmatrix}
+\mathrm{G}(p)
\hspace{24pt},
\end{equation}
$\mathrm{F}(p)$ désignant la matrice carrée $a^{-1}\hspace{1pt}\mathrm{diag}(\mathrm{pas}(\mathrm{L}_{1};p),\ldots,\mathrm{pas}(\mathrm{L}_{4};p))$ et $\mathrm{G}(p)$ la matrice colonne $a^{-1}\hspace{1pt}\mathrm{B}^{\mathrm{inf}}(p)$ avec
$$
\hspace{24pt}
\mathrm{B}^{\mathrm{inf}}(p)
\hspace{4pt}:=\hspace{4pt}
\begin{bmatrix}
\mathrm{n}^{\mathrm{inf}}_{p}(\mathrm{L}_{1})\hspace{1pt}\mathrm{pas}(\mathrm{L}_{1};p)-\mathrm{N}^{1}_{p}(\mathrm{L}_{1},\mathrm{Leech}) \\
\vspace{-10pt} \\
\mathrm{n}^{\mathrm{inf}}_{p}(\mathrm{L}_{2})\hspace{1pt}\mathrm{pas}(\mathrm{L}_{2};p)-\mathrm{N}^{1}_{p}(\mathrm{L}_{2},\mathrm{Leech}) \\
\vspace{-10pt} \\
\mathrm{n}^{\mathrm{inf}}_{p}(\mathrm{L}_{3})\hspace{1pt}\mathrm{pas}(\mathrm{L}_{3};p)-\mathrm{N}^{1}_{p}(\mathrm{L}_{3},\mathrm{Leech}) \\
\vspace{-10pt} \\
\mathrm{n}^{\mathrm{inf}}_{p}(\mathrm{L}_{4})\hspace{1pt}\mathrm{pas}(\mathrm{L}_{4};p)-\mathrm{N}^{1}_{p}(\mathrm{L}_{4},\mathrm{Leech})
\end{bmatrix}
\hspace{24pt}.
$$
La matrice $\mathrm{F}(p)$ est à coefficients entiers pour $p\not\equiv 1\bmod{23}$. Expliquons pourquoi. Comme dans la remarque qui suit III.4.3.3 on écrit $\mathrm{pas}(\mathrm{L}_{i};p)=\mathrm{pas}_{1}(\mathrm{L}_{i})\hspace{1pt}\mathrm{pas}_{2}(\mathrm{L}_{i};p)$. On constate que la matrice
$$
a^{-1}\hspace{2pt}\mathrm{diag}(\mathrm{pas}_{1}(\mathrm{L}_{1}),\ldots,\mathrm{pas}_{1}(\mathrm{L}_{4}))\hspace{2pt}\mathrm{diag}(23,1,1,1)
$$ (qui est indépendante de $p$) est à coefficients entiers~; or, par définition même, $23$ divise l'entier $\mathrm{pas}_{2}(\mathrm{E}_{24};p)$ pour $p\not\equiv 1\bmod{23}$.  On observera que si $\mathrm{F}(p)$ est à coefficients entiers alors, compte tenu de \eqref{systemeaffine} il est de même pour $\mathrm{G}(p)$. Le fait que $\mathrm{F}(p)$, $\mathrm{G}(p)$ sont à coefficients entiers et que les $\theta_{r}(p)$, $r\geq 7$, sont entiers conduit aux congruences évoquées plus haut.

\bigskip
(3) Soient $\theta_{r}(p;X_{1},X_{2},X_{3},X_{4})$, $7\leq r\leq 10$, les quatre polynômes affines de~$\mathbb{Q}[X_{1},X_{2},X_{3},X_{4}]$ définis par
\footnotesize
$$
\begin{bmatrix}
\theta_{7}(p;X_{1},X_{2},X_{3},X_{4}) \\ \theta_{8}(p;X_{1},X_{2},X_{3},X_{4}) \\  \theta_{9}(p;X_{1},X_{2},X_{3},X_{4}) \\ \theta_{10}(p;X_{1},X_{2},X_{3},X_{4})
\end{bmatrix}
:=\mathrm{diag}(p(p+1),p+1,p+1,p+1)^{-1}
(\mathrm{F}(p)
\begin{bmatrix} X_{1} \\ X_{2} \\ X_{3} \\ X_{4}
\end{bmatrix}
+\mathrm{G}(p))
$$
\normalsize
(la notation est étrange\ldots mais transparente, on observera que la condition (2) équivaut au fait que $\theta_{r}(p;x_{1},x_{2},x_{3},x_{4})$ est entier pour $7\leq r\leq 10$). Soit $\mathrm{T}_{p}(X_{1},X_{2},X_{3},X_{4})$ la matrice obtenue en substituant dans \eqref{tpformel} aux $\theta_{r}(p)$, $r\leq 7$, les $\theta_{r}(p;X_{1},X_{2},X_{3},X_{4})$~; $\mathrm{T}_{p}(X_{1},X_{2},X_{3},X_{4})$ est donc une matrice à $24$ lignes et $24$ colonnes dont les coefficients sont des polynômes affines de $\mathbb{Q}[X_{1},X_{2},X_{3},X_{4}]$ telle que l'on a $\mathrm{T}_{p}(x_{1},x_{2},x_{3},x_{4})=\mathrm{T}_{p}$. On note $\mathrm{N}_{p}(\mathrm{L}_{i},\mathrm{L}_{j};X_{1},X_{2},X_{3},X_{4})$ le coefficient d'indice $(j,i)$ de $\mathrm{T}_{p}(X_{1},X_{2},X_{3},X_{4})$ et on pose, pour $i\leq 23$,
$$
\mathrm{n}_{p}(\mathrm{L}_{i};X_{1},X_{2},X_{3},X_{4})
\hspace{4pt}:=\hspace{4pt}
\frac{\mathrm{N}_{p}(\mathrm{L}_{i},\mathrm{Leech};X_{1},X_{2},X_{3},X_{4})}{\mathrm{pas}(\mathrm{L}_{i};p)}
$$
(on a donc, par construction, $\mathrm{n}_{p}(\mathrm{L}_{i};X_{1},X_{2},X_{3},X_{4})=X_{i}+\mathrm{n}^{\mathrm{inf}}_{p}(\mathrm{L}_{i})$ pour $i\leq 4$).  On a par exemple 
$$
\mathrm{n}_{p}(\mathrm{L}_{5};X_{1},X_{2},X_{3},X_{4})
\hspace{4pt}=\hspace{4pt}
\sum_{i=1}^{4}\hspace{4pt}\alpha_{i}\hspace{1pt}\gamma_{i}(p)\hspace{2pt}X_{i}
+\beta(p)
$$
avec
\begin{multline*}
\hspace{12pt}
(\alpha_{1},\alpha_{2},\alpha_{3},\alpha_{4})=(\hspace{1pt}-1472\hspace{1pt},\hspace{1pt}-\frac{41}{8}\hspace{1pt},\hspace{1pt}\frac{ 119}{16}\hspace{1pt},\hspace{1pt}\frac{281}{256}\hspace{1pt}) \\
\gamma_{i}(p)=\frac{\mathrm{p.g.c.d.}(p-1,24\hspace{1pt}\mathrm{h}(\mathrm{L}_{5}))}{\mathrm{p.g.c.d.}(p-1,24\hspace{1pt}\mathrm{h}(\mathrm{L}_{i}))}
\hspace{12pt},\hspace{12pt}
\beta(p)=\frac{\mathrm{N}^{1}_{p}(\mathrm{L}_{5},\mathrm{Leech})}{\mathrm{pas}(\mathrm{L}_{i};p)}
\hspace{12pt}.
\end{multline*}
D'après III.4.3.3, $\mathrm{n}_{p}(\mathrm{L}_{5};x_{1},x_{2},x_{3},x_{4})=\mathrm{n}_{p}(\mathrm{L}_{5})$ est entier~; on en déduit que le quadruplet $(x_{1},x_{2},x_{3},x_{p})$ vérifie une congruence affine modulo un entier explicite dépendant de $p$, disons $\mathrm{m}(p)$. Précisons un peu. On constate que le dénominateur commun des rationnels $\alpha_{i}\hspace{1pt}\gamma_{i}(p)$ (mis sous forme irréductible) est $\mathrm{m}(p)$, la fonction $p\mapsto\mathrm{m}(p)$ pour $61\leq p\leq 113$ étant donnée par le tableau suivant (on a rajouté les valeurs de $\mathrm{m}(p)$ pour $p=109$ et $p=113$ en vue d'une future application)
\scriptsize
\begin{center}
\hspace{4pt}
\begin{tabular}
{|c|c|c|c|c|c|c|c|c|c|c|c|c|c|c|c|c|}
\hline
$p$ & 61 & 67 & 71 & 73 & 79 & 83 & 89 & 97 & 101 & 103 & 107 & 109 & 113 \\
\hline
$\mathrm{m}$ & 1280 &  256 &  1280 &  768 &  256 &  256 &  256 &  128 &  6400 &  256 &  256 &  768 & 128 \\
\hline
\end{tabular}
\hspace{4pt}.
\end{center}
\normalsize
On en déduit que $\mathrm{m}(p)\beta(p)$ est entier et la congruence évoquée plus haut.

\medskip
On observe que la valuation $2$-adique de $\mathrm{v}_{2}(\gamma_{i}(p))$ est nulle pour $i=1,2,3$ et que l'on a $\mathrm{v}_{2}(\gamma_{4}(p))=0$ pour $\mathrm{v}_{2}(p-1)\leq 3$, $\mathrm{v}_{2}(\gamma_{4}(p))=1$ pour $\mathrm{v}_{2}(p-1)\geq 4$. On déduit, en contemplant le quadruplet $(\alpha_{1},\alpha_{2},\alpha_{3},\alpha_{4})$, que la congruence que l'on vient de dégager détermine la classe de $x_{4}$ modulo $16$ pour $p\not=97,113$ (et $61\leq p\leq 113$) et modulo $8$ pour $p=97,113$. La valeur de $x_{4}$ modulo $16$ pour $61\leq p\leq 113$ est donnée par le tableau suivant
\scriptsize
\begin{center}
\hspace{4pt}
\begin{tabular}
{|c|c|c|c|c|c|c|c|c|c|c|c|c|c|c|c|c|}
\hline
$p$ & 61 & 67 & 71 & 73 & 79 & 83 & 89 & 97 & 101 & 103 & 107 & 109 & 113 \\
\hline
$x_{4}\bmod{16}$ & 3 &  0 &  2 &  11 &  7 &  3 &  2 &  0 ou 8 &  5 &  4 &  5 &  0 & 0 ou 8 \\
\hline
\end{tabular}
\hspace{4pt}.
\end{center}
\normalsize
Cette forme faible de la contrainte (3) sera notée ($3_{4}$) ci-après.

\bigskip
(4) On doit avoir $\vert\theta_{7}(p;x_{1},x_{2},x_{3},x_{4})\vert\leq 4\hspace{1pt}p^{\frac{19}{2}}$ et $\vert\theta_{r}(p;x_{1},x_{2},x_{3},x_{4})\vert\leq 4\hspace{1pt}p^{\frac{21}{2}}$ pour $r=8,9,10$.

\bigskip
\textit{Le cas $p=61$}

\medskip
On a $\underline{\mathrm{e}}_{61}=(0,0,0,5)$. La contrainte  ($3_{4}$) détermine le quadruplet cherché ~: $(x_{1},x_{2},x_{3},x_{4})=(0,0,0,3)$.

\bigskip
\textit{Le cas $p=67$}

\medskip
On a $\underline{\mathrm{e}}_{67}=(0,0,0,1)$. La contrainte  ($3_{4}$) détermine le quadruplet cherché ~: $(x_{1},x_{2},x_{3},x_{4})=(0,0,0,0)$.

\bigskip
\textit{Le cas $p=71$}

\medskip
On a $\underline{\mathrm{e}}_{71}=(0,0,0,6)$.  La contrainte  ($3_{4}$) détermine le quadruplet cherché ~: $(x_{1},x_{2},x_{3},x_{4})=(0,0,0,2)$.

\bigskip
\textit{Le cas $p=73$}

\medskip
On a $\underline{\mathrm{e}}_{73}=(0,6,9,20)$. La contrainte (2) donne notamment $(x_{1},x_{2},x_{3},x_{4})\equiv (0,3,5,11)\bmod{37}$. Le quadruplet cherché est $(0,3,5,11)$.
\vfill\eject

\bigskip
\textit{Le cas $p=79$}

\medskip
On a $\underline{\mathrm{e}}_{79}=(0,0,1,12)$.  La contrainte (2) donne $2\hspace{1pt}x_{3}+x_{4}\equiv 9\bmod{79}$~;\linebreak l'encadrement $0\leq 2\hspace{1pt}x_{3}+x_{4}\leq 14$ force $2\hspace{1pt}x_{3}+x_{4}=9$.  La contrainte  ($3_{4}$) détermine alors le quadruplet cherché ~: $(x_{1},x_{2},x_{3},x_{4})=(0,0,1,7)$.

\bigskip
\textit{Le cas $p=83$}

\medskip
On a $\underline{\mathrm{e}}_{83}=(0,0,0,6)$. La contrainte  ($3_{4}$) détermine le quadruplet cherché ~: $(x_{1},x_{2},x_{3},x_{4})=(0,0,0,3)$.

\bigskip
\textit{Le cas $p=89$}

\medskip
On a $\underline{\mathrm{e}}_{89}=(0,6,10,67)$. La contrainte (2) donne $-2\hspace{1pt}x_{2}-7\hspace{1pt}x_{3}+x_{4}+7\equiv 0\bmod{89}$~; l'encadrement $-75\leq -2\hspace{1pt}x_{2}-7\hspace{1pt}x_{3}+x_{4}-7\leq 74$ force $x_{4}=2\hspace{1pt}x_{2}+7\hspace{1pt}x_{3}-7$. La contrainte (3) donne $10\hspace{1pt}x_{2}+3\hspace{1pt}x_{3}-45\equiv 0\bmod{256}$ ; le même argument que précédemment montre que l'on a  $10\hspace{1pt}x_{2}+3\hspace{1pt}x_{3}-45=0$. En particulier $x_{2}$ est divisible par $3$. L'encadrement $0\leq x_{3}\leq 10$ force alors $x_{2}=3$. Le quadruplet cherché est $(0,3,5,34)$.

\bigskip
\textit{Le cas $p=97$}

\medskip
On a $\underline{\mathrm{e}}_{97}=(0,117,187,548)$. L'ordinateur dit que les quadruplets satisfaisant (1), (2), et (3) sont $(0, 22, 63, 432)$ et $(0, 71, 105, 272)$ (on peut aider l'ordinateur en observant que $(2)$ implique  $(x_{2},x_{3})\equiv (1,0)\bmod{7}$ et que ($3_{4}$) dit que $x_{4}$ est divisible par $8$). Le premier ne passe pas le ``test de Ramanujan'' (la contrainte (4))~: $(x_{1},x_{2},x_{3},x_{4})=(0, 71, 105, 272)$.

\bigskip
\textit{Le cas $p=101$}

\medskip
On a $\underline{\mathrm{e}}_{101}=(0, 78, 124, 3643)$. L'ordinateur dit que le seul quadruplet satisfaisant (1), (2) et (3) est $(0, 30, 63, 2149)$.

\bigskip
\textit{Le cas $p=103$}

\medskip
On a $\underline{\mathrm{e}}_{103}=(0, 29, 46, 273)$. L'ordinateur dit que les quadruplets satisfaisant (1), (2) et (3) sont $(0, 7, 46, 196)$, $(0, 15, 27, 148)$ et $(0, 23, 8, 100)$. Le premier et le troisième ne passent pas le test de Ramanujan : $(x_{1},x_{2},x_{3},x_{4})=(0, 15, 27, 148)$.

\bigskip
\textit{Le cas $p=107$}

\medskip
On a $\underline{\mathrm{e}}_{107}=(0, 14, 23, 141)$. L'ordinateur dit que le seul quadruplet satisfaisant (1), (2) et (3) est $(0, 7, 10, 53)$.
\vfill\eject

\bigskip
\textit{Les cas $p=109$ et $p=113$}

\bigskip
Une conséquence inattendue du théorème \ref{recapcong}, que nous démontrerons dans le prochain paragraphe en faisant appel à la théorie des représentations galoisiennes, est que le quadruplet $(x_{1},x_{2},x_{3},x_{4})$, introduit lors de l'étude des cas $61\leq p\leq 107$, vérifie des congruences affines modulo des diviseurs de $p+1$ dont certaines peuvent être ``indépendantes'' des congruences de la contrainte~(2). Ces contraites supplémentaires permettent de déterminer le quadruplet $(x_{1},x_{2},x_{3},x_{4})$ pour $p=109$ et $p=113$.

\bigskip
\textit{Le cas $p=109$}

\medskip
On a $\underline{\mathrm{e}}_{109}=(0, 337, 538, 1049)$. L'ordinateur dit qu'il existe $208$  quadruplets satisfaisant les contraintes (1), (2) et (3) et que parmi ceux-ci  $12$ satisfont la contrainte (4). On note $\underline{x}^{(1)},\underline{x}^{(2)},\ldots,\underline{x}^{(12)}$ ces $12$ quadruplets d'entiers naturels~; on suppose, pour fixer les idées, que l'on a $\underline{x}^{(1)}<\underline{x}^{(2)}<\ldots<\underline{x}^{(12)}$ pour l'ordre lexicographique. Pour $k=1,2,\ldots,12$, on vérifie que $\mathrm{T}_{109}(\underline{x}^{(k)})$ est à coefficients entiers et que $\mathrm{n}_{109}(\mathrm{L}_{i};\underline{x}^{(k)})$ est aussi entier pour $6\leq i\leq 23$ (nous avons tout fait pour que $\mathrm{n}_{109}(\mathrm{L}_{i};\underline{x}^{(k)})$ soit entier pour $1\leq i\leq 5$).

\medskip
Ce qui précède montre que la méthode qui nous a permis de déterminer $\tau_{6,8}(p)$, $\tau_{8,8}(p)$, $\tau_{12,6}(p)$ et $\tau_{4,10}(p)$ pour $p\leq 107$ est en échec pour $p=109$.

\medskip
On parvient cependant à lever l'indétermination évoquée ci-dessus de la façon suivante. On contemple la suite de quadruplets d'entiers relatifs
$$
\hspace{24pt}
(\hspace{2pt}\theta_{7}(109;\underline{x}^{(k)})\hspace{2pt},\hspace{2pt}\theta_{8}(109;\underline{x}^{(k)})\hspace{2pt},\hspace{2pt}\theta_{9}(109;\underline{x}^{(k)})\hspace{2pt},\hspace{2pt}\theta_{10}(109;\underline{x}^{(k)})\hspace{2pt})_{k=1,2,\ldots,12}
\hspace{24pt}.
$$
Par chance, un seul de ces quadruplets, à savoir le $5$-ième, vérifie la congruence imposée par le point (12) du théorème \ref{recapcong} (le quadruplet en question est d'ailleurs celui qui était le plus probable ``au sens de Sato-Tate''). On observera que le point (12) du théorème \ref{recapcong}  est une congruence modulo $11$ et que $11$ divise $109+1$~!

\medskip
On peut paraphraser ainsi ce qui précède. En plus des contraintes (1), (2), (3) et (4), les entiers naturels $x_{i}:=\mathrm{n}_{109}(\mathrm{L}_{i})-\mathrm{n}^{\mathrm{inf}}_{109}(\mathrm{L}_{i})$, $i=1,2,3,4$, sont soumis à la contrainte, disons (2-supp), imposée par le point (12) du théorème \ref{recapcong}. Il n'est pas difficile d'expliciter (2-supp), il s'agit de la congruence
$$
\hspace{24pt}
x_{1} + x_{2} + 6\hspace{1pt}x_{3} + x_{4} + 2
\hspace{4pt}\equiv\hspace{4pt}
0
\hspace{8pt}
\bmod{11}
\hspace{24pt}.
$$
Un seul quadruplet $(x_{1},x_{2},x_{3},x_{4})$ satisfait les contraintes (1), (2), (2-supp), (3) et (4)~: le quadruplet $(0,138,284,576)$.
\vfill\eject

\bigskip
\textit{Le cas $p=113$}

\medskip
On a $\underline{\mathrm{e}}_{113}=(0,227,361,1058)$. On considère cette fois les points (3) et (4) du théorème \ref{recapcong} (parce que $19$ divise $113+1$). Ces deux points imposent au quadruplet $(x_{1},x_{2},x_{3},x_{4})$ une contrainte, disons (2-supp), constituée de deux congruences affines modulo $19$~; on constate  que celles-ci sont indépendantes des deux congruences affines modulo $19$ qui apparaissent dans la contrainte~(2). Du coup, le quadruplet $(x_{1},x_{2},x_{3},x_{4})$ est déterminé modulo~$19$~; on trouve~:
$$
(x_{1},x_{2},x_{3},x_{4})
\hspace{4pt}\equiv\hspace{4pt}
(0,6,3,16)
\hspace{8pt}
\bmod{19}
$$
(la présence du $0$ en première position au second membre est rassurante~!). L'ordinateur montre ensuite qu'il n'existe que deux quadruplets satisfaisant les contraintes (1), (2), (2-supp) et (3)~: $(0,120,155,396)$ et $(0,177,326,244)$. Le second ne passe pas le test de Ramanujan. Les entiers $\tau_{6,8}(113)$, $\tau_{8,8}(113)$, $\tau_{12,6}(113)$ et $\tau_{4,10}(113)$ et l'endomorphisme $\mathrm{T}_{113}$ de $\mathbb{Z}[\mathrm{X}_{24}]$ sont donc déter\-minés.

\medskip
On trouve par exemple
$$
\mathrm{N}_{113}(\mathrm{L}_{12},\mathrm{L}_{21})
=
633323838523478069636624166862873752207360000
$$
(c'est le plus grand nombre de $p$-voisins, $p$ premier, parmi ceux que nous avons pu calculer).

\bigskip
\textit{Le cas $p=127$}

\medskip
Les méthodes qui nous ont permis de déterminer $\tau_{6,8}(p)$, $\tau_{8,8}(p)$, $\tau_{12,6}(p)$ et $\tau_{4,10}(p)$ pour $p\leq 113$ sont en échec pour $p=127$. Nous expliquons pourquoi ci-après (en petits caractères).

\medskip
\footnotesize
-- L'ordinateur dit qu'il existe $3329$  quadruplets satisfaisant les contraintes (1), (2), (3) et (4), disons $\underline{x}^{(1)},\underline{x}^{(2)},\ldots,\underline{x}^{(3329)}$. Pour $k=1,2,\ldots,3329$, on vérifie que $\mathrm{T}_{127}(\underline{x}^{(k)})$ est à coefficients entiers et que $\mathrm{n}_{127}(\mathrm{L}_{i};\underline{x}^{(k)})$ est aussi entier pour $6\leq i\leq 23$.

\medskip
-- Le seul nombre premier qui divise $127+1$ est $2$ et l'on constate que les entiers $\theta_{r}(127; \underline{x}^{(k)})$, $7\leq r\leq 10$, $1\leq k\leq 3329$, sont tous pairs, en accord avec la congruence (12) du théorème~\ref{recapcong}. En fait, en contemplant les quatre polynômes affines $\theta_{r}(127;X_{1},X_{2},X_{3},X_{4})$, $7\leq r\leq 10$, on constate que l'on a $\theta_{7}(127;\underline{x})\equiv134400\bmod{2^{18}}$,  $\theta_{8}(127;\underline{x})\equiv3840\bmod{2^{13}}$, $\theta_{9}(127;\underline{x})\equiv-3840\bmod{2^{13}}$ et $\theta_{10}(127;\underline{x})\equiv256\bmod{2^{10}}$ pour tout $\underline{x}$ dans $\mathbb{Z}^{4}$.
\normalsize

\subsection{Détermination des $\tau_{j,k}(p^{2})$ pour $p\leq 29$}\label{calcfinvpgenre2-2}

On  note, ci-dessous, respectivement $\varpi_{1},\varpi_{2}, \dots,\varpi_{10}$ les représentations automorphes $\Delta_{11}$, $\Delta_{15}$, $\Delta_{17}$, $\Delta_{19}$, $\Delta_{21}$, $\mathrm{Sym}^{2}\Delta_{11}$, $\Delta_{19,7}$, $\Delta_{21,9}$, $\Delta_{21,13}$ et $\Delta_{21,5}$ (voir \S\ref{ikedabocherer}, \S VI.\ref{exconjal}, définition \ref{definitiondeltawv}). La représentation $\varpi_{r}$ est donc dans $\Pi_{\mathrm{cusp}}(\mathrm{PGL}_{2})$ pour $r\leq 5$, dans $\Pi_{\mathrm{cusp}}(\mathrm{PGL}_{3})$ pour $r=6$ et dans $\Pi_\mathrm{cusp}(\mathrm{PGL}_{4})$ pour $r\geq 7$.

\medskip
Soit $p$ un nombre premier, les entiers $\theta_{r}(p)$, $1\leq r\leq 10$, que nous avons introduits au paragraphe 2 vérifient par définition la relation
$$
\hspace{24pt}
\theta_{r}(p)
\hspace{4pt}=\hspace{4pt}
p^{\frac{\mathrm{w}(\varpi_{r})}{2}}\hspace{1pt}\mathop{\mathrm{trace}}\mathrm{c}_p(\varpi_{r})
\hspace{24pt}.
$$
(Voir \S VI.\ref{parlanpig}, \S VIII.\ref{paralgreg}~; on rappelle que $\mathrm{w}(\pi)$ désigne le poids motivique d'une représentation automorphe algébrique $\pi$ dans $\Pi_{\mathrm{cusp}}(\mathrm{PGL}_n)$, ici l'indice $w$ pour
$\pi=\Delta_{w}$ ou $\Delta_{w,v}$.)

\medskip
On pose pareillement
$$
\hspace{24pt}
\theta_{r}(p^{2})
\hspace{4pt}:=\hspace{4pt}
p^{\mathrm{w}(\varpi_{r})}\hspace{1pt}\mathop{\mathrm{trace}}\hspace{2pt}({\mathrm{c}_p(\varpi_{r})}^{2})
\hspace{24pt}.
$$
Pour $r\geq 7$, la définition ci-dessus est en accord avec celle que nous avons choisie en IX.\ref{definitiontaujk}~:  $\theta_{7}(p^{2})=\tau_{6,8}(p^{2})$, $\theta_{8}(p^{2})=\tau_{8,8}(p^{2})$, $\theta_{9}(p^{2})=\tau_{12,6}(p^{2})$ et
$\theta_{9}(p^{2})=\tau_{4,10}(p^{2})$.

\medskip
Par contre, on prendra garde que, pour $r=1,2,3,4,5$, $\theta_{i}(p^{2})$ n'est pas la valeur en $p^{2}$ des fonctions arithmétiques $\tau_{12}, \tau_{16},\tau_{18}, \tau_{20},\tau_{22}$~; on a en fait $\theta_{1}(p^{2})=\tau_{12}(p)^2-2\hspace{1pt}p^{11}$, $\theta_{2}(p^{2})=\tau_{16}(p)^2-2\hspace{1pt}p^{15}$, $\theta_{3}(p^{2})=\tau_{18}(p)^2-2\hspace{1pt}p^{17}$, $\theta_{4}(p^{2})=\tau_{20}(p)^2-2\hspace{1pt}p^{19}$ et $\theta_{5}(p^{2})=\tau_{22}(p)^2-2\hspace{1pt}p^{21}$. Enfin, on vérifie aisément que l'on a
$\theta_{6}(p^{2})=\tau_{12}(p)^4-4\hspace{1pt}p^{11}\hspace{1pt}\tau_{12}(p)^2 + 3\hspace{1pt}p^{22}$.

\medskip
Soit $\mathrm{V}_{\mathrm{St}}$ la représentation standard de $\mathrm{SO}_{n}(\mathbb{C})$. Pour $p$ un nombre premier, on note $\mathrm{T}_{p}^{\psi^2}$ l'opérateur de Hecke dans $\mathrm{H}_p(\mathrm{O}_n)$ défini, \textit{via} l'isomorphisme de Satake, par la formule
$$
p^{2-n}\hspace{4pt}\mathrm{Sat}(\mathrm{T}_p^{\psi^2})
\hspace{4pt}=\hspace{4pt}
\mathop{\psi^{2}}\hspace{2pt}[\mathrm{V}_{\mathrm{St}}]
\hspace{4pt}:=\hspace{4pt}
[\mathrm{V}_{\mathrm{St}}\otimes\mathrm{V}_{\mathrm{St}}]-2\hspace{2pt} [\Lambda^{2}\hspace{1pt}\mathrm{V}_{\mathrm{St}}]
$$
(voir \S VI.\ref{isomsatake}). L'intégralité de cet opérateur découle du lemme
\ref{satakecoefcar} ; en fait on a dans le cas présent
\begin{equation}\label{tpsi2}
\mathrm{T}_{p}^{\psi^{2}}
\hspace{4pt}=\hspace{4pt}
\mathrm{T}_{p}^{2}
- 2\hspace{1pt}p\hspace{1pt}\mathrm{T}_{p,p} - 2\hspace{1pt}p\hspace{2pt}(\hspace{2pt}\sum_{i=0}^{\frac{n}{2}-2}\hspace{1pt}p^{2i}+p^{\frac{n}{2}-2}\hspace{2pt})
\end{equation}
d'après les formules VI.\eqref{formulesattp} et VI.\eqref{formulesattpp}.

\medskip
On suppose maintenant $n=24$ et on note $\lambda_{j}^{\psi^2}(p)$ la valeur propre de $\mathrm{T}_p^{\psi^2}$ sur le vecteur $\mathrm{v}_{j}$ de $\mathbb{Z}[\mathrm{X}_{24}]$, $1\leq j\leq 24$~; on a tout fait pour avoir
$$
\hspace{24pt}
\lambda_{j}^{\psi^2}(p)
\hspace{4pt}=\hspace{4pt}
p^{22}\hspace{2pt}\mathop{\mathrm{trace}}\hspace{2pt}\mathrm{St}(\mathrm{c}_p(\pi_{j})^{2})
\hspace{4pt}=\hspace{4pt}
\mathrm{C}_{j,0}(p^{2}) + \sum_{r=1}^{10}\hspace{4pt}\mathrm{C}_{j,r}(p^{2})\hspace{2pt}\theta_{r}(p^{2})
\hspace{24pt},
$$
les polyn\^omes $\mathrm{C}_{j,r}$ de $\mathbb{Z}[X]$ étant ceux introduits au paragraphe 2.
\vfill\eject

\medskip
On revient au cas général. La formule
$$
(p+1)\hspace{1pt}\mathrm{T}_{p,p}
\hspace{4pt}=\hspace{4pt}
\mathrm{T}_p^2-\mathrm{T}_{p^{2}}-\frac{(p^{\frac{n}{2}}-1)(p^{\hspace{1pt}\frac{n}{2}-1}+1)}{(p-1)}
$$
de l'exemple VI.\ref{calcultp2} et la formule \eqref{tpsi2}, montrent que l'opérateur $\mathrm{T}_{p^{2}}$ s'exprime en fonction de $\mathrm{T}_{p}^{\psi^{2}}$ et $\mathrm{T}_{p}^{2}$~:
\begin{equation}\label{tp2}
\hspace{24pt}
\mathrm{T}_{p^{2}}
\hspace{4pt}=\hspace{4pt}
\frac{p+1}{2\hspace{1pt}p}\hspace{4pt}\mathrm{T}_{p}^{\psi^{2}}
+
\frac{p-1}{2\hspace{1pt}p}\hspace{4pt}\mathrm{T}_{p}^{2}\hspace{4pt}
-p^{n-2}+p^{\frac{n}{2}-2}
\hspace{24pt}.
\end{equation}
Comme l'on a déterminé les $\tau_{j,k}(p)$ pour $p\leq 113$ en \ref{calcfinvpgenre2-1}, les considérations ci-dessus montrent que la détermination des $\tau_{j,k}(p^2)$ avec $p \leq 29$ se ramène à celle des nombres de voisins $\mathrm{N}_{p^{2}}(\mathrm{L}_{i},\mathrm{Leech})$ pour $i\leq 4$ et $p\leq 29$ (on rappelle que l'on a posé $\mathrm{E}_{24}=\mathrm{L}_{1}$, $\mathrm{E}_{16}\oplus\mathrm{E}_{8}=\mathrm{L}_{2}$, $\mathrm{E}_{8}\oplus\mathrm{E}_{8}\oplus\mathrm{E}_{8}=\mathrm{L}_{3}$ et $\mathrm{A}^{+}_{24}=\mathrm{L}_{4}$)~; en effet, le système qui exprime les $\mathrm{N}_{p^{2}}(\mathrm{L}_{i},\mathrm{Leech})$ en fonction des ``inconnues''  $\tau_{j,k}(p^2)$ est toujours un système de Cramer parce que la matrice $\mathrm{a}$ de \eqref{mata} est inversible. On détermine ces nombres de voisins ci-après.

\bigskip
\textit{Les cas $p\leq 3$}

\medskip
On a $\mathrm{N}_{p^{2}}(\mathrm{L}_{i},\mathrm{Leech})=0$ pour $i\leq 4$ et $p=2,3$ d'après la proposition III.4.1.1.

\bigskip
\textit{Le cas $p=5$}

\medskip
On a encore $\mathrm{N}_{25}(\mathrm{L}_{i},\mathrm{Leech})=0$ pour $i\leq 3$ d'après la proposition III.4.1.1~; le point (d) du théorème III.4.2.10 donne la valeur de $\mathrm{N}_{25}(\mathrm{L}_{4},\mathrm{Leech})$.

\bigskip
\textit{Le cas $p=7$}

\medskip
On adapte la méthode utilisée précédemment pour déterminer $\mathrm{N}_{p}(\mathrm{L}_{i},\mathrm{Leech})$, $i=1,2,3,4$, dans les cas $7\leq p\leq 59$.

\medskip
La proposition III.4.3.1, dont on reprend les notations, doit être modifiée comme suit~:

\bigskip
\textbf{Proposition 8.2.1} {\em Soient $L$ un réseau de Niemeier avec racines et $p$ un nombre premier qui ne divise pas l'indice de $Q$ dans $L$~; soit $S$ le stabilisateur, pour l'action de $W$, d'un élément de $\mathrm{P}^{\mathrm{r\acute{e}g}}_{L}(\mathbb{Z}/p^{2})$.

\medskip
{\em (a)} Le groupe $S$ s'identifie à un sous-groupe de $(\mathbb{Z}/p^{2})^{\times}$.

\medskip
{\em (b)} Soit $S^{p}$ l'image de $S$ par l'endomorphisme $x\mapsto x^{p}$ du groupe $(\mathbb{Z}/p^{2})^{\times}$~;  l'action de $S^{p}$ sur $R$ (induite par celle de $W$) est libre.}

\bigskip
Du coup, le scholie-définition III.4.3.3 doit être modifié comme suit~:

\bigskip
\textbf{Scholie-Définition 8.2.2} {\em Soient $L$ un réseau de Niemeier avec racines et $p$ un nombre premier~; on note $\mathrm{pas}\hspace{1pt}(L;p^{2})$ l'entier défini par
$$
\hspace{24pt}
\mathrm{pas}\hspace{1pt}(L;p^{2})
\hspace{4pt}:=\hspace{4pt}
\frac{\vert\mathrm{W}(L)\vert}{\mathrm{p.g.c.d.}\hspace{1pt}(\hspace{1pt}p\hspace{1pt}(p-1)\hspace{1pt},\hspace{1pt}24\hspace{1pt}p\hspace{1pt}\mathrm{h}(L)\hspace{1pt},\hspace{1pt}\vert\mathrm{W}(L)\vert\hspace{1pt})}
\hspace{24pt}.
$$
Si $p$ ne divise pas l'indice de $Q$ dans $L$, alors $\mathrm{N}_{p^{2}}(L,\mathrm{Leech})$ est divisible par $\mathrm{pas}\hspace{1pt}(L;p^{2})$~; on note dans ce cas $\mathrm{n}_{p^{2}}(L)$ l'entier défini par
$$
\hspace{24pt}
\mathrm{N}_{p^{2}}(L,\mathrm{Leech})
\hspace{4pt}=\hspace{4pt}
\mathrm{n}_{p^{2}}(L)\hspace{4pt}
\mathrm{pas}\hspace{1pt}(L;p^{2})
\hspace{24pt}.
$$}

\medskip
On procède ensuite comme lorsque l'on a déterminé $\mathrm{N}_{p}(\mathrm{L}_{i},\mathrm{Leech})$ pour $i=1,2,3,4$ et $7\leq p\leq 59$. \textit{Mutatis mutandis}, on définit des entiers $\mathrm{n}_{p^{2}}^{\mathrm{inf}}(L)$ et $\mathrm{n}_{p^{2}}^{\mathrm{sup}}(L)$, ``aisément calculables'', tels que l'on a l'encadrement
$$
\hspace{24pt}
\mathrm{n}_{p^{2}}^{\mathrm{inf}}(L)
\hspace{4pt}\leq\hspace{4pt}
\mathrm{n}_{p^{2}}(L)
\hspace{4pt}\leq\hspace{4pt}
\mathrm{n}_{p^{2}}^{\mathrm{sup}}(L)
\hspace{24pt}.
$$
Précisons un peu. Pour obtenir cet encadrement on contemple l'expression \eqref{tp2} de $\mathrm{T}_{p^{2}}$ en fonction de $\mathrm{T}_{p}$ et $\mathrm{T}_{p}^{\psi^{2}}$ obtenue plus haut, on utilise la détermination de $\mathrm{T}_{p}$ pour $p\leq 113$ et les inégalités de Ramanujan pour les $\tau_{j,k}(p^{2})$, à savoir $\vert\tau_{6,8}(p^{2})\vert\leq 4\hspace{1pt}\hspace{1pt}p^{19}$, $\vert\tau_{8,8}(p^{2})\vert\leq 4\hspace{1pt}\hspace{1pt}p^{21}$, $\vert\tau_{12,6}(p^{2})\vert\leq 4\hspace{1pt}\hspace{1pt}p^{21}$ et $\vert\tau_{4,10}(p^{2})\vert\leq 4\hspace{1pt}\hspace{1pt}p^{21}$ (pour des inégalités de Ramanujan plus fines voir \eqref{ramfin}).

\bigskip
On constate que l'on a $\mathrm{n}_{49}^{\mathrm{inf}}(\mathrm{L}_{i})=\mathrm{n}_{49}^{\mathrm{sup}}(\mathrm{L}_{i})$ pour $i=1,2,3,4$~; on en déduit la détermination de $\mathrm{N}_{49}(\mathrm{L}_{i},\mathrm{Leech})$ pour $i=1,2,3,4$.

\bigskip
\textit{Remarque.} On constate que l'on a $\mathrm{n}_{9}(\mathrm{L}_{16})=1$~;  l'énoncé 8.2.2 est donc,  en un certain sens, optimal.

\bigskip
\textit{Les cas $11\leq p\leq 29$}

\medskip
On adapte cette fois la méthode utilisée pour déterminer $\mathrm{N}_{p}(\mathrm{L}_{i},\mathrm{Leech})$, $i=1,2,3,4$, dans les cas $61\leq p\leq 107$.

\medskip
On pose $x_{i}=\mathrm{n}_{p^{2}}(\mathrm{L}_{i})-\mathrm{n}_{p^{2}}^{\mathrm{inf}}(\mathrm{L}_{i})$~; il s'agit à nouveau de déterminer le quadruplet d'entiers naturels $(x_{1},x_{2},x_{3},x_{4})$.

\medskip
On pose $\mathrm{e}_{p^{2}}(\mathrm{L}_{i}):=\mathrm{n}_{p^{2}}^{\mathrm{sup}}(\mathrm{L}_{i})-\mathrm{n}_{p^{2}}^{\mathrm{inf}}(\mathrm{L}_{i})$, $1\leq i\leq 23$, et
$$
\hspace{24pt}
\underline{\mathrm{e}}_{p^{2}}
\hspace{4pt}:=\hspace{4pt}
(\mathrm{e}_{p^{2}}(\mathrm{L}_{1}),\mathrm{e}_{p^{2}}(\mathrm{L}_{2}),\mathrm{e}_{p^{2}}(\mathrm{L}_{3}),\mathrm{e}_{p^{2}}(\mathrm{L}_{4}))
\hspace{24pt}.
$$
On a donc, par définition même, les inégalités
$$
\hspace{8pt}
x_{1}\leq\mathrm{e}_{p^{2}}(\mathrm{L}_{1})
\hspace{8pt},\hspace{8pt}
x_{2}\leq\mathrm{e}_{p^{2}}(\mathrm{L}_{2})
\hspace{8pt},\hspace{8pt}
x_{3}\leq\mathrm{e}_{p^{2}}(\mathrm{L}_{3})
\hspace{8pt},\hspace{8pt}
x_{4}\leq\mathrm{e}_{p^{2}}(\mathrm{L}_{4})\hspace{8pt}.
\leqno{(1)}
$$

\bigskip
\textit{Le cas $p=11$}

\medskip
On calcule tout d'abord $\underline{\mathrm{e}}_{11^{2}}$~; on trouve $\underline{\mathrm{e}}_{11^{2}}=(1,1868,270,17436)$.

\medskip
On exprime ensuite les entiers $\tau_{j,k}(11^{2})$ en fonction des ``inconnues'' $x_{1},x_{2},x_{3},x_{4}$~; on trouve (le calcul utilise la détermination que nous avons déjà faite de $\mathrm{T}_{11}$)~:
$$
\hspace{24pt}
\begin{bmatrix}
\tau_{6,8}(11^{2}) \\ \tau_{8,8}(11^{2}) \\ \tau_{12,6}(11^{2}) \\ \tau_{4,10}(11^{2}) \\
\end{bmatrix}
\hspace{4pt}=\hspace{4pt}
\frac{A}{61}\hspace{4pt}
\begin{bmatrix}
x_{1} \\ x_{2} \\ x_{3} \\ x_{4} \\
\end{bmatrix}
\hspace{4pt}+\hspace{4pt}\frac{B}{61}
\hspace{24pt},
$$
$A$ et $B$ désignant deux matrices explicites à coefficients entiers, respectivement de taille $(4,4)$ et $(4,1)$ (l'apparition du nombre premier $61$ est due au fait que $61$ divise $11^2+1$).

\medskip
On constate que la réduction modulo $61$ de la matrice $A$ est inversible~; on peut donc calculer la réduction modulo $61$ du quadruplet $(x_{1},x_{2},x_{3},x_{4})$, on trouve~:
$$
(x_{1},x_{2},x_{3},x_{4})
\hspace{4pt}\equiv\hspace{4pt}
(1,52,25,15)
\hspace{8pt}\bmod{61}
\leqno{(2)}
$$
(on observera que cette congruence force déjà $x_{1}=1$).

\medskip
On exprime également les entiers $\mathrm{n}_{11^{2}}(\mathrm{L}_{i})$, $5\leq i\leq 23$, en fonction de $x_{1},x_{2},x_{3},x_{4}$, on trouve~:
$$
\mathrm{n}_{11^{2}}(\mathrm{L}_{i})
\hspace{4pt}=\hspace{4pt}
a_{i,1}\hspace{1pt}x_{1}+a_{i,2}\hspace{1pt}x_{2}+a_{i,3}\hspace{1pt}x_{3}+a_{i,4}\hspace{1pt}x_{4}+b_{i}
$$
avec $a_{i,1},a_{i,2},a_{i,3},a_{i,4},b_{i}$ des nombres rationnels. Cette égalité montre que le quadruplet $(x_{1},x_{2},x_{3},x_{4})$ satisfait une certaine congruence affine (éventuellement triviale) modulo le p.p.c.m. des dénominateurs des $a_{i,j}$, $1\leq j\leq 4$~; on note (3) l'ensemble de ces nouvelles congruences.

\medskip
L'ordinateur dit que les seuls quadruplets vérifiant (1), (2) et (3) sont~:
$$
\hspace{12pt}
(1,662,269,6481)\hspace{12pt},\hspace{12pt}(1,1333,147,6481)\hspace{12pt},\hspace{12pt}
(1,1333, 208,17217)
\hspace{12pt}.
$$

\medskip
On constate enfin que le premier et le troisième quadruplet ne passent pas le test de Ramanujan. Les $\tau_{j,k}(11^{2})$ sont déterminés.

\bigskip
\textit{Méthode alternative}

\medskip
D'après la proposition IX.\ref{intcongtaujk}, le coefficient de $t^{2}$ dans le polynôme caractéristique $\mathop{\mathrm{d\acute{e}t}}(t-p^{\frac{\mathrm{w}(\varpi_i)}{2}}\hspace{1pt}\mathrm{c}_p(\varpi_{i}))$, pour $i=7,8,9,10$ est respectivement divisible par $p^{6},p^{6},p^{4},p^{8}$, c'est-à-dire que l'on a la congruence
\begin{equation}\label{cristalline}
\hspace{24pt}
\tau_{j,k}(p)^{2}
\hspace{4pt}\equiv\hspace{4pt}
\tau_{j,k}(p^{2})
\hspace{8pt}\bmod{2\hspace{1pt}p^{\hspace{1pt}k-2}}
\hspace{24pt}.
\end{equation}
On pose $\varepsilon_{j,k}(p):=\frac{\tau_{j,k}(p)^{2}-\tau_{j,k}(p^{2})}{2}$ et on exprime ces  $\varepsilon_{j,k}(p)$ en fonction du quadruplet d'entiers naturels $(x_{i})_{1\leq i\leq 4}:=(\mathrm{n}_{p^{2}}(\mathrm{L}_{i})-\mathrm{n}_{p^{2}}^{\mathrm{inf}}(\mathrm{L}_{i}))$ (on suppose que les $\tau_{j,k}(p)$ sont connus, ce qui est le cas pour les nombres premiers qui nous intéressent ici)~; on obtient une expression de la forme
$$
\begin{bmatrix}
\varepsilon_{6,8}(p) \\ \varepsilon_{8,8}(p) \\ \varepsilon_{12,6}(p) \\ \varepsilon_{4,10}(p) \\
\end{bmatrix}
\hspace{4pt}=\hspace{4pt}
\mathrm{E}(p)
\begin{bmatrix} x_{1} \\ x_{2} \\ x_{3} \\ x_{4}
\end{bmatrix}
+\mathrm{H}(p)
$$
avec $\mathrm{E}(p)$  et $\mathrm{H}(p)$ deux matrices à coefficients rationnels respectivement de taille $(4,4)$ et $(4,1)$. Le fait que les $\varepsilon_{j,k}(p)$ sont entiers et satisfont les congruences \eqref{cristalline} impose une contrainte au quadruplet $(x_{1},x_{2},x_{3},x_{4})$ que nous notons (2-bis).

\medskip
Pour $p=11$, on constate qu'il existe un seul quadruplet qui satisfait les contraintes (1) et (2-bis) à savoir $(1,1333,147,6481)$.

\bigskip
\textit{Les cas $p=13$, $p=17$ et $p=19$}

\medskip
On a~:

\smallskip
--\hspace{8pt}$\underline{\mathrm{e}}_{13^{2}}=(655,121728,14943,1135678)$~;

\smallskip
--\hspace{8pt}$\underline{\mathrm{e}}_{17^{2}}=(536541,5855913,9346120,46438144)$~;

\smallskip
--\hspace{8pt}$\underline{\mathrm{e}}_{19^{2}}=(2884703,84510145,134879385,4993470088)$.

\smallskip
On constate à nouveau que dans les trois cas il existe un seul quadruplet $(x_{1},x_{2},x_{3},x_{4})$ qui satisfait les contraintes (1) et (2-bis) à savoir~:

\smallskip
--\hspace{8pt} $(453,50943,3642,439453)$ pour $p=13$~;

\smallskip
--\hspace{8pt} $(217661,1571118,3271290,261210371)$ pour $p=17$~;

\smallskip
--\hspace{8pt} $(964326,29790571,55543719,3055506804)$ pour $p=19$.

\bigskip
\textit{Le cas $p=23$}

\medskip
On a $\underline{\mathrm{e}}_{23^{2}}=(93365728,753181406,1202088152,161617609778)$. Cette fois il existe deux quadruplets qui satisfont (1) et (2-bis). Le seul à passer le test de Ramanujan est $(52157635,398996852,418588772, 78467649933)$.

\bigskip
\textit{Le cas $p=29$}

\medskip
On a $\underline{\mathrm{e}}_{29^{2}}=(1662796593,308516971151,492397438725,2878328193860)$. Il existe $156$ quadruplets qui satisfont (1) et (2-bis). Parmi ceux-ci, il n'y en a que $6$ qui passent le test de Ramanujan.

\medskip
Pour effectuer la sélection ultime on procède comme dans le cas $p=11$, première méthode. On exprime $\mathrm{n}_{29^{2}}(\mathrm{L}_{5})$ en fonction de $(x_{1},x_{2},x_{3},x_{4})$ et on observe que l'intégralité de $\mathrm{n}_{29^{2}}(\mathrm{L}_{5})$ implique que le quadruplet $(x_{1},x_{2},x_{3},x_{4})$ satisfait une certaine congruence modulo $256$, disons (3). On constate que (3) n'est satisfaite que par un seul des $6$ quadruplets ci-dessus à savoir
$$
(773950187, 87165709281, 106617389411, 1454026724829)
$$
(en fait ce quadruplet est le seul à satisfaire à la fois (1) et (3)).

\bigskip
\textit{Remarque.} La détermination des $\tau_{j,k}(p^{2})$ pour $p\leq 29$ permet d'expliciter l'opérateur de Hecke $\mathrm{T}_{p^{2}}:\mathbb{Z}[\mathrm{X}_{24}]\to\mathbb{Z}[\mathrm{X}_{24}]$
pour $p\leq 29$. On trouve par exemple
\begin{multline*}
\mathrm{N}_{29^{2}}(\mathrm{E}_{24},\mathrm{L}_{21})
\hspace{4pt}=\hspace{4pt} \\
9787847431870605615736000813350868753051894303124387738419200000
\end{multline*}
(approximativement $0.98\times 10^{64}$, record battu~!).

\bigskip
\textit{Le cas $p=31$}

\medskip
La méthode que nous avons employée pour $11\leq p\leq 29$ est en échec pour $p=31$. Nous expliquons pourquoi ci-après (en petits caractères).

\medskip
\footnotesize
Cette méthode peut être décrite de la façon suivante. On pose  $\mathrm{n}_{p^{2}}(\mathrm{L}_{i})=\mathrm{n}_{p^{2}}^{\mathrm{inf}}(\mathrm{L}_{i})+x_{i}$, avec $x_{i}\in\mathbb{Z}$, pour $i=1,2,3,4$~; on pose $\underline{x}=(x_{1},x_{2},x_{3},x_{4})$, $\underline{x}$ est donc \textit{a priori} un élément de $\mathbb{Z}^{4}\subset\mathbb{R}^{4}$.

\medskip
Les inégalités de Ramanujan vérifiées par les $\tau_{j,k}(p^{2})$ disent que $\underline{x}$ appartient à un  parallélotope, disons $\mathrm{Par}_{p}$, de l'espace affine $\mathbb{R}^{4}$~; la définition des entiers $\mathrm{n}_{p^{2}}^{\mathrm{inf}}(\mathrm{L}_{i})$ et $\mathrm{n}_{p^{2}}^{\mathrm{sup}}(\mathrm{L}_{i})$ est telle que la condition $\underline{x}\in\mathrm{Par}_{p}$ implique l'encadrement $0\leq x_{i}\leq\mathrm{e}_{p^{2}}(\mathrm{L}_{i})$. Le fait que les $\varepsilon_{j,k}(p)$ sont des entiers, que ces entiers vérifient les congruences \eqref{cristalline} et enfin que les $\mathrm{n}_{p^{2}}(\mathrm{L}_{i})$ sont des entiers pour $5\leq i\leq23$ (ces $\mathrm{n}_{p^{2}}(\mathrm{L}_{i})$ s'expriment comme des fonctions affines de $\underline{x}$ à coefficients rationnels) font que $\underline{x}$ appartient à un translaté, disons $\Gamma_{p}^{\mathrm{aff}}$, d'un réseau, disons $\Gamma_{p}$, de l'espace vectoriel $\mathbb{R}^{4}$, contenant $\mathbb{Z}^{4}$.

\medskip
On considère le quotient
$$
\phi(p)
\hspace{4pt}:=\hspace{4pt}
\frac{\mathop{\mathrm{volume}}\hspace{2pt}(\hspace{1pt}\mathrm{Par}_{p}\hspace{1pt})}{\mathop{\mathrm{covolume}}\hspace{2pt}(\hspace{1pt}\Gamma_{p}\hspace{1pt})}
$$
(volume et covolume pour la mesure de Lebesgue). La différence essentielle entre les cas $p=29$ et $p=31$ est la suivante~: on a $\phi(29)\approx 0.02409$ et  $\phi(31)\approx 31918.2436$. Dans le premier cas nous avons pu montrer que l'intersection  $\Gamma_{29}^{\mathrm{aff}}\cap\mathrm{Par}_{29}$ contient un seul point et déterminer ce point.  Dans le second, le calcul de $\phi(31)$ indique, heuristiquement,  que le nombre de points de $\Gamma_{31}^{\mathrm{aff}}\cap\mathrm{Par}_{31}$ est de l'ordre $32000$~; en fait ce nombre de points est~$31995$.

\medskip
Le lecteur diligent objectera que nous avons été un peu paresseux en ce qui concerne les inégalités de Ramanujan. En effet, nous avons simplement utilisé le fait que, pour\linebreak $7\leq r\leq 10$, $\theta_{r}(p^{2})$ est somme de quatre nombres complexes de module~$p^{\mathrm{w}(\varpi_{r})}$ (on rappelle que $\theta_{7},\theta_{8},\theta_{9},\theta_{10}$ est une notation alternative pour $\tau_{6,8},\tau_{8,8},\tau_{12,6},\tau_{4,10}$). Or nous connaissons ici les $\theta_{r}(p)$, $7\leq r\leq 10$~; le fait que les racines dans $\mathbb{C}$ du polynôme caractéristique $\mathop{\mathrm{d\acute{e}t}}(t-p^{\frac{\mathrm{w}(\varpi_{r})}{2}}\hspace{1pt}\mathrm{c}_p(\varpi_{r}))$ sont de module $p^{\frac{\mathrm{w}(\varpi_{r})}{2}}$ équivaut aux inégalités
\begin{equation}\label{ramfin}
\hspace{24pt}
-4\hspace{2pt}p^{\mathrm{w}(\varpi_{r})}+\frac{\hspace{3pt}\theta_{r}(p)^{2}}{2}
\hspace{4pt}\leq\hspace{4pt}
\theta_{r}(p^{2})
\hspace{4pt}\leq\hspace{4pt}
(\hspace{1pt}2\hspace{2pt}p^{\frac{\mathrm{w}(\varpi_{r})}{2}}-\vert\theta_{r}(p)\vert\hspace{1pt})^{\hspace{1pt}2}
\hspace{24pt}.
\end{equation}
Ces inégalités montrent que le point $\underline{x}$ appartient à un parallélotope, disons $\mathrm{Par}_{p}^{\mathrm{slim}}$, contenu dans $\mathrm{Par}_{p}$. On a
$$
\hspace{24pt}
\frac{\mathop{\mathrm{volume}}\hspace{2pt}(\hspace{1pt}\mathrm{Par}_{p}^{\mathrm{slim}}\hspace{1pt})}{\mathop{\mathrm{volume}}\hspace{2pt}(\hspace{1pt}\mathrm{Par}_{p}\hspace{1pt})}
\hspace{4pt}=\hspace{4pt}
\prod_{r=7}^{10}\hspace{4pt}
(\hspace{2pt}1-\frac{\vert\theta_{r}(p)\vert}{4p^{\frac{\mathrm{w}(\varpi_{r})}{2}}}\hspace{2pt})^{\hspace{1pt}2}
\hspace{24pt}.
$$
Pour $p=31$, ce rapport  est approximativement $0.2115$ et le cardinal de $\Gamma_{31}^{\mathrm{aff}}\cap\mathrm{Par}_{31}^{\mathrm{slim}}$ est~$6735$.
\normalsize


\section{Congruences à la Harder} \label{parfinharder}

\medskip
Ce paragraphe se compose de trois parties.

\medskip
Dans la première, très élémentaire, on exploite l'observation suivante~: le seul fait que les endomorphismes $\mathrm{T}_{p}$ de $\mathbb{Z}[\mathrm{X}_{24}]$, $p$ premier, aient (après extension des scalaires à $\mathbb{Q}$) une base de vecteurs propres commune, à savoir celle de~$\mathrm{T}_{2}$, entraîne que les $\lambda_{j}(p)$ vérifient de nombreuses congruences. On obtient par exemple la congruence suivante~:
$$
\hspace{24pt}
(p+1)\hspace{4pt}(\tau_{4,10}(p)-\tau_{22}(p)-p^{13}-p^{8})
\hspace{4pt}\equiv\hspace{4pt}
0\hspace{8pt}\bmod{41}
\hspace{24pt}.
$$
Dans la deuxième partie, plus subtile, ``on divise par $p+1$'' certaines des congruences obtenues dans la première en faisant notamment appel à la théorie des représentations galoisiennes. On démontre par exemple que l'on a la congruence
$$
\tau_{4,10}(p)
\hspace{4pt}\equiv\hspace{4pt}
\tau_{22}(p)+p^{13}+p^{8}
\hspace{8pt}\bmod{41}
$$
conjecturée par Günter Harder \cite{harder}. Dans la troisième partie on analyse la forme que peut prendre \textit{a priori} une décomposition en irréductibles de la représentation $\ell$-adique résiduelle associée à un $\tau_{j,k}$. On déduit de cette analyse, et des calculs que l'on a faits au paragraphe 3, que certaines d'entre elles sont irréductibles ce qui explique pourquoi les  $\tau_{j,k}$ correspondants n'apparaissent pas dans les congruences dégagées dans la deuxième partie.

\vspace{0,75cm}
\textsc{Sur certaines congruences vérifiées par les $\lambda_{j}(p)$}

\medskip
On contemple à nouveau la formule \eqref{tpformel}~:
$$
\hspace{24pt}
\mathrm{T}_{p}
\hspace{4pt}=\hspace{4pt}
\mathrm{V}\hspace{4pt}\mathrm{diag}(\lambda_{1}(p),\lambda_{2}(p),\ldots,\lambda_{24}(p))\hspace{4pt}\mathrm{V}^{-1}
\hspace{24pt}.
$$
La matrice $\mathrm{V}$ ci-dessus est à coefficients entiers mais ce n'est pas le cas de la matrice~$\mathrm{V}^{-1}$. En effet  \texttt{PARI} nous dit que l'on a
\begin{multline*}
\hspace{8pt}
\vert\mathop{\mathrm{d\acute{e}t}}\mathrm{V}\vert
\hspace{4pt}=\hspace{4pt}
2^{\hspace{1pt}220}\hspace{2pt}.\hspace{2pt}3^{\hspace{1pt}85}\hspace{2pt}.\hspace{2pt}5^{\hspace{1pt}35}\hspace{2pt}.\hspace{2pt}7 ^{\hspace{1pt}23}\hspace{2pt}.\hspace{2pt}11^{\hspace{1pt}9}\hspace{2pt}.\hspace{2pt}13^{\hspace{1pt}10}\hspace{2pt}.\hspace{2pt}17^{\hspace{1pt}3}\hspace{2pt}.\hspace{2pt}19^{\hspace{1pt}2}\hspace{2pt}.\hspace{2pt}23^{\hspace{1pt}2}\hspace{2pt}.\hspace{2pt}\\41\hspace{2pt}.\hspace{2pt}131\hspace{2pt}.\hspace{2pt}283^{\hspace{1pt}2}\hspace{2pt}.\hspace{2pt}593\hspace{2pt}.\hspace{2pt}617^{\hspace{1pt}2}\hspace{2pt}.\hspace{2pt}691^{\hspace{1pt}10}\hspace{2pt}.\hspace{2pt}3617^{\hspace{1pt}4}\hspace{2pt}.\hspace{2pt}43867^{\hspace{1pt}3}
\hspace{8pt};
\end{multline*}
\texttt{PARI} nous dit également que le plus petit des entiers $d>0$ tels que $d\hspace{2pt}\mathrm{V}^{-1}$ est à coefficients entiers est
\begin{multline*}
\hspace{24pt}
\mathrm{D}
\hspace{4pt}:=\hspace{4pt}
2^{\hspace{1pt}21}\hspace{2pt}.\hspace{2pt}3^{\hspace{1pt}10}\hspace{2pt}.\hspace{2pt}5^{\hspace{1pt}5}\hspace{2pt}.\hspace{2pt}7 ^{\hspace{1pt}2}\hspace{2pt}.\hspace{2pt}11^{\hspace{1pt}3}\hspace{2pt}.\hspace{2pt}13^{\hspace{1pt}2}\hspace{2pt}.\hspace{2pt}17\hspace{2pt}.\hspace{2pt}19\hspace{2pt}.\hspace{2pt}23^{\hspace{1pt}2}\hspace{2pt}.\hspace{2pt}\\41\hspace{2pt}.\hspace{2pt}131\hspace{2pt}.\hspace{2pt}283\hspace{2pt}.\hspace{2pt}593\hspace{2pt}.\hspace{2pt}617\hspace{2pt}.\hspace{2pt}691^{\hspace{1pt}2}\hspace{2pt}.\hspace{2pt}3617\hspace{2pt}.\hspace{2pt}43867
\hspace{24pt}.
\end{multline*}
On voit donc que pour que la matrice $\mathrm{T}_{p}$ soit à coefficient entiers, il faut que de nombreuses congruences modulo des diviseurs de $\mathrm{D}$,  impliquant les valeurs propres $\lambda_{j}(p)$, soient vérifiées. La matrice $\mathrm{V}$ en main (merci Nebe-Venkov), l'obtention de ces congruences ne fait intervenir que la théorie des modules sur les anneaux principaux. L'expression \eqref{ljpformel} des $\lambda_{j}(p)$ en fonction de $\tau_{12}(p)$, $\tau_{16}(p)$, $\tau_{18}(p)$, $\tau_{20}(p)$, $\tau_{22}(p)$, $\tau_{6,8}(p)$, $\tau_{8,8}(p)$, $\tau_{12,6}(p)$ et $\tau_{4,10}(p)$ donne alors des congruences concernant ces fonctions arithmétiques.

\bigskip
La théorie des modules sur les anneaux principaux nous dit qu'il existe deux matrices $\mathrm{R}$ et $\mathrm{S}$ de $\mathrm{GL}_{24}(\mathbb{Z})$, et des entiers strictement positifs, $\mathrm{d}_{1},\mathrm{d}_{2},\ldots,\mathrm{d}_{24}$, avec $\mathrm{d}_{j}$ divisant $\mathrm{d}_{i}$ pour $j>i$, tels que l'on a
$$
\mathrm{V}
\hspace{4pt}=\hspace{4pt}
\mathrm{R}
\hspace{4pt}\mathrm{diag}(\mathrm{d}_{1},\mathrm{d}_{2},\ldots,\mathrm{d}_{24})\hspace{4pt}\mathrm{S}^{-1}
$$
(on observera que l'on a $\mathrm{d}_{1}=\mathrm{D}$ et $\prod_{i}\mathrm{d}_{i}=\vert\mathop{\mathrm{d\acute{e}t}}\mathrm{V}\vert$).

\medskip
Les deux conditions suivantes sont équivalentes~:

\medskip
-- La matrice $\mathrm{V}\hspace{2pt}\mathrm{diag}(\lambda_{1}(p),\lambda_{2}(p),\ldots,\lambda_{24}(p))\hspace{2pt}\mathrm{V}^{-1}$ est à coefficient entiers.

\medskip
-- La matrice
$$
\mathrm{diag}(\mathrm{d}_{1},\ldots,\mathrm{d}_{24})\hspace{2pt}\mathrm{S}^{-1}\hspace{2pt}\mathrm{diag}(\lambda_{1}(p),\lambda_{2}(p),\ldots,\lambda_{24}(p))\hspace{2pt}\mathrm{S}\hspace{2pt}
\mathrm{diag}(\mathrm{d}_{1}^{-1},\ldots,\mathrm{d}_{24}^{-1})
$$
est à coefficient entiers.

\medskip
Soit $k$ un entier avec $1\leq k\leq 24$~;  on pose
$$
\mathrm{E}_{k}
\hspace{4pt}:=\hspace{4pt}
\mathrm{S}^{-1}\hspace{2pt}\mathrm{diag}(\delta_{k,1},\delta_{k,2},\ldots,\delta_{k,24})\hspace{2pt}\mathrm{S}
$$
($\delta_{-,-}$ est le symbole de Kronecker) et on note $\mathrm{e}_{i,j,k}$ le coefficient d'indice $(i,j)$ de la matrice $\mathrm{E}_{k}$ (les $\mathrm{e}_{i,j,k}$ sont des entiers ``universels'', déterminés par $\mathrm{V}$ et un choix du couple $(\mathrm{R},\mathrm{S})$). La seconde condition ci-dessus est encore équivalente à la suivante~:
\vfill\eject

\medskip
-- Pour tout couple $(i,j)$ avec $i>j$ on a la congruence
\begin{equation}\label{lkpcong}
\hspace{24pt}
\sum_{k=1}^{24}
\hspace{4pt}
\mathrm{e}_{i,j,k}
\hspace{2pt}
\lambda_{k}(p)
\hspace{4pt}\equiv\hspace{4pt}
0
\hspace{8pt}\bmod{\frac{\mathrm{d}_{j}}{\mathrm{d}_{i}}}
\hspace{24pt}.
\end{equation}

\bigskip
Conceptualisons un peu ce qui précède.

\medskip
On observe tout d'abord que la seule propriété de $\mathrm{T}_{p}$ que nous ayons utilisée ci-dessus est que les $\mathrm{v}_{j}$, $1\leq j\leq 24$, sont tous des vecteurs propres de~$\mathrm{T}_{p}$. Soient $U$ un endomorphisme de $\mathbb{Z}[\mathrm{X}_{24}]$ vérifiant cette propriété et $\lambda_{j}(U)$, $1\leq j\leq 24$,  l'entier défini par l'égalité $U(\mathrm{v}_{j})=\lambda_{j}(U)\hspace{2pt}\mathrm{v}_{j}$~; alors pour tout couple $(i,j)$ avec $i>j$ on a la congruence
\begin{equation}\label{lkUcong}
\hspace{24pt}
\sum_{k=1}^{24}
\hspace{4pt}
\mathrm{e}_{i,j,k}
\hspace{2pt}
\lambda_{k}(U)
\hspace{4pt}\equiv\hspace{4pt}
0
\hspace{8pt}\bmod{\frac{\mathrm{d}_{j}}{\mathrm{d}_{i}}}
\hspace{24pt}.
\end{equation}

\medskip
On note $\mathrm{C}$ le sous-anneau de $\mathrm{End}_{\mathbb{Z}}(\mathbb{Z}[\mathrm{X}_{24}])$ constitué des endomorphismes~$U$ considérés ci-dessus~; les applications $\lambda_{j}:\mathrm{C}\to\mathbb{Z}\hspace{2pt},\hspace{2pt}U\mapsto\lambda_{j}(U)$ sont des homomorphismes d'anneaux dont le produit, disons $\underline{\lambda}$,
$$
\mathrm{C}\to\mathbb{Z}^{24}
\hspace{24pt},\hspace{24pt}
U\mapsto (\lambda_{1}(U),\lambda_{2}(U),\ldots,\lambda_{24}(U))
$$
est un homomorphisme d'anneaux injectif (ce qui montre en particulier que l'anneau $\mathrm{C}$ est commutatif). L'image de $\underline{\lambda}$ est le sous-anneau de $\mathbb{Z}^{24}$ constitué des $24$-uples $(x_{1},x_{2},\ldots,x_{24})$ vérifiant
$$
\sum_{k=1}^{24}
\hspace{4pt}
\mathrm{e}_{i,j,k}
\hspace{2pt}
x_{k}
\hspace{4pt}\equiv\hspace{4pt}
0
\hspace{8pt}\bmod{\frac{\mathrm{d}_{j}}{\mathrm{d}_{i}}}
$$
pour tout couple $(i,j)$ avec $i>j$.

\bigskip
\textit{Remarques}

\medskip
1) Puisque les coefficients de $\mathrm{V}$ sont premiers entre eux on a $\mathrm{d}_{24}=1$ (en fait on a $\mathrm{d}_{j}=1$ pour $j\geq 21$).

\medskip
2) Posons $\underline{\mathrm{d}}:=(\mathrm{d}_{1},\mathrm{d}_{2},\ldots,\mathrm{d}_{24})$ et notons $\Gamma(\underline{\mathrm{d}})$ le sous-groupe de $\mathrm{GL}_{24}(\mathbb{Z})$ intersection dans $\mathrm{GL}_{24}(\mathbb{Q})$ de $\mathrm{GL}_{24}(\mathbb{Z})$ et $\mathrm{diag}(\underline{\mathrm{d}})\hspace{2pt}\mathrm{GL}_{24}(\mathbb{Z})\hspace{2pt}\mathrm{diag}(\underline{\mathrm{d}})^{-1}$. On consta\-te que la classe de $\mathrm{S}$ dans l'ensemble fini $\mathrm{GL}_{24}(\mathbb{Z})/\Gamma(\underline{\mathrm{d}})$ ne dépend que de $\mathrm{V}$ et que l'on peut définir le sous-anneau $\mathrm{C}$ de $\mathbb{Z}^{24}$ en termes du $24$-uple $\underline{\mathrm{d}}$ et de cette classe.

\medskip
3) Soit $\ell$ un nombre premier, alors les homomorphismes d'anneaux $\mathrm{C}\to\mathbb{Z}/\ell$, considérés comme des éléments  du $\mathbb{Z}/\ell$-espace vectoriel $\mathrm{Hom}_{\mathbb{Z}}(\mathrm{C},\mathbb{Z}/\ell)$, sont linéairement indépendants (``indépendance des caractères''). Si la valuation $\ell$-adique de $\mathrm{D}$ est $1$, en clair si $\ell$ appartient à la liste
$$
\hspace{24pt}
\{17, 19, 41, 131, 283, 593, 617, 3617, 43867\}
\hspace{24pt},
$$
alors $\mathbb{Z}_{(\ell)}\otimes_{\mathbb{Z}}\mathop{\mathrm{coker}}\underline{\lambda}$ est annulé par la multiplication par $\ell$. L'observation précédente montre dans ce cas qu'il existe une relation d'équivalence sur $\{1,2,\ldots,24\}$, disons $\mathcal{R}_{\ell}$, uniquement déterminée, telle que $\mathbb{Z}_{(\ell)}\otimes_{\mathbb{Z}}\mathrm{C}$ est le sous-anneau de $\mathbb{Z}_{(\ell)}^{24}$ constitué des $24$-uples $(x_{1},x_{2},\ldots,x_{24})$ vérifiant les congruences $x_{i}\equiv x_{j}\bmod{\ell}$ pour $i\hspace{1pt}\mathcal{R}_{\ell}\hspace{1pt}j$. Nous verrons par exemple ci-après que $\mathbb{Z}_{(41)}\otimes_{\mathbb{Z}}\mathrm{C}$ est le sous-anneau de $\mathbb{Z}_{(41)}^{24}$ constitué des $24$-uples $(x_{1},x_{2},\ldots,x_{24})$ vérifiant $x_{18}\equiv x_{21}\bmod{41}$.

\bigskip
On peut déterminer en toute généralité la classe d'isomorphisme du $(\mathbb{Z}/\mathrm{D})$-module $\mathop{\mathrm{coker}}\underline{\lambda}$ à l'aide des ``routines'' d'algèbre linéaire de \texttt{PARI} (\texttt{mathnf}, \texttt{mathnfmod} et \texttt{matsnf}). On a effectué le calcul en considérant $\mathrm{C}$ comme le sous-module de $\mathbb{Z}^{24}$ constitué des $24$-uples $(x_{1},x_{2},\ldots,x_{24})$ tels que la matrice $\sum_{j}x_{j}\hspace{1pt}\mathrm{Proj}_{j}$ est à coefficients entiers (la notation $\mathrm{Proj}_{j}$ est introduite dans la démonstration de \ref{tpapprox}). On donne le résultat ci-après~:

\bigskip
\begin{prop}\label{isocoker} Soit $\ell$ un nombre premier qui divise $\mathrm{D}$, en clair un élément de la liste
$$
\hspace{18pt}
\{2, 3, 5, 7, 11, 13, 17, 19, 23, 41, 131, 283, 593, 617, 691, 3617, 43867\}
\hspace{18pt};
$$
on a des isomorphismes de la forme
$$
\hspace{24pt}
\mathbb{Z}_{(\ell)}\otimes_{\mathbb{Z}}\mathop{\mathrm{coker}}\hspace{1pt}(\hspace{1pt}\underline{\lambda}:\mathrm{C}\to\mathbb{Z}^{24}\hspace{1pt})
\hspace{4pt}\simeq\hspace{4pt}
\mathbb{Z}/\ell^{\mathrm{e}_{\ell,1}}\times\mathbb{Z}/\ell^{\mathrm{e}_{\ell,1}}\times\ldots\times\mathbb{Z}/\ell^{\mathrm{e}_{\ell,\mathrm{r}_{\ell}}}
\hspace{24pt},
$$
$\mathrm{e}_{\ell}:= (\mathrm{e}_{\ell,1},\mathrm{e}_{\ell,2},\ldots,\mathrm{e}_{\ell,\mathrm{r}_{\ell}})$ désignant la suite finie décroissante d'entiers naturels non nuls explicitée ci-dessous~:

\medskip
--\hspace{8pt}$\mathrm{e}_{2}=(21, 19, 17, 17, 15, 15, 14, 14, 12, 12, 11, 10, 9, 9, 9, 8, 8, 7, 6, 6, 3, 1)$~;

\smallskip
--\hspace{8pt}$\mathrm{e}_{3}=(10, 9, 9, 7, 7, 6, 5, 5, 5, 5, 5, 5, 5, 5, 4, 4, 4, 3, 3, 3, 2, 1)$~;

\smallskip
--\hspace{8pt}$\mathrm{e}_{5}=(5, 5, 3, 3, 3, 2, 2, 2, 2, 2, 2, 2, 2, 2, 2, 2, 2, 1, 1, 1, 1)$~;

\smallskip
--\hspace{8pt}$\mathrm{e}_{7}=(2, 2, 2, 2, 2, 2, 2, 1, 1, 1, 1, 1, 1, 1, 1, 1, 1, 1)$~;

\smallskip
--\hspace{8pt}$\mathrm{e}_{11}=(3, 2, 1, 1, 1, 1, 1)$~;

\smallskip
--\hspace{8pt}$\mathrm{e}_{13}=(2, 2, 1, 1, 1, 1, 1, 1, 1)$~;

\smallskip
--\hspace{8pt}$\mathrm{e}_{17}=(1, 1, 1)$~;

\smallskip
--\hspace{8pt}$\mathrm{e}_{19}=(1, 1)$~;

\smallskip
--\hspace{8pt}$\mathrm{e}_{23}=(2)$~;

\smallskip
--\hspace{8pt}$\mathrm{e}_{41}=(1)$~;

\smallskip
--\hspace{8pt}$\mathrm{e}_{131}=(1)$~;

\smallskip
--\hspace{8pt}$\mathrm{e}_{283}=(1,1)$~;

\smallskip
--\hspace{8pt}$\mathrm{e}_{593}=(1)$~;

\smallskip
--\hspace{8pt}$\mathrm{e}_{617}=(1,1)$~;

\smallskip
--\hspace{8pt}$\mathrm{e}_{691}=(2, 1, 1, 1, 1, 1, 1, 1, 1)$~;

\smallskip
--\hspace{8pt}$\mathrm{e}_{3617}=(1, 1, 1,1)$~;

\smallskip
--\hspace{8pt}$\mathrm{e}_{43867}=(1, 1, 1)$.
\end{prop}

\bigskip
On décrit maintenant deux généralisations des congruences \eqref{lkUcong} (et \eqref{lkpcong}).

\medskip
1) Les congruences \eqref{lkUcong} font, \textit{a priori}, intervenir simultanément les $24$ valeurs propres $\lambda_{j}(U)$. Soit $J$ un sous-ensemble arbitraire de $\{1,2,\ldots,24\}$~; nous décrivons ci-après un algorithme à celui qui conduit à \eqref{lkUcong} pour obtenir des congruences faisant intervenir seulement les $\lambda_{j}(U)$ avec $j$ dans $J$.

\medskip
On note respectivement $\mathrm{M}_{J}$ et $\mathrm{L}_{J}$, le sous-module de $\mathbb{Z}[\mathrm{X}_{24}]$ engendré par les~$\mathrm{v}_{j}$ avec $j\in J$ et l'intersection dans $\mathbb{Q}[\mathrm{X}_{24}]$ de $\mathbb{Q}\otimes_{\mathbb{Z}}\mathrm{M}_{J}$ et $\mathbb{Z}[\mathrm{X}_{24}]$. Toujours d'après la théorie des modules sur les anneaux principaux, il existe~:

\smallskip
-- une $J\times J$ matrice $\mathrm{S}_{J}=[\mathrm{s}_{J,i,j}]_{(i,j)\in J\times J}$, à coefficients entiers, inversible avec $\mathrm{S}_{J}^{-1}$ également à coefficients entiers,

\smallskip
-- des entiers strictement positifs $\mathrm{d}_{J,i}$, $i\in J$, avec $\mathrm{d}_{J,j}$ divisant $\mathrm{d}_{J,i}$ pour $i<j$,

\smallskip
tels que l'ensemble
$$
\{\hspace{2pt}\frac{1}{d_{J,j}}\sum_{i\in J}\hspace{2pt}\mathrm{s}_{J,i,j}\hspace{2pt}\mathrm{v}_{i}\hspace{2pt}\}_{j\in J}
$$
est une base de $\mathrm{L}_{J}$ (si bien que le quotient $\mathrm{L}_{J}/\mathrm{M}_{J}$ est isomorphe à la somme directe $\bigoplus_{j\in J}\mathbb{Z}/\mathrm{d}_{J,j}$).

\medskip
Soit $(i,j,k)$ un élément de $J\times J\times J$~; on note $\mathrm{e}_{J,i,j,k}$ le coefficient d'indice $(i,j)$ de la matrice $\mathrm{S}_{J}^{-1}\hspace{2pt}\mathrm{diag}((\delta_{k,i})_{i\in J})\hspace{2pt}\mathrm{S}_{J}$. Pour tout couple $(i,j)$ avec $i>j$ on a la congruence
\begin{equation}\label{lkUcongJ}
\sum_{k\in J}
\hspace{4pt}
\mathrm{e}_{J,i,j,k}
\hspace{2pt}
\lambda_{k}(U)
\hspace{4pt}\equiv\hspace{4pt}
0
\hspace{8pt}\bmod{\frac{\mathrm{d}_{J,j}}{\mathrm{d}_{J,i}}}
\end{equation}
(on observera que la somme $\sum_{k\in J}\hspace{1pt}\mathrm{e}_{J,i,j,k}$ est nulle).

\medskip
2) Soit $m$ un diviseur de $\mathrm{D}$~; on note $\mathbb{Z}_{(m)}$ le localisé de $\mathbb{Z}$ obtenu en inversant les éléments premiers à $m$. On obtient des congruences modulo des diviseurs d'une puissance de $m$ (et de $\mathrm{D}$) en remplaçant dans le point 1 ci-dessus l'anneau principal $\mathbb{Z}$ par l'anneau principal $\mathbb{Z}_{(m)}$.

\bigskip
On considère enfin un cas particulier de ce qui précède.

\medskip
On note $\mathcal{V}$ l'ensemble $\{\mathrm{v}_{1},\mathrm{v}_{2},\ldots,\mathrm{v}_{24}\}$. Soit $\mathcal{W}$ un sous-ensemble de $\mathcal{V}$~; on note $\mathrm{J}(\mathcal{W})$ le sous-ensemble de $\{1,2,\ldots,24\}$ constitué des $j$ avec $\mathrm{v}_{j}\in\mathcal{W}$. 

\medskip
Soient $m$ un diviseur de $\mathrm{D}$ et $\rho_{m}:\mathbb{Z}[\mathrm{X}_{24}]\to(\mathbb{Z}/m)[\mathrm{X}_{24}]$ l'homomorphisme de réduction modulo $m$. Soit $\mathcal{W}$ un sous-ensemble de $\mathcal{V}$. Supposons tout d'abord $m$ premier. Si le système de vecteurs $\rho_{m}(\mathcal{W})$ est lié et si $\mathcal{W}$ est minimal parmi les sous-ensembles de $\mathcal{V}$ avec cette propriété alors, pour tout $v$ dans $\mathcal{W}$, $\rho_{m}(\mathcal{W}-\{v\})$ est une base du sous-espace vectoriel engendré par $\rho_{m}(\mathcal{W})$. Supposons maintenant que $m$ est un diviseur quelconque de $\mathrm{D}$, nous dirons plus généralement que $\mathcal{W}$ est {\em $m$-lié minimal} si le sous-module de $(\mathbb{Z}/m)[\mathrm{X}_{24}]$ engendré par $\rho_{m}(\mathcal{W})$ est un $(\mathbb{Z}/m)$-module libre et si $\rho_{m}(\mathcal{W}-\{v\})$ en est une base pour tout $v$ dans $\mathcal{W}$. Si $\mathcal{W}$ est $m$-lié minimal alors le $\mathbb{Z}_{(m)}$-module $\mathbb{Z}_{(m)}\otimes_{\mathbb{Z}}(\mathrm{L}_{\mathrm{J}(\mathcal{W})}/\mathrm{M}_{\mathrm{J}(\mathcal{W})})$ est isomorphe à $\mathbb{Z}/\widetilde{m}$ avec $\widetilde{m}$ un multiple de $m$ qui divise une puissance de $m$. Le point 2 ci-dessus fournit des congruences modulo $\widetilde{m}$ et \textit{a fortiori} modulo $m$~:

\bigskip
\begin{prop}\label{minimal}
Soit $U$ un endomorphisme de $\mathbb{Z}[\mathrm{X}_{24}]$ qui admet les $\mathrm{v}_{j}$, $1\leq j\leq 24$, pour vecteurs propres. Soient $m$ un diviseur de $\mathrm{D}$ et $\mathcal{W}$ un sous-ensemble $m$-lié minimal de $\mathcal{V}$. Alors on a les congruences
$$
\lambda_{i}(U)
\hspace{4pt}\equiv\hspace{4pt}
\lambda_{j}(U)
\hspace{8pt}\bmod{m}
$$
pour tous $i$ et $j$ dans $\mathrm{J}(\mathcal{W})$.
\end{prop}

\medskip
Pour le confort du lecteur nous donnons une démonstration \textit{ab initio} de l'énoncé ci-dessus~:

\medskip
\textit{Démonstration.} Par définition, on dispose d'une relation de dépendance de la forme $\sum_{v\in\mathcal{W}}\hspace{2pt}\mu_{v}\hspace{2pt}\rho_{m}(v)=0$ avec $\mu_{v}\in(\mathbb{Z}/m)^{\times}$ et si $\sum_{v\in\mathcal{W}}\hspace{2pt}\mu'_{v}\hspace{2pt}\rho_{m}(v)=0$ est une autre relation de dépendance alors on a $\frac{\mu'_{v}}{\mu_{v}}=\frac{\mu'_{w}}{\mu_{w}}$ pour tous $v$ et $w$ dans~$\mathcal{W}$. La proposition en résulte.
\hfill$\square$

\bigskip
Il est clair que le cardinal d'un sous-ensemble $m$-lié minimal de $\mathcal{V}$ est supérieur ou égal à $2$. La proposition suivante, que l'on vérifie par inspection, dit que l'on a souvent l'égalité, au moins si $m$ est premier.

\medskip
\begin{prop}\label{minimal2}
Soient $\ell$ un diviseur premier de $\mathrm{D}$ différent de $3,5,7,11$ et $\mathcal{W}$ un sous-ensemble $\ell$-lié minimal de $\mathcal{V}$. Alors le cardinal de $\mathcal{W}$ est $2$.
\end{prop}

\medskip
\textit{Remarque.} Soit $\mathrm{P}_{\mathrm{X}_{24}}$ le $\mathbb{Z}$-schéma dont les $A$-points ($A$ anneau commutatif unitaire) sont les facteurs directs de rang $1$ de $A[\mathrm{X}_{24}]$ ($\mathrm{P}_{\mathrm{X}_{24}}$ est donc un avatar de l'espace projectif $\mathrm{P}^{23}$). On dispose d'un sous-ensemble canonique de $\mathrm{P}_{\mathrm{X}_{24}}(\mathbb{Q})=\mathrm{P}_{\mathrm{X}_{24}}(\mathbb{Z})$, à savoir l'ensemble des classes des $\mathrm{v}_{j}$~; on le note~$[\mathcal{V}]$. La proposition ci-dessus dit que cet ensemble à $24$ éléments est très loin d'être ``générique''. En effet, elle montre que pour les $\ell$ qui apparaissent dans son énoncé, les points de $\rho_{\ell}([\mathcal{V}])$ sont ``projectivement indépendants''.

\vspace{0,75cm}
\textsc{Exemples}

\bigskip
\textit{Le cas $m=43867$}

\medskip
Les sous-ensembles $m$-liés minimaux de $\mathcal{V}$ sont $\{\mathrm{v}_{1},\mathrm{v}_{11}\}$, $\{\mathrm{v}_{2},\mathrm{v}_{8}\}$ et $\{\mathrm{v}_{3},\mathrm{v}_{6}\}$. En prenant $U=\mathrm{T}_{p}$, $p$ premier, dans la proposition \ref{minimal} on obtient les congru\-ences modulo $43867$~:
\begin{equation}\label{cong43867-1}
\hspace{6pt}
\lambda_{1}(p)\hspace{2pt}\equiv\hspace{2pt}\lambda_{11}(p)
\hspace{18pt},\hspace{18pt}
\lambda_{2}(p)\hspace{2pt}\equiv\hspace{2pt}\lambda_{8}(p)
\hspace{18pt},\hspace{18pt}
\lambda_{3}(p)\hspace{2pt}\equiv\hspace{2pt}\lambda_{6}(p)
\hspace{18pt}.
\end{equation}
On on constate que l'on a

\medskip
--\hspace{24pt}$\lambda_{11}(p)-\lambda_{1}(p)
\hspace{4pt}=\hspace{4pt}
(p^{5} + p^{4} + p^{3} + p^{2} + p + 1)\hspace{4pt}(\tau_{18}(p)-p^{17}-1)$
\hspace{24pt},

\medskip
--\hspace{24pt}$\lambda_{8}(p)-\lambda_{2}(p)
\hspace{4pt}=\hspace{4pt}
(p^{4} + p^{3} + p^{2} + p )\hspace{4pt}(\tau_{18}(p)-p^{17}-1)$
\hspace{24pt},

\medskip
--\hspace{24pt}$\lambda_{6}(p)-\lambda_{3}(p)
\hspace{4pt}=\hspace{4pt}
(p^{3} + p^{2} )\hspace{4pt}(\tau_{18}(p)-p^{17}-1)$
\hspace{24pt}.

\medskip
Comme le $\mathrm{p.g.c.d.}$ des polynômes $X^{5} + X^{4} + X^{3} + X^{2} + X + 1$, $X^{4} + X^{3} + X^{2} + X $ et $X^{3} + X^{2}$ est $X+1$, on voit que \eqref{cong43867-1} implique la congruence
\begin{equation}\label{cong43867-2}
\hspace{24pt}
(p+1)\hspace{4pt}(\tau_{18}(p)-p^{17}-1)
\hspace{4pt}\equiv\hspace{4pt}
0\hspace{8pt}\bmod{43867}
\hspace{24pt}.
\end{equation}
Cette congruence est plus faible que la congruence bien connue (voir par exemple \cite{swinnertondyer})
\begin{equation}\label{cong43867-3}
\tau_{18}(p)-p^{17}-1
\hspace{4pt}\equiv\hspace{4pt}
0\hspace{8pt}\bmod{43867}
\end{equation}
(on notera quand même que \eqref{cong43867-2} implique \eqref{cong43867-3} pour $p\not\equiv -1\bmod{43867}$~!). Nous expliquerons cependant, au cours de la démonstration du théorème \ref{recapcong}, comment l'intervention de la théorie des représentations galoisiennes et une version plus élaborée de la proposition \ref{minimal} (Proposition \ref{congcaract}) permettent d'obtenir \eqref{cong43867-3} (ce qui fournit une démonstration bien compliquée de cette congruence~!).

\medskip
\textit{Remarque.} La proposition \ref{isocoker} dit que l'on a $\mathbb{Z}_{(43867)}\otimes_{\mathbb{Z}}\mathop{\mathrm{coker}}\underline{\lambda}
\approx(\mathbb{Z}/43867)^{3}$~; ce qui précède montre en fait que $\mathbb{Z}_{(43867)}\otimes_{\mathbb{Z}}\mathrm{C}$ est le sous-anneau de $\mathbb{Z}_{(43867)}^{24}$ constitué des $24$-uples $(x_{1},x_{2},\ldots,x_{24})$ vérifiant $x_{1}\equiv x_{11}$, $x_{2}\equiv x_{8}$ et $x_{3}\equiv x_{6}$ modulo $43867$.
\vfill\eject

\bigskip
\textit{Le cas $m=3617$}

\medskip
Les sous-ensembles $m$-liés minimaux de $\mathcal{V}$ sont $\{\mathrm{v}_{1},\mathrm{v}_{13}\}$, $\{\mathrm{v}_{2},\mathrm{v}_{12}\}$, $\{\mathrm{v}_{3},\mathrm{v}_{9}\}$ et $\{\mathrm{v}_{4},\mathrm{v}_{7}\}$ . On obtient cette fois les congruences modulo $3617$~:
\begin{equation*}\label{cong3617-1}
\hspace{8pt}
\lambda_{1}(p)\hspace{2pt}\equiv\hspace{2pt}\lambda_{13}(p)
\hspace{4pt},\hspace{4pt}
\lambda_{2}(p)\hspace{2pt}\equiv\hspace{2pt}\lambda_{12}(p)
\hspace{4pt},\hspace{4pt}
\lambda_{3}(p)\hspace{2pt}\equiv\hspace{2pt}\lambda_{9}(p)
\hspace{4pt},\hspace{4pt}
\lambda_{4}(p)\hspace{2pt}\equiv\hspace{2pt}\lambda_{7}(p)
\hspace{8pt}.
\end{equation*}
Comme précédemment, on constate que les congruences ci-dessus impliquent la congruence
\begin{equation*}\label{cong3617-2}
\hspace{24pt}
(p+1)\hspace{4pt}(\tau_{16}(p)-p^{15}-1)
\hspace{4pt}\equiv\hspace{4pt}
0\hspace{8pt}\bmod{3617}
\hspace{24pt}.
\end{equation*}

\bigskip
\textit{Le cas $m=691$}

\medskip
Les sous-ensembles $m$-liés minimaux sont de cardinal $2$ et au nombre de $12$. La  considération des deux ensembles $m$-liés minimaux $\{\mathrm{v}_{1},\mathrm{v}_{24}\}$ et $\{\mathrm{v}_{2},\mathrm{v}_{23}\}$ conduit à la congruence
\begin{equation*}\label{cong691}
\hspace{24pt}
(p+1)\hspace{4pt}(\tau_{12}(p)-p^{11}-1)
\hspace{4pt}\equiv\hspace{4pt}
0\hspace{8pt}\bmod{691}
\hspace{24pt}.
\end{equation*}

\bigskip
\textit{Le cas $m\hspace{1pt}=\hspace{1pt}283\hspace{1pt}.\hspace{1pt}617$}

\medskip
Les sous-ensembles $m$-liés minimaux sont  $\{\mathrm{v}_{1},\mathrm{v}_{5}\}$ et $\{\mathrm{v}_{2},\mathrm{v}_{4}\}$. On obtient cette fois la congruence
\begin{equation*}\label{cong283.617}
\hspace{24pt}
(p+1)\hspace{4pt}(\tau_{20}(p)-p^{19}-1)
\hspace{4pt}\equiv\hspace{4pt}
0\hspace{8pt}\bmod{283\hspace{1pt}.\hspace{1pt}617}
\hspace{24pt}.
\end{equation*}

\bigskip
\textit{Le cas $m\hspace{1pt}=\hspace{1pt}131\hspace{1pt}.\hspace{1pt}593$}

\medskip
Le seul sous-ensemble $m$-lié minimal est dans ce cas $\{\mathrm{v}_{1},\mathrm{v}_{3}\}$. Comme l'on a l'égalité $\lambda_{3}(p)-\lambda_{1}(p)=(p+1)(\tau_{22}(p)-p^{21}-1)$, on obtient la congruence
\begin{equation*}\label{cong131.593}
\hspace{24pt}
(p+1)\hspace{4pt}(\tau_{22}(p)-p^{21}-1)
\hspace{4pt}\equiv\hspace{4pt}
0\hspace{8pt}\bmod{131\hspace{1pt}.\hspace{1pt}593}
\hspace{24pt}.
\end{equation*}

\bigskip
\textit{Le cas $m=41$}

\medskip
Le seul sous-ensemble $m$-lié minimal est $\{\mathrm{v}_{18},\mathrm{v}_{21}\}$. Comme l'on a l'égalité $\lambda_{21}(p)-\lambda_{18}(p)=(p+1)(\tau_{4,10}(p)-\tau_{22}(p)-p^{13}-p^{8})$, on obtient la congruence
\begin{equation}\label{cong41-1}
\hspace{12pt}
(p+1)\hspace{4pt}(\tau_{4,10}(p)-\tau_{22}(p)-p^{13}-p^{8})
\hspace{4pt}\equiv\hspace{4pt}
0\hspace{8pt}\bmod{41}
\hspace{24pt}.
\end{equation}
Comme nous l'avons déjà dit, nous verrons que la théorie des représentations galoisiennes permet de montrer que l'on a en fait
\begin{equation}\label{cong41-2}
\hspace{18pt}
\tau_{4,10}(p)-\tau_{22}(p)-p^{13}-p^{8}
\hspace{4pt}\equiv\hspace{4pt}
0\hspace{8pt}\bmod{41}
\hspace{24pt},
\end{equation}
congruence conjecturée par G. Harder \cite{harder}.
\vfill\eject

\medskip
Le sous-ensemble $\{\mathrm{v}_{18},\mathrm{v}_{21}\}$ est $m$-lié minimal pour $m\hspace{1pt}=\hspace{1pt}2^{4}\hspace{1pt}.\hspace{1pt}3\hspace{1pt}.\hspace{1pt}41$ si bien que la congruence \eqref{cong41-1} se raffine en
\begin{equation*}\label{cong41-3}
(p+1)\hspace{4pt}(\tau_{4,10}(p)-\tau_{22}(p)-p^{13}-p^{8})
\hspace{4pt}\equiv\hspace{4pt}
0\hspace{8pt}\bmod{2^{4}\hspace{1pt}.\hspace{1pt}3\hspace{1pt}.\hspace{1pt}41}
\hspace{24pt}.
\end{equation*}

\medskip
\textit{Remarque.} Compte tenu de la proposition \ref{isocoker}, ce qui précède montre que $\mathbb{Z}_{(41)}\otimes_{\mathbb{Z}}\mathrm{C}$ est le sous-anneau de $\mathbb{Z}_{(41)}^{24}$ constitué des $24$-uples $(x_{1},x_{2},\ldots,x_{24})$ vérifiant $x_{18}\equiv x_{21}\bmod{41}$.

\bigskip
\textit{Le cas $m=23$}

\medskip
Le seul sous-ensemble $m$-lié minimal est $\{\mathrm{v}_{13},\mathrm{v}_{15}\}$~; en fait  $\{\mathrm{v}_{13},\mathrm{v}_{15}\}$ est $\widetilde{m}$-lié minimal avec $\widetilde{m}=23^{2}$. On obtient cette fois la congruence
\begin{equation}
\hspace{18pt}\label{cong23-1}
(p+1)\hspace{4pt}(\hspace{2pt}\tau_{8,8}(p)-(p^{6}+1)\hspace{1pt}\tau_{16}(p)\hspace{2pt})
\hspace{4pt}\equiv\hspace{4pt}
0\hspace{8pt}\bmod{23^{2}}
\hspace{24pt}.
\end{equation}
Là encore nous verrons plus tard que l'on a en fait
\begin{equation*}\label{cong23-2}
\hspace{24pt}
\tau_{8,8}(p)-(p^{6}+1)\hspace{1pt}\tau_{16}(p)
\hspace{4pt}\equiv\hspace{4pt}
0\hspace{8pt}\bmod{23^{2}}
\hspace{24pt}.
\end{equation*}

\medskip
Par le même argument que précédemment, la congruence \eqref{cong23-1} se raffine en
\begin{equation*}\label{cong23-3}
(p+1)\hspace{4pt}(\hspace{2pt}\tau_{8,8}(p)-(p^{6}+1)\hspace{1pt}\tau_{16}(p)\hspace{2pt})
\hspace{4pt}\equiv\hspace{4pt}
0\hspace{8pt}\bmod{2^{3}\hspace{1pt}.\hspace{1pt}3^{2}\hspace{1pt}.\hspace{1pt}23^{2}}
\hspace{24pt}.
\end{equation*}

\medskip
\textit{Remarque.} Compte tenu de la proposition \ref{isocoker}, ce qui précède montre que $\mathbb{Z}_{(23)}\otimes_{\mathbb{Z}}\mathrm{C}$ est le sous-anneau de $\mathbb{Z}_{(23)}^{24}$ constitué des $24$-uples $(x_{1},x_{2},\ldots,x_{24})$ vérifiant $x_{13}\equiv x_{15}\bmod{23^{2}}$.

\bigskip
\textit{Le cas $m=19$}

\medskip
Les sous-ensembles $m$-liés minimaux sont $\{\mathrm{v}_{9},\mathrm{v}_{10}\}$ et $\{\mathrm{v}_{21},\mathrm{v}_{22}\}$.

\medskip
La considération de $\{\mathrm{v}_{9},\mathrm{v}_{10}\}$ donne
\begin{equation}\label{cong19-1}
\hspace{6pt}
(p+1)\hspace{4pt}(\hspace{2pt}\tau_{12,6}(p)-(p^{4}+p^{2})\hspace{1pt}\tau_{16}(p)+p^{2}\hspace{1pt}\tau_{18}(p)-\tau_{22}(p)\hspace{2pt})
\equiv
0\hspace{8pt}\bmod{19}
\hspace{12pt}.
\end{equation}
Comme $\{\mathrm{v}_{9},\mathrm{v}_{10}\}$ est $m$-lié minimal pour $m\hspace{1pt}=\hspace{1pt}2^{4}\hspace{1pt}.\hspace{1pt}19$ on a aussi
\begin{equation*}\label{cong19-2}
\hspace{6pt}
(p+1)\hspace{4pt}(\hspace{2pt}\tau_{12,6}(p)-(p^{4}+p^{2})\hspace{1pt}\tau_{16}(p)+p^{2}\hspace{1pt}\tau_{18}(p)-\tau_{22}(p)\hspace{2pt})
\equiv
0\hspace{8pt}\bmod{2^{4}\hspace{1pt}.\hspace{1pt}19}
\hspace{6pt}.
\end{equation*}

\medskip
La considération de $\{\mathrm{v}_{21},\mathrm{v}_{22}\}$ donne
\begin{equation*}\label{cong19-3}
\hspace{12pt}
(p+1)\hspace{4pt}(\hspace{2pt}\tau_{4,10}(p)-(p^{8}+p^{2})\hspace{1pt}\tau_{12}(p)+p^{2}\hspace{1pt}\tau_{18}(p)-\tau_{22}(p)\hspace{2pt})
\equiv0\hspace{8pt}\bmod{19}
\hspace{12pt}.
\end{equation*}
Comme $\{\mathrm{v}_{21},\mathrm{v}_{22}\}$ est $m$-lié minimal pour $m\hspace{1pt}=\hspace{1pt}2^{4}\hspace{1pt}.3^{2}\hspace{1pt}.19$ on a aussi
\begin{equation*}\label{cong19-4}
(p+1)\hspace{4pt}(\hspace{2pt}\tau_{4,10}(p)-(p^{8}+p^{2})\hspace{1pt}\tau_{12}(p)+p^{2}\hspace{1pt}\tau_{18}(p)-\tau_{22}(p)\hspace{2pt})
\equiv 0
\bmod{2^{4}\hspace{1pt}.\hspace{1pt}3^{2}\hspace{1pt}.\hspace{1pt}19}
\hspace{4pt}.
\end{equation*}

\bigskip
\textit{Le cas $m=17$}

\medskip
Les sous-ensembles $m$-liés minimaux sont $\{\mathrm{v}_{5},\mathrm{v}_{9}\}$, $\{\mathrm{v}_{15},\mathrm{v}_{17}\}$ et$\{\mathrm{v}_{19},\mathrm{v}_{20}\}$. Les congruences associées respectivement à ces paires sont~:
\begin{equation}\label{cong17-1}
(p+1)\hspace{1pt}((p^{4}+p^{2})\hspace{1pt}\tau_{16}(p)-(p^{2}+1)\hspace{1pt}\tau_{20}(p)+\tau_{22}(p)-p^{17}-p^{4})\equiv 0
\hspace{4pt}\bmod{17}
\hspace{4pt};
\end{equation}
\begin{equation}\label{cong17-2}
(p+1)(\tau_{8,8}(p)-(p^{6}+p^{4})\tau_{12}(p)+(p^{4}+p^{2})\tau_{16}(p)-(p^{2}+1)\tau_{20}(p))\equiv 0
\bmod{17}
\hspace{4pt};
\end{equation}
\begin{equation}\label{cong17-3}
(p+1)\hspace{1pt}(\tau_{6,8}(p)-(p^{6}+p^{2})\tau_{12}(p)+p^{2}\hspace{1pt}\tau_{16}(p)-\tau_{20}(p))\equiv 0
\hspace{8pt}\bmod{17}
\hspace{4pt}.
\end{equation}

\bigskip
\textit{Le cas $m=13$}

\medskip
Les sous-ensembles $m$-liés minimaux sont de cardinal $2$ et au nombre de $12$. La considération des ensembles $m$-liés minimaux $\{\mathrm{v}_{6},\mathrm{v}_{10}\}$, $\{\mathrm{v}_{9},\mathrm{v}_{15}\}$, $\{\mathrm{v}_{10},\mathrm{v}_{11}\}$, $\{\mathrm{v}_{15},\mathrm{v}_{17}\}$ et $\{\mathrm{v}_{15},\mathrm{v}_{18}\}$ conduit respectivement aux congruences suivantes~:
\begin{equation}\label{cong13-1}
\hspace{8pt}
(p+1)\hspace{1pt}(\tau_{12,6}(p)-\tau_{22}(p)-p^{17} - p^{4})
\hspace{4pt}\equiv\hspace{4pt}
0
\hspace{8pt}
\bmod{13}
\hspace{8pt};
\end{equation}
\begin{equation}\label{cong13-2}
\hspace{8pt}
(p+1)\hspace{1pt}(\tau_{8,8}(p)-\tau_{22}(p)-p^{15}-p^{6})
\hspace{4pt}\equiv\hspace{4pt}
0
\hspace{8pt}
\bmod{13}
\hspace{8pt};
\end{equation}
\begin{equation}\label{cong13-3}
\hspace{8pt}
(p+1)\hspace{1pt}(\tau_{12,6}(p)-(p^{4}+1)\hspace{1pt}\tau_{18}(p))
\hspace{4pt}\equiv\hspace{4pt}
0
\hspace{8pt}
\bmod{13}
\hspace{8pt};
\end{equation}
\begin{equation*}\label{cong13-4}
\hspace{2pt}
(p+1)(\tau_{8,8}(p)-(p^{6}+p^{4})\tau_{12}(p)+(p^{4}+p^{2})\tau_{16}(p)-(p^{2}+1)\tau_{20}(p))\equiv 0
\bmod{13}
\hspace{2pt};
\end{equation*}
\begin{equation*}\label{cong13-5}
(p+1)\hspace{1pt}(\tau_{8,8}(p)-(p^{6}+p^{4})\hspace{1pt}\tau_{12}(p)+(p^{4}+p^{2})\hspace{1pt}\tau_{16}(p)-p^{2}\hspace{1pt}\tau_{18}(p)-\tau_{22}(p))
\equiv
0
\bmod{13}.
\end{equation*}
En fait, $\{\mathrm{v}_{10},\mathrm{v}_{11}\}$ est $m$-lié pour $m=2^{5}\hspace{1pt}.\hspace{1pt}7\hspace{1pt}.\hspace{1pt}13$ si bien que la congruence \eqref{cong13-3} se raffine en
\begin{equation}\label{cong13-6}
\hspace{24pt}
(p+1)\hspace{1pt}(\tau_{12,6}(p)-(p^{4}+1)\hspace{1pt}\tau_{18}(p))
\hspace{4pt}\equiv\hspace{4pt}
0
\hspace{8pt}
\bmod{2^{5}\hspace{1pt}.\hspace{1pt}7\hspace{1pt}.\hspace{1pt}13}
\hspace{24pt}.
\end{equation}

\bigskip
\textit{Le cas $m=11$}

\medskip
Les sous-ensembles $m$-liés minimaux sont $\{\mathrm{v}_{5},\mathrm{v}_{13}\}$, $\{\mathrm{v}_{10},\mathrm{v}_{15}\}$, $\{\mathrm{v}_{14},\mathrm{v}_{16}\}$,\linebreak $\{\mathrm{v}_{14},\mathrm{v}_{19}\}$, $\{\mathrm{v}_{16},\mathrm{v}_{19}\}$, $\{\mathrm{v}_{17},\mathrm{v}_{21}\}$ et $\{\mathrm{v}_{7},\mathrm{v}_{8},\mathrm{v}_{12}\}$ (observer que ce dernier a $3$ éléments~!).

\medskip
En prenant $\mathcal{W}=\{\mathrm{v}_{7},\mathrm{v}_{8},\mathrm{v}_{12}\}$ dans la proposition \ref{minimal} on obtient
$$
\hspace{24pt}
\lambda_{7}(p)\hspace{2pt}\equiv\hspace{2pt}\lambda_{8}(p)\hspace{2pt}\equiv\hspace{2pt}\lambda_{12}(p)
\hspace{8pt}
\bmod{11}
\hspace{24pt}.
$$
En prenant $\mathcal{W}=\{\mathrm{v}_{14},\mathrm{v}_{16}\}$, $\mathcal{W}=\{\mathrm{v}_{14},\mathrm{v}_{19}\}$ et  $\mathcal{W}=\{\mathrm{v}_{17},\mathrm{v}_{21}\}$ dans cette proposition, on obtient respectivement
$$
\hspace{2pt}
\lambda_{14}(p)\equiv\lambda_{16}(p)\bmod{11}
\hspace{3pt};\hspace{3pt}
\lambda_{14}(p)\equiv\lambda_{19}(p)\bmod{11^{2}}
\hspace{3pt};\hspace{3pt}
\lambda_{17}(p)\equiv\lambda_{21}(p)\bmod{11^{2}}
\hspace{2pt}.
$$
La deuxième congruence ci-dessus s'écrit encore
$$
\hspace{24pt}
p\hspace{1pt}(p+1)\hspace{1pt}(\tau_{6,8}(p)-\tau_{20}(p)-p^{13}-p^{6})
\hspace{4pt}\equiv\hspace{4pt} 0
\hspace{8pt}
\bmod{11^{2}}
\hspace{24pt},
$$
congruence qui par inspection du cas $p=11$ implique la suivante
\begin{equation}\label{cong11-1}
\hspace{24pt}
(p+1)\hspace{1pt}(\tau_{6,8}(p)-\tau_{20}(p)-p^{13}-p^{6})
\hspace{4pt}\equiv\hspace{4pt} 0
\hspace{8pt}
\bmod{11^{2}}
\hspace{24pt}.
\end{equation}
La troisième s'écrit encore
\begin{equation}\label{cong11-2}
\hspace{12pt}
(p+1)\hspace{1pt}(\tau_{4,10}(p)-(p^{2}+1)\hspace{1pt}\tau_{20}(p)+p^{2}\hspace{1pt}\tau_{18}(p)-p^{13}-p^{8})
\hspace{2pt}\equiv\hspace{2pt} 0
\hspace{4pt}
\bmod{11^{2}}
\hspace{12pt}.
\end{equation}

\bigskip
\textit{Exemples de spécialisations de la congruence \eqref{lkUcongJ} qui échappent à la proposition \ref{minimal}}

\bigskip
-- Comme les deux ensembles $\{\mathrm{v}_{9},\mathrm{v}_{15}\}$, $\{\mathrm{v}_{15},\mathrm{v}_{17}\}$ et $\{\mathrm{v}_{17},\mathrm{v}_{18}\}$ sont $13$-liés minimaux, on a $\lambda_{9}(p)\equiv\lambda_{15}(p)\bmod{13}$, $\lambda_{15}(p)\equiv\lambda_{17}(p)\bmod{13}$ et $\lambda_{17}(p)\equiv\lambda_{18}(p)\bmod{13}$ (la première et la deuxième congruences ont été exploitées ci-dessus, la troisième ne l'a pas été car elle ne fait pas intervenir les $\tau_{j,k}$). On va analyser ce que donne la congruence \eqref{lkUcongJ} pour $J=\{9,15,17,18\}$ (et $(i,j)=(4,1))$.

\medskip
On obtient (merci \texttt{PARI})~:
\begin{multline*}
2407302\hspace{2pt}\lambda_{9}(p) -513085\hspace{2pt}\lambda_{15}(p)-482792\hspace{2pt}\lambda_{17}(p)-1411425\hspace{2pt}\lambda_{18}(p)  \\
\hspace{4pt}\equiv\hspace{4pt}
0
\hspace{8pt}
\bmod{2^{8}\hspace{1pt}.\hspace{1pt}3^{2}\hspace{1pt}.\hspace{1pt}13^{2}\hspace{1pt}.\hspace{1pt}17}
\end{multline*}
et \textit{a fortiori}
$$
\hspace{24pt}
\lambda_{9}(p)+64\hspace{2pt}\lambda_{15}(p)-89\hspace{2pt}\lambda_{17}(p)+24\hspace{2pt}\lambda_{18}(p) 
\hspace{4pt}\equiv\hspace{4pt}
0
\hspace{8pt}
\bmod{13^{2}}
\hspace{24pt}.
$$

On a
\begin{multline*}
\hspace{8pt}
\lambda_{9}(p)+64\hspace{2pt}\lambda_{15}(p)-89\hspace{2pt}\lambda_{17}(p)+24\hspace{2pt}\lambda_{18}(p) 
\hspace{4pt}=\hspace{4pt} \\
\lambda_{9}(p)-\hspace{2pt}\lambda_{15}(p)+2\hspace{2pt}\lambda_{17}(p)-2\hspace{2pt}\lambda_{18}(p) \\
-65\hspace{2pt}(\lambda_{17}(p)-\lambda_{15}(p))
+26\hspace{2pt}(\lambda_{18}(p)-\lambda_{17}(p))
\hspace{8pt}.
\end{multline*}
Comme les deux différences $\lambda_{17}(p)-\lambda_{15}(p)$ et $\lambda_{18}(p)-\lambda_{17}(p)$ sont divisibles par $13$, on obtient au bout du compte la congruence
$$
\lambda_{9}(p)-\hspace{2pt}\lambda_{15}(p)+2\hspace{2pt}\lambda_{17}(p)-2\hspace{2pt}\lambda_{18}(p)
\hspace{4pt}\equiv\hspace{4pt}
0
\hspace{8pt}
\bmod{13^{2}}
$$
soit encore
\begin{equation*}\label{cong13-7}
(p+1)
(\tau_{8,8}(p)+2\hspace{1pt}p^{2}\hspace{1pt}\tau_{18}(p)-2\hspace{1pt}(p^{2}+1)\hspace{1pt}\tau_{20}(p)+\tau_{22}(p)-p^{15}-p^{6})
\equiv
0
\bmod{13^{2}}.
\end{equation*}

\medskip
\textit{Remarque.} L'ordinateur dit que l'on a en fait
$$
\lambda_{9}(p)-\hspace{2pt}\lambda_{15}(p)+2\hspace{2pt}\lambda_{17}(p)-2\hspace{2pt}\lambda_{18}(p)
\hspace{4pt}\equiv\hspace{4pt}
0
\hspace{8pt}
\bmod{2^{4}\hspace{1pt}.\hspace{1pt}3^{2}\hspace{1pt}.\hspace{1pt}5\hspace{1pt}.\hspace{1pt}13^{2}}
$$
pour $p\leq 113$ (rappelons que nous avons calculé les $\tau_{j,k}(p)$ pour $p\leq 113$). On peut montrer que l'on a
$$
\lambda_{9}(p)-\hspace{2pt}\lambda_{15}(p)+2\hspace{2pt}\lambda_{17}(p)-2\hspace{2pt}\lambda_{18}(p)
\hspace{4pt}\equiv\hspace{4pt}
0
\hspace{8pt}
\bmod{5}
$$
pour tout $p$ de la façon ci-après. On constate que $\{\mathrm{v}_{9},\mathrm{v}_{14},\mathrm{v}_{15}\}$ et $\{\mathrm{v}_{6},\mathrm{v}_{17},\mathrm{v}_{18}\}$ sont respectivement $25$-lié minimal et $5$-lié minimal~; on a donc en particulier $\lambda_{9}(p)\equiv\lambda_{15}(p)\bmod{25}$ et $\lambda_{17}(p)\equiv\lambda_{18}(p)\bmod{5}$.

\bigskip
-- On a vu plus haut que l'on a les congruences $\lambda_{14}(p)\equiv\lambda_{19}(p)\bmod{11^{2}}$ et $\lambda_{14}(p)\equiv\lambda_{16}(p)\bmod{11}$~;  on a donc aussi
$$
\hspace{24pt}
\lambda_{19}(p)-\lambda_{14}(p)+22\hspace{2pt}(\lambda_{16}(p)-\lambda_{14}(p))
\hspace{4pt}\equiv\hspace{4pt}
0
\hspace{8pt}\bmod{11^{2}}
\hspace{24pt}.
$$
La congruence \eqref{lkUcongJ} pour $J=\{14,16,19\}$ (et $(i,j)=(3,1)$) permet de montrer, par la même méthode que ci-dessus, que l'on a en fait
\begin{equation*}\label{cong11-3}
\hspace{24pt}
\lambda_{19}(p)-\lambda_{14}(p)+22\hspace{2pt}(\lambda_{16}(p)-\lambda_{14}(p))
\hspace{4pt}\equiv\hspace{4pt}
0
\hspace{8pt}\bmod{11^{3}}
\hspace{24pt}.
\end{equation*}

\bigskip
-- En prenant $J=\{1,2,23,24\}$ (et $(i,j)=(4,1)$) on obtient cette fois
$$
\hspace{24pt}
\lambda_{1}(p)-\hspace{2pt}\lambda_{2}(p)+2\hspace{2pt}\lambda_{23}(p)-2\hspace{2pt}\lambda_{24}(p)
\hspace{4pt}\equiv\hspace{4pt}
0
\hspace{8pt}
\bmod{691}^{2}
\hspace{24pt};
$$
cette congruence n'est pas très surprenante car le premier membre est égal à $(\tau(p)-p^{11}-1)^{2}$~!
\vfill\eject

\vspace{0,75cm}
\textsc{Où l'on explique comment la théorie des représentations galoisiennes permet de ``diviser par $p+1$'' certaines des  congruences précédentes}

\bigskip
\begin{thmv}\label{recapcong}
Pour tout nombre premier $p$, les congruences suivantes sont vérifiées~:

\smallskip
{\em (1)}\hspace{8pt}$\tau_{4,10}(p)\hspace{4pt}\equiv\hspace{4pt}\tau_{22}(p)+p^{13}+p^{8}
\hspace{8pt}\bmod{41}$\hspace{4pt}{\em (conjecture de Harder \cite{harder})};

\smallskip
{\em (2)}\hspace{8pt}$\tau_{8,8}(p)\hspace{4pt}\equiv\hspace{4pt}(p^{6}+1)\hspace{1pt}\tau_{16}(p)
\hspace{8pt}\bmod{23^{2}}$\hspace{4pt};


\smallskip
{\em (3)}\hspace{8pt}$\tau_{12,6}\hspace{4pt}\equiv\hspace{4pt}(p^{4}+p^{2})\hspace{1pt}\tau_{16}(p)
\hspace{8pt}\bmod{19}$\hspace{4pt};


\smallskip
{\em (4)}\hspace{8pt}$\tau_{4,10}(p)\hspace{4pt}\equiv\hspace{4pt}(p^{8}+p^{2})\hspace{1pt}\tau_{12}(p)\hspace{8pt}\bmod{19}$\hspace{4pt};


\smallskip
{\em (5)}\hspace{8pt}$\tau_{6,8}(p)\hspace{4pt}\equiv\hspace{4pt}(p^{6}+p^{2})\hspace{1pt}\tau_{12}(p)\hspace{8pt}\bmod{17}$\hspace{4pt};

\smallskip
{\em (6)}\hspace{8pt}$\tau_{8,8}(p)\hspace{4pt}\equiv\hspace{4pt}(p^{6}+p^{4})\hspace{1pt}\tau_{12}(p)\hspace{8pt}\bmod{17}$\hspace{4pt};


\smallskip
{\em (7)}\hspace{8pt}$\tau_{8,8}(p)\hspace{4pt}\equiv\hspace{4pt}p^{8} + p^{6} + p^{3} + p
\hspace{8pt}\bmod{13}$\hspace{4pt};


\smallskip
{\em (8)}\hspace{8pt}$\tau_{12,6}(p)\hspace{4pt}\equiv\hspace{4pt}p^{8}+p^{5}+p^{4}+p\hspace{8pt}\bmod{13}$\hspace{4pt};



\smallskip
{\em (9)}\hspace{8pt}$\tau_{6,8}(p)\hspace{4pt}\equiv\hspace{4pt}p^{8} + p^{6} + p^{3} + p
\hspace{8pt}
\bmod{11}$~;

\smallskip
{\em (10)}\hspace{8pt}$\tau_{6,8}(p)\hspace{4pt}\equiv\hspace{4pt} \tau_{20}(p)+p^{13}+p^{6}
\hspace{8pt}
\bmod{11^{2}}$~;

\smallskip
{\em (11)}\hspace{8pt}$\tau_{4,10}(p)\hspace{4pt}\equiv\hspace{4pt}p^{10}+p^{8}+p^{3}+p\hspace{8pt}
\bmod{11}$~;


\smallskip
{\em (12)}\hspace{8pt}$\tau_{8,8}(p)\hspace{4pt}\equiv\hspace{4pt}\tau_{12,6}(p)\hspace{8pt}
\bmod{11}$~;

\smallskip
{\em (13)}\hspace{8pt}$\tau_{12,6}(p)\hspace{4pt}\equiv\hspace{4pt}
p^{5}+ p^{4}+ p^{2}+ p
\hspace{8pt}
\bmod{7}$~;


\smallskip
{\em (14)}\hspace{8pt}$p\hspace{1pt}\tau_{6,8}(p)\hspace{4pt}\equiv\hspace{4pt}\tau_{8,8}(p)\hspace{4pt}\equiv\hspace{4pt}\tau_{4,10}(p)\hspace{8pt}
\bmod{7}$~;

\smallskip
{\em (15)}\hspace{8pt}$\tau_{8,8}(p)\hspace{4pt}\equiv\hspace{4pt}
2\hspace{1pt}(p^{3}+p^{2})\hspace{8pt}
\bmod{5}$~;

\smallskip
{\em (16)}\hspace{8pt}$\tau_{6,8,}(p)\hspace{4pt}\equiv\hspace{4pt}
\tau_{12,6}(p)\hspace{4pt}\equiv\hspace{4pt}
\tau_{4,10}(p)\hspace{4pt}\equiv\hspace{4pt}
p^{4}+p^{3}+p^{2}+p
\hspace{8pt}
\bmod{5}$~;

\smallskip
{\em (17)}\hspace{8pt}$\tau_{j,k}(p)\hspace{4pt}\equiv\hspace{4pt}
2\hspace{1pt}(p^{2}+p)\hspace{8pt}
\bmod{3}$~;

\smallskip
{\em (18)}\hspace{8pt}$\tau_{j,k}(p)\hspace{4pt}\equiv\hspace{4pt}
0\hspace{8pt}
\bmod{2}$.

\end{thmv}

\vspace{0,75cm}
\textit{Démonstration de la congruence (1)}

\bigskip
Nous avons montré plus haut que l'on a $\lambda_{18}(p)\equiv\lambda_{21}(p)\bmod{41}$ à l'aide de la proposition \ref{minimal} et constaté en invoquant le théorème \ref{thm24} (le résultat principal du mémoire~!) que l'on a
$$
\hspace{24pt}
\lambda_{21}(p)-\lambda_{18}(p)
\hspace{4pt}=\hspace{4pt}
(p+1)\hspace{2pt}(\hspace{2pt}\tau_{4,10}(p)-(\tau_{22}(p)+p^{13}+p^{8})\hspace{2pt})
\hspace{24pt}.
$$
Nous en avons déduit la congruence \eqref{cong41-1} qu'il nous faut maintenant ``diviser par $p+1$''. Pour cela nous allons faire intervenir les $24$ représentations galoisiennes $\ell$-adiques $\rho_{i,\ell}:\mathrm{Gal}(\overline{\mathbb{Q}}/\mathbb{Q})\to\mathrm{GL}_{24}(\overline{\mathbb{Q}}_{\ell})$ semi-simples, non ramifiées hors de $\ell$, introduites dans le premier paragraphe (avec ici $\ell=41$), caractérisées par les égalités $\lambda_{i}(p)=\mathop{\mathrm{trace}}\rho_{i,\ell}(\mathrm{Frob}_{p})$ pour tout $p\not=\ell$.
\vfill\eject 

\bigskip
On a vu aussi, dans ce même paragraphe, que le polynôme caractéristique de $\rho_{i,\ell}(\mathrm{Frob}_{p})$ est à coefficients entiers (et indépendant de $\ell$) et qu'il existe une représentation continue et semi-simple $\overline{\rho}_{i,\ell}:\mathrm{Gal}(\overline{\mathbb{Q}}/\mathbb{Q})\to\mathrm{GL}_{24}(\mathbb{F}_{\ell})$, unique à isomorphisme près, qui est non ramifiée hors de $\ell$ et telle que le polynôme caractéristique de $\overline{\rho}_{i,\ell}(\mathrm{Frob}_{p})$ est la réduction modulo $\ell$ du polynôme caractéristique de $\rho_{i,\ell}(\mathrm{Frob}_{p})$. Le lemme \ref{satakecoefcar} et la proposition \ref{minimal} impliquent l'énoncé suivant~:

\medskip
\begin{prop}\label{congcaract} Soient $m$ un diviseur de $\mathrm{D}$, $\ell$ un diviseur premier de $m$, et $\mathcal{W}$ un sous-ensemble $m$-lié minimal de $\mathcal{V}$. Alors on a la congruence
$$
\mathop{\mathrm{d\acute{e}t}}(t - \rho_{i,\ell}(\gamma))
\hspace{4pt}\equiv\hspace{4pt}
\mathop{\mathrm{d\acute{e}t}}(t-\rho_{j,\ell}(\gamma))
\hspace{8pt}\bmod{m\hspace{2pt}\mathbb{Z}_{\ell}}
$$
pour tous $i, j$ avec $\mathrm{v}_{i},\mathrm{v}_{j}$ dans $\mathcal{W}$ et tout $\gamma$ dans $\mathrm{Gal}(\overline{\mathbb{Q}}/\mathbb{Q})$. En particulier, les représentations $\overline{\rho}_{i,\ell}$ et $\overline{\rho}_{j,\ell}$ sont isomorphes.
\end{prop}

\medskip
En prenant $m=\ell=41$ et $\mathcal{W}=\{\mathrm{v}_{18},\mathrm{v}_{21}\}$ dans la proposition ci-dessus on obtient  $\overline{\rho}_{18,41}\simeq\overline{\rho}_{21,41}$. Pour obtenir une équation dont l'aspect soit similaire à celui de \eqref{cong41-1} on introduit le formalisme ci-après.

\medskip
Soit $\ell$ un nombre premier~; on note $\mathrm{A}_{\ell}$ l'anneau de Grothendieck des représentations continues de dimension finie, à coefficients dans $\mathbb{F}_{\ell}$, de $\mathrm{Gal}(\overline{\mathbb{Q}}/\mathbb{Q})$, disons non ramifiées hors de $\ell$. On a donc, par définition même, $\overline{\rho}_{21,41}-\overline{\rho}_{18,41}=0$ dans $\mathrm{A}_{41}$.

\medskip
Les représentations $\ell$-adiques notées $\rho_{\Delta_{w},\ell}$, $w=11,15,17,19,21$, et $\rho_{\Delta_{w,v},\ell}$, $(w,v)=(19,7),(21,9),(21,13),(21,5)$, dans la démonstration du théorème \ref{thmgal24}, seront ici notées respectivement $\mathrm{r}_{i;\ell}$, $i=12,16,18,20,22$, et $\mathrm{r}_{j,k;\ell}$, $(j,k)=(6,8),(8,8),(12,6),(4,10)$. Avec cette notation on a $\tau_{i}(p)=\mathop{\mathrm{trace}}\mathrm{r}_{i;\ell}(\mathrm{Frob}_{p})$ et $\tau_{j,k}(p)=\mathop{\mathrm{trace}}\mathrm{r}_{j,k;\ell}(\mathrm{Frob}_{p})$ pour tout $p\not=\ell$.

\medskip
Quelques rappels~:

\smallskip
-- Les représentations $\mathrm{r}_{i;\ell}$ et $\mathrm{r}_{j,k;\ell}$ sont respectivement de dimension $2$ et $4$.

\smallskip
-- Les représentations $\mathrm{r}_{i;\ell}$ peuvent être définies sur $\mathbb{Z}_{\ell}$. Il est probable que les représentations $\mathrm{r}_{j,k;\ell}$ puissent l'être aussi (voir la remarque \ref{remcoeffrhoil})~; nous utiliserons ci-après qu'elles peuvent être définies sur la clôture intégrale de $\mathbb{Z}_{\ell}$ dans une extension finie de $\mathbb{Q}_{\ell}$ (voir la démonstration du corollaire \ref{thmgal24red}).

\smallskip
-- On note $\omega_{\ell}:\mathrm{Gal}(\overline{\mathbb{Q}}/\mathbb{Q})\to\mathbb{Z}_{\ell}^{\times}$ l'homomorphisme (la représentation $\ell$-adique de dimension $1$) défini par l'action de $\mathrm{Gal}(\overline{\mathbb{Q}}/\mathbb{Q})$ sur les racines $\ell^{\alpha}$-ième de l'unité, $\alpha\geq 1$.

\smallskip
-- On a $\mathop{\mathrm{d\acute{e}t}}\mathrm{r}_{i;\ell}=\omega_{\ell}^{i-1}$, $\mathop{\mathrm{d\acute{e}t}}\mathrm{r}_{6,8;\ell}=\omega_{\ell}^{38}$ et $\mathop{\mathrm{d\acute{e}t}}\mathrm{r}_{j,k;\ell}=\omega_{\ell}^{42}$ pour $(j,k)=(8,8),(12,6),(4,10)$.

\smallskip
-- On note $\overline{\mathrm{r}}_{i;\ell}:\mathrm{Gal}(\overline{\mathbb{Q}}/\mathbb{Q})\to\mathrm{GL}_{2}(\mathbb{F}_{\ell})$, $\overline{\mathrm{r}}_{j,k;\ell}:\mathrm{Gal}(\overline{\mathbb{Q}}/\mathbb{Q})\to\mathrm{GL}_{4}(\mathbb{F}_{\ell})$ et $\overline{\omega}_{\ell}:\mathrm{Gal}(\overline{\mathbb{Q}}/\mathbb{Q})\to\mathbb{F}_{\ell}^{\times}$ les représentations résiduelles respectivement associées aux représentations $\ell$-adiques  $\mathrm{r}_{i;\ell}$, $\mathrm{r}_{j,k;\ell}$ (voir le corollaire \ref{thmgal24red}) et $\omega_{\ell}$~; $\overline{\omega}_{\ell}$ s'identifie à l'homomorphisme défini par l'action de $\mathrm{Gal}(\overline{\mathbb{Q}}/\mathbb{Q})$ sur les racines $\ell$-ième de l'unité, homomorphisme que l'on notera aussi $\chi_{\ell}$.

\medskip
Par la suite le nombre premier $\ell$ sera fixé,  aussi ferons-nous disparaître, sauf dans le cas de $\mathrm{A}_{\ell}$, l'indice $\ell$ de la notation.

\bigskip
Revenons maintenant à la démonstration de la congruence (1).

\medskip
On prend $\ell=41$. Les égalités
\begin{multline*}
\rho_{18}
\hspace{4pt}=\hspace{4pt}
(\omega^{14}\oplus\omega^{13} \oplus \omega^{12} \oplus 2\hspace{2pt}\omega^{11} \oplus \omega^{10} \oplus \omega^{9} \oplus \omega^{8})
\hspace{4pt}\oplus\hspace{4pt} \\
(\omega^{7}\oplus \omega^{6}\oplus \omega^{5} \oplus \omega^{4})\otimes\mathrm{r}_{12} 
\hspace{4pt}\oplus\hspace{4pt}
(\omega^{3} \oplus \omega^{2})\otimes\mathrm{r}_{18}
\hspace{4pt}\oplus\hspace{4pt}
(\omega\oplus1)\otimes\mathrm{r}_{18}
\end{multline*}
et
\begin{multline*}
\rho_{21}
\hspace{4pt}=\hspace{4pt}
(\omega^{12} \oplus 2\hspace{2pt}\omega^{11} \oplus \omega^{10})
\hspace{4pt}\oplus\hspace{4pt} \\
(\omega^{7}\oplus \omega^{6}\oplus \omega^{5} \oplus \omega^{4})\otimes\mathrm{r}_{12} 
\hspace{4pt}\oplus\hspace{4pt}
(\omega^{3} \oplus \omega^{2})\otimes\mathrm{r}_{18}
\hspace{4pt}\oplus\hspace{4pt}
(\omega\oplus1)\otimes\mathrm{r}_{4,10}
\end{multline*}
impliquent que l'on a dans l'anneau de Grothendieck $\mathrm{A}_{41}$ l'égalité
$$
\overline{\rho}_{21}-\overline{\rho}_{18}
\hspace{4pt}=\hspace{4pt}
(\chi+1)\hspace{4pt}(\hspace{2pt}\overline{\mathrm{r}}_{4,10}-(\overline{\mathrm{r}}_{18}+\chi^{13}+\chi^{8})\hspace{2pt})
$$
et l'équation promise
\begin{equation}\label{cong41-1-gal}
(\chi+1)\hspace{4pt}(\hspace{2pt}\overline{\mathrm{r}}_{4,10}-(\overline{\mathrm{r}}_{18}+\chi^{13}+\chi^{8})\hspace{2pt})
\hspace{4pt}=\hspace{4pt}
0
\end{equation}
qui est le ``pendant galoisien'' de \eqref{cong41-1}. Pour ``diviser par $\chi+1$'' cette équation on utilise la proposition \ref{divisionharder} ci-dessous. Avant d'énoncer cette proposition, il nous faut faire quelques observations et introduire encore quelques notations.

\bigskip
Soit $\ell$ un nombre premier.

\medskip
Soit $\rho$ une représentation continue de dimension finie, à coefficients dans $\mathbb{F}_{\ell}$, de $\mathrm{Gal}(\overline{\mathbb{Q}}/\mathbb{Q})$, non ramifiée hors de $\ell$.  L'application $\rho\mapsto\dim\rho$ induit un homomorphisme d'anneaux que l'on note encore $\dim:\mathrm{A}_{\ell}\to\mathbb{Z}$. Le théorème de Knonecker-Weber montre que le déterminant de $\rho$ est une puissance de~$\chi$. L'application $\rho\mapsto\mathop{\mathrm{d\acute{e}t}}\rho$ induit une application $\mathrm{A}_{\ell}\to\mathrm{C}_{\chi}$, $\mathrm{C}_{\chi}$ désignant le sous-groupe de $\mathrm{A}_{\ell}^{\times}$ engendré par $\chi$, que l'on note encore $\mathop{\mathrm{d\acute{e}t}}$~; puisque $\chi$ est d'ordre $\ell-1$ le groupe $\mathrm{C}_{\chi}$ est canoniquement isomorphe à  $\mathbb{Z}/(\ell-1)$. On vérifie que l'on a $\mathop{\mathrm{d\acute{e}t}}(x+y)=\mathop{\mathrm{d\acute{e}t}}(x)\mathop{\mathrm{d\acute{e}t}}(y)$ et $\mathop{\mathrm{d\acute{e}t}}(xy)=\mathop{\mathrm{d\acute{e}t}}(x)^{\dim y}\mathop{\mathrm{d\acute{e}t}}(y)^{\dim x}$ pour tous $x$ et $y$ dans $\mathrm{A}_{\ell}$.

\medskip
Le groupe abélien sous-jacent à l'anneau commutatif $\mathrm{A}_{\ell}$ est le groupe abélien libre engendré par l'ensemble, disons $\mathcal{S}$, des classes d'isomorphisme des représentations simples. Soit $H=\sum_{S\in\mathcal{S}}n_{S}\hspace{1pt}S$, avec $n_{S}\in\mathbb{Z}$,  un élément de~$\mathrm{A}_{\ell}$~; on pose
$$
\hspace{24pt}
\Vert H\Vert
\hspace{4pt}=\hspace{4pt}
\sum_{S\in\mathcal{S}}
\vert n_{S}\vert\dim S
\hspace{24pt}.
$$
L'application $\mathrm{A}_{\ell}\to\mathbb{N}\hspace{2pt},\hspace{2pt}H\mapsto\Vert H\Vert$ est une ``norme'', en clair les propriétés suivantes sont satisfaites~:

\smallskip
--\hspace{8pt}$H=0\iff\Vert H\Vert=0$~;

\smallskip
--\hspace{8pt}$\Vert n\hspace{1pt}H\Vert=\vert n\vert\hspace{1pt}\Vert H\Vert$ pour tout $n$ dans $\mathbb{Z}$~;

\smallskip
--\hspace{8pt}$\Vert H_{1}+H_{2}\Vert\leq\Vert H_{1}\Vert+\Vert H_{2}\Vert$ pour tous $H_{1}$ et $H_{2}$ dans $\mathrm{A}_{\ell}$.

\smallskip
Soient $\rho_{+}$ et $\rho_{-}$ deux représentations de $\mathrm{Gal}(\overline{\mathbb{Q}}/\mathbb{Q})$ à coefficients dans $\mathbb{F}_{\ell}$~; on observera que l'on a les égalités $\Vert\rho_{+}\Vert=\dim\rho_{+}$, $\Vert\rho_{-}\Vert=\dim\rho_{-}$ et l'inégalité $\Vert\rho_{+}-\rho_{+}\Vert\leq\dim\rho_{+}+\dim\rho_{-}$.

\medskip
Nous en arrivons enfin à l'énoncé que nous avions en vue~:

\medskip
\begin{prop}\label{divisionharder}
Soit $\ell\not=2$ un nombre premier~; soit $H$ un élément de $\mathrm{A}_{\ell}$. Si l'on a $(\chi+1)\hspace{1pt}H=0$ alors l'entier $\Vert H\Vert$ est divisible par $\ell-1$. Si l'on a en outre $\mathop{\mathrm{d\acute{e}t}}H=1$ alors l'entier $\Vert H\Vert$ est divisible par $2\hspace{1pt}(\ell-1)$.
\end{prop}

\medskip
\textit{Démonstration.} L'action évidente du groupe $\mathrm{C}_{\chi}$ sur le groupe abélien sous-jacent à $\mathrm{A}_{\ell}$, $(\chi^{k},x)\mapsto\chi^{k}x$, préserve le sous-ensemble $\mathcal{S}$ introduit ci-dessus. Soit $S$ un élément de $\mathcal{S}$~; on note $\Omega(S)$ l'orbite de $S$ sous l'action de $\mathrm{C}_{\chi}$, $\mathrm{m}(S)$ le cardinal de cette orbite (en clair, $\mathrm{m}(S)$ est le plus petit entier $k\geq 1$ tel que l'on a $\chi^{k}\hspace{1pt}S=S$) et $\mathbb{Z}[\Omega(S)]$ le sous-groupe (abélien libre) du groupe abélien sous-jacent à $\mathrm{A}_{\ell}$ qu'elle engendre. On a donc une décomposition en somme directe du groupe abélien sous-jacent à $\mathrm{A}_{\ell}$
$$
\hspace{24pt}
\mathrm{A}_{\ell}
\hspace{4pt}=\hspace{4pt}
\bigoplus_{S\in\mathcal{S}_{0}}
\mathbb{Z}[\Omega(S)]
\hspace{24pt},
$$
$\mathcal{S}_{0}\subset\mathcal{S}$ désignant un système de représentants pour l'ensemble quotient $\mathrm{C}_{\chi}\backslash\mathcal{S}$. Cette décomposition est compatible avec l'action de $\mathrm{C}_{\chi}$~; en particulier chaque facteur est envoyé dans lui-même par la multiplication par $\chi+1$.

\medskip
\begin{prop}\label{noyauharder}
Soit $S$ un élément de $\mathcal{S}$. Soient $\Omega(S)$ l'orbite de $S$ sous l'action de $\mathrm{C}_{\chi}$ et $\mathrm{m}(S)$ le cardinal de cette orbite, c'est-à-dire le plus entier $k\geq 1$ tel que l'on a $\chi^{k}\hspace{1pt}S=S$.

\smallskip
{\em (a)} L'entier $\mathrm{m}(S)$ divise $\ell-1$ et $\ell-1$ divise $\mathrm{m}(S)\dim S$.

\smallskip
{\em (b)} Le noyau de l'endomorphisme du groupe abélien $\mathbb{Z}[\Omega(S)]$ induit par la multiplication par $\chi+1$ est trivial si $\mathrm{m}(S)$ est impair et est engendré (comme groupe abélien) par
$$
(1-\chi+\chi^{2}-\chi^{3}+\ldots-\chi^{\hspace{1pt}\mathrm{m}(S)-1})\hspace{2pt}S
$$
si $\mathrm{m}(S)$ est pair.
\end{prop}

\medskip
\textit{Démonstration.} Le seul point qui n'est pas tout à fait évident est la seconde partie de (a). Pour s'en convaincre, observer que l'on a $\chi^{\mathrm{m}(S)}S=S$ et $\mathop{\mathrm{d\acute{e}t}}(\chi^{\mathrm{m}(S)}S)=\chi^{\mathrm{m}(S)\dim S}\mathop{\mathrm{d\acute{e}t}}S$.
\hfill$\square$

\bigskip
On reprend la démonstration de la proposition \ref{divisionharder}. Soit $\mathcal{S}_{0,0}$ le sous-ensemble de $\mathcal{S}_{0}$ constitué des $S$ avec $\mathrm{m}(S)$ pair~; la proposition ci-dessus montre que si l'on a $(\chi+1)H=0$ alors il existe des entiers relatifs~$n_{S}$, $S$ parcourant $\mathcal{S}_{0,0}$, tels que l'on a
\begin{equation}\label{noyauharderbis}
\hspace{24pt}
H
\hspace{4pt}=\hspace{4pt}
\sum_{S\in\mathcal{S}_{0,0}}
n_{S}\hspace{2pt}(1-\chi+\chi^{2}-\chi^{3}+\ldots-\chi^{\hspace{1pt}\mathrm{m}(S)-1})\hspace{2pt}S
\hspace{24pt}.
\end{equation}
Par définition même de $\Vert H\Vert$ on a
$$
\hspace{24pt}
\Vert H\Vert
\hspace{4pt}=\hspace{4pt}
\sum_{S\in\mathcal{S}_{0,0}}
\vert n_{S}\vert\hspace{2pt}\mathrm{m}(S)\dim S
\hspace{24pt}.
$$

La seconde partie du point (a) de la proposition \ref{noyauharder} dit que tous les produits $\mathrm{m}(S)\dim S$ ci-dessus sont divisibles par $\ell-1$, ce qui démontre la première partie de la proposition \ref{divisionharder}.

\medskip
Passons à la démonstration de la seconde partie de la proposition \ref{divisionharder}. L'égalité \eqref{noyauharderbis} implique la suivante
$$
\mathop{\mathrm{d\acute{e}t}}H
\hspace{4pt}=\hspace{4pt}
\chi^{-\frac{1}{2}
\sum_{S\in\mathcal{S}_{0,0}}
n_{S}\hspace{2pt}\mathrm{m}(S)\dim S}
$$
(observer que l'on a $\dim(1-\chi+\chi^{2}-\chi^{3}+\ldots-\chi^{\hspace{1pt}\mathrm{m}(S)-1})=0$ et $\mathop{\mathrm{d\acute{e}t}}(1-\chi+\chi^{2}-\chi^{3}+\ldots-\chi^{\hspace{1pt}\mathrm{m}(S)-1})=\chi^{-\frac{\mathrm{m}(S)}{2}}$), si bien que l'égalité $\mathop{\mathrm{d\acute{e}t}}H=1$ est équivalente à la congruence
$$
\hspace{24pt}
\sum_{S\in\mathcal{S}_{0,0}}
n_{S}\hspace{2pt}\mathrm{m}(S)\dim S
\hspace{4pt}\equiv\hspace{4pt}
0
\hspace{8pt}\bmod{2\hspace{1pt}(\ell-1)}
\hspace{24pt}.
$$
Comme $\vert n_{S}\vert$ et $n_{S}$ ont même parité et que tous les $\mathrm{m}(S)\dim S$ sont divisibles par $\ell-1$, les égalités $(\chi+1)H=0$ et $\mathop{\mathrm{d\acute{e}t}}H=1$ entraînent bien que $\Vert H\Vert$ est divisible par $2\hspace{1pt}(\ell-1)$.
\hfill
$\square\square$

\bigskip
\textit{Remarque.} L'égalité $(\chi+1)H=0$ implique $\dim H=0$ et $(\mathop{\mathrm{d\acute{e}t}}H)^{\hspace{1pt}2}=1$. Cette implication est ``optimale'' (pour $\ell\not=2$). Pour s'en convaincre prendre $H=H_{0}:=1-\chi+\chi^{2}-\chi^{3}+\ldots-\chi^{\hspace{1pt}\ell-2}$ et constater que l'on a $\mathop{\mathrm{d\acute{e}t}}H_{0}=\chi^{\frac{\ell-1}{2}}$. On observera également que l'on a $\Vert H_{0}\Vert=\ell-1$, ce qui montre que la première partie de la proposition \ref{divisionharder} est optimale~;  en outre, on a $(\chi+1)\hspace{1pt}(2\hspace{1pt}H_{0})=0$,  $\mathop{\mathrm{d\acute{e}t}}(2\hspace{1pt}H_{0})=1$ et  $\Vert 2\hspace{1pt}H_{0}\Vert=2\hspace{1pt}(\ell-1)$, ce qui montre la seconde partie de cette proposition est aussi optimale.

\bigskip
\textit{Démonstration de la congruence (1) à l'aide de \ref{divisionharder}.}

\medskip
On prend $\ell=41$ et $H=\overline{\mathrm{r}}_{4,10}-(\overline{\mathrm{r}}_{18}+\chi^{13}+\chi^{8})$. On a $\Vert H\Vert\leq 8$. Comme l'on a $(\chi+1)H=0$ (équation \eqref{cong41-1-gal}), la proposition \ref{divisionharder} dit que $\Vert H\Vert$ est divisible par $40$ (et même $80$ si l'on observe que l'on a $\mathop{\mathrm{d\acute{e}t}}H=1$). On en déduit $\Vert H\Vert=0$ et $H=0$. En évaluant les représentations $\overline{\mathrm{r}}_{4,10}$ et $\overline{\mathrm{r}}_{18}\oplus\chi^{13}\oplus\chi^{8}$,\linebreak ``en la classe de conjugaison $\mathrm{Frob}_{p}$'', on obtient la congruence (1) pour $p\not=41$. Par ailleurs la congruence \eqref{cong41-1} implique trivialement la congruence (1) pour $p=41$.\hfill$\square\square\square$

\vspace{0,75cm}
\textit{Démonstration de la congruence (2)}

\bigskip
On fixe $\ell=23$.

\medskip
En prenant $m=\ell=23$ et $\mathcal{W}=\{\mathrm{v}_{13},\mathrm{v}_{15}\}$ dans la proposition \ref{congcaract}, on obtient comme précédemment l'isomorphisme de représentations galoisiennes
\begin{equation}\label{hypothesedulemmeartinien}
\overline{\mathrm{r}}_{8,8}
\hspace{4pt}\simeq\hspace{4pt}
(\chi^{6}\oplus 1)\otimes\overline{\mathrm{r}}_{16}
\end{equation}
et la congruence
$$
\hspace{24pt}
\tau_{8,8}(p)
\hspace{4pt}\equiv\hspace{4pt}
(p^{6}+1)\hspace{1pt}\tau_{16}(p)
\hspace{8pt}\bmod{23}
\hspace{24pt}.
$$
Pour obtenir la congruence (2) (qui raffine à la fois la congruence ci-dessus et la congruence \eqref{cong23-1}) on utilise le lemme ci-dessous :

\medskip
\begin{lemme}\label{lemmeartinien}
Soient $B$ un anneau local artinien de corps résiduel $k$, $G$ un
groupe, $V_{1}$, $V_{2}$, $W_{1}$ et $W_2$ des $B[G]$-modules que l'on suppose libres de dimension finie comme $B$-modules.  On fait les hypothèses suivantes~:

\medskip
\begin{itemize}
\item[(i)] Pour $i=1,2$, les semi-simplifiés des $k[G]$-modules $k\otimes_{B}V_{i}$ et $k\otimes_{B}W_{i}$ sont isomorphes.

\smallskip
\item[(ii)] Les $k[G]$-modules $k\otimes_{B}V_{1}$ et $k\otimes_{B}V_{2}$ n'ont
aucun facteur de Jordan-Hölder en commun.

\smallskip
\item[(iii)] Pour tout $g$ dans $G$, on a $\mathop{\mathrm{d\acute{e}t}}(t-g_{\vert V_1\oplus V_2})=\mathop{\mathrm{d\acute{e}t}}(t -g_{\vert W_1\oplus W_2})$. 
\end{itemize}

\medskip
Alors pour $i=1,2$ et pour tout $g$, on a $\mathop{\mathrm{d\acute{e}t}}(t-g_{\vert V_i}) = \mathop{\mathrm{d\acute{e}t}}(t-g_{\vert W_i})$.
\end{lemme}
\vfill\eject

\bigskip
\textit{Démonstration.} Soient $U$ le $B[G]$-module $V_{1}\oplus V_{2}\oplus W_{1 }\oplus W_{2}$ et $R$ la $B$-algèbre image de $B[G]$ dans $\mathrm{End}_{B}(U)$.  Soit $J$ le radical de Jacobson de $R$.  Le $B$-module sous-jacent à $R$
étant de type fini, $J$ est le plus grand idéal bilatère nilpotent de $R$. C'est aussi le noyau de l'homomorphisme naturel de $R$ vers les endomorphismes du semi-simplifié de $k\otimes_{B}U$.  En particulier, on a $\mathfrak{m}R\subset J$, $\mathfrak{m}$ désignant l'idéal maximal de $B$, et $R/J$ est une
$k$-algèbre semi-simple de dimension finie.

\medskip
La théorie d'Artin-Wedderburn appliquée à $R/J$ et l'hypothèse (ii)
montrent que l'on peut trouver un idempotent $f $ dans $R/J$ tel que $f$ agisse
par l'identité sur le semi-simplifié de $k\otimes_{B}V_{1}$ et par $0$ sur
celui de $k\otimes_{B}V_{2}$.  Comme $J$ est nilpotent, cet idempotent se relève
en un idempotent $e$ dans $R$.  Ceci assure que $e$ agit par $0$  
sur $V_{2}$ (car il agit ainsi sur tous ses facteurs de Jordan-Hölder) et   
par l'identité sur $V_{1}$ et $W_{1}$ (pour la même raison).

\medskip
L'identité classique d'Amitsur, exprimant les coefficients du polynôme
caractéristique d'une somme de deux matrices comme une fonction
universelle des coefficients des polynômes caractéristiques de certains
monômes en ces deux matrices, montre que l'égalité de déterminants 
de l'hypothèse (iii) entraîne plus généralement $\mathop{\mathrm{d\acute{e}t}} (t-r_{\vert V_{1}\oplus V_{2}}) = \mathop{\mathrm{d\acute{e}t}}(t -r_{\vert W_{1}\oplus W_{2}})$ pour tout $r$ dans $R$. Le théorème en résulte~: soit $g$ un élément de $G$~; pour $i=1$ (resp.  $i=2$) on spécialise cette identité à $r=g\hspace{1pt}e$ (resp.  $r= g\hspace{1pt}(1-e)$).
\hfill$\square$

\bigskip
\textit{Démonstration de la congruence (2) à l'aide de \ref{lemmeartinien}.}

\medskip
On spécialise le lemme en question.

\medskip
(Rappel~: Le nombre premier $\ell$ est fixé égal à $23$, les notations $\mathrm{r}_{8,8}$, $\mathrm{r}_{16}$, $\omega$, $\rho_{13}$, $\rho_{15}$ et $\chi$ qui apparaissent ci-dessous sont les abréviations de $\mathrm{r}_{8,8;23}$, $\mathrm{r}_{16;23}$, $\omega_{23}$, $\rho_{13,23}$, $\rho_{15,23}$ et $\chi_{23}$.)

\medskip
La représentation $\mathrm{r}_{8,8}$ peut être définie sur la clôture intégrale, disons $\mathcal{D}$, de~$\mathbb{Z}_{23}$ dans une extension finie de $\mathbb{Q}_{23}$. La représentation $\mathrm{r}_{16}$ peut être définie sur $\mathbb{Z}_{23}$ et \textit{a fortiori} sur $\mathcal{D}$~; de même la représentation $\omega$ est définie sur $\mathbb{Z}_{23}$ et \textit{a fortiori} sur $\mathcal{D}$.

\medskip
-- On prend pour $B$ l'anneau quotient $\mathcal{D}/23^{2}$. L'anneau $B$ est local, son corps résiduel $k$ est un corps fini de caractéristique $23$~; $B$ est artinien (il est fini~!).

\medskip
-- On prend pour $G$ le groupe de Galois $\mathrm{Gal}(\overline{\mathbb{Q}}/\mathbb{Q})$.

\medskip
-- On prend pour $V_{1}$, $V_{2}$, $W_{1}$ et $W_{2}$, le $B$-module $B^{4}$ muni de l'action linéaire du groupe $G$ donnée respectivement par les représentations $(\omega^{6}\oplus 1)\otimes\mathrm{r}_{16}$, $\omega\otimes(\omega^{6}\oplus 1)\otimes\mathrm{r}_{16}$, $\mathrm{r}_{8,8}$ et $\omega\otimes\mathrm{r}_{8,8}$.

\medskip
-- L'hypothèse (i) est donnée par l'isomorphisme \eqref{hypothesedulemmeartinien} et l'extension des scalaires de $\mathbb{F}_{23}$ à $k$.

\medskip
-- Il n'est pas difficile de vérifier l'hypothèse (ii). La représentation (résiduelle modulo $23$) $\overline{r}_{16}$ est simple \cite{swinnertondyer} (pour s'en convaincre \textit{ab initio}, observer que l'on a $47\equiv 1\bmod{23}$ et $\tau_{16}(47)\not\equiv 2\bmod{23}$)~; il en résulte, toujours d'après Kronecker-Weber, que $k\otimes_{\mathbb{F}_{23}}\overline{r}_{16}$ est encore simple. Les facteurs de Jordan-Hölder de $k\otimes_{B}V_{1}$ (resp. $k\otimes_{B}V_{2}$), sont  donc $k\otimes_{\mathbb{F}_{23}}\overline{r}_{16}$ et $k\otimes_{\mathbb{F}_{23}}(\chi^{6}\hspace{2pt}\overline{r}_{16})$ (resp. $k\otimes_{\mathbb{F}_{23}}(\chi\hspace{2pt}\overline{r}_{16})$ et $k\otimes_{\mathbb{F}_{23}}(\chi^{7}\hspace{2pt}\overline{r}_{16})$). On conclut en observant que le déterminant de $\overline{r}_{16}$, $\chi^{6}\hspace{2pt}\overline{r}_{16}$, $\chi\hspace{2pt}\overline{r}_{16}$ et $\chi^{7}\hspace{2pt}\overline{r}_{16}$ est respectivement $\chi^{15}$, $\chi^{5}$, $\chi^{17}$ et $\chi^{7}$.

\medskip
- L'hypothèse (iii) est impliquée par la proposition \ref{congcaract} et le fait que l'on a $\rho_{13}=v\oplus\sigma$ et $\rho_{15}=w\oplus\sigma$, avec $v=(\omega\oplus 1)\otimes(\omega^{6}\oplus 1)\otimes\mathrm{r}_{16}$, $w=(\omega\oplus 1)\otimes\mathrm{r}_{8,8}$ et $\sigma$ une représentation $23$-adique de dimension $16$.

\bigskip
La conclusion du lemme dit que ``les polynômes caractéristiques en $\mathrm{Frob}_{p}$'', $p\not=23$, des représentations $23$-adiques $(\omega^{6}\oplus1)\otimes\mathrm{r}_{16}$ et $\mathrm{r}_{8,8}$ sont congrus modulo~$23^{2}$. \textit{A fortiori} ``les traces en $\mathrm{Frob}_{p}$'', $p\not=23$, sont congrues modulo~$23^{2}$, en d'autres termes la congruence (2) est satisfaite pour $p\not=23$. Le cas $p=23$ découle trivialement de la congruence \eqref{cong23-1}.
\hfill$\square\square$

\vspace{0,75cm}
\textit{Démonstration de la congruence (3)}

\bigskip
On fixe $\ell=19$.

\medskip
 En prenant $m=\ell=19$ et $\mathcal{W}=\{\mathrm{v}_{9},\mathrm{v}_{10}\}$ dans la proposition \ref{congcaract} et en utilisant la proposition \ref{divisionharder} on obtient l'égalité
 $$
 \overline{\mathrm{r}}_{12,6}-(\chi^{4}+\chi^{2})\hspace{2pt}\overline{\mathrm{r}}_{16}+\chi^{2}\hspace{2pt}\overline{\mathrm{r}}_{18}-\overline{\mathrm{r}}_{22}
 \hspace{4pt}=\hspace{4pt}
 0
 $$
 dans l'anneau de Grothendieck $\mathrm{A}_{19}$, soit encore l'isomorphisme de représentations
\begin{equation}\label{gal19-1}
\hspace{24pt}
 \overline{\mathrm{r}}_{12,6}\oplus\chi^{2}\hspace{2pt}\overline{\mathrm{r}}_{18}
 \hspace{4pt}\simeq\hspace{4pt}
 (\chi^{4}\oplus\chi^{2})\hspace{2pt}\overline{\mathrm{r}}_{16}\oplus\overline{\mathrm{r}}_{22}
\hspace{24pt}.
\end{equation}
Comme les représentations $\overline{\mathrm{r}}_{16}$ et $\overline{\mathrm{r}}_{22}$ sont simples (on peut s'en convaincre en observant que l'on a $\tau_{16}(5)\not\equiv\tau_{16}(43)\bmod{19}$ et $\tau_{22}(5)\not\equiv\tau_{22}(43)\bmod{19}$ alors que l'on a $5\equiv 43\bmod{19}$), la représentation $\chi^{2}\hspace{2pt}\overline{\mathrm{r}}_{18}$ est nécessairement isomorphe à l'une des représentations  $\chi^{4}\hspace{2pt}\overline{\mathrm{r}}_{16}$,  $\chi^{2}\hspace{2pt}\overline{\mathrm{r}}_{16}$ ou  $\overline{\mathrm{r}}_{22}$. Le calcul des déterminants montre que la seule possibilité est
\begin{equation}\label{gal19-2}
\hspace{24pt}
\chi^{2}\hspace{2pt}\overline{\mathrm{r}}_{18}
\hspace{4pt}\simeq\hspace{4pt}
\overline{\mathrm{r}}_{22}
\hspace{24pt}.
\end{equation}
Comme les représentations $ \overline{\mathrm{r}}_{12,6}$ et $(\chi^{4}\oplus\chi^{2})\hspace{2pt}\overline{\mathrm{r}}_{16}$ sont semi-simples, les isomorphismes \eqref{gal19-1} et \eqref{gal19-2} entraînent
$$
\hspace{24pt}
\overline{\mathrm{r}}_{12,6}
\hspace{4pt}\simeq\hspace{4pt}
(\chi^{4}\oplus\chi^{2})\hspace{2pt}\overline{\mathrm{r}}_{16}
\hspace{24pt}.
$$
Cet isomorphisme implique la congruence (3) pour $p\not=19$. On vérifie que cette congruence est aussi satisfaite pour $p=19$ (on peut éviter d'utiliser le calcul de $\tau_{12,6}(19)$ en observant  que les congruences $19^{2}\tau_{18}(19)\equiv\tau_{20}(19)$ et \eqref{cong19-1} impliquent (3) pour $p=19$).
\hfill$\square$

\vspace{0,75cm}
\textit{Démonstration des congruences (4), (5) et (6)}

\bigskip
La démonstration de la congruence (4) est la même que celle de la con\-gru\-ence~(3)~; elle est d'ailleurs plus rapide si l'on utilise \eqref{gal19-2}.

\bigskip
Passons à (5) et (6). Evidemment on fixe $\ell=17$.

\medskip
-- En utilisant le ``pendant galoisien'' de \eqref{cong17-3} et en copiant la démonstration de la congruence (3), on obtient les isomorphismes de représentations
\begin{equation}\label{gal17-1}
\chi^{2}\hspace{2pt}\overline{\mathrm{r}}_{16}
\hspace{4pt}\simeq\hspace{4pt}
\overline{\mathrm{r}}_{20}
\end{equation}
et
$$
\hspace{24pt}
\overline{\mathrm{r}}_{6,8}
\hspace{4pt}\simeq\hspace{4pt}
(\chi^{6}\oplus\chi^{2})\otimes\overline{\mathrm{r}}_{12}
\hspace{24pt}.
$$
Cet isomorphisme implique la congruence (5) pour $p\not=17$~; le cas $p=17$ peut se régler comme dans la cas de la congruence (3).

\medskip
-- Le pendant galoisien de \eqref{cong17-2} est l'équation suivante dans $\mathrm{A}_{17}$~:
$$
\hspace{24pt}
(\chi+1)\hspace{2pt}(\overline{\mathrm{r}}_{8,8}-(\chi^{6}+\chi^{4})\hspace{1pt}\overline{\mathrm{r}}_{12}+(\chi^{4}+\chi^{2})\hspace{1pt}\overline{\mathrm{r}}_{16}-(\chi^{2}+1)\hspace{1pt}\overline{\mathrm{r}}_{20})
\hspace{4pt}=\hspace{4pt}
 0
\hspace{24pt}.
$$
Cette equation et l'isomorphisme \eqref{gal17-1} entraînent
$$
\hspace{24pt}
(\chi+1)\hspace{2pt}(\overline{\mathrm{r}}_{8,8}-(\chi^{6}+\chi^{4})\hspace{1pt}\overline{\mathrm{r}}_{12})
\hspace{4pt}=\hspace{4pt}
 0
\hspace{24pt}.
$$
En invoquant la proposition \ref{divisionharder} on obtient l'équation
$$
\hspace{24pt}
\overline{\mathrm{r}}_{8,8}-(\chi^{6}+\chi^{4})\hspace{1pt}\overline{\mathrm{r}}_{12}
\hspace{4pt}=\hspace{4pt}
 0
\hspace{24pt},
$$
soit encore l'isomorphisme de représentations
$$
\overline{\mathrm{r}}_{8,8}
\hspace{4pt}\simeq\hspace{4pt}
(\chi^{6}\oplus\chi^{4})\otimes\overline{\mathrm{r}}_{12}
$$
(les deux membres sont semi-simples). Cet isomorphisme implique la con\-gru\-ence (6) pour $p\not=17$~; le cas $p=17$ peut se régler comme précédemment.

\medskip
\textit{Remarque.} Le pendant galoisien de \eqref{cong17-1} est l'équation suivante dans $\mathrm{A}_{17}$~:
$$
\hspace{24pt}
(\chi+1)\hspace{2pt}((\chi^{4}+\chi^{2})\hspace{1pt}\overline{\mathrm{r}}_{16}-(\chi^{2}+1)\hspace{1pt}\overline{\mathrm{r}}_{20}+\overline{\mathrm{r}}_{22}-\chi^{17}-\chi^{4})
\hspace{4pt}=\hspace{4pt}
0
\hspace{24pt}.
$$
Compte tenu de \eqref{gal17-1}, celle-ci donne l'équation
$$
\hspace{24pt}
(\chi+1)\hspace{2pt}(\overline{\mathrm{r}}_{22}-\chi^{17}-\chi^{4})
\hspace{4pt}=\hspace{4pt}
0
\hspace{24pt}.
$$
En invoquant à nouveau la proposition \ref{divisionharder} on obtient l'équation
$$
\hspace{24pt}
\overline{\mathrm{r}}_{22}-\chi^{17}-\chi^{4}
\hspace{4pt}=\hspace{4pt}
0
\hspace{24pt},
$$
soit encore l'isomorphisme de représentations
$$
\hspace{24pt}
\overline{\mathrm{r}}_{22}
\hspace{4pt}\simeq\hspace{4pt}
\chi^{17}\oplus\chi^{4}
\hspace{4pt}=\hspace{4pt}
\chi\oplus\chi^{4}
\hspace{24pt}.
$$
L'isomorphisme $\overline{\mathrm{r}}_{22}\simeq\chi\oplus\chi^{4}$ ci-dessus est l'un des isomorphismes donnés par Swinnerton-Dyer dans \cite{swinnertondyer}~; nous faisons un usage systématique de ce type d'isomorphisme ci-après.

\vspace{0,75cm}
\textit{Démonstration des congruences (7), (8), (9), (11) et (13)}

\bigskip
Le point  de départ de ces démonstrations est respectivement~:

\smallskip
-- l'équation dans $\mathrm{A}_{13}$ qui est le pendant galoisien de la congruence \eqref{cong13-2}~;

\smallskip
-- l'équation dans $\mathrm{A}_{13}$ qui est le pendant galoisien de la congruence \eqref{cong13-1} ou de la congruence \eqref{cong13-3}~;

\smallskip
-- l'équation dans $\mathrm{A}_{11}$ qui est le pendant galoisien de la réduction modulo $11$ de la congruence \eqref{cong11-1} (qui est une congruence modulo $11^{2}$)~;

\smallskip
-- l'équation dans $\mathrm{A}_{11}$ qui est le pendant galoisien de la réduction modulo $11$ de la congruence \eqref{cong11-2} (qui est une congruence modulo $11^{2}$)~;

\smallskip
-- l'équation dans $\mathrm{A}_{7}$ qui est le pendant galoisien de la réduction modulo $7$ de la congruence \eqref{cong13-6} (qui est une congruence modulo $2^{5}\hspace{1pt}.\hspace{1pt}7\hspace{1pt}.\hspace{1pt}13$ dont la réduction modulo $13$ est la congruence \eqref{cong13-3} mentionnée ci-dessus).

\medskip
Par la méthode que nous avons répétitivement employée ci-dessus, on exprime les $\overline{r}_{j,k}$ qui nous intéressent en fonction de certains  $\overline{r}_{i}$ et de $\chi$. Tous les $\overline{r}_{i}$ qui apparaissent peuvent à leur tour s'exprimer en fonction de $\chi$ grâce aux isomorphismes de Swinnerton-Dyer \cite{swinnertondyer}. Au bout du compte, on obtient des isomorphismes de la forme $\overline{r}_{j,k}\simeq\chi^{a_{1}}\oplus\chi^{a_{2}}\oplus\chi^{a_{3}}\oplus\chi^{a_{4}}$ qui conduisent aux congruences (7), (8), (9), (11) et (13).

\medskip
Traitons par exemple le cas de la congruence (13).

\medskip
On fixe $\ell=7$.   En prenant $m=\ell=7$ et $\mathcal{W}=\{\mathrm{v}_{10},\mathrm{v}_{11}\}$ dans la proposition \ref{congcaract} on obtient l'équation suivante dans $\mathrm{A}_{7}$~:
$$
\hspace{24pt}
(\chi+1)\hspace{1pt}(\overline{\mathrm{r}}_{12,6}-(\chi^{4}+1)\hspace{1pt}\overline{\mathrm{r}}_{18})
\hspace{4pt}=\hspace{4pt}
0
\hspace{24pt}.
$$
Puisque les représentations $\overline{\mathrm{r}}_{12,6}$ et $(\chi^{4}\oplus 1)\otimes\overline{\mathrm{r}}_{18}$ ont même déterminant (à savoir $\chi^{42}$) et que l'on a l'inégalité $\Vert\overline{\mathrm{r}}_{12,6}-(\chi^{4}+1)\hspace{1pt}\overline{\mathrm{r}}_{18}\Vert\leq 8$, la seconde partie de la proposition \ref{divisionharder} montre que l'on a en fait l'équation
$$
\hspace{24pt}
\overline{\mathrm{r}}_{12,6}-(\chi^{4}+1)\hspace{1pt}\overline{\mathrm{r}}_{18}
\hspace{4pt}=\hspace{4pt}
0
\hspace{24pt},
$$
ou ce qui revient au même l'isomorphisme de représentations
$$
\hspace{24pt}
\overline{\mathrm{r}}_{12,6}
\hspace{4pt}\simeq\hspace{4pt}
(\chi^{4}\oplus 1)\otimes\overline{\mathrm{r}}_{18}
\hspace{24pt}.
$$
Or Swinnerton-Dyer nous dit que l'on a $\overline{\mathrm{r}}_{18}\simeq\chi\oplus\chi^{4}$ si bien que l'on obtient au bout du compte
$$
\hspace{24pt}
\overline{\mathrm{r}}_{12,6}
\hspace{4pt}\simeq\hspace{4pt}
\chi^{5}\oplus\chi^{4}\oplus\chi^{2}\oplus\chi\hspace{24pt}.
$$
Cet isomorphisme donne la congruence (13) pour $p\not=7$~; le cas $p=7$ peut se déduire de \eqref{cong13-6}.

\vspace{0,75cm}
\textit{Démonstration de la congruence (10)}

\bigskip
Elle est semblable à celle la congruence (2).

\medskip
On fixe $\ell=11$.

\medskip
En prenant $m=\ell=11$ et $\mathcal{W}=\{\mathrm{v}_{14},\mathrm{v}_{19}\}$ dans la proposition \ref{congcaract} on obtient l'équation suivante dans $\mathrm{A}_{11}$~:
$$
\hspace{24pt}
\chi\hspace{1pt}(\chi+1)\hspace{2pt}(\overline{\mathrm{r}}_{6,8}-\overline{\mathrm{r}}_{20}-\chi^{13}-\chi^{6})
\hspace{24pt}.
$$
On ``divise cette équation par $\chi\hspace{1pt}(\chi+1)$'' en observant que $\chi$ est inversible et en utilisant la proposition \ref{divisionharder}~; on obtient ainsi l'isomorphisme de représentations galoisiennes
\begin{equation*}\label{gal11}
\overline{\mathrm{r}}_{6,8}
\hspace{4pt}\simeq\hspace{4pt}
\overline{\mathrm{r}}_{20}\oplus\chi^{13}\oplus\chi^{6}
\end{equation*}
et la congruence
\begin{equation}\label{divcong11}
\hspace{24pt}
\tau_{6,8}(p)
\hspace{4pt}\equiv\hspace{4pt}
\tau_{20}(p)+p^{13}+p^{6}
\hspace{8pt}\bmod{11}
\hspace{24pt}.
\end{equation}
Cette congruence se transforme en la congruence (9) en utilisant la congruence $\tau_{20}(p)\equiv p^{8}+p\bmod{11}$ de \cite{swinnertondyer}.

\medskip
Or nous avons vu, dans la première partie de ce paragraphe, que l'on a, pour tout nombre premier $p$, la congruence \eqref{cong11-1} suivante 
$$
\hspace{24pt}
(p+1)\hspace{1pt}(\tau_{6,8}(p)-\tau_{20}(p)-p^{13}-p^{6})
\hspace{4pt}\equiv\hspace{4pt} 0
\hspace{8pt}
\bmod{11^{2}}
\hspace{24pt}.
$$
Pour ``diviser cette congruence par $p+1$'' on emploie la méthode qui nous a permis d'obtenir la congruence (2) à partir de la congruence \eqref{cong23-1} (à savoir la mise en oeuvre du lemme \ref{lemmeartinien}).

\medskip
\textit{Remarque.} On ne peut avoir de congruence de la forme
$$
\tau_{6,8}(p)
\hspace{4pt}\equiv\hspace{4pt}
p^{a_{1}}+p^{a_{2}}+p^{a_{3}}+p^{a_{4}}
\hspace{8pt}
\bmod{11^{2}}
$$
pour tout nombre premier $p$, disons avec $p\not=11$. En effet, compte tenu de~(10), on aurait $\tau_{20}(p)\equiv p^{a_{1}}+p^{a_{2}}+p^{a_{3}}+p^{a_{4}}-p^{13}-p^{6}\bmod{11^{2}}$. Or cette congruence est en défaut pour le plus petit nombre premier $p$ avec $p\equiv 1\bmod{11^{2}}$, à savoir $p=727$~: $\tau_{20}(727)\equiv  68\not\equiv 2\bmod{11^{2}}$.

\medskip
\textit{Remarque.} Nous avons vu dans la première partie de ce paragraphe que l'on~a, pour tout nombre premier $p$, la congruence \eqref{cong11-2} suivante~:
$$
\hspace{24pt}
(p+1)\hspace{1pt}(\tau_{4,10}(p)-(p^{2}+1)\hspace{1pt}\tau_{20}(p)+p^{2}\hspace{1pt}\tau_{18}(p)-p^{13}-p^{8})
\hspace{4pt}\equiv\hspace{4pt}
0
\hspace{8pt}
\bmod{11^{2}}
\hspace{24pt}.
$$
Le pendant galoisien de la  réduction modulo $11$ de cette congruence, est l'équation suivante dans $\mathrm{A}_{11}$~:
$$
\hspace{24pt}
(\chi+1)\hspace{1pt}(\overline{\mathrm{r}}_{4,10}-(\chi^{2}+1)\hspace{1pt}\overline{\mathrm{r}}_{20}+\chi^{2}\hspace{1pt}\overline{\mathrm{r}}_{18}-\chi^{13}-\chi^{8})
\hspace{4pt}=\hspace{4pt}
0
\hspace{24pt}.
$$
En utilisant les isomorphismes $\overline{\mathrm{r}}_{20}\simeq\chi^{8}\oplus\chi$ et $\overline{\mathrm{r}}_{18}\simeq\chi^{6}\oplus\chi$ de \cite{swinnertondyer}, on obtient l'équation
$$
(\chi+1)\hspace{1pt}(\overline{\mathrm{r}}_{4,10}-\chi^{8}-\chi^{3}-\chi-1)
\hspace{4pt}=\hspace{4pt}
0
$$
qui après ``division par $p+1$'' conduit à la congruence (11). Mais cette fois le lemme \ref{lemmeartinien} ne permet pas de ``diviser par $p+1$'' la congruence \eqref{cong11-2}  car l'hypothèse concernant les facteurs de Jordan-Hölder n'est pas satisfaite. On constate cependant que la congruence
$$
\tau_{4,10}(p)
\hspace{4pt}\equiv\hspace{4pt}
(p^{2}+1)\hspace{1pt}\tau_{20}(p)-p^{2}\hspace{1pt}\tau_{18}(p)+p^{13}+p^{8}
\hspace{8pt}
\bmod{11^{2}}
$$
est vérifiée pour $p\leq 113$ (on rappelle que l'on a déterminé les $\tau_{j,k}(p)$ pour $p\leq 113$)~; on observera que ceci n'est vraiment une information que pour les nombres premiers $p\leq 113$ avec $p+1\equiv 0\bmod{11}$, à savoir $43$ et $109$~!

\vspace{0,75cm}
\textit{Démonstration des congruences (12) et (14)}

\bigskip
Les démonstrations de ces trois congruences sont toutes du même type (c'est pourquoi nous les avons mises dans le même chapeau). Démontrons par exemple que l'on a, pour tout nombre premier $p$, la congruence
\begin{equation}\label{newcong7}
\hspace{24pt}
p\hspace{1pt}\tau_{6,8}(p)
\hspace{4pt}\equiv\hspace{4pt}
\tau_{8,8}(p)
\hspace{8pt}
\bmod{7}
\hspace{24pt}.
\end{equation}

\medskip
On fixe $\ell=7$.

\medskip
En prenant $m=\ell=7$ et $\mathcal{W}=\{\mathrm{v}_{15},\mathrm{v}_{19}\}$ dans la proposition \ref{congcaract} on obtient l'équation suivante dans $\mathrm{A}_{7}$~:
\begin{multline}\label{gal7-1}
\hspace{24pt}
(\chi+1)\hspace{1pt}\overline{\mathrm{r}}_{8,8}-\chi(\chi+1)\hspace{1pt}\overline{\mathrm{r}}_{6,8} \\ +(\chi^{5}+\chi^{2})\hspace{1pt}\overline{\mathrm{r}}_{16}-(\chi^{5}+1)\hspace{1pt}\overline{\mathrm{r}}_{12}-\mathop{\mathrm{Sym}^{2}}\overline{\mathrm{r}}_{12}+\chi^{5}+2\hspace{1pt}\chi^{2}
\hspace{4pt}=\hspace{4pt}
0
\hspace{24pt}.
\end{multline}
Or on constate que les isomorphismes $\overline{\mathrm{r}}_{16}\simeq\chi^{2}+\chi$ et $\overline{\mathrm{r}}_{12}\simeq\chi^{4}+\chi$ de \cite{swinnertondyer} entraînent
$$
(\chi^{5}+\chi^{2})\hspace{1pt}\overline{\mathrm{r}}_{16}-(\chi^{5}+1)\hspace{1pt}\overline{\mathrm{r}}_{12}-\mathop{\mathrm{Sym}^{2}}\overline{\mathrm{r}}_{12}+\chi^{5}+2\hspace{1pt}\chi^{2}
\hspace{4pt}=\hspace{4pt}
0
$$
si bien que l'équation \eqref{gal7-1} se simplifie en la suivante~:
\begin{equation}\label{gal7-2}
\hspace{24pt}
(\chi+1)\hspace{1pt}\overline{\mathrm{r}}_{8,8}-\chi(\chi+1)\hspace{1pt}\overline{\mathrm{r}}_{6,8}
\hspace{4pt}=\hspace{4pt}
0
\hspace{24pt}.
\end{equation}
Posons $H=\overline{\mathrm{r}}_{8,8}-\chi\hspace{1pt}\overline{\mathrm{r}}_{6,8}$~; l'équation ci-dessus dit que l'on a $(\chi+1)\hspace{1pt}H=0$ et on constate que l'on a $\mathop{\mathrm{d\acute{e}t}}H=1$. On peut donc appliquer la seconde partie de la proposition \ref{divisionharder}~: $\Vert H\Vert$ est divisible par $12$. Comme l'on a \textit{a priori} $\Vert H\Vert\leq 8$, on en déduit l'égalité $H=0$, l'isomorphisme de représentations
\begin{equation}\label{gal7-3}
\hspace{24pt}
\chi\hspace{2pt}\overline{\mathrm{r}}_{6,8}
\hspace{4pt}\simeq\hspace{4pt}
\overline{\mathrm{r}}_{8,8}
\end{equation}
et la congruence \eqref{newcong7} pour $p\not=7$. Le cas $p=7$ est laissé au lecteur.

\vspace{0,75cm}
\textit{Démonstration des congruences (15), (16), (17) et (18)}

\bigskip
A nouveau, les démonstrations de ces congruences sont toutes du même type. Nous donnons quelques détails sur la démonstration des congruences (16)  et nous bornons à indiquer les modifications essentielles à apporter à celle-ci pour avoir les démonstrations de (15), (17) et (18).

\medskip
On fixe $\ell=5$.

\medskip
On considère le sous-ensemble $\mathcal{W}:=\{\mathrm{v}_{10},\mathrm{v}_{17},\mathrm{v}_{19},\mathrm{v}_{21}\}$ de $\mathcal{V}$~; on vérifie que~$\mathcal{W}$ est $5$-lié minimal. On en déduit, toujours grâce à la proposition \ref{congcaract}, que les quatre représentations $\overline{\rho}_{10}$, $\overline{\rho}_{17}$, $\overline{\rho}_{19}$ et $\overline{\rho}_{21}$ sont deux à deux isomorphes. L'isomorphisme $\overline{\rho}_{10}\simeq\overline{\rho}_{17}$ donne l'équation suivante dans $\mathrm{A}_{5}$~:
\begin{multline}\label{gal5-1}
\hspace{24pt}
(\chi+1)\hspace{1pt}\overline{\mathrm{r}}_{12,6}
\hspace{4pt}=\hspace{4pt} 
(\chi^{3}+\chi^{2}+\chi+1)\hspace{1pt}\overline{\mathrm{r}}_{20}-(\chi^{3}+\chi^{2})\hspace{1pt}\overline{\mathrm{r}}_{18} \\ 
+(\chi^{3}+\chi^{2}+\chi+1)\hspace{1pt}\overline{\mathrm{r}}_{12}-2\hspace{1pt}\chi^{3}-\chi^{2}-1
\hspace{24pt}.
\end{multline}
Compte tenu de \cite{swinnertondyer}, cette équation devient
\begin{equation}\label{gal5-2}
\hspace{24pt}
(1+\chi)\hspace{1pt}\overline{\mathrm{r}}_{12,6}
\hspace{4pt}=\hspace{4pt} 
2\hspace{1pt}
(1+\chi+\chi^{2}+\chi^{3})
\hspace{24pt}.
\end{equation}
Cette équation montre que la représentation (semi-simple) $(1\oplus\chi)\otimes\overline{\mathrm{r}}_{12,6}$ est somme directe de puissance de $\chi$~; il en résulte qu'il en est de même pour~$\overline{\mathrm{r}}_{12,6}$. On en déduit que l'on a dans $\mathrm{A}_{5}$ une équation de la forme
\begin{equation}\label{gal5-3}
\overline{\mathrm{r}}_{12,6}
\hspace{4pt}=\hspace{4pt} 
a_{0}+a_{1}\hspace{1pt}\chi+a_{2}\hspace{1pt}\chi^{2}+a_{3}\hspace{1pt}\chi^{3}
\end{equation}
avec $a_{k}$, $k=0,1,2,3$, des entiers vérifiant $a_{k}\geq 0$ et $a_{0}+a_{1}+a_{2}+a_{3}=\nolinebreak 1$. L'équation \ref{gal5-2} peut être réécrite de la façon suivante~:
$$
\hspace{24pt}
(1+\chi)\hspace{1pt}(\overline{\mathrm{r}}_{12,6}-(1+\chi+\chi^{2}+\chi^{3}))
\hspace{4pt}=\hspace{4pt} 
0
\hspace{24pt};
$$
cette réécriture et le point (b) de la proposition \ref{noyauharder} (prendre $S=\chi$) montrent qu'il existe un entier relatif $n$ tel que l'on a
$$
\hspace{24pt}
\overline{\mathrm{r}}_{12,6}
\hspace{4pt}=\hspace{4pt} 
1+\chi+\chi^{2}+\chi^{3}
+
n\hspace{1pt}(1-\chi+\chi^{2}-\chi^{3})
\hspace{24pt}.
$$
Les inégalités $a_{k}\geq 0$ montrent que l'on a $\vert n\vert\leq 1$. Le calcul du déterminant des deux membres montre que l'on a $n\equiv 0\bmod{2}$. On a donc $n=0$ et un isomorphisme de représentations
\begin{equation*}\label{gal5-4}
\hspace{24pt}
\overline{\mathrm{r}}_{12,6}
\hspace{4pt}\simeq\hspace{4pt} 
1\oplus\chi\oplus\chi^{2}\oplus\chi^{3}
\hspace{24pt}.
\end{equation*}
Cet isomorphisme implique la congruence $\tau_{12,6}(p)\equiv 1+p+p^{2}+p^{3}\bmod{5}$, ou ce qui revient au même $\tau_{12,6}(p)\equiv p^{4}+p^{3}+p^{2}+p\bmod{5}$, pour $p\not=5$. Le cas $p=5$ est laissé au lecteur.

\medskip
Les isomorphismes $\overline{\rho}_{19}\simeq\overline{\rho}_{17}$ et $\overline{\rho}_{21}\simeq\overline{\rho}_{17}$ conduisent pareillement aux congruences $\tau_{6,8}(p)\equiv p^{4}+p^{3}+p^{2}+p\bmod{5}$ et $\tau_{4,10}(p)\equiv p^{4}+p^{3}+p^{2}+p\bmod{5}$. (On observera que l'isomorphisme $\overline{\rho}_{19}\simeq\overline{\rho}_{17}$ conduit naturellement à l'isomorphisme $\chi\hspace{1pt}\overline{\mathrm{r}}_{6,8}\simeq 1\oplus\chi\oplus\chi^{2}\oplus\chi^{3}$ mais que l'on a $\chi^{-1}\hspace{1pt}(1\oplus\chi\oplus\chi^{2}\oplus\chi^{3})\cong 1\oplus\chi\oplus\chi^{2}\oplus\chi^{3}$.)

\bigskip
La congruence (15) s'obtient quant à elle en prenant par exemple $\mathcal{W}=\{\mathrm{v}_{3},\mathrm{v}_{15}\}$.

\bigskip
Passons enfin aux congruences (17) et (18). On peut les démontrer en prenant  $\ell=2,3$, $m=6$ et $\mathcal{W}=\{\mathrm{v}_{6},\mathrm{v}_{10}\},\{\mathrm{v}_{6},\mathrm{v}_{15}\},\{\mathrm{v}_{6},\mathrm{v}_{19}\},\{\mathrm{v}_{6},\mathrm{v}_{21}\}$ dans la proposition \ref{congcaract}.

\medskip
On rappelle que quand le nombre premier $\ell$ est $2$ ou $3$ les isomorphismes de \cite{swinnertondyer} qui nous concernent prennent une forme particulièrement simple : $\overline{\mathrm{r}}_{i}\simeq 1\oplus 1$ pour $\ell=2$ et $\overline{\mathrm{r}}_{i}\simeq 1\oplus\chi$ pour $\ell=3$.

\medskip
Pour $\ell=2$, on obtient par les méthodes précédentes
\begin{equation}\label{gal2}
\hspace{24pt}
\overline{\mathrm{r}}_{j,k}
\hspace{4pt}\simeq\hspace{4pt}
1\oplus 1\oplus 1\oplus 1
\hspace{24pt}.
\end{equation}
Pour $\ell=3$, on trouve que les représentations $\overline{\mathrm{r}}_{j,k}$ sont chacune isomorphe à l'une des trois représentations suivantes, $1\oplus 1\oplus 1\oplus 1$, $1\oplus 1\oplus\chi\oplus\chi$ ou $\chi\oplus\chi\oplus\chi\oplus\chi$. On lève l'ambiguité en utilisant le fait que la représentation duale $\overline{\mathrm{r}}_{j,k}^{*}$ est isomorphe à la représentation $\chi\hspace{1pt}\overline{\mathrm{r}}_{j,k}$ (voir le début de la remarque \ref{pairingsp4})~:
\begin{equation}\label{gal3}
\overline{\mathrm{r}}_{j,k}
\hspace{4pt}\simeq\hspace{4pt}
1\oplus 1\oplus\chi\oplus\chi
\hspace{24pt}.
\end{equation}

\vspace{0,75cm}
\textsc{Sur la décomposition en irréductibles des $\overline{\mathrm{r}}_{j,k;\ell}$}

\bigskip
Comme nous venons de le voir, les congruences (12) et (14) du théorème  \ref{recapcong} sont conséquences d'isomorphismes entre certaines représentations de la forme $\overline{\mathrm{r}}_{j,k;\ell}$ ou $\chi\otimes\overline{\mathrm{r}}_{j,k;\ell}$. Chacune des autres congruences de ce théorème, à l'exception de (2) et (10), est la manifestation d'une propriété de réducti\-bilité d'une représentation $\overline{\mathrm{r}}_{j,k;\ell}$. La proposition \ref{conjserrered} ci-dessous, sans doute bien connue, décrit de façon exhaustive les différentes possibilités pour une réduction d'une représentation de ce type ; elle explique en partie la structure des congruences que nous avons dégagées.

\medskip
Soit $\ell$ un nombre premier. Soit $\kappa$ un entier, on note $\mathrm{R}_{\kappa,\ell}$ l'ensemble (fini) des classes d'isomorphisme de représentations irréductibles $\mathrm{Gal}(\overline{\mathbb{Q}}/\mathbb{Q})\to\mathrm{GL}_2(\overline{\F}_\ell)$ de la forme $\rhob_{\pi,\iota}$, où $\pi\in \Pi_{\mathrm{cusp}}(\mathrm{PGL}_{2})$ est la représentation automorphe engendrée par un élément de $\mathrm{S}_{\kappa}(\mathrm{SL}_2(\mathbb{Z}))$ propre pour les opérateurs
de Hecke. On rappelle que l'on note $\nu:\mathrm{GSp}_{2g} \rightarrow\mathbb{G}_{\mathrm{m}}$ l'homomorphisme ``facteur de similitude'' (voir \S II.1). Enfin, si $S$ est une $\overline{\mathbb{F}}_{\ell}$-représentation irréductible de dimension
finie de $\mathrm{Gal}(\overline{\mathbb{Q}}/\mathbb{Q})$, on note $\mathrm{m}(S)$ le plus petit entier $m\geq 1$ tel que l'on a $\chi^{m}\otimes S \simeq S$ (on rappelle que l'on a posé $\chi=\overline{\omega_\ell}$) ; cette notation est en accord avec celle que l'on a introduite dans la proposition~\ref{noyauharder}.

\begin{prop}\label{conjserrered} Soient $\ell$ un nombre premier impair et
$$
r : \mathrm{Gal}(\overline{\mathbb{Q}}/\mathbb{Q}) \rightarrow 
\mathrm{GSp}_4(\overline{\F}_\ell)
$$
une représentation semi-simple, continue, et non ramifiée hors de $\ell$. Notons $w$ l'élément de $\mathbb{Z}/(\ell-1)$ tel que l'on a $\nu\circ r =\chi^{w}$, et supposons $w\equiv 1\bmod 2$.
\vfill\eject

\medskip
Alors on est dans un, et un seul, des cas suivants~:

\ps\begin{itemize}\ps
\item[(i)] il existe $a$ et $b$ dans $\mathbb{Z}/(\ell-1)$ tels que l'on a $r \simeq
\chi^{a}\oplus\chi^{b}\oplus\chi^{w-a}\oplus\chi^{w-b}$~;\ps\ps
\item[(ii)] il existe $\kappa\leq\ell+1$, une représentation $\rho$ dans $\mathrm{R}_{\kappa;\ell}$ et $a,b$ dans $\mathbb{Z}/(\ell-1)$ avec $2a+\kappa-1 \equiv w \bmod{\mathrm{m}(\rho)}$, tels que l'on a $r\simeq(\chi^{a}\otimes\rho)\oplus
\chi^{b}\oplus\chi^{w-b}$~; \ps\ps
\item[(iii)$_1$] il existe $\kappa\leq\ell+1$, une représentation $\rho$ dans $\mathrm{R}_{\kappa,\ell}$ et $a$ dans $\mathbb{Z}/(\ell-1)$ avec $2a+\kappa-1\not\equiv w \bmod{\mathrm{m}(\rho)}$, tels que l'on a $r \simeq(\chi^{a}\oplus\chi^{w-a-\kappa+1})\otimes\rho$~;\ps\ps
\item[(iii)$_2$] il existe une $\overline{\F}_\ell$-représentation $\rho$, irréductible
de dimension $2$, de déterminant $\mathop{\mathrm{d\acute{e}t}}\rho = \chi^a$ avec $a\in 2\mathbb{Z}$ et $a\not \equiv w \bmod \mathrm{m}(\rho)$, telle que l'on a $r\simeq (1\oplus\chi^{w-a})\otimes\rho$~;\ps\ps
\item[(iv)] pour $i=1,2$, il existe $\kappa_{i}\leq\ell+1$, une représentation $\rho_{i}$ dans $\mathrm{R}_{\kappa_i;\ell}$ et $a_{i}$ dans $\mathbb{Z}/(\ell-1)$
avec $2a_i+\kappa_{i}-1\equiv w \bmod{\mathrm{m}(\rho_{i})}$, tels que l'on a $r\simeq(\chi^{a_{1}}\otimes\rho_{1})\oplus(\chi^{a_{2}} \otimes \rho_2)$~;\ps\ps
\item[(v)] $r$ est irréductible.\ps\ps
\end{itemize}
\end{prop}

\medskip
\textit{Démonstration.} Soit $V$ est une $\overline{\mathbb{F}}_{\ell}$-représentation de dimension finie de $G=\mathrm{Gal}(\overline{\mathbb{Q}}/\mathbb{Q})$~; on note $\mathrm{i}(V)$ la représentation $V^{*}\otimes\chi^{w}$. Il est clair que l'on $V\cong\mathrm{i}(\mathrm{i}(V))$, et que l'application $V\mapsto\mathrm{i}(V)$ définit une auto-équivalence, exacte et contravariante, de la catégorie des $\overline{\mathbb{F}}_{\ell}$-représentations de dimension finie de $G$. En particulier, l'ensemble fini $\mathrm{J}(r)$ des facteurs de Jordan-Hölder de la représentation $r$ est stable par l'involution
$\mathrm{i}$.

\medskip
Observons que les éléments de $\mathrm{J}(r)$ sont de dimension $1$ ou $2$. En effet, supposons $r\simeq V\oplus W$ avec $V$ irréductible de dimension $3$ (et $W$ de dimension~$1$). La forme bilinéaire alternée non dégénérée dont est muni, par définition, le $\overline{\mathbb{F}}_{\ell}$-espace vectoriel sous-jacent à $r$, fournit un isomorphisme $G$-équivariant naturel $r\to\mathrm{i}(r)$ qui induit un isomorphisme $V\to\mathrm{i}(V)$. Il en résulte en particulier que la restriction à $V$ de la forme bilinéaire alternée évoquée ci-dessus est non dégénérée~; ce qui est absurde car la dimension de $V$ est impaire.

\medskip
Comme nous l'avons déjà dit, le théorème de Kronecker-Weber assure que les seuls homomorphismes $G\to\overline{\mathbb{F}}_\ell^\times$, supposés continus et non ramifiés hors de $\ell$, sont les puissances de $\chi$. Comme les entiers $\ell$ et $w$ sont impairs, on a $\mathrm{i}(\chi^{a})=\chi^{w-a}\neq \chi^a$ pour tout $a$ dans $\mathbb{Z}$, de sorte qu'aucune représentation de dimension $1$ de $G$ n'est ``fixée'' par $\mathrm{i}$. En particulier, $r$ est somme de $4$ caractères si, et seulement si, nous sommes dans le cas (i) de l'énoncé.

\medskip
Soit $V$ une $\overline{\mathbb{F}}_\ell$-représentation irréductible de
dimension $2$ de $G$ supposée continue et non ramifiée hors de $\ell$. 
Rappelons que l'on dit que $V$ est {\it impaire} si la classe de conjugaison
de $G$ constituée des conjugaisons complexes admet les valeurs propres $1$
et $-1$ dans la représentation $V$.  Il revient au m\^eme de dire que
$\mathop{\mathrm{d\acute{e}t}}V=\chi^{s}$ avec $s \equiv 1 \bmod 2$. Si $V$ est
impaire, la conjecture de Serre ``en niveau 1'', démontrée par Khare
\cite{khare}, assure qu'il existe $a$ dans $\mathbb{Z}/(\ell-1)$, un entier $\kappa\leq
\ell+1$ et une représentation $\rho$ dans $\mathrm{R}_{\kappa;\ell}$, tels que l'on a $V\simeq\chi^{a}\otimes\rho$. L'égalité $\mathop{\mathrm{d\acute{e}t}}\rho=\chi^{\kappa-1}$ montre de plus que l'on a $V \simeq \mathrm{i}(V)$ si et seulement si l'on a $2a+\kappa-1\equiv w\bmod{\mathrm{m}(\rho)}$.

\medskip
Supposons qu'il existe $V$ dans $\mathrm{J}(r)$ de dimension $2$ avec $\mathrm{i}(V)\not\simeq V$. Dans ce cas, on a $r \simeq V\oplus\mathrm{i}(V)$ et $\mathrm{i}(V) \simeq V\otimes\chi^{w}\hspace{2pt}(\mathop{\mathrm{d\acute{e}t}}V)^{-1}$. On est donc dans
le cas (iii)$_{1}$ si $V$ est impaire, et dans le cas (iii)$_{2}$ sinon. On observera que $V$ est impaire si et seulement si $\mathrm{i}(V)$ l'est, de sorte que les cas (iii)$_{1}$ et (iii)$_{2}$ s'excluent mutuellement.

\medskip
On peut donc supposer que toute représentation $V$ de $\mathrm{J}(r)$ est soit de dimension $1$, soit de dimension $2$ avec $\mathrm{i}(V)\simeq V$. Notons que dans ce dernier cas, $V$ est automatiquement impaire (car $w\equiv 1\bmod{2}$). D'après le théorème de Khare, $V$ est donc de la forme $\chi^{a}\otimes\rho$ avec $\rho\in 
\mathrm{R}_{\kappa;\ell}$ et $\kappa\leq\ell+1$, et de plus $2a+\kappa-1 \equiv w\bmod{\mathrm{m}(\rho)}$. On est donc dans les cas (ii) et (iv), selon respectivement
que $\mathrm{J}(r)$ contient une représentation de dimension $1$ ou deux
représentations de dimension $2$. Dans le cas restant, $r$ est irréductible. 
\hfill$\square$

\bigskip
Ce résultat et la table~\ref{taujk} permettent en retour de démontrer l'inexistence de certaines congruences. Donnons quelques exemples pour terminer. 

\begin{propv}\label{irredrhobjk} La représentation $\overline{\mathrm{r}}_{j,k;\ell}$ est irréductible (sur $\overline{\mathbb{F}}_{\ell}$) dans chacun des cas suivants~:

\smallskip
--\hspace{8pt}$(j,k)=(6,8)$ et $\ell=7,13,19$~;

\smallskip
--\hspace{8pt}$(j,k)=(8,8)$ et $\ell=7,11,19$~;

\smallskip
--\hspace{8pt}$(j,k)=(12,6)$ et $\ell=11,17$~;

\smallskip
--\hspace{8pt}$(j,k)=(4,10)$ et $\ell=7,13,17$.

\medskip
De plus, on a dans chacun de ces cas $\mathrm{m}(\overline{\mathrm{r}}_{j,k;\ell})=\ell-1$,
sauf pour $(j,k)=(6,8)$ et $\ell=13$, auquel cas on a simplement $\mathrm{m}(\overline{\mathrm{r}}_{6,8;13})\equiv 0\bmod{6}$.
\end{propv}

\medskip
\textit{Démonstration.} La représentation $\overline{\mathrm{r}}_{j,k;\ell}$, vue sur
$\overline{\F}_\ell$, satisfait les hypothèses de la proposition \ref{conjserrered} (avec $w=j+2k-3$), d'après la remarque \ref{pairingsp4}.  Appliquant cette proposition, il s'agit donc d'exclure pour chacun des triplets $(j,k,\ell)$ de l'énoncé ci-dessus, la possibilité d'une décomposition de la forme (i) à (iv).  La table \ref{taujk} est suffisamment fournie pour qu'il y ait plusieurs manières de procéder. Donnons quelques recettes simples~:
\vfill\eject

\bigskip
{\em \textsc{Critère 1.}  Si  $\overline{\rm r}_{j,k;\ell}$ est dans l'un des cas (i), (ii), (iii)$_1$ et (iv), et si $p$ est un nombre premier avec $p\equiv 1\bmod{\ell}$, on a respectivement~: 

\smallskip
\begin{itemize}\item[--] $\tau_{j,k}(p) \equiv 4\bmod{\ell}$ dans le cas (i)~;

\smallskip
\item[--] $\tau_{j,k}(p) \equiv \tau_{\kappa}(p)+2\bmod{\ell}$ dans le cas (ii)~;

\smallskip
\item[--] $\tau_{j,k}(p)\equiv 2\hspace{1pt}\tau_{\kappa}(p)\bmod{\ell}$ dans le cas
(iii)$_1$~;

\smallskip
\item[--] et $\tau_{j,k}(p)\equiv\tau_{\kappa_{1}}(p)+\tau_{\kappa_{2}}(p)\bmod{\ell}$
dans le cas (iv).
\end{itemize}

\bigskip
\textsc{Critère 2.} Si $\overline{\rm r}_{j,k;\ell}$ est dans le cas (iii)$_2$, alors
$\tau_{j,k}(p) \equiv 0 \bmod \ell$ pour tout nombre premier $p$ avec $p\equiv -1 \bmod
\ell$.}

\bigskip
-- Supposons d'abord $\ell=7$. On a $\mathrm{S}_\kappa(\mathrm{SL}_2(\mathbb{Z}))=0$ pour
tout $\kappa \leq \ell+1$, puis $\mathrm{R}_{\kappa,7}=\emptyset$ pour $\kappa\leq 8$.  Il suffit donc d'éliminer les cas (i) et (iii)$_2$.  On a $29\equiv 1\bmod{7}$ et l'on tire de la table \ref{taujk} les congruences
$$
\hspace{24pt}
\tau_{6,8}(29)\equiv\tau_{4,10}(29) \equiv \tau_{8,8}(29)\equiv
0\bmod{7}
\hspace{24pt},
$$
ce qui élimine le cas (i) par le critère 1.  On élimine de même le cas (iii)$_2$ à l'aide du critère 2~: on a $13\equiv -1\bmod{7}$, $\tau_{6,8}(13) \equiv 6\bmod{7}$ et $\tau_{4,10}(13)\equiv\tau_{8,8}(13)\equiv 1 \bmod{7}$.

\bigskip
-- Supposons $\ell=11$. On élimine les cas (i) et (iii)$_{2}$ comme précédemment, en observant que l'on a $23 \equiv 1\bmod{11}$ et
$\tau_{8,8}(23)\equiv\tau_{12,6}(23) \equiv 0\bmod{11}$, puis $43\equiv
-1 \bmod 11$ et $\tau_{8,8}(43) \equiv \tau_{12,6}(43) \equiv 6 \bmod{11}$. L'unique entier $\kappa\leq \ell+1$ tel que $\mathrm{S}_\kappa(\mathrm{SL}_2(\mathbb{Z}))
\neq 0$ est $\kappa=12$, et l'on a $\tau_{12}(23) \equiv -1\bmod{11}$.  On conclut à l'irréductibilité de $\overline{\mathrm{r}}_{j,k;11}$ pour $(j,k)=(8,8)$ et $(12,6)$ en observant $\tau_{j,k}(23) \not\equiv 1,-2 \bmod{11}$ (critère 1).  

\bigskip
-- Le cas $\ell=13$ est similaire car $\mathrm{S}_{14}(\mathrm{SL}_2(\mathbb{Z}))=0$. Le
critère 1 s'applique car $53 \equiv 1\bmod{13}$, $\tau_{12}(53) \equiv -3\bmod{13}$ et $\tau_{6,8}(53) \equiv \tau_{4,10}(53) \equiv 3 \bmod{13}$. Le critère 2 s'applique aussi car $103 \equiv -1 \bmod{13}$, $\tau_{6,8}(103)\equiv 11 \bmod{13}$ et $\tau_{4,10}(103) \equiv 5\bmod{13}$.

\bigskip
-- Dans le cas $\ell=17$, on conclut encore par les critères $1$ et $2$ grâce aux congruences suivantes : $103 \equiv 1\bmod{17}$, $\tau_{12}(103) \equiv 2 \bmod 17$, $\tau_{16}(103) \equiv 6 \bmod{17}$, $\tau_{18}(103)\equiv 8\bmod{17}$, $\tau_{4,10}(103)\equiv\tau_{12,6}(103)\equiv 1\bmod{17}$, et $67\equiv -1\bmod{17}$, $\tau_{4,10}(67)\equiv 8\bmod{67}$ et $\tau_{12,6}(67) \equiv 12\bmod{67}$.

\bigskip
-- Dans le cas $\ell=19$, le plus petit nombre premier $\equiv 1\bmod 19$ est $191>113$, qui est en dehors de la table \ref{taujk}.  En revanche, le critère 2 élimine encore le cas (iii)$_{2}$ car on a $37\equiv -1\bmod{19}$, $\tau_{6,8}(37) \equiv 4\bmod{19}$ et $\tau_{8,8}(37)\equiv 8\bmod{19}$.  Dégageons un autre critère.  Celui-ci  est basé sur l'observation suivante~: pour $\ell \leq 19$, si $\overline{\mathrm{r}}_{j,k;\ell}$ n'est pas dans le cas (iii)$_{2}$, alors tous ses facteurs de Jordan-Hölder sont définis sur $\mathbb{F}_{\ell}$.  Cela résulte en effet de la proposition \ref{conjserrered} et du fait que $\dim\mathrm{S}_\kappa(\mathrm{SL}_2(\mathbb{Z})) \leq 1$ pour $\kappa\leq\ell+1\leq 20$.

\bigskip
{\em \textsc{Critère 3.} Supposons que l'on a $\ell\leq 19$ et que $\overline{\mathrm{r}}_{j,k;\ell}$ n'est pas dans le cas (iii)$_{2}$. Supposons de plus qu'il existe un nombre premier $p\neq\ell$ tel que le polynôme $\mathrm{P}_{p}(t):=\mathop{\mathrm{d\acute{e}t}}(t-\mathrm{r}_{j,k;\ell}(\mathrm{ Frob}_{p}))$ de $\mathbb{Z}[t]$ soit irréductible modulo $\ell$. Alors la représentation $\overline{\mathrm{ r}}_{j,k;\ell}$ est irréductible.}

\bigskip
Comme nous avons déterminé les $\tau_{j,k}(p)$ pour $p\leq 113$ (table \ref{taujk}) et les $\tau_{j,k}(p^{2})$ pour $p\leq 29$ (table \ref{taujkp2}), la formule \eqref{polcarsp4} montre que nous disposons des polynômes $\mathrm{P}_{p}(t)$ pour $p\leq 29$~; dans le cas $\ell=19$, le critère ci-dessus est vérifié, pour $(j,k)=(6,8)$ et $p=3$, pour $(j,k)=(8,8)$ et $p=13$.

\bigskip
Il ne reste plus qu'à justifier la dernière assertion de la proposition \ref{irredrhobjk} portant sur les $\mathrm{m}(\overline{\mathrm{r}}_{j,k;\ell})$. On utilise l'observation suivante : soit $p\not=\ell$ un nombre premier avec $\tau_{j,k}(p)\not\equiv 0 \bmod{\ell}$, alors $\mathrm{m}(\overline{\mathrm{r}}_{j,k;\ell})$ est divisible par l'ordre de $p$ dans $(\mathbb{Z}/\ell)^{\times}$.

\bigskip
Par exemple, le nombre premier $3$ engendre $(\mathbb{Z}/7)^{\times}$ et on vérifie les congru\-ences $3\hspace{1pt}\tau_{6,8}(3)\equiv\tau_{4,10}(3)\equiv \tau_{8,8}(3) \equiv 4\bmod{7}$. Ceci montre $\mathrm{m}(\overline{r}_{6,8;7})=\mathrm{m}(\overline{r}_{8,8;7})=\mathrm{m}(\overline{\mathrm{r}}_{4,10;7})=6$. Les autres cas sont similaires. Pour $\ell=13$, on utilise que le nombre premier $2$ engendre $(\mathbb{Z}/\ell)^{\times}$, et la congruence $\tau_{4,10}(2)\equiv 10\bmod{13}$, de sorte que l'on a $\mathrm{m}(\overline{r}_{4,10;13})=12$. En revanche, on a $\tau_{6,8}(2)\equiv 0\bmod{13}$. On montre $\mathrm{m}(\overline{\mathrm{r}}_{6,8;13})\equiv 0 \bmod{6}$ à l'aide de la congruence $\tau_{6,8}(17)\equiv 7\bmod{13}$.
\hfill$\square$

\bigskip
Il est évident que les méthodes \textit{ad hoc} employées ci-dessus sont assez grossières, et qu'il est possible d'étudier les décompositions éventuelles des représentations $\overline{\mathrm{r}}_{j,k;\ell}$ pour des caractéristiques $\ell>19$.  Nous reportons cette étude, ainsi que la question plus intéressante de déterminer les images des $\overline{\mathrm{r}}_{j,k;\ell}$, à un travail ultérieur.  Pour nous autocongratuler, observons que les triplets $(j,k;\ell)$ intervenant dans l'énoncé de la proposition \ref{irredrhobjk} sont exactement ceux n'intervenant pas dans l'énoncé du théorème \ref{recapcong}, lorsque l'on a $\ell\leq 19$.

\begin{remarque}\label{remm6813} {\rm Supposons que $\overline{\mathrm{r}}_{j,k;\ell}$ est irréductible sur $\overline{\mathbb{F}}_\ell$ (et donc, en particulier, que l'on a $\ell>5$). Les assertions suivantes sont équivalentes :

\smallskip
{\em
(i) pour tout premier $p$ non carré modulo $\ell$, on a
$\tau_{j,k}(p) \equiv 0\bmod\ell$~;

\smallskip
(ii) $\mathrm{m}(\overline{\mathrm{r}}_{j,k;\ell})$ divise $\frac{\ell-1}{2}$~;

\smallskip
(iii) $\overline{r}_{j,k;\ell}$ est induite d'une représentation irréductible de dimension $2$, à coefficients dans $\F_{\ell^2}$, du groupe de Galois absolu de
$\mathbb{Q}(\sqrt{\ell^{*}})$ avec $\ell^{*}=(-1)^{\frac{\ell-1}{2}}\ell$.}

\bigskip
Dans le cas $(j,k)=(6,8)$ et $\ell=13$, une inspection de la table \ref{taujk} montre que la congruence de l'assertion (i) vaut pour tout premier $p \leq 113$. Il est tentant de conjecturer qu'elle est toujours satisfaite, autrement dit que l'on a $\mathrm{m}(\overline{r}_{6,8;13}) = 6$.}
\end{remarque}

\appendix

\titleformat{\chapter}
{\bf \large}
{\Alph{chapter}.}
{1 em}
{}

\chapter{Le réseau de Barnes-Wall et les séries thêta de Siegel des réseaux unimodulaires pairs de dimension $16$}
\label{appendicekneser}

\begin{center}
(d'après Martin Kneser \cite{kneser})
\end{center}

\parindent=0cm

\bigskip
Le réseau en question est un remarquable réseau de dimension $16$ découvert par Barnes et Wall en 1959 \cite{BW}. Il apparaît à de nombreuses reprises dans \cite{conwaysloane} où il est noté $\Lambda_{16}$ ou $\mathrm{BW}_{16}$ (Conway et Sloane mentionnent malicieusement qu'il a été redécouvert par beaucoup d'auteurs). Barnes et Wall définissent en fait dans \cite{BW} une suite de réseaux $(\Lambda_{2^{n}})_{n\in\mathbb{N}-\{0\}}$, $\Lambda_{2^{n}}$ étant de dimension $2^{n}$~;  les trois premiers réseaux de cette suite sont respectivement isomorphes à $\mathrm{I}_{2}$, $\mathrm{D}_{4}$ et $\mathrm{E}_{8}$. Pour une présentation simple et élégante des réseaux de Barnes-Wall nous conseillons \cite{NRS1} et \cite{NRS2}, on pourra également consulter \cite{BE}~; le réseau $\mathrm{U}$ qui apparaît dans la référence \cite{kneser} est un avatar de $\Lambda_{16}$. La définition de $\Lambda_{16}$ que nous donnons ci-après suffira à notre bonheur.

\medskip
Soit $I$ un lagrangien du $\mathrm{q}$-espace vectoriel $\mathbb{F}_{2}\otimes_{\mathbb{Z}}\mathrm{E}_{8}$ (on a donc $\mathrm{q}(I)=0$ et $\dim_{\mathbb{F}_{2}}I=4$) . On note $\Lambda_{16}$ le sous-module de $\mathrm{E}_{8}\oplus\mathrm{E}_{8}$ constitué des couples $(x_{1},x_{2})$ d'éléments de $\mathrm{E}_{8}$ dont les réductions modulo $2$, disons $\bar{x}_{1}$ et $\bar{x}_{2}$, vérifient $\bar{x}_{1}+\bar{x}_{2}\in I$~;  $\Lambda_{16}$ peut donc être vu comme un réseau entier (au sens quadratique) dans $\mathbb{Q}\otimes_{\mathbb{Z}}(\mathrm{E}_{8}\oplus\mathrm{E}_{8})$.

\bigskip
\textbf{Proposition A.1.} {\em Le réseau $\Lambda_{16}$ vérifie les propriétés suivantes~:
\begin{itemize}
\item[{\em (a)}] On a $\mathrm{q}(x)\geq 2$ (ou encore $x.x\geq 4$) pour tout $x$ dans $\Lambda_{16}-\{0\}$.
\item[{\em (b)}] On a $\xi.\xi\in\mathbb{Z}$ pour tout $\xi$ dans le réseau dual $\Lambda_{16}^{\sharp}$ et $\Lambda_{16}^{\sharp}$ muni de la forme quadratique $\xi\mapsto\xi.\xi$ est isomorphe (comme $\widetilde{\mathrm{q}}$-module) à $\Lambda_{16}$.
\item[{\em (c)}] Le module sous-jacent au $\mathrm{qe}$-module $\mathop{\mathrm{r\acute{e}s}}\Lambda_{16}$ est annulé par $2$ et le\linebreak  $\mathrm{qe}$-module $\mathop{\mathrm{r\acute{e}s}}\Lambda_{16}$ est isomorphe, via le plongement canonique de $\mathbb{F}_{2}$ dans $\mathbb{Q}/\mathbb{Z}$, au $\mathbb{F}_{2}$-$\mathrm{q}$-espace vectoriel hyperbolique $\mathrm{H}(I)$.
\end{itemize}}

\bigskip
\textit{Démonstration de (a).} Soit $x=(x_{1},x_{2})$ un élément de $\Lambda_{16}-\{0\}$. Si $x_{1}$ et $x_{2}$ sont non nuls alors on a $\mathrm{q}(x)=\mathrm{q}(x_{1})+\mathrm{q}(x_{2})\geq 1+1$. Si $x_{i}$ est nul alors $x_{3-i}$ appartient à $I$ (et est non nul) et $\mathrm{q}(x_{3-i})$ est pair.

\smallskip
\textit{Démonstration de (b).} On constate que le réseau $\Lambda_{16}^{\sharp}$ est le sous-module de $\mathbb{Q}\otimes_{\mathbb{Z}}(\mathrm{E}_{8}\oplus\mathrm{E}_{8})$ constitué des éléments $\xi=(\xi_{1},\xi_{2})$ tels que $(\xi_{1}+\xi_{2},\xi_{1}-\xi_{2})$ appartient à $\Lambda_{16}$. Or on a l'identité $\mathrm{q}(\xi_{1}+\xi_{2})+\mathrm{q}(\xi_{1}-\xi_{2})=2(\mathrm{q}(\xi_{1})+\mathrm{q}(\xi_{2}))$.

\smallskip
\textit{Démonstration de (c).} Soit $\Delta I$ le sous-espace vectoriel de $\mathbb{F}_{2}\otimes_{\mathbb{Z}}(\mathrm{E}_{8}\oplus\mathrm{E}_{8})=(\mathbb{F}_{2}\otimes_{\mathbb{Z}}\mathrm{E}_{8})\oplus(\mathbb{F}_{2}\otimes_{\mathbb{Z}}\mathrm{E}_{8})$ image diagonale de $I$~; on observe que $\Lambda_{16}$ est le sous-module de $\mathrm{E}_{8}\oplus\mathrm{E}_{8}$ image réciproque de $(\Delta I)^{\perp}$ par l'homomorphisme  $\mathrm{E}_{8}\oplus\mathrm{E}_{8}\to\mathbb{F}_{2}\otimes_{\mathbb{Z}}(\mathrm{E}_{8}\oplus\mathrm{E}_{8})$. Cette observation faite, la propriété (c) est une manifestation du phénomène général décrit ci-dessous.

\smallskip
Soit $L$ un $\mathrm{q}$-module sur $\mathbb{Z}$. Soient $p$ un nombre premier et $J$ un sous-espace vectoriel de $\mathbb{F}_{p}\otimes_{\mathbb{Z}}L$ avec $\mathrm{q}(J)=0$~; soit $M$ le sous-module de $L$ constitué des éléments dont la réduction modulo $p$ est orthogonale à $J$. Alors le groupe abélien $\mathop{\mathrm{r\acute{e}s}}M$ est annulé par $p$ et le $\mathrm{qe}$-module $\mathop{\mathrm{r\acute{e}s}}M$ est isomorphe, \textit{via} le plongement canonique de $\mathbb{F}_{p}$ dans $\mathbb{Q}/\mathbb{Z}$, au $\mathbb{F}_{p}$-$\mathrm{q}$-espace vectoriel hyperbolique $\mathrm{H}(J)$ ($\cong\mathrm{H}(J^{*})\cong\mathrm{H}(L/M)$).
\hfill$\square$

\bigskip
\textbf{Corollaire A.2.} {\em On a $\xi.\xi\geq 2$ pour tout $\xi$ dans $\Lambda_{16}^{\sharp}-\{0\}$.}

\bigskip
On explique maintenant comment déduire des propriétés du réseau $\Lambda_{16}$, en suivant la stratégie de Kneser, l'égalité des séries thêta
$$
\vartheta^{(g)}_{\mathrm{E}_{8}\oplus\mathrm{E}_{8}}
=
\vartheta^{(g)}_{\mathrm{E}_{16}}
$$
pour $g\leq 3$ (résultat dû à Witt pour $g\leq 2$ \cite{witt}.

\medskip
Cette égalité peut être reformulée en termes de représentations de formes quadratiques entières par $\mathrm{E}_{8}\oplus\mathrm{E}_{8}$ et $\mathrm{E}_{16}$. Précisons la terminologie. Soient $L$ un réseau unimodulaire pair et $G$ un $\mathbb{Z}$-module libre de dimension finie muni d'une forme quadratique à valeurs entières (compte tenu de ce qui va suivre, on peut supposer que ces valeurs sont positives ou nulles)~; une {\em représentation} de $G$ par $L$ est un homomorphisme $f:G\to L$ avec $\mathrm{q}(f(x))=\mathrm{q}(x)$ pour tout $x$ dans $G$. On note $\mathrm{Rep}(G,L)$ l'ensemble des représentations de $G$ par $L$~; il est clair que cet ensemble est fini. Voici la reformulation promise~:

\bigskip
\textbf{Théorème A.3.} {\em Soit $G$ un $\mathbb{Z}$-module libre de dimension finie muni d'une forme quadratique à valeurs entières. S l'on a $\dim G\leq 3$ alors les ensembles $\mathrm{Rep}(G,\mathrm{E}_{8}\oplus\mathrm{E}_{8})$ et $\mathrm{Rep}(G,\mathrm{E}_{16})$ ont même cardinal.}

\bigskip
\textit{Démonstration.} Le point-clé est l'observation ci-après.

\medskip
Soit $\gamma:\Lambda_{16}^{\sharp}\to\mathop{\mathrm{r\acute{e}s}}\Lambda_{16}$ l'homomorphisme de passage au quotient. On note $\mathcal{I}$ l'ensemble (fini) des lagrangiens de $\mathop{\mathrm{r\acute{e}s}}\Lambda_{16}$ et $\mathcal{L}$ l'ensemble des réseaux unimodulaires pairs contenant $\Lambda_{16}$ (et contenus dans $\Lambda_{16}^{\sharp}$)~; on rappelle que l'application $\mathcal{I}\to\mathcal{L}\hspace{2pt},\hspace{2pt}I\mapsto\gamma^{-1}(I)$ est une bijection (compatible avec les relations d'inclusion).
\vfill\eject

\medskip
Soit $J$ un sous-module de $\mathop{\mathrm{r\acute{e}s}}\Lambda_{16}$ avec $\mathrm{q}(J)=0$ et $\dim_{\mathbb{F}_{2}}J=3$. Soit $M$ le réseau $\gamma^{-1}(J)$~; on a $\mathop{\mathrm{r\acute{e}s}}M\cong J^{\perp}/J\simeq\mathrm{H}(\mathbb{Z}/2)$ (voir II.1.1). Soit $B$ le réseau image réciproque, par l'homomorphisme $M^{\sharp}\to\mathop{\mathrm{r\acute{e}s}}M$, de la ``droite" non isotrope (non isotrope au sens quadratique mais isotrope au sens bilinéaire)~;  $B$ est un réseau unimodulaire impair (voir  la discussion ``$2$-voisins, le point de vue de Borcherds" qui suit III.1.8). Comme $B$ est contenu dans $\Lambda_{16}^{\sharp}$, le corollaire A.2 implique que l'on a $x.x\geq 2$ pour tout $x$ dans $B-\{0\}$. Le scholie III.3.3.2 montre enfin que $B$ est isomorphe à $\mathrm{Bor}_{16}$ (définition III.3.3.2) et que les deux réseaux unimodulaires pairs, images réciproques des deux ``droites" isotropes de $\mathop{\mathrm{r\acute{e}s}}M$, sont non isomorphes (on observera incidemment que ce scholie montre aussi que la classe d'isomorphisme de $M$ est indépendante du choix de $J$).

\medskip
L'observation ci-dessus conduit au lemme ci-dessous. On appelle encore repré\-sentation de $G$ par $\Lambda_{16}^{\sharp}$ un homomorphisme $f:G\to\Lambda_{16}^{\sharp}$ avec $\mathrm{q}(f(x))=\mathrm{q}(x)$ pour tout $x$ dans $G$.

\bigskip
\textbf{Lemme A.4.} {\em Soit $f$ une représentation de $G$ par $\Lambda_{16}^{\sharp}$ avec $\dim G\leq 3$. Soit $\mathcal{L}(f)$ le sous-ensemble de $\mathcal{L}$ constitué des réseaux $L$ contenant $f(G)$~; soient $\mathcal{L}_{1}(f)$ et  $\mathcal{L}_{2}(f)$ les sous-ensembles de $\mathcal{L}(f)$ constitués des réseaux $L$ respectivement isomorphes à $\mathrm{E}_{8}\oplus\mathrm{E}_{8}$ et $\mathrm{E}_{16}$. Alors $\mathcal{L}_{1}(f)$ et  $\mathcal{L}_{2}(f)$ ont même cardinal.}

\bigskip
\textit{Démonstration.} On note $\mathcal{I}(f)$, $\mathcal{I}_{1}(f)$ et $\mathcal{I}_{2}(f)$ les sous-ensembles de $\mathcal{I}$ images des sous-ensembles $\mathcal{L}(f)$, $\mathcal{L}_{1}(f)$ et $\mathcal{L}_{2}(f)$ de $\mathcal{L}$ par la bijection inverse de la bijection $I\mapsto\gamma^{-1}(I)$. On note $\mathcal{J}$ l'ensemble des $J$ considérés plus haut et $\mathcal{J}(f)$ le sous-ensemble de $\mathcal{J}$ constitué des $J$ avec $J\supset (\gamma\circ f)(G)$. On note $\mathcal{K}(f)$ (resp. $\mathcal{K}_{1}(f)$, resp. $\mathcal{K}_{2}(f)$) le sous-ensemble de $\mathcal{I}(f)\times\mathcal{J}(f)$ (resp. $\mathcal{I}_{1}(f)\times\mathcal{J}(f)$, resp. $\mathcal{I}_{2}(f)\times\mathcal{J}(f)$) constitué des couples $(I,J)$ avec $I\supset J$. On note enfin $\pi_{\mathcal{J}}:\mathcal{K}(f)\to\mathcal{J}(f)$ et $\pi_{\mathcal{I}}:\mathcal{K}(f)\to\mathcal{I}(f)$ les applications $(I,J)\mapsto J$ et $ (I,J)\mapsto I$. Il est clair que $\pi_{\mathcal{J}}$ est surjective et que ses fibres ont toutes $2$ éléments. Pareillement $\pi_{\mathcal{I}}$ est surjective et ses fibres ont toutes $2^{(4-\delta(f))}-1$ éléments, $\delta(f)\leq 3$ désignant la dimension du $\mathbb{F}_{2}$-espace vectoriel $(\gamma\circ f)(G)$.  Compte tenu de ce qui précède, les applications $\mathcal{K}_{1}(f)\to\mathcal{J}(f)$ et $\mathcal{K}_{2}(f)\to\mathcal{J}(f)$, induites par $\pi_{\mathcal{J}}$, sont toujours surjectives. Comme l'on a  $\vert\mathcal{K}_{1}(f)\vert+\vert\mathcal{K}_{2}(f)\vert=\vert\mathcal{K}(f)\vert=2\vert\mathcal{J}(f)\vert$ (la notation $\vert-\vert$ désigne ici le cardinal d'un ensemble fini), on voit que l'on a $\vert\mathcal{K}_{1}(f)\vert=\vert\mathcal{J}(f)\vert$ et $\vert\mathcal{K}_{2}(f)\vert=\vert\mathcal{J}(f)\vert$. Par définition on a $\mathcal{K}_{1}(f)=\pi_{I}^{-1}(\mathcal{I}_{1}(f))$ et $\mathcal{K}_{2}(f)=\pi_{I}^{-1}(\mathcal{I}_{2}(f))$, d'où l'égalité $\vert\mathcal{I}_{1}(f)\vert=\vert\mathcal{I}_{2}(f)\vert$.
\hfill$\square$

\bigskip
\textit{Fin de la démonstration du théorème.} On note $\mathcal{L}_{1}$ et $\mathcal{L}_{2}$ les sous-ensembles de $\mathcal{L}$ constitués des réseaux respectivement isomorphes à $\mathrm{E}_{8}\oplus\mathrm{E}_{8}$ et $\mathrm{E}_{16}$~; on observe que l'on a $\mathcal{L}_{1}=\mathcal{L}_{1}(0)$, $\mathcal{L}_{2}=\mathcal{L}_{2}(0)$ ($0$ désignant l'unique représentation de $0$ par $\Lambda_{16}^{\sharp}$) et donc $\vert\mathcal{L}_{1}\vert=\vert\mathcal{L}_{2}\vert$. On note encore $\mathrm{Rep}(G,\Lambda_{16}^{\sharp})$ l'ensemble (fini) des représentations de $G$ par $\Lambda_{16}^{\sharp}$, on pose $\mathrm{r}_{1}(G)=\vert\mathrm{Rep}(G,\mathrm{E}_{8}\oplus\mathrm{E}_{8})\vert$,\linebreak $\mathrm{r}_{2}(G)=\vert\mathrm{Rep}(G,\mathrm{E}_{16})\vert$ et l'on considère les sous-ensembles de $\mathcal{L}_{i}\times\mathrm{Rep}(G,\Lambda_{16}^{\sharp})$, disons $\mathcal{R}_{i}(G)$ ($i=1,2$), constitués des couples $(L,f)$ avec $L\supset f(G)$. En projetant sur chacun des deux facteurs du produit $\mathcal{L}_{i}\times\mathrm{Rep}(G,\Lambda_{16}^{\sharp})$, on constate que l'on a 
$$
\hspace{24pt}
\vert\mathcal{R}_{i}(G)\vert=\vert\mathcal{L}_{i}\vert\hspace{2pt}\mathrm{r}_{i}(G)
\hspace{24pt}\text{et}\hspace{24pt}
\vert\mathcal{R}_{i}(G)\vert=\sum_{f\in\mathrm{Rep}(G,\Lambda_{16}^{\sharp})} \vert\mathcal{L}_{i}(f)\vert
\hspace{24pt}.
$$
Le lemme A.4 implique bien l'égalité $\mathrm{r}_{1}(G)=\mathrm{r}_{2}(G)$.
\hfill$\square\square$

\parindent=0cm

\chapter{\textbf{Formes quadratiques et voisins en dimension impaire}}\label{appendixb}

\bigskip
On présente dans cet appendice le pendant , ``en dimension impaire'', d'une partie de la théorie que nous avons développée, ``en dimension paire'', dans les chapitres II et III.

\vspace{0,75cm}
\textbf{1.} \textbf{Généralités concernant les formes quadratiques sur un module projectif de rang constant impair}

\bigskip
Soit $A$ un anneau commutatif unitaire. Au chapitre II, nous avons défini un $\mathrm{q}$-module sur $A$ comme un $A$-module projectif de type fini $L$ muni d'une forme quadratique $\mathrm{q}:L\to A$ telle que la forme bilinéaire symétrique associée est non dégénéré. Si $2$ est non inversible dans $A$ et si $L$ est de rang constant, alors cette non dégénéréscence force ce rang à être pair (considérer un homomorphisme de $A$ dans un corps $k$ de caractéristique $2$ et observer que la forme bilinéaire symétrique associée à la forme quadratique dont $k\otimes_{A}L$ est muni, est alternée). Si $L$ est de rang constant impair il est classique de dire que $\mathrm{q}$ est non dégénérée si la dégénérescence de la forme bilinéaire symétrique associée est ``minimale''. Nous précisons cette définion ci-après. Notre présentation met en avant la notion de demi-déterminant (voir par exemple \cite{knus}~); pour une présentation plus sophistiquée voir \cite[Expos\'e XII]{sga7}.

\bigskip
Soit $k\geq 1$ un entier~; on note respectivement
$$
\psi_{L}^{k}:\Lambda^{k}L\to \Lambda^{k-1}L\otimes L
\hspace{24pt},\hspace{24pt}
\phi_{L}^{k}:\Lambda^{k-1}L\otimes L\to\Lambda^{k}L
$$
les homomorphismes ``coproduit'' et ``produit'' induits par la structure d'algè\-bre de Hopf de l'algèbre extérieure $\Lambda L$ (les produits tensoriels et l'algèbre extérieure sont sur $A$).

\bigskip
\textbf{Lemme 1.1.} {\em Soit $L$ un $A$-module muni d'une forme quadratique $q:L\to A$~; soit $b:L\times L\to A$ la forme bilinéaire symétrique associée. La forme bilinéaire symétrique associée à la forme quadratique $(\Lambda^{k-1}b\otimes q)\circ\psi_{L}$ est $k\hspace{1pt}\Lambda^{k}b$.}

\bigskip
\textit{Démonstration.} Elle résulte de ce que l'homomorphisme composé $\phi_{L}^{k}\circ\psi_{L}^{k}$ est $k$-fois l'identité de $\Lambda^{k}L$.
\hfill$\square$

\bigskip
Soit $\beta$ une forme bilinéaire symétrique, on note $\mathrm{qd}(\beta)$ la forme quadratique $x\mapsto\beta(x,x)$~; on observe que la forme bilinéaire qui lui est associée est $2\hspace{1pt}\beta$.

\bigskip
La démonstration des énoncés 1.2, 1.3 et 1.4 ci-après est immédiate.

\bigskip
\textbf{Proposition-Définition 1.2} (puissances extérieures impaires d'une forme quadratique)\textbf{.} {\em Soient $L$ un $A$-module muni d'une forme quadratique $q$ et $k$ un entier impair. On pose
$$
\hspace{24pt}
\Lambda^{k}q=(\Lambda^{k-1}b\otimes q)\circ\psi_{L}^{k}-\frac{k-1}{2}\hspace{3pt}\mathrm{qd}\hspace{1pt}(\Lambda^{k}b)
\hspace{24pt},
$$
$b$ désignant la forme bilinéaire symétrique associée à $q$ ($\Lambda^{k}q$ est donc une forme quadratique sur $\Lambda^{k}L$).

\medskip
La forme bilinéaire symétrique associée à $\Lambda^{k}q$ est $\Lambda^{k}b$. On appelle $\Lambda^{k}q$  la {\em $k$-ième puissance extérieure} de $q$.}

\bigskip
\textbf{Proposition 1.3.} {\em Soit $L$ un $A$-module projectif de type fini. Soit $\mathcal{B}(L)$\linebreak (resp. $\mathcal{Q}(L)$) le $A$-module constitué des formes bilinéaires symétriques\linebreak $L\times L\to A$ (resp. des formes quadratiques $L\to A$). Si $L$ est de rang~$1$ alors l'homomorphisme de $A$-modules $\mathrm{qd}:\mathcal{B}(L)\to\mathcal{Q}(L)$ est un isomorphisme.}

\bigskip
\textbf{Proposition-Définition 1.4.} {\em Soit $L$ un $A$-module projectif de type fini de rang constant impair, disons $n$, muni d'une forme quadratique $q$~; on note~$b$ la forme bilinéaire symétrique associée à $q$. Les deux formes bilinéaires symétriques $\Lambda^{n}b$ et $\mathrm{qd}^{-1}(\Lambda^{n}q)$, dont le $A$-module $\Lambda^{n}L$, projectif de rang $1$, est muni, sont reliées par l'égalité
$$
\hspace{24pt}
\Lambda^{n}b
\hspace{4pt}=\hspace{4pt}
2\hspace{2pt}\mathrm{qd}^{-1}(\Lambda^{n}q)
\hspace{24pt}.
$$
La forme bilinéaire symétrique $\mathrm{qd}^{-1}(\Lambda^{n}q)$ (ou le $A$-module projectif $\Lambda^{n}L$ muni de cette forme) s'appelle le {\em demi-déterminant} de $q$ (on rappelle que la forme bilinéaire symétrique $\Lambda^{n}b$ s'appelle le déterminant de $b$, voir II.1). Le demi-déterminant de $q$ est noté $\frac{1}{2}\text{-}\mathop{\mathrm{d\acute{e}t}}q$ ou $\frac{1}{2}\text{-}\mathop{\mathrm{d\acute{e}t}}L$.

\smallskip
(La terminologie et la notation sont évidemment justifiées par le fait que l'on a $\mathop{\mathrm{d\acute{e}t}}L=2\hspace{2pt}(\frac{1}{2}\text{-}\mathop{\mathrm{d\acute{e}t}}L)$~!).}

\bigskip
\textit{Exemple.} Soit $a$ un élément de $A$. On rappelle que l'on note $\langle a\rangle$ le $A$-module $A$ muni de la forme bilinéaire symétrique $(x,y)\mapsto a\hspace{1pt}xy$~; on note donc $\mathrm{qd}(\langle a\rangle)$ le $A$-module $A$ muni de la forme quadratique $x\mapsto a\hspace{1pt}x^{2}$. Nous avons tout fait pour que le demi-déterminant de $\mathrm{qd}(\langle a\rangle)$ soit $\langle a\rangle$.
\vfill\eject

\bigskip
\textbf{Proposition 1.5.} {\em Soient $P$ et $L$ deux $A$-modules projectifs de type fini, respectivement  de rang constant pair et impair,  munis de forme quadratique. On a un isomorphisme canonique de $A$-modules projectifs de rang $1$ munis de forme bilinéaire symétrique
$$
\textstyle
\hspace{24pt}
\frac{1}{2}\text{-}\mathop{\mathrm{d\acute{e}t}}(P\oplus L)
\hspace{4pt}\cong\hspace{4pt}
\mathop{\mathrm{d\acute{e}t}}P
\otimes
\frac{1}{2}\text{-}\mathop{\mathrm{d\acute{e}t}}L
\hspace{24pt}.
$$}

\textit{Démonstration.} Soient respectivement $q_{P}$, $b_{P}$, $q_{L}$, $b_{L}$ les formes quadratiques et bilinéaires dont $P$ et $L$ sont munis~; soit $m$ le rang de $P$ et $n$ celui de $L$. Il faut montrer que l'on a $\Lambda^{m+n}(q_{P}\oplus q_{L})=\Lambda^{m}b_{P}\otimes\Lambda^{n}q_{L}$. En utilisant la naturalité des puissances extérieures des formes bilinéaires (resp. puissances extérieures impaires des formes quadratiques) on se ramène au ``cas universel''. Dans ce cas l'anneau $A$ est un anneau de polynômes, à coefficients dans~$\mathbb{Z}$, en $\frac{m(m+1)}{2}+\frac{n(n+1)}{2}$ indéterminées, et $2$ n'est pas diviseur de $0$ ; or on sait que l'on a $\mathop{\mathrm{d\acute{e}t}}(P\oplus L)\cong\mathop{\mathrm{d\acute{e}t}}P
\otimes\mathop{\mathrm{d\acute{e}t}}L$.
\hfill$\square$

\bigskip
\textbf{Définition 1.6.} {\em Soit $L$ un $A$-module projectif de type fini de rang constant impair muni d'une forme quadratique $q$. Nous dirons que $q$ est {\em minimalement dégénérée} si la forme bilinéaire symétrique $\frac{1}{2}\text{-}\mathop{\mathrm{d\acute{e}t}}q$ est non dégénérée. Pour faire court un $A$-module projectif de type fini de rang constant impair muni d'une forme quadratique minimalement dégénérée, sera appelé un {\em $\mathrm{q}\text{-}\mathrm{i}$-module} sur $A$.}

\bigskip
\textit{Exemple-Remarque.} Soient $A$ un anneau, $P$ un $\mathrm{q}$-module sur $A$ de rang constant pair et $u$ un élément de $A^{\times}$. La proposition 1.5 montre que la somme orthogonale $P\oplus\mathrm{qd}(\langle u\rangle)$ est un $\mathrm{q}\text{-}\mathrm{i}$-module sur $A$. La proposition 1.2 de \cite[Expos\'e XII]{sga7} dit que, localement pour la topologie étale, tout $\mathrm{q}\text{-}\mathrm{i}$-module est de ce type avec en outre $P$ hyperbolique.

\bigskip
\textsc{Groupes classiques (suite)}

\bigskip
Soit $L$ un $\mathrm{q}\text{-}\mathrm{i}$-module sur $A$, disons de rang $n$~; comme au chapitre~II, un endomorphisme $\alpha$ du $A$-module sous-jacent à $L$ est dit {\em orthogonal} s'il préserve la forme quadratique. La naturalité des puissances extérieures des formes quadratiques montre que $\Lambda^{n}\alpha$ est un endomorphisme orthogonal du $\mathrm{b}$-module $\frac{1}{2}\text{-}\mathop{\mathrm{d\acute{e}t}}L$ et donc que les endomorphismes $\Lambda^{n}\alpha$ et $\alpha$ sont des automorphismes. Les endomorphismes orthogonaux forment un groupe pour la composition que l'on appelle le {\em groupe orthogonal} de $L$ et que l'on note~$\mathrm{O}(L)$. Le foncteur $R\mapsto\mathrm{O}(R\otimes_{A}L)$, défini sur la catégorie des $A$-algèbres commutatives et à valeurs dans la catégorie des groupes, est un $A$-schéma en groupes que l'on note~$\mathrm{O}_{L}$. Compte tenu de ce l'on a vu plus haut, l'homomorphisme composé $\mathrm{O}_{L}\to\mathrm{GL}_{L}\overset{\mathrm{d\acute{e}t}}{\to}\mathbb{G}_{\mathrm{m}}$, induit un homomorphisme $\mathrm{d\acute{e}t}:\mathrm{O}_{L}\to\mu_{2}$ (observer que $\mathrm{O}_{\frac{1}{2}\text{-}\mathop{\mathrm{d\acute{e}t}}L}$ s'identifie à $\mu_{2}$). On note $\mathrm{SO}_{L}$ le noyau de ce dernier homomorphisme.

\medskip
Si $L$ est de rang $1$, alors le groupe $\mathrm{O}_{L}$ s'identifie encore $\mu_{2}$, si bien que l'on ne peut espérer que $\mathrm{O}_{L}$ soit lisse sur $A$, en toute généralité. C'est cependant le cas pour $\mathrm{SO}_{L}$~:

\bigskip
\textbf{Proposition 1.7.} {\em Pour tout $\mathrm{q\text{-}i}$-module $L$ sur un anneau $A$, le $A$-schéma en groupes $\mathrm{SO}_{L}$ est lisse sur $A$.}

\medskip
\footnotesize
\textit{Démonstration.} Puisque la propriété que l'on veut vérifier est locale pour la topologie étale, on peut supposer, d'après la proposition 1.2 de \cite[Expos\'e XII]{sga7} déjà citée plus haut, que l'on a $L=\mathrm{H}(A^{n})\oplus\mathrm{qd}(\langle u\rangle)$ avec $u$ dans $A^{\times}$. Comme les deux $\mathrm{q}$-modules $\mathrm{H}(A^{n})$ et $\langle u\rangle\otimes\mathrm{H}(A^{n})$ sont isomorphes on peut supposer en outre $u=1$ et donc au bout du compte $A=\mathbb{Z}$ et $L=\mathrm{H}(\mathbb{Z}^{n})\oplus\mathrm{A}_{1}$~; le groupe $\mathrm{SO}_{L}$ est alors le groupe que l'on a noté $\mathrm{SO}_{n+1,n}$ en \ref{enonceparamst}. Le fait que $\mathrm{SO}_{n+1,n}$ est lisse sur $\mathbb{Z}$ est bien connu. Nous montrons ci-après que cette propriété peut être vue comme une conséquence de II.1.5~; cette démonstration (très indirecte~!) est dans l'esprit du paragraphe 2.

\smallskip
On pose $P=\mathrm{H}(\mathbb{Z}^{n})\oplus\mathrm{H}(\mathbb{Z})$. On note $(e_{1},e_{2})$ la base canonique facteur $\mathrm{H}(\mathbb{Z})$ et on pose $e=e_{1}+e_{2}$ et $f=e_{1}-e_{2}$~; on a donc $\mathrm{q}(e)=1$, $\mathrm{q}(f)=-1$ et $e.f=0$. On constate que $L$ s'identifie (avec sa forme quadratique) à l'orthogonal de $f$.

\smallskip
Soit $\mathcal{C}$ la quadrique affine d'équation $\mathrm{q}=-1$~; on observe que ce $\mathbb{Z}$-schéma est lisse sur $\mathbb{Z}$. On note $\mathrm{O}_{P,f}$ le sous-groupe du groupe $\mathrm{O}_{P}$ (qui est lisse sur $\mathbb{Z}$ d'après II.1.5), stabilisateur de $f$, pour l'action évidente de $\mathrm{O}_{P}$ sur $\mathcal{C}$~;  le ``calcul différentiel'' montre que le groupe $\mathrm{O}_{P,f}$ est lisse sur $\mathbb{Z}$. L'égalité $L=f^{\perp}$ fournit un homomorphisme de schémas en groupes, disons $\omega:\mathrm{O}_{P,f}\to\mathrm{O}_{L}$. Les énoncés 1.8 et 1.9 ci-après concernent cet homomorphisme~; le second implique 1.7. La démonstration du premier est laissée au lecteur.

\medskip
\textbf{Proposition 1.8} {\em Le diagramme
$$
\begin{CD}
\mathrm{O}_{P,f}@>\omega>>\mathrm{O}_{L} \\
@V\widetilde{\mathop{\mathrm{d\acute{e}t}}}VV
@V\mathop{\mathrm{d\acute{e}t}}VV \\
\mathbb{Z}/2@>>>\mu_{2}
\end{CD}
$$
dans lequel la flèche notée $\widetilde{\mathop{\mathrm{d\acute{e}t}}}$ est la restriction de l'homomorphisme $\widetilde{\mathop{\mathrm{d\acute{e}t}}}:\mathrm{O}_{P}\to\mathbb{Z}/2$, est commutatif.}

\medskip
On note $\mathrm{SO}_{P,f}$ le noyau de l'homomorphisme $\widetilde{\mathop{\mathrm{d\acute{e}t}}}:\mathrm{O}_{P,f}\to\mathbb{Z}/2$. La proposition ci-dessus montre que l'homomorphisme $\omega$ induit un homomorphisme, disons $\omega_{\mathrm{S}}$, de $\mathrm{SO}_{P,f}$ dans $\mathrm{SO}_{L}$.

\medskip
\textbf{Proposition 1.9.} {\em L'homomorphisme de $\mathbb{Z}$-schémas en groupes
$$
\omega_{\mathrm{S}}:
\mathrm{SO}_{P,f}
\longrightarrow
\mathrm{SO}_{L}
$$
est un isomorphisme.}

\medskip
\textit{Démonstration.} Soit $A$ un anneau commutatif unitaire.

\smallskip
On se convainc tout d'abord l'injectivité de l'homomorphisme $\mathrm{SO}_{P,f}(A)\to\mathrm{SO}_{L}(A)$. Pour cela on contemple le diagramme commutatif de $\mathbb{Z}$-schémas en groupes
$$
\begin{CD}
\mathrm{O}_{\mathrm{H}(\mathbb{Z}),f}@>>>\mathrm{O}_{P,f}@>\widetilde{\mathop{\mathrm{d\acute{e}t}}}>>\mathbb{Z}/2 \\
@V\omega VV @V\omega VV @VVV \\
\mathrm{O}_{\mathrm{A}_{1}} @>>> \mathrm{O}_{L}
@>\mathop{\mathrm{d\acute{e}t}}>>  \mu_{2}
\end{CD} 
$$
dans lequel les deux flèches horizontales de gauche sont les inclusions évidentes. Comme un élément de $\mathrm{O}_{P}(A)$ dont la restriction à $A\otimes_{\mathbb{Z}}\mathrm{H}(\mathbb{Z}^{n})$ est l'identité s'identifie à un élément de $\mathrm{O}_{\mathrm{H}(\mathbb{Z})}(A)$, cette contemplation montre que les noyaux des $\mathrm{O}_{P,f}(A)\to\mathrm{O}_{L}(A)$ et $\mathrm{O}_{\mathrm{H}(\mathbb{Z}),f}(A)\to\mathrm{O}_{\mathrm{A}_{1}}(A)$ coïncident~; on conclut en observant que la composée des deux flèches horizontales du haut est un isomorphisme.

\smallskip
On se convainc ensuite de ce que l'homomorphisme $\mathrm{SO}_{P,f}(A)\to\mathrm{SO}_{L}(A)$ est surjectif. Soit $\alpha$ un élément de $\mathrm{O}_{L}(A)$. On écrit $A\otimes_{\mathbb{Z}}L=\mathrm{H}(A^{n})\oplus A\hspace{1pt}e$ et on pose $M=\alpha(\mathrm{H}(A^{n}))$. Soit $M^{\perp}$ l'orthogonal dans $P$ de $M$ vu comme un sous-module de $A\otimes_{\mathbb{Z}}P$~; $M^{\perp}$ possède les propriétés suivantes~:

\smallskip
-- $M^{\perp}$ est un $A$-module projectif de rang $2$,

\smallskip
-- la restriction à $M^{\perp}$ de la forme quadratique de $A\otimes_{\mathbb{Z}}P$ est non-dégénérée,

\smallskip
-- le $\mathrm{q}$-module $A\otimes_{\mathbb{Z}}P$ est isomorphe à la somme orthogonale $M\oplus M^{\perp}$,

\smallskip
-- le discriminant $\Delta(M^{\perp})$ est trivial (voir \cite[Expos\'e XII, 1.11]{sga7}),

\smallskip
-- $\alpha(e)$ et $f$ appartiennent à $M^{\perp}$.

\smallskip
La proposition suivante, dont la démonstration est laissée au lecteur, montre que le\linebreak $\mathrm{q}$-module $M^{\perp}$ est isomorphe à $\mathrm{H}(A)$ (observer que l'on a $\mathrm{q}(\alpha(e))=1$)~:

\bigskip
\textbf{Proposition 1.10.} {\em Soit $N$ un $\mathrm{q}$-module sur $A$ de rang $2$. Les propriétés suivantes sont équivalentes~:
\begin{itemize}
\item [(i)] $N$ est isomorphe au $\mathrm{q}$-module hyperbolique $\mathrm{H}(A)$~;
\item[(ii)] le discriminant $\Delta(M)$ est trivial et il existe $e$ dans $M$ avec $\mathrm{q}(e)=1$.
\end{itemize}}

\medskip
Soit $\gamma$ un automorphisme du $\mathrm{q}$-module $A\otimes_{\mathbb{Z}}P=\mathrm{H}(A^{n})\oplus\mathrm{H}(A)$ induit par l'isomorphisme $\alpha:\mathrm{H}(A^{n})\to M$, un isomorphisme $\beta:\mathrm{H}(A)\to M^{\perp}$ et l'isomorphisme $M\oplus M^{\perp}\to\nolinebreak P$. Comme le groupe orthogonal $\mathrm{O}_{\mathrm{H}(\mathbb{Z})}(A)$ agit transitivement sur l'ensemble des $x$ avec $\mathrm{q}(x)=\nolinebreak 1$, on peut supposer $\beta(f)=f$ et donc $\gamma\in\mathrm{O}_{P,f}(A)$~; comme l'homomomorphisme $\widetilde{\mathop{\mathrm{d\acute{e}t}}}:\mathrm{O}_{\mathrm{H}(\mathbb{Z}),f}\to\mathbb{Z}/2$ est un isomorphisme on peut supposer en outre $\gamma\in\mathrm{SO}_{P,f}(A)$. Par construction $\alpha^{-1}\circ\omega_{\mathrm{S}}(\gamma)$ est alors un élément de $\mathrm{SO}_{L}(A)$ qui est l'identité sur $\mathrm{H}(A^{n})$~; or un tel élément est l'identité.
\hfill$\square$

\medskip
\textit{Remarque.} L'homomorphisme $\mathrm{O}_{\mathrm{H}(\mathbb{Z}),f}\to\mathrm{O}_{P,f}$ s'identifie à un homomorphisme $\mathbb{Z}/2\to\mathrm{O}_{P,f}$ qui est une section ``centrale'' de l'homomorphisme $\widetilde{\mathop{\mathrm{d\acute{e}t}}}:\mathrm{O}_{P,f}\to\mathbb{Z}/2$~; il en résulte que le groupe $\mathrm{O}_{P,f}$ est canoniquement isomorphe au produit $\mathrm{SO}_{P,f}\times\mathbb{Z}/2$. Pareillement le groupe $\mathrm{O}_{L}$ est canoniquement isomorphe au produit $\mathrm{SO}_{L}\times\mu_{2}$ et l'homomorphisme $\omega:\mathrm{O}_{P,f}\to\mathrm{O}_{L}$ s'identifie au produit de l'isomorphisme $\omega_{S}$ et de l'homomorphisme canonique $\mathbb{Z}/2\to\mu_{2}$.

\medskip
Nous achevons ce paragraphe en dégageant l'énoncé suivant dont le lecteur n'aura pas de mal à décoder les notations~:

\medskip
\textbf{Scholie 1.11.} {\em Soit $A$ un anneau comutatif unitaire. Soit $P$ un $\mathrm{q}$-module sur $A$, de rang constant pair, muni d'un élément $e$ avec $\mathrm{q}(e)\in A^{\times}$~; soit $L$ l'orthogonal de $e$ dans $P$. Alors~:

\smallskip
-- $L$ est un $\mathrm{qi}$-module sur $A$,

\smallskip
-- le $A$-groupe $\mathrm{SO}_{L}$ s'identifie au $A$-groupe $\mathrm{SO}_{P,e}$,

\smallskip
-- les $A$-groupes  $\mathrm{O}_{P,e}$ et $\mathrm{O}_{L}$ s'identifient respectivement au produits $\mathrm{SO}_{L}\times\mathbb{Z}/2$ et $\mathrm{SO}_{L}\times\mu_{2}$,

\smallskip
-- l'homomorphisme canonique $\mathrm{O}_{P,e}\to\mathrm{O}_{L}$ s'identifie à l'homomorphisme induit par\linebreak l'homomorphisme $\mathbb{Z}/2\to\mu_{2}$.}

\normalsize
\vspace{0,75cm}
\textbf{2. Sur les $\mathrm{q}\text{-}\mathrm{i}$-modules sur $\mathbb{Z}$}

\bigskip
Un  $\mathrm{q}\text{-}\mathrm{i}$-module sur $\mathbb{Z}$ n'est rien d'autre qu'un $\mathbb{Z}$-module libre $L$, de dimension finie impaire, muni d'une forme bilinéaire symétrique paire (c'est-à-dire telle que l'on a $x.x$ pair pour tout $x$ dans $L$) avec $\vert\hspace{-1.75pt}\mathop{\mathrm{d\acute{e}t}}L\vert=2$. En fait, le point (a) de la proposition ci-dessous montre que les propriétés de la forme bilinéaire entraînent celle de la dimension.

\bigskip
\textbf{Proposition 2.1} (classification des $\mathbb{Z}$-$\mathrm{q}\text{-}\mathrm{i}$-modules)\textbf{.}

{\em 
\medskip
{\em (a)} Soit $L$ un $\mathbb{Z}$-module libre de dimension finie, muni d'une forme bilinéaire symétrique paire, avec $\vert\hspace{-1.75pt}\mathop{\mathrm{d\acute{e}t}}L\vert=2$.  Alors la signature de $L$ vérifie la congru\-ence $\tau(L)\equiv\pm 1\pmod{8}$ (si bien que la dimension de $L$ est impaire).

\medskip
{\em (b)} Soit $P$ un $\mathbb{Z}$-module libre de dimension finie, muni d'une forme bilinéaire symétrique paire, avec $\vert\hspace{-1.75pt}\mathop{\mathrm{d\acute{e}t}}L\vert=1$, et d'un élément $e$ avec $e.e=2\hspace{1pt}\epsilon$. Soit $L$ l'orthogonal de $e$, muni de la forme bilinéaire symétrique paire restriction de celle de $P$. Alors est $L$ est un $\mathbb{Z}$-module libre de dimension finie ($P/L$ étant libre de dimension $1$) avec $\mathop{\mathrm{d\acute{e}t}}L=2\hspace{1pt}\epsilon\mathop{\mathrm{d\acute{e}t}}P$.

\medskip
{\em (c)} Soit $(n,\epsilon)$ un élément de $\mathbb{N}\times\{\pm 1\}$ avec $n$ impair~; soit $\mathrm{QI}_{n,\epsilon}$ l'ensemble des classes d'isomorphisme de $\mathbb{Z}$-$\mathrm{q\text{-}i}$-modules $L$, avec $\dim L=n$ et $\tau(L)\equiv\nolinebreak\epsilon\pmod{8}$. Soit $(n,\epsilon)$ un élément de $\mathbb{N}\times\{\pm 1\}$ avec $n$ pair~; soit $\mathrm{QR}_{n,\epsilon}$\linebreak l'ensemble des classes d'isomorphisme de $\mathbb{Z}$-$\mathrm{q}$-modules munis d'un élément~$e$ avec $e.e=2\hspace{1pt}\epsilon$. Alors l'application
$$
(P,e)
\hspace{4pt}\mapsto\hspace{4pt}
e^{\perp}
$$
induit une bijection de $\mathrm{QR}_{n,\epsilon}$ sur $\mathrm{QI}_{n-1,-\epsilon}$.}

\bigskip
\textit{Démonstration du point (a) de 2.1.} On considère $L$ comme un $\widetilde{\mathrm{q}}$-module. Le groupe sous-jacent à son résidu $\mathop{\mathrm{r\acute{e}s}}L$ s'identifie à $\mathbb{Z}/2$ et sa forme quadratique d'enlacement vérifie $\mathrm{q}(\bar{1})=\frac{\epsilon}{4}$ avec $\epsilon=\pm{1}$. On est donc amené à introduire la somme orthogonale de $\widetilde{\mathrm{q}}$-modules $L\oplus\langle -\epsilon\rangle\otimes\mathrm{A}_{1}$ ($\mathrm{A}_{1}=\mathrm{Q}(\mathbf{A}_{1})$ n'est rien d'autre que le $\mathbb{Z}$-module $\mathbb{Z}$ muni de la forme quadratique $x\mapsto x^{2}$) dont le résidu s'identifie à la somme orthogonale de $\mathrm{qe}$-modules $\mathop{\mathrm{r\acute{e}s}}L\oplus\langle -1\rangle\otimes\mathop{\mathrm{r\acute{e}s}}L$. Ce $\mathrm{qe}$-module possède un unique lagrangien à savoir la diagonale~; on note $P$ le $\mathrm{q}$-module qui correspond à ce lagrangien \textit{via} II.1.1. La signature de $P$ étant divisible par $8$, on a bien la congruence $\tau(L)\equiv\epsilon\pmod{8}$.

\bigskip
\textit{Démonstration du point (b) de 2.1.} Elle résulte par exemple du point (c) de l'énoncé suivant dont il sera commode de disposer dans ce mémoire. La vérification de cet énoncé est laissée au lecteur.

\bigskip
\textbf{Proposition 2.2.} {\em Soit $A$ un anneau de Dedekind. Soient $L$ un $\mathrm{q}$-module sur~$A$, $M$ un sous-module et $M^{\perp}$ son orthogonal.

\smallskip
On suppose que $M$ est facteur direct dans $L$ (autrement dit que le quotient $L/M$ est sans torsion, on observera que $M^{\perp}$ est facteur direct dans
$L$ pour tout $M$).

\smallskip
Alors~:

\smallskip
{\em (a)} On a $(M^{\perp})^{\perp}=M$.

\smallskip
On suppose en outre que la restriction à $M$ de la forme bilinéaire de $L$ est non singulière, en clair que l'homomorphisme induit $M\to\mathrm{Hom}_{A}(M,A)$ est injectif~; $M$ est donc un $\widetilde{\mathrm{q}}$-module sur $A$.

\smallskip
Alors~:

\smallskip
{\em (b)} La restriction à $M^{\perp}$ de  cette forme
bilinéaire est aussi non singulière.

\smallskip
{\em (c)} L'homomorphisme de $A$-modules canonique
$M\oplus M^{\perp}\to L$ est injectif et l'on a une suite exacte de $A$-modules
$$
\hspace{24pt}
0
\longrightarrow M\oplus M^{\perp}
\longrightarrow L
\longrightarrow\mathop{\mathrm{r\acute{e}s}}M\longrightarrow 0
\hspace{24pt},
$$
dans laquelle l'homomorphisme $L\to\mathop{\mathrm{r\acute{e}s}}M$ est le composé de
l'isomorphisme $L\to\mathrm{Hom}_{A}(L,A)$ induit par
la forme bilinéaire de $L$, et des homomorphismes canoniques $\mathrm{Hom}_{A}(L,A)\to
\mathrm{Hom}_{A}(M,A)$, $\mathrm{Hom}_{A}(M,A)\to
\mathop{\mathrm{r\acute{e}s}}M$.

\smallskip
{\em (d)} Les isomorphismes de $A$-modules,
$\mathop{\mathrm{r\acute{e}s}}M\cong L/(M\oplus M^{\perp})$ et  $\mathop{\mathrm{r\acute{e}s}}
M^{\perp}\cong L/(M\oplus M^{\perp})$
(observer que $M$ et $M^{\perp}$
jouent des r\^{o}les symétriques), induisent un isomorphisme de $A$-modules
$\varphi:\mathrm{r\acute{e}s}
\hspace{1pt}M\to\mathop{\mathrm{r\acute{e}s}}M^{\perp}$
tel que l'on a\linebreak $\mathrm{q}(\varphi(\xi))=-\mathrm{q}(\xi)$ pour tout $\xi$ dans $\mathop{\mathrm{r\acute{e}s}}M$. En d'autres termes on a un isomorphisme canonique de $A$-$\mathrm{qe}$-modules
$$
\hspace{24pt}
\mathrm{r\acute{e}s}
\hspace{1pt}M^{\perp}\cong
\langle -1\rangle
\hspace{-1pt}\otimes\hspace{-1pt}
\mathop{\mathrm{r\acute{e}s}}M
\hspace{24pt}.
$$

\smallskip
{\em (e)} Le $\mathrm{q}$-module $L$
correspond via II.1.1 et
l'isomorphisme de $\mathrm{qe}$-modules
$\mathop{\mathrm{r\acute{e}s}}(M
\oplus M^{\perp})
\cong\mathop{\mathrm{r\acute{e}s}}M
\oplus\mathrm{r\acute{e}s}
\hspace{1pt}M^{\perp}$
au lagrangien de
$\mathop{\mathrm{r\acute{e}s}}M
\oplus\mathrm{r\acute{e}s}
\hspace{1pt}M^{\perp}$ qu'est le graphe
de $\varphi$.}

\bigskip
\textit{Remarques}

\smallskip
-- Il existe une ``version bilinéaire'' de la proposition 2.2 dans laquelle les $\mathrm{q}$-modules (resp. $\widetilde{\mathrm{q}}$-modules, resp. $\mathrm{qe}$-modules) sont remplacés par des\linebreak $\mathrm{b}$-modules (resp. $\widetilde{\mathrm{b}}$-modules, resp. $\mathrm{e}$-modules).

\smallskip
-- Il existe aussi, dans un cas particulier, une ``version bilinéaire-quadratique" de la proposition 2.2. Précisons un peu. Soit $L$ un $\mathrm{b}$-module impair sur $\mathbb{Z}$~; soit $u$ un vecteur de Wu de $L$ que l'on suppose indivisible et non isotrope. Par construction, l'orthogonal $u^{\perp}$ de $u$ dans $L$ est un $\widetilde{\mathrm{b}}$-module pair, en d'autres termes un $\widetilde{\mathrm{q}}$-module. La proposition 2.2 dit que $\mathop{\mathrm{r\acute{e}s}}u^{\perp}$, en tant que $\mathrm{e}$-module, est isomorphe à $\mathbb{Z}/u.u$ muni de la forme d'enlacement défini par $\bar{1}.\bar{1}=-\frac{1}{u.u}$. On vérifie que la forme quadratique d'enlacement est quant à elle définie par
$$
\hspace{24pt}
\mathrm{q}(\bar{1})
\hspace{4pt}=\hspace{4pt}
\frac{1}{2}\hspace{2pt}(1-\frac{1}{u.u})
\hspace{24pt}.
$$

\bigskip
\textit{Démonstration du point (c) de 2.1.} On note $\omega_{n,\epsilon}:\mathrm{QR}_{n,\epsilon}\to\mathrm{QI}_{n-1,-\epsilon}$ l'application induite par l'application $(P,e)\mapsto e^{\perp}$. Soit $(n,\epsilon)$ un élément de $\mathbb{N}\times\{\pm 1\}$ avec $n$ impair~; le $\mathrm{q}$-module $P$ qui apparaît dans la démonstration du point (a) de 2.1 est par construction muni d'un élément e avec $e.e=-2\hspace{1pt}\epsilon$. L'unicité du lagrangien invoquée dans cette construction montre que l'application $L\mapsto(P,e)$ induit une application $\mathrm{QI}_{n,\epsilon}\to\mathrm{QR}_{n+1,-\epsilon}$ que l'on note $\pi_{n,\epsilon}$. On constate que les deux applications $\omega_{n,\epsilon}$ et $\pi_{n-1,-\epsilon}$ sont inverses l'une de l'autre.
\hfill$\square$

\bigskip
\textbf{Scholie 2.3} {\em Soit $L$ un $\mathrm{q\text{-}i}$-modules sur $\mathbb{Z}$~; soit $\epsilon$ l'élément de $\{\pm 1\}$ défini par $\tau(L)\equiv\epsilon\pmod{8}$. Alors le $\mathrm{qe}$ module $\mathop{\mathrm{r\acute{e}s}}L$ est isomorphe à $\mathbb{Z}/2$ muni de la forme quadratique d'enlacement définie par $\mathrm{q}(\bar{1})=\frac{\epsilon}{4}$.}

\bigskip
La proposition 2.1 et le scholie II.2.1 conduisent également au scholie suivant~:

\bigskip
\textbf{Scholie 2.4.} {\em Soient $L_{1}$ et $L_{2}$ deux $\mathrm{q\text{-}i}$-modules sur $\mathbb{Z}$. Les deux conditions suivantes sont équivalentes~:
\begin{itemize}
\item [(i)] les deux $\mathrm{b}$-espaces vectoriels sur $\mathbb{Q}$, $\mathbb{Q}\otimes_{\mathbb{Z}}L_{1}$ et $\mathbb{Q}\otimes_{\mathbb{Z}}L_{2}$, sont iso\-morphes~;
\item [(ii)] les deux $\mathrm{b}$-espaces vectoriels sur $\mathbb{R}$, $\mathbb{R}\otimes_{\mathbb{Z}}L_{1}$ et $\mathbb{R}\otimes_{\mathbb{Z}}L_{2}$, sont isomorphes.
\end{itemize}}

\vspace{0,75cm}
\textsc{Genre d'un $\mathrm{q\text{-}i}$-module sur $\mathbb{Z}$}

\bigskip
La méthode de démonstration de la proposition 2.1 fournit l'énoncé suivant (l'intertitre fait référence au point (b))~:
\vfill\eject

\bigskip
\textbf{Proposition 2.5.} {\em Soit $L$ un $\mathrm{q\text{-}i}$-module sur $\mathbb{Z}$ de dimension $2n+1$ et de déterminant $2\hspace{1pt}\epsilon$ avec $\epsilon=\pm1$~; soit $p$ un nombre premier.

\medskip
{\em (a)} Le $\mathrm{q\text{-}i}$-espace vectoriel $\mathbb{F}_{p}\otimes_{\mathbb{Z}}L$ est isomorphe à  $\mathrm{H}(\mathbb{F}_{p}^{n})\oplus\langle (-1)^{n}\hspace{1pt}\epsilon\rangle\otimes\mathrm{A}_{1}$.

\medskip
{\em (b)} Le $\mathrm{q\text{-}i}$-module $\mathbb{Z}_{p}\otimes_{\mathbb{Z}}L$ est isomorphe à $\mathrm{H}(\mathbb{Z}_{p}^{n})\oplus\langle (-1)^{n}\hspace{1pt}\epsilon\rangle\otimes\mathrm{A}_{1}$.}

\bigskip
\textit{Remarques.} Pour $p\not=2$, un $\mathrm{q\text{-}i}$-espace vectoriel sur $\mathbb{F}_{p}$ (resp. un $\mathrm{q\text{-}i}$-module sur $\mathbb{Z}_{p}$) n'est rien d'autre qu'un $\mathrm{b}$-espace vectoriel (resp. $\mathrm{b}$-module) de dimension impaire. Soient $A$ un anneau (commutatif unitaire) et $u$ un élément de $A^{\times}$, alors on a $\langle u\rangle\otimes\mathrm{H}(A^{n})\cong\mathrm{H}(A^{n})$~; on peut donc remplacer dans l'énoncé 2.4 $\mathrm{H}(\mathbb{F}_{p}^{n})\oplus\langle (-1)^{n}\hspace{1pt}\epsilon\rangle\otimes\mathrm{A}_{1}$ par $\langle (-1)^{n}\hspace{1pt}\epsilon\rangle\otimes(\mathrm{H}(\mathbb{F}_{p}^{n})\oplus\mathrm{A}_{1})$ et $\mathrm{H}(\mathbb{Z}_{p}^{n})\oplus\langle (-1)^{n}\hspace{1pt}\epsilon\rangle\otimes\mathrm{A}_{1}$ par $\langle (-1)^{n}\hspace{1pt}\epsilon\rangle\otimes(\mathrm{H}(\mathbb{Z}_{p}^{n})\oplus\mathrm{A}_{1})$.

\vspace{0,75cm}
\textsc{Le cas défini positif}

\bigskip
Soit $L$ un  $\mathrm{q\text{-}i}$-module sur $\mathbb{Z}$ avec $\mathbb{R}\otimes_{\mathbb{Z}}L$ défini positif. Compte tenu de ce que nous avons vu plus haut, un tel $L$ n'est rien d'autre qu'un réseau (entier) pair de déterminant $2$. Nous abandonnons ci-après la terminologie ``$\mathrm{q\text{-}i}$-module sur $\mathbb{Z}$ défini positif'' (qui est loin d'être classique~!) pour la terminologie ``réseau pair de déterminant $2$''.

\bigskip
Soit $L$ un réseau pair de déterminant $2$~; le point (a) de 2.1 montre que l'on a la congruence $\dim L\equiv\pm 1\pmod{8}$.

\bigskip
On étudie d'abord le cas $\dim L\equiv-1\pmod{8}$. Le point (c) de la proposition 2.1 se spécialise de la façon suivante~:

\bigskip
\textbf{Proposition 2.6.} {\em Soit $n>0$ un entier avec $n\equiv-1\pmod{8}$~; soit $\mathrm{X}_{n}$ l'ensemble des classes d'isomorphisme de réseaux pairs $L$, avec $\dim L=n$ et $\mathop{\mathrm{d\acute{e}t}}L=2$. Soit $n>0$ un entier avec $n\equiv 0\pmod{8}$ pair~; soit $\mathrm{X}_{n}^{\mathrm{A}_{1}}$ l'ensemble des classes d'isomorphisme de réseaux unimodulaires pairs $P$ de dimension $n$ munis d'un élément~$e$ avec $e.e=2$ (en d'autres termes d'une racine). Alors l'application
$$
(P;e)
\hspace{4pt}\mapsto\hspace{4pt}
e^{\perp}
$$
induit une bijection de $\mathrm{X}_{n}^{\mathrm{A}_{1}}$ sur $\mathrm{X}_{n-1}$.}

\bigskip
(La justification de la notation $\mathrm{X}_{n}^{\mathrm{A}_{1}}$ est la suivante~: la donnée d'une racine de $P$ est équivalente à celle d'une représentation de $\mathrm{A}_{1}$ par $P$. Elle est le pendant de la notation $\mathrm{X}_{n}^{\mathrm{E}_{7}}$ qui sera introduite plus bas.)
\vfill\eject

\bigskip
\textit{Exemples}

\medskip
\textit{Détermination de $\mathrm{X}_{7}$.} Comme $\mathrm{X}_{8}$ n'a qu'un élément, à savoir la classe de~$\mathrm{E}_{8}$, et que le groupe Weyl de $\mathbf{E}_{8}$ (qui coïncide avec le groupe orthogonal de~$\mathrm{E}_{8}$) agit transitivement sur l'ensemble des racines, l'ensemble $\mathrm{X}_{7}$ n'a qu'un élé\-ment, à savoir la classe de l'orthogonal d'une racine dans $\mathrm{E}_{8}$. Cet orthogonal est noté $\mathrm{E}_{7}$, notation en accord avec celle que nous avons adoptée au chapitre~II~: $\mathrm{E}_{7}=\mathrm{Q}(\mathbf{E}_{7})$.

\medskip
\textit{Détermination de $\mathrm{X}_{15}$}. Comme le groupe de Weyl de $\mathbf{D}_{16}$ (qui coïncide avec le groupe orthogonal de $\mathrm{E}_{16}$) agit transitivement sur l'ensemble des racines et qu'il en est de même pour le groupe orthogonal de $\mathrm{E}_{8}\oplus\mathrm{E}_{8}$, l'ensemble $\mathrm{X}_{15}$ a deux éléments~:

\smallskip
- la classe de l'orthogonal d'une racine dans $\mathrm{E}_{16}$, disons $\mathrm{E}_{15}$~;

\smallskip
- la classe de $\mathrm{E}_{7}\oplus\mathrm{E}_{8}$.

\bigskip
\textit{Détermination de $\mathrm{X}_{23}$}

\medskip
Compte tenu de la proposition 2.6 et du théorème II.3.17, cette détermination est conséquence de l'observation suivante~:

\medskip
Soit $L$ un réseau unimodulaire pair de dimension $24$ avec racines. Soient $e_{1}$ et $e_{2}$ deux racines de $L$ ; soit $R_{i}$, $i=1,2$, la composante irréductible du système de racines $\mathrm{R}(L)$ à laquelle appartient $e_{i}$. Alors les deux conditions suivantes sont équivalentes~:

\smallskip
(i) Les deux systèmes de racines $R_{1}$ et $R_{2}$ sont isomorphes.

\smallskip
(ii) Il existe un élément $\alpha$ du groupe orthogonal $\mathrm{O}(L)$ tel que l'on a $\alpha(e_{1})=e_{2}$.

\smallskip
L'implication (ii)$\Rightarrow$(i) est évidente. L'implication (i)$\Rightarrow$(ii) se vérifie par inspection. Nous détaillons un peu cette vérification ci-après.

\medskip
On considère la décomposition en composantes irréductibles du système de racines $\mathrm{R}(L)$~:
$$
\hspace{24pt}
\mathrm{R}(L)
\hspace{4pt}\simeq\hspace{4pt}
\coprod_{R\in\mathcal{R}}\hspace{2pt}\mathrm{m}(R)\hspace{1pt}R
\hspace{24pt},
$$
$\mathcal{R}$ désignant, ci-dessus,  l'ensemble des classes d'isomorphisme de système de racines irréductibles de type ADE, $\mathrm{m}:\mathcal{R}\to\mathbb{N}$ une application à support fini et $\mathrm{m}(R)\hspace{1pt}R$ la réunion disjointe de $\mathrm{m}(R)$ copies de $R$. On rappelle que l'on note $\mathrm{A}(\mathrm{R}(L))$ le groupe orthogonal du réseau $\mathrm{Q}(\mathrm{R}(L))$ et que le groupe de Weyl $\mathrm{W}(\mathrm{R}(L))$ est  un sous-groupe distingué de $\mathrm{A}(\mathrm{R}(L))$. On rappelle également que l'on note $\mathrm{G}(\mathrm{R}(L))$ le groupe quotient $\mathrm{A}(\mathrm{R}(L))/\mathrm{W}(\mathrm{R}(L))$ et que l'on a un isomorphisme de groupes (canonique en un sens manifeste)~:
$$
\mathrm{G}(\mathrm{R}(L))
\hspace{2pt}\cong\hspace{2pt}
\prod_{R\in\mathcal{R}}\hspace{4pt}
(\hspace{2pt}\mathrm{G}(R)^{\mathrm{m}(R)}\rtimes\mathfrak{S}_{\mathrm{m}(R)}\hspace{2pt})
\hspace{2pt}=\hspace{2pt}
(\hspace{2pt}\prod_{R\in\mathcal{R}}\hspace{2pt}\mathrm{G}(R)^{\mathrm{m}(R)}\hspace{2pt})
\rtimes
(\hspace{2pt}\prod_{R\in\mathcal{R}}\hspace{2pt}\mathfrak{S}_{\mathrm{m}(R)}\hspace{2pt})
$$
($\mathfrak{S}_{\mathrm{m}(R)}$ désignant le groupe symétrique d'ordre $\mathrm{m}(R)$ qui opère de façon évidente sur le groupe $\mathrm{G}(R)^{\mathrm{m}(R)}$). On paraphrase  enfin le point (b) du scholie II.3.15~: le groupe orthogonal $\mathrm{O}(L)$ est le sous-groupe de $\mathrm{A}(\mathrm{R}(L))$ image inverse par l'homomorphisme $\mathrm{A}(\mathrm{R}(L))\to\mathrm{G}(\mathrm{R}(L))$  du sous-groupe qui stabilise le lagrangien $L/\mathrm{Q}(\mathrm{R}(L))$ du $\mathrm{qe}$-module $\mathop{\mathrm{r\acute{e}s}}\mathrm{Q}(\mathrm{R}(L))$. Dans \cite{erokhin2}, V.~A.~Erokhin explicite au cas par cas ce stabilisateur, qu'il note $\mathrm{H}(L)$, pour les $23$ classes d'isomorphisme de réseaux unimodulaires pairs de dimension $24$ avec racines. Il est clair que l'on  a une suite exacte de groupes canonique
$$
\hspace{24pt}
1\to\mathrm{H}_{1}(L)\to\mathrm{H}(L)\to\mathrm{H}_{2}(L)\to 1
\hspace{24pt},
$$
$\mathrm{H}_{1}(L)$ s'identifiant à un sous-groupe du produit $\prod_{R\in\mathcal{R}}\mathrm{G}(R)^{\mathrm{m}(R)}$ et $\mathrm{H}_{2}(L)$ à un sous-groupe du produit $\prod_{R\in\mathcal{R}}\mathfrak{S}_{\mathrm{m}(R)}$, à savoir l'image de la restriction à $\mathrm{H}(L)$ de l'homomorphisme canonique $\mathrm{G}(\mathrm{R}(L))\to\prod_{R\in\mathcal{R}}\mathfrak{S}_{\mathrm{m}(R)}$ (bien que cela ne soit pas explicitement précisé, les groupes $\mathrm{H}_{1}(L)$  et $\mathrm{H}_{2}(L)$ sont les groupes $G_{1}$ et $G_{2}$ dont les cardinaux apparaissent dans les colonnes 5 et 6 de \cite[Ch. 16, Table 16.1]{conwaysloane}).

\medskip
L'implication (i)$\Rightarrow$(ii) qui nous occupe résulte du fait que pour tout $R$ dans $\mathcal{R}$ l'image de $\mathrm{H}_{2}(L)$ dans $\mathfrak{S}_{\mathrm{m}(R)}$ est un sous-groupe transitif, ce que l'on vérifie en parcourant l'inventaire de \cite{erokhin2} (évidemment les seuls $R$ à considérer sont ceux pour lesquels on a $\mathrm{m}(R)\geq 2$~!). Explicitons, à titre d'exemple, ces sous-groupes transitifs, pour les 5 premiers systèmes de racines de l'inventaire en question~:

\medskip
1) $\mathrm{R}(L)=24\mathbf{A}_{1}$

\smallskip
L'image de $\mathrm{H}_{2}(L)$ dans $\mathfrak{S}_{24}$ est le groupe de Mathieu $\mathrm{M}_{24}$.

\medskip
2) $\mathrm{R}(L)=12\mathbf{A}_{2}$

\smallskip
L'image de $\mathrm{H}_{2}(L)$ dans $\mathfrak{S}_{12}$ est le groupe de Mathieu $\mathrm{M}_{12}$.

\medskip
3) $\mathrm{R}(L)=8\mathbf{A}_{3}$

\smallskip
Il existe une bijection de l'ensemble $\{1,2,\ldots,8\}$ sur l'ensemble sous-jacent à $\mathbb{F}_{2}^{3}$, vu comme un espace affine de dimension $3$ sur $\mathbb{F}_{2}$, qui induit un isomorphisme de $\mathrm{H}_{2}(L)$ sur le sous-groupe des transformations affines.

\medskip
4) $\mathrm{R}(L)=6\mathbf{A}_{4}$

\smallskip
Il existe une bijection de l'ensemble $\{1,2,\ldots,6\}$ sur l'ensemble sous-jacent à $\mathbf{P}^{1}(\mathbb{F}_{5})$ qui induit un isomorphisme de $\mathrm{H}_{2}(L)$ sur le sous-groupe des transformations projectives. (Attention : la liste de générateurs donnée dans \cite{erokhin2} pour $\mathrm{H}(L)$ est incomplète.)

\smallskip
\footnotesize
Plus précisément, la suite exacte $1\to\mathrm{H}_{1}(L)\to\mathrm{H}(L)\to\mathrm{H}_{2}(L)\to 1$ est dans ce cas isomorphe à la suite exacte $1\to\mathbb{F}_{5}^{\times}/\{\pm 1\}\to\mathrm{GL}_{2}(\mathbb{F}_{5})/\{\pm\mathrm{Id}\}\to\mathrm{PGL}_{2}(\mathbb{F}_{5})\to 1$.
\normalsize

\medskip
5) $\mathrm{R}(L)=4\mathbf{A}_{6}$

\smallskip
L'image de $\mathrm{H}_{2}(L)$ dans $\mathfrak{S}_{4}$ est le sous-groupe alterné  $\mathfrak{A}_{4}$.

\medskip
On constate donc au bout du compte que l'ensemble  $\mathrm{X}_{23}$ s'identifie au sous-ensemble du produit $\mathrm{X}_{24}\times\mathcal{R}$ constitué des couples $(x,r)$ tels que $r$ est la classe d'isomorphisme d'une composante irréductible du système de racines $\mathrm{R}(x)$ (l'abus de notation est véniel). Le cardinal de $\mathrm{X}_{23}$ s'obtient en contemplant la deuxième colonne de \cite[Ch. 16,  Table 16.1]{conwaysloane} (Table \ref{listeniemeier} de notre mémoire)~:
$$
\vert\mathrm{X}_{23}\vert
\hspace{4pt}=\hspace{4pt}
32
\scriptstyle{\hspace{4pt}=\hspace{4pt}1+2+1+1+1+2+2+2+1+1+3+1+2+1+1+2+1+2+1+1+1+1+1+0}
$$
(ce cardinal est évidemment à comparer au nombre des représentations de la table \ref{tableso23}).

\bigskip
On étudie maintenant le cas $\dim L\equiv 1\pmod{8}$.

\medskip
\textbf{Proposition 2.7.} {\em Soit $n>0$ un entier avec $n\equiv 1\pmod{8}$~; soit $\mathrm{X}_{n}$ l'ensemble des classes d'isomorphisme de réseaux pairs $L$, avec $\dim L=n$ et $\mathop{\mathrm{d\acute{e}t}}L=2$. Soit $n>0$ un entier avec $n\equiv 0\pmod{8}$ pair~; soit $\mathrm{X}_{n}^{\mathrm{E}_{7}}$ l'ensemble des classes d'isomorphisme de réseaux unimodulaires pairs $P$ de dimension $n$ munis d'un homomorphisme $f:\mathrm{E}_{7}\to P$ avec $\mathrm{q}(f(x))=\mathrm{q}(x)$ pour tout $x$ dans $\mathrm{E}_{7}$ (en d'autres termes d'une représentation de $\mathrm{E}_{7}$ par $P$).

\medskip
Alors l'application
$$
(P;f)
\hspace{4pt}\mapsto\hspace{4pt}
(f(\mathrm{E}_{7}))^{\perp}
$$
induit une bijection de $\mathrm{X}_{n}^{\mathrm{E}_{7}}$ sur $\mathrm{X}_{n-7}$.}

\bigskip
\textit{Démonstration.} C'est une variante de la démonstration du point (c) de 2.1. On considère cette fois la somme orthogonale $L\oplus\mathrm{E}_{7}$. D'après le scholie 2.3, le résidu de ce que $\widetilde{\mathrm{q}}$-module est encore isomorphe à $\mathop{\mathrm{r\acute{e}s}}L\oplus\langle -1\rangle\otimes\mathop{\mathrm{r\acute{e}s}}L$. On achève \textit{mutatis mutandis}.
\hfill$\square$

\medskip
\textit{Notation-remarque.} Soient $G$ et $L$ deux réseaux entiers~; on rappelle qu'une représentation de $G$ par $L$ est un homomorphisme de $\mathbb{Z}$-modules $f:G\to L$ avec $f(x).f(y)=x.y$ pour tous $x$ et $y$ dans $G$ et que l'ensemble (fini) de ces $f$ est noté $\mathrm{Rep}(G,L)$. On note  $\overline{\mathrm{Rep}}(G,L)$ le quotient $\mathrm{Rep}(G,L)\slash\mathrm{O}(G)$ de l'action à droite du groupe orthogonal $\mathrm{O}(G)$ sur $\mathrm{Rep}(G,L)$ (action qui est libre)~; $\overline{\mathrm{Rep}}(G,L)$ peut être vu comme l'ensemble des sous-modules de $L$ isomorphes, comme réseaux entiers, à $G$.

\bigskip
\textit{Exemples}

\medskip
\textit{Détermination de $\mathrm{X}_{17}$.} On constate par inspection que les  seuls systèmes de racines irréductibles de type ADE qui contiennent $\mathbf{E}_{7}$ sont $\mathbf{E}_{7}$ et $\mathbf{E}_{8}$ (utiliser par exemple  \cite[Ch. VI, \S1, Proposition 24]{bourbaki}). On en déduit, toujours par inspection, que les seuls systèmes de racines apparaissant dans la classification de Niemeier qui contiennent $\mathbf{E}_{7}$ sont $\mathbf{R}_{2}=\mathbf{D}_{16}\coprod\mathbf{E}_{8}$, $\mathbf{R}_{3}=\mathbf{E}_{8}\coprod\mathbf{E}_{8}\coprod\mathbf{E}_{8}$,  $\mathbf{R}_{6}=\mathbf{A}_{17}\coprod\mathbf{E}_{7}$ et $\mathbf{R}_{7}=\mathbf{D}_{10}\coprod\mathbf{E}_{7}\coprod\mathbf{E}_{7}$. Soit $\mathrm{P}_{i}$, $i\in\{2,3,6,7\}$, le réseau unimodulaire pair avec $\mathrm{R}(\mathrm{P}_{i})\approx\mathbf{R}_{i}$~; on vérifie que le groupe orthogonal $\mathrm{O}(\mathrm{P}_{i})$ agit transitivement sur l'ensemble $\overline{\mathrm{Rep}}(\mathrm{E}_{7},\mathrm{P}_{i})$ dans les quatre cas. La vérification est immédiate dans les trois premiers~; pour le quatrième on utilise l'observation faite à la fin de la cinquième illustration que nous avons donnée au chapitre II de la proposition 3.9. La proposition 2.7 montre donc que l'ensemble $\mathrm{X}_{17}$ a $4$ éléments. Précisons lourdement. On choisit un sous-réseau de $\mathrm{P}_{i}$ isomorphe à $\mathrm{E}_{7}$ et on note $\mathrm{L}_{i}$ l'orthogonal de ce sous-réseau~; on a $\mathrm{X}_{17}=\{[\mathrm{L}_{2}],[\mathrm{L}_{3}],[\mathrm{L}_{6}],[\mathrm{L}_{7}]\}$ (la notation $[L]$ désignant la classe d'isomorphisme d'un réseau pair $L$ de dimension $17$ et de déterminant~$2$).

\medskip
Il n'est pas difficile de donner une définition \textit{ab initio} des réseaux $\mathrm{L}_{i}$~:

\smallskip
-- $\mathrm{L}_{2}=\mathrm{E}_{16}\oplus\mathrm{A}_{1}$~;

\smallskip
-- $\mathrm{L}_{3}=\mathrm{E}_{8}\oplus\mathrm{E}_{8}\oplus\mathrm{A}_{1}$;

\smallskip
-- $\mathrm{L}_{6}=\mathrm{A}_{17}^{+}$ (le $\mathrm{qe}$-module $\mathop{\mathrm{r\acute{e}s}}\mathrm{A}_{17}$ est isomorphe à $\mathbb{Z}/18$ avec $\mathrm{q}(\bar{k})=\frac{17\hspace{1pt}k^{2}}{36}$, $\mathrm{A}_{17}^{+}$ est le réseau pair correspondant \textit{via} le point (b) de II.1.1 au sous-module engendré par $\bar{6}$)~;

\smallskip
-- $\mathrm{L}_{7}=(\mathrm{D}_{10}\oplus\mathrm{E}_{7})^{+}$  (le $\mathrm{qe}$-module $\mathop{\mathrm{r\acute{e}s}}(\mathrm{D}_{10}\oplus\mathrm{E}_{7})$ contient deux sous-modules non triviaux isotropes, $(\mathrm{D}_{10}\oplus\mathrm{E}_{7})^{+}$ est le réseau pair correspondant \textit{via} le point (b) de II.1.1 à l'un ou l'autre de ces deux sous-modules).

\medskip
\textit{Détermination de $\mathrm{X}_{9}$.} La proposition 2.7 montre que $\mathrm{X}_{9}$ a un seul élément à savoir la classe du réseau $\mathrm{E}_{8}\oplus\mathrm{A}_{1}$.

\bigskip
\textit{Sur la détermination de  $\mathrm{X}_{25}$.} L'isomorphisme $\mathrm{X}_{25}\cong\mathrm{X}_{32}^{\mathrm{E}_{7}}$ de la proposition 2.7 n'est pas très utile pour déterminer $\mathrm{X}_{25}$ car $\mathrm{X}_{32}$ n'a pas encore été déterminé. Par contre le point (c) de la proposition 2.1 peut être mis en oeuvre. Le point en question dit que $\mathrm{X}_{25}$ est en bijection avec l'ensemble des classes d'isomorphisme des couples $(P;e)$, $P$ désignant un $\mathrm{q}$-module sur $\mathbb{Z}$ de dimension $26$ et de signature $24$ et $e$ un élément de $P$ avec $e.e=-2$. Or tous les $P$ de ce type sont isomorphes d'après le théorème II.2.7. On peut choisir comme représentant de cette classe d'isomorphisme le $\mathrm{q}$-module $\mathrm{II}_{25,1}$ (notation de \cite{conwaysloane}\cite{borcherdsthese}) qui peut être vu comme le réseau du $\mathrm{q}$-espace vectoriel $\mathbb{Q}^{26}$ muni de la forme quadratique $\frac{1}{2}(\sum_{i=1}^{25}x_{i}^{2}-x_{26}^{2})$, engendré par le sous-module de $\mathbb{Z}^{26}$ constitué des $(x_{1},x_{2},\ldots,x_{26})$ avec $\sum_{i=1}^{26}x_{i}$ pair et le vecteur $\frac{1}{2}(1,1,\ldots,1)$. Soit $\mathrm{Y}$ l'ensemble des éléments $e$ de $\mathrm{II}_{25,1}$ avec $e.e=-2$~; Borcherds décrit dans \cite{borcherdsthese} un algorithme pour déterminer l'ensemble $\mathrm{O}(\mathrm{II}_{25,1})\backslash\mathrm{Y}$ qui d'après ce qui précède est en bijection avec $\mathrm{X}_{25}$ ($\vert\mathrm{X}_{25}\vert=121$).
\vfill\eject

\vspace{0,75cm}
\textbf{3. Théorie des $p$-voisins pour les $\mathrm{q\text{-}i}$-modules sur $\mathbb{Z}$}

\bigskip
Voici le pendant (au moins pour l'anneau de Dedekind $\mathbb{Z}$) du point (a) de la proposition III.1.1. 

\bigskip
\textbf{Proposition 3.1.} {\em Soit $V$ un $\mathrm{q}$-espace vectoriel sur $\mathbb{Q}$~; soient $L_{1}$ et $L_{2}$ deux réseaux entiers (au sens quadratique) de $V$ d'indice $2$ dans leurs duaux (les réseaux $L_{1}$ et $L_{2}$ sont donc en particulier deux $\mathrm{q\text{-}i}$-modules et $L_{1}\cap L_{2}$ un $\widetilde{\mathrm{q}}$-module sur~$\mathbb{Z}$).

\medskip
On pose $I_{1}=L_{1}/(L_{1}\cap L_{2})$, $I_{2}=L_{2}/(L_{1}\cap L_{2})$ et $R=(L_{1}^{\sharp}\cap L_{2}^{\sharp})/(L_{1}\cap L_{2})$.

\medskip
{\em (a)} On a $L_{1}\cap L_{2}^{\sharp}=L_{1}\cap L_{2}$ et $L_{2}\cap L_{1}^{\sharp}=L_{1}\cap L_{2}$.

\medskip
{\em (b)} Les deux homomorphismes canoniques $R\to L_{i}^{\sharp}/L_{i}$, $i=1,2$, sont des isomorphismes.

\smallskip
{\em (c)} Les trois inclusions, de $L_{1}$, $L_{2}$ et $L_{1}^{\sharp}\cap L_{2}^{\sharp}$, dans $(L_{1}\cap L_{2})^{\sharp}$ induisent un isomorphisme canonique de groupes abéliens
$$
I_{1}\hspace{2pt}\oplus\hspace{2pt}I_{2}
\hspace{2pt}\oplus\hspace{2pt}R
\hspace{4pt}\cong\hspace{4pt}
\mathop{\mathrm{r\acute{e}s}}(L_{1}\cap L_{2})
$$
(qui permet d'identifier la source et le but).

\medskip
{\em (d)} L'accouplement $I_{1}\times I_{2}\to\mathbb{Q}/\mathbb{Z}$ induit par la forme d'enlacement de $\mathop{\mathrm{r\acute{e}s}}(L_{1}\cap L_{2})$ est non dégénéré. Pour cette forme les deux sous-modules $I_{1}\oplus I_{2}$ et $R$ sont orthogonaux et canoniquement isomorphes, comme $\mathrm{qe}$-modules, respectivement à $\mathrm{H}(I_{1})$ et $\mathop{\mathrm{r\acute{e}s}}L_{1}$, si bien que le $\mathrm{qe}$-module $\mathop{\mathrm{r\acute{e}s}}(L_{1}\cap L_{2})$ est canoniquement isomorphe à la somme orthogonale $\mathrm{H}(I_{1})\oplus\mathop{\mathrm{r\acute{e}s}}L_{1}$.}

\bigskip
\textit{Démonstration.} On vérifie le point (a) et l'isomorphisme de $\mathrm{qe}$-modules\linebreak $R\cong\mathop{\mathrm{r\acute{e}s}}L_{1}$~; la vérification du reste de l'énoncé est laissée au lecteur.

\medskip
-- Soit $\mathcal{L}(V)$ l'ensemble des réseaux entiers de $V$ ordonné par inclusion~; on observe qu'un réseau entier de $V$ qui est d'indice $2$ dans son dual est un élément maximal de $\mathcal{L}(V)$ (en fait, tous les éléments maximaux sont de ce type). Le point (a) résulte de cette observation. On considère le réseau\linebreak $L_{1}+(L_{1}^{\sharp}\cap L_{2})$. Ce réseau appartient à $\mathcal{L}(V)$ et contient $L_{1}$~; on a donc l'égalité $L_{1}=L_{1}+(L_{1}^{\sharp}\cap L_{2})$ qui implique $L_{1}\cap L_{2}^{\sharp}=L_{1}\cap L_{2}$.

\medskip
-- Le sous-module $L_{1}$ de $(L_{1}\cap L_{2})^{\sharp}$ correspond \textit{via} le point (b) de la proposition II.1.1 au sous-module isotrope $I_{1}$ de $\mathop{\mathrm{r\acute{e}s}}(L_{1}\cap L_{2})$. Le point (c) de cette même proposition montre que l'on a $\mathop{\mathrm{r\acute{e}s}}L_{1}\cong I_{1}^{\perp}/I_{1}$. Or  $I_{1}^{\perp}/I_{1}$ s'identifie comme $\mathrm{qe}$-module à $R$.
\hfill$\square$

\bigskip
La vérification de l'énoncé ci-dessous est immédiate.

\bigskip
\textbf{Proposition-Définition 3.2.} {\em Soit $V$ un $\mathrm{q}$-espace vectoriel sur $\mathbb{Q}$~; soient $L_{1}$ et $L_{2}$ deux réseaux entiers (au sens quadratique) de $V$ d'indice $2$ dans leurs duaux.

\medskip
Soit $p$ un nombre premier~; les conditions suivantes sont équivalentes~:
\begin{itemize}
\item [(i)] $L_{1}\cap L_{2}$ est d'indice $p$ dans $L_{1}$~;
\item [(ii)] $L_{1}\cap L_{2}$ est d'indice $p$ dans $L_{2}$.
\end{itemize}

\medskip
Si ces conditions sont vérifiées alors on dit que $L_{1}$ et $L_{2}$ sont {\em $p$-voisins} (ou que $L_{2}$ est un {\em $p$-voisin} de $L_{1}$). Dans ce cas $L_{1}/(L_{1}\cap L_{2})$ et $L_{2}/(L_{1}\cap L_{2})$  sont les seuls sous-modules isotropes non
triviaux de $\mathop{\mathrm{r\acute{e}s}}(L_{1}\cap L_{2})$.}

\bigskip
On fixe maintenant un $\mathrm{q\text{-}i}$-module $L$ sur $\mathbb{Z}$ et on analyse l'ensemble des $p$-voisins de $L$ dans $\mathbb{Q}\otimes_{\mathbb{Z}}L$~; dans ce contexte, un $p$-voisin de $L$ est un réseau entier $L'$ dans $\mathbb{Q}\otimes_{\mathbb{Z}}L$ avec $L'$ d'indice $2$ dans $L'^{\sharp}$ et $L\cap L'$ d'indice $p$ dans~$L$. L'analyse en question se déroule comme au chapitre III~:

\medskip
On pose $M=L\cap L'$. D'après ce qui précède~:

\smallskip
-- Le réseau $pL'$ est contenu dans $M$.

\smallskip
-- L'image de l'homomorphisme composé $pL'\subset M\subset L\to L/pL$ est une droite isotrope, disons $c$, de $L/pL$ muni de sa structure de $\mathrm{q\text{-}i}$-espace vectoriel sur $\mathbb{Z}/p$. 

\smallskip
-- Le réseau $M$ est l'image inverse par  l'homomorphisme $L\to L/pL$ de $c^{\hspace{1pt}\perp}$, $c^{\perp}$ désignant le sous-espace vectoriel de $L/pL$ orthogonal de la droite $c$.

\smallskip
-- Le réseau $L'$ est l'image inverse par l'homomorphisme $M^{\sharp}\to\mathop{\mathrm{r\acute{e}s}}M$ de l'unique sous-module non trivial, isotrope pour la forme quadratique d'enlacement, distinct de $L/M$.

\medskip
Réciproquement~:

\bigskip
\textbf{Proposition 3.3.} {\em Soit $c$ une droite isotrope de $L/pL$~; soit $M$ le sous-module de $L$ image inverse par  l'homomorphisme $L\to L/pL$ de $c^{\hspace{1pt}\perp}$. Alors~:

\smallskip
{\em (a)} Le $\mathrm{qe}$-module $\mathop{\mathrm{r\acute{e}s}}M$ est isomorphe à $\mathrm{H}(\mathbb{Z}/p)\oplus\mathop{\mathrm{r\acute{e}s}}L$.

\smallskip
{\em (b)} L'image inverse par l'homomorphisme $M^{\sharp}\to\mathop{\mathrm{r\acute{e}s}}M$ de l'unique sous-module non trivial, isotrope pour la forme quadratique d'enlacement, distinct de $L/M$, est un $p$-voisin $L'$ de $L$ avec $L\cap L'=M$.}

\bigskip
\textit{Démonstration du point (a).} Elle est très semblable à celle de III.1.4~:

\smallskip
Soit $u$ un élément de $L$ dont la classe modulo $p$ engendre la droite $c$~; puisque cette droite est isotrope on a $\mathrm{q}(u)\equiv 0\bmod{p}$. Soit $b$ la forme bilinéaire symétrique associée à la forme quadratique de $L/pL$. Si $p$ est différent de~$2$ alors $b$ est non dégénérée si bien qu'il existe un un élément $v$ de $L$ avec $u.v\equiv 1\bmod{p}$. C'est encore le cas pour $p=2$. Pour s'en convaincre il suffit de montrer que la classe modulo $2$ de $u$ n'appartient pas à $\ker b$. Comme~$L$ est un $\mathrm{q\text{-}i}$-module sur $\mathbb{Z}$, $L/2L$ est naturellement un $\mathrm{q\text{-}i}$-espace vectoriel sur $\mathbb{Z}/2$~; le point (a) de 2.5 montre que $\ker b$ est de dimension $1$ et que la restriction de la forme quadratique à $\ker b$ est non triviale (en fait, il s'agit là d'un phénomène général pour les  $\mathrm{q\text{-}i}$-espaces vectoriels sur un corps de caractéristique $2$, voir par exemple la remarque qui suit 1.6).

\smallskip
On constate que $v$ et $\frac{u}{p}$ appartiennent à $M^{\sharp}$~; on pose $w=\frac{u}{p}-\frac{\mathrm{q}(u)}{p}\hspace{1pt}v$. On observe que l'on a dans $\mathbb{Q}/\mathbb{Z}$ les égalités $\mathrm{q}(v)=0$, $\mathrm{q}(w)=0$ et $v.w=\frac{1}{p}$. Soit $H$ le sous-module de $\mathop{\mathrm{r\acute{e}s}}M$ engendré par les classes de $v$ et $w$ (ou celles de $v$ et $\frac{u}{p}$)~; l'observation précédente montre que la restriction de la forme d'enlacement à $H$ est non dégénérée, que $H$ est un $\mathbb{Z}/p$-espace vectoriel de dimension $2$ dont une base est constituée des classes de $v$ et $w$, que $H$ est isomorphe comme $\mathrm{qe}$-module  à $\mathrm{H}(\mathbb{Z}/p)$ et enfin que le $\mathrm{qe}$-module $\mathop{\mathrm{r\acute{e}s}}M$ est isomorphe à la somme orthogonale $H\oplus H^{\perp}$. Soit $I$ le sous-module $L/M$ de $\mathop{\mathrm{r\acute{e}s}}M$ (le sous-module engendré par $v$)~; $I$ est isotrope et le $\mathrm{qe}$-module $I^{\perp}/I$ s'identifie à $H^{\perp}$. Le point (c) de la proposition II.1.1 montre donc que $H^{\perp}$ s'identifie à~$\mathop{\mathrm{r\acute{e}s}}L$.
\hfill$\square$

\medskip
On constate au bout du compte que l'énoncé III.1.5 concernant les $\mathrm{q}$-modules sur $\mathbb{Z}$ reste valable, mot pour mot,  pour les $\mathrm{q\text{-}i}$-modules sur $\mathbb{Z}$. Précisons. Soit $L$ un $\mathrm{q\text{-}i}$-module sur $\mathbb{Z}$. On note $\mathrm{C}_{L}$ la quadrique projective~;  $\mathrm{C}_{L}$ est encore lisse sur $\mathbb{Z}$ (c'est d'ailleurs le critère élégant que choisit Deligne dans \cite[Expos\'e XII]{sga7} pour caractériser les formes quadratiques qu'il qualifie de non dégénérées). Soit $c$ un point de $\mathrm{C}_{L}(\mathbb{Z}/p)$~; on note $\mathrm{vois}_{d}(L;c)$ le réseau $L'$ de $\mathbb{Q}\otimes_{\mathbb{Z}}L$ associé à $c$ par le procédé décrit dans le point (b) de 3.3 et $\mathrm{Vois}_{d}(L)$ l'ensemble des $p$-voisins de $L$ dans $\mathbb{Q}\otimes_{\mathbb{Z}}L$.

\bigskip
\textbf{Proposition 3.4.} {\em L'application
$$
\mathrm{C}_{L}(\mathbb{Z}/d)\to\mathrm{Vois}_{d}(L)
\hspace{24pt},\hspace{24pt}
c\mapsto\mathrm{vois}_{d}(L;c)
$$
est une bijection.}

\vspace{0,75cm}
\textbf{4. Théorie des $p$-voisins pour les réseaux pairs de déterminant~$2$}

\bigskip
Soit $\mathrm{E}$ l'un des deux réseaux pairs de déterminant $2$, $\mathrm{A}_{1}$ ou $\mathrm{E}_{7}$~; dans la suite de ce paragraphe, sauf mention expresse du contraire, $\mathrm{E}$ est fixé.

\medskip
Soit $n>0$ un entier avec $n\equiv-\dim\mathrm{E}\pmod{8}$.

\medskip
On note $\mathrm{U}_{n}$ le $\mathbb{Q}$-espace vectoriel $\mathbb{Q}^{n}$ muni de la forme quadratique $$
\hspace{24pt}
(x_{1},x_{2},\ldots,x_{n})
\hspace{4pt}\mapsto\hspace{4pt}
\frac{1}{2}\hspace{2pt}\sum_{i=1}^{n-1}x_{i}^{2}
\hspace{2pt}+\hspace{2pt}
x_{n}^{2}
\hspace{24pt}.
$$
On pose $\mathrm{V}_{n}=\mathrm{U}_{n}\oplus(\mathbb{Q}\otimes_{\mathbb{Z}}\mathrm{E})$ ($\mathrm{V}_{n}$ est isomorphe, comme $\mathrm{q}$-espace vectoriel, à $\mathbb{Q}^{n+\dim\mathrm{E}}$  muni de la forme quadratique $(x_{1},x_{2},\ldots,x_{n+\dim\mathrm{E}})\mapsto\frac{1}{2}\sum_{i}x_{i}^{2}$).

\medskip
Soit $\mathcal{X}(\mathrm{U}_{n})$) (resp. $\mathcal{X}(\mathrm{V}_{n})$) l'ensemble des réseaux pairs de déterminant $2$ de $\mathrm{U}_{n}$ (resp. des réseaux unimodulaires pairs de $\mathrm{V}_{n}$). On a donc $\mathrm{X}_{n}=\mathrm{O}(\mathrm{U}_{n})\backslash\mathcal{X}(\mathrm{U}_{n})=\mathrm{SO}(\mathrm{U}_{n})\backslash\mathcal{X}(\mathrm{U}_{n})$ (resp. $\mathrm{X}_{n+\dim E}=\mathrm{O}(\mathrm{V}_{n})\backslash\mathcal{X}(\mathrm{V}_{n})$). Soient $\mathcal{X}(\mathrm{V}_{n};\mathrm{E})$ le sous-ensemble de $\mathcal{X}(\mathrm{V}_{n})$ constitué des $P$ contenant $\mathrm{E}$ et $\mathrm{O}(\mathrm{V}_{n};\mathrm{E})$ le sous-groupe de $\mathrm{O}(\mathrm{V}_{n})$ constitué des éléments qui sont l'identité sur $\mathrm{E}$. Avec ces notations, on peut paraphraser les propositions 2.6 et 2.7 de la façon suivante~:

\bigskip
\textbf{Proposition 4.1.} {\em L'application
$$
\mathcal{X}(\mathrm{V}_{n};\mathrm{E})\to\mathcal{X}(\mathrm{U}_{n})
\hspace{24pt},\hspace{24pt}
P\mapsto P\cap\mathrm{U}_{n}
$$
est une bijection équivariante relativement à l'isomorphisme de groupes
$$
\hspace{24pt}
\mathrm{O}(\mathrm{V}_{n};\mathrm{E})\overset{\cong}{\to}\mathrm{O}(\mathrm{U}_{n})
\hspace{24pt}.
$$}

\medskip
La spécialisation au cas défini positif du paragraphe 3 conduit quant à elle à l'énoncé suivant~:

\bigskip
\textbf{Proposition 4.2.} {\em Soient $P_{1}$ et $P_{2}$ deux réseaux unimodulaires pairs de $\mathrm{V}_{n}$ contenant $\mathrm{E}$~;  soient respectivement, $L_{1}$ et $L_{2}$ les deux résaux pairs de déter\-minant $2$ de $\mathrm{U}_{n}$, $P_{1}\cap\mathrm{U}_{n}$ et $P_{2}\cap\mathrm{U}_{n}$. Soit $p$ un nombre premier, les deux conditions suivantes sont équivalentes~:
\begin{itemize}
\item [(i)] $P_{1}$ et $P_{2}$ sont $p$-voisins~;
\item [(ii)] $L_{1}$ et $L_{2}$ sont $p$-voisins.
\end{itemize}}

\bigskip
\textit{Démonstration de (i)$\Rightarrow$(ii).} Par définition même, l'homomorphisme canonique $L_{i}/(L_{1}\cap L_{2})\to P_{i}/(P_{1}\cap P_{2})$, $i=1,2$, est injectif. D'après 3.2 on a l'alternative suivante~: $L_{1}\cap L_{2}$ d'indice $p$ dans $L_{1}$ et $L_{2}$ ou $L_{1}=L_{2}$. Or l'application $P\mapsto P\cap\mathrm{U}_{n}$ est injective.
\hfill$\square$

\bigskip
\textit{Démonstration de (ii)$\Rightarrow$(i).} On considère le réseau $N$ de $\mathrm{V}_{n}$ suivant~:
$$
\hspace{24pt}
N
\hspace{4pt}:=\hspace{4pt}
(L_{1}\oplus\mathrm{E})\cap(L_{2}\oplus\mathrm{E})
\hspace{4pt}=\hspace{4pt}
(L_{1}\cap L_{2})\oplus\mathrm{E}
\hspace{24pt}.
$$
D'après le point (d) de 3.1 et le point (d) de 2.2, le $\mathrm{qe}$-module $\mathop{\mathrm{r\acute{e}s}}N$ est canoniquement  isomorphe à la somme orthogonale $$
\mathrm{H}(L_{1}/(L_{1}\cap L_{2}))
\hspace{4pt}\oplus\hspace{4pt}
\langle-1\rangle\otimes\mathop{\mathrm{r\acute{e}s}}\mathrm{E}
\hspace{4pt}\oplus\hspace{4pt}
\mathop{\mathrm{r\acute{e}s}}\mathrm{E}
$$
et l'on vérifie que $P_{1}\cap P_{2}$ est le sous-module de $N^{\sharp}$ correspondant \textit{via} le point (b) de la proposition II.1.1 à la ``diagonale'' du facteur $\langle-1\rangle\otimes\mathop{\mathrm{r\acute{e}s}}\mathrm{E}\oplus\mathop{\mathrm{r\acute{e}s}}\mathrm{E}$.\linebreak Le point (c) de cette même proposition montre que l'on a $\mathop{\mathrm{r\acute{e}s}}(P_{1}\cap P_{2})\cong\mathrm{H}(L_{1}/(L_{1}\cap L_{2}))$. On conclut à l'aide du point (a) de III.1.1.
\hfill$\square$

\bigskip
Soit comme ci-dessus $n>0$ un entier avec $n\equiv-\dim\mathrm{E}\pmod{8}$~; soit $p$ un nombre premier. Les opérateurs de Hecke $\mathop{\mathrm{T}_{p}}:\mathbb{Z}[\mathrm{X}_{n}]\to\mathbb{Z}[\mathrm{X}_{n}]$ sont définis comme dans le cas $n\equiv 0\pmod{8}$~:
$$
\hspace{24pt}
\mathop{\mathrm{T}_{p}}\hspace{1pt}[L]
\hspace{4pt}:=\hspace{4pt}
\sum_{L'\in\mathrm{Vois}_{p}(L)}[L']
\hspace{24pt},
$$
pour tout réseau pair $L$ de déterminant $2$ et de dimension $n$. Soient $L$ et $L'$  deux réseaux pairs de déterminant $2$ et de dimension $n$, l'entier $\mathrm{N}_{p}(L,L')$ est défini pareillement comme la $[L']$-coordonnée de $\mathop{\mathrm{T}_{p}}\hspace{1pt}[L]$.

\bigskip
Le point (a) de 2.5 montre que toutes les quadriques $\mathrm{C}_{L}(\mathbb{Z}/p)$ ont le même cardinal, disons encore $\mathrm{c}_{n}(p)$~; on a cette fois
$$
\hspace{24pt}
\mathrm{c}_{n}(p)=\sum_{m=0}^{n-2}\hspace{4pt}p^{m}
\hspace{24pt}.
$$
L'énoncé III.2.2 reste valable, mot pour mot~:

\bigskip
\textbf{Proposition 4.3.} {\em Soit  $p$ un nombre premier. On a
$$
\sum_{y\in\mathrm{X}_{n}}
\mathrm{N}_{p}(x,y)
\hspace{4pt}=\hspace{4pt}
\mathrm{c}_{n}(p)
$$
pour tout $x$ dans $\mathrm{X}_{n}$.}

\bigskip
On dégage ci-après l'analogue de  la proposition III.1.10.

\smallskip
-- On note $\mathrm{B}_{n}(p)$ l'ensemble des classes d'isomorphisme des $\widetilde{\mathrm{q}}$-modules $M$ sur~$\mathbb{Z}$, avec $\dim M=n$, $\mathbb{R}\otimes_{\mathbb{Z}}M>0$ et $\mathop{\mathrm{r\acute{e}s}}M\simeq\mathrm{H}(\mathbb{Z}/p)\oplus\langle-1\rangle\otimes\mathop{\mathrm{r\acute{e}s}}\mathrm{E}$~; on observera que tout automorphisme du $\mathrm{qe}$-module $\mathrm{H}(\mathbb{Z}/p)\oplus\langle-1\rangle\otimes\mathop{\mathrm{r\acute{e}s}}\mathrm{E}$ est l'identité sur le facteur $\langle-1\rangle\otimes\mathop{\mathrm{r\acute{e}s}}\mathrm{E}$ si bien que le $\mathrm{qe}$-module est muni canoniquement d'un facteur direct  $\langle-1\rangle\otimes\mathop{\mathrm{r\acute{e}s}}\mathrm{E}$ .

\smallskip
-- On note $\widetilde{\mathrm{B}}_{n}(p)$ l'ensemble des classes d'isomorphisme de couples  $(M;\omega)$ avec $M$ comme précédemment et $\omega$ une bijection de l'ensemble des sous-modules isotropes non-trivaux de $\mathop{\mathrm{r\acute{e}s}}M$ sur l'ensemble $\{1,2\}$. Par définition $\widetilde{\mathrm{B}}_{n}(p)$ est muni d'une action à gauche du groupe symétrique $\mathfrak{S}_{2}$ et l'ensemble quotient $\mathfrak{S}_{2}\backslash\widetilde{\mathrm{B}}_{n}(p)$ s'identifie à $\mathrm{B}_{n}(p)$.

\smallskip
-- Soit $(M;\omega)$ comme ci-dessus, on note $\mathrm{d}_{i}(M;\omega)$, $i=1,2$, l'image réciproque par la surjection $M^{\sharp}\to\mathop{\mathrm{r\acute{e}s}}M$ de $\omega^{-1}(i)$~; $\mathrm{d}_{1}(M;\omega)$ et $\mathrm{d}_{2}(M;\omega)$ sont deux réseaux pairs de déterminant $2$ (de dimension $n$), $p$-voisins dans $\mathbb{Q}\otimes_{\mathbb{Z}}M$. En passant aux classes d'isomorphisme on obtient deux applications que l'on note encore $\mathrm{d}_{1}$ et $\mathrm{d}_{2}$ de $\widetilde{\mathrm{B}}_{n}(p)$ dans $\mathrm{X}_{n}$.

\smallskip
Nous pouvons énoncer~:

\bigskip
\textbf{Proposition 4.4.} {\em Soit $p$ un nombre premier~; soient $x_{1}$ et $x_{2}$ deux éléments de $\mathrm{X}_{n}$. On a~:
$$
\mathrm{N}_{p}(x_{1},x_{2})
\hspace{4pt}=\hspace{4pt}
\sum_{\beta\hspace{2pt}\in\hspace{2pt}\mathrm{d}_{1}^{-1}(x_{1})\hspace{1pt}\cap\hspace{1pt}\mathrm{d}_{2}^{-1}(x_{2})}\hspace{8pt}
\frac{\vert\mathrm{O}(x_{1})\vert}{\vert\mathrm{O}(\beta)\vert}
$$
(avec $\vert\mathrm{O}(\beta)\vert=\vert\mathrm{Aut}(M;\omega)\vert$, $(M;\omega)$ représentant $\beta$).}

\bigskip
On reformule maintenant, au moins dans un cas particulier,  l'énoncé ci-dessus en faisant intervenir les réseaux unimodulaires pairs associés aux réseaux pairs de déterminant $2$.

\smallskip
Soit $L$ un réseau pair de déterminant $2$, de dimension $n$. On pose comme précédemment $P=(L\oplus\mathrm{E})^{+}$. Par construction, $P$ est canoniquement muni d'une représentation $\mathrm{i}:\mathrm{E}\to P$. On note $\rho(L)$ (resp. $\bar{\rho}(L)$) le cardinal de l'orbite de $\mathrm{i}$ (resp. $\mathrm{i}(L)$) sous l'action du groupe $\mathrm{O}(P)$. Comme nous l'avons déjà dit (au moins dans le cas $\mathrm{E}=\mathrm{E}_{7}$) on a $\rho(L)=\bar{\rho}(L)\hspace{1pt}\vert\mathrm{O}(\mathrm{E})\vert$. On dispose donc de deux fonctions, $\rho$ et $\bar{\rho}$, de $\mathrm{X}_{n}$ dans $\mathrm{N}-\{0\}$ (en fait il sera plus commode de travailler avec $\rho$ dans le cas $\mathrm{E}=\mathrm{A}_{1}$ et avec $\bar{\rho}$ dans le cas $\mathrm{E}=\mathrm{E}_{7}$).

\smallskip
Soit $M$ un $\widetilde{\mathrm{q}}$-module $M$ sur $\mathbb{Z}$, avec $\dim M=n$, $\mathbb{R}\otimes_{\mathbb{Z}}M>0$ et $\mathop{\mathrm{r\acute{e}s}}M\simeq\mathrm{H}(\mathbb{Z}/p)\oplus\langle-1\rangle\otimes\mathop{\mathrm{r\acute{e}s}}\mathrm{E}$. On pose $R=(M\oplus\mathrm{E})^{+}$. Décodons. Il résulte de l'observation faite plus haut que le $\mathrm{qe}$-module $\mathop{\mathrm{r\acute{e}s}}(M\oplus\mathrm{E})$ est muni canoniquement d'un facteur direct $\langle-1\rangle\otimes\mathop{\mathrm{r\acute{e}s}}\mathrm{E}\oplus\mathop{\mathrm{r\acute{e}s}}\mathrm{E}$~; $R$ est le réseau correspondant \textit{via} le point (b) de la proposition II.1.1 à la ``diagonale'' de facteur direct et le point (c) de cette même proposition montre que l'on a $\mathop{\mathrm{r\acute{e}s}}R\simeq\mathrm{H}(\mathbb{Z}/p)$. A nouveau $R$ est canoniquement muni d'une représentation $\mathrm{i}:\mathrm{E}\to R$ et l'on note encore $\rho(M)$ (resp. $\bar{\rho}(M)$) le cardinal de l'orbite de $\mathrm{i}$ (resp. $\mathrm{i}(M)$) sous l'action du groupe $\mathrm{O}(R)$. On dispose donc aussi de deux fonctions, $\rho$ et $\bar{\rho}$, de $\mathrm{B}_{n}(p)$ dans $\mathrm{N}-\{0\}$ avec $\rho(-)=\bar{\rho}(-)\hspace{1pt}\vert\mathrm{O}(\mathrm{E})\vert$.

\smallskip
Ce qui précède montre que l'on dispose d'applications naturelles de $\mathrm{X}_{n}$ (resp. $\mathrm{B}_{n}(p)$) dans $\mathrm{X}_{n+\dim\mathrm{E}}$ (resp. $\mathrm{B}_{n+\dim\mathrm{E}}\hspace{1pt}(p)$)~; ces applications sont notées $\pi$.

\medskip
Soit enfin $\mathrm{B}_{n}^{0}(p)$ le sous-ensemble de $\mathrm{B}_{n}(p)$ constitué des classes d'isomorphisme des $M$ comme ci-dessus, tels que les deux réseaux pairs de déter\-minant~$2$, $L_{1}$ et $L_{2}$, associés à $M$, sont non isomorphes. On note $\mathrm{e}$ l'application de $\mathrm{B}_{n}^{0}(p)$ dans l'ensemble des paires  d'éléments de $\mathrm{X}_{n}$, qui envoie $[M]$ sur la paire $\{[L_{1}],[L_{2}]\}$.

\bigskip
\textbf{Proposition 4.5.} {\em Soient $x_{1}$ et $x_{2}$ deux éléments distincts de $\mathrm{X}_{n}$. On a~:
$$
\hspace{4pt}
\mathrm{N}_{p}(x_{1},x_{2})
=\hspace{-4pt}
\sum_{b\hspace{2pt}\in\hspace{2pt}\mathrm{e}^{-1}(\{x_{1},x_{2}\})}
\frac{\vert\mathrm{O}(\pi(x_{1}))\vert}{\vert\mathrm{O}(\pi(b))\vert}\frac{\rho(b)}{\rho(x_{1})} 
=\hspace{-4pt}
\sum_{b\hspace{2pt}\in\hspace{2pt}\mathrm{e}^{-1}(\{x_{1},x_{2}\})}
\frac{\vert\mathrm{O}(\pi(x_{1}))\vert}{\vert\mathrm{O}(\pi(b))\vert}\frac{\bar{\rho}(b)}{\bar{\rho}(x_{1})}
\hspace{4pt}.
$$}

\bigskip
\textit{Démonstration.} On vérifie que l'hypothèse $x_{1}\not=x_{2}$ permet de transformer l'énoncé 4.4 en l'énoncé 4.5. On représente l'élément $\beta$ de $\widetilde{\mathrm{B}}_{n}(p)$ qui apparaît en 4.4 par un couple $(M;\omega)$~; on pose $L_{1}=\mathrm{d}_{1}(M;\omega)$, $P_{1}=(L_{1}\oplus\mathrm{E})^{+}$ et $R=(M\oplus\mathrm{E})^{+}$. On a~:
$$
\hspace{2pt}
\frac{\vert\mathrm{O}(L_{1})\vert}{\vert\mathrm{O}(M;\omega)\vert}
=
[\mathrm{O}(M):\mathrm{O}(M;\omega)]\hspace{2pt}\frac{\vert\mathrm{O}(L_{1})\vert}{\vert\mathrm{O}(M)\vert}
=
[\mathrm{O}(M):\mathrm{O}(M;\omega)]\hspace{2pt}\frac{\vert\mathrm{O}(P_{1})\vert}{\vert\mathrm{O}(R)\vert}
\hspace{2pt}
\frac{\rho(M)}{\rho(L_{1})}
\hspace{2pt}.
$$
Or l'hypothèse $x_{1}\not=x_{2}$ implique $\mathrm{O}(M;\omega)=\mathrm{O}(M)$~; en effet, dans le cas contraire il existe un élément de $\mathrm{O}(M)$ qui échange les deux sous-modules istropes non triviaux de $\mathop{\mathrm{r\acute{e}s}}\mathrm{M}$ et l'on a $x_{1}=x_{2}$.
\hfill$\square$

\vspace{0,75cm}
\textsc{$2$-voisins, le point de vue de Borcherds (suite)}

\medskip
Les propositions 2.6 et 2.7 disent que l'application $\pi:\mathrm{X}_{n}\to\mathrm{X}_{n+\dim\mathrm{E}}$ considérée ci-dessus est la composée d'une bijection $\mathrm{X}_{n}\cong\mathrm{X}_{n+\dim\mathrm{E}}^{\mathrm{E}}$ et d'une ``application d'oubli'' $\mathrm{X}_{n+\dim\mathrm{E}}^{\mathrm{E}}\to\mathrm{X}_{n+\dim\mathrm{E}}$. Pareillement, l'application $\pi:\mathrm{B}_{n}(p)\to\mathrm{B}_{n+\dim\mathrm{E}}\hspace{1pt}(p)$ est la composée d'une bijection $\mathrm{B}_{n}(p)\cong\mathrm{B}_{n+\dim\mathrm{E}}^{\mathrm{E}}\hspace{1pt}(p)$ et d'une application d'oubli $\mathrm{B}_{n+\dim\mathrm{E}}^{\mathrm{E}}\hspace{1pt}(p)\to\mathrm{B}_{n+\dim\mathrm{E}}\hspace{1pt}(p)$~; la définition de l'ensemble $\mathrm{B}_{n+\dim\mathrm{E}}^{\mathrm{E}}\hspace{1pt}(p)$ est sans surprise~: $\mathrm{B}_{n+\dim\mathrm{E}}^{\mathrm{E}}\hspace{1pt}(p)$ est l'ensemble des classes d'isomorphisme des couples $(R;f)$, $R$ et $f$ désignant respectivement un réseau pair avec  $\mathop{\mathrm{r\acute{e}s}}R\simeq\mathrm{H}(\mathbb{Z}/p)$ et une représentation de $\mathrm{E}$ dans $R$.

\medskip
On a vu au chapitre III (en suivant Borcherds) que l'ensemble $\mathrm{B}_{n+\dim\mathrm{E}}\hspace{1pt}(2)$ s'identifie à l'ensemble, que nous avons noté $\mathrm{B}_{n+\dim\mathrm{E}}$, des classes d'isomorphisme des réseaux unimodulaires impairs de dimension $n+\dim\mathrm{E}$. Pareillement, l'ensemble $\mathrm{B}_{n+\dim\mathrm{E}}^{\mathrm{E}}\hspace{1pt}(2)$ s'identifie à l'ensemble des classes d'isomorphisme des réseaux unimodulaires impairs de dimension $n+\dim\mathrm{E}$ munis d'une représentation de $\mathrm{E}$, ensemble que nous noterons $\mathrm{B}_{n+\dim\mathrm{E}}^{\mathrm{E}}$. On dispose aussi d'une identification à la Borcherds $\mathrm{B}_{n}(2)\cong\mathrm{B}_{n}$, $\mathrm{B}_{n}$ désignant l'ensemble des réseaux impairs de déterminant $2$~: à un réseau pair $M$ avec $\mathop{\mathrm{r\acute{e}s}}M\simeq\mathrm{H}(\mathbb{Z}/2)\oplus\langle-1\rangle\otimes\mathop{\mathrm{r\acute{e}s}}\mathrm{E}$ on associe le réseau correspondant \textit{via} II.1.1 (version bilinéaire) à l'unique sous-module non trivial de $\mathop{\mathrm{r\acute{e}s}}M$ qui est isotrope au sens bilinéaire mais non quadratique.

\medskip
Soit $(Q;f)$ un réseau unimodulaire impair de dimension $n+\dim\mathrm{E}$ muni d'une représententation $f:\mathrm{E}\to Q$. Le lecteur vérifiera que la correspondance $Q\mapsto M$ peut être décrite des deux façons suivantes~:

\smallskip
-- On considère le sous-module $R$ d'indice $2$ de $Q$ constitué des $x$ avec $x.x$ pair. On observe que l'on a $f(\mathrm{E})\subset R$ ; $M$ est l'orthogonal de $f(\mathrm{E})$ dans $R$.

\smallskip
-- On considère l'orthogonal de $f(\mathrm{E})$ dans $Q$, disons $\Lambda$. On observe que $\Lambda$ est impair (si ce n'était le cas, on aurait $R=\Lambda\oplus f(\mathrm{E})$, égalité interdite parce que le résidu bilinéaire du second membre n'est pas hyperbolique)~; $M$ est le sous-module de $\Lambda$ constitué des $x$ avec $x.x$ pair.

\medskip
Enfin l'application $\rho:\mathrm{B}_{n+\dim\mathrm{E}}^{\mathrm{E}}\to\mathbb{N}-\{0\}$ (resp $\bar{\rho}:\mathrm{B}_{n+\dim\mathrm{E}}^{\mathrm{E}}\to\mathbb{N}-\{0\}$) induite par l'identification $\mathrm{B}_{n+\dim\mathrm{E}}^{\mathrm{E}}\hspace{1pt}(2)\cong\mathrm{B}_{n+\dim\mathrm{E}}^{\mathrm{E}}$ associe à la classe du couple $(Q;f)$ le cardinal de l'orbite de $\mathrm{f}$ (resp. $\mathrm{f}(L)$) sous l'action du groupe $\mathrm{O}(Q)$.

\bigskip
On a observé au chapitre III, avec Nebe et Venkov, que si un réseau unimodulaire impair  $L$, de dimension divisible par $8$, représente $1$, alors les deux réseaux unimodulaires pairs $2$-voisins de $L$ sont isomorphes  (Corollaire III.1.16). Nous terminons ce paragraphe en dégageant l'énoncé technique ci-dessous que l'on peut voir comme le pendant, dans le contexte présent, de cette observation.

\bigskip
\textbf{Proposition 4.6.} {\em Soit $n>0$ un entier avec $n\equiv-\dim\mathrm{E}\pmod{8}$. Soient $L_{1}$ et $L_{2}$ deux réseaux pairs de déterminant $2$ dans un $\mathrm{q}$-espace vectoriel $U$ de dimension $n$~; on pose $V=U\oplus(\mathbb{Q}\otimes_{\mathbb{Z}}\mathrm{E})$. On suppose que $L_{1}$ et $L_{2}$ sont $2$-voisins dans $U$ et non isomorphes. On note $P_{i}$, $i=1,2$, le réseau unimodulaire pair dans $V$, contenant $\mathrm{E}$, associé à $L_{i}$ ; on note $Q$ le réseau unimodulaire impair dans $V$, contenant $\mathrm{E}$, dont les $2$ voisins pairs sont $P_{1}$ et~$P_{2}$.

\medskip
{\em (a)} Dans le cas $\mathrm{E}=\mathrm{E}_{7}$ le réseau $Q$ ne représente pas $1$.

\medskip
{\em (b)} Dans le cas $\mathrm{E}=\mathrm{A}_{1}$ on a l'alternative suivante~:

\smallskip
{\em (b.1)} Le réseau $Q$ ne représente pas $1$.

\smallskip
{\em (b.2)} Les réseaux $P_{1}$ et $P_{2}$ sont isomorphes et il existe un isomorphisme de réseaux $\phi:Q\cong\mathrm{I}_{2}\oplus Q'$ avec $\mathrm{I}_{2}$ contenant $\phi(\mathrm{A}_{1})$ et $Q'$ ne représentant pas $1$.}

\bigskip
\textit{Démonstration.}  On note $Q^{1}$ le sous-module de $Q$ engendré par les éléments $x$ avec $x.x=1$ et $Q'$ son orthogonal~; on a donc une décomposition en somme orthogonale $Q=Q^{1}\oplus Q'$, le réseau $Q^{1}$ est isomorphe à $\mathrm{I}_{m}$ (avec  $m=\dim_{\mathbb{Z}}Q^{1}$) et l'ensemble $\mathrm{R}(Q)$ des racines de $Q$ est isomorphe à la réunion disjointe $\mathrm{R}(\mathrm{I}_{m})\coprod\mathrm{R}(Q')=\mathrm{R}(\mathrm{D}_{m})\coprod\mathrm{R}(Q')$.

\medskip
Cas $\mathrm{E}=\mathrm{E}_{7}$. D'après ce qui précède $\mathrm{E}_{7}$ est contenu dans $Q'$ ce qui force $m=0$. En effet, s'il existe $e$ dans $Q$ avec $e.e=1$ alors la symétrie orthogonale $\mathrm{s}_{e}$ est l'identité sur $\mathrm{E}_{7}$ et échange $P_{1}$ et $P_{2}$ (voir III.1.16), si bien que $L_{1}$ et $L_{2}$ sont isomorphes.
\hfill$\square$

\medskip
Cas $\mathrm{E}=\mathrm{A}_{1}$ et $m\not=0$. Les réseaux $P_{1}$ et $P_{2}$ sont isomorphes d'après III.1.16. Soit $\alpha$ une racine de $\mathrm{A}_{1}$, disons ``la'' racine positive. On a nécessairement $\alpha\in\mathrm{Q}^{1}$ ou $\alpha\in\mathrm{Q}'$. Le même argument que précédemment interdit $\alpha\in\mathrm{Q}'$. On a donc $\alpha\in\mathrm{Q}^{1}$ ce qui entraîne $m\geq 2$. On a d'autre part $m<3$. En effet, si l'on a $m\geq 3$, alors il existe un élément $e$ dans $Q^{1}$ avec $e.e=1$ et $e.\alpha=0$ et l'on a encore $\mathrm{s}_{e}(\alpha)=\alpha$ et $\mathrm{s}_{e}(P_{1})=P_{2}$.
\hfill$\square$

\vspace{0,75cm}
\textbf{5. Exemples}

\bigskip
\textbf{5.1.} Détermination de $\mathrm{T}_{2}$ pour $n=17$

\bigskip
La matrice de l'opérateur de Hecke $\mathrm{T}_{2}:\mathbb{Z}[\mathrm{X}_{17}]\to\mathbb{Z}[\mathrm{X}_{17}]$, dans la base $(\mathrm{E}_{16}\oplus\mathrm{A}_{1},\mathrm{E}_{8}\oplus\mathrm{E}_{8}\oplus\mathrm{A}_{1},\mathrm{A}_{17}^{+},(\mathrm{D}_{10}\oplus\mathrm{E}_{7})^{+})$ est la suivante (on la note encore~$\mathrm{T}_{2}$)~:
$$
\hspace{24pt}
\mathrm{T}_{2}
\hspace{4pt}=\hspace{4pt}
\begin{bmatrix}
20265 & 18225 & 153 & 63 \\ 12870 & 14910 & 0 & 90 \\ 16384 & 0 & 21624 & 18432 \\ 16016 & 32400 & 43758 & 46950
\end{bmatrix}
\hspace{24pt}.
$$
On explique ci-dessous comment la théorie du paragraphe 4 conduit à cette égalité. On reprend les notations introduites à la fin du paragraphe 2.

\bigskip
Soient $i$ et $j$ deux éléments de l'ensemble $\{2,3,6,7\}$ avec $i\not=j$. La table de Borcherds \cite[Chap. 17]{conwaysloane} fournit la liste des classes d'isomorphisme des réseaux unimodulaires impairs de dimension $24$ dont les $2$-voisins unimodulaires pairs sont isomorphes à $\mathrm{P}_{i}$ et $\mathrm{P}_{j}$~:

\smallskip
-- Cette liste est vide pour $\{i,j\}=\{2,3\}$.

\smallskip
-- Cette liste contient un seul élément, disons $[\mathrm{Bor}_{i,j}]$, pour $\{i,j\}\not=\{2,3\}$.

\smallskip
La table de Borcherds montre aussi que dans le cas $\{i,j\}\not=\{2,3\}$ le groupe orthogonal $\mathrm{O}(\mathrm{Bor}_{i,j})$ agit transitivement sur l'ensemble $\overline{\mathrm{Rep}}(\mathrm{E}_{7},\mathrm{Bor}_{i,j})$.

\smallskip
On rappelle que le groupe $\mathrm{O}(\mathrm{P}_{i})$ agit transitivement sur l'ensemble $\overline{\mathrm{Rep}}(\mathrm{E}_{7},\mathrm{P}_{i})$ pour $i=2,3,6,7$.

\smallskip
La proposition III.3.3.1 donne
$$
\hspace{24pt}
\mathrm{N}_{2}(\mathrm{P}_{i},\mathrm{P}_{j})
\hspace{4pt}=\hspace{4pt}
0
\hspace{6pt}\text{pour}\hspace{6pt}\{i,j\}=\{2,3\}
\hspace{24pt},
$$
$$
\hspace{24pt}
\mathrm{N}_{2}(\mathrm{P}_{i},\mathrm{P}_{j})
\hspace{4pt}=\hspace{4pt}
\frac{\vert\mathrm{O}(\mathrm{P}_{i})\vert}{\vert\mathrm{O}(\mathrm{Bor}_{i,j})\vert}
\hspace{6pt}\text{pour}\hspace{6pt}\{i,j\}\not=\{2,3\}
\hspace{24pt}.
$$
La proposition 4.5 et la discussion qui suit cette proposition donnent donc
$$
\hspace{24pt}
\mathrm{N}_{2}(\mathrm{L}_{i},\mathrm{L}_{j})
\hspace{4pt}=\hspace{4pt}
0
\hspace{6pt}\text{pour}\hspace{6pt}\{i,j\}=\{2,3\}
\hspace{24pt},
$$
$$
\hspace{24pt}
\mathrm{N}_{2}(\mathrm{L}_{i},\mathrm{L}_{j})
\hspace{4pt}=\hspace{4pt}
\frac{\vert\overline{\mathrm{Rep}}(\mathrm{E}_{7},\mathrm{Bor}_{i,j})\vert}{\vert\overline{\mathrm{Rep}}(\mathrm{E}_{7},\mathrm{P}_{i})\vert}
\hspace{4pt}
\mathrm{N}_{2}(\mathrm{P}_{i},\mathrm{P}_{j})
\hspace{6pt}\text{pour}\hspace{6pt}\{i,j\}\not=\{2,3\}
\hspace{24pt}.
$$
Considérons par exemple le cas $\{i,j\}=\{7,2\}$. La table de Borcherds mentionnée ci-dessus montre que l'on a dans ce cas  $\vert\overline{\mathrm{Rep}}(\mathrm{E}_{7},\mathrm{Bor}_{7,2})\vert=1$~; en effet, Borcherds nous dit que l'ensemble des racines de $\mathrm{Bor}_{7,2}$ (qui a le numéro 150 dans sa table) est isomorphe à $\mathbf{D}_{10}\coprod\mathbf{E}_{7}\coprod\mathbf{D}_{6}\coprod\mathbf{A}_{1}$. On a d'autre part $\vert\overline{\mathrm{Rep}}(\mathrm{E}_{7},\mathrm{P}_{7})\vert=2$ puisque l'ensemble des racines de $\mathrm{P}_{7}$ est isomorphe à $\mathbf{D}_{10}\coprod\mathbf{E}_{7}\coprod\mathbf{E}_{7}$. On en déduit
$$
\hspace{24pt}
\mathrm{N}_{2}(\mathrm{L}_{7},\mathrm{L}_{2})
\hspace{4pt}=\hspace{4pt}
\frac{1}{2}
\hspace{4pt}\mathrm{N}_{2}(\mathrm{P}_{7},\mathrm{P}_{2})
\hspace{4pt}=\hspace{4pt}
63
\hspace{24pt}.
$$

\bigskip
\footnotesize
\textit{Remarque.} Le lecteur observera que l'on a $\mathrm{N}_{2}(\mathrm{E}_{16}\oplus\mathrm{A}_{1},\mathrm{E}_{8}\oplus\mathrm{E}_{8}\oplus\mathrm{A}_{1})=\mathrm{N}_{2}(\mathrm{E}_{16},\mathrm{E}_{8}\oplus\mathrm{E}_{8})$ et $\mathrm{N}_{2}(\mathrm{E}_{8}\oplus\mathrm{E}_{8}\oplus\mathrm{A}_{1},\mathrm{E}_{16}\oplus\mathrm{A}_{1})=\mathrm{N}_{2}(\mathrm{E}_{8}\oplus\mathrm{E}_{8},\mathrm{E}_{16})$\ldots et n'aura pas de mal à trouver une explication à ce phénomène.
\normalsize

\bigskip
Comme dans le cas $n=24$, on constate que les valeurs propres de l'opérateur de Hecke $\mathrm{T}_{2}:\mathbb{Z}[\mathrm{X}_{17}]\to\mathbb{Z}[\mathrm{X}_{17}]$ sont entières et simples. Du coup, si l'on admet la conjecture VIII.\ref{arthurstint}, le point (ii) du théorème IX.\ref{thm1517} détermine les opérateurs de Hecke $\mathrm{T}_{p}:\mathbb{Z}[\mathrm{X}_{17}]\to\mathbb{Z}[\mathrm{X}_{17}]$ pour tout nombre premier $p$. On obtient par exemple la formule suivante~:
\begin{multline*}
\frac{7}{286}
\hspace{4pt}
\mathrm{N}_{p}(\mathrm{E}_{16}\oplus\mathrm{A}_{1},\mathrm{E}_{8}\oplus\mathrm{E}_{8}\oplus\mathrm{A}_{1})
\hspace{4pt}=\hspace{4pt} \\
(5p^{4}+7p^{3}+7p^{2}+7p+5)\hspace{4pt}
\frac{p^{11}-\tau(p)+1}{691}-
26\hspace{4pt}\frac{p^{15}-\tau_{16}(p)+1}{3617}
\hspace{24pt},
\end{multline*}
$\tau_{16}(p)$ désignant ci-dessus le $p$-ième coefficient de Fourier de la forme modulaire (pour $\mathrm{SL}_{2}(\mathbb{Z})$) parabolique normalisée de poids $16$.

\bigskip
\textbf{5.2.} Détermination de $\mathrm{T}_{2}$ pour $n=15$

\bigskip
Le point (b) de la proposition 4.6 et les mêmes arguments que précédemment (en plus simple) donnent
$$
\hspace{8pt}
\frac{\mathrm{N}_{2}(\mathrm{E}_{15},\mathrm{E}_{7}\oplus\mathrm{E}_{8})}{\mathrm{N}_{2}(\mathrm{E}_{16},\mathrm{E}_{8}\oplus\mathrm{E}_{8})}
\hspace{4pt}=\hspace{4pt}
\frac{\vert\mathrm{R}(\mathrm{Bor}_{16})\vert}{\vert\mathrm{R}(\mathrm{E}_{16})\vert}
\hspace{8pt},\hspace{8pt}
\frac{\mathrm{N}_{2}(\mathrm{E}_{15},\mathrm{E}_{7}\oplus\mathrm{E}_{8})}{\mathrm{N}_{2}(\mathrm{E}_{16},\mathrm{E}_{8}\oplus\mathrm{E}_{8})}
\hspace{4pt}=\hspace{4pt}
\frac{\vert\mathrm{R}(\mathrm{Bor}_{16})\vert}{\vert\mathrm{R}(\mathrm{E}_{8}\oplus\mathrm{E}_{8}))\vert}
\hspace{8pt},
$$
$\mathrm{Bor}_{16}$ étant le réseau unimodulaire impair de dimension $16$, introduit en III.3.3.2, qui ``fait le pont'' entre $\mathrm{E}_{16}$ et $\mathrm{E}_{8}\oplus\mathrm{E}_{8}$. Comme l'on a $\mathrm{R}(\mathrm{Bor}_{16})=\mathbf{D}_{8}\coprod\mathbf{D}_{8}$ (Scholie III.3.3.2), on trouve
$$
\hspace{24pt}
\frac{\mathrm{N}_{2}(\mathrm{E}_{15},\mathrm{E}_{7}\oplus\mathrm{E}_{8})}{\mathrm{N}_{2}(\mathrm{E}_{16},\mathrm{E}_{8}\oplus\mathrm{E}_{8})}
\hspace{4pt}=\hspace{4pt}
\frac{7}{15}
\hspace{24pt},\hspace{24pt}
\frac{\mathrm{N}_{2}(\mathrm{E}_{7}\oplus\mathrm{E}_{8},\mathrm{E}_{15})}{\mathrm{N}_{2}(\mathrm{E}_{8}\oplus\mathrm{E}_{8},\mathrm{E}_{16})}
\hspace{4pt}=\hspace{4pt}
\frac{7}{15}
\hspace{24pt}.
\leqno{(*)}
$$
On en déduit la matrice de l'opérateur de Hecke $\mathrm{T}_{2}:\mathbb{Z}[\mathrm{X}_{15}]\to\mathbb{Z}[\mathrm{X}_{15}]$ dans la base $(\mathrm{E}_{15},\mathrm{E}_{7}\oplus\mathrm{E}_{8})$ (matrice que l'on note encore $\mathrm{T}_{2}$)~:
$$
\hspace{24pt}
\mathrm{T}_{2}
\hspace{4pt}=\hspace{4pt}
\begin{bmatrix}
10377 & 8505 \\ 6006 & 7878
\end{bmatrix}
\hspace{24pt}.
$$

\bigskip
Comme précédemment, si l'on admet la conjecture VIII.\ref{arthurstint}, le point (i) du théorème IX.\ref{thm1517} montre que les égalités  (*) ci-dessus se généralisent pour tout nombre premier $p$~:
$$
\hspace{21pt}
\frac{\mathrm{N}_{p}(\mathrm{E}_{15},\mathrm{E}_{7}\oplus\mathrm{E}_{8})}{\mathrm{N}_{p}(\mathrm{E}_{16},\mathrm{E}_{8}\oplus\mathrm{E}_{8})}
\hspace{4pt}=\hspace{4pt}
\frac{p^{3}-1}{p^{4}-1}
\hspace{24pt},\hspace{24pt}
\frac{\mathrm{N}_{p}(\mathrm{E}_{7}\oplus\mathrm{E}_{8},\mathrm{E}_{15})}{\mathrm{N}_{p}(\mathrm{E}_{8}\oplus\mathrm{E}_{8},\mathrm{E}_{16})}
\hspace{4pt}=\hspace{4pt}
\frac{p^{3}-1}{p^{4}-1}
\hspace{21pt}.
$$

\bigskip
\textbf{5.3.} Sur le problème de la détermination de $\mathrm{T}_{2}$ pour $n=23$

\medskip
On constate que l'application $\pi:\mathrm{X}_{15}\to\mathrm{X}_{16}$ est injective~; ceci implique que l'alternative (b.2) de 4.6 n'a pas lieu pour $n=15$ (cet argument a été en fait utilisé ci-dessus pour déterminer $\mathrm{T}_{2}$ pour $n=15$). Par contre l'application $\pi:\mathrm{X}_{23}\to\mathrm{X}_{24}$ n'est pas injective et nous allons voir que l'alternative (b.2) de 4.6 a bien lieu pour $n=23$.

\medskip
Explicitons par exemple $\pi^{-1}([\mathrm{E}_{16}\oplus\mathrm{E}_{8}])$. Le quotient $\mathrm{O}(\mathrm{E}_{16}\oplus\mathrm{E}_{8})\backslash\mathrm{R}(\mathrm{E}_{16}\oplus\mathrm{E}_{8})$ s'identifie à la réunion disjointe $\mathrm{O}(\mathrm{E}_{16})\backslash\mathrm{R}(\mathrm{E}_{16})\coprod\mathrm{O}(\mathrm{E}_{8})\backslash\mathrm{R}(\mathrm{E}_{8})$ (et est donc un ensemble à deux éléments)~; on en déduit, grâce à la proposition 2.6,  $\pi^{-1}([\mathrm{E}_{16}\oplus\mathrm{E}_{8}])=\{[\mathrm{E}_{15}\oplus\mathrm{E}_{8}],[\mathrm{E}_{16}\oplus\mathrm{E}_{7}]\}$.

\medskip
On montre ci-après que les deux éléments $[\mathrm{E}_{15}\oplus\mathrm{E}_{8}]$ et $[\mathrm{E}_{16}\oplus\mathrm{E}_{7}]$ de $\mathrm{X}_{23}$ sont $2$-voisins et que l'ensemble $\mathrm{e}^{-1}(\{[\mathrm{E}_{15}\oplus\mathrm{E}_{8}],[\mathrm{E}_{16}\oplus\mathrm{E}_{7}]\})$ (notation de~4.5), vu comme un sous-ensemble de $\mathrm{B}_{24}^{\mathrm{A}_{1}}$ (notation introduite dans la discussion du paragraphe 4 intitulée ``$2$-voisins, le point de vue de Borcherds (suite)''), est le singleton $\{[(\mathrm{I}_{2}\oplus\Lambda;\iota)]\}$,  $\Lambda$ désignant le réseau unimodulaire de dimension $22$ correspondant au lagrangien évident du $\mathrm{e}$-module $\mathop{\mathrm{r\acute{e}s}}(\mathrm{E}_{15}\oplus\mathrm{E}_{7})$ ($\Lambda$ pourrait être aussi noté $(\mathrm{E}_{15}\oplus\mathrm{E}_{7})^{+}$) et $\iota:\mathrm{A}_{1}\to\mathrm{I}_{2}\oplus\Lambda$ désignant la représentation induite par la représentation canonique $\mathrm{A}_{1}\to\mathrm{I}_{2}$.

\bigskip
On pose $S=\mathrm{A}_{1}\oplus\mathrm{A}_{1}\oplus\mathrm{E}_{15}\oplus\mathrm{E}_{7}$~; on note $\sigma$ l'élément évident de $\mathrm{O}(S)$ qui échange les deux facteurs $\mathrm{A}_{1}$.  On note $\varpi_{i}$ le générateur du résidu du\linebreak $i$-ème facteur de $S$. Le $\mathrm{qe}$-module $\mathop{\mathrm{r\acute{e}s}}S$ est donc un $\mathbb{Z}/2$-espace vectoriel de dimension $4$, de base $\{\varpi_{1},\varpi_{2},\varpi_{3},\varpi_{4}\}$, la forme quadratique d'enlacement étant définie par $\mathrm{q}(\varpi_{i})=\frac{1}{4}$ pour $i=1,2$,  $\mathrm{q}(\varpi_{i})=-\frac{1}{4}$ pour $i=3,4$ et $\varpi_{i}.\varpi_{j}=0$ pour $i\not=j$~; la structure de $\mathrm{qe}$-module de Venkov est quant à elle déterminée par par $\mathrm{qm}(\varpi_{i})=\frac{1}{4}$ pour $i=1,2$ et $\mathrm{qm}(\varpi_{i})=\frac{3}{4}$ pour $i=3,4$. Le $\mathrm{qe}$-module $\mathop{\mathrm{r\acute{e}s}}S$ possède deux lagrangiens~:

\smallskip
-- le sous-espace, disons $J_{1}$, engendré par $\varpi_{1}+\varpi_{3}$ et $\varpi_{2}+\varpi_{4}$,

\smallskip
-- le sous-espace, disons $J_{2}$, engendré par $\varpi_{1}+\varpi_{4}$ et $\varpi_{2}+\varpi_{3}$.

\smallskip
On observera que ces deux lagrangiens sont les graphes des deux isomorphismes de $\mathrm{qe}$-modules de $\mathop{\mathrm{r\acute{e}s}}(\mathrm{A}_{1}\oplus\mathrm{A}_{1})$ sur $\langle-1\rangle\otimes\mathop{\mathrm{r\acute{e}s}}(\mathrm{E}_{15}\oplus\mathrm{E}_{7})$ et qu'il sont échangés par $\sigma$. 

\smallskip
On note $P_{k}$, $k=1,2$, le réseau unimodulaire pair avec $S\subset P_{k}\subset S^{\sharp}$ et $P_{k}/S=J_{k}$. Il est clair que $P_{1}$ et $P_{2}$ sont tous deux isomorphes à $\mathrm{E}_{16}\oplus\mathrm{E}_{8}$ et qu'il sont échangés par $\mathbb{Q}\otimes_{\mathbb{Z}}\sigma$. 

\smallskip
On pose $K=J_{1}\cap J_{2}$~; $K$ est le sous-espace de dimension $1$ de $\mathop{\mathrm{r\acute{e}s}}S$ engendré par $\varpi_{1}+\varpi_{2}+\varpi_{3}+\varpi_{4}$. On note $R$ le réseau pair avec $S\subset R\subset S^{\sharp}$ et $R/S=K$~; il est clair que l'on a $R=P_{1}\cap P_{2}$. On identifie le $\mathrm{qe}$-module $\mathop{\mathrm{r\acute{e}s}}R$ à $K^{\perp}/K$ (Proposition II.1.1)~; on observe que $\mathop{\mathrm{r\acute{e}s}}R$ est engendré par les classes de $\varpi_{1}+\varpi_{3}$ et $\varpi_{1}+\varpi_{4}$, cette observation permet de se convaincre de ce que $\mathop{\mathrm{r\acute{e}s}}R$ est isomorphe à $\mathrm{H}(\mathbb{Z}/2)$. On constate donc que $P_{1}$ et $P_{2}$ sont $2$-voisins, disons dans $\mathbb{Q}\otimes_{\mathbb{Z}}S$.

\medskip
Soit $L_{k}$, $k=1,2$, l'orthogonal dans $P_{k}$ du premier facteur $\mathrm{A}_{1}$ de $S$. On constate que l'on a $L_{1}\simeq\mathrm{E}_{15}\oplus\mathrm{E}_{8}$ et $L_{2}\simeq\mathrm{E}_{16}\oplus\mathrm{E}_{7}$. La proposition 4.2 montre bien que $L_{1}$ et $L_{2}$ sont $2$-voisins (disons dans l'orthogonal, dans $\mathbb{Q}\otimes_{\mathbb{Z}}S$, du premier facteur $\mathrm{A}_{1}$ de $S$).

\medskip
Soit $J_{3}$ le sous-espace vectoriel de $\mathop{\mathrm{r\acute{e}s}}S$  engendré par $\varpi_{1}+\varpi_{2}$ et $\varpi_{3}+\varpi_{4}$~; $J_{3}$ est un ``lagrangien au sens bilinéaire'' de $\mathop{\mathrm{r\acute{e}s}}S$. Le $\mathrm{e}$-module $\mathop{\mathrm{r\acute{e}s}}(\mathrm{A}_{1}\oplus\mathrm{A}_{1})$ (resp. $\mathop{\mathrm{r\acute{e}s}}(\mathrm{E}_{15}\oplus\mathrm{E}_{7})$) possède un unique lagrangien, disons $J_{4}$ (resp. $J_{5}$)~; $J_{3}$ est la somme orthogonale de $J_{4}$ et $J_{5}$. Le réseau unimodulaire associé au couple  $(\mathrm{A}_{1}\oplus\mathrm{A}_{1};J_{4})$ est isomorphe à $\mathrm{I}_{2}$. Comme on l'a dit plus haut, on note $\Lambda$ le réseau unimodulaire de dimension $22$ associé au couple $(\mathrm{E}_{15}\oplus\mathrm{E}_{7};J_{5})$~; $\Lambda$ est le premier réseau de la table de Conway et Sloane  \cite[Chap. 16, Table 16.7, $\dim=22$]{conwaysloane}. On se convainc de ce que $\Lambda$ ne représente pas $1$ en invoquant la structure de $\mathrm{qe}$-module de Venkov de $\mathop{\mathrm{r\acute{e}s}}(\mathrm{E}_{15}\oplus\mathrm{E}_{7})$. Le réseau unimodulaire associé au couple $(S;J_{3})$ est donc isomorphe à $\mathrm{I}_{2}\oplus\Lambda$. On constate que l'on a $K\subset J_{3}$~; il en résulte que le lagrangien au sens bilinéaire de $\mathop{\mathrm{r\acute{e}s}}R$ est $J_{3}/K$ et que le réseau unimodulaire impair dont les $2$-voisins pairs sont $P_{1}$ et $P_{2}$ coïncide avec le réseau unimodulaire associé au couple $(S;J_{3})$.

\medskip
On se convainc enfin de l'égalité dans $\mathrm{B}_{24}^{\mathrm{A}_{1}}$
$$
\mathrm{e}^{-1}(\{[\mathrm{E}_{15}\oplus\mathrm{E}_{8}],[\mathrm{E}_{16}\oplus\mathrm{E}_{7}]\})
\hspace{4pt}=\hspace{4pt}
\{[(\mathrm{I}_{2}\oplus\Lambda;\iota)]\}
$$
en contemplant la première colonne de la table de Conway et Sloane évoquée ci-dessus.

\medskip
Compte tenu de cette égalité, on peut déterminer $\mathrm{N}_{2}(\mathrm{E}_{15}\oplus\mathrm{E}_{8},\mathrm{E}_{16}\oplus\mathrm{E}_{7})$ à l'aide par exemple de 4.4 :
$$
\hspace{4pt}
\mathrm{N}_{2}(\mathrm{E}_{15}\oplus\mathrm{E}_{8},\mathrm{E}_{16}\oplus\mathrm{E}_{7})
\hspace{2pt}=\hspace{2pt}
\frac{\vert\mathrm{O}(\mathrm{E}_{15}\oplus\mathrm{E}_{8})\vert}{\vert\mathrm{O}(\mathrm{A}_{1}\oplus\mathrm{E}_{15}\oplus\mathrm{E}_{7})\vert}
\hspace{2pt}=\hspace{2pt}
\frac{\vert\mathrm{O}(\mathrm{E}_{8})\vert}{\vert\mathrm{O}(\mathrm{A}_{1}\oplus\mathrm{E}_{7})\vert}
\hspace{2pt}=\hspace{2pt}
120
\hspace{4pt}.
$$

\bigskip
Pour terminer ce sous-paragraphe, reprenons pour $n=23$ l'observation faite par Nebe et Venkov pour $n=24$ ~: comme la somme $\sum_{y\in\mathrm{X}_{23}}\mathrm{N}_{2}(x,y)$ est connue pour tout $x$ (Proposition 4.3), il suffit pour déterminer l'opérateur de Hecke $\mathrm{T}_{2}:\mathbb{Z}[\mathrm{X}_{23}]\to\mathbb{Z}[\mathrm{X}_{23}]$ de calculer $\mathrm{N}_{2}(x,y)$ pour $x\not=y$. Cette observation et le point (b) de 4.6 (et l'exemple précédent) amènent à se poser la question suivante~:

\smallskip
Peut-on déterminer $\mathrm{T}_{2}$ pour $n=23$, à la Nebe-Venkov, en contemplant simplement la table de Borcherds \cite[Chap. 17]{conwaysloane} et de celle de Conway et Sloane \cite[Chap. 16, Table 16.7, $\dim=22$]{conwaysloane}~?

\chapter{Tables} 
\label{tableschap9}

\begin{table}[htp]

\renewcommand{\arraystretch}{1.5}

{\small

\hspace{-.6cm}
\begin{tabular}{|c|c||c|c|}
\hline
$g$ & $\psi(\pi,{\rm St})$ & $g$ & $\psi(\pi,{\rm St})$ \cr
\hline 

 &  & $7$ & ${\rm Sym}^2 \Delta_{11} \oplus \Delta_{17}[4]
\oplus \Delta_{11}[2]$ \cr \hline

$1$ & ${\rm Sym}^2 \Delta_{11}$ & $7$ & ${\rm Sym}^2 \Delta_{11} \oplus
\Delta_{15}[6]$ \cr
\hline
$2$ & $\Delta_{21}[2] \oplus [1]$ & $8$ & $\Delta_{15}[8] \oplus [1]$ \cr
\hline
$3$ & ${\rm Sym}^2 \Delta_{11} \oplus \Delta_{19}[2]$ & $8$ &
$\Delta_{21}[2] \oplus \Delta_{17}[2] \oplus \Delta_{11}[4] \oplus [1]$ \cr
\hline
$4$ &  $\Delta_{21}[2] \oplus \Delta_{17}[2] \oplus [1]$ & $8$ &
$\Delta_{19}[4] \oplus \Delta_{11}[4] \oplus [1]$\cr
\hline
$4$ & $\Delta_{19}[4] \oplus [1]$ & $8$ & $\Delta_{21,9}[2] \oplus
\Delta_{15}[4] \oplus [1]$\cr
\hline
$5$ & ${\rm Sym}^2 \Delta_{11} \oplus \Delta_{19}[2] \oplus
\Delta_{15}[2]$ & $9$ & ${\rm Sym}^2 \Delta_{11} \oplus \Delta_{19}[2]
\oplus \Delta_{11}[6]$\cr
\hline
$5$ & ${\rm Sym}^2 \Delta_{11} \oplus \Delta_{17}[4]$ & $9$ & ${\rm
Sym}^2 \Delta_{11} \oplus \Delta_{19,7}[2] \oplus \Delta_{15}[2] \oplus
\Delta_{11}[2]$\cr
\hline
$6$ & $\Delta_{17}[6] \oplus [1]$ & $10$ & $\Delta_{21}[2] \oplus
\Delta_{11}[8] \oplus [1]$\cr
\hline
$6$ & $\Delta_{21}[2] \oplus \Delta_{15}[4] \oplus [1]$ & $10$ &
$\Delta_{21,5}[2] \oplus \Delta_{17}[2] \oplus \Delta_{11}[4] \oplus [1]$\cr
\hline
$6$ & $\Delta_{21,13}[2] \oplus \Delta_{17}[2] \oplus [1]$ & $11$ &
${\rm Sym}^2 \Delta_{11} \oplus \Delta_{11}[10]$\cr
\hline
$7$ & ${\rm Sym}^2 \Delta_{11} \oplus \Delta_{19}[2] \oplus \Delta_{15}[2]
\oplus \Delta_{11}[2]$ & $12$ & $\Delta_{11}[12] \oplus [1]$ \cr
\hline
\end{tabular} \ps  
}
\caption{Param\`etres standards $\psi(\pi,{\rm St})$ des repr\'esentations
$\pi$ dans $\Pi_{\rm cusp}({\rm
Sp}_{2g})$ engendr\'ees par une forme de Siegel de poids $12$ pour ${\rm
Sp}_{2g}(\Z)$, en genre
$g\leq 12$. }
\label{tablek=12}
\end{table}

\begin{sidewaystable}[!htbp]
\centering
{\tiny
\begin{tabular}{c||c|c|c|c|c|c|}

\backslashbox{{\rm Gram}}{$(j,k)$} & $(0,10)$ & $(6,8)$ & $(0,12)$ & $(8,8)$ & $(12,6)$ & $(4,10)$ \cr   

\hline

& & & & & & \cr

$\left[\begin{array}{cc} 2 & -1 \\ -1 & 2\end{array}\right]$ & $1$ & $Y^2 X^2(X
-Y)^2$ & $1$ &
$Y^2 X^2 (X-Y)^2 (X^2 -Y X + Y^2) $ & $Y^4 X^4 (X-Y)^4$ & $(X^2-YX+Y^2)^2$ \cr

& & & & & & \cr

\hline 

& & & & & & \cr

$\left[\begin{array}{cc} 2 & 0 \\ 0 & 2\end{array}\right]$ & $-2$ & $-2
\, Y^2X^2 (X^2+Y^2)$ & $10$ &
$-2 \,Y^2 X^2 (X^4 - 5\, Y^2 X^2 + Y^4)$ & $-2 \,Y^4 X^4 (X-Y)^2 (X+Y)^2$ & $-2
\, (X^4-9 \,Y^2 X^2 +Y^4)$ \cr

& & & & & & \cr

\hline 

& & & & & & \cr

$\left[\begin{array}{cc} 2 & -1 \\ -1 & 4\end{array}\right]$ & $-16$ &
\shortstack{ $- 8 \,Y^2 ( X-2 \,Y ) ( X+Y )$ \\ $ \times ( 2\,X^2 - 2\, Y X + Y^2 )$ } &  $-88$ &
\shortstack{ $-8\,Y^2 (2\,X^6 - 6 \,Y X^5 + 14 \,Y^2 X^4 $ \\ $- 18\, Y^3 X^3 + 14
\, Y^4 X^2 - 6\,Y^5 X + 3\,Y^6)$} & \shortstack{$-16 \,Y^4
X^2(X-2\,Y)(X-Y)^2$ \\ $\times
(X+Y)(X^2-YX+Y^2)$} & 
\shortstack{$-8 (2 X^4 - 4\,Y X^3 + 21\,Y^2 X^2$ \\ $ - 19\,Y^3 X - 14\,Y^4)$}\cr

& & & & & & \cr

\hline 

& & & & & & \cr

$\left[\begin{array}{cc} 2 & 0 \\ 0 & 4\end{array}\right]$ & $36$ & $12
\,Y^2 (X^2+Y^2)(3\,X^2-2\,Y^2)$ & $-132$ &
$12 \,Y^2 (X^2-3 \,Y^2) (3 \,X^4 - 5\, Y^2 X^2 - Y^4) $ & \shortstack{$12
\,Y^4 X^2 (3 \, X^6 - 10 \, Y^2 X^4$ \\ $+ 3\, Y^4 X^2 - 4\,Y^4)$}  & $36 \,(X^4 - 7\,Y^2 X^2 -\, 7 Y^4)$ \cr

& & & & & & \cr

\hline

& & & & & & \cr
\shortstack{$\left[\begin{array}{cc} 2 & -1 \\ -1 & 6\end{array}\right]$\\$ $} &
\shortstack{$99$\\$ $\\$ $\\$ $} &
\shortstack{$3 \,Y^2 (33 \, X^4 - 66 \, Y X^3$ \\  $- 91 \, Y^2 X^2 + 124 \, Y^3 X - 44 \,
Y^4)$\\ $ $}
&
\shortstack{$1275$\\$ $\\$ $\\$ $} 
&
\shortstack{$3 \,Y^2 (33 \, X^6 - 99 \, Y X^5$ \\ $+ 410 \, Y^2 X^4 - 655 \, Y^3 X^3$ \\ $+ 343 \, Y^4 X^2 - 32 \, Y^5 X - 132 \, Y^6)$} 
& 
\shortstack{$3 \,Y^4 (33 \, X^8 - 132 \, Y X^7 + 142 \, Y^2 X^6$ 
\\ $+ 36 \, Y^3 X^5 - 207 \, Y^4 X^4 + 200 \, Y^5 X^3 $ 
\\ $+ 88 \, Y^6 \, X^2 - 160 \, Y^7 X + 80 \, Y^8)$}  
& 
\shortstack{$99 \, (X^4 - 2 \, Y X^3 + 23 \, Y^2 X^2$ \\ $- 22 \, Y^3 X - 7 \,
Y^4)$ \\ $ $} \cr

& & & & & & \cr

\hline 

& & & & & & \cr

$\left[\begin{array}{cc} 4 & -2 \\ -2 &  4\end{array}\right]$ & $240$ & $0$ & $2784$ &
$1344\, Y^2 X^2 (X-Y)^2 (X^2 -Y X + Y^2) $ & $-240 \, Y^4 X^4 (X-Y)^4$ &
$-1680 \, (X^2-YX+Y^2)^2$ \cr

& & & & & & \cr

\hline

& & & & & & \cr

\shortstack{$\left[\begin{array}{cc} 2 & 0 \\ 0 &  6 \end{array}\right]$\\ $
$} & \shortstack{$-272$ \\$ $ \\$ $}
& 
\shortstack{$-16\, Y^2 \, (X-Y) (X+Y)$\\ $\times (17 \, X^2 + 13 \, Y^2)$}
&
\shortstack{$736$ \\$ $ \\$ $}
&
\shortstack{$-16 \, Y^2 (X-Y) (X+Y) $\\$\times (17 \, X^4 - 62 \, Y^2 X^2 + 39\, Y^4)$}
& 
\shortstack{$-16 \, Y^4 (X-Y) (X+Y) $\\$\times (17 \, X^6 - 63 \, Y^2 X^4 -
6\, Y^4 X^2 - 20 \, Y^6)$}
&
\shortstack{$-16\, (17\, X^4 - 96 \, Y^2 X^2 - 144 \, Y^4)$\\$ $} \cr
& & & & & & \cr
\hline

\end{tabular}   }

\caption{\small Coefficients de Fourier en $\frac{1}{2}{\rm Gram}$ d'un g\'en\'erateur bien choisi de ${\rm S}_{j,k}$} 
\label{tablecoeffdjk}
\end{sidewaystable}

\begin{table}
\renewcommand{\arraystretch}{1.5}
{\scriptsize
\hspace{-1.3cm}{
\begin{tabular}{|c||c|c|c|c|}
\hline $p$   & $\tau_{6,8}(p)$ & $\tau_{8,8}(p)$ & $\tau_{12,6}(p)$ & $\tau_{4,10}(p)$ \\
\hline
\hline $2$ & $0$ & $1344$ & $-240$ & $-1680$ \\
\hline $3$ & $-27000$ & $-6408$ & $68040$ & $55080$ \\
\hline $5$ & $2843100$ & $-30774900$ & $14765100$ & $-7338900$ \\
\hline $7$ & $-107822000$ & $451366384$ & $-334972400$ & $609422800$ \\
\hline $11$ & $3760397784$ & $13030789224$ & $3580209624$ & $25358200824$ \\
\hline $13$ & $9952079500$ & $-328006712228$ & $91151149180$ & $-263384451140$ \\
\hline $17$ & $243132070500$ & $5520456217764$ & $-11025016477020$ & $-2146704955740$ \\
\hline $19$ & $595569231400$ & $-28220918878760$ & $-22060913325080$ & $43021727413960$ \\
\hline $23$ & $-6848349930000$ & $79689608755152$ & $195863810691120$ & $-233610984201360$ \\
\hline $29$ & $53451678149100$ & $-1105748270340$ & $-1743496339579620$ & $-545371828324260$ \\
\hline $31$ & $234734887975744$ & $1851264166857664$ & $1979302106496064$ & $830680103136064$ \\
\hline $37$ & $448712646713500$ & $22115741387845324$ & $-3685951226317460$ & $11555498201265580$ \\
\hline $41$ & $-1267141915544076$ & $-29442241674311916$ & $106065086529460884$ & $-56208480716702316$ \\
\hline $43$ & $-1828093644641000$ & $308109789751260712$ & $74859021001125400$ & $160336767963955000$ \\
\hline $47$ & $-6797312934516000$ & $43932618784857504$ & $156108802652634720$ & $-116311331328502560$ \\
\hline $53$ & $30226618925077500$ & $-1178253142902441108$ & $-1224706812408694260$ & $-1944489787072554420$ \\
\hline $59$ & $-51143734375273800$ & $-3366234739477561080$ & $6289866383536712760$ & $1843701997761637080$ \\
\hline $61$ & $7626516406720684$ & $-8962102322409921476$ & $4857626575164933724$ & $2376385974282228124$ \\
\hline $67$ & $-12252758021387000$ & $14381861853876396664$ & $10336923176891703880$ & $487223803841627560$ \\
\hline $71$ & $-225641741059730736$ & $40475791736823448944$ & $-39237199980379430256$ & $18272191888645387344$ \\
\hline $73$ & $486083162996216500$ & $-11604559187113183148$ & $9078939377243940820$ & $26899631446378070740$ \\
\hline $79$ & $1424574980940205600$ & $14996327278915320160$ & $71873557961577515680$ & $-80184572998399700960$ \\
\hline $83$ & $-1351980902639367000$ & $-154502893221792192168$ & $94316650925918995560$ & $157078549808482338120$ \\
\hline $89$ & $-1127953215815294700$ & $-49999331367987019020$ & $115915137334350529140$ & $22873692749841743220$ \\
\hline $97$ & $-2710671093611565500$ & $765838865005585444804$ & $894968190691418183620$ & $-219326787347594393660$ \\
\hline $101$ & $14595359522423307804$ & $-1274759541025862678196$ & $75745749887557044204$ & $867394381514415093804$ \\
\hline $103$ & $18796572299556586000$ & $1130145111856472690992$ & $229164380766640031440$ & $-657903326636255684720$ \\
\hline $107$ & $-23385476046562641000$ & $542230976527798722984$ & $-3571178446181577738600$ & $-395867979731685155400$ \\
\hline $109$ & $36219247764172458700$ & $-884687494456719863780$ & $-2024515635534667135940$ & $30287010492785677180$ \\
\hline $113$ & $-53733316769620465500$ & $705599831303150185572$ & $-4230007868022803115420$ & $1657202008073896578660$
\\
\hline
\end{tabular}
}
}

\caption{Quelques valeurs propres d'opérateurs de Hecke en genre $2$ : l'entier $\tau_{j,k}(p)$, pour $p$ premier et $p\leq 113$.}
\label{taujk}
\end{table}

\begin{table}
\renewcommand{\arraystretch}{1.5}
{\scriptsize
\hspace{1cm}{
\begin{tabular}{|c||c|c|}
\hline $p$   & $\tau_{6,8}(p^{2})$ & $\tau_{8,8}(p^{2})$ \\
\hline
\hline $2$ & 409600 & 348160  \\
\hline $3$ & 333371700 & 748312020  \\
\hline $5$ & -15923680827500 & -395299890927500  \\
\hline $7$ & -253514141409500 & -155544419215478300  \\
\hline $11$ & -75764187476725473836 & 19641545832571328136244  \\
\hline $13$ & -4843967045593944889100 & -596184280686941758305260  \\
\hline $17$ & 101161485715920379759300 & -208424259842935445790738620  \\
\hline $19$ & 2430966330762186234484084 & -1388004707990982166729991276  \\
\hline $23$ & -129889399810754988793919900 & -36435169742921431436190920540  \\
\hline $29$ & -7216762572241226809807993676 & -18636070203076686393140997747116  \\
\hline
\hline $p$   & $\tau_{12,6}(p^{2})$ & $\tau_{4,10}(p^{2})$ \\
\hline $2$ & 4276480 & -700160 \\
\hline $3$ & -8215290540 & 1854007380 \\
\hline $5$ & 722477627072500 & -904546757727500 \\
\hline $7$ & -1126868422025500700 & -391120313742441500 \\
\hline $11$ & -2263452414601610414156 & -18738678558496864257356 \\
\hline $13$ & -299941151717771094659180 & 323494600665947822387860 \\
\hline $17$ & -94260803115254202283241660 & 70477693184423227137834820 \\
\hline $19$ & -475514565037103383307581676 & -1048771276144665792567133676 \\
\hline $23$ & -505868492227965057753270620 & -93299515424177439346879450460 \\
\hline $29$ & -11097155072276494608459664937516 & -2847689414234249875206600521516 \\
\hline
\end{tabular}
}
}

\vspace{1cm}
\caption{Quelques valeurs propres d'op\'erateurs de Hecke en genre $2$ : l'entier $\tau_{j,k}(p^{2})$, pour $p$ premier et $p\leq 29$.}
\label{taujkp2}
\end{table}

\begin{table}[htp]
\renewcommand{\arraystretch}{1.5}

{\small 
\begin{center}
\begin{tabular}{|c|c|c|c|}

\hline $i$ & $\psi_i$ & $\lambda_i$ & $g_i$ \cr

\hline $1$ & $[23]\oplus [1]$ & $8390655$ & $0$ \cr

\hline $2$ & ${\rm Sym}^2 \Delta_{11} \oplus [21]$ & $4192830$ & $1$ \cr

\hline $3$ & $\Delta_{21}[2] \oplus [1] \oplus [19]$ & $2098332$ & $2$ \cr

\hline $4$ & ${\rm Sym}^2 \Delta_{11} \oplus \Delta_{19}[2] \oplus [17]$ & $1049832$
& $3$  \cr

\hline $5$ & $\Delta_{19}[4] \oplus [1] \oplus [15]$ & $533160$ & $4$ \cr 

\hline $6$ & $\Delta_{21}[2] \oplus \Delta_{17}[2] \oplus [1] \oplus [15]$ & $519120$
& $4$ \cr

\hline $7$ & ${\rm Sym}^2 \Delta_{11} \oplus \Delta_{19}[2] \oplus \Delta_{15}[2]\oplus [13]$ & $268560$
& $5$ \cr

\hline $8$ & ${\rm Sym}^2 \Delta_{11} \oplus \Delta_{17}[4] \oplus [13]$ & $244800$
& $5$ \cr

\hline $9$ & $\Delta_{21}[2] \oplus \Delta_{15}[4] \oplus [1] \oplus[11]$ & $145152$
& $6$ \cr

\hline $10$ & $\Delta_{21,13}[2] \oplus \Delta_{17}[2] \oplus [1] \oplus [11]$ & $126000$
& $6$ \cr

\hline $11$ & $\Delta_{17}[6] \oplus [1] \oplus [11]$ & $99792$ & $6$ \cr 

\hline $12$ & ${\rm Sym}^2 \Delta_{11} \oplus \Delta_{15}[6] \oplus [9]$ & $91152$
& $7$ \cr

\hline $14$ & ${\rm Sym}^2 \Delta_{11} \oplus \Delta_{19}[2] \oplus \Delta_{15}[2] \oplus \Delta_{11}[2] \oplus [9]$ & $69552$
& $7$ \cr

\hline $16$ & ${\rm Sym}^2 \Delta_{11} \oplus \Delta_{17}[4] \oplus \Delta_{11}[2] \oplus [9]$ & $45792$
& $7$ \cr

\hline $13$ & $\Delta_{15}[8] \oplus [1] \oplus [7]$ & $89640$ & $8$ \cr

\hline $15$ & $\Delta_{21,9}[2] \oplus \Delta_{15}[4] \oplus [1] \oplus [7]$ & $51552$
& $8$ \cr

\hline $17$ & $\Delta_{19}[4] \oplus \Delta_{11}[4] \oplus [1] \oplus [7]$ & $35640$
& $8$ \cr

\hline $18$ & $\Delta_{21}[2] \oplus \Delta_{17}[2] \oplus \Delta_{11}[4] \oplus [1] \oplus [7]$ & $21600$
& $8$ \cr

\hline $19$ & ${\rm Sym}^2 \Delta_{11} \oplus \Delta_{19,7}[2] \oplus \Delta_{15}[2] \oplus \Delta_{11}[2] \oplus [5]$
& $17280$ & $9$ \cr

\hline $20$ & ${\rm Sym}^2 \Delta_{11} \oplus \Delta_{19}[2] \oplus \Delta_{11}[6] \oplus [5]$ & $5040$
& $9$ \cr

\hline $21$ & $\Delta_{21,5}[2] \oplus \Delta_{17}[2] \oplus \Delta_{11}[4] \oplus [1] \oplus [3]$
& $-7920$ & $10$ \cr

\hline $22$ & $\Delta_{21}[2] \oplus \Delta_{11}[8] \oplus [1] \oplus [3]$ & $-16128$
& $10$\cr

\hline $23$ & ${\rm Sym}^2 \Delta_{11} \oplus \Delta_{11}[10] \oplus [1]$ & $-48528$ &
$11$ \cr

\hline $24$ & $\Delta_{11}[12]$ & $-98280$ & $12$ \cr
\hline
\end{tabular} \end{center}
}
\caption{Param\`etres standards $\psi(\pi,{\rm St})$ des $24$ repr\'esentations
$\pi$ dans $\Pi_{\rm disc}({\rm O}_{24})$ telles que $\pi_\infty$ est
triviale, rang\'ees par degr\'e croissant.}
\label{table24bis}
\end{table}

\newpage
${}$
\vspace{3 cm}
\begin{table}[htp]

{\scriptsize \renewcommand{\arraystretch}{1.8} \medskip
\begin{center}
\begin{tabular}{|c|c||c|c||c|c|}
\hline  $(m_1,\dots,m_8)$ & $\dim V_\lambda^\Gamma$ & $(m_1,\dots,m_8)$ & $\dim V_\lambda^\Gamma$ & $(m_1,\dots,m_8)$ & $\dim V_\lambda^\Gamma$ \\

\hline $(0, 0, 0, 0, 0, 0, 0, 0)$ & $1$ & $(4, 4, 0, 0, 0, 0, 0, 0)$ & $1$ & $(4, 4, 4, 4, 2, 2, 0, 0)$ & $1$ \\

\hline $(2, 2, 0, 0, 0, 0, 0, 0)$ & $1$ & $(4, 4, 2, 2, 0, 0, 0, 0)$ & $1$ & $(4, 4, 4, 4, 2, 2, 2, 2)$ & $1$ \\

\hline $(2, 2, 2, 2, 0, 0, 0, 0)$ & $1$ & $(4, 4, 2, 2, 2, 2, 0, 0)$ & $1$ & $(4, 4, 4, 4, 4, 0, 0, 0)$ & $1$ \\

\hline $(2, 2, 2, 2, 2, 2, 0, 0)$ & $1$ & $(4, 4, 2, 2, 2, 2, 2, 2)$ & $1$ & $(4, 4, 4, 4, 4, 2, 2, 0)$ & $1$ \\

\hline $(2, 2, 2, 2, 2, 2, 2, 2)$ & $1$ & $(4, 4, 4, 0, 0, 0, 0, 0)$ & $1$ & $(4, 4, 4, 4, 4, 4, 0, 0)$ & $1$ \\

\hline $(4, 0, 0, 0, 0, 0, 0, 0)$ & $1$ & $(4, 4, 4, 2, 2, 0, 0, 0)$ & $1$ & $(4, 4, 4, 4, 4, 4, 2, 2)$ & $1$ \\

\hline $(4, 2, 2, 0, 0, 0, 0, 0)$ & $1$ & $(4, 4, 4, 2, 2, 2, 2, 0)$ & $1$ & $(4, 4, 4, 4, 4, 4, 4, 0)$ & $1$ \\

\hline $(4, 2, 2, 2, 2, 0, 0, 0)$ & $1$ & $(4, 4, 4, 4, 0, 0, 0, 0)$ & $2$ & $(4, 4, 4, 4, 4, 4, 4, 4)$ & $2$ \\

\hline $(4, 2, 2, 2, 2, 2, 2, 0)$ & $1$ &  &  &  & \\

\hline

\end{tabular} \end{center}} \caption{Les $8$-uples d'entiers
$(m_1,\dots,m_8)$,
v\'erifiant $4 \geq m_1 \geq m_2 \geq \dots \geq m_8 \geq 0$, tels que $V_\lambda^\Gamma
\neq 0$, o\`u $V_\lambda$ est la repr\'esentation irr\'eductible de ${\rm
SO}(\R^{16})$ de plus haut poids $\lambda = \sum_{i=1}^8 m_i \varepsilon_i$
et $\Gamma={\rm SO}({\rm E}_8 \oplus {\rm E}_8)$.} \label{tableSO16}
\end{table}

\newpage
${}$
\vspace{2 cm}

\begin{table}[htp]
\renewcommand{\arraystretch}{1.5}
{\footnotesize
\begin{tabular}{cc}

{{}} $[22]$ &  $\Delta_{11}[11]$\cr
{{}} $\Delta_{15}[7]\oplus[8]$ &  
 $\Delta_{17}[5]\oplus[12]$\cr
{{}} $\Delta_{19}[3]\oplus[16]$ & 
 $\Delta_{21}\oplus[20]$\cr  
{{}} $\Delta_{17}[5]\oplus\Delta_{11}\oplus[10]$ & 
 $\Delta_{19}[3]\oplus\Delta_{11}[5]\oplus[6]$\cr
{{}} $\Delta_{19}[3]\oplus\Delta_{15}\oplus[14]$ &   
 $\Delta_{21}\oplus\Delta_{11}[9]\oplus[2]$\cr 
{{}} $\Delta_{21}\oplus\Delta_{15}[5]\oplus[10]$ & 
 $\Delta_{21}\oplus\Delta_{17}[3]\oplus[14]$\cr
{{}} $\Delta_{21}\oplus\Delta_{19}\oplus[18]$ &    
 $\Delta_{21,9}\oplus\Delta_{15}[5]\oplus[8]$\cr
{{}} $\Delta_{21,13}\oplus\Delta_{17}[3]\oplus[12]$ & 
 $\Delta_{19}[3]\oplus\Delta_{15}\oplus\Delta_{11}[3]\oplus[8]$\cr
{{}} $\Delta_{21}\oplus\Delta_{17}[3]\oplus\Delta_{11}[3]\oplus[8]$ & 
 $\Delta_{21}\oplus\Delta_{19}\oplus\Delta_{11}[7]\oplus[4]$\cr   
{{}} $\Delta_{21}\oplus\Delta_{19}\oplus\Delta_{15}[3]\oplus[12]$ & 
 $\Delta_{21}\oplus\Delta_{19}\oplus\Delta_{17}\oplus[16]$\cr   
{{}} $\Delta_{21,13}\oplus\Delta_{17}[3]\oplus\Delta_{11}\oplus[10]$ & 
 $\Delta_{21}\oplus\Delta_{19}\oplus\Delta_{15}[3]\oplus\Delta_{11}\oplus[10]$\cr
{{}} $\Delta_{21}\oplus\Delta_{19}\oplus\Delta_{17}\oplus\Delta_{11}[5]\oplus[6]$ &  
 $\Delta_{21}\oplus\Delta_{19}\oplus\Delta_{17}\oplus\Delta_{15}\oplus[14]$\cr  
{{}} $\Delta_{21,5}\oplus\Delta_{19}\oplus\Delta_{17}\oplus\Delta_{11}[5]\oplus[4]$ & 
 $\Delta_{21,9}\oplus\Delta_{19}\oplus\Delta_{15}[3]\oplus\Delta_{11}\oplus[8]$\cr 
{{}} $\Delta_{21,9}\oplus\Delta_{19,7}\oplus\Delta_{15}[3]\oplus\Delta_{11}\oplus[6]$ & 
 $\Delta_{21,13}\oplus\Delta_{19}\oplus\Delta_{17}\oplus\Delta_{15}\oplus[12]$\cr
{{}} $\Delta_{21}\oplus\Delta_{19}\oplus\Delta_{17}\oplus\Delta_{15}\oplus\Delta_{11}[3]\oplus[8]$ & 
 $\Delta_{21}\oplus\Delta_{19,7}\oplus\Delta_{17}\oplus\Delta_{15}\oplus\Delta_{11}[3]\oplus[6]$\cr
{{}} $\Delta_{21,5}\oplus\Delta_{19,7}\oplus\Delta_{17}\oplus\Delta_{15}\oplus\Delta_{11}[3]\oplus[4]$ & 
 $\Delta_{21,13}\oplus\Delta_{19}\oplus\Delta_{17}\oplus\Delta_{15}\oplus\Delta_{11}\oplus[10].$
\end{tabular} \ps
}
\caption{Param\`etres standards des $32$ repr\'esentations $\pi$ dans $\Pi_{\rm cusp}({\rm
SO}_{23})$ telles que $\pi_\infty=1$, admettant la conjecture VIII.\ref{conjaj2}. }
\label{tableso23}
\end{table}

\end{document}